\PassOptionsToPackage{pdfpagelabels=true}{hyperref}

\documentclass{amsbook}
\usepackage{amssymb}

\usepackage{enumitem, url}
\usepackage[pagebackref]{hyperref}
\makeatletter
\def\th@plain{\thm@headfont{\bfseries}\itshape}
\def\th@definition{\thm@headfont{\bfseries}\normalfont}
\makeatother
\newtheorem{theorem}{Theorem}[section]
\newtheorem{lem}[theorem]{Lemma}
\newtheorem{prop}[theorem]{Proposition}
\newtheorem{cor}[theorem]{Corollary}
\newtheorem{refo}[theorem]{Reformulation}

\theoremstyle{definition}
\newtheorem{defn}[theorem]{Definition}
\newtheorem{exe}[theorem]{Example}
\newtheorem{rem}[theorem]{Remark}
\newtheorem{notat}[theorem]{Notation}
\newtheorem{constr}[theorem]{Construction}
\newtheorem{ques}[theorem]{Question}
\newtheorem{conv}[theorem]{Convention}

\numberwithin{section}{chapter}
\numberwithin{equation}{chapter}

\newcommand{\N}{\mathbf N}
\newcommand{\Z}{\mathbf Z}
\newcommand{\Q}{\mathbf Q}
\newcommand{\R}{\mathbf R}
\newcommand{\C}{\mathbf C}
\newcommand{\K}{\mathbf K}
\newcommand{\A}{\mathbf A}
\newcommand{\F}{\mathbf F}
\newcommand{\T}{\mathbf T}
\newcommand{\G}{\mathbf G}
\newcommand{\Un}{\mathbf 1}
\newcommand{\So}{\mathbf{Sol}}
\newcommand{\Hi}{\mathcal H}
\newcommand{\Li}{\mathcal L}
\newcommand{\Ki}{\mathcal K}
\newcommand{\U}{\mathcal U}
\newcommand{\AC}{\mathcal A}
\newcommand{\M}{\mathcal M}
\newcommand{\Norm}{\mathcal N}
\newcommand{\OO}{\mathcal O}
\newcommand{\ENorm}{\mathcal E}
\newcommand{\Ideal}{\mathcal J}
\newcommand{\Ri}{\mathcal R}
\newcommand{\Pri}{\mathrm{Prim}}
\newcommand{\QD}{\mathrm{QD}}
\newcommand{\Max}{\mathrm{Max}}
\newcommand{\FC}{\mathrm{FC}}
\newcommand{\GL}{\mathrm{GL}}
\newcommand{\PGL}{\mathrm{PGL}}
\newcommand{\SL}{\mathrm{SL}}
\newcommand{\PSL}{\mathrm{PSL}}
\newcommand{\SO}{\mathrm{SO}}

\newcommand{\Sym}{\mathrm{Sym}}

\newcommand{\vN}{\mathrm{vN}}
\newcommand{\Tr}{\mathrm{Tr}}
\newcommand{\Aff}{\mathrm{Aff}}
\newcommand{\Ind}{\mathrm{Ind}}
\newcommand{\Res}{\mathrm{Res}}
\newcommand{\Proj}{\mathrm{Proj}}
\newcommand{\Pker}{\mathrm{Pker}}
\newcommand{\Hom}{\mathrm{Hom}}
\newcommand{\Comm}{\mathrm{Comm}}
\newcommand{\BS}{\mathrm{BS}} 
\newcommand{\Per}{\mathrm{Per}}
\newcommand{\Aut}{\mathrm{Aut}}
\newcommand{\Bohr}{\mathrm{Bohr}}
\newcommand{\Prof}{\mathrm{Prof}}
\newcommand{\Sub}{\mathrm{Sub}}
\newcommand{\Char}{\mathrm{Char}}

\usepackage{imakeidx}
\makeindex

\begin{document}

% \frontmatter
% hidden on December 7

\title{
\Huge
Unitary representations of groups,
\\
duals, and characters
% \normalsize
% \\
% .
% \\
% December 13, 2019
% \\
% Work in progress 
% \\
% % ---
% % Comments are welcome
% Changes needed, both math and LaTeX
}

\author{Bachir Bekka and Pierre de la Harpe}

%    author one information
% \author{Bachir Bekka}
\address{
IRMAR,
UMR-CNRS 6625,
Universit\'e de Rennes 1,
Campus Beaulieu,
F-35042 Rennes Cedex,
France.
}
\curraddr{}
\email{bachir.bekka@univ-rennes1.fr}

%    author two information
% \author{Pierre de la Harpe}
\address{
Section de math\'ematiques,
Universit\'e de Gen\`eve,
2--4 rue du Li\`evre,
C.P. 64,
CH--1211 Gen\`eve 4,
Suisse.
}
\curraddr{}
\email{Pierre.delaHarpe@unige.ch}

%    If any version of the Mathematics Subject Classification other
%    than the 2010 edition appears, then you have an old version
%    of the AMS-LaTeX collection and need to upgrade.  Download from
%    http://www.ams.org/tex/amslatex.html .


\thanks{
The first author has been supported in part 
by the ANR (French Agence Nationale de la Recherche)
through the projects Labex Lebesgue (ANR-11-LABX-0020-01)
and GAMME (ANR-14-CE25-0004).
\vskip15cm
December 16, 2019}

\maketitle

%    Dedication.  If the dedication is longer than a line or two,
%    remove the centering instructions and the line break.
%\cleardoublepage
%\thispagestyle{empty}
%    If this book uses the documentclass stml-l or mmono-s, change
%    13.5pc to 10.5pc.
%\vspace*{13.5pc}
%\begin{center}
%  Dedication text (use \\[2pt] for line break if necessary)
%\end{center}
%\cleardoublepage

%    Change page number to 7 if a dedication is present.

% \setcounter{page}{5}
% hidden on December 7

\tableofcontents

%    Include unnumbered chapters (preface, acknowledgments, etc.) here.
%-----------------------------------------------------------------------------
% Beginning of preface.tex
%-----------------------------------------------------------------------------
%
% AMS-LaTeX 1.2 sample file for a monograph, based on amsbook.cls.
% This is a data file input by chapter.tex.
%%%%%%%%%%%%%%%%%%%%%%%%%%%%%%%%%%%%

\chapter*{Foreword}

This is an expository book on unitary representations of topological groups,
and of several dual spaces,
which are spaces of such representations up to some equivalence.
The most important notions are defined for topological groups,
but a special attention is paid to the case of discrete groups.
\par

The unitary dual of a group $G$ is the space of equivalence classes
of its irreducible unitary representations;
it is both a topological space and a Borel space.
It is useful in important situations, but not of much use in many other situations,
for example because, for an infinite discrete group,
it is usually not countably separated as a Borel space.
\par

We discuss various candidates for being a dual space of a discrete group.
The two most important ones are the primitive dual,
the space of weak equivalence classes of unitary irreducible representations,
and the normal quasi-dual,
the space of quasi-equivalence classes of traceable factor representations.
The last space is parametrized by ``characters", which can be finite of infinite. 
The space of finite characters is an object particularly easy to describe:
it is in natural bijection with the space
of indecomposable central positive definite functions on the group,
that we call the Thoma dual.
We discuss other related spaces,
such as the space of finite dimensional unitary representations.
\par

All these spaces are much better behaved than the unitary dual.
They deserve to be studied on their own, and related to each other.
One may try to classify them for specific groups.
We illustrate the theory by systematically working out a series of specific examples:
Heisenberg groups,
affine groups of infinite fields,
solvable Baumslag--Solitar groups,
lamplighter groups,
and general linear groups.
\par

Operator algebras play an important role in the discussion,
in particular von Neumann algebras associated to a unitary representation
and C*-algebras associated to a locally compact group.

\chapter*{Introduction}
% Chapter 0
\label{ChapterIntroduction}
% \markboth{INTRODUCTION}{INTRODUCTION}

\emph{
Unitary group representations, that is,
representations of groups by linear isometries of a Hilbert space,
play an important role in a large variety of subjects,
including number theory, geometry, probability theory and theoretical physics. 
Given a group $G$, a major problem is to describe the most general representation of $G$
in terms of a dual space, that is, a set of ``elementary" representations of $G$,
up to a suitable equivalence relation.
The most common dual space of $G$ is the unitary dual $\widehat G$.
While the choice of $\widehat G$ as dual space is adequate for some classes of groups $G$
(as compact groups, locally compact abelian groups, or semi-simple Lie groups),
we will see that for other groups, especially for infinite discrete groups,
it is necessary to consider dual spaces different from $\widehat G$.
}

\section*
{Representation theory of finite groups}
 % Section 0.A Representation theory of finite groups
\label{S-Duals-Finite}

The theory of \textbf{linear representations} of finite groups
is a classical subject of fundamental importance.
It was founded by Frobenius in the late 1890s
and further developed by Burnside and Schur, among others (see below).
Let $G$ be a finite group and let $\pi$ be a representation of $G$
on a finite-dimensional complex vector space $V$.
There is a decomposition $V = V_1 \oplus \cdots \oplus V_n$
in a direct sum of \textbf{irreducible} invariant subspaces $V_i$.
Such a decomposition being usually not unique,
a more canonical one can be obtained as follows.
Ler $\{\pi_1, \cdots, \pi_k\}$
be the set of equivalence classes of irreducible representations of $G$.
The space $V$ admits a \emph{unique} decomposition $V = W_1 \oplus \cdots \oplus W_k$
in a sum of invariant subspaces $W_j$
such that the restriction of~$\pi$ to $W_j$
is equivalent to the direct sum of $n_j$ copies of $\pi_j$,
for some non-negative integer $n_j$;
the $W_j$ are called the \textbf{isotypical components} of $\pi$.
% where $n \pi = \underbrace{\pi \oplus \cdots \oplus \pi}_{n \hskip.2cm \text{times}}$.
As a result,
one faces two main problems in the representation theory of $G$:
\begin{itemize}
\item[(A)]
determine the \textbf{dual space} $\widehat G$ of $G$,
that is, the set of equivalence classes of irreducible representations of $G$;
\item [(B)]
for a given representation $\pi$ of $G$,
determine the multiplicities $n_\sigma$ in the isotypical decomposition
$\pi \simeq \bigoplus_{\sigma \in \widehat G} \hskip.1cm n_\sigma \sigma$,
where we use $\simeq$ to indicate equivalent representations.
\end{itemize}
\par

One way to approach Problem (A), already initiated by Frobenius,
is through the determination of the set of \textbf{characters}
of irreducible representations of $G$.
A representation $\pi$ of $G$ gives rise to its character
$$
\chi_\pi \, \colon \, G \to \C, \hskip.2cm g \mapsto \Tr(\pi(g)) ,
$$
where $\Tr$ stands for ``trace'';
note that $\chi_\pi$ is a class function.
(Recall that a class function or central function on a group $G$
is a function $f$ on $G$ which is constant on conjugacy classes:
$f(hgh^{-1}) = f(g)$ for all $g, h \in G$.)
\index{Central function = class function}
\index{Class function = central function}
Every representation is determined up to equivalence by its character.
The map $\sigma \mapsto \chi_\sigma$ sets up a bijection 
$$
\widehat G \, \xrightarrow{\approx} \, E(G)
$$
between $\widehat G$ and an appropriate set $E(G)$ of class functions on $G$.
By $\approx$, we indicate a bijection between sets,
or a homeomorphisms between topological spaces.
\index{$a5$@$\approx$ isomorphism! of sets (bijection)
or topological spaces (homeomorphism)}
\par

In view of Schur's orthogonality relations,
also Problem (B) can be solved through the use of characters.
\par

Following Frobenius, Molien, Burnside, E.~Noether, and others,
the representation theory of~$G$ can be viewed
as the study of modules over the \textbf{group algebra} 
$$
\C[G] \, = \, \Big\{\sum_{g \in G} c_g g \hskip.1cm \Big\vert \hskip.1cm c_g \in \C \Big\} .
$$
Indeed, every representation $\pi \, \colon G \to \GL(V)$
extends linearly to a representation
$$
\widetilde \pi \, \colon \, \C[G] \to \mathrm{End}(V) ,
\hskip.5cm \sum_{g \in G} c_g g \mapsto \sum_{g \in G} c_g \pi(g) .
$$
The map $\pi \to \widetilde \pi$ sets up a bijection between the representations of $G$
and the modules over $\C[G]$, preserving equivalence and irreducibility.
\par

The character $\chi_\pi$ of a representation $\pi \, \colon G \to \GL(V)$
has a natural extension $\widetilde \chi_\pi$ to $\C[G]$ given by
$$
\widetilde \chi_\pi \, \colon \, \C[G] \to \C, \hskip.2cm x \mapsto \Tr(\widetilde \pi (x)) .
$$
The map $\chi_\sigma \mapsto \widetilde \chi_\sigma$ sets up a bijection between 
$E(G)$ and an appropriate set $\Char(G)$ of linear functionals $\widetilde \chi$ 
on $\C[G]$ which satisfies the trace property: 
$$
\widetilde \chi (xy) \, = \, \widetilde \chi (yx)
\hskip.5cm \text{for all} \hskip.2cm
x,y \in \C[G].
$$
\par
 
The dual space $\widehat G$ can be described in terms of the ideal theory 
of the ring $\C[G]$ as follows:
the map $\sigma \to \ker(\widetilde \sigma)$ sets up a bijection 
$$
\widehat G \xrightarrow{\sim} \Pri(G)
$$
between $\widehat G$ and the space of \textbf{primitive ideals} of $\C[G]$.
(Recall that a primitive ideal in a ring $R$ is the annihilator of a simple $R$-module.)
% Observe that the set of primitive ideals of our special ring $R = \C[G]$
% coincides with the set of maximal two-sided ideals of $\C[G]$.

\section*
{Representations of compact groups and abelian groups} 
% Section 0.B Representations of compact groups and abelian groups
\label{S-RepsCompAbelian}
 
For a topological group $G$, representations of interest in this book
are \textbf{unitary representations} that is,
homomorphisms $\pi \, \colon G \to \U(\Hi)$ from $G$
to the unitary group of
a (not necessarily finite dimensional) complex Hilbert space $\Hi$,
which are continuous in an appropriate sense (see Section~\ref{S-DefUnitD}).
We will deal mostly with topological groups which are \emph{locally compact},
and also sometimes second-countable locally compact.
Such a group $G$ admits a \textbf{Haar measure},
that is, a translation-invariant positive measure $\mu_G$ on the Borel subsets of $G$;
second-countable locally compact groups
are characterized among standard Borel groups
by the fact that they admit a translation-invariant measure (see \cite{Mack--57}).
A locally compact group $G$ has a distinguished unitary representation,
the regular representation $\lambda_G$ on $L^2(G) = L^2(G, \mu_G)$,
which plays a central role in the theory --- as it already does for finite groups.
The \text{unitary dual} $\widehat G$ of $G$
is the set of equivalence classes of irreducible unitary representations of $G$. 

\subsection*{Compact groups}

Most of Frobenius' theory
extends to the case of a compact group $G$,
as shown by Schur and Weyl between 1924 and 1927.
First, every representation of $G$ is equivalent to a unitary representation. Then:
\begin{itemize}
\item[(I)]
every irreducible representation of $G$ is finite dimensional;
\item[(II)]
the map $\sigma \mapsto \chi_\sigma$ sets up a bijection between $\widehat G$
and the set of irreducible characters of $G$, which is a set of class functions on $G$;
\item[(III)]
every representation $\pi \, \colon G \to \U(\Hi)$ admits a
canonical isotypical decomposition
$\pi \simeq \bigoplus_{\sigma \in \widehat G} \hskip.1cm n_\sigma \sigma$;
\item[(IV)]
every unitary representation $\pi$ of $G$ ``extends" uniquely to a $*$-representation
$\widetilde \pi$ of the group algebra $L^1(G) = L^1(G, \mu_G)$;
\item [(V)]
the map $\rho \mapsto \ker(\widetilde \rho)$ sets up a bijection 
$$
\widehat G \, \xrightarrow{\sim} \, \Pri(G)
$$
between $\widehat G$ and the space of primitive (or maximal) ideals of $L^1(G)$.
\end{itemize}
The crucial tool for this extension is the use of averaging procedures
which are made possible by the fact that the Haar measure $\mu_G$ on $G$
is \emph{finite} when (and only when) $G$ is compact.
One should add to the previous list one more item,
a key step to find a complete system of irreducible representations of $G$:
\begin{itemize}
\item[(VI)]
every irreducible representation of $G$ is equivalent to a subrepresentation of 
the regular representation $\lambda_G$.
\end{itemize}

\subsection*{Locally compact abelian groups}

Let now $G$ be a locally compact abelian group.
Then $\widehat G$ coincides with the set of unitary characters $G \to \T$
and is itself a locally compact abelian group, for a natural topology.
By Pontrjagin and van Kampen duality theory, the dual group of 
$\widehat G$ can naturally be identified, as a topological group, with $G$.
This justifies the term ``dual'' used for $\widehat G$.
When $G$ is not locally compact abelian, the same word is used for $\widehat G$
even though this justification does not carry over.
\par

Let $\pi \, \colon G \to \U(\Hi)$ be a unitary representation of $G$.
In general, $\pi$ does not admit a direct sum decomposition into irreducible representations;
for example, when the group is infinite and discrete, its left regular representation
does not contain any irreducible subrepresentation.
% Proposition 7.A.3 a locally compact abelian group
% contains an irreducible subrepresentation
% (i.e., a subrepresentation of dimension $1$, because $G$ is abelian)
% if and only if the group is compact.
However, $\pi$ can be analyzed by means of
a unique projection-valued measure on $\widehat G$
or a unique collection of measure classes on $\widehat G$ with associated multiplicities
(see Chapter~\ref{Chapter-AbelianGroups}).
In case of the regular representation $\pi = \lambda_G$ on $L^2(G)$, 
such a decomposition generalizes the Plancherel formula
of classical Fourier analysis on $\R^n$ or $\T^n$.
\par

The set $E(G)$ coincides obviously with $\widehat G$ and here also
the map $\chi \to \ker(\widetilde \chi)$ is a bijection between $\widehat G$ 
and the space of primitive (or maximal) ideals of $L^1(G)$,
where $\widetilde \chi \, \colon L^1(G) \to \C$ is the multiplicative linear functional
canonically associated to $\chi$.
\par
 
As a result, we see that the representation theory of $G$
has the features (I)-(V) of Frobenius theory listed above, appropriately interpreted.
\par

Even when $G$ is non-compact, so that no unitary character of $G$
is equivalent to a subrepresentation of the regular representation $\lambda_G$,
the projection-valued measure on $\widehat G$
associated to $\lambda_G$ has $\widehat G$ as support.
Using a notion to be introduced in Section~\ref{SectionWC+FellTop},
it follows that every $\chi \in \widehat G$ is weakly contained in $\lambda_G$
and this is to be viewed as an analogue of feature (VI) in the list above.
 
\section*
{Primitive dual and locally compact groups of type~I}
% Section 0.C Primitive dual and locally compact groups of type~I
\label{S-PrimDualTypeI}
 
We turn now to the case of a \emph{general} locally compact group $G$.
By the Gel'fand--Raikov theorem, $G$ has plenty of irreducible unitary representations:
the unitary dual $\widehat G$ separates the points of $G$ (Theorem~\ref{Gel'fandRaikov}).
Observe that, when $G$ is not compact, $\widehat G$
will usually contain infinite dimensional representations
(for more on this, see Section~\ref{SectionFdrep}).
\par

As in the case of compact or abelian groups,
the unitary representations are in bijective correspondence
with the $*$-representations of the group algebra $L^1(G)$;
however, a more suitable algebra
is the \textbf{maximal C*-algebra} $C^*_{\rm max}(G)$ of $G$
which is the enveloping C*-algebra of $L^1(G)$.
Every unitary representation $\pi$ of $G$ ``extends"
to a representation $\widetilde \pi$ of $C^*_{\rm max}(G)$.
The \textbf{primitive ideal space} of $C^*_{\rm max}(G)$,
or \textbf{primitive dual} of $G$, is 
$$
\Pri(G) \, = \, \{\ker(\widetilde \pi) \mid \pi \in \widehat G\},
$$ 
and we have a canonical surjective map 
$$
\kappa^{\rm d}_{\rm prim} \, \colon \, \widehat G \, \twoheadrightarrow \, \Pri(G) ,
\hskip.2cm
\pi \mapsto \ker(\widetilde \pi) .
$$
The space $\Pri(G)$ is equipped with a natural topology, the Jacobson topology 
(see Section~\ref{SectionPrimC*}).
\par

There exists a description of $\Pri(G)$ at the group level:
it is the space of \textbf{weak equivalence classes}
of irreducible unitary representations of $G$.
Thus $\Pri(G)$ is a well-defined space for any topological group.
See Section \ref{SectionWC+FellTop} for the notion of weak equivalence,
and Sections~\ref{PrimIdealSpace} \& \ref{C*algLCgroup} for $\Pri(G)$.
\par

Besides $\widehat G$, the space $\Pri(G)$
is the most important alternative dual space of a group $G$.

\vskip.2cm
 
Let $\pi \, \colon G \to \U(\Hi)$ be a unitary representation of $G$ in a separable Hilbert space.
There is a way $\pi$ can be decomposed in irreducible parts.
As already mentioned in the previous section,
decomposition in direct sums do not suffice (unless $G$ is compact),
but there is a notion of \textbf{direct integral} of unitary representations,
and
$$
\pi \, = \, \int_{\widehat G}^{\oplus} n_\rho \hskip.1cm \rho \hskip.1cm d\mu(\rho) ,
$$
where $\mu$ is a Borel measure on $\widehat G$ (equipped with a natural Borel structure),
and every irreducible representation $\rho \in \widehat G$
occurs with a multiplicity $n_\rho$.
% (see Section~\ref{SectionDecomposingIrreps} = 1.G).
However, for some groups $G$, an unexpected pathology appears: 
there exist unitary representations $\pi$ of $G$ having two direct integral decompositions
$$
\pi
\, = \, \int_\Omega^\oplus \rho_\omega d\mu(\omega)
\, = \, \int_\Psi^\oplus \sigma_\psi d\nu(\psi),
$$
where $\Omega, \Psi$ are \emph{disjoint} Borel subsets of $\widehat G$,
and $\mu$ [respectively $\nu$] a measure on $\Omega$ [respectively $\Psi$],
and $\rho_\omega$ and $\sigma_\psi$ irreducible representations of $G$ 
for all $\omega \in \Omega$ and all $\psi \in \Psi$.
See Section~\ref{SectionDecomposingIrreps}.
\par

It turns out that the pathology mentioned above 
is related to the existence of factor representations of several types.
Let $\pi \, \colon G \to \U(\Hi)$ be a unitary representation.
The bicommutant $\pi(G)''$ of $\pi(G)$ is a von Neumann subalgebra of $\Li(\Hi)$
and $\pi$ is said to be \textbf{factorial}, or a \textbf{factor representation},
if $\pi(G)''$ is a factor, 
that is, if the centre of $\pi(G)''$ consists of scalar multiples of the identity only.
\par

Examples of factor representations
include representations of the form $n \pi$,
where $\pi$ is an irreducible unitary representation and $n$ a cardinal.
Groups without factor representations of another type
constitute an important class defined as follows:
a locally compact group $G$ is a \textbf{group of type I}
if every factor representation of $G$ is of the form $n \pi$
for some $\pi \in \widehat G$ and some cardinal $n$.
\par
 
The class of second-countable locally compact groups $G$ of type I
is characterized by a number of equivalent properties:
\begin{itemize}
\item
for every unitary representation $\pi$ of $G$,
the direct integral decomposition
$\pi = \int_{\widehat G} n_\rho \hskip.1cm \rho \hskip.1cm d\mu(\rho)$
mentioned above is unique in an appropriate sense
(Section~\ref{SectionDecomposingRepTypeI});
\item
the dual $\widehat G$, equipped with its natural Borel structure, is countably separated;
\item
the canonical map
$\kappa^{\rm d}_{\rm prim} \, \colon \widehat G \twoheadrightarrow \Pri(G)$
is injective, and therefore bijective.
\end{itemize}
The last two characterizations are part of Glimm's theorem (see Section \ref{SectionGlimm}).
\par

The class of locally compact groups of type I includes, among others,
abelian groups, compact groups, 
and Lie groups with finitely many connected components which are
either semisimple or nilpotent (see Section~\ref{SectionTypeI}).
As a rule, discrete groups are \emph{not} of type I:
more precisely, a discrete group $G$ is of type I
if and only if $G$ is virtually abelian (Thoma's theorem~\ref{discretTypes}).

\vskip.2cm

Let $G$ be a second-countable locally compact group.
Every unitary representation $\pi$ of $G$
has a canonical direct integral decomposition $\int_X \sigma_x d\mu(x)$
into \emph{factor} representations $\sigma_x$;
when $G$ is of type~I,
there are multiplicities $n_x$ and irreducible representations $\rho_x$
such that $\sigma_x = n_x \rho_x$ for all $x \in X$,
and therefore a canonical decomposition $\pi_x \simeq \int_X n_x \rho_x d\mu(x)$
into \emph{irreducible} representations $\rho_x$.
When $G$ is not of type I, factor representations can be of other types,
and this is a reason why direct integral decomposition in irreducible representations
are not canonical.

\vskip.2cm

Let $G$ be a second-countable locally compact group which is not of type I.
Since $\widehat G$ is not countably separated,
there is no systematic procedure 
to encode the equivalence classes of irreducible unitary representations of $G$.
By contrast, for the Borel structure associated to the Jacobson topology,
$\Pri(G)$ is a standard Borel space,
by a theorem of Effros (Subsection~\ref{SS:Subsection-Prim(A)}).
The space $\Pri(G)$ is therefore much better behaved than $\widehat G$,
and this justifies it as an alternative dual space for a locally compact group $G$.

\vskip.2cm

In Chapter \ref{ChapterTypeI}, there is a definition of being of type I
which applies to topological groups which are not necessarily locally compact.

\section*
{The normal quasi-dual of a locally compact group}
 % Section 0.D The normal quasi-dual of a locally compact group
 \label{S-NormalQS}

Let $G$ be a locally compact group.
If $G$ is not of type I, then $G$ admits factor representations that are not of type I.
It is convenient to introduce one more notion of equivalence:
two unitary representations $\pi_1, \pi_2$ of $G$ are \textbf{quasi-equivalent}
if there are cardinals $n_1, n_2$ such that $n_1 \pi_1$ and $n_2 \pi_2$ are equivalent
(more on this notion in Section \ref{Sectioncomppres})
The \textbf{quasi-dual} $\QD(G)$ of $G$ is the set of
quasi-equivalence classes of factor representations of $G$.
It has a natural topology, and there is a natural inclusion
$$
\kappa^d_{\rm qd} \, \colon \, \widehat G \, \hookrightarrow \, \QD(G)
$$
of the dual as a subspace of the quasi-dual. 
See Section \ref{SectionQuasidual}.
\par

When $G$ is $\sigma$-compact,
given a factor representation $\pi$ of $G$,
it is a fact (see Section \ref{C*kernelsprimitive})
that the $C^*$-kernel of $\pi$ is a primitive ideal of $C^*_{\rm max}(G)$,
that is, $\ker(\widetilde \pi) \in \Pri(G)$.
Consequently, there is a natural extension
$$
\kappa^{\rm qd}_{\rm prim} \, \colon \, \QD(G) \twoheadrightarrow \Pri(G)
$$
of the map $\kappa^{\rm d}_{\rm prim} \, \colon \widehat G \twoheadrightarrow \Pri(G)$.
\par

When $G$ is of type I, the inclusion $\widehat G \hookrightarrow \QD(G)$
is a surjective homeomorphism.
When $G$ is not of type I, a description of $\QD(G)$ is out of reach,
possibly even more than for $\widehat G$.
But one may try to describe $\Pri(G)$
as the collection of $C^*$-kernels of \emph{suitable} factor representations of $G$.
This leads to the definition of $\QD(G)_{\rm norm}$, as follows.

\vskip.2cm
 
Let $G$ be a
% second-countable 
locally compact group.
A unitary representation $\pi$ of $G$ is \textbf{normal} if it is factorial
and if the factor $\pi(G)''$ admits a (not necessarily finite) trace $\tau$
such that there exists a positive element $x \in C^*_{\rm max}(G)$
for which $0 < \tau(x) < \infty$.
The \textbf{normal quasi-dual} $\QD(G)_{\rm norm}$ is the set 
of normal factor representations of $G$, up to quasi-equivalence
(see Chapter~\ref{ChapterCharacters}).
It inherits from $\QD(G)$ a topology and the related Borel structure.
When $G$ is second-countable, the Borel space $\QD(G)_{\rm norm}$
is standard (Theorem~\ref{ThHalpern}).
\par

When $G$ is $\sigma$-compact, the restriction
of the map $\kappa^{\rm qd}_{\rm prim}$ defined above
provides a map
$$
\kappa^{\rm norm}_{\rm prim} \, \colon \, \QD(G)_{\rm norm} \to \Pri(G),
\hskip.2cm \pi \mapsto \ker(\widetilde \pi).
$$
When $G$ is a connected Lie group, a remarkable result of Puk\'anszky
shows that $\kappa^{\rm norm}_{\rm prim}$ is a bijection
(Chapter \ref{ChapterCharacters}).
This result indicates that, for such a group,
normal factor representations form the natural class of unitary representations
one should seek to classify.
It should be mentioned that
not every irreducible unitary representation of $G$ is normal when $G$ is not of type I
(Corollary~\ref{Cor-IrredNormalRep}). 
\par

For a general locally compact group $G$,
the map $\kappa^{norm}_{prim}$ can fail to be injective or surjective
(see Section~\ref{Section-CharactersPrimitive}).
Nevertheless, we think that $\QD(G)_{\rm norm}$ is an interesting object
which deserves to be studied as a dual space of $G$ in its own right.

\section*
{Characters of a locally compact group}
 % Section 0.E Characters of a locally compact group
\label{S-CharactersLCGroup}

Let $G$ be a second-countable locally compact group.
The normal quasi-dual $\QD(G)_{\rm norm}$ admits a parametrization
in terms of the set $\Char(G)$ of characters of $G$:
a normal factor representation $\pi$ of $G$ is uniquely determined, up to quasi-equivalence,
by its \textbf{character} $\chi_\pi$ which is
a (not necessarily finite) trace on the maximal C*-algebra $C^*_{\rm max}(G)$ of $G$
(see Section~\ref{Section-NormRepChar}).
\par

Assume that $G$ is of type I.
Then all the dual spaces we have considered above may be canonically identified
$$
\widehat G \, \xrightarrow{\approx} \,
\Pri(G) \, \xleftarrow{\approx} \,
\QD(G)_{\rm norm} \, \xrightarrow{\approx} \,
{\mathrm Char}(G).
$$
Nevertheless, something is gained:
the possibility to parametrize the unitary dual $\widehat G$
by means of the space $\Char(G)$ of characters,
in analogy to the case of finite or compact groups.
This programme has been carried out by Harish--Chandra for a connected semisimple Lie group,
in which case the set of characters can be identified
with a set of appropriate central distributions on the group \cite{Hari--54}. 

\section*
{Thoma's dual space}
 % Section 0.F Thoma's dual space
\label{S-ThomaDualSpace}

Let $G$ be a second-countable locally compact group.
The set $\Char(G)$ of characters of $G$ consists of two parts:
\begin{itemize}
\item
the set $E(G)$ of finite characters of $G$;
\item
the set $\Char(G) \smallsetminus E(G)$ of infinite characters.
\end{itemize}
The set $E(G)$ admits a simple description:
it can be viewed as the set of continuous functions 
$\varphi \, \colon G \to \C$ which are central, of positive type, normalized by $\varphi(e) = 1$,
and extremal with these properties.
(See Section \ref{S-FPosType} for functions of positive type,
and Section~\ref{SThoma'sdual} for the Thoma dual.)
\par

In general, the set $\Char(G) \smallsetminus E(G)$ is non empty
and seems difficult to describe
(see Section~\ref{Section-ExamplesInfiniteChar}, for the construction of a few
infinite characters of the Baumslag--Solitar group $\BS(1, p)$ and of the lamplighter group).
By contrast, $E(G)$ can be determined for several discrete groups
(see Example~\ref{exThomaDiscrete} and Chapter~\ref{ChapterThomadualExamples}).
Moreover, for some classes of discrete groups $G$ such as groups of polynomial growth,
all characters are finite (that is, $\Char(G) = E(G)$) and the canonical map 
$E(G) \twoheadrightarrow \Pri(G)$ is a bijection
(Theorem \ref{Th-HoweKa}).
% (see Proposition~\ref{Prop-NilGr}) = 11.D.2 ????
\par

Thoma made a crucial use of the space $E(G)$
in his characterization of discrete groups which are of type I
and suggested that $E(G)$ should be viewed
as a dual space of discrete groups \cite{Thom--64a, Thom--64b}.
We will refer to $E(G)$ as the \textbf{Thoma's dual} of $G$.
Again, Thoma's dual consists of two parts:
\begin{itemize}
\item
the set $E(G)_{\rm fd}$ of almost periodic functions in $E(G)$; 
this set $E(G)_{\rm fd}$ paramet\-rizes the space $\widehat G_{\rm fd}$
of equivalence classes of finite dimensional irreducible unitary representations of $G$;
\item
the set $E(G) \smallsetminus E(G)_{\rm fd}$,
of which each element gives rise to a homomorphism
of $G$ in the unitary group $\U(\mathcal M)$ of a factor $\mathcal M$ of type II$_1$. 
\end{itemize}

For a discrete group $G$ which is not of type I, the set $E(G) \smallsetminus E(G)_{\rm fd}$ 
is non empty (Theorem~\ref{ThGlimm2}).
By contrast, when $G$ is connected, we have $E(G) = E(G)_{\rm fd}$,
by a result of Kadison and Singer (\cite{KaSi--52}; see also \cite{SevN--50}).
In particular, if $G$ is a semisimple connected Lie group without compact factor,
$E(G)$ consists only of the trivial character $1_G$.

\section*
{Historical comments on dual spaces of groups}
% Section 0.G Historical comments on duals spaces of groups
\label{S-Historic}
 
\subsection*{Unitary dual}

The study of linear representations of finite groups was initiated
by Frobenius in a series of articles between 1896 and 1900 (see Vol.\ III of \cite{Frob--Col}).
This work was continued, among others, by Burnside and Schur
who established in particular the fact
that such representations are equivalent to unitary representations.
\par
 
Using invariant integration, Schur and Weyl extended between 1924 and 1927
most of Frobenius' representation theory from finite groups to compact groups
(see in particular \cite{Schu--24} and \cite{Weyl--24}, as well as \cite{Bore--86}).
This culminated with the Peter-Weyl theorem \cite{PeWe--27}
and the classification of the unitary dual 
of simple compact Lie groups (see \cite{Weyl--Col}).
\par

The proof of the existence of a Haar measure \cite{Haar--33}
opened up the possibility to a theory of unitary representations for general locally compact groups.
A strong impetus for the development of the theory
was provided by applications to quantum theory.
This motivated Stone, von Neumann, and Weyl
to determine the unitary dual of the real Heisenberg group
(see \cite{Ston--30} and \cite{vNeu--31}),
and Wigner to study the unitary dual of the Poincar\'e group $\SO(3,1) \ltimes \R^4$
(see \cite{Wign--39}, as well as \cite{Rose--04}).
\par

A duality theory for locally compact abelian groups
was developed around 1935 by Pontrjagin and van Kampen,
and exposed in the influential books by Pontrjagin and Weil,
\cite{Pont--39} and \cite{Weil--40}.
The classical spectral theorem for a one-parameter family of unitary operators
(which amounts to the classification of unitary representations of $\R$)
was generalized to locally compact abelian groups
by Segal, Naimark, Ambrose and Godement
(see SNAG's theorem in Section~\ref{SNAG}).
\par

Gel'fand and Raikov showed in \cite{GeRa--43} that every locally compact group $G$ 
has enough irreducible unitary representations to separate the points of $G$,
triggering a systematic study of unitary representations for such groups.
First major results on classifications of unitary duals did concern
the affine group $\Aff(\R)$ of $\R$ \cite{GeNa--47a}, 
the special linear group $\SL(2,\C)$ \cite{GeNa--47b},
and $\SL(2,\R)$ \cite{Barg--47}.
Further developments include Mackey's analysis of induced representations 
for extensions of locally compact groups \cite{Mack--52, Mack--58},
Harish--Chandra's monumental work on the unitary dual of reductive Lie groups
\cite{Hari--Col},
and Kirillov's orbit method for nilpotent Lie groups \cite{Kiri--62}.
For a history of the subject of unitary group representations
see \cite{Mack--76}.

\subsection*{Operator algebras and primitive duals}

As already mentioned above,
the representation theory of a finite group $G$
amounts to the study of modules over the group algebra $\C[G]$,
which is a semisimple ring by Maschke's theorem.
This approach goes back to the first works on the subject
by Frobenius, Molien, and Burnside, in the last years of 19th century;
see \cite{Hawk--71, Hawk--72, Hawk--74}.
%
\iffalse
The importance for representations of ideals of group algebras
was later emphasized by E.~Noether;
see her lecture notes \cite{Noet--29},
and the discussion in \cite{Roqu--06}.
\fi
%
\par

Gel'fand developed in \cite{Gelf--41} a general theory of Banach algebras
and, together with Naimark, laid down the foundations for the theory of C*-algebras \cite{GeNa--43}.
As already mentioned above when discussing the primitive dual,
the unitary representation theory of a locally compact group $G$
can be viewed as the representation theory of a specific C*-algebra,
namely the maximal C*-algebra $C^*_{\rm max}(G)$ of $G$.
\par

It was soon realized by Godement, Mautner, Segal, and others
that the theory of von Neumann algebras,
which was developed by Murray and von Neumann
in a series of papers published in the 30s
(see \cite[Volume III, Rings of operators]{vNeu--Col}),
is of fundamental importance for the study of unitary representations,
in particular in questions concerning existence and uniqueness of decomposition
of such representations as direct integrals of irreducible unitary representations.
\par

Given a separable C*-algebra $A$, its dual or spectrum $\widehat A$ is defined
as the set of irreducible representations of $A$ in complex Hilbert spaces,
modulo unitary equivalence.
As mentioned above, Glimm showed that $\widehat A$,
equipped with a natural Borel structure,
is countably separated exactly when $A$ is of type I \cite{Glim--61a}.
\par

For an algebra $A$, Jacobson introduced in \cite{Jaco--45}
the set $\Pri(A)$ of primitive ideals of $A$ as a dual space of $A$
and endowed it with a topology which nowadays bears his name.
For a C*-algebra $A$, the topological space $\Pri(A)$ is a quotient of $\widehat A$.
Its topology was studied by Fell \cite{Fell--60a}.
In many situations, $\Pri(A)$ has much better topological properties than $\widehat A$.
In particular, for a locally compact group~$G$,
Fell gives a description of the topology of $\Pri(G) = \Pri(C^*_{\rm max}(G))$
in terms of weak containment of irreducible representations of $G$.
The relevance of $\Pri(A)$ for the description of a type I C*-algebra $A$
as algebra of sections of a bundle of elementary algebras
was demonstrated by Kaplansky \cite{Kapl--51b}.
This study was further pursued in Fell~\cite{Fell--61}.
As already mentioned, Effros \cite{Effr--63} proved the important result
that $\Pri(A)$ is always a standard Borel space when $A$ is separable.
All these facts demonstrate the importance
of the study of $\Pri(G)$ for a locally compact group $G$.

\subsection*{Normal quasi-dual space and characters}

Characters of finite groups were considered by Frobenius
even before he defined linear representations of such groups
(see Vol.\ III of \cite{Frob--Col}).
A famous formula for the characters of irreducible representations of compact Lie groups
was given by Weyl (see Vol.\ II of \cite{Weyl--Col}).
Gel'fand and Naimark described the characters of irreducible unitary representations 
of $G$ as conjugation-invariant distributions,
for $G = \SL_2(\C)$ in \cite{GeNa--47b} and 
other classical groups in \cite{GeNa--57}.
A deep and far-reaching study of characters of irreducible unitary representations 
of a general semisimple Lie group was undertaken by Harish--Chandra~\cite{Hari--56}.
\par

Godement introduced and studied a notion of
normal factor representations of a unimodular locally compact group in~\cite{Gode--54}.
This notion was extended and formulated in the context of C*-algebras
(or even involutive Banach algebras) by Guichardet in \cite{Guic--63}.
The normal quasi-dual space $\QD(A)_{\rm norm}$ of a C*-algebra $A$.
is defined as the set of normal factor representations of $A$,
modulo quasi-equivalence of representations.
Answering of question of Dixmier \cite[7.5.4]{Dixm--C*}),
Halpern showed in \cite{Halp--75} that the normal quasi-dual space $\QD(A)_{\rm norm}$
is a standard Borel space for any separable C*-algebra $A$.
\par

Let $G$ be a connected Lie group.
As already mentioned mentioned above when discussing the normal dual,
Puk\'anszky \cite{Puka--74} showed that the natural map
from $\QD(G)_{\rm norm} = \QD(C^*_{\rm max}(G))_{\rm norm}$
to $\Pri(G)$ is a bijection.
A detailed analysis of $\QD(G)_{\rm norm}$ was performed by Guichardet \cite{Guic--63}
for the Baumslag--Solitar group $G = \BS(2,1)$. 
\par

Vershik and Kerov gave in \cite{VeKe--81} a description of $\QD(G)_{\rm norm}$ 
for the infinite symmetric group $G = \Sym_{\rm fin}(\N)$.
\par

The subset $\QD(G)_{\rm fin}$ of $\QD(G)_{\rm norm}$
consisting of factor representations of finite type was studied by Godement \cite{Gode--51b}
for a locally compact unimodular group~$G$.
He showed in particular that the family of such representations
is point separating if and only if $G$ is a so-called SIN-group
(that is, if $G$ has arbitrary small conjugation-invariant neighbourhoods).
\par
 
Thoma's dual
(that is, the set $E(G)$ of characters of representations from $\QD(G)_{\rm fin}$)
was systematically studied by Thoma in the case of a discrete group $G$
and explicitly determined for some concrete examples
\cite{Thom--64a, Thom--64b, Thom--64c}).

\subsection*{Other dual spaces} 

Given a locally compact groups $G$, one may consider other dual spaces of $G$,
besides the ones discussed above, as for instance
the space $\Max (C^*_{\rm max}(G))$ of maximal two-sided ideals
of $C^*_{\rm max}(G)$
and the space $\Pri_*(L^1(G))$ of kernels of irreducible $*$-representations of $L^1(G)$.
There is an injection
$$
\Max(C^*_{\rm max}(G)) \hookrightarrow \Pri(G)
$$
and a surjective canonical map
$$
\Pri(G) \twoheadrightarrow \Pri_*(L^1(G)) .
$$
In general, these maps are not bijective (see Section~\ref{SectionVariants})
and so $\Max (C^*_{\rm max}(G))$ and $\Pri_*(L^1(G))$
provide only partial information on $\Pri(G)$.
\par

There is no duality involved in connection with the dual spaces we introduced so far, 
apart from the Pontrjagin--van Kampen duality in case of an abelian group~$G$.
So, the use of the word ``dual space" is questionnable;
it would be more accurate to speak of ``structure space"
when we refer to $\widehat G$ or $\Pri(G)$.
Nevertheless, we will stick to the more traditional ``dual space".
\par
 
There are genuine duality theories for locally compact groups,
such as the Tannaka--Krein duality for compact groups
(see for instance \cite[\S~30]{HeRo--70}),
and various duality theories for general locally compact groups,
developed by Stinespring, Kac, Ernest, Enock and Schwartz, among others
(for an overview, see \cite{EnSc--92}).
These duality theories involve, at least in the non-compact case,
objects such as Hopf algebras which are quite different from our dual spaces.

\section*{Overview}
% Section 0.H Overview
\label{S-Overview}

The first aim of this book is to introduce and discuss
the dual spaces introduced above for a group $G$:
\begin{itemize}
\item
the unitary dual $\widehat G$;
\item
the space $\widehat G_{\rm fd}$ of equivalence classes
of finite dimensional irreducible representations of $G$;
\item
the primitive dual $\Pri(G)$;
\item
the normal quasi-dual $\QD(G)_{\rm norm}$;
\item
the space of characters $\Char(G)$;
\item
and Thoma's dual $E(G)$.
\end{itemize}
As much as possible, basic definitions are formulated for topological groups,
even if most of the interesting statements and examples
are stated for locally compact groups, indeed often for second-countable ones.
Systematically, the definitions will be illustrated by a family of basic examples,
which are all discrete groups:
\begin{itemize}
\item[$\bullet$]
the infinite dihedral group $D_\infty$;
% (Section \ref{SectionInfDiGroup}) = 3.A
\item[$\bullet$]
two-step nilpotent discrete groups,
in particular the Heisenberg groups $H(R)$ over a commutative ring $R$,
and more specifically $H(\K)$ over an infinite field $\K$,
and $H(\Z)$ over the integers $\Z$;
% (Section \ref{Section-IrrRepTwoStepNil} = 3.B
\item[$\bullet$]
the affine group $\Aff(\K)$ of a field $\K$;
% (Section \ref{Section-IrrRepAff} = 3.C
\item[$\bullet$]
the solvable Baumslag--Solitar group $\BS(1, p)$ for a prime $p$;
% (Section \ref{Section-IrrRepBS} = 3.D
\item[$\bullet$]
the lamplighter group, which is the wreath product $\Z \wr (\Z / 2 \Z)$;
% (Section \ref{Section-IrrRepLamplighter} = 3.E
\item[$\bullet$]
and the general linear group $\GL_n(\K)$ over a field $\K$.
% (Section \ref{Section-IrrRepGLN} = 3.F
\end{itemize}
\par

Among other things, these worked out examples
(except $D_\infty$, the only group of type I in the list)
will demonstrate that, in general:
\begin{itemize}
\item[$\bullet$]
a complete description of $\widehat G$ is hopeless or useless;
\item[$\bullet$]
a description of the spaces $\widehat G_{\rm fd}$, $\Pri(G)$ and $E(G)$ is often possible;
\item[$\bullet$]
the map $\kappa^{\rm d}_{\rm prim} \, \colon \widehat G \twoheadrightarrow \Pri(G)$
is not injective, indeed has large fibers;
\item[$\bullet$]
a full description of the space $\Char(G) \smallsetminus E(G)$
seems to be difficult in general.
\end{itemize}
\par

Along the way, we establish results on the dual spaces
of locally compact groups $G$ which can be written as semi-direct product
of a countable discrete subgroup $H$ and an open abelian normal subgroup $N$
(observe that, apart from $\GL_n(\K)$, every group appearing above is of this form).
Let $G = H \ltimes N$ be such a group. We will
\begin{itemize}
\item[$\bullet$]
show that the description of $\widehat G$ amounts to the description
of $G$-quasi-invariant probability measures on $\widehat N$
and of certain cocycles $\widehat N \times H \to \U(\Ki)$
with values in the unitary group of a Hilbert space $\Ki$,
a usually hopeless task
(Theorem~\ref{Theo-AllRepSemiDirect2});
\item[$\bullet$]
give a description of the space $\widehat G_{\rm fd}$,
when the subgroup $H$ is moreover abelian
(Theorem~\ref{Theo-FiniteDimRepSemiDirect});
\item[$\bullet$]
construct factor representations of $G$ of various types 
(Theorem~\ref{Theo-FacRepSemiDirect});
\item[$\bullet$]
construct normal factor representations of $G$ and determine their characters
(Theorem~\ref {Theo-CharactersSemiDirect}).
\end{itemize}

\vskip.2cm

We discuss in detail some of the basic tools needed for the proofs;
in particular, we will give an account on 
\begin{itemize}
\item
functions of positive type and the Gel'fand--Naimark--Segal construction
(Section \ref{S-FPosType}, \ref{SectionPrimC*}, and \ref{Section-GNS-Traces});
\item
induced representations of a locally compact group from an open subgroup,
and the Mackey--Shoda irreducibility and equivalence criteria
(Section~\ref{Section-IrrIndRep});
\item
unitary representations of a locally compact abelian group,
in terms of multip\-licity-free representations,
and also in terms of projection-valued measures
(Chapter~\ref{Chapter-AbelianGroups});
\item
von Neumann algebras associated to a representation (Chapter \ref{ChapterTypeI});
\item
C*-algebras associated to a locally compact group (Chapter \ref{ChapterAlgLCgroup});
\item
traces on C*-algebras and the theory of Hilbert algebras (Chapter~\ref{ChapterCharacters});
\item
the group measure space construction of Murray and von Neumann
(Chapter~\ref{ChapterGroupMeasureSpace}).
\end{itemize}
Various other reminders on topology, measure theory and abelian harmonic analysis
are collected in an appendix to this book.

\vskip.2cm

Some comments are in order concerning operator algebras.
On the one hand, the theory of unitary representations of groups 
and that of operator algebras cannot be separated, 
neither conceptually nor historically.
On the other hand, we would find frustrating
to go through all the necessary prerequisites of operator algebras
before addressing group representations.
We have tried to propose a helpful selection of ``reminders''
on operator algebras, scattered through the chapters on group representations.
For what is not explained here in detail,
most citations refer to Dixmier's books \cite{Dixm--vN, Dixm--C*}.
For one shorter exposition, we can also recommend \cite[Chap.~14]{Wall--92}.
\begin{center}
\emph{It is the occasion to celebrate Jacques Dixmier 
\\
and the writing of his two books on operator algebras.
\\
More than 50 years after publication,
\\
they remain unrivalled references on the subject.}
\end{center}

%-----------------------------------------------------------------------
% End of preface
%-----------------------------------------------------------------------
\include{contents}

% \mainmatter
% hidden on December 7

%    Include main chapters here.
\chapter{Unitary dual and primitive dual}
% Chapter 1
\label{ChapterUnitaryDual}

\emph{
In Section \ref{S-DefUnitD},
we introduce unitary group representations and unitary duals.
The unitary dual $\widehat G$ of a topological group $G$
is the set of equivalence classes of its irreducible unitary representations;
we write simply ``dual'' for ``unitary dual''.
}
\par

\emph{
A unitary representation $\pi$ of a topological group $G$ is cyclic
if there exists a vector $\xi$ in the Hilbert space of $\pi$
such that $\pi(G)\xi$ generates the Hilbert space of~$\pi$.
There is a correspondence between equivalence classes
of pairs $(\pi, \xi)$, where $\pi$ is a unitary representation with cyclic vector $\xi$,
and a class of functions on $G$, called functions of positive type.
Under this correspondence, irreducible representations
are mapped into indecomposable functions of positive type and vice-versa.
This is the subject of Section \ref{S-FPosType} on the Gel'fand--Naimark--Segal construction.
}
\par

\emph{
In Section \ref{SectionWC+FellTop},
we introduce the so-called Fell topology on the dual $\widehat G$ of a topological group $G$
and show how it is related to the notion of weak containment between unitary representations. 
Section~\ref{SectionFellPropers} is devoted to
the study of the interplay between topological properties of $G$
and topological properties of its dual $\widehat G$, for a locally compact group $G$.
}
\par

\emph{
It turns out that, in general, the topological space $\widehat G$ does not even have
the weakest of all separation properties, the T$_0$ property. 
In Section \ref{PrimIdealSpace}, we introduce our second dual space,
the primitive dual $\Pri(G)$ of a topological group $G$,
as the set of weak equivalence classes of irreducible representations of $G$.
In a sense, $\Pri(G)$ is the largest T$_0$ quotient of $\widehat G$.
We will give later (in Chapter \ref{ChapterAlgLCgroup}) an alternative
description of $\Pri(G)$ in case $G$ is a locally compact group.
 }
\par

\emph{One of the main procedure for the construction
of unitary representations of a group $G$ is induction,
a way to build representations of $G$ out of representations of a subgroup $H$. 
In Section \ref{Section-IrrIndRep},
we give an account of induced representations in case the subgroup $H$ is open.
In particular, a treatment of the Mackey-Shoda criteria of irreducibility
and equivalence of induced representations is presented.
}
\par

\emph{
In Section \ref{SectionDecomposingIrreps},
we give a short introduction to the notions of
direct integrals of Hilbert spaces, direct integrals of unitary representations,
and direct integrals of the associated von Neumann algebras.
We discuss how every unitary representation
of a second-countable locally compact group is equivalent to
a direct integral of \emph{irreducible} representations.
Moreover, we show that such direct integral decompositions are not unique in general.
}
\par

\emph{
As a preparation for later chapters,
we give a characterization of the decomposable operators
between two direct integrals of Hilbert spaces
in the final Section~\ref{Section-DecomposableOperators}.
}

\section
{Definition of the unitary dual}
% Section 1.A
\label{S-DefUnitD}

\index{$h1$@$\Hi$ Hilbert space}
Let $\Hi$ be a Hilbert space;
in this book, Hilbert spaces are always assumed
to be \emph{complex}.
We denote by $\Li (\Hi)$ 
the involutive algebra of bounded linear operators on $\Hi$
and by
$$
\U(\Hi) \, = \, \{ x \in \Li (\Hi) \mid x^*x = xx^* = \mathrm{Id}_{\Hi} \}
$$
the group of unitary operators on $\Hi$.
Sometimes, the field $\C$ of complex numbers is viewed as a Hilbert space.
\index{Involutive algebra $\Li (\Hi)$}
\index{$h4$@$\Li (\Hi)$ bounded linear operators on $\Hi$}
\index{Unitary group $\U(\Hi)$}
\index{$h6$@$\U(\Hi)$ unitary group of $\Hi$}
\index{$m1$@$\C$ complex numbers}
\par

On $\Li (\Hi)$, the \textbf{strong topology}
is the locally convex topology defined by the set of seminorms
$x \mapsto \Vert x \xi \Vert$, for $\xi \in \Hi$,
and the \textbf{weak topology}
is the locally convex topology defined by the set of seminorms
$x \mapsto \vert \langle x \xi \mid \eta \rangle \vert$, for $\xi, \eta \in \Hi$.
The restrictions of these two topologies to $\U(\Hi)$ coincide,
and make $\U(\Hi)$ a topological group.
More on this in Appendix \ref{AppHspacesop}.
% this is easy to check,
% and follows for example from \cite[Chap.~I, \S~3, no~1]{Dixm--vN}.
\index{Strong topology on $\Li(\Hi)$ and $\U(\Hi)$}
 \index{Weak topology on $\Li(\Hi)$ and $\U(\Hi)$}
\par

Let $G$ be a topological group.
The unit element in $G$ is denoted most often by~$e$,
exceptionally by $0$ or $1$.
\index{$b1$@$G$ topological group}
\index{$b2$@$e \in G$ unit element}

\begin{defn}
% 1.A.1
\label{defunitaryrep}
A \textbf{unitary representation} $\pi$ of $G$ in a Hilbert space $\Hi$
is a homomorphism $\pi \, \colon G \to \U(\Hi)$
such that the associated map
$G \times \Hi \to \Hi, \hskip.1cm (g, \xi) \mapsto \pi(g)\xi$, is continuous;
we often write $\Hi_\pi$ rather than $\Hi$ for the Hilbert space of $\pi$,
and $(\pi, \Hi)$ for $\pi$.
\end{defn}

\begin{center}
From now on in this book, all representations of groups are unitary
\\
unless explicitly written otherwise,
\\
and we write ``\textbf{representation}'' for ``unitary representation''.
\end{center}
\index{Unitary representation = representation}
\index{Representation = unitary group representation}
Note also that our definition of ``representation'' \emph{includes} the continuity requirement.

\begin{prop}
% 1.A.2
\label{repcontinuoussepjoint}
Let $G$ be a topological group, $\Hi$ a Hilbert space,
and $\pi \,\colon G \to \U(\Hi)$ a homomorphism.
The following two conditions are equivalent:
\begin{enumerate}[label=(\roman*)]
\item\label{iDErepcontinuoussepjoint}
the map $G \times \Hi \to \Hi, \hskip.2cm (g, \xi)\to \pi(g)\xi$ is continuous,
i.e., $\pi$ is a representation of $G$ in $\Hi$;
\item\label{iiDErepcontinuoussepjoint}
the map $\pi$ is continuous with respect to the strong topology on $\mathcal{U}(\Hi)$.
\end{enumerate}
\end{prop}

\begin{proof}
One implication is obvious and the other follows from the inequality
$$
\Vert \pi(g)\xi - \pi(h)\eta \Vert \, \le \, \Vert \pi(g)\xi - \pi(h)\xi \Vert + \Vert \pi(h)\xi - \pi(h)\eta \Vert
\, = \, \Vert \pi(g)\xi - \pi(h)\xi \Vert + \Vert \xi - \eta \Vert
$$
for $g,h \in G$ and $\xi, \eta \in \Hi$.
\end{proof}

\index{Character! $1$@unit $1_G$}
\index{Unit representation or unit character $1_G$}
\index{Representation! unit}
\index{$b6$@$1_G$ unit representation of $G$}
\index{Subrepresentation}
We write $1_G$ for the 
\textbf{unit representation} of $G$ in $\C$,
defined by $1_G(g) = \mathrm{Id}_\C$ for all $g \in G$.

A $\pi(G)$-invariant closed subspace $\Ki$ of $\Hi$ gives rise to
a \textbf{subrepresentation} $G \to \U(\Ki)$,
assigning to $g \in G$ the restriction of $\pi(g)$ to $\Ki$.
\par

\index{Intertwiner! $1$@between representations of a group}
\index{$h8$@$\Li (\Hi_1, \Hi_2)$ bounded linear operators from $\Hi_1$ to $\Hi_2$}
\index{$h4$@$\Li (\Hi)$ bounded linear operators on $\Hi$}
\index{Commutant}
Let $(\pi_1, \Hi_1), (\pi_2, \Hi_2)$ be two representations of $G$.
A bounded linear operator $x \, \colon \Hi_1 \to \Hi_2$
is an \textbf{intertwiner} between $\pi_1$ and $\pi_2$
if $x \pi_1(g) = \pi_2(g) x$ for all $g \in G$.
Observe that the set $\Hom_G(\pi_1, \pi_2)$
of intertwiners between $\pi_1$ and $\pi_2$
is a linear subspace of the space $\Li (\Hi_1, \Hi_2)$
of bounded linear operators from $\Hi_1$ to $\Hi_2$.
For one representation $(\pi, \Hi)$,
the set $\Hom_G(\pi, \pi)$ is the commutant $\pi(G)'$;
it always contains the set $\C \mathrm{Id}_\Hi$
of scalar multiples of the identity of $\Hi$. 
More generally, the \textbf{commutant}
of a subset $S$ of $\Li (\Hi)$ is the subset
$$
S' \, = \, \{x \in \Li (\Hi) \mid xs = sx
\hskip.2cm \text{for all} \hskip.2cm 
s \in S \} .
$$

\begin{defn}
% 1.A.3
\label{Def-equivalentdisjoint}
Two representations $\pi_1$ and $\pi_2$ of a topological group $G$ are \textbf{equivalent},
and this is written $\pi_1 \simeq \pi_2$,
if there exists a unitary operator $V \, \colon \Hi_1 \to \Hi_2$
such that $\pi_2(g) = V \pi_1(g) V^*$ for all $g \in G$.
% \par
% 
The representations $\pi_1$ and $\pi_2$ are \textbf{disjoint}
% and this is written $\pi_1 \downspoon \pi_2$,
% iciici n'accepte pas $\pi_1 \mackey \pi_2$,
if there does not exist a subrepresentation of~$\pi_1$
equivalent to a subrepresentation of $\pi_2$.
% \par
%
(We come back to these definitions in Chapter \ref{ChapterTypeI}.)
\index{Equivalent! $1$@representations}
\index{$a1$@$\simeq$ equivalence of representations!}
\index{Representation! equivalent $\simeq$}
\index{Representation! disjoint}
\index{Disjoint! representations}
\end{defn}

\index{Direct sum of representations}
\index{Representation! direct sum}
The \textbf{direct sum} of a family $(\pi_\iota)_{\iota \in I}$ of representations of $G$
is the representation $\pi = \bigoplus_{\iota \in I} \pi_\iota$ of $G$
in the Hilbert direct sum $\bigoplus_{\iota \in I} \Hi_{\pi_\iota}$
defined by $\pi(g) = \bigoplus_{\iota \in I} \pi_\iota(g)$ for all $ g \in G$.

\begin{prop}
% 1.A.4
Let $\pi$ be a representation of a topological group $G$
such that $\Hi_\pi$ contains a $\pi(G)$-invariant closed subspace~$\Ki$.
\par

Then the orthogonal $\Ki^\perp$ is also a $\pi(G)$-invariant closed subspace of $\Hi$.
It follows that $\pi$ is the direct sum of the two subrepresentations of $G$
in the subspaces $\Ki$ and $\Ki^\perp$.
\end{prop}

\begin{proof}
Let $\xi \in \Ki^\perp$. For all $\eta \in \Ki$, we have
$$
\langle \pi(g) \xi \mid \eta \rangle \, = \, \langle \xi \mid \pi(g^{-1}) \eta \rangle
\, \in \, \langle \xi \mid \Ki \rangle \, = \, \{0\} ,
$$
hence $\pi(g) \xi \in \Ki^\perp$.
This shows that $\Ki^\perp$ is $\pi(G)$-invariant.
\end{proof}

\begin{lem}
% 1.A.5
\label{Prop-EquSubRep}
Let $(\pi_1, \Hi_1)$, $(\pi_2, \Hi_2)$ be two representations of a topological group $G$,
and $T \in \Hom_G(\pi_1, \pi_2)$.
Set $\Ki_1 = (\ker T)^{\perp}$ and let $\Ki_2$ denote the closure of the image of $T$. 
\par

Then $\Ki_1$ and $\Ki_2$ are closed invariant subspaces
of $\Hi_1$ and $\Hi_2$ respectively,
and the subrepresentation of $\pi_1$ defined by $\Ki_1$
is equivalent to the subrepresentation of $\pi_2$ defined by $\Ki_2$. 
\end{lem}

\begin{proof}
Since $T$ intertwines $\pi_1$ and $\pi_2$,
the space $\Ki_1$ is $\pi_1(G)$-invariant,
the space $T(\Ki_1)$ is $\pi_2(G)$-invariant,
and therefore $\Ki_2$ is $\pi_2(G)$-invariant.
\par

One checks immediately that 
$T^* \in \Li (\Hi_2, \Hi_1)$ intertwines $\pi_2$ and $\pi_1$,
and that $T^*T \in \Li (\Hi_1)$ intertwines $\pi_1$ with itself.
Since $\vert T \vert = (T^*T)^{1/2}$
is a limit in the strong operator topology of polynomials in $T^*T$
(see Section \ref{AppProjValMeas}),
$\vert T \vert$ also intertwines $\pi_1$ with itself.
\par

Consider the partial isometry $U$ of the polar decomposition $T = U \vert T \vert$;
recall that $\ker U = \ker T$, and that the restriction of $U$ to 
$\Ki_1$ is an isometry onto $\Ki_2$.
By definition of $U$
(i.e., $U\xi = 0$ if $\xi \in \ker T$ and $U \vert T \vert \xi = T \xi$ if $\xi \in (\ker T)^\perp$,
see Appendix \ref{AppProjValMeas}),
it follows that $U$ intertwines the restriction of $\pi_1$ to $\Ki_1$
and the restriction of $\pi_2$ to $\Ki_2$.
Therefore these two restrictions are equivalent.
\end{proof}

In the following proposition, the implications
\ref{iiDEProp-Prop-EquSubRep} $\Rightarrow$ \ref{iDEProp-Prop-EquSubRep}
and
\ref{iiiDEProp-Prop-EquSubRep} $\Rightarrow$ \ref{ivDEProp-Prop-EquSubRep}
are immediate consequences of Lemma~\ref{Prop-EquSubRep},
and the converse implications are trivial.

\begin{prop}
% 1.A.6
\label{Prop-Prop-EquSubRep}
Let $(\pi_1, \Hi_1)$ and $(\pi_2, \Hi_2)$ be two representations of a topological group $G$.
On the one hand, the following two properties are equivalent:
\begin{enumerate}[label=(\roman*)]
\item\label{iDEProp-Prop-EquSubRep}
$\pi_1$ and $\pi_2$ are equivalent;
\item\label{iiDEProp-Prop-EquSubRep}
$\Hom_G(\pi_1, \pi_2)$ contains an invertible operator from $\Hi_1$ to $\Hi_2$.
\end{enumerate}
On the other hand, the following two properties are equivalent:
\begin{enumerate}[label=(\roman*)]
\addtocounter{enumi}{2}
\item\label{iiiDEProp-Prop-EquSubRep}
$\pi_1$ and $\pi_2$ are disjoint;
\item\label{ivDEProp-Prop-EquSubRep}
$\Hom_G(\pi_1, \pi_2) = \{0\}$.
\end{enumerate}
\end{prop}

\begin{cor}
% 1.A.7
\label{Cor-Prop-EquSubRep}
Let $\pi$, and $\pi_i$ for $i \in I$, be representations of $G$.
Assume that $\pi$ and $\pi_i$ are disjoint for every $i \in I$.
\par

Then $\pi$ and $\bigoplus_{i \in I}\pi_i$ are disjoint.
\end{cor}

\begin{proof}
In view of Proposition~\ref{Prop-Prop-EquSubRep},
it suffices to show that $\Hom_G(\pi, \bigoplus_{i \in I}\pi_i) = \{0\}$.
\par

Let $\Hi_i$ be the Hilbert space of $\pi_i$ and $\Hi = \bigoplus_{i \in I} \Hi_i$.
Denote by $P_i$ the orthogonal projection from $\Hi$ to $\Hi_i$.
Let $T \in \Hom_G(\pi, \bigoplus_{i \in I}\pi_i)$. 
Then $P_iT \in \Hom_G(\pi, \pi_i)$ and hence $P_iT = 0$ for every $i \in I$.
Therefore, $T = 0$.
\end{proof}

The following proposition shows that
the equivalence of representations satisfies a form of the
\textbf{Schr\"oder--Bernstein property}.

\begin{prop}
% 1.A.8
\label{SchroderBernstein}
Let $\pi_1, \pi_2$ be two representations of $G$
in Hilbert spaces $\Hi, \Ki$ respectively.
Assume that $\pi_1$ is equivalent to a subrepresentation of $\pi_2$
and $\pi_2$ is equivalent to a subrepresentation of $\pi_1$.
\par

Then $\pi_1$ and $\pi_2$ are equivalent.
\end{prop}

\begin{proof}
By hypothesis, there exist
\begin{enumerate}
\item[--]
a $G$-invariant closed subspace $\Ki'$ of $\Ki$
and a $G$-equivariant partial isometry $S$ of initial space $\Hi$ and final space $\Ki'$,
\item[--]
a $G$-invariant closed subspace $\Hi'$ of $\Hi$
and a $G$-equivariant partial isometry $T$ of initial space $\Ki$ and final space $\Hi'$.
\end{enumerate}
Set $U = TS$ and $\Hi'' = U(\Hi)$. Then
$$
\Hi \, \supset \, \Hi' \, \supset \, \Hi''
\leqno{(\sharp)} 
$$
and $U$ is a $G$-equivariant partial isometry of initial space $\Hi$
and final space $\Hi''$. For all $m \ge 0$, set
$$
\Hi^{(2m)} \, = \, U^m(\Hi)
\hskip.5cm \text{and} \hskip.5cm
\Hi^{(2m+1)} \, = \, U^m(\Hi') .
$$
Observe that $\Hi^{(0)} = \Hi$, $\Hi^{(1)} = \Hi'$, and $\Hi^{(2)} = \Hi''$.
It follows from ($\sharp$) that
$$
\Hi^{(2m)} \, \supset \, \Hi^{(2m+1)} \, \supset \, \Hi^{(2m+2)}
\hskip.5cm \text{for all} \hskip.2cm
m \ge 0 .
$$
Set $\mathcal J = \bigcap_{n \ge 0} \Hi^{(n)}$.
We have orthogonal Hilbert sums
$$
\begin{aligned}
\Hi \, &= \, \mathcal J 
\oplus (\Hi^{(0)} \ominus \Hi^{(1)})
\oplus (\Hi^{(1)} \ominus \Hi^{(2)})
\oplus (\Hi^{(2)} \ominus \Hi^{(3)})
\oplus (\Hi^{(3)} \ominus \Hi^{(4)})
\oplus \cdots
\\
\Hi' \, &= \, \mathcal J 
\oplus (\Hi^{(1)} \ominus \Hi^{(2)})
\oplus (\Hi^{(2)} \ominus \Hi^{(3)})
\oplus (\Hi^{(3)} \ominus \Hi^{(4)})
\oplus (\Hi^{(4)} \ominus \Hi^{(5)})
\oplus \cdots ,
\end{aligned}
$$?
where $\Hi^{(i)} \ominus \Hi^{(j)}$ denotes the orthogonal complement of 
$\Hi^{(j)}$ in $\Hi^{(i)}$.
Since $U$ induces for all $m \ge 0$ a map
$(\Hi^{(2m)} \ominus \Hi^{(2m+1)}) \to (\Hi^{(2m+2)} \ominus \Hi^{(2m+3)})$
which is a $G$-equivariant surjective isometry, we have also
$$
\Hi' \, \simeq \, \mathcal J 
\oplus (\Hi^{(1)} \ominus \Hi^{(2)})
\oplus (\Hi^{(0)} \ominus \Hi^{(1)})
\oplus (\Hi^{(3)} \ominus \Hi^{(4)})
\oplus (\Hi^{(2)} \ominus \Hi^{(3)})
\oplus \cdots .
$$
It follows that the representations of $G$ in $\Hi$ and $\Hi'$ are equivalent.
Since the representations of $G$ in $\Hi'$ and $\Ki$ are equivalent
by hypothesis, this concludes the proof.
% [Modulo a mild translation, this is the proof 
% of Proposition 1 in \S~1 of Chapter III of \cite{Dixm--vN}.]
\end{proof}

\begin{defn}
% 1.A.9
\label{Def-IrredRep}
A representation $\pi$ of the topological $G$ is \textbf{irreducible}
if its Hilbert space $\Hi_\pi$ is not $\{0\}$
and if the only $\pi(G)$-invariant closed subspaces of $\Hi_\pi$
are $\{0\}$ and $\Hi_\pi$.
\index{Irreducible representation}
\index{Representation! irreducible}
\end{defn}

\begin{rem}
% 1.A.10
\label{Rem-IrredRep}
(1)
Irreducibility of a representation is a property of its equivalence class:
if $\pi$ is an irreducible representation of $G$,
then every representation which is equivalent to $\pi$ is irreducible.

\vskip.2cm

(2)
Representations of $G$ of dimension $1$ are irreducible.
They are equivalent to \textbf{unitary characters} of $G$,
that is to continuous homomorphisms $G \to \T$.
\par

Here and below, $\R$ denotes the field of real numbers,
and $\T$ the circle group,
viewed sometimes as $\{ z \in \C \mid \vert z \vert = 1\}$,
and sometimes as $\R / \Z$.
\index{Unitary character}
\index{Character! $2$@unitary character}
\index{$l 3$@$\R$ real numbers}
\index{$m2$@$\T$ circle group}
\end{rem}

In the following proposition,
Claim \ref{iDESchurEtc} is known as Schur's lemma.
For the easy proof, we refer to \cite[Appendix A]{BeHV--08}.
\index{Schur's lemma}

\begin{prop}
% 1.A.11
\label{SchurEtc}
Let $\pi, \pi_1, \pi_2$ be representations of a topological group $G$,
\begin{enumerate}[label=(\arabic*)]
\item\label{iDESchurEtc}
The representations $\pi$ is irreducible if and only if
$\Hom_G(\pi, \pi) = \C \mathrm{Id}_{\Hi_\pi}$.
\item\label{iiDESchurEtc}
Suppose that $\pi_1$ and $\pi_2$ are irreducible;
then:
\par

\hskip.5cm
either $\pi_1$ and $\pi_2$ are equivalent and 
$\dim_\C \Hom_G(\pi_1, \pi_2) = 1$,
\par

\hskip.5cm
or $\pi_1$ and $\pi_2$ are not equivalent and
$\Hom_G(\pi_1, \pi_2) = \{0\}$.
\end{enumerate} 
\end{prop}

The next proposition shows that the dimensions of the Hilbert spaces
in which irreducible representations of $G$ exist are bounded by some cardinal.
There is a reminder on the dimensions of Hilbert spaces in Appendix \ref{AppHspacesop}.

\begin{prop}
% 1.A.12
\label{Pro-IrrRepBoundedCard}
Let $G$ be a topological group.
Let $\aleph$ be a cardinal such that
there exists in $G$ a dense subset, say $X$, of cardinality at most $\aleph$.
\par

Then, for every irreducible representation $(\pi, \Hi)$ of $G$,
the dimension of $\Hi$ is bounded by $\aleph$.
\par

In particular, if $G$ is separable, then $\Hi$ is a separable Hilbert space.
\end{prop}

\begin{proof}
Let $\xi$ be a non-zero vector in $\Hi$ and let $\Ki$ be the closed
subspace of $\Hi$ generated by $\pi(G)\xi$. 
Since $X$ is dense in $G$ and since $X$ has cardinality $\aleph$, it is clear that 
$\dim \Ki \le \aleph$.
Since $\Ki$ is $G$-invariant and $\pi$ is irreducible, we have $\Ki = \Hi$.
\end{proof}

The following consequence of Proposition~\ref{Pro-IrrRepBoundedCard}
is crucial for the subject of this book.

\begin{cor}
% 1.A.13
The equivalence classes of irreducible representations
of a topological group $G$ constitute a set, say $\widehat G$.
\end{cor}

\begin{proof}
Let $\kappa$ be a cardinal which satisfies the condition
stated in Proposition~\ref{Pro-IrrRepBoundedCard}
and let $\Hi$ be a fixed Hilbert space of dimension $\kappa$.
The family $\mathcal R$ of irreducible representations of $G$ in closed linear subspaces of $\Hi$ is a set.
Proposition~\ref{Pro-IrrRepBoundedCard} shows that every irreducible representation of $G$
is equivalent to some representation of $\mathcal R$.
So, $\widehat{G}$ is the quotient space $\mathcal R / \simeq$ of $\mathcal R$ 
by the relation of equivalence of representations.
\end{proof}

\begin{defn}
% 1.A.14
\label{Def-UnitaryDual}
The \textbf{unitary dual} $\widehat G$ of $G$ is the space
of equivalence classes of irreducible representations of $G$.
\par

Most often below, we use simply ``\textbf{dual}'' for ``unitary dual''.
\par

When $G$ is abelian, its dual $\widehat G$
can be identified with the group $\Hom (G, \T)$ of continuous homomorphisms $G \to \T$.
\index{Unitary dual}
\index{$b7$@$\widehat G$ (unitary) dual of the topological group $G$}
\index{Dual! $1$@short for Unitary dual}
\end{defn}

The dual $\widehat G$ is a topological space for the Fell topology, 
of which we postpone a first definition until Section \ref{SectionWC+FellTop}.
For $G$ locally compact, two equivalent definitions of the Fell topology
are discussed in Sections \ref{SectionQuasidual} and \ref{C*algLCgroup}.
\par

Next, we examine the action of automorphisms of $G$ on the representation theory of $G$.
For actions, see Definition \ref{appendixaction}.

\begin{notat}
% 1.A.15
\label{notationaction}
\index{Action of $G$ on $X$, respectively $G \curvearrowright X$ and $X \curvearrowleft G$}
When we wish to emphasize that
an action of a group $G$ on a set $X$ is a left action $G \times X \to X$,
we refer to this action by the notation
$$
G \curvearrowright X .
$$
Similarly, for a right action $X \times G \to X$, we use
$$
X \curvearrowleft G.
$$
\end{notat}

\begin{rem}
% 1.A.16
\label{Rem-AutoRep}
Let $\Aut(G)$ be the group of bicontinuous automorphisms of $G$.
We consider $\Aut(G)$ as acting on $G$ on the left.

\vskip.2cm

\index{Conjugate representation} 
\index{Representation! conjugate}
(1)
Let $\theta \in \Aut(G)$. For a representation $(\pi, \Hi)$ of $G$,
the \textbf{conjugate representation} of $\pi$ by $\theta$ 
is the representation $\pi^\theta$ of $G$ which is defined on the same Hilbert space $\Hi$ by 
$$
\pi^\theta(g) \, = \, \pi(\theta(g))
\hskip.5cm \text{for} \hskip.2cm
g \in G.
$$
\par

Let $\aleph$ be a cardinal
(for example, the cardinal $\aleph_0$ of the set of integers).
For each cardinal $n \le \aleph$, let $\Hi_n$ be a Hilbert space of dimension $n$.
Denote by ${\rm Rep}_\aleph (G)$ the set of representations of $G$
in one of the Hilbert spaces $\Hi_n$, with $n \le \aleph$.
There is a natural \textbf{right} action
${\rm Rep}_\aleph (G) \times \Aut(G) \to {\rm Rep}_\aleph (G)$
defined by
$$
(\pi, \theta) \mapsto \pi^\theta
\hskip.5cm \text{for} \hskip.2cm
\pi \in {\rm Rep}_\aleph (G)
\hskip.2cm \text{and} \hskip.2cm
\theta \in \Aut(G).
$$
The subset ${\rm Irr}_\aleph (G)$ of ${\rm Rep}_\aleph (G)$
consisting of irreducible representations is invariant by $\Aut(G)$,
hence the action just defined restricts to a right action
of $\Aut(G)$ on ${\rm Irr}_\aleph (G)$.

\vskip.2cm

(2)
Let $\aleph$ be a cardinal such that
every irreducible representation of $G$ is of dimension at most $\aleph$
(see Proposition \ref{Pro-IrrRepBoundedCard}).
The dual $\widehat G$ can be identified with
the quotient of ${\rm Irr}_\aleph (G)$ be the relation of equivalence of representations,
i.e., $\widehat G = {\rm Irr}_\aleph (G) / {\simeq}$.
The action described in (1) induces a right action
$\widehat G \curvearrowleft \Aut(G)$ of $\Aut(G)$ on $\widehat G$.

\vskip.2cm

(3)
Let $\pi$ a representation of $G$ and $a \in G$. 
By a slight abuse of notation, we denote by $\pi^a$
the conjugate representation of $\pi$ by the inner automorphism defined by~$a$.
Observe that $\pi^a \, \colon g \mapsto \pi(aga^{-1})$ is equivalent to $\pi$;
indeed, the unitary operator $\pi(a)$ intertwines $\pi^a$ and $\pi$.

\vskip.2cm

(4)
We extend as follows the situation considered in Item (3).
Let $N$ be a closed normal subgroup of $G$ and 
$\rho$ a representation of $N$.
Every element $a \in G$ defines an automorphism of $N$;
the conjugate representation $\rho^a$ of $\rho$ by $a$
is the conjugate representation of $\rho$ by this automorphism,
so that
$$
\rho^a(n) \, = \, \rho( an a^{-1}),
\hskip.5cm \text{for all} \hskip.2cm
n \in N .
$$
\par

The map 
$$
\widehat N \times G \to \widehat N, 
\hskip.5cm
(\rho,a) \mapsto \rho^a 
$$
defines a \emph{right} action of $G$ on the dual space $\widehat N$ of $N$.
Since $\rho^a$ is equivalent to $\rho$ for $a \in N$,
this action factorizes to an action $\widehat N \curvearrowleft G/N$ of the quotient group $G/N$.
\end{rem}

\section[Functions of positive type and GNS construction]
{Functions of positive type and GNS construction}
% Section 1.B
\label{S-FPosType}

The material of this section can be found in many places,
for example in \cite[\S~13]{Dixm--C*} and \cite[Appendix~C]{BeHV--08}.
\par

Let $G$ be a topological group.

\begin{defn}
% 1.B.1
\label{Def-FPT}
A \textbf{function of positive type} on $G$
is a continuous function $\varphi \, \colon G \to \C$ such that
$$
\sum_{i, j = 1}^n c_i \overline{c_j} \varphi(g_j^{-1}g_i) \, \ge \, 0
\hskip.5cm \text{for all} \hskip.2cm 
n \ge 1, \hskip.2cm
c_1, \hdots, c_n \in \C ,
\hskip.2cm \text{and} \hskip.2cm
g_1, \hdots, g_n \in G.
$$
Many authors use ``\textbf{positive definite function}'' instead of
``function of positive type''.
\index{Positive type, for a function}
\index{Positive definite function, see functions of positive type!}
\index{Function of positive type}
\end{defn}

\begin{rem}
% 1.B.2
\label{Rem-FPT}
(1)
Let $\varphi$ be such a function of positive type on $G$.
For $n = 2$, $g_1 = e$, and $g_2 = g$, the definition implies 
that the matrix
% iciici petite matrice iciici
$\bigl(\begin{smallmatrix}
\varphi(e) & \varphi(g^{-1}) \\ \varphi(g) & \varphi(e)
\end{smallmatrix} \bigr)$
is positive Hermitian; it follows that
$$
\varphi(g^{-1}) \, = \, \overline{\varphi(g)}
\hskip.5cm \text{and} \hskip.5cm
\vert \varphi (g) \vert \le \varphi (e)
\hskip.5cm \text{for all} \hskip.2cm 
g \in G .
$$
Note that the functions $\check{\varphi}$ and $\Re{\varphi}$,
defined by $\check{\varphi}(g) = \varphi(g^{-1})$
and $(\Re(\varphi))(g) = \frac{1}{2}(\varphi(g) + \check{\varphi}(g))$
for all $g \in G$,
are also of positive type.

\vskip.2cm

(2)
Functions of positive type have been first studied
in classical harmonic analysis on groups like $\T$ or $\R$.
For an historical survey, we refer to \cite{Stew--76}.
\end{rem}

\index{$d1$@$P(G), P_{\le 1}(G), P_1(G)$ functions of positive type on $G$}
\index{Algebras! $1$@continuous functions, $C(X)$, $C^b(X)$, $\C^c(X)$, $C^0(X)$}
Denote by $P(G)$ the set of continuous functions $\varphi \, \colon G \to \C$ 
which are of positive type.
For $\varphi, \psi \in P(G)$ and $c \in \R_+$,
it follows from the definition that
$\varphi + \psi$, $c\varphi$, and $\overline \varphi$ are of positive type.
[For $\varphi \psi$, see Example \ref{3exeGNS}(4).]
Therefore $P(G)$ is a convex cone in the space $C^b(G)$
of bounded continuous complex-valued functions on $G$,
and is closed under complex conjugation.
Let $P_{\le 1}(G)$ the subset of functions in $P(G)$
such that $\sup_{g \in G}\vert \varphi(g) \vert = \varphi(e) \le 1$,
and $P_1(G)$ the subset of $P_{\le 1}(G)$
of functions normalized by $\varphi(e) = 1$.
The subsets $P_{\le 1}(G)$ and $P_1(G)$ are also convex.

\begin{defn}
% 1.B.3
\label{exposfunctdiagcoeff}
\index{Matrix coefficient}
\index{Diagonal matrix coefficient}
Let $\pi$ be a representation of $G$ in a Hilbert space $\Hi_\pi$.
A \textbf{matrix coefficient}
of $\pi$ is a continuous function $\varphi_{\pi, \xi, \eta}$ on $G$ defined by
$$
\varphi_{\pi, \xi, \eta} (g) \, = \, \langle \pi(g) \xi \mid \eta \rangle 
\hskip.5cm \text{for all} \hskip.2cm
g \in G ,
$$
for some $\xi, \eta \in \Hi_\pi$.
A \textbf{diagonal matrix coefficient} 
is a function of the form $\varphi_{\pi, \xi, \xi}$;
for simplicity, we rather write $\varphi_{\pi, \xi}$.
\par

Diagonal matrix coefficients are of positive type.
Indeed, for any positive integer~$n$
and $c_1, \hdots, c_n \in \C$, $g_1, \hdots, g_n \in G$, $\xi \in \Hi_\pi$, we have
$$
\sum_{i, j = 1}^n c_i \overline{c_j} \varphi_{\pi, \xi}(g_j^{-1}g_i) 
\, = \,
\Big\langle \sum_{i = 1}^n c_i \pi(g_i) \xi
\hskip.1cm \Big\vert \hskip.1cm
\sum_{j=1}^n c_j \pi(g_j) \xi \Big\rangle
\, \ge \, 0 .
$$
These coefficients $\varphi_{\pi, \xi}$ are the
\textbf{functions of positive type associated to}~$\pi$.
Note that $\varphi_{\pi, \xi}(e) = \Vert \xi \Vert^2$.
\index{Function of positive type! associated to a representation}
\end{defn}

Since two equivalent representations give rise to 
the same space of diagonal matrix coefficients,
functions of positive type are associated to equivalence classes of representations.
Below, we will often use (abusively) the same symbol for a representation
and its equivalence class.

\begin{defn}
% 1.B.4
\label{defcyclcrep}
\index{Cyclic representation}
\index{Representation! cyclic}
Let $G$ be a topological group
and $\pi$ a representation of $G$ in a Hilbert space $\Hi$.
A vector $\xi \in \Hi$ is \textbf{cyclic}
if the smallest closed subspace of $\Hi$
containing $\pi(G) \xi$ is $\Hi$ itself;
the representation $\pi$ is \textbf{cyclic} 
if it has a cyclic vector.
\end{defn}

\begin{constr}
% 1.B.5
\label{constructionGNS2}
\index{Gel'fand--Naimark--Segal representation!}
\index{GNS, short for Gel'fand--Naimark--Segal}
\index{Gel'fand--Naimark--Segal representation! $1$@with $\varphi \in P(G)$} 
The \textbf{Gel'fand--Naimark--Segal representation},
or shortly \textbf{GNS representation}, 
or \textbf{GNS construction},
shows that the diagonal matrix coefficients of Definition \ref{exposfunctdiagcoeff}
provide all examples of functions of positive type on a topological group $G$.
More precisely, the construction below shows that,
given $\varphi \in P(G)$, $\varphi \ne 0$,
there exist a representation $\pi_\varphi$ of $G$ in a Hilbert space $\Hi_\varphi$, 
and a cyclic vector $\xi_\varphi \in \Hi_\varphi$, 
such that $\varphi$ is a diagonal matrix coefficient of $\pi$, i.e., 
such that $\varphi (g) = \langle \pi_\varphi (g) \xi_\varphi \mid \xi_\varphi \rangle$
for all $g \in G$.

\vskip.2cm

\index{$j1$@$\C[G]$ group algebra of a group $G$}
Let $\C[G]$ denote the complex vector space of functions 
$f \, \colon G \to \C$ with finite support.
Then $\C[G]$ is a unital $*$-algebra, for the convolution product 
defined by 
$$
(f_1 \ast f_2) (x) \, = \, \sum_{y \in G} f_1(xy^{-1}) f_2(y)
\hskip.5cm \text{for} \hskip.2cm
f_1, f_2 \in \C[G], \hskip.1cm x \in G,
$$
and the involution given by
$$
f^*(x) \, = \, \overline{f(x^{-1})} 
\hskip.5cm \text{for} \hskip.2cm
f \in \C[G], \hskip.1cm x \in G .
$$
For $g \in G$, the characteristic function $\delta_g$ of $\{g\}$ is in $\C[G]$.
For $g, h \in G$, we have $\delta_g \delta_h = \delta_{gh}$
and $\delta_g^* = \delta_{g^{-1}}$.
\par

Let $\varphi$ be a non-zero function of positive type on $G$.
A positive Hermitian form $\Phi_\varphi$ is defined on $\C[G]$ by
$$
\Phi_\varphi( f_1, f_2) \, = \, 
\sum_{x,y \in G} f_1(x)\overline{ f_2(y)} \varphi(y^{-1} x)
\hskip.5cm \text{for all} \hskip.2cm 
f_1,f_2 \in \C[G].
$$
It follows from the Cauchy--Schwarz inequality 
that the set $J_\varphi$ of all $f \in \C[G]$ such that $\Phi_\varphi(f, f) = 0$ 
is a left ideal of $\C[G]$.
We write $[f]$ for the class in the vector space $\C[G]/J_\varphi$
of a function $f \in \C[G]$.
We define $\Hi_\varphi$ as the Hilbert space completion of $\C[G]/J_\varphi$ 
for the scalar product given on $\C[G]/J_\varphi$ by 
$$
\langle [f_1] \mid [f_2] \rangle \, = \, \Phi_\varphi( f_1, f_2)
\hskip.5cm \text{for all} \hskip.2cm 
f_1, f_2 \in \C[G] .
$$
\par

For $g \in G$, the operator on $\C[G]$ defined by 
$f \mapsto \big(x \mapsto f(g^{-1}x)\big)$
preserves the bilinear form $\Phi_\varphi$,
passes therefore to the quotient $\C[G]/J_\varphi$, 
and extends to a unitary operator on $\Hi_\varphi$, denoted by $\pi_\varphi(g)$. 
It is clear that $\pi_\varphi(gh) = \pi_\varphi(g) \pi_\varphi(h)$ for all $g, h \in G$.
\par

Let $f \in \C[G]$; set $\xi = [f] \in \C[G]/J_\varphi$. Then 
$$
\begin{aligned}
\langle \pi_\varphi(g) \xi \mid \xi \rangle \, 
&\, = \, \sum_{x,y \in G} f(x) \overline{f(y)} \varphi(y^{-1}gx) ,
\\
\langle \xi \mid \pi_\varphi(g) \xi \rangle \,
&\, = \, \sum_{x,y \in G} f(x) \overline{f(y)} \varphi(y^{-1} g^{-1} x) ,
\end{aligned}
$$
so that 
$$
\Vert \pi_\varphi(g)\xi - \xi \Vert^2 \, = \, 
2 \Vert \xi \Vert^2 
- \sum_{x,y \in G} f(x) \overline{f(y)}
\left( \varphi(y^{-1}gx) + \varphi(y^{-1} g^{-1} x) \right) .
$$
Since $\varphi$ is continuous and $f$ has finite support,
$$
\lim_{g \to e} \sum_{x,y \in G} f(x) \overline{f(y)} \varphi(y^{-1}gx)
\, = \,
\lim_{g \to e} \sum_{x,y \in G} f(x) \overline{f(y)} \varphi(y^{-1}g^{-1}x)
\, = \,
\Vert \xi \Vert^2 .
$$
It follows that $g \mapsto \pi_\varphi (g)\xi$ is continuous
for every $\xi \in \C[G]/J_\varphi$, and therefore for every $\xi \in \Hi_ \varphi$. 
This shows that the assignment $g \mapsto \pi_\varphi(g)$
is continuous, i.e., that it is a representation of $G$ in $\Hi_\varphi$.
Moreover, the vector $\xi_\varphi := [\delta_e]$ is a cyclic vector
of norm $\sqrt{ \varphi (e) }$, and we have 
$$
\varphi(g) \, = \, \langle \pi_\varphi(g) \xi_\varphi \mid \xi_\varphi \rangle 
\hskip.5cm \text{for all} \hskip.2cm 
g \in G.
$$
The \textbf{GNS triple} associated to $\varphi$
consists of the Hilbert space $\Hi_\varphi$,
the representation $\pi_\varphi$, 
and the cyclic vector $\xi_\varphi$.
\end{constr}

\begin{rem}
% 1.B.6
\label{GNSvariants}

The terminology ``GNS construction'',
which applies both to groups and C*-algebras (see \ref{constructionGNS9}),
refers to two seminal articles: \cite{GeNa--48} and \cite{Sega--47}.
\par

The GNS construction described above 
associates a representation of a group~$G$
to a function of positive type on the group.
It is closely related to the GNS construction
which associates a representation of a C*-algebra $A$
to a state on the C*-algebra.
See Subsection \ref{SS:Subsection-States}, and \cite[\S~2]{Dixm--C*}. 
\par

Construction \ref{constructionGNS2} uses the algebra $\C[G]$
of finitely supported functions, as in \cite{Thom--64a};
this will be suited for our treatment of traces in Section~\ref{Section-GNS-Traces}.
There are variants using, instead of $\C[G]$, 
a space of continuous functions on $G$, as in \cite[C.1--C.4]{BeHV--08},
or $L^1(G, \mu_G)$ when $G$ is locally compact
with Haar measure $\mu_G$, as in \cite[3.20]{Foll--16}.
% or an appropriate C*-algebra, as in Subsection \ref{SS:Subsection-States}.
\end{rem}

\begin{exe}
% 1.B.7
\label{3exeGNS}
(1)
On a topological group $G$,
the constant function of value~$1$ is of positive type,
and the corresponding GNS representation is the unit representation $1_G$.
\par

More generally, every unitary character $\chi \, \colon G \to \T$ is a function of positive type.
With the notation of \ref{constructionGNS2},
the ideal $J_\chi$ is the kernel of the algebra homomorphism
$\C[G] \to \C, \hskip.1cm f \mapsto \sum_{x \in G} \chi(x) f(x)$, of codimension one,
and the GNS representation $\pi_\chi$ of $G$ is equivalent to $\chi$.

\vskip.2cm

(2)
Let $H$ be an open subgroup of the topological group $G$,
and $\psi \, \colon G \to \C$ the characteristic function of $H$.
Then $\psi$ is of positive type,
and $\pi_\psi$ is equivalent to 
the quasi-regular representation of $G$ in $\ell^2(G/H)$,
with cyclic vector the characteristic function
of the point $H$ in $G/H$
(see \ref{quasiregularrep} for the definition).
\par

More generally, let $\chi \, \colon H \to \T$ be a unitary character,
and $\widetilde \chi \, \colon G \to \C$ the function mapping $g$ to $\chi (g)$
if $g \in H$ and to $0$ if $g \in G \smallsetminus H$.
Then $\widetilde \chi$ is of positive type,
and $\pi_{\widetilde \chi}$ is equivalent to the representation $\Ind_H^G \chi$,
as defined in \ref{Section-IrrIndRep} below.
(This is a particular case of what is proved
in Proposition \ref{diagcoeffinduced}.)

\vskip.2cm

(3)
In particular, for a discrete group $\Gamma$,
the characteristic function $\delta_e$ of the unit element
is of positive type, and $\pi_{\delta_e}$ is equivalent
to the left regular representation~$\lambda_\Gamma$.

\vskip.2cm

(4)
Let $\varphi, \psi$ be two functions of positive type on a topological group.
The function $\varphi + \psi$ is of positive type,
and the representation $\pi_{\varphi + \psi}$ 
is equivalent to a subrepresentation of the direct sum $\pi_\varphi \oplus \pi_\psi$.
The function $\varphi \psi$ is of positive type,
and the representation $\pi_{\varphi \psi}$ 
is equivalent to a subrepresentation of the tensor product $\pi_\varphi \otimes \pi_\psi$.

\vskip.2cm

(5)
Let $G$ be a \emph{locally compact abelian} group
and $\widehat G$ the dual group of $G$.
For a probability measure $\mu$ on the $\sigma$-algebra 
$\mathcal B (\widehat G)$ of Borel subsets of $\widehat G$,
let $\mathcal{F}(\mu)$ denote the Fourier--Stieltjes transform of $\mu$,
which is a bounded continuous function on $G$.
Bochner's theorem establishes that
the map $\mu \mapsto \mathcal{F}(\mu)$ is a bijection
from the set of probability measures on $\mathcal B (\widehat G)$
to the set $P_1(G)$ of normalized functions of positive type on $G$
(Theorem \ref{TheoremBochner}).
\par

In Construction \ref{defpimupourGab} below,
we associate to a probability measure $\mu$ on $\widehat G$
a canonical representation $\pi_\mu$ of $G$ on $L^2(\widehat G, \mu)$.
Let $1_{\widehat G}$ denote the constant function of value $1$ on $\widehat G$.
Then the matrix coefficient
$\langle \pi_{\mu}(g) 1_{\widehat G} \mid 1_{\widehat G} \rangle
= \int_{\widehat G} \chi(g) d\mu(\chi)$
is the inverse Fourier transform $\overline{\mathcal F}(\mu)$ of $\mu$,
hence $\langle \pi_{\mu}(g) 1_{\widehat G} \mid 1_{\widehat G} \rangle
= \mathcal F (\check{\mu})$,
where $\check{\mu}$ is the image of $\mu$ by the map
$G \to G, \hskip.1cm g \mapsto g^{-1}$
(see Appendix \ref{AppLCA+Pont}).
\end{exe}

\index{Equivalent! $2$@cyclic representations}
Let us say that two cyclic representations $\pi_j \, \colon G \to \U(\Hi_j)$, 
with cyclic vector $\xi_j \in \Hi_j$ ($j = 1, 2$),
are \textbf{equivalent}
if there exists a surjective isometry
$U \, \colon \Hi_1 \to \Hi_2$ such that $U \pi_1(g) = \pi_2(g) U$ for all $g \in G$,
and $U \xi_1 = \xi_2$.
In this case, it is clear that the functions of positive type 
associated to $\xi_1$ and $\xi_2$ coincide, that is,
$\langle \pi_1(g)\xi_1 \mid \xi_1 \rangle = \langle \pi_2(g)\xi_2 \mid \xi_2 \rangle$
for all $g \in G$.
The next result shows that the converse statement is also true.

\begin{prop}
% 1.B.8
\label{GNSbijP(G)cyclic}
Let $G$ be a topological group. 
The assignments
$$
\varphi \rightsquigarrow ( \pi_\varphi, \Hi_\varphi, \xi_\varphi)
$$
given by the GNS construction and
$$
(\pi, \Hi, \xi) \, \mapsto \, \big(\varphi_{\pi, \xi} \, \colon 
g \mapsto \langle \pi(g) \xi \mid \xi \rangle\big)
$$
induce bijections inverse to each other between
\begin{enumerate}[label=(\roman*)]
\item\label{iDEGNSbijP(G)cyclic}
the set $P_1(G)$ of normalized functions of positive type on $G$, and
\item\label{iiDEGNSbijP(G)cyclic}
the set of triples $(\pi, \Hi, \xi)$ 
where $\pi$ is a representation of $G$ in $\Hi$ 
and $\xi \in \Hi$ a unit cyclic vector, modulo equivalence. 
\end{enumerate}
\end{prop}

\begin{proof}
Let us verify that the assignment
of the function $\varphi_{\pi, \xi}$
to an equivalence class $(\pi, \Hi, \xi)$
is injective.
The other verifications are left to the reader.
\par

Let $(\pi_1, \Hi_1, \xi_1)$ and $(\pi_2, \Hi_2, \xi_2)$ be triples 
consisting of cyclic representations $\pi_1$ and $\pi_2$, 
with cyclic vectors $\xi_1$ and $\xi_2$ respectively.
Let $\varphi_1, \varphi_2$ be the functions on $G$ defined by
$$
\varphi_1(g) \, = \, \langle \pi_1(g)\xi_1 \mid \xi_1 \rangle
\hskip.5cm \text{and} \hskip.5cm
\varphi_2(g) \, = \, \langle \pi_2(g)\xi_2 \mid \xi_2 \rangle 
\hskip.5cm \text{for all} \hskip.2cm 
g \in G.
$$
Assume that $\varphi_1 = \varphi_2$, and denote this function by $\varphi$.
We have to show that the triples $(\pi_1, \Hi_1, \xi_1)$ and $(\pi_2, \Hi_2, \xi_2)$
are equivalent.
\par

For all $g_1, \hdots, g_n \in G$ and $c_1, \hdots, c_n \in \C$, we have
$$
\Big\Vert \sum_{i = 1}^n c_i \pi_2(g_i) \xi_2 \Big\Vert^2
\, = \,
\sum_{i, j = 1}^n c_i \overline{c_j} \varphi(g_j^{-1} g_i)
\, = \, 
\Big\Vert \sum_{i = 1}^n c_i \pi_ 1 (g_i) \xi_ 1 \Big\Vert^2 .
$$
It follows that the map 
$\sum_{i = 1}^n c_i \pi_1(g_i) \xi_1 \mapsto \sum_{i = 1}^n c_i \pi_2(g_i) \xi_ 2$ 
is well defined on the linear span of $\pi_1(G)\xi_1$ 
and extends to a Hilbert space isomorphism $U \, \colon \Hi_1 \to \Hi_ 2$.
One checks that $U\pi_1(g) = \pi_ 2(g)U$ for all $g \in G$, and $U \xi_1 = \xi_2$,
as was to be shown.
\end{proof}

\begin{rem}
% 1.B.9
\label{Rem-AutoPositiveType}
There is a natural \emph{right} action $P(G) \curvearrowleft \Aut(G)$
of the group $\Aut(G)$ of bicontinuous automorphisms of $G$
on the the set $P(G)$ of functions of positive type on $G$, given by 
$$
\varphi^\theta(g) = \varphi(\theta(g))
\hskip.5cm \text{for} \hskip.2cm
\theta \in \Aut(G), \hskip.1cm \varphi \in P(G), \hskip.1cm g \in G.
$$
Let $\varphi \in P(G)$ and $\theta \in \Aut(G)$.
If $(\pi, \Hi, \xi)$ is a GNS-triple associated to $\varphi$,
then $(\pi^\theta, \Hi, \xi)$ is a GNS-triple associated to $\varphi^\theta$,
where $\pi^\theta$ is the conjugate representation of $\pi$ by $\theta$
as defined in Remark~\ref{Rem-AutoRep}.
\end{rem}

In Proposition \ref{P1=cyclicmodeqIRR}, we will characterize
those $\varphi \in P(G)$ for which $\pi_\varphi$ is irreducible.
This requires a definition and two preliminary lemmas.
\par

\index{Indecomposable! $1$@for $\varphi \in P_{\le 1}(G)$}
Let ${\rm Extr} (P_{\le 1}(G))$ be the set of 
functions $\varphi$ in $P_{\le 1}(G)$ that are 
\textbf{indecomposable} in the following sense: 
if $\varphi = t \varphi_1+ (1 - t) \varphi_2$ 
for $\varphi_1, \varphi_2 \in P_{\le 1}(G)$ and $t \in \mathopen] 0,1 \mathclose[$, 
then $\varphi_1 = \varphi_2 = \varphi$.
Define similarly ${\rm Extr} (P_1(G))$.
\index{$d3$@${\rm Extr} (\cdot)$
indecomposable elements}
\par

As the following lemma shows, there is a simple relationship between the sets 
${\rm Extr} (P_{\le 1}(G))$ and ${\rm Extr} (P_1(G))$.

\begin{lem}
% 1.B.10
\label{Lem-ExtP}
For every topological group $G$, we have
$$
{\rm Extr} (P_{\le 1}(G))
\, = \, 
{\rm Extr} (P_1(G)) \cup \{0\}.
$$
\end{lem}

\begin{proof}
Let us first show that 
$
{\rm Extr} (P_1(G)) \cup \{0\} 
\subset 
{\rm Extr} (P_{\le 1}(G)).
$
\par

Let $\varphi \in {\rm Extr} (P_1(G))$; 
let $\varphi_1, \varphi_2 \in P_{\le 1}(G)$ and $t \in \mathopen] 0,1 \mathclose[$
be such that $t\varphi_1 + (1 - t)\varphi_2 = \varphi$.
Then $t\varphi_1(e) + (1 - t)\varphi_2(e) = 1$. 
Therefore $\varphi_1(e) = \varphi_2(e) = 1$,
that is, $\varphi_1, \varphi_2 \in P_1(G)$,
and hence $\varphi_1 = \varphi_2 = \varphi$.
Therefore $\varphi \in {\rm Extr} (P_{\le 1}(G))$.
\par

Let $\varphi_1, \varphi_2 \in P_{\le 1}(G)$ and $t \in \mathopen] 0,1 \mathclose[$
be such that $t\varphi_1 + (1 - t)\varphi_2 = 0$. 
Then $\varphi_1(e) = 0$ and $\varphi_2(e) = 0$,
that is, $\varphi_1 = 0$ and $\varphi_2 = 0$. 
Therefore $0 \in {\rm Extr} (P_{\le 1}(G))$. 

\vskip.2cm

For the opposite inclusion, consider
$\varphi \in {P}_{\le 1}(G)$ such that $\varphi \notin (P_1(G) \cup \{0\})$. 
For $t := \varphi(e)$, we have $0 < t < 1$, $\dfrac{1}{t}\varphi \in P_1(G)$,
and $\varphi = t(\dfrac{1}{t}\varphi) + (1 - t)0$.
It follows that $\varphi \notin {\rm Extr} (P_{\le 1}(G))$.
\end{proof}

Let $G$ be a topological group, and $\varphi \in P(G)$.
For $\psi \in P(G)$, write $\psi \le \varphi$ if $\varphi - \psi \in P(G)$. 
Let $T \in \pi_\varphi(G)'$ be such that $0 \le T \le I$.
Here $\pi_\varphi(G)'$ is the commutant of $\pi_\varphi(G)$
and $I$ stands for the identity operator on $\Hi_\varphi$.
Define $\varphi_T \, \colon G \mapsto \C$ by
\index{$h5$@$I = \mathrm{Id}_\Hi \in \Li (\Hi)$ identity operator on $\Hi$}
$$
\varphi_T(g) 
\, = \,
\langle \pi_\varphi(g)T \xi_ \varphi \mid \xi_ \varphi \rangle
\, = \,
\langle \pi_ \varphi(g)T^{1/2}\xi_ \varphi \mid T^{1/2}\xi_ \varphi \rangle
\hskip.5cm \text{for all} \hskip.2cm
g \in G .
$$
Then $\varphi_T \in P(G)$ and, since $\varphi - \varphi_T = \varphi_{I - T}$, 
we have $\varphi_T \le \varphi$.
(For $T \ge 0$ and $T^{1/2}$, there is a reminder in Appendix \ref{AppAlgC*}.)

\begin{lem}
% 1.B.11
\label{Lem-CommutantCyclic}
The map $T \mapsto \varphi_T$
is a bijective correspondence between 
the set of $T \in \pi_ \varphi(G)'$ such that $0 \le T \le I$
and the set of $\psi \in P(G)$ such that $\psi \le \varphi$. 
\end{lem}

\begin{proof}
Let $T, S \in \pi_ \varphi(G)'$ be such that $0 \le T, S \le I$ 
be such that $\varphi_T = \varphi_S$; then 
$$
\langle \pi_ \varphi(g) \xi_ \varphi \mid T \xi_\varphi \rangle
\, = \, 
\langle \pi_ \varphi(g)\xi_ \varphi \mid S \xi_ \varphi \rangle
\hskip.5cm \text{for all} \hskip.2cm 
g \in G .
$$
Since $\xi_ \varphi$ is cyclic for $\pi_ \varphi$, this implies that 
$T \xi_ \varphi = S \xi_ \varphi$.
Moreover, since $T, S \in \pi_ \varphi(G)'$, 
$$
T \pi_\varphi(g) \xi_\varphi
\, = \, \pi_ \varphi(g) T \xi_\varphi 
\, = \, \pi_ \varphi(g) S \xi_\varphi 
\, = \, S \pi_\varphi(g)\xi_\varphi
\hskip.5cm \text{for all} \hskip.2cm 
g \in G ,
$$
so that $T = S$. 
This shows that the map $T \mapsto \varphi_T$ is injective.

\vskip.2cm

Let $\psi \in P(G)$ be such that $\psi \le \varphi$.
Denote by $\Phi_\psi$ the associated positive Hermitian form on $\C[G]$.
For $f_1, f_2 \in \C[G]$, we have
$$
\vert \Phi_\psi(f_1, f_2) \vert^2 
\, \le \, \Phi_\psi(f_1, f_1) \Phi_\psi(f_2, f_2)
\, \le \, \Phi_\varphi(f_1,f_1) \Phi_\varphi(f_2,f_2)
\, = \, \Vert [f_1]\Vert^2 \Vert [f_2]\Vert^2.
$$
This shows that $\Phi_\psi$ induces 
a continuous positive Hermitian form on $\Hi_ \varphi$,
defined by $([f_1], [f_2]) \mapsto \Phi_\psi(f_1, f_2)$.
It follows that there exists $T \in \Li (\Hi_ \varphi)$, with $0 \le T \le I$, such that
$$
\Phi_\psi( f_1, f_2) \, = \, \langle [f_1] \mid T [f_2] \rangle
\hskip.5cm \text{for all} \hskip.2cm 
f_1,f_2 \in \C[G] .
$$
Therefore
$$
\psi (y^{-1}x)
\, = \, \Phi_\psi( \delta_x, \delta_y)
\, = \, \langle \pi_\varphi(x)\xi_ \varphi \mid T \pi_\varphi(y)\xi_ \varphi \rangle
\hskip.5cm \text{for all} \hskip.2cm 
x,y \in G .
\leqno(\sharp)
$$
\par

For $x, y, z \in G$, we have
$$
\begin{aligned}
 \langle \pi_\varphi(x)\xi_ \varphi &\mid T \pi_\varphi(z)\pi_\varphi(y)\xi_ \varphi \rangle
 \, = \, \psi((zy)^{-1}x)
 \\ 
 \, &= \,
 \psi(y^{-1}(z^{-1}x))
 \, = \, \langle \pi_\varphi(z^{-1}x)\xi_ \varphi \mid T \pi_\varphi(y)\xi_ \varphi \rangle
 \\
\, &= \,
\langle \pi_\varphi(x)\xi_\varphi \mid \pi_\varphi(z) T \pi_\varphi(y)\xi_ \varphi \rangle.
\end{aligned}
$$
Since $\xi_ \varphi$ is a cyclic vector, it follows that $T \in \pi_ \varphi(G)'$.
By ($\sharp$), we have $\psi = \varphi_T$.
This shows that the map $T \mapsto \varphi_T$ is surjective.
\end{proof}

\begin{prop}
% 1.B.12
\label{P1=cyclicmodeqIRR}
Let $G$ be a topological group.
Let $\varphi \in P_1(G)$;
denote by $(\Hi_\varphi, \pi_\varphi, \xi_\varphi)$
the corresponding GNS triple.
\par

The representation $\pi_\varphi$ is irreducible 
if and only if $\varphi \in {\rm Extr} (P_1(G))$.
\end{prop}

\begin{proof}
Assume first that $\pi_\varphi$ irreducible; 
let $\varphi_1, \varphi_2 \in P_1(G)$ and $t \in \mathopen] 0,1 \mathclose[$
be such that $\varphi = t\varphi_1 + (1 - t)\varphi_2$. 
Let $j \in \{1, 2 \}$; since $\varphi_j \le \varphi$, 
there exists by Lemma \ref{Lem-CommutantCyclic}
$T_j \in \pi_\varphi(G)'$ with $0 \le T_j \le I$
such that $\varphi_j = \varphi_{T_j}$.
By irreducibility of $\pi_\varphi$, we then have $T_j = \lambda_j \mathrm{Id}$ 
for some $\lambda_j \in \C$;
indeed, $T_j = \mathrm{Id}$, because $\varphi_j(e) = 1$. 
Therefore $\varphi_1 = \varphi_2 = \varphi$, 
and it follows that $\varphi \in {\rm Extr} (P_1(G))$. 

\vskip.2cm

Assume now that $\pi_\varphi$ is not irreducible;
there exists an orthogonal projection $P \in \pi_\varphi(G)'$ 
such that $0 \ne P \ne I$.
Set $t := \varphi_P(e) = \Vert P \xi_\varphi \Vert^2$.
We claim that $t \ne 0$. 
Indeed, assume that $t = 0$. 
Then $P\xi_\varphi = 0$ 
and hence $P \pi(g) \xi_\varphi = \pi(g) P\xi_\varphi = 0$ for all $g \in G$; 
since $\xi_\varphi$ is a cyclic vector, this implies $P = 0$,
in contradiction with our assumption.
Similarly, we have $t \ne 1$, because $I-P \ne 0$.
\par

Set $\varphi_1 = \frac{1}{t} \varphi_P$ and $\varphi_2 = \frac{1}{1 - t} \varphi_{I-P}$; 
then $\varphi_1, \varphi_2 \in P_1(G)$ and $\varphi = t \varphi_1 + (1 - t) \varphi_2$; 
so, $\varphi \notin {\rm Extr} (P_1(G))$.
\end{proof}

\subsection*{The case of locally compact groups}

Assume from now on that $G$ is a locally compact group.
Choose a left Haar measure $\mu_G$ on $G$,
i.e., on the $\sigma$-algebra $\mathcal B (G)$ of the Borel subsets of $G$.
\par

The convex cone $P(G)$ of continuous functions of positive type on $G$
and its subsets $P_{\le 1}(G)$, $P_1(G)$, ${\rm Extr} (P_{\le 1}(G))$,
${\rm Extr} (P_1(G))$
can be viewed as subspaces of $L^\infty(G, \mu_G)$.
We consider these sets as topological spaces,
for the topology induced by the weak$^*$-topology of $L^\infty(G, \mu_G)$.

\begin{theorem}
% 1.B.13
\label{Banach--Alaoglu--Krein--Milman--Raikov}
Let $G$ be a locally compact group.
Let the notation be as above.
\begin{enumerate}[label=(\arabic*)]
\item\label{1DEBanach--Alaoglu--Krein--Milman--Raikov}
The convex subset $P_{\le 1}(G)$ of $L^\infty(G, \mu_G)$ 
is compact for the weak$^*$-topology.
\item\label{2DEBanach--Alaoglu--Krein--Milman--Raikov}
The convex hull of ${\rm Extr} (P_{\le 1}(G))$
is weak$^*$-dense in $P_{\le 1}(G)$.
\item\label{3DEBanach--Alaoglu--Krein--Milman--Raikov}
The convex hull of ${\rm Extr} (P_1(G))$
is weak$^*$-dense in $P_1(G)$.
\item\label{4DEBanach--Alaoglu--Krein--Milman--Raikov}
On $P_1(G)$, the weak$^*$-topology 
and the topology of uniform convergence on compact subsets of $G$ coincide.
\item\label{5DEBanach--Alaoglu--Krein--Milman--Raikov}
For $\Gamma$ a discrete group, the subspace $P_1(\Gamma)$
is weak$^*$-closed in $\ell^\infty(\Gamma)$.
\end{enumerate}
\end{theorem}

\begin{proof}[Sketch of proof and reference]
Claim \ref{1DEBanach--Alaoglu--Krein--Milman--Raikov}
is a consequence of the Banach--Alaoglu theorem.
Claim \ref{2DEBanach--Alaoglu--Krein--Milman--Raikov}
is a consequence of the Krein--Milman theorem,
and Claim \ref{3DEBanach--Alaoglu--Krein--Milman--Raikov} follows.
Claim \ref{4DEBanach--Alaoglu--Krein--Milman--Raikov} is a theorem of Raikov.
It follows that the convex hull of ${\rm Extr} (P_1(G))$ is dense in $P_1(G)$
for the topology of uniform convergence on compact subsets of~$G$.
Claim \ref{5DEBanach--Alaoglu--Krein--Milman--Raikov} follows from
Claim \ref{1DEBanach--Alaoglu--Krein--Milman--Raikov}
because the evaluation map
$\ell^\infty(\Gamma) \to \C, \hskip.1cm \varphi \mapsto \varphi(e)$
is continuous for the weak$^*$-topology.
\par

For details, we refer to \cite[C.5]{BeHV--08}. 
\end{proof}

\begin{rem}
% 1.B.14
\label{bijftp}
(1)
Note that the space $P_1(G)$ need not be closed in $L^\infty(G, \mu_G)$
for the weak$^*$-topology when $G$ is not discrete.
For example, the weak$^*$-closure of $P_1(\T)$ in $L^\infty(\T, \mu_T)$ contains $0$.
Indeed, for $k \in \Z$, let $e_k \in P_1(\T)$ be the unitary character
defined on $\T = \{ z \in \C \mid \vert z \vert = 1 \}$ by $e_k(z) = z^k$.
By the Riemann--Lebesgue lemma, we have
$$
\lim_{k \to \infty} \int_\T f(z) e_k(z) d\mu_T(z) \, = \, 0
\hskip.5cm \text{for all} \hskip.2cm
f \in L^1(\T, \mu_T) ,
$$
i.e., the sequence $(e_k)_{k \ge 1}$ converges to $0 \in L^\infty(\T, \mu_T)$
for the weak$^*$-topology.
More generally, the Riemann--Lebesgue lemma holds 
for an arbitrary locally compact abelian group
(Proposition 3 in \cite[Chap.~II, \S~1, No 2]{BTS1--2}),
and the same argument as above shows that the closure of $P_1(G)$ contains $0$
for any locally compact abelian group $G$ such that $\widehat G$ is not compact.

\vskip.2cm

(2)
For a locally compact group $G$,
the bijection established in Proposition \ref{GNSbijP(G)cyclic}
fits in a larger picture;
indeed, there are canonical bijections
between the following sets:
\begin{enumerate}[label=(\roman*)]
\item\label{aDEbijftp}
$P_1(G)$, 
\item\label{bDEbijftp} 
the set of cyclic representations $(\pi, \Hi, \xi)$ of $G$
with $\Vert \xi \Vert = 1$ modulo equivalence,
\item\label{cDEbijftp}
the set of positive forms on the group algebra $L^1(G, \mu_G)$,
\item\label{dDEbijftp}
the set of positive forms on the maximal C*-algebra $C^*_{\rm max}(G)$,
as defined in Section \ref{C*algLCgroup}.
\end{enumerate}
We come back to this in Example \ref{exGNSpourC*}(2).

\vskip.2cm

(3)
A non-discrete topological group $G$
has ``discontinuous representations'', as can be checked as follows.
Denote by $G_{\rm disc}$ the group $G$ made discrete
(that is, viewed as a discrete group).
The Dirac delta function $\delta$ supported at the origin of $G_{\rm disc}$
is of positive type, and discontinuous on $G$.
If follows that the associated GNS representation $\pi_\delta \, \colon G_{\rm disc} \to \U(\Hi)$
is not continuous when viewed as a homomorphism $G \to \U(\Hi)$.
\index{$b5$@$G_{\rm disc}$ topological group made discrete}
\end{rem}

As we now show,
locally compact groups have a distinguished \emph{faithful}\ representation.

\begin{exe}
% 1.B.15
\label{Leftregrep}
\index{Left regular representation}
\index{Representation! left regular}
Let $G$ be a locally compact group
and $\mu_G$ a left Haar measure on $G$.
The \textbf{left regular representation}
$\lambda_G$ of $G$ is defined on $L^2(G, \mu_G)$ by 
$$
\lambda_G(g) f(x) \, = \, f(g^{-1}x) 
\hskip.5cm \text{for all} \hskip.2cm
g, x \in G
\hskip.2cm \text{and} \hskip.2cm
f \in L^2(G, \mu_G) .
$$
This representation is faithful.
Indeed, let $g \in G$, $g \ne e$.
Since $G$ is a locally compact space,
there exist a compact neighbourhood $V$ of $e$ in $G$ such that $g^{-1} \notin V$,
and a continuous function $f \, \colon G \to \mathopen[ 0,1 \mathclose]$
supported in $V$ (hence $f \in L^2(G, \mu_G)$) such that $f(e) = 1$.
Then $f(e) = 1 \ne 0 = f(g^{-1}) = (\lambda_G (g) f) (e)$,
hence $f \ne \lambda_G (g) f$ and $\lambda_G (g) \ne \mathrm{Id}_{L^2(G, \mu_G)}$;
this shows that $\lambda_G$ is faithful.
\par

\index{Right regular representation}
\index{Representation! right regular}
The \textbf{right regular representation} $\rho_G$ of $G$
is defined on $L^2(G, \mu_G)$ by 
$$
\rho_G(g) f(x) \, = \, \Delta_G(g)^{1/2} f(xg)
\hskip.5cm \text{for all} \hskip.2cm
g, x \in G, \, f \in L^2(G, \mu_G),
$$
where $\Delta_G$ is the modular function of $G$
(see Appendix \ref{AppLCG});
note that this uses the \emph{right} action of $G$ on itself,
and yet $\mu_g$ is still a \emph{left} Haar measure on $G$.
\par

The representations $\lambda_G$ and $\rho_G$ are equivalent.
Indeed, if $U \, \colon L^2(G, \mu_G) \to L^2(G, \mu_G)$ is defined by
$$
(U f) (g) \, = \, \Delta_G(g)^{-1/2} f(g^{-1})
\hskip.5cm \text{for all} \hskip.2cm
f \in L^2(G, \mu_G)
\hskip.2cm \text{and} \hskip.2cm
g \in G ,
$$
then $U$ is unitary and $U \lambda_G(g) = \rho_G(g) U$ for all $g \in G$.
More on this in \cite[A.4]{BeHV--08}.
\par

If $G \ne \{e\}$, the commutant $\lambda_G(G)'$ contains $\rho_G(G)$;
in particular, $\lambda_G(G)'$
contains operators which are not multiples of the identity.
It follows that the left-regular representation is not irreducible.
\end{exe}

\index{Gel'fand--Raikov theorem}
When $G$ is locally compact, the Gel'fand--Raikov theorem
of \cite{GeRa--43} establishes that $\widehat G$ is sufficiently large,
in the sense that it separates the points of $G$.
See also \cite[13.6.6]{Dixm--C*}).
The proof uses the faithfulness of the left regular representation.

\begin{theorem}[\textbf{Gel'fand--Raikov}]
% 1.B.16
\label{Gel'fandRaikov}
Let $G$ be a locally compact group.
\par

For every $g \in G$ with $g \ne 1$,
there exists an irreducible representation $\pi$ of $G$ in a Hilbert space $\Hi$
such that $\pi (g) \ne \mathrm{Id}_\Hi$.
\end{theorem}

\begin{proof}
Let $g \in G, g \ne 1$.
By Proposition \ref{P1=cyclicmodeqIRR}, it suffices to show
that there exists $\varphi \in {\rm Extr} (P_1(G))$
with $\varphi(g) \ne 1$.
\par

Since the regular representation $\lambda_G$ of $G$ is faithful,
$\lambda_G(g)$ is not the identity on $L^2(G, \mu_G)$; 
hence, there exists $f \in L^2(G, \mu_G)$ with $\Vert f \Vert =1$ 
and $\langle\lambda_G(g) f \mid f \rangle \ne 1$.
Consider the function $\psi$ on $G$ defined by 
$\psi(x) = \langle \lambda_G(x) f \mid f \rangle$ for $x \in G$.
Then $\psi \in P_1(G)$ and $\psi(g)\ne 1$.
Since the convex hull of ${\rm Extr} (P_{1}(G))$
is dense in $P_{1}(G)$ for the topology of uniform convergence on compact subsets
(Theorem \ref{Banach--Alaoglu--Krein--Milman--Raikov}),
it follows that there exists $\varphi \in {\rm Extr} (P_1(G))$
with $\varphi(g) \ne 1$.
\end{proof}

\begin{rem} 
% 1.B.17
\label{Rem-RepNonLC}
(1)
For second-countable locally compact groups, there is another proof
of the Gel'fand-Raikov theorem~\ref{Gel'fandRaikov},
using a direct integral decomposition of a faithful representation
(e.g., the left regular representation)
into irreducible representations \cite[Pages 109--110]{Mack--76}.

\vskip.2cm

(2)
The dual $\widehat G$ of a general topological group $G$
need not contain more than the unit representation $1_G$.
Indeed, there are topological groups $G$ such that, 
for every representation $\pi$ of $G$ in a Hilbert space $\Hi$ (irreducible or not),
$\pi(g) = \mathrm{Id}_\Hi$ for all $g \in G$.
Examples include
abelian topological groups described in \cite{Bana--83, Bana--91},
the group of orientation preserving homeomorphisms
of the closed unit interval $\mathopen[ 0,1 \mathclose]$
with the compact-open topology
(this is a special case of a result of \cite{Megr--01},
where isometric representations in reflexive Banach spaces are considered),
the group of isometries 
of the universal Urysohn metric space $\mathbb U_1$ of diameter one 
with the compact-open topology \cite[Corollary 1.4]{Pest--07},
and examples occurring in \cite[Corollary 4.10]{BYTs--16}.

\vskip.2cm

(3)
There are also topological groups $G$ that have faithful representations
but no irreducible representation distinct from $1_G$.
This is the case of the group of measurable functions 
from $\mathopen[ 0,1 \mathclose]$ to the unit circle,
with the topology induced by the weak topology 
on the von Neumann algebra $L^\infty( \mathopen[ 0,1 \mathclose] )$,
that makes it an abelian Polish group \cite{Glas--98} 
(see also \cite[Example C.5.10]{BeHV--08}).
\end{rem}

\section[Weak containment and Fell topology]
{Weak containment and Fell topology
for representations of topological groups}
% Section 1.C
\label{SectionWC+FellTop}

In this section, we introduce a topology 
on the dual space $\widehat G$ of a topological group~$G$ 
as well as the related notion of weak containment for
representations of $G$.

\subsection*{Weak containment of representations}

\begin{defn}
% 1.C.1
\label{DefWC}
\index{Representation! weakly contained $\preceq$}
\index{Weakly! contained (representation)}
\index{$a2$@$\preceq$ weak containment of representations}
Let $G$ be a topological group.
Let $\pi, \pi_2$ be two representations of $G$. 
Define $\pi_1$ to be \textbf{weakly contained} in $\pi_2$,
and write $\pi_1 \preceq \pi_2$,
if every function of positive type associated to $\pi_1$
can be approximated, uniformly on every compact subset of $G$,
by sums of functions of positive type associated to $\pi_2$.
\par

The relation $\pi_1 \preceq \pi_2$ only depends
on the equivalence classes of $\pi_1$ and $\pi_2$.
The relation of weak containment is transitive:
if $\pi_1 \preceq \pi_2$ and $\pi_2 \preceq \pi_3$, then $\pi_1 \preceq \pi_3$.
\par

(Note that there is a related \emph{but different} notion of
weak containment due to Zimmer,
see \cite[7.3.5]{Zimm--84}.)
\par

The representations $\pi_1$ and $\pi_2$ are \textbf{weakly equivalent},
written $\pi_1 \sim \pi_2$,
if $\pi_1 \preceq \pi_2$ and $\pi_2 \preceq \pi_1$.
\index{Representation! weakly equivalent $\sim$}
\index{Weakly! equivalent (representation)}
\index{$a3$@$\sim$ weak equivalence of representations}
% \nomenclature{$\sim$}{weak equivalence of representations} iciici
\end{defn}

\begin{rem}
% 1.C.2
\label{propweakcont}
Let $G$ be a topological group and $\pi_1, \pi_2$ two representations of~$G$.

\vskip.2cm

(1)
If $\pi_1$ is contained in $\pi_2$, then $\pi_1$ is weakly contained in $\pi_2$.
%\par
It follows that,
if $\pi_1$ and $\pi_2$ are equivalent, then they are weakly equivalent.

\vskip.2cm

(2)
Let $n_1,n_2$ be two cardinal numbers, not $0$;
then $n_1 \pi_1$ is weakly equivalent to~$\pi_1$,
and $\pi_1 \preceq \pi_2$ if and only if $n_1 \pi_1 \preceq n_2 \pi_2$.

\vskip.2cm

In some particular situations, weak containment implies containment; see Proposition \ref{PropOnWeAndEq}.
\end{rem}

For the notion of amenability and for the next theorem, we refer to \cite[Appendix G]{BeHV--08}.
\index{Amenable! group}

\begin{theorem}[\textbf{Hulanicki--Reiter}]
% 1.C.3
\label{HulanickiReiter}
For a locally compact group $G$, the following properties are equivalent:
\begin{enumerate}[label=(\roman*)]
\item\label{iDEHulanickiReiter}
the unit representation $1_G$ is weakly contained
in the regular representation $\lambda_G$;
% \item[(iv$_b$)]
% there exists a finite dimensional representation of $G$ weakly contained in $\lambda_G$,
\item\label{iiDEHulanickiReiter}
every representation of $G$ 
is weakly contained in $\lambda_G$;
\item\label{iiiDEHulanickiReiter}
$G$ is amenable.
\end{enumerate}
\end{theorem}

We will need several times the following criterion for weak containment. For a proof, see \cite[Lemma F.1.3]{BeHV--08}.

\begin{prop}
% 1.C.4 
\label{oneisenough}
Let $G$ be a topological group and $\pi_1, \pi_2$ two representations of $G$.
Assume that $\pi_1$ has a cyclic vector $\xi$; 
let $\varphi_\xi$ be the corresponding function of positive type:
$\varphi_\xi(g) = \langle \pi_1 (g) \xi \mid \xi \rangle$ for all $g \in G$.
The following conditions are equivalent:
\begin{enumerate}[label=(\roman*)]
\item\label{iDEoneisenough}
$\pi_1$ is weakly contained in $\pi_2$, i.e., 
for every function of positive type $\varphi$ associated to $\pi_1$,
for every $\varepsilon > 0$, and every compact subset $Q$ of $G$,
there exist functions of positive type $\psi_1, \hdots, \psi_k$ associated to $\pi_2$
such that
\hfill\par\noindent
$\sup_{g \in Q} \Big\vert \varphi(g) - \sum_{i = 1}^k \psi_i(g) \Big\vert < \varepsilon$.
\item\label{iiDEoneisenough}
For every $\varepsilon > 0$ and every compact subset $Q$ of $G$,
there exist functions of positive type $\psi_1, \hdots \psi_k$ associated to $\pi_2$
such that
\hfill\par\noindent
$\sup_{g \in Q} \vert \varphi_\xi (g) - \sum_{i = 1}^k \psi_i(g) \vert < \varepsilon$.
\end{enumerate}
\end{prop}

\subsection*{Fell topology}

\begin{defn}
% 1.C.5
\label{DefFellTop}
Let $G$ be a topological group.
Let $\mathcal R$ be a set of equivalence classes of representations of $G$
and let $\pi \in \mathcal R$.
For a finite sequence $\varphi_1, \hdots, \varphi_n$
of functions of positive type associated to $\pi$,
a compact subset $Q$ of $G$, and a real number $\varepsilon > 0$,
define $W(\pi; \varphi_1, \hdots, \varphi_n, Q, \varepsilon)$
to be the subset of those $\rho \in \mathcal R$ which have the following property:
there exist functions $\psi_1, \hdots, \psi_n$ from $G$ to $\C$,
each one a sum of functions of positive type associated to $\rho$, such that
\index{Fell topology}
$$
\vert \varphi_i (g) - \psi_i (g) \vert \, < \, \varepsilon
\hskip.5cm \text{for all} \hskip.2cm
i \in \{1, \hdots, n\}
\hskip.2cm \text{and} \hskip.2cm
g \in Q .
$$
The sets $W(\pi; \varphi_1, \hdots, \varphi_n, Q, \varepsilon)$ constitute
a base of a topology on $\mathcal R$ called the 
\textbf{Fell topology}.
\par

Among the most important cases for $\mathcal R$,
there are the dual $\widehat G$, 
the set of equivalence classes of cyclic representations of $G$,
and the set of equivalence classes of representations of $G$
in Hilbert spaces of dimension at most some cardinal $\aleph$
(for example in separable Hilbert spaces when $G$ is a second-countable group).
\end{defn}

In case $\mathcal R = \widehat G$,
the following proposition provides another description of the Fell topology.

\begin{prop}
% 1.C.6
\label{Pro-ConvFell}
Let $G$ be a topological group and $(\rho_\iota)_{\iota \in I}$ a net in $\widehat G$.
Let $\pi \in \widehat G$
and let $\varphi$ be a normalized function of positive type associated to $\pi$.
The following properties are equivalent:
\begin{enumerate}[label=(\roman*)]
\item\label{iDEPro-ConvFell}
the net $(\rho_\iota)_{\iota \in I}$ converges to $\pi$ in the Fell topology;
\item\label{iiDEPro-ConvFell}
there exists a net $(\psi_\iota)_{\iota \in I}$,
where each $\psi_\iota$ is a sum of functions of positive type associated to $\rho_\iota$,
such that $\lim_\iota \psi_\iota = \varphi$ uniformly on compact subsets of $G$.
\end{enumerate}
Suppose moreover that $G$ is locally compact.
Properties \ref{iDEPro-ConvFell} and \ref{iiDEPro-ConvFell} are equivalent to:
\begin{enumerate}[label=(\roman*)]\addtocounter{enumi}{2}
\item\label{iiiDEPro-ConvFell}
there exists a net $(\psi_\iota)_{\iota \in I}$
of functions of positive type associated to $\rho_\iota$
such that $\lim_\iota \psi_\iota = \varphi$ uniformly on compact subsets of $G$.
\end{enumerate}
\end{prop}

\begin{proof}
The equivalence of \ref{iDEPro-ConvFell} and \ref{iiDEPro-ConvFell}
is an immediate consequence of Proposition~\ref{oneisenough}.
For \ref{iiiDEPro-ConvFell}, we refer to \cite[Proposition F.1.4]{BeHV--08}.
\end{proof}

The next proposition gives a description of Fell's topology on $\widehat G$
in terms of weak containment.

\begin{prop}
% 1.C.7
\label{weakcont=Felltop}
Let $G$ be a locally compact group, $(\pi_\iota)_{\iota \in I}$ a net in $\widehat G$,
and $\pi \in \widehat G$. The following properties are equivalent:
\begin{enumerate}[label=(\roman*)]
\item\label{iDEweakcont=Felltop}
$(\pi_\iota)_{\iota \in I}$ converges to $\pi \in \widehat G$;
\item\label{iiDEweakcont=Felltop}
$\pi$ is weakly contained in $\bigoplus_{\jmath \in J} \pi_{\jmath}$
for every subnet $(\pi_{\jmath})_{\jmath \in J}$ of $(\pi_\iota)_{\iota \in I}$.
\end{enumerate}
\end{prop}

\begin{proof}
The fact that \ref{iDEweakcont=Felltop} implies \ref{iiDEweakcont=Felltop}
is a straightforward consequence of Proposition~\ref{Pro-ConvFell} and holds for an arbitrary topological group $G$.
The fact that \ref{iiDEweakcont=Felltop} implies \ref{iDEweakcont=Felltop}
follows from the definition of Fell's topology
in combination with Proposition F.1.4 in \cite{BeHV--08}.
\end{proof}

\begin{prop}
% 1.C.8
\label{dualquotientdeExtrP1}
Let $G$ be a locally compact group, ${\rm Extr} (P_1(G))$
the space of indecomposable normalized functions of positive type on $G$,
viewed as a subspace of $L^\infty(G)$ with the weak$^*$-topology,
and $\widehat G$ the dual of $G$, with the Fell topology.
\par

The surjective map
$p \, \colon {\rm Extr} (P_1(G)) \twoheadrightarrow \widehat G$
given by the Gelfand--Naimark construction
(see Construction \ref{constructionGNS2} and Proposition \ref{P1=cyclicmodeqIRR})
is continuous and open.
\end{prop}

\begin{proof}
Recall from Theorem \ref{Banach--Alaoglu--Krein--Milman--Raikov}
that the topology on ${\rm Extr} (P_1(G))$
coincides with the topology of uniform convergence on compact subsets of $G$. 
The continuity of $p$ is therefore an immediate consequence
of Proposition~\ref{Pro-ConvFell}.
\par

To show that $p$ is open, let $U$ be an open subset of ${\rm Extr} (P_1(G))$.
Assume, by contradiction, that $p(U)$ is not open.
Then there exists $\pi\in p(U)$
and a net $(\pi_\iota)_{\iota \in I}$ in $\widehat{G} \smallsetminus p(U)$
such that $\lim_\iota \pi_i = \pi$.
Let $\varphi \in U$ be a normalized function of positive type associated to $\pi$.
By Proposition~\ref{Pro-ConvFell},
there exist functions of positive type $\varphi_\iota$ associated to $\pi_\iota$
such that $\lim_\iota \varphi_\iota = \varphi$.
We have $\varphi_\iota \notin p^{-1}(p(U))$ and therefore $\varphi_\iota\notin U$.
Since $\varphi \in U$ and $U$ is open, this is a contradiction.
\end{proof}

\begin{rem}
% 1.C.9
\label{Rem-FellTopLCOtherDef}
(1)
For a locally compact group $G$, weak containment
can be equivalently defined in terms of C*-algebras
(Proposition \ref{tradweqGC*}, see also Proposition~\ref{Pro-WeakContainmentOperatorNorm}),
and there are two other equivalent definitions of the Fell topology on $\widehat G$.
See \cite[\S\S~3 \& 18]{Dixm--C*},
as well as our Sections \ref{SectionQuasidual} and \ref{C*algLCgroup}.

\vskip.2cm

(2)
The terminology ``Fell topology'' for topological groups and C*-algebras
refers to articles published by Fell
in the early 60s, in particular to \cite{Fell--60a}.
For groups, this topology appears earlier
in Godement's thesis \cite[Section 16]{Gode--48}.
\end{rem}

\subsection*{Functorial properties of Fell's topology}

We examine how Fell's topology behaves with respect to continuous homomorphisms.

\begin{prop}
% 1.C.10
\label{Fell-Homomorphism}
Let $G$ and $H$ be topological groups.
\begin{enumerate}[label=(\arabic*)]
\item\label{iDEFell-Homomorphism}
Let $\alpha \, \colon G \to H$ be a continuous homomorphism with \emph{dense} image.
The induced map 
$$
\widehat \alpha \, \colon \widehat H \to \widehat G,
\hskip.2cm
\pi \mapsto \pi \circ \alpha
$$
is injective and continuous for the Fell topologies on $\widehat H$ an $\widehat G$.
\item\label{iiDEFell-Homomorphism}
Let $\theta$ be a bicontinuous automorphism of the topological group $G$.
The induced map $\widehat \theta \, \colon \widehat G \to \widehat G$ is a homeomorphism.
\item\label{iiiDEFell-Homomorphism}
Let $G$ be a topological group, $N$ be a closed normal subgroup of $G$,
and $p \, \colon G \twoheadrightarrow G/N$ the canonical projection.
Assume that, for every compact subset $\overline{Q}$ of $G/N$,
there exists a compact subset $Q$ of $G$ with $p(Q)=\overline{Q}$.
Then $\widehat p$ is a homeomorphism between $\widehat{G/N}$ 
 and a closed subspace of $\widehat G$. 
 \end{enumerate}
\end{prop}

\begin{proof}
Observe that \ref{iiDEFell-Homomorphism}
is an immediate consequence of \ref{iDEFell-Homomorphism}.

\vskip.2cm

To show \ref{iDEFell-Homomorphism}, let $\pi \in \widehat H$.
Then $\pi \circ \alpha$ is a representation of $G$,
since $\alpha$ is a continuous homomorphism;
moreover, $\pi \circ \alpha$ is irreducible, since $\alpha$ has dense image.
So, $\widehat \alpha \, \colon \widehat H \to \widehat G$ is well defined.
The injectivity of $\widehat \alpha$
is a consequence of the density of the image of $\alpha$.
The fact that $\widehat \alpha$ is continuous with respect to the Fell topologies
is straightforward to check. 

\vskip.2cm

To show \ref{iiiDEFell-Homomorphism},
observe that the image of $\widehat p$ coincides with the subset 
$$
X \, = \, \{\pi \in \widehat G \mid \pi |_N = \mathrm{Id} \}
$$
of $\widehat G$.
It is clear that $X$ is a closed subset of $\widehat G$.
The continuity of the inverse map
$\widehat{p}^{-1} \, \colon X \to \widehat{G/N}$ follows easily from the fact that,
for every compact subset $\overline{Q}$ of $G/N$,
there exists a compact subset $Q$ of $G$ with $p(Q) = \overline{Q}$.
\end{proof}

\begin{rem}
% 1.C.11
\label{Rem-Fell-Homomorphism}
(1)
In Proposition~\ref{Fell-Homomorphism}~\ref{iDEFell-Homomorphism},
the image of $\widehat{\alpha}$ is not necessarily closed,
as shown by the following two examples.
\par

Consider first the injection
$\alpha \, \colon \R \to \T^2, \hskip.2cm t \mapsto (e^{2 \pi i t}, e^{2 \pi i \theta t})$,
where $\theta$ is an irrational number.
The dual $\widehat \R$ of $\R$ consists of the unitary characters
$\chi_k : \R \to \T, \hskip.2cm t \mapsto e^{2 \pi k t}$, for all $k \in \R$,
and the image $\widehat \alpha (\widehat{\T^2})$ is 
$\{ \chi_{m + \theta n} \in \widehat \R \mid m,n \in \Z \}$,
which is a dense countable subgroup of $\widehat \R$.
\par

Consider next the canonical injection $\alpha \, \colon \Q \to \R$,
where $\Q$ is equipped with the discrete topology.
The image $\widehat \alpha(\widehat \R)$ is not closed in $\widehat \Q$.
Indeed, otherwise, $\widehat \alpha(\widehat \R)$
would be a locally compact group and, hence (see \cite[Theorem 5.29]{HeRo--63}), 
the continuous homomorphism $\widehat \alpha$
would be a homeomorphism between $\widehat \R$ and $\widehat \alpha(\widehat \R)$;
since $\widehat \Q$ and therefore $\widehat \alpha(\widehat \R)$ is compact,
this would imply that $\R$, which is homeomorphic to $\widehat \R$, is compact.

\vskip.2cm

(2)
The condition on the pair $(G,N)$ stated
in Proposition~\ref{Fell-Homomorphism}~\ref{iiiDEFell-Homomorphism}
is always satisfied if $G$ is a locally compact group
(see Lemma B.1.1 in \cite{BeHV--08}),
or if $G$ is a complete metrizable topological group
(see Chap.~9, \S~2, Proposition 2 in \cite{BTG5--10}).
However, there are examples of topological groups 
for which this condition is not satisfied (see Remark 1.7.9 in \cite{BeHV--08}).
\end{rem}

\index{Commutator $[a,b] := a^{-1}b^{-1}ab$}
\index{Commutator subgroup = derived group}
\index{Derived group = commutator subgroup}
\index{$b4$@$G_{\rm ab}$ abelianization of the group $G$}
\begin{exe}
% 1.C.12
\label{Example-CommutatorSubgroup}
Let $G$ be a topological group. Let $[G,G]$ denote the derived group of $G$,
also called the commutator subgroup of $G$,
generated by its commutators $[g, h] := g^{-1}h^{-1}gh$, with $g, h \in G$.
Let $\overline{[G,G]}$ be its closure,
and $G_{\rm ab} := G / \overline{[G,G]}$ the Hausdorff abelianized group of $G$.
Assume that $G$ is either locally compact 
or metrizable complete.
\par

Then $\widehat{G_{\rm ab}}$ can be identified with
the closed subspace $\{\pi \in \widehat G \mid \dim \pi = 1\}$ of $\widehat G$.
In particular, every representation of $G$ of dimension $1$
defines a closed point in $\widehat G$.
\end{exe}

\section{Topological properties of the dual of a group}
% Section 1.D
\label{SectionFellPropers}

We review how topological properties of a locally compact group $G$
and of its dual space $\widehat G$ are related.

\subsection*{Properties for general locally compact groups}

Recall that a topological space $X$ is \textbf{quasi-compact}
if every open covering of $X$ contains a finite covering;
equivalently, if every net $(x_\iota)_{\iota \in I}$ 
in $X$ has a convergent subnet $(x_{\jmath})_{\jmath \in J}$.
A topological space $X$ is \textbf{locally quasi-compact} if every 
point in $X$ has a quasi-compact neighbourhood.

\begin{prop}
% 1.D.1
\label{Pro-DualQuasiCompact}
Let $G$ be a locally compact group.
\par

Then $\widehat G$, equipped with Fell's topology, is locally quasi-compact.
\end{prop}

\begin{proof}
Let $\pi \in \widehat G$.
Fix a normalized function of positive type $\varphi$ associated to $\pi$.
By continuity of $\varphi$,
there exists a compact neighbourhood $Q$ of the unit element in $G$
such that $\sup_{g \in Q} \vert \varphi(g) - 1 \vert \le 1/4$. 
\par

Let $\mu_G$ be a Haar measure on $G$.
Choose $f \in C^c(G)$ with support contained in $Q$
and with the following properties:
$f \ge 0$, $\int_G f(g) d\mu_G(g) = 1$, and $f = f^*$.
Observe that $\int_G \psi(g)f(g) d\mu_G(g) \in \R$
for every function of positive type $\psi$ on $G$; indeed:
$$
\begin{aligned}
\overline{ \int_G \psi(g) f(g) d\mu_G(g) }
\, &= \, \int_G \psi(g^{-1}) f(g) d\mu_G(g)
\, = \, \int_G \psi(g) f(g^{-1}) \Delta_G(g^{-1}) d\mu_G(g)
\\
\, &= \, \int_G \psi(g) f^*(g) d\mu_G(g)
\ = \, \int_G \psi(g) f(g) d\mu_G(g) 
\end{aligned}
$$
(for the definitions of $f^*$ and $\Delta_G$, 
and for some of the equalities above, see Appendix~\ref{AppLCG}).
\par

Define $U$ to be the set of $\rho$ in $\widehat G$ for which there exists
an associated normalized function of positive type $\psi$ such that 
$$
\left\vert \int_G (\varphi(g)- \psi(g))f(g) d\mu_G(g) \right\vert \, \le \, 1/4.
$$
The proof consists in showing that $U$ is a quasi-compact neighbourhood of $\pi$.

\vskip.2cm

It is clear that $\pi \in U$.
Let $W = W(\pi; \varphi, Q, 1/4)$ be the set of $\rho$ in $\widehat G$
for which there exists a normalized function of positive type $\psi$ associated to $\rho$
such that 
$$
\vert \varphi(g) - \psi (g) \vert \, \le \, 1/4
\hskip.5cm \text{for all} \hskip.2cm
g \in Q ;
$$
then $W$ is a neighbourhood of $\pi$ for Fell's topology 
which is obviously contained in~$U$.
It follows that $U$ is a neighbourhood of $\pi$.
\par

Let $(\rho_\iota)_{\iota \in I}$ be a net in $U$.
To show that $U$ is quasi-compact,
we have to show that this net has a subnet which converges to some element of $U$.
\par

For every $\iota \in I$, let $\psi_\iota$ be a normalized function of positive type
associated to $\rho_\iota$ such that 
$$
\left\vert \int_G (\varphi(g) - \psi_\iota(g)) f(g) d\mu_G(g) \right\vert \, \le \, 1/4.
$$
Since $P_{\le 1}(G)$ is compact in the weak$^*$-topology,
there exists a subnet $(\psi_{\jmath})_{\jmath \in J}$
which is weak$^*$-convergent to some function $\psi_0 \in P_{\le 1}(G)$.
It is clear that 
$$
\left\vert \int_G (\varphi(g) - \psi_0(g)) f(g) d\mu_G(g) \right\vert \le 1/4.
\leqno(*)
$$
Moreover, since $\sup_{g \in Q} \vert\varphi(g) - 1 \vert \le 1/4$
and $\int_Q f(g) d\mu_G(g) = 1$, we have
$$
\left\vert \int_G \varphi(g)f(g) d\mu_G(g) \right\vert
\, \ge \, \Big\vert \int_Q f(g) d\mu_G(g) \Big\vert 
-
\Big\vert \int_Q ( 1 - \varphi(g) ) f(g) d\mu_G(g) \Big\vert
\, \ge \, 1 - 1/4.
$$
It follows therefore from $(*)$ that $\psi_0 \ne 0$. 
\par

Let $\rho_0$ be the representation of $G$ associated to $\psi_0$ by the GNS construction.
Then $(\rho_{\jmath})_{\jmath \in J}$ converges to $\rho_0$.
Set 
$$
X \, := \, \left\{ \sigma \in \widehat G \mid \sigma \preceq \rho_0 \right\}
$$
(it is the so-called \emph{support} of the representation $\rho_0$).
Let $\sigma \in X$. 
Then $\sigma \preceq \bigoplus_{\jmath' \in J'} \rho_{\jmath'}$
for every subnet $(\rho_{\jmath'})_{\jmath' \in J'}$ of $(\rho_{\jmath})_{\jmath \in J}$. 
Therefore, $(\rho_{\jmath})_{\jmath \in J}$ converges to $\sigma$, 
by Proposition~\ref{weakcont=Felltop}.
It remains to show that $X \cap U \ne \emptyset$.
\par

Assume, by contradiction, that $X \cap U = \emptyset$.
Let $P_X$ be the set of functions in $P_1(G)$ associated
to representations from $X$. Then
$$
\left\vert \int_G (\varphi(g) - \psi(g)) f(g) d\mu_G(g) \right\vert
\, > \, 1/4 
\hskip.5cm \text{for all} \hskip.2cm
\psi \in P_X.
$$
Let $\mathcal{C}$ be the closure of the convex hull of $P_X$ for the weak$^*$-topology.
Then $ P_X$ is the set of extreme points of $\mathcal{C}$ and we have 
$$
\left\vert \int_G (\varphi(g) - \psi(g)) f(g) d\mu_G(g) \right\vert \, \ge \, 1/4 
\hskip.5cm \text{for all} \hskip.2cm
\psi \in \mathcal{C}.
$$
Consider the continuous affine map 
$\Phi \, \colon \mathcal{C} \to \R$ defined by 
$$
\Phi(\psi) \, = \, \int_G (\varphi(g) - \psi(g)) f(g) d\mu_G(g)
\hskip.5cm \text{for} \hskip.2cm
\psi \in \mathcal{C}.
$$
The image of $\Phi$ is a compact interval $\mathopen[ a,b \mathclose]$ of $\R$
with either $a \ge 1/4$ or $b \le -1/4$. 
\par

Let $x \in \{a,b\}$.
Then $\Phi^{-1}(x)$ is a compact convex subset of $\mathcal{C}$.
Every extreme point of $\Phi^{-1}(x)$ must be an extreme point of $\mathcal{C}$,
that is, a point in $P_X$.
Consequently, $\vert x \vert > 1/4$.
This shows that $a > 1/4$ or $b < -1/4$.
Setting $\delta = a$ if $a > 1/4$ and $\delta = -b$ if 
$b < -1/4$, we have therefore
$$
\left\vert \int_G (\varphi(g) - \psi(g)) f(g) d\mu_G(g) \right\vert \, \ge \, \delta \, > \, 1/4
\hskip.5cm \text{for all} \hskip.2cm
\psi \in \mathcal{C}.
$$
\par

Now, $\rho_0$ is weakly equivalent to $\bigoplus_{\sigma \in X} \sigma$
(see Proposition F.2.7 in \cite{BeHV--08});
hence, $\psi_0$ belongs to $\mathcal{C}$.
It follows that 
$$
\left\vert \int_G (\varphi(g) - \psi_0(g))f(g) d\mu_G(g) \right\vert
\, \ge \, \delta \, > \, 1/4 ,
$$
and this is in contradiction with $(*)$.
\end{proof}

\begin{rem}
% 1.D.2
\label{Rem-LocQuasiCompact}
Proposition~\ref{Pro-DualQuasiCompact}
can be deduced from a far more general result
concerning the spectrum of an arbitrary C*-algebra;
see \cite[3.3.8]{Dixm--C*}.
\end{rem}

We will often deal with locally compact groups which are second-countable.
The following proposition shows that this property is inherited by their duals.
(For second-countable locally compact abelian groups,
see also Proposition~\ref{Prop-DualAbelianSecondCountGroup}.)

\begin{prop}
% 1.D.3
\label{Prop-DualSecondCountGroup}
Let $G$ be a second-countable locally compact group.
\par

Then $\widehat G$ is second-countable.
\end{prop}

\begin{proof}
Let $\mu_G$ be a Haar measure on $G$. Since $G$ is second-countable,
$L^1(G, \mu_G)$ is separable and hence the unit ball of $L^\infty(G, \mu_G)$ 
is second-countable for the weak$^*$-topology.
It follows that the subspace ${\rm Extr} (P_1(G))$ of extreme points in 
$P_1(G)$ is second-countable for the weak$^*$-topology.
Proposition~\ref{dualquotientdeExtrP1} shows 
therefore that $\widehat G$ is second-countable.
\end{proof}

\begin{rem}
% 1.D.4
\label{Rem-DualSecondCountable}
Proposition~\ref{Prop-DualSecondCountGroup}
is a special case of Corollary~\ref {Cor-PrimSeparable}
which states that the spectrum of any separable C*-algebra is second-countable.
\end{rem}

\begin{rem}
% 1.D.5
\label{Rem-TopologyDual}
It will be shown in Proposition~\ref{P(A)Baire}
that the dual of a locally compact group is a Baire space.
\end{rem}

\subsection*{Properties for Abelian groups}

Let $G$ be a locally compact abelian group. 
The Fell topology on $\widehat G = \Hom (G, \T)$
is simply the topology of uniform convergence of unitary characters 
on compact subsets of $G$.

\begin{prop}
% 1.D.6
\label{Prop-LocCompGroups}
Let $G$ be a locally compact abelian group.
\par

Then $\widehat G$ is locally compact.
\end{prop}

\begin{proof}
By Proposition~\ref{Pro-DualQuasiCompact}, $\widehat G$ is locally quasi-compact.
Since a net of complex-valued functions can have at most
one limit point, $\widehat G$ is a Hausdorff topological space.
\end{proof}

The next proposition shows that there is a perfect duality between discrete and compact 
abelian groups.

\begin{prop}
% 1.D.7
\label{Prop-AbelianDiscreteCompactGroups}
Let $G$ be a locally compact abelian group.
The following conditions are equivalent:
\begin{enumerate}[label=(\roman*)]
\item\label{iDEProp-AbelianDiscreteCompactGroups}
$G$ is discrete (respectively, compact);
\item\label{iiDEProp-AbelianDiscreteCompactGroups}
$\widehat G$ is compact (respectively, discrete).
\end{enumerate}
\end{prop}

\begin{proof}
It is clear that $\widehat G $ is compact if $G$ is discrete.
% c'est essentiellement Tikhonov !!!
\par

Assume now that $G$ is compact. To show that $\widehat G$ is discrete,
it suffices to show that the unit character $1_G$ is open in $\widehat G$.
For any $z \in \T \smallsetminus \{1\}$,
there exists $n \in \Z$ such that $\vert z^n - 1 \vert \ge \sqrt{2}$.
This implies that,
for every $\chi \in \widehat G \smallsetminus \{ 1_G \}$, we have 
$$
\sup_{g \in G} \vert \chi(g) - 1_G(g) \vert \, \ge \sqrt{2} ,
$$
i.e., $\{ 1_G \}$ is open in $\widehat G$.
\end{proof}

\subsection*{Properties for compact groups}

We show that the dual of a compact group is a discrete space.

\begin{prop}
% 1.D.8
\label{Prop-FellCompGroups}
Let $G$ be a compact group.
\par

Then $\widehat G$, equipped with Fell's topology, is discrete.
\end{prop}

\begin{proof} 
Let $(\pi_\iota)_\iota$ be a net in $\widehat G$ converging to $\pi \in \widehat G$.
 Let $\varphi$ be a normalized function
of positive type associated to $\pi$.
There exists a net $(\psi_\iota)_\iota$
of functions of positive type $\psi_\iota$ associated to $\pi_\iota$ such that 
$\lim_\iota \psi_\iota = \varphi$ uniformly on $G$
(see Proposition~\ref{Pro-ConvFell}).
\par

Let $\mu_G$ be a Haar measure on $G$.
Then 
$$
\lim_\iota \int_G \varphi(x) \overline{\psi_\iota(x)} d\mu_G(x)
\, = \, \int_G \vert \varphi(x) \vert^2 d\mu_G(x) \, \ne \, 0
$$
and, hence, there exists $\iota_0$ such that 
$$
\int_G \varphi(x) \overline{\psi_\iota(x)} d\mu_G(x) \, \ne \, 0
\hskip.5cm \text{for all} \hskip.2cm
\iota \ge \iota_0.
$$
It follows from Schur's orthogonality relations 
(see, for example, \cite[Theorem 5.8]{Foll--16})
that $\pi$ is equivalent to $\pi_\iota$ for all $\iota \ge \iota_0$.
\end{proof}

\begin{rem}
% 1.D.9
\label{Rem-DiscretDual}
(1)
The converse of Proposition~\ref{Prop-FellCompGroups} holds:
let $G$ be a locally compact group such that $\widehat G$ is discrete;
then $G$ is compact.
See \cite[Theorem 1]{Shte--71},
and \cite[Theorem 3.4]{Bagg--72} for second-countable $G$;
see also \cite[Theorem 7.6]{Wang--75}.

\vskip.2cm

(2)
By Propositions \ref{Prop-DualSecondCountGroup} and \ref{Prop-FellCompGroups},
the dual of a second-countable compact group is discrete and countable.
However, there are non-compact locally compact groups with countable dual.
The following example appears in \cite[4.5, 4.6]{Bagg--72},
where it is attributed to Fell, and in \cite[\S~4]{Shte--71}.
For more on locally compact groups with a countable dual, see \cite[Section 6]{Wang--75}.
\end{rem}

\begin{exe}[\textbf{A non-compact locally compact group with countable dual}]
% 1.D.10
\label{Example-CompactNonDiscrete}
Let $p$ be a prime, $\Q_p$ the field of $p$-adic numbers, 
$\Z_p$ its ring of integers, and $\Z_p^\times$ the multiplicative group 
of invertible elements in $\Z_p$.
Consider the natural action of $\Z_p^\times$ by automorphisms of $\Q_p$
given by multiplication.
Let $M_p := \Z_p^\times \ltimes \Q_p$ be the $p$-adic motion group.
We claim that the dual space of $\widehat{M_p}$ is countable.
\index{$m3$@$\Q_p$ $p$-adic numbers}
\index{$m4$@$\Z_p$ $p$-adic integers}
\par

Indeed, identify $\widehat{\Q_p}$ with $\Q_p$ (see Section~\ref{Section-IrrRepBS}).
Under this identification, the dual action of $\Z_p^\times$ on $\widehat{\Q_p}$
corresponds to its natural action on $\Q_p$.
The $\Z_p^\times$-orbits in $\Q_p$ are $\{0\}$
and the sets $O_j$ of $p$-adic numbers of valuation $j \in \Z$.
Moreover, the stabilizer in $\Z_p^\times$
of every point in $\Q_p \smallsetminus \{0\}$ is trivial.
\par

Choose, for every $j \in \Z$, a representative $\chi_j$ for the orbit $O_j$
and let $\pi_j$ be the induced representation $\Ind_{\Q_p}^{M_p} \chi_j$
(compare Section~\ref{Section-IrrIndRep}).
Using the ``Mackey machine" (see Remark \ref{MackMach}),
one can show that
$$
\widehat{M_p} \, = \, \widehat{\Z_p^\times} \cup \{\pi_j \mid j \in \Z\}.
$$
Since $\widehat{\Z_p^\times}$ is countable, $\widehat{M_p}$ is countable.
\par

Observe that $\widehat{M_p}$ is not discrete:
in fact, every $\chi \in \widehat{\Z_p^\times}$
belongs to the closure of $ \{\pi_j \mid j\in \Z\}$.
\end{exe}

\begin{rem}
% 1.D.11
\index{Kazhdan Property (T)}
Duals spaces of groups may have discrete subsets that reflect important properties of the group.
In particular, a locally compact group $G$ has \textbf{Kazhdan Property (T)}
if the class of the unit representation $1_G$ is isolated in $\widehat G$;
this property has shown to be a very significant one \cite{Kazh--67}, \cite{BeHV--08}.
\end{rem}

\subsection*{Properties for discrete groups}

We show that the dual of a discrete group is quasi-compact.

\begin{prop}
% 1.D.12
\label{Pro-DualDiscreteGroup}
Let $\Gamma$ be a discrete group.
\par

Then $\widehat \Gamma$ is a quasi-compact topological space.
\end{prop}

\begin{proof}
Let $(\pi_\iota)_{\iota \in I}$ be a net in $\widehat \Gamma$.
Choose, for every $\iota \in I$,
a normalized function of positive type $\varphi_\iota$ associated to $\pi_\iota$.
Since $P_1(\Gamma)$ is compact
(Theorem \ref{Banach--Alaoglu--Krein--Milman--Raikov}),
there exists a subnet $(\varphi_{\jmath})_{\jmath \in J}$
which converges pointwise on $\Gamma$
to some normalized function of positive type $\varphi$.
Let $\rho$ be the representation associated to $\varphi$ by the GNS construction.
Choose $\pi \in \widehat G$ with $\pi \preceq \rho$.
Then $\pi \preceq \bigoplus_{\jmath' \in J'} \pi_{\jmath'}$ for every subnet $(\pi_{\jmath'})_{\jmath' \in J'}$
of $(\pi_{\jmath})_{\jmath \in J}$.
Therefore $(\pi_{\jmath})_{\jmath \in J}$ converges to $\pi$, 
by Proposition~\ref{weakcont=Felltop}.
\end{proof}

\begin{rem}
% 1.D.13
\label{dualnondenpourgpediscretinfini}
(1)
It will be shown in Proposition \ref{GpeDenDualNonden} that
the dual of a countable infinite discrete group is uncountable.

\vskip.2cm

(2)
As the following example shows, a non discrete locally compact group
may have a quasi-compact dual space.
Another example appears in \cite[4.7]{Bagg--72} and \cite[\S~4]{Shte--71}.
\end{rem}

\begin{exe}[\textbf{A non-discrete locally compact group with quasi-compact dual}]
% 1.D.14
\label{Exe-CompactNonDiscrete}
Let $\K$ be an infinite field.
We view its additive group $(\K, +)$ as a discrete abelian group.
The multiplicative group $H := \K^\times$ acts naturally on $\K$.
Set $N := \widehat \K$ and consider the dual action of $H$ on the infinite compact group $N$.
Let $G$ be the locally compact group $H \ltimes N$,
where $H$ is equipped with the discrete topology.
Observe that $G$ is non discrete and that $N$ is a 
compact and open normal subgroup of $G$.
We claim that $\widehat G$ is quasi-compact.
\par

Indeed, by Pontrjagin duality (see Theorem~\ref{PontDuality17}),
$\widehat N$ can be canonically identified with $\K$,
the action of $H$ on $\widehat N$ corresponding to its natural action on $\K$.
The $H$-orbits on $\widehat N \cong \K$ are $\{1_N\}$
and $\widehat N \smallsetminus \{1_N\}$;
moreover, $H$ acts freely on $\widehat N \smallsetminus \{1_N\}$.
\par

Choose $\chi_0 \in \widehat N \smallsetminus \{1_N\}$
and let $\pi_0$ be the induced representation $\Ind_{N}^G \chi_0$.
Then $\pi_0$ is irreducible (see Theorem~\ref{Theo-IrredInducedRep}) and,
using the Mackey machine
as in Example~\ref{Example-CompactNonDiscrete}, one can show that
$$
\widehat G \, = \, \widehat H \cup \{\pi_0\}.
$$
Since $\widehat H$ is compact,
it is obvious that the topological space $\widehat G$ is quasi-compact. 
\end{exe}

\subsection*{Separation properties of the dual space of a group}

Let $G$ be a locally compact group which is either abelian or compact.
Propositions \ref{Prop-AbelianDiscreteCompactGroups}
and \ref{Prop-FellCompGroups} show in particular that
the topological space $\widehat G$ has the T$_2$ separation property
(that is, $\widehat G$ is a Hausdorff topological space),
hence it has the weaker properties T$_1$ and~T$_0$.
For a reminder on the properties T$_0$, T$_1$, and T$_2$,
see Appendix \ref{AppTop}.
As we now see, the dual space $\widehat G$ of a locally compact group $G$
fails in general to have these separation properties.
\par

\index{LC for ``locally compact''}
We use sometimes \textbf{LC group} as a shorthand for ``locally compact group''.
The statement of the following theorem involves notions
(group of type I, CCR group) to be defined in Chapter \ref{ChapterTypeI}
and its proof uses results stated in Chapter~\ref {ChapterAlgLCgroup}.

\begin{theorem}
% 1.D.15
\label{Theo-SeparationPropertiesDual}
Let $G$ be a second-countable locally compact group. 
Let $\widehat G$ be its dual, equipped with the Fell topology.
\begin{enumerate}[label=(\arabic*)]
\item\label{iDETheo-SeparationPropertiesDual}
$\widehat G$ has property T$_0$ if and only if $G$ is of type I.
\item\label{iiDETheo-SeparationPropertiesDual}
$\widehat G$ has property T$_1$ if and only if $G$ is a CCR group.
\end{enumerate}
\end{theorem}

\begin{proof}
\ref{iDETheo-SeparationPropertiesDual}
Assume that $G$ of type I
and let $\pi$ and $\sigma$ be two distinct points in $\widehat G$.
It follows from Glimm theorem (Theorem~\ref{ThGlimm} below)
that $\pi$ and $\sigma$ are not weakly equivalent.
We may assume that, say, $\pi$ is not weakly contained in $\rho$.
It follows from the definition of the Fell topology
that there exists a neighbourhood $U$ of $\pi$ with $\rho\notin U$.
This shows that $\widehat G$ is a T$_0$-space.
\par

The converse statement is also a part of Glimm theorem;
we prefer to give an independent proof.
Assume that $G$ is not of type I.
By Corollary~\ref{Cor-FacRepPrimitiveG-NonTypI} below,
there exist two (in fact, uncountably many) distinct points $\pi$ and $\rho$ in $\widehat G$
which are weakly equivalent.
Then, every neighbourhood of $\pi$ contains $\rho$
and every neighbourhood of $\rho$ contains $\pi$.
This shows that $\widehat G$ is not a T$_0$-space.
\par

\vskip.2cm

\ref{iiDETheo-SeparationPropertiesDual}
Assume that $\widehat G$ has property T$_1$.
Then, by \ref{iDETheo-SeparationPropertiesDual}, $G$ is of type I.
Let $\pi$ be a point in $\widehat G$.
By Glimm theorem, $\pi(C^*_{\rm max}(G))$ contains
the ideal $\Ki(\Hi)$ of compact operators on the Hilbert space $\Hi$ of $\pi$,
where $C^*_{\rm max}(G)$ denotes the maximal C*-algebra of $G$,
to be introduced in Chapter~\ref{ChapterAlgLCgroup}.
We have to show that $\pi(C^*_{\rm max}(G))$ coincides with $\Ki(\Hi)$.
\par

Assume, by contradiction that $\Ki(\Hi)$ is a proper ideal of $\pi(C^*_{\rm max}(G))$. 
Then $I := \pi^{-1}(\Ki(\Hi))$ is a proper ideal of $C^*_{\rm max}(G)$.
So, the quotient $A := C^*_{\rm max}(G)/I$ is a non-zero C*-algebra.
Let $\rho$ be a (non-zero) irreducible representation of $A$.
Then $\rho \circ p$ is an irreducible representation of $C^*_{\rm max}(G)$,
where $p \,\colon C^*_{\rm max}(G) \twoheadrightarrow A$ is the canonical homomorphism.
Denote again by $\rho$ the representation of $G$
corresponding to $\rho \circ p$ (see Section~\ref{C*algLCgroup}).
Since $\textnormal{C*ker}(\rho)$ obviously contains $\textnormal{C*ker}(\pi)$,
the representation $\rho$ is weakly contained in $\pi$
(see Proposition~\ref{tradweqGC*}). 
However, $\{ \pi \}$ is a closed point in $\widehat G$, since $\widehat G$ is a T$_1$-space.
Therefore, $\rho$ is weakly equivalent to $\pi$;
this means that $\textnormal{C*ker}(\rho) = \textnormal{C*ker}(\pi)$. 
Since $I \subset \textnormal{C*ker}(\rho)$,
it follows that $I \subset \textnormal{C*ker}(\pi)$.
This is a contradiction, since $\Ki(\Hi)$ is contained in $\pi(C^*_{\rm max}(G))$.
\par

Conversely, assume that $G$ is a CCR group. Let $\pi \in \widehat G$.
Then $\pi$ is a CCR representation and hence $\{ \pi \}$ is a closed point in $\widehat G$,
by Corollary~\ref{Cor-RepCCR} below.
\end{proof}

\begin{rem}
% 1.D.16
\label{closedpointsdual}
(1)
Let $G$ be a locally compact group.
Every one-dimensional representation of $G$ defines a closed point of $\widehat G$,
as observed in Example \ref{Example-CommutatorSubgroup}.
In particular $\widehat G$ always contains a closed point, $1_G$.
Compare with Remark \ref{Rem-TopPropPrimitiveSpace}.
\par

More generally,
every finite dimensional irreducible representation $\pi$ of $G$
defines a point in the dual $\widehat G$ which is closed.
Indeed, such a representations is CCR (see Remark \ref{Rem-CCR-Rep})
and therefore $\{\pi\}$ is closed in $\widehat G$ (Corollary \ref{Cor-RepCCR}).

\vskip.2cm

(2)
The converse statement to (1) holds for countable groups:
let $\Gamma$ be a countable group and $\pi$ an irreducible representation of $\Gamma$
such that $\{\pi\}$ is a closed point in $\widehat \Gamma$. Then $\pi$ is finite dimensional.
\par

Indeed, the spectrum (see Section~\ref{SectionPrimC*}) of 
the C*-algebra $\mathcal A$ generated by $\pi(G)$ 
consists only of the singleton $\{\pi\}$.
This implies that $\mathcal A$ is isomorphic
to the algebra of compact operators on a Hilbert space (see \cite[4.7.3]{Dixm--C*});
since $\mathcal A$ is unital, it follows that $\mathcal A$ is finite dimensional.
As $\pi$ is irreducible, the claim is proved.

\vskip.2cm

(3)
Let $\Gamma$ be a discrete group. The following five conditions are equivalent
(type~I groups and CCR groups are defined in Chapter \ref{ChapterTypeI}):
\begin{enumerate}[label=(\roman*)]
\item\label{iDEclosedpointsdual}
$\Gamma$ has an abelian subgroup of finite index,
\item\label{iiDEclosedpointsdual}
the group $\Gamma$ is type I,
\item\label{iiiDEclosedpointsdual}
the dual space $\widehat \Gamma$ is T$_0$,
\item\label{ivDEclosedpointsdual}
the group $\Gamma$ is CCR,
\item\label{vDEclosedpointsdual}
the dual space $\widehat \Gamma$ is T$_1$.
\end{enumerate}
The equivalence
\ref{iDEclosedpointsdual} $\Longleftrightarrow$ \ref{iiDEclosedpointsdual}
is Thoma Theorem \ref{discretTypes-Thoma}.
Equivalences \ref{iiDEclosedpointsdual} $\Longleftrightarrow$ \ref{iiiDEclosedpointsdual}
and \ref{ivDEclosedpointsdual} $\Longleftrightarrow$ \ref{vDEclosedpointsdual}
are those of Theorem \ref{Theo-SeparationPropertiesDual}.
The implication \ref{iDEclosedpointsdual} $\Rightarrow$ \ref{ivDEclosedpointsdual}
follows from Theorem \ref{Theo-MooreGroups}~\ref{4DElistegroupesMoore} quoted below.
The implication \ref{vDEclosedpointsdual} $\Rightarrow$ \ref{iiiDEclosedpointsdual}
is obvious.
\end{rem}

\begin{rem} 
% 1.D.17
We do not know of a characterization
of the class ${\mathcal{HD}}$ of LC groups with a Hausdorff dual space.
Here is a sample of results concerning this class of groups,
\par

--- Every central-by-compact group belongs to ${\mathcal{HD}}$
(see \cite[Corollary 5.2]{GrMo--71b}).
A LC group $G$ is called \textbf{central-by-compact}
if $G/Z$ is compact, where $Z$ is the centre of $G$.
The class of central-by-compact groups includes
direct products of a compact group and an abelian LC group;
the next example shows that this inclusion is strict.
\par

--- The set $\Z^2 \times \{ \pm 1 \}$ with the multiplication defined by
$(m_1, m_2, \varepsilon) (n_1, n_2, \delta)
= (m_1+n_1, m_2+n_2, \varepsilon \delta (-1)^{m_1n_2})$
is a discrete group in $\mathcal{HD}$ with centre
$\{ (m_1, m_2, 1) \in \Z^2 \times \{ \pm 1 \} \mid m_1 \equiv m_2 \equiv 0 \pmod{2} \}$
of index $4$.
The group is not a direct product of this centre by a finite group.
This example appears in \cite[Page 701]{Palm--78}.
\par

--- A countable discrete group belongs to ${\mathcal{HD}}$
if and only if it is central-by-compact (see \cite[Proposition~6]{Raeb--82}).
\par

--- A connected LC group $G$ belongs to ${\mathcal{HD}}$ if and only if
$G$ is compact-by-abelian (see \cite{BaSu--81}).
\par

--- Let $G$ be a semi-direct product $K \ltimes A$,
with $K$ a compact group and $A$ a locally compact abelian group.
A necessary and sufficient condition 
for $G$ to belong to ${\mathcal{HD}}$ is given in \cite[Section 9]{Bagg--68};
see also \cite[Proposition 5.62]{KaTa--13}.
\end{rem}

\section
{Primitive dual of a topological group}
% Section 1.E
\label{PrimIdealSpace}
As was seen in Theorem~\ref{Theo-SeparationPropertiesDual},
 the dual space $\widehat G$ of a LC group $G$ does not have the 
 T$_0$ separation property in general.
Given an arbitrary topological $X$, one can associate to $X$
a canonical topological space with the T$_0$ separation property.
When applied to $\widehat G$, this construction leads to the following definition
(see Reformulation \ref{refdefprim} below).

 \begin{defn}
 % 1.E.1
 \label{Def-PrimitiveDual}
 Let $G$ be a topological space.
The \textbf{primitive dual} $\Pri(G)$ 
is the space of weak equivalence classes of irreducible representations of $G$,
with the topology defined as follows.
The natural surjective map
\begin{equation}
\label{candualprim}
\tag{$\kappa$1}
\kappa^{\rm d}_{\rm prim} \, \colon \, \widehat G \twoheadrightarrow \Pri(G) 
\end{equation}
associates to an equivalence class of irreducible representations
its weak equivalence class
($\kappa$ stands for $\kappa \alpha \nu o \nu \iota \kappa o \sigma$,
i.e., for canonical).
The topology on $\Pri(G)$ is defined as the quotient by (\ref{candualprim}) 
of the Fell topology on the dual $\widehat G$.
\index{Primitive dual! $1$@of a topological group}
\index{$c1$@$\Pri(G)$ primitive dual of the topological group $G$}
\index{Dual! $2$@surjection on primitive dual}
\end{defn}

\begin{rem}
% 1.E.2
For a locally compact group $G$,
we will give in \ref{C*algLCgroup}
an equivalent definition of the topological space $\Pri(G)$
as the space of primitive ideals of the maximal C*-algebra of $G$,
equipped with the Jacobson topology.
When $G$ is moreover $\sigma$-compact, 
we will extend the map $\kappa^{\rm d}_{\rm prim}$
to the quasi-dual $\QD(G)$ of $G$
defined in Section \ref{SectionQuasidual}.
\end{rem}

We reformulate the definition of $\Pri(G)$ and its topology 
in terms of the notion of universal Kolmogorov quotient of a topological space,
i.e., of the largest T$_0$ quotient of a topological space (see Appendix \ref{AppTop}).

\begin{refo}
% 1.E.3
\label{refdefprim}
Let $G$ be a topological group.
\par

The universal Kolmogorov quotient of the dual $\widehat G$
is the pair made of the primitive dual $\Pri(G)$
and the natural surjective map $\widehat G \twoheadrightarrow \Pri(G)$.
\end{refo}

Observe that the map $\widehat G \twoheadrightarrow \Pri(G)$
is in general far from being injective.
Indeed, when $G$ is a second-countable LC group,
this map is injective (and therefore bijective) if and only if $G$ is of type I.
More on this in Section \ref{noninjectivity}.

\begin{rem}
% 1.E.4
\label{Rem-TopPropPrimitiveSpace}
(1)
Let $G$ be a topological group.
By definition (see \ref{refdefprim}), the primitive dual $\Pri(G)$ is a T$_0$-space.
% voir aussi [Palmer, page 691] 

\vskip.2cm

(2)
We will determine in Chapter \ref{ChapterPrimExa}
the primitive duals of some discrete groups.
In most cases, these duals are not T$_1$ spaces.

\vskip.2cm

(3)
Here are a few results concerning the T$_1$ separation property for primitive duals
for some classes of groups.

--- Let $\Gamma$ be a finitely generated solvable group.
Then $\Pri(\Gamma)$ is a T$_1$ space if and only if
$\Gamma$ is virtually nilpotent (see \cite[Theorem 5]{MoRo--76}).

--- Discrete nilpotent groups have T$_1$ primitive duals (see \cite{Pogu--82}).

--- Points in primitive duals of LC-groups need not be locally closed
(i.e., need not be closed in any neighbourhood of themselves).
For example, the primitive dual of a non-abelian free group of finite rank
does not have any locally closed point 
(remark following Proposition 6.1 in \cite{MoRo--76}).
\index{Free group}

--- If $G$ is a connected LC group, all points of $\Pri(G)$ are locally closed.
Moreover, $\Pri(G)$ is a T$_1$ space if and only if
$G$ is a projective limit of connected Lie groups,
each of which is type $R$ on its radical
\cite[Theorem 1]{MoRo--76}.
(A connected Lie group $G$ is of type $R$ on its radical 
if the adjoint action of $G$ on the radical of its Lie algebra is distal.)
\par

We refer to Page 726 of the review article by Palmer \cite{Palm--78}
for more on locally compact groups $G$ for which $\Pri(G)$ is a T$_1$-space.

\vskip.2cm

(4)
Concerning the Hausdorff property T$_2$ for primitive duals, we quote only the following result.
Let $G$ be a LC group $G$ in which every conjugacy class has a compact closure;
then the space $\Pri(G)$ is Hausdorff.
See \cite{LiMo--74} for $\sigma$-compact groups,
and \cite{Kani--79} in general.

\vskip.2cm

(5)
Let $G$ be a locally compact group.
We will show in Proposition~\ref{P(A)Baire} that $\Pri{(G)}$ is a Baire space.
\end{rem}

We will continue our discussion of primitive duals,
both for locally compact groups and for C*-algebras,
in Chapter \ref{ChapterGroupMeasureSpace}.

\section
[Induced representations]
{Induced representations, irreducibility, and equivalence}
% Section 1.F
\label{Section-IrrIndRep}

The construction of irreducible representations of a group 
is often based on induction of representations, 
which is a procedure to build representations of a group 
out of representations of its subgroups.
\par

To avoid technical problems, we will restrict ourselves 
in this section to the case of \emph{open} subgroups.
In fact, given a LC group $G$, a \emph{closed} subgroup $H$ of $G$,
and a representation $\sigma$ of $H$,
an induced representation $\Ind_H^G \sigma$ of $G$ can always be defined. 
For this more general construction, 
see \cite{Mack--52}, \cite{Mack--58}, \cite{Foll--16}, 
\cite{BeHV--08} and \cite{KaTa--13}.
\par

In the particular case of inducing a \emph{one-dimensional} representation $\sigma$
of an open subgroup $H$ of $G$,
we will state and prove sufficient conditions under which $\Ind_H^G \sigma$ is irreducible,
and under which two such representations are not equivalent.

\vskip.2cm

\index{Transversal $T$, and $G = \bigsqcup_{t \in T} tH$}
Let $G$ be a topological group, and $H$ an open subgroup of $G$.
A \textbf{left transversal} for $H$ in $G$ is a subset $T$ of $G$
such that $G$ is the disjoint union $\bigsqcup_{t \in T} tH$.
It is sometimes convenient, but it is not necessary, to choose $e \in T$
as a representative of the class $H$ in $H \backslash G$.
\par

Similarly, a \textbf{right transversal} for $H$ in $G$ is a subset $T$ of $G$
such that $G = \bigsqcup_{t \in T} Ht$.
Observe that, if $T$ is a left transversal for $H$ in $G$,
then the set $T^{-1}$ of inverses of elements in $T$
is a right transversal for $H$ in $G$.

\begin{defn}
% 1.F.1
\label{openeasierthanclosed}

Let $G$ be a topological group,
$H$ an open subgroup,
and $(\sigma, \Ki_\sigma)$ a representation of $H$.
An \textbf{induced representation} of $\sigma$ from $H$ to $G$
is a representation $\pi$ of $G$ on some Hilbert space $\Hi$
that satisfies the following conditions:
\begin{enumerate}[label=(\arabic*)]
\item\label{iDEopeneasierthanclosed}
there exists a closed subspace $\Ki$ of $\Hi$ that is invariant by $\pi(H)$;
\item\label{iiDEopeneasierthanclosed}
the representation $\Res^G_H (\pi) \vert_{\Ki} \, \colon H \to \U(\Ki)$ 
obtained by restriction of $\pi$ to $H$ and $\Ki$ is equivalent to $\sigma$;
\item\label{iiiDEopeneasierthanclosed}
the space of $\pi$ is the orthogonal Hilbert sum $\Hi = \bigoplus_{t \in T} \pi(t) (\Ki)$,
for some left transversal $T$ of $H$ in $G$.
\end{enumerate}
\index{Induced representation}
\index{Representation! induced}
\end{defn}

An induced representation $\pi$ as above
is often denoted by $\Ind_H^G \sigma$.

\begin{rem}
% 1.F.2
\label{Rem-DefindRep}

Note that, in view of \ref{iDEopeneasierthanclosed}, 
$\pi(t) (\Ki)$ only depends on the left coset in $G$ of $t \in T$ 
and, hence, if \ref{iiiDEopeneasierthanclosed} holds for one left transversal $T$ for $H$,
it holds for any such transversal.
Moreover, $G$ acts in the obvious way transitively on the set
$$
\{ \Li \mid \Li = \pi(t)(\Ki) \hskip.2cm \text{for some} \hskip.2cm t \in T \}
\, = \,
\{ \Li \mid \Li = \pi(g)(\Ki) \hskip.2cm \text{for some} \hskip.2cm g \in G \}
$$
and the stabilizer of $\Ki$ is $H$.
\par

Let us check that 
\begin{enumerate}[label=(\arabic*)]
\addtocounter{enumi}{3}
\item\label{ivDEopeneasierthanclosed}
an induced representation of $\sigma$
is \emph{unique} up to equivalence, 
\end{enumerate}
that is, 
if $\pi_1, \pi_2$ are representations of $G$ which satisfy
Conditions \ref{iDEopeneasierthanclosed} to \ref{iiiDEopeneasierthanclosed}, 
then $\pi_1$ and $\pi_2$ are equivalent.
\par

Indeed, for $j = 1, 2$, 
there exist by \ref{iDEopeneasierthanclosed} and \ref{iiDEopeneasierthanclosed}
a closed $\pi_j(H)$-invariant subspace $\Ki_j$ 
such that $\Res^G_H (\pi_j) \vert_{\Ki_j}$ 
is equivalent to $\sigma$. 
Therefore there exists a surjective isometry $V \, \colon \Ki_1 \to \Ki_2$ which intertwines 
$\Res^G_H (\pi_1) \vert_{\Ki_1}$ 
and $\Res^G_H (\pi_2) \vert_{\Ki_2}$.
We claim that the map $U \, \colon \Hi_1 \to \Hi_2$ given by 
$$
U \pi_1(g)\eta \, = \, \pi_2(g) V \eta
\hskip.5cm \text{for} \hskip.2cm
g \in G, \, \eta \in \Ki_1
$$
is well-defined; indeed, if $\pi_1(g) \eta = \pi_1(g') \eta'$
for $g' \in G$ and $\eta' \in \Ki_1$, then $g' = gh$ and $\eta' = \pi_1(h^{-1}) \eta$
for some $h \in H$ and hence
$$
\pi_2(g') V(\eta') \, = \, \pi_2(gh) V (\pi_1(h^{-1}) \eta)
\, = \, \pi_2(gh) \pi_2(h^{-1}) V (\eta) \, = \, \pi_2(g) V(\eta).
$$
Then $U$ is an isometry, and it is surjective because its image
contains $\pi2(t)(\Ki_2)$ for all $t \in T$.
Finally, $U$ intertwines $\pi_1$ and $\pi_2$;
indeed, let $t \in T$ and $\xi \in \pi_1(t)(\Ki_1)$,
say $\xi = \pi_1(t)\eta$ for $\eta \in \Ki_1$; for $g \in G$,
$$
U \pi_1(g) \xi \, = \, U \pi_1(gt) \eta \, = \,
\pi_2(g) \pi_2(t) V \eta \, = \, \pi_2(g) U \pi_1(t) \eta \, = \
\pi_2(g) U \xi .
$$

\vskip.2cm

Some properties of $\Ind_H^G \sigma$ are straightforward consequences of the definition.
For example:
\begin{enumerate}[label=(\arabic*)]
\addtocounter{enumi}{4}
\item\label{vDEopeneasierthanclosed}
The dimension of $\Ind_H^G \sigma$ is $\vert G/H \vert \dim (\sigma)$.
In particular $\Ind_H^G \sigma$ is finite dimensional if and only if
$H$ of finite index in $G$ and $\sigma$ is finite dimensional.
\item\label{viDEopeneasierthanclosed}
If $N$ is a normal subgroup of $G$ contained in $H$
and such that $\sigma(n) = \mathrm{Id}_\Ki$ for all $n \in N$,
then $\left(\Ind_H^G \sigma \right) (n) = \mathrm{Id}_\Hi$ for all $n \in N$.
\end{enumerate}
\end{rem}

\begin{rem}
% 1.F.3
\label{remdefind}
The characterization of induced representations by 
Conditions \ref{iDEopeneasierthanclosed} to \ref{iiiDEopeneasierthanclosed} 
of \ref{openeasierthanclosed}
is valid for open subgroups only.
In the more general case of induction from closed subgroups,
these conditions need not be satisfied.
Indeed, consider for instance the affine 
group of the real line $\Aff(\R) = \R^\times \ltimes \R$,
with its natural locally compact topology (see Section~\ref{Section-IrrRepAff}).
Let $\pi$ be the induced representation $\Ind_N^{\Aff(\R)} \chi$,
for the normal subgroup $N = \{ (1,b) \mid b \in \R \}$
and the unitary character of $N$ given by $\chi(1,b)= e^{ i b}$.
Then (see Remark~\ref{AffKLie})
$\pi$ can be realized on $L^2(\R^\times, dt/\vert t \vert)$ by the formula
$$
\pi(a,b) \xi (t)
\, = \,
\exp \left( -2 \pi i bt \right) \xi(at)
$$
for all $\xi \in L^2(\R^\times, dt / \vert t \vert)$.
Using Fourier transform on $N \approx \R$,
we see that the restriction of $\pi$ to $N$
is equivalent to the regular representation $\lambda_N$;
in particular, $\pi|_N$ contains no finite dimensional subrepresentation.
\par

In the general case of inducing from closed subgroups of LC groups,
there is however a characterization of induced representations in terms of
``transitive systems of imprimitivity''.
We refer to \cite[\S~6]{Mack--58}, as well as to \cite[Chap.~6]{Foll--16}.
\end{rem}

\begin{constr}
% 1.F.4
\label{constructionInd}
Let $G$ be a topological group, $H$ an \emph{open} subgroup of $G$, 
and $(\sigma, \Ki)$ a representation of $H$.
In in the literature, there are several models
for the induced representation $\Ind_H^G \sigma$
as defined in \ref{openeasierthanclosed}.
We are going to describe three of them.

\vskip.2cm

(1)
Consider the Hilbert space $\Hi$ of maps $f \, \colon G \to \Ki$ with the following properties
\begin{enumerate}[label=(\alph*)]
\item\label{aDEconstructionInd}
$f(hx) = \sigma(h) f(x)$ for all $h \in H, \hskip.1cm x \in G$;
\item\label{bDEconstructionInd}
$\sum_{Hx \in H \backslash G} \Vert f(x) \Vert^2 < \infty$.
\end{enumerate}
Observe that, in view of \ref{aDEconstructionInd}
and the unitarity of $\sigma$, 
the number $\Vert f(x) \Vert^2$ only depends on the coset $Hx$;
therefore, the sum in \ref{bDEconstructionInd},
which is a sum over one $x$ for each coset in $H \backslash G$,
is well defined.
Observe also that, for $g \in G$ and $f \in \Hi$,
the function $G \to \Ki, \hskip.2cm x \mapsto f(xg)$
is again in $\Hi$.
The induced representation
$\pi := \Ind_H^G \sigma$
is the representation of $G$ defined in $\Hi$ by
$$
(\pi(g) f) (x) \, = \, f(xg) 
\hskip.5cm \text{for all} \hskip.2cm 
g \in G, \hskip.1cm f \in \Hi
\hskip.2cm \text{and} \hskip.2cm 
x \in G.
$$
Using another model for $\pi$, we will verify below that
the homomorphism 
$$
G \to \U(\Hi), \hskip.2cm g \mapsto \pi(g)
$$
is indeed continuous.
\par

Let us check that $\pi$ as defined above is indeed 
an induced representation of $\sigma$
in the sense of Definition \ref{openeasierthanclosed}.
Let $T$ be a left transversal for $H$ in $G$;
denote by $t_0$ the element in $T \cap H$.
For every $t \in T$, let $\Hi_t$ be the space
of functions $f \, \colon G \to \Ki$ with support contained in $Ht^{-1}$ 
(in particular of functions with support in $H$ when $t = t_0$).
Then $\Hi_{t_0}$ is a closed $\pi(H)$-invariant subspace of $\Hi$.
Moreover $\Hi_t$ is closed in $\Hi$ for all $t \in T$,
and $\Hi = \bigoplus_{t \in T} \Hi_t$.
This shows that Conditions \ref{iDEopeneasierthanclosed}
and \ref{iiiDEopeneasierthanclosed}
of Definition \ref{openeasierthanclosed} are satisfied.
For $\xi \in \Ki$, define a function $f_\xi \, \colon G \to \Ki$
by $f_\xi(g) = \sigma(t_0^{-1}g) \xi$ if $g \in H$
and $f_\xi(g) = 0$ if $g \in G \smallsetminus H$.
Then $f_\xi \in \Hi_{t_0}$, and $\xi \mapsto f_\xi$
defines an isometry $U$ of $\Ki$ onto $\Hi_{t_0}$;
moreover, $U \sigma(h) \xi = \pi(h) U \xi$
for all $h \in H$ and $\xi \in \Ki$, 
so that Condition \ref{iiDEopeneasierthanclosed}
of Definition \ref{openeasierthanclosed}
is also satisfied, as was to be shown.

\vskip.2cm

(2)
It will be convenient for us to work with another model
for induced representations, as follows.
Choose a right transversal $T$,
so that we have a disjoint union $G = \bigsqcup_{t \in T} Ht$.
For every $g \in G$ and $t \in T$, 
there are uniquely defined elements
$$
\alpha(t, g) \, \in \, H
\hskip.2cm \text{and} \hskip.2cm
t \cdot g \, \in \, T
\hskip.2cm \text{such that} \hskip.2cm
tg \, = \, \alpha(t, g) (t \cdot g) .
\leqno(\ast)
$$
It is routine to check that
the map $\alpha \, \colon T \times G \to H$ satisfies
$$
\alpha(t, g_1 g_2) \, = \, \alpha(t, g_1) \alpha(t \cdot g_1, g_2)
\hskip.5cm \text{for all} \hskip.2cm
t \in T
\hskip.2cm \text{and} \hskip.2cm
g_1, g_2 \in G ,
$$
i.e., it is a cocycle, see (\ref{eqq/RNcocR}) in Section \ref{Section-AllIrredRep} below,
\index{Cocycle!}
and that
$$
T \times G \to T, \hskip.2cm (t, g) \mapsto t \cdot g
$$
is an action of $G$ on $T$ \emph{on the right},
which is equivalent to the natural action of $G$ on $H \backslash G$ 
given by right multiplication.
\par

Let $\ell^2(T, \Ki)$ denote the Hilbert space 
of functions from $T$ to $\Ki$ such that
\hfill\par\noindent
$\sum_{t \in T} \Vert f(t) \Vert^2 < \infty$, with scalar product
$\langle f_1 \mid f_2 \rangle := \sum_{t \in T} \langle f_1(t) \mid f_2(t) \rangle_{\Ki}$.
The induced representation $\pi = \Ind_H^G \sigma$ of $G$
is defined on $\ell^2(T, \Ki)$ by
\index{$h2$@$\ell^2(T, \Ki)$ Hilbert space}
$$
(\pi(g) f) (t) \, = \, \sigma (\alpha(t, g)) (f(t \cdot g)) 
\hskip.5cm \text{for all} \hskip.2cm 
g \in G, \hskip.1cm f \in \ell^2(T, \Ki)
\hskip.2cm \text{and} \hskip.2cm 
t \in T.
$$
We leave it to the reader to check that $\pi$ satisfies
Conditions \ref{iDEopeneasierthanclosed} to \ref{iiiDEopeneasierthanclosed}
of Definition \ref{openeasierthanclosed}.
\par

We now show that $\pi$ is indeed a representation of $G$,
i.e., that the map
$$
G \to \ell^2(T, \Ki), \hskip.2cm g \mapsto \pi(g)f
$$
is continuous for every $f \in \ell^2(T, \Ki)$.
\par

For $t \in T$ and $\xi \in \Ki$,
let $f_{t, \xi}$ be the element of $\ell^2(T, \Ki)$
defined by $f_{t, \xi} (t) = \xi$ and $f_{t, \xi} (s) = 0$ for $s \in T, s \ne t$.
By definition of $\pi$, we have
$$
\pi(g) f_{t, \xi} \, = \, f_{t \cdot g^{-1}, \sigma( \alpha(t \cdot g^{-1}, g)) \xi}
$$
for $g \in G$.
Since the linear span of $\{f_{t, \xi} \mid t \in T, \hskip.1cm \xi \in \Ki \}$
is dense in $\ell^2(T, \Ki)$, it suffices to prove that 
the map $g \mapsto \pi(g)f_{t, \xi}$ is continuous 
for $t \in T$ and $\xi \in \Ki$.
Moreover, as $\pi$ is a homomorphism, 
it suffices to prove the continuity of this map at the unit element $e$.
\par

Consider now a fixed pair of $t \in T$ and $\xi \in \Ki$, and $\varepsilon > 0$.
Since $\sigma$ is a representation of $H$, 
there exists a neighbourhood $U$ of $e$ in $H$ such that 
$$
\Vert \sigma(h) \xi - \xi \Vert \, < \, \varepsilon 
\hskip.5cm \text{for all} \hskip.2cm 
h \in U.
$$
Since $H$ is open, $U$ is a neighbourhood of $e$ in $G$ and we can find
a neighbourhood $V$ of $e$ such that 
$V \subset U$ and $t^{-1}Vt \subset U$.
For $g \in t^{-1}Vt$, we have $tg \in V t \subset Ut \subset Ht$,
so that $t \cdot g = t$ and $\alpha(t, g) \in U$,
hence $\pi(g)f_{t, \xi} = f_{t, \sigma(\alpha(t, g))\xi}$.
It follows that
$$
\Vert \pi(g)f_{t, \xi} - f_{t, \xi} \Vert \, = \,
\Vert \sigma(\alpha(t, g)) \xi - \xi \Vert \, < \, \varepsilon
$$
for every $g \in t^{-1}Vt$.

\vskip.2cm

(3)
Occasionally
(see Construction \ref{constructionInd+mu}),
we will prefer to realize $\Ind_H^G \sigma$ on $\ell^2(T, \Ki)$,
where $T$ is now a \emph{left} transversal for $H$ in $G$.
For $g \in G$ and $t \in T$, define $g \cdot t \in T$ and $\beta(g, t) \in H$
by the equality $gt = (g \cdot t)\beta(g, t)$.
The map $\beta \, \colon G \times T \to H$ satisfies the cocycle relation
$$
\beta(g_1g_2, t) \, = \, \beta(g_1, g_2 \cdot t) \beta(g_2, t)
\hskip.5cm \text{for all} \hskip.2cm
g_1, g_2 \in G
\hskip.2cm \text{and} \hskip.2cm
t \in T ,
$$
see (\ref{eqq/RNcocL}) in Section \ref{Section-AllIrredRep} below.
\par

The induced representation $\pi = \Ind_H^G \sigma$ is defined by
$$
(\pi(g) f) (t) \, = \, \sigma \left( \beta(g^{-1}, t)^{-1} \right) \left( f (g^{-1} \cdot t) \right) 
\hskip.5cm \text{for all} \hskip.2cm 
g \in G, \hskip.1cm f \in \ell^2(T, \Ki)
\hskip.2cm \text{and} \hskip.2cm 
t \in T.
$$
\end{constr}

\begin{exe}
% 1.F.5
\label{quasiregularrep}
The \textbf{quasi-regular representation} of a topological group $G$
associated to an open subgroup $H$
is the induced representation $\Ind_H^G 1_H$.
\par

It can be realized as the natural representation of $G$ on $\ell^2(H \backslash G)$
 given by right translations,
or equivalently as the natural representation of $G$ on $\ell^2(G/H)$
given by left translations,
and then usually denoted by $\lambda_{G/H}$.
The intertwining isometry $S$ from $\ell^2(H \backslash G)$ onto $\ell^2(G/H)$
is defined by $(S f) (xH) = f(H x^{-1})$ for all $f \in \ell^2(H \backslash G)$ and $x \in G$.
\index{Quasi-regular representation}
\index{Representation! quasi-regular}
\par

When $G$ is locally compact, as already mentioned above,
there is also a quasi-regular representation $\Ind_H^G 1_H$
for every \emph{closed} subgroup $H$ of $G$.
We refer to \cite[Appendix B.1]{BeHV--08},
or to Definition \ref{defquasiregG/K}
for the particular case of a compact subgroup $H$ in $G$.
\end{exe}

We state here the following properties of the induction process for open subgroups,
but they hold just as well for closed subgroups of locally compact groups.
See, e.g., \cite[Appendix E]{BeHV--08}.

\begin{prop}
% 1.F.6
\label{InductionQqPropr}
Let $G$ be a topological group, $H$ an open subgroup of $G$,
and $K$ an open subgroup of $H$.
Let $\pi$ be a representation of $G$;
let $\sigma_1, \sigma_2, \hdots$ be representations of $H$,
and $\tau$ a representation of $K$.
\begin{enumerate}[label=(\arabic*)]
\item\label{iDEInductionQqPropr}
If $\sigma_1$ and $\sigma_2$ are equivalent representations of $H$,
then $\Ind_H^G \sigma_1$ and $\Ind_H^G \sigma_2$
are equivalent representations of $G$.
\item\label{iiDEInductionQqPropr} 
$\Ind_H^G (\bigoplus_{i \ge 1} \sigma_i)$ and
$\bigoplus_{i \ge 1} \big( \Ind_H^G \sigma_i \big)$
are equivalent.
\item\label{iiiDEInductionQqPropr} 
$\Ind_H^G (\Ind_K^H \tau)$
and $\Ind_K^G \tau$ are equivalent.
\item\label{ivDEInductionQqPropr}
$\Ind_H^G (\Res^G_H \pi)$ is equivalent to
the tensor product $\pi \otimes \lambda_{G/H}$
of $\pi$ with the quasi-regular representation $\lambda_{G/H}$.
\end{enumerate}
\end{prop}

Next, we identify the conjugate representation of an induced representation
by a bicontinuous automorphism of $G$.
\par

If $(\sigma, \Ki)$ is a representation of a subgroup $H$ of $G$ and $\theta \in \Aut(G)$,
we denote by $\sigma^\theta$ the representation
of the subgroup $\theta^{-1}(H)$ on $\Ki$ defined by 
$$
\sigma^\theta(x) \, = \, \sigma (\theta(x))
\hskip.5cm \text{for all} \hskip.2cm
x \in \theta^{-1}(H).
$$
Observe that, in case $H = G$,
this representation is the conjugate representation $\sigma^\theta$ of $G$
discussed in Remark~\ref{Rem-AutoRep}.
\index{Conjugate representation} 
\index{Representation! conjugate}

\begin{prop}
% 1.F.7
\label{PropAutondRep}
Let $G$ be a topological group, 
$H$ an open subgroup,
and $\sigma$ a representation of $H$ in a Hilbert space $\Ki$. 
Let $\theta \in \Aut(G)$.
\par

The representation $( \Ind_H^G \sigma )^\theta$ is equivalent to
the representation $\Ind_{\theta^{-1}(H)}^G( \sigma^\theta)$.
\end{prop}

\begin{proof}
Let $T$ be a right transversal for $H$ in $G$.
Then $S := \theta^{-1}(T)$ is a right transversal for $\theta^{-1}(H)$ in $G$.
\par

For $t \in T$ and $g \in G$,
denote now by $\alpha_T(t,g) \in H$ and $t \cdot_T g \in T$
the elements defined as in ($\ast$) of \ref{constructionInd}(2),
and similarly for $\alpha_S(s,g) \in \theta^{-1}(H)$ and $s \cdot_S g \in S$ when $s \in S$.
\par

Let $g \in G$ and $t \in T$. On the one hand, we have
$$
\theta( \theta^{-1}(t) g) \, = \, t\theta(g) \, = \, \alpha_T(t, \theta(g)) (t \cdot_T \theta(g))
$$
and, on the other hand
$$
\theta( \theta^{-1}(t) g) 
\, = \, \theta\left( \alpha_S(\theta^{-1}(t), g) (\theta^{-1}(t)\cdot_S g)\right)
\, = \, \theta( \alpha_S(\theta^{-1}(t), g)) \theta( \theta^{-1}(t)\cdot_S g).
$$ 
It follows that 
$$
\theta(\alpha_S(\theta^{-1}(t), g)) \, = \, \alpha_T(t, \theta(g))
\hskip.5cm \text{and} \hskip.5cm
\theta(\theta^{-1}(t)\cdot_S g) \, = \, t \cdot_T \theta(g).
\leqno(\ast \ast)
$$
\par

We realize $\Ind_H^G \sigma$ in $\ell^2(T, \Ki)$ 
and $\Ind_{\theta^{-1}(H)}^G(\sigma^\theta)$ in $\ell^2(S, \Ki)$.
Let 
$$
V \, \colon \, \ell^2(T, \Ki) \to \ell^2(S, \Ki)
$$
be the surjective isometry defined by
$$
(Vf) (\theta^{-1}(t)) \, = \, f(t)
\hskip.5cm \text{for} \hskip.2cm
f \in \ell^2(T, \Ki), \, t \in T.
$$
In view of ($\ast \ast$), we have
$$
\begin{aligned}
\Big( V ((\Ind_H^G \sigma)^\theta) (g) f \Big) (\theta^{-1}(t)) 
\, &= \, \Big( V ((\Ind_H^G \sigma) (\theta(g)) f \Big) (\theta^{-1}(t)) 
\\
\, &= \, \Big( (\Ind_H^G \sigma) (\theta(g) ) f \Big) (t) 
\\
\, &= \, \sigma (\alpha_T(t, \theta(g)) f (t \cdot_T \theta(g))
\\
\, &= \, \sigma(\theta(\alpha_S(\theta^{-1}(t),g))) f(\theta(\theta^{-1}(t)\cdot_S g))
\\
\, &= \, \sigma^\theta (\alpha_S(\theta^{-1}(t),g)) (Vf) ( \theta^{-1}(t) \cdot_S g)
\\
\, &= \, \Big( (\Ind_{\theta^{-1}(H)}^G (\sigma^\theta)) (g) Vf \Big) (\theta^{-1}(t)),
\end{aligned} 
$$
for $f \in \ell^2(T, \Ki)$, $t \in T$ and $g \in G$.
This shows that $V$ intertwines
$(\Ind_H^G \sigma)^\theta$ and $\Ind_{\theta^{-1}(H)}^G(\sigma^\theta)$.
\end{proof}

We will need several times the following elementary proposition.
Let $N$ be a normal subgroup of $G$ and $(\sigma, \Ki)$ a representation of $N$.
Let $a \in G$. Recall (see Remark \ref{Rem-AutoRep}) that $\sigma^a$ 
is the conjugate representation of $\sigma$ by $a$,
defined on $\Ki$ by
$$
\sigma^a(n) \, = \, \sigma(ana^{-1})
\hskip.5cm \text{for all} \hskip.2cm
n \in N .
\leqno(\sharp) 
$$

\begin{prop}
% 1.F.8
\label{PropConjIndRep}
Let $G$ be a topological group, 
$N$ an open normal subgroup,
and $\sigma$ a representation of $N$ in a Hilbert space $\Ki$. 
Let $a \in G$.
\begin{enumerate}[label=(\arabic*)]
\item\label{iDEPropConjIndRep}
The representations $\Ind_N^G(\sigma^a)$ and $\Ind_N^G \sigma$ 
are equivalent.
\item\label{iiDEPropConjIndRep}
The restriction of $\Ind_N^G \sigma$ to $N$ is equivalent
to the direct sum $\bigoplus_{t \in T} \sigma^{t^{-1}}$,
where $T \subset G$ is a right transversal for $N$ in $G$.
\par
In particular, if $N$ is contained in the centre of $G$,
the restriction of $\Ind_N^G \sigma$ to $N$ is equivalent to a multiple of $\sigma$.
\end{enumerate}
\end{prop}

\begin{proof}
\ref{iDEPropConjIndRep}
By Proposition~\ref{PropAutondRep}, 
the representations $(\Ind_N^G \sigma)^a$ and $\Ind_N^G(\sigma^a)$ are equivalent.
Since $(\Ind_N^G \sigma)^a$ and $\Ind_N^G(\sigma)$ are equivalent, the claim is proved.

\vskip.2cm

\ref{iiDEPropConjIndRep}
Set $\pi := \Ind_N^G \sigma$ and $\Hi := \ell^2(T, \Ki)$.
By the defining property of induced representations
(Definition~\ref{openeasierthanclosed}),
there exists a $\pi(N)$-invariant closed subspace $\Hi_e$ of $\Hi$
such that $\pi \vert_N$ restricted to $\Hi_e$ is equivalent to $\sigma$
and such that
$$
\Hi \, = \, \bigoplus_{t \in T} \pi(t) (\Hi_e).
$$
For every $n \in N$, $t \in T$, and $\xi \in \Hi_e$, we have 
$$
\pi(n) (\pi(t)\xi) \, = \, \pi(t) (\pi(t^{-1}n t) \xi).
$$
This shows that $\pi(t) (\Hi_e)$ is $\pi(N)$-invariant
and that $\pi \vert_N$ restricted to $\pi(t)\Hi_e$
is equivalent to $\sigma^{t^{-1}}$, for every $t \in T$.
\par

If $N$ is contained in the centre of $G$,
then $\sigma^{g} = \sigma$ for every $g \in G$
and hence $\pi \vert_N$ is equivalent to a multiple of $\sigma$.
\end{proof}

The following proposition has already been referred to, 
in Example \ref{3exeGNS}.
It will be needed again,
in Chapters \ref{ChapterPrimExa} and \ref{ChapterThomadualExamples}.
\par

Given a group $G$, a subgroup $H$, 
and a complex-valued functions $\varphi$ on $H$,
the \textbf{trivial extension} of $\varphi$ to $G$
is the function $\widetilde \varphi$ on $G$ defined by
$$
\widetilde \varphi(g) \, = \,
\begin{cases}
\varphi(g)& \text{if} \hskip.2cm g \in H\\
0& \text{if} \hskip.2cm
g \in G \smallsetminus H.
\end{cases}
$$
Observe that, if $G$ is a topological group, $H$ an open subgroup,
and $\varphi$ a continuous function, then $\widetilde \varphi$ is continuous.
\index{Trivial extension of a function defined on a subgroup}

\begin{prop}
% 1.F.9
\label{diagcoeffinduced}
Let $G$ be a topological group, $H$ an open subgroup,
$\varphi \, \colon H \to \C$ a function of positive type, $\varphi \ne 0$, 
and $\widetilde \varphi \, \colon G \to \C$ its trivial extension.
\par

Then $\widetilde \varphi$ is a function of positive type on $G$.
\par

Denote by $\sigma$ the GNS representation of $H$ associated to $\varphi$.
The GNS representation $\pi$ of $G$ associated to $\widetilde \varphi$
is equivalent to the induced representation $\Ind_H^G \sigma$.
In particular, $\Ind_H^G \sigma$ is cyclic
and has a cyclic vector $\xi$ such that
$$
\widetilde \varphi (g) \, = \,
\langle \big( \Ind_H^G \sigma \big) (g) \xi \mid \xi \rangle
\hskip.5cm \text{for all} \hskip.2cm 
g \in G.
$$
\end{prop}

\noindent
[\emph{Note.} When $G$ is locally compact,
there is a version of this proposition for $H$ closed (rather than $H$ open) in $G$,
for which we refer to \cite{Blat--63}, or \cite[Theorem 6.13]{Foll--16}.
% Lemma 2.31 and Theorem 2.32 in \cite{KaTa--13}.
This involves measures of positive type (rather than functions)
and modulus functions on $G$ and $H$.]

\begin{proof}
Denote by $(\Ki, \sigma, \eta)$ the GNS triple associated to $\varphi$,
so that $\varphi(h) = \langle \sigma(h) \eta \mid \eta \rangle$ for all $h \in H$.
Let $\pi := \Ind_H^G \sigma$ be 
the induced representation of $\sigma$ from $H$ to $G$,
as in the second model of Construction \ref{constructionInd}.
Recall that this involves a right transversal $T$,
so that $G = \bigsqcup_{t \in T} Ht$,
that $\pi(G)$ acts on $\Hi = \ell^2(T, \Ki)$, and that we have
$$
(\pi(g) f) (t) \, = \, \sigma(\alpha(t, g)) f(t \cdot g)
\hskip.5cm \text{for all} \hskip.2cm
g \in G, \hskip.2cm f \in \Hi,
\hskip.2cm \text{and} \hskip.2cm
t \in T ,
$$
with $\alpha(t, g) \in H$ and $t \cdot g \in T$.
We assume moreover, for simplicity, that $e \in T$.
\par

Let $\xi \in \Hi$ be the map defined by $\xi(e) = \eta$
and $\xi(t) = 0$ for all $t \in T \smallsetminus \{e\}$. We have
$$
\begin{aligned}
\langle \pi(h) \xi \mid \xi \rangle 
\, &\, = \, 
\sum_{t \in T} \langle \sigma(\alpha(t, h)) (\xi(t \cdot h)) \mid \xi(t) \rangle \, = \,
\langle \sigma(\alpha(e, h)) (\xi(e \cdot h)) \mid \xi(e) \rangle 
\\
\, &\, = \, 
\langle \sigma(h) \eta \mid \eta \rangle \, = \,
\varphi(h) 
\hskip2cm \text{for all} \hskip.2cm 
h \in H,
\\
\langle \pi(g) \xi \mid \xi \rangle
\, &\, = \,
0 
\hskip5.2cm \text{for all} \hskip.2cm 
g \in G \smallsetminus H .
\end{aligned}
$$
Therefore $\widetilde \varphi (g) = \langle \pi(g)\xi \mid \xi \rangle$ for all $g \in G$,
i.e., $\widetilde \varphi$ is the function of positive type
associated to $\pi$ and $\xi$;
in particular, $\widetilde \varphi$ is of positive type.
\par

To show that the GNS representation associated to $\widetilde \varphi$
is equivalent to $\pi$,
it suffices by Proposition \ref{GNSbijP(G)cyclic}
to check that $\xi$ is a cyclic vector for $\pi$.
\par

For every $t \in T$ and $\kappa \in \Ki$, let $\xi_t^\kappa \in \Hi$
be defined by $\xi_t^\kappa(t) = \kappa$
and $\xi_t^\kappa(s) = 0$ when $s \in T \smallsetminus \{t\}$;
observe that $\xi_e^\eta = \xi$.
The linear span $X$
of $\{\xi_t^\kappa \in \Hi \mid t \in T, \kappa \in \Ki \}$
is dense in $\Hi$. 
On the one hand, since $\eta$ is cyclic for $\sigma$,
the linear span of the set 
$$
\{ \xi_t^{\sigma(h)\eta} \in \Hi \mid t \in T, \hskip.1cm h \in H \}
$$
is dense in $X$, and hence in $\Hi$.
On the other hand,
for $t \in T$ and $h \in H$, we have
$$
\pi(t^{-1}h) \xi_e^\eta \, = \, \xi_t^{\sigma(h)\eta} ;
$$
indeed, $\alpha(t, t^{-1}h) = h$ and, for $s \in T$,
we have $s \cdot t^{-1}h = e$ if and only if $s = t$.
This shows that $\xi_e^\eta = \xi$ is cyclic for $\pi$.
\end{proof}

\index{Monomial representation}
\index{Representation! monomial}
From now on, we restrict ourselves to \textbf{monomial representations}, 
that is, to representations of the form $\Ind_H^G \chi$, 
where $\chi \, \colon H \to \T$ is a unitary character of $H$.
We use Model (2) in Construction \ref{constructionInd},
so that he corresponding induced representation $\pi = \Ind_{H}^G \chi$
is realized on $\ell^2(T)$ and we have
$$
(\pi(g) f) (t) \, = \, \chi (\alpha(t, g)) f(t \cdot g)
\hskip.5cm \text{for all} \hskip.2cm
g \in G, \hskip.1cm f \in \ell^2(T), 
\hskip.2cm \text{and} \hskip.2cm 
t \in T .
$$
We assume moreover that $e \in T$,
and we denote by $(\delta_t)_{t \in T}$ 
the orthonormal basis of $\ell^2(T)$
defined by $\delta_t(t) = 1$ and $\delta_t(s) = 0$
when $s \in T$, $s \ne t$.
Observe that we have
$$
\pi(t^{-1}) \delta_e \, = \, \delta_t
\hskip.5cm \text{for all} \hskip.2cm
t \in T .
$$
In particular $\delta_e$ is a cyclic vector for~$\pi$. 
\par

For $i = 1, 2$, let $H_i$ be an open subgroup of $G$ 
and $\chi_i$ a unitary character of $H_i$. 
Choose a right transversal $T_i$ for $H_i$ in $G$, with $e \in T_i$. 
For $g \in G$ and $t \in T_i$, 
we have as in $(\ast)$ of \ref{constructionInd}
uniquely defined elements $\alpha_i(t,g) \in H_i$ and $t \cdot_i g \in T_i$ 
such that $tg = \alpha_i(t,g) (t \cdot_i g)$.
\par

Our study of the monomial representations $\pi_i = \Ind_{H_i}^G \chi_i$,
their irreducibility and their equivalence, is based on the following lemma. 

\begin{lem}
% 1.6.10
\label{Lem-IrredInducedRep}
We keep the notation above.
Let $S \in \Hom_G(\pi_1, \pi_2)$. Set $f := S \delta_e \in \ell^2(T_2)$.
\begin{enumerate}[label=(\arabic*)]
\item\label{iDELem-IrredInducedRep}
We have $f = 0$ if and only if $S = 0$.
\item\label{iiDELem-IrredInducedRep}
Assume that $\pi_1 = \pi_2$; 
then $f$ is a scalar multiple of $\delta_e$
if and only if $S$ is a scalar operator.
\item\label{iiiDELem-IrredInducedRep}
For every $t \in T_2$ and $h \in H_1$, we have 
$$
\chi_2(\alpha_2(t, h)) f(t \cdot_2 h) \, = \, \chi_1(h) f(t).
$$
\item\label{ivDELem-IrredInducedRep}
Let $t \in T_2$ be such that
the $H_1$-orbit $\{t\cdot_2 h \mid h \in H_1 \}$ in $T_2$ is infinite; 
then $f(t) = 0$.
\item\label{vDELem-IrredInducedRep}
Let $t \in T_2$ and $h \in H_1$ be such that $f(t) \ne 0$ and $t\cdot_2h = t$;
then $tht^{-1} \in H_2$ and $\chi_1(h) = \chi_2(tht^{-1})$.
\end{enumerate}
\end{lem}

\begin{proof}
Since $S \in \Hom_G(\pi_1, \pi_2)$,
Items \ref{iDELem-IrredInducedRep} and \ref{iiDELem-IrredInducedRep}
follow from the fact, mentioned above, 
that $\delta_e$ is a cyclic vector for $\pi_1$.

\vskip.2cm

\ref{iiiDELem-IrredInducedRep} 
Let $h \in H_1$.
Since $\alpha_1(e, h) = h$ and $e \cdot_1 h = e$, 
and $s \cdot_1 h \ne e$ for $s \in T_1 \smallsetminus \{e\}$, we have
$$
\begin{aligned}
(\pi_1(h)\delta_e) (e)
\, &= \, \chi_1(\alpha_1e, h)) \delta_e(e \cdot_1 h) = \chi_1(h)
\\
(\pi_1(h)\delta_e) (s)
\, &= \, \chi_1(\alpha_1(s, h)) \delta_e( s \cdot_1 h) \, = \, 0
\hskip.5cm \text{if} \hskip.2cm
s \ne e ,
\end{aligned}
$$
so that $\pi_1(h) \delta_e = \chi_1(h) \delta_e$.
\par

Let $t \in T_2$. We have
$$
(\pi_2(h) f) (t)
\, = \, (\pi_2(h) S \delta_e) (t)
\, = \, (S \pi_1(h) \delta_e) (t)
\, = \, \chi_1(h) (S \delta_e) (t)
\, = \, \chi_1(h) f(t)
$$
and
$$
(\pi_2(h) f) (t) \, = \, \chi_2(\alpha_2(t, h)) f(t\cdot_2 h).
$$
The claim follows.

\vskip.2cm

\ref{ivDELem-IrredInducedRep} 
For $t \in T_2$ and $h \in H_1$, 
it follows from \ref{iiiDELem-IrredInducedRep} that
$\vert f (t\cdot_2 h) \vert = \vert f(t) \vert$,
i.e., that $\vert f \vert$ is constant on the $H_1$-orbit of $t$ in $T_2$. 
Since $f \in \ell^2(T_2)$, this implies that $f(t) = 0$ if this orbit is infinite.

\vskip.2cm

\ref{vDELem-IrredInducedRep}
Let $t \in T_2$ and $h \in H_1$ be such that $t\cdot_2h = t$ and $f(t) \ne 0$.
By \ref{iiiDELem-IrredInducedRep}, we have
$$
\chi_2(\alpha_2(t, h)) \, = \, \chi_1(h).
$$
Moreover, since 
$$
th \, = \, \alpha_2(t, h) (t \cdot_2 h) \, = \, \alpha_2(t, h) t ,
$$
we have $tht^{-1} = \alpha_2(t, h) \in H_2$
and $\chi_1(h) = \chi_2(\alpha_2(t, h)) = \chi_2(tht^{-1})$.
\end{proof}

If $\chi$ is a unitary character of the subgroup $H$ and $g \in G$, 
we denote by $\chi^g$ the unitary character of $g^{-1} H g$ defined by 
$$
\chi^g(h) \, = \, \chi(ghg^{-1}) 
\hskip.5cm \text{for all} \hskip.2cm
h \in g^{-1} H g
$$
(this extends a notation of Remark \ref{Rem-AutoRep}
to a more general situation).
\par

\index{Commensurator}
Recall that the \textbf{commensurator} of a subgroup $H$ of $G$
is the set $\Comm_G (H)$,of those elements $g \in G$ such that 
$H \cap g^{-1} Hg$ is of finite index in both $H$ and $g^{-1} Hg$. 
Equivalent definitions are recorded in Proposition \ref{defScommensurator}.
\par

We are now ready to formulate the irreducibility and equivalence criteria 
for monomial representations.

\begin{theorem}[\textbf{Irreducibility of monomial representations}]
% 1.F.11
\label{Theo-IrredInducedRep}
Let $G$ be a topological group, $H$ an open subgroup of $G$, 
and $\chi$ a unitary character of $H$.
Assume that
\begin{enumerate}
\item[$\bullet$]
for every $g \in \Comm_G (H) \smallsetminus H$,
the restrictions of $\chi$ and $\chi^g$
to the subgroup $H \cap g^{-1} H g$ do not coincide. 
\end{enumerate} 
Then the induced representation $\Ind_H^G \chi$ of $G$ is irreducible.
\end{theorem}

\begin{proof}
Let $T$ be a right transversal for $H$ in $G$, with $e \in T$. 
\par

We will show the contraposition, 
and therefore assume that $\pi = \Ind_H^G \chi$ is not irreducible. 
We claim that there exists $t \in T \smallsetminus \{e\}$ 
such that $t \in \Comm_G (H)$ and $\chi = \chi^t$ on $H \cap t^{-1} H t$.
This will prove the theorem because $t \notin H$.
\par

Since $\pi$ is not irreducible, there exists $S \in \Hom_G(\pi, \pi) = \pi(G)'$ 
which is not a scalar operator.
On the one hand, the support of $f := S \delta_e \in \ell^2(T)$ is not reduced to $\{e\}$, 
by Lemma~\ref{Lem-IrredInducedRep}~\ref{iiDELem-IrredInducedRep}.
On the other hand, by Item \ref{ivDELem-IrredInducedRep} of the same lemma, 
we have $f(s) = 0$ for every $s \in T$ with an infinite $H$-orbit. 
Therefore, there exists $t \in T \smallsetminus \{e\}$ with a finite
$H$-orbit such that $f(t) \ne 0$. 
\par

Let $t^* \in T$ be such that $Ht^{-1} = Ht^*$, that is, $t^* = e \cdot t^{-1}$.
We have
$$
\begin{aligned}
(\pi(t)\delta_e) (t^*)
\, &= \, \chi(\alpha(t^*,t)) \delta_e(t^* \cdot t) = \chi(\alpha(t^*,t)) ,
\\
(\pi(t)\delta_e) (s)
\, &= \, \chi(\alpha(s,t)) \delta_t( s \cdot t) \, = \, 0
\hskip.5cm \text{if} \hskip.2cm
s \ne t^* ,
\end{aligned}
$$
and therefore $\pi(t) \delta_e = \chi(\alpha(t^*, t)) \delta_{t^*}$.
\par

Consider the adjoint operator $S^*$ of $S$, which is also in $\pi(G)'$, 
and set $f' := S^*\delta_e$.
We claim that $f'(t^*) \ne 0$.
\par

Indeed, since $\pi (t^{-1})\delta_e = \delta_t$ 
and $\pi (t)\delta_e = \chi(\alpha(t^*,t)) \delta_{t^*}$, we have
$$
\begin{aligned}
f'(t^*) 
\, &= \, \langle f' \mid \delta_{t^*} \rangle
\, = \, \langle S^*\delta_e \mid \delta_{t^*} \rangle
\, = \, \langle \delta_e \mid S \delta_{t^*} \rangle
\\
\, &= \, \langle \delta_e \mid S \chi(\alpha(t^*, t))^{-1} \pi(t) \delta_e \rangle
\, = \, \chi(\alpha(t^*,t))\langle \delta_e \mid S \pi(t) \delta_e \rangle
\\
\, &= \, \chi(\alpha(t^*,t))\langle \delta_e \mid \pi(t) S \delta_e \rangle
\, = \, \chi(\alpha(t^*,t))\langle \pi(t^{-1}) \delta_e \mid f \rangle
\\
\, &= \, \chi(\alpha(t^*,t)) \overline{\langle f \mid \delta_t \rangle}
\, = \, \chi(\alpha(t^*,t))\overline{f(t)},
\end{aligned}
$$
and in particular $f'(t^*) \ne 0$, as claimed.
\par

Lemma \ref{Lem-IrredInducedRep}~\ref{ivDELem-IrredInducedRep}
applied to $S^*$ shows that $t^*$ has a finite $H$-orbit. 
Therefore $t \in \Comm_G (H)$,
by Proposition \ref{defScommensurator}~\ref{6DEdefScommensurator}.
\par
 
Let $h \in H \cap t^{-1}Ht$. As $f(t) \ne 0$ and $t \cdot h = t$, 
Lemma~\ref{Lem-IrredInducedRep}~\ref{vDELem-IrredInducedRep}, 
applied with $H_1 = H_2 = H$ and $\chi_1 = \chi_2 = \chi$, 
shows that $\chi(h) = \chi(tht^{-1}) = \chi^t(h)$. 
So, $\chi$ and $\chi^t$ coincide on $H \cap t^{-1}Ht$.
\end{proof}

\begin{rem}
% 1.F.12
\label{Rem-IrredInducedRepConverse}
The sufficient condition for the irreducibility of 
$\Ind_H^G \chi$ stated in Theorem \ref{Theo-IrredInducedRep} 
is also necessary; see Theorem 6' in \cite{Mack--51}.
Corollary \ref{Cor-NormalSubgIrr} below provides
a proof of this fact in the special case where $H$ is a normal subgroup.
\end{rem}

The following two corollaries particularize Theorem~\ref{Theo-IrredInducedRep}
to two opposite cases, self-commensurating subgroups and normal subgroups.

\begin{cor}
% 1.F.13
\label{Cor-QuasiReg}
Let $H$ be an open subgroup of $G$ such that $\Comm_G (H) = H$.
\par

Then $\Ind_H^G \chi$ is irreducible for every unitary character $\chi$ of $H$.
In particular, the quasi-regular representation $\lambda_{G/H}$ is irreducible.
\index{Quasi-regular representation}
\index{Representation! quasi-regular}
\end{cor}

\begin{exe}
% 1.F.14
\label{Exa-IrredRepFreeGroup}
(1)
Let $F$ be the free group on two generators $a$ and $b$.
Let $A$ be the subgroup generated by $a$; it is an infinite cyclic group.
One easily checks that $\Comm_F (A) = A$.
\par

Therefore $\Ind_A^F \chi$ is irreducible for every unitary character $\chi$ of $A$.

\vskip.2cm

(2)
Let $\Gamma$ be a group of permutations of an infinite set $X$
and let $\pi$ denote the corresponding permutation representation of $\Gamma$ on $\ell^2(X)$,
defined by $(\pi(\gamma) f)(x) = f(\gamma^{-1} x)$ for all
$\gamma \in \Gamma$, $f \in \ell^2(X)$, and $x \in X$.
If the action is $2$-transitive, then the representation $\pi$ is irreducible.
\index{Action! $2$-transitive}
\par

Indeed, if $x_0$ denote some point of $X$ and $H$ its isotropy subgroup,
the hypothesis of $2$-transitivity implies that
$H$ is its own commensurator in $\Gamma$ (see Proposition A.5.4),
hence $\pi = \lambda_{\Gamma/H} = \Ind_H^\Gamma 1_H$ is irreducible by Corollary 1.F.12.
\par

More generally, $\Ind_H^\Gamma \rho$ is irreducible
for any irreducible finite dimensional representation $\rho$ of $H$.
This is a particular case of a result in \cite{Mack--51},
and appears also in several other references, for example in \cite{Obat--89}.
\par

In case $\Gamma$ acts $2$-transitively on a finite set $X$,
the corresponding permutation representation is the direct sum of two representations,
the unit representation $1_\Gamma$
and another irreducible representation of dimension $\vert X \vert - 1$.
This goes back to Burnside, see \cite[\S~250]{Burn--11}.
\end{exe}

\begin{cor}
% 1.F.15
\label{Cor-NormalSubgIrr}
Let $N$ be an open normal subgroup of $G$
and $\chi$ a unitary character of $N$. The following conditions are equivalent:
\begin{enumerate}[label=(\roman*)]
\item\label{iDECor-NormalSubgIrr}
$\chi^g \ne \chi$ for every $g \in G \smallsetminus N$;
\item\label{iiDECor-NormalSubgIrr}
the induced representation $\Ind_N^G \chi$ is irreducible.
\end{enumerate}
\end{cor}

\begin{proof}
The fact that \ref{iDECor-NormalSubgIrr} implies \ref{iiDECor-NormalSubgIrr} 
being a particular case of Theorem~\ref{Theo-IrredInducedRep},
we only have to show that \ref{iiDECor-NormalSubgIrr} implies \ref{iDECor-NormalSubgIrr}.
\par

Assume, by contraposition, that we have $\chi^{g_0} = \chi$
for some $g_0 \in G \smallsetminus N$.
It will be convenient to use here Model (1)
in Construction \ref{constructionInd} for $\pi = \Ind_N^G \chi$; 
so, $\pi$ is realized in the space $\Hi$ of the functions $f \, \colon G \to \C$
such that $f(ng) = \chi(n) f(g)$ for all $g \in G$, $n \in N$
and such that $\vert f \vert \in \ell^2(N \backslash G)$.
For $g \in G$, we have $\pi(g)f(x) = f(xg)$ for all $x \in G$.
\par
 
For $f \in \Hi$, define $Tf \, \colon G \to \C$ by $Tf(g) = f(g_0g)$.
Then $Tf \in \Hi$ and $T \, \colon \Hi \to \Hi$ is a unitary operator.
Indeed, for all $g \in G$, $n \in N$, we have 
$$
\begin{aligned}
Tf(ng) 
\, &= \, f(g_0 n g_0^{-1} g_0 g)
\, = \, \chi^{g_0}(n) f(g_0 g)
\\
\, &= \, \chi(n) Tf(g)
\end{aligned}
$$
and 
$$
\begin{aligned}
\Vert Tf \Vert^2
\, &= \, \sum_{Nx \in N \backslash G} \vert Tf(x) \vert^2
\, = \, \sum_{Nx \in N \backslash G} \vert f(g_0x) \vert^2
\\
\, &= \, \sum_{Nx \in N \backslash G} \vert f(x) \vert^2
\, = \, \Vert f \Vert^2.
\end{aligned}
$$
Moreover, $T$ intertwines $\pi$ with itself. Indeed, for all $g \in G$, we have
$$
\begin{aligned}
\pi(g) (Tf) (x)
\, &= \, Tf (x g)
\, = \, f(g_0 xg)
\\
\, &= \, (\pi(g)f) (g_0x)
\, = \, T (\pi(g)f) (x)
\end{aligned}
$$
for all $x \in G$. Since $T$ is obviously not a scalar operator,
it follows that $\pi$ is not irreducible.
\end{proof}

Corollary~\ref{Cor-NormalSubgIrr} 
will be applied in Chapter \ref{Chapter-ExamplesIndIrrRep}
to produce irreducible representations for various groups.
\par

We turn now to the question of equivalence of monomial representations.

\begin{theorem}[\textbf{Non-equivalence of monomial representations}]
% 1.F.16
\label{Theo-EquiInducedRep}
Let $G$ be a topological group, $H_1, H_2$ open subgroups of $G$, 
and $\chi_1, \chi_2$ unitary characters of $H_1, H_2$ respectively. 
Assume that
\begin{enumerate}
\item[$\bullet$]
for every $g \in G$ such that $g^{-1} H_2 g \cap H_1$ 
has finite index in $g^{-1} H_2 g$ and in $H_1$,
the restrictions of $\chi_2^g$ and $\chi_1$ 
to the subgroup $g^{-1} H_2 g \cap H_1$ do not coincide.
\end{enumerate}
Then $\Ind_{H_1}^G \chi_1$ and $\Ind_{H_2}^G \chi_2$ are not equivalent.
\end{theorem}

\begin{proof}
We use Model (2) of Construction \ref{constructionInd}.
For $i = 1, 2$, let $T_i$ be a right transversal for $H_i$ in $G$,
with $e \in T_i$.
\par

We will again show the contraposition, and therefore assume that 
$\pi_1 := \Ind_{H_1}^G \chi_1$ and $\pi_2 := \Ind_{H_2}^G \chi_2$ are equivalent. 
We claim that there exists $t \in T_2$ such that
$t^{-1} H_2 t \cap H_1$ has finite index in $t^{-1} H_2 t$ and in $H_1$,
and such that $\chi_2^t$ and $\chi_1$ coincide on $t^{-1} H_2 t \cap H_1$. 
This will prove the theorem.
\par

Since $\pi_1$ and $\pi_2$ are equivalent, 
there exists $S \ne 0$ in $\Hom_G(\pi_1, \pi_2)$.
Set $f := S \delta_e \in \ell^2(T_2)$.
\par

On the one hand, we have $f \ne 0$, 
by Lemma \ref{Lem-IrredInducedRep}~\ref{iDELem-IrredInducedRep}.
On the other hand, by Item \ref{ivDELem-IrredInducedRep} of the same lemma,
we have $f(t) = 0$ for every $t \in T_2$ with an infinite $H_1$-orbit. 
Therefore there exists $t \in T_2$
such that $f(t) \ne 0$, with a finite $H_1$-orbit.
In particular, the stabilizer $t^{-1} H_2 t \cap H_1$ of $t$ in $H_1$ has finite index in $H_1$.
\par

Set $t^* = e \cdot_1 t^{-1} \in T_1$.
Consider the adjoint operator $S^*$, which is in $\Hom_G(\pi_2, \pi_1)$, 
and set $f' := S^*\delta_e \in \ell^2(T_1)$.
As in the proof of Theorem~\ref{Theo-IrredInducedRep}, we have
$$
f'(t^*) \, = \, \chi_1(\alpha_1(t^*,t))\overline{f(t)}
$$
and hence $f'(t^*) \ne 0$.
\par

Lemma~\ref{Lem-IrredInducedRep}~\ref{ivDELem-IrredInducedRep} 
applied now to $S^*$ 
shows that $t^*$ has a finite $H_2$-orbit in $T_1$. 
Since $t^* = e \cdot_1 t^{-1}$, the stabilizer of $t^*$ in $H_2$
is $t H_1 t^{-1} \cap H_2$. 
Therefore $t H_1 t^{-1} \cap H_2$ has finite index in $H_2$
and so $t^{-1} H_2 t \cap H_1$ has finite index in $t^{-1} H_2 t$.
\par
 
Let $h \in t^{-1} H_2 t \cap H_1$. 
As $f(t) \ne 0$ and $t\cdot_2 h = t$ ,
Lemma~\ref{Lem-IrredInducedRep}~\ref{vDELem-IrredInducedRep}, 
shows that $\chi_1(h) = \chi_2(tht^{-1})$; 
hence, $\chi_2^t$ and $\chi_1$ coincide on $t^{-1}H_2t \cap H_1$.
\end{proof}

\begin{rem}
% 1.F.17
\label{Rem-Mackey--Shoda}

Theorems \ref{Theo-IrredInducedRep} and \ref{Theo-EquiInducedRep}
are known as the Mackey--Shoda criteria 
for irreducibility and equivalence for monomial representations.
The terminology refers to \cite{Shod--33}, for the case finite groups, 
and \cite{Mack--51}, for the general case.
\par

These criteria are also necessary conditions
for irreducibility and equivalence, respectively
(compare Corollaries~\ref{Cor-NormalSubgIrr}, \ref{Cor-SelfCommSubgEquiv},
and \ref{Cor-NormalSubgEquiv}).
Moreover, these criteria remain valid 
when the inducing representations
are finite-dimensional \cite{Klep--61, Corw--75}.
However, the irreducibility criterion 
may fail when the inducing representation
is infinite-dimensional \cite{BeCu--03}.
\end{rem}

We explicit now two consequences of Theorem~\ref{Theo-EquiInducedRep},
first for self-commensurating subgroups and then for normal subgroups.
\par

\index{Commensurate subgroups}
Recall that two subgroups $H$ and $L$ of a group $G$ are \textbf{commensurate}
if $H \cap L$ has finite index in both $H$ and $L$.
Two commensurate subgroups $H, L$ of $G$ have the same commensurator:
$\Comm_G(H) = \Comm_G(L)$;
see Proposition \ref{defScommensurator}.

\begin{cor}[\textbf{Equivalence of monomial irreducible representations}]
% 1.F.18
\label{Cor-SelfCommSubgEquiv}
Let $G$ be a topological group, $H_1, H_2$ two open subgroups,
$\chi_1$ a unitary character of $H_1$, and $\chi_2$ a unitary character of $H_2$.
Suppose that $\Comm_G(H_1) = H_1$ and $\Comm_G(H_2) = H_2$. 
The following properties are equivalent:
\begin{enumerate}[label=(\roman*)]
\item\label{iDECor-SelfCommSubgEquiv}
the irreducible representations $\Ind_{H_1}^G \chi_1$ and $\Ind_{H_2}^G \chi_2$
are equivalent;
\item\label{iiDECor-SelfCommSubgEquiv}
there exists $g \in G$ such that $g^{-1} H_2 g = H_1$ and $\chi_2^g = \chi_1$.
\end{enumerate}
\end{cor}

\begin{proof}
Assume that $g^{-1} H_2 g = H_1$ and $\chi_2^g = \chi_1$ for some $g \in G$.
Then $\Ind_{H_1}^G \chi_1$ is equivalent to $(\Ind_{H_2}^G \chi_2)^g$
(see Proposition~\ref{PropAutondRep}) and hence to $\Ind_{H_2}^G \chi_2$.
\par

Conversely,
assume that $\Ind_{H_1}^G \chi_1$ and $\Ind_{H_2}^G \chi_2$ are equivalent.
Then, by Theorem \ref{Theo-EquiInducedRep},
there exists $g \in G$ such that $g^{-1} H_2 g \cap H_1$ 
has finite index in $g^{-1} H_2 g$ and in $H_1$,
and such that $\chi_2^g$ and $\chi_1$ coincide on $g^{-1} H_2 g \cap H_1$.
\par

In particular, $H_1$ and $g^{-1} H_2 g $ are commensurate.
Since $\Comm_G(H_1) = H_1$ and 
$$
\Comm_G(g^{-1} H_2 g) \, = \, g^{-1} \Comm_G(H_2) g \, = \, g^{-1} H_2 g,
$$
it follows that $g^{-1} H_2 g = H_1$ and $\chi_2^g = \chi_1$.
\end{proof}

Note that there exist triples $(G, H_1, H_2)$,
where $G$ is a group and $H_1, H_2$ subgroups which are not conjugate in $G$,
such that the representations $\Ind_{H_1}^G 1_{H_1}$ and $\Ind_{H_2}^G 1_{H_2}$
are equivalent (and not irreducible).
For examples with $G$ finite, we quote \cite{Suna--86} and \cite{Broo--88}.

\vskip.2cm

The following corollary is a direct consequence of
Theorem \ref{Theo-EquiInducedRep} and Proposition~\ref{PropConjIndRep}.

\begin{cor}
% 1.F.19
\label{Cor-NormalSubgEquiv}
Let $N$ be an open normal subgroup of $G$ 
and $\chi_1, \chi_2$ unitary characters of $N$.
\par

Then $\Ind_N^G \chi_1$ and $\Ind_N^G \chi_2$ 
are non-equivalent representations of $G$
if and only if $\chi_1^g \ne \chi_2$ for every $g \in G$.
\end{cor}

We will later need (Proof of Theorem~\ref{Theo-PrimIdealTwoStepNilpotent})
the following weak containment result involving induced representations.
\par

\index{Amenable! homogeneous space}
Given an open subgroup $H$ of a topological group $G$,
we say that $G/H$ is an \textbf{amenable homogeneous space} in the sense of Eymard,
if the unit representation~$1_G$ is weakly contained
in the quasi-regular representation $(\lambda_{G/H}, \ell^2(G/H))$;
for more on this notion, see \cite{Eyma--72}.
Observe that, when $H$ is normal in $G$,
the amenability of the homogeneous space $G/H$ is equivalent
to the amenability of the (discrete)
group $G/H$ \cite[Theorem G.3.2]{BeHV--08}.

\begin{prop}
% 1.F.20
\label{Prop-WeakContainmentInduced}
Let $G$ be a topological group, $H$ an open subgroup of $G$
and $(\pi, \Hi)$ a representation of $G$.
\begin{enumerate}[label=(\arabic*)]
\item\label{iDEProp-WeakContainmentInduced}
Assume that $G/H$ is amenable.
Then $\pi$ is weakly contained in the induced representation $\Ind_H^G (\pi \vert_H)$.
\item\label{iiDEProp-WeakContainmentInduced}
Assume that $H$ is normal in $G$ and that $G/H$ is abelian. 
Then $\Ind_H^G (\pi \vert_H)$ is weakly equivalent to the direct sum 
$$
\bigoplus_{\chi \in \widehat{G/H}} \hskip.1cm \pi \chi
$$
where $\pi \chi \simeq \pi \otimes \chi$ is the representation of $G$ on $\Hi$
given by $(\pi \chi) (g) = \pi(g) \chi(gH)$ for all $g \in G$.
\end{enumerate}
\end{prop}

\begin{proof} 
Recall that $\Ind_H^G (\pi \vert_H)$ is equivalent
to the tensor product $\pi \otimes \lambda_{G/H}$
(Proposition~\ref{InductionQqPropr}~\ref{ivDEInductionQqPropr}).

\vskip.2cm

\ref{iDEProp-WeakContainmentInduced}
Since $1_G \prec \lambda_{G/H}$ by hypothesis,
$$
\pi \, = \, \pi \otimes 1_G
\, \prec \, \pi \otimes \lambda_{G/H}
\, \simeq \, \Ind_H^G (\pi \vert_H)
$$
by continuity of the tensor product operation \cite[Proposition F.3.2]{BeHV--08}.

\vskip.2cm

\ref{iiDEProp-WeakContainmentInduced}
Since $G/H$ is an abelian group, 
$\lambda_{G/H}$ is weakly equivalent to $\bigoplus_{\chi \in \widehat{G/H}} \chi$,
see \cite[Proposition F.2.7]{BeHV--08}. 
It follows that
$\Ind_H^G(\pi \vert_H) \, \simeq \, \pi \otimes \lambda_{G/H}$
is weakly equivalent to
$\pi \otimes \Big( \bigoplus_{\chi \in \widehat{G/H}} \chi \Big) \, \simeq \,
\bigoplus_{\chi \in \widehat{G/H}} \pi \chi$,
again by continuity of the tensor product.
\end{proof}

\section
[Irreducible decompositions]
{On decomposing representations
% \\
into irreducible representations}
% Section 1.G
\label{SectionDecomposingIrreps}

For a second-countable locally compact group $G$,
arbitrary representations can be decomposed into irreducible representations.
This appears in several situations:
\begin{enumerate}
\item[--]
Decomposing representations into irreducible components
is of fundamental importance already in the classical Fourier analysis, 
when $G$ is a torus $\T^n$ or the translation group $\R^n$ of an Euclidean space.
\item[--]
The theory of unitary representations of groups,
and in particular decompositions into irreducible components,
provide an organizing principle for the study of special functions,
which appear as solutions
of many partial differential equations of mathematical physics \cite{Vile--68}.
\item[--]
In case $G$ is of type I, its regular representation
can be decomposed into irreducible representations
using the so-called Plancherel measure on $\widehat G$ \cite[18.8]{Dixm--C*}.
\item[--]
Decomposition into irreducible representations relates to spectral analysis.
For example, if $\Gamma$ is a discrete subgroup of $\SL_2(\R)$
such that $\Gamma \backslash \SL_2(\R)$ is compact,
understanding the irreducible subrepresentations
of the natural representation of $\SL_2(\R)$ in $L^2(\Gamma \backslash \SL_2(\R)$ 
is essentially equivalent to understanding the spectrum
of the Laplace--Beltrami operator on the upper half-plane;
see \cite[Chap.~I, \S~4]{GGPS--69} and \cite[in particular Theorem 2.7.1]{Bump--97}.
\end{enumerate}
\par

For compact groups, decompositions in direct sums are appropriate:
every representation of such a group is a direct sum of irreducible representations.
But, for non-compact groups, including $\Z$ and $\R$ (see Example~\ref{Exa-DirIntRegRep}),
it is in general necessary to use the more general notion
of direct integral decomposition.
\par

The theory of direct integral decompositions
was founded by von Neumann in the late 1930s \cite{vNeu--49},
and made explicit for group representations in the early 1950s 
\cite[Theorem 1.2]{Maut--50a}, \cite[Chapitre IV]{Gode--51a}.
We give a brief account of this theory 
and refer to \cite[Chap.~II]{Dixm--vN} and \cite[\S~8]{Dixm--C*}
for a comprehensive exposition.
% \cite[\S~7]{Foll--16} for a short introduction).

\subsection*{Direct integrals of representations}

\index{Field of Hilbert spaces over a measure space}
Let $(X, \mathcal B, \mu)$ be a measure space,
where $\mu$ is a $\sigma$-finite measure on $(X, \mathcal B)$.
A \textbf{field of Hilbert spaces over} $X$ is a family $(\Hi_x)_{x \in X}$,
where $\Hi_x$ is a Hilbert space for each $x \in X$. 
Elements in $\prod_{x \in X} \Hi_x $ are called \textbf{vector fields over} $X$.
\par

\index{Fundamental family of measurable vector fields over a Borel space}
A \textbf{fundamental family of measurable vector fields} for a field $(\Hi_x)_{x \in X}$ over $X$
is a sequence $(e^n)_{n \in \N}$ of vector fields over $X$,
where $e^n = (e^n_x)_{x \in X}$ and $e^n_x \in \Hi_x$ for all $x \in X$,
with the following properties:
\begin{enumerate}
\item[$\bullet$]
the map $x \mapsto \langle e^m_x \mid e^n_x \rangle$ is measurable,
for every $m, n \in \N$,
\item[$\bullet$]
the linear span of $\{e^n_x \hskip.2cm \colon \hskip.2cm n \in \N\}$
is dense in $\Hi_x$, for every $x \in X$.
\end{enumerate}
Note that, for such a fundamental family to exist, 
the second property implies that the Hilbert space $\Hi_x$ is separable
for every $x \in X$.
\par

Fix a fundamental family of measurable vector fields $(e^n)_{n \in \N}$.
A vector field $\xi \in \prod_{x \in X} \Hi_x$ is said to be
a \textbf{measurable vector field} if the functions
$x \mapsto \langle \xi_x \mid e^n_x \rangle$
are measurable for all $n \in \N$.
The set $M$ of measurable vector fields is a linear subspace of $\prod_{x \in X} \Hi_x$. 
Moreover, for $\xi, \eta \in M$, the function
$X \to \C$, $x \mapsto \langle \xi_x \mid \eta_x \rangle$
is measurable.
The pair $\left( (\Hi_x)_{x \in X}, M \right)$ is called
a \textbf{measurable field of Hilbert spaces} over $X$.
\index{Measurable field! $1$@measurable vector field over a Borel space}
\index{Measurable field! $2$@measurable field of Hilbert spaces over a Borel space} 
\par

\index{Square-integrable vector field over a Borel space}
Two measurable vector fields are equivalent if they are equal $\mu$-almost everywhere.
From now on, we write (abusively) ``measurable vector fields''
instead of ``equivalence classes of measurable vector fields'' for this equivalence relation.
A measurable vector field $\xi$ is a \textbf{square-integrable vector field} if 
$$
\int_X \Vert \xi_x \Vert^2 d\mu(x) \, < \, \infty .
$$
Equipped with the obvious inner product,
the set $\Hi$ of all square-integrable vector fields is a Hilbert space called
the \textbf{direct integral} of the field $(\Hi_x)_{x \in X}$ of Hilbert spaces over $X$
and denoted by
$$
\Hi \, = \, \int_X^\oplus \Hi_x d\mu(x).
$$
\index{Direct integral of Hilbert spaces }

\begin{exe}
% 1.G.1
\label{Exa-DirectIntHilbertSpace}
(1)
Let $X$ be a countable set and $\mu$ the counting measure.
Let $(\Hi_x)_{x \in X}$ be a field of separable Hilbert spaces over $X$.
An appropriate fundamental family of measurable vector fields in this case
is a family $(e^n)_{n \in \N}$ such that
$(e^n_x)_{n \in \N}$ contains an orthonormal basis of $\Hi_x$ for all $x \in X$.
Then every vector field is measurable and 
$$
\int^\oplus_X \Hi_x d\mu(x) \, = \, \bigoplus_{x \in X} \Hi_x.
$$

\vskip.2cm

(2)
Let $\Ki$ be a fixed separable Hilbert space. 
Set $\Hi_x = \Ki$ for every $x \in X$. 
Choose an orthonormal basis $(\xi_n)_{n \ge 1}$ of $\Ki$ and define
a fundamental family of measurable vector fields $(e^n)_{n \ge 1}$ over $X$
by $e^n_x = \xi_n$ for every $x \in X$ and ${n \ge 1}$. 
The measurable vector fields over $X$ are the mappings $X \to \Ki$ 
which are measurable with respect to the Borel structure on $\Ki$
given by the weak topology.
Then $\int^\oplus_X \Hi_x d\mu(x)$ coincides with
the Hilbert space $L^2(X, \mu, \Ki)$
of square-integrable measurable maps $ X \to \Ki$,
where two such functions are identified if they agree $\mu$-almost everywhere,
equipped with the scalar product 
$$
\langle F_1 \mid F_2 \rangle \, := \, 
\int_{X} \langle F_1(x) \mid F_2(x) \rangle d\mu(x)
\hskip.5cm \text{for all} \hskip.2cm
F_1, F_2 \in L^2(X, \mu, \Ki) .
$$
The space $\int^\oplus_X \Hi_x d\mu(x)$ is called a \textbf{constant field} of Hilbert spaces.
In the special case of a set $X$ with the counting measure,
this space coincides with the space denoted by $\ell^2(X, \Ki)$
in Construction \ref{constructionInd}(2).
\index{$h3$@$L^2(X, \mu, \Ki), L^2(\widehat G, \mu, \Ki)$ Hilbert space}
\end{exe}

The following proposition is often used to reduce questions involving
direct integrals of Hilbert spaces to the situation of constant fields of Hilbert spaces.
For the proof, we refer to Chap. II, \S~1, Proposition 2 in \cite{Dixm--vN}
or Proposition 7.21 in \cite{Foll--16}.

\begin{prop}
% 1.G.2
\label{Pro-DecConstantFieldHilbertSpaces}
Let $\Hi \, = \, \int_X^\oplus \Hi_x d\mu(x)$ be a direct integral of Hilbert spaces
over the measure space $(X,\mu)$.
For every $n =1, 2, \dots, \infty$,
let 
$$
X_n \, = \, \{ x \in X \mid \dim \Hi_x = n \},
$$
and let $\Ki_n$ be a fixed Hilbert space of dimension $n$.
Then, for every $n$, the following holds:
\begin{enumerate}[label=(\roman*)]
\item[(i)]
$X_n$ is measurable;
\item[(ii)]
the direct integral $\int_{X_n}^\oplus \Hi_x d\mu(x)$ of Hilbert spaces
is isomorphic to $L^2(X, \Ki_n)$ over $X_n$, that is,
there exists a family $(U_x)_{x \in X_n}$ of Hilbert space isomorphisms
$U_x \, \colon \Hi_x \to \Ki_n$ such that
$\xi = (\xi_x)_{x \in X_n}$ is a measurable vector field over $X_n$
if and only if $x \mapsto U_x \xi_x$ is a measurable map from $X_n$ to $\Ki_n$. 
\end{enumerate}
As a consequence, we have an isomorphism of Hilbert spaces over $X$
$$
\Hi \, = \, \oplus_n \int_{X_n}^\oplus \Hi_x d\mu(x)
\approx \oplus_n L^2(X_n, \mu, \Ki_n) .
$$
\end{prop}

\begin{defn}
% 1.G.3
\label{defn:decomposableop}
Let $\Hi = \int_X^\oplus \Hi_x d\mu(x)$ be a direct integral of Hilbert spaces
over a measure space $(X, \mu)$. 
For every $x \in X$, let $T_x \in \Li (\Hi_x)$.
Then $x \mapsto T_x$ is a \textbf{measurable field of operators} over $X$ if the map
$x \mapsto \langle T_x \xi_x \mid \eta_x \rangle$
is measurable for all $\xi, \eta \in \Hi$.
Assume that $x \mapsto T_x$ is $\mu$-essentially bounded, that is,
$x \mapsto \Vert T_x \Vert$ is a $\mu$-essentially bounded function;
then a bounded operator $T \, \colon \Hi \to \Hi$ is defined by
$$
(T \xi)_x \, := \, T_x (\xi_x)
\hskip.5cm \text{for} \hskip.2cm
\xi \in \Hi, \hskip.1cm x \in X.
$$
and we write 
$$
T \, = \, \int^\oplus_X T_x d\mu(x).
$$
In this case, the norm $\Vert T \Vert$ is the $\mu$-essential supremum
of the function $x \mapsto \Vert T_x \Vert$ (see \cite[Chap.~II, \S~2, no~3]{Dixm--vN}).
Bounded operators on $\Hi$ of this form
are called \textbf{decomposable} over $X$.
This notion will be sightly generalized in Section \ref{Section-DecomposableOperators}.
\index{Decomposable operator over a Borel space}
\end{defn}

\begin{defn}
% 1.G.4
\label{defn:diagonalisableop}
To every $\varphi \in L^\infty(X, \mu)$
is associated a decomposable operator $m(\varphi)$ 
on $\Hi = \int_X^\oplus \Hi_x d\mu(x)$ defined by
$$
m(\varphi) \, := \, \int^\oplus_X \varphi(x){\mathrm{Id}}_{\Hi_x} d\mu(x).
$$
Thus, $m(\varphi)$ is given by
$$
(m(\varphi) \xi)_x \, = \, \varphi(x) \xi_x
\hskip.5cm \text{for all} \hskip.2cm 
\xi \in \Hi, \hskip.1cm x \in X.
$$
An operator in $\Li (\Hi)$ of the form $m(\varphi)$ for 
$\varphi \in L^\infty(X, \mu)$ is called \textbf{diagonalisable}.
The set $\{m(\varphi)\mid \varphi \in L^\infty(X, \mu)\}$ of diagonalisable
operators on $\Hi$ is obviously a commutative selfadjoint subalgebra of $\Li (\Hi)$.
\index{Diagonalisable operator}
\end{defn}

Let $G$ be a second-countable locally compact group. 
For every $x \in X$, let $\pi_x$ be a representation of $G$ on $\Hi_x$.
Assume that $x \mapsto \pi_x$ is a \textbf{mesurable field of representations} of $G$,
that is,
there is a measurable field $x \mapsto \Hi_x$ of Hilbert spaces over $X$
and, for every $g \in G$,
a measurable field $x \mapsto \pi_x(g)$ of operators on the $\Hi_x$~'s over $X$
such that $g \mapsto \pi_x(g)$ is a representation of $G$ for all $x \in X$. 
Then 
$$
\pi(g) \, := \, \int^\oplus_X \pi_x(g) d\mu(x)
$$
is a unitary operator on $\Hi$ or every $g \in G$.
It is clear that $g \mapsto \pi(g)$ is a homomorphism
from $G$ to the unitary group of $\Hi = \int^\oplus_X \Hi_x d\mu(x)$.
It can be shown \cite[Proposition 18.7.4]{Dixm--C*} that the map $g \mapsto \pi(g)$ 
is strongly continuous, so that $\pi$ is a representation of $G$ on $\Hi$.
We write
$$
\pi \, = \, \int^\oplus_X \pi_x d\mu(x)
$$
and $\pi$ is called the \textbf{direct integral}
of the family $(\pi_x)_{x \in X}$.
\index{Direct integral of representations}

\vskip.2cm

Let $H$ be a closed subgroup of $G$
and $\sigma_1, \sigma_2$ representations of $H$.
The following proposition generalizes the fact (see Proposition~\ref{InductionQqPropr})
that the direct sum $\Ind_H^G \sigma_1 \oplus \Ind_H^G \sigma_2$ 
is equivalent to the induced representation $\Ind_H^G ( \sigma_1 \oplus \sigma_2)$.

\begin{prop}
% 1.G.5
\label{Prop-DirectIntIndRep}
Let $H$ be a closed subgroup of $G$ and $\sigma = \int^\oplus_X \sigma_x d\mu(x)$
a direct integral of a measurable field $(\sigma_x)_{x \in X}$
of representations of $H$ over $X$.
\par

Then $\Ind_H^G \sigma$ is equivalent to
the direct integral $\int^\oplus_X (\Ind_H^G \sigma_x) d\mu(x)$
of the family $(\Ind_H^G \sigma_x)_{x \in X}$.
\end{prop}

For a proof, we refer to \cite[Theorem 10.1]{Mack--52}.

\subsection*{Direct integral decompositions of representations}

We are going to describe how a representation
of a second-countable locally compact group $G$
admits decompositions as direct integrals of representations, in several ways.
The basic tool for this is the following result;
for a proof, see Theorem~8.3.2 and Section~18.7 in \cite{Dixm--C*},
or Theorem 7.29 in \cite{Foll--16}.
\par

Let $(\pi, \Hi)$ be a representation of a group $G$.
The \textbf{commutant} of $\pi(G)$, that is, 
$$
\pi(G)' \, = \, \{ T \in \Li (\Hi) \mid T \pi(g) = \pi(g)T
\hskip.2cm \text{for all} \hskip.2cm
g \in G \},
$$
is a von Neumann subalgebra of $\Li (\Hi)$.

\begin{theorem}
% 1.G.6
\label{Theo-IntDecUniRep1} 
Let $G$ be a second-countable locally compact group
and let $\pi$ be a representation of $G$ on a separable Hilbert space $\Hi$.
Let $\mathcal A$ be an abelian von Neumann subalgebra of the commutant $\pi(G)'$.
\par

Then there exist a standard Borel space $X$,
a $\sigma$-finite measure $\mu$ on $X$,
a measurable field of Hilbert spaces $x \mapsto \Hi_x$ over $X$, 
a measurable field of representations $x \mapsto \pi_x$ of $G$ in the $\Hi_x$~'s over $X$,
and an isomorphism of Hilbert spaces $U \, \colon \Hi \to \int_X^\oplus \Hi_x d\mu(x)$,
with the following properties:
\begin{enumerate}[label=(\arabic*)]
\item\label{iDETheo-IntDecUniRep1}
$
U \pi(g) U^{-1} \, = \, \int^\oplus_X \pi_x(g) d\mu(x)
\hskip.5cm \text{for all} \hskip.2cm
g \in G ;
$
\item\label{iiDETheo-IntDecUniRep1}
$U \mathcal A U^{-1}$ coincides with the algebra of diagonalisable operators on 
$\int_X^\oplus \Hi_x d\mu(x)$.
\end{enumerate}
\end{theorem}

There are two important choices of the algebra $\mathcal A$
in Theorem~\ref{Theo-IntDecUniRep1}
which lead to two integral decompositions with particular properties:
\begin{enumerate}
\item[$\bullet$]
$\mathcal A$ is a maximal abelian subalgebra of $\pi(G)'$;
\item[$\bullet$]
$\mathcal A$ is the centre of $\pi(G)'$.
\end{enumerate}
In the first case, the representations $\pi_x$ in Theorem~\ref{Theo-IntDecUniRep1}
are irreducible for $\mu$-almost every $x \in X$.
More precisely, we have the following Proposition \ref{Pro-IntDecUniMaxAbelian};
for a proof, see \cite[Chap.~II, \S~3, no~1, corollaire 1]{Dixm--vN},
and also \cite[8.5.1]{Dixm--C*}.
(In the second case, when $\mathcal A$ is the centre of $\pi(G)'$, 
the $\pi_x$'s in Theorem~\ref{Theo-IntDecUniRep1}
are factor representations for $\mu$-almost every $x \in X$ and are mutually disjoint.
This leads to the central decomposition of $\pi$, to be treated later;
see Definition \ref{Def-FactorialRep} for factor representations
and Theorem~\ref{thmDirectIntFact} for the decomposition result.)

\begin{prop}
% 1.G.7
\label{Pro-IntDecUniMaxAbelian} 
Let $G$ be a second-countable locally compact group.
Consider a standard Borel space $X$,
a $\sigma$-finite measure $\mu$ on $X$,
a measurable field $x \mapsto \Hi_x$ of Hilbert spaces over $X$, and 
a measurable field $x \mapsto \pi_x$ of representations of $G$ in the $\Hi_x$~'s over $X$.
Set 
$$
\Hi \, = \, \int_X^\oplus \Hi_x d\mu(x)
\hskip.5cm \text{and} \hskip.5cm
\pi \, = \, \int_X^\oplus \pi_x d\mu(x) .
$$
Let $\mathcal A$ be the algebra of diagonalisable operators in $\Li (\Hi)$.
The following properties are equivalent:
\begin{enumerate}[label=(\roman*)]
\item\label{iDEPro-IntDecUniMaxAbelian}
$\mathcal A$ is a maximal abelian subalgebra of $\pi(G)'$;
\item\label{iiDEPro-IntDecUniMaxAbelian}
$\pi_x$ is irreducible for $\mu$-almost every $x \in X$.
\end{enumerate}
\end{prop}

Since, by Zorn's lemma, maximal abelian von Neumann subalgebras of $\pi(G)'$ always exist, 
the following result is an immediate consequence of Theorem~\ref{Theo-IntDecUniRep1}
and Proposition~\ref{Pro-IntDecUniMaxAbelian}.

\begin{theorem}
% 1.G.8
\label{thmDirectIntIrreps}
Let $G$ be a second-countable locally compact group
and let $\pi$ be a representation of $G$ on a separable Hilbert space $\Hi$.
\par

Then there exist a standard Borel space $X$,
a $\sigma$-finite measure $\mu$ on $X$,
a measurable field $x \mapsto \Hi_x$ of Hilbert spaces over $X$,
and a measurable field $x \mapsto \pi_x$ of \emph{irreducible} representations 
of $G$ in the $\Hi_x$'s,
such that $\pi$ is equivalent to the direct
integral $\int^\oplus_X \pi_x d\mu(x)$.
\index{Decomposition of representations! $1$@into irreducible representations}
\end{theorem}

\begin{exe}
% 1.G.9
\label{Exa-DirIntRegRep}
Let $G$ be a second-countable locally compact abelian group. 
Choose the normalizations of the Haar measures $\mu_G$ on $G$
and $\mu_{\widehat G}$ on the dual group $\widehat G$
so that Plancherel's Theorem holds (Theorem~\ref{PlancherelTh}).
Denote by $\pi$ the representation of $G$ on $L^2(\widehat G, \mu_{\widehat G})$
defined by
$$
(\pi(g)\xi) (\chi) \, = \, \chi(g)\xi(\chi)
\hskip.5cm \text{for all} \hskip.2cm
\xi \in L^2(\widehat G, \mu_{\widehat G}), \hskip.1cm
\chi \in \widehat G, \hskip.1cm
g \in G .
$$
Then the Fourier transform
$$
\mathcal{F} \, \colon L^2(G, \mu_G) \, \to \, L^2(\widehat G, \mu_{\widehat G})
$$
is an isomorphism of Hilbert spaces which intertwines
the regular representation $\lambda_G$ on $L^2(G, \mu_G)$
of Example \ref{Leftregrep} and the representation $\pi$:
$$
\mathcal F \lambda_G (g) \, = \, \pi(g) \mathcal F
\hskip.5cm \text{for all} \hskip.2cm
g \in G.
$$
For every $\chi \in \widehat G$, set $\Hi_\chi = \C$;
then
$$
L^2(\widehat G, \mu_{\widehat G}) \, = \, \int^\oplus_{\widehat G} \Hi_\chi d\mu (\chi)
\hskip.5cm \text{and} \hskip.5cm
\pi \, = \, \int^\oplus_{\widehat G} \ \chi d\mu(\chi),
$$
so that $\lambda_G$ is equivalent to 
$\int^\oplus_{\widehat G} \chi d\mu(\chi)$.
\end{exe}

\subsection*{Non uniqueness of irreducible decompositions}

Let $\pi$ be a representation of a second-countable locally compact group $G$
on a separable Hilbert space.
By Theorem~\ref{thmDirectIntIrreps},
$\pi$ is equivalent to a direct integral $\int^\oplus_X \pi_x d\mu(x)$
of irreducible representations $\pi_x$ over a standard Borel space $X$.
It is clear that no uniqueness can be expected
for the choice of the measure space $(X, \mu)$.
For instance, if a measure $\nu$ on $X$ is equivalent to $\mu$,
then $\int^\oplus_X \pi_x d\nu(x)$ is equivalent to $\int^\oplus_X \pi_x d\mu(x)$,
as is easily checked.
\par

One may ask the following more reasonable question.
Assume that $\int^\oplus_{Y} \sigma_y d\nu(y)$
is another direct integral decomposition $\pi$ into irreducible representations
over a standard Borel space $Y$; 
does there exist measurable subsets $N, N'$ of $X$ and $Y$
such that $\mu(N) = \nu(N') = 0$ and 
$$
\{\pi_x \mid x \in X \smallsetminus N\}
\, = \,
\{\sigma_y \mid y \in Y \smallsetminus N' \},
$$
where these two last spaces are viewed as subspaces of $\widehat G$ ? 
\par

This question has a positive answer for some groups,
those of type I; see Section \ref{SectionTypeI} for this notion.
In fact, for groups of type I, one has a refinement of Theorem~\ref{thmDirectIntIrreps}
with a stronger uniqueness statement (Theorem~\ref{thmDirectIntIrreps+} below).
On the contrary, for groups which are not of type I,
the question has always a negative answer, as the following result of Dixmier shows.

\begin{theorem}
% 1.G.10
\label{Theo-NonUniqueIntDecIrrUniRep} 
Let $G$ be a second-countable locally compact group which is not of type I.
\par

There exist a representation $\pi$ of $G$ in a separable Hilbert space
and two direct integral decompositions into irreducible representations
$$
\pi \, = \, \int^\oplus_X \pi_x d\mu(x)
\, = \, \int^\oplus_{Y} \sigma_y d\mu(y)
$$
such that $\pi_x$ and $\sigma_y$ are non equivalent for all $(x,y) \in X \times Y$.
\end{theorem}

The proof of Theorem~\ref{Theo-NonUniqueIntDecIrrUniRep}
can be found \cite[Corollaire~2]{Dixm--64b},
where a more general result is proved
about representations of non type I separable C*-algebras.
We verify below this non-uniqueness property for two concrete examples of representations;
the first example appeared in \cite{Yosh--51} and the second one in \cite{Mack--51}.

\begin{exe}
% 1.G.11
\label{Exa-NonUniqueDecIrrRep}
(1)
Let $F$ be the free group on two generators $a$ and $b$.
Let $A$ be the subgroup generated by $a$ and $B$ the subgroup generated by $b$.
\par

For the regular representation $\lambda_F$ of $F$, we have (by induction in stages)
$$
\lambda_F \, \simeq \, \Ind_{\{e\}}^F 1_{\{e\}}\, \simeq \,
\Ind_{A}^F (\Ind_{\{e\}}^A 1_{\{e\}}) \, \simeq \, \Ind_{A}^F (\lambda_A)
$$
and similarly
$$
\lambda_F \, \simeq \, \Ind_{B}^F (\lambda_B).
$$
We have 
$$
\lambda_A \simeq \int_{\widehat A}^\oplus \chi d\mu (\chi), 
$$
where $\mu$ is a Haar measure on $\widehat A$;
see Example~\ref{Exa-DirIntRegRep}.
Therefore, by Proposition~\ref{Prop-DirectIntIndRep}, we obtain a first integral decomposition 
$$
\lambda_F\simeq \int_{\widehat A}(\Ind_A^F \chi) d\mu (\chi).
$$
Similarly, we have a second integral decomposition
$$
\lambda_F\simeq \int_{\widehat B}(\Ind_B^F \psi) d\nu (\psi),
$$
where $\nu$ is a Haar measure on $\widehat B$.
\par

Now, the induced representations $\Ind_A^F \chi$ and $\Ind_B^F \psi$
are irreducible for every $\chi \in \widehat A$ and $\psi \in \widehat B$
(Example~\ref{Exa-IrredRepFreeGroup}).
Moreover, since $g^{-1} B g \cap A$ has infinite index in $A$
(in fact, $g^{-1} B g \cap A = \{e\}$) for every $g \in F$,
the representations $\Ind_A^F \chi$ and $\Ind_B^F \psi$
are not equivalent for $\chi \in \widehat A$ and $\psi \in \widehat B$,
by Theorem~\ref{Theo-EquiInducedRep}.

\vskip.2cm

\index{$l 2$@$\Q$ rational numbers}
(2)
Let $\Gamma$ be the affine group of the rational field $\Q$, that is, 
$$
\Gamma \, = \, \left(\begin{matrix}
\Q^\times & \Q \\ 0 & 1
\end{matrix} \right)
$$
(see Section~\ref{Section-IrrRepAff}). 
Consider
$$
\text{the normal subgroup} \hskip.2cm
N \, = \, \left(\begin{matrix}
1& \Q \\ 0 & 1
\end{matrix} \right)
\hskip.2cm \text{and the subgroup} \hskip.2cm
H \, = \, \left(\begin{matrix}
\Q^\times &0 \\ 0 & 1
\end{matrix} \right)
$$
of $\Gamma$.
As in (1), we have two direct integral decompositions of the regular representation
$\lambda_\Gamma$ of $\Gamma$:
$$
\lambda_\Gamma \simeq \int_{\widehat N}(\Ind_N^\Gamma \chi) d\mu (\chi)
\hskip.5cm \text{and} \hskip.5cm
\lambda_\Gamma \simeq \int_{\widehat H}(\Ind_H^\Gamma \psi) d\nu (\psi).
$$
It follows from the Mackey--Shoda criteria (compare with Theorem~\ref{Prop-IrredRepAffine})
that the induced representations $\Ind_N^\Gamma \chi$ and $\Ind_H^\Gamma \psi$
are irreducible and non equivalent for $\chi \in \widehat N$ and $\psi \in \widehat H$. 
\index{Affine group! $3$@$\Aff(\Q)$}
\index{Free group}
\end{exe}

\section
{Decomposable operators on fields of Hilbert spaces} 
% Section 1.H
\label{Section-DecomposableOperators}

In this section, we introduce some notions and results
concerning diagonalisable and decomposable operators on fields of Hilbert spaces,
going back to \cite{vNeu--49} and \cite{Gode--51a}.
Our exposition will closely follow \S~2 in Chapter II of \cite{Dixm--vN}.
\par

The tools we introduce will be needed for Section \ref{Section-CanDecRepAbGr},
in particular for Theorem~\ref{Thm-DecRepAbelianGroups}
on the canonical decomposition of representations of abelian groups,
in our construction in Section~\ref{Section-AllIrredRep}
of representations of some semi-direct products $H \ltimes N$,
as well as in Section~\ref{SectionMSC} and Chapter~\ref{Chap:NormalInfiniteRep}.

\vskip.2cm

\index{Measurable field! $3$@measurable field of operators over a Borel space} 
We first generalize the notion of decomposable operator introduced
in Definition \ref{defn:decomposableop}.
Let $(X, \mathcal B)$ be a Borel space, equipped with a $\sigma$-finite measure $\mu$.
Let $x \mapsto \Hi_x$ and $x \mapsto \Ki_x$
be measurable fields of Hilbert spaces over $(X, \mu)$. Set 
$$
\Hi \, = \, \int_X^\oplus \Hi_x d\mu(x)
\hskip.2cm \text{and} \hskip.2cm
\Ki \, = \, \int_X^\oplus \Ki_x d\mu(x).
$$
Let $(T_x)_x$ be a field of operators, $T_x \in \Li (\Hi_x, \Ki_x)$ for each $x \in X$,
such that the map
$$
x \mapsto \langle T_x \xi_x \mid \eta_x \rangle
$$
is measurable for all $\xi \in \Hi, \eta \in \Ki$,
and such that $x \mapsto \Vert T_x \Vert$ is a $\mu$-essentially bounded function.
Let $T \, \colon \Hi \to \Ki$ be the bounded operator defined by
$$
(T \xi)_x \, = \, T_x (\xi_x)
\hskip.5cm \text{for} \hskip.2cm
\xi \in \Hi, \hskip.1cm x \in X;
$$
we write 
$$
T \, = \, \int^\oplus_X T_x d\mu(x).
$$
The norm of $T$ coincides with the $\mu$-essential supremum
of the function $x \mapsto \Vert T_x \Vert$ (see \cite[Chap.~II, \S~2, no~3]{Dixm--vN}).
Bounded operators $\Hi \to \Ki$ obtained in this way
are said to be \textbf{decomposable} over $X$.
\index{Decomposable operator over a Borel space}
\par

Let $\varphi \in L^\infty(X, \mu)$.
Let $m_1(\varphi) \in \Li (\Hi)$ and $m_2(\varphi) \in \Li (\Ki)$
be the associated diagonalisable operators on respectively $\Hi$ and $\Ki$
(see Definition \ref{defn:diagonalisableop}).
It is obvious that 
$$
m_2(\varphi) T \, = \, Tm_1(\varphi)
$$
for every decomposable operator $T \in \Li (\Hi, \Ki)$.
The main result of this section is that this property characterizes the decomposable operators.

\begin{theorem}
% 1.H.1
\label{caractdecop-DirectIntegral}
Let 
$$
\Hi \, = \, \int_X^\oplus \Hi_x d\mu(x)
\hskip.2cm \text{and} \hskip.2cm
\Ki \, = \, \int_X^\oplus \Ki_x d\mu(x)
$$
be two direct integrals of Hilbert spaces over a Borel space $(X, \mathcal B)$,
equipped with a $\sigma$-finite measure $\mu$. Let $T \in \Li (\Hi, \Ki)$ be such that 
$$
m_2(\varphi) T= Tm_1(\varphi)
\hskip.2cm \text{for all} \hskip.2cm
\varphi \in L^\infty(X, \mu).
$$
\par

Then $T$ is a decomposable operator.
\end{theorem}

We will give below the proof of Theorem~\ref{caractdecop-DirectIntegral}
in case $x \mapsto \Hi_x$ and $x \mapsto \Ki_x$
are constant fields of Hilbert spaces over $(X, \mu)$,
that is (see Example~\ref{Exa-DirectIntHilbertSpace}),
in case $\int_X^\oplus \Hi_x d\mu(x) = L^2(X, \mu, \Ki_1)$ and 
$\int_X^\oplus \Ki_x d\mu(x) = L^2(X, \mu, \Ki_2)$,
for fixed Hilbert spaces $\Ki_1$ and $\Ki_2$.
The proof in the general case follows
with the use of Proposition \ref{Pro-DecConstantFieldHilbertSpaces}; 
see Th\'eor\`eme 1, Chap.~II, \S~2 in \cite{Dixm--vN}.
\par

Let $\Ki$ be a Hilbert space.
To every $\varphi \in L^\infty(X, \mu)$ is associated
an operator $m(\varphi)$ on $L^2(X, \mu, \Ki)$ of multiplication by $\varphi$;
see Definition \ref{defn:diagonalisableop}.
The set 
$$
\mathcal A \, := \, \{m(\varphi) \mid \varphi \in L^\infty(X, \mu)\}
$$
of diagonalizable operators on $L^2(X, \mu, \Ki)$
is an abelian $*$-subalgebra of $\Li (L^2(X, \mu, \Ki))$.
The following result shows in particular that $\mathcal A$ is a von Neumann algebra.
%(It is Proposition 7 in \cite[Chap.~II, \S~2, no~4]{Dixm--vN}.)
% [Dixm--vN, page 163].

\begin{prop}
% 1.H.2
\label{Prop-AbelianVN}
Let $(X, \mathcal B)$ be a Borel space 
equipped with a $\sigma$-finite measure $\mu$,
and let $\Ki$ be a Hilbert space.
\begin{enumerate}[label=(\arabic*)]
\item\label{iDEProp-AbelianVN}
Let $L^\infty (X, \mu)$ be equipped with the weak$^*$-topology and 
$\Li (L^2(X, \mu, \Ki))$ with the weak operator topology. The map 
$$
L^\infty (X, \mu) \to \Li (L^2(X, \mu, \Ki)), \hskip.2cm
\varphi \mapsto m(\varphi)
$$
is continuous, of image $\mathcal A$.
\item\label{iiDEProp-AbelianVN}
The subalgebra $\mathcal A$ of diagonalisable operators
is closed in $\Li (L^2(X, \mu, \Ki))$ for the weak operator topology
and is therefore an abelian von Neumann algebra.
\end{enumerate}
\end{prop}

\begin{proof}
\ref{iDEProp-AbelianVN}
For $\varphi \in L^\infty(X, \mu)$ and $F_1, F_2 \in L^2(X, \mu, \Ki)$,
we have
$$
\langle m(\varphi) F_1 \mid F_2 \rangle
\, = \, \int_X \varphi(x) \langle F_1(x) \mid F_2(x) \rangle \hskip.1cm d\mu(x).
$$
Since the function $x \mapsto \langle F_1(x) \mid F_2(x)\rangle$
belongs to $L^1(X, \mu)$,
this shows that, for fixed $F_1, F_2 \in L^2(X, \mu, \Ki)$, the map
$$
\varphi \, \mapsto \, \langle m(\varphi) F_1 \mid F_2 \rangle
$$
is continuous, when $L^\infty (X, \mu)$ is equipped with the weak$^*$-topology.
 
\vskip.2cm

\ref{iiDEProp-AbelianVN}
The unit ball of $L^\infty (X, \mu)$ is weak$^*$-compact.
It follows from \ref{iDEProp-AbelianVN} that the unit ball of $\mathcal A$
is compact in the weak operator topology.
This implies that $\mathcal A$ is closed
in the weak operator topology (see \cite[Chap.~I, \S~3, Th.\ 2]{Dixm--vN}).
\end{proof}

For $\varphi \in L^2(X, \mu)$ and $\xi \in \Ki$,
let $\varphi \otimes \xi$ denote the vector-valued function in $L^2(X, \mu, \Ki)$
defined by 
$$
\varphi \otimes \xi(x) \, = \, \varphi(x) \xi
\hskip.5cm \text{for all} \hskip.2cm 
x \in X.
$$

\begin{lem}
% 1.H.3
\label{Lem-DensHilbert}
Let $(X, \mathcal B)$ and $\mu$ be as above.
Assume that $\Ki$ is a separable Hilbert space and 
let $\{\xi_i \mid i \ge 1 \}$ be a countable total subset of $\Ki$.
\par

For every total subset $D$ of $L^2(X, \mu)$, the set 
$$
\{\varphi \otimes \xi_i \mid \varphi \in D, \hskip.2cm i \ge 1 \}
$$
is a total subset of $L^2(X, \mu, \Ki)$.
\end{lem}

\begin{proof}
Let $F \in L^2(X, \mu, \Ki)$ be such that 
$$
\langle \varphi \otimes \xi_i \mid F \rangle \, = \, 0
\hskip.5cm \text{for all} \hskip.2cm
\varphi \in D, \hskip.2cm i \ge 1.
$$
It suffices to show that $F=0$ in $L^2(X, \mu, \Ki)$.
\par

Let $i \ge 1$. Then 
$$
\int_X \varphi(x) \hskip.1cm \langle \xi_i \mid F(x) \rangle \hskip.1cm d\mu(x)
\, = \, 0
\hskip.5cm \text{for all} \hskip.2cm
\varphi \in D.
$$
Observe that, by Cauchy--Schwarz inequality,
the function $x \mapsto \langle \xi_i \mid F(x) \rangle$ belongs to $L^2(X, \mu)$. 
Since $D$ is a total subset of $L^2(X, \mu)$, 
there exists a subset $N_i \in \mathcal B$
such that $\mu(N_i) = 0$ and $\langle \xi_i \mid F(x) \rangle = 0$
for every $x \in X \smallsetminus N_i$.
Let $N = \bigcup_{i \ge 1} N_i$. Then $\mu(N) = 0$ and we have 
$$
\langle \xi_i \mid F(x) \rangle \, = \, 0 
\hskip.5cm \text{for all} \hskip.2cm 
x \in X \smallsetminus N, \hskip.1cm i \ge 1.
$$
Since $\{\xi_i \mid i \ge 1 \}$ is a total subset of $\Ki$,
we have therefore $F(x) = 0$ for every $x \in X \smallsetminus N$.
\end{proof} 

We now state and prove Theorem~\ref{caractdecop-DirectIntegral}
for constant fields of Hilbert spaces. 

Let $\Ki_1, \Ki_2$ be Hilbert spaces. 
Let $T \, \colon \, X \to \Li (\Ki_1, \Ki_2), \hskip.2cm x \mapsto T(x)$
be a map with the following properties:
\begin{itemize}
\setlength\itemsep{0em}
\item[$\bullet$]
$T$ is measurable, that is, 
$x \mapsto \langle T(x) \xi \mid \eta \rangle$
is measurable for all $\xi \in \Ki_1, \eta \in \Ki_2$;
\item[$\bullet$] $x \mapsto \Vert T(x) \Vert$ belongs to $L^\infty(X, \mu)$.
\end{itemize} 
The operator $\widetilde T$, defined by 
$$
(\widetilde T F) (x) \, = \, T(x) F(x)
\hskip.5cm \text{for all} \hskip.2cm
F \in L^2(X, \mu, \Ki_1)
\hskip.2cm \text{and} \hskip.2cm
x \in X,
$$
is a decomposable operator from $L^2(X, \mu, \Ki_1)$ to $L^2(X, \mu, \Ki_2)$,
with norm equal to the $\mu$-essential supremum of $x \mapsto \Vert T(x) \Vert$.
Every decomposable operator from $L^2(X, \mu, \Ki_1)$ to $L^2(X, \mu \Ki_2)$
is of this form. 

\begin{theorem}
% 1.H.4
\label{caractdecop}
Let $(X, \mathcal B)$ be a Borel space, equipped with a $\sigma$-finite measure~$\mu$,
and let $\Ki_1, \Ki_2$ be \emph{separable} Hilbert spaces.
Let 
$$
\widetilde T \, \colon \, L^2(X, \mu, \Ki_1) \to L^2(X, \mu, \Ki_2)
$$
be a bounded linear operator such that
$$
\widetilde T m_1(\varphi) \, = \, m_2(\varphi) \widetilde T
\hskip.5cm \text{for all} \hskip.2cm
\varphi \in L^\infty(X, \mu).
$$
\par

Then $\widetilde T$ is a decomposable operator.
\end{theorem}
% dans \cite{Dixm--vN}, c'est au Chap.~II, no~5, Th\'eor\`eme 1, 
% page 164 de mon \'edition.

\begin{proof}
Let $\nu$ be a probability measure on $(X, \mathcal B)$ 
which is equivalent to $\mu$ 
(see Proposition \ref{PropEquivProbaMeasure}).
For $i = 1, 2$, the map 
$U_i \, \colon L^2(X, \mu, \Ki_i) \to L^2(X, \nu, \Ki_i)$,
given by 
$$
U_i(F) \, = \, \left(\frac{d\mu}{d\nu}\right)^{1/2} 
\hskip.5cm \text{for all} \hskip.2cm
F \in L^2(X, \mu, \Ki_i),
$$
is an isomorphism of Hilbert spaces. Set
$$
\widetilde S \, = \, U_2 \widetilde T U_1^{-1}
\, \colon L^2(X, \nu, \Ki_1) \to L^2(X, \mu, \Ki_2).
$$
It is clear that $\widetilde S m_1(\varphi) \, = \, m_2(\varphi) \widetilde S$
for every $\varphi \in L^\infty(X, \nu)$ and that $\widetilde S$ is a decomposable
operator if and only if $\widetilde T$ is a decomposable operator.
\par

Therefore, we can assume without loss of generality that $\mu$ is a 
probability measure on $X$. 
\par

Since $\Ki_1$ is separable,
there exists a countable total subset $\{ \xi_j \mid j \ge 1 \}$ of $\Ki_1$. 
For every $\varphi \in L^\infty(X, \mu)$, we have
$$
\widetilde T (\varphi \otimes \xi_j)
\, = \, \widetilde T \left( m_1(\varphi) (\mathbf{1}_X \otimes \xi_j) \right)
\, = \, m_2(\varphi) \widetilde T (\mathbf{1}_X \otimes \xi_j).
$$
Since $L^\infty(X, \mu)$ is dense in $L^2(X, \mu)$
and in view of Lemma~\ref{Lem-DensHilbert},
this shows that $\widetilde T$ is determined by its values on the set
$$
\{ \mathbf{1}_X \otimes \xi_j \mid j \ge 1 \}.
$$
For every $j \ge 1$, set 
$$
F_j \, := \, \mathbf{1}_X \otimes \xi_j \in L^2(X, \mu, \Ki_1)
\hskip.5cm \text{and} \hskip.5cm 
G_j \, := \, \widetilde T(F_j) \in L^2(X, \mu, \Ki_2).
$$
\par

Let $Q$ be a countable dense subset of $\C$.
For $n \ge 1$ and for every $\lambda = (\lambda_1, \hdots, \lambda_n) \in Q^n$,
set 
$$
F_\lambda \, := \, 
\sum_{j = 1}^n \lambda_j F_j
\hskip.5cm \text{and} \hskip.5cm
G_\lambda \, := \, \sum_{j = 1}^n \lambda_j G_j.
$$
Let $n \ge 1$ and $\lambda \in Q^n$.
For every $\varphi \in L^\infty(X, \mu)$, we have 
$$
\left( \widetilde T m_1(\varphi) \right) F_\lambda
\, = \, \left( m_2(\varphi) \widetilde T \right) F_\lambda
\, = \, m_2(\varphi)G_\lambda
$$
and therefore 
$$
\begin{aligned}
\int_X \vert \varphi(x) \vert^2 \Vert G_\lambda(x) \Vert^2 d\mu(x)
\, &= \, \Vert m_2(\varphi)G_\lambda \Vert^2
\, = \, \Vert \left( \widetilde T m_1(\varphi) \right) F_\lambda \Vert^2
\\
\, &\le \, \Vert \widetilde T \Vert^2 \Vert m_1(\varphi)F_\lambda\Vert^2
\, = \, \Vert \widetilde T \Vert^2
\int_X \vert \varphi(x) \vert^2 \Vert F_\lambda(x) \Vert^2 d\mu(x) .
\end{aligned}
$$
Since this inequality holds for every $\varphi \in L^\infty(X, \mu)$,
there exists a subset $N_\lambda \in \mathcal B$ with $\mu(N_\lambda) = 0$
such that 
$$
\Vert G_\lambda(x) \Vert 
\, \le \, \Vert \widetilde T \Vert \Vert F_\lambda(x) \Vert 
$$
for all $x \in X \smallsetminus N_\lambda$.
\par

Set $N := \bigcup_{\lambda} N_\lambda$,
where $\lambda$ runs over the countable set $\bigcup_{n \ge 1} Q^n$.
Then $\mu(N) = 0$ and 
\begin{equation}
\label{eqq/dem5B3}
\tag{$*$}
\Big\Vert \sum_{j = 1}^n \lambda_j G_j(x) \Big\Vert
\, \le \,
\Vert \widetilde T \Vert \hskip.1cm \Big\Vert \sum_{j = 1}^n \lambda_j F_j(x) \Big\Vert 
\end{equation}
for all $x \in X \smallsetminus N$, all $n \ge 1$,
and all $\lambda_1, \hdots, \lambda_n \in Q$.
\par

Let $x \in X \smallsetminus N$.
Since $F_j(x) = \xi_j$ and since $\{\xi_j \mid j \ge 1 \}$ is a total subset of $\Ki_1$,
Inequality (\ref{eqq/dem5B3})
implies that we can define a bounded linear map
$T(x) \, \colon \Ki_1 \to \Ki_2$ such that 
$T(x) \xi_j = G_j(x)$ for every $j \ge 1$,
and then we have $\Vert T(x) \Vert \le \Vert \widetilde T \Vert$.
\par

Setting $T(x) = 0$ for $x \in N$, we obtain a measurable map 
$$
T \, \colon \, X \to \Li (\Ki_1, \Ki_2), \hskip.2cm x \mapsto T(x)
$$
such that 
$$
\Vert T(x) \Vert \, \le \, \Vert \widetilde T \Vert
\hskip.5cm \text{for $\mu$-almost all} \hskip.2cm
x \in X .
$$
Let $\widetilde T' \, \colon L^2(X, \mu, \Ki_1) \to L^2(X, \mu, \Ki_2)$
be the decomposable operator defined by the map $T$.
For every $j \ge 1$,
we have $\widetilde T' (F_j) = G_j$ by definition of $F_j$ and $\widetilde T'$,
hence $\widetilde T'( F_j) = \widetilde T(F_j)$, that is, 
$$
\widetilde T' (\mathbf{1}_X \otimes \xi_j)
\, = \, \widetilde T (\mathbf{1}_X \otimes \xi_j).
$$
Since, as observed above, $\widetilde T'$ and $\widetilde T$
are determined by their values
on the set $\{\mathbf{1}_X \otimes \xi_j \mid j \ge 1 \}$,
we have $\widetilde T = \widetilde T'$. 
\end{proof}

We draw an immediate consequence of Theorem~\ref{caractdecop}.
Recall from Proposition~\ref{Prop-AbelianVN} that,
for a Hilbert space $\Ki$ and a measure space $(X, \mu)$,
the algebra $\mathcal A$
of diagonalisable operators in $\Li (L^2(X, \mu, \Ki))$ is a von Neumann algebra.
%(It is Proposition 7 in \cite[Chap.~II, \S~2, no~4]{Dixm--vN}.)
% [Dixm--vN, page 163].

\begin{cor}
% 1.H.5
\label{Cor-AbelianVN}
Let $(X, \mathcal B)$ be a Borel space 
equipped with a $\sigma$-finite measure $\mu$.
Let $\Ki$ be a separable Hilbert space
and $\mathcal A$ the von Neumann algebra of diagonalisable operators
in $\Li (L^2(X, \mu, \Ki))$. 
\par
The commutant $\mathcal A'$ of $\mathcal A$ in $\Li (L^2(X, \mu, \Ki))$ is
an abelian algebra if and only if $\Ki$ is one-dimensional.
In particular, the von Neumann algebra $\mathcal A$
is \emph{maximal abelian}, that is, $\mathcal A = {\mathcal A}'$,
if and only if $\Ki$ is one-dimensional.
\end{cor}

\begin{proof}
Assume that $\Ki$ is not one-dimensional.
Viewing every operator in $\Li (\Ki)$
as a constant valued decomposable operator in $\Li (L^2(X, \mu, \Ki))$,
we see that $\Li (\Ki)$ is isomorphic to a subalgebra of ${\mathcal A}'$.
Therefore ${\mathcal A}'$ is not abelian. 
\par

Assume that $\Ki$ is one-dimensional.
Then the algebra of decomposable operators
in $\Li (L^2(X, \mu, \Ki)) = \Li (L^2(X, \mu))$
coincides with the algebra $\mathcal A$
of diagonalisable operators.
By Theorem~\ref{caractdecop}, we have therefore 
${\mathcal A}' = \mathcal A$. 
\end{proof} 

\section
{Direct integrals of von Neumann algebras}
% Section 1.I
\label{S: DirectInt-VN}

Given a direct integral $\pi=\int^\oplus_X \pi_x d\mu(x)$ of 
representations of a group $G$, we are interested in the relationship 
between the von Neumann algebra $\pi(G)''$ 
and the family of von Neumann algebras $(\pi_x(G)'')_{x \in X}$.
\par

Let $\Hi = \int_X^\oplus \Hi_x d\mu(x)$ be a direct integral of Hilbert spaces
over a measure space $(X, \mu)$, as in Section~\ref{SectionDecomposingIrreps}.
For every $x \in X$, let ${\mathcal M}_x$ be a von Neumann subalgebra
of $ \Li (\Hi_x)$.
We will say that $x \mapsto {\mathcal M}_x$ is a field of von Neumann algebras over $X$.
For general facts about von Neumann algebras, see Appendix~\ref{AppAlgvN}.

\begin{defn}
% 1.I 1
\label{Def-Field-vN}
The field $x \mapsto {\mathcal M}_x$ of von Neumann algebras over $X$
is \textbf{measurable} if there exists a sequence 
$$
x \mapsto T^{(1)}_x, \hskip.2cm x \mapsto T^{(2)}_x, \hskip.2cm \hdots
$$
of measurable fields of operators over $X$ such that,
for $\mu$-almost every $x \in X$,
the von Neumann algebra ${\mathcal M}_x$ is generated by 
the sequence $T^{(1)}_x, T^{(2)}_x, \dots$. 
\index{Measurable field of von Neumann algebras}
\end{defn}

\begin{exe}
% 1.I.2
\label{Exa-measuableFieldvN}
(1)
The field $x \mapsto \C {\rm Id}_x$ 
is obviously a measurable field of von Neumann algebras.

The field $x \mapsto \Li (\Hi_x)$ is also a measurable field of von Neumann algebras.
Indeed, in view of Proposition~\ref{Pro-DecConstantFieldHilbertSpaces},
we may assume that $x \mapsto \Hi_x$ is a constant field
corresponding to a fixed separable Hilbert space $\Ki$.
Let $(S_n)_{n \ge 1}$ be a sequence generating $\Li (\Ki)$.
It suffices to take the sequence of constant fields $(x \mapsto S_n)_{n \ge 1}$
over $X$.

\vskip.2cm

(2)
Let $G$ be a second countable locally compact group and let 
$$
\pi \, = \, \int^\oplus_X \pi_x d\mu(x)
$$
be a direct integral of representations of $G$.
Then $x \mapsto \pi(G)_x''$ is a measurable field of von Neumann algebras.
Indeed, let $(g_n)_{n \ge 1}$ be a dense sequence in $G$.
For every $x \in X$, the sequence $(\pi_x(g_n))_{n \ge 1}$ generates $\pi_x(G)''$.
\end{exe}

In the following proposition, we state some crucial properties
of measurable fields of von Neumann algebras;
for the proof, see Chap. II, \S~2, Proposition 1 and Lemme 1 in \cite{Dixm--vN}. 

\begin{prop}
% 1.I.3
\label{Lemma-FieldvN}
Let $x \mapsto {\mathcal M}_x$ be a measurable field
of von Neumann algebras over $X$.
\begin{enumerate}[label=(\arabic*)]
\item
The set $\mathcal{M}$ of all decomposable operators
$T \, = \, \int^\oplus_X T_x d\mu(x)$ with $T_x \in {\mathcal M}_x$
for $\mu$-almost every $x \in X$ is a von Neumann algebra of $\Li(\Hi)$.
\item
Assume that $X$ is a standard Borel space.
Then the field $x \mapsto {\mathcal M}_x'$ is measurable,
where the von Neumann algebra ${\mathcal M}_x'$ is the commutant of 
$ {\mathcal M}_x$ in $\Li(\Hi_x)$.
\end{enumerate}
\end{prop} 

Next, we introduce the notion of a decomposable von Neumann algebra.

\begin{defn}
% 1.I.4
\label{Def-DecomposableField-vN}
(1)
Let $x \mapsto {\mathcal M}_x$
be a measurable field of von Neumann algebras over $X$.
The von Neumann algebra of all decomposable 
$T = \int^\oplus_X T_x d\mu(x)$ with $T_x \in {\mathcal M}_x$
for $\mu$-almost every $x \in X$ is called
the \textbf{direct integral of the von Neumann algebras} $({\mathcal M})_{x \in X}$
and denoted by $\int^\oplus_X {\mathcal M}_x d\mu(x)$.

\vskip.2cm

(2)
A von Neumann algebra $\mathcal M$ of $\Li(\Hi)$ is \textbf{decomposable} 
if there exists a measurable field $x \mapsto {\mathcal M}_x$
of von Neumann algebras of $X$ such that 
$$
\mathcal M \, = \, \int^\oplus_X {\mathcal M}_x d\mu(x).
$$
\index{Direct integral of von Neumann algebras}
\index{Decomposable von Neumann algebra}
\end{defn}

\begin{exe}
% 1.I.5
\label{Exa-measurableFieldvN}
(1)
The decomposable von Neumann algebra $\int^\oplus_X \C {\rm Id}_x d\mu(x)$ 
is the algebra $\mathcal A$ of diagonalizable operators.
The decomposable von Neumann algebra $\int^\oplus_X {\Li(\Hi_x)}_x d\mu(x)$ 
is the algebra of decomposable operators in $\Li(\Hi)$.

\vskip.2cm

(2)
Let $\pi = \int^\oplus_X \pi_x d\mu(x)$ be a direct integral
of representations of a second countable locally compact group $G$.
The measurable field $x \mapsto \pi(G)_x''$ of von Neumann algebras
is not necessarily decomposable.
For an obvious counterexample, consider the case
where $x \mapsto \Hi_x$ is a constant field given by a Hilbert space $\Ki$ 
and $\pi_X$ is the identity representation of $G$ in $\Ki$ for every $x \in X$. 
\end{exe}

The commutant of a decomposable von Neumann algebra
admits a neat description.
% for the proof, see Chap. II, \S 3, Th\'eor\`eme 4 in \cite{Dixm--vN}. 

\begin{theorem}
% 1.I.6
\label{Theo-DecomposableField-vN}
Assume that $X$ is a standard Borel space
and let $\mathcal M = \int^\oplus_X {\mathcal M}_x d\mu(x)$
be a decomposable von Neumann algebra.
\par

Then $\mathcal M' \, =\int^\oplus_X {\mathcal M}_x' d\mu(x)$
\end{theorem}

\begin{proof}
Let $(T^{(n)})_{n \ge 1}$ be a sequence in $\mathcal M$ which generates $\mathcal M$;
upon enlarging this sequence, we may assume that
${\mathcal M}_x$ is generated by $(T^{(n)}_x)_{n \ge 1}$,
for $\mu$-almost every $x \in X$.
\par

Let $T = \int^\oplus_X T_x d\mu(x)$ be a decomposable operator.
Then $T \in \mathcal M'$ if and only if
$$
TT^{(n)} \, = \, T^{(n)}T
\hskip.5cm \text{and} \hskip.5cm
T{T^{(n)}}^* \, = \, {T^{(n)}}^*T
$$
for every $n \ge 1$. Hence, $T \in \mathcal M'$ if and only if
$$
T_xT^{(n)}_x \, = \, T^{(n)}_xT_x
\hskip.5cm \text{and} \hskip.5cm
T_x{T^{(n)}_x}^* \, = \, {T^{(n)}_x}^*T_x
$$
for every $n \ge 1$ and $\mu$-almost every $x \in X$.
This shows that $T \in \mathcal M'$
if and only if $T_x \in {\mathcal M}_x'$ for $\mu$-almost every $x \in X$. 
Since $\mathcal M'$ is decomposable (Proposition~\ref{Lemma-FieldvN}),
it follows that 
$$
\mathcal M' \, = \, \int^\oplus_X {\mathcal M}_x' d\mu(x).
$$
\end{proof}

We give a necessary and sufficient
for a von Neumann subalgebra of $\Li(\Hi)$ to be decomposable.

\begin{theorem}
% 1.I.7
\label{Theo-CharDec}
Assume that $X$ is a standard Borel space
and let $\mathcal M$ be a von Neumann subalgebra of $\Li(\Hi)$.
The following conditions are equivalent.
\begin{enumerate}[label=(\roman*)]
\item\label{iDETheo-CharDec}
$\mathcal M$ is decomposable over $X$.
\item\label{iiDETheo-CharDec}
We have $\mathcal A \subset \mathcal M \subset \mathcal B$,
where $\mathcal A, \mathcal B$ are respectively the algebras of diagonalizable
and decomposable operators on $\Hi = \int^\oplus_X {\Hi}_x d\mu(x)$.
\end{enumerate}
\end{theorem}

\begin{proof}
It is obvious that \ref{iDETheo-CharDec} implies \ref{iiDETheo-CharDec}.
Assume that \ref{iiDETheo-CharDec} holds.
Recall from Theorem~\ref{caractdecop-DirectIntegral}
that $\mathcal B= \mathcal A'$
and hence $\mathcal B'= \mathcal A$.
Since $\mathcal A \subset \mathcal M$, we have therefore
$\mathcal M' \subset \mathcal B$. 
So, both $\mathcal M$ and $\mathcal M'$ consist of decomposable operators.
\par

Let $(T^{(n)})_{n \ge 1}$ and $(S^{(n)})_{n \ge 1}$ be sequences
generating $\mathcal M$ and $\mathcal M'$ respectively. Write
$$
T^{(n)} \, = \, \int^\oplus_X T_x^{(n)} d\mu(x)
\hskip.2cm \text{and} \hskip.2cm
S^{(n)} \, = \, \int^\oplus_X S_x^{(n)} d\mu(x).
$$
For every $x \in X$, let $\mathcal M_x$
denote the von Neumann algebra generated by $(T^{(n)}_x)_{n \ge 1}$. 
We claim that $\mathcal M = \int^\oplus_X {\mathcal M}_x d\mu(x)$.
\par

Observe first that $S_x^{(n)}$ and ${S^{(n)}_x}^*$
commute with $T_x^{(m)}$ and ${T^{(m)}_x}^*$
for every $n, m \ge 1$ and for $\mu$-almost every $x \in X$.
\par

Let $T = \int^\oplus_X T_x d\mu(x)$ with $T_x \in {\mathcal M}_x$
for $\mu$-almost every $x \in X$. Then, for every $n \ge 1$,
the operator $T_x$ commutes with $S_x^{(n)}$ and ${S^{(n)}_x}^*$
for $\mu$-almost every $x \in X$;
hence, $T$ commutes with $S^{(n)}$ and ${S^{(n)}}^*$ for every $n \ge 1$. 
Therefore, $T \in (\mathcal M')' = \mathcal M$.
So, we have shown that 
$$
\int^\oplus_X {\mathcal M}_x d\mu(x) \, \subset \, \mathcal M.
$$
\par

By Proposition~\ref{Lemma-FieldvN}, $x \mapsto {\mathcal M}_x'$
is a measurable field of von Neumann algebras;
moreover, by Theorem~\ref{Theo-DecomposableField-vN},
the commutant of $\int^\oplus_X {\mathcal M}_x d\mu(x)$
is $ \int^\oplus_X {\mathcal M}_x'd\mu(x)$.
So, to prove that $\mathcal M$
is contained in $\int^\oplus_X {\mathcal M}_x d\mu(x)$,
it suffices to show that 
$$
\int^\oplus_X {\mathcal M}_x'd\mu(x) \, \subset \, \mathcal M'.
$$ 
Let $S = \int^\oplus_X S_x d\mu(x)$ with $S_x \in {\mathcal M}_x'$
for $\mu$-almost every $x \in X$. Then, for every $n \ge 1$,
the operator $S_x$ commutes with $T_x^{(n)}$ and ${T^{(n)}_x}^*$
for $\mu$-almost every $x \in X$.
It follows that $S$ commutes with $T^{(n)}$ and ${T^{(n)}}^*$
for every $n \ge 1$.
Therefore, $S \in \mathcal M'$.
\end{proof}

Recall from Example~\ref{Exa-measurableFieldvN} 
that, if $\int^\oplus_X \pi_x d\mu(x)$ is a direct integral
of representations of a second countable locally compact group $G$,
then $x \mapsto \pi_x(G)''$ is a measurable field of von Neumann algebras.
The following corollary is therefore an immediate consequence
of Theorem~\ref{Theo-CharDec} and Theorem~\ref{Theo-DecomposableField-vN}.

\begin{cor}
% 1.I.8 
\label{Cor-Theo-CharDec}
Assume that $X$ is a standard Borel space and let 
$\pi=\int^\oplus_X \pi_x d\mu(x)$ be a direct integral of 
representations of a second countable locally compact group $G$
on $\Hi = \int^\oplus_X {\Hi}_x d\mu(x)$.
The following conditions are equivalent.
\begin{enumerate}[label=(\roman*)]
\item
$\pi(G)''=\int^\oplus_X \pi_x(G)''d\mu(x)$.
\item
$\pi(G)''$ contains the algebra $\mathcal A$
of diagonalizable operators on $\Hi$.
\end{enumerate}
Moreover, when these conditions hold, we have 
$$
\pi(G)' \, = \, \int^\oplus_X \pi_x(G)'d\mu(x).
$$
\end{cor}

\begin{rem}
% 1.I.9
\label{Rem-CentralDec}
We will apply Corollary~\ref{Cor-Theo-CharDec} 
in the case where $\int^\oplus_X \pi_x d\mu(x)$
is the central decomposition of a representation $\pi$ of $G$
(see Section~\ref{SectionDecomposingFact}).
\end{rem}

%-----------------------------------------------------------------------
% End of chapter 1
%-----------------------------------------------------------------------

\chapter[Representations of LC abelian groups]
{Representations of locally compact abelian groups}
% Chapter 2
\label{Chapter-AbelianGroups}

\emph{In this chapter, we discuss representations
of a second-countable locally compact abelian group $G$,
irreducible or not.
In Section~\ref{Section-CanRepAbGr}, we introduce
the canonical representation $\pi_\mu$ of $G$
attached to a probability measure $\mu$ on $\widehat G$,
and show that every representation of $G$
is equivalent to a direct sum of such representations (Corollary~\ref{Cor-RepAbCyclic}).
}
\par

\emph{
This result is refined in Section~\ref{Section-CanDecRepAbGr},
where we establish that a representation of $G$ in a separable Hilbert space
is equivalent to a representation determined by
a set of multiplicities $I \subset \overline{\N^*}$
and a sequence $(\mu_m)_{m \in I}$
of mutually singular probability measures on $\widehat G$
(Theorem~\ref{Thm-DecRepAbelianGroups}).
This refinement is a much deeper result; 
on the one hand, it can be seen as a generalization
of the multiplicity theorem of Hahn--Hellinger,
which classifies selfadjoint operators up to unitary equivalence
(see Remark~\ref{Rem-CanonicalDec-Type I});
on the other hand, it admits an extension to representations of type I groups
(Theorem \ref{thmDirectIntIrreps+}).
}
\par

\emph{
The SNAG Theorem, which is the theme of Section~\ref{SNAG},
shows that representations of $G$ are in a bijective correspondence
with projection-valued measures on $\widehat G$. 
In Section \ref{SecCwcpvm} and \ref{SecCandecpm},
we reformulate containment and weak containment properties
as well as the canonical decomposition result of representations
of a second-countable locally compact abelian group $G$
in terms of projection-valued measures on~$\widehat G$.
}

\emph{
These facts are useful to understand irreducible representations
of some non-abelian groups having abelian normal subgroups
(Section \ref{SNAG-NonAbelian})
}

\section
{Canonical representations of abelian groups}
 % Section 2.A
\label{Section-CanRepAbGr}

\index{LCA for ``locally compact abelian''}
Let $G$ be a locally compact abelian group, for short a \textbf{LCA group}.
Let $\widehat G$ its dual group, which is itself a LCA group.
From now on in this chapter and until Section~\ref{SecCandecpm}, 
$$
\text{we assume that $G$ is a \emph{second-countable} locally compact abelian group.}
\leqno{(*)}
$$
For some reminder on LCA groups, see Appendix \ref{AppLCA+Pont}.
\par

In this chapter, a ``measure'' on $G$ or on $\widehat G$
is always meant to be a measure on the Borel $\sigma$-algebra,
$\mathcal B (G)$ or $\mathcal B (\widehat G)$ respectively,
which is finite on compact subsets.
In particular, because of ($*$), such a measure is a Radon measure
(see Theorem \ref{regmeas2ndc}).

\begin{constr}
% 2.A.1
\label{defpimupourGab}
Let $\mu$ be a Radon measure on $\widehat G$.
Note that the Hilbert space $L^2(\widehat G, \mu)$ is separable
since $\widehat G$ is second-countable
(see Proposition~\ref{Prop-DualAbelianSecondCountGroup}
and Theorem~\ref{regmeas2ndc}).
\par

The \textbf{canonical representation} of $G$ on $L^2(\widehat G, \mu)$
is the representation $\pi_\mu$ defined by
\begin{equation}
\label{eqq/exemplepimu}
\tag{$\sharp$}
(\pi_\mu(g)f) (\chi) \, = \, \chi(g) f(\chi)
\hskip.5cm \text{for all} \hskip.2cm
g \in G, \hskip.1cm f \in L^2(\widehat G, \mu), 
\hskip.2cm \text{and} \hskip.2cm
\chi \in \widehat G .
\end{equation}
\index{Canonical representation $\pi_\mu$ of a LCA group $G$
on $L^2(\widehat G, \mu$)}
\index{Representation! canonical (for a LCA group)}
For a bounded continuous function $\varphi$ on $\widehat G$,
we will often make no notational distinction between $\varphi$
and its equivalence class in $L^\infty(\widehat G, \mu)$
and denote as in Definition \ref{defn:diagonalisableop}
by $m(\varphi)$ the multiplication operator 
$$
m(\varphi) \, \colon \, L^2(\widehat G, \mu) \to L^2(\widehat G, \mu),
\hskip.2cm 
f \mapsto \varphi f
$$
associated to $\varphi$.
In particular, for $g \in G$, the map $\widehat g$
defined on $\widehat G$ by $\widehat g (\chi) = \chi(g)$ for $\chi \in \widehat G$
is a unitary character on $\widehat G$ and we have
$$
\pi_\mu(g) \, = \, m(\widehat g).
$$
\par

Let us check that $\pi_\mu$ is indeed a representation of $G$.
It is clear that $\pi_\mu(g)$ is a unitary operator on $L^2(\widehat G, \mu)$
for every $g \in G$
and that the map $g \mapsto \pi_\mu(g)$ is a group homomorphism.
The issue is to check that the map
$\pi_\mu \, \colon G \to \U (L^2(\widehat G, \mu))$ is continuous.
\par

Let $f \in L^2(\widehat G, \mu)$. By polarization it suffices to show that
$g \mapsto \langle \pi_\mu(g)f \mid f \rangle$
is a continuous function on $G$.
For every $g \in G$, we have
$$
\langle \pi_\mu(g) f \mid f \rangle
\, = \, \int_{\widehat G} \chi(g) \vert f(\chi) \vert^2 d\mu(\chi)
\, = \, \int_{\widehat G} \widehat g (\chi) \vert f(\chi) \vert^2 d\mu(\chi)
\, = \, \overline{\mathcal{F}}(\mu_f) (g),
$$
where $\mu_f$ is the finite positive measure on $\widehat G$
defined by $d\mu_f(\chi) = \vert f(\chi) \vert^2d\mu(\chi)$
and $\overline{\mathcal{F}}(\mu_f)$ its inverse Fourier-Stieltjes transform
(see Appendix~\ref{AppLCA+Pont}).
Since $\overline{\mathcal{F}}(\mu_f)$ is a continuous function on $G$,
this concludes the proof.
\end{constr}

\begin{prop}
% 2.A.2
\label{Prop:equivcanrepLCA}
Let $\mu, \nu$ two Radon measures on $\widehat G$ which are equivalent to each other.
\par

Then the canonical representations $\pi_\mu, \pi_\nu$ are equivalent to each other.
\end{prop}

About the converse, see Proposition
\ref{Prop-CanRepAb-Subrepr}~\ref{iiDEProp-CanRepAb-Subrepr}.
 
\begin{proof}
Consider the Radon--Nikodym derivative $\frac{d\nu}{d\mu}$
of $\nu$ with respect to $\mu$.
Define a linear operator
$$
T \, \colon \, L^2(\widehat G, \mu) \to L^2(\widehat G, \nu),
$$
by 
$$
T(f) (\chi) \, = \, f(\chi) \sqrt{ \frac{d\mu}{d\nu}(\chi) }
\hskip.5cm \text{for all} \hskip.2cm
f \in L^2(\widehat G, \mu), \hskip.1cm \chi \in \widehat G.
$$
The map $T$ is invertible;
indeed, $T^{-1}$ is given by the analogous formula,
with $\frac{d\nu}{d\mu}(\chi)$ instead of $\frac{d\mu}{d\nu}(\chi)$.
The map $T$ is also isometric:
for every $f \in L^2(\widehat G, \mu)$, we have
$$
\Vert T(f) \Vert^2
\, = \, \int_{\widehat G} \vert f(\chi)\vert^2 \frac{d\mu}{d\nu}(\chi) d\nu(\chi)
\, = \, \int_{\widehat G} \vert f(\chi)\vert^2 d\mu(\chi)
\, = \, \Vert f \Vert^2.
$$
Moreover, $T$ intertwines the representations $\pi_\mu$ and $\pi_\nu$:
$$
\begin{aligned}
\big( T\pi_\mu(g) (f) \big) (\chi)
\, &= \, \big( \pi_\mu(g) (f) \big) (\chi) \sqrt{ \frac{d\mu}{d\nu}(\chi) }
\, = \, \chi(g) f(\chi) \sqrt{ \frac{d\mu}{d\nu}(\chi) }
\\
\, &= \, \chi(g) (T (f)) (\chi)
\, = \, \big( \pi_\nu(g) (T (f)) \big) (\chi)
\end{aligned}
$$
for $g \in G$ and $\chi \in \widehat G$.
\end{proof}

\begin{exe}
% 2.A.3
\label{expimupourGab}
(1)
If $\mu_{\widehat G}$ is a Haar measure on $\widehat G$,
then $\pi_{\mu_{\widehat G}}$ is equivalent to the regular representation $\lambda_G$ of~$G$.
More precisely,
let $\mu_G$ be the Haar measure on $G$ for which
the Fourier transform
$\mathcal F \, \colon L^2(G, \mu_G) \to L^2(\widehat G, \mu_{\widehat G})$ 
is unitary;
as already noted in Example \ref{Exa-DirIntRegRep},
we have $\pi_{\mu_{\widehat G}} (g) \mathcal F = \mathcal F \lambda_G (g)$ for all $g \in G$.

\vskip.2cm

(2)
Let $\chi_0 \in \widehat G$,
and $\mu$ the Dirac measure at $\chi_0$.
Then $L^2(\widehat G, \mu) = \C$ and $\pi_\mu = \chi_0$.
\end{exe}

\begin{rem}
% 2.A.4
\label{rempimupourGab}

Let $G$ be s above a second-countable LCA group.
In the rest of this chapter,
$$
\begin{aligned}
&
\text{we will consider canonical representations $\pi_\mu$ of $G$}
\\
& \hskip1.5cm
\text{for \emph{probability} measures on $\widehat G$ only.}
\end{aligned}
\leqno{(**)}
$$
This is because we are interested in representations of $G$ up to equivalence only,
and because every Radon measure on $\widehat G$
is equivalent to a probability Radon measure on $\widehat G$
(Proposition \ref{PropEquivProbaMeasure}).
\end{rem}

As we will see, representations of the form $\pi_\mu$
are the building blocks of the most general representations of $G$,
in the sense that every representation of $G$ in a separable Hilbert space
is equivalent to a direct sum $\bigoplus_{k} \pi_{\mu_k}$.
\par

Our first goal is to show that,
given a sequence $(\mu_k)_k$ of probability measures on $\widehat G$
singular with each other,
the direct sum $\bigoplus_{k} \pi_{\mu_k}$ is itself a representation of the form $\pi_\mu$,
for some probability measures $\mu$ on $\widehat G$.

\vskip.2cm

\begin{prop}
% 2.A.5
\label{Prop-RepAbDirectSum}
Let $(\mu_k)_{k \ge 1}$ be a (possibly finite) sequence
of probability measures on $\widehat G$,
singular with each other.
Let $(c_k)_{k \ge 1}$ be a sequence of positive real numbers 
with $\sum_{k \ge 1} c_k = 1$;
set $\mu = \sum_{k \ge 1} c_k \mu_k$.
\par

Then $\pi_\mu$ is equivalent to $\bigoplus_{k \ge 1} \pi_{\mu_k}$.
\end{prop}

\begin{proof}
Since the $\mu_k$~'s are mutually singular probability measures,
we can find a sequence $(A_k)_{k \ge 1}$ of pairwise disjoint subsets
in $\mathcal B (\widehat G)$ such that $\mu_k(A_k) = 1$ for every $k \ge 1$. 
Let 
$$
T \, \colon \, L^2(\widehat G, \mu) \to \bigoplus_{k \ge 1} L^2(\widehat G, \mu_k),
\hskip.2cm
f \mapsto \bigoplus_{k \ge 1} \sqrt{c_k} \Un_{A_k} f.
$$
Observe that, since every $\mu_k$ is absolutely continuous with respect to $\mu$,
the map $T$ is well-defined.
Denoting by $\Vert \cdot \Vert_\mu$ and 
$\Vert \cdot \Vert_{\mu_k}$ the norms on 
$L^2(\widehat G, \mu)$ and $ L^2(\widehat G, \mu_k)$, we have
$$
\begin{aligned}
\Vert f \Vert^2_\mu
\, &= \, \int_{\widehat G}
\vert f(\chi) \vert^2 d\mu(\chi)
\, = \, \sum_{k \ge 1} \int_{A_k} \vert \sqrt{c_k} \Un_{A_k}(\chi) f(\chi) \vert^2 d\mu_k(\chi)
\\
\, &= \, \sum_{k \ge 1} \Vert \sqrt{c_k} \Un_{A_k} f \Vert^2_{\mu_k}
\, = \, \Vert Tf \Vert^2 ,
\end{aligned}
$$
hence $T$ is an isometry. 
\par

Moreover, it is straightforward that $T$ intertwines
$\pi_\mu$ and $\bigoplus_{k \ge 1} \pi_{\mu_k}$.
\par

Finally, $T$ is surjective;
more precisely, one checks that the operator given by
$$
\bigoplus_{k \ge 1} L^2(\widehat G, \mu_k) \to L^2(\widehat G, \mu),
\hskip.2cm
\bigoplus_{k \ge 1} f_k \mapsto \sum_{k \ge 1} \frac{1}{\sqrt{c_k}} \Un_{A_k} f_k
$$
is bounded, and is the inverse $T^{-1}$ of $T$.
\end{proof}

We now show that the representations of the form $\pi_\mu$
are exactly the cyclic representations of $G$ (up to equivalence). 

\begin{prop}
% 2.A.6
\label{Prop-RepAbCyclic}
Let $G$ be as above a second-countable LCA group.
\begin{enumerate}[label=(\arabic*)]
\item\label{iDEProp-RepAbCyclic}
Let $\mu$ be a probability measure on $\widehat G$.
Then $\pi_\mu$ is a cyclic representation, with cyclic vector $\Un_{\widehat G}$. 
\item\label{iiDEProp-RepAbCyclic}
Every cyclic representation of $G$ is equivalent to a representation $\pi_\mu$
for a probability measure $\mu$ on $\widehat G$.
\end{enumerate}
\end{prop}

\begin{proof}
\ref{iDEProp-RepAbCyclic}
Since 
$$
\{ \pi_\mu(g)\Un_{\widehat G} \mid g \in G \} \, = \, \{ \widehat g \mid g \in G \},
$$
the claim follows from
Lemma~\ref{Lemma-DensityFourierTransform}~\ref{1DELemma-DensityFourierTransform},
case $p = 2$.

\vskip.2cm

\ref{iiDEProp-RepAbCyclic}
Let $(\pi, \Hi)$ be a cyclic representation of $G$
and let $\xi \in \Hi$ be a cyclic unit vector for $\pi$.
The function 
$$
\varphi \, \colon \, G \to \C , \hskip.2cm g \mapsto \langle \pi(g) \xi \mid \xi \rangle
$$
is a normalized function of positive type on $G$ (see~\ref{S-FPosType}). 
Therefore, by Bochner's theorem \ref{TheoremBochner},
there exists a probability measure $\mu$ on $\widehat G$ such that 
$\varphi = \overline{\mathcal{F}}(\mu)$. 
\par

We claim that $\pi$ is equivalent to $\pi_\mu$. Indeed, we have 
$$
\langle\pi_\mu(g)\Un_{\widehat G} \mid \Un_{\widehat G} \rangle 
\, = \, \int_{\widehat G} \widehat g (\chi) d\mu(\chi)
\, = \, \overline{\mathcal{F}}(\mu)
\, = \, \varphi(g),
\hskip.5cm \text{for every} \hskip.2cm
g \in G ,
$$
i.e., $\varphi$ coincides with the function of positive type
$\langle \pi_\mu(\cdot) \Un_{\widehat G} \mid \Un_{\widehat G} \rangle$.
Since $\xi$ and $\Un_{\widehat G}$ are cyclic vectors
for respectively $\pi$ and $\pi_\mu$,
the claim follows from Proposition~\ref{GNSbijP(G)cyclic}.
\end{proof}

\begin{cor}
% 2.A.7
\label{Cor-RepAbCyclic}
Let $\pi$ a representation of $G$ on a separable Hilbert space $\Hi$.
\par

There exists a (possibly finite) sequence $(\mu_k)_{k \ge 1} $
of probability measures on $\widehat G$
such that $\pi$ is equivalent to the direct sum $\bigoplus_{k \ge 1}\pi_{\mu_k}$. 
\end{cor}

\begin{proof}
Since the corollary is trivial when $\Hi = \{0\}$,
we assume from now on that $\Hi$ is not $0$-dimensional.
By Proposition \ref{Prop-RepAbCyclic}~\ref{iiDEProp-RepAbCyclic},
it suffices to show that $\Hi$ can be written as
a direct sum $\bigoplus_{k \ge 1} \Ki_k$
of $\pi(G)$-invariant closed subspaces $\Ki_k$
such that the restriction of $\pi$ to $\Ki_k$ is cyclic. 
\par

Let $\{\xi_n \mid n \ge 1 \}$ be a countable total subset of unit vectors in $\Hi$.
We construct by induction on $n \ge 1$
a strictly increasing sequence $(k_n)_{n \ge 1}$ of positive integers,
and a sequence $(\Ki_n)_{n \ge 1}$ 
of $\pi(G)$-invariant cyclic and mutually orthogonal subspaces of $\Hi$
such that, for all $n \ge 1$, we have
$\xi_k \in \bigoplus_{i = 1}^n \Ki_i$ for all $k \le k_n$.
[We leave it to the reader to formulate the slight changes necessary below
in case the construction provides a finite sequence $(k_n)_{n \ge 1}$,
possibly with the last $k_n$ the extended positive integer $\infty$.]
\par

Let $\Ki_1$ be the closed subspace of $\Hi$ generated by $\pi(G)\xi_1$;
set $k_1 = 1$.
Then $\Ki_1$ is $\pi(G)$-invariant, it contains $\xi_1$, 
and the restriction of $\pi$ to $\Ki_1$ is cyclic.
\par

Assume now that, for $n \ge 1$,
sequences $(k_1 < \hdots < k_n)$ and $(\Ki_1, \hdots, \Ki_n)$ have been constructed,
with the desired properties. 
Let $\eta_{k_n+1}$ be the projection of $\xi_{k_n+1}$ on
the orthogonal complement of $\bigoplus_{i = 1}^n \Ki_i$. 
Two cases may occur.
\begin{enumerate}
\item[$\bullet$]
$\eta_{k_n+1} = 0$. Then $\xi_{k_n+1} \in \bigoplus_{i = 1}^n \Ki_i$
and we may replace $k_n$ by $k_n+1$. 
\item[$\bullet$]
$\eta_{k_n+1} \ne 0$. Let then $\Ki_{n+1}$ be the closed subspace of $\Hi$
generated by $\pi(G)\eta_{k_n+1}$.
Then $\Ki_{n+1}$ is $\pi(G)$-invariant, cyclic,
and orthogonal to $\bigoplus_{i = 1}^n \Ki_i$.
The projection of $\xi_{k_n+1}$ on the orthogonal complement
of $\bigoplus_{i = 1}^n \Ki_i$, which is $\eta_{k_n+1}$,
as well as the projection of $\xi_{k_n+1}$ on $\bigoplus_{i = 1}^n \Ki_i$,
belong both to $\bigoplus_{i = 1}^{n+1} \Ki_i$.
Therefore, $\xi_{k_n+1} \in \bigoplus_{i = 1}^{n+1} \Ki_i$. 
\end{enumerate}
This provides the construction of the desired sequences.
\par
 
Since $\xi_k \in \bigoplus_{i \ge 1} \Ki_i$ for every $k \ge 1$,
and since $\{\xi_k \mid k \ge 1 \}$ is total in $\Hi$, we have
$$
\Hi \, = \, \bigoplus_{i \ge 1} \Ki_i 
$$
and the corollary follows.

\vskip.2cm

\emph{Note.}
Using Zorn's lemma,
% [Sunders, 3.4.9]
it is easy to show that any representation of a topological group in any Hilbert space
is a direct sum of cyclic representations \cite[2.2.7]{Dixm--C*}.
\end{proof}

\begin{rem}
% 2.A.8
\label{Rem-RepAbCyclic}
(1)
Concerning the equivalence classes of the probability measures $\mu_k$ 
of Corollary~\ref{Cor-RepAbCyclic},
no uniqueness can be expected.
For instance, let $\mu_1$ and $\mu_2$
be probability measures on $\widehat G$ with disjoint supports
and let $\mu = (\mu_1 + \mu_2)/2$;
then $\pi_\mu$ is equivalent to $\pi_{\mu_1} \oplus \pi_{\mu_2}$,
by Proposition~\ref{Prop-RepAbDirectSum}.
\par

We will give in Theorem~\ref{Thm-DecRepAbelianGroups}
a refinement of Corollary~\ref{Cor-RepAbCyclic}, with a uniqueness property.
For this, we will need to consider versions ``with multiplicities"
of the canonical representations $\pi_\mu$
of Construction \ref{defpimupourGab}.

\vskip.2cm

(2)
In the special case $G = \R$, Corollary \ref{Cor-RepAbCyclic}
is equivalent to the ``multiplication operator form"
of the spectral theorem for a (possibly unbounded) selfadjoint operator
on a separable Hilbert space (see \cite{Halm--51}).
In fact, Corollary \ref{Cor-RepAbCyclic} 
as well as its refined version (Theorem~\ref{Thm-DecRepAbelianGroups} below)
can both be seen as particular cases 
of a theorem for representations of abelian C*-algebras,
for which we refer to Chapter 2 of \cite{Arve--76},
or Section 3.5 of \cite{Sund--97}.
\end{rem}

Next, we investigate when two representations of the form $\pi_\mu$ have a common subrepresentation.

\begin{prop}
% 2.A.9
\label{Prop-CanRepAb-Subrepr}
Let $\mu_1, \mu_2$ be probability measures on $\widehat G$.
\begin{enumerate}[label=(\arabic*)]
\item\label{iDEProp-CanRepAb-Subrepr}
There exists a non-zero bounded operator
$L^2(\widehat G, \mu_1) \to L^2(\widehat G, \mu_2)$ 
intertwining $ \pi_{\mu_1}$ and $\pi_{\mu_2}$ if and only if 
$\mu_1$ and $\mu_2$ are not singular with each other.
\par
In particular, if $\mu_1$ and $\mu_2$ are singular with each other, then 
the representations $ \pi_{\mu_1}$ and $ \pi_{\mu_2}$ are disjoint.
\item\label{iiDEProp-CanRepAb-Subrepr}
The representations $\pi_{\mu_1}$ and $\pi_{\mu_2}$ are equivalent if and only if 
the measures $\mu_1$ and $\mu_2$ are equivalent.
\end{enumerate}
\end{prop}

\begin{proof}
\ref{iDEProp-CanRepAb-Subrepr}
Assume that there exists a non-zero bounded operator
$$
T \, \colon \, L^2(\widehat G, \mu_1) \to L^2(\widehat G, \mu_2)
$$
which intertwines $\pi_{\mu_1}$ and $\pi_{\mu_2}$. 
Then $\Hi_1 := (\ker T)^\perp$ is $\pi_{\mu_1}(G)$-invariant,
the closure $\Hi_2$ of the image of $T$ is $\pi_{\mu_2}(G)$-invariant, and 
there exists an isometric bijective linear map $U \, \colon \Hi_1 \to \Hi_2$ 
intertwining $ \pi_{\mu_1}$ and $ \pi_{\mu_2}$
(see Lemma~\ref{Prop-EquSubRep}).
\par

Let $P \, \colon L^2(\widehat G, \mu_1) \twoheadrightarrow \Hi_1$ be the orthogonal projection;
set $f_1 := P(\Un_{\widehat G})$.
By Proposition~\ref{Prop-RepAbCyclic},
$\Un_{\widehat G}$ is a cyclic vector for $\pi_{\mu_1}$.
Since $P$ intertwines $\pi_{\mu_1}$ with itself,
$f_1$ is a cyclic vector for the restriction of $\pi_{\mu_1}$ to $\Hi_1$.
\par

For $i = 1, 2$ and $\varphi \in C^b(\widehat G)$,
denote by $m_i(\varphi)$
the multiplication operator by $\varphi$ on $L^2(\widehat G, \mu_i)$.
We have 
$$
U m_1(\widehat g) \, = \, m_2(\widehat g) U
\hskip.5cm \text{for all} \hskip.2cm
g \in G.
$$
Set
$$
f_2 \, := \, U(f_1) \in L^2(\widehat G, \mu_2).
$$
For every $g \in G$, we have
$$
U(\widehat g f_1)
\, = \, U(m_1(\widehat g ) (f_1))
\, = \, m_2(\widehat g ) (U(f_1 ))
\, = \, \widehat g f_2.
$$
\par

Consider the non-zero finite measures $\nu_1$ and $\nu_2$
on $\widehat G$ given by 
$$
d\nu_1(\chi) \, := \, \vert f_1(\chi)\vert^2 d\mu_1(\chi)
\hskip.5cm \text{and} \hskip.5cm
d\nu_2(\chi) \, := \, \vert f_2(\chi)\vert^2 d\mu_2(\chi).
$$
For $i = 1, 2$,
denote by $\langle \cdot \mid \cdot \rangle_{L^2(\widehat G, \mu_i)}$
the scalar product on $L^2(\widehat G, \mu_i)$.
For every $g \in G$, we have
$$
\overline{\mathcal{F}}(\nu_1) (g)
\, = \, \int_{\widehat G} \widehat g (\chi) \vert f_1(\chi)\vert^2 d\mu_1(\chi)
\, = \, \langle \widehat g f_1 \mid f_1 \rangle_{L^2(\widehat G, \mu_1)}.
$$
Since $U$ is an isometry, it follows that
$$
\begin{aligned}
\overline{\mathcal{F}}(\nu_1) (g)
\, &= \, \langle \widehat g f_1 \mid f_1 \rangle_{L^2(\widehat G, \mu_1)}
\, = \, \langle U(\widehat g f_1) \mid U(f_1) \rangle_{L^2(\widehat G, \mu_2)}
\\
\, &= \, \langle \widehat g f_2 \mid f_2 \rangle_{L^2(\widehat G, \mu_2)}
\, = \, \int_{\widehat G} \widehat g (\chi) \vert f_2(\chi)\vert^2 d\mu_2(\chi)
\, = \, \overline{\mathcal{F}}(\nu_2) (g) .
\end{aligned}
$$
Therefore we have
$$
\overline{\mathcal{F}}(\nu_1) \, = \, \overline{\mathcal{F}}(\nu_2) ,
$$
hence $\nu_1 = \nu_2$,
by injectivity of the Fourier--Stieltjes transform (\ref{Pro-InjectivityFourier}).
\par

Since $\nu_1$ is absolutely continuous with respect to $\mu_1$,
and $\nu_2$ is absolutely continuous with respect to $\mu_2$,
it follows that $\mu_1$ and $\mu_2$ are not singular with each other.

\vskip.2cm

Conversely, assume that $\mu_1$ and $\mu_2$ are not singular with each other.
There exists a probability measure $\mu$ on $\widehat G$ 
and non-negative functions $\psi_1$ in $L^1(\widehat G, \mu_1)$
and $\psi_2$ in $ L^1(\widehat G, \mu_2)$
such that 
$$
d\mu(\chi) \, = \, \psi_1(\chi) d\mu_1(\chi)
\hskip.5cm \text{and} \hskip.5cm
d\mu(\chi) \, = \, \psi_2(\chi) d\mu_2(\chi).
$$
Upon replacing $\psi_1$ by $\psi_1 / \max_{\chi \in \widehat G} \psi_1(\chi)$, 
we can assume that $\psi_1 \le \Un_{\widehat G}$. 
\par

Define a linear operator
$$
T \, \colon \, L^2(\widehat G, \mu_1) \to L^2(\widehat G, \mu_2),
$$
by 
$$
T(f) (\chi) \, = \, f(\chi) \psi_2(\chi)^{1/2}
\hskip.5cm \text{for all} \hskip.2cm
f \in L^2(\widehat G, \mu_1), \hskip.1cm \chi \in \widehat G.
$$
Then $T$ is a well defined and bounded;
indeed, for every $f \in C^b (\widehat G)$, we have
$$
\begin{aligned}
\Vert T (f) \Vert^2
\, &= \, \int_{\widehat G} \vert f(\chi) \vert^2 \psi_2(\chi) d\mu_2(\chi)
\, = \, \int_{\widehat G} \vert f(\chi) \vert^2 d\mu(\chi)
\\
\, &= \, \int_{\widehat G} \vert f(\chi)\vert^2 \psi_1(\chi)d\mu_1(\chi)
\, \le \, \int_{\widehat G} \vert f(\chi)\vert^2 d\mu_1(\chi)
\, = \, \Vert f \Vert^2.
\end{aligned}
$$
Moreover, $T$ intertwines the representations $\pi_{\mu_1}$ and $\pi_{\mu_2}$:
$$
\begin{aligned}
\big( T \pi_{\mu_1}(g) (f) \big) (\chi)
\, &= \, \big( \pi_{\mu_1}(g) (f) \big) (\chi) \psi_2(\chi)^{1/2}
\, = \, \widehat g (\chi) f(\chi) \psi_2(\chi)^{1/2} 
\\
\, &= \, \widehat g (\chi) (T (f)) (\chi)
\, = \, \big( \pi_{\mu_2}(g) T (f) \big) (\chi)
\end{aligned}
$$
for $g \in G$, $f \in L^2(\widehat G, \mu_1)$, and $\chi \in \widehat G$.
Finally, $T \ne 0$, since $\mu \ne 0$, and hence $\psi_2 \ne 0$.

\vskip.2cm

\ref{iiDEProp-CanRepAb-Subrepr}
If $\mu_1$ and $\mu_2$ are equivalent,
then $\pi_{\mu_1}$ and $\pi_{\mu_2}$ are equivalent,
as was shown in Proposition \ref{Prop:equivcanrepLCA}.

\vskip.2cm
 
Conversely, assume that there exists a bijective isometry
$$
T \, \colon \, L^2(\widehat G, \mu_1) \to L^2(\widehat G, \mu_2)
$$
which intertwines $\pi_{\mu_1}$ and $\pi_{\mu_2}$.
Then $\ker T = \{0\}$ and the first part of the proof of \ref{iDEProp-CanRepAb-Subrepr}
shows that $\mu_1$ is absolutely continuous with respect to $\mu_2$.
The same proof for $T^{-1}$ in place of $T$
shows that $\mu_2$ is absolutely continuous with respect to $\mu_1$.
It follows that $\mu_1$ and $\mu_2$ are equivalent.
\end{proof}

We denote by $\overline{\N^*}$ the set $\{1, 2, \hdots, \infty \}$
of extended positive integers.
\index{$k9$@$\overline{\N^*} = \{1, 2, \hdots, \infty \}$ extended positive integers}

\begin{constr}
% 2.A.10
\label{defpimupourGab++}
Let $\mu$ be a probability measure on $\widehat G$
and $n \in \overline{\N^*}$ an extended positive integer.
We associate to the pair $(\mu, n)$ a representation $\pi_\mu^{(n)}$ of $G$,
defined as follows.
\par

\index{$h3$@$L^2(X, \mu, \Ki), L^2(\widehat G, \mu, \Ki)$ Hilbert space}
Let $\Ki$ be a Hilbert space of dimension $n$. 
Let $L^2(\widehat G, \mu, \Ki)$ denote the Hilbert space
of $L^2$-functions from $\widehat G$ to $\Ki$
(already introduced in Example~\ref{Exa-DirectIntHilbertSpace}).
A representation $\pi_\mu^{(n)}$ of $G$ on $L^2(\widehat G, \mu, \Ki)$
is defined by
$$
(\pi_\mu^{(n)}(g)F) (\chi) \, = \, \chi(g) F(\chi) 
\hskip.5cm \text{for all} \hskip.2cm
g \in G, \hskip.1cm F \in L^2(\widehat G, \mu, \Ki), 
\hskip.2cm \text{and} \hskip.2cm
\chi \in \widehat G .
$$
Note that $\pi_\mu^{(n)}$ is equivalent to the multiple 
$$
n \pi_\mu \, = \,
\underbrace{\pi_\mu \oplus \cdots \oplus \pi_\mu}_{n-\text{times}}
$$
of the canonical representation $\pi_\mu$ 
of Construction \ref{defpimupourGab}.
\end{constr}
 
Given a representation $\pi_\mu^{(n)}$,
we will identify in Section \ref{Section-CanDecRepAbGr}
the commuting algebra
$$
\{ T \in \Li (L^2(\widehat G, \mu, \Ki))
\mid 
T \pi_\mu^{(n)}(g) = \pi_\mu^{(n)}(g) T
\hskip.2cm \text{for all} \hskip.2cm
g \in G \}
$$
and, more generally, the space of intertwining operators between two such representations.

\section[Canonical decomposition of representations]
{Canonical decomposition of representations of abelian groups}
% Section 2.B
\label{Section-CanDecRepAbGr}

Let $G$ be a second-countable locally compact abelian group.
Refining Corollary~\ref{Cor-RepAbCyclic},
we will give in this section a canonical decomposition
for every representation of $G$ on a separable Hilbert space,
in terms of the representations $\pi_\mu^{(n)}$
defined in \ref{defpimupourGab++}.
\par

Given a probability measure $\mu$ on $\widehat G$,
an extended positive integer $n \in \overline{\N^*}$,
and a Hilbert space $\Ki$ of dimension $n$,
recall that the representation $\pi_\mu^{(n)}$
is defined on $L^2(\widehat G, \mu, \Ki)$ by 
$$
(\pi_\mu^{(n)}(g)F) (\chi) \, = \, \chi(g) F(\chi)
\hskip.5cm \text{for} \hskip.2cm
g \in G, \hskip.1cm F \in L^2(\widehat G, \mu, \Ki), 
\hskip.2cm \text{and} \hskip.2cm
\chi \in \widehat G .
$$
For $g \in G$, we have a unitary character
$\widehat g \, \colon \widehat G \to \T$ defined by 
$\widehat g (\chi) = \chi(g)$ for all $\chi \in \widehat G$
and a diagonalisable operator $m(\widehat g)$ on $L^2(\widehat G, \mu, \Ki)$,
as in Definition \ref{defn:diagonalisableop}.
Then we can write
$$
\pi_\mu^{(n)}(g) \, = \, m(\widehat g)
\hskip.5cm \text{for all} \hskip.2cm
g \in G.
$$

\vskip.2cm

We draw from Section~\ref{Section-DecomposableOperators}
some consequences on common subrepresentations
between representations of the form $\pi_\mu^{(n)}$.
The following proposition generalizes parts of Proposition~\ref{Prop-CanRepAb-Subrepr}.

\begin{prop}
% 2.B.1
\label{Prop-RepAb-Subrep} 
Consider two probability measures $\mu_1, \mu_2$ on $\widehat G$,
two extended positive integers $n_1,n_2$ in $\overline{\N^*}$,
and two Hilbert spaces $\Ki_1, \Ki_2$,
respectively of dimension $n_1, n_2$.
\begin{enumerate}[label=(\arabic*)]
\item\label{iDEProp-RepAb-Subrep} 
Assume that $\mu_1, \mu_2$ are singular with each other.
Then the representations $\pi_{\mu_1}^{(n_1)}$ and $\pi_{\mu_2}^{(n_2)}$ are disjoint:
there exists no non-zero bounded operator
between $L^2(\widehat G, \mu_1, \Ki_1)$ and $L^2(\widehat G, \mu_2, \Ki_2)$
which intertwines $\pi_\mu^{(n_1)}$ and $\pi_\mu^{(n_2)}$.
\item\label{iiDEProp-RepAb-Subrep}
Assume that the representations $\pi_{\mu_1}^{(n_1)}$ and $\pi_{\mu_2}^{(n_2)}$
are equivalent.
Then the measures $\mu_1$ and $\mu_2$ are equivalent.
\end{enumerate}
\end{prop}

\begin{proof}
\ref{iDEProp-RepAb-Subrep}
Assume, by contradiction, that there exists a non-zero bounded operator
$$
\widetilde T \, \colon \, L^2(\widehat G, \mu_1, \Ki_1) \to L^2(\widehat G, \mu_2, \Ki_2)
$$
intertwining $\pi_{\mu_1}^{(n_1)}$ and $\pi_{\mu_2}^{(n_2)}$.
Recall that the representation $\pi_{\mu_i}^{(n_i)}$
is equivalent to $n_i \pi_{\mu_i}$, for $i = 1, 2$.
\par

There exists a $\pi_{\mu_1}^{(n_1)}(G)$-invariant closed subspace $\Hi_1$
of $L^2(\widehat G, \mu_1, \Ki_1)$ 
such that the restriction of $\pi_{\mu_1}^{(n_1)}$ to $\Hi_1$
is equivalent to $\pi_{\mu_1} = \pi_{\mu_1}^{(1)}$,
and the restriction of $\widetilde T$ to $\Hi_1$ is not zero.
There exists also a $\pi_{\mu_2}^{(n_2)}(G)$-invariant closed subspace $\Hi_2$ 
of $L^2(\widehat G, \mu_2, \Ki_2)$ such that
the restriction of $\pi_{\mu_2}^{(n_2)}$ to $\Hi_2$
is equivalent to $\pi_{\mu_2} = \pi_{\mu_2}^{(1)}$,
and such that the map
$$
P \circ (\widetilde T \vert_{\Hi_1}) \, \colon \, \Hi_1 \to \Hi_2
$$
is not zero,
where $P \, \colon L^2(\widehat G, \mu_2, \Ki_2) \twoheadrightarrow \Hi_2$
is the orthogonal projection.
Since $P \circ \widetilde T$ intertwines
$\pi_{\mu_1}^{(n_1)}$ and $\pi_{\mu_2}^{(n_2)}$
and since the restrictions of $\pi_{\mu_1}^{(n_1)}$ and $\pi_{\mu_2}^{(n_2)}$
to $\Hi_1$ and 
$\Hi_2$ are equivalent to $\pi_{\mu_1}$ and $\pi_{\mu_2}$, 
this contradicts
Proposition \ref{Prop-CanRepAb-Subrepr}~\ref{iDEProp-CanRepAb-Subrepr}.

\vskip.2cm

\ref{iiDEProp-RepAb-Subrep}
The proof is an extension of the proof of
Proposition~\ref{Prop-CanRepAb-Subrepr}~\ref{iiDEProp-CanRepAb-Subrepr}.
Let 
$$
\widetilde T \, \colon \, L^2(\widehat G, \mu_1, \Ki_1) \to L^2(\widehat G, \mu_2, \Ki_2)
$$
be a bijective linear isometry intertwining $\pi_{\mu_1}^{(n_1)}$ and $\pi_{\mu_2}^{(n_2)}$.
\par

For $\varphi \in L^\infty(\widehat G, \mu_i)$,
in particular for $\varphi \in C^b (\widehat G)$, and $\xi \in \Ki_i$, 
let $\varphi \otimes \xi$
denote the vector-valued function in $L^2(\widehat G, \mu_i, \Ki_i)$
defined as in Section~\ref{Section-DecomposableOperators}
by $(\varphi \otimes \xi) (\chi) = \varphi(\chi) \xi$.
Fix a unit vector $\xi \in \Ki_1$ and set
$$
F \, := \, \widetilde T (\Un_{\widehat G}\otimes \xi) \in L^2(\widehat G, \mu_2, \Ki_2)
$$
For every $g \in G$, we have
$$
\widetilde T (\widehat g \otimes \xi) \, = \,
\widetilde T (m_1(\widehat g) (\Un_{\widehat G} \otimes \xi))
\, = \, m_2(\widehat g) (\widetilde T (\Un_{\widehat G} \otimes \xi))
\, = \, m_2(\widehat g ) F.
$$
On the one hand, we have
$$
\begin{aligned}
\overline{\mathcal{F}}(\mu_1) (g)
\, &= \, \int_{\widehat G} \widehat g (\chi)d\mu_1(\chi)
\, = \, \int_{\widehat G} \widehat g (\chi) \Vert \xi \Vert^2 d\mu_1(\chi)
\\
\, &= \, \langle m_1(\widehat g ) (\Un_{\widehat G} \otimes \xi)
\mid \Un_{\widehat G} \otimes \xi \rangle
\hskip.5cm \text{for all} \hskip.2cm
g \in G .
\end{aligned}
$$
On the other hand, consider the bounded Borel measure $\mu$ on $\widehat G$, 
given by 
$$
d\mu(\chi) \, = \, \Vert F(\chi) \Vert^2 d\mu_2(\chi).
$$
Then, since $\widetilde T$ is an isometry, we have
$$
\begin{aligned}
\overline{\mathcal{F}}(\mu) (g)
\, &= \, \int_{\widehat G} \widehat g (\chi) \Vert F(\chi) \Vert^2 d\mu_2(\chi)
\, = \, \langle m_2(\widehat g) F \mid F \rangle
\\
\, &= \, \langle m_2(\widehat g) \widetilde T (\Un_{\widehat G} \otimes \xi) 
\mid \widetilde T (\Un_{\widehat G} \otimes \xi) \rangle
\\
\, &= \, \langle \widetilde T (m_1(\widehat g) (\Un_{\widehat G} \otimes \xi))
\mid \widetilde T (\Un_{\widehat G} \otimes \xi) \rangle
\\
\, &= \, \langle m_1(\widehat g) (\Un_{\widehat G} \otimes \xi))
\mid \Un_{\widehat G} \otimes \xi \rangle
\hskip.5cm \text{for all} \hskip.2cm
g \in G .
\end{aligned}
$$
Therefore, 
$$
\overline{\mathcal{F}}(\mu_1) \, = \, \overline{\mathcal{F}}(\mu).
$$
By injectivity of the Fourier--Stieltjes transform, it follows that 
$$
d\mu_1(\chi) \, = \, d \mu(\chi) \, = \, \Vert F(\chi) \Vert^2 d\mu_2(\chi).
$$
So, $\mu_1$ is absolutely continuous with respect to $\mu_2$. 
\par

The same argument, applied to $T^{-1}$ in place of $T$, shows that 
$\mu_2$ is absolutely continuous with respect to $\mu_1$. 
\end{proof}

Next, we identify the intertwining operators
between two representations of the form 
 $ \pi_\mu^{(n_1)}$ and $ \pi_\mu^{(n_2)}$
 for a given probability measure $\mu$ on $\widehat G$.
 
\begin{prop}
% 2.B.2
\label{Prop-RepAb-Equiv} 
Consider a probability measure $\mu$ on $\widehat G$,
two extended positive integers $n_1,n_2$ in $\overline{\N^*}$,
and two Hilbert spaces $\Ki_1, \Ki_2$ respectively of dimension $n_1, n_2$.
\par

The space of bounded linear operators
$L^2(\widehat G, \mu, \Ki_1) \to L^2(\widehat G, \mu, \Ki_2)$ 
intertwining $\pi_\mu^{(n_1)}$ and $\pi_\mu^{(n_2)}$
consists of the decomposable operators.
In particular, $\pi_\mu^{(n_1)}$ and $\pi_\mu^{(n_2)}$ are equivalent
if and only if $n_1 = n_2$.
\end{prop}

\begin{proof}
For $i = 1, 2$ and $\varphi \in L^{\infty}(\widehat G, \mu)$,
denote by $m_i(\varphi)$
the corresponding multiplication operator on $L^2(\widehat G, \mu, \Ki_i)$.

\vskip.2cm

Suppose first that $\widetilde T$ is decomposable:
$(\widetilde T F) (\chi) = T(\chi) F(\chi)$
for some $\mu$-essentially bounded measurable map
$T \, \colon \widehat G \to \Li (\Ki_1, \Ki_2)$
and for all $F \in L^2(\widehat G, \mu, \Ki_1)$.
Then
$$
(\widetilde T \pi_\mu^{(n_1)}(g) F) (\chi) \, = \,
T(\chi) \chi(g) F(\chi) \, = \,
\chi(g) T(\chi) F(\chi) \, = \,
(\pi_\mu^{(n_2)}(g) \widetilde T F) (\chi)
$$
for all $g \in G$ and $F \in L^2(\widehat G, \mu, \Ki_1)$,
so that $\widetilde T$ intertwines
$\pi_\mu^{(n_1)}$ and $\pi_\mu^{(n_2)}$.

\vskip.2cm

Conversely, suppose now that $\widetilde T$
intertwines $\pi_\mu^{(n_1)}$ and $\pi_\mu^{(n_2)}$,
i.e., that 
$$
\widetilde T m_1(\widehat g) \, = \, m_2(\widehat g) \widetilde T
\hskip.5cm \text{for all} \hskip.2cm
g \in G.
$$
We also have
$$
\widetilde T m_1(\varphi) \, = \, m_2(\varphi) \widetilde T
$$
for all trigonometric polynomial $\varphi \in \mathrm{Trig}(\widehat G)$,
by linearity,
and indeed for all $\varphi \in L^\infty(\widehat G, \mu)$
because $\mathrm{Trig}(\widehat G)$
is weak$^*$-dense in $L^\infty(\widehat G, \mu)$; see
Lemma~\ref{Lemma-DensityFourierTransform}~\ref{2DELemma-DensityFourierTransform}.
It follows from Theorem \ref{caractdecop}
that $\widetilde T$ is decomposable.
\end{proof}

\begin{rem}
% 2.B.3
\label{Rem-Prop-RepAb-Equiv} 
In the particular case $n_1 = n_2 = 1$, Proposition~\ref{Prop-RepAb-Equiv}
shows that $\pi_\mu(G)' = \Hom_G(\pi_\mu, \pi_\mu)$ coincides with
the algebra $\{ m(\varphi) \mid \varphi \in L^\infty(G, \mu) \}$
of diagonalisable operators.
In particular, $\pi_\mu(G)'$ is abelian
and so, with a terminology to be introduced in Section~\ref{Sectioncomppres},
the canonical representation $\pi_\mu$ is \textbf{multiplicity-free}
(see Theorem~\ref{Thm-CanRepAbMultFree}).
\index{Multiplicity-free representation}
\end{rem}

We will need the following two corollaries of Proposition~\ref{Prop-RepAb-Equiv}.
Given a probability measure $\mu$ on $\widehat G$, 
a Hilbert space $\Ki$, and a Borel subset $B$ in $\widehat G$,
we denote by $E_\mu(B)$ the projection operator on 
$L^2(\widehat G, \mu, \Ki))$ defined by 
$$
(E_\mu(B)F) (\chi) \, = \, \Un_B(\chi) F(\chi)
\hskip.5cm \text{for all} \hskip.2cm
F \in L^2(\widehat G, \mu, \Ki), \hskip.1cm
\chi \in \widehat G.
$$

\begin{cor}
% 2.B.4
\label{Cor-Prop-RepAb-Equiv1} 
Let $\mu$ a probability measure on $\widehat G$,
and $\Ki$ a Hilbert space of some dimension $n \in \overline{\N^*}$. 
For a closed subspace $\Hi$ of $L^2(\widehat G, \mu, \Ki)$,
the following properties are equivalent:
\begin{enumerate}[label=(\roman*)]
\item\label{iDECor-Prop-RepAb-Equiv1}
$\Hi$ is invariant under the commutant $\pi_\mu^{(n)}(G)'$ of $\pi_\mu^{(n)}(G)$;
\item\label{iiDECor-Prop-RepAb-Equiv1}
$\Hi$ is the range of a projection $E_\mu(B)$ for $B \in \mathcal B (\widehat G)$.
\end{enumerate}
\end{cor}

\begin{proof}
Assume that $\Hi$ is the range of a projection $E_\mu(B)$ for $B \in \mathcal B (\widehat G)$.
It is clear that $E_\mu(B)$ commutes with 
every decomposable operator on $L^2(\widehat G, \mu, \Ki)$, that is, 
with every element in $\pi_\mu^{(n)}(G)'$, by Proposition~\ref{Prop-RepAb-Equiv}.
Therefore, $\Hi$ is invariant under $\pi_\mu^{(n)}(G)'$.

\vskip.2cm

Conversely, assume that $\Hi$ is invariant under $\pi_\mu^{(n)}(G)'$.
The orthogonal projection $P$ from $L^2(\widehat G, \mu, \Ki)$ to $\Hi$
belongs therefore to the von Neumann algebra $\pi_\mu^{(n)}(G)''$.
Proposition~\ref{Prop-AbelianVN} and the proof of Proposition~\ref{Prop-RepAb-Equiv} 
show that $\pi_\mu^{(n)}(G)''$ coincides
with the algebra $\{m(\varphi)\mid \varphi \in L^\infty(\widehat G, \mu)\}$
of diagonalisable operators on $L^2(\widehat G, \mu, \Ki)$.
Therefore, $P = m(\Un_B) = E_\mu(B)$ for some $B \in \mathcal B (\widehat G)$.
\end{proof}

\begin{notat}
% 2.B.5
\label{notationKn}
For the rest of this section, we denote for each $n \in \overline{\N^*}$
by $\Ki_n$ a Hilbert space of dimension $n$, fixed once for all.
\end{notat}

\begin{cor}
% 2.B.6
\label{Cor-Prop-RepAb-Equiv2} 
Let $I$ a subset of $\overline{\N^*}$,
and $(\mu_n)_{n \in I}$ a sequence
of mutually singular probability measures on $\widehat G$.
Consider the representation
$$
\pi \, = \, \bigoplus_{n \in I} \pi_{\mu_n}^{(n)}
\hskip.5cm \text{of $G$ in the Hilbert space} \hskip.5cm
\Hi \, = \, \bigoplus_{n \in I} L^2(\widehat G, \mu_n, \Ki_n ) .
$$
\par

For every $n \in I$, the orthogonal projection
$P_n$ of $\Hi$ onto $L^2(\widehat G, \mu_n, \Ki_n)$ 
belongs to the centre of the von Neumann algebra $\pi(G)'$.
\end{cor}

\begin{proof}
Since $L^2(\widehat G, \mu_n, \Ki_n)$ is invariant under $\pi(G)$,
we have $P_n \in \pi(G)'$ for every $n \in I$.
\par

Let $T \in \pi(G)'$. Then, $P_k \circ T \circ P_\ell \in \pi(G)'$
for every $k, \ell \in I$. 
For $k \ne \ell $, the measures $\mu_k, \mu_\ell$ are singular with each other
and it follows therefore from Proposition~\ref{Prop-RepAb-Subrep} 
that $P_k \circ T \circ P_\ell = 0$. Therefore
$$
TP_n \, = \, \Big( \sum_{m \in I} P_m \Big) TP_n \, = \,
P_n T P_n \, = \, P_nT \Big( \sum_{m \in I} P_m \Big) \, = \, P_n T,
$$
for every $n \in I$. This shows that $P_n \in \pi(G)''$ for every $n \in I$.
\end{proof}

We will use in the sequel the following notation.
Let $\mu$ be a probability measure on a Borel space $(X, \mathcal B)$
and let $B \in \mathcal B$ be such that $\mu(B) > 0$.
We denote by $\mu \vert B$ the probability measure on $X$
defined by
$$
(\mu \vert B) (C) \, = \, \frac{ \mu(C \cap B)}{\mu(B)} 
\hskip.5cm \text{for every } \hskip.2cm
C \in \mathcal B .
$$
Observe that $\mu \vert B$ is absolutely continuous with respect to $\mu$
and that the Radon--Nikodym derivative $\frac{d (\mu \vert B)}{d\mu}$
vanishes outside $B$.
When $X$ is a disjoint union
$\bigsqcup_{n \ge 1} B_n$ of measurable subsets $B_n \in \mathcal B$,
we have
$$
\mu \, = \, \sum_{n \ge 1, \mu(B_n) > 0} \mu(B_n) \hskip.1cm \mu \vert B_n.
$$

\begin{lem}
% 2.B.7
\label{Lem-RestrMeasure}
Let $\mu$ be a probability measure on $(X, \mathcal B)$
and $f \, \colon X \to \R_+$ a measurable function such that $\int_X f d\mu = 1$.
Set $B_f = \{x \in X \mid f(x) > 0 \}$.
\par

Then $\mu(B_f) > 0$ and $\mu\vert B_f$ is equivalent
to the probability measure $\nu = f \mu$ on $X$.
\end{lem}

\begin{proof}
We have $\nu (B_f) = \int_{B_f} f(x) d\mu(x) = 1$, and therefore $\mu(B_f) > 0$.
For $C \in \mathcal B$, we have
$$
\nu(C) \, = \, \int_{X} \mathbf{1}_C(x) f(x) d\mu(x)
\, = \, \int_{C \cap B_f} f(x) d\mu(x).
$$
Since $f > 0$ on $B_f$, it follows that $\nu(C) = 0$ if and only if $\mu(C \cap B_f) = 0$.
\end{proof}

The main result in this chapter is the following theorem,
which shows that every representation of $G$ in a separable Hilbert space 
admits a decomposition in terms of representations of the form $\pi_\mu^{(n)}$,
\emph{with uniqueness up to the appropriate equivalence.}

\begin{theorem}[\textbf{Canonical decomposition of representations of abelian groups}]
% 2.B.8
\label{Thm-DecRepAbelianGroups}
Let $G$ be a second-countable locally compact abelian group
and $\pi$ a representation of $G$ in a separable Hilbert space $\Hi$.
\par

There exist a set $I$ of extended positive integers
and a sequence $(\mu_n)_{n \in I}$ of probability measures on $\widehat G$
with the following properties:
\begin{enumerate}[label=(\arabic*)]
\item\label{iDEThm-DecRepAbelianGroups}
the measures $\mu_n$ are singular with each other;
\item\label{iiDEThm-DecRepAbelianGroups}
the representation $\pi$ is equivalent to the direct sum
$\bigoplus_{n \in I} \pi_{\mu_n}^{(n)}$.
\end{enumerate}
\par

Moreover, this decomposition is unique in the following sense:
assume that there exist a set $I'$ of extended positive integers
and a sequence $(\mu'_n)_{n \in I'}$ of probability measures on $\widehat G$
which are mutually singular with each other,
such that $\pi$ is equivalent to the direct sum $\bigoplus_{n \in I'} \pi_{\mu'_n}^{(n)}$.
Then $I = I'$,
and $\mu_n$ and $\mu'_n$ are equivalent for every $n \in I$.
\end{theorem}

\begin{proof}
Since the theorem is trivial when $\Hi = \{0\}$,
we assume from now on that $\Hi$ is not $0$-dimensional.
By Corollary~\ref{Cor-RepAbCyclic},
there exist a (possibly finite) sequence $(\lambda_i)_{i \ge 1}$
of probability measures on $\widehat G$ such that 
$\pi$ is equivalent to $\bigoplus_{i \ge 1} \pi_{\lambda_i}$. 
These probability measures are not necessarily mutually singular; 
the strategy is to replace them by probability measures $\mu_n$
which are mutually singular, at the cost of 
considering representations $\pi_{\mu_n}^{(n)}$ with multiplicities $n$. 
Our proof of the existence part of Theorem~\ref{Thm-DecRepAbelianGroups}
is patterned after that of Part (1) of Theorem 3.5.3 in \cite{Sund--97}.

\vskip.2cm

$\bullet$ {\it First step.} 
Choose a sequence of positive real numbers $(c_i)_{i \ge 1}$
such that $\sum_{i \ge 1} c_i = 1$
and define a probability measure $\mu$ on $\widehat G$ by
$$
\mu \, := \, \sum_{i \ge 1} c_i \lambda_i .
$$
For every $i \ge 1$, note that $\lambda_i$ is absolutely continuous with respect to $\mu$
and set
$$
B_i \, = \, \Big\{ \chi \in \widehat G
\hskip.2cm \Big\vert \hskip.2cm
\frac{d\lambda_i}{d\mu} (\chi) > 0 \Big\}
\, \in \, \mathcal B (\widehat G) .
$$
Note that $\mu (B_i) > 0$,
and that the probability measure $\mu \vert B_i$
is equivalent to $\lambda_i = \frac{d\lambda_i}{d\mu} \mu$
(Lemma \ref{Lem-RestrMeasure}).
Since the representation $\pi$ is equivalent to $\bigoplus_{i \ge 1} \pi_{\lambda_i}$,
it is also equivalent to $\bigoplus_{i \ge 1} \pi_{\mu \vert B_i}$
(see Proposition \ref{Prop:equivcanrepLCA}).
\par

Set
$$
B \, := \, \bigcup_{i \ge 1} B_i \, \in \, \mathcal B (\widehat G).
$$
Observe that $\mu(B) = 1$. 
For $n \in \overline{\N^*}$ and $i \ge 1$, set
$$
A_n \, := \, \Big\{\chi \in B
\hskip.2cm \Big\vert \hskip.2cm
\sum_{j \ge 1} \Un_{B_j}(\chi) = n \Big\} 
$$
and
$$
B_{i,n} \, := \, B_i \cap A_n .
$$
In words, $A_n$ [respectively $B_{i,n}$]
is the set of $\chi \in B$ [respectively $\chi \in B_i$]
which belong to exactly $n$ of the sets $B_j$.
The $A_n$'s are clearly disjoint, and so are the $B_{i, n}$'s for fixed $i \ge 1$,
so that we have disjoint unions
$$
\displaystyle
B \, = \, \bigsqcup_{n \in \overline{\N^*}} A_n
\leqno{(1)}
$$
and
$$
\displaystyle
B_i \, = \, \bigsqcup_{n \in \overline{\N^*}} B_{i,n}
\hskip.5cm \text{for all} \hskip.2cm
i \ge 1 .
\leqno{(2)}
$$
(Note that one may have $A_n = \emptyset$ for several values of $n$:
for example, if all $B_i$ coincide, then all but one of the $A_n$~'s are empty;
of course, the non-empty $A_n$ are pairwise disjoint.
Similarly, one may have $B_{i,n} = \emptyset$.)
If $n \ne n'$, then $B_{i,n}$ and $B_{i,n'}$ are disjoint;
however, if $i \ne i'$, the sets $B_{i,n}$ and $B_{i',n}$ need not be disjoint.
\par

Let $i \ge 1$ and $n \in \overline{\N^*}$.
For $\ell \in \N$ such that $0 \le \ell \le n$ if $n < \infty$
(and no condition on $\ell$ if $n = \infty$), let
$$
B_{i,n, \ell} \, := \, \Big\{ \chi \in B_{i,n}
\hskip.2cm \Big\vert \hskip.2cm
\sum_{j = 1}^n \Un_{B_j}(\chi) = \ell \Big\}
$$
be the set of these $\chi \in B_{i,n}$ which belong to exactly $\ell$
of $B_1, B_2, \hdots, B_n$. We have disjoint unions
$$
\displaystyle
B_{i,n} \hskip.2cm
\, = \, \bigsqcup_{\ell = 0}^n B_{i,n, \ell}
\leqno{(3)}
$$
and, for every $\ell \in \N$ such that $0 \le \ell \le n$ if $n < \infty$
(and no condition if $n = \infty$), 
$$
\displaystyle
A_n
\, = \, \bigsqcup_{i \ge 1} B_{i, n, \ell} .
\leqno{(4)}
$$
Let us check (4).
By definition, the union $\bigcup_{i \ge 1} B_{i, n, \ell}$ is contained
in $A_n$ for every $\ell \le n$.
To show the reverse inclusion, let $\chi \in A_n$.
Then the set $\{j \ge 1 \mid \chi \in B_j \}$ consists of $n$ indices,
say $i_1 < i_2 < \cdots$, 
a finite sequence with last index $i_n$ if $n < \infty$
and an infinite sequence if $n = \infty$.
Fix an integer $\ell$ with $1 \le \ell \le n$ in case $n < \infty$ and $1 \le \ell$ in case $n = \infty$.
Then $\chi \in B_{i_\ell}$ and hence $\chi \in B_{i_\ell, n, \ell}$;
moreover $\chi \notin B_{i, n, \ell}$ for every $i \ne i_\ell$.
Therefore $A_n = \bigcup_{i \ge 1} B_{i, n, \ell}$;
moreover this is a disjoint union, hence (4) holds.

\vskip.2cm

$\bullet$ {\it Second step.}
For a Borel set $C \subset \widehat G$ such that $\mu(C) = 0$,
we set $\mu \vert C = 0$.
Whenever a representation $\pi_{\mu \vert C}$
for such a $\mu$-null set $C$ appears below,
it stands for the representation on the $0$-dimensional Hilbert space.
\par

Let $I$ be the (possibly finite) sequence
of those extended positive integers $n$
such that $\mu(A_n) > 0$.
For all $n \in I$, define
$$
\mu_n \, := \, \mu \vert A_{n} .
$$
Fix $i \ge 1$ and $n \in I$. Define
$$
\mu_{i,n} \, := \, \mu \vert B_{i,n}.
$$
It follows from (2) and from Proposition~\ref{Prop-RepAbDirectSum}
that the representation $\pi_{\mu \vert B_i}$
is equivalent to $\bigoplus_{n\in I} \pi_{\mu_{i,n}}$. 
Consider now $\ell \in \N$ such that $0 \le \ell \le n$ if $n < \infty$
(and no condition if $n = \infty$); define
$$
\mu_{i,n, \ell} \, := \, \mu\vert {B_{i,n, \ell}}.
$$
It follows from (3) and again from Proposition~\ref{Prop-RepAbDirectSum}
that $\pi_{\mu_{i,n}}$ is equivalent $\bigoplus_{\ell \le n} \pi_{\mu_{i,n, \ell}}$.
Therefore,
$$
\pi_{\mu \vert B_i}
\hskip.5cm \text{is equivalent to} \hskip.5cm 
\bigoplus_{n\in I}\bigoplus_{\ell \le n} \pi_{\mu_{i, n, \ell}},
$$
for every $i \ge 1$.
(When $n = \infty$, the notation $\ell \le n$
should be understood as $\ell \in \N$.)
Since $\pi$ is equivalent to $\bigoplus_{i \ge 1} \pi_{\mu\vert{B_{i}}}$,
it follows that 
$$
\pi
\hskip.5cm \text{is equivalent to} \hskip.5cm
\bigoplus_{i \ge 1}\bigoplus_{n\in I}\bigoplus_{\ell \le n} \pi_{\mu_{i, n, \ell}}
\, = \, 
\bigoplus_{n \in I}\bigoplus_{\ell \le n} \bigoplus_{i \ge 1}\pi_{\mu_{i, n, \ell}} .
$$
By (4) and Proposition~\ref{Prop-RepAbDirectSum},
$\bigoplus_{i \ge 1}\pi_{\mu_{i, n, \ell}}$ is equivalent to
$\pi_{\mu_{n}}$ for \emph{every} $\ell \le n$. 
Therefore, $\bigoplus_{\ell \le n} \bigoplus_{i \ge 1}\pi_{\mu_{i, n, \ell}}$
is equivalent to $\pi_{\mu_n}^{(n)}$.
It follows that
$$
\pi
\hskip.5cm \text{is equivalent to} \hskip.5cm
\bigoplus_{n\in I}\pi_{\mu_n}^{(n)}.
$$
We have thus proved the existence part of the theorem.

\vskip.2cm

%%Retablissement de la version originale du ``Third step"
$\bullet$ {\it Third step.}
We prove the uniqueness property.
By the existence part, we may assume that $\pi$
is the representation $\bigoplus_{n \in I} \pi_{\mu_n}^{(n)}$ on the Hilbert space
$$
\Hi \, = \, \bigoplus_{n \in I} L^2(\widehat G, \mu_n, \Ki_n),
$$
where $I$ is a subset of $ \overline{\N^*}$,
$(\mu_n)_{n \in I}$ 
is a sequence of mutually singular probability measures on $\widehat G$,
and $\Ki_n$ is the Hilbert space of dimension $n$ chosen in \ref{notationKn}.
Set $\Hi_n := L^2(\widehat G, \mu_n, \Ki_n)$ for $n\in I$.
\par

For a subset $J$ of $\overline{\N^*}$, let 
$$
\Hi \, = \, \bigoplus_{m \in J} \Hi_m',
$$
be an orthogonal decomposition of $\Hi$ into $\pi(G)$-invariant closed subspaces $\Hi_m'$
such that the restriction of $\pi$ to $\Hi_m'$ is equivalent to the representation 
$\pi_{\mu_m'}^{(m)}$ on $L^2(\widehat G, \mu'_m, \Ki_m)$,
for a sequence of mutually singular probability measures $(\mu_m')_{m \in J}$
on $\widehat G$.
We claim that $J = I$ and that $\Hi_n' = \Hi_n$ for every $n \in I$.
Once established, this will finish the proof.
\par

For $n \in I$ and $m \in J$, denote by $P_n$ and $P'_m$
the orthogonal projections from $\Hi$ to $\Hi_n$ and $\Hi_m'$.
\par

Fix $n \in I$ and $m \in J$ with $P'_m \circ P_n \ne 0$.
By Corollary~\ref{Cor-Prop-RepAb-Equiv2},
$P'_m \circ P_n$ belongs to the centre of 
the von Neumann algebra $\pi(G)''$.
Therefore the orthogonal complement $\Ki$ of $\ker (P'_m \circ P_n)$
and the closure $\Ki'$ of the image of $P'_m \circ P_n$
are non-zero closed subspaces of $\Hi_n$ and $\Hi_m'$ respectively
which are invariant under $\pi(G)$,
and the representations of $G$ on $\Ki$ and $\Ki'$ are equivalent
(see Proposition~\ref{Prop-EquSubRep},
and compare with the proof of Proposition \ref{Prop-CanRepAb-Subrepr}). 
However, by Corollary~\ref{Cor-Prop-RepAb-Equiv1}, 
$\Ki$ and $\Ki'$ are the range of
projections $E_{\mu_n}(B)$ and $E_{\mu_m'}(B')$
for some Borel subsets $B$ and $B'$ of $\widehat G$. 
We have $\mu_n(B)>0$ and $\mu_m'(B')>0$ since $\Ki$ and $\Ki'$ are non-zero.
This shows that the representations
of $G$ on $\Ki$ and $\Ki'$ defined by $\pi$ 
are equivalent to the representations
$\pi_{\mu_n\vert B}^{(n)}$ and $\pi_{\mu_m' \vert {B'}}^{(m)}$ respectively. 
It follows from Proposition~\ref{Prop-RepAb-Subrep}
that $\mu_n\vert B$ and $\mu_m' \vert {B'}$
are equivalent and from Proposition~\ref{Prop-RepAb-Equiv} that $n = m$.
\par
 
So, we have proved that 
$$
P'_m \circ P_n \, = \, P_n \circ P'_m \, = \, 0
\hskip.5cm \text{whenever} \hskip.2cm
n \ne m.
$$
Since $\mathrm{Id}_\Hi = \sum_{n \in I} P_n = \sum_{m \in J } P'_m$,
it follows that $I = J$
and that $\Hi_n = \Hi_n'$ for every $n \in I$.
\end{proof}

\begin{rem}
% 2.B.9
\label{Rem-CanonicalDec-Type I}
(1)
Theorem~\ref{Thm-DecRepAbelianGroups} gives a complete set of invariants
for equivalence classes of representations of $G$ in separable Hilbert spaces:
a subset $I$ of $\overline{\N^*}$ of multiplicities
and a sequence $([\mu_n])_{n \in I}$ of equivalence classes
of mutually singular probability measures on $\widehat G$.

\vskip.2cm

(2)
In the special case $G = \R$, Theorem~\ref{Thm-DecRepAbelianGroups}
is equivalent to the classical Hahn--Hellinger theorem
which gives a complete set of invariants
for the unitary equivalence classes of selfadjoint operators on separable Hilbert spaces.
See \cite{Halm--51};
see also \cite{Sund--97} for a proof using a C*-algebra approach.

\vskip.2cm

(3)
Anticipating on a notion to be introduced in Chapter~\ref{ChapterTypeI},
we observe that a LCA group is a group of type I (see Theorem \ref{explesTypeI}).
Theorem~\ref{Thm-DecRepAbelianGroups} is a special case
of a more general decomposition result for representations of such groups
(Theorem~\ref{thmDirectIntIrreps+}).
\end{rem}

The case for which the set $I$ of extended positive integers
in Theorem \ref{Thm-DecRepAbelianGroups} consists of a single number
will be of special interest for us (see Proposition \ref{Prop-RestNormalSub-Bis}).

\begin{defn}
% 2.B.10
\label{Def-HomogenousRep}
A representation $\pi$ of a second-countable LCA group $G$
in a separable Hilbert space $\Hi$
is said to be \textbf{homogeneous}, or more precisely \textbf{$n$-homogeneous}, 
if the set $I$ of extended positive integers attached to $\pi$ by Theorem
\ref{Thm-DecRepAbelianGroups}~\ref{iDEThm-DecRepAbelianGroups}
consists of a single number $n \in \overline{\N^*}$;
equivalently: if $\pi$ is equivalent to the representation $\pi_\mu^{(n)}$,
for some probability measure $\mu$ on $\widehat G$
and some extended positive integer $n \in \overline{\N^*}$.
\index{Representation! homogeneous}
\index{Homogeneous representation}
\end{defn}

\begin{rem}
% 2.B.11
\label{Rem-HomogenousRep}
A representation $\pi$ of $G$ is homogeneous if and only if 
its commutant $\pi(G)'$ is a homogenous von Neumann algebra,
in the sense of D\'efinition~1 in \cite[Chap.~III, \S~3]{Dixm--vN}.
\end{rem}

\section
{The SNAG Theorem}
% Section 2.C
\label{SNAG}

Let $G$ be a second-countable locally compact abelian group.
In Sections \ref{Section-CanRepAbGr} and \ref{Section-CanDecRepAbGr},
we described the representations of $G$
in terms of probability measures on $\widehat G$ and multiplicities.
As with the spectral theorem for a selfadjoint operator on a Hilbert space
(see Remark~\ref{Rem-RepAbCyclic}),
it is convenient to give another description of representations of $G$,
in terms of projection-valued measures;
for these, see the reminder in \ref{AppProjValMeas}.

\begin{exe}[\textbf{Projection-valued measures on $\widehat G$}]
% 2.C.1
\label{Ex-StandardProjectionValuedMeasure}
Let $\mathcal B (\widehat G)$ the $\sigma$-algebra
of Borel subsets of the dual group $\widehat G$.
Let $\mu$ be a probability measure on $\widehat G$.

\vskip.2cm

(1)
For every $B \in \mathcal B (\widehat G)$,
define a projection $E_\mu(B)$ on $L^2(\widehat G, \mu)$ by 
$$
(E_\mu(B)f) (\chi) \, = \, \Un_B(\chi) f(\chi)
\hskip.5cm \text{for all} \hskip.2cm
f \in L^2(\widehat G, \mu), \hskip.1cm \chi \in \widehat G.
$$
The map
$$
E_\mu \, \colon \, B \mapsto E_\mu(B)
$$
is easily checked to be a projection-valued measure on $\widehat G$
and $L^2(\widehat G, \mu)$;
see Appendix \ref{AppProjValMeas}.
\par

Let $\pi_\mu$ be the canonical representation
of $G$ on $L^2(\widehat G, \mu)$, as in Construction \ref{defpimupourGab};
recall that, for all $g \in G$ and $f_1, f_2 \in L^2(\widehat G, \mu)$, we have
$$
\langle \pi_\mu(g) f_1 \mid f_2 \rangle \, = \,
\int_{\widehat G} (\pi_\mu(g) f_1) (\chi) \overline{f_2(\chi)} d\mu(\chi) \, = \,
\int_{\widehat G} \chi(g) f_1(\chi) \overline{f_1(\chi)} d\mu(\chi) .
$$
Therefore
$$
\pi_\mu(g) \, = \, \int_{\widehat G} \chi(g) dE_\mu (\chi)
\hskip.5cm \text{for all} \hskip.2cm
g \in G .
$$
The measure $\mu_{f_1, f_2}$ on $\widehat G$
defined by $E_\mu$ and $f_1, f_2 \in L^2(\widehat G, \mu)$
(see Appendix \ref{AppProjValMeas})
is given by
$$
\mu_{f_1,f_2}(B) \, = \, \int_{B} f_1(\chi)\overline{f_2}(\chi) d\mu(\chi)
\hskip.5cm \text{for all} \hskip.2cm
B \in \mathcal B (\widehat G) ,
$$
i.e., by $\mu_{f_1, f_2} = f_1 \overline{f_2} \mu$.
Note that $\mu_{f_1, f_2}$ is a complex measure.

\vskip.2cm

(2)
More generally, let $n \in \overline{\N^*}$
and $\Ki$ a Hilbert space of dimension $n$. 
Let $\pi_\mu^{(n)}$ be the associated representation
of $G$ on $L^2(\widehat G, \mu, \Ki)$,
as in Construction \ref{defpimupourGab++}.
A projection-valued measure
$$
E_\mu^{(n)} \, \colon \, \mathcal B \to \Proj(L^2(\widehat G, \mu, \Ki))
$$
is defined by the same formula,
$$
(E_\mu^{(n)} (B)F) (\chi) \, = \, \Un_B(\chi) F(\chi)
\hskip.5cm \text{for all} \hskip.2cm
B \in \mathcal B, \hskip.1cm F \in L^2(\widehat G, \mu, \Ki), \hskip.1cm
\chi \in \widehat G .
$$
One shows as above that
$$
\pi_\mu^{(n)}(g) \, = \, \int_{\widehat G} \chi(g) dE_\mu^{(n)} (\chi)
\hskip.5cm \text{for all} \hskip.2cm
g \in G.
$$
Observe that, for $B \in \mathcal B (\widehat G)$,
we have $E_\mu^{(n)}(B) = 0$ if and only if $\mu(B) = 0$.
\end{exe}

The following result shows in particular that
there is a bijective correspondence
between representations of the LCA group $G$
and projection-valued measures on $\widehat G$.
The name ``SNAG Theorem'' refers to Stone, Naimark, Ambrose, and Godement
\cite{Ston--30, Naim--43, Ambr--44, Gode--44}.

%For a discussion and a proof of the following theorem, see, for example,
%\cite[Theorem 4.45]{Foll--16} and \cite[Section D.3]{BeHV--08}.

% Dictionary
% Snag = An unexpected or hidden obstacle or drawback.
% ?there's one small snag?

\begin{theorem}[\textbf{SNAG Theorem}]
% 2.C.2
\label{thmSNAG}
Let $G$ be a second-countable locally compact abelian group and $\Hi$ a Hilbert space.
\begin{enumerate}[label=(\arabic*)]
\item\label{iDEthmSNAG}
Let $\pi$ be a representation of $G$ in $\Hi$.
There exists a unique projection-valued measure $E$
on $\widehat G$ and $\Hi$
such that
$$
\pi(g) \, = \, \int_{\widehat G} \chi(g) dE (\chi)
\hskip.5cm \text{for all} \hskip.2cm
g \in G.
$$
For every $\xi, \eta \in \Hi$, we have
$$
\langle \pi(g) \xi \mid \eta \rangle
\, = \, \overline{\mathcal{F}}(\mu_{\xi, \eta}) (g)
\hskip.5cm \text{for all} \hskip.2cm
g \in G,
$$
where $\mu_{\xi, \eta}$ is the complex Borel measure on $\widehat G$
defined by $E$ and $\xi, \eta$.
\item\label{iiDEthmSNAG}
Let $G$, $(\pi, \Hi)$ and $E$ be as in \ref{iDEthmSNAG}.
The commutant $\pi(G)'$ of $\pi(G)$ coincides
with the commutant of the set $\{E(B) \mid B \in \mathcal B (\widehat G)\}$.
In particular, $E(B)$ belongs to the von Neumann algebra $\pi(G)''$ 
and the range of $E(B)$ is $\pi(G)$-invariant,
for every $B \in \mathcal B (\widehat G)$.
\item\label{iiiDEthmSNAG}
If $E$ is a projection-valued measure on $\widehat G$ and $\Hi$,
the formula
$$
\pi_E(g) \, := \, \int_{\widehat G} \chi(g) dE (\chi)
\hskip.5cm \text{for all} \hskip.2cm
g \in G
$$
defines a representation of $G$ in $\Hi$.
\end{enumerate}
\end{theorem}

\begin{proof}
\ref{iDEthmSNAG}
We first prove the uniqueness of $E$.
Let $E, E' \, \colon \mathcal B (\widehat G) \to \Proj(\Hi)$ 
be two projection-valued measures such that 
$$
\pi(g) 
\, = \, \int_{\widehat G} \chi(g) dE (\chi)
\, = \, \int_{\widehat G} \chi(g) dE' (\chi)
\hskip.5cm \text{for all} \hskip.2cm
g \in G.
$$
Let $\xi, \eta \in \Hi$.
Denote by $\mu_{\xi, \eta}$ and $\mu'_{\xi, \eta}$
the associated bounded measures on 
$\widehat G$ defined by $E$ and $E'$, respectively.
For every $\xi, \eta \in \Hi$, we have
$$
\begin{aligned}
\langle \pi(g) \xi \mid \eta\rangle
\, &= \,
\int_{\widehat G} \widehat g (\chi)d\mu_{\xi, \eta}(\chi)
\, = \,
\overline{\mathcal{F}}(\mu_{\xi, \eta}) (g)
\\
\, &= \,
\int_{\widehat G} \widehat g (\chi)d\mu'_{\xi, \eta}(\chi)
\, = \,
\overline{\mathcal{F}}(\mu'_{\xi, \eta}) (g)
\end{aligned}
$$
for all $g \in G$.
By injectivity of the Fourier--Stieltjes transform,
this implies $\mu_{\xi, \eta} = \mu'_{\xi, \eta}$.
Since this holds for all $\xi, \eta \in \Hi$, it follows that $E = E'$.

\vskip.2cm

Next, we show the existence of a projection-valued measure
with the desired property. 
\par

By Corollary~\ref{Cor-RepAbCyclic} there exists a family $(\mu_i)_{i \in I}$
of probability measures on $\mathcal B (\widehat G)$
such that $\pi$ is equivalent to the direct sum $\bigoplus_{i \in I}\pi_{\mu_i}$
(for the case of $\Hi$ not separable,
see the note at the end of the proof of Corollary \ref{Cor-RepAbCyclic}).
So, it suffices to prove the existence of a projection-valued measure
$E \, \colon \mathcal B (\widehat G) \to \Proj(\Hi)$
with the desired property in the case of
$$
\Hi \, = \, \bigoplus_{i \in I} L^2(\widehat G, \mu_i)
\hskip.5cm \text{and} \hskip.5cm
\pi \, = \, \bigoplus_{i \in I}\pi_{\mu_i}.
$$
\par
 
For every $i \in I$, let 
$$
E_{\mu_i} \, \colon \, \mathcal B (\widehat G) \to \Proj (L^2(\widehat G, \mu_i))
$$
be the projection-valued measure
as in Example~\ref{Ex-StandardProjectionValuedMeasure}.
As noted there, we have
$$
\pi_{\mu_i}(g) \, = \, \int_{\widehat G} \chi(g) dE_{\mu_i} (\chi)
\hskip.5cm \text{for all} \hskip.2cm
g \in G.
$$
Let $E \, \colon \mathcal B (\widehat G) \to \Proj (\Hi)$
be the direct sum of the $E_{\mu_i}$~'s:
$$
E(B) \, = \, \bigoplus_{i \in I} E_{\mu_i}(B)
\hskip.5cm \text{for all} \hskip.2cm
B \in \mathcal B (\widehat G).
$$
Then $E$ is a projection-valued measure and 
$$
\pi(g) \, = \, 
\bigoplus_{i \in I} \int_{\widehat G} \chi(g) dE_{\mu_i} (\chi)
\, = \,
\int_{\widehat G} \chi(g) d(\bigoplus_{i \in I}E_{\mu_i}) (\chi)
\, = \,
\int_{\widehat G} \chi(g) dE (\chi),
$$
for all $g \in G$. 
 
\vskip.2cm

\ref{iiDEthmSNAG}
Let $T \in \Li (\Hi)$.
For $g \in G$ and $\xi, \eta \in \Hi$, we have, by \ref{iDEthmSNAG},
$$
\langle \pi (g) T\xi \mid \eta \rangle
\, = \, \overline{\mathcal{F}}(\mu_{T\xi, \eta}) (g)
$$
and 
$$
\langle T \pi(g)\xi \mid \eta\rangle
\, = \, \langle \pi(g)\xi \mid T^*\eta\rangle
\, = \, \overline{\mathcal{F}}(\mu_{\xi, T^*\eta}) (g) .
$$
Therefore $T \in \pi(G)'$ if and only if
$\overline{\mathcal{F}}(\mu_{T\xi, \eta}) (g) =
\overline{\mathcal{F}}(\mu_{\xi, T^*\eta}) (g)$
for every $g \in G$ and $\xi, \eta \in \Hi$.
By injectivity of the Fourier--Stieltjes transform, 
we have $T \in \pi(G)'$
if and only if $\mu_{T\xi, \eta} = \mu_{\xi, T^*\eta}$
for every $\xi, \eta \in \Hi$,
hence if and only if
$$
\langle E(B) T\xi \mid \eta\rangle
\, = \,
\langle E(B) \xi \mid T^*\eta\rangle
\, = \,
\langle TE(B) \xi \mid \eta\rangle
$$
for all $B \in \mathcal B (\widehat G)$ and $\xi, \eta \in \Hi$.
Therefore, $T \in \pi(G)'$ if and only $TE(B) = E(B)T$
for every $B \in \mathcal B (\widehat G)$.

\vskip.2cm

\ref{iiiDEthmSNAG}
Let $E \, \colon \mathcal B (\widehat G) \to \Proj(\Hi)$ be a projection-valued measure.
The map 
$$
f \, \mapsto \, \int_{\widehat G} f(\chi) dE (\chi)
$$
is a $*$-homomorphism from ${\rm Bor}^b(\widehat G)$ to $\Li (\Hi)$;
for the $*$-algebra ${\rm Bor}^b(\widehat G)$ of bounded Borel functions
from $\widehat G$ to $\C$
and for this $*$-homomorphism, see Construction \ref{pvmfc}.
In particular, 
$$
\pi_E(g) \, := \, \int_{\widehat G} \widehat g (\chi) dE (\chi) = \int_{\widehat G} \chi(g) dE (\chi)
$$
defines a unitary operator on $\Hi$ for every $g \in G$ and 
$$
\pi_E(g_1)\pi_E(g_2)
\, = \, \int_{\widehat G} \chi(g_1) \chi(g_2) dE (\chi)
\, = \, \int_{\widehat G} \chi(g_1g_1) dE (\chi)
\, = \, \pi_E(g_1g_2)
$$
for all $g_1, g_2 \in G$.
So, $g \mapsto \pi_E(g)$ is a homomorphism from $G$ to the unitary group of $\Hi$. 
For $\xi, \eta \in \Hi$, we have (as in \ref{iDEthmSNAG})
$$
\langle \pi_E(g) \xi \mid \eta \rangle
\, = \, \overline{\mathcal{F}}(\mu_{\xi, \eta}) (g)
\hskip.5cm \text{for all} \hskip.2cm
g \in G.
$$
Since $\overline{\mathcal{F}}(\mu_{\xi, \eta})$ is a continuous on $\widehat G$,
it follows that $g \mapsto \pi_E(g)$ is a representation of~$G$.
\end{proof}

Next, we show how the projection-valued measures
associated to a representation of $G$ and to its conjugate by an automorphism
are related.
\par

Recall from Remark~\ref{Rem-AutoRep}
that $\Aut(G)$ acts on the right by automorphisms of~$\widehat G$:
$$
(\chi, \theta) \mapsto \chi^\theta
\hskip.5cm \text{for} \hskip.2cm
\chi \in \widehat G, \hskip.1cm \theta \in \Aut(G),
$$
with $\chi^\theta = \chi \circ \theta$.

\begin{cor}
% 2.C.3
\label{Cor-SNAG-Auto}
Let $(\pi, \Hi)$ a representation of $G$, and $\theta \in \Aut(G)$.
Let $E \, \colon \mathcal B (\widehat G) \to \Proj(\Hi)$
be the projection-valued measure associated to $\pi$.
\par

Then the projection-valued measure associated to the conjugate representation $\pi^\theta$
is the map $E^\theta \, \colon \mathcal B (\widehat G) \to \Proj(\Hi)$ defined by 
$$
E^{\theta}(B) \, = \, E(B^{\theta^{-1}})
\hskip.5cm \text{for} \hskip.2cm
B \in \mathcal B (\widehat G),
$$
where $B^{\theta^{-1}} = \{\chi^{\theta^{-1}}\mid \chi \in B\}$.
\end{cor}

\begin{proof}
One checks immediately that $E^{\theta}$
is indeed a projection-valued measure on $\widehat G$ and $\Hi$.
\par

Let $\xi, \eta \in \Hi$;
let $\mu_{\xi, \eta}$ and $\nu_{\xi, \eta}$ be the associated bounded measures
on $\widehat G$ defined by $E$ and $E^\theta$.
Denote by $\theta_*(\mu_{\xi, \eta})$ the measure on $\widehat G$
which is the image of $\mu_{\xi, \eta}$ under the map 
$$
\widehat G \to \widehat G, \hskip.2cm \chi \mapsto \chi^\theta.
$$
\par

For every $B \in \mathcal B (\widehat G)$, we have
$$
\nu_{\xi, \eta}(B)\, = \, \langle E^{\theta}(B) \xi \mid \eta \rangle
\, = \, \langle E(B^{\theta^{-1}}) \xi \mid \eta \rangle
\, = \, \mu_{\xi, \eta} (B^{\theta^{-1}})
\, = \theta_*(\mu_{\xi, \eta}) (B).
$$
Therefore, we have
$$
\nu_{\xi, \eta}= \theta_*(\mu_{\xi, \eta}).
$$
It follows that, for every $g \in G$, we have
$$
\begin{aligned}
\langle \pi^\theta (g) \xi \mid \eta\rangle
\, &= \, \langle \pi (\theta(g)) \xi \mid \eta\rangle
\, = \, \int_{\widehat G} \widehat{\theta(g)} (\chi)d\mu_{\xi, \eta}(\chi)
\\
\, &= \, \int_{\widehat G} \widehat g (\chi^\theta)d\mu_{\xi, \eta}(\chi)
\, = \, \int_{\widehat G} \widehat g (\chi) d\theta_*(\mu_{\xi, \eta}) (\chi)
\, = \, \int_{\widehat G} \widehat g (\chi) d\nu_{\xi, \eta}(\chi).
\end{aligned}
$$
Therefore, we have
$$
\pi^\theta (g) = \int_{\widehat G} \chi(g) dE^\theta (\chi)
\hskip.5cm \text{for all} \hskip.2cm
g \in G
$$ 
and the uniqueness part of the SNAG Theorem shows that 
$E^\theta$ is the projection-valued measure associated to $\pi^\theta$.
\end{proof}

\section[Containment and weak containment]
{Containment and weak containment
% \\
in terms of projection-valued measures}
% Section 2.D
\label{SecCwcpvm}

Let $G$ be a second-countable LCA group,
$\Hi$ a Hilbert space,
$\pi$ a representation of $G$ in $\Hi$,
and $E$ the corresponding projection-valued measure on $\widehat G$ and $\Hi$.
Containment and weak containment of unitary characters of $G$ in~$\pi$
can be expressed neatly in terms of the atoms and of the support of $E$
(atoms and supports are defined in Appendix \ref{AppProjValMeas}).

\begin{prop}
% 2.D.1
\label{Prop-ContAbelian}
Let $\Hi$ a Hilbert space,
$\pi$ a representation of $G$ in $\Hi$, 
and $E$ the projection-valued measure associated to $\pi$
as in the SNAG theorem \ref{thmSNAG}.
Let $\chi_0 \in \widehat G$. Set
$$
\Hi^{\chi_0} \, = \, \left\{ \xi \in \Hi \mid \pi(g) \xi = \chi_0(g) \xi 
\hskip.2cm \text{for all} \hskip.2cm 
g \in G \right\}.
$$
(a)
The following properties are equivalent:
\begin{enumerate}[label=(\roman*)]
\item\label{iDEProp-ContAbelian} 
$\chi_0$ is contained in $\pi$, i.e., $\Hi^{\chi_0} \ne \{0\}$;
\item\label{iiDEProp-ContAbelian}
$\chi_0$ is an atom of $E$;
\end{enumerate}
when they hold, $E(\{\chi_0\})$ is the orthogonal projection on $\Hi^{\chi_0}$.

\vskip.2cm

\noindent
(b)
The following properties are equivalent:
\begin{enumerate}[label=(\roman*)]
\addtocounter{enumi}{2}
\item\label{iiiDEProp-ContAbelian} 
$\chi_0$ is weakly contained in $\pi$;
\item\label{ivDEProp-ContAbelian}
$\chi_0$ belongs to the support of $E$.
\end{enumerate}
\end{prop}

\begin{proof}
(a)
Assume that $\chi_0$ is contained in $\pi$.
% ; then the space $\Hi^{\chi_0}$ is non-zero. 
Let $\xi$ be a unit vector in $\Hi^{\chi_0}$ and $\mu_\xi$
the probability measure on $\widehat G$
associated to $E$ and $\xi$. For every $g \in G$, we have
$$
\chi_0(g) \, = \, \langle \pi(g) \xi \mid \xi \rangle
\, = \, \int_{\widehat G} \chi (g )d\mu_\xi(\chi).
$$
Since the points in the unit circle are extremal points of the unit disc, 
it follows that $\mu_\xi$ is the Dirac measure $\delta_{\chi_0}$; so
$$
\langle E(\{\chi_0\})\xi \mid \xi \rangle \, = \, \mu_\xi (\{\chi_0\}) \, = \, 1.
$$
Therefore $E(\{\chi_0\})\xi = \xi$, that is, $\xi$ is in the range of $E(\{\chi_0\})$,
and in particular $E(\{\chi_0\}) \ne \{0\}$.

\vskip.2cm

Assume now that $E(\{\chi_0\}) \ne 0$. 
Let $\xi$ be a unit vector in the range of the projection $E(\{\chi_0\})$. 
Since $\langle E(\{\chi_0\}) \xi \mid \xi \rangle = 1$,
the probability measure on $\widehat G$ associated to $E$ and $\xi$
is the Dirac measure $\delta_{\chi_0}$.
Therefore, for every $g \in G$, we have
$$
\langle \pi(g)\xi \mid \xi \rangle \, = \, 
\int_{\widehat G} \chi (g) d\delta_{\chi_0}(\chi) \, = \,
\chi_0(g).
$$
By the equality case of the Cauchy--Schwarz inequality, we have
$$
\pi(g)\xi \, = \, \chi_0(g) \xi
\hskip.5cm \text{for all} \hskip.2cm
g \in G,
$$
that is, $\xi \in \Hi^{\chi_0}$, and so $\Hi^{\chi_0} \ne \{0\}$.
\par

Observe that we have shown in the course of the proof 
that $\Hi^{\chi_0}$ coincides with the range of $E(\{\chi_0\})$. 

\vskip.2cm

(b) 
Assume that $\chi_0$ is weakly contained in $\pi$. 
Then there exists a sequence $(\xi_i)_i$ of unit vectors in $\Hi$ such that 
$$
\lim_i \langle \pi(g) \xi_i \mid \xi_i \rangle \, = \, \chi_0(g)
$$
uniformly on compact subsets of $G$ (see \cite[F.1.4]{BeHV--08}).
Denote by $\mu_i$ the probability measure on $\widehat G$ 
associated to $E$ and $\xi_i$
(with the notation of Appendix \ref{AppProjValMeas}, $\mu_i$ should be $\mu_{\xi_i}$).
We have therefore
$$
\lim_i {\overline{\mathcal{F}}}(\mu_i) (g) \, = \,
\lim_i \int_{\widehat G} \chi(g) d\mu_i(\chi) \, = \,
\chi_0(g)
\leqno{(*)}
$$
uniformly on compact subsets of $G$. 

%\marginpar{Passage pas clair}
Let $U$ be an arbitrary compact neighbourhood of $\chi_0$.
We claim that ${\lim_i \mu_i(U) = 1}$;
once proved, it will follow that
$$
\mu_i(U) \, = \, \langle E(U)\xi_i \mid \xi_i \rangle \, \ne \, 0
$$
for $i$ large enough and therefore $E(U) \ne 0$,
showing that $\chi_0$ is in the support of $E$.
\par

Denote by $\mu_G$ a Haar measure on $G$.
As is well-known (see Theorem (31.34) in \cite{HeRo--70}),
there exists a function $\varphi \, \colon G \to \mathopen[ 0,1 \mathclose]$ in $L^1(G, \mu_G)$
such that
% there exists a function $\varphi \in L^1(G, \mu_G)$ such that
% $\overline{\mathcal F}(\varphi) \le \Un_U$ and 
$\overline{\mathcal F}(\varphi) (\chi_0) = 1$
and $\overline{\mathcal F}(\varphi) (\chi) = 0$ for all $\chi \in G \smallsetminus U$.
\par

For every probability measure $\mu$ on $\widehat G$, we have,
by Fubini's theorem
(which can be applied, because $\widehat G$ is second-countable,
and therefore $\mu$ and $\mu_G$ are $\sigma$-finite measures)
$$
\begin{aligned}
\int_{\widehat G} \overline{\mathcal{F}}(\varphi) (\chi) d\mu(\chi)
\, &= \,
\int_{\widehat G} \int_{G}\chi(g) \varphi(g) d\mu_G(g) d\mu(\chi)
\\
\, &= \,
 \int_{G}\varphi(g) \left(\int_{\widehat G} \chi(g) d\mu(\chi)\right) d\mu_G(g)
 \\
\, &= \,
\int_{G}\varphi(g) \overline{\mathcal{F}}(\mu) (g) d\mu_G(g).
\end{aligned}
$$
As 
$$
\vert \varphi(g)\overline{\mathcal{F}}(\mu_i) (g) \vert \, \le \,
\vert \varphi(g) \vert
\hskip.5cm \text{for all} \hskip.2cm
g \in G,
$$
it follows from Lebesgue dominated convergence theorem that 
$$
\begin{aligned}
\lim_i \int_{\widehat G} \overline{\mathcal F}(\varphi) (\chi) d\mu_i(\chi)
\, &= \, \int_{\widehat G} \lim_i \varphi(g) \overline{\mathcal F} (\mu_i) (g) d\mu_G(g)
\\
\, &= \, \int_{G} \varphi(g) \chi_0(g) d\mu_G(g)
\, = \,
\overline{\mathcal{F}}(\varphi) (\chi_0)
\, = \,
1.
\end{aligned}
$$
Since $\overline{\mathcal{F}}(\varphi) \le \Un_U$ and hence
$$
\int_{\widehat G} \overline{\mathcal{F}}(\varphi) (\chi) d\mu_i(\chi)
\, \le \, 
\int_{\widehat G} \Un_U(\chi) d\mu_i(\chi)
\, = \,
\mu_i(U) \, \le \, 1,
$$
it follows that $\lim_i \mu_i(U) = 1$ and the claim is proved.

\vskip.2cm

Conversely, assume that $E(U) \ne 0$ for every neighbourhood $U$ of $\chi_0$.
Let $(U_i)_i$ be a sequence of decreasing open neighbourhoods of $\chi_0$ such that 
$$
\bigcap_i U_i \, = \, \{\chi_0\}.
$$
For every $i$, let $\xi_i$ be a unit vector in the range of the projection $E(U_i)$.
Since 
$$
E(B)\xi_i \, = \, E(B)E(U_i)\xi_i \, = \, E(B \cap U_i)\xi_i \, = \, 0
$$
for every $B \in \mathcal B (\widehat G)$ such that $B \cap U_i = \emptyset$,
the support of the corresponding probability measure 
$\mu_i := \mu_{\xi_i}$ on $\widehat G$ 
is contained in $U_i$. 
\par

Let $K$ be a compact subset of $G$ and $\varepsilon > 0$.
Then the family $\{\widehat g \mid g \in K \} \subset C(\widehat G)$ is equicontinuous
and we have therefore, for $i$ large enough, 
$$
\vert \widehat g (\chi) - \widehat g (\chi_0) \vert \, \le \, \varepsilon
\hskip.5cm \text{for all} \hskip.2cm
\chi \in U_i
\hskip.2cm \text{and} \hskip.2cm
g \in K .
$$
Since the support $\mu_i$
is contained in $U_i$, it follows that, for $i$ large enough, 
$$
\left\vert \int_{\widehat G} \chi(g) d\mu_i(\chi) - \chi_0(g) \right\vert
\, = \,
\left\vert \int_{\widehat G} \widehat g (\chi) d\mu_i(\chi) -\widehat g (\chi_0) \right\vert
\, \le \,
\int_{U_i} \vert \widehat g (\chi)-\widehat g (\chi_0) \vert d\mu_i(\chi)
\, \le \,
\varepsilon
$$
for all $g \in K$. 
This shows that 
$$
\lim_i \langle \pi(g) \xi_i \mid \xi_i \rangle \, = \,
\lim_i \int_{\widehat G} \chi(g) d\mu_i(\chi) \, = \, \chi_0(g)
$$
uniformly on $K$.
Therefore, $\chi_0$ is weakly contained in $\pi$.
\end{proof}

For later reference, 
% proof of Theorems 4.A.1 and 10.D.1
% c'est \`a dire \ref{Prop-RepGrDiedral} et \ref{Theo-PrimIdealBS}.
we note the following consequence of 
Proposition~\ref{Prop-ContAbelian}.

\begin{cor} 
% 2.D.2
\label{Cor-ContAbelian}
Let $G$, $\Hi$, $\pi$, and $E$ be as in Proposition \ref{Prop-ContAbelian}.
Assume that the support of $E$ is a finite subset 
$\{\chi_1, \hdots, \chi_n\}$ of $\widehat G$. 
For $i \in \{1, \hdots, n\}$, define $\Hi^{\chi_i}$ as in \ref{Prop-ContAbelian}.
\par

Then $\Hi^{\chi_i} \ne \{0\}$ for every $i \in\{1, \hdots, n\}$,
and $\Hi$ decomposes as the direct sum 
$$
\Hi \, = \, \Hi^{\chi_1} \oplus \cdots \oplus \Hi^{\chi_n}.
$$
In particular, 
$\pi$ is equivalent to $k_1 \chi_1 \oplus \cdots \oplus k_n\chi_n$, 
for multiplicities $k_1, \hdots, k_n \in \{1, 2, 3, \hdots, \infty \}$.
\end{cor}

\begin{proof}
The complement $U$ of $\{\chi_1, \hdots, \chi_n\}$ in $\widehat G$ 
is an open set with $E(U) = 0$. 
For $i \in \{1, \hdots, n\}$, we have $E(\{ \chi_i \}) \ne 0$, 
since $\chi_i$ is in the support of $E$.
Moreover, by the proof of Proposition~\ref{Prop-ContAbelian},
$E(\{ \chi_i \})$ is the orthogonal projection on $\Hi^{\chi_i}$.
Since 
$$
\mathrm{Id}_{\Hi} \, = \, 
E(U)+ E(\{\chi_1, \hdots, \chi_n\}) \, = \, 
E(\{\chi_1 \}) \oplus \cdots \oplus E(\{\chi_n\}),
$$
the claim follows.
\end{proof}

\begin{rem}
% 2.D.3
\label{Rem-DecRepAbelianGroups}
Let $\pi$ be a representation of $G$.
By the existence and the uniqueness part 
of Theorem~\ref{Thm-DecRepAbelianGroups},
there exists a \emph{unique} sequence $([\mu_k])_k$
of equivalence classes of probability measures $\mu_k$ on $\widehat G$
attached to $\pi$.
\par

For every $k$, denote by $A_k$ the set of atoms
and by $S_k$ the support of $\mu_k$.
Obviously, $A_k$ and $S_k$ only depend on the class of $\mu_k$.
Let $E$ be the projection-valued measure associated to $\pi$. 
As is easily checked, the set of atoms of $E$ is $\bigcup_k A_k$
and the support of $E$ is the closure of $\bigcup_k S_k$.
\end{rem}

\section
[Canonical decompositions]
{Canonical decomposition of projection-valued measures}
% Section 2.E
\label{SecCandecpm}

Let $G$ be a second-countable locally compact abelian group.
Theorem~\ref{Thm-DecRepAbelianGroups} shows
that every representation of $G$ on a separable Hilbert space
admits a canonical decomposition as a direct sum of representations
of the form $\pi_{\mu_k}^{(n_k)}$.
Using the SNAG Theorem, we will reformulate this result
in terms of projection-valued measures:
every projection-valued measure on $\widehat G$ is equivalent, in a canonical way,
to a direct sum of projection-valued measures 
associated to representations of the form $\pi_{\mu_k}^{(n_k)}$.
\par

We define equivalence of projection-valued measures in the obvious way:
two projection-valued measures 
$$
E \, \colon \, \mathcal B (\widehat G) \to \Proj(\Hi)
\hskip.5cm \text{and} \hskip.5cm
E' \, \colon \, \mathcal B (\widehat G) \to \Proj(\Hi')
$$
are \textbf{equivalent} if there exists a Hilbert space isomorphism
$T \, \colon \Hi \to \Hi'$ such that 
$$
TE(B) \, = \, E'(B) T
\hskip.5cm \text{for all} \hskip.2cm
B \in \mathcal B (\widehat G).
$$
By the SNAG Theorem \ref{thmSNAG},
this is the case if and only if the associated representations $\pi_E$ and $\pi_{E'}$ are equivalent.
\index{Projection-valued measure! equivalent}
\index{Equivalent! projection-valued measures}

\begin{theorem}
[\textbf{Canonical decomposition of projection-valued measures of abelian groups}]
% 2.E.1
\label{Theo-DecPVM}
Let $G$ be a second-countable locally compact abelian group. Let
$$
E \, \colon \, \mathcal B (\widehat G) \to \Proj(\Hi)
$$
a projection-valued measure on $\widehat G$ and a separable Hilbert space $\Hi$.
Denote by $\M \subset \Li (\Hi)$ the von Neumann algebra 
$\{ E(B) \mid B \in \widehat G \}''$ generated by the image of $E$.
\par

There exist a set $I \subset \overline{\N^*}$ of extended positive integers,
a sequence $(P_n)_{n \in I}$ of mutually orthogonal projections in the centre of $\M$,
and a sequence $(\mu_n)_{n \in I}$
of mutually singular probability measures on $\widehat G$
with the following properties:
\begin{enumerate}[label=(\arabic*)]
\item\label{iDETheo-DecPVM}
$E = \bigoplus_{n \in I} E_n$,
where $E_n$ is the projection-valued measure $E_n = P_n \circ E$.
\item\label{iiDETheo-DecPVM}
For every $n \in I$, the projection-valued measure $E_n$ is equivalent to the 
projection-valued measure $E_{\mu_n}^{(n)}$
as in Example~\ref{Ex-StandardProjectionValuedMeasure}.
\end{enumerate}
\par

Moreover, this decomposition is unique in the following sense:
assume that there exists a subset $J \subset \overline{\N^*}$,
a sequence $(Q_m)_{m \in J}$ of mutually orthogonal projections in the commutant $\M'$ of $\M$,
and a sequence $(\mu_m')_{m \in J}$ 
of mutually singular probability measures on $\widehat G$,
with properties analogous to \ref{iDETheo-DecPVM} and \ref{iiDETheo-DecPVM}.
Then $J = I$, and $Q_n = P_n$, 
and $\mu_n$ and $\mu_n'$
are equivalent for every $n \in I$.
\end{theorem}

\begin{proof}
Let $\pi$ be the representation of $G$ in $\Hi$ defined by $E$.
In view of Theorem~\ref{Thm-DecRepAbelianGroups}, we can assume that 
$$
\pi \, = \, \bigoplus_{n \in I} \pi_{\mu_n}^{(n)}
$$
for a subset $I$ of $\overline{\N^*}$
and a sequence $(\mu_n)_{n \in I}$
of mutually singular probability measures on $\widehat G$.
So, we have
$$
\Hi \, = \, \bigoplus_{n \in I} L^2(\widehat G, \mu_n, \Ki_n),
$$
where $\Ki_n$ is a Hilbert space of dimension $n$.
\par

For every $n \in I$, denote by $P_n$ the orthogonal projection
of $\Hi$ onto $L^2(\widehat G, \mu_n, \Ki_n)$.
The mutually orthogonal projections $P_n$ belong to the centre of $\M= \pi(G)''$,
by Corollary~\ref{Cor-Prop-RepAb-Equiv2}.
\par

We claim that the sequence $(P_n)_{n \in I}$ satisfies 
Properties \ref{iDETheo-DecPVM} and \ref{iiDETheo-DecPVM} of the theorem.

\vskip.2cm

\ref{iDETheo-DecPVM}
Since $\sum_{n \in I} P_n = I$, we have
$E = \bigoplus_{n \in I} E_n$ for $E_n = P_n \circ E$.

\vskip.2cm

\ref{iiDETheo-DecPVM}
It is clear that $E_n$ coincides with the projection-valued measure $E_{\mu_n}^{(n)}$
associated to $\pi_{\mu_n}^{(n)}$, for every $n\in I$.

\vskip.2cm

It remains to check the uniqueness assertion.
Let $J \subset \overline{\N^*}$,
 $(Q_m)_{m \in J}$ a sequence of mutually orthogonal projections in the centre of $\M$,
and $(\mu_m')_{m \in I}$ a sequence
of mutually singular probability measures on $\widehat G$,
with properties \ref{iDETheo-DecPVM} and \ref{iiDETheo-DecPVM}.
As $\pi$ is determined by $E$,
it follows from the uniqueness part of Theorem~\ref{Thm-DecRepAbelianGroups},
that $J = I$ and that $\mu_n'$ is equivalent to $\mu_n$ for every $n \in I$.
\par
 
For every $n\in I$, the range of $Q_n$ is $\pi(G)$-invariant, since $Q_n\in \pi(G)'$;
moreover, the subrepresentation of $\pi$ defined by the range of $Q_n$
is equivalent to $\pi_{\mu_n}^{(n)}$, since $\mu_n'$ is equivalent to $\mu_n$.
 Since $\pi_{\mu_n}^{(n)}$ and $\pi_{\mu_m}^{(m)}$
are disjoint (see Proposition~\ref{Prop-RepAb-Subrep}),
we have therefore $Q_m P_n = 0$ for all $m, n \in I$ with $m\ne n$.
As $\sum_{m\in I} Q_m = \sum_{n\in I} P_n = I $, it follows that $Q_n = P_n$ for every $n\in I$.
\end{proof}

\section
{Applying the SNAG Theorem to non-abelian groups}
% Section 2.F
\label{SNAG-NonAbelian}

We will be interested in representations of groups $G$
which are non-abelian but which have a ``large'' abelian normal subgroup $N$
(this holds for the examples of Sections
\ref{SectionInfDiGroup} to \ref{Section-IrrRepLamplighter}).
We will analyze a given representation of $G$
by applying the SNAG theorem
and the canonical decomposition Theorem~\ref{Theo-DecPVM}
to its restriction to~$N$. 
This procedure is the first step of the so-called Mackey machine;
see Remark \ref{MackMach}.
\par

In order to state our next proposition,
we need one more notion.
A representation $(\pi, \Hi)$ of a group $G$ is said to be \textbf{factorial},
or a \textbf{factor representation}, 
if the von Neumann algebra $\pi(G)'' \subset \Li (\Hi)$ 
generated by $\pi(G)$ is a factor, 
that is, if the centre of $\pi(G)''$ consists of
the scalar operators on $\Hi$ only.
We come back to this notion in Chapter~\ref{ChapterTypeI}; see in particular
Definition \ref{Def-FactorialRep} and Proposition \ref{Pro-FactorRepVN}.
\index{Representation! factor}
\index{Factor representation}
\index{von Neumann algebra! $1$@factor}
\index{Factor}
\par

Let $N$ be closed abelian normal subgroup of $G$.
Recall from Remark~\ref{Rem-AutoRep}
the natural right action $\widehat N \curvearrowleft G$,
defined by $\chi^g (n) = \chi(gng^{-1})$ for $\chi \in \widehat N$, $g \in G$, and $n \in N$.
For a subset $B$ of $\widehat N$,
we write $B^{g}$ for the image of $B$ under the action of $g \in G$,
that is, $B^g = \{\chi^g \mid \chi \in B\}$.

\begin{prop}
% 2.F.1
\label{Prop-RestNormalSub}
Let $G$ be a locally compact group, 
$N$ a closed second-countable \emph{abelian} normal subgroup of $G$,
and $(\pi, \Hi)$ a representation of $G$.
Let $E$ be the projection-valued measure on $\widehat N$ associated to 
the restriction $\pi \vert_N$ of $\pi$ to~$N$.
\begin{enumerate}[label=(\arabic*)]
\item\label{iDEProp-RestNormalSub}
For every $g \in G$, we have
$$
\pi(g)E(B)\pi(g)^{-1} \, = \, E(B^{g^{-1}})
\hskip.5cm \text{for all} \hskip.2cm 
B \in \mathcal B (\widehat N).
$$
\item\label{iiDEProp-RestNormalSub}
Assume that $\pi$ is factorial.
Then the support of $E$
coincides with the closure of a $G$-orbit in $\widehat N$.
\end{enumerate}
\end{prop}

\begin{proof} 
\ref{iDEProp-RestNormalSub}
Set $\sigma := \pi \vert_N$. 
For $g \in G$, let
$$
E^g \, \colon \, \mathcal B (\widehat N) \to \Proj (\Hi)
$$
be the projection-valued measure
associated to the conjugate representation $\sigma^g$ of~$\sigma$. 
\par

On the one hand, since 
$$
\sigma^g(n) \, = \, \pi(g) \sigma(n) \pi(g^{-1}) 
\hskip.5cm \text{for all} \hskip.2cm 
n \in N,
$$
we have 
$$
E^g(B) \, = \, \pi(g)E(B) \pi(g^{-1}) 
\hskip.5cm \text{for all} \hskip.2cm
B \in \mathcal B (\widehat N).
$$
On the other hand, by Corollary~\ref{Cor-SNAG-Auto}, we have
$$
E^g(B) \, = \, E(B^{g^{-1}})
\hskip.5cm \text{for all} \hskip.2cm 
B \in \mathcal B (\widehat N).
$$
Claim \ref{iDEProp-RestNormalSub}
follows from the uniqueness part of the SNAG Theorem.

\vskip.2cm

\ref{iiDEProp-RestNormalSub}
The support $S$ of $E$ is a closed subset of $\widehat N$,
which is $G$-invariant, by \ref{iDEProp-RestNormalSub}.
Since $N$ is a second-countable LCA group, $\widehat N$ is also second-countable
(Proposition~\ref{Prop-DualAbelianSecondCountGroup}).
Let $(U_i)_{i \in \N}$ be a countable basis for the topology of $\widehat N$.
Denote by $I \subset \N$ the subset of integers $i$
such that $U_i \cap S \ne \emptyset$.
\par

Let $i \in I$. Then
$$
U_i^G \, := \, \{ \chi^g \mid \chi \in U_i, \hskip.1cm g \in G \}
$$
is a $G$-invariant open subset of $\widehat N$. 
For the projection $P_i = E(U_i^G)$,
we have $P_i \ne 0$, since $U_i^G$ is an open set 
with a non-empty intersection with $S$. 
\par

On the one hand, by the SNAG theorem,
$P_i$ belongs to the von Neumann algebra $\pi(N)''$ generated by $N$; 
in particular, $P_i \in \pi(G)''$.
On the other hand, since $U_i^G$ is $G$-invariant, 
it follows from \ref{iDEProp-RestNormalSub}
that $P_i$ commutes with $\pi(g)$
for every $g \in G$, that is, $P_i \in \pi(G)'$.
Therefore, 
$$
P_i \, \in \, \pi(G)'' \cap \pi(G)',
$$
that is, $P_i$ belongs to the centre of the von Neumann algebra $\pi(G)''$.
Since $\pi(G)''$ is a factor and $P_i \ne 0$,
it follows that $P_i = I$.
We have therefore, for every $i \in I$,
$$
E(\widehat N \smallsetminus U_i^G) \, = \, I - P_i \, = \, 0 .
$$
Therefore, since $I$ is a countable set
and by the properties of a projection-valued measure,
we have
$$
E \Big( \bigcup _{i \in I} (\widehat N \smallsetminus U_i^G) \Big) \, = \, 0.
$$
It follows that 
$E \big( \bigcap _{i \in I} U_i^G \big) = I$ and therefore 
$$
E \Big( S \cap \big( \bigcap _{i \in I} U_i^G \big) \Big)
\, = \, 
E(S) E \Big( \bigcap _{i \in I} U_i^G \Big) \, = \, I ,
$$
because $E(S) = I$.
In particular, we have
$$
S \cap \big( \bigcap _{i \in I} U_i^G \big) \, \ne \, \emptyset .
$$
\par

Let $\chi_0 \in S \cap \big( \bigcap _{i \in I} U_i^G \big)$.
We claim that the $G$-orbit of $\chi_0$ is dense in $S$.
Indeed, let $\chi \in S$ and let $U$ be an open neighbourhood of $\chi$.
Then, $U_i \subset U$ for some $i \in I$.
Since $\chi_0 \in U_i^G$,
the $G$-orbit of $\chi_0$ has a non-empty intersection with 
$U_i$ and hence with $U$. 
\end{proof}

\begin{rem}
% 2.F.2
\label{Rem-RestNormalSub}
The formula of \ref{iDEProp-RestNormalSub} in
Proposition \ref{Prop-RestNormalSub} 
shows that $B \mapsto E(B)$ is a system of imprimitivity for $(\pi, \Hi)$
based on the $G$-space $\widehat N$, in the sense of \cite{Mack--58}.
\end{rem}

We will need later (proof of Theorem~\ref{Theo-AllRepSemiDirect2})
the following refinement of
Proposition~\ref{Prop-RestNormalSub}~\ref{iiDEProp-RestNormalSub}.

\begin{prop}
% 2.F.3
\label{Prop-RestNormalSub-Bis}
Let $G$ be a locally compact group,
$N$ a closed second-countable abelian normal subgroup of $G$,
and $\pi$ a representation of $G$ in a separable Hilbert space $\Hi$.
Let $I$ be the subset of $\overline{\N^*}$
and let $(\mu_n)_{n \in I}$ be a sequence
of probability measures on $\widehat N$
attached to the restriction $\pi \vert_N$ of $\pi$ to $N$,
as in Theorem~\ref{Thm-DecRepAbelianGroups}.
\begin{enumerate}[label=(\arabic*)]
\item\label{iProp-RestNormalSub-Bis}
Every probability measure $\mu_n$ is quasi-invariant under $G$.
\item\label{iiProp-RestNormalSub-Bis}
Assume that $\pi$ is factorial.
Then $\pi \vert_N$ is a homogeneous representation of~$N$,
that is, $\pi \vert_N$ is equivalent to the representation 
$(\pi_\mu^{(n)}, L^2(\widehat N, \mu, \Ki))$
for a probability measure $\mu$ on $\widehat N$ which is quasi-invariant by $G$
and a Hilbert space $\Ki$ of some dimension $n \in \overline{\N^*}$.
Moreover, $\mu$ is ergodic under the $G$-action.
\end{enumerate}
\end{prop}

\begin{proof}
\ref{iProp-RestNormalSub-Bis}
Set $\sigma = \pi \vert_N$ and let $E \, \colon \mathcal B (\widehat N) \to \Proj (\Hi)$
be the projection-valued measure associated to $\sigma$.
Let $(P_n)_{n \in I}$ be the sequence of orthogonal projections in the centre of $\sigma(N)''$
attached to $E$ as in Theorem~\ref{Theo-DecPVM},
so that $E = \bigoplus_{n \in I} P_n \circ E$.
\par

Fix $g \in G$
and denote by $E^g$ the projection-valued measure 
associated to the conjugate representation $\sigma^g$.
\par

On the one hand, we have 
$$
E^g(B) \, = \, \pi(g)E(B)\pi(g^{-1})
\hskip.5cm \text{for all} \hskip.2cm
B \in \mathcal B (\widehat N).
$$
It follows that $I$, $(\mu_n)_{n \in I}$ and $(\pi(g) P_n \pi(g^{-1}))_{n \in I}$
are respectively a subset of $\overline{\N^*}$,
a sequence of mutually singular probability measures on $\widehat N$,
and a sequence of mutually orthogonal projections in the centre of $\sigma^g(N)''$,
satisfying together Properties \ref{iDETheo-DecPVM} and \ref{iiDETheo-DecPVM}
of Theorem~\ref{Theo-DecPVM}
for the projection-valued measure $E^g$. 
\par

On the other hand, by Proposition~\ref{Prop-RestNormalSub}, we have
$$
E^g(B) \, = \, E(B^{g^{-1}})
\hskip.5cm \text{for all} \hskip.2cm
B \in \mathcal B (\widehat N).
$$
In particular, this implies (see Theorem (\ref{thmSNAG}(ii))
that $\sigma^g(N)'' = \sigma(N)''$.
\par

For every $n \in I$, consider the conjugate representation
$\left( \big({\pi_{\mu_n}^{(n)}}\big)^g, L^2(\widehat N, \mu_n, \Ki_n) \right)$;
we have
$$
\left( \big({\pi_{\mu_n}^{(n)}}\big)^g (x) F \right) (\chi)
\, = \, \left( \pi_{\mu_n}^{(n)}(g x g^{-1}) F \right) (\chi)
\, = \, \chi(g x g^{-1}) F(\chi)
\, = \, \chi^{g}(x)F(\chi),
$$
for all $x \in N, \chi \in \widehat N$ and $F \in L^2(\widehat N, \mu_n, \Ki_n)$.
It follows that $\big({\pi_{\mu_n}^{(n)}}\big)^g$
is (equivalent to) the representation of $N$
corresponding to the projection-valued measure $P_n \circ E^g$.
Since $\big({\pi_{\mu_n}^{(n)}}\big)^g$ is equivalent to 
$\pi_{g_*(\mu_n)}^{(n)}$,
it follows that $I$, $(g_*(\mu_n))_{n \in I}$ and $(P_n)_{n \in I}$
are another subset of $\overline{\N^*}$,
a sequence of mutually singular probability measures on $\widehat N$,
and a sequence of mutually orthogonal projections in the centre of $\sigma(N)'' = \sigma^g(N)''$,
satisfying together Properties \ref{iDETheo-DecPVM} and \ref{iiDETheo-DecPVM}
of Theorem \ref{Theo-DecPVM}
for the projection-valued measure $E^g$.
\par

It follows from the uniqueness part of Theorem~\ref{Theo-DecPVM} that:
\begin{itemize}
\setlength\itemsep{0em}
\item
$g_*(\mu_n)$ and $\mu_n$ are equivalent for every $n \in I$; 
\item
$P_n = \pi(g) P_n \pi(g^{-1})$ for every $n \in I$.
\end{itemize}
Since this holds for every $g \in G$,
and since $P_n$ belongs to $\pi(N)''$ and hence to $\pi(G)''$, 
this implies that 
\begin{itemize}
\setlength\itemsep{0em}
\item
every measure $\mu_n$ is quasi-invariant by $G$;
\item every projection $P_n$ belongs belongs to the centre of $\pi(G)''$.
\end{itemize}

\vskip.2cm
 
\ref{iiProp-RestNormalSub-Bis}
Assume now that $\pi$ is factorial.
Since every $P_n$ belongs to the centre of $\pi(G)''$,
and $\sum_{n \in I} P_n = I$,
it follows that the subset $I$ of $\overline{\N^*}$
consists of a single element $n \in \overline{\N^*}$,
that is, $\pi$ is $n$-homogeneous for this $n$.
\par
 
Denote now by $\mu$ the probability measure on $\widehat N$ associated to $n$. 
Let $B \in \mathcal B (\widehat N)$ be $G$-invariant.
The projection $E(B)$ belongs to $\pi(G)''$,
and also to $\pi(G)'$, because
$$
E(B) \, = \, E(B^g) \, = \, \pi(g)E(B)\pi(g^{-1})
\hskip.5cm \text{for all} \hskip.2cm
g \in G.
$$
It follows that $E(B)$ belongs to the centre of $\pi(G)''$,
and we have $E(B) = 0$ or $E(B) = I$.
This implies that $\mu(B) = 0$ or $\mu(B) = 1$,
and shows that $\mu$ is ergodic.
\end{proof}

\begin{cor}
% 2.F.4
\label{Cor-RestNormalSub-Bis}
Let $G$ be a locally compact group,
$N$ a closed second-countable abelian normal subgroup of $G$,
and $\pi$ a \emph{factor} representation of $G$ in a separable Hilbert space $\Hi$.
Let $E$ be the projection-valued measure associated to $\pi \vert_N$.
\begin{enumerate}[label=(\arabic*)]
\item\label{iCor-RestNormalSub-Bis}
The equivalence class of the ergodic $G$-quasi-invariant probability measure 
$\mu$ on $\widehat N$ attached to $\pi \vert_N$
as in Proposition~\ref{Prop-RestNormalSub-Bis} is unique. 
More specifically, the sets $B \in \mathcal B (\widehat N)$ 
such that $\mu(B) = 0$ for any such measure $\mu$ 
are exactly the sets $B \in \mathcal B (\widehat N)$ 
such that $E(B) = 0$.
\item\label{iiCor-RestNormalSub-Bis}
The support (respectively, the set of atoms)
of every probability measure $\mu$ on $\widehat N$ attached to $\pi \vert_N$
as in Proposition~\ref{Prop-RestNormalSub-Bis}
coincides with the support (respectively, the set of atoms) of $E$.
\end{enumerate}
\end{cor}

\begin{proof}
Let $\mu$ be an ergodic $G$-quasi-invariant probability measure 
$\widehat N$ attached to $\pi \vert_N$.

\vskip.2cm

\ref{iCor-RestNormalSub-Bis}
By Proposition~\ref{Prop-RestNormalSub-Bis}~\ref{iiProp-RestNormalSub-Bis},
we can assume that $\Hi = L^2(\widehat N, \mu, \Ki)$
for a Hilbert space $\Ki$ of some dimension $n \in \overline{\N^*}$
and that $\pi \vert_N$ coincides with the canonical representation 
$\pi_\mu^{(n)}$ of $N$ in $L^2(\widehat N, \mu, \Ki)$.
As we have already observed in the beginning of Section~\ref{SecCandecpm},
concerning the projection-valued measure
$E_\mu^{(n)}$ associated to $\pi_\mu^{(n)}$,
we have, for every $B \in \mathcal B (\widehat N)$~:
$$
E_\mu^{(n)} (B) \, = \, 0 
\, \Longleftrightarrow \,
\mu(B) \, = \, 0.
$$
This shows Item \ref{iCor-RestNormalSub-Bis}

\vskip.2cm

Item \ref{iiCor-RestNormalSub-Bis}
is an immediate consequence of Item \ref{iCor-RestNormalSub-Bis}.
\end{proof}

\begin{rem}[\textbf{Mackey machine}]
% 2.F.5
\label{MackMach}
\index{Mackey machine}
We give some indications on the \textbf{Mackey machine},
which is a procedure for producing and analyzing
representations of a locally compact group $G$
in terms of the representations of a closed normal subgroup $N$
and of subgroups of $G/N$.
The terminology refers to \cite{Mack--58}.
For details and proofs, see also \cite[Section 6.4]{Foll--16}
or/and \cite[Section 4.3]{KaTa--13}.
% voir Theorem 4.29 page 161
The procedure can be adapted to a variety of situations,
with various degrees of sophistication.
\par

Here, we only address the easiest situation,
and assume that the three following properties hold.
\begin{enumerate}[label=(\alph*)]
\item
The group $G$ is a \emph{semi-direct product} $G = H \ltimes N$
of a closed subgroup $H$ of $G$
with an \emph{abelian} closed normal subgroup $N$ of $G$.
\item
The groups $H$ and $N$ are second-countable, hence so is $G$.
\item
The semi-direct product is \emph{regular},
i.e., there exists a Borel subset $B \subset \widehat N$
such that $B$ intersects each $H$-orbit in $\widehat N$ in exactly one point,
where the right action $\widehat N \curvearrowleft H$
is the restriction to $H$ of the natural action of $G$ on the dual of $N$
(this is not the orthodox notion of a regular semi-direct product,
but our Condition (c) implies the usually defined property).
\end{enumerate}
The machine produces irreducible representations, indeed all of them, as follows.
For $\chi \in \widehat N$, denote by $\mathcal O_\chi$ its $H$-orbit
and by $H_\chi = \{ h \in H \mid \chi^h = \chi \}$ its isotropy in~$H$.
For any representation $\sigma$ of $H_\chi$,
we have a representation $\sigma \chi$ defined by
$$
(\sigma \chi)(h,n) \ = \, \sigma(h) \chi(n)
\hskip.5cm \text{for all} \hskip.2cm
(h,n) \in H_\chi \ltimes N .
$$
Set $\pi_{\sigma, \chi} = \Ind_{H_\chi \ltimes N}^G (\sigma \chi)$. Then
\begin{enumerate}[label=(\arabic*)]
\item\label{iDEMackMach}
Let $\chi \in \widehat N$
and $\sigma \in \widehat{H_\chi}$ an irreducible representation of $H_\chi$;
the representation $\pi_{\sigma, \chi}$ of $G$ is irreducible.
\item\label{iiDEMackMach}
Every irreducible representation of $G$ is equivalent to a representation $\pi_{\sigma, \chi}$,
for some $\chi \in \widehat N$ and $\sigma \in \widehat{H_\chi}$.
\item\label{iiiDEMackMach}
Two irreducible representations of $G$ of the previous form,
say $\pi_{\sigma_1, \chi_1}$ and $\pi_{\sigma_2, \chi_2}$, are equivalent
if and only if $\chi_1, \chi_2$ are in the same $H$-orbit,
i.e., there exists $h \in H$ such that $\chi_2 = (\chi_1)^h$,
and the representations $x \mapsto \sigma_1(x)$
and $x \mapsto \sigma_2(h^{-1} x h)$ of $H_{\chi_1}$ are equivalent.
\end{enumerate}
There are refined versions of the machine,
to deal with locally compact groups $G$
having a closed normal subgroup $N$ which is of type I,
and such that the extension $N \to G \to G/N$ satisfies appropriate conditions;
but sophistication increases,
among other reasons because projective representations have to be considered.
\end{rem}

%-----------------------------------------------------------------------
% End of chapter 2
%-----------------------------------------------------------------------
\chapter[Examples of irreducible representations]
{Examples of irreducible representations}
% Chapter 3
\label{Chapter-ExamplesIndIrrRep}

\emph{Using the irreducibility and equivalence criteria 
from Section~\ref{Section-IrrIndRep},
we construct families of non-equivalent irreducible representations 
for various discrete groups: 
\begin{enumerate}
\item[$\bullet$]
the infinite dihedral group $D_\infty$
(Section \ref{SectionInfDiGroup}),
\item[$\bullet$]
two-step nilpotent discrete groups,
in particular Heisenberg groups over fields
and over the integers
(Section \ref{Section-IrrRepTwoStepNil}), 
\item[$\bullet$]
the affine group $\Aff(\K)$ over a field $\K$
(Section \ref{Section-IrrRepAff}), 
\item[$\bullet$]
the Baumslag--Solitar group $\BS(1, p)$ for a prime $p$
(Section \ref{Section-IrrRepBS}), 
\item[$\bullet$]
the lamplighter group $\Z \ltimes \bigoplus_{k \in \Z} \Z / 2 \Z$
(Section \ref{Section-IrrRepLamplighter}),
\item[$\bullet$]
and the general linear group $\GL_n(\K)$ over $\K$
(Section \ref{Section-IrrRepGLN}).
\end{enumerate}
}
\par

\emph{
For $D_\infty$, we obtain in this way
a complete description of the dual $\widehat{D_\infty}$.
Moreover, this space coincides with the primitive dual $\Pri(D_\infty)$,
and also with the quasi-dual $\QD(D_\infty)$
introduced in Section \ref{SectionQuasidual}.
}
\par

\emph{
By contrast, for any group $\Gamma$ of Sections
\ref{Section-IrrRepTwoStepNil} to \ref{Section-IrrRepGLN},
we will see that it is impossible to select in these families, in a measurable way,
a set of pairwise inequivalent irreducible representations,
and a fortiori impossible to parametrize in a measurable way
the dual $\widehat \Gamma$ of $\Gamma$.
This is a manifestation of the fact that these groups are not of type I
(a notion defined in Section~\ref{SectionTypeI}).
Moreover, as checked below (Section \ref{noninjectivity}),
each group of Sections \ref{Section-IrrRepTwoStepNil} to \ref{Section-IrrRepGLN}
has an uncountable family of representations which are pairwise non-equivalent
and nevertheless weakly equivalent to each other.
}
\par

\emph{In later chapters,
% Chapter \ref{Chap-AllIrrRed} = 5
we will indicate further constructions providing more irreducible representations
of the groups listed above.}
\par

\emph{The final Section \ref{ExamplesND} refers briefly
to the duals of some non-discrete groups.
}

\section
{Infinite dihedral group}
% Section 3.A
\label{SectionInfDiGroup}

\index{Dihedral group}
Let $D_\infty$ be the \textbf{infinite dihedral group},
viewed here as the semi-direct product $\Z / 2 \Z \ltimes \Z$,
for the non-trivial action of the two-element group $Z / 2 \Z = \{1, \varepsilon\}$ on $\Z$, 
that is, for the action given by $\varepsilon \cdot n = -n$ for all $n \in \Z$.
The group has a matrix form:
$D_\infty = \begin{pmatrix} \{\pm 1 \} & \Z \\ 0 & 1 \end{pmatrix}$.
It is well known that $D_\infty$ is also isomorphic to
the free product of two groups of order two.
\par

Observe that the commutator subgroup $[D_\infty, D_\infty]$ 
is the subgroup $2 \Z$ of $\Z$,
and that the abelianized group $(D_\infty)_{\rm ab} = D_\infty/[D_\infty, D_\infty]$ 
is isomorphic to the Vierergruppe $V_4 = \Z / 2 \Z \times \Z / 2 \Z$.
\par

In this section, we write $\Gamma$ for $D_\infty$ 
and $N$ for its normal subgroup $\Z$ of index $2$.
As usual, we identity $\widehat N$
with the circle group $\T = \{ z \in \C \mid \vert z \vert = 1\}$
through the map
$$
\T \, \overset{\approx}{\longrightarrow} \, \widehat N,
\hskip.2cm 
z \, \mapsto \, \chi_z ,
$$
where 
$$
\chi_z (n) \, = \, z^n
\hskip.5cm \text{for all} \hskip.2cm
n \in \Z ,
$$
and we view $\widehat N$ as a right $\Gamma$-space.
The action of $\varepsilon$ on $\widehat N$ is described by
$$
\chi_z^\varepsilon \, = \, \chi_{\overline z}^{\phantom{\varepsilon}}
\hskip.5cm \text{for all} \hskip.2cm
z \in \T.
$$
In particular, $\chi_1$ and $\chi_{-1}$ are $\Gamma$-fixed.
For $z \ne \pm 1$, the $\Gamma$-orbit of $\chi_z$ 
consists of the two-element set $\{ \chi_z, \chi_{\overline z}^{\phantom{a}} \}$. 
The space of $\Gamma$-orbits in $\widehat N$ can be identified
with the quotient space $\T/{\sim}$ for the equivalence relation
on $\T$ defined by $z \sim \overline z$.
\par

For every $z \in \T$, define a representation $\rho_z$ of $\Gamma$ by 
$$
\rho_z \, = \, \Ind_N^\Gamma \chi_z .
$$
It is of dimension $2$. Observe that $\rho_{\overline z}$ is equivalent to $\rho_z$,
since $\chi_{\overline z}$ is in the $\Gamma$-orbit of $\chi_z$
(see Proposition~\ref{PropConjIndRep}).

\begin{prop}
% 3.A.1
\label{Prop-RepGrDiedral}
Let $\Gamma = D_\infty$. 
The map 
$$
\Phi \, \colon \, \widehat{\Gamma_{\rm ab}} \sqcup
\big( (\T \smallsetminus \{1,-1 \})/{\sim} \big)
\hskip.2cm \longrightarrow \hskip.2cm
\widehat \Gamma,
$$
defined by $\Phi(\chi) = \chi$ for $\chi \in \widehat{\Gamma_{\rm ab}}$
and $\Phi(\{z, \overline z\}) = \rho_z$ for $z \in \T \smallsetminus \{1,-1 \}$, 
is a bijection.
\end{prop}

Note that $\widehat{\Gamma_{\rm ab}} = \widehat{V_4}$ is a four elements set,
and $(\T \smallsetminus \{1,-1 \})/{\sim}$ is in natural bijection
with the upper open half-circle
$\{ e^{i \theta} \in \T \mid 0 < \theta < \pi \}$.

\begin{proof}
The method of proof of Proposition~\ref{Prop-RepGrDiedral} 
is a special case of the Mackey machine,
mentioned in Remark \ref{MackMach}.

\vskip.2cm

$\bullet$ {\it First step.} 
We claim that the range of $\Phi$ is contained in $\widehat \Gamma$
and that $\Phi$ is injective.
\par

Indeed, for $z \in \T \smallsetminus \{1,-1 \}$,
the stabilizer of $\chi_z$ in $\Gamma$ is $N$,
so that the representation $\rho_z$ is irreducible, 
by the Mackey--Shoda irreducibility criterion
(Corollary \ref{Cor-NormalSubgIrr}).
Moreover, if $z,w \in \T$ are not in the same equivalence class, 
then $\rho_z$ and $\rho_w$ are not equivalent
by the Mackey--Shoda equivalence criterion
(Corollary~\ref{Cor-NormalSubgEquiv}).
This proves the non-trivial part of the first step.

\vskip.2cm

Let $(\rho, \Hi)$ be an irreducible representation of $\Gamma$. 
Assume that $\rho$ is not one-dimensional . To prove that $\Phi$ is surjective, 
we have to show that $\rho$ is equivalent to
$\rho_z$ for some $z \in \T \smallsetminus \{1,-1 \}$. 
We denote by $E$ the projection-valued measure on $\widehat N$
associated to $\rho \vert_N$, and by $S$ the support of $E$. 

\vskip.2cm

$\bullet$ {\it Second step.} 
We claim that $S$ consists of exactly 
one $\Gamma$-orbit in $\widehat N$.
\par

Indeed, by Proposition
\ref{Prop-RestNormalSub}~\ref{iiDEProp-RestNormalSub}, 
$S$ is the closure of a $\Gamma$-orbit in $\widehat N$.
Since the $\Gamma$-orbits are finite and hence closed, the claim follows.

\vskip.2cm

$\bullet$ {\it Third step.}
We claim that $S \ne \{\chi_1 \}$ and $S \ne \{\chi_{-1}\}$.
\par

Indeed, assume that $S = \{\chi_1 \}$ or $S = \{\chi_{-1}\}$. 
Then $\Hi = \Hi^{\chi_1}$ or $\Hi = \Hi^{\chi_{-1}}$, 
by Corollary~\ref{Cor-ContAbelian}. 
Observe that $\chi_1$ and $\chi_{-1}$ are trivial on the 
commutator subgroup $[\Gamma, \Gamma] = 2 \Z \subset N$. 
Therefore $\rho$ factorizes through $\Gamma_{\rm ab}$
and hence is one-dimensional.
This is a contradiction.

\vskip.2cm

$\bullet$ {\it Fourth step.} 
We have $\rho = \Ind_N^\Gamma \chi_z$ for some 
$z \in \T \smallsetminus \{-1,1 \}$.
\par

Indeed, by the second and the third step, we have $S = \{z, \overline z\}$ for some 
$z \in \T \smallsetminus \{-1,1 \}$. 
By Corollary~\ref{Cor-ContAbelian} again, we have 
$$
\Hi \, = \, \Hi^{\chi_z} \bigoplus \Hi^{\chi_{\overline z}}.
$$
The formula from Proposition 
\ref{Prop-RestNormalSub}~\ref{iDEProp-RestNormalSub} shows that 
$$
\Hi^{\chi_{\overline z}} \, = \, \rho(\varepsilon)\Hi^{\chi_z}.
$$
Since $\{e, \varepsilon\}$ is a transversal for $\Gamma/N$, we see that
$\rho = \Ind_N^\Gamma \sigma$, where the representation
$\sigma$ of $N$ is the multiple of $\chi_z$ 
defined on $\Hi^{\chi_z}$ (see Definition \ref{openeasierthanclosed}).
Since $\rho$ is irreducible, $\sigma$ is irreducible and hence $\sigma = \chi_z$
(that is, $\dim \Hi^{\chi_z} = 1$).
\end{proof}

Proposition~\ref{Prop-RepGrDiedral}
describes the dual $\widehat{D_\infty}$ of the infinite dihedral group as a set. 
We describe now this dual as a topological space.
\par

Set $C = \{e^{i \theta} \in \T \mid 0 \le \theta \le \pi \} \times \{1, -1 \}$;
it is a compact topological space,
homeomorphic to the the disjoint union of two closed intervals.
Define an equivalence relation on $C$
by $(e^{i \theta}, j) \sim (e^{i \theta'}, j')$ 
if either $0 < \theta = \theta' < \pi$,
or $\theta = \theta' \in \{0, \pi \}$ and $j = j'$;
we write $[ e^{i \theta}, j ] \in C / \sim$ for the equivalence class
of $(e^{i \theta}, j)$ in $C$.
Note that $C / \sim$ is a non-Hausdorff quasi-compact topological space,
that can be described as an interval having
pairs of non-separated points at each of its ends.
\par

For the infinite dihedral group, we choose the matrix model
$D_\infty = \begin{pmatrix} \{\pm 1 \} & \Z \\ 0 & 1 \end{pmatrix}$
and we will write $(\varepsilon, n)$ for 
$\begin{pmatrix} \varepsilon &n \\ 0 &1 \end{pmatrix}\in D_\infty$.
There are four unitary characters of $D_\infty$
denoted by $\rho_{0, 1}, \rho_{0, -1}, \rho_{\pi, 1}, \rho_{\pi, -1}$
and defined, for $(\varepsilon,n)\in D_\infty$, by
$$
\begin{aligned}
\rho_{0, 1}(\varepsilon, n) \, = \, 1 ,
\hskip1cm
&\rho_{0, -1}(\varepsilon, n) \, = \, \varepsilon ,
\\
\rho_{\pi, 1}(\varepsilon, n) \, = \, (-1)^n ,
\hskip1cm
&\rho_{\pi, -1}(\varepsilon, n) \, = \, (-1)^n\varepsilon .
\end{aligned}
$$
Let $z = e^{i \theta} \in \T$ be such that $0 \le \theta \le \pi$.
For $\Ind_N^\Gamma \chi_z$, we write now $\rho_\theta$,
rather than~$\rho_z$ as above.
According to Definition \ref{openeasierthanclosed},
the induced representation $\rho_\theta$
can be viewed as acting on the Hilbert space $\C^2$, and defined by
$$
\rho_\theta (1, n)\, = \,
\begin{pmatrix} e^{in\theta} & 0 \\ 0 & e^{-in\theta} \end{pmatrix}
\hskip.5cm \text{and} \hskip.5cm
\rho_\theta (-1, n) \, = \,
\begin{pmatrix} 0 & e^{in\theta} \\ e^{-in\theta} & 0 \end{pmatrix}
$$
for all $n \in \Z$.
The normalized function $\varphi_{\rho_\theta, s,t}$ of positive type
associated to $\rho_\theta$ and defined by a unit vector
$\xi = \begin{pmatrix} s \\ t \end{pmatrix} \in \C^2$ is given by
$$
\varphi_{\rho_\theta, s,t} (\varepsilon, n) \, = \,
\left\{ \begin{aligned}
e^{in\theta} \hskip.1cm \vert s \vert^2 + e^{-in\theta} \hskip.1cm \vert t \vert^2
\hskip2.9cm &\text{if} \hskip.2cm
\varepsilon = 1 ,
\\
e^{in\theta} \hskip.1cm \overline{s} \hskip.1cm t + e^{-in\theta} s \overline{t}
\, = \, 2 {\mathrm Re} (e^{in\theta} \overline{s} t )
\hskip.6cm &\text{if} \hskip.2cm
\varepsilon = -1 ,
\end{aligned} \right.
\leqno{({\rm PT})}
$$
for all $n \in \Z$
(see Definition \ref{exposfunctdiagcoeff}).
\par

Note that the representation $\rho_\theta$ is cyclic.
When $0 < \theta < \pi$, it is irreducible.
When $\theta = 0$, we have a direct sum $\rho_0 = \rho_{0,1} \oplus \rho_{0, -1}$
in a a direct sum of Hilbert spaces
$\C^2 = \C \begin{pmatrix} 1 / \sqrt 2 \\ 1 / \sqrt 2 \end{pmatrix}
\oplus \C \begin{pmatrix} \phantom{-}1 / \sqrt 2 \\ -1 / \sqrt 2 \end{pmatrix}$.
When $\theta = \pi$, we have $\rho_\pi = \rho_{\pi,1} \oplus \rho_{\pi, -1}$,
in the same decomposition of $\C^2$ in a direct sum
of two one-dimensional Hilbert spaces.

\begin{prop}
% 3.A.2
\label{PropRepGrDiedral}
Let the notation be as above.
Define a map $\Psi \, \colon C / \sim \hskip.1cm \to \widehat{D_\infty}$ by
$$
\begin{aligned}
\Psi([ e^{i \theta}, j ]) \, = \, \rho_\theta
\hskip1cm &\text{for} \hskip.2cm
[ e^{i \theta}, j ] \in C / \sim ,
\\
\Psi ( [1, 1] ) \, = \, \rho_{0, 1} , \hskip1cm
&\Psi ( [1, -1] ) \, = \, \rho_{0, -1} , 
\\
\Psi ( [-1, 1] ) \, = \, \rho_{\pi, 1} , \hskip1cm
&\Psi ( [-1, -1] ) \, = \, \rho_{\pi, -1} .
\end{aligned}
$$
Then $\Psi$ is a homeomorphism.
\end{prop}

\begin{proof}
By Proposition \ref{Prop-RepGrDiedral}, the map $\Psi$ is a bijection.

\vskip.2cm

$\bullet$ {\it First step.} 
The map $\Psi$ is continuous.
\par

Since the space $C / \sim$ is first countable, it suffices to consider sequences, rather than nets.
Consider a point $[e^{i \theta}, j]$ in $C / \sim$.
We have to check that, for any sequence $( [e^{i \theta_k}, j_k ] )_{k \ge 1}$
converging to $[ e^{i \theta}, j ]$ in $C / \sim$,
the sequence of representations $( \Psi( [e^{i \theta}, j] ) )_{k \ge 1}$
converges to $\Psi( [e^{i \theta}, j] )$ in $\widehat{D_\infty}$,
i.e., that,
for each $u, v \in \C$,
there exists a sequence $( \varphi_{\rho_{\theta_k}, u_k, v_k} )_{k \ge 1}$
of functions of positive type associated to the $\rho_{\theta_k}$~'s
such that
$$
\lim_{k \to \infty} \varphi_{\rho_{\theta_k}, u_k, v_k} \, = \, \varphi_{\rho_\theta, u, v}
\hskip.5cm
\text{(limit in the topology of simple convergence).}
\leqno{({\rm Lim })}
$$
\par

Assume first that $0 < \theta < \pi$.
Then $\Psi( [e^{i \theta}, j] ) = \rho_\theta$.
It is obvious from the explicit formulas (${\rm PT}$)
that $\lim_{k \to \infty} \varphi_{\rho_{\theta_k}, u, v} \, = \, \varphi_{\rho_\theta, u, v}$
for all $u, v \in \C$.
\par

Assume next that $\theta = 0$ and $j=-1$.
Then $\Psi( [1,-1] ) = \rho_{0, -1}$.
For any vector $\begin{pmatrix} \phantom{-}u \\ -u \end{pmatrix}$
in the space $\C \begin{pmatrix} \phantom{-}1 / \sqrt 2 \\ -1 / \sqrt 2 \end{pmatrix}$
of $\rho_{0, -1}$,
the function of positive type $\varphi_{\rho_0, u, -u}$
is a constant multiple of the unitary character $\rho_{0, -1}$.
It is again obvious from the explicit formulas (${\rm PT}$)
that $\lim_{k \to \infty} \varphi_{\rho_{\theta_k}, u, -u} \, = \, \varphi_{\rho_0, u, -u}$
for all $u \in \C$.
The arguments for the three other cases $(\theta, j) = (0, 1), (\pi, 1), (\pi, -1)$
are similar.

\vskip.2cm

$\bullet$ \emph{Second step.}
The map $\Psi^{-1} \, \colon \widehat{D_\infty} \to C / \sim$ is continuous.
\par

Let $A$ be a subset of $\widehat{D_\infty}$ and $\rho \in \overline{A}$.
We have to show that $\Psi^{-1}(\rho) \in \overline{ \Psi^{-1}(A) }$.
Two cases may occur.

\vskip.2cm

$\circ$ \emph{First case.}
$\rho$ is in the closure of $A \cap \{\rho_{0, 1}, \rho_{0, -1}, \rho_{\pi, 1}, \rho_{\pi, -1}\}$.
Then $\rho$ belongs to $A \cap \{\rho_{0, 1}, \rho_{0, -1}, \rho_{\pi, 1}, \rho_{\pi, -1}\}$,
since every unitary character of $D_\infty$ is a closed point in $\widehat{D_\infty}$
(see Corollary \ref{Example-CommutatorSubgroup}).
In particular, $\Psi^{-1}(\rho) \in \overline{\Psi^{-1}(A)}$. 

\vskip.2cm

$\circ$ \emph{Second case.}$
\rho$ is in the closure of $A \cap \{\rho_{\theta}\mid 0<\theta<\pi\}$.
\par
 
Let $\xi$ be a unit vector in the Hilbert space of $\rho$
and $\varphi_{\rho, \xi}$ the corresponding normalized function of positive type
associated to $\rho$.
Then there exists a sequence $(\theta_k)_{k\ge 1}$ in $(0, \pi)$ with $\rho_{\theta_k} \in A$
and sequences $(s_k)_{k \ge 1}, (t_k)_{k \ge 1}$ in $\C$
with $\vert s_k \vert^2 +\vert t_k \vert^2 =1$
such that 
$$
\lim_k \varphi_{\rho_{\theta_k}, s_k, t_k} (1,1) \, = \, \varphi_{\rho, \xi}(1,1).
$$
Therefore, we have (see the formulas ${\rm PT}$) 
$$
\lim_k \left( e^{i \theta_k} \vert s_k \vert^2 + e^{-i \theta_k} \vert t_k \vert^2 \right)
\, = \, \varphi_{\rho, \xi}(1,1).
\leqno{(*)}
$$
\par

Assume first that $\rho = \rho_{\theta}$ for $0 < \theta < \pi$. 
Choose $\xi = \begin{pmatrix} 1 \\ 0 \end{pmatrix}$.
Then $\varphi_{\rho, \xi}(1,1) = e^{i \theta}$ and it follows from $(*)$
and from the uniform convexity of the unit disc in $\R^2$
that $\lim_k \theta_k = \theta$. 
So, $\Psi^{-1}(\rho) = [ e^{i \theta}, 1 ] = \lim_k [ e^{i \theta_k}, 1 ]$
is in the closure of $\Psi^{-1}(A)$.
\par

Assume next that $\rho = \rho_{0, 1}$ or $\rho =\rho_{0, -1}$.
With $\xi=1$, we have $\varphi_{\rho, \xi}(1,1) = 1$ and therefore
$$
\lim_k \left( e^{i \theta_k} \vert s_k \vert^2 + e^{-i \theta_k} \vert t_k \vert^2 \right)
\, = \, 1.
$$
As before, this implies that $\lim_k \theta_k = 0$.
So, $\overline{\Psi^{-1}(A)}$ contains $[1, 1]$ and $[1, -1]$
and this shows that $\Psi^{-1}(\rho) \in \overline{\Psi^{-1}(A)}$.
\par 

Finally, assume that $\rho = \rho_{\pi, 1}$ or $\rho=\rho_{\pi, -1}$.
With $\xi = 1$, we have $\varphi_{\rho, \xi}(1,1) = -1$ and therefore
$$
\lim_k \left(e^{i\theta_k} \vert s_k \vert^2 + e^{-i\theta_k} \vert t_k \vert^2\right)
\, = \, -1.
$$
As before, this implies that $\lim_k \theta_k= \pi$.
So, $\overline{\Psi^{-1}(A)}$ contains $[-1, 1]$ and $[-1, -1]$
and this shows that $\Psi^{-1}(\rho)\in \overline{\Psi^{-1}(A)}$. 
\end{proof}

\section
{Two-step nilpotent groups}
% Section 3.2
\label{Section-IrrRepTwoStepNil}

\index{Two-step nilpotent group}
Let $\Gamma$ be a \textbf{two-step nilpotent group}, that is, 
$\Gamma$ is a group of which the derived group $[\Gamma, \Gamma]$
is contained in the centre $Z = Z(\Gamma)$.
Observe that $\Gamma$ could be abelian; we do not exclude this case,
as we sometimes have to consider abelian quotients of~$\Gamma$.
\par

Theorem~\ref{Theo-IrredTwoStepNil} provides
a procedure constructing irreducible representations of $\Gamma$,
which is reminiscent of Kirillov's orbit method for nilpotent connected real Lie groups 
\cite{Kiri--62}.
Our exposition is inspired in part by \cite{Howe--77b},
which provides a description of $\Pri(\Gamma)$
for finitely generated discrete nilpotent groups.
We will come back to representations of two-step nilpotent groups
in (Sub)sections \ref{SS:HeisGpsOverRings},
\ref{Sect-PrimIdealHeisenberg}, and \ref{ThomaDualTwoStepNil}.
% Ajouter aussi une r\'ef au Theorem \ref{Th-HoweKa} ???
% r\'eponse du 5 juillet 18 : non.

\begin{defn}
% 3.B.1
\label{defbihomomorphism}
\index{Bihomomorphism}
Let $\Delta$ be a group. A \textbf{bihomomorphism} from $\Delta$ to $\T$
is a map $\omega \, \colon \Delta \times \Delta \to \T$
such that the partial maps
$$
\omega_a \, \colon \, \Delta \, \to \, \T, \hskip.2cm \delta \mapsto \omega(a, \delta)
\hskip.5cm \text{and} \hskip.5cm
\omega^b \, \colon \, \Delta \, \to \, \T, \hskip.2cm \delta \mapsto \omega(\delta, b)
$$
ae homomorphisms.
For such a bihomomorphism $\omega$, the map
$$
\Delta \, \to \, \Hom (\Delta, \T) , \hskip.2cm a \mapsto \omega_a
$$
is a homomorphism;
when $\Delta$ is abelian, this is a homomorphism $\Delta \to \widehat \Delta$.
\end{defn}

\begin{lem}
% 3.B.2
\label{standard2stepnilp}
Let $\Gamma$ be a two-step nilpotent group. Let $Z$ be its centre.
\begin{enumerate}[label=(\arabic*)]
\item\label{iDEstandard2stepnilp}
Every maximal abelian subgroup of $\Gamma$ is normal in $\Gamma$
and contains its centre.
\item\label{iiDEstandard2stepnilp}
For every $\gamma \in \Gamma$, the maps
$$
\Gamma \to Z, \hskip.2cm a \mapsto [a, \gamma] 
\hskip.5cm \text{and} \hskip.5cm
\Gamma \to Z, \hskip.2cm a \mapsto [\gamma, a] = [a, \gamma]^{-1}
$$
are group homomorphisms.
\item\label{iiiDEstandard2stepnilp}
Let $\psi \in \widehat Z$ be a unitary character of $Z$. The map
$$
\omega \, \colon \, \Gamma/Z \times \Gamma/Z \to \T, 
\hskip.2cm 
(aZ, bZ) \mapsto \psi([ b, a])
$$
is a bihomomorphism from $\Gamma/Z$ to $\T$.
\item\label{ivDEstandard2stepnilp}
Let $\psi \in \widehat Z$ and $\omega$ be as in \ref{iiiDEstandard2stepnilp};
assume moreover that $\psi$ is faithful character. Then
$$
\Gamma/Z \, \longrightarrow \, \widehat{\Gamma/Z}, 
\hskip.2cm 
a Z \, \mapsto \, \omega_a
$$
is an injective homomorphism.
\end{enumerate}
\end{lem}

\begin{proof}
\ref{iDEstandard2stepnilp}
Maximal abelian subgroups exist in $\Gamma$, by Zorn's lemma;
let $N$ be such a subgroup.
Since $N$ is its own centralizer, $Z$ is contained in $N$.
Therefore $[\Gamma, N] \subset [N, N] \subset Z \subset N$.
It follows that $N$ is normal in $\Gamma$.

\vskip.2cm

\ref{iiDEstandard2stepnilp}
For all $a,b \in \Gamma$,
$$
\begin{aligned}
\phantom{x}
% tr\`es bizarre, mais \c ca ne va pas sans ce \phantom !!!!!
[ab, \gamma]
& \, = \,
b^{-1}a^{-1}\gamma^{-1} ab\gamma \, = \, 
b^{-1}(a^{-1}\gamma^{-1} a \gamma) \gamma^{-1}b\gamma
\\
& \, = \,
(a^{-1}\gamma^{-1} a \gamma) b^{-1} \gamma^{-1}b\gamma \, = \, 
[\gamma, a][\gamma, b] ,
\end{aligned}
$$
i.e., the map $\Gamma \to Z, \hskip.2cm a \mapsto [a, \gamma]$
is a group homomorphism.
The argument for $\Gamma \to Z, \hskip.2cm a \mapsto [\gamma, a]$
is similar.
\par
Observe that these homomorphisms factorize through $\Gamma/Z$.

\vskip.2cm

Claim \ref{iiiDEstandard2stepnilp} is a straightforward consequence
of \ref{iiDEstandard2stepnilp},
and Claim \ref{ivDEstandard2stepnilp}
is equally straightforward.
\end{proof}

We continue with a non-abelian two-step nilpotent group $\Gamma$,
its centre $Z$, and a maximal abelian subgroup $N$ of $\Gamma$.
We have a continuous homomorphism
$$
r \, \colon \, \widehat N \to \widehat Z, 
\hskip.2cm \chi \mapsto \chi \vert_Z,
$$
given by restriction to $Z$ of unitary characters of $N$.
By Pontrjagin theory, unitary characters extend from closed subgroups of LCA groups
(see, e.g., \cite[Chap.~II, \S~1, no~7, th.\ 4]{BTS1--2}).
In particular, the map $r$ is surjective.
Since $Z$ is central, the action of $\Gamma$ on $\widehat Z$ is trivial
and hence
$$
r(\chi^\gamma) \, = \, r(\chi)
\hskip.5cm \text{for all} \hskip.2cm 
\chi \in \widehat N, \gamma \in \Gamma.
$$
For $\psi \in \widehat Z$, set 
$$
\widehat N (\psi) \, := \, r^{-1}(\psi) \, = \,
\left\{ \chi \in \widehat N \hskip.1cm \big\vert \hskip.1cm
\chi \vert_Z = \psi \right\} .
$$
If we agree to identify $\widehat{N/Z}$ to a subgroup of $\widehat N$,
we also have for any $\chi_0 \in \widehat N (\psi)$
$$
\widehat N (\psi) \, = \,
\left\{\chi \in \widehat N \hskip.1cm \big\vert \hskip.1cm
\chi = \chi_0 \rho
\hskip.2cm \text{for some} \hskip.2cm
\rho \in \widehat{N/Z} \right\} .
$$
Observe that $\widehat N (\psi)$ is a $\Gamma$-invariant subset of $\widehat N$:
for $a \in \Gamma$ and $\chi \in \widehat N (\psi)$,
the conjugate unitary character $\chi^a$ is again in $\widehat N (\psi)$.

\begin{lem}
% 3.B.3
\label{Lem-FreeActDualHeis}
Let $\Gamma$ be a non abelian two-step nilpotent group,
$Z$ its centre,
$N$~a maximal abelian subgroup of $\Gamma$,
and $\psi \in \widehat Z$ a \emph{faithful} unitary character of $Z$;
let $a \in \Gamma$.
Let $\widehat N (\psi)$ be as above,
and, for $a \in \Gamma$, let $\omega_a$ be as in Definition \ref{defbihomomorphism}.
\begin{enumerate}[label=(\arabic*)]
\item\label{iDELem-FreeActDualHeis} 
Let $a \in \Gamma \smallsetminus N$. 
The restriction of $\omega_a$ to $N/Z$, that is the map 
$$
\omega_a \vert_{N/Z} \, \colon \, 
N/Z \to \T, \hskip.2cm n Z \mapsto \psi([n, a]) ,
$$
is a unitary character of $N/Z$ that is distinct from $1_{N/Z}$.
\item\label{iiDELem-FreeActDualHeis} 
Let $\chi \in \widehat N (\psi)$ and $a \in \Gamma$.
Then $\chi^{a^{-1}} = (\omega_a \vert_{N/Z})\chi \in \widehat N (\psi)$, 
where $\omega_a \vert_{N/Z}$ is viewed as a unitary character of $N$.
\item\label{iiiDELem-FreeActDualHeis} 
The group $\Gamma/N$ acts freely on $\widehat N (\psi)$.
\item\label{ivDELem-FreeActDualHeis} 
The $\Gamma$-orbit of every element $\chi \in \widehat N (\psi)$
is dense in $\widehat N (\psi)$.
\end{enumerate}
\end{lem}

\begin{proof}
\ref{iDELem-FreeActDualHeis} 
Since $\Gamma$ is not abelian, $N \subsetneqq \Gamma$.
Since $a \notin N$, there exists $n \in N$ such that $[n,a] \ne 1$.
Since $\psi$ is faithful, we have
$\left( \psi_a \vert_{N/Z} \right) (n Z) = \psi([n,a]) \ne 1$.
Therefore $\psi_a \vert_{N/Z} \ne 1_{N/Z}$.

\vskip.2cm

\ref{iiDELem-FreeActDualHeis} 
Since $\chi \vert_Z = \psi$, we have, for every $n \in N$,
$$
\chi^{a^{-1}}(n) \, = \, \chi(a^{-1} n a) \, = \, \chi(n[n,a]) \, = \,
\psi([n,a]) \chi(n) \, = \, \omega_a(n Z)\chi(n).
$$

\vskip.2cm

\ref{iiiDELem-FreeActDualHeis} 
Let $a \in \Gamma \smallsetminus N$ and $\chi \in \widehat{N_\psi}$. 
By \ref{iDELem-FreeActDualHeis}, we have $\omega_{a^{-1}} \vert_{N/Z} \ne 1_{N/Z}$.
Therefore $\chi^a \ne \chi$, by \ref{iiDELem-FreeActDualHeis}.

\vskip.2cm

\ref{ivDELem-FreeActDualHeis} 
Any two elements in $\widehat N (\psi)$ 
differ by a unitary character of $N$ that is trivial on~$Z$. 
In view of \ref{iiDELem-FreeActDualHeis}, 
it is enough to show that the subgroup
$$
H \, := \, \{ \rho \in \widehat{N/Z} \mid \rho = \omega_a \vert_{N/Z}
\hskip.2cm \text{for some} \hskip.2cm
a \in \Gamma \}
$$
is dense in $\widehat{N/Z}$. 
\par

By Pontrjagin duality, it is equivalent to show that 
$H^\perp = \{e\}$, where 
$$
H^\perp \, = \, \{nZ \in N/Z \mid \rho(nZ) = 1 
\hskip.2cm \text{for all} \hskip.2cm 
\rho \in H \}
$$
(see Corollaire 6 of Th\'eor\`eme 4 in \cite[Chap.~II, \S~1, no~7]{BTS1--2}).
\par

Let $n \in N$ be such that
$\omega_a (n Z) = 1$ for all $a \in \Gamma$,
i.e., such that $\psi([n,a]) = 1$ for all $a \in \Gamma$.
Then $n \in Z$, because $\psi$ is faithful,
i.e., $n Z$ is the neutral element of $N/Z$,
as was to be shown.
\end{proof}

\begin{defn}
% 3.B.4
\label{defcentralchar}
Let $G$ be a topological group; let $Z$ denote its centre.
Let $\pi$ be an irreducible representation of $G$.
Schur's lemma implies that, for $g \in G$,
there exists $\lambda_\pi(g) \in \T$
such that $\pi(g) = \lambda_\pi(g) \mathrm{Id}_{\Hi_\pi}$.
The map
$$
\lambda_\pi \, \colon \, Z \to \T , \hskip.2cm g \mapsto \lambda_\pi(g)
$$
is a unitary character of $Z$ called the \textbf{central character} of $\pi$.
It only depends on the equivalence class of $\pi$.
For more one this, see Section~\ref{Section-CentralCharacter}.
\index{Central character}
\index{Character! $3$@central}
\end{defn}

\begin{theorem}
% 3.B.5
\label{Theo-IrredTwoStepNil}
Let $\Gamma$ be a non-abelian two-step nilpotent discrete group with centre $Z$
and $N$ a maximal abelian subgroup of $\Gamma$.
Let $\psi$ be a \emph{faithful} unitary character of $Z$; 
set $\widehat N (\psi) = \{ \chi \in \widehat N \mid \chi \vert_{Z} = \psi \}$.
\begin{enumerate}[label=(\arabic*)]
\item\label{iDETheo-IrredTwoStepNil}
Let $\chi \in \widehat N (\psi)$.
Then $\Ind_N^\Gamma \chi$ is an irreducible representation of $\Gamma$,
with central character $\psi$.
\item\label{iiDETheo-IrredTwoStepNil}
Let $\chi_1, \chi_2 \in \widehat N (\psi)$ be 
such that $\chi_1^\gamma \ne \chi_2$ for every $\gamma \in \Gamma$. 
Then $\Ind_N^\Gamma \chi_1$ and $\Ind_N^\Gamma \chi_2$ 
are non-equivalent irreducible representations of $\Gamma$.
\end{enumerate}
\end{theorem}

\begin{proof}
\ref{iDETheo-IrredTwoStepNil} 
By Lemma \ref{Lem-FreeActDualHeis}~\ref{iiiDELem-FreeActDualHeis} 
we have $\chi^\gamma \ne \chi$ for every
$\gamma \in \Gamma \smallsetminus N$. 
Therefore $\Ind_N^\Gamma \chi$ is irreducible
by the Mackey--Shoda irreducibility criterion (Corollary~\ref{Cor-NormalSubgIrr}).
\par

It remains to show that $\psi$ is the central character of $\pi = \Ind_N^\Gamma \chi$.
Choose a right transversal $T$ for $N$ in $\Gamma$,
so that $\Gamma = \bigsqcup_{t \in T} Nt$.
For $n \in N$ and $t \in T$,
let $\alpha(t, n) \in N$ and $t \cdot n \in T$
be as in Construction \ref{constructionInd}(2),
so that $tn = \alpha(t, n) (t \cdot n)$.
In the particular case of $z \in Z$,
we have $\alpha(t, z) = z$ and $t \cdot z = t$.
It follows that
$$
(\pi(t) f) (t) \, = \, \chi(z) f(t) \, = \, \psi(z) f(t)
$$
for all $z \in Z$, $f \in \ell^2(T, \C) = \Hi_\pi$, and $t \in T$.
Therefore $\pi(z) = \psi(z)I$ for every $z \in Z$.

\vskip.2cm

\ref{iiDETheo-IrredTwoStepNil}
It follows from \ref{iDETheo-IrredTwoStepNil} that 
the representations $\Ind_N^\Gamma \chi_1$ and $\Ind_N^\Gamma \chi_2$
are irreducible and from the Mackey--Shoda equivalence criterion 
(Corollary~\ref{Cor-NormalSubgEquiv}) 
that they are not equivalent.
\end{proof}

\begin{rem}
% 3.B.6
\label{Rem-IrredTwoStepNil}
The central characters of the irreducible representations of $\Gamma$ 
in Theorem~\ref{Theo-IrredTwoStepNil} are faithful.
Irreducible representations of $\Gamma$ having non-faithful central characters 
can be constructed as follows.
(Note that a group need not contain any faithful central character;
this is for example the case of groups containing
non-cyclic finite central subgroups.)
\par

Let $\psi$ be a unitary character of the centre $Z$ of $\Gamma$,
not necessarily a faithful character.
Consider the quotient group $\Gamma / \ker \psi$,
its centre $Z(\Gamma / \ker \psi)$,
and the canonical projection $p \, \colon \Gamma \twoheadrightarrow \Gamma / \ker \psi$.
For some time, we view $\psi$
as a unitary character of $Z / \ker \psi$,
and extend it as a unitary character, say $\widetilde \psi$, 
of $Z(\Gamma / \ker \psi)$;
observe that the subgroup $\ker \widetilde \psi$ of $\Gamma / \ker \psi$
is central, and therefore normal, so that
$p^{-1}\left( \ker \widetilde \psi \right)$ is a normal subgroup of $\Gamma$.
Set $Z_\psi = p^{-1} \left( Z(\Gamma / \ker \psi) \right)$;
then $Z_\psi$ is a normal subgroup of $\Gamma$,
and $Z(\Gamma / \ker \psi) = Z_\psi / \ker \psi$.
\par

We view now $\widetilde \psi$ as a unitary character of $Z_\psi$,
and therefore $\ker \widetilde \psi$ as a subgroup of $\Gamma$ which is normal;
note that $\Gamma \supset Z_\psi \supset \ker \widetilde \psi$.
Consider the quotient group 
$\overline \Gamma = \Gamma / \ker \widetilde \psi$.
We claim that $Z_\psi / \ker \widetilde \psi$ is the centre of $\overline \Gamma$.
\par

Since $\ker \psi \subset \ker \widetilde \psi$,
the quotient $Z_\psi / \ker \widetilde \psi$
is in the centre of $\overline \Gamma$.
Conversely, let $\gamma \in \Gamma$
be such that $\gamma \ker \widetilde \psi$
is in the centre of $\overline \Gamma$;
then $\gamma x \in x \gamma \ker \widetilde \psi$,
i.e., $[\gamma, x] \in \ker \widetilde \psi$,
for all $x \in \Gamma$.
Since $[\gamma, x] \in [\Gamma, \Gamma] \subset Z$,
and $\psi \vert_Z = \widetilde \psi \vert_Z$,
it follows that $[\gamma, x] \in \ker \psi$ for all $x \in \Gamma$,
i.e., that $\gamma \in Z_\psi$;
hence the centre of $\overline \Gamma$ 
is contained in $Z_\psi / \ker \widetilde \psi$.
This proves the claim.
\par

Moreover, $\widetilde \psi$, viewed as a unitary character 
of the centre of $\overline \Gamma$, is faithful.
Now, Theorem~\ref{Theo-IrredTwoStepNil} applies to $\widetilde \psi$
and provides an irreducible representation $\pi$ of $\overline \Gamma$
having $\widetilde \psi$ as central character. 
Lifting $\pi$ to $\Gamma$, we obtain 
an irreducible representation of $\Gamma$ having $\psi$ as central character.
\end{rem}

\subsection*{Heisenberg group over a ring}

Next, we treat Heisenberg groups,
which are prominent examples of two-step nilpotent groups.
We will deal with them in a way
independent of Theorem~\ref{Theo-IrredTwoStepNil}
(which concerns faithful central characters)
and Remark~\ref{Rem-IrredTwoStepNil}
(which concerns non-faithful central characters).
\par

\index{Heisenberg group! $1$@$H(R)$ over a ring $R$}
Let $R$ be a commutative ring with unit.
Let $\Gamma = H(R)$ be the \textbf{Heisenberg group} over $R$, that is,
$$
\Gamma \, = \, 
\begin{pmatrix} 1 & R & R \\ 0 & 1 & R \\ 0 & 0 & 1 \end{pmatrix} .
$$
We identify $\Gamma$ with $R^3$, 
by writing $(a,b,c)$ for the matrix
$$
\begin{pmatrix}
1 & a & c \\
0 & 1 & b \\
0 & 0 & 1
\end{pmatrix}.
$$
The group law on $\Gamma$ is given by 
$$
(a,b,c) (a', b', c') \, = \, (a+a',b+b', c+c'+ab') 
\hskip.5cm \text{for} \hskip.2cm 
(a,b,c), \hskip.1cm (a', b', c') \in R^3.
$$
The centre of $\Gamma$ is $Z = \{(0,0, c) \mid c \in R\} \approx R$.
Since 
$$
[(a,b,c), (a', b', c')] \, = \, (0,0, ab'-a'b) \in Z 
\hskip.5cm \text{for} \hskip.2cm 
(a,b,c), (a', b', c') \in R^3,
$$
we see that $\Gamma$ is two-step nilpotent;
moreover, $[\Gamma, \Gamma] = Z$.
The bihomomorphism associated to a faithful unitary character $\psi$ of $Z$
as in Lemma \ref{standard2stepnilp}~\ref{iiiDEstandard2stepnilp}
can be viewed as the ``symplectic form" $\omega \, \colon R^2 \times R^2 \to \T$,
given by 
$$
\omega((a,b), (a',b')) \, = \, \psi(ab'-a'b) 
\hskip.5cm \text{for} \hskip.2cm
(a,b), (a', b') \in R^2.
$$
The subgroup of $\Gamma$
$$
N \, = \, \{(0,b, c) \ \in \Gamma \mid b,c \in R\} \, \approx \, R^2 
$$
is maximal abelian.
\par

Fix a unitary character $\psi$ of the centre $Z$ of $\Gamma$.
Let $I_\psi$ be the \textbf{ideal} of $R$ defined by 
$$
I_\psi \, := \, \{ a \in R \mid aR \subset \ker \psi \};
$$
in other words, $I_\psi$ is the largest ideal of $R$ contained in $\ker \psi$.
Define
$$
Z_\psi \, := \, \{(a, b, c) \in \Gamma \mid a, b \in I_\psi, \hskip.1cm c \in R\}
\hskip.2cm \text{and} \hskip.2cm
N_\psi \, := \, \{(a, b, c) \in \Gamma \mid a \in I_\psi, \hskip.1cm b,c \in R\} .
$$
As easily checked, $Z_\psi$ and $N_\psi$ are normal subgroups of $\Gamma$,
and $Z_\psi \subset N_\psi$. 
The subgroup $\ker \psi$ of $Z$ is central, and therefore normal, in $\Gamma$;
denote by $\overline \Gamma$ the quotient group $\Gamma / \ker \psi$
and by $p \, \colon \Gamma \twoheadrightarrow \overline \Gamma$
the canonical projection.

\begin{lem}
% 3.B.7
\label{Lem-IrredHeisRing}
For $\psi \in \widehat Z$.
Let $Z_\psi$, $N_\psi$, and $p \, \colon \Gamma \twoheadrightarrow \overline \Gamma$
be as above
\begin{enumerate}[label=(\arabic*)]
\item\label{iDELem-IrredHeisRing}
The subgroup $Z_\psi$ is the inverse image in $\Gamma$
of the centre $\overline Z$ of $\overline \Gamma$.
\item\label{iiDELem-IrredHeisRing}
The subgroup $N_\psi$ is the inverse image in $\Gamma$
of the subgroup $p(N_\psi)$,
which is a maximal abelian subgroup in $\overline \Gamma$.
\end{enumerate}
In particular, if $\psi$ is faithful, then $Z_\psi = Z$ and $N_\psi = N$.
\end{lem}

\begin{proof}
Let $\gamma = (a, b, c) \in \Gamma$.

\vskip.2cm

\ref{iDELem-IrredHeisRing}
We have $[\gamma, \gamma'] = (0, 0, ab' - a'b) \in \ker \psi$
for all $\gamma' = (a', b', c') \in \Gamma$
if and only if $ab' \in \ker \psi$ and $ba' \in \ker \psi$ for all $a',b' \in R$,
that is, 
if and only $a \in I_\psi$ and $b \in I_\psi$.
Therefore, $[\gamma, \gamma'] \in \ker \psi$ for all $\gamma' \in \Gamma$
if and only if $\gamma \in Z_\psi$. 
This shows that $Z_\psi = p^{-1} (\overline Z)$.

\vskip.2cm

\ref{iiDELem-IrredHeisRing}
We have $[\gamma, \gamma'] = (0,0, ab' - a'b) \in \ker \psi$
for all $\gamma' = (a', b', c') \in N_\psi$
if and only if $ab' \in \ker \psi$ for all $b' \in R$,
that is, if and only $a \in I_\psi$.
Therefore, $[\gamma, \gamma'] \in \ker \psi$ for all $\gamma' \in N_\psi$
if and only if $\gamma \in N_\psi$.
This shows that $\overline{N_\psi} := p(N_\psi)$
is a maximal abelian subgroup in $\overline \Gamma$ 
and that $N_\psi = p^{-1}(\overline{N_\psi})$.
\end{proof}

Let $\psi \in \widehat Z$;
view the ring $R$ and the ideal $I_\psi$ of $R$ defined above as abelian groups. 
For $(\alpha, \beta) \in \widehat{I_\psi} \times \widehat R$,
define $\chi_{\psi, \alpha, \beta} \, \colon N_\psi \to \T$ by the formula
$$
\chi_{\psi, \alpha, \beta}(a,b,c) \, = \, \alpha(a) \beta(b) \psi(c)
\hskip.5cm \text{for} \hskip.2cm
a \in I_\psi, b, c \in R.
$$
The following lemma gives a concrete description of the sets
$$
\widehat{Z_\psi} (\psi)
\, := \, \{ \chi \in \widehat{Z_\psi} \mid \chi \vert_{Z} = \psi \}
\hskip.2cm \text{and} \hskip.2cm
\widehat{N_\psi} (\psi)
\, := \, \{ \chi \in \widehat{N_\psi} \mid \chi \vert_{Z} = \psi \} .
$$

\begin{lem}
% 3.B.8
\label{Lem-IrredHeisRing2}
Let $\psi \in \widehat Z$. With the notation introduced above, the maps
$$
\begin{aligned}
&\widehat{I_\psi} \times \widehat{I_\psi} \, \to \, \widehat {Z_\psi} (\psi) ,
\hskip.2cm
&(\alpha, \beta) \, \mapsto \, \chi_{\psi, \alpha, \beta}
\\
&\widehat{I_\psi} \times \widehat R \, \to \, \widehat {N_\psi} (\psi) ,
\hskip.2cm
&(\alpha, \beta) \, \mapsto \, \chi_{\psi, \alpha, \beta}
\end{aligned}
$$
are bijections.
\end{lem}

\begin{proof}
We give the proof for $\widehat{N_\psi} (\psi)$.
That for $\widehat{Z_\psi} (\psi)$ is similar.
\par

Let $(\alpha, \beta) \in \widehat{I_\psi} \times \widehat R$.
Then $\chi_{\psi, \alpha, \beta}$ is a unitary character of $N_\psi$.
Indeed, let $\gamma = (a, b, c)$ and $\gamma' = (a', b', c')$ in $N_\psi$. 
Since $ab' \in \ker \psi$, we have
$$
\begin{aligned}
\chi_{\psi, \alpha, \beta}(\gamma \gamma')
\, &= \, \chi_{\psi, \alpha, \beta}(a+a',b+b', c+c'+ab')
\\
\, &= \, \alpha(a+a') \beta(b+b') \psi (c+c') \psi(ab')
\\
\, &= \, \alpha(a+a') \beta(b+b') \psi(c+c')
\\
\, &= \, \alpha(a) \alpha(a') \beta(b) \beta(b') \psi(c) \psi(c')
\\
\, &= \, \chi_{\psi, \alpha, \beta}(\gamma) \chi_{\psi, \alpha, \beta}(\gamma').
\end{aligned}
$$
It is obvious that $\chi_{\psi, \alpha, \beta}\vert_{Z} = \psi$.
\par 

Conversely, let $\chi$ be a unitary character of $N_\psi$
such that $\chi \vert_Z = \psi$. Let $a,a' \in I_\psi$ and $b, b' \in R$.
Since $\chi(0, 0, ab'-a'b) = \psi(ab'-a'b) = 1$, we have
$$
\begin{aligned}
\chi(a+a',b+b',0)
\, &= \, \chi ((a+a', b+b', 0) (0, 0, ab'-a'b))
\\
\, &= \, \chi (a+a', b+b', ab'-a'b)
\\
\, &= \, \chi((a,b,0) (a', b', 0))
\\
\, &= \, \chi(a, b, 0) \chi(a' ,b', 0).
\end{aligned}
$$
This shows that $(a,b) \mapsto \chi(a,b,0)$ is a unitary character
of the group $I_\psi \times R$,
which can therefore be written as $(a,b) \mapsto \alpha(a) \beta(b)$
for a uniquely defined pair 
$(\alpha, \beta) \in \widehat{I_\psi} \times \widehat R$.
Moreover, since $\chi \vert_Z = \psi$
and since $Z$ is contained in the centre of $N_\psi$, 
it follows that $\chi = \chi_{\psi, \alpha, \beta}$.
\end{proof}

We are ready to construct a family of irreducible representations of $H(R)$.

\begin{theorem}
% 3.B.9
\label{Theo-IrredHeisRing}
Let $R$ be a commutative ring with unit,
$\Gamma$ the Heisenberg group $H(R)$,
and $Z$ the centre of $\Gamma$.
Let $\psi \in \widehat Z$ be a unitary character;
let $N_\psi$ be the normal subgroup of $\Gamma$
defined just before Lemma \ref{Lem-IrredHeisRing}.
\begin{enumerate}[label=(\arabic*)]
\item\label{iDETheo-IrredHeisRing}
Let $\chi$ be a unitary character of $N_\psi$ such that $\chi \vert_Z = \psi$.
Then $\Ind_{N_\psi}^\Gamma \chi$ is an irreducible representation of $\Gamma$,
and its central character is $\psi$.
\item\label{iiDETheo-IrredHeisRing}
Let $\chi_1, \chi_2$ be unitary characters of $N_\psi$
such that $\chi_1 \vert_Z = \chi_2 \vert_Z = \psi$,
and $\chi_1^\gamma \ne \chi_2$ for every $\gamma \in \Gamma$. 
Then $\Ind_{N_\psi}^\Gamma \chi_1$ and $\Ind_{N_\psi}^\Gamma \chi_2$ 
are non-equivalent irreducible representations of $\Gamma$.
\item\label{iiiDETheo-IrredHeisRing}
Let $\psi' \in \widehat Z$, with $\psi' \ne \psi$.
Let $\chi \in \widehat{N_\psi}$ be such that $\chi \vert_Z = \psi$
and $\chi' \in \widehat{N_{\psi'}}$ be such that $\chi' \vert_Z = \psi'$.
Then the irreducible representations
$\Ind_{N_\psi}^\Gamma \chi$ and $\Ind_{N_{\psi'}}^\Gamma \chi'$ of $\Gamma$
are not equivalent.
\end{enumerate}
\end{theorem}

\begin{proof}
\ref{iDETheo-IrredHeisRing}
To show that $\Ind_{N_\psi}^\Gamma \chi$ is irreducible,
by the Mackey--Shoda criterion (Corollary~\ref{Cor-NormalSubgIrr}),
it suffices to check that $\chi^{\gamma^{-1}} \ne \chi$
for all $\gamma \in \Gamma \smallsetminus N_\psi$.
An argument similar to that of
\ref{Lem-FreeActDualHeis}~\ref{iiiDELem-FreeActDualHeis}
shows that this is the case.
\par

Indeed, let $\gamma \in \Gamma$
be such that $\chi^{\gamma^{-1}} = \chi$.
Since $\chi \vert_Z = \psi$, we have
$$
\chi(\gamma') \, = \, \chi^{\gamma^{-1}} (\gamma')
\, = \, \chi(\gamma'[\gamma', \gamma])
\, = \, \psi([\gamma', \gamma]) \chi(\gamma')
\hskip.5cm \text{for all} \hskip.2cm
\gamma' \in N_\psi .
$$
Therefore $[N_\psi, \gamma] \subset \ker \psi$,
and hence $\gamma \in N_\psi$
since the image of $N_\psi$ in $\Gamma/\ker \psi$ is a maximal abelian subgroup
(Lemma~\ref{Lem-IrredHeisRing}).
\par

For the claim concerning the central character,
see the proof of Theorem~\ref{Theo-IrredTwoStepNil}~\ref{iDETheo-IrredTwoStepNil}.

\vskip.2cm

\ref{iiDETheo-IrredHeisRing}
The claim is a direct consequence of the Mackey--Shoda equivalence criterion 
(Corollary~\ref{Cor-NormalSubgEquiv}).

\vskip.2cm

\ref{iiiDETheo-IrredHeisRing}
The representations
$\Ind_{N_\psi}^\Gamma \chi$ and $\Ind_{N_{\psi'}}^\Gamma \chi'$
are not equivalent because their central characters are not equal.
\end{proof}

In the situation of Theorem~\ref{Theo-IrredHeisRing},
the representations $\Ind_{N_\psi}^\Gamma \chi$
and $\Ind_{N_\psi}^\Gamma \chi^\gamma$
are equivalent for all $\gamma \in \Gamma$, by Proposition \ref{PropConjIndRep}.
The theorem can therefore be restated as follows.

\begin{cor}
% 3.B.10
\label{Cor-IrredHeisRingRestated}
Let the notation be as in Theorem \ref{Theo-IrredHeisRing}.
For $\psi \in \widehat Z$,
let the normal subgroup $N_\psi$ of $\Gamma$
and the subset $\widehat{N_\psi} (\psi) = \{ \chi \in \widehat{N_\psi} \mid \chi \vert_{Z} = \psi \}$
of the dual of $N_\psi$
be as above.
Define 
$$
\displaystyle
X \, := \, \bigsqcup_{\psi \in \widehat Z} \{ \psi \} \times \widehat{N_\psi} (\psi) ,
$$
a set equipped with the right $\Gamma$-action 
$$
X \times \Gamma \, \to \, X,
\hskip.2cm
((\psi, \chi), \gamma) \, \mapsto \, (\psi, \chi^\gamma).
$$
\par

Then $\Ind_{N_\psi}^\Gamma \chi$ is irreducible for all $(\psi, \chi) \in X$,
and the map 
$$
X \, \to \, \widehat \Gamma, 
\hskip.2cm
\chi \, \mapsto \, \Ind_{N_\psi}^\Gamma \chi
$$
factorizes to an injective map from the space of $\Gamma$-orbits in 
$X$ to the dual $\widehat \Gamma$ of $\Gamma$.
\end{cor}

We will consider two particular cases:
$R = \K$, a field, and $R = \Z$, the ring of integers.

\subsection*{Heisenberg group over a field}

\index{Heisenberg group! $2$@$H(\K)$ over a field $\K$}
Let $\K$ be a field, $\Gamma = H(\K)$,
and $Z$ the centre of $\Gamma$.
Fix a unitary character $\psi \in \widehat Z \approx \widehat \K$;
let $I_\psi$ be defined as before Lemma \ref{Lem-IrredHeisRing}.
Two cases may occur:
\begin{enumerate}
\item[$\bullet$]
If $\psi = 1_Z$, then $I_\psi = \K$ and hence $N_\psi = \Gamma$.
\item[$\bullet$]
If $\psi \ne 1_Z$, then $I_\psi = \{0\}$ and hence $N_\psi$
coincides with the maximal abelian subgroup 
$$
N \, = \, \{(0, b, c) \mid b, c \in \K\} \approx \K^2.
$$
\end{enumerate}
Corollary~\ref{Cor-IrredHeisRingRestated} takes therefore the following form.

\begin{cor}
% 3.B.11
\label{Cor-IrredHeisField}
Let $\K$ be a field, $\Gamma = H(\K)$, with centre $Z$,
and $p \, \colon \Gamma \twoheadrightarrow \Gamma/Z$ the canonical projection.
Let $N$ be the maximal abelian subgroup
$\{(0,b, c) \in \Gamma \mid b,c \in \K\}$.
For $\psi \in \widehat Z$, set
$\widehat N (\psi) = \{ \chi \in \widehat N \mid \chi \vert_Z = \psi \}$.
Define
$$
X \, := \, \widehat{\Gamma/Z} \hskip.2cm \bigsqcup \hskip.2cm
\Big( \bigsqcup_{\psi \in \widehat Z \smallsetminus \{1_Z\}}
\{ \psi \} \times \widehat N (\psi) \Big) .
$$
\par

Then $\Ind_{N}^\Gamma \chi$ is irreducible
for all $(\psi, \chi) \in \bigsqcup_{\psi \in \widehat Z \smallsetminus \{1_Z\}}
\{ \psi \} \times \widehat N (\psi)$.
The map $X \to \widehat \Gamma$, defined by
$$
\chi \, \mapsto \, \chi \circ p
\hskip.5cm \text{for} \hskip.2cm
\chi \in \widehat{\Gamma/Z}, 
$$
and
$$
(\psi, \chi) \, \mapsto \, \Ind_{N}^\Gamma \chi
\hskip.5cm \text{for} \hskip.2cm
\psi \in \widehat Z \smallsetminus \{1_Z\}, \hskip.1cm \chi \in \widehat N_\psi
$$
factorizes to an injective map from the space of $\Gamma$-orbits in $X$
to the dual $\widehat \Gamma$ of $\Gamma$.
\end{cor}

\begin{rem} 
% 3.B.12
Let $(\psi, \chi) \in \bigsqcup_{\psi \in \widehat Z \smallsetminus \{1_Z\}}
\{ \psi \} \times \widehat N (\psi)$.
The irreducible representation $\Ind_{N}^\Gamma \chi$
is finite dimensional if $\K$ (and therefore $\Gamma$) is finite,
and infinite dimensional if $\K$ is infinite.
\par

When $\K$ is finite,
the map $X \to \widehat \Gamma$ in Corollary~\ref{Cor-IrredHeisField}
is surjective, onto $\widehat{H(\K)}$;
see Corollary \ref{Cor-FinDimRepHeis-finiteFields}.
When $\K$ is infinite, this is far from being true;
see Remark \ref{Rem-MoreIrredRep}(1).
\end{rem}

\begin{rem}
% 3.B.13
\label{dualsoffields}
Fields, considered as discrete abelian groups for the addition,
have duals which are compact abelian groups;
these duals appear above in Corollary \ref{Cor-IrredHeisField},
since $Z \approx \K$,
and again below,
in Remarks \ref{Rem-IrredRepAffine}
and \ref{Rem-MoreIrredRep}(2) concerning $\Aff (\K)$.
% and in Lemma \ref{Lem-FreeActDualGL_n} 8.E.1
They can be described as follows.
\par

The dual $\widehat \Q$ of $\Q$
is isomorphic to the quotient $\mathbf{A}_\Q / \Q$ 
of the rational adele ring by its cocompact discrete subgroup $\Q$.
There is a nice exposition of this fact in \cite{Conr}.
\par

Let $\K$ be a field of characteristic $0$, viewed as a discrete abelian group.
Since $\K$ is a vector space over $\Q$,
it is a direct sum of copies of $\Q$,
hence its dual $\widehat \K$ is isomorphic to a direct product
of copies of $\mathbf{A}_\Q / \Q$.
\par

Let $p$ be a prime.
The dual of the finite field $\F_p$ of order $p$ is
isomorphic (non-canonically) to $\F_p$ itself.
\par

Let $\K$ be a field of characteristic $p$, viewed as a discrete abelian group.
Since $\K$ is a direct sum of copies of $\F_p$,
the dual $\widehat \K$ is isomorphic to a direct product of copies of $\F_p$.
% Voir aussi \cite[25.4, 25.5 \& 25.6]{HeRo--63}.
\end{rem}

\subsection*{Heisenberg group over the integers}

\index{Heisenberg group! $3$@$H(\Z)$}
Let $\Gamma = H(\Z)$ be the Heisenberg group over $\Z$.
The irreducible representations of this group which appear below
can also be found in \cite[Pages 210--218]{Foll--16}.
Further irreducible representations of $H(\Z)$
will be constructed in Chapter~\ref{Chap-AllIrrRed}
(see Remark~\ref{Rem-MoreIrredRep}).
\par

The centre of $\Gamma$ is 
$Z = \{ (0, 0, c) \in \Gamma \mid c \in \Z \} \approx \Z$
and we have
$$
\widehat Z \, = \,
\{\psi_\theta \mid \theta \in \mathopen[ 0,1 \mathclose[ \}
\, \approx \, \T ,
$$
where, for $\theta \in \mathopen[ 0,1 \mathclose[$, the character $\psi_\theta$ is defined by 
$$
\psi_\theta (0,0, c) \, = \, e^{2 \pi i \theta c}
\hskip.5cm \text{for all} \hskip.2cm
(0, 0, c) \in Z .
$$
We denote by $\widehat Z_\infty$ the subset of $\widehat Z$
of elements of infinite order, i.e., of faithful characters of $\Z$, i.e., 
$$
\widehat Z_\infty \, = \, \{ \psi_\theta \mid \theta \in \mathopen] 0,1 \mathclose]
\hskip.2cm \text{and} \hskip.2cm
\theta \notin \Q \} .
$$
\par

Fix $\psi = \psi_\theta \in \widehat Z_\infty$.
Since $\psi$ is faithful, the corresponding ideal is $I_{\psi_\theta} = \{0\}$
(notation as before Theorem~\ref{Theo-IrredHeisRing}),
so that $N_\psi$ coincides with
the maximal abelian subgroup $N = \{(0, b, c) \mid b, c \in \Z\}$. 
Define $\widehat N (\psi)$ as in Theorem \ref{Theo-IrredTwoStepNil},
i.e., $\widehat N (\psi) = \{ \chi \in \widehat N \mid \chi \vert_{Z} = \psi \}$.
Then 
$$
\widehat N (\psi) \, = \, 
\{ \chi_{ \theta, \beta} \mid \beta \in \mathopen[ 0,1 \mathclose[ \} 
\, \approx \, \T,
$$
where, for $\theta, \beta \in \mathopen[ 0,1 \mathclose[$,
the character $\chi_{ \theta, \beta} \in \widehat N$ is defined by 
$$
\chi_{\theta, \beta} (0,b, c) \, = \, e^{2 \pi i (\beta b + \theta c)}
\hskip.5cm \text{for all} \hskip.2cm
(0, b, c) \in N .
$$
\par

Consider $\chi = \chi_{ \theta, \beta } \in \widehat N (\psi)$. 
For every $\gamma = (a,b,c) \in \Gamma$
and $\gamma' = (0,b', c') \in N$, we have
$$
\begin{aligned}
\chi^\gamma(\gamma') 
\, &= \, \chi(\gamma \gamma' \gamma^{-1})
\, = \, \chi(\gamma') \chi([\gamma', \gamma^{-1}])
\, = \, \chi(\gamma') \psi([\gamma', \gamma])
\\
\, &= \, e^{ 2 \pi i (\beta b' + \theta c') } e^{ 2 \pi i \theta (ab') }
\, = \, e^{ 2 \pi i [ (\beta + a \theta)b' + \theta c')] }
\\
\, &= \, \chi_{\theta, \beta + a \theta} (0,b', c').
\end{aligned}
$$
Therefore the $\Gamma$-orbit of $\chi_{\theta, \beta}$ is 
$$
\{ \chi_{\theta, \beta + a \theta} \in \widehat N \mid a \in \Z\}.
$$
Observe that $\{\beta + a \theta \pmod 1 \mid a \in \Z\}$
is the orbit of $\beta$ under the irrational rotation $R_\theta \, \colon \T \to \T$.
\par

The following result is a direct consequence of Corollary~\ref{Cor-IrredHeisRingRestated}.

\begin{cor}
% 3.B.14
\label{Cor-IrredRepHeisInteger}
Consider the Heisenberg group $\Gamma = H(\Z)$ over $\Z$,
its centre $Z = \{(0 ,0, c) \in \Gamma \mid c \in \Z\}$,
and its maximal abelian subgroup $N = \{(0,b, c) \in \Gamma \mid b,c \in \Z\}$. 
For each $\theta \in \mathopen[ 0,1 \mathclose[$ irrational,
consider as above
\par
the unitary character $\psi_\theta \in \widehat Z_\infty$,
\par
the subset
$\widehat N (\psi_\theta)
= \{\chi_{ \theta, \beta} \in \widehat N \mid \beta \in \mathopen[ 0,1 \mathclose[ \}
\approx \T$
of the character group $\widehat N$,
\par
the irrational rotation 
$R_\theta \, \colon \mathopen[ 0,1 \mathclose[ \to \mathopen[ 0,1 \mathclose[$,
$\beta \mapsto \beta + \theta \pmod 1$,
\par
and the orbit space $\mathcal O_\theta = \mathopen[ 0,1 \mathclose[ /R_\theta$.
\par\noindent
Then:
\begin{enumerate}[label=(\arabic*)]
\item\label{iDECor-IrredRepHeisInteger}
For each $\beta \in \mathopen[ 0,1 \mathclose[$
and irrational $\theta \in \mathopen[ 0,1 \mathclose[$,
the representation $\Ind_N^\Gamma \chi_{\theta, \beta}$ of $\Gamma$
is irreducible, infinite dimensional, with central character $\psi_\theta$.
\item\label{iiDECor-IrredRepHeisInteger}
For each irrational $\theta \in \mathopen[ 0,1 \mathclose[$,
the map
$\mathopen[ 0,1 \mathclose[ \mapsto \widehat \Gamma,
\hskip.2cm 
\beta \mapsto \Ind_N^\Gamma \chi_{ \theta, \beta}$
factorizes to an injective map 
from $\mathcal O_\theta$ into $\widehat \Gamma$.
\item\label{iiiDECor-IrredRepHeisInteger}
For $\theta, \theta' \in \mathopen[ 0,1 \mathclose[$
irrational and distinct,
the images of the corresponding maps in \ref{iiDECor-IrredRepHeisInteger}
are disjoint, so that these maps provide an injective map
$$
\bigsqcup_{\theta \in \mathopen[ 0,1 \mathclose[, \theta \notin \Q} 
\mathcal O_\theta \, \hookrightarrow \, \widehat \Gamma .
$$
\end{enumerate}
\end{cor}

\begin{rem}
% 3.B.15
\label{Rem-Cor-IrredHeisIntegers}
(1)
For an irrational $\theta \in \mathopen[ 0,1 \mathclose[$
and for $\beta \in \mathopen[ 0,1 \mathclose[$, 
let $\pi_{\theta, \beta} := \Ind_N^\Gamma \chi_{ \theta, \beta}$
be the irreducible representation 
of $\Gamma$ of Corollary \ref{Cor-IrredRepHeisInteger}. 
\par

We choose $H = \{(a, 0, 0) \in \Gamma \mid a \in \Z\} \approx \Z$
as transversal for $N \backslash \Gamma$.
Construction~\ref{constructionInd}(2), of which we keep the notation here,
shows that the induced representation $\pi_{\theta, \beta}$
can be realized on the Hilbert space $\ell^2(\Z)$ as follows.
For $t = (k, 0, 0)$ in the transversal and $g = (a, b, c) \in \Gamma$,
we have
$$
(k, 0, 0) (a, b, c) \, = \, (k + a, b, c + kb) \, = \,
(0, b, c + kb) (k + a, 0, 0) ,
$$
hence
$$
\alpha(t, g) \, = \, (0, b, c + kb)
\hskip.5cm \text{and} \hskip.5cm
t \cdot g \, = \, (k + a, 0, 0) .
$$
Thus
$$
(\pi_{ \theta, \beta}(a,b,c) f) (k) \,= \, \chi_{\theta, \beta}(0, c, c + kb) f(a + k)
\, = \, e^{2 \pi i \theta c} e^{2 \pi i (\beta + k\theta) b} f(a + k)
$$
for $f \in \ell^2(\Z)$, $(a, b, c) \in \Gamma$, and $k \in \Z$. 

\vskip.2cm

(2)
The map of \ref{Cor-IrredRepHeisInteger}~\ref{iiiDECor-IrredRepHeisInteger}
above is far from being surjective.
See Remarks~\ref{Rem-MoreIrredRep}, and also \cite{Brow--73a}.

\vskip.2cm

(3)
Given a unitary character $\psi_\theta$
of the centre $Z$ of $\Gamma = H(\Z)$ with irrational $\theta$,
Corollary \ref{Cor-IrredRepHeisInteger} shows
that every $\Gamma$-orbit in $\widehat N (\psi_\theta)$ provides 
one equivalence class of irreducible representations of $\Gamma$ 
with central character $\psi_\theta$. 
Moreover, these orbits are in a one-to-one correspondence 
with the orbits of the irrational rotation $R_\theta$ on $\T$. 
Since there are uncountably many $R_\theta$-orbits,
there are uncountably many non-equivalent irreducible representations
of $\Gamma$ with central character $\psi_\theta$.
\end{rem}

We will show in Subsection~\ref{SS:BorelStructureQD} that
the dual space $\widehat \Gamma$ of $\Gamma$
is equipped with a natural Borel structure, called the Mackey--Borel structure.
We now deduce from Corollary~\ref{Cor-IrredRepHeisInteger}
that the Borel space $\widehat \Gamma$ is not countably separated.
This is a manifestation of the fact that $\Gamma$ is not type I 
(see Section~\ref{SectionTypeI} for this notion)
and is a particular instance of Glimm theorem 0f Chapter~\ref{ChapterAlgLCgroup}.

\begin{cor}
% 3.B.16
\label{Cor-Cor-IrredRepHeisInteger}
Let $\Gamma = H(\Z)$ be the Heisenberg group over $\Z$.
\par

The Borel space $\widehat \Gamma$ is not countably separated.
\end{cor}

\begin{proof}
Fix an irrational number $\theta \in \mathopen[ 0,1 \mathclose[$.
It is well-known that the rotation $R_\theta$ is ergodic
with respect to the Lebesgue measure $\mu$ on $\mathopen[ 0,1 \mathclose[$
(see for instance \cite[Chapter I, Example 1.4.ii]{BeMa--00}).
Since $\mu$ is not supported by an $R_\theta$-orbit,
this implies that the orbit space $\mathcal O_\theta = \mathopen[ 0,1 \mathclose[ /R_\theta$
is not countably separated for its canonical quotient Borel structure
(see \cite[Proposition 2.1.10]{Zimm--84}).
In order to show that $\widehat \Gamma$ is not countably separated,
it suffices therefore to show that the map 
$$
\mathopen[ 0,1 \mathclose[ \, \to \, \widehat \Gamma,
\hskip.2cm
\beta \, \mapsto \, \Ind_N^\Gamma \chi_{ \theta, \beta}
$$
of Corollary~\ref{Cor-IrredRepHeisInteger}~\ref{iiDECor-IrredRepHeisInteger} is measurable. 
\par

As shown in Remark~\ref{Rem-Cor-IrredHeisIntegers},
the representations $\pi_{\theta, \beta} = \Ind_N^\Gamma \chi_{ \theta, \beta}$
are realized in $\ell^2(\Z)$, for all $\beta \in \mathopen[ 0,1 \mathclose[$.
By definition of the Mackey--Borel structure on $\widehat \Gamma$,
see Subsection~\ref{SS:BorelStructureQD},
it suffices to prove that the functions 
$\beta \mapsto \langle \pi_{\theta, \beta}(a, b, c) f_1 \mid f_2 \rangle$
are Borel measurable,
for all $(a, b, c) \in \Gamma$ and $f_1, f_2 \in \ell^2(\Z)$.
\par

Let $(a, b, c) \in \Gamma$ and $f_1, f_2 \in \ell^2(\Z)$.
As measurability is inherited by linear combinations and pointwise limits,
we can assume that $f_1 = \delta_{k_1}$ and $f_2 = \delta_{k_2}$ for $k_1, k_2 \in \Z$.
In this case, by the formula of Remark~\ref{Rem-Cor-IrredHeisIntegers}
for $\pi_{ \theta, \beta}$, we have
$$
\langle \pi_{\theta, \beta}(a,b,c) \delta_{k_1} \mid \delta_{k_2} \rangle
\, = \,
\begin{cases}
e^{2 \pi i \theta c} e^{2 \pi i (\beta + (k_1 - a) \theta) b}
& \text{if} \hskip.2cm k_2 = k_1-a
\\
0
& \text{if} \hskip.2cm k_2 \ne k_1-a.
\end{cases}
$$
Clearly, this shows that
$\beta \mapsto \langle \pi_{\theta, \beta}(a, b, c) \delta_{k_1} \mid \delta_{k_2} \rangle$
is measurable.
\end{proof}

We now turn to the construction of irreducible representations of 
$\Gamma = H(\Z)$ with a non faithful central character. 
\par

For this, fix a character $\psi = \psi_{\theta} \in \widehat Z$
for $\theta \in \mathopen[ 0,1 \mathclose[$ \textbf{rational}.
Write
$$
\theta \, = \, \frac{p}{n}
\hskip.5cm \text{where} \hskip.5cm
n \in \N^*, \hskip.1cm p \in \{0, 1, \hdots, n - 1 \}
\hskip.2cm \text{and} \hskip.2cm
(p, n) = 1 ,
$$
so that $\psi$ is of order $n$.
Then (with the notation as before Theorem~\ref{Theo-IrredHeisRing}), we have 
$$
I_{\psi_{\theta}} \, = \, \ker \psi_\theta \, = \, n\Z ,
$$
and $N_\psi$ coincides with the normal subgroup 
$$
\Gamma(n) \, = \,
\{(a,b, c) \in \Gamma \mid a \in n\Z, \hskip.1cm b \in \Z, \hskip.1cm c \in \Z\}
$$
of $\Gamma$.
\par

We parametrize the dual of $n\Z$ by $\mathopen[ 0,1/n \mathclose[$:
every element of $\widehat{n\Z}$ is of the form $n\Z \ni a \mapsto e^{2 \pi i \alpha a} \in \T$ 
for a unique $\alpha \in \mathopen[ 0,1/n \mathclose[$. 
\par 

By Lemma~\ref{Lem-IrredHeisRing2},
every unitary character $\chi$ of $\Gamma(n)$ such that $\chi \vert_Z = \psi_\theta$
is of the form $\chi_{\theta, \alpha, \beta}$
for a unique pair $(\alpha, \beta) \in \mathopen[ 0,1/n \mathclose[ \times \mathopen[ 0,1 \mathclose[$, where 
$$
\chi_{\theta, \alpha, \beta}(a,b,c) \, = \, e^{2 \pi i (\alpha a + \beta b + \theta c)} 
\hskip.5cm \text{for} \hskip.2cm
a \in n\Z, \hskip.1cm b \in \Z, \hskip.1cm c \in \Z.
$$
For every $\gamma' = (a',b',c') \in \Gamma$
and $\gamma = (a, b, c) \in \Gamma(n)$, we have $\theta a b' \in \Z$ and
$$
\begin{aligned}
\chi_{\theta, \alpha, \beta}^{\gamma'}(\gamma)
\, &= \chi_{\theta, \alpha, \beta}(\gamma) \psi_{\theta}([\gamma, (\gamma')^{-1}]
\\
\, &= \, e^{2 \pi i (\alpha a + \beta b + \theta c)} e^{2 \pi i \theta (a'b - ab')}
\\
\, &= \, e^{2 \pi i (\alpha a + \beta b + \theta c)} e^{2 \pi i \theta (a'b)}
\\
\, &= \, e^{2 \pi i [\alpha a + (\beta + a' \theta)b + \theta c)]}
\\
\, &= \, \chi_{ \theta, \alpha, \beta + a' \theta} (\gamma) .
\end{aligned}
$$
Therefore the $\Gamma$-orbit of $\chi_{\theta, \alpha, \beta}$ is 
$$
\{ \chi_{ \theta, \alpha, \beta + a' \theta} \in \widehat N \mid a' \in \Z \}
\, = \,
\{ \chi_{ \theta, \alpha, \beta + \frac{k}{n}} \in \widehat N \mid k \in \Z \} 
$$
and so contains a character $\chi_{ \theta, \alpha, \beta'}$ for a unique pair $(\alpha, \beta')$
in $\mathopen[ 0,1/n \mathclose[ \times \mathopen[ 0,1/n \mathclose[$.
\par

In view of the previous remarks,
the following result is a direct consequence of
Corollary~\ref{Cor-IrredHeisRingRestated}.

\begin{cor}
% 3.B.17
\label{Cor-IrredFiniHeisIntegers}
Let $n$ be a positive integer,
$\psi = \psi_{\theta}$ a unitary character of $Z$ of order $n$,
and $(\alpha, \beta) \in
\mathopen[ 0,1/n \mathclose[ \times \mathopen[ 0,1/n\mathclose[$.
Let $\Gamma(n)$ and $\chi_{\theta, \alpha, \beta}$ be defined as above.
\par

The induced representation 
$$
\pi_{\theta, \alpha, \beta} \, := \,
\Ind_{\Gamma(n)}^\Gamma \chi_{\theta, \alpha, \beta}
$$
is an irreducible representation of $\Gamma=H(\Z)$ of dimension $n$.
\par

Moreover, for another triple
$\theta' = \frac{p'}{n'}$ and $(\alpha', \beta') \in
\mathopen[ 0,1/n' \mathclose[ \times \mathopen[ 0,1/n' \mathclose[$
of the same kind,
the representations
$\pi_{\theta, \alpha, \beta}$ and $\pi_{\theta', \alpha', \beta'}$ 
are equivalent
if and only if $(\theta, \alpha, \beta) = (\theta', \alpha', \beta')$.
\end{cor}

\begin{rem}
% 3.B.18
\label{Rem-IrredFiniHeisIntegers} 
We will see later (Corollary~\ref{Cor-FinDimRepHeis-Integers}) that 
\textbf{every} finite dimensional irreducible representation of $\Gamma=H(\Z)$
is equivalent to one of the representations $\pi_{\theta, \alpha, \beta}$ 
of Corollary~\ref{Cor-IrredFiniHeisIntegers}.
\end{rem}

\begin{rem}
% 3.B.19
\label{dualHeisLie}
\index{Heisenberg group! $4$@$H(\R)$}
Consider the real Heisenberg group $H(\R)$,
now with its natural connected and locally compact topology.
Up to isomorphism, it is the only non-abelian nilpotent
connected simply connected real Lie group of dimension~$3$.
In $H(\R)$, we consider
the abelian normal closed subgroup
$$
N \, = \,
\begin{pmatrix} 1 & 0 & \R \\ 0 & 1 & \R \\ 0 & 0 & 1 \end{pmatrix}
\, \subset \, 
H(\R) \, = \,
\begin{pmatrix} 1 & \R & \R \\ 0 & 1 & \R \\ 0 & 0 & 1 \end{pmatrix} .
$$
Write again $(x, y, z)$
for $\begin{pmatrix} 1 & x & z \\ 0 & 1 & y \\ 0 & 0 & 1 \end{pmatrix}$.
The unitary characters of $H(\R)$ and $N$ are
$$
\begin{aligned}
&\chi_{u, v} \, \colon \, H(\R) \to \T , \hskip.2cm &(x, y, z) \mapsto w^{2 \pi i (ux + vy)}
\hskip.5cm \text{for} \hskip.2cm
u, v \in \R ,
\\
&\chi_{v,w} \, \colon \, N \to \T , \hskip.2cm &(0, y, z) \mapsto w^{2 \pi i (vy + wz)}
\hskip.5cm \text{for} \hskip.2cm
v, w \in \R .
\end{aligned}
$$
\par

It is a version of the Stone--von Neumann theorem
\cite{Ston--30, vNeu--31}
% \cite{Mack--49}
that the dual of $H(\R)$ can be seen as a disjoint union:
$$
\widehat{H(\R)} \, = \, \R^2 \sqcup \R^\times ,
$$
where
\begin{enumerate}
\item[$\cdot$]
$\R^2$ is identified with the space of unitary characters of $H(\R)$, i.e., 
the dual of the abelian group $H(\R)_{\rm ab} = H(\R)/[H(\R),H(\R)] \approx \R^2$;
\item[$\cdot$]
$\R^\times$ is identified with
the space of equivalence classes of infinite-dimensional representations
$\Ind_N^{H(\R)} \chi_{0,w}$, for $w \in \R^\times$.
\end{enumerate}
For the Fell topology, the induced topology 
on either $\R^2$ or $\R^\times$ is the usual one,
but a sequence of points in $\R^\times$ that tends to $0$ in the usual sense
tends to all points of $\R^2$ in $\R^2 \sqcup \R^\times = \widehat{H(\R)}$.
\par

One method to obtain these results
is to use the Mackey machine.
The group $H(\R)$ acts naturally on the dual of $N$, that is on
$$
\widehat N \, = \, \{ \chi_{v,w} \mid (v,w) \in \R^2 \} .
$$
The set of orbits splits in two parts:
\begin{enumerate}
\item[$\cdot$]
the set of one-point orbits, i.e., of orbits $\{\chi_{v,0}\}$, where $v \in \R$,
\item[$\cdot$]
the set of line-orbits $\mathcal O_w = \{\chi_{t,w} \mid t \in \R \}$, where $w \in \R^\times$.
\end{enumerate}
The one-point orbits provide the unitary characters of $H(\R)$,
and the line-orbits provide 
the infinite-dimensional irreducible representations of $H(\R)$.
Details can be found in \cite{Fell--62}; see also \cite[Section 6.7]{Foll--16}.

\vskip.2cm 

More generally, let $G$ be a nilpotent, simply connected, connected, real Lie group
with Lie algebra $\mathfrak g$; 
the natural action of $G$ on the dual space $\mathfrak g^*$ 
is called the coadjoint representation.
Kirillov has constructed a map
% \marginpar{More on Kirillov orbit method, to come}
from the space of co-adjoint orbits $G \backslash \mathfrak g^*$
to the dual $\widehat G$, 
has shown that it is a continuous bijection,
and has conjectured that it is a homeomorphism \cite{Kiri--62};
the conjecture has been established by Brown \cite{Brow--73b}.
This has been generalized to soluble, simply connected, connected, real Lie groups
for which the exponential map is surjective \cite{LeLu--94}.
Indeed, the orbit method has shown to be a fruitful inspiration
for very large classes of groups and various problems \cite{Kiri--99}.
\end{rem}

\section
{The affine group of a field}
% Section 3.C
\label{Section-IrrRepAff}

\index{Affine group! $4$@$\Aff(\K)$ of a field $\K$} 
Let $\K$ be field
and let $\Aff(\K)$ be the \textbf{affine group} over $\K$,
that is, the group of affine transformations of $\K$;
in matrix form:
$$
\begin{pmatrix} \K^\times & \K \\ 0 & 1 \end{pmatrix}
\, \approx \, \K^\times \ltimes \K .
$$
The semi-direct product refers to the left action $\K^\times \curvearrowright \K$
defined by $(a, t) \mapsto at$ for all $a \in \K^\times$ and $t \in \K$.
The derived group and the abelianized group,
$$
[\Aff(\K), \Aff(\K)] \, = \, 
\begin{pmatrix} 1 & \K \\ 0 & 1 \end{pmatrix}
\, \approx \, \K
\hskip.5cm \text{and} \hskip.5cm
\Aff(\K) / [\Aff(\K), \Aff(\K)]
\, \approx \, \K^\times ,
$$
are both abelian, and we have a short exact sequence
$$
\K \, \approx \, \begin{pmatrix} 1 & \K \\ 0 & 1 \end{pmatrix}
 \, \hookrightarrow \,
\Aff(\K) = \begin{pmatrix} \K^\times & \K \\ 0 & 1 \end{pmatrix}
\, \twoheadrightarrow \,
\Aff(\K) / [\Aff(\K), \Aff(\K)] \approx \K^\times .
$$
\par

The left action $\K^\times \curvearrowright \K$
extends to a left action $\Aff(\K) \curvearrowright \K$
defined by $\left( \begin{pmatrix} a & s \\ 0 & 1 \end{pmatrix} , t \right) \mapsto at$,
with kernel $[\Aff(\K) , \Aff(\K)]$.
By duality, there is a corresponding right action 
$\widehat \K \curvearrowleft \K^\times, \hskip.1cm (\chi, a) \mapsto \chi^a$,
defined by
$$
\chi^a(t) \, = \, \chi(at)
\hskip.5cm \text{for all} \hskip.2cm 
\chi \in \widehat \K , \hskip.1cm a \in \K^\times , t \in \K .
$$
As above, it extends to a right action $\widehat \K \curvearrowleft \Aff(\K)$
with kernel $[\Aff(\K) , \Aff(\K)]$.
\par

For convenience until the end of this section,
we often write $\Gamma$ for the group $\Aff(\K)$, and $N$ for its derived group.
We identify $N$ with $\K$ and $\Gamma/N$ with $\K^\times$.

\begin{lem}
% 3.C.1
\label{Lem-FreeActDualAff}
We keep the notation above.
\begin{enumerate}[label=(\arabic*)]
\item\label{iDELem-FreeActDualAff}
The group $\Gamma/N$ acts freely on $\widehat N \smallsetminus \{1_N\}$.
\item\label{iiDELem-FreeActDualAff}
If $\K$ is finite, the action of $\Gamma$ on $\widehat N \smallsetminus \{1_N\}$
is transitive.
\item\label{iiiDELem-FreeActDualAff}
If $\K$ is infinite, the $\Gamma$-orbit of every $\chi \in \widehat N \smallsetminus \{1_N\}$ 
is dense in $\widehat N$.
\end{enumerate}
\end{lem}

\begin{proof}
\ref{iDELem-FreeActDualAff} 
Let $\chi \in \widehat N$, $\chi \ne 1_N$, and $a \in \K^\times$.
Since $\chi^a(t) = \chi(at)$ or all $t \in \K$,
we have $\chi^a = \chi$ if and only $(a - 1)t \in \ker (\chi)$ for all $t \in \K$.
Since $\ker (\chi) \ne \K$, this holds if and only if $a = 1$.
Therefore $\Gamma/N$ acts freely on $\widehat N \smallsetminus \{1_N\}$.

\vskip.2cm

\ref{iiDELem-FreeActDualAff}
Let $\chi \in \widehat N$, $\chi \ne 1_N$.
If $\K$ is finite of cardinality $q$, 
it follows from \ref{iDELem-FreeActDualAff} that the $\Gamma$-orbit of $\chi$
has cardinality $q-1$
and the action of $\Gamma$ on $\widehat N \smallsetminus \{1_N\}$
is therefore transitive.

\vskip.2cm

\ref{iiiDELem-FreeActDualAff}
Write $\chi^0$ for $1_N$,
so that $\chi^a$ is now defined for all $a \in \K$.
Then
$$
H \, := \, \{ \rho \in \widehat N \mid \rho = \chi^a 
\hskip.2cm \text{for some} \hskip.2cm 
a \in \K \}
$$
is a subgroup of $\widehat N$.
Indeed, we have, for all $a,b,t \in \K$,
$$
(\chi^a \chi^b) (t) \, = \, \chi(at) \chi(bt) \, = \, \chi ((a + b)t) \, = \, \chi^{a + b}(t).
$$
\par

Assume, by contradiction, that the $\Gamma$-orbit 
$$
\mathcal O \, := \, \{ \rho \in \widehat N \mid \rho = \chi^a 
\hskip.2cm \text{for some} \hskip.2cm 
a \in \K^\times \}
$$
is not dense in $\widehat N$. 
We have $H = \mathcal O \cup \{1_N\}$;
since $N = \K$ is infinite, 
its Pontrjagin dual $\widehat N$ is not discrete.
It follows that $H$ is not dense in $\widehat N$,
so that the quotient group $\widehat N / \overline H$ 
is a nontrivial LCA group.
By Pontrjagin duality 
% (see \cite[Theorem 24.2]{HeRo--63}), 
there exists $t_0 \in \K \smallsetminus \{0\}$ 
such that $\chi^a(t_0) = 1$ for all $\chi^a \in H$, 
that is, such that $\chi(at_0) = 1$ for all $a \in \K$.
This implies that $\chi(t) = 1$ for all $t \in \K$, 
which contradicts the fact that $\chi \ne 1_N$.
\end{proof}

\begin{theorem}
% 3.C.2
\label{Prop-IrredRepAffine}
Let $\Gamma = \Aff(\K)$ the group of affine transformations of a field $\K$
and $N \approx \K$ the derived group of $\Gamma$, as above.
\begin{enumerate}[label=(\arabic*)]
\item\label{iDEProp-IrredRepAffine}
For every $\chi \in \widehat \K \smallsetminus \{1_\K\}$,
the induced representation $\Ind_{N}^\Gamma \chi$ is irreducible.
\item\label{iiDEProp-IrredRepAffine}
Let $\chi_1, \chi_2 \in \widehat \K \smallsetminus \{1_\K\}$ be 
such that $\chi_1^\gamma \ne \chi_2$ for every $\gamma \in \Gamma$. 
Then $\Ind_{N}^\Gamma \chi_1$ and $\Ind_{N}^\Gamma \chi_2$
are non-equivalent irreducible representations of $\Gamma$.
\end{enumerate}
\end{theorem}

\noindent
Concerning \ref{iiDEProp-IrredRepAffine} above,
observe that, for $\chi \in \widehat K$ and $\gamma = (b,t) \in \Gamma$,
we have $\chi^\gamma = \chi^b$.

\begin{proof}
Claim \ref{iDEProp-IrredRepAffine} follows from
the Mackey--Shoda criterion for irreducibility of monomial representations
(Corollary~\ref{Cor-NormalSubgIrr})
and from Lemma~\ref{Lem-FreeActDualAff}~\ref{iDELem-FreeActDualAff}.
Claim \ref{iiDEProp-IrredRepAffine} follows from
the Mackey--Shoda criterion for equivalence of monomial representations
(Corollary \ref{Cor-NormalSubgEquiv}).
\end{proof}

Theorem \ref{Prop-IrredRepAffine} can be restated as follows.

\begin{prop}
% 3.C.3
\label{Prop-IrredRepAffine-bis}
Let $\Gamma = \Aff(\K)$ and $N \approx \K$
be as in Theorem \ref{Prop-IrredRepAffine}.
The map 
$$
\widehat \K \to \widehat \Gamma, \hskip.2cm \chi \mapsto 
\begin{cases}
\Ind_N^\Gamma \chi
&\text{if} \hskip.2cm \chi \in \widehat \K \smallsetminus \{1_\K\} 
\\
1_\Gamma
& \text{if} \hskip.2cm \chi =1_\K 
\end{cases}
$$
factorizes to an injective map from the space of $\Gamma$-orbits in $\widehat \K$
to the dual $\widehat \Gamma$ of $\Gamma$.
\end{prop}

\begin{rem}
% 3.C.4
\label{Rem-IrredRepAffine-Fini}
Let $\chi \in \widehat \K\smallsetminus \{1_\K\}$.
Following Construction \ref{constructionInd}(2),
the induced representation $\Ind_N^\Gamma \chi$
can be realized as follows.
The subgroup
$$
\left\{ \begin{pmatrix} a & 0 \\ 0 & 1 \end{pmatrix} \in \Aff(\K)
\hskip.1cm \Big\vert \hskip.1cm
a \in \K^\times \right\}
\, \approx \, \K^\times
$$
is a right transversal for $N$ in $\Gamma$, 
i.e.,
$$
\Gamma \, = \, \bigsqcup_{a \in \K^\times}
N \begin{pmatrix} a & 0 \\ 0 & 1 \end{pmatrix} .
$$
Therefore the induced representation $\pi_\chi := \Ind_N^\Gamma \chi$
acts on $\ell^2(\K^\times)$, and is given by
$$
(\pi_{\chi}(a,b) f) (t) \, = \, \chi(tb) f(at) 
\hskip.5cm \text{for all} \hskip.2cm
(a,b) \in \Aff(\K), \hskip.1cm
f \in \ell^2(\K^\times), 
\hskip.2cm \text{and} \hskip.2cm
t \in \K^\times .
$$
Note that, when $\K$ is finite, $\pi_\chi$ is the unique
irreducible representation of $\Aff(\K)$ of dimension greater than one
up to equivalence (Corollary~\ref{Cor-FinDimRepAffine}).
\end{rem}

\begin{rem}
% 3.C.5
\label{Rem-IrredRepAffine} 
Theorem \ref{Prop-IrredRepAffine} shows that 
every $\K^\times$-orbit in $\widehat \K$ provides 
one equivalence class of irreducible representations 
of $\Gamma = \Aff(\K)$.
\par

Assume from now on that $\K$ is infinite.
Let $\widehat \K$ be equipped with its normalized Haar measure $\mu$.
By uniqueness of $\mu$ as translation-invariant probability measure,
the action of $\K^\times$ on $\widehat \K$ preserves $\mu$.
Moreover this action is ergodic, 
since $\{0\}$ is the only finite orbit
for the action of $\K^\times$ on $\widehat{\widehat \K} \approx \K$
(see \cite[Chapter I, Proposition 1.5]{BeMa--00}).
\par

As in Corollary~\ref{Cor-Cor-IrredRepHeisInteger},
we conclude that the Borel space
$\widehat \Gamma$ is not countably separated.
This is another illustration of Glimm theorem of Chapter~\ref{ChapterAlgLCgroup}.
\end{rem}

\begin{rem}
% 3.C.6
\label{AffKLie}
Let $\K$ be a topological field.
The affine group $\Aff(\K) = \K^\times \ltimes \K$
has a topology as a subspace of the product $\K \times \K$
which makes it a topological group.
When $\K$ is moreover a local field, $\Aff(\K)$ is a
second-countable locally compact group,
of which the dual can be described as follows.
The results can be obtained by using the Mackey machine.

\vskip.2cm

\index{$l 5$@$\R_+ = \mathopen[ 0, \infty \mathclose[$
non-negative real numbers}
\index{$l 6$@$\R^\times_+ = \mathopen] 0, \infty \mathclose[$
positive real numbers}
\index{Affine group! $2$@$\Aff(\R)_0$, orientation preserving group}
Consider first the \textbf{orientation preserving affine group}
of the real line 
$$
G = \Aff(\R)_0 = \begin{pmatrix} \R^\times_+ & \R \\ 0 & 1 \end{pmatrix} ,
$$
which is a connected locally compact group;
here $\R^\times_+ := \{a \in \R \mid a > 0 \}$.
The abelian normal closed subgroup of interest is the translation group
$N = \begin{pmatrix} 1 & \R \\ 0 & 1 \end{pmatrix} \approx \R$.
The dual of $\Aff(\R)_0$ can be described as
$$
\widehat{\Aff(\R)_0} \, = \, 
\R \sqcup \{\pi_{-1}\} \sqcup \{\pi_{1}\} \, = \, 
\{ \chi_t \}_{t \in \R} \sqcup \{\pi_{-1}\} \sqcup \{\pi_{1}\} .
$$
The three parts correspond to the three orbits
of the natural action of $\Aff(\R)_0$ 
[equivalently of $\R^\times_+$]
on $\widehat N \approx \R$:
the orbit $\{0\}$ corresponds to the family of unitary characters
$$
\chi_t \, \colon \, \Aff(\R)_0 \to \T, \hskip.2cm
\begin{pmatrix} a &b \\ 0 & 1 \end{pmatrix} \mapsto a^{it} ,
\hskip.5cm \text{where} \hskip.2cm
t \in \R ,
$$
the orbits $\mathopen] -\infty,0 \mathclose[$
and $\mathopen] 0, \infty \mathclose[$
correspond to $\pi_\varepsilon = \Ind_N^{\Aff(\R)_0} \rho_\varepsilon$,
where $\varepsilon \in \{-1,1 \}$ and
$$
\rho_\varepsilon \, \colon \, N \to \T, \hskip.2cm
\begin{pmatrix} 1 &b \\ 0 & 1 \end{pmatrix} 
\mapsto e^{\varepsilon i b} .
$$
For the Fell topology, the induced topology on $\R$ is the usual one,
the closure of $\{\pi_1 \}$ is $\R \sqcup \{\pi_1 \}$,
and the closure of $\{\pi_{-1}\}$ is $\R \sqcup \{\pi_{-1}\}$.
We refer to \cite[Page 132]{Mack--52} or \cite[Page 263]{Fell--62},
or \cite[\S~6.7]{Foll--16}. 
\par

\index{Affine group! $1$@$\Aff(\R)$}
The affine group
$\Aff(\R) = \begin{pmatrix} \R^\times & \R \\ 0 & 1 \end{pmatrix}$
has two connected components, one being $\Aff(\R) _0$.
The abelian normal closed subgroup to consider
is the same translation group $N$ as above.
The action of $\Aff(\R)$ [equivalently of $\R^\times$]
has now \emph{two} orbits.
There exists an infinite-dimensional irreducible representation
$\pi$ of $\Aff(\R)$, unique up to equivalence.
The dual of $\Aff(\R)$ can be described as
$$
\widehat{\Aff(\R)} \, = \, 
\widehat{\R^\times} \sqcup \{ \pi \} \, = \,
(\R \times \{0,1 \}) \sqcup \{ \pi \} .
$$
For $\R^\times = \{a \in \R \mid a \ne 0 \}$,
the dual $\widehat{\R^\times}$ can be identified to $\R \times \{0, -1 \}$;
a point $(t, \delta) \in \R \times \{0,1 \}$ is the equivalence class
of the one-dimensional representation
$$
\begin{pmatrix} a & b \\ 0 & 1 \end{pmatrix} 
\mapsto \vert a \vert^{it} {\rm sign}(a)^\delta .
$$
The representation $\pi$ is that induced from $\rho_{1}$
(or equivalently from $\rho_{-1}$);
equivalently, $\pi$ is the representation of $\Aff(\R)$ 
on $L^2(\R^\times, dt/\vert t \vert)$ defined by
$$
\left( \pi \begin{pmatrix} a & b \\ 0 & 1 \end{pmatrix} \right) \xi (t)
\, = \,
\exp \left( -2 \pi i bt \right) \xi(at)
$$
for all $\begin{pmatrix} a & b \\ 0 & 1 \end{pmatrix} \in \Aff(\R)$,
$\xi \in L^2(\R^\times, dt / \vert t \vert)$,
and $t \in \R^\times$.
See \cite[Page 224]{GGPS--69}, where the setting is more generally that
of a non-discrete locally compact field $\mathbf K$, rather than $\R$ as here.
For the Fell topology,
the subspace $\R \times \{0, -1 \}$ of $\widehat{\Aff(\R)}$
has its usual topology
and the closure of $\{ \pi \}$ is the whole of $\widehat{\Aff(\R)}$.
\par

It is known that the regular representation $\lambda$ of $\Aff(\R)$
is equivalent to a multiple of the infinite-dimensional irreducible representation $\pi$
\cite[Example 1, Page 510]{KlLi--72}.
% pour le groupe de Heisenberg, voir page 493
In other words, anticipating on the terminology of Chapter \ref{ChapterTypeI},
$\lambda$ is a factor representation of type I$_\infty$, 
equivalent to $\infty \pi$, and therefore weakly equivalent to $\pi$.

\vskip.2cm

A similar description holds for the dual of the affine group
$\Aff(\K) = \begin{pmatrix} \K^\times & \K \\ 0 & 1 \end{pmatrix}$ 
of a non-Archimedean locally compact field $\K$.
% \cite[F.2.6]{BeHV--08}.
Its dual is the union of the Pontrjagin dual of $\K^\times$
and \emph{one} infinite-dimensional representation~$\pi$,
dense in the dual for the Fell topology.
\end{rem}

\section
{Solvable Baumslag--Solitar groups $\BS(1, p)$}
% Section 3.4
\label{Section-IrrRepBS}

\index{Baumslag--Solitar group $\BS(1, p)$}
Let $p$ be a prime.
The corresponding \textbf{Baumslag--Solitar group} is the solvable group 
$\BS(1, p)$ defined by the presentation

\begin{equation}
\label{eqq/BS}
\tag{BS}
\BS(1, p) \, = \, \langle t, x \mid txt^{-1} = x^p \rangle.
\end{equation}
Let
$$
\Z[1/p] \, = \, \{ kp^n \mid k, n \in \Z \}
$$
be the subring of $\Q$ generated by $1/p$. 
The group $\BS(1, p)$ is isomorphic to the group of matrices
$$
\begin{pmatrix} p^\Z & \Z[1/p] \\ 0 & 1 \end{pmatrix}
\, = \,
\left\{\begin{pmatrix} p^n & b \\ 0 & 1 \end{pmatrix} \in \GL_2(\Z[1/p])
\hskip.1cm \Big\vert \hskip.1cm
n \in \Z, b \in \Z[1/p] \right\}.
$$
There is a decomposition in semi-direct product
$$
\BS(1, p) \, = \, A \ltimes N \, \approx \, \Z \ltimes \Z[1/p] ,
$$
where
$$
A \, = \, \left\{ \begin{pmatrix} p^n & 0 \\ 0 & 1 \end{pmatrix}
\hskip.1cm \Big\vert \hskip.1cm 
n \in \Z \right\} \, \approx \, \Z 
\hskip.5cm \text{and} \hskip.5cm
N \, = \, \left\{ \begin{pmatrix} 1& b \\ 0 & 1 \end{pmatrix}
\hskip.1cm \Big\vert \hskip.1cm
b \in \Z[1/p] \right\} \, \approx \, \Z[1/p] ,
$$
and where elements of $A \approx \Z$ 
act on $N \approx \Z[1/p]$ by multiplication by powers of $p$.
In particular, the group $\BS(1, p)$ is solvable.
\par

Two remarks are in order.
The first is that the presentation (\ref{eqq/BS})
defines a solvable Baumslag--Solitar groups for any $p \in \Z$,
not only for $p$ prime.
Though part of what follows carries over to $\BS(1, p)$ for $p \in \Z$,
we will stick to the examples with $p$ prime, for simplicity.
The second is that non-solvable Baumslag--Solitar groups,
despite their importance from many points of view,
do not appear in this book.

\subsection*{Dual of the group $\Z[1/p]$ and $p$-adic solenoid}

We need to determine the dual group $\widehat N$ of $N \approx \Z[1/p]$.
For this, we introduce the field $\Q_p$ of $p$-adic numbers.
On the one hand,
we identify the ring $\Z[1/p]$ with its image in the product $\Q_p \times \R$,
by the diagonal embedding: 
$$
\Z[1/p] \, = \, \{ (a,a) \in \Q_p \times \R \mid a \in \Z[1/p] \} .
$$
As is well-known, $\Z[1/p]$ is a cocompact discrete subring of $\Q_p \times \R$;
see, e.g., \cite[Example 5.C.10(2)]{CoHa--16}.
On the other hand, we introduce the subgroup
$$
\Delta \, := \, \{ (a,-a) \in \Q_p \times \R \mid a \in \Z[1/p] \}
\, \approx \, \Z[1/p] .
$$
\par

Recall that $\widehat\R$ is isomorphic to $\R$ through the map
\begin{equation}
\label{eqq/dualR}
\tag{dual $\R$}
\R \, \to \, \widehat \R, \hskip.2cm y \mapsto e_y
\end{equation}
where $e_y$ is the unitary character
defined by $e_y(t) = e^{2 \pi i yt}$ for all $t \in \R$.
\par

We will identify $\widehat{\Q_p}$ with $\Q_p$, 
in the following way.
Every $s \ne 0$ in $\Q_p$ can be uniquely written as 
$$
s \, = \, \sum_{j=m}^{\infty}a_jp^j
$$
for integers
$m \in {\mathbf Z}$ and $a_j \in \{0, \hdots, p-1 \}$ for $j \ge m$, with $a_m \ne 0$.
The ``fractional part" $\{s\}$ of $s \in \Q_p$ is the number
in $\mathopen[ 0,1 \mathclose[ \cap \Z[1/p]$ defined by 
$$
\{s\} \, = \,
\begin{cases}
\sum_{j=m}^{-1}a_jp^j & \text{if} \hskip.2cm m < 0 
\\
0 & \text{if} \hskip.2cm m \ge 0
\end{cases}
$$
if $s \ne 0$, and $\{ s \} = 0$ if $s = 0$.
Then $s = \{s\} + [s]$, where the ``integer part" $[s] = s - \{s\}$
belongs to the ring $\Z_p$ of $p$-adic integers.
Observe that, for $s \in \Q_p$, we have $s - \{s\} \in \Z$
if and only if $s \in \Z[1/p]$. 
A special unitary character $\chi$ of $\Q_p$ is defined by 
$$
\chi \colon \Q_p \to \T ,
\hskip.2cm 
s \mapsto e^{2 \pi i \{s\}}.
$$
Note that $\ker (\chi) = \Z_p$, so that $\chi$ is locally constant,
and in particular continuous.
For $x \in \Q_p$, define $\chi_x \in \widehat{\Q_p}$ by 
$\chi_x(s) = \chi(xs)$.
The map
\begin{equation}
\label{eqq/dualQp}
\tag{dual $\Q_p$}
\Q_p \to \widehat{\Q_p}, \hskip.2cm x \mapsto \chi_x
\end{equation}
is a topological group isomorphism
(see \cite[Chapitre 2, \S~1, No 9]{BTS1--2}).
% (see \cite[D.4.5]{BeHV--08}).
\par

We have therefore an isomorphism
$$
\Phi \, \colon \, 
\Q_p \times \R \to \widehat{\Q_p \times \R},
\hskip.2cm (x,y) \mapsto \chi_x e_y.
$$
Observe that $\Z$ acts on $\Q_p \times \R$
by multiplications by powers of $p$,
that we have
$$
\Phi(p^nx, p^ny) \, = \, \Phi(x,y)^{p^n}
\hskip.5cm \text{for all} \hskip.2cm
n \in \Z 
\hskip.2cm \text{and} \hskip.2cm 
(x,y) \in \Q_p \times \R ,
$$
and that the closed subring $\Delta$ of $\Q_p \times \R$
is invariant by this action.
\par

For the dual of the subgroup $\Z[1/p]$ of $\Q_p \times \R$,
we have by Pontrjagin duality
$$
\widehat{\Z[1/p]} \, \approx \, \widehat{\Q_p \times \R} / \Z[1/p]^\perp,
\leqno{(*)}
$$
where
$$
\Z[1/p]^\perp \, := \, 
\left\{
\chi \in \widehat{\Q_p \times \R} \hskip.1cm \big\vert \hskip.1cm \chi(a) = 1
\hskip.2cm \text{for all} \hskip.2cm 
a \in \Z[1/p]
\right\}.
$$

\begin{lem}
% 3.28
\label{Lem-DualAdicInteger}
By restriction to $\Delta$, the isomorphism
$\Phi \, \colon \Q_p \times \R \to \widehat{\Q_p \times \R}$
provides an isomorphism $\Delta \to \Z[1/p]^\perp$.
\end{lem}

\begin{proof}
Let $(x,y) \in \Q_p \times \R$. 
We have to prove that 
$$
(x,y) \in \Delta
\, \Longleftrightarrow \,
\Phi(x,y) \in \Z[1/p]^\perp .
$$
\par

Assume that $(x,y) \in \Delta$, that is, $x \in \Z[1/p]$ and $y = -x$.
Then, for every $a \in \Z[1/p]$, we have 
$$
\Phi(x,y) (a,a)
\, = \, \chi_x(a) e_{-x}(a) \, = \, e^{2 \pi i \{xa\}}e^{-2 \pi i xa}
\, = \, e^{2 \pi i (\{xa\}-xa)} \, = \, 1,
$$
since $\{xa\} - xa \in \Z$. Therefore, $\Phi(x,y) \in \Z[1/p]^\perp$.
\par

Conversely, assume that $\Phi(x,y) \in \Z[1/p]^\perp$, i.e.,
$$
\Phi(x,y) (a,a) \, = \, e^{2 \pi i (\{xa\} + ya)} \, = \, 1
\hskip.5cm \text{for all} \hskip.2cm
a \in \Z[1/p] ,
$$
i.e., $\{xa\} + ya \in \Z$ for every $a \in \Z[1/p]$.
Note that this implies $y \in \Z[1/p]$. 
Since $xa - \{xa\} \in \Z_p$ for all $a \in \Q_p$, 
we have
$$
(x + y)a \, = \, (xa - \{xa\}) + \{xa\} + ya \, \in \, \Z_p
\hskip.5cm \text{for all} \hskip.2cm
a \in \Z[1/p] .
$$
We have in particular
$$
\frac{x + y}{p^n} \, \in \, \Z_p
\hskip.5cm \text{for all} \hskip.2cm
n \in \N ,
$$
and therefore $x + y = 0$. 
It follows that $x \in \Z[1/p]$, and $(x,y) \in \Delta$.
\end{proof}

\index{$m5$@$p$-adic solenoid $\So_p$}
\index{Solenoid $\So_p$}
The compact group
$$
\So_p \, := \, (\Q_p \times \R) / \Delta,
$$
is known as the \textbf{$p$-adic solenoid}.
Recall that a solenoid is a topological group $G$ such that
there exists a continuous homomorphism $\R \to G$ with dense image.
A locally compact solenoid is abelian,
and either topologically isomorphic to $\R$ or compact.
More on this particular solenoid $\So_p$ in Appendix to Chapter I of \cite{Robe--00}.
\par

The following proposition shows that
$\So_p$ can be identified with the dual of $N \approx \Z[1/p]$.
It is a consequence of Lemma \ref{Lem-DualAdicInteger} and of (*);
see also Exercice~9 in \cite[Chapitre 2, \S~2]{BTS1--2}.
\par

We write $[(x,y)]$ for the class in $\So_p$ of $(x,y) \in \Q_p \times \R$,
and $[\chi_x e_y]$ for the class in $\widehat{\Q_p \times \R} / \Z[1/p]^\perp$
of $\chi_x e_y \in \widehat{\Q_p \times \R}$.

\begin{prop}
% 3.D.2
\label{dualsolenoid}
For a prime $p$, the map
\begin{equation}
\label{eqq/dualZ1/p}
\tag{dual $\Z[1/p]$}
\left\{
\begin{aligned}
\So_p \hskip.2cm &\overset{\approx}{\to} \hskip.2cm
\widehat{\Z[1/p]} \approx \widehat{\Q_p \times \R} / \Z[1/p]^\perp
\\
[(x,y)] \hskip.2cm &\mapsto
\hskip1cm [\chi_x e_y]
\end{aligned}
\right.
\end{equation}
is an isomorphism from the $p$-adic solenoid
to the dual of the group $\Z[1/p]$.
\end{prop}

The action of $\Z$ on $\widehat{\Z[1/p]}$ by multiplications by powers of $p$
corresponds to the action of $\Z$ on $\So_p$
for which the generator $1 \in \Z$ acts by the map
$$
T_p \, \colon \, \left\{
\begin{aligned}
\So_p \hskip.2cm &\to \hskip.2cm \So_p
\\
[(x,y)] \hskip.2cm &\hskip.2cm \mapsto [(px, py)] 
\end{aligned}
\right. .
$$
We will need to identify the set $\Per(T_p)$ of periodic points of $T_p$ in $\So_p$.

\begin{lem}
% 3.30
\label{Lem-OrbitsDualAdicInteger}
The set $\Per(T_p)$ of $T_p$-periodic points in $\So_p$ is the countable subset 
$$
\Per(T_p) \, = \, \left\{ \left[\left( \frac{a}{p^n - 1}, \frac{-a}{p^n - 1} \right)\right]
\in \So_p
\hskip.1cm \bigg\vert \hskip.1cm
n \in \N^*, a \in \Z[1/p]\right\}
$$
of $\So_p$.
\end{lem}

\begin{proof}
For $(x,y) \in \Q_p \times \R$ and $n \in \N^*$, we have
$T_p^n [(x,y)] = [(x,y)]$ if and only if
$(p^n - 1)x \in \Z[1/p]$, $(p^n - 1)y \in \Z[1/p]$, and $(p^n - 1) (x + y) = 0$,
if and only if
$x = a/(p^n - 1)$ and $y = -a/(p^n - 1)$ for some $a \in \Z[1/p]$.
\end{proof}

\subsection*{Some irreducible representations of the group $\BS(1, p)$}

We return to the Baumslag--Solitar group
$$
\BS(1, p) \, = \, A \times N ,
\hskip.2cm \text{where} \hskip.2cm
A = \begin{pmatrix} p^\Z & 0 \\ 0 & 1 \end{pmatrix} \approx \Z
\hskip.2cm \text{and} \hskip.2cm
N = \begin{pmatrix} 1 & \Z[1/p] \\ 0 & 1 \end{pmatrix} \approx \Z[1/p] ,
$$
and where the semi-direct product refers to the left action
$A \curvearrowright N, \hskip.1cm (p^n, t) \mapsto p^n t$.
The dual right action
$\widehat N \curvearrowleft A, \hskip.1cm (\chi, p^n) \mapsto \chi^{p^n}$
is defined by $\chi^{p^n}(t) = \chi(p^n t)$
for all $\chi \in \widehat N$, $n \in \Z$, and $t \in N$.
By the isomorphism of Proposition \ref{dualsolenoid}, this action is equivalent to
$$
\So_p \curvearrowleft A, \hskip.2cm 
\big( [(x,y)], p^n \big) \mapsto [(p^n x, p^n y)] = T_p^n([(x,y)])
$$
for all $[(x,y)] \in \So_p$ and $n \in \Z$.
Moreover, the set
$$
\widehat N_{\rm per} \, = \, \{ \chi \in \widehat N \mid
\hskip.2cm \text{the orbit} \hskip.2cm
\widehat N \cdot A \hskip.2cm \text{is finite} \}
$$
corresponds to the set $\Per(T_p)$ of $T_p$-periodic points in $\So_p$.

\begin{prop}
% 3.D.4
\label{Prop-IrredRepBS}
Let $p$ be a prime,
$\Gamma = \BS(1, p) = A \ltimes N$
the corresponding solvable Baumslag--Solitar group,
and $\widehat N_{\rm per} \subset \widehat N$ be as above.
\begin{enumerate}[label=(\arabic*)]
\item\label{iDECor-IrredRepBS}
For every $\chi \in \widehat N \smallsetminus \widehat N_{\rm per}$,
the induced representation $\pi_{\chi} := \Ind_{N}^\Gamma \chi$ is irreducible.
\item\label{iiDECor-IrredRepBS}
Let $\chi_1, \chi_2 \in \widehat N \smallsetminus \widehat N_{\rm per}$
be such that $\chi_1^\gamma \ne \chi_2$ for every $\gamma \in A$. 
Then $\pi_{\chi_1}$ and $\pi_{ \chi_2}$
are non-equivalent irreducible representations of $\Gamma$.
\end{enumerate}
\end{prop}

\begin{proof}
Since $A$ obviously acts freely
on $\widehat N \smallsetminus \widehat N_{\rm per}$,
Items \ref{iDECor-IrredRepBS} and \ref{iiDECor-IrredRepBS}
follow from Mackey--Shoda criteria for irreducibility 
and equivalence of monomial representations
(Corollaries~\ref{Cor-NormalSubgIrr} and \ref{Cor-NormalSubgEquiv}).
\end{proof}

Proposition \ref{Prop-IrredRepBS} can be restated as follows.

\begin{prop}
% 3.D.5
\label{Prop-IrredRepBS-bis}
For $s \in \So_p$,
denote by $\chi_{s}$ the unitary character of $N$ under the identification
$\widehat N = \widehat{\Z[1/p]} \approx \So_p$
of Proposition \ref{dualsolenoid}.
The map 
$$
\So_p \smallsetminus \Per(T_p) \, \to \, \widehat \Gamma, 
\hskip.2cm
s \, \mapsto \, \Ind_N^\Gamma \chi_{s}
$$
factorizes to an injective map from the space of non-periodic $T_p$-orbits in $\So_p$
into the dual $\widehat \Gamma$ of $\Gamma = \BS(1, p)$.
\end{prop}

\begin{rem}
% 3.D.6
\label{Rem-Prop-IrredBS-bis} 
Proposition~\ref{Prop-IrredRepBS-bis} shows that 
every non-periodic $T_p$-orbit in $\So_p$ provides 
one equivalence class of irreducible representations 
of $\BS(1, p)$.
\par

Equip the compact group $\So_p$ with the normalized Haar measure $\mu$.
Then, by uniqueness of $\mu$ as a translation-invariant probability measure,
$\mu$ is $T_p$-invariant.
Moreover, the action of $T_p$ on $(\So_p, \mu)$ is ergodic.
Indeed, the character group of $\So_p$ can be identified with $\Z[1/p]$,
with the dual action of $T_p$ on $\Z[1/p]$ corresponding to multiplication by $p$.
Since $\{0\}$ is the only $T_p$-fixed point in $\Z[1/p]$,
the transformation $T_p$ of $\So_p$ is ergodic
(see \cite[Chap I, Proposition 1.5]{BeMa--00}).
\par

As in Corollary \ref{Cor-Cor-IrredRepHeisInteger}
and Remark~\ref{Rem-IrredRepAffine},
we conclude that the Borel space
$\widehat{\BS(1, p)}$ is not countably separated.
\end{rem}

We proceed now with the construction
of irreducible representations of $\BS(1, p)$
associated to the periodic orbits in $\So_p$. 
It will turn out that these representations are finite-dimensional.
We will see later (Corollary~\ref{Cor-FinDimRepBS})
that they exhaust the space $\widehat{\BS(1, p)}_{\rm fd}$
of finite-dimensional irreducible representations of $\BS(1, p)$.
\par

For $n \in \N^*$, we denote by $\So_p(n)$ the elements in $\So_p$
with $T_p$-period $n$.
For $k \in \Z$ and $b \in \Z[1/p]$, we write $(k, b)$ for the element
$\begin{pmatrix} p^k & b \\ 0 & 1 \end{pmatrix}$ of $\BS(1, p)$.
Fix an integer $n \in \N^*$ and a point $s \in \So_p(n)$.
The stabilizer of $s$ in $A$ is the subgroup
$$
A(n) \, = \, 
\left\{ (k,0) \in \BS(1, p) \mid k \in n\Z \right\} \, \approx \, n\Z .
$$
For $\Gamma = \BS(1, p)$, the set
$$
\Gamma(n) \, := \, A(n) \ltimes N \, = \,
\left\{ (k, b) \in \Gamma \mid k \in n\Z, \hskip.1cm b \in \Z[1/p] \right\}
$$
is a normal subgroup of index $n$ in $\Gamma$.
As in Section~\ref{Section-IrrRepTwoStepNil},
we parametrize the dual of $n\Z$ by $\mathopen[0,1/n \mathclose[$;
for $\theta \in \mathopen[0,1/n \mathclose[$,
we define $\chi_{s, \theta} \, \colon \Gamma(n) \to \T$ by 
$$
\chi_{s, \theta} (k, b) \, = \, \chi_\theta(k)\chi_s(b)
\hskip.5cm \text{for all} \hskip.2cm
(k, b) \in \Gamma(n).
$$
Since $\chi_s$ is $A(n)$-invariant, 
$\chi_{s, \theta}$ is a unitary character of $\Gamma(n)$; 
obviously, we have $\chi_{s, \theta} \vert_N = \chi_s$. 
Conversely, it is clear that every unitary character $\chi$ of $\Gamma(n)$ 
such that $\chi \vert_N = \chi_s$ is of the form $\chi_{s, \theta}$
for a unique $\theta \in \mathopen[0, 1/n \mathclose[$.

\begin{prop}
% 3.D.7
\label{Prop-IrredRepBS-FiniteDim}
Let $\Gamma = \BS(1, p)$ be the Baumslag--Solitar group.
Fix an integer $n \ge 1$.
\begin{enumerate}[label=(\arabic*)]
\item\label{iDEProp-IrredRepBS-FiniteDim}
Let $s \in \So_p(n)$ and $\theta \in \mathopen[0,1/n \mathclose[$.
Then 
$$
\pi_{s, \theta} \, := \, \Ind_{\Gamma(n)}^\Gamma \chi_{s, \theta}
$$
is an irreducible representation of $\Gamma$ of dimension $n$,
where $\Gamma(n)$ is the normal subgroup of $\Gamma$ defined above.
\item\label{iiDEProp-IrredRepBS-FiniteDim}
Let $(s, \theta), (s', \theta') \in \So_p(n) \times \mathopen[0,1/n \mathclose[$.
The representations $\pi_{s, \theta}$ and $\pi_{s', \theta'}$ are equivalent
if and only if $s$ and $s'$ are in the same $T_p$-orbit and $\theta = \theta'$.
\end{enumerate}
\end{prop}

\begin{proof}
\ref{iDEProp-IrredRepBS-FiniteDim}
It is clear that $\chi_{s, \theta}^\gamma \ne \chi_{s, \theta}$
for all $\gamma \in \Gamma \smallsetminus \Gamma(n)$. 
Therefore $\pi_{s, \theta}$
is irreducible, by Mackey--Shoda criterion 
(Corollary~\ref{Cor-NormalSubgIrr}).
Moreover, since $\Gamma(n)$ has index $n$ in $\Gamma$ 
and since $\chi_{s, \theta}$ is one-dimensional, 
the induced representation $\pi_{s, \theta}$ is finite-dimensional,
of dimension $n$.

\vskip.2cm

\ref{iiDEProp-IrredRepBS-FiniteDim}
Let $(s, \theta), (s', \theta') \in \So_p(n) \times \mathopen[0,1/n \mathclose[$.
The unitary characters $\chi_{s, \theta}$ and $\chi_{s', \theta'}$ of $\Gamma(n)$
are conjugate under $\Gamma$
if and only if $s$ and $s'$ are in the same $T_p$-orbit and $\theta = \theta'$.
Therefore, the claims follows from
the Mackey--Shoda criterion for equivalence of monomial representations 
(Corollary~\ref{Cor-NormalSubgEquiv}). 
\end{proof}

\section
{Lamplighter group}
% Section 3.E
\label{Section-IrrRepLamplighter}

\index{Lamplighter group}
The \textbf{lamplighter group} is the solvable group 
$\Gamma$ defined by the presentation
$$
\Gamma \, = \, \langle t, x \mid x^2 = 1, \hskip.1cm [t^kxt^{-k}, t^lxt^{-l}] = 1 
\hskip.5cm \text{for all} \hskip.2cm 
k, l \in \Z \rangle .
$$
The group $\Gamma$ can also be described as the wreath product 
$$
\Z \wr (\Z / 2 \Z) = A \ltimes N ,
\hskip.5cm \text{where} \hskip.5cm
A = \Z
\hskip.2cm \text{and} \hskip.2cm
N = \bigoplus_{k \in \Z} \Z / 2 \Z ,
$$
where the action of $\Z$ on $\bigoplus_{k \in \Z} \Z / 2 \Z$
defining the semi-direct product is given by shifting the coordinates.
\par

The dual group of $N$ can be identified
with the compact group
$$
X \, := \, \prod_{k \in \Z} \{0,1 \} .
$$
Under this identification, the action of $\Z$ on $\widehat N$ is given by 
the shift transformation 
$$
T \, \colon \, \prod_{k \in \Z} \{0,1 \} \to \prod_{k \in \Z} \{0,1 \}, 
\hskip.2cm
(x_n)_{n \in \Z} \mapsto (x_{n + 1})_{n \in \Z}.
$$
Moreover, the normalized Haar measure on $\widehat N = X$
is the measure $\mu = \otimes^{\Z} \nu$,
where $\nu$ is the uniform probability measure on $\{0,1 \}$.
\par

Observe that $\mu$ is $T$-invariant
(as it should be, since $T$ is an automorphism of $X$).
The action of $T$ on $(X, \mu)$ is ergodic,
since $\{0\}$ is the only $T$-fixed point in $\widehat X = N$
(see \cite[Chap I, Proposition 1.5]{BeMa--00}).
\par

For an integer $n \ge 1$, the set $X(n)$ of elements in $X$ with $T$-period $n$
consists of the sequences $(x_m)_{m \in \Z}\in \{0,1 \}^\Z$ with 
$$
x_{j + kn} \, = \, x_j
\hskip.5cm \text{for all} \hskip.2cm
j = 0, \cdots, n - 1
\hskip.2cm \text{and} \hskip.2cm
k \in \Z.
$$
So, $X(n)$ has exactly $2^n$ elements and the set 
$$
\Per(T) \, = \, \bigcup_{n \ge 1} X(n)
$$
of periodic points of $T$ is a countable dense subset of $X$. 
\par

%(see \cite[Chap I, Proposition 1.5]{BeMa--00}).
We can construct a family of irreducible representations
of the lamplighter group associated to the non-periodic $T$-orbits, 
as we did in Section~\ref{Section-IrrRepBS}
for the Baumslag--Solitar group $\BS(1, p)$. 
The proof of Proposition \ref{Prop-IrredRepBS}
carries over \emph{mutatis mutandis} and yields the following result.

\begin{prop}
% 3.E.1
\label{Prop-IrredRepLamplighter-bis}
Let $\Gamma = A \ltimes N = \Z \ltimes \bigoplus_{k \in \Z} \Z / 2 \Z$
and $X = \prod_{k \in \Z} \{0,1 \}$ be as above.
For $x \in X$,
denote by $\chi_{x}$ the corresponding unitary character
of $N$ under the identification $\widehat N \approx X$.
\par

For every $x \in X \smallsetminus \Per(T)$,
the representation $\Ind_N^\Gamma \chi_{x}$ is irreducible. The map
$$
X \smallsetminus \Per(T) \, \to \, \widehat \Gamma, 
\hskip.2cm
x \, \mapsto \, \Ind_N^\Gamma \chi_{x}
$$
factorizes to an injective map from the space of non-periodic $T$-orbits in $X$
into the dual $\widehat \Gamma$ of $\Gamma$.
\end{prop}

As for the Baumslag--Solitar group $\BS(1, p)$,
we will also construct irreducible representations of $\Gamma = A \ltimes N$
associated to periodic orbits in $X$.
It will turn out that these representations are finite-dimensional
and we will see later (Corollary~\ref{Cor-FinDimRepLamplighter})
that they exhaust the space $\widehat \Gamma_{\rm fd}$
of finite-dimensional irreducible representations of~$\Gamma$.
\par

Fix an integer $n \ge 1$.
Let $x \in X(n)$ be an element of $X$ of $T$-period $n$,
and $\chi_x \in \widehat N$ the corresponding unitary character of $N$.
The stabilizer of $\chi_x$ in $A$ is the subgroup $A(n) := n\Z$ of $A = \Z$.
Then $\Gamma(n) := A(n) \ltimes N$ is a normal subgroup of index $n$ in $\Gamma$.
As in Section~\ref{Section-IrrRepBS}, we parametrize $\widehat{A(n)}$
by the interval $\mathopen[0,1/n \mathclose[$:
for every $\theta \in \mathopen[0,1/n \mathclose[$,
we have a unitary character $\chi_{x, \theta}$ of $\Gamma(n)$ defined by 
$$
\chi_{x, \theta} (k, b) \, = \, \chi_\theta(k)\chi_x(b)
\hskip.5cm \text{for all} \hskip.2cm
(k, b) \in \Gamma(n) = n\Z \ltimes N.
$$
The following proposition is proved along the same lines as the corresponding 
result (Proposition~\ref{Prop-IrredRepBS-FiniteDim}) for the Baumslag--Solitar group.

\begin{prop}
% 3.E.2
\label{Prop-IrredRepLamplighter-FiniteDim}
Let $\Gamma = A \ltimes N$ be the lamplighter group; set $X = \widehat N$. 
Fix an integer $n \ge 1$;
let $X(n)$ be the subset of $X$
and $\Gamma(n)$ the subgroup of index $n$ of $\Gamma$ defined above.
\begin{enumerate}[label=(\arabic*)]
\item\label{iDEProp-IrredRepLamplighter-FiniteDim}
Let $x \in X(n)$ and $\theta \in \mathopen[0,1/n \mathclose[$,
and $\chi_{x, \theta} \in \widehat{\Gamma(n)}$ as above.
Then 
$$
\pi_{x, \theta} \, := \, \Ind_{\Gamma(n)}^\Gamma \chi_{x, \theta}
$$
is an irreducible representation of $\Gamma$ of dimension $n$.
\item\label{iiDEProp-IrredRepLamplighter-FiniteDim}
Let $(x, \theta), (x', \theta') \in X(n) \times \mathopen[0, 1/n \mathclose[$.
The representations
$\pi_{x, \theta}$ and $\pi_{x', \theta'}$ are equivalent if and only if 
$x$ and $x'$ are in the same $T$-orbit and $\theta = \theta'$.
\end{enumerate}
\end{prop}

\section
{General linear groups}
% Section 3.F
\label{Section-IrrRepGLN}

Let $\K$ be a field and $n$ an integer, $n \ge 2$.
Let $\Gamma = \GL_n(\K)$ be the \textbf{general linear group} over $\K$.
\index{General linear group! $\GL_n(\K)$ with $\K$ a field}
\par

Basic representations of $\Gamma$ are the 
\textbf{principal series representations},
which are defined as follows. Let
\index{Principal series! of $\GL_n(\K)$}
$$
B \, = \, \begin{pmatrix} \K^*& \K\phantom{^*} & \cdots & \K\phantom{^*} 
\\ 
0&\K^* & \cdots & \K\phantom{^*}
\\
\vdots & \vdots & \ddots & \vdots
\\
0 & 0 & \cdots & \K^*
\end{pmatrix}
$$
be the subgroup of upper-triangular matrices in $\Gamma$.
Then $B = A\ltimes N$ is the semi-direct product of the subgroup 
$$
A \, = \, \left\{\, \begin{pmatrix} a_1 & 0 & \cdots &0 
\\ 
0 &a_2 & \cdots & 0
\\
\vdots & \vdots & \ddots & \vdots
\\
0 & 0 & \cdots & a_n
\end{pmatrix} 
\hskip.2cm \Bigg\vert \hskip.2cm 
a_1, \hdots, a_n \in \K^\times \right\} \approx (\K^\times)^n
$$
of diagonal matrices with the subgroup 
$$
N \, = \, \begin{pmatrix} 1& \K & \cdots & \K 
\\ 
0 & 1 & \cdots & \K
\\
\vdots & \ddots & \ddots & \vdots
\\
0 & 0 & \cdots & 1
 \end{pmatrix}
$$
of upper-triangular matrices in $\Gamma$
with $1$'s on the diagonal.

A fundamental fact about $B$ is the following result.

\begin{lem}
% 3.F.1
\label{Lem-CommBorelSub}
Assume that $\K$ is infinite.
\par

Then $B$ coincides with its own commensurator $\Comm_\Gamma (B)$. 
\end{lem}
\index{Commensurator}

\begin{proof}
Let $X$ denote the set of complete flags in $\K^n$
and $x_0$ the complete flag
$$
V_1 \, = \, \langle e_1 \rangle 
\subset V_2 \, = \, \langle e_1, e_2 \rangle 
\subset \cdots 
\subset V_n \, = \, \langle e_1, \cdots, e_{n} \rangle,
$$
where $(e_1, e_2, \hdots, e_n)$ is the canonical basis of $\K^n$.
\par

The group $\Gamma = GL_n(\K)$ acts transitively on $X$ and $B$ is the stabilizer of $x_0$.
The claim is that $x_0$ is the only point in $X$ with a finite $B$-orbit;
equivalently, the claim is that, for every subgroup $H$ of $B$ of finite index,
$x_0$ is the unique $H$-fixed point.
(On commensurators, see Proposition \ref{defScommensurator}.)
\par

Let $H$ be a subgroup of finite index in $B$.
It suffices to prove that $V_1, \hdots, V_n$
are the unique non-zero $H$-invariant subspaces of $\K^n$.
\par

Let $V$ be a non-zero $H$-invariant subspace of $\K^n$.
Let $j \in \{1, \hdots, n\}$ be minimal such that $V \subset V_j$.
There exists $v \in V$ such that $v = \sum_{k=1}^j \alpha_k e_k$ 
for scalars $\alpha_k \in \K$, with $\alpha_j \ne 0$.
If $j = 1$, then $V = V_j$ and the claim is proved. 
So, we can assume that $j \ge 2$.
\par

Let $i \in \{1, \hdots, j-1 \}$. 
Consider the subgroup $B_{i, j} \cong \K$ of $B$
consisting of the elementary matrices $E_{i, j} (\alpha)$, for $\alpha \in \K$. 
(Thus, $E_{i, j} (\alpha)$ is the matrix $(\alpha_{k, l})_{1 \le k, l \le n}$ 
with $\alpha_{k, k} = 1$ for all $k$, 
with $\alpha_{i, j} = \alpha$ and with $\alpha_{k, l} = 0$ otherwise.)
Then $H_{i, j} := B_{i, j} \cap H$ is a finite index subgroup of $B_{i, j}$.
Observe that $H_{i, j}$ is infinite, since $\K$ is infinite; 
in particular, $H_{i, j}$ is not the trivial group 
and so, there exists $\alpha \ne 0$ such that $E_{i, j} (\alpha) \in H$.
On the one hand,
$$
E_{i, j}(\alpha)v - v \in V,
$$
since $V$ is $H$-invariant.
On the other hand,
$$
E_{i, j}(\alpha)v - v
\, = \, E_{i, j} (\alpha) \Big( \sum_{k=1}^j \alpha_k e_k \Big) - \sum_{k=1}^j \alpha_k e_k
\, = \, \alpha \alpha_j e_i.
$$
Since $\alpha \alpha_j \ne 0$, this shows that $e_i \in V$ for all $i \in \{1, \hdots, j-1 \}$.
Therefore we also have $e_j \in V$, since $v \in V$.
Therefore, $V = V_j$ and the claim is proved.
\end{proof}

Let $\chi = (\chi_1, \hdots, \chi_n) \in (\K^\times)^n$.
Define a representation of $B$, again denoted by $\chi$,
by lifting the unitary character of $A$ defined by $\chi$ to $B$, that is,
$$
\chi
\begin{pmatrix} 
a_1 & \ast & \cdots & \ast 
\\ 
0 & a_n & \cdots & \ast
\\
\vdots & \vdots & \ddots & \vdots
\\
0 & 0 & \cdots & a_n
 \end{pmatrix}
 \, = \, \chi_1(a_1) \cdots \chi_n(a_n).
$$
The following result is a consequence of Lemma~\ref{Lem-CommBorelSub} 
in combination with the Mackey--Shoda criteria 
(Theorems~\ref{Theo-IrredInducedRep} and \ref{Theo-EquiInducedRep}).

\begin{theorem}
% 3.F.2
\label{Theo-RepIrrGLn}
Let $\Gamma = GL_n(\K)$ for an infinite field $\K$.
\begin{enumerate}[label=(\arabic*)]
\item\label{iDETheo-RepIrrGLn}
For every $\chi \in (\K^\times)^n$, the representation $\Ind_B^\Gamma \chi$ is irreducible.
\item\label{iiDETheo-RepIrrGLn}
For $\chi, \chi' \in (\K^\times)^n$ with $\chi \ne \chi'$,
the representations $\Ind_B^\Gamma \chi$ 
and $\Ind_B^\Gamma \chi'$ are not equivalent.
\end{enumerate}
\end{theorem}

\section
{Some non-discrete examples}
% Section 3.G
\label{ExamplesND}

There are some well-studied examples of topological groups, non-discrete ones,
of which the dual has been determined, either as a set,
or more rarely as a topological space.
\par

We have already described
the dual of the real Heisenberg group $H(\R)$ in Remark \ref{dualHeisLie},
and the dual of the affine group $\Aff(\R)$ in Remark \ref{AffKLie}.
\par

For the simple complex Lie group $G = \SL_2(\C)$,
the dual as a set has been determined as soon as 1947,
independently by Bargmann, Gel'fand and Naimark, and Harish--Chandra
\cite{Barg--47, GeNa--47b, Hari--47}.
They showed that $\widehat G$ consists of three parts,
the first is the so-called principal series
(a disjoint union of a half-line with an infinite number of copies of $\R$),
the second the supplementary series
(an open interval)
and the third a point (the representation $1_G$).
The Fell topology on $\widehat G$ has been described by Fell
\cite{Fell--61}.
See also \cite{Fell--60a} for the dual of $\SL_n(\C)$.
\par

For $G = \SL_2(\R)$, we refer to \cite[Section 7.6]{Foll--16}
for a short description of the dual and its topology, without proofs;
see also \cite{Vale--84}.
The determination of the dual as a set was obtained in \cite{Barg--47},
and the Fell topology in \cite{Mili--71}.
For $\SO(n) \ltimes \R^n$, we refer to \cite[Example 5.59]{KaTa--13}.
\par

For semisimple groups in general,
there are descriptions of the Fell topology on parts of the dual, 
in particular on parts defined by the irreducible principal series,
or by the discrete series.
For the simple linear groups of so-called split-rank $1$,
there is also a description of the Fell topology on the reduced dual,
i.e., on that part of the dual consisting of irreducible representations
that are weakly contained in the regular representation.
See \cite{Lips--70}.

\begin{exe}[\textbf{some groups which are not locally compact}]
% 3.G.1
\label{exUinfty}
Consider now the following non-LC topological groups.
Let $\Hi$ be an infinite-dimensional separable Hilbert space.
\par

Let first $\U_\infty(\Hi)$ be the group of unitary operators $x$ on $\Hi$
such that $x-\mathrm{Id}_\Hi$ is compact, with the topology inherited from
the distance defined by $d(x, y) = \Vert x - y \Vert$.
The dual of $\U_\infty(\Hi)$ is countable and discrete \cite{Kiri--73, Ol's--78}.
Moreover, every representation of $\U_\infty(\Hi)$
is a direct sum of irreducible subrepresentations \cite[Theorem 1.11]{Ol's--78}.
\par

Representations of infinite-dimensional unitary groups have been studied
earlier in \cite{Sega--57}.

Let then $\U(\Hi)_{\rm str}$ denote the unitary group of $\Hi$
with the strong topology.
The inclusion $i \, \colon \U_\infty(\Hi) \to \U(\Hi)_{\rm str}$ is continuous with dense image,
so that there is an embedding $\pi \mapsto \pi \circ i$
from the dual of $\U(\Hi)_{\rm str}$ into that of $\U_\infty(\Hi)$.
Moreover, the explicit description by Kirillov and Ol'shanskii of the dual of $\U_\infty(\Hi)$
shows that every irreducible representation of $\U_\infty(\Hi)$
is of the form $\pi \circ i$ for some irreducible representation $\pi$ of $\\U(\Hi)_{\rm str}$;
it follows that the restriction of irreducible representations
provides a homeomorphism from the dual of $\U(\Hi)_{\rm str}$ 
onto that of $\U_\infty(\Hi)$.
If $\rho$ is a representation of $\U(\Hi)_{\rm str}$, 
its restriction to $\U_\infty(\Hi)$ is continuous,
and therefore is a direct sum of irreducible representations,
so that $\rho$ itself is a direct sum of irreducible representations.
\par

With the terminology of Chapter \ref{ChapterTypeI},
this shows that the groups
$\U_\infty(\Hi)$ and $\U(\Hi)_{\rm str}$ are of type I
(see Theorem~\ref{explesTypeI-bis}).
\end{exe}

%-----------------------------------------------------------------------
% End of chapter 3
%-----------------------------------------------------------------------

\chapter[Finite dimensional representations]
{Finite dimensional irreducible representations}
% Chapter 4
\label{Chapter-AllFiniteDimensionalRep}

\emph{This chapter is dedicated to the study
of the space $\widehat G_{\rm fd}$ of equivalence classes
of finite dimensional irreducible representation of a group $G$.
We describe $\widehat G_{\rm fd}$
for a locally compact group which is a semi-direct product $G = H \ltimes N$, where
\begin{enumerate}
\item[$\bullet$]
$H$ is an abelian discrete subgroup of $G$;
\item[$\bullet$]
$N$ is a locally compact abelian normal closed subgroup of $G$.
\end{enumerate}
It turns out that every finite dimensional irreducible representation $\pi$ of $G$
is equivalent to a representation of the form $ \Ind_{M}^G \chi$,
where $M$ is a finite index subgroup of $G$ which contains $N$
and $\chi$ a unitary character of $M$ (Theorem~\ref{Theo-FiniteDimRepSemiDirect}).
}
\par

\emph{
This implies that the finite dimensional irreducible representations 
constructed in Chapter \ref{Chapter-ExamplesIndIrrRep}
exhaust the space $\widehat \Gamma_{\rm fd}$
for the groups $\Gamma$ considered there, i.e.,
for the Heisenberg group over a commutative ring $R$,
the affine group $\Aff(\K)$ over a field $\K$,
the Baumslag--Solitar group $\BS(1, p)$,
and the lamplighter group
--- the case of $\GL_n(\K)$ over an infinite field $\K$ is treated separately
(Section \ref{Section-FiniteDimRepForSomeGroups}).
For all these groups (apart from the special case of finite groups),
$\widehat \Gamma_{\rm fd}$ is a much smaller space that $\widehat \Gamma$.
}
\par

\emph{
The set $\widehat \Gamma_{\rm fd}$
separates the points of $\Gamma$ for some of these groups:
for the Heisenberg group $H(\Z)$ over the integers, 
and for the Baumslag--Solitar group $\BS(1, p)$).
By contrast, $\widehat \Gamma_{\rm fd}$
consists only of one-dimensional representations for other groups:
for $H(\K)$ or $\GL_n(\K)$ over an infinite field $\K$.
}
\par

\emph{There are several interesting classes of groups $G$
defined in terms of the space $\widehat G_{\rm fd}$.
We first survey results concerning the so-called Moore groups,
which are those groups $G$ such that
$\widehat G_{\rm fd} = \widehat G$
(\ref{SS:MooreGr}).
Maximally almost periodic groups
are topological groups $G$
such that $\widehat G_{\rm fd}$ separates the points of $G$
(\ref{SS:MAP}),
and minimally almost periodic groups are these $G$
such that $\widehat G_{\rm fd} = \{ 1_G \}$
(\ref{SS:map}).
These classes can be characterized in terms of properties
of the Bohr compactification $\Bohr(G)$ of $G$,
a compact group which encodes the space $\widehat G_{\rm fd}$}
(\ref{SS:Profinite}).
\par

\emph{We discuss the class of discrete groups $\Gamma$
with the property that every finite dimensional representation of $\Gamma$
has finite image; these are the groups $\Gamma$
for which the Bohr compactification $\Bohr(\Gamma)$
coincides with another compactification of $\Gamma$,
the profinite completion $\Prof(\Gamma)$ of $\Gamma$
(Proposition~\ref{Pro-IsoBohrProfini}).
A prominent example of a group with this property
is the group $\SL_n(\Z)$ for $n \ge 3$ (Corollary~\ref{Cor-SLnZ-BohrProf}).}

\section[Some semi-direct products]
{Finite dimensional irreducible representations of some semi-direct products}
% Section 4.A
\label{Section-FinDimRepSemiDirect}

\begin{conv}
% [semi-direct product]
% 4.A.1
\label{convsemidirect}
We need to specify the definition we will adopt for a semi-direct product.
Let $H$ be a group acting on the left by automorphisms on another group~$N$;
let $\varphi \, \colon H \to \Aut (N)$ denote the corresponding group homomorphism.
The \textbf{semi-direct product} $G = H \ltimes N$ is defined
as the group with underlying set $H\times N$ and with the product
\index{Semi-direct product}
$$
(h_1, n_1) (h_2, n_2) = (h_1 h_2, n_1 \varphi (h_1) (n_2))
\hskip.5cm \text{for} \hskip.2cm
(h_1, n_1), (h_2, n_2) \in G.
$$
The obvious maps $H \to G$ and $N \to G$ are injective group homomorphisms;
we can therefore identify $H$ with the subgroup $H \times \{e\}$ of $G$
and $N$ with the normal subgroup $\{e\} \times N$ of $G$.
When $\varphi $ is understood, we also write $h \cdot n$ for $\varphi (h) (n)$.
Observe that we way write $nh = (h,n)$ but that, in general, $hn \ne (h,n)$;
indeed, we have $hn = (h, h \cdot n)$.
\par

The following identities will be often used below;
we use the identifications noted above.
For $h \in H, n' \in N$, we have $h n' h^{-1} = h \cdot n'$.
When $N$ is moreover abelian, for $g = (h,n) \in G$ and $n' \in N$,
we have $g n' g^{-1} = h n' h^{-1} = h \cdot n'$.
\par

Beware that, elsewhere and for example in \cite{BA1--3},
the semi-direct product of $H$ and $N$ is defined as the set $N \times H$,
with the same law as above, and is then denoted by $N \rtimes H$. 

\vskip.2cm

Suppose moreover that $H$ and $N$ are locally compact groups
and that $\varphi$ is a continuous homomorphism from $H$
to the group of bicontinuous automorphisms of $N$.
Then $H \ltimes N$ is a locally compact group for the product topology on $H \times N$.
\end{conv}

\index{$b8$@$\widehat G_{\rm fd}$ finite-dimensional part of $\widehat G$}
\index{Dual! $3$@finite-dimensional part}
\index{Finite-dimensional part of the dual}
For a topological group $G$, denote by $\widehat G_{\rm fd}$
the subset of its dual containing those equivalence classes of irreducible representations
which are finite-dimensional.
\par

Given a closed normal subgroup $N$ of $G$,
recall from Section~\ref{Section-IrrIndRep}
that $G$ acts (on the right) on the dual space $\widehat N$ by $(\pi,g) \to \pi^g$,
where $\pi^g$ is given $\pi^g(n) = \pi(gng^{-1})$
for $\pi \in \widehat N$, $g \in G$, and $n \in N$.
\par

Assume from now on that $G = H \ltimes N$ is a semi-direct product
of an abelian discrete subgroup $H$
and a locally compact abelian normal closed subgroup $N$.
For $\chi \in {\widehat N}$, let
$$
H_\chi \, = \, \{ h \in H \mid \chi^h = \chi \}
$$
be the stabilizer of $\chi$ in $H$.
Observe that $G_\chi = H_\chi \ltimes N$ is the stabilizer of $\chi$ in $G$.

\begin{lem}
% 4.A.2
\label{Lem-FiniteDim}
Consider $G = H \ltimes N$ and $\chi \in {\widehat N}$,
with $H_\chi \subset H$ and $G_\chi \subset G$, as above.
\begin{enumerate}[label=(\arabic*)]
\item\label{1DELem-FiniteDim}
For every $\alpha \in \widehat{H_\chi}$, the map 
$$
\alpha \otimes \chi \, \colon \, H_\chi \ltimes N \to \T,
\hskip.2cm 
(h, n) \mapsto \alpha(h)\chi(n).
$$
is a unitary character of $G_\chi$ which extends $\chi$.
\item\label{2DELem-FiniteDim}
Let $\pi$ be an irreducible representation of $G_\chi$
such that $\pi \vert_N$ is a multiple of~$\chi$.
Then $\pi$ is a unitary character of $G_\chi$
and there exists $\alpha \in \widehat{H_\chi}$
such that $\pi = \alpha \otimes \chi$.
\end{enumerate}
\end{lem}

\begin{proof}
\ref{1DELem-FiniteDim}
It is obvious that $\alpha \otimes \chi$ is continuous and that it extends $\chi$.
We have to check that 
$\alpha \otimes \chi$ is a homomorphism.
This is indeed the case: for $h_1, h_2 \in H_\chi$ and 
$n_1,n_2 \in N$, we have
$$
\begin{aligned}
(\alpha \otimes \chi) \big( (h_1, n_1) (h_2, n_2) \big)
\, &= \,
(\alpha \otimes \chi) \big( (h_1 h_2, e)(e, n_1 (h_1 \cdot n_2)) \big)
\\
\, &= \, \alpha(h_1) \alpha(h_2) \chi(n_1) \chi( h_1 \cdot n_2 )
\\
\, &= \, \alpha(h_1) \alpha(h_2) \chi(n_1) \chi(n_2)
\\
\, &= \, (\alpha \otimes \chi) (h_1, n_1) \hskip.1cm (\alpha \otimes \chi) (h_2, n_2) .
\end{aligned}
$$

\vskip.2cm

\ref{2DELem-FiniteDim}
For $g_1 = (h_1, n_1), g_2 = (h_2, n_2) \in G_\chi$, we have
$$
\begin{aligned}
\pi(g_1g_2)
\, &= \, \pi \big( (h_1h_2, e) (e, n_1 (h_1 \cdot n_2) ) \big)
\\
\, &= \, \pi(h_1) \pi(h_2) \chi(n_1) \chi( h_1 \cdot n_2)
\\
\, &= \, \pi(h_1) \pi(h_2) \chi(n_1) \chi(n_2)
\\
\, &= \, \pi(h_1h_2) \chi(n_1n_2) 
\end{aligned}
$$
and hence, since $H$ is abelian, $\pi(g_1g_2) = \pi(g_2 g_1)$.
This shows that $\pi$ is trivial on the commutator subgroup of $G_\chi$,
that is $\pi$ is a unitary character of $G_\chi$.
Moreover, $\pi = \alpha\otimes \chi$ for $\alpha := \pi \vert_{H_\chi} \in \widehat{H_\chi}$. 
\end{proof}

Set
$$
{\widehat N}_{\rm per} \, = \, 
\{\chi \in \widehat N \mid \text{the orbit} \hskip.2cm \chi \cdot H
\hskip.2cm \text{is finite} \} 
$$
(equivalently, the orbit $\chi \cdot G$ is finite).
Thus, ${\widehat N}_{\rm per}$ is the set $\chi \in \widehat N$
such that $H_\chi$ is a subgroup of finite index in $H$. 
Observe that, since $H$ is abelian, we have 
$$
H_\chi \, = \, H_{\chi^h}
\hskip.5cm \text{for all} \hskip.5cm
h \in H.
$$
Set
$$
X_{\rm fd} \, = \,
\bigsqcup_{\chi \in {\widehat N}_{\rm per}} \{ \chi \} \times \widehat{H_\chi}.
$$
There is a natural action $X_{\rm fd} \curvearrowleft H$ of $H$
on $X_{\rm fd}$ given by 
$$
(X_{\rm fd},H) \, \to \, X_{\rm fd},
\hskip.5cm
((\chi, \alpha), h) \, \mapsto \, (\chi^h, \alpha).
$$

\begin{theorem}
% 4.A.3
\label{Theo-FiniteDimRepSemiDirect}
Let $G = H \ltimes N$ be the semi-direct product
of a discrete abelian subgroup $H$
and a locally compact abelian normal subgroup $N$ of $G$.
Let $X_{\rm fd}$ be as above.
\par

The map $\Phi \, \colon X_{\rm fd} \to \widehat G_{\rm fd}$, defined by 
$$
\Phi \, \colon \, (\chi, \alpha) \, \mapsto \, \Ind_{H_\chi \ltimes N}^G (\alpha \otimes \chi)
$$
induces a bijection
between the space of $H$-orbits in $X_{\rm fd}$
and the space $\widehat G_{\rm fd}$
of finite-dimensional irreducible representations of $G$.
\end{theorem}

\begin{proof}
$\bullet$ \emph{First step.} 
We claim that the range of $\Phi$ is contained in $\widehat G_{\rm fd}$.
\par

Indeed, let $\chi \in {\widehat N}_{\rm per}$ and $ \alpha \in \widehat H_\chi$. 
Then 
$$
(\alpha \otimes \chi)^h \, = \, \alpha \otimes \chi^h \, \ne \, \alpha \otimes \chi
\hskip.5cm \text{for all} \hskip.2cm
h \in H\smallsetminus H_\chi .
$$
Therefore $\Ind_{H_\chi \ltimes N}^G (\alpha \otimes \chi)$
is irreducible, by the Mackey--Shoda criterion 
(Corollary~\ref{Cor-NormalSubgIrr}).
Moreover, since $H_\chi$ has finite index in $\Gamma$ 
and since $\alpha \otimes \chi$ is finite-dimensional, 
the induced representation $\Ind_{H_\chi \ltimes N}^G (\alpha \otimes \chi)$ is finite-dimensional.

\vskip.2cm

$\bullet$ \emph{Second step.} 
We claim that $\Phi$ factorizes to an injective map
from the space of $H$-orbits in $X_{\rm fd}$
to $\widehat G_{\rm fd}$.
\par

Indeed, let $(\chi, \alpha)$ and $ (\chi', \alpha')$ in $X_{\rm fd}$.
\par

Assume that $(\chi, \alpha)$ and $ (\chi', \alpha')$ are in the same $H$-orbit. 
Then $\chi$ and $\chi'$ are conjugate under $H$ and $\alpha = \alpha'$.
Therefore, $H_\chi = H_{\chi'}$ and the unitary characters
$\alpha\otimes \chi$ and $\alpha' \otimes \chi'$ of $H_\chi \ltimes N$
are conjugate under $H$. 
It follows that $\Ind_{H_\chi \ltimes N}^G (\alpha \otimes \chi)$
and $\Ind_{H_\chi \ltimes N}^G (\alpha' \otimes \chi')$
are equivalent (Proposition~\ref{PropConjIndRep}).
\par

Assume now that $(\chi, \alpha)$ and $ (\chi', \alpha')$ are not in the same $H$-orbit.
Two cases may occur:
\par

$\circ$ \emph{First case:}
$\chi$ and $\chi'$ are conjugate under $H$.
Then $H_\chi=H_{\chi'}$ and $\alpha\ne \alpha'$.
Therefore the unitary characters $\alpha \otimes \chi$ and $\alpha' \otimes \chi'$
of $H_\chi \ltimes N$ are not in the same $H$-orbit. 
It follows from the Mackey--Shoda criterion for equivalence of monomial representations 
(Corollary~\ref{Cor-NormalSubgEquiv}) 
that $\Ind_{H_\chi \ltimes N}^G (\alpha \otimes \chi)$
and $\Ind_{H_\chi \ltimes N}^G (\alpha' \otimes \chi')$ are not equivalent.
\par

$\circ$ \emph{Second case:}
$\chi$ and $\chi'$ are not conjugate under $H$.
By Proposition~\ref{PropConjIndRep}~\ref{iiDEPropConjIndRep},
the restriction of $\Ind_{H_\chi \ltimes N}^G (\alpha \otimes \chi)$ to $H_\chi \ltimes N$
is equivalent to a direct sum of conjugates of $\alpha \otimes \chi$;
similarly, the restriction of $\Ind_{H_{\chi'}\ltimes N}^G (\alpha' \otimes \chi')$
to $H_{\chi'}\ltimes N$ is equivalent
to a direct sum of conjugates of $\alpha' \otimes \chi'$.
Therefore, the restriction of $\Ind_{H_\chi \ltimes N}^G (\alpha \otimes \chi)$ to $N$
is equivalent to a direct sum of conjugates of $\chi$
and the restriction of $\Ind_{H_{\chi'} \ltimes N}^G (\alpha' \otimes \chi')$ to $N$
is equivalent to a direct sum of conjugates of $\chi'$. 
Since $\chi$ and $\chi'$ are not conjugate,
it follows that $\Ind_{H_\chi \ltimes N}^G (\alpha \otimes \chi) $
and $\Ind_{H_\chi' \ltimes N}^G (\alpha' \otimes \chi')$ cannot be equivalent.

\vskip.2cm

$\bullet$ \emph{Third step.}
We claim that $\Phi$ is surjective.
\par

Indeed, let $(\pi, \Hi)$ be a finite-dimensional irreducible representation
of $G$.
The restriction $\pi \vert_N$ is a finite-dimensional representation of the abelian
normal subgroup $N$. So, there exists a finite subset $S$ of 
$\widehat N$ such that $\Hi = \bigoplus_{\chi \in S} \Hi^{\chi}$,
where 
$$
\Hi^{\chi} \, := \, \left\{ \xi \in \Hi \mid \pi(n) \xi = \chi(n) \xi 
\hskip.5cm \text{for all} \hskip.2cm
n \in N \right\}
$$
is non-zero for every $\chi \in S$.
\par

Since $N$ is a normal subgroup of $G$, we have 
$$
\pi(g) \Hi^{\chi} \, = \, \Hi^{\chi^{g^{-1}}} 
\hskip.5cm \text{for all} \hskip.2cm
g \in G, \, \chi \in S.
\leqno{(*)}
$$
Therefore $S$ is $H$-invariant.
Since $\pi$ is irreducible, it follows that $S$ is a single $H$-orbit.
Since $S$ is finite, 
$S$ is the $H$-orbit of some $\chi \in {\widehat N}_{\rm per}$.
\par

By $(*)$, we have 
$$
\pi(h) \Hi^{\chi} \, = \, \Hi^{\chi}
$$
for every $h$ in the stabilizer $H_\chi$ of $\chi$ in $H$.
\par

Let $T$ be a transversal for the right coset space $(H_\chi \ltimes N)\backslash G$,
with $e \in T$.
Since $S$ is the $G$-orbit of $\chi$, we have 
$$
\Hi \, = \, \bigoplus_{t \in T} \pi(t)\Hi^{\chi}.
$$
This shows that $\pi$ is equivalent to the induced representation 
$\Ind_{H_\chi \ltimes N}^G \sigma$, where $\sigma$ is 
the subrepresentation of $\pi \vert_{H_\chi \ltimes N}$
defined on the $(H_\chi \ltimes N)$-invariant subspace
$\Hi^{\chi}$ (see Definition~\ref{openeasierthanclosed}).
\par

We claim that $\sigma = \alpha \otimes \chi $ for some $\alpha \in \widehat{H_\chi}$.
Indeed, $\sigma$ is an irreducible representation 
of $G_\chi = H_\chi \ltimes N$, since $\Ind_{H_\chi \ltimes N}^G \sigma \simeq \pi$
is irreducible (see Proposition~\ref{InductionQqPropr}.iii).
Moreover, $\sigma\vert_N$ is a multiple of $\chi$.
So, the claim follows from Lemma~\ref{Lem-FiniteDim}
and we conclude that $\pi$ is equivalent to
$\Ind_{H_\chi \ltimes N}^G(\alpha \otimes \chi) $.
\end{proof}

\section[Examples]
{All finite dimensional irreducible
% \\
representations for some groups}
% Section 4.B
\label{Section-FiniteDimRepForSomeGroups}

We are going to apply Theorem~\ref{Theo-FiniteDimRepSemiDirect}
to Heisenberg groups, affine groups,
solvable Baumslag--Solitar groups, and the lamplighter group. 

\subsection
{Heisenberg groups over rings}
% 4.B.a
\label{SS:HeisGpsOverRings}

Let $H(R)$ be the Heisenberg group over a unital commutative ring $R$,
as in Section~\ref{Section-IrrRepTwoStepNil}.
Recall that we may identify $H(R)$ with $R^3$, 
equipped with the group law 
$$
(a,b,c) (a', b', c') \, = \, (a + a', b + b', c + c' + ab').
$$
The centre of $H(R)$ is $Z = \{ (0,0,c) \mid c \in R \} \approx R$,
and $N = \{ (0,b,c) \mid b, c \in R \} \approx R^2$
is a maximal abelian subgroup of $H(R)$.
Observe that $H(R)$ is a semi-direct product 
$H \ltimes N$ for $H = \{ (a,0,0) \mid a \in R \} \approx R$.
\par

Let $\chi$ be a unitary character of $N$. Then $\chi = \chi_{\psi, \beta}$
for a unique pair $(\psi, \beta) \in (\widehat R)^2$,
where $\chi_{\psi, \beta}$ is defined by 
$$
\chi_{\psi, \beta}(0,b,c) \, = \, \beta(b) \psi(c)
\hskip.5cm \text{for} \hskip.2cm
b, c \in R.
$$ 
For $h = (a,0,0) \in H$, we have
$$
\chi_{\psi, \beta}^{h}(0,b,c) \, = \, \beta(b) \psi(ab) \psi(c)
\, = \, \chi_{\psi, \beta \psi^a} (0, b, c)
\hskip.5cm \text{for} \hskip.2cm
b, c \in R ,
$$
where $\psi^a \in \widehat R$ is defined by $\psi^a(b) = \psi(ab)$ for $b \in R$.
It follows that the $H$-orbit of $\chi_{\psi, \beta}$ is
$$
\{\chi_{\psi, \beta \psi^a} \mid a \in R \},
$$
and that the stabilizer of $\chi_{\psi, \beta}$, which only depends on $\psi$, is 
$$
H_\psi \, = \, \{ (a, 0, 0) \mid a \in I_\psi \},
$$
where $I_\psi$ is the ideal 
$$
I_\psi \, = \, \{ a \in R \mid aR \subset \ker \psi \}
$$
of $R$, already considered in Section~\ref{Section-IrrRepTwoStepNil}.
\par

Denote by ${\widehat R}_{\rm per}$ the set of $\psi \in \widehat R$
such that $I_\chi$ has \textbf{finite index} in $R$.
We see that the set ${\widehat N}_{\rm per}$
of unitary characters of $N$ with a finite $H$-orbit is 
$$
{\widehat N}_{\rm per} \, = \, 
\{\chi_{\psi, \beta} \mid
\psi \in {\widehat R}_{\rm per}, \hskip.1cm \beta \in \widehat R \}.
$$
\par

Recall from Lemma~\ref{Lem-FiniteDim} that, for $\alpha \in \widehat{I_\psi}$,
the formula
$$
\alpha \otimes \chi_{\psi, \beta} \, \colon \, (a,b,c) \mapsto \alpha(a) \beta(b) \psi(c)
$$
defines a unitary character of $H_\psi \ltimes N$.
Set
$$
X_{\rm fd} \, = \,
\bigsqcup_{\psi \in {\widehat R}_{\rm per}} \{ \psi \} \times \{( \alpha, \beta)
\mid \alpha \in \widehat{I_\psi}, \hskip.1cm \beta \in \widehat R \},
$$
with the $R$-action $X_{\rm fd} \curvearrowleft R$ on $X_{\rm fd}$ given by 
$$
X_{\rm fd} \times R \, \to \, X_{\rm fd},
\hskip.2cm
((\psi, \alpha, \beta), a) \, \mapsto \, (\psi, \alpha, \beta \psi^a).
$$
\par

In view of these remarks, the following result
is a direct consequence of Theorem~\ref{Theo-FiniteDimRepSemiDirect}. 

\begin{cor}
% 4.B.1
\label{Cor-FiniteDimRepFinDimRepHeis-Ring}
Let $\Gamma$ be the Heisenberg group $H(R)$ over a unital commutative ring $R$.
We keep the notation above.
The map $ X_{\rm fd} \to \widehat \Gamma_{\rm fd}$, defined by 
$$
(\psi, \alpha, \beta) \, \mapsto \,
\Ind_{H_\psi \ltimes N}^\Gamma (\alpha \otimes \chi_{\psi, \beta}),
$$
induces a bijection between the space of $R$-orbits in $X_{\rm fd}$
and the space $\widehat \Gamma_{\rm fd}$
of equivalence classes of finite-dimensional irreducible representations of $\Gamma$.
\end{cor}

$\bullet$
First particular case: $R = \K$ is a field.
As already noted in Section \ref{Section-IrrRepTwoStepNil},
for $\psi \in \widehat \K \smallsetminus \{1_\K\}$,
we have $I_\psi = \{0\}$ and hence $H_\psi = \{0\}$;
for $\psi = 1_\K$, we have $I_\psi = \K$ and hence $H_\psi = H$.
\par

We will need the following lemma.

\begin{lem}
% 4.B.2
\label{Lem-FinDimRepHeis-finiteFields}
Let $\K$ be a field and $\psi \in \widehat \K \smallsetminus \{1_\K\}$.
\par
Then $\{\psi^a \mid a \in \K\}$ is a dense subgroup of $\widehat \K$,
where $\psi^a$ is defined by $\psi^a(x) = \psi(ax)$.
In particular, $\{\psi^a \mid a \in \K\} = \widehat \K$ when $\K$ is finite.
\end{lem}

\begin{proof}
It is clear that $L := \{\psi^a \mid a \in \K\}$ is a subgroup of $\widehat \K$.
By Pontrjagin duality, it suffices to show that $L^\perp = \{0\}$,
where $L^\perp = \bigcap_{\chi \in L} \ker \chi$.
\par

Let $x \in L^\perp$.
Then $\psi(ax) = 1$ for all $a \in \K$ and hence $x = 0$;
indeed, otherwise $\{ax \mid a \in \K\} = \K$
and this would contradict the assumption that $\psi \ne 1_\K$.
\end{proof}

$\circ$
Assume first that $\K$ is infinite.
Then $\widehat{\K}_{\rm per} = \{1_\K\}$
and the following result is a direct consequence of
Corollary~\ref{Cor-FiniteDimRepFinDimRepHeis-Ring}.

\begin{cor}
% 4.B.3
\label{Cor-FinDimRepHeis-InfiniteFields}
Let $\Gamma = H(\K)$ for an infinite field $\K$.
\par

The unitary characters of $\Gamma$
are the only finite-dimensional irreducible representations of $\Gamma$,
that is, the map 
$$
\widehat{\Gamma/Z} \, \to \, \Gamma_{\rm fd},
\hskip.5cm
\chi \, \mapsto \, \chi \circ p
$$
is a bijection,
where $Z$ denotes the centre of $\Gamma$
and $p \, \colon \Gamma \twoheadrightarrow \Gamma/Z$ the canonical projection.
\end{cor}

$\circ$
Assume now that $\K$ is finite. Let $\psi \in \widehat \K \smallsetminus \{1_\K\}$.
Recall that, for $\beta_0 \in \widehat \K$,
the $H$-orbit of $\chi_0 = \chi_{\psi, \beta_0}$ is
$$
\{\chi_{\psi, \beta_0 \psi^a} \mid a \in \K\}.
$$
Therefore the $H$-orbit of $\chi_0$
coincides with $\{\chi_{\psi, \beta} \mid \beta \in \widehat \K \}$,
by Lemma~\ref{Lem-FinDimRepHeis-finiteFields}.
In particular, $\chi_{\psi, 1_\K}$ belongs to the $H$-orbit of $\chi_0$.
\par

In view of these remarks,
the following result is again a direct consequence
of Corollary~\ref{Cor-FiniteDimRepFinDimRepHeis-Ring}.
 
\begin{cor}
% 4.B.4
\label{Cor-FinDimRepHeis-finiteFields}
Let $\Gamma = H(\K)$ for a finite field $\K$.
Let $Z$ denote the centre of $\Gamma$
and $p \, \colon \Gamma \twoheadrightarrow \Gamma/Z$ the canonical projection.
\par

The map 
$$
\left( \widehat \K \smallsetminus \{1_\K\} \right)
\sqcup
\widehat{\Gamma/Z} \, \to \, \widehat \Gamma,
$$
defined by 
$$
\psi \, \mapsto \, \Ind_{N}^\Gamma \chi_{\psi, 1_\K}
\hskip.5cm \text{for} \hskip.2cm
\psi \in \widehat \K \smallsetminus \{1_\K\}
$$
and 
$$
\chi \, \mapsto \, \chi \circ p
\hskip.5cm \text{for} \hskip.5cm
\chi \in \widehat{\Gamma/Z},
$$
is a bijection.
\index{Heisenberg group! $2$@$H(\K)$ over a field $\K$}
\end{cor}

\begin{rem}
% 4.B.5
\label{Rem-RepHeisenbergFini}
Let $\Gamma = H(\K)$ for a finite field $\K$
and let $\psi \in \widehat Z \smallsetminus \{1_Z\}$.
Using the model of Construction~\ref{constructionInd}(2),
we give an explicit formula for the irreducible 
representation $\pi_\psi := \Ind_{N}^\Gamma \chi_{\psi, 1_\K}$ 
of $\Gamma$ of Corollary~\ref{Cor-FinDimRepHeis-finiteFields}:
for every $(a, b, c)$ in $\Gamma$,
the operator $\pi_\psi(a,b,c)$ acts on $\ell^2(\K)$ by 
$$
(\pi_\psi(a,b,c) f) (x) \, = \, \psi(c) \psi(bx) f(x + a)
\hskip.5cm \text{for} \hskip.2cm
f \in \ell^2(\K), \hskip.1cm x \in \K.
$$
Observe that a consequence of Corollary~\ref{Cor-FinDimRepHeis-finiteFields}
is that any two irreducible representations of $\Gamma$
with central character $\psi$ are equivalent;
this result may be viewed as an analog for Heisenberg groups over finite fields
of the classical Stone--von Neumann theorem on the uniqueness of the commutation relations (see \cite[Theorem 6.50]{Foll--16}, and Remark \ref{dualHeisLie}).
\end{rem}

$\bullet$
Second particular case: $R = \Z$ is the ring of rational integers.
Recall that $\widehat Z \approx \widehat \Z =
\{\psi_\theta \mid \theta \in \mathopen[ 0,1 \mathclose[ \}$,
with $\psi_\theta(n) = e^{2 \pi i n \theta}$ for $n \in \Z$.
\par

Let $\psi = \psi_\theta$. We have $I_\psi = \{0\}$
if $\psi$ has infinite order, that is, if $\theta$ is irrational.
We have $I_\psi = n\Z$ if $\psi$ has order $n \ge 1$, that is, 
if $\theta = p/n$ for $p \in \{0, 1, \hdots, n - 1 \}$, with $p$ and $n$ coprime.
As a consequence, denoting by $\widehat{Z}_n$
the elements of order $n \ge 1$ in $\widehat Z$, we have
$$
\widehat{\Z}_{\rm per} \, = \, \bigsqcup_{n \ge 1} \widehat{Z}_n.
$$
\par

Fix $n \ge 1$ and $\psi = \psi_\theta \in \widehat{Z}_n$. We have 
$$
H_\psi \, = \, \{(a,0,0)\mid a \in n\Z\}.
$$
So, the stabilizer $H_\psi \ltimes N$ of $\psi$ in $\Gamma = H(\Z)$
only depends on $n$ and coincides with the normal subgroup
$$
\Gamma(n) \, = \, \{(a,b,c)\mid a \in n\Z, b, c \in \Z\},
$$
which was already introduced in Section~\ref{Section-IrrRepTwoStepNil}.
Every unitary character $\chi$ of $\Gamma(n)$ with $\chi \vert_Z = \psi_\theta$
is of the form $\chi_{\theta, \alpha, \beta}$ for a unique pair
$(\alpha, \beta) \in \mathopen[ 0,1/n \mathclose[ \times \mathopen[ 0,1 \mathclose[$, where 
$$
\chi_{\theta, \alpha, \beta}(a,b,c) \, = \, e^{2 \pi i (\alpha a + \beta b + \theta c)} 
\hskip.5cm \text{for} \hskip.2cm
a \in n\Z, \hskip.1cm b \in \Z, \hskip.1cm c \in \Z.
$$
The $H$-orbit of $\chi_{\theta, \alpha, \beta}$ is 
$$
\{ \chi_{\theta, \alpha, \beta + a \theta} \in \widehat N \mid a \in \Z \}
\, = \,
\{ \chi_{\theta, \alpha, \beta + \frac{k}{n}} \in \widehat N \mid k \in \Z\} 
$$
and so contains a character $\chi_{ \theta, \alpha, \beta'}$ for a unique
$\beta' \in \mathopen[ 0,1/n \mathclose[$.
Set
$$
X \, = \,
\bigsqcup_{n \ge 1}
\bigsqcup_{\psi \in \widehat{Z}_n} 
\{ \psi \} \times \mathopen[ 0,1/n \mathclose[^2 .
$$

\par

The following result is therefore a direct consequence of
Corollary~\ref{Cor-FiniteDimRepFinDimRepHeis-Ring}.

\begin{cor}
% 4.B.6
\label{Cor-FinDimRepHeis-Integers}
Let $\Gamma = H(\Z)$.
We keep the notation above.
\par

For each $n \ge 1$ and $(\psi, \alpha, \beta) \in X$,
with $\psi = \psi_\theta \in \widehat \Z_n$, 
the representation $\Ind_{\Gamma(n)}^\Gamma \chi_{\theta, \alpha, \beta}$
of $\Gamma$ is irreducible of dimension $n$.
\par

The map $X \to \widehat \Gamma_{\rm fd}$
defined by 
$
(\psi, \alpha, \beta) \mapsto \Ind_{\Gamma(n)}^\Gamma \chi_{\theta, \alpha, \beta}
$
is a bijection.
\end{cor}

This corollary shows that the set $\widehat{H(\Z)}_{\rm fd}$
can be seen as a disjoint union
parametrized by $\mathopen[ 0,1 \mathclose[ \cap \Q$
of two-dimensional tori.

\subsection
{Affine group of arbitrary fields}
% 4.B.b

Let $\Gamma = \Aff(\K)$ be the affine group of a field $\K$,
as in Section~\ref{Section-IrrRepAff}.
Recall that $\Gamma = \K^\times \ltimes \K$,
for the action of $\K^\times$ on $\K$ defined by
$(a, t) \mapsto at$ for all $a \in \K^\times$ and $t \in \K$.
We set
$$
H \, = \, \K^\times
\hskip.5cm \text{and} \hskip.5cm
N \, = \, \K .
$$

\begin{cor}
% 4.B.7
\label{Cor-FinDimRepAffine}
Let $\Gamma = \Aff(\K)$ the group of affine transformations of a field $\K$. 
\begin{enumerate}[label=(\arabic*)]
\item\label{1DECor-FinDimRepAffine}
Assume that $\K$ is infinite. The unitary characters of $\Gamma$
are the only finite-dimensional irreducible representations of $\Gamma$,
that is, the map 
$$
\widehat{\Gamma/N} \, \to \, \widehat \Gamma_{\rm fd},
\hskip.2cm \chi \, \mapsto \, \chi \circ p
$$
is a bijection,
where $p \, \colon \Gamma \twoheadrightarrow \Gamma/N$ is the canonical projection.
\item\label{2DECor-FinDimRepAffine}
Assume that $\K$ is finite of cardinality $q$.
Fix any $\chi \in \widehat \K \smallsetminus \{1_\K\}$.
Up to equivalence, the irreducible representations of
$\Gamma$ are the representations of dimension one
and the representation $\Ind_N^\Gamma \chi$ of dimension $q-1$.
\end{enumerate}
\end{cor}

\begin{proof}
For $\chi \in \widehat \K \smallsetminus \{ 1_\K \}$, we have $H_\chi = \{ 1 \}$
by Lemma~\ref{Lem-FreeActDualAff}~\ref{iDELem-FreeActDualAff}.

\vskip.2cm

In particular, if $\K$ is infinite, then ${\widehat N}_{\rm per} = \{1_N\}$ 
and \ref{1DECor-FinDimRepAffine} is a direct consequence
of Theorem~\ref{Theo-FiniteDimRepSemiDirect}.

\vskip.2cm

\ref{2DECor-FinDimRepAffine}
If $\K$ is finite of cardinality $q$,
then $\Gamma$ has $\vert \Gamma / [\Gamma, \Gamma] \vert = q-1$
representations of degree $1$,
and one representation of degree $q-1$
denoted by $\Ind_N^\Gamma \chi$ (where $\chi \ne 1_\K$)
in Proposition \ref{Prop-IrredRepAffine-bis}.
Since
$$
(q-1)^2 + (q-1)\times 1^2 \, = \, (q-1)q \, = \, \vert \Gamma \vert ,
$$
any irreducible representation of $\Gamma$
is equivalent to one of these (see \cite[Proposition 15.4.1]{Dixm--C*}).
\end{proof}

\subsection
{Solvable Bausmslag--Solitar group}
% 4.B.c

Let $p$ be a prime.
Recall from Section~\ref{Section-IrrRepBS}
that the Baumslag--Solitar $\Gamma = \BS(1, p)$ is the semi-direct product
$$
\BS(1, p) \, = \, A \times N ,
\hskip.2cm \text{where} \hskip.2cm
A = \begin{pmatrix} p^\Z & 0 \\ 0 & 1 \end{pmatrix} \approx \Z
\hskip.2cm \text{and} \hskip.2cm
N = \begin{pmatrix} 1 & \Z[1/p] \\ 0 & 1 \end{pmatrix} \approx \Z[1/p] ,
$$
and the generator $\begin{pmatrix} p & 0 \\ 0 & 1 \end{pmatrix}$ of $A$
acts on $N$ by multiplication by $p$.
We often denote by $(k, b)$ the element
$\begin{pmatrix} p^k & b \\ 0 & 1 \end{pmatrix}$ of $\Gamma$,
where $k \in \Z$ and $b \in \Z[1/p]$. 
\index{Baumslag--Solitar group $\BS(1, p)$}
\par

Recall also that we can identify $\widehat N$ with the $p$-adic solenoid 
$$
\So_p \, = \, (\Q_p \times \R) / \Delta
\hskip.5cm \text{for} \hskip.2cm 
\Delta \, = \, \{(a,-a)\mid a \in \Z[1/p]\};
$$
the action of $A$ on $\widehat N$ corresponds to 
the action of $\Z$ on $\So_p$
for which the generator $1 \in \Z$ acts by
the transformation $T_p \, \colon \So_p \to \So_p$
of multiplication by $p$.
For $s \in \So_p$, denote by $\chi_{s}$ the corresponding unitary character of $N$
under the identification $\widehat N \approx \So_p$.
The set $\widehat N_{\rm per}$
of characters of $N$ with finite $A$-orbit coincides with 
$$
\{\chi_s\mid s\in \Per(T_p)\},
$$
where $\Per(T_p)$ is the set of $T_p$-periodic points in $\So_p$. 
\par

For $n \in \N^*$, let $s \in \So_p(n)$,
where $\So_p(n)$ is the set of elements in $\So_p$ with $T_p$-period $n$.
The stabilizer of $\chi_s$ in $A$ is 
$$
A(n) \, = \, 
\left\{(k,0) \in \BS(1, p) \mid k \in n\Z \right\} \, \approx \, n\Z .
$$
and so the stabilizer of $\chi_s$ in $\Gamma$ is the normal subgroup
$$
\Gamma(n) \, := \, A(n) \ltimes N \, = \,
\left\{ (k, b) \in \Gamma \mid k \in n\Z, \hskip.1cm b \in \Z[1/p] \right\}
$$
of index $n$ in $\Gamma$.
We parametrize the dual of $n\Z$ by $\mathopen[ 0,1/n \mathclose[$
through the map 
$$
\mathopen[ 0,1/n \mathclose[ \, \to \, \widehat{n\Z}, 
\hskip.2cm 
\theta \, \mapsto \, \chi_\theta,
$$
with $\chi_\theta$ defined by $\chi_\theta(k) = e^{2 \pi i \theta k}$
for all $k \in n\Z$
\par 

Every unitary character $\chi$ of $\Gamma(n)$ with $\chi \vert_N = \chi_s$ 
is of the form $\chi_{s, \theta}$ for a unique $\theta \in \mathopen[ 0,1/n \mathclose[$, 
where $\chi_{s, \theta} \, \colon \Gamma(n) \to \T$ is defined by 
$$
\chi_{s, \theta} (k, b) \, = \, \chi_\theta(k)\chi_s(b)
\hskip.5cm \text{for all} \hskip.2cm
(k, b) \in \Gamma(n).
$$
\par

We obtain the following direct consequence of Theorem~\ref{Theo-FiniteDimRepSemiDirect}.

\begin{cor}
% 4.B.8
\label{Cor-FinDimRepBS}
Let $p$ be a prime and
$\Gamma = \BS(1, p)$ be corresponding Baumslag--Solitar group.
Set
$$
X_{\rm fd} \, = \, 
\bigsqcup_{n \ge 1} \hskip.1cm
\big( \So_p(n) \times \mathopen[ 0,1/n \mathclose[ \big)
$$
and let $\widetilde{T_p} \, \colon X_{\rm fd} \to X_{\rm fd}$ be defined by
$\widetilde{T_p} (s, \theta) = (T_ps, \theta)$.
\par

The map $X_{\rm fd} \to \widehat \Gamma_{\rm fd}$, defined by 
$$
(s, \theta) \, \mapsto \, \Ind_{\Gamma(n)}^\Gamma \chi_{s, \theta} 
\hskip.5cm \text{for} \hskip.2cm
n \ge 1, \hskip.1cm (s, \theta) \in \So_p(n) \times \mathopen[ 0,1/n \mathclose[,
$$
induces a bijection
between the space of $\widetilde{T_p}$-orbits in $X_{\rm fd}$
and the space $\widehat \Gamma_{\rm fd}$
of equivalence classes of finite-dimensional irreducible representations of $\Gamma$.
\end{cor}

\subsection
{Lamplighter group}
% 4.B.d

Recall from Section~\ref{Section-IrrRepLamplighter}
that the lamplighter group is the semi-direct product $\Gamma = A \ltimes N$
of $A = \Z$ with $N = \bigoplus_{k \in \Z} \Z / 2 \Z$,
where the action of $\Z$ on $\bigoplus_{k \in \Z} \Z / 2 \Z$
is given by shifting the coordinates.
Recall also that $\widehat N$ can be identified
with $X = \prod_{k \in \Z} \{0,1 \}$,
the dual action of $\Z$ on $\widehat N$ being given
by the shift transformation $T$ on $X$.
\par

For $x \in X$, we denote by $\chi_x$ the corresponding element in $\widehat N$
under the identification $\widehat N \approx X$.
The set $\widehat N_{\rm per}$ of characters of $N$ with finite $A$-orbits coincides with 
$$
\{\chi_x \mid x \in \Per(T)\},
$$
where $\Per(T)$ is the set of $T$-periodic points in $X$. 
\par

For $n \in \N^*$, let $x \in X(n)$,
where $X(n)$ is the set of elements in $X$ with $T$-period~$n$.
The stabilizer of $\chi_x$ in $A$ is $n\Z$ 
and so the stabilizer of $\chi_x$ is the normal subgroup
$$
\Gamma(n) \, = \, n\Z \ltimes N
$$
of index $n$ in $\Gamma$.
\par

Parametrizing the dual of $n\Z$ by $\mathopen[ 0,1/n \mathclose[$ as before,
every unitary character $\chi$ of $\Gamma(n)$ with $\chi \vert_N = \chi_x$ 
is of the form $\chi_{x, \theta}$ for a unique $\theta \in \mathopen[ 0,1/n \mathclose[$, 
where $\chi_{x, \theta} \, \colon \Gamma(n) \to \T$ is defined by 
$$
\chi_{x, \theta} (k, b) \, = \, \chi_\theta(k)\chi_x(b)
\hskip.5cm \text{for all} \hskip.2cm
(k, b) \in \Gamma(n).
$$
\par

We obtain the following direct consequence of Theorem~\ref{Theo-FiniteDimRepSemiDirect}.

\begin{cor}
% 4.B.9
\label{Cor-FinDimRepLamplighter}
Let $\Gamma = A \ltimes N$ be the lamplighter group,
Set
$$
X_{\rm fd} \, = \, 
\bigsqcup_{n \ge 1} \hskip.1cm
\big( X(n) \times \mathopen[ 0,1/n \mathclose[ \big)
$$
and let $\widetilde T \, \colon X_{\rm fd} \to X_{\rm fd}$ be defined by
$\widetilde T(x, \theta) = (T x, \theta)$.
\par

The map $ X_{\rm fd} \to \widehat \Gamma_{\rm fd}$, defined by 
$$
(s, \theta)\, \mapsto \, \Ind_{\Gamma(n)}^\Gamma \chi_{x, \theta} 
\hskip.5cm \text{for} \hskip.2cm
n \ge 1, \hskip.1cm (x, \theta) \in X(n) \times \mathopen[ 0,1/n \mathclose[,
$$
induces a bijection
between the space of $\widetilde T$-orbits in $X_{\rm fd}$
and the space $\widehat \Gamma_{\rm fd}$
of finite-dimensional irreducible representations of $\Gamma$.
\end{cor}

\section[Classes of groups in terms of fd representations]
{Classes of groups
% \\
in terms of fd representations}
% Section 4.C
\label{SectionFdrep}

Given a topological group $G$, one may ask
how large is the set $\widehat G_{\rm fd}$ compared to~$\widehat G$.
More specifically, consider the following properties that $G$ may possess:
\begin{enumerate}
\item[$\bullet$]
every irreducible representation of $G$ is finite dimensional, that is, 
$\widehat G = \widehat G_{\rm fd}$;
\item[$\bullet$]
$G$ has sufficiently many finite dimensional representations,
in the sense that $\widehat G_{\rm fd}$ separates the points of $G$;
\item[$\bullet$]
$G$ has only trivial finite dimensional representations, that is,
$\widehat G_{\rm fd} = \{ 1_G \}$;
\item[$\bullet$]
the set $\widehat G_{\rm fd}$ is dense in $\widehat G$.
\end{enumerate}
We will discuss in the following subsections
the classes of LC groups singled out by either one of the three first properties.
The class of discrete groups with the fourth property
will be treated in Section \ref{SS:Faithful}.

\subsection
{Moore groups}
% subsection 4.C.a
\label{SS:MooreGr}

\index{Moore group}
A locally compact group $G$ is a \textbf{Moore group} if all irreducible representations 
of $G$ are finite-dimensional. 
LCA groups and compact groups are clearly Moore groups.
Other obvious examples are products $K \times A$,
where $K$ is a compact group and $A$ a LCA group.
More generally,
using the Mackey machine (see Remark~\ref{MackMach}),
it can be shown that a central-by-compact LC group is a Moore group.
As indicated by the following result,
the class of Moore groups is a rather restricted class of LC groups.
\par

For the next theorem, a Lie group is a LC group in which
the connected component of the identity is open, and is a connected Lie group.

\begin{theorem}
% 4.C.1
\label{Theo-MooreGroups}
\begin{enumerate}[label=(\arabic*)]
\item\label{1DElistegroupesMoore}
A Lie group is a Moore group if and only if it has an open subgroup 
of finite index that is central-by-compact.
% (i) = Theorem 2 in \cite{Moor--72}
\item\label{2DElistegroupesMoore}
A connected locally compact group is a Moore group
if and only if it is a Cartesian product $K \times \R^n$, with $K$ compact and $n \ge 0$.
% Pour $G$ almost connected, Moore = MAP,
% donc $G = K \times R^n$ par Freudenthal-Weil
\item\label{3DElistegroupesMoore}
A locally compact group is a Moore group if and only if it is a projective limit
of Moore Lie groups. 
% (iii) = Theorem 3 in \cite{Moor--72}
In particular: 
\item\label{4DElistegroupesMoore}
A discrete group is a Moore group
if and only if it has an abelian subgroup of finite index.
% (iv) = Thoma, voir Moore page 402 (c'est 2e page).
\item\label{5DElistegroupesMoore}
The class of Moore groups is stable under subgroups, quotient groups,
inverse limits, and finite extensions.
\end{enumerate}
\end{theorem}

\begin{proof}[References for the proof]
See \cite{Moor--72}, as well as \cite{Robe--69}.
The article by Moore was available as a preprint
before the article by Robertson was submitted. 
\end{proof}

\begin{rem}
% 4.C.2
\label{Rem-Theo-MooreGroups}
(1)
Note that \ref{1DElistegroupesMoore} in Theorem~\ref{Theo-MooreGroups}
\emph{does not} carry over to locally compact groups:
the example in \cite[Section 5]{Moor--72} is a Moore group
without any central-by-compact finite index subgroup.

\vskip.2cm

(2)
Note also that an abelian-by-compact locally compact group need not be a Moore group.
For example, the Euclidean group $\SO(2) \ltimes \R^2$
has a well-known family of infinite dimensional irreducible representations
indexed by the orbits of $\SO(2)$ in $\widetilde{\R^2} \smallsetminus \{0\}$,
i.e., indexed by $\R^\times_+$;
see, e.g., \cite[Pages 195--196]{Mack--76}.

\vskip.2cm

(3)
Concerning \ref{1DElistegroupesMoore} and \ref{2DElistegroupesMoore}
in Theorem~\ref{Theo-MooreGroups},
note that a connected locally compact group $G$
is central-by-compact if and only if
$G$ is isomorphic to a direct product $K \times \R^n$
for some compact group $K$ and some integer $n \ge 0$
(see \cite[16.4.6]{Dixm--C*} and Example \ref{Ex-MAP-Groups} below).
\end{rem}

\subsection[Maximally almost periodic groups and Bohr compactification]
{Maximally almost periodic groups
\\
and Bohr compactification}
% subsection 4.C.b
\label{SS:MAP}

\index{Maximally almost periodic group, or MAP group}
A topological group $G$ is \textbf{maximally almost periodic}, or \textbf{MAP},
if its finite-dimensional representations separate its points,
i.e., if for any $g \ne e$ in $G$ 
there exists a finite-dimensional representation $\pi$ of $G$
in some Hilbert space $\Hi$ such that $\pi(g) \ne \mathrm{Id}_\Hi$.
This class is strictly larger than the class of Moore groups.
It was introduced by von Neumann \cite{vNeu--34};
see also \cite[Chap.~VII]{Weil--40} and \cite[\S~16]{Dixm--C*}.
\par

We will put this class in perspective with Proposition \ref{Pro-AP-MAP}.
For this, we need the notion of Bohr compactification.

\begin{theorem}
% 4.C.3
\label{Theo-Bohr}
Let $G$ be a topological group.
\par

There exists a compact group $K$
and a continuous homomorphism $\beta \, \colon G \to K$ with dense image,
with the following universal property:
for every compact group $L$ and every continuous homomorphism $\alpha \, \colon G \to L$,
there exists a continuous homomorphism $\alpha' \, \colon K \to L$
such that $\alpha = \alpha' \circ \beta$.
\par

The pair $(K, \alpha)$ is unique in the following sense:
let $(K', \beta')$ be a pair consisting of a compact group 
and a continuous homomorphism $\beta' \, \colon G \to K'$ with dense image
satisfying the same universal property.
Then there exists an isomorphism $\alpha \, \colon K \to K'$ of topological groups
such that $\beta' = \alpha \circ \beta$.
\end{theorem}

\begin{proof}
$\bullet$ {\it Existence part.} 
Here is one way to prove the existence of $K$.
Choose a family $(\pi_i, \Hi_i)_{i \in I}$ of representatives
for the set of equivalence classes 
% in $\widehat G_{\rm fd}$.
of finite dimensional representations of $G$.
Let 
$$
\beta \, \colon \, G \to \prod_{i \in I} \U(\Hi_i),
\hskip.5cm
g \mapsto \bigoplus_{i \in I} \pi_i(g) .
$$
Let $K$ be the closure of $\beta(G)$ in the compact group $\prod_{i \in I} \U(\Hi_i)$.
Then $\beta$ is a continuous homomorphism and $K$ is a compact group.
\par

Let $L$ be a compact group and $\alpha \, \colon G \to L$ a continuous homomorphism. 
Let $(\sigma_j, \Ki_j)_{j \in J}$ be a family of representatives
for the set of equivalence classes
% in $\widehat L$.
of finite dimensional representations of $L$.
By Peter--Weyl theory,
% every $\Ki_j$ is finite dimensional and
the continuous homomorphism 
$$
\prod_{j \in J} \sigma_j \, \colon \, L \to \prod_{j \in J} \U(\Ki_j)
$$
is injective.
So, $L$ is topologically isomorphic to its image
in the compact group $\widetilde L := \prod_{j \in J} \U(\Ki_j)$.
We may therefore assume that $L$ is a closed subgroup of $\widetilde L$.
\par

Let $j \in J$ and let $p_j \, \colon \widetilde L \twoheadrightarrow \U(\Ki_j)$
be the canonical projection.
The composition $p_j \circ \alpha \, \colon G \to \U(\Ki_j)$ is 
% an irreducible 
a finite dimensional representation of $G$.
Therefore $p_j \circ \alpha$ is equivalent to $\pi_{i(j)}$ for a unique $i(j) \in I$.
Let $V_j \, \colon \Hi_{i(j)} \to \Ki_j$ be a unitary operator
intertwining $\pi_{i(j)}$ and $p_j \circ \alpha$, so that 
$$
V_j \pi_{i(j)} (g) V_j ^{-1} \, = \, p_j \circ \alpha(g)
\hskip.5cm \text{for all} \hskip.2cm
g \in G.
$$
\par

Let $\alpha' \, \colon K \to \widetilde L$ be the continuous homomorphism defined by 
$$
\alpha' \big( (U_i)_{i \in I} \big) \, = \, \big( V_j U_{i(j)} V_j ^{-1} \big)_{j \in J} 
$$
for $(U_i)_{i \in I} \in K \subset \prod_{i \in I} \U(\Hi_i)$.
Then
$$
\alpha'(\beta(g)) \, = \, \big( V_j \pi_{i(j)}(g) V_j ^{-1} \big)_{j \in J} \, = \,
\big( p_j \circ \alpha(g) \big)_{j \in J} \, = \, \alpha(g)
$$
for every $g \in G$.
Therefore $\alpha = \alpha' \circ \beta$.

\vskip.2cm

$\bullet$ {\it Uniqueness part.}
Let $(K', \beta')$ be a pair consisting of a compact group 
and a continuous homomorphism $\beta' \, \colon G \to K'$ with dense image
and satisfying the same universal property as $(K, \beta)$.
Then, by the universal property, there exist continuous homomorphisms
$\alpha \, \colon K \to K'$ and $\alpha' \, \colon K' \to K$
such that $\beta' = \alpha \circ \beta$ and $\beta = \alpha' \circ \beta'$.
Then
$$
\alpha'(\alpha(\beta(g))) \, = \, \alpha'(\beta'(g)) \, = \, \beta(g)
\hskip.5cm \text{for all} \hskip.2cm
g \in G
$$
and hence $\alpha' \circ \alpha = {\mathrm{Id}}_K$, since $\beta(G)$ is dense in $K$.
Similarly, we have $\alpha \circ \alpha' = {\mathrm{Id}}_{K'}$
and this shows that $\alpha$ is an isomorphism.
\end{proof}

\index{Bohr compactification}
The compact group associated to $G$ as in Theorem~\ref{Theo-Bohr}
is denoted by $\Bohr(G)$ and the corresponding homomorphism 
by $\beta_B \, \colon G \to \Bohr(G)$; 
the group $\Bohr(G)$, or more precisely 
the pair $(\Bohr(G), \beta_B)$,
is called the \textbf{Bohr compactification} of $G$.
It follows from the definition that $\Bohr(\cdot)$ is a covariant functor:
every continuous group homomorphism $G_1 \to G_2$
gives rise to a continuous homomorphism
$\Bohr(G_1) \to \Bohr(G_2)$.

\vskip.2cm

We are going to describe the Bohr compactification of a LCA group.
For a topological group $G$, we denote by $G_{\rm disc}$
the group $G$ made discrete (that is, viewed as a discrete group).
\par

Let $G$ be a LCA group.
Denote by $\Gamma = \widehat G_{\rm disc}$
the dual of $G$, viewed as a discrete group.
Observe that $\widehat \Gamma$ is a compact group,
and that we have an injective continuous homomorphism
$i \, \colon G \to \widehat \Gamma, \hskip.1cm g \mapsto \widehat g$,
defined by $\widehat g (\chi) = \chi(g)$ for all $g \in G$ and $\chi \in \Gamma$.

\begin{prop}
% 4.C.4
\label{Prop-BohrAbelian}
Let $G$ be a locally compact abelian group
and let $\Gamma = {\widehat G}_{\rm disc}$ be the dual group of $G$ made discrete.
\par

The Bohr compactification of $G$ coincides with the compact group $\widehat \Gamma$,
together with the embedding $i \, \colon G \to \widehat \Gamma$ as above.
\end{prop}

\begin{proof}
Observe first that $i(G)$ is dense in $\widehat \Gamma$.
Indeed, by Pontrjagin duality,
every element in $\widehat{\widehat \Gamma}$ is of the form
$\hat{\chi}$ for some $\chi \in \Gamma = {\widehat G}_{\rm disc}$. 
Therefore the annihilator of $i(G)$ in $\widehat{\widehat \Gamma}$
consists of the unit character $1_G$. 
Therefore, $i(G)$ is dense in $\widehat \Gamma$, by Pontrjagin duality.
\par

To check that the pair $(\widehat \Gamma, i)$
has the universal property of a Bohr compactification,
let $K$ be a compact group and $\alpha \, \colon G \to K$ a continuous homomorphism.
The closure $L$ of $\alpha(G)$ in $K$ is a compact abelian group.
\par

Let $\alpha^t \, \colon \widehat L \to \widehat G$ be the continuous homomorphism
obtained by duality from $\alpha \, \colon G \to L$.
Observe that, since $L$ is compact, we can view $\alpha^t$ as a homomorphism
from the discrete groups $\widehat L$
to the discrete group $\Gamma= {\widehat G}_{\rm disc}$.
Identifying $L$ with $\widehat{\widehat L}$,
we can view the dual homomorphism of $\alpha^t$
as a continuous homomorphism $\alpha' \, \colon \widehat \Gamma \to L$.
For every $g \in G$ and $\chi \in \widehat L$, we have 
$$
\chi(\alpha'( i (g))) \, = \, i (g) (\alpha^t(\chi)) \, = \, \chi(\alpha(i(g))).
$$
This shows that $\alpha' \circ i = \alpha$. 
\end{proof}

Let $G$ be a topological group.
The kernel of the continuous homomorphism
$\beta_B \, \colon G \to \Bohr(G)$
is the \textbf{von Neumann kernel} $N_G$ of $G$~;
it is a topologically characteristic closed subgroup of $G$.
There is a natural bijection from the subset $\widehat G_{\rm fd}$
of $\widehat G$ introduced above
onto the dual $\widehat{\Bohr(G)}$
of the Bohr compactification of $G$.
(Observe however that
the Fell topology on $\widehat{\Bohr(G)}$ is discrete,
but that $\widehat G_{\rm fd}$ in general
is not a discrete subspace of $\widehat G$.)
\index{von Neumann kernel of a topological group}
\par

Next, we characterize the functions on $G$
which ``extend" to continuous functions on $\Bohr(G)$.
For a proof, we refer to \cite[16.2.1]{Dixm--C*}.

\begin{theorem}
% 4.C.5
\label{Theo-AP}
Let $G$ be a topological group. 
For a continuous function $f \, \colon G \to \C$, the following properties are equivalent:
\begin{enumerate}[label=(\roman*)]
\item\label{iDETheo-AP}
there exists a continuous function $f' \, \colon \Bohr(G) \to \C$
such that $f = f' \circ \beta_B$;
\item\label{iiDETheo-AP}
the set of left translates $\{ g \mapsto f(a^{-1}g) \}_{a \in G}$
is relatively compact in the Banach space 
$C^b(G)$ of bounded continuous functions on $G$;
\item\label{iiiDETheo-AP}
$f$ is a uniform limit of linear combinations of coefficients of
finite-dimen\-sional irreducible representations of $G$.
\end{enumerate}
\end{theorem}

A function on the topological group $G$
which satisfies the equivalent properties of Theorem~\ref{Theo-AP}
is called \textbf{almost periodic},
and the space $AP(G)$ of these is a closed subalgebra
of the Banach algebra $C^b(G)$,
for the norm $f \mapsto \sup_{g \in G} \vert f(g) \vert$ of uniform convergence on $G$. 
\index{Almost periodic function}
\par

Maximally almost periodic groups can be characterized in terms 
of almost periodic functions.

\begin{prop}
% 4.C.6
\label{Pro-AP-MAP}
For a topological group $G$, the following conditions are equivalent
%(see \cite[\S~16]{Dixm--C*}):
\begin{enumerate}[label=(\roman*)]
\item\label{iDEPro-AP-MAP}
$G$ is MAP;
\item\label{iiDEPro-AP-MAP}
there exist a compact group $K$
and an injective continuous homomorphism $G \to K$;
\item\label{iiiDEPro-AP-MAP}
the homomorphism $\beta_G \, \colon G \to \Bohr(G)$
of $G$ in its Bohr compactification is injective,
that is, the von Neumann kernel of $G$ is trivial;
\item\label{ivDEPro-AP-MAP}
almost periodic functions separate the points of $G$, 
i.e., for $g,g' \in G$ with $g \ne g'$, there exists $f \in AP(G)$ such that $f(g) \ne f(g')$.
\end{enumerate}
\end{prop}

\begin{proof}
Assume that $G$ is MAP and let $(\pi_i, \Hi_i)_{i \in I}$ be a family
of representatives for the equivalence classes in $\widehat G_{\rm fd}$.
The continuous homomorphism
$$
\prod_{i \in I} \pi_i \, \colon \, G \to K \, := \, \prod_{i \in I} \U(\Hi_i)
$$
is injective and this shows that \ref{iDEPro-AP-MAP} implies \ref{iiDEPro-AP-MAP}. 
\par

The fact that \ref{iiDEPro-AP-MAP} implies \ref{iiiDEPro-AP-MAP}
follows from the universal property
of the Bohr compactification $(\Bohr, \beta_G)$ of $G$. 
\par

The characterization \ref{iDETheo-AP} in Theorem~\ref{Theo-AP}
of functions in $AP(G)$
shows that \ref{iiiDEPro-AP-MAP} implies \ref{ivDEPro-AP-MAP}. 
\par

Assume that \ref{ivDEPro-AP-MAP} holds.
It follows from Theorem~\ref{Theo-AP}~\ref{iiiDETheo-AP}
that finite-dimen\-sional irreducible representations of $G$ separate the points of $G$.
Therefore $G$ is MAP. 
\end{proof}

\begin{exe}
% 4.C.7
\label{Ex-MAP-Groups}
(1)
A Moore group is clearly MAP.
In particular, locally compact abelian groups and compact compact groups are MAP.

\vskip.2cm

\index{SIN group}
(2)
For a connected locally compact group $G$,
the following five properties are equivalent:
\begin{enumerate}[label=(\roman*)]
\item\label{iDEEx-MAP-Groups}
$G$ is a direct product $K \times \R^n$,
with $K$ a compact group and $n \ge 0$;
\item\label{iiDEEx-MAP-Groups}
$G$ is MAP;
\item\label{iiiDEEx-MAP-Groups}
$G$ is a Moore group;
\item\label{ivDEEx-MAP-Groups}
the quotient of $G$ by its centre is compact;
\item\label{vDEEx-MAP-Groups}
every neighborhood of the identity in $G$ contains a compact
neighborhood which is invariant under all inner automorphisms,
i.e., $G$ is a SIN group
(for more on SIN groups, see Section \ref{S:FaithfulFinTypRep}.
\end{enumerate}
For a proof, we refer to \cite[16.4.6]{Dixm--C*}.
The equivalence of \ref{iDEEx-MAP-Groups} and \ref{iiDEEx-MAP-Groups}
is due to Freudenthal for second-countable connected LC groups \cite{Freu--36}
and to Weil for connected LC groups \cite[\S~32]{Weil--40}.
For more equivalences, see Diagram 3 in \cite{Palm--78}.
\par

The equivalence of \ref{iiDEEx-MAP-Groups} and \ref{iiiDEEx-MAP-Groups}
carries over to connected-by-compact locally compact groups \cite[Theorem 2.18]{GrMo--71a}.

\vskip.2cm

(3)
For every topological group $G$,
the quotient $G/N_G$ by its von Neumann kernel is MAP.

\vskip.2cm

(4)
Let $G$ be a MAP group; every subgroup of $G$ is MAP.
Let $G_1, G_2$ be topological groups; the direct product $G_1 \times G_2$ is MAP
if and only if $G_1$ and $G_2$ are MAP.
Let $G$ be a LC group and $N$ a closed normal subgroup
which is either compact or equal to the centre of $G$;
if $G$ is MAP, so is $G/N$ \cite[Proposition 1]{LeRo--68}.
Every locally compact MAP group is unimodular \cite[Theorem 2]{LeRo--68}.
\end{exe} 

\index{Residually finite group}
There is an interesting class of MAP groups among discrete groups. 
Let $\Gamma$ be a \textbf{residually finite} discrete group, 
i.e., a group with the following property:
for every $\gamma \in \Gamma \smallsetminus \{e\}$, 
there exists a finite group $F$ and a homomorphism $\varphi \, \colon \Gamma \to F$
such that $\varphi(\gamma) \ne e$.
\par

It is clear that residually finite groups are MAP.
The converse does not hold, 
as witnessed by the additive group $\Q$ of rational numbers,
and indeed by any non-trivial divisible abelian group, 
since these are MAP and are not residually finite.
\par

However, among finely generated groups,
the MAP groups are exactly the residually finite ones.

\begin{prop}
% 4.C.8
\label{Pro-MAP-FinResidual}
For a finitely generated group $\Gamma$,
the following properties are equivalent: 
\begin{enumerate}[label=(\roman*)]
\item\label{iDEPro-MAP-FinResidual}
$\Gamma$ is MAP;
\item\label{iiDEPro-MAP-FinResidual}
$\Gamma$ is residually finite.
\end{enumerate}
\end{prop}

\begin{proof}
Indeed, suppose $\Gamma$ is MAP.
% \footnote{Voir aussi Harpe-Robertson-Valette, prop.~4.}
Then $\Gamma$ is a subgroup of a compact group $\Sigma$.
For $\gamma \in \Gamma$, there exists a representation $\pi \, \colon \Sigma \to \U(n)$
such that $\pi(\gamma) \ne \mathrm{Id}_{\C^n}$.
Since finitely generated groups that are linear
[in particular finitely generated groups that are subgroups of $\U(n)$ for some $n$]
are residually finite by Mal'cev Theorem \cite{Mal'c--40}, 
there exists a finite group $F$ and a homomorphism $\rho \, \colon \pi(\Gamma) \to F$
such that $\rho(\pi(\gamma)) \ne e$; hence $\Gamma$ is residually finite.
As mentioned above, the converse implication is obvious.
\end{proof}

\begin{exe}
% 4.C.9
\label{Ex-DiscreteMAP-NonMAP}
(1)
Finitely generated groups which are residually finite, and therefore MAP,
include linear groups (see the proof of Proposition~\ref{Pro-MAP-FinResidual})
and polycyclic groups, in particular nilpotent groups
(a result of K.\ Hirsch for which we refer to \cite[5.14.17]{Robi--96}).

\vskip.2cm

(2)
Discrete groups which are \emph{not} MAP include the following groups,
where $\K$ stands for an infinite field:
\begin{enumerate}
\item[$\cdot$]
Heisenberg groups $H(\K)$,
of which the von Neumann kernel coincides with the centre, 
by Corollary~\ref{Cor-FinDimRepHeis-InfiniteFields};
\index{Heisenberg group! $2$@$H(\K)$ over a field $\K$}
\item[$\cdot$]
affine groups $\Aff (\K)$,
for which the von Neumann kernel coincides
with the translation group $\K$,
by Corollary~\ref{Cor-FinDimRepAffine}
\index{Affine group! $4$@$\Aff(\K)$ of an infinite field $\K$} 
\item[$\cdot$]
general linear groups $\GL_n(\K)$ for $n \ge 2$,
of which the von Neumann kernel is $\SL_n(\K)$,
%by the argument of Example \ref{SLnotMAP} = 4.20
and special linear groups $\SL_n(\K)$,
by Corollary~\ref{Cor-FinDimRepGLn} and Corollary~\ref{Cor-FinDimRepGLn}
below;
\item[$\cdot$]
and of course the m.a.p.\ groups defined below
(except the one-element group).
\end{enumerate}
\index{General linear group! $\GL_n(\K)$ with $\K$ a field}
\end{exe}

\subsection
{Bohr compactification and profinite completion}
% subsection 4.C.c
\label{SS:Profinite}

\index{Profinite group}
A compact group $G$ is \textbf{profinite} if $G$ is topologically isomorphic
to the projective limit of finite groups.
Let us recall this last notion.
\par

\index{Projective limit of finite groups}
Let $(F_i)_{i \in I}$ be a family of finite groups indexed by a directed set $I$,
with respect to a collection of homomorphisms $p_{i, j} \, \colon F_i \to F_j$ for $j \le i$
such that $p_{i, i} = {\mathrm{Id}}_{F_i}$ and $p_{i, j} \circ p_{j, k} = p_{i, k}$
for $k \le j \le i$.
The associated \textbf{projective limit}
is the closed subgroup $\varprojlim_i F_i$
of the compact group $\prod_{i \in I} F_i$ defined by 
$$
\varprojlim_i F_i \, = \, \left\{ (g_i)_{i \in I} \mid g_j = p_{i, j}(g_i)
\hskip.2cm \text{for all} \hskip.2cm
i, j \in I, j \le i \right\}.
$$
Among compact groups, the profinite groups admit the following nice characterization
(for more on this class of groups, see \cite{Wils--98}).

\begin{prop}
% 4.C.10
\label{Prop-CharProfiniteGr}
For a compact group $G$, the following properties are equivalent:
\begin{enumerate}[label=(\roman*)]
\item\label{iDEProp-CharProfiniteGr}
$G$ is profinite; 
\item\label{iiDEProp-CharProfiniteGr}
$G$ is topologically isomorphic to a closed subgroup of a product of finite groups;
\item\label{iiiDEProp-CharProfiniteGr}
$G$ is totally disconnected.
\end{enumerate}
\end{prop}

\begin{proof}
It is obvious that \ref{iDEProp-CharProfiniteGr} implies \ref{iiDEProp-CharProfiniteGr}.

\vskip.2cm

Assume that $G$ is a closed subgroup of a product
$K = \prod_{i \in I} F_i$ of finite groups.
Then $G$ is compact since $K$ is compact.
For every $i \in I$, the canonical projection $p_i \, \colon K \twoheadrightarrow F_i$ is continuous.
Therefore the image of the connected component $G_0$ of the identity of $G$ under $p_i$
is connected, and therefore reduced to the identity, since $F_i$ is finite.
It follows that $G_0$ is trivial, that is, $G$ is totally disconnected.
This shows that \ref{iiDEProp-CharProfiniteGr} implies \ref{iiiDEProp-CharProfiniteGr}.

\vskip.2cm

To show that \ref{iiiDEProp-CharProfiniteGr} implies \ref{iDEProp-CharProfiniteGr},
assume that $G$ is totally disconnected. 
Let $\mathcal N$ be the family of open normal subgroups of $G$.
\par

For every $N \in \mathcal N$, the quotient $G/N$ is a finite group,
since $G$ is compact and $G/N$ is discrete.
Let $K = \varprojlim_N G/N$ be the projective limit
of the family $(G/N)_{N \in \mathcal N}$,
with respect to the canonical homomorphisms 
$$
p_{MN} \, \colon \, G/M \to G/N
\hskip.5cm \text{if} \hskip.2cm
M \subset N.
$$
We claim that the natural homomorphism
$$
\Phi \, \colon \, G \to K = \varprojlim_N G/N
$$
is a topological isomorphism. Since every quotient homomorphism $G \to G/N$
is continuous, the map $\Phi$ is continuous. 
Therefore, since $G$ is compact, it suffices to show that $\Phi$ is bijective. 
\par

Let us show that $\Phi$ is injective.
It is well-known that the connected component of a LC group $H$ coincides
with the intersection of all open subgroups of $H$
(\cite[Chap.~III, \S~4, no~6]{BTG1--4}).
Since $G$ is totally disconnected, it follows that 
$$
\ker(\Phi) \, = \, \bigcap_{N \in \mathcal N} N \, = \, \{e\}.
$$
Therefore $\Phi$ is injective. 
\par

To show that $\Phi$ is surjective, let $(g_N N)_{N \in \mathcal N} \in K$.
We claim that the intersection $\bigcap_{N \in \mathcal N} g_N N$
of cosets in $G$ is non empty. 
\par

Indeed, let $\mathcal{F}$ be a finite subset of $\mathcal N$.
Then, $M := \bigcap_{N \in \mathcal{F}} N$ is a normal open subgroup of $G$.
Since $(g_N N)_{N \in \mathcal N} \in \varprojlim_N G/N$,
we have $g_N = p_{M, N}(g_M)$ for every $N \in \mathcal{F}$.
Therefore $g_M \in \bigcap_{N \in \mathcal{F}} g_N N$.
We have shown that $\bigcap_{N \in \mathcal{F}} g_N N$ is non empty
for every finite subset $\mathcal{F}$ of $\mathcal N$. 
Observe that every coset $g_N N$ is closed, since open subgroups are closed.
As $G$ is compact, it follows form the finite intersection property
that $\bigcap_{N \in \mathcal N} g_N N$ is non empty.
For $g \in \bigcap_{N \in \mathcal N} g_N N$,
we have $\Phi(g) = (g_N N)_{N \in \mathcal N}$.
Therefore $\Phi$ is surjective.
\end{proof}

For a discrete group $\Gamma$,
we discuss now the relationship of ${\Bohr(G)}$
with another completion of $\Gamma$, the profinite completion of $\Gamma$.

\begin{theorem}
% 4.C.11
\label{Theo-ProFinCompletion}
Let $\Gamma$ be a discrete group.
\par

There exists a profinite group $G$ and a homomorphism $\beta \, \colon \Gamma \to G$
with dense image, with the following universal property:
for every profinite group $K$ 
and every homomorphism $\alpha \, \colon \Gamma \to K$,
there exists a continuous homomorphism 
$\alpha' \, \colon G \to K$
such that $\alpha = \alpha' \circ \beta$.
\par

The pair $(G, \beta)$ is unique in the following sense:
let $(G', \beta')$ be a pair consisting of a profinite group
and a homomorphism $\beta' \, \colon \Gamma \to G'$ with dense image
satisfying the same universal property.
Then there exists an isomorphism $\alpha \, \colon G \to G'$ of topological groups
such that $\beta' = \alpha \circ \beta$.
\end{theorem}

\begin{proof}
Here is one way to construct $G$. 
Let $G$ be the closure of the image of $\Gamma$
in the compact group $\prod_{N} \Gamma/N$,
where $N$ runs over the collection 
of all normal subgroups of finite index in $\Gamma$. 
Let $\beta \, \colon \Gamma\to G$ be the natural map. 
\par

One checks as in the proof of Theorem~\ref{Theo-Bohr} that 
$(G, \beta)$ has the desired universal property
and that it is the unique pair with this property. 
\end{proof}

The profinite group associated to $\Gamma$ as in Theorem~\ref{Theo-ProFinCompletion}
is denoted by $\Prof(\Gamma)$
and the corresponding homomorphism by
$\beta_P \, \colon \Gamma \to \Prof(\Gamma)$.
The pair $(\Prof(\Gamma), \beta_P)$,
is called the \textbf{profinite completion} of $\Gamma$. 
\index{Profinite completion}

\begin{rem} 
% 4.C.12
\label{Rem-ProFinCompletion}
(1)
As shown by the universal property, 
$\Prof(\Gamma)$ may also be realized 
as the projective limit $\varprojlim \Gamma/N$ of the family $(\Gamma/N)_N$,
where the $N$~'s are subgroups of finite index in $\Gamma$, 
with respect to the canonical maps 
$\Gamma/N \twoheadrightarrow \Gamma/M$ when $N \subset M$;
and the homomorphism $\beta$ is defined in terms
of the projections $\Gamma \twoheadrightarrow \Gamma/N$.
\par

Alternatively, one may define $\Prof(\Gamma)$ as 
the largest Hausdorff quotient of the completion of $\Gamma$
in the group topology for which 
the collection of all finite index normal subgroups in $\Gamma$
is a fundamental system of open neighborhoods of the identity.

\vskip.2cm

(2)
Note that $\Prof(\cdot)$ is a covariant functor:
every group homomorphism $\Gamma_1 \to \Gamma_2$
gives rise to a continuous homomorphism
$\Prof(\Gamma_1) \to \Prof(\Gamma_2)$.
For properties of this functor, see for example \cite[Section 3.2]{RiZa--10}.

\vskip.2cm

(3)
Observe that the group $\Gamma$ is residually finite if and only if
the homomorphism $\beta_P \, \colon \Gamma \to \Prof(\Gamma)$
is injective.
\end{rem}

Next, we determine the profinite completion of $\Z$. 
(In the literature, this completion is often denoted by $\widehat \Z$,
a notation which conflicts with our use of $\widehat \Gamma$ 
as the dual of a group $\Gamma$.)

\begin{exe}
% 4.C.13
\label{Ex-ProfZ}
We may realize the profinite completion $\Prof(\Z)$ of $\Z$ 
as the projective limit 
$$
\Prof(\Z) \, = \, \varprojlim_n \Z/n\Z
$$
with respect to the canonical maps 
$\Z/n\Z \twoheadrightarrow \Z/m\Z$ when $m \vert n$,
together with the obvious injection 
$$
\Z \, \hookrightarrow \, \varprojlim_{n} \Z/n\Z.
$$
\par

We can also realize $\Prof(\Z)$
as the product $\prod_{p \in P} \Z_p$ over the set $P$ of all primes
of the groups $\Z_p$ of $p$-adic integers,
together with the natural injection $\beta \, \colon \Z \hookrightarrow \prod_{p \in P} \Z_p$. 
Indeed, to check this, it suffices to prove that the pair $(\prod_{p \in P}\Z_p, \beta)$
satisfies the universal property of Theorem~\ref{Theo-ProFinCompletion}. 
\par

First, observe that $\prod_{p \in P} \Z_p$ is a totally disconnected compact group
and hence a profinite group. 
Next, $\beta(\Z)$ is dense in $\prod_{p \in P} \Z_p$.
Indeed, let $a = (a_p)_p \in \prod_{p \in P} \Z_p$.
A basis of neighbourhoods of $a$ is given by sets of the form 
$$
U_{p_1, \dots, p_k}^{e_1, \dots, e_k}(a)
\, = \,
\left\{ (x_p)_{p \in P} \in \prod_{p \in P} \Z_p
\hskip.1cm \bigg\vert \hskip.1cm
x_{p_i}-a_{p_i} \in p_i^{e_i} \Z_{p_i}
\hskip.2cm \text{for} \hskip.2cm
i = 1, \dots, k \right\},
$$
for pairwise distinct primes $p_1, \dots , p_k$ and integers $e_1, \dots, e_k\ge 1$. 
Given such a neighbourhood, for every $i = 1, \dots, k$, we can write
$$
a_{p_i} \, = \, m_{i} + b_{p_i}
\hskip.5cm \text{with} \hskip.2cm
m_i \in \Z
\hskip.2cm \text{and} \hskip.2cm
b_{p_i} \in p_i^{e_i} \Z_{p_i}.
$$
By the Chinese Remainder Theorem, there exists an integer $m \in \Z$ such that 
$$
m - m_i \, \in \, p_i^{e_i} \Z
\hskip.5cm \text{for every} \hskip.2cm
i = 1, \dots, k.
$$
It is clear that
$\beta(m) \in U_{p_1, \dots, p_k}^{e_1, \dots, e_k}(a)$.
\par

Let $K$ a profinite group and $\alpha \, \colon \Z \to K$ a homomorphism.
We claim that there exists a continuous homomorphism 
$\alpha' \, \colon \prod_{p \in P} \Z_p \to K$
such that $\alpha = \alpha' \circ \beta$.
Since $K$ is a subgroup of a product of finite groups, we may assume that 
$K = \prod_{i \in I} F_i$ for a family $(F_i)_{i \in I}$ of finite groups
and $\alpha = (\alpha_i)_{i \in I}$ for homomorphisms $\alpha_i \, \colon \Z \to F_i$.
Therefore it suffices to show that, for every $i \in I$, there exists a continuous homomorphism 
$\alpha'_i \, \colon \prod_{p \in P} \Z_p \to F_i$ such that $\alpha_i = \alpha'_i \circ \beta$.
\par

Fix $i \in I$ and let $n \ge 1$ be such that $\ker \alpha_i = n \Z$. 
Write $n = p_1^{e_1} \cdots p_k^{e_k}$ for pairwise distinct primes $p_1, \dots , p_k$
and for integers $e_1, \dots, e_k \ge 1$. Then 
$$
\prod_{j = 1}^k \Z_{p_j} / p_j^{e_j}\Z_{p_j} 
\, \approx \,
\prod_{j = 1}^k \Z/p_j^{e_j}\Z 
\, \approx \,
\Z / n\Z .
$$
More precisely, let
$$
f_n \, \colon \, \prod_{p \in P} \Z_p \twoheadrightarrow \prod_{j=1}^k \Z_{p_j} / p_j^{e_j} \Z_{p_j}.
$$
be the canonical projection.
Then 
$$
f_n \circ \beta \, \colon \, \Z \to \prod_{j=1}^k \Z_{p_j} / p_i^{e_j}\Z_{p_j}
$$
is a surjective homomorphism with kernel $n\Z$.
Therefore the map $\alpha_i' \, \colon \prod_{p \in P} \Z_p \to F_i$, given by 
$$
\alpha'_i (x) \, = \, \alpha_i (m) ,
$$
for $x \in \prod_{p \in P} \Z_p$
and for any $m \in \Z$ such that $f_n(\beta(m)) = f_n(x)$,
is a well defined homomorphism.
Moreover, we have
$$
\alpha_i'(\beta(m)) \, = \, \alpha_i(m)
% \, = \, \alpha_i
\hskip.5cm \text{for all} \hskip.2cm
m \in \Z,
$$
that is, $\alpha_i= \alpha_i' \circ \beta$.
\par

Observe that, by uniqueness of the Bohr compactification,
it follows in particular that there is an isomorphism of topological groups
$$
\Phi \, \colon \, \varprojlim_n \Z/n\Z \to \prod_{p \in P} \Z_p,
$$
which is uniquely determined by the conditions
$$
\Phi(\beta_1(m)) \, = \, \beta_2(m)
\hskip.5cm \text{for} \hskip.2cm
m \in \Z ,
$$
where, $\beta_1 \, \colon \Z \hookrightarrow \varprojlim_n \Z/n\Z$
and $\beta_2 \, \colon \Z \hookrightarrow \prod_{p \in P} \Z_p$
denote the canonical injections.~$\square$
\end{exe}

Let $\Gamma$ be a discrete group.
Observe that, by the universal property, there is a continuous homomorphism 
$c^B_P \, \colon \Bohr(\Gamma) \to \Prof(\Gamma)$
such that $\beta_P = c^B_P \circ \beta_B$; it follows that 
$c^B_P$ has dense image and is therefore surjective, 
since $c^B_P $ is continuous and since $\Bohr(\Gamma)$ is compact.
\par

In general, the surjective homomorphism $c^B_P$ is not injective;
for example, when $\Gamma = \Z$ is infinite cyclic, we have
(see Proposition~\ref{Prop-BohrAbelian})
$$
\Bohr(\Z) \, \approx \, 
\big( ({\widehat \Z})_{\rm disc} \big) \widehat{\phantom{G}} \, \approx \,
\big( (\R / \Z)_{\rm disc} \big) \widehat{\phantom{G}}
$$
and so $\Bohr(\Z)$ is not second-countable,
whereas (see Example~\ref{Ex-ProfZ})
$$
\Prof(\Z) \, \approx \, \prod_{p \in P} \Z_p
$$
is a second-countable compact abelian group.
\par

One may ask which groups $\Gamma$ have the 
property that the homomorphism
$c^B_P$ is an isomorphism.
This property admits equivalent reformulations.

\begin{prop}
% 4.C.14
\label{Pro-IsoBohrProfini}
Let $\Gamma$ be a discrete group. 
The following properties are equivalent:
\begin{enumerate}[label=(\roman*)]
\item\label{iDEPro-IsoBohrProfini}
the continuous surjective homomorphism
$c^B_P \, \colon \Bohr(\Gamma)
\twoheadrightarrow \Prof(\Gamma)$
is an isomorphism;
\item\label{iiDEPro-IsoBohrProfini}
the image of every finite-dimensional representation
$\Gamma \to \U(n)$ is finite;
\item\label{iiiDEPro-IsoBohrProfini}
the compact group $\Bohr(\Gamma)$ is totally disconnected.
\end{enumerate}
\end{prop}

\begin{proof}
Assume that \ref{iDEPro-IsoBohrProfini} holds.
Let $\pi \, \colon \Gamma \to \U(n)$ be a finite-dimensional representation of $\Gamma$. 
By the universal property of $\Bohr(\Gamma)$,
there is a continuous homomorphism 
$\pi' \, \colon \Bohr(\Gamma) \to \U(n)$
such that $\pi = \pi' \circ \beta_B$.
The map 
$$
f \, = \, \pi' \circ ({c^B_P})^{-1} \, \colon \, \Prof(\Gamma) \to \U(n)
$$
is a continuous homomorphism. 
Therefore there exists a neighbourhood $U$ of $e$
in $\Prof(\Gamma)$ such that $\Vert f(g) - I \Vert < 1$
for every $g \in U$, 
where $\Vert \cdot \Vert$ is the operator norm 
on the ambient algebra $\Li (\C^n)$ of $\U(n)$,
and $I$ stands for $\mathrm{Id}_{\C^n}$.
There exists a normal subgroup $N$ of finite index in
$\Prof(\Gamma)$ such that $N\subset U$.
We have then
$$
\Vert f(g^n) - I \Vert \, = \, \Vert f(g)^n - I \Vert < 1
$$
for every $g \in N$ and every $n \in \Z$. 
It follows that $1$ is the only eigenvalue of the unitary matrix $f(g)$ 
and hence that $f(g) = I$ for every $g \in N$. 
So, $f$ is trivial on $N$ and hence has finite image. 
It follows that $\pi = f \circ c^B_P \circ \beta_B$ has finite image. 
Therefore \ref{iiDEPro-IsoBohrProfini} holds.

\vskip.2cm

Assume that \ref{iiDEPro-IsoBohrProfini} holds.
Since $\Bohr(\Gamma)$
can be realized as the image of the homomorphism
$$
\prod_{\pi \in \widehat G_{\rm fd}} \pi \, \colon \,
G \to \prod_{\pi \in \widehat G_{\rm fd}} \U(\Hi_\pi) ,
$$
$\Bohr(\Gamma)$ is 
a closed subgroup of a Cartesian product of finite groups,
and therefore is totally disconnected. This proves \ref{iiiDEPro-IsoBohrProfini}.

\vskip.2cm

Assume that \ref{iiiDEPro-IsoBohrProfini} holds. 
Then $\Bohr(\Gamma)$ is profinite. 
Therefore by the universal property of $\Prof(\Gamma)$,
there exists a continuous map 
$c^P_B \, \colon \Prof(\Gamma) \to \Bohr(\Gamma)$
such that $c^P_B \circ \beta_P = \beta_B$.
By density of $\beta_P(\Gamma)$ in $\Prof(\Gamma)$ 
and of $\beta_B(\Gamma)$ in $\Bohr(\Gamma)$,
we see that $c^P_B$ is the inverse map to $c^B_P$.
\end{proof}

We characterize abelian groups for which $c^B_P$ is an isomorphism.

\begin{cor}
% 4.C.15
\label{Cor-AbelianBohrProf}
Let $\Gamma$ be a discrete \emph{abelian} group.
The following properties are equivalent:
\begin{enumerate}[label=(\roman*)]
\item\label{iDECor-AbelianBohrProf}
the map $c^B_P \, \colon \Bohr(\Gamma)
\twoheadrightarrow \Prof(\Gamma)$
is an isomorphism;
\item\label{iiDECor-AbelianBohrProf}
the group $\Gamma$ is a direct sum $\bigoplus_{i \in I} C_i$
of finite cyclic groups $C_i$ of bounded orders. 
\end{enumerate}
\end{cor}

\begin{proof}
Indeed, by Pr\"ufer theorem (see \cite[4.35]{Robi--96}), 
the abelian group $\Gamma$ is periodic of bounded exponent
(that is, there exists an integer $N \ge 1$ such that $\gamma^N = e$
for every $\gamma \in \Gamma$)
if and only if $\Gamma$ is direct sum of finite cyclic groups of bounded orders.
So, it suffices to show that $c^B_P$ is an isomorphism
if and only if $\Gamma$ is periodic of bounded exponent.

\vskip.2cm

Assume that $\Gamma$ is periodic of bounded exponent $N$.
Then, $\chi(\gamma)^N = \chi(\gamma^N) = 1$
for all $\gamma \in \Gamma, \chi \in \widehat \Gamma$.
So, every $\chi \in \widehat \Gamma$ has finite image
and hence $c^B_P$ is an isomorphism, by Proposition \ref{Pro-IsoBohrProfini}.

\vskip.2cm

Conversely, assume that $c^B_P$ is an isomorphism,
that is, every $\chi \in \widehat \Gamma$ has finite image.
So, $G_{\rm disc}$ is a periodic group,
where $G$ is the compact group $\widehat \Gamma$.
Then $G = \bigcup_{n \ge 1}G_n$,
where $G_n$ is the closed subset $\{\chi \in G \mid \chi^n = 1 \}$ of $G$. 
By Baire theorem, there exists a non-empty open subset $U$ of $G$
which is contained in $G_n$ for some $n \ge 1$. 
Then $UU^{-1}$ is an open neighbourhood of $e$ contained in $G_n$.
The subgroup $H$ generated by $UU^{-1}$ is contained in $G_n$.
Moreover, $H$ is open in $G$ and has therefore finite index,
since $G$ is compact. It follows that $G = G_N$ for $N = n \vert G/H \vert $,
that is, $G_{\rm disc}$ is periodic of exponent $N$.
As $\chi(\gamma^N) = \chi^N(\gamma) = 1$
for all $\gamma \in \Gamma$ and $\chi \in \widehat \Gamma = G$, 
this implies that $\Gamma$ is periodic of exponent $N$.
\end{proof}

Observe that we recover from Corollary~\ref{Cor-AbelianBohrProf}
the following known fact (see \cite[4.35]{Rudi--62}):
if $\Gamma$ is a discrete abelian group of finite exponent,
then $\Bohr(\Gamma)$ is totally disconnected.

\vskip.2cm

Next, for $\Gamma= \SL_n(\Z)$ and $n \ge 3$,
we show that the map $c^B_P$ is an isomorphism.
Combining this result with the positive solution of the congruence subgroup theorem,
we will give a precise description of $\Bohr(\Gamma)$. 

\begin{cor}
% 4.C.16
\label{Cor-SLnZ-BohrProf}
Let $n$ be an integer with $n \ge 3$.
\begin{enumerate}[label=(\arabic*)]
\item\label{1DECor-SLnZ-BohrProf}
The map $c^B_P \, \colon \Bohr(\SL_n(\Z))
\twoheadrightarrow \Prof(\SL_n(\Z))$
is an isomorphism.
\item\label{2DECor-SLnZ-BohrProf}
The Bohr compactification $\Bohr(\SL_n(\Z))$ coincides
with $\SL_n(\Prof(\Z))$, 
together with the homomorphism $\SL_n(\Z) \to \SL_n(\Prof(\Z))$
induced by the natural embedding 
$$
\Z \to \Prof(\Z) \approx \prod_{p \in P} \Z_p.
$$
\end{enumerate}
\end{cor}

\begin{proof}
\ref{1DECor-SLnZ-BohrProf}
Consideration of the homomorphisms of reduction modulo $N$, for $N \ge 1$,
shows that $\Gamma = \SL_n(\Z)$ is residually finite for every $n \ge 2$.
We claim that every finite-dimensional representation of $\Gamma$ has finite image,
when $n \ge 3$.
\par

To prove the claim, we assume first that $n = 3$.
Let $\pi \, \colon \Gamma \to \U(d)$ be 
a finite-dimensional representation of $\Gamma = \SL_3(\Z)$.
Consider the restriction of $\pi$ to the subgroup
\index{Heisenberg group! $3$@$H(\Z)$}
$$
H \, := \, \begin{pmatrix}
1 & \Z & \Z \\
0 & 1 & \Z \\
0 & 0 & 1
\end{pmatrix} \subset \Gamma,
$$
which is a copy of the Heisenberg group $H(\Z)$ over the integers.
By Corollary \ref{Cor-FinDimRepHeis-Integers},
there exists an integer $N \ge 1$ such that $\pi$ is trivial on the elementary matrix 
$$
E_{1,3}(N) \, = \, \begin{pmatrix}
1 & 0 & N
\\
0 & 1 & 0
\\
0 & 0 & 1
\end{pmatrix}. 
$$
Now, the smallest normal subgroup $\Delta_N$ of $\Gamma$
containing $E_{1,3} (N)$ has finite index in $\Gamma$;
see Corollary 3 in \cite[\S~17.2]{Hump--80}
for a proof of this elementary but non-trivial fact,
which indeed holds for any $n \ge 3$.
Since $\pi$ factorizes through the quotient $\Gamma/\Delta_N$, 
the claim is proved for $n = 3$.
\par

We assume now that $n \ge 3$.
Let $\sigma$ denote the restriction of $\pi$
to the upper-left corner subgroup $\SL_3(\Z)$ of $\Gamma = \SL_n(\Z)$.
By the previous argument, there exists $N \ge 1$ such that
$\sigma(\gamma^N) = I$ for all $\gamma \in \SL_3(\Z)$.
Therefore $\ker (\pi)$ contains $E_{1,3}(N)$,
and hence the smallest normal subgroup
of $\Gamma$ containing $E_{1,3}(N)$,
which is a group of finite index in $\Gamma$.
This concludes the proof of the claim for all $n \ge 3$.
\par

Therefore the Bohr compactification $\Bohr(\SL_n(\Z))$ 
is isomorphic to the profinite completion of $\SL_n(\Z)$.

\vskip.2cm

\ref{2DECor-SLnZ-BohrProf}
By the congruence subgroup theorem \cite{BaLS--64, Menn--65},
the homomorphism
$$
\Prof(\SL_n(\Z)) \to \SL_n(\Prof(\Z))
$$
defined by the universal property is an isomorphism for $n \ge 3$.
The claim follows from this fact and from Item \ref{1DECor-SLnZ-BohrProf}.
\end{proof}

\begin{rem}
% 4.C.17
\label{Exa-TotallyDisconnectedBohr}
(1)
In contrast to the result in Corollary~\ref{Cor-SLnZ-BohrProf}, the group $\SL_2(\Z)$ has finite-dimensional representations with infinite images.
Indeed, there is an observation of Hausdorff at the heart
of the Hausdorff--Banach--Tarski paradox:
the rotation group $\SO(3)$ contains subgroups isomorphic
to the free product $(\Z / 2 \Z) \ast (\Z / 3 \Z)$,
i.e., to the quotient $\PSL_2(\Z)$ of $\SL_2(\Z)$ by its centre of order $2$.
This is in the appendix to \S~X.1 in \cite{Haus--14};
see also, for example, \cite{OsAd--76}.
\par

In particular, $\SL_2(\Z)$ has infinite three-dimensional representations,
and the epimorphism 
$\Bohr(\SL_2(\Z)) \twoheadrightarrow \Prof(\SL_2(\Z))$
has a non-trivial kernel.

\vskip.2cm

(2)
Let $\Gamma$ be a group that has a quotient
which is infinite, abelian, and not a direct sum of cyclic groups of bounded orders
(see Corollary~\ref{Cor-AbelianBohrProf} ); for example, let $\Gamma$ be a group that has a quotient
isomorphic to one of $\Z, \K, \K^\times$,
with $\K$ an infinite field of characteristic $0$.
Then $\Gamma$ has an infinite one-dimensional representations,
so that, again, the epimorphism 
$\Bohr(\Gamma) \twoheadrightarrow \Prof(\Gamma)$
has a non-trivial kernel.
\par

This applies to many of our favourite examples:
to $H(\K)$, $H(\Z)$, $\BS(1, p)$, $\Aff(\K)$, and $\GL_n(\K)$,
with again $\K$ an infinite field of characteristic $0$.
\end{rem}

\subsection
{Minimally almost periodic groups}
% subsection 4.C.d
\label{SS:map}

\index{Minimally almost periodic group, or m.a.p.\ group}
A topological group $G$ is \textbf{minimally almost periodic}, 
or \textbf{m.a.p.},
if every finite-dimensional representation of $G$ is the identity, that is, 
if $\widehat G_{\rm fd} = \{ 1_G \}$. 
\par

The straightforward proof of the following characterizations of 
m.a.p.\ groups is left to the reader
(compare with the chacterizations of MAP groups in Proposition~\ref{Pro-AP-MAP}).

\begin{prop}
% 4.C.18
\label{Pro-map}
For a topological group $G$, the following properties are equivalent: 
\begin{enumerate}[label=(\roman*)]
\item\label{iDEPro-map}
$G$ is a m.a.p.\ group; 
\item\label{iiDEPro-map}
every continuous homomorphism from $G$ to a compact group
is the identity;
\item\label{iiiDEPro-map}
the Bohr compactification of $G$ is the group with one element;
\item\label{ivDEPro-map}
every almost periodic continuous function on $G$ is constant.
\end{enumerate}
\end{prop}

The following characterization of finitely generated group
which are m.a.p.\ should be compared with Proposition~\ref{Pro-MAP-FinResidual}.

\begin{prop}
% 4.C.19
\label{Pro-map-FinGen}
For a finitely generated group $\Gamma$, the following 
properties are equivalent: 
\begin{enumerate}[label=(\roman*)]
\item\label{iDEPro-map-FinGen}
$\Gamma$ is m.a.p.;
\item\label{iiDEPro-map-FinGen}
$\Gamma$ has no other finite quotient than the one-element group.
\end{enumerate}
\end{prop}

\begin{proof}
Indeed, Proposition~\ref{Pro-map}~\ref{iiDEPro-map}
shows that \ref{iDEPro-map-FinGen} implies \ref{iiDEPro-map-FinGen}. 
To show the converse, assume $\Gamma$ is not m.a.p.
Therefore there exists a non-trivial representation $\pi \, \colon \Gamma \to \U(n)$.
Since finitely generated groups which are linear
are residually finite (Mal'cev Theorem), 
there exists a finite group $F$ with more than one element
and a surjective homomorphism $\rho \, \colon \pi(\Gamma) \to F$.
Then $F$ is a non-trivial finite quotient of $\Gamma$.
Therefore, \ref{iiDEPro-map-FinGen} implies \ref{iDEPro-map-FinGen}.
\end{proof}

The following result in the case of $\K = \R$ appears
in the original article of von Neumann \cite[Page 483]{vNeu--34}.
See also Example $\mathcal G^{(\gamma)}$ in Section 5 of \cite{vNWi--40}.

\begin{prop}
% 4.C.20
\label{SLnotMAP}
Let $\K$ be an infinite field and $n$ an integer, $n \ge 2$.
\par

The special linear group $\SL_n(\K)$ is minimally almost periodic.
\end{prop}

\begin{proof}
We have to show that,
for every integer $d \ge 1$ and
every finite-dimensional representation $\pi \, \colon \SL_n(\K) \to \U(d)$,
we have $\ker (\pi) = \SL_n(\K)$.
Since proper normal subgroups of $\SL_n(\K)$ are central 
(see for example Theorem 4.9 in \cite{Arti--57}),
it suffices to show that $\ker (\pi)$ contains a matrix $\gamma$
which is not central in $\SL_n(\K)$.
Moreover, if $\alpha \, \colon \SL_2(\K) \to \SL_n(\K)$ is any injective homomorphism,
$\pi \circ \alpha$ is a finite dimensional representation of $\SL_2(\K)$;
so it suffices to prove the claim in the case $n = 2$. 
\par

For $x \in \K$ and $a \in \K^\times$, consider the elements
$$
\gamma_x \, = \, \begin{pmatrix} 1 & x \\ 0 & 1 \end{pmatrix} 
\hskip.5cm \text{and} \hskip.5cm
\delta_a \, = \, \begin{pmatrix} a & 0 \\ 0 & a^{-1} \end{pmatrix} 
$$
of $\SL_2(\K)$.
Observe that the $\gamma_x$~'s commute with each other and that 
$$
\delta_a \gamma_x \delta_a ^{-1} \, = \, 
\gamma_{a^2x} .
$$
It follows that $\{\pi(\gamma_{a^2}) \mid a \in \K^\times\}$ 
is a family of conjugate and commuting matrices in $\U(d)$.
Therefore, the $\pi(\gamma_{a^2})$~'s have all the same set 
$\{\lambda_1, \dots, \lambda_d\}$ of eigenvalues and, moreover, we can find a
common basis in $\C^d$ for which they are represented by diagonal matrices. 
However, there are only finitely many diagonal matrices with 
diagonal entries from the set $\{\lambda_1, \dots, \lambda_d\}$ and the set 
$\{a^2 \mid a \in \K^\times\}$ is infinite, since $\K$ is infinite.
It follows that there exist $a,b \in \K^\times$ with $a^2 \ne b^2$ such that
$\pi( \gamma_{a^2} ) = \pi( \gamma_{b^2})$.
The element 
$$
\gamma_{a^2}\gamma_{b^2}^{-1} \, = \, \gamma_{a^2-b^{2}} \in \SL_n(\K)
$$
belongs therefore to the kernel of $\pi$ and is clearly not central.
\end{proof}

\begin{rem}
% 4.C.21
(1)
The result of Proposition \ref{SLnotMAP} can be extended to other groups, as follows:
let $\Gamma$ be a group generated by its subgroups isomorphic to 
either $\SL_2(\K)$ or $\PSL_2(\K)$, for some infinite field $\K$; 
then every finite-dimensional representation of $\Gamma$
is a multiple of the identity $1_\Gamma$.
\par

For $\K = \R$, this extended result and structure theory
imply that every connected real Lie group
which is semisimple, without compact factors and with finite centre,
is minimally almost periodic, as a discrete group, 
and therefore \emph{a fortiori} as a Lie group (see \cite[Remark 3]{Vale--86}).
In fact, an even stronger result is true \cite{SevN--50}:
such a group has no non-trivial representation
in the unitary group of a factor of finite type (see Section~\ref{SectionTypesI+IILC}).

\vskip.2cm

(2)
By contrast, a second-countable locally compact group 
that is solvable and non-compact
is not minimally almost periodic \cite[Theorems 3.4 \& 4.1]{Schm--84}.
\end{rem}

We can now easily determine all finite dimensional representations 
of $\GL_n(\K)$ for an infinite field.

\begin{cor}
% 4.C.22
\label{Cor-FinDimRepGLn}
Let $\K$ be an infinite field and $n$ an integer, $n \ge 2$.
\par

Every irreducible finite dimensional representation of $\Gamma = \GL_n(\K)$
is a unitary character of $\Gamma$; more precisely, the map
$$
\widehat \K \to \widehat \Gamma_{\rm fd}, \hskip.2cm \chi \mapsto \varphi_\chi,
$$
where $\varphi_\chi(\gamma) = \chi(\det (\gamma))$ for $\gamma \in \Gamma$, 
is a bijection
\end{cor}

\begin{proof}
Let $\pi$ be an irreducible finite dimensional representation of $\Gamma$.
By proposition~\ref{SLnotMAP}, $\SL_n(\K)\subset \ker \pi$. 
Since $[\Gamma, \Gamma] = \SL_n(\K) = \ker \det$, the claim follows.
\end{proof}

\begin{exe}
% 4.C.23
\label{OtherExplesmap}
There are several other known classes of m.a.p.\ groups:
The group ${\rm Alt} (\N)$
of even finitely supported permutations of the integers (a simple group),
more generally every countable group which is the union
of a strictly increasing sequence of simple subgroups,
as well as every countable infinite group 
which is finitely generated and simple \cite[Examples 2.12 \& 2.13]{BCRZ--16}.
Also, there exists a countable infinite m.a.p.\ group
that is orderable, locally solvable, and perfect;
for this and other examples, we refer to \cite{GrKO--15}.
% \cite[Corollary 3.2]{GrKO--15}
\end{exe}

%-----------------------------------------------------------------------
% End of chapter 4
%-----------------------------------------------------------------------
\chapter[Describing all irreducible representations]
{Describing all irreducible representations of some semi-direct products}
% Chapter 5
\label{Chap-AllIrrRed}

\emph{Let $\Gamma$ be a countable discrete group.
The procedure of induction of representations
from subgroups of $\Gamma$,
as described in Section~\ref{Section-IrrIndRep},
does produce irreducible representations of $\Gamma$, 
but usually does not produce all of them.
We will now describe a general construction of representations
which produces} all \emph{irreducible representations, up to equivalence,
of a second-countable LC group $G$ of the following form:
$G = H \ltimes N$ is a semi-direct product, where
\begin{enumerate}
\item[$\bullet$]
$H$ is a discrete subgroup of $G$;
\item[$\bullet$]
$N$ is a locally compact abelian normal closed subgroup of $G$.
\end{enumerate}
In fact, this procedure will give a description
of an arbitrary} factor \emph{representation of $G$, up to equivalence
(Theorem~\ref{Theo-AllRepSemiDirect2}).
The description of a factor representation or an irreducible representation of $G$
involves an equivalence class of probability measures $\mu$ on $\widehat N$
and a cocycle $\widehat N \times H \to \U(\Ki)$ over $(\widehat N, \mu)$
with values in the unitary group of some Hilbert space $\Ki$. 
}
\par

\emph{
One should mention that, even when $H$ is abelian,
the procedure we are going to describe does not allow in general
to establish a concrete list of all irreducible representations of $G$.
See the comments in Remark~\ref{Rem-AllRepSemiDirect}.
}
\par

\emph{
This procedure is an extension of the Mackey machine
to the context of ergodic non-transitive actions
and can be cast in the more general framework of groupoid representations
\cite{Rams--76}.
}

\section
{Constructing some irreducible representations}
% Section 5.A
\label{Section-MoreIrredRep}

Let $G = H \ltimes N$ be a semi-direct product
of a countable \emph{discrete} group $H$
with a second-countable locally compact abelian normal subgroup $N$;
note that $N$ need not be discrete.
In this section, we exhibit a family of representations of $G$
attached to $H$-quasi-invariant measures
on the second-countable LCA group $\widehat N$.

\begin{defn}
% 5.A.1
\label{defnRN}
Let $(X, \mathcal B)$ be a Borel space,
$\mu$ a $\sigma$-finite measure on $(X, \mathcal B)$,
and $H \curvearrowright X$, or
$$
H \times X \to X, \hskip.2cm (h, x) \mapsto hx ,
$$
a \textbf{left} action of a group $H$ on $(X, \mathcal B)$, preserving the class of $\mu$.
For $h \in H$, denote by $(h^{-1})_*(\mu)$ the image of $\mu$ by the action of $h^{-1}$,
and consider the Radon--Nikodym derivative 
$$
d (h,x) \, := \, \frac{d (h^{-1})_*(\mu)} {d\mu} (x)
\hskip.5cm \text{for all} \hskip.2cm
x \in X
$$
of $(h^{-1})_*(\mu)$ with respect to $\mu$.
Then (see Proposition \ref{cocycleRNtaumu})
we have the cocycle relation
$$
\label{eqq/RNcocL}
% \tag{RN}
d (h_1h_2, x) \, = \, d (h_1, h_2x) d (h_2, x)
\hskip.5cm \text{for all} \hskip.2cm
h_1, h_2 \in H, \hskip.1cm x \in X .
$$
The map $d \, \colon H \times X \to \R_+^\times$
is called the \textbf{Radon--Nikodym cocycle}
of the action $H \curvearrowright X$.
\index{Radon--Nikodym! cocycle}
\index{Cocycle! $1$@Radon--Nikodym}

\vskip.2cm

In case of a \textbf{right} group action $X \curvearrowleft H$,
the cocycle is defined in the same way, and precisely as follows.
For $h \in H$, we write $\mu_*(h^{-1})$ for the image of $\mu$ by the action of $h^{-1}$,
and we define the Radon--Nikodym cocycle $c$ by
$$
c(x, h) \, := \, \frac{ d\mu_*(h^{-1}) }{d \mu} (x)
\hskip.5cm \text{for all} \hskip.2cm
x \in X .
$$
Then we have the cocycle relation, which reads now
\begin{equation}
\label{eqq/RNcocR}
\tag{RN}
c(x, h_1 h_2) \, = \, c(x, h_1) c(xh_1, h_2)
\hskip.5cm \text{for all} \hskip.2cm
x \in X, \hskip.1cm h_1, h_2 \in H.
\end{equation}
\small
[Indeed, let $X \times H \to X$ be a right action
giving rise to the Radon--Nikodym map $c \, \colon X \times H \to \R_+^\times$;
there is a left action $H \times X \to X$ defined by $(h, x) \mapsto hx := xh^{-1}$,
giving rise to a cocycle $d$, and $c(x, h) = d(x, h^{-1})$;
the cocycle relation for $c$ follows from this last equality
and from the cocycle relation for $d$.]
\normalsize
\end{defn}

\begin{constr}
% 5.A.2
\label{construcrep1}
Let $G = H \ltimes N$ be a semi-direct product
of a countable discrete subgroup $H$
with a second-countable locally compact abelian normal subgroup $N$,
with respect to a left action $H \curvearrowright N$ of $H$ on $N$;
see Convention \ref{convsemidirect}.
Recall that we identify the groups $H$ and $N$ to subgroups of $G$.
For $h \in H$ and $n \in N$, we have
$
(h,e) (e,n) (h,e)^{-1} = (h, h \cdot n) (h^{-1}, e) = (e, h \cdot n) ,
$
so that the left action of $H$ on $N$ can then be written as
$H \times N \to N, \hskip.1cm (h, n) \mapsto hnh^{-1}$;
it extends naturally to the $G$-action
$G \times N \to N, \hskip.1cm (g,n) \mapsto ghg^{-1}$.
\par

This left action $H \curvearrowright N$ corresponds to a right action
$\widehat N \curvearrowleft H, \hskip.1cm (\chi, h) \mapsto \chi^h$,
defined by
$$
\chi^h(n) \, = \, \chi(hnh^{-1})
\hskip.5cm \text{for all} \hskip.2cm
\chi \in \widehat N, \hskip.1cm h \in H, \hskip.1cm n \in N .
$$
The latter extends naturally to a right action of $\widehat N \curvearrowleft G$,
defined by $\chi^g(n) = \chi(gng^{-1})$
for all $\chi \in \widehat N$, $g \in G$, and $n \in N$.
\par

Let $\mu$ be a probability measure on $\widehat N$
which is quasi-invariant by $G$ (equivalently by $H$).
For every $g = (h, n) \in G$,
we define an operator $\pi_\mu(g)$
on the Hilbert space $L^2(\widehat N, \mu)$ by
$$
(\pi_\mu(g)f) (\chi) \, = \, \chi (n) c(\chi, h)^{1/2} f(\chi^{h}) 
\hskip.5cm \text{for all} \hskip.2cm
f \in L^2(\widehat N, \mu), \hskip.1cm \chi \in \widehat N .
$$
By definition of the Radon--Nikodym derivative,
$\pi_\mu(g)$ is a unitary operator on $L^2(\widehat N, \mu)$;
the cocycle identity (\ref{eqq/RNcocR}) implies that 
$$
G \to \U (L^2(\widehat N, \mu)), \hskip.2cm g \mapsto \pi_\mu(g)
$$
is a group homomorphism.
Moreover, since $\widehat N$ is second-countable, 
the Hilbert space $L^2(\widehat N, \mu)$ is separable
(Theorem \ref{regmeas2ndc}~\ref{2DEregmeas2ndc})
and hence $\pi_\mu$ is continuous,
when $\U (L^2(\widehat N, \mu))$ is equipped with the strong operator topology
(see Proposition A.6.1 in \cite{BeHV--08}).
So, $(\pi_\mu, L^2(\widehat N, \mu))$ is a representation of $G$. 
\end{constr}

\begin{rem}
% 5.A.3
\label{Rem-construcrep1}
Let $G = H \ltimes N, \mu$ and $\pi_\mu$ be as in Construction~\ref{construcrep1}.

\vskip.2cm

(1)
The representation $\pi_\mu$ is an extension to $G$
of the canonical representation of the LCA group $N$ associated to $\mu$
considered in Section~\ref{Section-CanRepAbGr}
(and also denoted by $\pi_\mu$ there). 
\par

In case $H = \{e\}$ and $G = N$,
the representation $\pi_\mu$
is precisely that of Construction \ref{defpimupourGab}.

\vskip.2cm

(2)
The restriction of $\pi_\mu$ to $H$ is the Koopman representation of $H$
associated to the action $(\widehat N, \mu) \curvearrowleft H$.
See Section~\ref{SectionMSC},
where Koopman representations are introduced in a more general context.

\vskip,2cm

(3)
Assume that $\mu$ is \emph{ergodic} under the action of $H$. 
Then $\pi_\mu$ is an \emph{irreducible} representation of $G$.
This is a particular case of Corollary \ref{Cor-AllRepSemiDirect1} below.
\end{rem}

\section
{Constructing all irreducible representations} 
% Section 5.B
\label{Section-AllIrredRep}

The construction of the representations $\pi_\mu$
in Section~\ref{Section-MoreIrredRep}
can be considerably generalized 
by adding an operator-valued cocycle over $(\widehat N, \mu)$
as a further ingredient.

\begin{defn}
% 5.B.1
\label{defcocycle}
Let $\mu$ be a $\sigma$-finite measure
on a Borel space $(X, \mathcal B)$.
Let $H$ be a topological group equipped with a \textbf{right} group action
$(X, \mathcal B) \curvearrowleft H$ which is measurable
and which preserves the class of $\mu$.
\par

\index{Cocycle!}
Let $L$ be is a topological group.
A \textbf{cocycle $\alpha$ over $(X, \mu)$ with values in $L$} 
is a map $\alpha \, \colon X \times H \to L$ which satisfies the identity
\begin{equation}
\label{eqq/cocycle}
\tag{Coc}
\alpha(x, h_1h_2) \, = \, \alpha(x, h_1) \alpha(xh_1, h_2),
\end{equation}
for $\mu$-almost all $x \in X$ and all $h_1, h_2 \in H$,
and which is measurable when $H$ and $L$ are equipped with the 
Borel structures underlying their respective topologies.
Note that (\ref{eqq/cocycle}) implies that $\alpha(x, e) = e$
and $\alpha(x, h)^{-1} = \alpha(xh, h^{-1})$
for $\mu$-almost all $x \in X$ and all $h \in H$.
\par

Observe that, if $\mu_1$ and $\mu_2$ are two $\sigma$-finite measures
on $X$ which are quasi-invariant by $H$
and absolutely continuous with each other,
a map $\alpha \, \colon X \times H \to L$
is a cocycle over $(X, \mu_1)$
if and only if it is a cocycle over $(X, \mu_2)$.
Observe also that, if $\alpha_1 \, \colon X \times H \to L$ is a cocycle over $(X, \mu)$
and $\varphi \, \colon X \to L$ is a measurable map,
then the map $\alpha_2 \, \colon X \times H \to L$, defined by 
\begin{equation}
\label{eqq/cohomolog}
\tag{Coh}
\alpha_2(x, h) \, = \, \varphi(x)\alpha_1(x, h) \varphi(xh)^{-1},
\end{equation}
for all $x \in X$ and $h \in H$,
is also a cocycle over $(X, \mu)$.
\par

Two cocycles $\alpha_1, \alpha_2 \, \colon X \times H \to L$ over $(X, \mu)$ 
are \textbf{cohomologous}, or \textbf{equivalent},
if there exists a measurable map $\varphi \, \colon X \to L$ such that 
(\ref{eqq/cohomolog}) holds
for $\mu$-almost all $x \in X$ and all $h \in H$. 
\index{Cocycle! $2$@cohomologous}
\par

A cocycle $\alpha \, \colon X \times H \to L$ over $(X, \mu)$ 
which is cohomologous to the identity cocycle $I \, \colon (x, h) \mapsto e$
is called a \textbf{coboundary}. Thus, $\alpha$ is a coboundary 
if and only if there exists a measurable map $\varphi \, \colon X \to L$ such that 
$$
\alpha(x, h) \, = \, \varphi(x)\varphi(xh)^{-1}
$$
for $\mu$-almost every $x \in X$ and all $h \in H$.
\end{defn}

\begin{rem}
% 5.B.2
Let $(X, \mathcal B)$, $\mu$, $H$, $(X, \mathcal B) \curvearrowleft H$, and $L$ be as above.
We describe a construction which provides a justification for
the defining relations (\ref{eqq/cocycle}) and (\ref{eqq/cohomolog}).
\par

Let $(Y, \mathcal C)$ be a measure space.
Denote by $\mathcal M (X, Y)$ the space of measurable maps from $X$ to $Y$,
modulo the equivalence relation of equality $\mu$-almost everywhere.
There is a natural left action $H \curvearrowright \mathcal M (X, Y)$ defined by
$(h \ast f) (x) = f(xh)$ for all $h \in H$, $f \in \mathcal M (X, Y)$, and $x \in X$.
\par

Assume that, moreover, $L$ acts on $(Y, \mathcal C)$,
say on the left for notational convenience.
Let $\alpha \, \colon X \times H \to L$ be a cocycle.
There is a \textbf{twisted action} of $H$ on $\mathcal M (X, Y)$ defined by
$$
(h \ast_\alpha f) (x) \, = \, \alpha(x, h) f(xh)
\hskip.5cm \text{for all} \hskip.2cm
h \in H, \hskip.1cm f \in \mathcal M (X, Y), \hskip.1cm x \in X .
$$
The cocycle identity (\ref{eqq/cocycle}) is precisely what is needed
for this to be an action, i.e., 
for $h_1 \ast_\alpha (h_2 \ast_\alpha f) = (h_1h_2) \ast_\alpha f$
for all $h_1, h_2 \in H$ and $f \in \mathcal M (X, Y)$.

\par

If $\alpha$ and $\beta$ are two cocycles over $(X, \mu)$ with values in $L$
which are cohomologous,
then, because of (\ref{eqq/cohomolog}), 
the corresponding twisted actions $\ast_\alpha$ and $\ast_\beta$ are equivalent
in an appropriate sense (see \cite[Section 4.2]{Zimm--84}).
\end{rem}

\begin{exe}
% 5.B.3
\label{examplecocycle}

(1)
Let $H$ be a subgroup of a group $G$
and $T$ a right transversal, so that $G = \bigsqcup_{t \in T} Ht$.
With the notation of Construction \ref{constructionInd}(2), the map
$$
T \times G \to H, \hskip.2cm (t,g) \mapsto \alpha(t,g)
$$
is a cocycle over $T$ (with the counting measure) with values in $H$.

\vskip.2cm

(2)
Consider a Borel space $(X, \mathcal B)$,
a $\sigma$-finite measure $\mu$ on it,
and a group $H$ acting on $X$ on the right
by measurable transformations preserving the class of $\mu$.
Denote by $\mu_*(h)$ the direct image of $\mu$ by the action of $h$.
The Radon--Nikodym cocycle $c_\mu \, \colon X \times H \to \R^\times_+$
is given by
$$
c(x, h) \, = \, \frac{ d\mu_*(h^{-1}) }{d \mu} (x) .
$$
as in Definition \ref{defnRN}.
Let $\nu$ be a another $\sigma$-finite measure on $(X, \mathcal B)$
which is equivalent to $\mu$.
Then $\nu = \varphi \mu$ for a measurable function $\varphi \, \colon X \to \R^\times_+$
which is locally $\mu$-integrable.
We have
$$
\begin{aligned}
c_\nu (x, h)
\, &= \, \frac{ d\nu_*(h^{-1}) }{d \mu} (x)
\, = \, \frac{ d\nu_*(h^{-1}) }{ d\mu_*(h^{-1}) }(x)
\hskip.1cm \frac{ d\mu_*(h^{-1}) }{ d\mu }(x)
\hskip.1cm \frac{d\mu}{d\nu}(x)
\\
\,& = \, \varphi(xh) \hskip.1cm c_\mu(x, h) \hskip.1cm \frac{1}{\varphi(x)} ,
\end{aligned}
$$
that is,
$$
c_\mu(x, h) \, = \, \varphi(x)\hskip.1cm c_\nu(x, h)\hskip.1cm \varphi(xh)^{-1}
\hskip.5cm \text{for all} \hskip.2cm
x \in X
\hskip.2cm \text{and} \hskip.2cm
h \in H .
$$
This shows that the cocycles $c_\nu$ and $c_\mu$ are cohomologous.
\par

In particular, there exists a $\sigma$-finite measure $\nu$ on $(X, \mathcal B)$
which is equivalent to $\mu$ and $H$-invariant,
i.e., and such that $c_\nu(c, h) = 1$ for all $x \in X$ and $h \in H$,
if and only if $c_\mu$ is a coboundary,
that is, if and only if $c_\mu$ is of the form
$$
c_\mu(x, h) \, = \, \frac{\varphi(x)}{\varphi(xh)},
$$
for some measurable and locally $\mu$-integrable
function $\varphi \, \colon X \to \R^\times_+$.

\vskip.2cm

(3)
Let $G = H \ltimes N$ be a semi-direct product 
with $H$ countable discrete,
$N$ locally compact abelian,
$\mu$ a probability measure on $\widehat N$
which is quasi-invariant by $G$ (equivalently: by $H$),
and $\pi_\mu \, \colon G \to \U (L^2(\widehat N, \mu))$
the resulting representations,
as in Construction \ref{construcrep1}.
The map
$$
\alpha \, \colon \, \widehat N \times G \to \T , \hskip.2cm
(\chi, (h, n)) \mapsto \chi (n)
$$
is a cocycle, which satisfies the identity (\ref{eqq/cocycle}).
Note that, in this construction, the definition of the representation $\pi_\mu$
involves the product of two cocycles,
one is $\alpha$, the other the Radon-Nikodym cocycle.
They will both appear again in the definition of the representation
of Construction \ref{construcrep2}.

\vskip.2cm

(4) 
Let $H, L$ be topological groups,
$f \, \colon H \to L$ a measurable group homomorphism,
and $(X, \mathcal B, \mu)$ a $\sigma$-finite measure space
on which $H$ acts by measurable transformations preserving the class of $\mu$.
Then 
$$
\alpha \, \colon \, X \times H \to L , \hskip.2cm (x, h) \mapsto f(h)
$$
is a cocycle over $(X, \mu)$ with values in $L$.
\par

More generally, let $H, L, M$ be topological groups, $\alpha \, \colon X \times H \to L$
a cocycle over $(X, \mu)$ with values in $L$,
and $f \, \colon L \to M$ a measurable group homomorphism. Then 
$$
f \circ \alpha \, \colon \, X \times H \to M, \hskip.2cm (x, h) \mapsto f(\alpha(x, h))
$$
is a cocycle over $(X, \mu)$ with values in $M$.

\vskip.2cm

(5)
Let $(X, \mathcal B, \mu)$ be a $\sigma$-finite measure space
and $T \, \colon X \to X$ a Borel isomorphism preserving the class of $\mu$.
Consider $(X, \mathcal B, \mu)$ as a $\Z$-space,
for the generator $1$ of $\Z$ acting by $T$
(yet we continue to write $T^ix$ the action of a power $T^i$ of $T$ on $x$).
Let $H$ be a topological group and
$f \, \colon X \to H$ a measurable map.
Define a measurable map $\alpha_f \, \colon X \times \Z \to H$ by 
$$
\alpha_f(x, n) \, = \,
\begin{cases}
\prod_{i = 0}^{n - 1} f(T^i x)
&\text{if} \hskip.2cm n > 0 ,
\\
e &\text{if} \hskip.2cm n = 0 ,
\\
\prod_{i = 1}^{\vert n \vert} f(T^{-i}x)^{-1}
&\text{if} \hskip.2cm n < 0 .
\end{cases}
$$
One checks that $\alpha_f$ is a cocycle.
\par

Every cocycle $\alpha \, \colon X \times \Z \to H$ is of the form $\alpha_f$
for some measurable map $f \, \colon X \to H$.
Indeed, if $f \, \colon X \to H$ is defined by $f(x) = \alpha (x,1)$ for all $x \in X$,
then $\alpha = \alpha_f$, by the cocycle relation (\ref{eqq/cocycle}).
\par

A cocycle $\alpha = \alpha_f$ is a coboundary
if and only if there exists a measurable map $\varphi \, \colon X \to H$
such that 
$$
f(x) \, = \, \varphi(x)^{-1}\varphi (Tx)
$$
for $\mu$-almost every $x \in X$.

\vskip.2cm

(6) Let $X = \T$ be the circle,
equipped with the normalized Lebesgue measure~$\mu$.
Consider the action of $H = \Z$ given on $\T$
by an irrational rotation $R_\theta$.
Let
$$
f \, \colon \, \T \to \T, \hskip.2cm e^{2 \pi i t} \mapsto e^{2 \pi i t}
$$
denote the identity map.
We claim that the corresponding cocycle $\alpha_f \, \colon \T \times \Z \to \T$
with values in the circle group $\T$ is a not a coboundary.
\par
 
Indeed, assume, by contradiction, that there exists a measurable function
$\varphi \, \colon \T \to \T$ such that 
$$
f(e^{2 \pi i t}) \, = \, \varphi(e^{2 \pi i t})^{-1} \varphi (R_\theta(e^{2 \pi i t})),
$$
that is, 
$$
e^{2 \pi i t} \varphi(e^{2 \pi i t}) \, = \, \varphi (e^{2 \pi i (\theta + t)})
$$
for almost every $t \in \R$. 
For $n \in \Z$, let 
$$
c_n(\varphi) \, = \, \int_0^{1} \varphi(e^{2 \pi i t}) e^{-2 \pi i nt} dt
$$
be the $n$-th Fourier coefficient of $\varphi \in L^\infty(\T, \mu)$. 
The $n$-th Fourier coefficient
of the function $t\mapsto \varphi (e^{2 \pi i (\theta + t)})$
is $e^{-2 \pi i n\theta} c_n(\varphi)$
and the $n$-th Fourier coefficient
of the function $t\mapsto e^{2 \pi i t} \varphi(e^{2 \pi i t})$
is $c_{n - 1}(\varphi)$.
Therefore we have
$$
c_n(\varphi) \, = \, e^{2 \pi i n\theta} c_{n - 1}(\varphi)
\hskip.5cm \text{for all} \hskip.2cm 
n \in \Z.
$$
This implies that $(\vert c_n(\varphi)\vert)_{n \in \Z}$ is a constant sequence.
However, $(\vert c_n(\varphi)\vert)_{n \in \Z} \in \ell^2(\Z)$, 
since $\varphi \in L^2(\T, \mu)$. 
It follows that $c_n(\varphi) = 0$ for all $n \in \Z$,
that is, $\varphi = 0$.
This is impossible, since
$f(e^{2 \pi i t}) = \varphi(e^{2 \pi i t})^{-1}\varphi (R_\theta(e^{2 \pi i t}))$
and $f \ne 0$.

\vskip.2cm

More generally, one can show along the same lines that,
for $k \in \Z$, $k \ne 0$, 
the cocycles $\alpha_{f_k} \, \colon \T \times \Z \to \T$, 
defined by the functions 
$$
f_k \, \colon \, \T \to \T, \hskip.2cm e^{2 \pi i t} \mapsto e^{2 \pi i kt},
$$
are pairwise non-cohomologous.
For a more general result, see \cite[Lemma 2.2]{Furs--61}.
\end{exe}

\begin{constr}
% 5.B.4
\label{construcrep2}
Let us go back to a locally compact group $G = H \ltimes N$,
as in Section \ref{Section-MoreIrredRep}.
Let $\mu$ be a probability measure on $\widehat N$,
which is quasi-invariant by the action of $G$
(equivalently: by the action of $H$). 
We will be interested in cocycles over $(\widehat N, \mu)$
with values in the unitary group of a Hilbert space.
Recall that we have the $\R_+$-valued Radon--Nikodym cocycle $c$
defined by $c(\chi, h) \, := \, \frac{ d\mu_*(h^{-1}) }{d \mu} (\chi)$
for all $\chi \in \widehat N$ and $h \in H$.
\par

Let $\Ki$ be a Hilbert space. The unitary group $\U(\Ki)$,
equipped with the strong operator topology, is a topological group.
Let $\alpha \, \colon \widehat N \times H \to \U(\Ki)$
be a cocycle over $(\widehat N, \mu)$ with values in $\U(\Ki)$.
For every $g = (h, n) \in G$,
we define an operator $\pi_{\mu, \alpha}(g)$
on the Hilbert space $L^2(\widehat N, \mu, \Ki)$ by 
$$
(\pi_{\mu, \alpha}(h, n)F) (\chi)
\, = \,
\chi (n) \hskip.1cm c(\chi, h)^{1/2} \hskip.1cm \alpha(\chi, h) \hskip.1cm (F(\chi^{h}))
$$
for all $F \in L^2(\widehat N, \mu, \Ki)$ and $\chi \in \widehat N$.
Using the cocycle identities (\ref{eqq/RNcocR}) and (\ref{eqq/cocycle}),
one checks as in Section~\ref{Section-MoreIrredRep} that
$$
H \ltimes N \to \U (L^2(\widehat N, \mu, \Ki)), \hskip.2cm
(h, n) \mapsto \pi_{\mu, \alpha}(h, n)
$$
is a representation of $H \ltimes N$ on $L^2(\widehat N, \mu, \Ki)$.
\end{constr}

\begin{rem}
% 5.B.5
\label{Rem-RepCoc}
Observe that $\pi_{\mu, \alpha}$ is an extension to $G$
of the canonical representation $\pi_\mu^{(k)}$
of the LCA group $N$ associated to $(\mu, \Ki)$
and considered in Sections~\ref{Section-CanRepAbGr}
and \ref{Section-CanDecRepAbGr}, where $k$ is the dimension of $\Ki$. 
The results obtained in these sections about the representations $\pi_\mu^{(k)}$ 
will be useful in our study of the representations $\pi_{\mu, \alpha}$.
\end{rem}

As the next lemma shows, the representation $\pi_{\mu, \alpha}$
depends, up to equivalence, on the \emph{class} of $\mu$ only.

\begin{lem}
% 5.B.6
\label{Lem-EquiMeasuresEquiRep}
Let $\mu_1$ and $\mu_2$ be two probability measures
on $\widehat N$ which are quasi-invariant by $H$
and absolutely continuous with each other.
Let $\alpha \, \colon \widehat N \times H \to \U(\Ki)$ be a cocycle over 
$(\widehat N, \mu_1)$ and therefore also a cocycle over $(\widehat N, \mu_2)$,
with values in the unitary group of a Hilbert space $\Ki$.
\par

Then the representations $\pi_{\mu_1, \alpha}$ and $\pi_{\mu_2, \alpha}$
are equivalent.
\end{lem}

\begin{proof}
We have $\mu_1 = \varphi \mu_2$,
where $\varphi := \dfrac{d\mu_1}{d\mu_2}$
is the Radon--Nikodym derivative of $\mu_1$ with respect to $\mu_2$.
The operator 
$$
T \, \colon \, L^2(\widehat N, \mu_1, \Ki) \to L^2(\widehat N, \mu_2, \Ki),
$$
defined by 
$$
T(F) (\chi) \, = \, \varphi(\chi)^{1/2} F(\chi)
\hskip.3cm \text{for all} \hskip.2cm
F \in L^2(\widehat N, \mu_1, \Ki), \hskip.1cm \chi \in \widehat N
$$
is an isomorphism of Hilbert spaces.
Moreover, $T$ intertwines the representations
$\pi_{\mu_1, \alpha}$ and $\pi_{\mu_2, \alpha}$;
indeed, if
$c_{\mu_i}(\chi, h) = \dfrac{ d(\mu_i)_*(h_i^{-1}) } { d\mu_i }(x)$ for $i = 1, 2$,
we have (see Example \ref{examplecocycle}(2)) 
$$
c_{\mu_2}(\chi, h) \, = \,
\frac{ \varphi(\chi) }{ \varphi(\chi^h) } c_{\mu_1} (\chi, h)
$$
and hence
$$
\begin{aligned}
\big( T \pi_{\mu_1, \alpha}(h, n) (F) \big) (\chi)
\, &= \, \varphi(\chi)^{1/2} \big( \pi_{\mu_1, \alpha}(h, n) (F) \big) (\chi)
\\
\, &= \, \big( \varphi(\chi) c_{\mu_1} (\chi, h))^{1/2} \chi (n) \alpha(\chi, n) \big( F(\chi^{h}) \big)
\\
\, &= \, \big( \varphi(\chi^h) c_{\mu_2} (\chi, h))^{1/2} \chi (n) \alpha(\chi, n) \big( F(\chi^{h}) \big)
\\
\, &= \, \chi (n) c_{\mu_2} (\chi, h)^{1/2} \alpha(\chi, n) \big( T(F) (\chi^h) \big)
\\
\, &= \, \big( \pi_{\mu_2, \alpha}(h, n) T (F) \big) (\chi)
\end{aligned}
$$
for $h \in H, \hskip.1cm n \in N, \hskip.1cm F \in L^2(\widehat N, \mu_1, \Ki)$,
and $\chi \in \widehat N$.
It follows that 
$\pi_{\mu_1, \alpha}$ and $\pi_{\mu_2, \alpha}$ are equivalent.
\end{proof}

We proceed now towards Theorem \ref{Theo-AllRepSemiDirect1},
concerning irreducibility and equivalence
of representations of the form $\pi_{\mu, \alpha}$.

\begin{defn}
% 5.B.7
\label{defintertwinercocycle}
As above, let $\mu$ be a probability measure on $\widehat N$
which is quasi-invariant by the action of $H$.
We need to introduce the notion
of an intertwiner between cocycles over $(\widehat N, \mu)$,
taking their values in the unitary groups of possibly different Hilbert spaces.
We will also define \emph{irreducible} cocycles with values in such groups.
\par

Let
$$
\alpha_1 \, \colon \, \widehat N \times H \to \U(\Ki_1)
\hskip.5cm \text{and} \hskip.2cm
\alpha_2 \, \colon \, \widehat N \times H \to \U(\Ki_2)
$$
be measurable cocycles over $(\widehat N, \mu)$
with values in the unitary groups
of Hilbert spaces $\Ki_1$ and $\Ki_2$. 
\index{Intertwiner! $2$@between cocycles on a group}
An \textbf{intertwiner} between $\alpha_1$ and $\alpha_2$
is a decomposable operator 
$$
\widetilde T \, \colon \, L^2(\widehat N, \mu, \Ki_1) \, \to \, L^2(\widehat N, \mu, \Ki_2) ,
$$
associated to an essentially bounded measurable map
$$
T \, \colon \, \widehat N \to \Li (\Ki_1, \Ki_2),
\hskip.2cm \chi \mapsto T(\chi)
$$
as in Section \ref{Section-DecomposableOperators}, i.e., defined by
$$
(\widetilde T F) (\chi) \, = \, T(\chi) F(\chi)
\hskip.5cm \text{for all} \hskip.2cm
F \in L^2(\widetilde N, \mu, \Ki_1)
\hskip.2cm \text{and} \hskip.2cm
\chi \in \widetilde N ,
$$
such that the relation
$$
\alpha_2(\chi, h) T(\chi^{h}) \, = \,T(\chi)\alpha_1(\chi, h) 
$$
holds for all $h \in H$ and $\mu$-almost all $\chi \in \widetilde N$.
\par

Though $\alpha_1$ and $\alpha_2$ do not have values in the same unitary group,
\index{Cocycle! $2$@cohomologous}
we extend (slightly) the terminology defined earlier in this section:
the cocycles $\alpha_1$ and $\alpha_2$ are
\textbf{cohomologous} (or \textbf{equivalent})
if there exists an intertwiner $\widetilde T$ such that
$T(\chi) \, \colon \Ki_1 \to \Ki_2$ is an isomorphism of Hilbert spaces
for $\mu$-almost every $\chi \in \widehat N$.
\par

A measurable cocycle $\alpha \, \colon \widehat N \times H \to \U(\Ki)$
over $(\widehat N, \mu)$ is said to be \textbf{irreducible}
if every intertwiner $L^2(\widehat N, \mu, \Ki) \to L^2(\widehat N, \mu, \Ki)$
between $\alpha$ and itself belongs to the subalgebra 
$$
\mathcal A \, = \, \{m(\varphi) \mid \varphi \in L^\infty(\widehat N, \mu)\}
$$
of diagonalisable operators in $\Li (L^2(\widehat N, \mu, \Ki))$.
\index{Cocycle! $3$@irreducible}
\end{defn}

Given a probability measure $\mu$ on $\widehat N$
quasi-invariant by $H$,
the following lemma shows that the intertwiners
between two cocycles $\alpha_1, \alpha_2$ over $(\widehat N, \mu)$
as in Construction \ref{construcrep2}
coincide with the intertwiners
between the associated representations $\pi_{\mu, \alpha_1}, \pi_{\mu, \alpha_2}$.

\begin{lem}
% 5.B.8
\label{Lem-AllRepSemiDirect}
Let $G = H \ltimes N$ be a semi-direct product of a countable discrete subgroup $H$
with a second-countable locally compact abelian normal subgroup $N$.
Let $\mu$ be a probability measure on $\widehat N$
which is quasi-invariant by the $H$-action.
For $i = 1, 2$, let $\Ki_i$ be a separable Hilbert space,
and let $\alpha_i \, \colon \widehat N \times H \to \U(\Ki_i)$
be a cocycle over $(\widehat N, \mu)$.
\par

For an operator
$\widetilde T \, \colon L^2(\widehat N, \mu, \Ki_1) \to L^2(\widehat N, \mu, \Ki_2)$,
the following properties are equivalent:
\begin{enumerate}[label=(\roman*)]
\item\label{iDELem-AllRepSemiDirect}
$\widetilde T$ is an intertwiner 
between the cocycles $\alpha_1$ and $\alpha_2$;
\item\label{iiDELem-AllRepSemiDirect}
$\widetilde T$ is an intertwiner between the representations
$\pi_{\mu, \alpha_1}$ and $\pi_{\mu, \alpha_2}$.
\end{enumerate}
\end{lem}

\begin{proof}
Set $\Hi_i = L^2(\widehat N, \mu_i, \Ki_i)$, for $i = 1, 2$.
For $\varphi \in L^{\infty}(\widehat N, \mu_i)$, denote by $m_i(\varphi)$
the operator on $\Hi_i$ given by multiplication by $\varphi$,
as in Section~\ref{Section-DecomposableOperators}.
Observe that, for $n \in N$, we have 
$$
\pi_{\mu_i, \alpha_i}(n) \, = \, m_i(\widehat n).
$$
Recall that $\widehat n$ denotes the function in $C^b(\widehat N)$
defined by $\widehat n (\chi) = \chi(n)$ for $\chi \in \widehat N$.

\vskip.2cm

\ref{iDELem-AllRepSemiDirect} $\Rightarrow$ \ref{iiDELem-AllRepSemiDirect}
Let $\widetilde T \, \colon \Hi_1 \to \Hi_2$
be an intertwiner between $\alpha_1$ and $\alpha_2$,
associated to an essentially bounded measurable map
$T \, \colon \widehat N \to \Li (\Ki_1, \Ki_2)$.
As a decomposable operator, $\widetilde T$ is such that
$$
\widetilde T m_1(\widehat n) \, = \, m_2(\widehat n)\widetilde T
\hskip.5cm \text{for all} \hskip.2cm
n \in N
$$
(see Section \ref{Section-DecomposableOperators}),
so that $\widetilde T$ intertwines
$\pi_{\mu, \alpha_1} \vert_N$ and $\pi_{\mu, \alpha_2} \vert_N$.
\par

Let $h \in H$. We have
$$
\begin{aligned}
&
\big( \pi_{\mu, \alpha_2}(h, e) \widetilde T (F) \big) (\chi)
\, = \, c(\chi, h)^{1/2} \alpha_2 (\chi, h) \big( \widetilde T (F) (\chi^{h}) \big)
\\
& \hskip1cm
\, = \, c(\chi, h)^{1/2} \alpha_2(\chi, h) \big( T(\chi^h) (F(\chi^{h}) \big)
\, = \, c(\chi, h)^{1/2} \big( T(\chi) \alpha_1(\chi, h) \big) (F(\chi^h))
\\
& \hskip1cm
\, = \, T(\chi) \big( (\pi_{\mu, \alpha_1}(h, e)F) (\chi) \big)
\, = \, \big( \widetilde T \pi_{\mu, \alpha_1}(h, e) (F) \big) (\chi)
\end{aligned}
$$
for all $F \in \Hi_1$ and $\mu$-almost all $\chi \in \widehat N$.
Therefore 
$$
\pi_{\mu, \alpha_2}(h, e) \widetilde T \, = \,
\widetilde T \pi_{\mu, \alpha_1}(h, e)
\hskip.5cm \text{for all} \hskip.2cm
h \in H .
$$
It follows that $\widetilde T$ intertwines
$\pi_{\mu, \alpha_1}$ and $\pi_{\mu, \alpha_2}$.

\vskip.2cm

\ref{iiDELem-AllRepSemiDirect} $\Rightarrow$ \ref{iDELem-AllRepSemiDirect}
Conversely, let $\widetilde T \in \Li (\Hi_1, \Hi_2)$
be an intertwiner between $\pi_{\mu, \alpha_1}$ and $\pi_{\mu, \alpha_2}$.
Since $\widetilde T$ intertwines
$\pi_{\mu, \alpha_1} \vert_N$ and $\pi_{\mu, \alpha_2} \vert_N$,
and since $\pi_{\mu, \alpha_1} \vert_N$ and $\pi_{\mu, \alpha_2} \vert_N$
coincide with the canonical representations
$\pi_\mu^{(n_1)}$ and $\pi_\mu^{(n_2)}$ from Section~\ref{Section-CanDecRepAbGr},
where $n_1$ and $n_2$ are the dimensions of $\Ki_1$ and $\Ki_2$, 
it follows from Proposition~\ref{Prop-RepAb-Equiv} 
that $\widetilde T$ is a decomposable operator,
associated to an essentially bounded measurable map
$$
\widehat N \, \to \, \Li (\Ki_1, \Ki_2), \hskip.2cm \chi \mapsto T(\chi) .
$$
Let $h \in H$.
Since $\pi_{\mu, \alpha_2}(h, e) \widetilde T = \widetilde T \pi_{\mu, \alpha_1}(h, e)$,
we have as above
$$
\begin{aligned}
&
c(\chi, h)^{1/2} T(\chi) \alpha_1(\chi, h) (F(\chi^{h}))
\, = \, T(\chi) \Big( (\pi_{\mu, \alpha_1}(h, e)F) (\chi) \Big)
\\
& \hskip1cm
\, = \, \Big( \widetilde T \pi_{\mu, \alpha_1}(h, e) (F) \Big) (\chi)
\, = \, \big( \pi_{\mu, \alpha_2}(h, e) \widetilde T (F) \big) (\chi)
\\
& \hskip1cm
\, = \, c(\chi, h)^{1/2} \alpha_2(\chi, h) (\widetilde T (F)) (\chi^{h})
\, = \, c(\chi, h)^{1/2} \alpha_2(\chi, h)T(\chi^h) (F(\chi^{h}))
\end{aligned}
$$
for all $F \in L^2(\widehat N, \mu, \Ki_1)$
and $\mu$-almost all $\chi \in \widehat N$.
It follows that 
$$
T(\chi)\alpha_1(\chi, h) \, = \, \alpha_2(\chi, h)T(\chi^h)
$$
for every $h\in H$ and $\mu$-almost every $\chi \in \widehat N$;
hence, $\widetilde T$ is an intertwiner between $\alpha_1$ and $\alpha_2$.
\end{proof} 
 
The first main result in this chapter concerns the irreducibility
of representations of the form $\pi_{\mu, \alpha}$ and equivalence
between such representations.

\begin{theorem}
% 5.B.9
\label{Theo-AllRepSemiDirect1}
Let $G = H \ltimes N$ be a semi-direct product
of a countable discrete subgroup $H$ with a 
second-countable locally compact abelian normal subgroup~$N$.
\par

Let $\mu, \mu_1, \mu_2$ be probability measures on $\widehat N$
which are quasi-invariant by the $H$-action,
$\Ki, \Ki_1, \Ki_2$ separable Hilbert spaces,
$\alpha, \alpha_1, \alpha_2$ cocycles over $\widehat N \times H$
with values in $\U(\Ki), \U(\Ki_1), \U(\Ki_2)$ respectively,
and $\pi_{\mu, \alpha}, \pi_{\mu_1, \alpha_1}, \pi_{\mu_2, \alpha_2}$
the corresponding representations, as in Construction \ref{construcrep2}.
\begin{enumerate}[label=(\arabic*)]
\item\label{1DETheo-AllRepSemiDirect1}
Suppose that $\mu$ is ergodic for the $H$-action.
Then the representation $\pi_{\mu, \alpha}$ is irreducible
if and only if the cocycle $\alpha$ is irreducible.
\item\label{2DETheo-AllRepSemiDirect1}
If the representations $\pi_{\mu_1, \alpha_1}$ and $\pi_{\mu_2, \alpha_2}$
are equivalent,
then the measures $\mu_1$ and $\mu_2$ are equivalent.
\item\label{3DETheo-AllRepSemiDirect1}
Suppose that $\mu_1 = \mu_2 = \mu$.
Then $\pi_{\mu, \alpha_1}$ and $\pi_{\mu, \alpha_2}$ are equivalent
if and only if the cocycles $\alpha_1, \alpha_2$ are cohomologous.
\end{enumerate}
\end{theorem}

\begin{proof}
\ref{1DETheo-AllRepSemiDirect1}
Set $\Hi = L^2(\widehat N, \mu, \Ki)$. 
Assume that $\alpha$ is irreducible
and let $\widetilde T \in \Li (\Hi)$ be an intertwining operator for $\pi$ with itself. 
It follows from Lemma~\ref{Lem-AllRepSemiDirect}
that $\widetilde T$ is an intertwiner between $\alpha$ and itself.
Since $\alpha$ is irreducible, $\widetilde T$ is a diagonalisable operator:
there exists $\varphi \in L^\infty(\widehat N, \mu)$ such that
$$
\widetilde T \, = \, m(\varphi).
$$
The fact $\widetilde T$ is an intertwiner between $\alpha$ and itself
means that
$$
\varphi(\chi^h)\alpha(\chi, h) \, = \, \varphi(\chi) \alpha(\chi, h),
$$
that is, since $\alpha(\chi, h)$ is a non-zero operator,
$$
\varphi(\chi^{h}) \, = \, \varphi(\chi),
$$
for every $h \in H$ and $\mu$-almost every $\chi \in \widehat N$.
Since the $H$-action on $\widehat N$ is ergodic,
it follows that $\varphi$ is constant $\mu$-almost everywhere
(see Propositions \ref{ergodicity1=2} and \ref{prop-ergodicity}),
and hence $\widetilde T$ is a scalar operator.
This shows that $\pi_{\mu, \alpha}$ is irreducible.
\par

Conversely, assume that $\pi_{\mu, \alpha}$ is irreducible
and let $\widetilde T$ be an intertwiner between $\alpha$ and itself.
It follows from Lemma~\ref{Lem-AllRepSemiDirect} that $\widetilde T$ 
intertwines $\pi$ with itself.
Therefore $\widetilde T$ is scalar operator, by irreducibility of $\pi$.
This implies that $\alpha$ is irreducible.

\vskip.2cm

\ref{2DETheo-AllRepSemiDirect1}
Assume that $\pi_{\mu_1, \alpha_1}$ and $\pi_{\mu_2, \alpha_2}$ are equivalent.
Then, in particular, $\pi_{\mu_1, \alpha_1} \vert_N$ and $\pi_{\mu_2, \alpha_2} \vert_N$ are equivalent.
Therefore $\mu_1$ and $\mu_2$ are equivalent, by Proposition~\ref{Prop-RepAb-Subrep}.

\vskip.2cm

\ref{3DETheo-AllRepSemiDirect1}
Assume that $\pi_{\mu, \alpha_1}$ and $\pi_{\mu, \alpha_2}$ are equivalent.
There exists a Hilbert space isomorphism
$\widetilde T \, \colon \Hi_1 \to \Hi_2$ 
which intertwines $\pi_{\mu, \alpha_1}$ and $\pi_{\mu, \alpha_2}$,
where $\Hi_i = L^2(\widehat N, \mu, \Ki_i)$ for $i = 1, 2$.
It follows from Lemma~\ref{Lem-AllRepSemiDirect}
that $\widetilde T$ is an intertwiner between $\alpha_1$ and $\alpha_2$.
In particular, $\widetilde T$ is decomposable. 
Since $\widetilde T$ is an isomorphism between $\Hi_1$ and $\Hi_2$,
we see that $T(\chi)$ is an isomorphism between $\Ki_1$ and $\Ki_2$ 
for $\mu$-almost every $\chi \in \widehat N$,
where $\chi \mapsto T(\chi)$ is the field of operators from $\Ki_1$ to $\Ki_2$
associated to $\widetilde T$.
This shows that $\alpha_1$ and $\alpha_2$ are cohomologous.
\par

Similarly, if $\alpha_1$ and $\alpha_2$ are cohomologous,
then Lemma~\ref{Lem-AllRepSemiDirect}
implies that $\pi_{\mu, \alpha_1}$ and $\pi_{\mu, \alpha_2}$ are equivalent.
\end{proof}

\begin{rem}
% 5.B.10
\label{Rem-notergnotirr}
The ergodicity assumption in
\ref{Theo-AllRepSemiDirect1}~\ref{1DETheo-AllRepSemiDirect1} is necessary.
More precisely:
let $\mu$ be a probability measure on $\widehat N$
which is quasi-invariant by the $H$-action, and not ergodic,
and let $\alpha \, \colon \widehat N \times H \to \U(\Ki)$
be a cocycle over $(\widehat N, \mu)$;
then $\pi_{\mu, \alpha}$ is not irreducible. 
\par

Indeed, we can write $\mu = (\mu_1 + \mu_2)/2$
for $H$-quasi-invariant and mutually singular probability
measures $\mu_1$ and $\mu_2$.
We have a decomposition 
$$
L^2(\widehat N, \mu, \Ki) \, = \, 
L^2(\widehat N, \mu_1, \Ki) \oplus L^2(\widehat N, \mu_2, \Ki),
$$
into $\pi_{\mu, \alpha}(G)$-invariant non-trivial subspaces.
\end{rem}

Since every cocycle with values in the circle group $\T$ is irreducible,
the following corollary is a straightforward consequence
of Theorem~\ref{Theo-AllRepSemiDirect1}.

\begin{cor}
% 5.B.11
\label{Cor-AllRepSemiDirect1}
Let $G = H \ltimes N$ be as in Theorem~\ref{Theo-AllRepSemiDirect1}.
Let $\mu$ be a probability measure on $\widehat N$
which is quasi-invariant and ergodic by the $H$-action 
and let $\alpha \, \colon \widehat N \times H \to \T$
be a cocycle over $(\widehat N, \mu)$ with values in the circle group.
\begin{enumerate}[label=(\arabic*)]
\item\label{1DECor-AllRepSemiDirect1}
The representation $\pi_{\mu, \alpha}$ of $G$ on $L^2(\widehat N, \mu)$ is irreducible.
\par
Moreover, if $\beta \, \colon \widehat N \times H \to \T$
is another cocycle over $(\widehat N, \mu)$ which is not cohomologous to $\alpha$,
then $\pi_{\mu, \alpha}$ and $\pi_{\mu, \beta}$ are not equivalent.
\item\label{2DECor-AllRepSemiDirect1}
Let $\nu$ be another probability measure on $\widehat N$
which is quasi-invariant and ergodic by the $H$-action,
and not equivalent to $\mu$.
For any cocycle $\beta \, \colon \widehat N \times H \to \T$ over $(\widehat N, \nu)$,
the irreducible representations $\pi_{\mu, \alpha}$ and $\pi_{\nu, \beta}$
are not equivalent.
\end{enumerate}
\end{cor}

\begin{exe}
% 5.B.12
\label{Example-IrredRepDiscreteHeisenberg}
Let $\Gamma = H(\Z)$ be the Heisenberg group over $\Z$.
Recall that $\Gamma$ is the semi-direct product $H \ltimes N$
and that its centre is $Z$,
where 
$$
\begin{aligned}
H \, &= \, \{(a, 0, 0) \in \Gamma \mid a \in \Z\} ,
\\
N \, &= \, \{(0, b, c) \in \Gamma \mid b, c \in \Z\} ,
\\
Z \, &= \, \{(0, b, c) \in \Gamma \mid c \in \Z\} .
\end{aligned}
$$
Fix $\theta \in \mathopen[ 0,1 \mathclose[$ irrational. Set
$$
\omega \, := \, e^{2 \pi i \theta} \in \T ,
$$
so that the irrational rotation $R_\theta \, \colon \T \to \T$ is given by 
$z \mapsto z \omega$.
Consider the character
$$
\psi_\theta \, \colon \, Z \to \T ,
\hskip.5cm
(0, 0, c) \mapsto \omega^c = e^{2 \pi i \theta c} .
$$
The subset $\widehat N (\psi_\theta)$ of $\widehat N$
is the set of unitary characters $\chi$ of $N$ with restriction $\psi_\theta$ on $Z$.
For each $z \in \T$, we have a character
$$
\chi_z \in \widehat N (\psi_\theta),
\hskip.5cm
\chi_z(0, b, c) = z^b \omega^c
\hskip.5cm \text{for all} \hskip.2cm
b,c \in \Z .
$$
The spaces $\T$ and $\widehat N (\psi_\theta)$
can be identified by means of the map
$$
\T \, \to \, \widehat N (\psi_\theta), \hskip.2cm z \, \mapsto \, \chi_z .
$$
The right action
$\widehat N (\psi_\theta) \curvearrowleft H, \hskip.1cm (\chi_z, h) \mapsto \chi_z^h$
is given by
$$
\chi_z^h \, = \, \chi_{z \omega^a}
\, = \, \chi_{R_\theta^a z} \, \in \, \widehat N (\psi_\theta)
\hskip.5cm \text{for all} \hskip.2cm 
h = (a, 0, 0) \in H
\hskip.2cm \text{and} \hskip.2cm
z \in \T .
$$
Indeed, we have
$$
\begin{aligned}
\chi_z^h (0, b, c) 
\, &= \, \chi_z( h (0, b, c) h^{-1} )
\, = \, \chi_z(0, b, c + ab)
\\
\, &= \, z^b \omega^{c + ab}
\, = \, (z \omega^a)^b \omega^c
\, = \, \chi_{z \omega^a} (0, b, c) .
\end{aligned}
$$
Let $\mu$ be the image of the Lebesgue measure $\lambda$ on $\T$
under the map $z \mapsto \chi_z$;
we view $\mu$ as a measure on $\widehat N$ which vanish
on Borel sets disjoint from $\widehat N (\psi_\theta)$,
and we identify the probability space $(\widehat N, \mu)$ with $(\T, \lambda)$.
Then $\mu$ is an $H$-invariant ergodic probability measure on $\widehat N$,
since $\lambda$ is an $R_\theta$-invariant ergodic probability measure on $\T$.
The representation $\pi_\mu$ of $\Gamma$ associated to $\mu$
as in Construction~\ref{construcrep1} is given on $L^2(\T, \lambda)$ by the formula
$$
(\pi_\mu(a, b, c)f) (z) \, = \, \chi_z(0, b, c) f(z \omega^a)
\, = \, z^b \omega^c f(z \omega^a)
$$
for all $f \in L^2(\T, \lambda)$, $(a, b, c) \in \Gamma$, and $z \in \T$.
By Corollary~\ref{Cor-AllRepSemiDirect1}~\ref{1DECor-AllRepSemiDirect1},
the representation $\pi_\mu$ is irreducible.
\end{exe}

\begin{exe}
% 5.B.13
\label{Example-IrredRepDiscreteHeisenbergBIS}
In the situation of the previous example,
consider moreover the cocycle $\alpha \, \colon \T \times \Z \to \T$
over $(\widehat N, \mu)$
associated to the identity map $f \, \colon \T \to \T$,
as in Example~\ref{examplecocycle}(5).
Then $\alpha$ is given by
$$
\alpha (z, a) \, = \,
\begin{cases}
\prod_{i = 0}^{a - 1} R_\theta^i(z) = z^{a} \omega^{a(a - 1)/2} 
&\text{if} \hskip.2cm a > 0
\\
1
&\text{if} \hskip.2cm a = 0
\\
\prod_{i = 1}^{\vert a \vert} R_\theta^{-i}(z)^{-1} = z^{a} \omega^{a(a - 1)/2}
&\text{if} \hskip.2cm a < 0 .
\end{cases}
$$
The representation $\pi_{\mu, \alpha}$ of $\Gamma$
on $L^2(\T, \lambda) = L^2(\widehat N, \mu)$
associated to $\mu$ and $\alpha$ is given by
$$
\begin{aligned}
(\pi_{\mu, \alpha}(a, b, c) f) (z)
\, &= \,
\chi_z (0, b, c) \alpha(z,a) f(z \omega^a) \, = \,
z^b \omega^c \alpha(z,a) f(z \omega^a)
\\
\, &= \,
z^{a + b} \omega^{a(a - 1)/2 + c} f(z \omega^a) 
\end{aligned}
$$
for all $(a, b, c) \in \Gamma$, $f \in L^2(\T, \lambda)$, and $z \in \T$.
The representation $\pi_{\mu, \alpha}$ is irreducible,
again by Corollary~\ref{Cor-AllRepSemiDirect1}.
It is not equivalent to $\pi_\mu$,
since $\alpha$ is not a coboundary
(Example \ref{examplecocycle}(6)).
In Remark~\ref{Rem-MoreIrredRep},
we will also observe that $\pi_{\mu, \alpha}$ is not equivalent
to any one of the irreducible representations of $\Gamma$ 
of Section~\ref{Section-IrrRepTwoStepNil}.
\par

More generally, for an integer $i \ne 0$, let $\alpha_i \, \colon \T \times \Z \to \T$
be the cocycle over $(\widehat N, \mu)$ associated to
the map $\T \to \T, \hskip.1cm z \mapsto z^i$.
The associated representations $\pi_{\mu, \alpha_i}$ of $H(\Z)$
are irreducible and are mutually non equivalent, since the $\alpha_i$~'s
are mutually non cohomologous.
Moreover, no $\pi_{\mu, \alpha_i}$ is equivalent
to any one of the irreducible representations
of $\Gamma$ of Section~\ref{Section-IrrRepTwoStepNil}.
\par

For more representations of $H(\Z)$,
see Remark \ref{Rem-MoreIrredRep} (1) below. 
\end{exe}

As second main result in this chapter,
we establish the remarkable fact that every factor representation
of a semi-direct product $G = H \ltimes N$ as above
is equivalent to a representation $\pi_{\mu, \alpha}$
of Construction \ref{construcrep2}.

\begin{theorem}
% 5.B.14
\label{Theo-AllRepSemiDirect2}
Let $G = H \ltimes N$ be a semi-direct product of a countable discrete subgroup $H$
with a second-countable locally compact abelian normal subgroup~$N$.
Let $\pi$ be a \emph{factor} representation of $G$ in a separable Hilbert space $\Hi$.
\par

Then $\pi$ is equivalent to a representation of the form $\pi_{\mu, \alpha}$,
for a probability measure $\mu$ on $\widehat N$
which is quasi-invariant and ergodic by the $H$-action,
and a cocycle $\alpha \, \colon \widehat N \times H \to \U(\Ki)$
for a separable Hilbert space $\Ki$.
\end{theorem}

Recall from Theorem~\ref{Theo-AllRepSemiDirect1}
that $\alpha$ is irreducible if $\pi$ is irreducible.

\begin{proof}
By Proposition~\ref{Prop-RestNormalSub-Bis},
there exists a probability measure $\mu$ on $\widehat N$
which is quasi-invariant and ergodic by the $H$-action
and a Hilbert space $\Ki$ such that 
the restriction $\pi \vert_N$ of $\pi$ to $N$ is equivalent
to the representation $\pi_{\mu, I} \vert_N$
on $L^2(\widehat N, \mu, \Ki)$,
where $I$ is the trivial cocycle on $(\widehat N, \mu)$,
defined by $I(\chi, h) = \mathrm{Id}_{\Ki}$ for all $\chi \in \widehat N$ and $h \in H$.
By definition, $\pi_{\mu, I}$ is given by 
$$
(\pi_{\mu, I}(h, n)F) (\chi) \, = \, \chi (n) c(\chi, h)^{1/2} F(\chi^{h})
$$
for all $(h, n) \in H \ltimes N$ and $F \in L^2(\widehat N, \mu, \Ki)$.
We can therefore assume without loss of generality that 
$\Hi = L^2(\widehat N, \mu, \Ki)$ and $\pi \vert_N = \pi_{\mu, I} \vert_N$.
\par

Set 
$$
\widetilde{T_h} \, := \, \pi(h) \pi_{\mu, I}(h^{-1}) \in \U(\Hi) 
\hskip.5cm \text{for} \hskip.2cm 
h \in H .
$$
For every $n \in N$, we have $\pi(h^{-1}nh) = \pi_{\mu, I}(h^{-1}nh)$
since $h^{-1}nh\in N$, and therefore
$$
\begin{aligned}
\widetilde{T_h} \pi_{\mu, I}(n)
\, &= \, \pi(h)\pi_{\mu, I}(h^{-1})\pi_{\mu, I}(n)
\, = \, \pi(h)\pi_{\mu, I}(h^{-1}nh) \pi_{\mu, I}(h^{-1})
\\
\, &= \, \pi(h)\pi(h^{-1}nh) \pi_{\mu, I}(h^{-1})
\, = \, \pi(n)\pi(h) \pi_{\mu, I}(h^{-1})
\, = \, \pi(n) \widetilde{T_h}
\\
\, &= \, \pi_{\mu, I}(n) \widetilde{T_h}.
\end{aligned}
$$
It follows from Proposition~\ref{Prop-RepAb-Equiv}
that $\widetilde{T_h}$ is a decomposable operator.
Let
$$
T_h \, \colon \, \widehat N \to \Li (\Ki), \hskip.2cm \chi \mapsto T_h(\chi)
$$
be an essentially bounded measurable map
to which $\widetilde{T_h}$ is associated.
For every $h \in H$ and $\chi \in \widehat N$,
we have $T_h(\chi) \in \U(\Ki)$ since $\widetilde{T_h} \in \U(\Hi)$. 
\par

Define 
$$
\alpha \, \colon \, \widehat N \times H \to \U(\Ki),
\hskip.2cm (\chi, h) \mapsto T_h(\chi) .
$$
Observe that $\alpha$ is measurable.
We claim that $\alpha$ is a cocycle over $(\widehat N, \mu)$.
Indeed, let $h_1, h_2 \in H$. 
Since $\pi$ and $\pi_{\mu, I}$ are representations of $G$, we have
$$
\begin{aligned}
\widetilde{T_{h_1h_2}}
\, &= \, \pi(h_1) \pi(h_2)\pi_{\mu, I}(h_2^{-1}) \pi_{\mu, I}(h_1^{-1})
\\
\, &= \, \pi(h_1)\pi_{\mu, I}(h_1^{-1})\pi_{\mu, I}(h_1)\pi(h_2)
\pi_{\mu, I}(h_2^{-1}) \pi_{\mu, I}(h_1^{-1})
\\
\, &= \, \widetilde{T_{h_1}} \pi_{\mu, I}(h_1) \widetilde{T_{h_2}} \pi_{\mu, I}(h_1^{-1}).
\end{aligned}
$$
It follows that, for every $F \in L^2(\widehat N, \mu, \Ki)$
and $\chi \in \widehat N$, we have
$$
\begin{aligned}
T_{h_1h_2}(\chi) (F(\chi))
\, &= \,
(\widetilde{T_{h_1h_2}} F) (\chi)
\\
\, &= \,
\Big( \widetilde{T_{h_1}}\pi_{\mu, I}(h_1) \widetilde{T_{h_2}}\pi_{\mu, I}(h_1^{-1}) F \Big) (\chi)
\\
\, &= \,
T_{h_1} (\chi) \Big( \pi_{\mu, I}(h_1) \widetilde{T_{h_2}}\pi_{\mu, I}(h_1^{-1}) F \Big) (\chi)
\\
\, &= \,
c(\chi, h_1)^{1/2} T_{h_1}(\chi) \Big( (\widetilde{T_{h_2}}\pi_{\mu, I}(h_1^{-1}) F) (\chi^{h_1}) \Big)
\\
\, &= \,
c(\chi, h_1)^{1/2} T_{h_1}(\chi) \Big( T_{h_2}(\chi^{h_1}) (\pi_{\mu, I}(h_1^{-1}) F) (\chi^{h_1}) \Big)
\\
\, &= \,
c(\chi, h_1)^{1/2} c(\chi^{h_1}, h_1^{-1})^{1/2} T_{h_1}(\chi) \Big( T_{h_2}(\chi^{h_1}) (F(\chi)) \Big)
\\
\, &= \,
c(\chi, h_1h_1^{-1})^{1/2} T_{h_1}(\chi) \Big( T_{h_2}(\chi^{h_1}) (F(\chi)) \Big)
\\
\, &= \,
T_{h_1}(\chi)T_{h_2}(\chi^{h_1}) (F(\chi)).
\end{aligned}
$$
This implies that, for $\mu$-almost every $\chi \in \widehat N$, we have
$$
T_{h_1h_2}(\chi) \, = \, T_{h_1}(\chi)T_{h_2}(\chi^{h_1})
$$
and therefore
$$
\alpha( \chi, h_1h_2) \, = \, T_{h_1h_2}(\chi)
\, = \, T_{h_1}(\chi)T_{h_2}(\chi^{h_1})
\, = \, \alpha(\chi, h_1) \alpha(\chi^{h_1}, h_2).
$$
Therefore $\alpha$ satisfies the cocycle relation. 
\par

Next, we claim that $\pi$ coincides
with the representation $\pi_{\mu, \alpha}$ on $L^2(\widehat N, \mu, \Ki)$.
Indeed, since $\pi(h) = \widetilde{T_{h}} \pi_{\mu, I}(h)$, we have,
for every $h \in H$ and $F \in L^2(\widehat N, \mu, \Ki)$, 
$$
\begin{aligned}
(\pi(h) F) (\chi)
\, &= \, T_h(\chi) ((\pi_{\mu, I}(h) F) (\chi))
\, = \,c(\chi, h)^{1/2}T_h(\chi) (F(\chi^{h}))
\\
\, &= \,c(\chi, h)^{1/2}\alpha(\chi, h) (F(\chi^{h}))
\, = \,(\pi_{\mu, \alpha}(h)F) (\chi).
\end{aligned}
$$
Therefore $\pi \vert_H = \pi_{\mu, \alpha} \vert_H$.
Together with $\pi \vert_N = \pi_{\mu, \alpha} \vert_N$,
this implies that $\pi = \pi_{\mu, \alpha}$.
\end{proof}

\begin{rem}
% 5.B.15
\label{Rem-AllRepSemiDirect}
(1)
Let $G = H \ltimes N$ be as in Theorem~\ref{Theo-AllRepSemiDirect2}.
In view of Theorems~\ref{Theo-AllRepSemiDirect1} and \ref{Theo-AllRepSemiDirect2},
the description of the dual $\widehat G$ of $G$
is reduced to the following two problems.
For every extended positive integer $n$ in $\{1, \hdots, +\infty \}$,
we choose a Hilbert space $\Hi_n$ of dimension $n$.
\begin{itemize}
\setlength\itemsep{0em}
\item[(a)]
Find a collection $\mathcal P$ of probability measures on $\widehat N$
which are quasi-invariant and ergodic for the $H$-action,
with the following properties:
\begin{enumerate}
\item[--]
every probability measure on $\widehat N$
which is quasi-invariant and ergodic by the $H$-action
is equivalent to a measure in $\mathcal P$,
\item[--]
two measures in $\mathcal P$ are not
equivalent to each other.
\end{enumerate}
\item[(b)]
For every $\mu$ in $\mathcal P$,
and every extended positive integer $n$,
% $n \in \{1, \hdots, +\infty \}$,
find a collection $\mathcal C_{\mu,n}$ of irreducible cocycles
$\alpha \, \colon H \times \widehat N \to \U(\Hi_n)$ over $(\widehat N, \mu)$, 
with the following properties:
\begin{enumerate}
\item[--]
every such cocycle is cohomologous to a cocycle in $\mathcal C_{\mu,n}$,
\item[--]
two cocycles in $\mathcal C_{\mu,n}$ are not cohomologous.
\end{enumerate}
\end{itemize}
In general, as shown by the examples discussed
in Remark~\ref{Rem-MoreIrredRep} below,
both Problem~(a) and Problem (b) seem intractable.
\par

For instance, assume that $G$ is a discrete amenable group
and that $N$ is infinite.
Then the normalized Haar measure on $\widehat N$
is a non-atomic $H$-invariant probability measure
and it follows from Corollary~\ref{Cor-Theo-Schmidt} below
that there exist uncountably many non-equivalent probability measures
on $\widehat N$ which are quasi-invariant and ergodic by the $H$-action.
Their classification seems to be out of reach,
even in classical situations such as the one given by an irrational rotation 
of the circle equipped with Lebesgue measure
(Remark~\ref{Rem-MoreIrredRep}).
\par

Similarly, the problem of classifying the cocycles over $(\widehat N, \mu)$
for a fixed probability measure $\mu$ seems to be intractable.

\vskip.2cm

(2)
Theorem~\ref{Theo-AllRepSemiDirect2} is known to several experts.
It is a special case of results on groupoid representations in \cite{Rams--76}
(see especially \S~10 there).
It is also mentioned in \cite[\S~3]{Suth--83}
in case $H$ acts freely on $\widehat N \smallsetminus \{1 \}$.
\end{rem}

\section
{Identifying induced representations} 
% Section 5.C
\label{Section-IdentifyIrredRep}

As in Sections~\ref{Section-MoreIrredRep} and \ref{Section-AllIrredRep},
we consider a locally compact group of the form ${G = H \ltimes N}$,
where $N$ is a second-countable locally compact abelian normal subgroup
and $H$ a countable discrete subgroup of $G$.
\par

Among the irreducible representations $\pi_{\mu, \alpha}$ of $G$
of Theorem~\ref{Theo-AllRepSemiDirect1},
we will identify the ones which are obtained by inducing a unitary character of $N$.
As a result, the representations constructed in 
Sections~\ref{Section-IrrRepTwoStepNil}, \ref{Section-IrrRepAff},
and \ref{Section-IrrRepBS} 
correspond to representations of the form $\pi_{\mu, \alpha}$
for very special probability measures $\mu$ and trivial cocycles $\alpha$. 

\begin{prop}
% 5.C.1
\label{Prop-NonInduced}
Let $G = H\ltimes N$ be as above,
% Theorems \ref{Theo-AllRepSemiDirect1} and \ref{Theo-AllRepSemiDirect2},
$\mu$ a probability measure on~$\widehat N$
which is quasi-invariant and ergodic by the $H$-action,
$\Ki$ a Hilbert space, 
and $\alpha \, \colon H \times \widehat N \to \U(\Ki)$
an irreducible cocycle over $(\widehat N, \mu)$.
Let $\chi_0 \in \widehat N$
and let $X$ the $H$-orbit of $\chi_0$ in $\widehat N$. 
The following properties are equivalent:
\begin{enumerate}[label=(\roman*)]
\item\label{iDEProp-NonInduced}
the representation $\pi_{\mu, \alpha}$ is equivalent to
the induced representation $\Ind_N^G \chi_0$;
\item\label{iiDEProp-NonInduced}
the stabilizer of $\chi_0$ in $H$ is trivial,
$\mu$ is atomic and $\mu(X) = 1$,
the Hilbert space $\Ki$ is one-dimensional,
and $\alpha$ is cohomologous to the trivial cocycle
$H \times \widetilde N \to \T$ of constant value $1$.
\end{enumerate}
\end{prop}

\begin{proof}
We assume first that \ref{iiDEProp-NonInduced} holds
Then $\pi_{\mu, \alpha}$ is equivalent to $\pi_\mu = \pi_{\mu, I}$,
by Theorem~\ref{Theo-AllRepSemiDirect1}.
We have to show that $\pi := \Ind_N^G \chi_0$
is equivalent to $\pi_\mu$.
\par

We choose $H$ as transversal for $N \backslash G$,
and consider therefore the disjoint union $G = \bigsqcup_{h \in H} Nh$.
Let $(h_1, n_1) \in G$ and $h \in H$;
if we use the notation of Construction~\ref{constructionInd}(2) for
$$
(h, e)(h_1, n_1) \, = \, \alpha(n, (h_1, n_1)) (h \cdot (h_1, n_1)) \, \in \, N H
$$
and $h_1 \cdot_G n_1$ for the action of $H$ on $N$
defining the semi-direct product as in Convention \ref{convsemidirect},
we have
$$
\alpha( h, (h_1,n_1) ) \, = \, h \cdot_G n_1 \, = \, h n_1 h^{-1}
\hskip.5cm \text{and} \hskip.5cm
h \cdot (h_1, n_1) \, = \, hh_1 ,
$$
and in particular $\chi_0 (\alpha( h, (h_1, n_1) ) = \chi_0^h (n_1)$.
The induced representation $\pi$ is therefore realized on $\ell^2(H)$ by 
$$
(\pi(h_1,n_1)f) (h) \, = \, \chi_0^{h}(n_1)f(h h_1)
\hskip.5cm \text{for all} \hskip.2cm
f \in \ell^2(H), \hskip.1cm h_1, h \in H, \hskip.1cm n_1 \in N .
$$
\par

Since $\mu$ is atomic with $\mu(X) = 1$, 
we have $L^2(\widehat N, \mu) = \ell^2(X, \mu)$.
For $h \in H$ and $\chi \in \widehat N$, we have
$$
\mu(\{\chi^{h}\}) \, = \, c(\chi, h) \mu(\{ \chi \}),
$$
by definition of the Radon--Nikodym derivative
(see Theorem \ref{Radon--NikodymTheorem}).
We define a bijective linear map $T \, \colon \ell^2(X, \mu) \to \ell^2(H)$ by 
$$
(Tf) (h) \, = \, f(\chi_0^{h}) \mu(\{\chi_0^{h}\})^{1/2}
\, = \, f(\chi_0^{h}) c(\chi_0, h)^{1/2} \mu(\{\chi_0\})^{1/2}
$$
for all $f \in \ell^2(X, \mu)$ and $h \in H$.
\par

Then $T$ is unitary; indeed, since $\chi_0$ has trivial stabilizer in $H$
and since $X$ is the $H$-orbit of $\chi_0$, we have
$$
\Vert Tf \Vert ^2
\, = \, \sum_{h \in H} \vert (Tf) (h) \vert^2
\, = \, \sum_{h \in H} \vert f(\chi_0^{h}) \vert^2 \mu(\{\chi_0^{h}\})
\, = \, \sum_{\chi \in X} \vert f(\chi) \vert^2 \mu(\{ \chi \})
\, = \, \Vert f \Vert ^2,
$$
for every $f \in \ell^2(X, \mu)$. 
Moreover, for $(h_1,n_1) \in G$, $f \in \ell^2(X, \mu)$, and $h \in H$,
we have
$$
\begin{aligned}
T(\pi_\mu(h_1,n_1) f) (h)
\, &= \, (\pi_\mu(h_1,n_1) f) (\chi_0^{h}) \mu(\{\chi_0^{h}\})^{1/2}
\\
\, &= \, \chi_0^{h} (n_1) c(\chi_0^{h}, h_1)^{1/2}f(\chi_0^{hh_1}) \mu(\{\chi_0^{h}\})^{1/2}
\\
\, &= \, \chi_0^{h} (n_1) \frac{ \mu(\{\chi_0^{hh_1}\})^{1/2} }{ \mu(\{\chi_0^{h}\})^{1/2} }
 f(\chi_0^{hh_1}) \mu(\{\chi_0^{h}\})^{1/2}
\\
\, &= \, \chi_0^{h} (n_1)\mu(\{\chi_0^{hh_1}\})^{1/2} f(\chi_0^{hh_1})
\\
\, &= \, \chi_0^{h} (n_1) Tf(hh_1)
\\
\, &= \, (\pi(h_1,n_1)Tf) (h).
\end{aligned}
$$
This shows that $T$ intertwines $\pi_\mu$ and $\pi$;
hence, $\pi_\mu$ and $\pi = \Ind_N^G \chi_0$ are equivalent,
i.e., \ref{iDEProp-NonInduced} holds.

\vskip.2cm

Next, we show that \ref{iDEProp-NonInduced} implies \ref{iiDEProp-NonInduced}.
Assume that $\pi_{\mu, \alpha}$ and $\Ind_N^G \chi_0$ are equivalent.
In particular, $\Ind_N^G \chi_0$ is irreducible,
since $\pi_{\mu, \alpha}$ is irreducible by Theorem~\ref{Theo-AllRepSemiDirect1}.
It follows then from the Mackey--Shoda Criterion (Corollary~\ref{Cor-NormalSubgIrr})
that the stabilizer of $\chi_0$ in $H$ is trivial. 
\par

We write $\pi$ for $\Ind_N^G \chi_0$.
By Proposition~\ref{PropConjIndRep},
$\pi \vert_N$ is equivalent to the representation 
$$
\bigoplus_{h \in H} \chi_0^{h} \, = \, \bigoplus_{\chi \in X} \chi.
$$
So, $\pi_{\mu, \alpha} \vert_N$ is equivalent to $\bigoplus_{\chi \in X} \chi$.

The measure class attached to $\pi \vert_N$
as in Proposition~\ref{Prop-RestNormalSub-Bis} coincides with the class of $\mu$. 
Since $\pi_{\mu, \alpha} \vert_N$ is equivalent to $\bigoplus_{\chi \in X} \chi$,
it follows that $X$ is contained in the set of atoms of $\mu$.
Therefore, we have $\mu(X) = 1$, by the ergodicity of $\mu$.
The fact (proved above) that
\ref{iiDEProp-NonInduced} implies \ref{iDEProp-NonInduced}
shows then that $\pi_\mu = \pi_{\mu, I}$
is equivalent to the induced representation $\pi$. 
\par

As a consequence, $\pi_{\mu, \alpha}$ is equivalent to $\pi_{\mu, I}$.
It follows from Theorem~\ref{Theo-AllRepSemiDirect1}
that $\alpha$ is cohomologous to the trivial cocycle with values in $\T$
and hence that $\Ki$ is of dimension one.
\end{proof}

\begin{rem}
% 5.C.2
\label{Rem-MoreIrredRep}
Using Theorem~\ref{Theo-AllRepSemiDirect1},
we can construct uncountably many irreducible representations
of the discrete countable groups $\Gamma = H \ltimes N$
appearing in Sections \ref{Section-IrrRepTwoStepNil}
to \ref{Section-IrrRepLamplighter}, 
which are not equivalent to the representations constructed there by induction.

\vskip.2cm

(1)
Let $\Gamma = H(\Z) = H \ltimes N$ be the Heisenberg group over $\Z$, where 
$$
H \, = \, \{(a, 0, 0) \in \Gamma \mid a \in \Z\}
\hskip.5cm \text{and} \hskip.5cm
N \, = \, \{(0, b, c) \in \Gamma \mid b, c \in \Z\} 
$$
as before.
Fix $\theta \in \mathopen[ 0,1 \mathclose[$ irrational
and let as before $\widehat N (\psi_\theta)$ be the $H$-invariant set in $\widehat N$
of unitary characters $\chi$ of $N$
such that $\chi(0, 0, c) = e^{2 \pi i \theta c}$ for all $c \in \Z$. 
Then $\widehat N (\psi_\theta)$ is homeomorphic to $\T$,
and the action of $H$ on $\widehat N (\psi_\theta)$ corresponds to
the action of the irrational rotation $R_\theta$ on the circle group $\T$
(see Corollary~\ref{Cor-IrredRepHeisInteger}
and Example \ref{Example-IrredRepDiscreteHeisenberg}).
\par

We constructed in Example \ref{Example-IrredRepDiscreteHeisenbergBIS}
a countably infinite family of irreducible representations of $\Gamma$
which are non-equivalent to those appearing
in Corollaries \ref{Cor-IrredRepHeisInteger} and \ref{Cor-IrredFiniHeisIntegers}.
We now show that
there are \emph{uncountably many} such representations. 

\vskip.2cm

It is known that there are uncountably many pairwise non-equivalent
$R_\theta$-quasi-invariant and ergodic probability measures on $\T$
which are non-atomic. (See \cite{Kean--71} and \cite{KaWe--72};
this follows also from Corollary~\ref{Cor-Theo-Schmidt} below.)
\par

In fact, a more precise result is true, as follows from \cite[Theorem 4.11]{Krie--71}:
% Theorem (4.11) in \cite{Krie--71}: 
since the Lebesgue measure
is the unique $R_\theta$-invariant probability measure on $\T$,
the set of equivalence classes of
$R_\theta$-quasi-invariant and ergodic probability measures on $\T$,
equipped with a suitable Borel structure, is not countably separated.
This shows that the classification of equivalence classes
of $R_\theta$-quasi-invariant and ergodic probability measures on $\T$ is hopeless.
\par

In particular, Theorem~\ref{Theo-AllRepSemiDirect1} provides us
with uncountably many irreducible representations of $\Gamma$
which are non-equivalent to those appearing in 
Corollary~\ref{Cor-IrredRepHeisInteger} and Corollary~\ref{Cor-IrredFiniHeisIntegers}.

\vskip.2cm

We can construct further irreducible representations of $\Gamma$.
Indeed, as shown in Example~\ref{examplecocycle}(6), 
there are several pairwise non-cohomologous measurable cocycles 
$$
\alpha \, \colon \, \T \times \Z \to \T
$$
for $R_\theta$ over $(\T, \mu)$, where $\mu$ is the Lebesgue measure on $\T$
(see also \cite{Furs--61}, \cite{Veec--69}, and \cite{Brow--73a}).
By Corollary~\ref{Cor-AllRepSemiDirect1},
every such cocycle $\alpha$ defines an irreducible representation $\pi_{\mu, \alpha}$
as in Theorem~\ref{Theo-AllRepSemiDirect1}.
More generally, non-cohomologous irreducible cocycles 
$$
\alpha \, \colon \, \T \times \Z \to \U(n)
$$
with values in the unitary group $\U(n)$ exist for every $n \ge 1$ (see \cite{BaMe--86}).
In this way, we obtain new representations $\pi_{\mu, \alpha}$ of $\Gamma$.
This fact is also mentioned in Theorem 6.67 and the remarks on Page 218 of \cite{Foll--16}.

\vskip.2cm

The previous remarks also show that Theorem 10 in \cite{GoSV--16} is erroneous,
where it is claimed that every irreducible representation of $H(\Z)$
is equivalent to a representation of the form $\pi_\mu = \pi_{\mu, I}$,
for an ergodic quasi-invariant measure $\mu$ on $\widehat N$.

\vskip.2cm

(2)
Let $\Gamma = \Aff(\K) = \K^\times \ltimes \K$ be the affine group
of a countable infinite field~$\K$.
The normalized Haar measure on the compact group $\widehat \K$
is a $\K^\times$-invariant and ergodic probability measure.
Therefore, by Theorem~\ref{Theo-Schmidt}
%\cite[Corollary 10.7]{Schm--77a},
there exist uncountably many pairwise non-equivalent
$\sigma$-finite $\K^\times $-invariant and ergodic measures on $\widehat \K$
which are non-atomic.
It follows that there exist uncountably many pairwise non-equivalent 
$\K^\times $-quasi-invariant and ergodic probability measures on $\widehat \K$
which are non-atomic.
As above, we obtain uncountably many irreducible representations of $\Gamma$
which are non-equivalent to those of Theorem \ref{Prop-IrredRepAffine};

\vskip.2cm

(3)
Similarly, we obtain from Theorem~\ref{Theo-AllRepSemiDirect1}
uncountably many irreducible representations
of the Baumslag--Solitar group $\Gamma = \BS(1, p)$ and of the lamplighter group
which are non-equivalent to those given in
Section~\ref{Section-IrrRepBS} and Section~\ref{Section-IrrRepLamplighter}.
\end{rem}

\section[Infinite non-atomic ergodic invariant measures]
{On the existence of infinite non-atomic
% \\
measures which are invariant and ergodic}
% Section 5.D
\label{Section-TheoSchmidt}

Let $\Gamma$ be a countable \emph{amenable} group,
$(X, \mathcal B)$ a standard Borel space,
$\mu$ a non-atomic probability measure on $(X, \mathcal B)$,
and $\Gamma \curvearrowright (X, \mathcal B, \mu)$ an action
for which the measure $\mu$ is invariant and ergodic.
Following \cite{Schm--77b},
we will show that there are uncountably many
infinite $\sigma$-finite non-atomic measures on $(X, \mathcal B)$
which are invariant and ergodic under the action of $\Gamma$.
\par

\index{Orbit equivalence of group actions} 
For this, we need to recall the notion of orbit equivalence
for measure preserving group actions on probability spaces.
Let $\Gamma_1$ and $\Gamma_2$ be two countable groups
acting respectively on standard Borel spaces
$(X_1, \mathcal B_1)$ and $(X_2, \mathcal B_2)$
equipped with non-atomic probability measures $\mu_1$ and $\mu_2$
which are invariant and ergodic.
The actions $\Gamma_1 \curvearrowright (X_1, \mathcal B_1, \mu_1)$
and $\Gamma_2 \curvearrowright (X_2, \mathcal B_2, \mu_2)$
 are \textbf{orbit equivalent} if there exists a Borel isomorphism 
$\Phi \, \colon X_1 \to X_2$ such that $\Phi_*(\mu_1) = \mu_2$ and such that 
$$
\Phi(\Gamma_1 \cdot x) \, = \, \Gamma_2 \cdot \Phi(x)
$$
for $\mu_1$-almost every $x \in X_1$.
\par

The following theorem,
due to Dye and Ornstein--Weiss \cite{Dye--59, OrWe--80}, 
is the fundamental result in the theory of orbit equivalence of group actions.
\par

\begin{theorem}[\textbf{Dye and Ornstein--Weiss}]
% 5.D.1
\label{Theo-Dye} 
Let $\Gamma_1$ and $\Gamma_2$ be two amenable countable groups.
\par

Then any two non-atomic probability measure preserving ergodic actions
$\Gamma_1 \curvearrowright (X_1, \mathcal B_1, \mu_1)$
and $\Gamma_2 \curvearrowright (X_2, \mathcal B_2, \mu_2)$
are orbit equivalent.
\end{theorem}

Theorem~\ref{Theo-Dye} will be the main ingredient
in the proof of the following result from \cite{Schm--77b}.

\begin{theorem}[\textbf{K.~Schmidt}]
% 5.D.2
\label{Theo-Schmidt} 
Let $\Gamma$ be an amenable countable group,
$(X, \mathcal B)$ a standard Borel space,
$\mu$ a non-atomic probability measure on $(X, \mathcal B)$,
and $\Gamma \curvearrowright (X, \mathcal B)$ a measurable action.
Assume that $\mu$ is $\Gamma$-invariant and ergodic.
\par

Then there exist uncountably many mutually non-equivalent 
$\sigma$-finite and non-atomic infinite measures on $(X, \mathcal B)$
which are $\Gamma$-invariant and ergodic.
\end{theorem}

\begin{proof}
Set $\Lambda := \Gamma\times \Gamma$ and $Y := X \times X$.
Consider the product action of $\Lambda$ on $(Y, \mathcal B \otimes \mathcal B)$: 
$$
(\gamma_1 , \gamma_2) (x_1, x_2) \, = \, (\gamma_1 \cdot x_1, \gamma_2 \cdot x_2)
\hskip.5cm \text{for all} \hskip.2cm 
(\gamma_1, \gamma_2) \in \Lambda, \hskip.1cm (x_1,x_2) \in Y.
$$
This action preserves the probability measure $\nu := \mu \otimes \mu$
on $(Y, \mathcal B \otimes \mathcal B)$
and is obviously ergodic. 
\par

Since $\mu$ is non-atomic and ergodic,
it is well known that $\mu$-almost every $x \in X$
has an infinite orbit (see Proposition~\ref{Prop-Aperiodic}).
Fix $x_0 \in X$ with an infinite $\Gamma$-orbit. Set 
$$
Y_{x_0} \, := \,
\{ (\gamma \cdot x_0, x) \in X \times X \mid \gamma \in \Gamma, x \in X \} \, = \,
\bigsqcup_{t \in \mathcal{O}} \{t\} \times X, 
$$
where $\mathcal{O}$ is the $\Gamma$-orbit of $x_0$.
Observe that $Y_{x_0}$ is a $\Lambda$-invariant Borel subset of $Y$.
Define a Borel measure $\nu_{x_0}$ on $Y$ by 
$$
\int_Y f d\nu_{x_0} \, = \, \sum_{t \in \mathcal{O}} \int_X f(t, x) d\mu(x)
$$
for every positive Borel function $f$ on $Y$. 
Then:
\begin{itemize}
\setlength\itemsep{0em}
\item[$\bullet$]
$\nu_{x_0}(Y \smallsetminus Y_{x_0}) = 0$.
\item[$\bullet$]
$\nu_{x_0}$ is infinite, since 
$$
\int_Y \mathbf{1}_Y d\mu_{x_0}
\, = \, \sum_{t \in \mathcal{O}} \int_X \mathbf{1}_X(x) d\mu(x)
\, = \, \sum_{ t\in \mathcal{O}} 1
\, = \, \infty .
$$
\item[$\bullet$]
$\nu_{x_0}$ is $\sigma$-finite. Indeed, $\mathcal{O}$ is countable, 
$$
Y \, = \, \Big(Y \smallsetminus Y_{x_0} \Big)
\bigsqcup 
\Big( \bigsqcup_{t \in \mathcal{O}} \{t\} \times X \Big) ,
$$
we have $\nu_{x_0}(Y \smallsetminus Y_{x_0}) = 0$,
and $\nu_{x_0} (\{t\} \times X) = 1$ for every $t \in \mathcal{O}$.
\item[$\bullet$]
$\nu_{x_0}$ is non-atomic, since $\mu$ is non-atomic.
\item[$\bullet$]
$\nu_{x_0}$ is $\Lambda$-invariant.
Indeed, let $f$ be a positive Borel function on $Y$
and $(\gamma_1, \gamma_2)$ in $\Lambda$.
For the translate $_{(\gamma_1, \gamma_2)} f$ of $f$ by $(\gamma_1, \gamma_2)$,
we have, by $\Gamma$-invariance of $\mu$,
$$
\begin{aligned}
\int_{Y}{ _{(\gamma_1, \gamma_2)}} f d\nu_{x_0}
\, &= \, \sum_{t\in\mathcal{O}} \int_X{_{(\gamma_1, \gamma_2)}} f(t, x) d\mu(x)
\, = \, \sum_{t\in\mathcal{O}} \int_X f(\gamma_1^{-1} \cdot t, \gamma_2^{-1}\cdot x) d\mu(x)
\\
\, &= \, \sum_{t\in\mathcal{O}}\int_X f( \gamma_1^{-1} \cdot t, x) d\mu(x)
\, = \, \sum_{t\in\mathcal{O}}\int_X f(t, x) d\mu(x)
\, = \, \int_Y fd\nu_{x_0} .
\end{aligned}
$$
\item[$\bullet$]
$\nu_{x_0}$ is ergodic.
Indeed, let $B$ be a $\Lambda$-invariant measurable subset of $Y$.
For every $t \in \mathcal{O}$, the $t$-section
$$
B_{t} \, := \, \{x \in X \mid (t, x) \in B\}
$$
is a $\Gamma$-invariant measurable subset of $X$, since $B$ is $\Lambda$-invariant.
By ergodicity of $\mu$, it follows that
either $\mu(B_{t}) = 1$ or $\mu(B_{t}) = 0$ for every $t \in \mathcal{O}$.
Moreover, for every $\gamma \in \Gamma$, we have
$B_{\gamma \cdot t} = ((\gamma, e)^{-1} B)_{t}$
and hence
$$
B_{\gamma \cdot t} \, = \, B_{t} ,
$$
by $\Lambda$-invariance of $B$.
It follows that we have either $\mu(B_{t}) = 1$ for all $t \in \mathcal{O}$
or $\mu(B_{t}) = 0$ for all $t \in \mathcal{O}$.
Since $B = \bigsqcup_{t \in \mathcal{O}} \{t\}\times B_t$, we have in the first case
$$
\nu_{x_0} (Y \smallsetminus B)
\, = \, \sum_{t \in \mathcal{O}} \mu(X \smallsetminus B_t)
\, = \, 0
$$
and in the second case
$$
\nu_{x_0} (B)
\, = \, \sum_{t\in\mathcal{O}} \mu(B_{t})
\, = \, 0.
$$
\end{itemize} 
\par

By Theorem~\ref{Theo-Dye}, there exists a Borel isomorphism 
$\Phi \, \colon Y \to X$ such that $\Phi_*(\nu) = \mu$ and 
$$
\Phi (\Lambda \cdot y) \, = \, \Gamma \cdot \Phi(y)
$$
for $\nu$-almost every $y \in Y$.
Let $y = (x_0, x) \in Y$ be such that $\Phi(\Lambda\cdot y) = \Gamma\cdot \Phi(y)$. 
The image 
$$
\mu_{x_0} \, := \, \Phi_*(\nu_{x_0})
$$
of $\nu_{x_0}$ is a measure on $(X, \mathcal B)$,
which is clearly infinite, $\sigma$-finite and non-atomic,
since $\nu_{x_0}$ has these properties.
Moreover, the following hold:
\begin{itemize}
\setlength\itemsep{0em}
\item[$\bullet$]
$\mu_{x_0}$ is $\Gamma$-invariant.
Indeed, let $A$ be a Borel subset of $X$ and $\gamma \in \Gamma$.
Then $\mu_{x_0}(A) = \nu_{x_0}(B)$, where
$B := \Phi^{-1}(A)$. Since 
$$
\Phi(\Lambda\cdot B) \, = \, \Gamma\cdot \Phi(B) \, = \, \Gamma\cdot A,
$$
there exists a partition $B = \bigsqcup_{i} B_i$ in measurable subsets $B_i$ such that 
$$
\Phi^{-1}(\gamma \cdot A) \, = \, \bigsqcup_{i} \lambda_i B_i
$$
for some $\lambda_i \in \Lambda$.
Then, using the $\Lambda$-invariance of $\nu_{x_0}$, we have
$$
\begin{aligned}
\mu_{x_0}(\gamma A)
\, &= \, \nu_{x_0} (\Phi^{-1} (\gamma A))
\, = \,
\sum_i \nu_{x_0} (\lambda_i B_i)
\\
\, &= \,
\sum_i \nu_{x_0} (B_i)
\, = \,
\nu_{x_0}(B)
\, = \,
\mu_{x_0} (A).
\end{aligned}
$$
\item[$\bullet$]
$\mu_{x_0}$ is ergodic.
Indeed let $A$ be $\Gamma$-invariant Borel subset of $X$.
Then $B = \Phi^{-1}(A)$ is $\Lambda$-invariant, since 
$$
\Lambda\cdot B \, = \, \Phi^{-1}(\Gamma\cdot A) \, = \, \Phi^{-1}(A) \, = \, B.
$$
Therefore $\nu_{x_0}(Y \smallsetminus B) = 0$ or $\nu_{x_0}(B) = 0$,
that is, $\mu_{x_0}(X \smallsetminus A) = 0$ or $\mu_{x_0}( A) = 0$.
\end{itemize}
\par

Finally, let $x_0, x_0' \in X$ with distinct $\Gamma$-orbits.
Then, with the notation as above, we have
$$
\nu_{x_0}(Y_{x_0'}) \, = \, 0
\hskip.5cm \text{and} \hskip.5cm
\nu_{x_0'}(Y \smallsetminus Y_{x_0'}) \, = \, 0 ,
$$
so that the two measures $\nu_{x_0}$ and $\nu_{x_0'}$ are mutually singular.
It follows that $\mu_{x_0}$ and $\mu_{x_0'}$ are mutually singular.
Since there are uncountably many $\Gamma$-orbits in $X$,
we have constructed uncountably many
non-equivalent measures on $X$ with the required properties.
\end{proof}

In the situation of Theorem~\ref{Theo-Schmidt}, one is sometimes interested 
in the existence of an uncountable family of \emph{probability} measures 
on $X$ with similar properties (see Remark~\ref{Rem-MoreIrredRep}).
As we now show, this remains true
at the cost of replacing the invariance of the measures by their quasi-invariance. 
Observe that no better result can be expected as there might be
a \emph{unique} probability measure on $X$ which is invariant under the $\Gamma$-action;
these are called uniquely ergodic actions.
For example, the action of $\Z$ given by an irrational rotation on the circle
is uniquely ergodic (see \cite[Example 3.10]{BeMa--00}).

%The following consequence of Theorem~\ref{Theo-Schmidt}
% was used in Remark~\ref{Rem-MoreIrredRep}
\begin{cor}
% 5.D.3
\label{Cor-Theo-Schmidt} 
Let $\Gamma$ and $(X, \mathcal B, \mu)$ be as in Theorem~\ref{Theo-Schmidt}.
\par

Then there exist uncountably many mutually non-equivalent non-atomic
\emph{probability} measures on $(X, \mathcal B)$
which are $\Gamma$-quasi-invariant and ergodic.
\end{cor}

\begin{proof}
By Theorem~\ref{Theo-Schmidt}, there exists an uncountable family
$(\mu_\iota)_{\iota \in I}$ of non-equivalent 
$\sigma$-finite and non-atomic infinite measures on $(X, \mathcal B)$
which are $\Gamma$-invariant and ergodic.
By Lemma~\ref{PropEquivProbaMeasure},
for every $\iota \in I$, there exists a probability measure $\nu_\iota$
on $(X, \mathcal B)$ which is equivalent to $\mu_\iota$.
It is clear that every $\nu_\iota$ is quasi-invariant and ergodic by $\Gamma$
and that the $\nu_\iota$~'s are mutually non-equivalent. 
\end{proof}

\begin{rem}
% 5.D.4
\label{Rem-Schmidt}
For various strengthenings of Theorem~\ref{Theo-Schmidt}
and other results around the existence of infinite $\sigma$-finite invariant measures,
see the lecture notes \cite{Schm--77a}.
\end{rem}

%-----------------------------------------------------------------------
% End of chapter 5
%-----------------------------------------------------------------------
\chapter[Groups of type I]{Groups of type I}
% Chapter 6
\label{ChapterTypeI}

\emph{
Let $\pi$ be a representation of a second-countable LC group $G$ in a separable Hilbert space.
As seen in Section~\ref{SectionDecomposingIrreps},
$\pi$ can always be decomposed into a direct integral of irreducible representations,
but two such decompositions of $\pi$ may be completely different
(see Theorem~\ref{Theo-NonUniqueIntDecIrrUniRep}).
We will define in this chapter the class of type I representations of $G$,
for which a direct integral decomposition is essentially unique.
}
\par

\emph{
Following Mackey (see \cite{Mack--76} and \cite{Mack--78}),
we define in Section~\ref{Sectioncomppres} type~I representations
by means of the notion of subordination of representations,
which is weaker than containment,
and the corresponding relation of quasi-equival\-ence.
One should think of quasi-equivalent representations as representations
which have the same ``kinds" of subrepresentations.
Representations with only one ``kind" of subrepresentations are called factor representations;
more precisely, a representation is factorial if it cannot be 
decomposed as a direct sum of two disjoint representations.
Next, we introduce multiplicity-free representations:
a representation $\pi$ is multiplicity-free
% multiplicity-free in \cite{Dixm--C*} 5.4.4 et 5.4.5
if, for every decomposition $\pi = \pi_1 \oplus \pi_2$,
the subrepresentations $\pi_1$ and $\pi_2$ are disjoint.
Type I representations are defined as the representations
which are quasi-equivalent to a multiplicity-free representation.
}
\par

\emph{
For a representation $\pi$ of a topological group $G$ on a Hilbert space $\Hi$,
the commutant $\pi(G)'$ and the bicommutant $\pi(G)''$ of $\pi(G)$ in $\Li (\Hi)$
are von Neumann subalgebras of $\Li (\Hi)$.
This allows us to reformulate in Section~\ref{S:QE-Factor-VNAlgebras} 
all the notions defined in Section~\ref{Sectioncomppres}
%(subordination, quasi-equivalent representation, multiplicity-free representations,
% factor representations, and type I representations)
in terms of von Neumann algebras.
When $G$ is locally compact and second-countable,
every multiplicity-free representation of $G$
has an essentially unique decomposition as direct integral of irreducible representations;
moreover, every type I representation of $G$ has a canonical decomposition
as sum of multiplicity-free representations.
We also discuss multiplicity-free representations of abelian LC groups,
and multiplicity-free representations associated to Gel'fand pairs.
}
\par

\emph{
We introduce in Section~\ref{SectionQuasidual}
the quasi-dual $\QD(G)$ of a topological group $G$,
as the space of quasi-equivalence classes of factor representations of $G$;
it has a natural Borel structure, called the Mackey--Borel structure of $\QD(G)$.
Quasi-equivalence classes of factor representations of type I
constitute a subset of the quasi-dual which can be identified with the dual,
so that $\widehat G$ can be seen as a subspace of $\QD(G)$.
When $G$ is locally compact and second-countable,
every representation of $G$ has a canonical decomposition (the so-called central decomposition)
as a direct integral of factor representations,
to which is associated a class of measures on $\QD(G)$.
}
\par

\emph{
A topological (not necessarily LC) group $G$ is defined to be of type I
if all representations of $G$ are of type I;
equivalently, if every factor representation of $G$ is a multiple of an irreducible one
(Theorem~\ref{Theo-TypIFactorRep}).
For a second-countable locally compact group which is of type I,
every representation of $G$ has a canonical decomposition
as a direct integral of irreducible representations,
given by a family of mutually singular measures on $\widehat G$
and a multiplicity function (Theorem~\ref{thmDirectIntIrreps+}).
 }
\par

\emph{
We introduce in Section~\ref{S:ClassTypeI-GCRGroups}
the class of GCR (or postliminal) groups and show that every such group is of type I
(Theorem~\ref{Theo-GCR-Group}).
We discuss in detail the important subclass of CCR (or liminal) groups.
We give a complete proof of the fact that $\SL_2(\R)$ belongs to this subclass
and we indicate how these arguments can be extended to general 
semisimple connected real Lie groups. 
}

\section
{Comparing representations, quasi-equivalence}
% Section 6.A
\label{Sectioncomppres}

Let $G$ be a topological group.

\subsection
{Disjointness and quasi-equivalence}
% subsection 6.A.a
\label{SS:QE}

For a representation $\pi$ of $G$ and a cardinal $n$,
recall that $n \pi$ denotes the direct sum of~$n$ representations equivalent to~$\pi$.
When $G$ is second-countable, 
it suffices below to consider countable cardinals
$n \in \{1, 2, 3, \hdots, \infty \}$, where $\infty$ stands for $\aleph_0$.
For larger groups, larger cardinals are necessary; 
nevertheless, we will use below the ambiguous notation ``$\infty$''. 

\begin{defn}
% 6.A.1
\label{Def-Subordinate}
Let $\pi_1, \pi_2$ be two representations of $G$ in Hilbert spaces of positive dimensions.

\vskip.2cm

Recall from Section \ref{S-DefUnitD}
that $\pi_1, \pi_2$ are \textbf{disjoint}
if there does not exists a pair of
non-trivial subrepresentations $\rho_1 \le \pi_1$, $\rho_2 \le \pi_2$
such that $\rho_1$ and $\rho_2$ are equivalent.
\index{Representation! disjoint}
\index{Disjoint! representations}

\vskip.2cm

The representation $\pi_1$ is \textbf{subordinate} to $\pi_2$ 
if every non-trivial subrepresentation of $\pi_1$
contains a non-trivial subrepresentation 
equivalent to a subrepresentation of $\pi_2$,
i.e., if no non-trivial subrepresentation of $\pi_1$ is disjoint from $\pi_2$.
\index{Representation! subordinate}
\index{Subordinate representation}

\vskip.2cm

The representations $\pi_1$ and $\pi_2$ are 
\textbf{quasi-equivalent},
and we write $\pi_1 \approx \pi_2$,
if each one is subordinate to the other
\index{Representation! quasi-equivalent $\approx$}
\index{Quasi-equivalent! representations}
\index{$a4$@$\approx$ quasi-equivalence of representations}
\end{defn}

Here are a few straightforward immediate consequences of the definitions.

\begin{prop}
% 6.A.2
\label{Prop-EqQuasiEq-Immediate}
Let $\pi_1, \pi_2$ be two representations of~$G$.
\begin{enumerate}[label=(\arabic*)]
\item\label{1DEProp-EqQuasiEq-Immediate}
If $\pi_1$ is contained in $\pi_2$, then $\pi_1$ is subordinate to $\pi_2$.
If $\pi_1$ and $\pi_2$ are equivalent, then they are quasi-equivalent.
\item\label{2DEProp-EqQuasiEq-Immediate}
If $\pi_1$ is irreducible, then $\pi_1$ is subordinate to $\pi_2$
if and only if it is contained in $\pi_2$.
\item\label{3DEProp-EqQuasiEq-Immediate}
If $\pi_1$ and $\pi_2$ are irreducible, then they are either equivalent or disjoint.
Two irreducible representations are quasi-equivalent
if and only if they are equivalent.
\end{enumerate}
\end{prop}

\begin{rem}
% 6.A.3
\label{Les3equivrepr}
Beware that there are three distinct equivalence relations defined for representations.
Two representations $\pi_j \, \colon G \to \U(\Hi_j)$, $j = 1, 2$,
are defined to be
\begin{enumerate}
\item[$\bullet$]
equivalent, written $\pi_1 \simeq \pi_2$,
if there exists a unitary operator $V$ from $\Hi_1$ onto $\Hi_2$
such that $\pi_2(g) = V\pi_1(g)V^*$ for all $g \in G$,
as defined in Section \ref{S-DefUnitD};
\item[$\bullet$]
weakly equivalent, written $\pi_1 \sim \pi_2$,
if $\pi_1 \preceq \pi_2$ and $\pi_2 \preceq \pi_1$,
as defined in Section \ref{SectionWC+FellTop};
equivalently if $\textnormal{C*ker}(\pi_1) = \textnormal{C*ker}(\pi_2)$,
see Section \ref{C*algLCgroup};
\item[$\bullet$]
quasi-equivalent, written $\pi_1 \approx \pi_2$,
as defined above.
\end{enumerate}
\index{Schr\"oder--Bernstein property for representations}
Each of these three notions
has the appropriate \textbf{Schr\"oder--Bernstein property}.
More precisely, for equivalence,
if $\pi_1$ is a subrepresentation of $\pi_2$ and $\pi_2$ a subrepresentation of $\pi_1$,
then $\pi_1$ and $\pi_2$ are equivalent, by Proposition \ref{SchroderBernstein}.
Similarly, for weak equivalence:
if $\pi_1 \preceq \pi_2$ and $\pi_2 \preceq \pi_1$,
then $\pi_1 \sim \pi_2$, by definition.
Similarly, for quasi-equivalence:
if $\pi_1$ is subordinate to $\pi_2$ and $\pi_2$ subordinate to $\pi_1$,
then $\pi_1 \approx \pi_2$, by definition.
\end{rem}

Let us compare equivalence and quasi-equivalence,
as we did in Remark \ref{propweakcont} for equivalence and weak equivalence.
We do this first in a particularly simple situation, 
and then more generally in Proposition \ref{eqetsubordinate}.
\par

Let $\pi, \rho$ be two representations of $G$.
Assume that $\pi$ and $\rho$ can be written as direct sums of irreducible subrepresentations
$$
\pi \, = \, \bigoplus_{i \in I} \pi_i
\hskip.5cm \text{and} \hskip.5cm 
\rho \, = \, \bigoplus_{j \in J} \rho_j
$$
(this is always possible in case $G$ is compact).
The representations $\pi$ and $\rho$ are disjoint
if and only if $\pi_i$ and $\rho_j$ are not equivalent
for all $(i, j) \in I \times J$.
They are quasi-equivalent if and only if
every $\pi_i$ is equivalent to some $\rho_j$
and every $\rho_j$ is equivalent to some $\pi_i$.

\begin{prop}
% 6.A.4
\label{eqetsubordinate}
Let $\pi_1, \pi_2, \sigma_1$ be three representations of $G$.
\begin{enumerate}[label=(\arabic*)]
\item\label{1DEeqetsubordinate}
$\sigma_1$ is subordinate to $\pi_1$ if and only 
there exists a cardinal $n$
such that $\sigma_1$ is equivalent to a subrepresentation of $n \pi_1$.
\item\label{2DEeqetsubordinate}
$\pi_1$ and $\pi_2$ are quasi-equivalent if and only if
there exists a cardinal $n$
such that $n \pi_1$ and $n \pi_2$ are equivalent.
\end{enumerate}
\end{prop}

\begin{proof}
We only write the proof for representations in separable Hilbert spaces,
and we will therefore deal with countable cardinals only.
The non separable case is similar, upon using come cardinal arithmetic.

\vskip.2cm

To show Claim \ref{1DEeqetsubordinate},
assume first that $\sigma_1$ is equivalent to a subrepresentation of $n \pi_1$.
Write 
$$
n \pi_1 \, = \, \bigoplus_i \rho_i ,
$$
where the $\rho_i$~'s are representations equivalent to $\pi_1$.
Let $P_i$ be the orthogonal projection
of the space of $n \pi_1$ onto the space of $\rho_i$.
Let $\tau$ be a subrepresentation of $\sigma_1$.
By hypothesis, there exists $T \ne 0$ in $\Hom_G(\tau, n \pi_1)$.
Therefore there exists $i$ such that $P_iT \ne 0$,
so that $\Hom_G(\tau, \rho_i) \ne \{0\}$, and therefore $\Hom_G(\tau, \pi_1) \ne \{0\}$,
i.e., $\tau$ and $\pi_1$ are not disjoint (Proposition~\ref{Prop-Prop-EquSubRep}).
It follows that $\sigma_1$ is subordinate to $\pi_1$.
\par

Assume now that $\sigma_1$ is subordinate to $\pi_1$.
Denote by $\Hi$ the Hilbert space of~$\sigma_1$.
Consider the family $\mathcal F$ of $G$-invariant closed subspaces of $\Hi$
such that the corresponding subrepresentation of $\sigma_1$
is equivalent to a subrepresentation of $\pi_1$.
By Zorn's lemma, there exists a maximal subfamily $(\Ki_j)_j$
of mutually orthogonal subspaces $\Ki_j \in \mathcal F$.
Set $\Ki \, = \, \bigoplus_j \Ki_j$.
\par

We claim that $\Ki = \Hi$.
Indeed, assume by contradiction that this is not the case.
Then the orthogonal complement $\Ki^\perp$ of $\Ki$ is non-zero and $G$-invariant.
The subrepresentation of $\sigma_1$ defined by $\Ki^\perp$
and the representation $\pi_1$ are disjoint.
This contradicts the fact that $\sigma_1$ is subordinate to $\pi_1$.
\par

From the equality $\Ki = \Hi$,
it follows that $\sigma_1$ is equivalent to a subrepresentation of~$\infty \pi_1$.

\vskip.2cm

To show Claim \ref{2DEeqetsubordinate},
assume first that $\infty \pi_1$ and $\infty \pi_2$ are equivalent.
Then, by Claim \ref{1DEeqetsubordinate},
$\pi_1$ is subordinate to $\pi_2$ and $\pi_2$ is subordinate to $\pi_1$.
This means that $\pi_1$ and $\pi_2$ are quasi-equivalent.
\par

Assume now that $\pi_1$ and $\pi_2$ are quasi-equivalent.
By Claim \ref{1DEeqetsubordinate},
$\pi_1$ is equivalent to a subrepresentation of $\infty \pi_2$,
hence $\infty \pi_1$ is equivalent to a subrepresentation of $\infty \pi_2 \simeq \infty (\infty \pi_2)$;
similarly, $\infty \pi_2$ is equivalent to a subrepresentation of $\infty \pi_1$.
It follows from Proposition~\ref{SchroderBernstein},
that $\infty \pi_1$ and $\infty \pi_2$ are equivalent.
\end{proof}

As a consequence of Proposition~\ref{eqetsubordinate} and Remark~\ref{propweakcont},
quasi-equivalence and weak-equivalence compare as follows.

\begin{cor}
% 6.A.5
\label{QE-WEQ}
Let $\pi_1, \pi_2$ be two representations of~$G$.
\par

If $\pi_1$ is subordinate to $\pi_2$, then $\pi_1$ is weakly contained in $\pi_2$.
In particular, if $\pi_1$ and $\pi_2$ are quasi-equivalent, then they are weakly equivalent.
\end{cor}

\begin{rem}
% 6.A.6
\label{Rem-eqetquasieq}
The converse statement of the last claim of Corollary \ref{QE-WEQ} fails in general,
that is, weakly equivalent representations need not be quasi-equivalent.
Indeed, consider for example the additive group $\R$,
its regular representation $\lambda_\R$,
and the direct sum 
$$
\pi \, = \, \bigoplus_{q \in Q} \chi_q,
$$
where $Q$ is a countable dense subset of $\R$,
and $\chi_q$ denotes the unitary character
$\chi_q \, \colon \R \to \T, \hskip.1cm x \mapsto e^{i qx}$.
The representations $\lambda_\R$ and $\pi$ are weakly equivalent
(more on this in, e.g, \cite[F.2.7 \& F.6.4]{BeHV--08}), but they are not quasi-equivalent.
Indeed, otherwise $\chi_q$ would be equivalent to a subrepresentation of $\lambda_\R$
and this is absurd,
because $\lambda_\R$ does not have any irreducible subrepresentation. 
\end{rem}

\subsection
{Factor representations}
% 6.A.b
\label{SS:FactRep}

We introduce now a class of representations which are called factorial
for a reason which will become clear in the next section
(see Proposition~\ref{Pro-FactorRepVN}).
Another terminology for the same notion is ``primary representation'';
it is used, for example, in \cite{Mack--76}.

\begin{defn}
% 6.A.7
\label{Def-FactorialRep}
A representation $\pi$ of $G$ is a \textbf{factor representation},
or is \textbf{factorial},
if, for every non-trivial decomposition $\pi = \pi_1 \oplus \pi_2$,
the representations $\pi_1$ and $\pi_2$ are not disjoint.
\par

(As a consequence of Proposition \ref{Cor-FactorRepSubRep} below,
$\pi$ is factorial if and only if, for every non-trivial decomposition $\pi = \pi_1 \oplus \pi_2$,
the representations $\pi_1$ and $\pi_2$ are quasi-equivalent.)
\index{Representation! factor}
\index{Factor representation}
\end{defn}

\begin{exe}
% 6.A.8
\label{Exa-FacRepTrivial}
It follows from the definition that irreducible representations are factorial.
\par

More generally,
let $\pi'$ be an irreducible representation of $G$ and let $n \ge 1$ be a cardinal.
Then $\pi := n \pi'$ is factorial.
Indeed, let $\pi = \pi_1 \oplus \pi_2$ be a non trivial direct sum decomposition.
Then, for $i = 1$ and $i = 2$, we have $\Hom_G(\pi_i, \pi')\ne \{0\}$ by Corollary~\ref{Cor-Prop-EquSubRep};
since $\pi'$ is irreducible, it follows from 
Lemma~\ref{Prop-EquSubRep} that $\pi_1$ and $\pi_2$ contain subrepresentations
which are equivalent to $\pi'$.
It follows that $\pi_1$ and $\pi_2$ are not disjoint.
\end{exe}

The following proposition characterizes the factor representations
appearing in Example~\ref{Exa-FacRepTrivial}.

\begin{prop}
% 6.A.9
\label{Prop-QE-IrredRepFactRep}
Let $\pi$ be representation of $G$.
The following properties are equivalent:
\begin{enumerate}[label=(\roman*)]
\item\label{iDEProp-QE-IrredRepFactRep}
$\pi$ is a factor representation which contains an irreducible representation of $G$;
\item\label{iiDEProp-QE-IrredRepFactRep}
$\pi$ is a multiple of an irreducible representation of $G$.
\end{enumerate}
\end{prop}

There is one more equivalent condition in Proposition \ref{Prop-TypeIFac}. 

\begin{proof}
The fact that \ref{iiDEProp-QE-IrredRepFactRep} implies \ref{iDEProp-QE-IrredRepFactRep}
was shown in Example~\ref{Exa-FacRepTrivial}.
To show the converse,
assume that $\pi$ contains an irreducible representation $\pi'$.
\par

Let $\Hi$ be the Hilbert space of $\pi$.
There exists a maximal family $(\Ki_j)_{j \in J} $ of $G$-invariant
and mutually orthogonal closed subspaces $\Ki_i$ of $\Hi$
such that, for every $j \in J$, the subrepresentation $\pi_j$ of $\pi$ defined by $\Ki_j$
is equivalent to $\pi'$
(compare with the proof of Proposition~\ref{eqetsubordinate}~\ref{1DEeqetsubordinate}). 
Set $\Ki = \bigoplus_{j \in J} \Ki_j$.
Then the subrepresentation of $\pi$ defined by $\Ki^\perp$ is disjoint from every $\pi_j$,
hence disjoint from $\bigoplus_{j \in J} \pi_j$, by Corollary~\ref{Cor-Prop-EquSubRep}.
Since $\pi$ is factorial, we must have $\Ki^\perp = \{0\}$, that is, $\Ki = \Hi$.
This shows that $\pi$ is a multiple of $\pi'$.
\end{proof}

For examples of factor representations which are not multiples of irreducible representations, 
see below Section~\ref{S:ICCGroupsNontypeI} and Chapter~\ref{Chap:NormalInfiniteRep}.
\par

The following result shows that factoriality is inherited by subrepresentations.

\begin{prop}
% 6.A.10
\label{Pro-SubrepFactRep}
Let $\pi$ be a factor representation of $G$,
and let $\pi_1, \pi_2$ be non trivial subrepresentations of $\pi$
acting on mutually orthogonal subspaces.
\par

Then $\pi_1$ and $\pi_2$ are not disjoint.
In particular, every non trivial subrepresentation of $\pi$ is factorial.
\end{prop}

\begin{proof}
Let $\Hi$ be the Hilbert space of $\pi$
and $\Hi_1, \Hi_2$ the mutually orthogonal subspaces of $\Hi$ defining $\pi_1$ and $\pi_2$. 
\par

Assume, by contradiction, that 
$\pi_1$ and $\pi_2$ are disjoint. Set 
$$
\Hi' \, = \, \Hi_1 \oplus \Hi_2.
$$ 
There exists a maximal family $(\Ki_j)_{j \in J}$ of $G$-invariant
and mutually orthogonal closed subspaces $\Ki_j$ of $\Hi'^\perp$ 
such that, for every $j \in J$, the subrepresentation $\pi_j$ of $\pi$ defined by $\Ki_j$
is disjoint from $\pi_1$ 
(compare with the proof of
Proposition~\ref{eqetsubordinate}~\ref{1DEeqetsubordinate}).
\par

Let $\Hi_0$ be the orthogonal complement of $\bigoplus_{j \in J} \Ki_j$ in $\Hi'^\perp$
and denote by $\pi_0$ the subrepresentation of $\pi$ defined by $\Hi_0$.
Then we have
$$
\pi \, = \, (\pi_1 \oplus \pi_2) \oplus (\pi_0 \oplus \bigoplus_{j \in J} \pi_j)
\, = \, (\pi_1 \oplus \pi_0) \oplus (\pi_2 \oplus \bigoplus_{j \in J} \pi_j).
$$
On the one hand, $\pi_2 \oplus \bigoplus_{j \in J} \pi_j$ is disjoint from $\pi_1$
by Corollary~\ref{Cor-Prop-EquSubRep}.
On the other hand, no subrepresentation of $\pi_0$ is disjoint from $\pi_1$,
by maximality of the family $(\Ki_j)_{j \in J}$;
it follows that $\pi_2 \oplus \bigoplus_{j \in J} \pi_i$ is disjoint from $\pi_0$.
Therefore $\pi_2 \oplus \bigoplus_{j \in J} \pi_j$ is disjoint from $\pi_1 \oplus \pi_0$.
This contradicts the fact that $\pi$ is factorial.
\end{proof}

The next proposition shows how two factor representations
can be compared.
We will see later (Corollary~\ref{Cor-FactorRepSubRep})
that any two non disjoint factor representations are quasi-equivalent.

\begin{prop}
% 6.A.11
\label{Pro-CompFactRep}
Let $\pi$ and $\pi'$ be factor representations of $G$.
\par
Then one of the following relations holds:
$\pi$ and $\pi'$ are disjoint, 
or $\pi$ is contained in $\pi'$,
or $\pi'$ is contained in $\pi$.
\end{prop}

\begin{proof}
Assume that $\pi$ and $\pi'$ are not disjoint. We have to show that either
$\pi$ is contained in $\pi'$ or $\pi'$ is contained in $\pi$.
\par

Let $\Hi$ and $\Hi'$ be respectively the Hilbert spaces of $\pi$ and $\pi'$. 
There exists a family $(\Ki_i, \Ki'_i)_{i \in I}$ of pairs
of $G$-invariant closed subspaces $\Ki_i$ and $\Ki_i'$ of $\Hi$ and $\Hi'$ respectively
with the following three properties: 
\begin{enumerate}
\item[$\bullet$]
the $\Ki_i$'s are mutually orthogonal and the $\Ki'_i$'s are mutually orthogonal;
\item[$\bullet$]
for every $i \in I$, the subrepresentation $\pi_i$ of $\pi$ defined by $\Ki_i$
is equivalent to the subrepresentation $\pi'_i$ of $\pi'$ defined by $\Ki'_i$;
\item[$\bullet$]
$\big( \bigoplus_{i \in I} \Ki_i, \bigoplus_{i \in I} \Ki'_i \big)$
is a maximal pair of subspaces of $(\Hi, \Hi')$
with the two previous properties.
\end{enumerate}
Let $\Hi_0$ and $\Hi'_0$ be the orthogonal complements
of $\bigoplus_{i \in I} \Ki_i$ and $\bigoplus_{i \in I} \Ki'_i$.
Denote by $\pi_0$ and $\pi'_0$ the subrepresentations
of $\pi$ and $\pi'$ defined by $\Hi_0$ and $\Hi'_0$.
We have 
$$
\pi \, = \, \pi_0 \oplus \bigoplus_{i \in I} \pi_i
\hskip.5cm \text{and} \hskip.5cm
\pi' \, = \, \pi'_0 \oplus \bigoplus_{i \in I} \pi'_i .
$$
Moreover, $\bigoplus_{i \in I} \pi_i$ is equivalent to $\bigoplus_{i \in I} \pi'_i$
and, by maximality of the family $(\Ki_i, \Ki'_i)_{i \in I}$,
the representations $\pi_0$ and $\pi'_0$ are disjoint.
\par

Two cases may occur.
\begin{enumerate}
\item[$\bullet$]
The subrepresentation $\pi_0$ is trivial. Then $\pi$ is contained in $\pi'$.
\item[$\bullet$]
The subrepresentation $\pi_0$ is non trivial. 
We claim that $\pi'_0$ is trivial, that is, $\pi'$ is contained in $\pi$.
\end{enumerate}
Indeed, assume, by contradiction that $\pi'_0$ is non trivial.
Since $\pi$ is factorial, $\pi_0$ and $\bigoplus_{i \in I} \pi_i$ are not disjoint.
Therefore $\pi_0$ and $\bigoplus_{i \in I} \pi'_i$ are not disjoint;
so, $\bigoplus_{i \in I} \pi'_i$ contains a non trivial subrepresentation $\sigma'$
which is equivalent to a subrepresentation of $\pi_0$.
By Proposition~\ref{Pro-SubrepFactRep}, $\sigma'$ and $\pi'_0$ are not disjoint,
since $\pi'_0$ is non trivial.
It follows that $\pi_0$ and $\pi'_0$ are not disjoint and this is a contradiction.
The claim is established.
\end{proof}

\subsection
{Multiplicity-free and type I representations}
% subsection 6.A.c
\label{SS:TypIRep}

In this subsection, we introduce type I representations
and their building blocks, the multiplicity-free representations.

\begin{defn}
% 6.A.12
\label{Def-MultFreeTypeI Rep}
A representation $\pi$ of $G$ is \textbf{multiplicity-free}
% multiplicity-free in \cite{Dixm--C*} 5.4.4 et 5.4.5
if, for every non-trivial decomposition $\pi = \pi_1 \oplus \pi_2$,
the representation $\pi_1$ and $\pi_2$ are disjoint.
\index{Representation! multiplicity-free}
\index{Multiplicity-free representation}
\par

The representation $\pi$ of $G$ is \textbf{type I}, or of \textbf{discrete type}, 
if it is quasi-equivalent to a multiplicity-free representation.
\index{Representation! type I}
\index{Type I! $1$@representation}
\par
 
It follows from the definitions that the type I property is stable by quasi-equival\-ence:
if $\pi$ is a representation of $G$ of type I,
then so is every representation which is quasi-equivalent to $\pi$.
\end{defn}

The next remark indicates a few straightforward examples
of multiplicity-free and type I representations.
For more interesting examples, see Subsection \ref{SS:ExaMultFreeRep}.
% see Theorem~\ref{Thm-CanRepAbMultFree},
% Corollary~\ref{Cor-Thm-CanRepAbMultFree},
% and Theorem~\ref{Thm-MultFreeGel'fandPair}.

\begin{rem}
% 6.A.13
\label{Exa-MutiplicityFreeAndTypeI Trivial}
(1)
Let $(\pi)_{i \in I}$ be a family of mutually non-equivalent irreducible representations of $G$.
Then $\pi \, \colon= \bigoplus_{i \in I} \pi_i$ is a multiplicity-free representation.
\par

Indeed, since the $\pi_i$'s are mutually non-equivalent and irreducible,
we have by Proposition~\ref{Prop-Prop-EquSubRep}, 
Schur's lemma and Lemma \ref{Lem-MatrixCommutant} below
$$
\Hom_G(\pi, \pi)
\, = \, \bigoplus_{i \in I} \Hom_G(\pi_i, \pi_i)
\, = \, \bigoplus_{i \in I} \C \mathrm{Id}_{\Hi_i},
$$
where $\Hi_i$ denotes the Hilbert space of $\pi_i$. In particular, $\Hom_G(\pi, \pi)$
is an abelian algebra and the claim follows from Proposition~\ref{Pro-MultFreeVN}
below.

\vskip.2cm
 
(2)
Let $(n_i)_{i \in I}$ be a family of cardinals
and $(\pi)_{i \in I}$ a family of mutually non-equivalent irreducible representations of $G$.
Then $\pi := \bigoplus_{i \in I} n_i\pi_i$ is a type I representation.
\par

Indeed, $\pi' := \bigoplus_{i \in I} \pi_i$ is multiplicity-free by (1),
and $\pi$ is quasi-equivalent to $\pi'$, by Proposition~\ref{eqetsubordinate}.

\vskip.2cm

(3)
Let $G$ be a compact group. Then every representation of $G$ is of type I.
\par

Indeed, as is well-known, every representation of such a group is a direct sum of 
irreducible representations (see, for instance, \cite[Chap.~5]{Robe--83})
and hence the claim follows from (2).

\vskip.2cm

(4)
A factor representation of $G$ is multiplicity-free if and only if it is irreducible
(Proposition \ref{Pro-SubrepFactRep}).
We will characterize below factor representations which are of type I
(Proposition~\ref{Prop-TypeIFac}).

\vskip.2cm

(5)
Let $\pi$ be a multiplicity-free representation and $\sigma$ a subrepresentation of $\pi$.
Then $\sigma$ is multiplicity-free.
\par

Indeed, let $\sigma = \sigma_1 \oplus \sigma_2$
be a non-trivial decomposition of $\sigma$. Then, we have a non-trivial decomposition
$$
\pi \, = \, \sigma_1 \oplus (\sigma_2 \oplus \sigma'),
$$
where $\sigma'$ is the subrepresentation of $\pi$
on the orthogonal complement of the subspace corresponding to $\sigma$. 
Therefore $\sigma_1$ and $\sigma_2 \oplus \sigma'$ are disjoint.
This implies that $\sigma_1$ and $\sigma_2$ are disjoint.

\vskip.2cm

(6)
Let $\pi_1, \pi_2, \sigma_1$ be three representations of~$G$.
If $\sigma_1$ is quasi-equivalent to a subrepresentation of $\pi_1$
and $\pi_1$ is quasi-equivalent to $\pi_2$,
then there exists a subrepresentation $\sigma_2$ of $\pi_2$
such that $\sigma_1$ and $\sigma_2$ are quasi-equivalent.
\par

To show this, it suffices to consider the situation in which
$\sigma_1$ is a subrepresentation of $\pi_1$.
Let $\Hi_1, \Hi_2$ be the Hilbert spaces of $\pi_1, \pi_2$
and let $\Ki_1$ be the subspace of $\Hi_1$ defining $\sigma_1$.
Let~$\Li$ be the closed subspace of $\Ki_1$
generated by the family of $G$-invariant closed subspaces of $\Ki_1$
such that the corresponding subrepresentation of $\sigma_1$
is equivalent to a subrepresentation of $\pi_2$.
We claim that $\Li = \Ki_1$.
Indeed, assume by contradiction that this is not the case.
Then the orthogonal complement $\Li^\perp$ of $\Li$ in $\Ki_1$ is non-zero and $G$-invariant.
The subrepresentation of $\sigma_1$ defined by $\Li^\perp$
is disjoint from $\pi_2$.
This contradicts the fact that $\pi_1$ is subordinate to $\pi_2$.

Let now $\Ki_2$ be the closed subspace of $\Hi_2$
generated by the family of $G$-invariant closed subspaces of $\Hi_2$
such that the corresponding subrepresentation of $\pi_2$
is equivalent to a subrepresentation of $\sigma_1$.
Let $\sigma_2$ be the subrepresentation of $\pi_2$ defined by $\Ki_2$.
It is clear that $\sigma_1$ is subordinate to $\sigma_2$ and that 
$\sigma_2$ is subordinate to $\sigma_1$.
 
\vskip.2cm

(7)
Let $\pi$ be a type I representation of $G$
and $\sigma$ a subrepresentation of $\pi$.
Then $\sigma$ is type I.
\par

Indeed, $\pi$ is quasi-equivalent to a multiplicity-free representation 
$\pi'$ and so $\sigma$ is quasi-equivalent to a subrepresentation $\sigma'$ of $\pi'$.
Since $\sigma'$ is multiplicity-free by (5), the claim follows from (6).
\end{rem}

We now show that a given quasi-equivalence class
contains at most one multip\-licity-free representation, up to equivalence.

\begin{cor}
% 6.A.14
\label{Eq-QE-MutiplicityFree}
Let $\pi_1, \pi_2$ be two representations of~$G$.
Assume that $\pi_1$ and $\pi_2$ are multiplicity-free and quasi-equivalent.
\par

Then $\pi_1$ and $\pi_2$ are equivalent.
\end{cor}

\begin{proof}
In view of Proposition~\ref{SchroderBernstein}, it suffices 
to show that $\pi_1$ is equivalent to a subrepresentation of $\pi_2$.
\par

As $\pi_1$ is subordinate to $\pi_2$, by the proof of Proposition~\ref{eqetsubordinate},
we have an orthogonal decomposition 
$$
\Hi \, = \, \bigoplus_{i \in I} \Hi_i
$$
of the Hilbert space $\Hi$ of $\pi_1$ into $G$-invariant closed subspaces $\Hi_i$
such that the subrepresentation $\sigma_i$ of $\pi_1$ corresponding to $\Hi_i$
is equivalent to a subrepresentation $\rho_i$ of $\pi_2$;
so, for every $i$, there exists a bijective isometric intertwiner 
$T_i \in \Hom_{G} (\sigma_i, \rho_i)$.
\par

Denote by $P_i$ the orthogonal projection from the Hilbert space $\Ki$ of $\pi_2$
on the subspace $\Ki_i$ of $\rho_i$.
On the one hand, we have for all $i, j$ 
$$
T_j ^{-1}P_jT_i \, \in \, \Hom_{G} (\sigma_{i}, \sigma_{j}).
$$
On the other hand, since $\pi_1$ is multiplicity-free, the $\sigma_i$~'s are mutually disjoint.
If $i \ne j$, it follows from Proposition~\ref{Prop-Prop-EquSubRep} that $T_j ^{-1}P_jT_i = 0$,
hence $P_j = 0$ on the range of $T_i$.
This shows that the subspaces $\Ki_i$ are mutually orthogonal.
\par

Let $T \, \colon \Hi \to \Ki$ be defined by
$$
T\vert_{\Hi_i} \, = \, T_i
\hskip.5cm \text{for every} \hskip.3cm
i .
$$
It is clear that $T$ is an isometric $G$-equivariant embedding,
showing that $\pi_1$ is equivalent to a subrepresentation of $\pi_2$.
\end{proof}

\section[Group representations and von Neumann algebras]
{Group representations and von Neumann algebras}
% Section 6.B
\label{S:QE-Factor-VNAlgebras}

As in the previous section, $G$ denotes a topological group.
\par

The following proposition is a basic tool for the reformulation of the notions
defined in Section \ref{Sectioncomppres} in terms of von Neumann algebras.

\begin{prop}
% 6.B.1
\label{Prop-Commutant}
Let $(\pi, \Hi)$ be a representation of $G$.
\begin{enumerate}[label=(\arabic*)]
\item\label{1DEProp-Commutant}
The commutant $\pi(G)' = \Hom_G(\pi, \pi)$ and the bicommutant $\pi(G)''$ 
are von Neumann subalgebras of $\Li (\Hi)$.
Moreover, $\pi(G)''$ is the von Neumann algebra generated by $\pi(G)$.
\item\label{2DEProp-Commutant}
Let $\Ki$ be a closed linear subspace of $\Hi$,
and $P \in \Li (\Hi)$ the orthogonal projection on $\Ki$.
Then $\Ki$ is $G$-invariant if and only if $P \in \pi(G)'$.
\end{enumerate}
\end{prop}

\begin{proof}
Since the representation $\pi$ is unitary, $\pi(G)^* = \pi(G)$,
and Claim \ref{1DEProp-Commutant} follows from
the von Neumann bicommutant theorem;
see Appendix \ref{AppAlgvN}.

\vskip.2cm

For \ref{2DEProp-Commutant},
it is clear that $\Ki$ is $G$-invariant if $P \in \pi(G)'$. Conversely, 
assume that $\Ki$ is $G$-invariant. Then $\pi(g) P=P\pi(g)P$ for all $g \in G$.
Therefore 
$$
P\pi(g) \, = \, (\pi(g^{-1}) P)^* \, = \, (P\pi(g^{-1})P)^* \, = \, P\pi(g)P
$$
and hence $\pi(g) P = P\pi(g)$ for all $g \in G$.
\end{proof}

Let $\Hi$ be a Hilbert space which is a Hilbert direct sum
of a family $(\Hi)_{i \in I}$ of closed subspaces indexed by a set $I$.
For each $i \in I$, denote by $P_i \, \colon \Hi \twoheadrightarrow \Hi_i$
the orthogonal projection.
An operator $T \in \Li(\Hi)$ can be described by a matrix $(T_{i, j})_{i, j \in I}$,
where $T_{i, j} \in \Li( \Hi_j, \Hi_i)$ is defined by
$$
T_{i, j} (\xi) \,:= \, P_i T(\xi)
\hskip.5cm \text{for all} \hskip.2cm
i, j \in I
\hskip.2cm \text{and} \hskip.2cm
\xi \in \Hi_j,
$$
and
$$
T \Big( (\xi_i)_{i \in I} \Big) \, = \, \Big( \sum_{j \in I} P_i T \xi_j \Big)_{i \in I}
\hskip.5cm \text{for all} \hskip.2cm
(\xi_i)_{i \in I} \in \Hi = \bigoplus_{i \in I} \Hi_i .
$$

\begin{lem}
% 6.B.2
\label{Lem-MatrixCommutant}
Let $(\pi_i, \Hi_i)_{i \in I}$ be a family of representations of $G$,
and let $\pi = \bigoplus_{i \in I} \pi_i$ be their direct sum,
in the Hilbert space $\Hi$ which is the direct sum $\bigoplus_{i \in I} \Hi_i$.
Let $T \in \Li (\Hi)$ and let $(T_{i, j})_{i, j \in I}$ be the matrix of $T$ 
associated to the decomposition $\Hi = \bigoplus_{i \in I} \Hi_i$.
The following properties are equivalent:
\begin{enumerate}[label=(\roman*)]
\item\label{iDELem-MatrixCommutant}
$T \in \pi(G)' = \Hom_G(\pi, \pi)$;
\item\label{iiDELem-MatrixCommutant}
$T_{i, j} \in \Hom_G(\pi_j, \pi_i)$ for every $i, j \in I$.
\end{enumerate}
In particular, when the representation $\pi_i$ are pairwise disjoint,
Property \ref{iDELem-MatrixCommutant} is equivalent to
\begin{enumerate}[label=(\roman*)]
\addtocounter{enumi}{2}
\item\label{iiiDELem-MatrixCommutant}
$T_i \in \pi_i(G)'$ and $T_{i, j} = 0$ for all $i, j \in I$ such that $i \ne j$.
\end{enumerate}
\end{lem}

Let $\Hi$ be a Hilbert space,
$\Ki$ a closed subspace,
and $P$ the orthogonal projection of $\Hi$ onto $\Ki$.
Let $\mathcal{M}$ be a von Neumann subalgebra of $\Li (\Hi)$.
In case $P \in \mathcal M'$ or $P \in \mathcal M$, we denote by $\mathcal M_P$ 
the von Neumann subalgebra $\{P T \vert_\Ki \mid T \in \mathcal M\}$ of $\Li (\Ki)$. 
Observe that $\mathcal M_P$ is isomorphic to the subalgebra $P\mathcal M P$ of $\Li (\Hi)$.
\par

When $P \in \mathcal M'$, the subspaces $\Ki$ and $\Ki^\perp$
are invariant under $\mathcal M$,
and $\mathcal M_P$ is the induced von Neumann algebra as in Appendix \ref{AppAlgvN}.
\par

We record the following useful consequence of Lemma~\ref{Lem-MatrixCommutant},
for $I$ consisting of two elements.

\begin{prop}
% 6.B.3
\label{Prop-MatrixCommutant}
Let $(\pi, \Hi)$ be a representation of $G$ and set $\mathcal M := \pi(G)''$.
Let $P \in \mathcal M'$ be a projection. Then 
$$
(\mathcal M_P)' \, = \, (\mathcal M')_P.
$$
\end{prop}

\begin{proof}
For $T \in (\mathcal M_P)'$,
let $\widetilde{T}$ be the bounded operator on $\Hi$
with matrix $(\widetilde{T}_{i, j})_{1 \le i, j \le 2}$
associated to the decomposition $\Hi = \Ki \oplus \Ki^\perp$
and given by $\widetilde{T}_{1, 1} = T$ and $T_{i, j} = 0$ otherwise.
Lemma~\ref{Lem-MatrixCommutant} shows that $\widetilde{T} \in \mathcal M'$.
Therefore $T \in (\mathcal M')_P$.
\par

Conversely, let $T \in (\mathcal M')_P$.
Then there exists $\widetilde{T} \in \mathcal M '$
with matrix $(\widetilde{T}_{i, j})_{1 \le i, j \le 2}$
associated to the decomposition $\Hi = \Ki \oplus \Ki^\perp$
such that $\widetilde{T}_{1, 1} = P \widetilde T P = T$.
For all $S \in \mathcal M$, we have
$T PSP = P \widetilde T PSP = P \widetilde T SP =
P S \widetilde T PP = PSP \widetilde T P = PSP T$.
It follows that
% We have $P \widetilde{T} P \in \mathcal M'$ and $T = P \widetilde{T} P \vert_\Ki$.
% Therefore 
$T \in (\mathcal M_P)'$.
\end{proof}

\subsection
{Disjointness and quasi-equivalence, again}
% subsection 6.B.a
\label{SS:QE+}
\index{Representation! quasi-equivalent $\approx$}
\index{Quasi-equivalent! representations}
\index{$a4$@$\approx$ quasi-equivalence of representations}

Let $(\pi, \Hi)$ be a representation of $G$
and $\Ki_1$ a $G$-invariant closed subspace of $\Hi$.
Denote by $\pi_1$ and $\pi_2$ the corresponding subrepresentations
of $\pi$ defined by $\Ki_1$ and $\Ki_2 := \Ki_1^\perp$.
Then $\pi(G)$ is contained in $\pi_1(G)\oplus \pi_2(G) $,
which is itself contained in $\Li (\Ki_1)\oplus \Li (\Ki_2)\subset \Li (\Hi)$.
This implies the inclusions
$$
\pi_1(G)' \oplus \pi_2(G)' \subset \pi(G)'
\leqno{(a.1)} 
$$
and
$$
\pi(G)'' \subset \pi_1(G)'' \oplus \pi_2(G)''.
\leqno{(a.2)} 
$$
We now show that equality in the inclusions (a.1) and (a.2)
corresponds to disjointness of $\pi_1$ and $\pi_2$.

\begin{prop}
% 6.B.4
\label{Pro-DisjointRepVN}
Let $(\pi, \Hi)$ be a representation of $G$
which is the direct sum of two subrepresentations $\pi_1, \pi_2$
on orthogonal subspaces $\Ki_1, \Ki_2$ of $\Hi$, such that $\Hi = \Ki_1 \oplus \Ki_2$.
Let $P_1$ be the orthogonal projection from $\Hi$ onto $\Ki_1$.
The following properties are equivalent:
\begin{enumerate}[label=(\roman*)]
\item\label{iDEPro-DisjointRepVN}
$\pi_1$ and $\pi_2$ are disjoint;
\item\label{iiDEPro-DisjointRepVN}
the von Neumann algebras $\pi(G)'$ and $\pi_1(G)' \oplus \pi_2(G)'$ coincide;
\item\label{iiiDEPro-DisjointRepVN}
the von Neumann algebras $\pi(G)''$ and $\pi_1(G)'' \oplus \pi_2(G)''$ coincide;
\item\label{ivDEPro-DisjointRepVN}
$P_1$ belongs to the centre of $\pi(G)'$, that is, belongs to $\pi(G)' \cap \pi(G)''$.
\end{enumerate}
\end{prop}

\begin{proof}
Assume that $\pi_1$ and $\pi_2$ are disjoint. 
Let $T \in \pi(G)'$.
Then, by Lemma~\ref{Lem-MatrixCommutant} and Proposition \ref{Prop-Prop-EquSubRep},
the matrix $(T_{i, j})_{1 \le i, j \le 2}$ of $T$
associated to the decomposition $\Hi = \Ki_1 \oplus \Ki_2$ is diagonal;
and $T_{i, i} \in \pi_i(G)'$ for $i = 1, 2$.
Therefore $\pi(G)' \subset \pi_1(G)' \oplus \pi_2(G)'$.
By (a.1), $\pi(G)' = \pi_1(G)' \oplus \pi_2(G)'$.
This shows that \ref{iDEPro-DisjointRepVN} implies \ref{iiDEPro-DisjointRepVN}.

\vskip.2cm

Assume that $\pi_1(G)' \oplus \pi_2(G)' = \pi(G)'$.
Then $\pi_1(G)'' \oplus \pi_2(G)'' \, \subset \, \pi(G)''$,
hence, $\pi(G)'' = \pi_1(G)'' \oplus \pi_2(G)''$ by (a.2).
This shows that \ref{iiDEPro-DisjointRepVN} implies \ref{iiiDEPro-DisjointRepVN}. 

\vskip.2cm

Assume that $\pi_1(G)'' \oplus \pi_2(G)'' = \pi(G)''$.
For the matrix $(T_{i, j})_{1 \le i, j \le 2}$ of $P_1$
associated to the decomposition $\Hi = \Ki_1 \oplus \Ki_2$,
we have $T_{1, 1} = \mathrm{Id}_{\Ki_1} \in \pi_1(G)''$ and $T_{i, j} = 0$ otherwise.
Therefore $P_1 \in \pi(G)''$.
Since $P_1 \in \pi(G)'$, it follows that $P_1$ belongs to the centre of $\pi(G)'$. 
This shows that \ref{iiiDEPro-DisjointRepVN} implies in \ref{ivDEPro-DisjointRepVN}.

\vskip.2cm

Assume that $\pi_1$ and $\pi_2$ are not disjoint. 
Then, by Proposition \ref{Prop-Prop-EquSubRep},
there exists a non-zero element $S$ in $\Hom_G(\pi_1, \pi_2)$.
Let $T \in \Li (\Hi)$ be the operator with matrix $(T_{i, j})_{1 \le i, j \le 2}$
associated to the decomposition $\Hi = \Ki_1 \oplus \Ki_2$
defined by $T_{2, 1} = S$ and $T_{i, j} = 0$ otherwise.
Then $T \in \pi(G)'$, by Lemma~\ref{Lem-MatrixCommutant}.
It is clear that $TP_1 \ne P_1T$.
So, $P_1$ does not belong to the centre of $\pi(G)'$.
It follows that \ref{ivDEPro-DisjointRepVN} implies \ref{iDEPro-DisjointRepVN}.
\end{proof}

The following proposition gives a characterization
of quasi-equivalence of representations in terms of von Neumann algebras.
For a reminder on induced von Neumann algebras such as $\pi_2(G)''_E$ below,
see Appendix \ref{AppAlgvN}.

\begin{prop}
% 6.B.5
\label{Pro-SubordinateVN}
Let $\pi_1$ and $\pi_2$ be representations of $G$.
The following properties are equivalent:
\begin{enumerate}[label=(\roman*)]
\item\label{iDEPro-SubordinateVN}
$\pi_1$ is subordinate to $\pi_2$;
\item\label{iiDEPro-SubordinateVN}
there exist a central projection $E \in \pi_2(G)''$ and an isomorphism 
$$
\Phi \, \colon \, \pi_2(G)''_{E} \to \pi_1(G)''
$$
such that $ \Phi(E\pi_2(g)E) = \pi_1(g)$ for all $g \in G$.
\end{enumerate}
\par

In particular, $\pi_1$ and $\pi_2$ are quasi-equivalent
if and only if there exists an isomorphism $\Phi \, \colon \pi_2(G)'' \to \pi_1(G)''$
such that $ \Phi(\pi_2(g)) = \pi_1(g) $ for all $g \in G$.
\end{prop}

\begin{proof}
\ref{iDEPro-SubordinateVN} $\Rightarrow$ \ref{iiDEPro-SubordinateVN}
Assume that $\pi_1$ is subordinate to $\pi_2$.
Then $\pi_1$ is equivalent to a subrepresentation of $n \pi_2$ for some cardinal $n$,
by Proposition~\ref{eqetsubordinate}.
Denote by $\Ki$ a Hilbert space of dimension $n$.
\par

Let $\Hi_1$ and $\Hi_2$ be the Hilbert spaces of $\pi_1$ and $\pi_2$.
We realize $\pi := n \pi_2$ on the Hilbert space $\Hi := \Hi_2 \otimes \Ki$ by 
$$
\pi(g) \, = \, \pi_2(g) \otimes \mathrm{Id}_{\Ki}
\hskip.5cm \text{for all} \hskip.2cm
g \in G .
$$
Observe that
$$
\pi_2(G)'' \to \pi(G)'', \hskip.2cm T \mapsto T \otimes \mathrm{Id}_{\Ki},
$$
is an isomorphism of von Neumann algebras
mapping $\pi_2(g)$ to $\pi(g)$ for all $g \in G$.
Therefore, it suffices to show the following claim:
there exists a central projection $E \in \pi(G)''$
and an isomorphism $\Phi \, \colon \pi(G)''_E \to \pi_1(G)''$
such that $ \Phi(E \pi(g) E) = \pi_1(g)$ for all $g \in G$.
\par

We view $\pi_1$ as a subrepresentation of $\pi$.
Let $P \in \pi(G)'$ be the orthogonal projection onto the subspace defining $\pi_1$
and let $E$ be the central support of $P$
(for the notion of central support, see Definition \ref{defcentralsupportoperator}).
Since $\pi(G)''_{P} = \pi_1(G)''$ and since the induction map 
$\pi(G)''_{E} \to \pi(G)''_{P}$ is an isomorphism
(Proposition \ref{inductionisiso}
applied to the von Neumann algebra $\mathcal M = \pi(G)''_{E}$)),
the claim follows.

\vskip.2cm

\ref{iiDEPro-SubordinateVN} $\Rightarrow$ \ref{iDEPro-SubordinateVN}
Assume now that there exist a central projection $E \in \pi_2(G)''$
and an isomorphism $\Phi \, \colon E \pi_2(G)''E \to \pi_1(G)''$
such that $\Phi(E\pi_2(g)E) = \pi_1(g)$ for all $g \in G$.
Upon replacing $\pi_2$ by its subrepresentation $g \mapsto E \pi_2(g)$
defined on the range of $E$,
we may assume that $E$ is the identity operator on $\Hi_2$.
\par

By the structure of isomorphisms of von Neumann algebras
(Proposition \ref{isovNalgebra}), there exist
\par

$\bullet$
a Hilbert space $\Ki$,
\par

$\bullet$
a projection $P \in (\pi_2(G) \otimes \mathrm{Id}_{\Ki})' \subset \Li (\Hi_2 \otimes \Ki)$,
with range $\Hi_2' \subset \Hi_2 \otimes \Ki$,
\par

$\bullet$
and an isomorphism of Hilbert spaces $U \, \colon \Hi_2' \to \Hi_1$,
\par

\noindent
such that 
$$
\Phi(T) \, = \, U \circ (T \otimes \mathrm{Id}_{\Ki})_P \circ U^{-1}
\hskip.5cm \text{for every} \hskip.2cm
T \in \pi_2(G)''.
$$
(Observe that $\Hi_2'$ is invariant under $\pi_2(G) \otimes \mathrm{Id}_{\Ki}$,
since $P \in (\pi_2(G) \otimes \mathrm{Id}_{\Ki})'$.)
Therefore, we have 
$$
\pi_1(g) \, = \, \Phi(\pi_2(g)) \, = \, U \circ (\pi_2(g) \otimes \mathrm{Id}_{\Ki}) \circ U^{-1}
\hskip.5cm \text{for every} \hskip.2cm
g \in G .
$$
This shows that $\pi_1$ is equivalent to a subrepresentation of $n \pi_2$,
where $n$ is the cardinal of $\Ki$.
Therefore $\pi_1$ is subordinate to $\pi_2$, by Proposition~\ref{eqetsubordinate}.

\vskip.2cm

It remains to prove the last statement in Proposition~\ref{Pro-SubordinateVN}.
Assume that there exists an isomorphism $\Phi \, \colon \pi_2(G)'' \to \pi_1(G)''$
such that $ \Phi(\pi_2(g)) = \pi_1(g)$ for all $g \in G$;
then, by what we have just seen, $\pi_1$ and $\pi_2$ are quasi-equivalent.
\par

Conversely, assume that $\pi_1$ and $\pi_2$ are quasi-equivalent.
Then there exist a central projection $E \in \pi_2(G)''$
and an isomorphism $\Phi \, \colon \pi_2(G)''_{E} \to \pi_1(G)''$
such that $ \Phi(E \pi_2(g) E) = \pi_1(g)$ for all $g \in G$.
We claim that $E = I$.
Indeed, let $\pi', \pi''$ be the subrepresentations of $\pi_2$ defined by $E$ and by $I - E$. 
On the one hand, $\pi'$ is quasi-equivalent to $\pi_1$, by what we have proved above.
On the other hand, $\pi'$ and $\pi''$ are disjoint, by Proposition~\ref{Pro-DisjointRepVN}.
Therefore $\pi_1$ and $\pi'' $ are disjoint.
Since $\pi_2=  \pi' \oplus \pi''$ is quasi-equivalent to $\pi_1$, it follows that $\pi''$ is trivial,
that is, $E = I$.
\end{proof}

\subsection
{Factor representations, again}
% subsection 6.B.b
\label{SS:FactRep+}

\index{Representation! factor}
\index{Factor representation}
Factor representations can be characterized in terms of von Neumann algebras:

\begin{prop}
% 6.B.6
\label{Pro-FactorRepVN}
For a representation $\pi$ of~$G$,
the following properties are equivalent:
\begin{enumerate}[label=(\roman*)]
\item\label{iDEPro-FactorRepVN}
$\pi$ is a factor representation;
\item\label{iiDEPro-FactorRepVN}
$\pi(G)''$ is a factor.
\end{enumerate}
\end{prop}

\begin{proof}
Let $\Hi$ be the Hilbert space of $\pi$.
Assume that $\pi$ is a factor representation
and let $P$ a non-zero projection in the centre $\mathcal Z$ of $\pi(G)''$.
By Proposition \ref{Prop-Commutant}, $\pi$ is the direct sum
of the subrepresentations $\pi_1$ and $\pi_2$
defined on the ranges of $P$ and $\mathrm{Id}_{\Hi} - P$.
By Proposition~\ref{Pro-DisjointRepVN}, the subrepresentations $\pi_1$ and $\pi_2$ are disjoint,
hence $\pi_2$ is trivial, that is, $P = \mathrm{Id}_{\Hi}$.
Therefore $0$ and $\mathrm{Id}_{\Hi}$ are the only projections in $\mathcal Z$.
It follows from Proposition~\ref{Prop-ProjvN} that $\pi(G)''$ is a factor.

\vskip.2cm

Conversely, assume that $\pi(G)''$ is a factor
and let $\pi = \pi_1 \oplus \pi_2$ be a direct sum decomposition
in two non trivial subrepresentations.
Since the centre of $\pi(G)'$ consists only of the scalar multiples of $\mathrm{Id}_{\Hi}$,
the orthogonal projection onto the subspace defining $\pi_1$
does not belong to the centre of $\pi(G)'$.
It follows from Proposition~\ref{Pro-DisjointRepVN} that $\pi_1$ and $\pi_2$ are not disjoint.
Therefore $\pi$ is a factor representation.
\end{proof}

As a first consequence of Proposition~\ref{Pro-FactorRepVN},
we show that factoriality of a representation is a property of its quasi-equivalence class.

\begin{cor}
% 6.B.7
\label{Cor-FactorRepQE}
Let $\pi$ be a factor representation of $G$.
\par

Then every representation which is quasi-equivalent to $\pi$ is factorial.
In particular, every multiple $n \pi$ of $\pi$ is factorial.
\end{cor}

\begin{proof}
Let $\sigma$ be a representation which is quasi-equivalent to $\pi$.
By Proposition~\ref{Pro-SubordinateVN}, 
the von Neumann algebras $\sigma(G)''$ and $\pi(G)''$ are isomorphic.
The claim follows now from Proposition~\ref{Pro-FactorRepVN}.
\end{proof}

A second consequence of Proposition~\ref{Pro-FactorRepVN}
is that subrepresentations of a factor representation belong to the same quasi-equivalence class.

\begin{cor}
% 6.B.8
\label{Cor-FactorRepSubRep}
Let $\pi$ be a factor representation of $G$.
\par

Then every non trivial subrepresentation of $\pi$ is quasi-equivalent to $\pi$.
\end{cor}

\begin{proof}
Let $\Hi$ be the Hilbert space of $\pi$
and let $\sigma$ be a non trivial subrepresentation of $\pi$.
By Proposition~\ref{Pro-FactorRepVN},
the only non-zero projection in the centre of $\pi(G)''$ is $\mathrm{Id}_{\Hi}$.
Since $\sigma$ is subordinate to $\pi$,
it follows from Proposition~\ref{Pro-SubordinateVN}
that there exists an isomorphism $\Phi \, \colon \pi(G)'' \to \sigma(G)''$
such that $ \Phi(\pi(g)) = \sigma(g)$ for all $g \in G$,
hence that $\sigma$ is quasi-equivalent to $\pi$.
\end{proof}

The following consequence of
Proposition~\ref{Pro-FactorRepVN} and Corollary~\ref{Cor-FactorRepSubRep}
should be compared with the analogous statement for irreducible representations
in Proposition~\ref{Pro-IrrRepBoundedCard}.

\begin{cor}
% 6.B.9
\label{FacSeparable}
Let $X$ a dense subset of $G$ of cardinal~$\aleph$.
\par

Every factor representation of $G$ is quasi-equivalent to a factor representation
in a Hilbert space of dimension at most $\aleph$.
\par

In particular, if $G$ is separable, every factor representation 
of $G$ is quasi-equivalent to a factor representation in a separable Hilbert space.
\end{cor}

\begin{proof}
Let $(\pi, \Hi)$ be a factor representation of $G$.
Let $\xi$ be a non-zero vector in $\Hi$ and let $\Ki$ be the closed
subspace of $\Hi$ generated by $\pi(G)\xi$.
Since $\pi(X)\xi$ is dense in $\pi(G)\xi$ and since $X$ has cardinal $\aleph$,
we have $\dim \Ki \le \aleph$.
Moreover, $\Ki$ is $G$-invariant and 
so defines a subrepresentation $\sigma$ of $\pi$.
By
% Proposition~\ref{Pro-FactorRepVN} and
Corollary~\ref{Cor-FactorRepSubRep},
$\sigma$ is a factor representation and is quasi-equivalent to $\pi$.
\end{proof}

\subsection
{Multiplicity-free and type I representations, again}
% subsection 6.B.c
\label{SS:TypIRep+}
\index{Representation! type I}
\index{Type I! $1$@representation}

Multiplicity-free representations can also be characterized in terms of von Neumann algebras.

\index{Representation! multiplicity-free}
\index{Multiplicity-free representation}
\begin{prop}
% 6.B.10
\label{Pro-MultFreeVN}
For a representation $\pi$ of $G$,
the following properties are equivalent:
\begin{enumerate}[label=(\roman*)]
\item\label{iDEPro-MultFreeVN}
$\pi$ is multiplicity-free;
\item\label{iiDEPro-MultFreeVN}
$\pi(G)'$ is an abelian von Neumann algebra.
\end{enumerate}
\end{prop}

\begin{proof}
Assume that $\pi$ is multiplicity-free. Then, by Proposition~\ref{Pro-DisjointRepVN},
every projection contained in $\pi(G)'$ belongs to the centre of $\pi(G)'$, 
In view of Proposition~\ref{Prop-ProjvN}, it follows that $\pi(G)'$ is abelian.

\vskip.2cm
 
Conversely, assume that $\pi(G)'$ is abelian.
Let $\pi = \pi_1 \oplus \pi_2$ be a decomposition of $\pi$
into subrepresentations $\pi_1$ and $\pi_2$.
Since the orthogonal projection on the subspace defining $\pi_1$ belongs to $\pi(G)'$,
it follows from Proposition~\ref{Pro-DisjointRepVN} that $\pi_1$ and $\pi_2$ are disjoint.
\end{proof}

The direct sum of multiplicity-free representations is not multiplicity-free in general;
for instance, $\pi \oplus \pi$ is not multiplicity free for any representation $\pi$.
However, the following corollary of Proposition~\ref{Pro-MultFreeVN} holds.

\begin{cor}
% 6.B.11
\label{Cor-DirectSumMultFreeRep}
Let $ (\pi_i)_{i \in I} $ a family of \emph{pairwise disjoint} multiplicity-free representations of $G$.
\par

Then $ \bigoplus_{i \in I} \pi_i$ is multiplicity-free.
\end{cor}

\begin{proof}
Set $\pi := \bigoplus_{i \in I} \pi_i$. 
Since the $\pi_i$'s are pairwise disjoint, we have
$$
\pi(G)' \, = \, \bigoplus_{i \in I} \pi_i(G)'
$$
(see the proof of Proposition~\ref{Pro-DisjointRepVN}).
Therefore $\pi(G)'$ is abelian, and $\pi$ is multiplicity-free
by Proposition~\ref{Pro-MultFreeVN}.
\end{proof}

The following characterization of type I representations in terms of von Neumann algebras
is a direct consequence of
Proposition~\ref{Pro-SubordinateVN} and Proposition~\ref{Pro-MultFreeVN}.

\begin{prop}
% 6.B.12
\label{Pro-TypeIVN}
For a representation $\pi$ of $G$,
the following properties are equivalent:
\begin{enumerate}[label=(\roman*)]
\item\label{iDEPro-TypeIVN}
$\pi$ is of type I;
\item\label{iiDEPro-TypeIVN}
there exist a representation $\sigma$ of $G$ with abelian commutant $\sigma(G)'$ 
and an isomorphism $\Phi \, \colon \pi(G)'' \to \sigma(G)''$
such that $ \Phi(\pi(g)) = \sigma(g) $ for all $g \in G$.
\end{enumerate}
\end{prop}

Next, we show that the property of being type I is inherited by direct sums of representations.

\begin{prop}
% 6.B.13
\label{Prop-DirectSumTypeIRep}
(1)
Let $(\pi, \Hi)$ be a representation of $G$
and $P$ a projection in $\pi(G)'$ with central support $E$.
Assume that the subrepresentation of $\pi$ defined on the range of $P$ is of type I.
\par

Then the subrepresentation of $\pi$ defined on the range of $E$ is of type I.

\vskip.2cm

(2)
Let $ (\pi_i)_{i \in I} $ be a family of type I representations of $G$.
\par

Then $ \bigoplus_{i \in I} \pi_i$ is a type I representation.
\end{prop}

\begin{proof}
(1)
By Proposition~\ref{Pro-TypeIVN},
there exist a representation $\sigma$ of $G$ with abelian commutant $\sigma(G)'$ 
and an isomorphism $\Phi \, \colon \pi(G)''_P \to \sigma(G)''$
such that $\Phi(P\pi(g)P) = \sigma(g)$ for all $g \in G$.
The induction map $\Psi \, \colon \pi(G)''_E \to \pi(G)''_{P}$ is an isomorphism.
Therefore
$$
\Phi \circ \Psi \, \colon \, \pi(G)''_E \to \sigma(G)''
$$
is an isomorphism and we have 
$$
\Phi \circ \Psi(E\pi(g)E) \, = \, \Phi(P\pi(g) P) \, = \, \sigma(g)
\hskip.5cm \text{for all} \hskip.2cm
g \in G.
$$
Therefore, again by Proposition~\ref{Pro-TypeIVN}, the subrepresentation of $\pi$ defined on 
the range of $E$ is of type I.

\vskip.2cm

(2)
Set 
$$
\Hi \, := \, \bigoplus \Hi_i
\hskip.5cm \text{and} \hskip.5cm 
\pi \, := \, \bigoplus_{i \in I} \pi_i,
$$
where $\Hi_i$ is the Hilbert space of $\pi_i$. For $i \in I$, denote by 
$P_i \in \pi(G)'$ the orthogonal projection from $\Hi$ onto $\Hi_i$
and by $E_i$ the central support of $P_i$.
\par

By Zorn's lemma, we can find a maximal family $(F_j)_{j \in J}$
of \emph{mutually orthogonal} projections $F_j$
in the centre $\mathcal Z$ of $\pi(G)''$
with the property that the subrepresentation of $\pi$ defined on the range of $F_j$ is of type I.
Set 
$$
F \, := \, \bigoplus_{j \in J} F_j .
$$
Then $F$ is a projection in $\mathcal Z$.
We claim that $F = \mathrm{Id}_{\Hi}$. 
\par
 
Indeed, assume, by contradiction, that $\mathrm{Id}_{\Hi} - F \ne 0$.
Since $\bigoplus_{i \in I} P_i = \mathrm{Id}_{\Hi}$, there exists $i \in I$ such that 
$(\mathrm{Id}_{\Hi} - F) P_i \ne 0$ and therefore $(\mathrm{Id}_{\Hi} - F)E_i \ne 0$.
The range $\Ki'$ of $(\mathrm{Id}_{\Hi} - F)E_i$ is invariant under $\pi(G)$ and $\pi(G)'$,
since $F$ and $E_i$ belong to $\mathcal Z$.
Therefore the orthogonal projection $F'$ on $\Ki'$ is a non-zero projection in $\mathcal Z$. 
Let $\pi'$ be the subrepresentation of $\pi$ defined on the range of $F'$.
\par
 
On the one hand, $\Ki'$ is contained in the range of $E_i$;
so, $\pi'$ is contained in the subrepresentation $\pi'_i$ of $\pi$ defined on the range of $E_i$.
Since $\pi_i$ is of type I, it follows from (1) that $\pi'_i$ is of type I.
Therefore $\pi'$ is of type I (see Remark~\ref{Exa-MutiplicityFreeAndTypeI Trivial}.6).
On the other hand, $\Ki'$ is contained in the range of $\mathrm{Id}_{\Hi} - F$
and so $F'$ is orthogonal to $F_j$ for every $j \in J$.
This is a contradiction with the maximality of the family $(F_j)_{j \in J}$. 
Therefore, we have $\bigoplus_{j \in J} F_j = \mathrm{Id}_{\Hi}$, as claimed.
\par
 
For every $j \in J$, denote by $\widetilde \pi_j$ the subrepresentation of $\pi$
defined on the range of $F_j$.
Then $\pi = \bigoplus_{j \in J} \widetilde \pi_j$. Moreover, since the $F_j$'s are mutually orthogonal 
projections in the centre of $\pi(G)''$, the representations $\widetilde \pi_j$'s are pairwise disjoint
and we have 
$$
\pi(G)'' \, = \, \bigoplus_{j \in J} \widetilde \pi_j(G)''
$$
(compare with the proof of Proposition~\ref{Pro-DisjointRepVN}).

For every $j \in J$, 
there exist a representation $\sigma_j$ of $G$ with abelian commutant $\sigma_j(G)'$ 
and an isomorphism $\Phi_j \, \colon \widetilde \pi_j(G)'' \to \sigma_j(G)''$
such that $\Phi_j(\widetilde \pi_j(g)) = \sigma_j(g) $ for all $g \in G$ (Proposition~\ref{Pro-TypeIVN}).
Then
$$
\Phi \, := \, \bigoplus_{j \in J} \Phi_j \, \colon \,
\pi(G)'' \to \bigoplus_{j \in J}\sigma_j(G)''
$$
is an isomorphism of von Neumann algebras
such that $ \Phi_j(\pi(g)) = \bigoplus_{j \in J}\sigma_j(g)$ for all $g \in G$.
Therefore (Proposition~\ref{Pro-SubordinateVN}), $\pi$ is quasi-equivalent
to $\sigma := \bigoplus_{j \in J}\sigma_j$. 
The $\sigma_j$~'s are pairwise disjoint, since $\sigma_j$ is quasi-equivalent 
to $\widetilde \pi_j$ and since the $\widetilde \pi_j$'s are pairwise disjoint.
Therefore $\sigma$ is multiplicity-free, by Corollary~\ref{Cor-DirectSumMultFreeRep}.
Therefore, $\pi$ is of type I.
\end{proof}

We now show that the representations
appearing in Proposition~\ref{Prop-QE-IrredRepFactRep}
are exactly the factor representations of type $I$.

\begin{prop}
% 6.B.14
\label{Prop-TypeIFac}
Let $\pi$ be a factor representation of $G$. The following properties are equivalent:
\begin{enumerate}[label=(\roman*)]
\item\label{iDEProp-TypeIFac}
$\pi$ is of type I;
\item\label{iiDEProp-TypeIFac}
$\pi$ is a multiple of an irreducible representation of $G$.
\end{enumerate}
When $\pi$ has these properties, there exists an irreducible representation $\sigma$,
unique up to equivalence, such that
every representation of $G$ which is quasi-equivalent to $\pi$
is equivalent to a multiple of $\sigma$.
\end{prop}

\begin{proof}
That \ref{iiDEProp-TypeIFac} implies \ref{iDEProp-TypeIFac}
is a particular case of Remark \ref{Exa-MutiplicityFreeAndTypeI Trivial} (2).
\par

For the converse, assume that $\pi$ is type I, i.e.,
that $\pi$ is quasi-equivalent to a multiplicity-free representation $\sigma$ of $G$.
By Corollary~\ref{Cor-FactorRepQE}, $\sigma$ is a factor representation.
Proposition~\ref {Pro-MultFreeVN} implies that $\sigma(G)'$ is one dimensional;
this means that $\sigma$ is an irreducible representation.
Therefore, $\sigma$ is equivalent to a subrepresentation of $\pi$.
by Proposition~\ref{Prop-EqQuasiEq-Immediate}.
Therefore $\pi$ is a multiple of $\sigma$,
by Proposition~\ref{Prop-QE-IrredRepFactRep}.
\par
 
To show the last claim, let $\pi'$ be a representation of $G$
which is quasi-equivalent to $\pi$.
Then $\pi'$ is a factor representation of type I
and the equivalence of \ref{iDEProp-TypeIFac} and \ref{iiDEProp-TypeIFac},
applied to $\pi'$,
shows that $\pi'$ is a multiple of an irreducible representation $\sigma'$.
Since $\sigma$ and $\sigma'$ are quasi-equivalent and irreducible, 
$\sigma$ and $\sigma'$ are equivalent.
\end{proof}

It is worth summarizing the results
we obtained so far about the quasi-equivalence class of a factor representation
Proposition~\ref{Pro-CompFactRep}, Corollary~\ref {Cor-FactorRepQE},
Corollary~\ref{Cor-FactorRepSubRep}, and Proposition~\ref{Prop-TypeIFac}).

\begin{theorem}
% 6.B.15
\label{propfact}
Let $\pi$ be a factor representation of a topological group $G$.
\begin{enumerate}[label=(\arabic*)]
\item\label{1DEpropfact}
Every representation which is quasi-equivalent to $\pi$ is factorial.
\item\label{2DEpropfact}
Every non trivial subrepresentation of $\pi$ is factorial and quasi-equivalent to $\pi$.
\item\label{3DEpropfact}
Let $\pi_1$ be a representation in the quasi-equivalence class of $\pi$;
then, either $\pi_1$ is a subrepresentation of $\pi$ or $\pi$ is a subrepresentation of $\pi_1$.
\item\label{4DEpropfact}
Let $\sigma$ be a factor representation which is not quasi-equivalent to $\pi$;
then $\sigma$ and $\pi$ are disjoint.
\item\label{5DEpropfact}
The representation $\pi$ is of type I if and only if
$\pi$ is a multiple of an irreducible representation $\sigma$ of $G$;
when $\pi$ is of type I, $\sigma$ is uniquely defined up to equivalence by $\pi$.
\end{enumerate}
\end{theorem}

\begin{rem}
% 6.B.16
\label{rem-propfact}
\noindent
(1)
In view of Claim \ref{5DEpropfact} of Theorem \ref{propfact},
we may refine the type I classification of factor representations as follows:
given a cardinal $n \ge 1$, a factor representation is said to be of type I$_n$
if $\pi$ is quasi-equivalent to an irreducible representation of dimension~$n$.

\vskip.2cm

(2)
As we will see later, 
see Section~\ref{S:ICCGroupsNontypeI} and Chapter~\ref{Chap:NormalInfiniteRep},
there are factor representations which are not of type I. 
They are classified in various other types,
via the Murray--von Neumann classification of factors in types
(see Section~\ref{SectionvN+rep}),
It is worth mentioning that, starting from Theorem~\ref{propfact},
one may independently develop a ``multiplicity theory''
for factor representations and recover the Murray--von Neumann type classification,
as shown in \cite[Chap.~1]{Mack--76}.
\end{rem}

The following result shows that a type I representation
has a unique direct sum decomposition
into representations of well-defined multiplicities.
For the particular case of a second-countable locally compact abelian group,
see the proof of Theorem~\ref{Thm-DecRepAbelianGroups}.
For the general case, we quote \cite[Chap.~III, \S~3, no~1]{Dixm--vN}
(where $\Hi$ need not be separable, 
and the set $I$ of cardinals need not be inside $\overline{\N^*}$);
see also \cite[5.4.9]{Dixm--C*}.
% et \cite[A.50]{Dixm--C*}.

\begin{theorem}
% 6.B.17
\label{Theo-DecTypeI-MultiplicityFree}
Let $\pi$ be a type I representation of $G$
in a separable Hilbert space~$\Hi$.
\par

There exist a set $I \subset \overline{\N^*}$ of extended positive integers
and an orthogonal decomposition $\Hi = \bigoplus_{m\in I} \Hi_m$ into
$G$-invariant closed subspaces $\Hi_m$ with the the following properties:
\begin{enumerate}[label=(\arabic*)]
\item\label{1DETheo-DecTypeI-MultiplicityFree}
for $m, n \in I$ with $m \ne n$, the subrepresentations of $\pi$ defined in $\Hi_m$
and $\Hi_n$ are disjoint;
\item\label{2DETheo-DecTypeI-MultiplicityFree}
for every $m\in I$, the subrepresentation of $\pi$ defined in $\Hi_m$
is equivalent to $m \pi_m$ for a multiplicity-free representation $\pi_m$ of $G$.
\end{enumerate}

Moreover, the orthogonal projection on $\Hi_m$ belongs to the centre of $\pi(G)'$ for every $m \in I$
and this decomposition is unique in the following sense:
if we have $\Hi = \bigoplus_{n \in J} \Ki_n$
for a set $J$ of extended integers and for $G$-invariant closed subspaces $\Ki_n$
with properties \ref{1DETheo-DecTypeI-MultiplicityFree} and \ref{2DETheo-DecTypeI-MultiplicityFree}, 
then $I = J$ and $\Hi_m = \Ki_m$ for all $m \in I$.
 \end{theorem}

For multiplicity-free representations,
we have a canonical decomposition in the sense of the following theorem:

\begin{theorem}
% 6.B.18
\label{Thm-DecDirectIntegral-MultiplicityFree}
Let $G$ be a second-countable locally compact group
and $\pi$ a multiplicity-free representation of $G$
in a separable Hilbert space $\Hi$. Let 
$$
\pi \, \simeq \, \int^\oplus_X \pi_x d\mu(x)
\, \simeq \, \int^\oplus_Y \sigma_y d\nu(y)
$$
be two integral decompositions of $\pi$ over standard Borel spaces $X, Y$
equipped with $\sigma$-finite measures $\mu, \nu$, 
and measurable fields $x \to \pi_x$, $y \to \sigma_x$ of irreducible representations.
\par

Then (after appropriate modifications on subsets of measure zero)
there exists a Borel isomorphism
$f \, \colon X \to Y$ such that $f_*(\mu)$ is equivalent to $\nu$
and such that 
$\sigma_{f(x)}$ and $\pi_x$ are equivalent 
for $\mu$-almost every $x \in X$.
\end{theorem}

\begin{proof}
Let $x \twoheadrightarrow \Hi_x$ and $y \twoheadrightarrow \Ki_y$
be the measurable fields of Hilbert spaces
defining $\int^\oplus_X \pi_x d\mu(x)$ and $\int^\oplus_Y \sigma_y d\nu(y)$.
We can and will identify $\Hi$ with $\int^\oplus_X \Hi_x d\mu(x)$
and $\pi$ with $\int^\oplus_X \pi_x d\mu(x)$.
Set
$$
\Ki \, := \, \int^\oplus_Y \Ki_y d\nu(y)
\hskip.2cm \text{and} \hskip.2cm
\sigma \, := \, \int^\oplus_Y \sigma_y d\nu(y)
$$
and let $U \, \colon \Hi \to \Ki$ be an isomorphism of Hilbert spaces
such that $U \pi(g) = \sigma(g)U$ for every $g \in G$.
\par

Let $\mathcal A$ and $\mathcal B$ be
the algebras of diagonalisable operators in $\Li (\Hi)$ and $\Li (\Ki)$ respectively.
Recall that we have isomorphisms
$L^\infty(X, \mu) \to \mathcal A, \hskip.1cm \varphi \mapsto m(\varphi)$ and
$L^\infty(Y, \nu) \to \mathcal A, \hskip.1cm \psi \mapsto m(\psi)$.
Since the $\pi_x$'s and the $\sigma_y$'s are irreducible, 
$\mathcal A$ and $\mathcal B$ are maximal abelian subalgebras
of $\pi(G)'$ and $\sigma(G)'$,
by Proposition~\ref{Pro-IntDecUniMaxAbelian}.
As $\pi$ and $\sigma$ are multiplicity-free,
$\pi(G)'$ and $\sigma(G)' = U \pi(G)'U^{-1}$ are abelian
and hence coincide with $\mathcal A$ and $\mathcal B$, respectively.
Therefore 
$$
m(\varphi) \mapsto U m(\varphi) U^{-1}
$$
is an isomorphism from $\mathcal A$ onto $\mathcal B$,
In this way, we obtain an isomorphism of $*$-algebras 
$$
F \, \colon \, L^\infty (X, \mu) \to L^\infty (Y, \nu) 
$$
such that $U m(\varphi) U^{-1} = m(F(\varphi))$.
By Theorem~\ref{Theo-VonNeumann},
upon modifying $X$ and $Y$ on subsets of measure zero,
we can assume that
there exists a Borel isomorphism $f \, \colon X \to Y$ such that 
$f_*(\mu)$ is equivalent to $\nu$ and such that 
$F(\varphi) (f(x)) = \varphi(x)$ for all $\varphi \in L^\infty (X, \mu)$ and $x \in X$.
It follows that, without loss of generality,
we can identify the measure space $(Y, \nu)$
with $(X, \mu)$ and $f$ with $\mathrm{Id}_X$;
so, $F = \mathrm{Id}_{L^\infty(X, \mu)}$ and we have 
$$
U m(\varphi) U^{-1} \, = \, m(\varphi)
\hskip.2cm \text{for every} \hskip.2cm
\varphi \in L^\infty (X, \mu).
$$
It follows from Proposition~\ref{caractdecop-DirectIntegral}
that $U$ is a decomposable operator;
so, there exists a measurable field $x \to U_x$
of Hilbert space isomorphisms $U_x \, \colon \Hi_x \to \Ki_x$
such that
$$
(U \xi)_x \, = \, U_x (\xi_x)
\hskip.5cm \text{for} \hskip.2cm
\xi \in \Hi
\hskip.2cm \text{and} \hskip.2cm
x \in X.
$$
Since $U \pi(g) = \sigma(g) U$, it follows that $U_x \pi_x(g) = \sigma_x(g) U_x$
for every $g \in G$ and $\mu$-almost every $x \in X$. 
\end{proof}

\subsection
{Examples of multiplicity-free representations of abelian groups}
% subsection 6.B.d
\label{SS:ExaMultFreeRep}
\index{Representation! multiplicity-free}
\index{Multiplicity-free representation}

\begin{exe}
% 6.B.19
\label{Exa-RegRep}
The left regular representation $\lambda_G$ of LC group $G$
is multipli\-city-free if and only if the group $G$ is abelian.
\par

To show this, it suffices by Proposition~\ref{Pro-MultFreeVN}
to show that $\lambda_G(G)'$ is abelian if and only if $G$ is abelian.
Assume first that $\lambda_G(G)'$ is abelian.
Since $\rho_G(G)$ is obviously contained in $\lambda_G(G)'$
and since $\rho_G$ is faithful, it follows that $G$ is abelian.
For the converse, we use the general fact that the algebra $\lambda_G(G)'$ 
coincides with the von Neumann algebra $\rho_G(G)''$ generated by $\rho_G(G)$;
see Theorem~\ref{Theo-RegRepHilbertAlg} and Remark \ref{Rem-Theo-RegRepHilbertAlg}.
% where this fact is proved for every unimodular group $G$),
When $G$ is abelian, this implies that $\lambda_G(G)'$ is abelian.
\end{exe}

Let $G$ be a second-countable LC abelian group.
Recall from Construction~\ref{defpimupourGab++} that,
given a probability measure $\mu$ on $\widehat G$,
an extended positive integer $n \in \overline{\N^*}$,
and a Hilbert space $\Ki$ of dimension $n$,
we associate the representation $\pi_\mu^{(n)}$ of $G$
defined on $L^2(\widehat G, \mu, \Ki)$ by 
$$
(\pi_\mu^{(n)}(g)F) (\chi) \, = \, \chi(g) F(\chi)
\hskip.5cm \text{for} \hskip.2cm
g \in G, \hskip.1cm F \in L^2(\widehat G, \mu, \Ki), 
\hskip.2cm \text{and} \hskip.2cm
\chi \in \widehat G .
$$
When $n = 1$, the representation $\pi_\mu^{(1)}$
is the \textbf{canonical representation} $\pi_\mu$
of Construction~\ref{defpimupourGab} associated to $\mu$.

\begin{theorem}
% 6.B.20
\label{Thm-CanRepAbMultFree}
Let $G$ be a second-countable locally compact abelian group
and $\pi$ a representation of $G$ in a separable Hilbert space. 
The following properties are equivalent:
\begin{enumerate}[label=(\roman*)]
\item\label{iDEThm-CanRepAbMultFree}
$\pi$ is multiplicity-free;
\item\label{iiDEThm-CanRepAbMultFree}
$\pi$ is equivalent to a canonical representation $\pi_\mu$
for some probability measure on the Borel subsets of $\widehat G$.
\end{enumerate}
\end{theorem}

\begin{proof}
We know that the commutant $\pi_\mu(G)'$ is abelian 
(Proposition~\ref{Prop-RepAb-Equiv}
and Remark~\ref{Rem-Prop-RepAb-Equiv}),
hence $\pi_\mu$ is multiplicity-free, by Proposition~\ref{Pro-MultFreeVN}.
So, \ref{iiDEThm-CanRepAbMultFree} implies~\ref{iDEThm-CanRepAbMultFree}.

\vskip.2cm

For the converse, consider a representation $\pi$ of $G$
in a separable Hilbert space.
By Theorem~\ref{Thm-DecRepAbelianGroups}, we may assume that 
$$
\pi \, = \, \bigoplus_{n \in I} \pi_{\mu_n}^{(n)},
$$
where $I$ is a set of extended positive integers,
$(\mu_n)_{n \in I}$ a sequence of probability measures on $\widehat G$,
and $\pi_{\mu_n}^{(m)}$ for each $n \in I$
the representation on $L^2(\widehat G, \mu_n, \Ki_n)$ defined above.
\par

On the one hand, by Proposition~\ref{Prop-RepAb-Equiv},
the commutant $\pi_{\mu_n}^{(n)}(G)'$ of $\pi_{\mu_n}^{(n)}(G)$
inside $\Li (L^2(\widehat G, \mu_n, \Ki_n))$ consists of
all decomposable operators in $\Li (L^2(\widehat G, \mu_n, \Ki_n))$.
On the other hand,
observe that there is an obvious injection of $\pi_{\mu_n}^{(n)}(G)'$ into $\pi(G)'$.
\par

Assume now that $\pi$ is multiplicity-free,
so that $\pi(G)'$ is abelian by Proposition \ref{Pro-MultFreeVN},
and therefore $\pi_\mu^{(n)}(G)'$ is abelian;
this is possible only if $n = 1$.
It follows that $I = \{ 1 \}$, that is, $\pi = \pi_{\mu_1}$.
So, \ref{iDEThm-CanRepAbMultFree} implies \ref{iiDEThm-CanRepAbMultFree}.
\end{proof}

The following corollary is a direct consequence of Theorem~\ref{Thm-DecRepAbelianGroups}
and Theorem~\ref{Thm-CanRepAbMultFree}.

\begin{cor}
% 6.B.21
\label{Cor-Thm-CanRepAbMultFree}
Let $G$ be a second-countable locally compact abelian group.
\par

Every representation $\pi$ of $G$ in a separable Hilbert space is of type I.
\end{cor}

 In a terminology to be introduced in Section~\ref{SectionTypeI},
Corollary~\ref{Cor-Thm-CanRepAbMultFree} states that
a second-countable LC abelian group is of type I.
This result extends to every topological abelian group,
as we will mention again below (Theorem \ref{explesTypeI-bis}).

\subsection
{Examples of multiplicity-free representations:
quasi-regular representations of Gel'fand pairs}
% subsection 6.B.d

Let $G$ be a LC group and $K$ a compact subgroup. 
Let $C^c(K\backslash G/K)$ be the space of
complex-valued continuous functions of compact support on $G$
which are left and right-invariant under translation by elements from $K$, that is, 
$$
C^c (K \backslash G / K) \, = \,
\{ \varphi \in C^c(G) \mid
\varphi(kgk') = \varphi(g)
\hskip.2cm \text{for all} \hskip.2cm
g \in G, \hskip.1cm k,k' \in K \}.
$$
Recall that the convolution of two functions $f_1, f_2 \in C^c(G)$ is 
the function $f_1 \ast f_2 \in C^c(G)$ defined by
$$
f_1 \ast f_2(g) \, = \, \int_G f_1(x) f_2(x^{-1} g) d\mu_G(x)
\, = \, \int_G f_1(gx) f_2(x^{-1}) d\mu_G(x)
\hskip.5cm \text{for} \hskip.2cm
g \in G,
$$
where $\mu_G$ is a left Haar measure on $G$.
(Note that this product $\ast$ does depend on the choice of the Haar measure $\mu_G$.)
Then $C^c(G)$ is an algebra for the convolution product
and $C^c (K \backslash G / K) $ is a subalgebra of $C^c(G)$.
\par

There is a linear projection $\Phi \, \colon C^c(G) \to C^c (K \backslash G / K)$ given by 
$$
\Phi(f) (g) = \int_K \int_K f(k_1g k_2) \hskip.1cm d\mu_K(k_1) \hskip.1cm d\mu_K(k_2)
\hskip.5cm \text{for} \hskip.2cm
f \in C^c(G), \ g \in G,
$$
where $\mu_K$ is the normalized Haar measure on $K$;
in other words, $\Phi(f) = \mu_K \ast f\ast \mu_K$ for $f \in C^c(G)$.

\begin{defn}
% 6.B.22
\label{Def-Gel'fandPair}
A \textbf{Gel'fand pair} is a pair $(G, K)$,
where $G$ is a locally compact group and $K$ a compact subgroup,
such that the convolution algebra $C^c (K \backslash G / K)$ is commutative.
\index{Gel'fand pair} 
\end{defn}

Here are two first straightforward examples of Gel'fand pairs.

\begin{exe}
% 6.B.23
\label{Exa-Gel'fandPair}
(1)
Let $G$ be a LC abelian group. It is obvious that 
$(G, \{e\})$ is a Gel'fand pair.

\vskip.2cm

\index{Motion group}
(2)
Let $G$ be a \textbf{motion group}, that is, $G = K \ltimes A$ is the semi-direct product
of a compact subgroup $K$ and a closed abelian normal subgroup $A$. 
Then $(G, K)$ is a Gel'fand pair. Indeed, the algebra $C^c (K \backslash G / K)$ 
is canonically isomorphic to the convolution algebra $C^c(A)$,
which is obviously commutative.
\end{exe}

More interesting examples of Gel'fand pairs
are provided by the following result from \cite{Gelf--50}.
\par

In the sequel, for a function $f$ on a group $G$,
we denote by $\check{f}$ the function on $G$ 
defined by $\check{f}(g) = f(g^{-1})$ for $g \in G$.
 
\begin{theorem}[Gel'fand]
% 6.B.24
\label{Theo-Gel'fandPairSemissimple}
Let $G$ be a locally compact group and $K$ a compact subgroup.
Assume that there exists a continuous involutive automorphism $\theta$ of $G$
with the following property:
\begin{enumerate}
\item[]
$\theta(g^{-1})$ belongs to the double coset $KgK$, for every $g \in G$.
\end{enumerate}
Then $(G, K)$ is a Gel'fand pair.
\end{theorem}

\begin{proof}
The image $\theta_*(\mu_G)$ of the left Haar measure $\mu_G$ by $\theta$
is obviously a left-invariant Radon measure on $G$.
Therefore $\theta_*(\mu_G) = \lambda \mu_G$ for some $\lambda>0$.
Since $\theta^2 = I$, we have $\lambda = 1$,
that is, $\mu_G$ is invariant under $\theta$.
It follows that $\theta$ defines an automorphism $f \mapsto f^\theta$
of the convolution algebra $C^c(G)$, given by
$$
f^\theta (g) \, = \, f(\theta(g))
\hskip.5cm \text{for all} \hskip.2cm
f \in C^c(G), \hskip.1cm g \in G.
$$
Let $f_1, f_2 \in C^c (G)$. Then, on the one hand, we have
$$
f_1^\theta \ast f_2^\theta \, = \, (f_1 \ast f_2)^\theta.
\leqno{(*)}
$$
On the other hand, we have
$$
\begin{aligned}
(\check{f_1} \ast \check{f_2}) (g)
\, &= \, \int_G f_1(x^{-1}) f_2(g^{-1}x) d\mu_G(x)
\\
\, &= \, \int_G f_2(y) f_1(y^{-1} g^{-1}) d\mu_G(y)
\, = \, f_2 \ast f_1(g^{-1})
\end{aligned}
$$
for every $g \in G$, that is, 
$$
\check{f_1} \ast \check{f_2} \, = \, (f_2 \ast f_1)\check{}.
\leqno{(**)}
$$
Now, let $f \in C^c (K \backslash G / K)$; it follows from the assumption that
$$
f^\theta(g) \, = \, f(\theta(g)) \, = \, f(g^{-1})
\hskip.5cm \text{for all} \hskip.2cm
g \in G,
$$
that is,
$$
f^\theta \, = \, \check{f} .
\leqno{(***)}
$$
Relations $(*), (**)$ and $(***)$ imply that 
$$
f_1 \ast f_2
\, \overset{(**)}{=} \, (\check{f_2} \ast \check{f_1})\check{}
\, \overset{(***)}{=} \, (f_2^\theta \ast f_1^\theta)^\theta
\, \overset{(*)}{=} \, f_2 \ast f_1
\hskip.5cm \text{for all} \hskip.2cm
f_1, f_2 \in C^c (K \backslash G / K) ,
$$
and this concludes the proof.
\end{proof}
 
\begin{exe}
% 6.B.25
\label{Exa-Gel'fandPairSemissimple}
(1)
Let $G$ be a LC group. Assume that there exists an action of $G$
by isometries on a metric space $(X,d)$ with the following properties:
\begin{itemize}
\setlength\itemsep{0em}
\item
the action $G \curvearrowright (X,d)$
is transitive on equidistant pairs of points, that is,
for every $(x,y), (x',y')$ in $X \times X$ such that $d(x,y) = d(x',y')$
there exists $g \in G$ such that $x' = gx$ and $y' = gy$;
\item
the stabilizer $K$ in $G$ of a point $x_0\in X$ is compact.
\end{itemize} 
Then $(G, K)$ is a Gel'fand pair.
Indeed, let us check that $\theta = \mathrm{Id}_G$ has the property stated
in Theorem~\ref{Theo-Gel'fandPairSemissimple}.
Let $g \in G$. Then 
$$
d(x_0, g x_0) \, = \, d(g^{-1}x_0, x_0) \, = \, d(x_0,g^{-1}x_0)
$$
and so there exists $k\in K$ with $g^{-1} x_0 = kgx_0$ and therefore
$g^{-1} \in K g K$. 
\par

The following examples are classical.
Let $n$ be a positive integer.
The Gel'fand pair $(\SO(n), \SO(n-1))$ arises
from the natural action of $\SO(n)$ on the sphere ${\mathbf S}^{n-1}$.
The Gel'fand pair $(\SO(n) \ltimes \R^n, \SO(n)$ arises
from the natural action of $(\SO(n) \ltimes \R^n$ on the Euclidean space $\R^n$.
For $n \ge 2$, the Gel'fand pair $(\SO(n,1), \SO(n))$ arises
from the natural action of $\SO(n,1)$ on the hyperbolic space ${\mathbf H}^{n}(\R)$.
\par

Let $T_d$ be a regular tree of some some degree $d \ge 3$.
Let $\Aut (T_d)$ be its automorphism group, with the topology of the pointwise convergence,
which makes it a locally compact group.
Then the action of $\Aut (T_d)$ on $T_d$ is transitive on equidistant pairs of points.
If $K$ denotes the stabilizer of some vertex of $T_d$,
it follows that $(\Aut (T_d), K)$ is a Gelfand pair.

\vskip.2cm

(2)
Let $G$ be a non-compact semisimple connected real Lie group with finite centre 
and let $K$ be a maximal compact subgroup of $G$.
There exists an involutive automorphism $\theta$ of $G$,
called \emph{Cartan involution}, 
such that $K$ is the set of $\theta$-fixed elements in $G$.
Let $\mathfrak g$ be the Lie algebra of $G$ and let 
$$
\mathfrak p \, = \, \{ X \in \mathfrak g \mid d\theta_e(X) = -X \},
$$
where $d\theta_e \, \colon \mathfrak g \to \mathfrak g$
is the differential of $\theta$ at the group unit.
We have the so-called \textbf{Cartan decomposition}
$$
G \, = \, K \exp\mathfrak p
$$
of $G$ (see \cite[Theorem 6.31]{Knap--02},
and our Lemma~\ref{Lem-CartanDecSL2} for the special case $G = \SL_2(\R)$).
\par

We claim that $(G, K)$ is a Gel'fand pair.
To show this, let us check that $\theta$ has the property
stated in Theorem~\ref{Theo-Gel'fandPairSemissimple}.
Let $g = k \exp (X)$ with $k \in K$ and $X \in \mathfrak p$.
Then $g^{-1} = \exp (-X) k^{-1} $ and so 
$$
\begin{aligned}
\theta(g^{-1})
\, &= \, \theta(\exp (-X)) k^{-1}
\, = \, \exp (-d\theta_e(X)) k^{-1}
\, = \, \exp (X) k^{-1}
\\
\, &= \, k^{-1} g k^{-1}
\, \in \, K gK.
\end{aligned}
$$
\par
An example of such a pair $(G, K)$ is the pair $(SL_n(\R), \SO(n))$, for $n \ge 2$.

\vskip.2cm

(3)
Let $\K$ be a non-Archimedean local field,
$\Bbb G$ be an algebraic group defined over $\K$,
and $G = \Bbb G(\K)$ the group of $\K$-points of $\Bbb G$,
which is a locally compact group.
Assume that $\Bbb G$ is reductive.
Then $G$ has maximal compact subgroups which are good,
in the sense of Bruhat--Tits;
let $K$ be one of them.
Then $(G, K)$ is a Gelfand pair \cite[4.4.9]{BrTi--72}.
% page 82 (79) de l'article
\par

For example, $(\SL_d(\Q_p), \SL_d(\Z_p))$ is a Gelfand pair
for all primes $p$ and all integers $d \ge 2$.
\end{exe}

\begin{rem}
% 6.B.26
In the case of $(G, K) = (\SL_2(\R), \SO(2))$,
we will see in Proposition~\ref{Prop-CommutativeAlgebraSpherical} below
that $C^c (K \backslash G / K)$ is part of a whole family
of commutative algebras attached to $(G, K)$.
\end{rem} 

Let $G$ be a LC group with left Haar measure $\mu_G$.
Observe that the modular function $\Delta_G$ of $G$ is determined by the equality 
$$
\int_G f(x^{-1})d\mu_G(x) \, = \, \int_G \Delta_G(x^{-1}) f(x)d\mu_G(x)
\hskip.5cm \text{for all} \hskip.2cm
f \in C_c(G)
$$
(see \cite[Lemma A.3.4]{BeHV--08}).

\begin{lem}
% 6.B.27
\label{Lem-UnimodularGel'fandPair}
Let $G$ be a locally compact group containing a compact subgroup $K$
such that $(G, K)$ is a Gel'fand pair.
\par

Then $G$ is unimodular.
\end{lem}
% avant : cit\'e de \cite[Chap.~IV, \S~3, Theorem~3.1]{Helg--84}.

\begin{proof}.
It suffices to show that 
$$
\int_G f(x^{-1})d\mu_G(x) \, = \, \int_G f(x)d\mu_G(x)
\hskip.5cm \text{for all} \hskip.2cm
f \in C_c(G) .
\leqno{(+)}
$$
We claim first that Equality $(+)$ holds for $f \in C^c (K \backslash G / K)$.
Indeed, choose $h\in C^c (K \backslash G / K)$ with $h(x) = h(x^{-1}) = 1$
for every $x \in \mathrm{supp} (f)$.
Then, since $ C^c (K \backslash G / K)$ is commutative,
we have $f \ast h = h \ast f$ and hence
$$
\begin{aligned}
\int_G f(x^{-1})d\mu_G(x)
&\, = \, \int_G f(x^{-1}) h(x) d\mu_G(x)
\, = \, (f\ast h) (e) = (h\ast f) (e)
\\
&\, = \, \int_G h(x^{-1}) f(x) d\mu_G(x
\, = \, \int_G f(x)d\mu_G(x).
\end{aligned}
$$
Next, we prove Equality $(+)$ for arbitrary $f \in C_c(G)$.
\par

Since $K$ is compact, we have $\Delta_G\vert_K = 1$.
Therefore $\mu_G$ is right-invariant by $K$.
It follows that, for every $f \in C_c(G)$, we have
$$
\begin{aligned}
\int_G (\mu_K \ast f\ast \mu_K ) (x) \mu_G(x)
&\, = \, \int_G \left(\int_K \int_K f(k_1x k_2) d\mu_K(k_1)d\mu_K(k_2)\right) d\mu_G(x)
\\
&\, = \, \int_G f(x)d\mu_ G(x),
\end{aligned} 
$$
where $\mu_K$ is the normalized Haar measure on $K$.
Since $\mu_K \ast f\ast \mu_K$ belongs to $C^c (K \backslash G / K)$
and since $\mu_K \ast \check{f} \ast \mu_K = \check{\mu_K \ast f \ast \mu_K}$,
Equality $(+)$ holds for $f$.
\end{proof}

\begin{defn}
% 6.B.28
\label{defquasiregG/K}
Let $G$ be a unimodular LC group and $K$ a compact subgroup. 
Then $G/K$ carries a $G$-invariant Borel measure $\mu_{G/K}$, which 
is unique up to a scalar multiple (see \cite[Corollary B.1.7]{BeHV--08}).
We denote by $\lambda_{G/K}$
the \textbf{quasi-regular representation} $\lambda_{G/K}$ of $G$ on $L^2(G/K, \mu_{G/K})$,
that is, the natural representation of $G$
given on $L^2(G/K, \mu_{G/K})$ by left translations.
\par

Compare with Example \ref{quasiregularrep}.
\index{Quasi-regular representation}
\index{Representation! quasi-regular}
\end{defn}

Let $(\pi, \Hi)$ be a representation of a LC group $G$.
Then there is an associated representation $C^c(G) \to \Li (\Hi)$,
of the convolution algebra $C^c(G)$, denoted by $\pi$ again, defined by
$$
\pi(f) \, = \, \int_G f(g) \pi(g) d\mu_G(g) 
\hskip.5cm \text{for all} \hskip.2cm 
f \in C^c(G)
$$
More generally, let $M^b(G)$ be the Banach $*$-algebra of complex Radon measures on $G$.
Then a $*$-representation $M^b(G) \to \Li (\Hi)$, again denoted by $\pi$, is defined by
$$
\pi(\mu) \, = \, \int_G \pi(g) d\mu(g)
\hskip.5cm \text{for all} \hskip.2cm 
\mu \in M^b(G);
$$
see Section~\ref{C*algLCgroup} below.
The von Neumann algebras generated by $\pi(G)$, by $\pi(C^c(G))$, and by $M^b(G)$
all coincide (see Proposition~\ref{vNpiGallsame}).
\par

The following criterion for Gel'fand pairs is well-known (see \cite[Theorem 9.7.1]{Wolf--07}).

\begin{theorem}
% 6.B.29
\label{Thm-MultFreeGel'fandPair}
Let $G$ be locally compact group and $K$ a compact subgroup.
The following properties are equivalent:
\begin{enumerate}[label=(\roman*)]
\item\label{iDEThm-MultFreeGel'fandPair}
$(G, K)$ is a Gel'fand pair;
\item\label{iiDEThm-MultFreeGel'fandPair}
the quasi-regular representation 
$\lambda_{G/K}$ of $G$ on $L^2(G/K))$ is multiplicity-free. 
\end{enumerate}
\end{theorem}

\index{Quasi-regular representation}
\index{Representation! quasi-regular}

\begin{proof}
We view $L^2(G/K)$ as the closed subspace of 
$L^2(G)$ consisting of the functions of $L^2(G)$
which are invariant under right translation by elements in $K$.
Observe that $L^2(G/K)$ is invariant under the regular representation 
$\lambda_G$ and that the quasi-regular $\lambda_{G/K}$ coincides
with the subrepresentation of $\lambda_G$ defined by $L^2(G/K)$.
\par

The orthogonal projection $P$ from $L^2(G)$ onto $L^2(G/K)$
is given $P = \lambda_G (\mu_K)$, that is, 
$$
P(f) (gK) \, = \, \int_K f(gk) d\mu_K(k),
$$
for $f \in L^2(G)$ and $\mu_G$-almost every $g \in G$,
where $\mu_G$ is a Haar measure on $G$
and $\mu_K$ the normalized Haar measure on $K$.
\par

Let $\mathcal M$ be the von Neumann algebra 
generated by $\lambda_{G/K}(C^c(G))$ in $\Li (L^2(G/K))$.
On the one hand, we have
$$
\mathcal M \, = \, P\Li (G) P,
$$ 
where $\Li (G) \, = \, \lambda_G(C^c(G))'' \subset \Li (L^2(G))$
is the von Neumann algebra of $G$.
On the other hand, by Lemma~\ref{Lem-UnimodularGel'fandPair}, $G$ is unimodular and so 
$$
\Li (G) ' \, = \, \Ri (G),
$$
where $\Ri (G) = \rho(G)''$ is the von Neumann algebra
generated by the right regular representation $\rho_G$ of $G$
(see Theorem~\ref{Theo-RegRepHilbertAlg} below).
By Proposition~\ref{Prop-MatrixCommutant}, we have
$$
\mathcal M' \, = \, P\Ri (G)P ,
$$
and it follows that $P \rho_G(C^c(G)) P$ is dense in $\mathcal M'$
for the strong operator topology.
Since
$$
\{P \rho_G(f) P\mid f \in C^c(G)\}
\, = \, \{\rho_{G}( \mu_K \ast f \ast \mu_K) \mid f \in C^c(G)\}
\, = \, \rho_G(C^c (K \backslash G / K))
$$
and since $\rho_G$ is faithful on $C^c(G)$, 
we see that $C^c (K \backslash G / K)$ is commutative
if and only if $\mathcal M'$ is commutative,
hence by Proposition~\ref{Pro-MultFreeVN} if and only if $\lambda_{G/K}$ is multiplicity-free.
\end{proof}

\section
{The quasi-dual of a topological group}
% Section 6.3
\label{SectionQuasidual}

In this section, we introduce the quasi-dual of a topological group $G$.
We will see that this new dual space as well as the unitary dual of $G$
carry natural Borel structures.
In particular, these Borel structures will allow us
to perform direct integrals decompositions of representations
over the quasi-dual or the dual of $G$.

\subsection
{The quasi-dual}
% subsection 6.C.a
\label{SS:QuasiDual}

Let $G$ be a topological group.

\begin{defn}
% 6.C.1
\label{Def-QuasiDual}
The \textbf{quasi-dual} $\QD(G)$ of $G$ is 
the space of quasi-equivalence classes of factor representations of $G$.
\index{Quasi-dual}
\index{$c2$@$\QD(G)$ quasi-dual of the topological group $G$}
\end{defn}

\begin{defn}
% 6.C.2
\label{Def-TypIPartQD}
Every irreducible representation of $G$ is factorial, indeed factorial of type~I,
and two irreducible representations of $G$ are quasi-equivalent
if and only if they are equivalent (Proposition~\ref{Prop-EqQuasiEq-Immediate}). 
Consequently, there is a well defined injection 
$\kappa^{\rm d}_{\rm qd}$ from the dual $\widehat G$ into the quasi-dual $\QD(G)$,
mapping the equivalence class of an irreducible representation of $G$
to its quasi-equivalence class.
\par

The image of the map $\kappa^{\rm d}_{\rm qd}$,
which is the subspace of $\QD(G)$ consisting of the quasi-equivalence classes
of factor representations of type I, 
is called the \textbf{type I part of the quasi-dual} and is denoted by $\QD(G)_{\rm I}$.
\end{defn}

\index{Quasi-dual! $1$@type I part}
\index{Type I! $2$@part of the quasi-dual}
The map $\kappa^{\rm d}_{\rm qd}$ can be written
as the composition of a bijection and an inclusion:
\begin{equation}
\label{eqq/d/qd}
\tag{$\kappa$2}
\kappa^{\rm d}_{\rm qd} \, \colon \, \widehat G \overset{\text{bij}}{\longrightarrow} 
\QD(G)_{\rm I} 
\subset \QD(G) .
\end{equation}
As will be seen in Section~\ref{S:ICCGroupsNontypeI},
the map $\kappa^{\rm d}_{\rm qd}$ is not surjective in general.

\subsection[Topologies and Borel structures on duals and quasi-duals]
{Topologies and Borel structures on duals and quasi-duals}
% subsection 6.C.b
\label{SS:BorelStructureQD}

Let $\aleph$ be a cardinal.
For every cardinal $n \le \aleph$, 
let $\Hi_n$ be a fixed Hilbert space of dimension $n$.
Let $\mathcal R_n$ be a set of representations of $G$ in $\Hi_n$;
the main examples of interest of sets $\mathcal R_n$ are
\begin{itemize}
\setlength\itemsep{0em}
\item
the set ${\rm Rep}_n(G)$ of all representations of $G$ in $\Hi_n$,
\item
the subset ${\rm Fac}_n(G)$ of all factor representations of $G$ in $\Hi_n$,
\item
and the subset ${\rm Irr}_n(G)$ of all irreducible representations of $G$ in $\Hi_n$.
\end{itemize}
We consider $\mathcal R_n$
with the topology of weak uniform convergence on compact subsets,
for which a net $\left(\pi_\nu \right)_\nu$ converges to $\pi$
if, for every $\xi, \eta \in \Hi_n$,
the net $\left( \langle \pi_\nu(g) \xi \mid \eta \rangle \right)_\nu$ 
converges to $\langle \pi(g) \xi \mid \eta \rangle$
uniformly on every compact subset of $G$.
We equip the disjoint union 
$$
\mathcal R_{\le \aleph} \, := \, \bigsqcup_{n \le \aleph} \mathcal R_n
$$
with the sum topology.
(All this as in \cite[18.1]{Dixm--C*}.)
\par

The space $\mathcal R_{\le \aleph}$
is also endowed with the associated Borel structure,
equivalently with the coarsest Borel structure for which the functions
$\pi \mapsto \langle \pi(g) \xi \mid \eta \rangle$ are Borel measurable,
for all $g \in G$ and $\xi, \eta \in \Hi_\pi$.

\begin{rem}
% 6.C.3
\label{Rem-Mackey topology}
A net $\left(\pi_\nu \right)_\nu$ in $\mathcal R_n$ converges to $\pi$
if and only if, for every $\xi \in \Hi_n$, we have
$$
\lim_\nu \Vert \pi_\nu(g) \xi -\pi(g) \xi \Vert \, = \, 0,
$$
uniformly on compact subsets of $G$.
Indeed, this follows from the relation
$$
\Vert \pi_\nu(g) \xi -\pi(g) \xi \Vert^2 \, = \,
2 \left( \Vert \xi \Vert^2 - \mathrm{Re} \langle \pi_\nu(g) \xi \mid \pi(g) \xi \rangle \right).
$$
\end{rem}

\begin{defn}
% 6.C.4
\label{DefFellTop2}
Let $\aleph$ be a cardinal such that
every factor representation of $G$ is quasi-equivalent to representation
from ${\rm Fac}_{\le \aleph}(G)$ (see Corollary~\ref{FacSeparable}).
Then
\begin{center}
\emph{the quotient space
\hskip.1cm ${\rm Fac} (G)_{\le \aleph}/ \hskip-.1cm \approx$ \hskip.1cm
can be identified with the quasi-dual $\QD(G)$}
\end{center}
(recall that $\approx$ denotes quasi-equivalence of representations).
The Fell topology on $\QD(G)$ is the quotient by the relation $\approx$
of the topology on ${\rm Fac}_n(G)$ defined above.
\par

Let $\aleph'$ be a cardinal such that
every irreducible representation of $G$ is equivalent to representation
from ${\rm Irr}_{\le \aleph'}(G)$
(see Proposition \ref{Pro-IrrRepBoundedCard},
and observe that $\aleph' = \aleph$ is a possible choice).
Then
\begin{center}
\emph{the quotient space
\hskip.1cm ${\rm Irr}_{\le \aleph'}(G) / \hskip-.1cm \simeq$ \hskip.1cm
can be identified with the dual $\widehat G$}
\end{center}
(recall that $\simeq $ denotes equivalence of representations).
The topology on $\widehat G$
which is the quotient by the relation $\simeq$
of the topology on ${\rm Irr}_{\le \aleph'}(G)$ defined above
is the Fell topology of Section \ref{SectionWC+FellTop}.
\index{Fell topology}
\par

Since irreducible representations are equivalent if and only if they are quasi-equivalent,
the identifications above show again that $\widehat G$ is naturally a subspace of $\QD(G)$.
\end{defn}

\begin{defn}
% 6.C.5
\label{Def-MackeyBorelStructure}
Let $\aleph, \aleph'$ be cardinals as in the previous definition.
Let ${\rm Fac}_{\le \aleph}(G)$ and ${\rm Irr}_{\le \aleph'}(G)$ and be endowed
with the Borel structure introduced above.
\par

The \textbf{Mackey--Borel structure} on the quasi-dual $\QD(G)$ of $G$ 
is by definition the quotient Borel structure on 
${\rm Fac}_{\le \aleph}(G) / \hskip-.1cm \approx$.
\par

The \textbf{Mackey--Borel structure} on the dual $\widehat G$ of $G$
is by definition the quotient Borel structure on 
${\rm Irr}_{\le \aleph'}(G) / \hskip-.1cm \simeq$.
\index{Mackey--Borel structure}
\end{defn}

\begin{rem}
% 6.C.6
\label{Rem-BorelStructure-Topology}
(1)
The Mackey--Borel structure on $\widehat G$
is always finer than the Borel structure associated to the Fell topology.
Indeed, this follows from the continuity of the map ${\rm Irr}_{\le \aleph'}(G) \to \widehat G$.
% see \cite[3.8.3]{Dixm--C*},
% where an hypothesis of separability is stated, but not used. 
% 19 dec 16 : Effectivement on peut se passer de la s\'eparabilit\'e de $A$ 
% (si on consid\`ere d'autres cardinaux).
These two Borel structures on $\widehat G$ need not coincide,
as will be seen below (Glimm Theorem \ref{ThGlimm}).
\par

The following trivial example shows an analogous phenomenon.
Define on the real line $\R$ an equivalence relation $\simeq$ by
$x \simeq y$ if $y-x \in \Q$.
On the quotient $\R / \hskip-.1cm \simeq$, the quotient topology is trivial,
with exactly two open subset $\emptyset$ and $\R / \hskip-.1cm \simeq$,
and consequently the associated Borel structure is trivial.
However, the quotient Borel structure is non-trivial:
every countable subset $B$ of $\R / \hskip-.1cm \simeq$ has an inverse image in $\R$
that is countable, in particular Borel, and therefore $B$ is Borel.
% Voir Berberian page 83-sur-84.

\vskip.2cm

(2)
Assume now that $G$ is a \emph{second-countable locally compact group}.
Then we may choose for $\aleph$ the cardinality $\aleph_0$ of countable infinite sets.
The spaces ${\rm Fac}_{\le \aleph_0}(G)$ and ${\rm Irr}_{\le \aleph_0}(G)$
are known to be Polish spaces \cite[3.7.1, 3.7.4]{Dixm--C*}.
\par

Moreover, $\QD(G)_{\rm I}$, the type I part of the quasi-dual,
is a Borel subspace of $\QD(G)$,
which is therefore Borel isomorphic to $\widehat G$ \cite[7.3.6]{Dixm--C*}.
\par
% Voir aussi \cite[Th\'eor\`eme 2]{Dixm--62} !!!

Also, every point of $\QD(G)$ is a Borel subset of $\QD(G)$; see \cite[7.2.4]{Dixm--C*}.
\end{rem}

\subsection
{Decomposition into factor representations}
% subsection 6.C.c
\label{SectionDecomposingFact}

Let now $G$ be a second-countable locally compact group.
\par

Recall that $G$ may have representations with two direct integral decompositions
involving \emph{irreducible} representations
belonging to two disjoint parts of $\widehat G$
(Theorem~\ref{Theo-NonUniqueIntDecIrrUniRep}).
\par

Given a representation $\pi$ of $G$,
one may consider direct integral decompositions of $\pi$
into \emph{factor} representations.
Of course, uniqueness of such decompositions may also fail,
as irreducible representations are factorial.
However, as we now show, one can distinguish among such decompositions 
a \emph{canonical} one.
\par

For the notion of a direct integral of von Neumann algebras
which is involved in the statement of the following result,
see Section~\ref{S: DirectInt-VN}.

\begin{theorem}
% 6.C.7
\label{thmDirectIntFact0}
Let $G$ be a second-countable locally compact group
and let $\pi$ be a representation of $G$ on a separable Hilbert space $\Hi$.
Let $\mathcal Z$ be the center of $\pi(G)'$.
\par

There exist a standard Borel space $X$,
a $\sigma$-finite measure $\mu$ on $X$,
a measurable field of Hilbert spaces $x \mapsto \Hi_x$ over $X$, 
a measurable field of representations $x \mapsto \pi_x$ of $G$
in the $\Hi_x$~'s over $X$,
and an isomorphism of Hilbert spaces
$U \, \colon \Hi \to \int_X^\oplus \Hi_x d\mu(x)$,
with the following properties:
\begin{enumerate}[label=(\arabic*)]
\item\label{iDEthmDirectIntFact0}
$
U \pi(g) U^{-1} \, = \, \int^\oplus_{X} \pi_x(g) d\mu(x)
\hskip.5cm \text{for all} \hskip.2cm
g \in G ;
$
\item\label{iiDEthmDirectIntFact0}
$U \mathcal Z U^{-1}$ coincides with
the algebra of diagonalisable operators on $\int_X^\oplus \Hi_x d\mu(x)$;
\item\label{iiiDEthmDirectIntFact0}
$\pi_x$ is factorial for $\mu$-almost every $x \in X$;
\item\label{ivDEthmDirectIntFact0}
there exists a measurable subset $N$ of $X$ with $\mu(N) = 0$,
such that $\pi_x$ and $\pi_y$ are disjoint
for every $x, y \in X \smallsetminus N$ with $x \ne y$.
\end{enumerate}
Moreover, we have
$$
U \pi(G)'' U^{-1} \, = \, \int^\oplus_{X} \pi_x(G)'' d\mu(x)
\hskip.5cm \text{and} \hskip.5cm
U \pi(G)' U^{-1} \, = \, \int^\oplus_{X} \pi_x(G)'d\mu(x) .
$$
\end{theorem}

\begin{proof}
The existence of $(X, \mu)$, $\int_X^\oplus \Hi_x d\mu(x)$,
$\int^\oplus_{X} \pi_x d\mu(x)$, and $U$,
as stated with properties \ref{iDEthmDirectIntFact0} and \ref{iiDEthmDirectIntFact0}
follows from Theorem~\ref{Theo-IntDecUniRep1}.
So, we may assume that 
$\Hi = \int_X^\oplus \Hi_x d\mu(x)$, $\pi = \int_X^\oplus \pi_x d\mu(x)$,
and $\mathcal Z$ is the algebra of diagonalizable operators on $\Hi$. 
\par

The fields $x \mapsto \pi_x(G)''$ and $x \mapsto \pi_x(G)'$
of von Neumann algebras are measurable
(see Example~\ref{Exa-measuableFieldvN}
and Lemma~\ref{Lemma-FieldvN}).
Moreover, by Corollary~\ref{Cor-Theo-CharDec}, we have 
$$
\pi(G)'' \, = \, \int^\oplus_{X} \pi_x(G)''d\mu(x)
\hskip.5cm \text{and} \hskip.5cm
\pi(G)' \, = \, \int^\oplus_{X} \pi_x(G)'d\mu(x) .
$$
So, the last statement in the theorem holds.
It remains to show that properties \ref{iiiDEthmDirectIntFact0} and \ref{ivDEthmDirectIntFact0} are satisfied.
\par

We claim that the field $x \mapsto \pi_x(G)'' \cap \pi_x(G)'$
of von Neumann algebras is measurable.
Indeed, by Lemma~\ref{Lemma-FieldvN},
it suffices to show that the field 
$$
x \mapsto \left( \pi_x(G)'' \cap \pi_x(G)' \right)'
$$
is measurable. 
This is the case, since the von Neumann algebra $(\pi_x(G)'' \cap \pi_x(G)')'$
is generated by $\pi_x(G)' \cup \pi_x(G)''$ for every $x\in X$
and the fields $x \mapsto \pi_x(G)'$ and $x \mapsto \pi_x(G)''$ are measurable. 
It follows that 
$$
\pi(G)'' \cap \pi(G)' \, = \, \int^\oplus_{X} (\pi_x(G)'' \cap \pi_x(G)') d\mu(x). 
$$
Since $\pi(G)'' \cap \pi(G)'$ is the center $\mathcal Z$ of $\pi(G)''$
and since $\mathcal Z$ coincides
with the algebra of diagonalisable operators on $\Hi$,
it follows that $\pi_x(G)'' \cap \pi_x(G)'' = \C{\rm I}_{\Hi_x}$
for $\mu$-almost every $x \in X$.
So, property \ref{iiiDEthmDirectIntFact0} holds.
\par

To show that property \ref{ivDEthmDirectIntFact0} holds,
let $T = \int^\oplus_{X} T_x d\mu(x) \in \pi(G)''$.
Since $T$ is a strong limit of a sequence from $\pi(C^c(G))$,
there exists a set $N(T)$ with $\mu(N(T)) = 0$
and a sequence $(f_n)_{n \ge 1}$ in $C^c(G)$
such that $(\pi_x(f_n))_{n \ge 1}$ converges strongly to $T_x$
for every $x \in X \smallsetminus N(T)$
(see Chap. II, \S~2, Proposition 4 in \cite{Dixm--vN}).
Let $x, y \in X \smallsetminus N(T)$ and $R \in \Hom_G(\pi_x, \pi_y)$.
Then 
$$
R \pi_x(f_n) \, = \, \pi_y(f_n) R
\hskip.5cm \text{for all } \hskip.2cm
n \ge 1
$$
and it follows that $R T_x = T_y R$.
\par
 
Let now $(X_n)_{n \ge 1}$ be a point separating sequence of Borel subsets of $X$.
For every $n \ge 1$,
let $P^{(n)} \in \mathcal Z$ be the diagonalizable operator
defined by $\Un_{X_n}$, that is, 
$$
P^{(n)} \, = \, \int^\oplus_{X} \Un_{X_n}(x) {\rm I}_{\Hi_x} d\mu(x).
$$
For $N := \bigcup_{n \ge 1} N (P^{(n)})$, we have $\mu(N) = 0$,
where $N(P^{(n)})$ is defined as above. 
\par

Let $x, y \in X \smallsetminus N$ with $x \ne y$
and let $R \in \Hom_G(\pi_x, \pi_y)$.
Then, by what was seen above, we have
$$
R P^{(n)}_x \, = \, P^{(n)}_y R
\hskip.5cm \text{for every } \hskip.2cm
n \ge 1.
$$
Since $(X_n)_{n \ge 1}$ is point separating, there exists $n\ge 1$ such that 
$\Un_{X_n}(x) = 1$ and $\Un_{X_n}(y) = 0$
and hence $P^{(n)}_x(x) = {\rm I}_{\Hi_x}$
and $P^{(n)}_x(y) = 0$.
This implies that $R = 0$, so that $\pi_x$ and $\pi_y$ are disjoint.
\end{proof}

Let $\pi$ be a representation of the second countable LC group $G$
and $\int^\oplus_{X} \pi_x(g) d\mu(x)$
a direct integral over a standard Borel space $X$ 
as in Theorem~\ref{thmDirectIntFact0}.
Denoting by $[\pi_x]$ the quasi-equivalence class of $\pi_x$,
the map $\Phi \,\colon X \to \QD(G), \hskip.2cm x \mapsto [\pi_x]$ 
is injective;
it can be seen that $\Phi$ is a Borel isomorphism
between $X$ and a Borel subset of $\QD(G)$.
In this way, we obtain a direct integral decomposition of $\pi$ over $\QD(G)$,
as in the following theorem;
for the details of the proof, we refer to \cite[8.4.2, 18.7.6]{Dixm--C*}. 

\vskip.2cm

For the notion of a standard measure on a Borel space which appears
in the following result, see Appendix~\ref{AppMeasureB}.

\begin{theorem}
% 6.C.7
\label{thmDirectIntFact}
Let $G$ be a second-countable locally compact group
and $\pi$ a representation of $G$ in a separable Hilbert space $\Hi$.
Let $\mathcal Z$ be the centre of the von Neumann algebra generated by $\pi(G)$. 
\index{Decomposition of representations! $3$@into factor representations} 
\par

There exist a standard measure $\mu$ on $\QD(G)$,
a measurable field of Hilbert spaces 
$s \to \Hi_s$ over $\QD(G)$,
a measurable field of factor representation 
$s \to (\rho_s, \Hi_s)$ over $\QD(G)$
such that the quasi-equivalence class of $\rho_s$ is $s$ for all $s \in \QD(G)$,
and an isomorphism of Hilbert spaces
$U \, \colon \Hi \to \int^\oplus_{\QD(G)} \Hi_s d\mu (s)$
with the following properties:
\begin{enumerate}[label=(\arabic*)]
\item\label{1DEthmDirectIntFact}
% $U \pi(g) U^{-1} = \mathlarger{\int}^\oplus_{\QD(G)} \rho_s(g) d\mu(s)$
$U \pi(g) U^{-1} = \int^\oplus_{\QD(G)} \rho_s(g) d\mu(s)$
for all $g \in G$;
\item\label{2DEthmDirectIntFact}
$U \mathcal Z U^{-1}$ coincides with the von Neumann algebra 
of diagonalisable operators on 
% $\mathlarger{\int}^\oplus_{\QD(G)} \Hi_s d\mu(s)$.
$\int^\oplus_{\QD(G)} \Hi_s d\mu(s)$.
\end{enumerate}
Moreover, this decomposition is unique in the following sense:
if $\pi$ is equivalent to a direct integral of factor representations
${\int}^\oplus_{\QD(G)} \rho_s d\nu(s)$ as above
with properties 
\ref{1DEthmDirectIntFact} and \ref{2DEthmDirectIntFact}, 
then $\mu$ and $\nu$ are equivalent,
and (after appropriate modifications on subsets of measure zero)
$\rho_s$ and $\sigma_s$ are equivalent 
for all $s \in \QD(G)$.
\end{theorem}

\begin{defn}
% 6.C.8
\label{Def-CentralDec}
Let $G$ be a second-countable locally compact group
and $\pi$ a representation of $G$ in a separable Hilbert space $\Hi$.
A decomposition 
$$
\pi \, \simeq \, \int^\oplus_{\QD(G)} \rho_s d\mu(s)
$$
as in Theorem~\ref{thmDirectIntFact}
is called a \textbf{central direct integral decomposition} of $\pi$
and $\mu$ is the \textbf{measure associated} to this decomposition.
The class of the measure $\mu$ on $\QD(G)$ is uniquely determined by $\pi$.
\index{Central direct integral decomposition}
\end{defn}

The relation of subordination of representations
admits a neat characterization in terms of measures on $\QD(G)$;
see \cite[8.4.4, 8.4.5]{Dixm--C*}. 

\begin{prop}
% 6.C.9
\label{Pro-CentralDecQE}
Let $G$ be a second-countable locally compact group.
Let $\pi$ and $\rho$ be representations of $G$ in separable Hilbert spaces.
Let $\mu$ and $\nu$ be the measures on $\QD(G)$
associated to central direct integral decompositions of $\pi$ and $\rho$.
The following properties are equivalent:
\begin{enumerate}[label=(\roman*)]
\item\label{iDEPro-CentralDecQE}
$\pi$ is quasi-equivalent to a subrepresentation of $\rho$;
\item\label{iiDEPro-CentralDecQE}
$\mu$ is absolutely continuous with respect to $\nu$.
\end{enumerate}
\end{prop}

The following corollary is an immediate consequence of Proposition~\ref{Pro-CentralDecQE}.

\begin{cor}
% 6.C.10
\label{Cor-CentralDecQE}
Let $G$, $\pi, \rho$ and $\mu, \nu$ be as in Proposition~\ref{Pro-CentralDecQE}.
\begin{enumerate}[label=(\arabic*)]
\item\label{1DECor-CentralDecQE}
The representations $\pi$ are $\rho$ are quasi-equivalent if and only if
the measures $\mu$ and $\nu$ are equivalent.
\item\label{2DECor-CentralDecQE}
The representations $\pi$ are $\rho$ are disjoint if and only if
the measures $\mu$ and $\nu$ are mutually singular.
\end{enumerate}
\end{cor}

\section
{Groups of type I}
% Section 6.D
\label{SectionTypeI}

We now ready to introduce the central notion of groups of type I.

\begin{defn}
% 6.D.1
\label{Def-TypeIGroup}
A topological group $G$ is \textbf{type I} if all its representations are type I.
\index{Type I! $3$@group}
\end{defn}

The following result,
which is a consequence of Proposition \ref{Prop-DirectSumTypeIRep},
implies that a separable group is of type I
if all its representations in separable Hilbert spaces are type I.

\begin{prop}
% 6.D.2
\label{Cor-Prop-DirectSumTypeIRep}
Let $G$ be a topological group
and $X$ a dense subset of $G$ of cardinal $\aleph$.
Assume that every representation of $G$
in a Hilbert space of dimension at most $\aleph$ is of type I.
\par

Then $G$ is of type I.
\end{prop}

\begin{proof}
Let $\pi$ be a representation of $G$ in a Hilbert space $\Hi$.
Then we can decompose $\pi$ as a direct sum $\bigoplus_{i \in I} \pi_i$
of subrepresentations $\pi_i$ in Hilbert spaces of dimension at most $\aleph$
(see the proof of Proposition~\ref{Pro-IrrRepBoundedCard} or that of Corollary~\ref{FacSeparable}).
Therefore the claim follows from Proposition~\ref{Prop-DirectSumTypeIRep}.
\end{proof}

\begin{rem}
% 6.D.3
\label{aaaa}
Let $G$ be a topological group of type I.
The natural map $\kappa^{\rm d}_{\rm qd}$ of Definition \ref{Def-TypIPartQD}
is a bijection from the dual $\widehat G$ onto the quasidual $\QD (G)$.
\end{rem}

\subsection
{Factor representations of type I groups}
% subsection 6.D.a
\label{SS:FactorTypeI}

The following theorem shows that in Definition~\ref{Def-TypeIGroup} it suffices to require
all factor representations of $G$ to be of type I, at least when $G$ is LC.
Recall that a factor representation is of type I
if and only if it is a multiple of an irreducible representation (Proposition~\ref{Prop-TypeIFac}).

\begin{theorem}
% 6.D.4
\label{Theo-TypIFactorRep}
Let $G$ be a second-countable locally compact group.
The following properties are equivalent:
\begin{enumerate}[label=(\roman*)]
\item\label{iDETheo-TypIFactorRep}
$G$ is of type I;
\item\label{iiDETheo-TypIFactorRep}
every factor representation of $G$ is of type I.
\end{enumerate}
\end{theorem}

\begin{proof}
We only have to show that \ref{iiDETheo-TypIFactorRep}
implies \ref{iDETheo-TypIFactorRep}. 
\par

Assume that every factor representation of $G$ is of type I.
Let $\pi$ be a representation of $G$; we have to show that $\pi$ is of type I.
Since $G$ is second-countable, and therefore separable,
it suffices by Proposition~\ref{Cor-Prop-DirectSumTypeIRep}
to consider the case of a representation $\pi$ in a separable Hilbert space $\Hi$.
\par

By Theorem~\ref{thmDirectIntFact0},
we may assume that $\Hi = \int_X^\oplus \Hi_x d\mu(x)$
and that $\pi = \int_X^\oplus \pi_x d\mu(x)$
for a measurable field $x \mapsto \Hi_x$ of Hilbert spaces
and a measurable field $x \mapsto \pi_x$ of representations
over a standard Borel space $X$ with the following properties:
\begin{enumerate}[label=(\roman*)]
\item[(i)]
$\pi_x$ is factorial for $\mu$-almost every $x\in X$;
\item[(ii)]
the center $\mathcal Z$ of $\pi(G)''$
is the algebra of diagonalizable operators on $\Hi$.
\end{enumerate}
For every $p =1, 2, \dots, \infty$, let 
$$
X_p \, = \, \{x \in X \mid \dim \Hi_x = p\},
$$
and let $\Ki_p$ be a fixed Hilbert space of dimension $p$.
By Proposition~\ref{Pro-DecConstantFieldHilbertSpaces}, we have
$$
\Hi \, = \, \oplus_{p} \int_{X_p}^\oplus \Hi_x d\mu(x)
\, \approx \, \oplus_{p} L^2(X_p, \mu, \Ki_p).
$$
So, $\pi$ is equivalent to a direct sum $\oplus_{p}\pi^{(p)}$,
where $\pi^{(p)}$ is a representation in $L^2(X_p, \mu, \Ki_p)$
which is a direct integral $\int_{X_p} \pi_x^{(p)} d\mu(x)$ over $X_p$
of representations in $\Ki_p$.
By Proposition~\ref{Prop-DirectSumTypeIRep},
we may therefore assume that $\Hi = L^2(X, \mu, \Ki)$
and that $\pi$ is a direct integral $\int_{X} \pi_x d\mu(x)$ over $X$
of representations $\pi_x$ in a fixed separable Hilbert space $\Ki$, 
with the properties (i) and (ii) above.
\par

Recall from Theorem~\ref{thmDirectIntFact0} that 
$x \mapsto \pi_x(G)''$ and $x \mapsto \pi_x(G)'$
are measurable fields of von Neumann algebras and that 
$$
\pi(G)'' \, = \, \int^\oplus_{X} \pi_x(G)''d\mu(x)
\hskip.2cm \text{and} \hskip.2cm
\pi(G)' \, = \, \int^\oplus_{X} \pi_x(G)'d\mu(x).
$$
Let $(T^{(n)})_{n \ge 1}$ and $(S^{(n)})_{n \ge 1}$
be sequences respectively in $\pi(G)''$ and $\pi(G)'$
such that $(T^{(n)}_x)_{n \ge 1}$ and $(S^{(n)}_x)_{n \ge 1}$
generate $\pi_x(G)''$ and $\pi_x(G)'$ for $\mu$-almost every $x$ in $X$.
\par

We may of course assume that $T^{(n)}$ 
(respectively $S^{(n)}$) has norm $1$
and that ${T^{(n)}}^* \in \{T^{(m)} \mid m \ge 1\}$
(respectively, ${S^{(n)}}^* \in \{S^{(m)} \mid m \ge 1\}$) for every $n \ge 1$.
\par

Since $\pi_x$ is factorial,
$\pi_x$ is a multiple of an irreducible representation of $G$
(Proposition~\ref{Prop-TypeIFac}),
for $\mu$-almost every $x\in X$;
in particular, there exists a non-zero projection $P_x \in \pi_x(G)'$
such that the reduced algebra $P_x \pi_x(G)'P_x$ is abelian.
We are going to show that we can choose the $P_x$'s
so that the field $x \mapsto P_x$ is measurable.
\par

The unit ball $\Li_1$ of $\Li(\Ki)$, equipped with the weak operator topology,
is metrizable and compact and hence a standard Borel space.
Let $X_1$ be a Borel subset of $X$ with $\mu(X \smallsetminus X_1) = 0$ 
and such that $\pi_x$ is factorial for every $x \in X_1$.
\par

Let $\Omega$ be the set of pairs $(x, S) \in X_1 \times \Li_1$
which satisfy the following conditions:
% with the following properties:
\begin{enumerate}
\item\label{1DETheo-TypIFactorRep}
$ST^{(n)}_x = T^{(n)}_xS$ for every $n \ge 1$,
\item\label{2DETheo-TypIFactorRep}
$S^2 = S = S^*$ (that is, $S$ is a projection),
\item\label{3DETheo-TypIFactorRep}
$(SS^{(n)}_x S) (SS^{(m)}_x S) = (SS^{(m)}_x S) (SS^{(n)}_xS)$
for all $m, n \ge 1$,
\item\label{4DETheo-TypIFactorRep}
$S \ne 0$.
\end{enumerate}
It is obvious that every one of these properties
defines a Borel subset of $ X_1 \times \Li_1$;
hence $\Omega$ is a Borel subset of $X_1\times \Li_1$.
Moreover, for every $x\in X_1$, the set 
$$
\{ S \in \Li_1 \mid (x, S) \in \Omega\}
$$
is not empty, by what we have seen above.
It follows from von Neumann selection theorem
(Theorem~\ref{Theo-vNSelection})
that there exists a Borel subset $X_2$ of $X_1$
with $\mu (X \smallsetminus X_2) = 0$
and a Borel map 
$$
X_2 \to \Li_1, \hskip.2cm x \mapsto P_x
$$
such that $(x, P_x) \in \Omega$ for every $x \in X_2$. 
\par

Let $x\in X_2$. Conditions
(\ref{1DETheo-TypIFactorRep}), (\ref{2DETheo-TypIFactorRep}),
(\ref{3DETheo-TypIFactorRep}, (\ref{3DETheo-TypIFactorRep})
imply, respectively,
that $P_x$ belongs to $\pi_x(G)'$,
that $P_x$ is a projection, 
that $P_x \pi_x(G)' P_x$ is abelian,
and that $P_x \ne 0$. 
\par

Set $P_x = 0$ for $x \in X \smallsetminus X_2$ and let 
$$
P \, := \, \int^\oplus_{X} P_xd\mu(x).
$$
Then $P$ is a projection in $\pi(G)'$ such that $P\pi(G)' P$ is abelian.
\par
 
We claim that the central support of $P$ is ${\rm Id}_\Hi$. 
Indeed, let $E$ be a central projection in $\pi(G)''$ with $P \le E$.
Since the center of $\pi(G)''$ is the algebra of diagonalizable operators on $\Hi$,
we have $E = \int_Y \Un_x {\rm Id}_\Ki$
for a Borel subset $Y$ of $X$.
As $P_x \ne 0$ for every $x \in X_2$ and 
$\mu(X \smallsetminus X_2) = 0$,
it follows that $\mu(X \smallsetminus Y) = 0$, that is, $E = {\rm Id}_\Hi$.
\par 

So, the induction map $\pi(G)'' \to P \pi(G)'' P$
is an isomorphism of von Neumann algebras (see Proposition~\ref{inductionisiso}).
The commutant of $P\pi(G)'' P$ coincides with $P\pi(G)' P$
(see Proposition~\ref{Prop-MatrixCommutant});
since $P\pi(G)' P$ is abelian,
it follows from Proposition~\ref{Pro-TypeIVN} that $\pi$ is of type I.
\end{proof}

\subsection
[Decomposition into irreducible representations]
{Decomposition into irreducible representations
\\
for groups of type I}
% subsection 6.D.b
\label{SectionDecomposingRepTypeI}

Recall that every representation of second-countable LC group on a separable Hilbert space
can be written, in general in several different ways, as a direct integral of irreducible representations
(Theorem~\ref{Theo-NonUniqueIntDecIrrUniRep}).
\par 
 
Let $G$ be a second-countable LC group of type I.
It follows from Theorem~\ref{Theo-DecTypeI-MultiplicityFree}
that every representation of $G$ is, in a canonical way,
a direct sum of multiplicity-free representations.
The following result gives a description of the multiplicity-free representations of $G$.

\begin{prop}
% 6.D.5
\label{Pro-Multiplity-FreeTypeI}
Let $G$ be a second-countable locally compact group of type I
and $\pi$ a multiplicity-free representation of $G$ in a separable Hilbert space.
\par

There exists a positive Borel measure $\mu$ on 
$\widehat G$ such that $\pi$ is equivalent to the direct integral 
${\int}^\oplus_{\widehat G} \sigma d\mu(\sigma)$.
Moreover, the measure $\mu$ is unique up to equivalence.
\end{prop}

\begin{proof}
Using the central decomposition of $\pi$ as in Theorem~\ref{thmDirectIntFact},
we can assume that
$$
\pi \, = \, {\int}^\oplus_{\QD(G)} \pi_x d\mu(x)
\, = \, {\int}^\oplus_{\widehat G} \pi_xd\mu(x).
$$
Since $G$ is type I, we have $\pi_x = n(x) \sigma_x$
for an integer $n(x)$ and an irreducible representation $\sigma_x$.
Let $\mathcal A$ be the subalgebra of $\Li (\Hi)$ 
consisting of the decomposable operators
of the form $T = {\int}^\oplus_{\widehat G} T_xd\mu(x)$,
where $T_x \in \pi_x(G)'$ for every $x$.
It is clear that $\mathcal A$ is contained in $\pi(G)'$.
Since $\pi$ is multiplicity-free, $\pi(G)'$ and hence $\mathcal A$ is abelian.
It follows that $\pi_x(G)'$ is abelian and hence that $n(x) = 1$ for $\mu$-almost
every $x$. So, $\pi_x$ is irreducible for $\mu$-almost every $x$.
\par

The uniqueness of the class of $\mu$ follows from Corollary~\ref{Cor-CentralDecQE}
\end{proof}

\begin{rem}
% 6.D.6
\label{Rem-Pro-Multiplity-FreeTypeI}
(1)
The converse statement of Proposition~\ref{Pro-Multiplity-FreeTypeI} is also true:
when $G$ be a a second-countable LC group of type I, every representation of the form 
${\int}^\oplus_{\widehat G} \pi_x d\mu(x)$ for a positive Borel measure $\mu$
on $\widehat G$ is multiplicity-free (see \cite[8.6.4]{Dixm--C*});
so, the map 
$$
\mu \mapsto {\int}^\oplus_{\widehat G} \pi_x d\mu(x)
$$
induces a bijection between equivalence classes of Borel measures on $\widehat G$ 
and equivalence classes of multiplicity-free representations of $G$.

\vskip.2cm

(2)
Even when $G$ is not of type I, recall (Theorem \ref{Thm-DecDirectIntegral-MultiplicityFree})
that a multiplicity-free representation of $G$ admits an essentially unique
integral decomposition into irreducible representations over a standard Borel space.
The proof of Proposition~\ref{Pro-Multiplity-FreeTypeI} could have been based on this fact.
\end{rem}

The following result is a refinement of Theorem \ref{thmDirectIntIrreps} 
for LC groups of type I.

\begin{theorem}
[\textbf{Canonical decomposition into irreducible representations}]
% 6.D.7
\label{thmDirectIntIrreps+}
Let $G$ be a second-countable locally compact group
and $\pi$ a representation of $G$
in a separable Hilbert space $\Hi$.
Assume that $G$ is of type I. 
\index{Decomposition of representations! 
$2$@into irreducible representations for type I groups}
\par

% \marginpar{Mieux $\infty$ ou $\aleph_0$ ?} tel du 14 oct : $\infty$ mieux.
There exist a set $I \subset \overline{\N^*}$ of extended positive integers
and a sequence $(\mu_m)_{m \in I}$ of positive measures on $\widehat G$
with the following properties:
\begin{enumerate}[label=(\arabic*)]
\item\label{1DEthmDirectIntIrreps+}
the measures $\mu_m$ are singular with each other;
\item\label{2DEthmDirectIntIrreps+}
the representation $\pi$ is equivalent to
$$
\bigoplus_{m \in I} m \hskip.1cm
\int_{\widehat G}^\oplus \sigma_x d\mu_m(x) .
$$
\end{enumerate}
Moreover, this decomposition is unique in the following sense:
if $\pi$ is equivalent to 
$\bigoplus_{n \in J} n \hskip.1cm 
\int_{\widehat G}^\oplus \sigma_x d\nu_n(x)$
% \mathlarger{\int}_{\widehat G}^\oplus \sigma_x d\nu_n(x)$
for a set $J$ of extended integers
and for mutually singular measures $\nu_n$ on $\widehat G$,
then $I = J$ and $\mu_m$ and $\nu_m$ are equivalent for all $m \in I$.
\end{theorem}

\begin{proof}
By Theorem~\ref{Theo-DecTypeI-MultiplicityFree},
there exist a set $I \subset \overline{\N^*}$ of extended positive integers
and an orthogonal decomposition $\pi = \bigoplus_{m\in I} \pi^{(m)}$ with the following properties:
\begin{itemize}
\setlength\itemsep{0em}
\item
for every $m\in I$, we have $\pi^{(m)} = m\sigma^{(m)}$
for a multiplicity-free representation $\sigma^{(m)}$ of $G$;
\item
$\pi^{(m)}$ and $\pi^{(n)}$ are disjoint for $m \ne n$.
\end{itemize}
By Proposition~\ref{Pro-Multiplity-FreeTypeI}, for every $m \in I$, 
there exists a positive Borel measure $\mu$ on 
$\widehat G$ such that $\sigma^{(m)}$ is equivalent to the direct integral 
${\int}^\oplus_{\widehat G} \sigma_x d\mu_m(x)$.
Therefore $\pi$ is equivalent to 
$$
\bigoplus_{m \in I} m \hskip.1cm
\int_{\widehat G}^\oplus \sigma_x d\mu_m(x) .
$$
Moreover, $\mu_m$ and $\mu_n$ are disjoint for $m \ne n$ by Corollary~\ref{Cor-CentralDecQE}.
\par

The uniqueness statement follows from
the uniqueness of the decomposition $\pi = \bigoplus_{m\in I} \pi^{(m)}$
of Theorem~\ref{Theo-DecTypeI-MultiplicityFree}.
\end{proof}

\begin{rem}
% 6.D.8
\label{Rem-thmDirectIntIrreps+}

(1)
In the special case of an abelian second-countable LC group,
Theorem~\ref{thmDirectIntIrreps+}
boils down to Theorem~\ref{Thm-DecRepAbelianGroups}.
%For type I representations of arbitrary LC groups, 
%Theorem \ref{thmDirectIntIrreps} can also be refined, but not to the same extent;
%we refer to \cite[\S~5.4]{Dixm--C*}.

\vskip.2cm

(2)
Theorem~\ref{thmDirectIntIrreps+} fails for \emph{every} LC group which is not of type I,
as noted in Theorem~\ref{Theo-NonUniqueIntDecIrrUniRep}.
\end{rem}

\section
{A class of groups of type I}
% Section 6.E
\label{S:ClassTypeI-GCRGroups}

We introduce the important class of the so-called GCR groups.
This is a class of type I locally compact groups; in fact, it will turn out 
later (Section~\ref{SectionGlimm}) that every second-countable LC group 
of type I is a GCR group.

\subsection
{On the algebra of compact operators}
% subsection 6.E.a
\label{SS:CompactOperators}

Let $\Hi$ a Hilbert space.
The subset $\Ki(\Hi)$ of $\Li (\Hi)$ consisting of compact operators
is a selfadjoint norm-closed two-sided ideal of $\Li (\Hi)$ and in particular a C*-algebra
(for more details on C*-algebras and their representations, see Section~\ref{SectionPrimC*}).
Observe that the identity representation $\mathrm{Id} \, \colon \Ki(\Hi) \to \Li (\Hi)$ is irreducible.
\par

We recall a few standard facts about the representation theory of $\Ki(\Hi)$;
for a proof, we refer to \cite[4.1]{Dixm--C*}.
 
\begin{prop}
% 6.E.1
\label{Pro-FactsAboutCompactOperators}
Consider the C*-algebra $\Ki(\Hi)$ of compact operators on a Hilbert space $\Hi$.
\begin{enumerate}[label=(\arabic*)]
\item\label{1DEPro-FactsAboutCompactOperators}
Every non-degenerate representation of $\Ki(\Hi)$
is equivalent to a multiple of the identity representation $\mathrm{Id} \, \colon \Ki(\Hi) \to \Li (\Hi)$;
\item\label{2DEPro-FactsAboutCompactOperators}
$\Ki(\Hi)$ is simple: the only closed two sided ideals of $\Ki(\Hi)$ are $\{0\}$ and $\Ki(\Hi)$;
\item\label{3DEPro-FactsAboutCompactOperators}
let $A$ be a C*-subalgebra of $\Ki(\Hi)$ such that $\mathrm{Id} \, \colon A \to \Li (\Hi)$ is irreducible;
then $A = \Ki(\Hi)$.
\end{enumerate}
\end{prop}

For the definition of ``non-degenerate'', see Section \ref{SectionPrimC*}.

\subsection
{GCR groups}
% 6.E.b
\label{SS:GCR groups}

Let $G$ be locally compact group, with left Haar measure $\mu_G$.
Let $(\pi, \Hi)$ be a representation of $G$.
Denote again by $\pi$ the associated representation $C^c(G) \to \Li (\Hi)$
of the convolution algebra; this has already appeared
just after Definition \ref{defquasiregG/K},
and will appear again in Section~\ref{C*algLCgroup} below.

\begin{defn}
% 6.E.2
\label{Def-GCR-Group}
A locally compact group $G$ is a \textbf{GCR group}, or \textbf{postliminal group},
if, for every irreducible representation $(\pi, \Hi)$ of $G$,
the norm closure of $\pi(C^c(G))$ in $\Li (\Hi)$
contains a \emph{non-zero} compact operator.
\index{GCR! $1$@group}
\end{defn} 

\begin{rem}
% 6.E.3
\label{Rem-GCR-Group}
Let $G$ be a locally compact group.

\vskip.2cm

(1)
Assume that $G$ is a GCR group and let $(\pi, \Hi)$ be an irreducible representation of $G$. 
Then $\Ki(\Hi)$ is contained in the norm closure of $\pi(C^c(G))$.
\par

Indeed, denote by $A$ the norm closure of $\pi(C^c(G))$ and set $J = A \cap \Ki(\Hi)$.
Then $J$ is a closed two-sided of the C*-subalgebra of $A$ of $\Li (\Hi)$,
and $J \ne \{0\}$ because $G$ is a GCR group.
As $\mathrm{Id} \, \colon A \to \Li (\Hi)$ is irreducible,
it follows easily that $\mathrm{Id} \, \colon J \to \Li (\Hi)$ is irreducible \cite[2.11.3]{Dixm--C*}.
Therefore we have $J = \Ki(\Hi)$, by Proposition~\ref{Pro-FactsAboutCompactOperators}.

\vskip.2cm

(2)
We introduce below the class of GCR C*-algebras.
% (Definition~\ref{Def-GCR-C*}).
It turns out that $G$ is a GCR group
if and only if its maximal C*-algebra as defined in Section~\ref{C*algLCgroup} is GCR,
at least when $G$ is second-countable;
see Section~\ref{SectionGlimm}.

\vskip.2cm

(3)
We will define in Section~\ref{Section-NormRepChar}
the notion of a \textbf{normal} representation of $G$. 
In view of (1) above, $G$ is a GCR group if and only if every irreducible
representation of $G$ is normal (Proposition~\ref{Pro-IrredNormalRep}).
\end{rem}

We will need at several places the following characterization of weak containment 
of representations of $G$ in terms of operator norms.
It will be discussed later:
see Proposition \ref{tradweqGC*} and Remark~\ref{Rem-tradweqGC*}.

\begin{prop}
% 6.E.4
\label{Pro-WeakContainmentOperatorNorm}
Let $G$ be a locally compact group and $\pi_1, \pi_2$ two representations of $G$.
The following properties are equivalent:
\begin{enumerate}[label=(\roman*)]
\item\label{iDEPro-WeakContainmentOperatorNorm}
$\pi_1$ is weakly contained in $\pi_2$;
\item\label{iiDEPro-WeakContainmentOperatorNorm}
for every $f \in C^c(G)$, we have $\Vert \pi_1(f) \Vert \le \Vert \pi_2(f) \Vert$. 
\end{enumerate}
\end{prop}

The importance of the class of GCR groups lies in the following result.

\begin{theorem}
% 6.E.5
\label{Theo-GCR-Group}
Let $G$ be a second-countable locally compact GCR group.
\par

Then $G$ is of type I.
\end{theorem} 

\begin{proof} 
In view of Theorem~\ref{Theo-TypIFactorRep},
it suffices to show that every factor representation of $G$ is of type I.
\par

Let $(\pi, \Hi_\pi)$ be a factor representation of $G$ on a Hilbert space $\Hi_\pi$.
Since $G$ is second-countable,
$\pi$ is weakly equivalent to an \emph{irreducible} representation $(\sigma, \Hi_\sigma)$ of $G$,
as we will show later (Corollary~\ref{Cor-FacRepPrimitiveG}). 
Let $A$ be the norm closure of $\sigma(C^c(G))$ in $\Li (\Hi_\sigma)$.
Since $G$ is a GCR group, we have 
$$
\Ki(\Hi_\sigma) \subset A,
$$ 
by Remark~\ref{Rem-GCR-Group}(1).
As $\pi$ is weakly equivalent to $\sigma$, we have
$$
\Vert \pi(f) \Vert \, = \, \Vert \sigma(f) \Vert
\hskip.2cm \text{for every} \hskip.2cm
f \in C^c(G) 
$$
(see Proposition~\ref{Pro-WeakContainmentOperatorNorm}).
This implies that the map 
$$
\sigma(C^c(G)) \to \pi(C^c(G)), \hskip.2cm \sigma(f) \mapsto \pi(f)
$$
is well-defined and extends to an isometric $*$-homomorphism
$\Phi \, \colon A \to \Li (\Hi_\pi)$,
that is, to a faithful representation $\Phi$ of $A$ on $\Hi_\pi$. 
\par
 
By Proposition~\ref{Pro-FactsAboutCompactOperators}, the restriction of $\Phi$ 
to $\Ki(\Hi_\sigma)$ is a multiple of the identity representation of $\Ki(\Hi_\sigma)$ on $\Hi_\sigma$;
this means that there exist a cardinal $n \ge 1$
and a Hilbert space isomorphism $U \, \colon \Hi_\pi \to \Hi_\sigma \otimes \Hi$
such that
$$
U \Phi(T) U^{-1} \, = \, T \otimes \mathrm{Id}_{\Hi}
\hskip.2cm \text{for every} \hskip.2cm
T \in \Ki(\Hi_\sigma),
$$
where $\Hi$ is a Hilbert space of dimension $n$. 
\par
 
Observe that $\Ki(\Hi_\sigma)$ is dense in $ \Li (\Hi_\sigma)$
and hence in $A$, for the strong operator topology.
It follows that we have
$$
U \Phi(T) U^{-1} \, = \, T \otimes \mathrm{Id}_{\Hi}
\hskip.2cm \text{for every} \hskip.2cm
T \in A,
$$
and therefore
$$
U \pi(f) U^{-1} \, = \, \sigma(f) \otimes \mathrm{Id}_{\Hi} 
\hskip.2cm \text{for every} \hskip.2cm
f \in C^c(G).
$$
Let $\mu_G$ be a Haar measure on $G$
and $(f_n)_n$ be an approximate identity for $L^1(G, \mu_G)$
consisting of functions $f_n \in C^c(G)$
(see Section \ref{AppAlgC*}).
Then, for every $g \in G$, we have
$$
\lim_n \pi({}_g f_n) \, = \, \pi(g)
\hskip.5cm \text{and} \hskip.5cm
\lim_n \sigma( {}_g f_n ) \, = \, \sigma(g)
$$
in the strong operator topology,
where ${}_gf_n$ is defined by ${}_gf_n(x) = f_n (g^{-1}x)$.
It follows that 
$$
U \pi(g) U^{-1} \, = \, \sigma(g) \otimes \mathrm{Id}_{\Hi}
\hskip.2cm \text{for every} \hskip.2cm
g \in G.
$$
This shows that $\pi$ is equivalent to $n\sigma$. Therefore $\pi$ is of type I.
\end{proof}

\begin{rem}
% 6.E.6
\label{Rem-Theo-GCR-Group}
Theorem~\ref{Theo-GCR-Group} admits a converse:
\begin{itemize}
\setlength\itemsep{0em}
\item
every second-countable locally compact of type I is a GCR group.
\end{itemize}
This is a deep result which is valid more generally in the context of separable C*-algebras
and which is part of Glimm theorem, to be discussed in Section~\ref{SectionGlimm}
\end{rem}

\subsection
{CCR groups}
% subsection 6.E.c
\label{CCR-groups}

We introduce an important subclass of the class of GCR groups.

\begin{defn}
% 6.E.7
\label{Def-CCR-Group}
Let $G$ be a locally compact group.
\par

A representation $(\pi, \Hi)$ of $G$ is said to be a \textbf{CCR representation}
if the algebra $\pi(C^c(G))$ is contained in $\Ki(\Hi)$.
\index{Liminal! $1$@representation}
\index{CCR! $1$@representation}
\index{Representation! CCR}
\par

The group $G$ is a \textbf{CCR group}, or \textbf{liminal group},
if every irreducible representation of $G$ is a CCR representation.
\index{CCR! $2$@group}
\index{Liminal! $1$@locally compact group}
\end{defn} 

Concerning the meaning of the acronyms GCR and CCR, see Remark~\ref{Rem-Def-GCR-C*}.

\begin{rem}
% 6.E.8
\label{Rem-CCR-Rep}
Let $G$ be a locally compact group.

\vskip.2cm

(1)
Every finite dimensional representation of $G$ is CCR.

\vskip.2cm

(2)
If $G$ is a discrete group and $\pi$ is an infinite dimensional representation of $G$,
then $\pi$ is not CCR: indeed, $\pi(e) = I$ is not a compact operator.

\vskip.2cm

(3)
If $G$ is a CCR group, then obviously $G$ is a GCR group.

\vskip.2cm

(4)
Let $(\pi, \Hi)$ be an \emph{irreducible} CCR representation of $G$.
Then $\pi(C^c(G))$ is dense in $\Ki(\Hi)$ in the operator norm.
This follows from Remark~\ref{Rem-GCR-Group}(1).

\vskip.2cm

(5)
Let $(\pi, \Hi)$ be an \emph{irreducible} CCR representation of $G$.
Then $\{\pi\}$ is a closed point of the dual space $\widehat G$,
equipped with the Fell topology as defined in \ref{S-DefUnitD}
(see Corollary~\ref{Cor-RepCCR} below).
 \end{rem}

As the next example shows, the class of GCR groups is strictly larger than the class of CCR groups.

\begin{exe}
% 6.E.9
\label{Exa-GCR-NonCCR}
Let $G = \Aff(\R)$ be the affine group over the real line.
We view $G$ as the space $\R^\times \times \R$,
equipped with the group law
given by $(a,b) (a',b') = (aa', ab'+b)$,
and with the natural topology.
It is a second-countable locally compact group,
which has a left Haar measure $da db / a^2$.
Recall from Remark~\ref{AffKLie} that, up to equivalence,
$G$ has a unique infinite dimensional irreducible representation $\pi$
defined on $L^2(\R^\times, dt / \vert t \vert)$ by 
$$
\pi(a,b) \xi(t) \, = \, e^{- 2 \pi i bt} \xi(at)
\hskip.5cm \text{for} \hskip.2cm
(a, b) \in G, \hskip.1cm \xi \in L^2(\R^\times, dt / \vert t \vert), \hskip.1cm t \in \R^\times.
$$
Let $f \in C^c(G)$. It can be checked (see \cite[Corollaire p. 164]{Khal--74})
that $\pi(f)$ is a compact operator if and only if
$$
\int_{\R} f(a,b) db \, = \, 0
\hskip.2cm \text{for every} \hskip.2cm
a \in \R^\times.
$$
In particular, $G$ is a GCR group but not a CCR group.
Since $\{\pi\}$ is not closed in $\widehat G$ as mentioned in Remark~\ref{AffKLie},
the fact that $G$ is not a CCR group
is also an immediate consequence of Corollary~\ref{Cor-RepCCR} below.
\end{exe}

Before we proceed with examples of CCR groups,
we recall the following facts about the representation theory of compact groups.
Let $K$ be a compact group. 
\begin{itemize}
\setlength\itemsep{0em}
\item
Every irreducible representation of $K$ is finite dimensional.
\item
For $\sigma\in \widehat{K}$, let $\chi_\sigma$ be the character
of $\sigma$, that is, the continuous function 
$$
K \to \C, \hskip.2cm x \mapsto \mathrm{Tr}(\sigma(x)).
$$
For $\sigma, \sigma' \in \widehat{K}$,
we have the following ``orthogonality relations":
$$
\chi_\sigma\ast \chi_{\sigma'} \, = \, 0
\hskip.5cm \text{if} \hskip.2cm 
\sigma \ne \sigma'
\hskip.5cm \text{and} \hskip.5cm
\chi_\sigma\ast \chi_{\sigma} \, = \, d_\sigma^{-1}\chi_{\sigma} ,
$$
where $\star$ denotes the convolution product
and $d_\sigma$ the dimension of $\sigma$.
\item
Every representation $\pi$ of $K$ decomposes as direct sum
$\pi = \bigoplus_{\sigma \in \widehat{K}} n_\sigma \sigma$ ;
more precisely, if $(\pi, \Hi)$ is a representation of $K$,
we have a decomposition 
$$
\Hi \, = \, \bigoplus_{\sigma \in \widehat{K}} \Hi_\sigma
$$
into isotypical components
($\Hi_\sigma$ is the sum of irreducible subspaces on which 
$\pi$ is equivalent to $\sigma$)
and the orthogonal projection $\Hi\to \Hi_\sigma$
is the operator $\pi(d_\sigma \overline{\chi_\sigma})$.
\end{itemize}
For all this, see for instance \cite[Chap.~15]{Dixm--C*}
%\cite[Theorem (27.44)]{HeRo--70} = AbstractHarmAn2
or \cite[Chap.~5]{Foll--16}.

\vskip.2cm

Important examples of CCR groups are groups which contain large compact subgroups
in the sense of the following definition which appears in \cite[p.641]{Mack--63}
(see also \cite[Chap.~4]{Warn--72}, where a stronger notion is considered).

\begin{defn}
% 6.E.10
\label{Def-LargeCompactSubgroup}
Let $G$ be a topological group.
A compact subgroup $K$ is a \textbf{large compact subgroup}
of $G$ if every $\pi \in \widehat G$ has finite $K$-multiplicities, that is, 
in the decomposition
$\pi \vert_K = \bigoplus_{\sigma \in \widehat{K}} n_\sigma \sigma$, 
we have $n_\sigma<\infty$ for every $\sigma\in \widehat{K}$.
\index{Large compact subgroup}
\end{defn} 

\begin{prop}
% 6.E.11
\label{Pro-LargeCompactSubgroupTypeI}
Let $G$ be a locally compact group which contains a large compact subgroup $K$.
\par

Then $G$ is a CCR group. In particular, $G$ is of type I.
\end{prop}

\begin{proof}
Let $(\pi, \Hi)$ be an irreducible representation of $G$.
Let
$$
\Hi \, = \, \bigoplus_{\sigma\in \widehat{K}} \Hi_\sigma
$$
be the decomposition of $\Hi$ into $K$-isotypical components.
Since $K$ is a large compact subgroup,
the restriction of $\pi$ to $K$ has finite multiplicities
and this means that $\Hi_\sigma$ is finite dimensional for every $\sigma\in \widehat{K}$;
equivalently, the operator $\pi(d_\sigma \overline{\chi_\sigma})$,
which is the orthogonal projection on $\Hi_\sigma$,
is of finite rank for every $\sigma\in \widehat{K}$.
\par

Let $\mu_G$ be a Haar measure on $G$.
Let $f \in C^c(G)$. Viewing $f$ as a vector in $L^2(G, \mu_G)$
and considering the restriction to $K$
of the left regular representation $\lambda_G$ of $G$,
we have
$$
f \, = \, \sum_{\sigma\in \widehat{K}} \lambda_G(\overline{\chi_\sigma}) (f)
\, = \, \sum_{\sigma\in \widehat{K}} d_\sigma \overline{\chi_\sigma} \ast f,
$$
where the sum converges in $L^2(G)$.
Since the support of $\overline{\chi_\sigma} \ast f$
is contained in the compact set $K\mathrm{supp} (f)$,
we also have
$f = \sum_{\sigma\in \widehat{K}} d_\sigma \overline{\chi_\sigma} \ast f$,
with a converging sum in $L^1(G, \mu_G)$.
Therefore, we have
$$
\pi(f)
\, = \, \pi \left(\sum_{\sigma\in \widehat{K}} d_\sigma \overline{\chi_\sigma} \ast f\right)
\, = \, \sum_{\sigma\in \widehat{K}} \pi(d_\sigma \overline{\chi_\sigma}) \pi(f),
$$
where the last sum converges in the operator norm.
Since $\pi(d_\sigma \overline{\chi_\sigma})$ is of finite rank,
it follows that $\pi(f)$ is a compact operator.
\end{proof}

Important examples of groups with a large compact subgroups are reductive groups
over local fields.

\begin{theorem}
% 6.E.12
\label{Thm-ReductiveCCR groups}
 Let $G$ be a reductive linear algebraic group over a local field.
Let $K$ be a maximal compact subgroup of $G$.
\par

Then $K$ is a large compact subgroup 
and so $G$ is a CCR group.
\end{theorem}

For a proof of Theorem~\ref{Thm-ReductiveCCR groups}, 
see \cite[Theorem 2]{Gode--52} or \cite[15.5.6]{Dixm--C*} in the real case
(see also \cite[Page 230]{Hari--53}),
and \cite[Theorem 1]{Bern--74} in the non-archimedean case.
We will give below a complete proof of Theorem~\ref{Thm-ReductiveCCR groups}
in the special case $G = \SL_2(\R)$.
In the case of a linear semisimple real Lie group,
the proof follows a similar strategy (see Remark~\ref{Rem-Theo-CCRforSL2}),
but is more involved, because $K$ is not necessarily abelian.

\begin{exe}
% 6.E.13
\label{Example-CCR groups}
(1)
Another class of examples of groups containing a large compact subgroup 
are motion groups (see Example~\ref{Exa-Gel'fandPair}).
More generally, assume that the LC compact group $G$ can be written
as the product $G = KA$ of a compact subgroup $K$ and a closed abelian subgroup $A$.
\par

Then $K$ is a large compact subgroup and so $G$ is a CCR group (see \cite[Theorem 5]{Gode--52}). 

\vskip.2cm

(2)
Every connected nilpotent Lie group is a CCR group (see \cite{Dixm--59, Kiri--62}).
For a characterization of connected and simply connected
solvable Lie groups that are CCR, see Chapter V in \cite{AuMo--66}.
\end{exe}

\subsection
{$\SL_2(\R)$ is a CCR group}
% subsection 6.E.d
\label{SS:SL2R-CCR groups}

Let $G = \SL_2(\R)$.
We are going to show that the maximal compact subgroup $K = \SO(2)$
is a large subgroup of $G$,
and therefore that $G$ is a CCR group
(see Proposition~\ref{Pro-LargeCompactSubgroupTypeI}).
\par

Recall that $K = \SO(2)$ is the group of matrices of the form
$$
r(\theta) \, = \, \left(\begin{matrix}
\phantom{-}\cos \theta & \sin \theta \\ -\sin \theta & \cos \theta
\end{matrix} \right)
$$
for $ \theta \in \R$.
The dual group of the abelian group $K$ is 
$$
\widehat{K} \, = \, \{ \chi_n \mid n \in \Z\} \, \cong \, \Z,
$$
where $\chi_n$ is defined by 
$$
\chi_n(r(\theta)) \, = \, e^{i n \theta}
\hskip.5cm \text{for every} \hskip.2cm
\theta \in \R.
$$
For $n \in \Z$, set 
$$
C^c(G)_{ \chi_n}
\, = \, \chi_{-n} \ast C^c(G) \ast \chi_{-n}
\, = \, \{ \chi_{-n} \ast f\ast \chi_{-n }\mid f \in C^c(G)\}.
$$
For a function $f \, \colon G \to \C$, define $f^*$ by $f^*(x) = \overline{f(x^{-1})}$.

\begin{lem}
% 6.E.14
\label{Lem-PropertiesCG-chi}
For every $n \in \Z$, the following holds:
\begin{enumerate}[label=(\arabic*)]
\item\label{1DELem-PropertiesCG-chi}
$C^c(G)_{\chi_n}$ is a subalgebra of the convolution algebra $C^c(G)$
and is invariant under the involution $f\mapsto f^*$;
\item\label{2DELem-PropertiesCG-chi}
$C^c(G)_{\chi_n}$ is the linear subspace of $C^c(G)$ consisting of the 
functions $f \in C^c(G)$ satisfying the equation
$$
f(r(\theta_1) x r(\theta_2))
\, = \, \chi_{n}(r(\theta_1))f(x)\chi_{n}(r(\theta_2))
\, = \, e^{i n \theta_1}f(x) e^{i n \theta_2}
$$
for all $\theta_1, \theta_2 \in \R$ and $x \in G$;
\end{enumerate}
\end{lem}

\begin{proof}
Item \ref{1DELem-PropertiesCG-chi} follows from the relations 
$$
\chi_{-n}\ast \chi_{-n} \, = \, \chi_{-n}
\hskip.2cm \text{and} \hskip.2cm
\chi_{-n}^* \, = \, \chi_{-n};
$$
Item \ref{2DELem-PropertiesCG-chi} can be checked by a direct computation.
\end{proof}

Observe that, in the case $n = 0$, the algebra $C^c(G)_{\chi_n}$ 
is the algebra $C^c(K \backslash G / K)$ of $K$-invariant functions in
$C^c(G)$ considered in Subsection~\ref{SS:ExaMultFreeRep}.
\par

The Cartan decomposition of $G$
(see Example~\ref{Exa-Gel'fandPairSemissimple})
can be further refined, as follows.

\begin{lem}
% 6.E.15
\label{Lem-CartanDecSL2}
Let $A$ be the subgroup of diagonal matrices with positive diagonal entries
in $G = \SL_2(\R)$.
\par

Then $G = KAK$.
\end{lem}

\begin{proof}
Let $x \in G$. Denote by $x^t$ the transpose of the matrix $x$.
Then $xx^t$ is a symmetric positive definite matrix.
So, there exists a symmetric positive definite matrix $y$ with $y^2 = x^t x$.
Then $xy^{-1} \in K$, since
$$
(xy^{-1})^t xy^{-1} \, = \, (y^{-1})^t x^t xy^{-1} \, = \, y^{-1} y^2 y^{-1} \, = \, I.
$$
So, $x = k_1 y$ for some $k_1 \in K$.
(This is as $x = u \vert x \vert$ in Appendix \ref{AppProjValMeas}.)
Moreover, since $y$ is symmetric and positive definite,
there exists $k_2 \in K$ such that $a := k_2^{-1} y k_2$ belongs to $A$.
Therefore $x = (k_1k_2) a (k_2^{-1})$.
\end{proof}

There is a corresponding Cartan decomposition $G = KAK$
% as in Lemma~\ref{Lem-CartanDecSL2},
for any semisimple connected real Lie group $G$.
See for example \cite[Theorem 7.39]{Knap--02}.
\par

Since $(G, K)$ is a Gel'fand pair (Example~\ref{Exa-Gel'fandPairSemissimple}),
the algebra $C^c(K \backslash G /K)$ is commutative.
We now generalize this property to all algebras $C^c(G)_{ \chi_n}$.

\begin{prop}
% 6.E.16
\label{Prop-CommutativeAlgebraSpherical}
The algebra $C^c(G)_{ \chi_n}$ is commutative for every $n \in \Z$.
\end{prop}

\begin{proof}
The proof bears similarities with the proof of
Gelf'and theorem~\ref{Theo-Gel'fandPairSemissimple}).
Following the proof of Theorem~1 in \cite[Chap.~II, \S~1]{Lang--85}, we consider 
\begin{itemize}
\setlength\itemsep{0em}
\item
the anti-automorphism $\tau$ of $G$ of order 2 given by $x \mapsto x^t$ and
\item
the automorphism $\sigma$ of $G$ of order 2 given by $x \mapsto x_0 x x_0^{-1}$, where 
$$
x_0 \, = \,
\left(\begin{matrix}
1 & \phantom{-}0 \\ 0 & -1
\end{matrix} \right) .
$$
\end{itemize}
Observe that $\tau$ and $\sigma$ act both as the identity
on the subgroup $A$ of diagonal matrices.
Moreover, $\tau$ and $\sigma$ coincide on $K$:
for $k = r(\theta)\in K$, we have 
$$
\tau(k) \, = \, k^t \, = \, k^{-1}
$$
and 
$$
\sigma(k) \, = \, x_0 r(\theta) x_0^{-1} \, = \, r(-\theta) \, = \, k^{-1}.
$$
\par

For $f \in C^c(G)$, define $f^\tau \in C^c(G)$ by 
$$
f^\tau (x) \, = \, f(\tau(x))
\hskip.5cm \text{for} \hskip.2cm
x \in G ,
$$
and define similarly $f^\sigma \in C^c(G)$. 
\par

The arguments of the proof of Theorem~\ref{Theo-Gel'fandPairSemissimple}
show that $\tau$ and $\sigma$ preserve the Haar measure $\mu_G$ on $G$
(recall that $G$ is unimodular).
These arguments also show that, for every $f_1, f_2 \in C^c(G)$, we have
$$
(f_1 \ast f_2)^\tau \, = \, f_2^\tau \ast f_1^\tau
$$
(since $\tau$ is an anti-automorphism of $G$) and 
$$
(f_1 \ast f_2)^\sigma \, = \, f_1^\sigma \ast f_2^\sigma
$$
(since $\sigma$ is an automorphism of $G$).
\par 
 
Let $f \in C^c(G)_{ \chi_n}$. We claim that $f^\tau = f^\sigma$.
Indeed, let $x \in G$; write $x = k_1 a k_2$ with $k_1, k_2 \in K$ and $a\in A$
(see Lemma~\ref{Lem-CartanDecSL2}). On the one hand, we have
$$
f^\tau(x)
\, = \, f(k_2^t a^t k_1^t)
\, = \, f(k_2^{-1} a k_1^{-1})
\, = \, \chi_n(k_2) \chi_n(k_1) f(a);
$$
on the other hand, we have
$$
f^\sigma(x)
\, = \, f(\sigma(k_1) \sigma(a)\sigma(k_2)
\, = \, f(k_1^{-1} a k_2^{-1})
\, = \, \chi_n(k_1) \chi_n(k_2) f(a).
$$
\par

Now, let $f_1, f_2 \in C^c(G)_{ \chi_n}$.
Then $(f_1 \ast f_2)^\tau = (f_1 \ast f_2)^\sigma$
and hence
$$
f_2^\tau \ast f_1^\tau \, = \, f_1^\sigma \ast f_2^\sigma
$$
that is,
$$
f_2^\sigma \ast f_1^\sigma \, = \, f_1^\sigma \ast f_2^\sigma
$$
Therefore, 
$$
f_2 \ast f_1
\, = \, (f_2^\sigma \ast f_1^\sigma)^\sigma
\, = \ (f_1^\sigma \ast f_2^\sigma)^\sigma
\, = \ f_1 \ast f_2 ,
$$
and this concludes the proof.
\end{proof}

We can now state the main result of this subsection.

\begin{theorem}
% 6.E.17
\label{Theo-CCRforSL2}
The subgroup $K = \SO(2)$ is a large subgroup of $G = \SL_2(\R)$.
\par

Consequently, $G = \SL_2(\R)$ is a CCR group and is therefore of type I.
\end{theorem}

\begin{proof}
Let $(\pi, \Hi)$ be an irreducible representation of $G$. Let 
$$
\Hi \, = \, \bigoplus_{n \in \Z} \Hi_{n}
$$
be the decomposition of $\Hi$ into $K$-isotypical components $\Hi_{n} := \Hi_{\chi_n}$.
We will show that $\Hi_n$ is finite-dimensional,
in fact at most one-dimensional, for every $n \in \Z$.

\vskip.2cm

Let $n \in \Z$ be fixed in the sequel.

\vskip.2cm

$\bullet$ {\it First step.}
We claim that $\Hi_n$ is invariant under 
$\pi(C^c(G)_{ \chi_n})$.
\par

Indeed, recall that the orthogonal projection 
$P_n \, \colon \Hi \to \Hi_n$ is the operator $\pi( \overline{\chi_n} ) = \pi(\chi_{-n})$.
Let $f \in C^c(G)_{ \chi_n}$. Then $f = \chi_n \ast f \ast \chi_n$ and hence
$$
\pi(f) \, = \, \pi (\chi_n \ast f \ast \chi_n)
\, = \, \pi (\chi_n) \pi(f) \pi (\chi_n)
\, = \, P_n \pi(f) P_n .
$$ 
This shows that $\Hi_n$ is invariant under $\pi(f)$.
\par
 
As a result, we have a representation 
$$
\pi_n \, \colon \, C^c(G)_{ \chi_n} \to \Li (\Hi_n),
\hskip.3cm
f \mapsto \pi(f)\vert_{\Hi_n} .
$$ 
Recall that $C^c(G)_{ \chi_n}$ is invariant under the involution $f \mapsto f^*$;
therefore, $\pi_n$ is a $*$-homomorphism.
 
\vskip.2cm

$\bullet$ {\it Second step.}
We claim that the commutant $\pi_n(C^c(G)_{ \chi_n})'$
of the selfadjoint subalgebra $\pi_n(C^c(G)_{ \chi_n})$ of $\Li (\Hi_n)$
consists only of scalars operators.
\par
 
Indeed, it suffices to show that the von Neumann algebra 
$\pi_n(C^c(G)_{\chi_n})''$ generated by $\pi_n(C^c(G)_{\chi_n})$ 
coincides with $\Li (\Hi_n)$. Let $T \in \Li (\Hi_n)$ and set
$$
\widetilde T \, := \,TP_n \in \Li (\Hi).
$$
Since $\pi$ is irreducible, the von Neumann algebra generated by $\pi(C^c(G))$ 
coincides with $\Li (\Hi)$.
Therefore there exists a net $(f_i)_{i \in I}$ in $C^c(G)$
such that $\lim_i \pi(f_i) = \widetilde T$ in the strong (or weak) operator topology.
Then $\chi_{-n} \ast f_i \ast \chi_{-n} \in C^c(G)_{ \chi_n}$ and 
$$
\lim_i \pi(\chi_{-n} \ast f_i\ast \chi_{-n})
\, = \, \lim_i (P_n \pi(f_i) P_n)
\, = \, P_n \lim_n \pi(f_i) P_n
\, = \, P_n \widetilde T P_n
\, = \, P_n T P_n.
$$
Therefore 
$$
\lim_ i \pi_n(\chi_{-n} \ast f_i\ast \chi_{-n} ) \, = \, T
$$
and the claim is proved.

\vskip.2cm

$\bullet$ {\it Third step.}
We claim that $\pi_n$ is an irreducible representation of $C^c(G)_{ \chi_n}$.
\par

Indeed, let $\Ki$ be a $\pi_n(C^c(G)_{ \chi_n})$-invariant closed subspace
of $\Hi_{n}$ and let 
$$
P \, \colon \, \Hi_n \to \Ki
$$
be the orthogonal projection.
Let $f \in C^c(G)_{\chi_n}$.
Then, $P \pi_n(f) P = \pi_n(f) P$ and therefore
$$
P\pi_n(f) \, = \, (\pi_n(f^*) P)^* \, = \, (P \pi_n(f^*) P)^* \, = \, P \pi_n(f) P,
$$
that is, $\pi_n(f) P = P\pi_n(f)$. So, $P \in \pi_n(C^c(G)_{ \chi_n})'$.
Therefore by the second step, $P$ is a scalar operator,
that is, $\Ki = \{0\}$ or $\Ki = \Hi_n$.

\vskip.2cm

$\bullet$ {\it Fourth step.}
The dimension of $\Hi_n$ is at most one.
\par

Indeed, $\pi_n(C^c(G)_{ \chi_n})'$
consists only of scalar operators, by the second step.
Since $\pi_n$ is irreducible (third step)
and $C^c(G)_{ \chi_n}$ is commutative (Proposition~\ref{Prop-CommutativeAlgebraSpherical}),
it follows that $\dim \Hi_n\le 1$.
\end{proof}

\begin{rem}
% 6.E.18
\label{Rem-Theo-CCRforSL2}
Parts of the proof of Theorem~\ref{Theo-CCRforSL2} are valid in a much greater generality. 
\par

Indeed, let $G$ be a unimodular LC compact and $K$ a compact subgroup of $G$.
For $\sigma \in \widehat{K}$, let $\psi_\sigma = d_\sigma \chi_\sigma$,
where $\chi_\sigma$ is the character of $\sigma$ and $d_\sigma$ its dimension.
Define
$$
C^c(G)_\sigma \, = \, \psi_\sigma \ast C^c(G) \ast \psi_\sigma.
$$
Then $C^c(G)_\sigma$ is a subalgebra of the convolution algebra $C^c(G)$
and is invariant under the involution $f\mapsto f^*$. 
\par

Let $(\pi, \Hi)$ be an irreducible representation of $G$ and let 
$$
\Hi \, = \, \bigoplus_{\sigma\in \widehat{K},} \Hi_{\sigma}
$$
be the decomposition of $\Hi$ into $K$-isotypical components. 
\par

The arguments of the three first steps in the proof of Theorem~\ref{Theo-CCRforSL2}
carry over, with the obvious changes, and show that the following holds:
\begin{itemize}
\setlength\itemsep{0em}
\item 
$\Hi_\sigma$ is invariant under $\pi(C^c(G)_{ \sigma})$
and so defines a $*$-representation $\pi_\sigma$ of $C^c(G)_{ \sigma}$;
\item
the commutant of $\pi_\sigma(C^c(G)_{\sigma})$ in $\Li (\Hi_\sigma)$
consists only of scalar operators;
\item
$(\pi_\sigma, \Hi_\sigma)$ is an irreducible representation of $C^c(G)_{ \sigma}$.
\end{itemize}
\par

Assume now that $G$ is a linear semisimple real Lie group
and $K$ a maximal compact subgroup of $G$.
Using the fact that $G$ has a decomposition $G = K S$,
where $S$ is a closed solvable subgroup (``Iwasawa decomposition"),
one shows (see \cite[15.5.5]{Dixm--C*} or \cite[Theorem 2]{Gode--52})
that the last step in the proof of Theorem~\ref{Theo-CCRforSL2}
remains valid in the following form:
\begin{itemize}
\setlength\itemsep{0em}
\item
the dimension of $\Hi_\sigma$ is at most $d_\sigma$.
\end{itemize}
\end{rem}

\subsection*{On the class of groups of type I}

We give a non exhaustive list of topological groups which are known to be of type I.

\begin{theorem}[Some locally compact groups of type I]
% 6.E.19
\label{explesTypeI}
The following locally compact groups are of type I:
\begin{enumerate}[label=(\arabic*)]
\item\label{iDEexplesTypeI}
Second-countable locally compact abelian groups.
\item\label{iiDEexplesTypeI}
Compact groups. 
\item\label{iiiDEexplesTypeI}
Nilpotent connected locally compact groups, 
and more generally connected-by-compact locally compact groups
whose solvable radical is actually nilpotent.
\item\label{ivDEexplesTypeI} Reductive algebraic groups over a local field.
% Harish--Chandra pour r\'eductifs, livre Wallach 3.4.10 et 14.6.10
% \footnote{A la page 275, Dixmier se r\'ef\`ere \`a UN AUTRE article de Harish --- de 1954.}
% \marginpar{R\'ef correcte~?}
\item\label{vDEexplesTypeI}
Connected real algebraic groups.
\item\label{viDEexplesTypeI}
Adelic reductive groups of the form $G = \mathbf{G}(\A)$,
where $\mathbf{G}$ is a reductive algebraic group defined over $\Q$
and $\A$ the adele ring of the field $\Q$.
\item\label{viiDEexplesTypeI}
Some solvable locally compact groups, e.g., 
all solvable connected real Lie groups that are of dimensions at most $4$,
and all solvable connected real Lie groups 
for which the exponential map is surjective.
\item\label{viiiDEexplesTypeI}
Full groups of automorphisms of Bruhat--Tits trees.
\end{enumerate}
\end{theorem}

\begin{proof}[References for proofs and comments]
Claim \ref{iDEexplesTypeI} follows from Corollary \ref{Cor-Thm-CanRepAbMultFree}
and Corollary~\ref{Cor-Prop-DirectSumTypeIRep};
Claim \ref{iiDEexplesTypeI} is straightforward.
\par

For \ref{iiiDEexplesTypeI}
see \cite{Dixm--59},
% \cite{Kiri--62},
and \cite{Lips--72}.
For \ref{ivDEexplesTypeI}, see Theorem~\ref{Thm-ReductiveCCR groups}.
% see \cite[Page 230]{Hari--53}
% (also \cite[15.5.6]{Dixm--C*}) and \cite{Dixm--57}.
% For \ref{vDEexplesTypeI},
% see \cite{Bern--74}.
\par

For Claim \ref{vDEexplesTypeI}, see \cite{Dixm--57}.
\par

Claim \ref{viDEexplesTypeI} follows from
\cite[Chap.~III, \S~3, no.~3]{GGPS--69},
combined with our \ref{ivDEexplesTypeI}. See also \cite[Section 6]{Moor--65}.
\par

In contrast to this and in relation with \ref{iiiDEexplesTypeI},
it is shown in \cite[beginning of Section 7]{Moor--65} that,
for a nilpotent algebraic group $\mathbf G$ defined over $\Q$,
the adelic group $\mathbf{G}(\A)$ is type I if and only if $\mathbf G$ is abelian.
\par

Concerning \ref{viiDEexplesTypeI},
see \cite{Take--57} for groups of which the exponential map is surjective.
A necessary and sufficient condition for a solvable connected simply connected Lie group
to be of type I is established in \cite{AuKo--71},
in terms of coadjoint orbits.
\par

For \ref{viiiDEexplesTypeI}, see \cite{Ol's--77} and \cite{Ol's--80}.
For other groups of automorphisms of trees, see \cite{HoRa--19}.
\end{proof}

The universal covering group of a connected (even solvable) Lie group of type I
need not be type I \cite{Dixm--61}.

\begin{theorem}[Some more topological groups of type I]
% 6.E.20
\label{explesTypeI-bis}
The following topological groups are type I:
\begin{enumerate}[label=(\arabic*)]
\item\label{iDEexplesTypeI-bis}
Abelian topological groups.
\item\label{iiDEexplesTypeI-bis}
The groups $\U_\infty(\Hi)$ and $\U(\Hi)_{\rm str}$,
as defined above in \ref{exUinfty}.
\item\label{iiiDEexplesTypeI-bis}
The group $\Sym (\N)$ of permutations of $\N$
with the topology of simple convergence.
\item\label{ivDEexplesTypeI-bis}
The group $\Aut(\Q, \le)$ of automorphisms of the ordered set
of the rationals with the topology of simple convergence.
\item\label{vDEexplesTypeI-bis}
The group ${\rm Homeo}(2^\N)$
of homeomorphisms of the Cantor space with the compact-open topology.
\end{enumerate}
\end{theorem}

\begin{proof}[References for proofs and comments]
Claim \ref{iDEexplesTypeI-bis} is a consequence of the (non obvious) fact
that abelian von Neumann algebras are of type I (see Proposition~\ref{Pro-TypeIVN-bis}).
\par

For Claim \ref{iiDEexplesTypeI-bis}, see
\cite{Kiri--73}, \cite{Ol's--78}, and \ref{ExamplesND} above.
For a discussion of the two latter articles and further results on representations of unitary groups,
see \cite{Neeb--14}.
\par

For \ref{iiiDEexplesTypeI-bis}, 
see \cite[Theorems 3 \& 4]{Lieb--72}.
More generally, consider a countable infinite set $X$,
the Polish group $\Sym(X)$ of permutations of $X$,
and a subgroup $G$ of $\Sym(X)$; assume that
\begin{enumerate}[label=(\alph*)]
\item\label{aDEexplesTypeI}
$G$ is oligomorphic,
i.e., for all $n \ge 1$
its diagonal action of $G$ on $X^n$ has a finite number of orbits,
\item\label{bDEexplesTypeI}
for every separable Hilbert space $\Ki$, 
every homomorphism $G \to \U(\Ki)$ is continuous.
\end{enumerate}
Then $G$ has countably many equivalence classes of irreducible representations,
and every representation of $G$ in a separable Hilbert space
is a direct sum of irreducible representations (see \cite{Tsan--12}).
In particular, $G$ is of type I.
This applies to the groups $\Aut(\Q)$ 
and ${\rm Homeo}(2^\N)$ of \ref{vDEexplesTypeI-bis}.
\end{proof}

We mention a list of operations under which the class of type I groups is stable.

\begin{prop}
% 6.E.21
\label{stableTypeI}
(Stability properties of type I groups)
\begin{enumerate}[label=(\arabic*)]
\item\label{iDEstableTypeI}
Open subgroups of type I second-countable locally compact groups are type I.
For a second-countable locally compact group $G$ and a closed subgroup $H$ of finite index,
$G$ is type I if and only if $H$ is type I.
\item\label{iiDEstableTypeI}
Quotients of type I groups are type I:
if $G$ is a type I topological group and $N$ a closed normal subgroup,
then $G/N$ is type I.
\item\label{iiiDEstableTypeI}
A Cartesian product of a finite number of type I groups is type I.
\item\label{ivDEstableTypeI}
Let $G = \varprojlim G_\alpha$ be a projective limit of locally compact groups
$G_\alpha = G/K_\alpha$,
with $K_\alpha$ a compact normal subgroup of $G$ for all $\alpha$;
then the locally compact group $G$ is type I
if and only if each $G_\alpha$ is type I.
\item\label{vDEstableTypeI}
Let $G$ be a topological group and $H$ a closed subgroup such that 
$G/H$ is locally compact 
and possesses a non-zero, finite, $G$-invariant Radon measure;
if $H$ is type I then $G$ is type I. 
\par

In particular, (type I)-by-compact locally compact-groups are type I, 
and a fortiori abelian-by-compact locally compact groups are type I.
\end{enumerate}
\end{prop}

\begin{proof}[References for proofs]
For \ref{iDEstableTypeI}, see Proposition 2.4 and Corollary 2.5 in \cite{Kall--73}.
Fact \ref{iiDEstableTypeI} is a straightforward consequence of the definitions.
For \ref{iiiDEstableTypeI}, see, e.g., \cite[Theorem 3.2]{Mack--76};
% page 128
if $G_1, \hdots, G_k$ are LC groups of type I,
the dual of the product $\prod_{i = 1} G_i$
is naturally homeomorphic to the Cartesian product $\prod_{i = 1}^k \widehat G_i$
of the duals of the $G_i$~'s.
Claim \ref{ivDEstableTypeI} follows from the definitions
and from Proposition~\ref{FacKK} below
(see also the remark following Proposition 2.2 in \cite{Moor--72}).
Claim \ref{vDEstableTypeI} is Theorem 1 in \cite{Kall--70}.
The particular case of \ref{vDEstableTypeI} for second-countable LC groups
is alternatively a consequence of \cite[Theorem 9.3]{Mack--58}.
\end{proof}

%-----------------------------------------------------------------------
% End of chapter 6
%-----------------------------------------------------------------------
\chapter{Non type I groups}
% Chapter 7
\label{NonTypeIGroups}

\emph{
There are large classes of groups which are not of type I.
Following Murray and von Neumann, we first show in Section~\ref{S:ICCGroupsNontypeI}
that the regular representation of a discrete icc group
is a factor representation which is not of type I.
}

\emph{
Recall from Proposition~\ref{Pro-TypeIVN}
that a type I representation $\pi$ of a group $G$
is characterized by a property of the von Neumann algebra $\pi(G)''$.
Von Neumann algebras are subdivided into various types,
the so-called types I$_{\rm f}$, I$_\infty$, II$_1$, II$_\infty$, and III.
We recall this type classification in Section \ref{SectionvN},
after a short reminder on traces on operator algebras. 
This leads in Section~\ref{SectionvN+rep}
to the definition of group representations of corresponding types. 
Historically, types were first defined by Murray and von Neumann in 1936,
before being imported in the theory of unitary group representations,
notably by Mautner, Godement, Segal, and Mackey
(see \cite{Maut--50a}, \cite{Maut--50b}, \cite{Gode--51a}, \cite{Sega--51}, \cite{Mack--53}).
}

\emph{
We discuss several groups that are not of type I.
In particular, we state Thoma's characterization of discrete groups of type I
(Theorem~\ref{discretTypes-Thoma}).
The regular representation $\lambda_\Gamma$ of a discrete group 
$\Gamma$ generates a von Neumann algebra of finite type.
It turns out that $\Gamma$ is of type I if and only if $\lambda_\Gamma$ is of type I
(Theorem~\ref{discretTypes}).
We indicate a class of examples of groups $\Gamma$ 
for which $\lambda_\Gamma$
has both a type I component and a type II component
(Theorem~\ref{Theo-Kaplansky}).
We mention in \ref{SectionTypesI+IILC} results 
about the type of the regular representation for classes of non discrete locally compact groups.
}

\emph{
In Section~\ref{S:IntDecFactorRep},
we discuss the direct integral decompositions
of a factor representation $\pi$ of a second-countable LC group $G$
into irreducible representations.
We will see that almost every irreducible representation
involved in such a direct integral decomposition of $\pi$
is weakly equivalent to $\pi$.
This last result has important consequences:
every factor representation of $G$ is weakly equivalent
to some irreducible representation of $G$,
and the quotient map $\widehat G \twoheadrightarrow \Pri(G)$
has a natural extension
to a map $\QD(G) \twoheadrightarrow \Pri(G)$;
moreover, when $G$ is not of type I, 
many fibers of the map $\widehat G \twoheadrightarrow \Pri(G)$ are uncountable.
% (Corollary~\ref{Cor-FacRepPrimitiveG-NonTypI}).
}

\section{A class of non type I groups}
% Section 7.A
\label{S:ICCGroupsNontypeI}

As we will see later (Theorem~\ref{discretTypes}),
most discrete groups are not of type I.
For the moment, we single out the following venerable result of Murray and von Neumann
\cite[Lemma 5.3.4]{MuvN--43};
see also Chap.~III, \S~7, no~6 in \cite{Dixm--vN}.
% pages 282--283.
For the definition and examples of icc groups, see Appendix~\ref{AppGpactions}.
For the definition of factors of type II$_1$, see Appendix~\ref{AppAlgvN}.

\begin{prop}
% 7.A.1
\label{iccfactorII1}
\index{Icc group}
Let $\Gamma$ be a discrete group, $\Gamma \ne \{e\}$,
and $\Li (\Gamma) := \lambda_\Gamma(\Gamma)''$
the von Neumann algebra of the left regular representation of $\Gamma$.
The following properties are equivalent:
\begin{enumerate}[label=(\roman*)]
\item\label{iDEiccfactorII1}
$\Gamma$ is an icc group;
\item\label{iiDEiccfactorII1}
$\Li (\Gamma)$ is a factor.
\end{enumerate}
When this is the case, $\Li (\Gamma)$ is a factor of type II$_1$.
\end{prop}

\begin{proof}
Let $\rho_\Gamma$ be the right regular representation $\rho_\Gamma$
of $\Gamma$ on $\ell^2(\Gamma)$.
Observe that
$\lambda_\Gamma(\gamma) \rho_\Gamma(\gamma')
= \rho_\Gamma(\gamma') \lambda_\Gamma(\gamma) $
for all $\gamma, \gamma' \in \Gamma$.
% (This is straightforward. We \emph{do not} use here Theorem \ref{Theo-RegRepHilbertAlg}
% according to which the von Neumann algebra
% generated by $\{ \rho_\Gamma(\gamma') \mid \gamma' \in \Gamma \}$
% is \emph{equal} to the commutant of $\Li (\Gamma)$.)
\par

Let $T \in \Li (\Gamma)$. Then
$$
\rho_\Gamma(\gamma) T \, = \, T \rho_\Gamma(\gamma)
\hskip.5cm \text{for all} \hskip.2cm
\gamma \in \Gamma .
$$
Suppose moreover that $T$ is in the centre of $\Li (\Gamma)$. Then
$$
\lambda_\Gamma(\gamma) T \, = \, T \lambda_\Gamma(\gamma)
\hskip.5cm \text{for all} \hskip.2cm
\gamma \in \Gamma .
$$
Set
$$
f \, := \, T\delta_e \in \ell^2(\Gamma).
$$
Then, for every $\gamma \in \Gamma$,
we have $\lambda_\Gamma(\gamma) \rho_\Gamma(\gamma)\delta_e \, = \, \delta_e$
and
$$
\lambda_\Gamma(\gamma)\rho_\Gamma(\gamma) f
\, = \, ( \lambda_\Gamma(\gamma)\rho_\Gamma(\gamma)T)\delta_e
\, = \, (T \lambda_\Gamma(\gamma) \rho_\Gamma(\gamma))\delta_e
\, = \, T\delta_e \, = \, f.
$$
Therefore $f$ is a function in $\ell^2(\Gamma)$ which is central.
It follows that the support of $f$ consists
of elements with finite conjugacy class.

\vskip.2cm

\ref{iDEiccfactorII1} $\Rightarrow$ \ref{iiDEiccfactorII1}
Assume that $\Gamma$ is icc. Then $f$ is supported by $\{e\}$, that is
$T\delta_e = f = c\delta_e$ for a scalar $c \in \C$.
For every $\gamma \in \Gamma$, we have
$$
T\delta_\gamma
\, = \,T \lambda_\Gamma(\gamma)\delta_e
\, = \, \lambda_\Gamma(\gamma) T\delta_e
\, = \, c \lambda_\Gamma(\gamma) \delta_e \, = \, c \delta_\gamma.
$$
Since the span of $\{\delta_\gamma \mid \gamma \in \Gamma\}$
is dense in $\ell^2(\Gamma)$, 
it follows that $T = cI$.
This shows that $\lambda_\Gamma(\Gamma)''$ is a factor.
Moreover, this factor is finite,
because it has a finite trace $T \mapsto \langle T \delta_e \mid \delta_e \rangle$,
and infinite dimensional,
because $\Gamma$ is infinite.
This means that $\lambda_\Gamma(\Gamma)''$ is a factor of type II$_1$.

\vskip.2cm

\ref{iiDEiccfactorII1} $\Rightarrow$ \ref{iDEiccfactorII1}
Assume now that $\Gamma$ is not icc, that is, there exists
a finite conjugacy class $C \subset \Gamma$ with $C \ne \{e\}$.
Consider the operator
$$
T \, := \, \lambda_\Gamma(\Un_C) \, = \, 
\sum_{\gamma \in C} \lambda_\Gamma(\gamma),
$$
which belongs to $\lambda_\Gamma(\Gamma)''$.
Then $T \in \lambda_\Gamma(\Gamma)'$;
indeed, for every $x \in \Gamma$, we have
$$
\begin{aligned}
T \lambda_\Gamma(x)
\, &= \, \sum_{\gamma \in C} \lambda_\Gamma(\gamma) \lambda_\Gamma(x)
\, = \, \lambda_\Gamma(x)\sum_{\gamma \in C} \lambda_\Gamma(x^{-1}\gamma x )
\\
\, &= \, \lambda_\Gamma(x)\sum_{\gamma \in C} \lambda_\Gamma(\gamma)
\, = \, \lambda_\Gamma(x) T.
\end{aligned}
$$
Therefore $T$ belongs to the centre
of $\lambda_\Gamma(\Gamma)''$.
Moreover, $T$ is not a multiple of the identity, since $T\delta_e = \Un_C$. 
So, $\Li (\Gamma)$ is not a factor.
\end{proof}

\begin{defn}
% 7.A.2
\label{defvNofGamma}
The \textbf{von Neumann algebra of a discrete group} $\Gamma$
is the von Neumann algebra $\Li(\Gamma) := \lambda_\Gamma(\Gamma)''$
of its left regular representation $\lambda_\Gamma \,\colon \Gamma \to \U (\ell^2(\Gamma))$.
\index{von Neumann algebra $\Li(\Gamma)$ of a discrete group $\Gamma$}
\end{defn}

A remarkable property of von Neumann algebras of discrete groups,
also observed by Murray and von Neumann,
is that they admit faithful finite normal traces,
as shown by the following proposition
(these terms are defined in Section~\ref{SectionvN}).
Recall that $(\delta_\gamma)_{\gamma \in \Gamma}$
denotes the standard basis of the Hilbert space $\ell^2(\Gamma)$,
and that $\delta_e$ is the basis element defined by the unity element $e \in \Gamma$;
we denote as in Appendix \ref{AppAlgC*} by $\Li (\Gamma)_+$
the cone of positive operators in $\Li (\Gamma)$.
 
\begin{prop}
% 7.A.3
\label{Pro-TraceDiscreteGroup}
Let $\Gamma$ be a discrete group. The linear functional 
$$
t \, \colon \, \Li (\Gamma) \to \C, \hskip.2cm T \mapsto \langle T \delta_e \mid \delta_e \rangle
$$
has the following properties:
\begin{enumerate}[label=(\arabic*)]
\item\label{iDEPro-TraceDiscreteGroup}
$t$ is continuous for the weak operator topology on $\Li (\Gamma)$;
\item\label{iiDEPro-TraceDiscreteGroup}
$t(\mathrm{Id}_{\ell^2(\Gamma)}) = 1$;
\item\label{iiiDEPro-TraceDiscreteGroup}
$t(ST) = t(TS)$ for all $S,T \in \Li (\Gamma)$;
\item\label{ivDEPro-TraceDiscreteGroup}
$t$ is faithful, that is, $t(T) > 0$ for every $T \in \Li (\Gamma)_+$, $T \ne 0$.
\end{enumerate}
\end{prop}
 
\begin{proof}
Properties \ref{iDEPro-TraceDiscreteGroup} and \ref{iiDEPro-TraceDiscreteGroup} are obvious.

\vskip.2cm

For $\gamma \in \Gamma$, we have
$t(\lambda_\Gamma(\gamma)) = 1$ if $\gamma = e$
and $t(\lambda_\Gamma(\gamma)) = 0$ otherwise.
This implies that 
$$
t(\lambda_\Gamma(\gamma_1)\lambda_\Gamma(\gamma_2))
\, = \,
t(\lambda_\Gamma(\gamma_1 \gamma_2))
\, = \, 
t(\lambda_\Gamma(\gamma_2 \gamma_1))
\, = \,
t(\lambda_\Gamma(\gamma_2)\lambda_\Gamma(\gamma_1))
$$
for all $\gamma_1, \gamma_2 \in \Gamma$.
By linearity and density, Property \ref{iiiDEPro-TraceDiscreteGroup} follows.

\vskip.2cm

To show Property \ref{ivDEPro-TraceDiscreteGroup},
let $T \in \Li (\Gamma)_+$ be such that $t(T) = 0$;
we have to show that $T = 0$.
Write $T = S^*S$ for $S \in \Li (\Gamma)$. 
Then 
$$
0 \, = \, t(T) \, = \, t(S^*S) \, = \, \langle S^*S \delta_e \mid \delta_e \rangle
\, = \, \Vert S \delta_e \Vert^2 ,
$$
that is, $S \delta_e = 0$.
Now, $\delta_e$ is a cyclic vector
for the right regular representation $\rho_\Gamma$ and 
$S \in \Li (\Gamma) \subset \rho_\Gamma(\Gamma)'$;
hence $S \xi = 0$ for every $\xi \in \ell^2(\Gamma)$.
It follows that $S = 0$, that is, $T = 0$.
 \end{proof} 

We are going to show that, for an icc group $\Gamma$,
the regular representation $\lambda_\Gamma$ is not of type I.
Our proof will depend on the following fact, which is of independent interest.

\begin{prop}
% 7.A.4
\label{Pro-NoIrredSubRegRep}
Let $\Gamma$ be a discrete group. Assume that $\lambda_\Gamma$
contains an irreducible subrepresentation.
\par

Then $\Gamma$ is finite.
\end{prop}

% Retablissement de la preuve de Janvier 2019
\begin{proof}
Assume that $\lambda_\Gamma$ contains an 
irreducible subrepresentation $\pi$.
Let $\Hi$ be the $\Gamma$-invariant closed subspace of $\ell^2(\Gamma)$
defining $\pi$.

Let $P \in \Li (\Gamma)'$ be the orthogonal projection from $\ell^2(\Gamma)$ onto $\Hi$.
The induced von Neumann algebra $\Li (\Gamma)_P = P \Li (\Gamma) P$,
viewed as a subalgebra of $\Li (\Hi)$, coincides with $\pi(\Gamma)''$
and hence with $\Li (\Hi)$, since $\pi$ is irreducible.
Let $E \in \Li (\Gamma) \cap \Li (\Gamma)'$ be the central support of $P$.
Recall from Subsection \ref{SS:QE+} that the induction map
$$
\Phi \, \colon \, \Li (\Gamma)_E \to \Li (\Gamma)_P
$$
is an isomorphism.
Consider the linear functional $\varphi \,\colon \Li (\Hi) \to \C$ defined by 
$$
\varphi(T) \, = \, \langle \Phi^{-1}(T) \delta_0 \mid \delta_0 \rangle
\hskip.2cm \text{for} \hskip.2cm
T \in \Li (\Hi) = \Li (\Gamma)_P .
$$
Then, by Proposition~\ref{Pro-TraceDiscreteGroup}, $\varphi$ has the following properties:
\begin{enumerate}[label=(\alph*)]
\item\label{aDEPro-NoIrredSubRegRep}
$\varphi$ is continuous for the weak operator topology on $\Li (\Hi)$;
\item\label{bDEPro-NoIrredSubRegRep}
$\varphi (ST)\, = \, \varphi (TS)$ for all $S,T \in \Li (\Hi)$.
\item\label{cDEPro-NoIrredSubRegRep}
$\varphi \ne 0$.
\end{enumerate}
It follows that $\Hi$ is finite dimensional.
\par

Indeed, by Property \ref{bDEPro-NoIrredSubRegRep},
we have $\varphi(p) = \varphi(q)$ for any two projections $p,q \in \Li (\Hi)$
with the same dimension,
as $p$ and $q$ are conjugate under the unitary group of $\Hi$.
Since $\mathrm{Id}_{\Hi} = \sum_i p_i$ in the weak operator topology
for one-dimensional projections $p_i$,
it follows from Properties \ref{aDEPro-NoIrredSubRegRep} and \ref{cDEPro-NoIrredSubRegRep}
that $\Hi$ is finite dimensional.
\par

Since $\pi$ is finite-dimensional, the trivial representation $1_\Gamma$
is contained in the tensor product $\pi \otimes \overline \pi$,
where $\overline \pi$ is the conjugate representation of $\pi$
\cite[A.1.13]{BeHV--08}. 
Therefore $1_\Gamma$ is contained in $\lambda_\Gamma \otimes \overline \lambda_\Gamma$,
which is equivalent to $\lambda_\Gamma$ (see \cite[Appendix E]{BeHV--08}).
So, $\ell^2(\Gamma)$ contains a non-zero invariant vector; this is only possible if 
$\Gamma$ is finite.
\end{proof}

The analogue for LC groups of Proposition \ref{Pro-NoIrredSubRegRep} does not hold.
Indeed, there are second-countable LC groups $G$ which are Fell groups, i.e., 
which are such that the regular representation $\lambda_G$
is a direct sum of irreducible representations, and which are non-compact.
We refer to \cite[Section IV]{Bagg--72} for the original example of Fell,
and to \cite{Robe--78} for some other examples.
\par

The following corollary is an immediate consequence of 
Propositions~\ref{iccfactorII1}, \ref{Pro-NoIrredSubRegRep}, and \ref{Prop-TypeIFac}.

\begin{cor}
% 7.A.5
\label{Cor-iccfactorII1}
Let $\Gamma$ be a discrete icc group.
\par

Then $\Gamma$ is not of type I.
\end{cor}

\section
{Operator algebras, traces and types}
% Section 7.B
\label{SectionvN}

The aim of the present section is to recall some of the terminology
from the theory of operator algebras relevant for the type classification 
of representations.

\subsection*{Traces on C*-algebras}

\index{$e2$@$A$ C*-algebra}
\index{Trace! $1$@on a C*-algebra}
Traces are important for both von Neumann algebras and C*-algebras,
and also because characters of group representations are defined in terms of traces
(see Chapter \ref{ChapterCharacters}).
Consider first a C*-algebra~$A$.
The cone of positive elements
is $A_+ := \{ x \in A \mid x = y^*y \hskip.2cm \text{for some} \hskip.2cm y \in A \}$.
For $x,y \in A_+$, the order $x \le y$ is defined by $y-x \in A_+$.
More on $A_+$ in Appendix \ref{AppAlgC*}.
\par

A \textbf{trace} on $A$ is a function $t \,\colon A_+ \to \mathopen[ 0, \infty \mathclose]$
such that 
\begin{enumerate}[label=(\arabic*)]
\item 
$t(x+y) = t(x) + t(y)$ for all $x,y \in A_+$;
\item 
$t(\lambda x) = \lambda t(x)$ 
for all $x,y \in A_+$ and $\lambda \in \mathopen[ 0, \infty \mathopen[$
\par\noindent
(it is agreed that the product $0 \infty$ is $0$);
\item 
$t(x^*x) = t(xx^*)$ for all $x \in A_+$.
\end{enumerate}
Such a trace is
\begin{enumerate}[label=(\arabic*)]
\addtocounter{enumi}{3}
\item
\index{Semi-finite trace}
\textbf{semi-finite} 
if $t(x) = \sup \{t(y) \mid y \in A_+, \hskip.1cm y \le x, \hskip.1cm t(y) < \infty \}$
for all $x \in A_+$;
\item
\textbf{finite}
\index{Finite trace}
if $t(x) < \infty$ for all $x \in A_+$;
\item
\index{Lower semi-continuous trace}
\textbf{lower semi-continuous}
if, for all $x \in A_+$ and $\lambda \in \R$ such that $\lambda < t(x)$,
there exists $\varepsilon > 0$ such that $\lambda < t(y)$
for all $y \in A_+$ such that $\Vert y - x \Vert < \varepsilon$;
% equivalently, if $t(x) \le \liminf_n t(x_n)$ for all $x \in A_+$
% and all sequences $(x_n)_{n \ge 1}$ in $A_+$ such that $x = \lim_n x_n$
% (for the equivalence: \cite[Chap.~IV, \S~6, no~2]{BTG1--4});
\item
\index{Faithful trace}
\textbf{faithful}
if, for $x \in A_+$, the equality $t(x) = 0$ implies $x = 0$.
\end{enumerate}
\index{Trace! $2$@domination}
Let $t_1, t_2$ be two traces on $A$;
then $t_1$ \textbf{dominates} $t_2$
if $t_1(x) \ge t_2(x)$ for all $x \in A_+$.

\begin{exe}
% 7.B.1.
\label{ExampleTraceC*}
(1)
Let $A$ be a C*-algebra, not $\{0\}$.
The function $t$ defined on $A_+$ by $t(0) = 0$ and $t(x) = \infty$ for all $x \ne 0$ in $A_+$
is a lower semi-continuous faithful trace that is not semi-finite.

\vskip.2cm

(2)
Let $c^0(\N)$ be the C*-algebra of sequences $(z_i)_{i \ge 1}$
of complex numbers such that $\lim_{i \to \infty} z_i = 0$,
equipped with the norm $\Vert (z_i)_{i \ge 1} \Vert_\infty := \sup_{i \ge 1} \vert z_i \vert$.
Its positive cone $c^{0}_+$ is the subset of sequences
in $c^0(\N)$ of non-negative real numbers.
Let $\mu = (\mu_i)_{i \ge 1}$ be a sequence of non-negative real numbers;
a trace $t_\mu$ on $c^0(\N)$ is defined by
$$
t_\mu \, \colon \, c^0(\N)_+ \to \mathopen[ 0, \infty \mathclose] ,
\hskip.5cm
(z_i)_{i \ge 1} \, \mapsto \, \sum_{i = 1}^\infty \mu_i z_i .
$$
This trace is clearly semi-finite and it is an exercise in calculus
to check that it is lower semi-continuous.
The trace $t_\mu$ is finite if and only if $\sum_{i = 1}^\infty \mu_i < \infty$,
and faithful if and only if $\mu_i > 0$ for all $i \ge 1$.
\par

It is easy to show that every semi-finite and semi-continuous trace on $c^0(\N)$
is of the form $t_\mu$ for a sequence $\mu$ as above.
\par

For a generalization of this example, see Proposition~\ref{Pro-TraceCommC*}.

\vskip.2cm

(3)
Let $t$ be the trace on $C^0(\R)$ defined on $f \ge 0$
by $t(f) = 0$ if the support of $f$ is compact and $t(f) = \infty$ otherwise. 
Then $t$ is neither lower semi-continuous nor semi-finite.

\vskip.2cm

(4)
Let $\mathcal K (\Hi)$ be the C*-algebra of compact operators 
on a separable infinite-dimensional Hilbert space $\Hi$.
There exist so-called ``Dixmier traces'' $t$ on $\mathcal K(\Hi)$,
not the zero trace,
which vanish on the usual trace-class operators \cite{Dixm--66}. 
They can be shown to be semi-finite, and not lower semi-continuous.

\vskip.2cm

(5) 
Let $G$ be a unimodular LC group,
$C^*_\lambda(G)$ the reduced C*-algebra of $G$
(see Chapter \ref{ChapterAlgLCgroup}),
and $C^c(G)$ the sub-$*$-algebra of continuous functions on $G$ with compact supports.
A semi-finite faithful trace $t$ is defined on $C^*_\lambda(G)$ by
$$
t (f \ast f^*) \, = \, (f \ast f^*) (e) \, = \, \Vert f \Vert_2^2
\hskip.5cm \text{for all} \hskip.2cm
f \in C^c(G),
$$
and suitably extended on elements of $C^*_\lambda(G)_+$;
see Example~\ref{Ex-HilbertAlgebra}(2).
\par

The trace $t$ is finite if and only if $G$ is discrete.
In this case, we write $\Gamma$ rather than $G$ for the group,
and $\C[\Gamma]$ rather than $C^c(\Gamma)$ for the group algebra.
(See also Example \ref{ExampleTracevN}(2).)

\vskip.2cm

(6)
Let $G$ be a locally compact group. 
There exists a non-zero finite trace on $C^*_\lambda(G)$
if and only if the amenable radical of $G$ is open \cite{FoSW, KeRa--17}.
\end{exe}

The following proposition gives a complete description
of lower semi-continuous and semi-finite traces on a commutative C*-algebra.
\index{Algebras! $1$@continuous functions, $C(X)$, $C^b(X)$, $\C^c(X)$, $C^0(X)$}

\begin{prop}
% 7.B.2
\label{Pro-TraceCommC*}
Let $X$ be a second countable locally compact space.
Let $C^0(X)$ be the abelian C*-algebra of continuous functions on $X$
which vanish at infinity.
Let $\mu$ be a Radon measure on an open subset $U$ of $X$. 
Define $t_{\mu} \,\colon C^0(X)_+ \to \mathopen[ 0, \infty \mathclose]$ by 
$$
t_{\mu} (f) \, = \, \int_U f(x) d\mu(x)
\hskip.5cm \text{for every} \hskip.2cm
f \in C^0(X)_+.
$$
\begin{enumerate}[label=(\arabic*)]
\item\label{iDEPro-TraceCommC*}
$t_\mu$ is a lower semi-continuous and semi-finite trace on $C^0(X)$;
\item\label{iiDEPro-TraceCommC*}
$t_\mu$ is finite if and only if $\mu$ is finite;
\item\label{iiiDEPro-TraceCommC*}
$t_\mu$ is faithful if and only if $\mu(V\cap U) \ne 0$ 
for every non-empty open subset $V$ of $X$.
\item\label{ivDEPro-TraceCommC*}
Let $t$ be a lower semi-continuous and semi-finite trace on $C^0(X)$. Then
$t=t_\mu$ for a Radon measure $\mu$ on an open subset $U$ of $X$.
\end{enumerate}
\end{prop}

\begin{proof}
\ref{iDEPro-TraceCommC*}
Let $(f_n)_{n \ge 1}$ be a sequence in $C^0(X)_+$ and $f \in C^0(X)_+$
be such that $\lim_n f_n = f$.
Then, by Fatou's lemma, we have
$$
t_\mu (f) \, = \, \int_U \liminf f_n(x) d\mu(x) \, \le \,
\liminf \int_U f_n(x) d\mu(x) \, = \, \liminf t_\mu (f_n),
$$
and the lower continuity of $t_\mu$ follows.
\par

Choose an increasing sequence $(\varphi_n)_{n \ge 1}$
of functions $\varphi_n \in C^c(X)_+$ with supports contained in $U$
and such that $\lim _n \varphi_n (x) = 1$ for every $x \in U$.
Let $f \in C^0(X)_+$. We have $\varphi_n f \le f$;
moreover, for every $x \in U$, the sequence $(\varphi_n(x)f(x))_{n \ge 1}$ 
is increasing and converges to $f(x)$.
Hence, by Lebesgue monotone convergence theorem, we have 
$$
\lim_n t_{\mu}(\varphi_n f) \, = \, \lim_n \int_U \varphi_n(x)f(x) d\mu(x)
\, = \, \int_U f (x) d\mu(x) \, = \, t_\mu (f).
$$
Since $\varphi_n f$ has compact support contained in $U$
and since $\mu$ is a Radon measure on~$U$,
we have $t_\mu (\varphi_n f) < +\infty$
for every $n \ge 1$. This shows that $t_\mu$ is semi-finite.

\vskip.2cm

Item \ref{iiDEPro-TraceCommC*} is obvious. is clear

\vskip.2cm

\ref{iiiDEPro-TraceCommC*}
Assume that $t_\mu$ is faithful
and let $V$ be a non-empty open subset of $X$.
Let $f \in C^c(X)_+, f \ne 0$, with support contained in $V$.
Then $t_\mu (f) \ne 0$ and it follows that $\mu(V \cap U) \ne 0$.
\par

Assume that $\mu(V \cap U) \ne 0$ 
for every non-empty open subset $V$ of $X$.
Let $f \in C^0(X)_+$ with $f \ne 0$.
The set $V_f = \{ x \in X \mid f(x) > 0 \}$ is open and non-empty;
so, $\mu(V_f\cap U) > 0$ and therefore
$$
t_\mu (f) \, = \, \int_{V_f \cap U} f (x) d\mu(x) \, > \, 0
$$
This shows that $t_\mu$ is faithful.

\vskip.2cm

\ref{ivDEPro-TraceCommC*}
Let $\mathfrak{m_t}$ be the ideal of definition of $t$,
that is, the linear span of
$$
\mathfrak{m}_{t,+} \, = \, \{ f \in C^0(X)_+ \mid t(f) \, < \, +\infty \}
$$
(see Section~\ref{Section-RepTraceC*}).
The closure $J$ of $\mathfrak{m_t}$ is a closed ideal in $C^0(X)$.
Therefore, there exists a closed subset $A$ of $X$ such that 
$$
J \, = \, \{ f \in C^0(X) \mid f \vert_A = 0 \}.
$$
Set $U := X \smallsetminus A$.
We identify $J$ with $C_0(U)$ 
and so $\mathfrak{m_t}$ with a dense ideal of $C^0(U)$.
\par

By Lemma~\ref{Lem-DenseIdeal}, $\mathfrak{m_t}$ contains $C^c(U)$.
Hence, $t \vert_{C^c(U)}$ is a positive linear form on $C^c(U)$ and there exists
therefore a Radon measure $\mu$ on $U$ such that
$$
t(f) \, = \, \int_U f(x) d\mu(x)
\hskip.5cm \text{for all} \hskip.2cm
f \in C^c(U).
$$
\par

We claim that $t = t_\mu$.
Indeed, let $(\varphi_n)_{n \ge 1}$ be an increasing sequence in $C^c(U)_+$
such that $\lim_n \varphi_n=1$ on $U$.
\par

Let $f \in C^0(X)_+$.
Then $\varphi_n f \in C^c(U)_+$ and $\varphi_n f \le f$ for every $n \ge 1$;
moreover we have $\lim_n \varphi_n f(x) = f(x)$ for every $x \in U$.
So, by Lebesgue monotone convergence theorem, we have
$$
\lim_n t(\varphi_nf) \, = \, \lim_n \int_U \varphi_n(x)f(x) d\mu(x) \, = \,
\int_U f(x) d\mu(x) \, = \, t_\mu(f).
$$
It suffices therefore to show that $\lim_n t(\varphi_n f) = t(f)$
for every $f \in C^0(X)_+$.
\par

Assume first that $f \in C^0(U)_+$.
Then $\lim_n \varphi_n f = f$ and hence, by the lower semi-continuity of $t$, we have
$$
\liminf t(\varphi_n f) \, \ge \, t(f);
$$
so, $\lim_n t( \varphi_n f) = t(f)$, since $t( \varphi_n f)) \le t(g)$.
\par

Next, let $f \in C^0(X)_+$ such that $f \notin C^0(U)_+$.
Then $t(f) = +\infty$.
Let $C >0$.
Since $t$ is semi-finite, there exists
$g \in \mathfrak{m}_{t,+}$ such that 
$$
g \, \le \, f
\hskip.5cm \text{and} \hskip.5cm
C < t(g).
$$
By what we have seen above, 
we have $\lim_n t(\varphi_n g) = t(g)$. 
Hence, there exists $n_0$ such that 
$$
C \, < \, t(\varphi_n g)
\hskip.5cm \text{for all} \hskip.2cm
n \ge n_0 .
$$
Since $\varphi_n g \le \varphi_n f$, we have therefore
$$
C \, < \, t(\varphi_n f)
\hskip.5cm \text{for all} \hskip.2cm
n \ge n_0 
$$
and so 
$$
\lim_n t(\varphi_nf) \, = \, +\infty \, = \, t(f).
$$
\end{proof}

\begin{rem}
% 7.B.3
\label{Rem-Pro-TraceCommC*}
The result in Proposition~\ref{Pro-TraceCommC*}
admits an appropriate generalization which holds for any C*-algebra of type I;
see Corollaire 2 of Th\'eor\`eme 2 in \cite{Dixm--63}.
\end{rem}

Every C*-algebra $\mathcal{A}$ has a unique minimal dense ideal, called
the Pedersen ideal of $\mathcal{A}$ (see \cite[\S~5.6]{Pede--79}).
The following elementary lemma identifies this ideal 
in the case where $\mathcal{A}$ is commutative.
The lemma has been used in the proof of Proposition~\ref{Pro-TraceCommC*}.

\begin{lem}
% 7.B.4
\label{Lem-DenseIdeal}
Let $X$ be a locally compact topological space.
\par

Every dense ideal of $C^0(X)$ contains $C^c(X)$.
\end{lem}

\begin{proof}
Let $J$ be a dense ideal in $C^0(X)$.
Let $f \in C^c(X)$; denote by $K$ the support of $f$.
Let $x \in K$. There exists $g_x \in J$ such that $g_x(x) \ne 0$
and we can therefore find a neighbourhood $U_x$ of $x$ such 
that $g_x > 0$ on $U_x$.
By compactness of $K$, there exists $x_1, \dots, x_n \in K$ such that
$$
K \, \subset \, \bigcup_{i = 1}^n U_{x_i}.
$$
Then $g := \sum_{i = 1}^n \vert g_{x_i} \vert^2$ belongs to $J$ and $g > 0$ on $K$. 
Define $h \,\colon X \to \C$ by 
$$
h(x) \, = \, \begin{cases}
\dfrac{f(x)}{g(x)} 
& \text{for $x \in K$,}
\\
0
& \text{for $x \notin K$.}
\end{cases}
$$
Then $h$ is continuous on $X$ with compact support and therefore
$f = hg$ belongs to~$J$.
\end{proof}

Traces on C*-algebras are addressed again below,
in Section \ref{Section-RepTraceC*}.

\subsection*{Traces on von Neumann algebras}

\index{Trace! $3$@on a von Neumann algebra}
\index{Normal trace}
Let now $\mathcal M$ be a von Neumann algebra;
in particular, it is a C*-algebra, and the definitions above apply.
A trace $t$ on $\mathcal M$ is \textbf{normal} if, 
for every increasing filtering subset $\mathcal F \subset \mathcal M_+$
with least upper bound $x \in \mathcal M_+$,
the value $t(x)$ is the least upper bound of 
$\{ t(y) \in \mathopen[ 0, \infty \mathclose] \mid y \in \mathcal F \}$.
\par

A trace on $\mathcal M$ is normal if and only if it is lower semi-continuous
with respect to the weak topology
% \cite[Remarque 2]{Dixm--63}.
\cite[Chap.~I, \S~6, exercice~1]{Dixm--vN}.

\begin{exe}
% 7.B.5
\label{ExampleTracevN}
(1)
If $\Hi$ is a separable Hilbert space,
the standard trace $\mathrm Tr$ on $\Li (\Hi)$ 
is semi-finite, normal, and faithful.
Recall that the standard trace is defined by
$\mathrm Tr (x) = \sum_{i \ge 1} \langle x e_i \mid e_i \rangle
\in \mathopen[ 0, \infty \mathclose]$ for $x \in \Li (\Hi)_+$,
where $(e_i)_{i \ge 1}$ is an orthonormal basis of $\Hi$;
the sum is independent of the choice of the basis.
\par

Every semi-finite normal trace on $\Li (\Hi)$
is proportional to ${\mathrm Tr}$;
see \cite[Chap.~I, \S~6, no~6]{Dixm--vN}.
\par

The standard trace on $\Li (\Hi)$
is finite if and only if $\Hi$ is finite-dimensional.
\index{$h4$@$\Li (\Hi)$ bounded linear operators on $\Hi$}

\vskip.2cm

If $\Hi$ is infinite-dimensional, the construction 
referred to in Example \ref{ExampleTraceC*}(4)
provides ``Dixmier traces'' on $\Li (\Hi)$
that are semi-finite and not normal \cite{Dixm--66}.

\vskip.2cm

(2)
Let $\Gamma$ be a discrete group 
and $\lambda_\Gamma$ the left regular representation
of $\Gamma$ on $\ell^2(\Gamma)$. The map
$$
t \, \colon \, T \mapsto \langle T \delta_e \mid \delta_e \rangle
$$
is a faithful finite normal trace on the von Neumann algebra 
$\Li (\Gamma) \subset \Li (\ell^2(\Gamma))$
generated by $\lambda_\Gamma(\Gamma)$,
as shown in Proposition~\ref{Pro-TraceDiscreteGroup}.
\index{Left regular representation}
\index{$j2$@$\Li (\Gamma)$ von Neumann algebra associated to
the left regular representation of a group $\Gamma$}
\par

Observe that $t$ coincides on the group algebra $\C[\Gamma]$ with
the trace of Example \ref{ExampleTraceC*}(5).

\vskip.2cm

(3)
On an abelian von Neumann algebra, finite traces need not be normal.
Indeed, by the Hahn--Banach theorem, 
there exist non-zero positive linear forms 
on $L^\infty(\mathopen[ 0,1 \mathclose], dt)$
which vanish on the closed subspace $C(\mathopen[ 0,1 \mathclose])$
of continuous functions.
Such a form is a finite trace on $L^\infty(\mathopen[ 0,1 \mathclose], dt)$;
it is not normal since $C(\mathopen[ 0,1 \mathclose])$
is weakly dense in $L^\infty(\mathopen[ 0,1 \mathclose], dt)$.
\end{exe}

\subsection*{Types of von Neumann algebras}

A von Neumann algebra $\mathcal M$ is
\begin{enumerate}[label=(\Alph*)]
\item\label{ADEdefvnA}
\textbf{type I} or \textbf{discrete}
if it is isomorphic to a von Neumann algebra with an abelian commutant;
\index{von Neumann algebra! $2$@type I} 
\index{Type I! $4$@von Neumann algebra}
% \cite[I.8.1]{Dixm--vN}
\item\label{BDEdefvnA}
\textbf{continuous} 
if there does not exist any projection $p \ne 0$ in
the centre $\mathcal M \cap \mathcal M'$ such that $p\mathcal M p$ is discrete;
\index{von Neumann algebra! $3$@continuous}
% \cite[I.8.1]{Dixm--vN}
\item\label{CDEdefvnA}
\textbf{finite} [respectively \textbf{semi-finite}] 
if, for every $x \in \mathcal M_+$, $x \ne 0$,
there exists a normal trace $t$ on $\mathcal M$
that is finite [respectively semi-finite] and such that $t(x) > 0$;
\index{von Neumann algebra! $4$@finite}
\index{von Neumann algebra! $5$@semi-finite}
\item\label{DDEdefvnA}
\textbf{properly infinite} [respectively \textbf{purely infinite}]
if there does not exist any non-zero finite normal trace on $\mathcal M$
[respectively any semi-finite normal trace on $\mathcal M$];
\index{Properly infinite! $1$@von Neumann algebra}
\index{Purely infinite! $1$@von Neumann algebra}
\index{von Neumann algebra! $9$@properly infinite}
\index{von Neumann algebra! $9$@purely infinite}
\item\label{EDEdefvnA}
\textbf{type II$_1$} if it is continuous and finite;
\index{von Neumann algebra! $6$@type II$_1$}
\index{Type II, III! $1$@von Neumann algebra}
\item\label{FDEdefvnA}
\textbf{type II$_\infty$} if it is continuous, semi-finite, and properly infinite;
\index{von Neumann algebra! $7$@type II$_\infty$}
\item\label{GDEdefvnA}
\textbf{type III}
if it is purely infinite, i.e., 
if the only normal semi-finite trace on $\mathcal M_+$ is the zero trace;
\index{von Neumann algebra! $8$@type III}
\end{enumerate}
For factors (see Section~\ref{SS:FactRep+}), the typology can be refined:
a factor of type I is of \textbf{type~I$_n$} 
if it is isomorphic to the algebra ${\rm M}_n(\C)$
of $n$-by-$n$ complex matrices, for some integer $n \ge 1$,
and of \textbf{type I$_\infty$} if it is infinite-dimensional, 
isomorphic to the algebra $\Li (\Hi)$ 
for some infinite-dimensional Hilbert space $\Hi$.
Moreover, type III can be subdivided in types III$_\lambda$, 
for $\lambda \in \mathopen[ 0,1 \mathclose]$;
this will appear briefly in Chapters
\ref{ChapterGroupMeasureSpace} and \ref{Chap:NormalInfiniteRep}.

\begin{rem}
% 7.B.6
\label{Rem-TypI-Semifinite}
(1)
Type I von Neumann algebras can also be characterized as follows:
a von Neumann algebra $\mathcal M$ is of type I
if and only if there exists a projection $p\in \mathcal M$ with central support $I$
such that $p\mathcal M p$ is abelian
\cite[Chap.~I, \S~8, Th\'eor\`eme 1]{Dixm--vN}.

\vskip.2cm

(2)
Finiteness and semi-finiteness have other characterizations in terms of traces:
A von Neumann algebra $\mathcal M$ is semi-finite
if and only if it has a \emph{faithful} semi-finite normal trace.
When $\mathcal M$ is countably decomposable, $\mathcal M$ is finite
if and only if it has a faithful finite normal trace \cite[Chap.~I, \S~6, no~7]{Dixm--vN}.
(A von Neumann algebra $\mathcal M$ is countably decomposable 
if every collection of pairwise orthogonal non-zero projections in $\mathcal M$ is countable.
For example, a von Neumann algebra that can act faithfully
on a separable Hilbert space is countably decomposable.)

\vskip.2cm

(3) The type of a von Neumann algebra $\mathcal M$
inside some algebra $\Li (\Hi)$ is closely related
to the type of its commutant $\mathcal M'$:
\begin{enumerate}
\item[--]
$\mathcal M$ is of type I if and only if $\mathcal M'$ is of type I
\cite[Chap.~I, \S~8, no~2]{Dixm--vN};
\item[--]
$\mathcal M$ is of type II if and only if $\mathcal M'$ is of type II
(in \cite[Chap.~I]{Dixm--vN}, see
\S~6, no~8, cor.~1 of prop.~13,
and \S~8, no~2, cor.~1 of th.\ 1);
% pages 104 and 124 resp.
\item[--]
$\mathcal M$ is of type III if and only if $\mathcal M'$ is of type III
\cite[Chap.~I, \S~6, no~8]{Dixm--vN}.
% \cite[Chap.~I, \S~6, no~8, cor.~3 of prop.~13]{Dixm--vN}
\end{enumerate}
\end{rem}

The following proposition reduces the study of general von Neumann algebras
to that of von Neumann algebras of one of the types introduced above.

\begin{prop}
% 7.B.7
\label{sommetypeavN}
Every von Neumann algebra $\mathcal M$ is canonically a product
$$
\mathcal M_{\rm I_f} \times \mathcal M_{{\rm I}_\infty} \times 
\mathcal M_{{\rm II}_1} \times \mathcal M_{{\rm II}_\infty} \times
\mathcal M_{\rm III}
$$
where the five terms in the product
are respectively of type I and finite, of type I and infinite,
of type II$_1$, of type II$_\infty$, and of type III
(and each can be absent).
\end{prop}

\begin{proof}[Reference for the proof]
In \cite[Chap.~I]{Dixm--vN}, see
\S~6, no~7, cor.~1 of prop.~8
and \S~8, no~1, cor.~1 of prop.~1.
\end{proof}

\section{Types of group representations}
% Section 7.C
\label{SectionvN+rep}

Let $G$ be a topological group.
Recall from Section~\ref{S:QE-Factor-VNAlgebras} that,
for every representation $(\pi, \Hi)$ of $G$ in a Hilbert space $\Hi_\pi$,
there are two associated von Neumann subalgebras of $\Li (\Hi_\pi)$,
which are commutant of each other: $\pi(G)''$ and $\pi(G)'$.
\par

We now relate type I representations introduced in Section~\ref{Sectioncomppres}
to type I von Neumann algebras.

\begin{prop}
% 7.C.1
\label{Pro-TypeIVN-bis}
For a representation $\pi$ of $G$,
the following properties are equivalent:
\begin{enumerate}[label=(\roman*)]
\item\label{iDEPro-TypeIVN-bis}
$\pi$ is of type I;
\item\label{iiDEPro-TypeIVN-bis}
the von Neumann algebra $\pi(G)''$ is of type I;
\item\label{iiiDEPro-TypeIVN-bis}
the von Neumann algebra $\pi(G)'$ is of type I.
\end{enumerate}
\end{prop}

\begin{proof}
The fact that \ref{iDEPro-TypeIVN-bis} implies \ref{iiDEPro-TypeIVN-bis}
follows from Proposition~\ref{Pro-TypeIVN} and the definition of a type I von Neumann algebra.
The fact that a von Neumann algebra and its commutant are together of type I or not
(see Remark~\ref{Rem-TypI-Semifinite}) shows that
\ref{iiDEPro-TypeIVN-bis} and \ref{iiiDEPro-TypeIVN-bis} are equivalent.
\par

It remains to show that \ref{iiiDEPro-TypeIVN-bis} implies \ref{iDEPro-TypeIVN-bis}.
Assume that $\pi(G)'$ is of type I.
Then there exists a projection $p \in \mathcal \pi(G)'$ with central support $I$
such that $p \mathcal \pi(G)' p$ is abelian (see Remark~\ref{Rem-TypI-Semifinite}).
Let $\sigma$ be the subrepresentation of $\pi$ on the range of~$p$. 
Since $\sigma(G)' \cong p\mathcal \pi(G)' p$, the representation $\sigma$ is multiplicity-free 
(see Proposition~\ref{Pro-MultFreeVN}).
As the central support of $p$ is $I$, it follows from Proposition~\ref{Prop-DirectSumTypeIRep}
that the representation $\pi$ is of type I. 
\end{proof}

We can now define various types of group representations.

\begin{defn}
% 7.C.2
\label{Def-TypesRep}
The representation $\pi$ of $G$ is said to be of type I, II, II$_1$, II$_\infty$, or III
if the von Neumann algebra $\pi(G)''$ is of type I, II, II$_1$, II$_\infty$, or III.
\par

The representation $\pi$ is said to be
of discrete, continuous, finite, semi-finite, properly infinite, or purely infinite type 
if the von Neumann algebra $\pi(G)''$ is
of discrete, continuous, finite, semi-finite, properly infinite, or purely infinite type.
\index{Representation! type II, II$_1$, II$_\infty$, III}
\index{Representation! finite type}
\index{Type II, III! $2$@representation}
\index{Properly infinite! $2$@representation}
\index{Purely infinite! $2$@representation}
\end{defn}

\begin{rem}
% 7.C.3
\label{Rem-Def-TypesRep}
Let $\pi$ be a representation of $G$. 

\vskip.2cm

(1)
It follows from the definitions
that every representation quasi-equivalent to $\pi$
is of the same type as $\pi$.
See Proposition~\ref{Pro-SubordinateVN}.

\vskip.2cm

(2)
Our notion of representation of finite type does \textbf{not} agree
with the notion of ``finite representation" from \cite[Page 30]{Mack--76},
which is not invariant by quasi-equivalence.

\vskip.2cm

(3)
$\pi$ is of continuous type if and only if, for every non-zero subrepresentation $\rho$ of $\pi$,
there are two equivalent subrepresentations $\rho', \rho''$ of $\rho$ in orthogonal invariant subspaces, such that $\rho' \oplus \rho'' = \rho$.

\vskip.2cm

(4)
If $\pi$ is a representation of type III on a separable Hilbert space,
then every representation that is quasi-equivalent to $\pi$ is equivalent to $\pi$
\cite[5.6.6]{Dixm--C*}.
\end{rem}

As for von Neumann algebras (Proposition~\ref{sommetypeavN}),
the study of general representations can be reduced to that 
of representations of specific types.

\begin{prop}
% 7.C.4
\label{sommetyperep}
Every representation $\pi$ of a topological group is canonically a direct sum
$$
\pi_{\rm I_f} \oplus \pi_{{\rm I}_\infty} \oplus 
\pi_{{\rm II}_1} \oplus \pi_{{\rm II}_\infty} \oplus \pi_{\rm III}
$$
where the representations in the direct sum
are respectively of type I and finite, of type I and infinite,
of type II$_1$, of type II$_\infty$, and of type III
(some of these can be in spaces of dimension $0$).
\end{prop}

With the notions we have introduced,
we can rephrase parts of Propositions~\ref{iccfactorII1} and \ref{Pro-TraceDiscreteGroup}
in the following way.

\begin{prop}
% 7.C.5
\label{Prop-TypRegRepDiscrete}
Let $\Gamma$ be a discrete group.
\begin{enumerate}[label=(\arabic*)]
\item\label{iDEProp-TypRegRepDiscrete}
The regular representation $\lambda_\Gamma$ is a representation of finite type.
\item\label{iiDEProp-TypRegRepDiscrete}
If $\Gamma \ne \{e\}$, the following three conditions are equivalent:
\begin{itemize}
\setlength\itemsep{0em}
\item[---]
$\Gamma$ is icc,
\item[---]
$\lambda_\Gamma$ is a factor representation,
\item[---]
$\lambda_\Gamma$ is a factor representation of type II$_1$.
\end{itemize}
\end{enumerate}
\end{prop}.

\section{Types of representations of discrete groups}
% Section 7.D
\label{SectionTypesI+IId}

All groups appearing in this section will be discrete groups.

\subsection*{Thoma's theorem}
% subsection 7.D.a

We have seen before
% (see Proposition~\ref{Prop-TypRegRepDiscrete})
% that an icc discrete group is not of type I.
that a discrete group need not be of type I.
In fact, we have the following characterization. 

\begin{theorem}[Thoma]
% 7.D.1
\label{discretTypes-Thoma}
Let $\Gamma$ be a discrete group.
The following properties are equivalent:
\begin{enumerate}[label=(\roman*)]
\item\label{iDEdiscretTypes-Thoma}
$\Gamma$ is type I;
\item\label{iiDEdiscretTypes-Thoma}
$\Gamma$ is virtually abelian, that is,
$\Gamma$ contains an abelian subgroup of finite index.
\end{enumerate}
\end{theorem}

For a proof, we refer to the original article \cite{Thom--68}
(see also \cite{Thom--64a} for the case of countable groups).
The deep implication is
\ref{iDEdiscretTypes-Thoma} $\implies$ \ref{iiDEdiscretTypes-Thoma}.

\subsection*{Type of the regular representation}
% subsection 7.D.b

In complement of Theorem \ref{discretTypes-Thoma},
here are other results about the type of the regular representation $\lambda_\Gamma$
of a discrete group $\Gamma$. 
We denote by ${\rm D}(\Gamma)$ the commutator group of $\Gamma$,
and by $\FC(\Gamma)$ its FC-centre, i.e.,
the normal subgroup of elements in $\Gamma$ 
with a finite conjugacy class (see Appendix \ref{AppGpactions}).

\begin{theorem}
% 7.D.2
\label{discretTypes}
Let $\Gamma$ be a group.
\begin{enumerate}[label=(\arabic*)]
\item\label{iDEdiscretTypes}
$\Gamma$ is of type I if and only if $\lambda_\Gamma$ is of type I.
\item\label{iiDEdiscretTypes}
$\lambda_\Gamma$ is of type II
if and only if either $[\Gamma \,\colon \FC(\Gamma)] = \infty$,
or $[\Gamma \,\colon \FC(\Gamma)] < \infty$ and ${\rm D}(\FC(\Gamma))$ is infinite.
\item\label{iiiDEdiscretTypes}
In particular, if $\Gamma$ is not virtually abelian,
$\FC(\Gamma)$ of finite index in $\Gamma$,
and ${\rm D}(\FC(\Gamma))$ is a finite group, then $\lambda_\Gamma$
has both a type I component and a type II component.
\item\label{ivDEdiscretTypes}
If $\Gamma$ is finitely generated, 
then $\lambda_\Gamma$ is either of type I or of type II.
\end{enumerate}
\end{theorem}

\noindent Note:
\emph{Since $\lambda_\Gamma$ is finite, one can replace
``type II'' by ``type II$_1$'' in \ref{iiDEdiscretTypes} to \ref{ivDEdiscretTypes}.}

\begin{proof}[References for proofs]
Parts \ref{iDEdiscretTypes} and \ref{iiDEdiscretTypes}
are S\"atze 2 and 1 in \cite{Kani--69}.
% (see also \cite{Smit--72}).
% $\FC(\Gamma)$ d'indice infini implique type II d\'ej\`a dans [Mautner, Duke, 1951]
Claim \ref{iiiDEdiscretTypes} follows.
\par

Claim \ref{ivDEdiscretTypes} follows from 
Theorem~\ref{discretTypes-Thoma} \& \ref{iiDEdiscretTypes},
and the following classical fact: a finitely generated group $\Gamma$
with $\FC(\Gamma)$ of finite index is virtually abelian
(see Proposition \ref{dueToSchur}),
and therefore of type I.
% See also \cite[13.11.14]{Dixm--C*}.
\end{proof}

\begin{rem}
% 7.D.3
(1)
Let $\Gamma$ be a countable group which is not virtually abelian.
Then $\Gamma$ has factor representations of type II$_1$, and of type III,
see Theorem \ref{ThGlimm2} below.
% D'apr\`es \cite{Guic--63}, $\BS(1, 2)$, un groupe icc, 
% a des repr\'esentations factorielles de type II$_\infty$
% In particular, the inclusion map (\ref{eqq/d/qd}) of the dual of $\Gamma$ 
% into the quasi-dual of $\Gamma$
% (see \ref{SectionQuasidual}) is not onto.

\vskip.2cm

(2)
Claim \ref{iDEdiscretTypes} of the theorem
\emph{does not} carry over to LC groups:
there are groups with left regular representation of type I that are not of type I
(see Section \ref{SectionTypesI+IILC}).
\index{Left regular representation}
\end{rem}

\begin{exe}
% 7.D.4
\label{typesmainex}
Except the infinite dihedral group,
our main examples of Chapter \ref{Chapter-ExamplesIndIrrRep}
are not virtually abelian,
and therefore not of type I.
We point out below to some explicit examples
of factor representations of type II$_1$ of these groups.

\vskip.2cm

\index{Heisenberg group}
(1)
The Heisenberg group $H(\Z)$ of Section \ref{Section-IrrRepTwoStepNil}
is not of type I.
For factor representations of $H(\Z)$ of types II and III,
we refer to Example \ref{Exa-FacRepSemiDirect}.

\vskip.2cm

(2)
Let $\Gamma$ be a icc group,
so that the left regular representation of $\Gamma$ is factorial of type II$_1$.
Assume moreover that $\Gamma$ is amenable (for example solvable).
Then $\Li (\Gamma)$ is \emph{the} hyperfinite factor of type II$_1$,
by the famous result of Connes on the uniqueness of
the hyperfinite factor of type II$_1$ \cite[see in particular Corollary 7.2]{Conn--76}.
This applies to
\begin{enumerate}
\item[--]
$\Aff(\K)$, the affine group over an infinite field $\K$,
as in Section \ref{Section-IrrRepAff},
\item[--]
the Baumslag--Solitar group $\BS(1,n)$ for every integer $n \ge 2$,
as in Section~\ref{Section-IrrRepBS},
\item[--]
the restricted wreath product $H \wr A$ whenever $H$ and $A$
are non-trivial abelian countable groups,
for example the lamplighter group $\Z \wr (\Z / 2 \Z)$
of Section \ref{Section-IrrRepLamplighter}.
\end{enumerate}
\index{Affine group! $4$@$\Aff(\K)$ of an infinite field $\K$} 
\index{Baumslag--Solitar group $\BS(1, p)$}
\index{Lamplighter group}

\vskip.2cm

\index{General linear group! $\GL_n(\K)$ with $\K$ a field}
(3)
For an infinite field $\K$ and an integer $n \ge 2$,
the general linear group $\Gamma = \GL_n(\K)$
of Section \ref{Section-IrrRepGLN} is not type I.
The quotient group $\PGL_n(\K) = \GL_n)\K) / \K^\times$ can be shown to be icc,
for example because it is Zariski dense
in an appropriate algebraic group with trivial centre
(compare \cite[Proposition 2]{BeHa--94}),
so that its left regular representation is of type II$_1$;
it follows that the quasi-regular representation
$\lambda_{\GL_n(\K) / \K^\times}$ of $\GL_n(\K)$
in $\ell^2(\PGL_n(\K))$ is factorial of type II$_1$.

\vskip.2cm

(4)
The question to know whether appropriate quasi-regular representations
provide examples of factor representation of types II and III
has attracted some attention
(quasi-regular representations have been defined just before \ref{Cor-QuasiReg}).
Here is one example:
\index{Quasi-regular representation}
\index{Representation! quasi-regular}
\par

For the integral translation subgroup
of the rational affine group
$$
T \, = \, \begin{pmatrix} 1 & \Z \\ 0 & 1 \end{pmatrix}
\hskip.2cm \subset \hskip.2cm
\Aff(\Q) = \begin{pmatrix} \Q^\times & \Q \\ 0 & 1 \end{pmatrix} ,
$$
the quasi-regular representation $\lambda_{\Aff(\Q)/T}$ is of type III
\cite[Proposition 3.6]{Bind--93}.
\end{exe}

\subsection*{Groups with regular representation of mixed type}
% subsection 7.D.c

Let $\Gamma$ be a discrete group which is not virtually abelian;
assume that $[\Gamma \,\colon \FC(\Gamma)] < \infty$
and that ${\rm D}(\FC(\Gamma)$ is finite.
It follows from Theorem \ref{discretTypes}
that the regular representation of $\Gamma$ must involve representations of different types. 
The following result of Kaplansky \cite[Theorem~2]{Kapl--51a} gives 
the explicit type decomposition of $\lambda_\Gamma$ for a class of examples.
Recall that $\Li (\Gamma)$ denotes the von Neumann algebra generated by $\lambda_\Gamma$.

\begin{theorem}
% 7.D.5
\label{Theo-Kaplansky}
Let $\Gamma$ be an infinite discrete group with centre $Z$.
Assume that $Z$ is finite of prime order $p$ and that ${\rm D}(\Gamma) = Z$.
\par
Then 
$$
\Li (\Gamma) \, = \, \mathcal A \oplus \bigoplus_{i = 1}^{p-1} \mathcal M_i,
$$
where $\mathcal A$ is an abelian von Neumann algebra.
and the $\mathcal M_i$'s are factors of type II$_1$.
\end{theorem}

\begin{proof}
Let $z_0$ be a generator for the cyclic group $Z$ of order $p$.
The unitary dual $\widehat Z$ can be identified with the group
$C_p$ of $p$-th roots of unity:
to every $\omega \in C_p$
corresponds the unitary character $\chi_\omega \,\colon Z \to \T$
defined by $\chi_\omega(z_0^j) = \omega^j$ for $j \in \{0, 1, \hdots, p-1\}$.
\par

We have an orthogonal decomposition of $\ell^2(\Gamma)$ under
the action of $\lambda_{\Gamma}(z_0)$,
given by 
$$
\ell^2(\Gamma\,) = \, \bigoplus_{\omega \in C_p} \Hi_\omega,
$$
where
$$
\Hi_\omega \, = \, \{ f \in \ell^2(\Gamma) \mid f(z_0x) = \omega f(x)
\hskip.2cm \text{for every} \hskip.2cm
x \in \Gamma\}.
$$
Let $\omega \in C_p$.
Observe that $\Hi_\omega$ is invariant under $\lambda_\Gamma(\Gamma)$
as well as under $\rho_\Gamma(\Gamma)$.
The corresponding orthogonal projection $P_\omega \,\colon \ell^2(\Gamma) \to \Hi_\omega$
belongs to the centre of the von Neumann algebra $\Li (\Gamma)$
and is given by 
$$
P_\omega (f) (x) \, = \, \dfrac{1}{p} \sum_{i = 0}^{p-1} \omega^{-i}f(z_0^i x)
\hskip.5cm \text{for} \hskip.2cm
f \in \ell^2(\Gamma)
\hskip.1cm \text{and} \hskip.1cm
x \in \Gamma.
$$
Denote by $\pi_\omega$ the restriction of $\lambda_\Gamma$ to $\Hi_\omega$.
The von Neumann algebra $\pi_\omega (\Gamma)''$
coincides with the von Neumann algebra $\Li (\Hi)_{P_\omega}$ induced by $P_\omega$.
It is clear that 
$$
\pi \, = \, \bigoplus_{\omega \in C_p} \pi_\omega
$$ 
and that 
$$
\Li (\Gamma) \, = \, \bigoplus_{\omega \in C_p} \pi_\omega (\Gamma)''.
$$

\vskip.2cm

\noindent
$\bullet$ \emph{First case.} Assume that $\omega = 1$.
We claim that $\pi_1 (\Gamma)''$ is abelian. 
\par

Indeed, $\Hi_1$ can be identified with $\ell^2(\Gamma/Z)$ and $\pi_1$ with 
$\lambda_{\Gamma/Z}$. The claim follows, since $\Gamma/Z$ is abelian.

\vskip.2cm

\noindent
$\bullet$ \emph{Second case.} Assume that $\omega \ne 1$. 
We claim that $\pi_\omega$ is factorial of type II$_1$.
\par

Indeed, set 
$$
f_\omega \, := \, P_\omega (\delta_e) \in \Hi_\omega .
$$ 
Then $f_\omega$ is a cyclic vector for the restriction of $\rho_\Gamma$ to $\Hi_\omega$,
since $\delta_e$ is a cyclic vector for $\rho_\Gamma$
and since $P_\omega \in \rho_\Gamma(\Gamma)'$.
Therefore, the linear map 
$$
\pi_\omega(\Gamma)'' \to \Hi_\omega, \hskip.5cm T \mapsto T(f_\omega)
$$
is injective.
\par

We first prove that $\pi_\omega(\Gamma)''$ is a factor. 
Let $T \in \Li (\Hi_\omega)$ be in the centre of $\pi_\omega(\Gamma)''$.
In order to show that $T$ is a scalar multiple of the identity,
it suffices check that the function
$$
f \, := \, T(f_\omega) \in \Hi_\omega
$$
is a scalar multiple of $f_\omega$.
\par

As in the proof of Proposition ~\ref{iccfactorII1},
$f$ is invariant under conjugation by elements of $\Gamma$. 
We claim that $f(x) = 0$ for every $x \in \Gamma \smallsetminus Z$.
\par

Indeed, assume that $f(x)\ne0$ for some $x \in \Gamma \smallsetminus Z$.
Then there exists $y \in \Gamma$ with $yxy^{-1}\ne x$. Since 
${\rm D}(\Gamma) = Z$, we have therefore $yxy^{-1}= z_0^i x$ for some
$i \in \{1, \dots, p-1 \}$. It follows that 
$$
f(x) \, = \, f(yxy^{-1}) \, = \, f(z_0^i x) \, = \, \omega^i f(x)
$$
and this is a contradiction, since $\omega^i \ne 1$.
\par

So, $f = 0$ on $\Gamma \smallsetminus Z$.
Since $f \in \Hi_\omega$, we have, on the one hand,
$$
f(z_0^i) \, = \, \omega^i f(z_0)
\hskip.3cm \text{for every} \hskip.2cm
i \in \{0, \dots, p-1 \};
$$
on the other hand, a straightforward computation shows that
$$
f_\omega(x) \, = \, 
\begin{cases} 
\dfrac{\omega^i}{p}& \text{if} \hskip.2cm x = z_0^i
\hskip.2cm \text{for} \hskip.2cm
i \in \{0, \dots, p-1 \}
\\
0 & \text{otherwise.}
\end{cases}
$$
Therefore $f$ is a scalar multiple of $f_\omega$. 
Therefore, $\pi_\omega(\Gamma)''$ is a factor.
\par

Next, we show that the factor $\pi_\omega(\Gamma)''$ is of finite type.
The canonical trace of $\Li (\Hi)$ (as in Proposition~\ref{Pro-TraceDiscreteGroup})
defines by restriction a faithful trace on 
$$
\pi_\omega(\Gamma)'' \, \cong \, P_\omega \Li (\Hi) P_\omega.
$$
Therefore $\pi_\omega(\Gamma)''$ is a finite factor. 
\par

It remains to show that $\pi_\omega(\Gamma)''$ is infinite dimensional.
For this, it suffices to check that the image of $\pi_\omega(\Gamma)''$
under the map 
$$
\pi_\omega(\Gamma)'' \to \ell^2(\Gamma), \hskip.2cm T \mapsto T(f_\omega)
$$
is infinite dimensional.
This is indeed the case: since $\Gamma$ is infinite and $Z$ is finite,
we can find an infinite sequence $(x_n)_{n \ge 1}$ of elements $x_n \in \Gamma$
with $x_m Z\ne x_n Z$ for all $m,n$ with $m \ne n$.
Set 
$$
f_n \, := \, \pi_\omega(x_n) (f_\omega)
\hskip.5cm \text{for} \hskip.2cm
n \ge 1.
$$ 
One checks that $\mathrm{supp}(f_n) = x_n Z$.
Therefore the family $(f_n)_{n \ge 1}$ is linearly independent. 
\end{proof}

\begin{exe}
% 7.D.6
\label{HeisenbergFp}
A specific example as in Theorem~\ref{Theo-Kaplansky}.
is the ``Heisenberg group'' $H_\infty(\F_p)$ defined as follows.
Let $p$ be a prime with $p \ge 3$. Denote by $\F_p$ the field of order $p$,
by $V = \F_p^{(\N)}$ a vector space over $\F_p$ of countable infinite dimension,
by $(x_i)_{i \in \N}$ a typical element of $V$, 
with $x_i \in \F_p$ and $x_i = 0$ for almost all $i \in \N$,
and by $\omega$ the symplectic form on $V \oplus V$, defined by
\index{Heisenberg group! $5$@$H_\infty(\F_p)$}
$$
\omega((x,y),(x',y')) \, = \, \sum_{i \in \N} (x_iy'_i -y_ix'_i)
\hskip.5cm \text{for} \hskip.2cm 
(x,y), (x',y') \in V \oplus V.
$$
Then $H_\infty(\F_p)$ is the group with underlying set $V \oplus V \oplus \F_p$
and with multiplication defined by
$$
(x,y,z) (x',y',z') \, = \, (x+x', y+y', z+z'+ \omega((x,y),(x',y')))
$$
for $(x,y,z), (x',y',z') \in H_\infty(\F_p)$.
%(In some terminology, $H_\infty(\F_p)$ is an infinitely generated extra-special $p$-group.)
The additive group of $\F_p$ is identified with the subgroup of $H_\infty(\F_p)$
of elements of the form $(0, 0, z)$.
It is straightforward to check (using that $p$ is odd) that
\begin{enumerate}[label=(\arabic*)]
\item
the commutator group and the centre of $H_\infty(\F_p)$ both coincide with $\F_p$,
\item
all conjugacy classes of $H_\infty(\F_p)$ are finite,
i.e., $\FC(\Gamma) = \Gamma$,
\item
$H_\infty(\F_p)$ is not virtually abelian.
\end{enumerate}
It follows from Theorem~\ref{Theo-Kaplansky} that 
$\lambda(H_\infty(\F_p))''$ 
is a direct sum of an abelian von Neumann algebra
and $p-1$ factors of type II$_1$.
\end{exe}

\section[Types of representations of LC groups]
{On types of representations of locally compact groups}
% Section 7.E
\label{SectionTypesI+IILC}
We quote below some known results about the possible types of the regular representation 
for some classes of locally compact groups $G$.

\vskip.2cm

(1)
If $G$ is unimodular, then $\lambda_G$ has no type III part \cite{Sega--50b}.
Also, if $G$ is connected, then $\lambda_G$ has no type III part;
this is a fact established in \cite{Dixm--69}, 
building up on previous results for second-countable groups and solvable Lie groups
from \cite{Gode--51b}, \cite{Sega--50a}, \cite[Lemma 7.4]{Maut--50b}, 
and Puk\'anszky (\cite[Section IV]{Puka--71} and a previous announcement).
% cit\'e par \cite{Mack--76}, page 63. 
\par

If $G$ is connected, then $G$ does not have any representation of type II$_1$ \cite{KaSi--52}
(see Corollary~\ref{Cor-Theo-KadisonSinger}).
\vskip.2cm

(2)
If $G$ is a LC group, $\lambda_G$ is finite
if and only if $G$ a SIN group (see Theorem~\ref{Theo-SIN-FiniteVN}).

\vskip.2cm

(3)
More generally, for $G$ locally compact,
$\lambda_G$ has a non-trivial finite part
if and only if there exists at least one invariant compact neighbourhood of $e$.
Moreover, when this holds,
there exists a compact normal subgroup $K$ of $G$,
say with quotient $Q := G/K$, such that the finite part of $\lambda_G(G)''$
is isomorphic to $\lambda_Q(Q)''$
\cite[Section 4]{Tayl--76}.
There are similar results for the part of $\lambda_G(G)''$ that is type I and finite
\cite[Section 5]{Tayl--76}.

\vskip.2cm

(4)
A tensor product $\mathcal M \otimes \mathcal N$ is a factor of type II$_\infty$
whenever $\mathcal M$ is a factor of type II$_1$ 
and $\mathcal N$ a factor of type I$_\infty$.
As we have already noted in Remark~\ref{AffKLie},
the left regular representation of the affine group $\Aff(\R)$
is of type I$_\infty$. 
From this, Proposition~\ref{Prop-TypRegRepDiscrete},
and the previous observation,
it follows that, for any icc discrete group $\Gamma$,
the left regular representation of the Cartesian product
$\Gamma \times \Aff(\R)$ is of type II$_\infty$.
\index{Affine group! $1$@$\Aff(\R)$}

\vskip.2cm

(5)
Nilpotent connected LC groups are of type I, 
as mentioned in Theorem~\ref{explesTypeI}
but connected solvable LC groups need not be.
Consider for example the Mautner group, 
that is the semi-direct product $M := \C^2 \rtimes \R$ 
for the action of $\R$ on $\C^2$
given by $(z,z')t = (e^{2 \pi i t}z, e^{2 \pi i \lambda t}z')$ 
for all $t \in \R$ and $(z,z') \in \C^2$, for some irrational $\lambda \in \R$.
% def plus gŽnŽrale du groupe de Mautner ? voir peut-\^etre Kaniuth-Taylor page 151
This group is not of type I, 
see \cite[Page 138, and Theorem 9 Page 110]{AuMo--66};
indeed its regular representation is of type II$_\infty$ \cite[Proposition 4.1]{Kaji--83}.
Observe that the Mautner group is unimodular.
As much as we know, Mautner himself has not published anything more on this group
than some allusions in \cite{Maut--50b}.
% Plus pr\'ecis\'ement : Mautner, Ann. Math. 52, 1950, bas de la page 529).
\index{Mautner group $M := \C^2 \rtimes \R$}

\vskip.2cm

(6)
Mackey has found a LC group that has the following properties:
it is solvable, simply connected, not unimodular, not type I,
but its regular representation is type I \cite{Mack--61};
see also \cite[Example 3.e]{KlLi--73}.

\vskip.2cm

(7)
For examples of groups with $\lambda_G$ of type III,
see below Section \ref{Section-VNRegular}.

\section[Non type I factor representations]
{Non type I factor representations and irreducible representations}
% Section 7.F
\label{S:IntDecFactorRep}

Let $G$ be a second-countable LC group and $\pi$ a representation of $G$ in a Hilbert space.
Recall from Section~\ref{SectionDecomposingIrreps} that $\pi$ can be written
as a direct integral of irreducible representations $\pi\simeq \int^\oplus_X \pi_x d\mu(x)$
and that such a decomposition is usually not unique.
However, if $\pi$ is a factor representation,
the following result shows that almost all $\pi_x$'s are weakly equivalent to $\pi$.
% (see Section \ref{SectionWC+FellTop} for this notion).
The proof we give, which apparently is due to Fell, is patterned
after the one sketched in \cite[Remarque on page 100]{Dixm--60a}.

\begin{prop}
% 7.F.1
\label{Pro-IntDecFactorRep}
Let $G$ be a second-countable locally compact group 
and $x \mapsto \pi_x$ a measurable field of representations 
of $G$ over a standard Borel space $X$
equipped with $\sigma$-finite measure $\mu$.
Assume that the direct integral $\pi := \int^\oplus_X \pi_x d\mu(x)$ is a factor representation.
\par

Then there exists a Borel subset $N$ of $X$ with $\mu(N) = 0$ such that 
$\pi_x$ is weakly equivalent to $\pi$, for every $x \in X \smallsetminus N$.
\end{prop}

\begin{proof}
% Let $x \mapsto \Hi_x$ be the measurable field of Hilbert spaces on $X$ defining $\pi$, so that 
% $$
% \Hi \, := \, \int^\oplus_X \Hi_x d\mu(x)
% $$ 
% is the Hilbert space of $\pi$. 
% \par
%
Let $E \subset X$ be a Borel subset with $\mu(E) > 0$. Let 
$$
\pi_E \, := \, \int_E \pi_x d\mu(x)
$$
be the direct integral of the $\pi_x$'s over $E$. 
It is clear that $\pi_E$ is a non trivial subrepresentation of $\pi$.
Therefore since $\pi$ is factorial, $\pi_E$ is quasi-equivalent to $\pi$ 
by Corollary~\ref{Cor-FactorRepSubRep}.
It follows that $\pi_E$ is weakly equivalent to $\pi$ (Corollary~\ref{QE-WEQ}).
Therefore, on the one hand, we have, by Proposition~\ref{Pro-WeakContainmentOperatorNorm},
$$
\Vert \pi_E(f) \Vert \, = \, \Vert \pi(f) \Vert
\hskip.2cm \text{for every} \hskip.2cm
f \in C^c(G).
$$
On the other hand, since $\pi_E(f)$ is a decomposable operator, $\Vert \pi_E(f)\Vert$
is the $\mu$-essential supremum over $E$ of the function $x \mapsto \Vert \pi_x(f)\Vert$
(see Section~\ref{Section-DecomposableOperators}).
\par

We have shown that $\Vert \pi(f) \Vert$
coincides with the $\mu$-essential supremum over $E$ of $x \mapsto \Vert \pi_x(f)\Vert$
for every $f \in C^c(G)$ and every Borel subset $E$ of $X$ with $\mu(E)>0$.
\par

Since $G$ is second-countable, we can find
a sequence $(f_n)_{n \ge 1}$ in $C^c(G)$ which is dense for the $L^1(G)$-norm.
For every pair of integers $n,m \ge 1$, set
$$
E_{n,m} \, := \, \left\{ x \in X \mid \Vert \pi_x(f_n)\Vert \le \Vert \pi(f_n)\Vert - \frac{1}{m}\right\}.
$$
Then $E_{n,m}$ is a measurable subset of $X$ and it follows from what we have just shown
that $\mu(E_{n,m}) = 0$. 
Let
$$
N= \bigcup_{n,m} E_{m,n}.
$$
 Then $\mu(N) = 0$ and for $x \in X \smallsetminus N$,
we have
$$
\Vert \pi_x(f_n)\Vert \, = \, \Vert \pi(f_n)\Vert
\hskip.2cm \text{for every} \hskip.2cm
n \ge 1.
$$
Therefore 
$\pi_x$ is weakly equivalent to $\pi$ for $x \in X \smallsetminus N$, 
by Proposition~\ref{Pro-WeakContainmentOperatorNorm} again.
\end{proof}

The following corollary of Proposition~\ref{Pro-IntDecFactorRep}
establishes an important link between the quasi-dual $\QD(G)$
and the primitive dual $\Pri(G)$ of $G$,
introduced respectively in Sections~\ref{SectionQuasidual}
and \ref{PrimIdealSpace}.
Recall that $ \Pri(G)$ is in canonical way a quotient of $\widehat G$.

\begin{cor}
% 7.F.2
\label{Cor-FacRepPrimitiveG}
Let $G$ be a second-countable locally compact group and $\pi$ a factor representation of $G$.
\par

Then there exists an irreducible representation $\pi_{\rm irr}$ of $G$
which is weakly equivalent to $\pi$.
The map $\QD(G) \, \twoheadrightarrow \, \Pri(G)$ defined by
$$
(\text{quasi-equivalence class of $\pi$}) \mapsto (\text{weak equivalence class of $\pi_{\rm irr}$)}
$$
extends the natural map $\kappa^{\rm d}_{\rm prim} \,\colon \widehat G \twoheadrightarrow \Pri(G)$
of Definition \ref{Def-PrimitiveDual}.
\end{cor}

\begin{proof}
By Corollary~\ref{FacSeparable} and Corollary~\ref{QE-WEQ},
we can assume that the Hilbert space of $\pi$ is separable.
Then $\pi$ admits a direct integral decomposition 
$$
\pi \, \simeq \, \int^\oplus_X \pi_x d\mu(x)
$$
into irreducible representations (see Theorem~\ref{thmDirectIntIrreps}).
By Proposition~\ref{Pro-IntDecFactorRep}, almost every $\pi_x$ can be used for $\pi_{\rm irr}$.
\end{proof}

\begin{rem}
% 7.F.3
(1)
Corollary~\ref{Cor-FacRepPrimitiveG} follows also from Corollary~\ref{Cor-FacRepPrimitive}.

\vskip.2cm

(2)
Corollary~\ref{Cor-FacRepPrimitiveG} holds more generally
for a $\sigma$-compact locally compact group $G$ (see Proposition~\ref{FacKK}).
\end{rem}

In the case of non type I factor representations, we can prove the following 
stronger version of Corollary~\ref{Cor-FacRepPrimitiveG}.

\begin{cor}
% 7.F.4
\label{Cor-FacRepPrimitiveG-NonTypI}
Let $G$ be a second-countable locally compact group and $\pi$ a factor representation of $G$.
Assume that $\pi$ is not of type I.
\par

Then exist uncountably many pairwise non equivalent irreducible 
representations of $G$ which are weakly equivalent to $\pi$.
\end{cor}

\begin{proof}
As in the proof of Corollary~\ref{Cor-FacRepPrimitiveG},
we can assume that the Hilbert space of $\pi$ is separable.
Let
$$
\pi \, \simeq \, \int^\oplus_X \pi_x d\mu(x)
$$
be a direct integral decomposition into irreducible representations
over a standard Borel space $X$ defined by
a measurable field $x \mapsto \Hi_x$ of Hilbert spaces.
By Proposition~\ref{Pro-IntDecFactorRep}, we can assume that 
the $\pi_x$'s are all weakly equivalent to $\pi$.
\par

Let $\sigma$ be an irreducible representation of $G$.
Set 
$$
E_{\sigma} \, := \, \{x \in X \mid \pi_x \simeq \sigma\}
$$
(recall that $\simeq$ denotes equivalence of representations).
\par

We claim that $E_\sigma$ is a Borel subset of $X$.
Indeed, let $p \in \{1, 2, \dots, +\infty\}$ be the dimension of $\sigma$, let 
$$
X_p \, = \, \{x \in X \mid \dim \Hi_x = p \},
$$
and let $\Ki$ be a fixed Hilbert space of dimension $p$.
By Proposition~\ref{Pro-DecConstantFieldHilbertSpaces}, 
$X_p$ is a Borel subset of $X$ and the subrepresentation of $\pi$ defined on
$$
\int_{X_p}^\oplus \Hi_x d\mu(x) \approx L^2(X_p, \mu, \Ki)
$$
is equivalent to a representation in $L^2(X_p, \mu, \Ki)$
which is a direct integral $\int_{X_p} \pi_x^{(p)} d\mu(x)$ over $X_p$
of representations $\pi_x^{(p)}$ in $\Ki$.
Let $\mathrm{Irr}_p(G)$ be the set of irreducible representations of $G$
in $\Ki$, equipped with its Mackey--Borel structure
as in Subsection~\ref{SS:BorelStructureQD}.
The set 
$$
\Omega_\sigma \, = \, \{\rho \in \mathrm{Irr}_p(G) \mid \rho \simeq \sigma\}
$$
consists of the representations $\rho \in \mathrm{Irr}_p(G)$
for which $\Hom_G(\rho, \sigma)$ has dimension $1$.
By \cite[Lemma 3.7.3]{Dixm--C*}, $\Omega_\sigma$ is a Borel subset of $\mathrm{Irr}_p(G)$.
Since 
$$
X_p \to \mathrm{Irr}_p(G), \hskip.2cm x \mapsto \pi_x^{(p)}
$$
is measurable and since
$$
E_\sigma \, = \, \{ x \in X_p \mid \pi_x^{(p)} \in \Omega_\sigma\},
$$
$E_\sigma$ is a Borel subset of $X$.
\par

Then $E_\sigma$ is a Borel subset of $X$.
We claim that $\mu(E_{\sigma}) = 0$. 
Once proved, this will imply the claim, since it will follow that 
$\mu(\bigcup_{\sigma\in C} E_{\sigma}) = 0$ for every countable subset $C$ of $\widehat G$.
\par

Assume by contradiction $\mu(E_{\sigma})>0$ for some $\sigma\in \widehat G$,
with Hilbert space $\Hi_\sigma$.
Consider the representation $(\pi_{E_\sigma}, \Hi_{E_\sigma})$, where
$$
\pi_{E_\sigma} \, := \, \int^\oplus_{E_\sigma} \pi_x d\mu(x) 
\hskip.2cm \text{and} \hskip.2cm
\Hi_{E_\sigma} \, := \, \int^\oplus_{E_\sigma} \Hi_x d\mu(x).
$$
For every $x \in E_\sigma$,
there exists an isomorphism of Hilbert spaces $U_x \,\colon \Hi_\sigma \to \Hi_x$
intertwining $\sigma$ and $\pi_x$. 
It follows that there exists an isomorphism of Hilbert spaces 
$$
\Hi_\sigma \otimes L^2(E_\sigma, \mu) \to \Hi_{E_\sigma}
$$
intertwining a multiple of $\sigma$ and $\pi_{E_\sigma}$
(see Th\'eor\`eme 2 in \cite[Chap.~II, \S~2]{Dixm--vN}).
Therefore $\pi_{E_\sigma}$ is equivalent to a multiple of $\sigma$ and is therefore of type I.
However, $\pi_{E_\sigma}$ is a non trivial subrepresentation of 
the factor representation $\pi$ and is therefore quasi-equivalent to $\pi$.
Therefore $\pi$ is of type I and this is a contradiction
\end{proof}

%-----------------------------------------------------------------------
% End of chapter 7
%-----------------------------------------------------------------------
\chapter[C*-algebras and LC groups]
{Representations of C*-algebras and of LC groups}
% Chapter 8
\label{ChapterAlgLCgroup}

\emph{Let $G$ be a second-countable LC group 
and $\widehat G$ its dual, equipped with the Mackey--Borel structure,
as defined in Section~\ref{SectionQuasidual}.
On the one hand, as we will see in this chapter (Glimm theorem~\ref{ThGlimm}),
$\widehat G$ is not countably separated,
unless $G$ is of type I, as defined in Section~\ref{SectionTypeI}.
On the other hand, a result of Effros (Theorem~\ref{EffrosTheoremG})
shows that the primitive dual $\Pri(G)$ is always a standard Borel space.
This shows that $\Pri(G)$ is much better behaved than $\widehat G$,
and one may hope to classify $\Pri(G)$ for some specific non type I groups. 
This will be done in Chapter \ref{ChapterPrimExa} for most of our test examples.
}
\par

\emph{Given a representation $\pi$ of a locally compact group $G$, 
it is useful to ``extend" $\pi$ to group algebras associated to $G$.
Such a group algebra, which is universal in an appropriate sense, 
is the maximal C*-algebra $C^*_{\rm max}(G)$ of $G$,
which is defined in Section~\ref{C*algLCgroup}.
The representations of the group $G$ correspond precisely
to the non-degenerate representations of the C*-algebra $C^*_{\rm max}(G)$.
The spectrum (or dual space) $\widehat A$
and the primitive dual (or primitive ideal space) $\Pri(A)$
are defined for a general C*-algebra $A$,
and agree with $\widehat G$ and $\Pri(G)$
in case $A = C^*_{\rm max}(G)$ for a locally compact group $G$.
}
\par

\emph{Section \ref{SectionPrimC*} is devoted to 
the study of $\Pri(A)$ for a separable C*-algebra $A$,
where $\Pri(A)$ is equipped with the Borel structure induced by the Jacobson topology.
In particular, we give a complete proof of Effros' theorem:
$\Pri(A)$ is a standard Borel space.
}
\par

\emph{
In Section~\ref{SS:Functioriality}, we discuss in some length 
the functorial properties of maximal C*-algebras
and reduced C*-algebras of locally compact groups.
The assignment $G \rightsquigarrow C^*_{\rm max}(G)$
for a general locally compact group $G$ does not have good functorial properties,
with respect to arbitrary continuous homomorphisms.
However, as we will see, $\Gamma \rightsquigarrow C^*_{\rm max}(\Gamma)$
is a covariant functor from the category of discrete groups to the category of C*-algebras.
In contrast, this fails for the assignment
$\Gamma \rightsquigarrow C^*_{\lambda}(\Gamma)$,
where $C^*_{\lambda}(\Gamma)$ denotes the reduced C*-algebra of the group $\Gamma$.
}
\par

\emph{
Several results about a second-countable LC group $G$
use the fact that $C^*_{\rm max}(G)$ is a separable C*-algebra.
In Section \ref{C*kernelsprimitive}, we show how
some of these results can be extended to $\sigma$-compact LC groups,
using as one tool the Kakutani--Kodaira Theorem~\ref{Kakutani--Kodaira}.
In particular, the fact that a factor representation of $G$
is weakly equivalent to an irreducible one
(which was previously proved in case $G$ is second-countable,
see Corollary~\ref{Cor-FacRepPrimitiveG})
is valid for any $\sigma$-compact LC group.
}
\par

\emph{
In Section \ref{Section-CentralCharacter},
we come back to the notion of central character of a group representations,
that has already appeared in Section \ref{Section-IrrRepTwoStepNil},
and we show that it can be associated to a primitive ideal. 
}
\par

\emph{
Section \ref{SectionGlimm} contains a statement of Glimm theorem,
which characterizes type I for second-countable LC groups $G$.
In particular, this class of groups coincides
with the class of GCR groups from Section~\ref{SS:GCR groups}.
Actually, Glimm theorem characterizes separable C*-algebras of type I
(see Remark~\ref{Rem-Def-GCR-C*} for this notion).
We prefer to state a formulation in terms of second-countable groups,
and to show how it extends to the case of $\sigma$-compact LC groups.
}

\emph{
In Section \ref{vNalgLCgroup},
we discuss alternative definitions of the von Neumann algebra $\pi(G)''$
associated to a representation $\pi$ of a locally compact group $G$
(see Section~\ref{S:QE-Factor-VNAlgebras}).
In Section \ref{SectionVariants},
we discuss briefly one way to replace the maximal C*-algebra
of a locally compact group $G$ by other group algebras.
}

\section
{C*-algebras}
% Section 8.A
\label{SectionPrimC*}

Our goal in this section is a result
about primitive duals of second-countable LC groups,
Theorem \ref{EffrosTheoremG},
which follows as in \cite{Effr--63}
from a general result on separable C*-algebras,
Theorem \ref{EffrosTheoremA}.
For this, we need to relate the primitive ideal space of a C*-algebra 
to its space of pure states. 

\subsection
{Spectrum and primitive ideal space}
% subsection 8.A.a
\label{SS:SpecEtPrimDeA}

A \textbf{complex *-algebra} is a complex algebra $A$ endowed with an involution,
which is a map $A \to A, \hskip.1cm x \mapsto x^*$ which is
semilinear ($(x+y)^* = x^* + y^*$ and $(\lambda x)^* = \overline \lambda x^*$
for all $x,y \in A$ and $\lambda \in \C$),
involutive ($(x^*)^* = x$ for all $x \in A$),
and such that $(xy)^* = y^*x^3$ for all $x,y \in A$.
A \textbf{representation} of $A$ in a Hilbert space $\Hi$
is a morphism of $*$-algebras $\pi \,\colon A \to \Li (\Hi)$.
For a family $(\pi_\iota \,\colon A \to \Li (\Hi_\iota))_{\iota \in I}$,
there is natural notion of direct sum
$\bigoplus_{\iota \in I} \pi_\iota \,\colon A \to \Li (\bigoplus_{\iota \in I} \Hi_\iota)$.
Two representations $\pi_1, \pi_2$ of $A$ on Hilbert spaces $\Hi_1, \Hi_2$
are \textbf{equivalent} if there exists a non-zero \textbf{intertwiner} between them, i.e.,
a surjective isometry $U \,\colon \Hi_1 \to \Hi_2$
such that $U\pi_1(x) = \pi_2(x)U$ for all $x \in A$.
\index{Intertwiner! $3$@between representations of a C*-algebra}
\par

\index{Trivial representation}
\index{Representation! trivial (for a C*-algebra)}
\index{Non-degenerate representation of a *-algebra}
\index{Representation! non-degenerate}
A representation $\pi$ of a complex *-algebra $A$ in a Hilbert space $\Hi$
is \textbf{trivial} if $\pi(A) = \{0\}$.
It is \textbf{non-degenerate}
if the closed linear span $\overline{\pi(A)\Hi}$ of
$$
\{ \pi(a)\xi \mid a \in A, \hskip.1cm \xi \in \Hi\}
$$
is $\Hi$ itself.
Every representation $\pi \,\colon A \to \Li (\Hi)$ is a direct sum 
$$
\pi \, = \, \pi_0 \oplus \pi_1
$$
of a trivial representation $\pi_0$ and a non-degenerate representation $\pi_1$;
indeed,
$$
\Hi_0 \, = \, \{\xi \in \Hi \mid \pi(A)\xi = \{0\}\}
\hskip.5cm \text{and} \hskip.5cm
\Hi_1= \overline{\pi(A)\Hi}
$$
are $\pi(A)$-invariant closed subspaces of $\Hi$,
the corresponding subrepresentations $\pi_0$ and $\pi_1$
are respectively trivial and non-degenerate,
and $\Hi = \Hi_0\oplus \Hi_1$.
\par 

\index{Representation! irreducible}
A representation $\pi$ of $A$ in $\Hi$ is \textbf{irreducible}
if there are exactly two $\pi(A)$-invariant closed subspaces of $\Hi$,
namely $\{0\}$ and $\Hi$.
Note that an irreducible representation is either trivial and in $\Hi \approx \C$, 
or non-degenerate.
\index{$a5$@$\approx$ isomorphism! of Hilbert spaces} 

\vskip.2cm

\index{Representation! of a C*-algebra}
A \textbf{C*-algebra} is a complex involutive algebra $A$
endowed with a complete norm such that
$\Vert xy \Vert \le \Vert x \Vert \Vert y \Vert$ and $\Vert x \Vert^2 = \Vert x^*x \Vert$
for all $x,y \in A$.
A representation $\pi \,\colon A \to \Li (\Hi)$ of a C*-algebra $A$
is automatically norm decreasing:
$\Vert \pi(x) \Vert \le \Vert x \Vert$ for all $x \in A$ \cite[1.3.7]{Dixm--C*};
moreover, when $\pi$ is faithful, i.e., when $\pi$ is injective,
then $\Vert \pi(x) \Vert = \Vert x \Vert$ for all $x \in A$ \cite[1.8.1]{Dixm--C*}.
A representation $\pi \,\colon A \to \Li (\Hi)$ of a C*-algebra $A$
is irreducible, as defined above, if and only if
there are exactly two $\pi(A)$ invariant subspaces of $\Hi$ (closed or not),
namely $\{0\}$ and $\Hi$ \cite[2.8.4]{Dixm--C*}.

\begin{defn}
% 8.A.1
\label{Def-DualAPrimA}

The \textbf{spectrum} $\widehat A$ of a C*-algebra $A$
it the space of equivalence classes
of non-trivial irreducible representations of $A$.
We will define a topology on $\widehat A$
after the introduction of $\Pri(A)$.
\index{Spectrum of a C*-algebra}
\end{defn}

\begin{defn}
% 8.A.2
\label{Def-PrimA}
A two-sided ideal $I$ in $A$ is \textbf{primitive} 
if there exists a non-trivial irreducible representation $\pi \,\colon A \to \Li_(\Hi_\pi)$
such that $I = \ker (\pi)$.
The \textbf{primitive ideal space}
$\Pri(A)$ of $A$ 
is the set of primitive ideals of $A$.
\index{Primitive ideal}
\index{Primitive ideal space! $2$@of a C*-algebra}
\par

Observe that a primitive ideal is closed and proper, by definition.
%Recall that a closed two-sised ideal of $A$ is necessarily selfadjoint \cite[1.8.2]{Dixm--C*}.
The canonical map
\begin{equation}
\label{eqq/d/spA}
\tag{$\kappa$1$_A$}
\kappa^{\rm d}_{\rm prim} \, \colon \, \widehat A \twoheadrightarrow \Pri(A) ,
\hskip.5cm \pi \mapsto \ker (\pi) 
\end{equation}
is surjective, by definition
(compare with the map (\ref{candualprim}) in Section \ref{PrimIdealSpace}).
\end{defn}

\begin{rem}
% 8.A.3
\label{Rem-PrimitiveIdealAlg}

(1)
The notion of primitive ideal space is defined more generally 
for associative algebras.
Consider a complex associative algebra $\mathcal A$.
A two-sided ideal $I$ of $\mathcal A$ is 
\textbf{primitive}
if there exists a complex vector space $V$
and an irreducible representation 
$\pi_{\rm alg} \,\colon \mathcal A \to {\rm End}(V)$,
in the sense of associative algebras,
such that $I = \ker(\pi_{\rm alg})$.
The primitive ideal space of $\mathcal A$
is the space $\Pri(\mathcal A)$ of all primitive ideals of $\mathcal A$,
with the Jacobson topology.
Recall that the space of an irreducible representation is not $\{0\}$, by definition,
so that a primitive ideal is always proper (i.e., it cannot be $\mathcal A$ itself).
By definition, the closure of a subset $\Sigma$ of $\Pri(\mathcal A)$ is as above
the set of all $J \in \Pri(\mathcal A)$ such that
$\bigcap_{I \in \Sigma} I \subset J$.

\vskip.2cm

(2)
For two-sided ideals in C*-algebras,
primitive ideals as defined in (1) are the same as
primitive ideals as in Definition~\ref{Def-PrimA}.
More precisely:
\par

Let $A$ be a C*-algebra.
Let $\pi_{\rm alg}$ be an irreducible representation of $A$ in a complex vector space, 
in the sense of associative algebras.
Then there exists a C*-algebra representation $\pi_*$ of $A$ in some Hilbert space $\Hi_\pi$
that is equivalent to $\pi_{\rm alg}$. 
Moreover, if $\pi_*, \pi'_*$ are two irreducible C*-algebra representations of $A$,
they are algebraically equivalent 
(i.e., there exists an invertible linear operator $T$ 
from the space of $\pi_*$ to that of $\pi'_*$
such that $\pi'_*(x) = T\pi_*(x)T^{-1}$ for all $x \in A$)
if and only if they are equivalent in the sense of C*-algebras
(i.e., there exists a unitary operator $V$ 
from the space of $\pi_*$ to that of $\pi'_*$ 
such that $\pi'_*(x) = V\pi_*(x)V^*$ for all $x \in A$)
\cite[2.9.6]{Dixm--C*}.
\end{rem}

We now introduce a topology on the primitive ideal space of $A$.

\begin{defn}
% 8.A.4
\label{Def-JacobsonTopA}
The space $\Pri(A)$ is furnished with the
\textbf{Jacobson topology}, also called the \textbf{hull-kernel topology},
defined as follows:
the closure of a subset $\Sigma$ of $\Pri(A)$ is the set of all $J \in \Pri(A)$ such that
$$
\bigcap_{I \in \Sigma} I \, \subset \, J. 
$$
\index{Jacobson topology}
\end{defn}

Here is an important result about the separation properties of the topological space $\Pri(A)$.
We will simply say \textbf{maximal ideal} of $A$ for an ideal
which is maximal among the closed and proper two-sided ideals of $A$.

\begin{prop}
% 8.A.5
\label{Prop-SepPrimA}
Let $A$ be a C*-algebra.
\begin{enumerate}[label=(\arabic*)]
\item\label{iDEProp-SepPrimA}
$\Pri(A)$ is a T$_0$-space.
\item\label{iiDEProp-SepPrimA}
Let $I \in \Pri(A)$. Then $I$ is a closed point in $\Pri(A)$
if and only if $I$ is a maximal ideal of $A$.
In particular, $\Pri(A)$ is a T$_1$-space
if and only if every primitive ideal of $A$ is a maximal ideal.
\end{enumerate}
\end{prop}

\begin{proof}
\ref{iDEProp-SepPrimA}
Let $I_1, I_2$ be two distinct primitive ideals in $A$. We have therefore, say,
$I_1 \not \subset I_2$. By definition of the Jacobson topology, the closure
of $I_1$ is 
$$
X \, = \, \{ I \in \Pri(A) \mid I_1 \subset I\}.
$$
Therefore, $\Pri(A) \smallsetminus X$ is a neighbourhood of $I_2$
which does not contains $I_1$.

\vskip.2cm

\ref{iiDEProp-SepPrimA}
It is clear that $I$ is closed if and only if $I$ is maximal among the primitive ideals of $A$
(see the proof of \ref{iDEProp-SepPrimA}).
We claim that if $I$ is maximal among the primitive ideals of $A$,
then $I$ is maximal among the closed and proper two-sided ideals of $A$.
\par

Indeed, let $J$ be a closed and proper two-sided ideal of $A$ with $I \subset J$.
Let $\pi$ be an irreducible representation of the quotient C*-algebra $A/J$. 
Then $\pi$ lifts to an irreducible representation, again denoted by $\pi$, of $A$
such that $J \subset \ker (\pi)$. Since $\ker(\pi) \in \Pri(A)$,
it follows that $I = \ker(\pi)$ and hence $I = J$.
\end{proof}
 
\begin{rem}
% 8.A.6
\label{Rem-JacobsonTopA}

(1)
The fact that $\Pri(A)$ is a T$_0$-space should be compared
with the definition of $\Pri(G)$ from Reformulation~\ref{refdefprim}
together with Proposition~\ref{tradweqGC*}.

\vskip.2cm

(2)
The space $\Pri(A)$ is a Baire space
(see Proposition~\ref{P(A)Baire} below).

\vskip.2cm

(3)
When $A$ is separable, $\Pri(A)$ is second-countable
(see Corollary \ref{Cor-PrimSeparable} below).

\vskip.2cm

(4)
When $A$ has a unit, $\Pri(A)$ is quasi-compact \cite[3.1.8]{Dixm--C*}.
The converse does not hold:
the C*-algebra of compact operators on a separable infinite-dimensional Hilbert space 
has no unit, and its primitive ideal space consists of one point
(see Proposition~\ref{Pro-FactsAboutCompactOperators}).
\end{rem}

We can now define a topology and the induced Borel structure on the spectrum of $A$.

\begin{defn}
% 8.A.7
\label{Def-TopSpectre}
The \textbf{topology on the spectrum} $\widehat A$
of a C*-algebra $A$ is defined as follows: 
a subset of $\widehat A$ is open if and only if 
its inverse image under the map
$$
\kappa^{\rm d}_{\rm prim} \, \colon \, \widehat A \twoheadrightarrow \Pri(A)
$$
is open in $\Pri(A)$.
This topology induces the \textbf{Borel structure on the spectrum}~$\widehat A$.
% \index
\end{defn}

\begin{rem}
% 8.A.8
\label{Rem-Def-TopSpectre}

(1)
There are several equivalent ways to define the topology on $\widehat A$;
see \cite[\S~3]{Dixm--C*}.

\vskip.2cm

(2)
It is obvious that the map
$\kappa^{\rm d}_{\rm prim} \,\colon \widehat A \twoheadrightarrow \Pri(A)$ 
is continuous and open.

\vskip.2cm

(3)
Using Proposition~\ref{Prop-SepPrimA},
it is an easy exercise to show that $\widehat A$ is a T$_0$-space
if and only if the map $\kappa^{\rm d}_{\rm prim}$ is a homeomorphism.

\vskip.2cm

(4)
In general, the map $\kappa^{\rm d}_{\rm prim}$ is far from being a bijection;
see \ref{noninjectivity} below.
\end{rem}

The spectrum $\widehat A$ and the primitive ideal space $\Pri(A)$
introduced above are analogues for a C*-algebra $A$
of the dual $\widehat G$ and the primitive dual $\Pri(G)$
for a topological group $G$.
Moreover, there is also an analogue of the quasi-dual $\QD(G)$,
the quasi-spectrum $\mathrm{QS}(A)$ of $A$
(see Remark~\ref{Rem-Cor-FacRepPrimitive}).
%We will not say more on this here, and rather refer to Dixmier,
%in particular \cite[7.3.6]{Dixm--C*} 
%and \cite[Cor.\ 3 of Th.\ 2]{Dixm--60a}.

\subsection
{Pure states and primitive ideal space}
% subsection 8.A.b
\label{SS:Subsection-States}

Let $A$ be C*-algebra. 
Let $A^*_+$ denote the set of positive linear forms on $A$;
recall that these forms are continuous,
i.e., that $A^*_+$ is a subset of the Banach space dual $A^*$ of $A$ \cite[2.1.8]{Dixm--C*}.
Denote by $A^*_{+, \le 1}$, respectively $A^*_{+, 1}$, 
the subset of $A^*_+$ of forms of norm at most $1$, respectively equal $1$;
a \textbf{state} on $A$ is an element of $A^*_{+, 1}$.
\index{State on a C*-algebra}
\index{$e4$@$A^*_+, A^*_{+, \le 1}, A^*_{+, 1}$ 
positive linear forms on a C*-algebra $A$}
\par

When $A$ has a unit, it is easy to check that $\Vert \varphi \Vert = \varphi(1)$
for all $\varphi \in A^*_+$ \cite[2.1.4]{Dixm--C*}, 
so that $A^*_{+, 1} = \{ \varphi \in A^*_+ \mid \varphi(1) = 1 \}$.

\begin{constr}
% 8.A.9
\label{constructionGNS9}
The \textbf{Gel'fand--Naimark--Segal construction} associates
to every $\varphi \in A^*_+ \smallsetminus \{0\}$ 
a GNS triple $(\pi_\varphi, \Hi_\varphi, \xi_\varphi)$, 
consisting of a representation $(\pi_\varphi, \Hi_\varphi)$ of $A$
and a cyclic vector $\xi_\varphi \in \Hi_\varphi$
such that $\varphi(x) = \langle \pi_\varphi(x)\xi_\varphi\, \mid \, \xi_\varphi \rangle$, 
in the following way.
\index{Gel'fand--Naimark--Segal representation!}
\index{Gel'fand--Naimark--Segal representation! $2$@with $\varphi \in A^*_+$}
\par

The set $J = \{x \in A \mid \varphi(x^*x) = 0\}$
is a closed left ideal of $A$.
A scalar product is defined on $A / J$ by 
$$
\langle x+J \mid y+J \rangle \, = \, \varphi(y^*x) 
\hskip.5cm \text{for all} \hskip.2cm 
x,y \in A.
$$
Let $\Hi_\varphi$ be the Hilbert space completion of $A / J$.
Then a representation $\pi_\varphi$ of $A$ on $\Hi_\varphi$ 
is defined on $A / J$ by 
$$
\pi_\varphi(a) (x+J) \, = \, ax +J 
\hskip.5cm \text{for all} \hskip.2cm 
a,x \in A.
$$
Moreover, $x+J \mapsto \varphi(x)$ is well-defined on $A/J$ 
and extends to a continuous linear form on 
$\Hi_\varphi$, again denoted by $\varphi$.
Therefore there exists a unique vector $\xi_\varphi \in \Hi_\varphi$ 
such that
$$
\varphi(x) \, = \, \langle x+J \mid \xi_\varphi \rangle
\hskip.5cm \text{for all} \hskip.2cm
x \in A .
$$
(When $A$ is unital $\xi_\varphi$ is simply the image of $1$ in $\Hi_\varphi$.) 
We have
$$
\begin{aligned}g
\langle y+J \mid \pi_\varphi(x)\xi_\varphi \rangle
\, &= \, \langle \pi_\varphi(x^*) (y+J) \mid \xi_\varphi \rangle
\, = \, \langle x^* y+J \mid \xi_\varphi \rangle
\\
\, &= \, \varphi(x^*y) \, = \, \langle y + J \mid x + J \rangle
\end{aligned}
$$
for all $x,y \in A$, and hence $\pi_\varphi(x) \xi_\varphi = x+J$ for every $x \in A$.
This implies that $\xi_\varphi$ is cyclic 
and $\varphi(x) = \langle \pi_\varphi(x) )\xi_\varphi \mid \xi_\varphi \rangle$
for every $x \in A$. 
\par

The GNS triple $(\pi_\varphi, \Hi_\varphi, \xi_\varphi)$ 
associated as above to $\varphi$ is unique up to unitary equivalence:
if $(\pi, \Hi)$ is a representation of $A$ with cyclic vector $\xi$
such that $\Vert \xi \Vert = \Vert \xi_\varphi \Vert$
and $\varphi(x) = \langle \pi(x)\xi \mid \xi \rangle$ for all $x \in A$, 
then there exists a unitary Hilbert space isomorphism $S \,\colon \Hi \to \Hi_\varphi$ 
intertwining $\pi$ and $\pi_\varphi$ and such that $S(\xi) = \xi_\varphi$. 
\end{constr} 

\begin{exe}
% 8.A.10
\label{exGNSpourC*}
(1)
Consider a compact space $X$ and the C*-algebra $C(X)$.
The space of states $C(X)^*_{+,1}$ can be identified with
the space of probability Radon measures on~$X$;
this is one form of the Riesz representation theorem,
see Theorem \ref{Riesz}.
For $\mu$ such a measure,
there are identifications of the Hilbert space $\Hi_\mu$ with $L^2(X, \mu)$,
of the cyclic vector $\xi_\mu$ with the constant function of value $1$ on $X$,
and of the operator $\pi_\mu(f)$ with the multiplication by $f$ on $L^2(X, \mu)$
for all $f \in C(X)$.
The representation $\pi_\mu$ is irreducible
if and only if it is of dimension $1$,
if and only if $\mu$ is a Dirac measure.

\vskip.2cm
 
(2)
In anticipation of Section \ref{C*algLCgroup}, 
consider a locally compact group $G$ with left Haar measure $\mu_G$,
its group algebra $L^1(G, \mu_G)$,
and its maximal C*-algebra $C^*_{\rm max}(G)$.
Every normalized function of positive type $\varphi \in P_1(G)$
defines a normalized positive linear form $\omega^1_\varphi$ on $L^1(G, \mu_G)$
by $\omega^1_\varphi(f) = \int_G \varphi(g) f(g) d\mu_G(g)$,
which extends canonically to a state
$\omega = \omega^*_\varphi \in C^*_{\rm max}(G)^*_{+,1}$ \cite[13.4.5, 2.7.5]{Dixm--C*}.
The GNS representation $\pi_{\varphi}$ of $G$
corresponds to the GNS representation $\pi_\omega$ of $C^*_{\rm max}(G)$.
\par

[Recall that a positive linear form on a C*-algebra
is automatically continuous \cite[2.1.8]{Dixm--C*};
also, every positive linear form on $L^1(G, \mu_G)$ is continuous
(this follows from a theorem of Varopoulos, see \cite[32.27]{HeRo--70}).]
% Also, every continuous positive linear form on $L^1(G)$
% has a unique positive linear extention on $C^*_{\rm max}(G)$,
% blablabla \cite[2.7.5]{Dixm--C*}.
\end{exe}

\index{Indecomposable! $2$@for $\varphi \in A^*_{+, \le 1}$}
\index{Pure state on a C*-algebra}
Let $A$ be a C*-algebra.
In the unit ball of the dual Banach space $A^*$,
the subsets $A^*_{+, \le 1}$ and $A^*_{+,1}$ are convex.
We denote the corresponding sets of indecomposable elements by
${\rm Extr} (A^*_{+, \le 1})$ and ${\rm P}(A)$;
the set ${\rm P}(A)$ is the set of \textbf{pure states} of $A$.
The set $A^*_{+, \le 1}$ is compact in the unit ball of $A^*$ for the weak$^*$-topology,
and we have ${\rm Extr} (A^*_{+, \le 1}) = \{0\} \sqcup {\rm P}(A)$
\cite[2.5.5]{Dixm--C*}.
Note that ${\rm Extr} (A^*_{+,1})$ is closed in ${\rm Extr} (A^*_{+, \le 1})$
when $A$ has a unit,
and that $0$ is in its closure when $A$ has no unit \cite[2.12.13]{Dixm--C*}.
\par

For $\varphi \in A^*_{+, \le 1} \smallsetminus \{0\}$, 
the associated GNS representation $(\pi_\varphi, \Hi_\varphi)$ is irreducible 
if and only if $\varphi \in {\rm P}(A)$ \cite[2.5.4]{Dixm--C*}.
We therefore have a map 
$$
{\rm P}(A) \twoheadrightarrow \widehat A,
\hskip.5cm
\varphi\mapsto \pi_\varphi,
$$
where $\widehat A$ is the spectrum of $A$ 
(as defined in \ref{SS:SpecEtPrimDeA}).
This map is surjective: the inverse image of $\pi \in \widehat A$ 
consists of the states of the form
$x \mapsto \langle \pi(x)\xi \mid \xi \rangle$, where $\xi$ is a unit vector
in the Hilbert space of $\pi$.
Composition of this map with the canonical surjection
$$
\kappa^{\rm d}_{\rm prim} \, \colon \, \widehat A \twoheadrightarrow \Pri(A),
\hskip.2cm \pi \mapsto \ker(\pi)
$$
of Definition \ref{Def-PrimA} yields a surjective map
$$
\Phi \, \colon \, {\rm P}(A) \twoheadrightarrow \Pri(A),
\hskip.2cm \varphi \mapsto \ker (\pi_\varphi).
$$

\begin{prop}
% 8.A.11
\label{Pro-EtatsPrim}
Let $A$ be a C*-algebra.
Let its pure states space ${\rm P}(A)$ be equipped with the weak$^*$-topology
and its primitive ideal space $\Pri(A)$ be equipped with the Jacobson topology.
\par

The maps $\Phi \,\colon {\rm P}(A) \twoheadrightarrow \Pri(A)$
and ${\rm P}(A) \twoheadrightarrow \widehat A$
defined above are continuous and open.
\end{prop}

\begin{proof}
The topology on $\Pri(A)$ may be described in terms of pure states as follows
(see \cite[3.4.10]{Dixm--C*}).
Let $S \subset \Pri(A)$ and $J \in \Pri(A)$; set 
$$
Q(S) \, = \, \Phi^{-1}(S) \subset {\rm P}(A).
$$
Then $J \in \overline S$ if and only if $\Phi^{-1}(J) \in \overline{Q(S)}$.
This shows that $\Phi$ is continuous.
\par

Let $U$ be an open subset of ${\rm P}(A)$; we claim that $V := \Phi(U)$ is open in $\Pri(A)$.
Indeed, assume, by contradiction, that $V$ is not open. Then,
some $J \in V$ is in the closure of $\Pri(A) \smallsetminus V$.
Let $\varphi \in U$ be such that $\Phi(\varphi) = J$.
Then $\varphi$ belongs to the closure of the set
$$
\Phi^{-1}(\Pri(A) \smallsetminus V) \, = \, {\rm P}(A) \smallsetminus \Phi^{-1}(V).
$$
Since ${\rm P}(A) \smallsetminus \Phi^{-1}(V) \subset {\rm P}(A) \smallsetminus U$
and since ${\rm P}(A) \smallsetminus U$ is closed, it follows that 
$\varphi \in {\rm P}(A) \smallsetminus U$ and this is a contradiction.
\par

The claim for ${\rm P}(A) \twoheadrightarrow \widehat A$
%(which is \cite[3.4.11]{Dixm--C*}) 
follows.
\end{proof}

The following proposition will play a crucial role in the sequel.

\begin{prop}
% 8.A.12
\label{Pro-EtatsStandard}
Let $A$ be a separable C*-algebra.
\par

Then ${\rm P}(A)$ is a Polish space.
\end{prop}

\begin{proof} 
Since $A$ is separable, the unit ball $A^*_{\le 1}$ of $A^*$ 
is a metrizable compact space and hence a Polish space.
Therefore, the same holds for the closed subspace $X := A^*_{+, \le 1}$
of $A^*_{\le 1}$.
\par

Let $d$ be a metric on $X$ defining its topology; 
set $X_0 = \{0\}$ and, for $n \ge 1$, let $X_n$ be the set of $\varphi \in X$ 
such that there exist $\varphi_1, \varphi_2 \in X$
with $\varphi = \frac{1}{2} (\varphi_1+ \varphi_2)$
and $d(\varphi_1, \varphi_2) \ge 1/n$.
Then every $X_n$ is closed in $X$ and 
$$
{\rm P}(A) \, = \,
X \smallsetminus \bigcup_{n \ge 0} X_n
\, = \, \bigcap_{n \ge 0} (X \smallsetminus X_n)
$$
is a G$_\delta$ set in the Polish space $X$. 
Therefore, ${\rm P}(A)$ is a Polish space (see Proposition \ref{subspacePolish}).
\end{proof}

Here is a first consequence of Proposition~\ref{Pro-EtatsStandard}.
\begin{cor}
% 8.A.13
\label{Cor-PrimSeparable}
Let $A$ be a separable C*-algebra.
\par

The primitive ideal space $\Pri(A)$ and the spectrum $\widehat{A}$
are second-countable topological spaces.
\end{cor}

\begin{proof}
By Proposition~\ref{Pro-EtatsStandard}, 
there exists a sequence $(U_n)_n$ of open sets 
generating the topology of ${\rm P}(A)$. 
Since $\Phi$ is open and continuous (Proposition~\ref{Pro-EtatsPrim}),
the sequence of open sets $(\Phi(U_n))_n$
generates the topology of $\Pri(A)$. 
So, $\Pri(A)$ is second-countable.

The claim for $\widehat A$ is proved similarly.
\end{proof}

Here is a second consequence of Proposition~\ref{Pro-EtatsStandard}.

\begin{prop}
% 8.A.14
\label{P(A)Baire}
Let $A$ be a C*-algebra.
\par

The primitive ideal space $\Pri(A)$ and the spectrum $\widehat{A}$ are Baire spaces.
\index{Baire space}
\index{Topological space! Baire}
\end{prop}

\begin{proof}
We provide a proof for the particular case of a separable C*-algebra $A$,
and refer to \cite[3.4.13]{Dixm--C*} for the general case.
\par

Let $(V_n)_n$ be a sequence of open dense subsets of $\Pri(A)$.
Since $\Phi$ is continuous and open, every $U_n := \Phi^{-1}(V_n)$ 
is open and dense in ${\rm P}(A)$. 
However, the space ${\rm P}(A)$, 
being a Polish space (Proposition \ref{Pro-EtatsStandard}),
is a Baire space.
Therefore $\bigcap_n U_n$ is dense in ${\rm P}(A)$; 
since $\Phi$ is continuous, it follows that 
$$
\bigcap_n V_n \, = \, \Phi \left( \bigcap_n U_n \right)
$$
is dense in $\Pri(A)$. So, $\Pri(A)$ is a Baire space.

The claim for $\widehat A$ is proved similarly.
\end{proof}

\subsection
{A primitive ideal space is a standard Borel space}
% subsection 8.A.c
\label{SS:Subsection-Prim(A)}

Let $A$ be a C*-algebra.
\par

\index{C*-semi-norm} 
A \textbf{C*-semi-norm} on $A$ is a semi-norm
$N \,\colon A \to \R_+$ such that 
\begin{enumerate}[label=(\arabic*)]
\item
$N(x) \le \Vert x \Vert$,
\item 
$N(xy) \le N(x)N(y)$,
\item 
$N(x^*x) = N(x)^2$,
\end{enumerate}
for all $x,y \in A$.
\par
Denote by $\Ideal(A)$ the set of closed two-sided ideals of $A$
and by $\mathcal N (A)$ the set of C*-semi-norms on $A$.
\index{$e5$@$\Ideal(A)$ set of closed two-sided ideals of the C*-algebra $A$}
\index{$e6$@$\mathcal N (A)$ C*-semi-norms on the C*-algebra $A$}

\begin{exe}
% 8.A.15
\label{Ex-C^*-seminorm}
(1)
The maximum of two C*-semi-norms is a C*-semi-norm
(but the minimum is not).

\vskip.2cm

(2)
Let $J \in \Ideal(A)$.
Then $A/J$, equipped with the quotient norm
$$
x+J \, \mapsto \, \Vert x+J \Vert := \inf_{y \in J} \Vert x+y \Vert ,
$$
is a C*-algebra and the map
$$
N_J \, \colon \, A \to \R_+ 
\hskip.5cm \text{defined by} \hskip.5cm
N_J(x) = \Vert x+J\Vert
\hskip.5cm \text{for all} \hskip.2cm
x \in A
$$
is a C*-semi-norm on $A$.
Observe that $\{x \in A\, \mid \, N_J(x) = 0\} = J$.
The next lemma shows that every C*-semi-norm on $A$ is of the form
$N_J$, for some closed two-sided ideal $J$.
\end{exe}

\begin{lem}
% 8.A.16
\label{Lem-C^*-seminorm}
We keep the notation above.
\begin{enumerate}[label=(\arabic*)]
\item\label{iDELem-C^*-seminorm}
The map $J \mapsto N_J$ is a one-to-one correspondence
from $\Ideal(A)$ onto $\mathcal N (A)$.
\item\label{iiDELem-C^*-seminorm}
For every $J_1,J_2 \in \Ideal(A)$, 
we have $N_{J_1 \cap J_2} = \max \{ N_{J_1}, N_{J_2} \}$.
\end{enumerate}
\end{lem}

\begin{proof}
\ref{iDELem-C^*-seminorm}
Let $N \in \mathcal N (A)$. Set $J := \{x \in A \mid N(x) = 0\}$.
Then $J \in \Ideal(A)$. 
Moreover, the map 
$\overline N \,\colon x+J \mapsto N(x)$ is a norm on the quotient $A/J$.
The completion $\overline A$ of $A/J$ with respect to $\overline N$ 
is a C*-algebra. 
The inclusion map from $A/J$, 
equipped with the quotient norm of $A$, to $\overline A$,
is an injective homomorphism between C*-algebras and is therefore an isometry.
Therefore $N = N_J$ and this shows that $J \mapsto N_J$ is surjective. 
\par

Let $J_1,J_2 \in \Ideal(A)$ be such that $N_{J_1} = N_{J_2}$. Then
$$
J_1 \, = \, \{x \in A\, \mid \, N_{J_1}(x) = 0\} \, = \, \{x \in A\, \mid \, N_{J_2}(x) = 0\} \, = \,J_2
$$
and this shows that $J \mapsto N_J$ is injective. 

\vskip.2cm

\ref{iiDELem-C^*-seminorm}
The $*$-algebra $A/J_1 \oplus A/J_2$, equipped with the norm 
$$
\Vert (x+J_1, y+ J_2) \Vert \, = \, \max \{ \Vert x+J_1 \Vert, \Vert y + J_2 \Vert \} 
$$
is a C*-algebra. The map 
$$
A/J_1 \cap J_2 \to A/J_1 \oplus A/J_2, 
\hskip.5cm 
x+ J_1 \cap J_2 \mapsto (x+J_1, x+ J_2)
$$
is an injective homomorphism of C*-algebras, hence it is an isometry. 
Therefore, we have 
$$
N_{J_1 \cap J_2}(x)
\, = \,
\Vert x+J_1 \cap J_2 \Vert
\, = \, \max \{ \Vert x + J_1 \Vert, \Vert x + J_2 \Vert \}
\, = \, \max \{ N_{J_1}(x), N_{J_2}(x) \},
$$
for all $x \in A$. 
\end{proof}

\index{Prime ideal}
Recall that a two-sided ideal $J$ in $A$ is \textbf{prime} 
if $J \ne A$ and if, for $J_1, J_2 \in \Ideal(A)$ such that $J_1 J_2 \subset J$,
at least one of the inclusions $J_1 \subset J$, $J_2 \subset J$ holds.
Here, we are only interested in closed prime ideals.
For examples of closed prime ideals, see Proposition \ref{Prop-PrimPrimitive}.
\par

Observe that, for $J_1,J_2 \in \Ideal(A)$, we have $J_1J_2 = J_1 \cap J_2$.
Indeed, it is obvious that $J_1J_2 \subset J_1 \cap J_2$. 
To show the reverse inclusion,
it suffices to show that $(J_1 \cap J_2)_+ \subset J_1J_2$, 
since the C*-algebra $J_1 \cap J_2$ is the linear span of 
$(J_1 \cap J_2)_+$. 
Let $x \in (J_1 \cap J_2)_+$. Then $x^{1/2} \in J_1 \cap J_2$
and hence $x = x^{1/2}x^{1/2} \in J_1 \cap J_2$.

\vskip.2cm

\index{Extremal semi-norm}
A C*-semi-norm $N \in \mathcal N (A)$ is \textbf{extremal} if, 
whenever $N_1, N_2 \in \mathcal N (A)$ are such that $N = \max \{ N_1, N_2 \}$, 
then $N_1 = N$ or $N_2 = N$.
We denote by $\ENorm(A)$ the set of extremal C*-semi-norms on $A$ 
that are different from $0$.
\index{$d3$@${\rm Extr} (\cdot)$ indecomposable elements}
\index{$e7$@$\ENorm(A)$ extremal C*-semi-norms on the C*-algebra $A$}
\par

The next lemma is Theorem 2.1 in \cite{Effr--63}.

\begin{lem}
% 8.A.17
\label{Lem-PrimeExtremal}
Let $J \in \Ideal(A)$.
The following properties are equivalent:
\begin{enumerate}[label=(\roman*)]
\item\label{iDELem-PrimeExtremal}
$J$ is a closed prime ideal;
\item\label{iiDELem-PrimeExtremal}
whenever $J = J_1 \cap J_2$ for $J_1,J_2 \in \Ideal(A)$ holds, 
then $J_1 = J$ or $J_2 = J$;
\item\label{iiiDELem-PrimeExtremal}
$N_J \in \ENorm(A)$.
\end{enumerate}
\end{lem}

\begin{proof}
Assume that $J$ is prime. If $J = J_1 \cap J_2$, then $J = J_1J_2$ 
and hence $J_1 = J$ or $J_2 = J$.
This shows that \ref{iDELem-PrimeExtremal} implies \ref{iiDELem-PrimeExtremal}.

\vskip.2cm

Let $J_1,J_2 \in \Ideal(A)$ be such that $J_1 \cap J_2 \subset J$.
Set $I_1 = \overline{J+J_1}$ and $I_2 = \overline{J+J_2}$.
Then $I_1, I_2 \in \Ideal(A)$ and $I_1 I_2 \subset J$; 
hence, we have
$$
I_1 \cap I_2 \, = \, I_1I_2 \subset J \subset I_1 \cap I_2,
$$
that is, $J = I_1 \cap I_2$. Therefore, $I_1 = J$ or $I_2 = J$, that is,
$J_1 \subset J$ or $J_2 \subset J$. So, $J$ is a prime ideal.
This shows that \ref{iiDELem-PrimeExtremal}implies \ref{iDELem-PrimeExtremal}.

\vskip.2cm

The C*-semi-norm $N_J$ is extremal if and only if, 
whenever $N_J = \max \{ N_{J_1}, N_{J_2} \}$ for $J_1, J_2 \in \Ideal(A)$ holds, 
then $J_1 = J$ or $J_2 = J$.
Since $\max \{ N_{J_1}, N_{J_2} \} = N_{J_1 \cap J_2}$
(Lemma~\ref{Lem-C^*-seminorm}), this shows that 
\ref{iiiDELem-PrimeExtremal} and \ref{iiDELem-PrimeExtremal} are equivalent.
\end{proof}

Observe that, for a subset $S$ of $\Pri(A)$, 
we have $\bigcap_{J \in S} J = \{0\}$ if and only if $S$ is dense in $\Pri(A)$.
In particular, if $S$ is a closed subset of $\Pri(A)$ 
such that $\bigcap_{J \in S} J = \{0\}$,
then $S = \Pri(A)$.
% indeed, see \cite[2.9.7]{Dixm--C*}). 
% $I$ is a closed two-sided ideal of $A$, then $I = \bigcap_J J$, 
% where the intersection is taken over all primitive ideals 
% $J$ which contain $I$ \cite[2.9.7]{Dixm--C*}.

\begin{prop}
% 8.A.18
\label{Prop-PrimPrimitive}
Let $A$ be a C*-algebra and $J$ a closed two-sided ideal in $A$.
\begin{enumerate}[label=(\arabic*)]
\item\label{iDEProp-PrimPrimitive}
If $J \in \Pri(A)$, then $J$ is a closed prime ideal.
\item\label{iiDEProp-PrimPrimitive}
Assume that $A$ is separable. If $J$ is a closed prime ideal, then $J \in \Pri(A)$. 
\end{enumerate}
\end{prop}

\begin{proof}
\ref{iDEProp-PrimPrimitive}
% voir aussi \cite[VII.3.1, page 191]{Davi--96}.
Let $(\pi, \Hi)$ be an irreducible representation of $A$ 
such that $J = \ker(\pi)$. 
Let $J_1,J_2 \in \Ideal(A)$ be such that $J_1J_2 \subset J$. 
Assume that $J_2$ is not contained in $J$, 
that is, $\pi(J_2) \ne \{0\}$. 
Then $\overline{\pi(J_2)\Hi}$ is a closed $\pi(A)$-invariant subspace 
and therefore $\overline{\pi(J_2)\Hi} = \Hi$,
by irreducibility of $\pi$.
\par

Since 
$$
\pi(J_1)\overline{\pi(J_2)\Hi} \, \subset \, 
\overline{\pi(J_1)\pi(J_2)\Hi} \, = \,
\overline{\pi(J_1J_2)\Hi} \, \subset \,
\overline{\pi(J)\Hi} = \{0\},
$$
it follows that $\pi(J_1) = \{0\}$, that is, $J_1 \subset J$.

\vskip.2cm

\ref{iiDEProp-PrimPrimitive}
It suffices to show that if $\{0\}$ is a prime ideal of $A$, then $\{0\}$ is a primitive. 
Indeed, assume this is proved and let $J$ be a closed prime ideal of $A$. 
Then $\{0\}$ is a prime ideal of the separable C*-algebra $A/J$ 
and hence $\{0\} \in \Pri(A/J)$; 
this means that $A/J$ has a faithful irreducible representation $\pi$. 
Then $\pi \circ p \in \widehat A$, where $p \,\colon A \to A/J$ is the canonical projection, 
and hence $J = \ker(\pi \circ p) \in \Pri(A)$.
\par
 
Assume that $\{0\}$ is a prime ideal of $A$. 
Since $\Pri(A)$ is second-countable (Corollary~\ref{Cor-PrimSeparable}),
there exists a sequence $(V_n)_n$ of non-empty open subsets 
generating the topology of $\Pri(A)$. 
\par

For $n$, set 
$$
I_n \, = \, \bigcap_{J \in V_n} J 
\hskip.5cm \text{ and } \hskip.5cm
J_n \, = \, \bigcap_{J \in \Pri(A) \smallsetminus V_n} J.
$$
Then $I_n$ and $J_n$ are closed two-sided ideals of $A$ and $I_n \cap J_n = \{0\}$.
Moreover, $J_n \ne \{0\}$, since $\Pri(A) \smallsetminus V_n$
is a closed proper subset of $\Pri(A)$.
The ideal $\{0\}$ being prime, it follows that $I_n = 0$.
This means that $V_n$ is a dense subset of $\Pri(A)$. 
Since $\Pri(A)$ is a Baire space
(Corollary~\ref{Cor-PrimSeparable}), 
$\bigcap_n V_n$ is therefore dense in $\Pri(A)$.
\par 
 
Let $J \in \bigcap_n V_n$. Then $\{J\}$ is dense in $\Pri(A)$, 
since the sequence $(V_n)_n$ generates the topology of $\Pri(A)$. 
Therefore, $J\subset I$ for every $I \in \Pri(A)$ and hence $J = 0$.
So, $\{0\}$ is a primitive ideal. 
\end{proof}

\begin{rem}
% 8.A.19
% Claim \ref{iiDEProp-PrimPrimitive} is Cor.~1 of Th.~2 in \cite{Dixm--60a}.
There are known C*-algebras (non separable ones)
in which the closed ideal $\{0\}$ is prime and not primitive \cite{Weav--03}.
\end{rem}

As a consequence of Proposition~\ref{Prop-PrimPrimitive},
we now show that kernels of factor representations of a separable C*-algebra
are primitive ideals.
A first proof of this result was already given in Corollary~\ref{Cor-FacRepPrimitiveG},
when $A$ is the maximal C*-algebra
of a second-countable LC group (see Section~\ref{C*algLCgroup}).
For another proof, see \cite[Corollaire 3]{Dixm--60a}.
\par

As for group representations (see Proposition~\ref{Pro-FactorRepVN}), 
a representation $\pi$ of the C*-algebra $A$
is said to be \textbf{factorial} (or a \textbf{factor representation})
if the von Neumann algebra $\pi(A)''$ generated by $\pi(A)$ is a factor.

\begin{cor}
% 8.A.20
\label{Cor-FacRepPrimitive}
Let $A$ be separable C*-algebra and $\pi$ a factor representation of $A$.
\par

Then $\ker \pi \in \Pri(A)$. 
\end{cor}

\begin{proof}
Let $\Hi$ be the Hilbert space of $\pi$ and set $J := \ker \pi$.
By Proposition~\ref{Prop-PrimPrimitive}~\ref{iiDEProp-PrimPrimitive}
it suffices to show that $J$ is a prime ideal.
\par

Let $J_1,J_2 \in \Ideal(A)$ be such that $J_1J_2 \subset J$. 
Assume that $J_2$ is not contained in $J$, that is, $\pi(J_2) \ne \{0\}$. 
\par

Then $\overline{\pi(J_2)\Hi}$ is a non-zero closed subspace of $\Hi$
which is invariant under both $\pi(A)$ and its commutant $\pi(A)'$.
It follows that the orthogonal projection 
$$
P \, \colon \, \Hi \to \overline{\pi(J_2)\Hi}
$$
is a non-zero projection in the centre of $\pi(A)''$.
Since $\pi$ is factorial, we have $P = \mathrm{Id}_{\Hi}$,
that is $\overline{\pi(J_2)\Hi} = \Hi$.
Therefore 
$$
\pi(J_1)\Hi \, = \,
\pi(J_1)\overline{\pi(J_2)\Hi} \, = \, 
\overline{\pi(J_1)\pi(J_2)\Hi} \, = \,
\overline{\pi(J_1J_2)\Hi} \subset \overline{\pi(J)\Hi} = \{0\},
$$
and so $\pi(J_1) = \{0\}$, that is, $J_1 \subset J$.
\end{proof}

\begin{rem}
% 8.A.21
\label{Rem-Cor-FacRepPrimitive}
In analogy with the group case (see Proposition~\ref{Pro-SubordinateVN}),
two factor representations $\pi_1$ and $\pi_2$ of $A$ are said to be \textbf{quasi-equivalent}
if there exists an isomorphism $\Phi \,\colon \pi_1(A)'' \to \pi_2(A)''$
such that $\Phi(\pi_1(a)) = \pi_2(a)$ for every $a\in A$.
The \textbf{quasi-spectrum} $\mathrm{QS}(A)$ of $A$ is the space
of quasi-equivalence classes of factor representations of $A$.
Observe that two quasi-equivalent representations of $A$ have the same kernels.
\par

When $A$ is separable, Corollary~\ref{Cor-FacRepPrimitive} shows that
the map $\pi\mapsto \ker \pi$ induces a map $\mathrm{QS}(A) \twoheadrightarrow \Pri(A)$.
Observe that this map is an extension to $\mathrm{QS}(A)$
of the canonical map $\widehat A \twoheadrightarrow \Pri(A)$.
\end{rem}

We equip the set $\mathcal N (A)$ of C*-semi-norms on $A$
with the topology of pointwise convergence on $A$.
Then $\mathcal N (A)$ is a compact space, 
since it consists of bounded functions on $A$. 
When $A$ is separable, $\Norm(A)$ is a compact metrizable space 
and hence a Polish space.

\begin{lem}
% 8.A.22
\label{Lem-CSeminormCompact}
Assume that the C*-algebra $A$ is separable.
\par

The set $\ENorm(A)$ of extremal C*-semi-norms on $A$ 
is a G$_\delta$ set in the Polish space $\Norm(A)$.
In particular, $\ENorm(A)$ is a standard Borel space.
\end{lem}
 
\begin{proof}
Let $(x_n)_{n \in \N}$ be a dense sequence in $A$. 
For $i, j \in \N$ with $i \ne j$ and $k \ge 1$, set 
$$
U_{i, j}^{(k)} \, := \, 
\left\{ (N, N') \in \Norm(A) \times \Norm(A) 
\, \Big\vert \, 
N(x_i) +\frac{1}{k} \le N'(x_i) \text{ and } N'(x_j) +\frac{1}{k} \le N(x_j) \right\}.
$$
Then $U_{i, j}^{(k)}$ is a closed and hence compact subset 
of $\Norm(A) \times \Norm(A)$. 
\par
 
The map 
$$
f \, \colon \, \Norm(A) \times \Norm(A) \to \Norm(A), 
\hskip.5cm
(N, N')\mapsto \max \{ N, N' \}
$$
is continuous and we have 
$$
f \left( \bigcup_{i, j \in \N, i \ne j, k \ge 1} U_{i, j}^{(k)} \right ) \cup \{0\}
\, = \, \Norm(A) \smallsetminus \ENorm(A).
$$
Since $\{0\}$ is closed and each $f(U_{i, j}^{(k)})$ is compact,
it follows that $\Norm(A) \smallsetminus \ENorm(A)$
is a countable union of closed subsets of $\Norm(A)$;
therefore, $\ENorm(A)$ is a G$_\delta$ set in $\Norm(A)$. 
\end{proof}

The next result is Theorem 2.1 in \cite{Effr--63}.

\begin{theorem}[\textbf{Effros}]
% 8.A.23
\label{EffrosTheoremA}
Let $A$ be separable C*-algebra.
\par

The primitive ideal space
$\Pri(A)$ of $A$ is a standard Borel space.
\end{theorem}

\begin{proof} 
Consider the injective map 
$$
F \, \colon \, \Pri(A) \hookrightarrow \Norm(A), 
\hskip.5cm
J \mapsto N_J .
$$
By Lemma~\ref{Lem-PrimeExtremal}
and Proposition~\ref{Prop-PrimPrimitive}, we have $F(\Pri(A)) = \ENorm(A)$.
In view of Lemma~\ref{Lem-CSeminormCompact}, it suffices to show
that $F \,\colon \Pri(A) \to \ENorm(A)$ is an isomorphism of Borel spaces.
\par

Let us show that $F$ is a Borel map. By the definition
of the topology of $\Norm(A)$, it suffices to show that $J \mapsto N_J(x)$
is a Borel function for every $x \in A$. 
We are going to prove that $J \mapsto N_J(x)$ is even lower semi-continuous. 
\par

To show this, let $\alpha \in \R$. 
Let $J_0 \in \Pri(A)$ be in the closure of the set 
$$
X_{\alpha} \, := \, \{J \in \Pri(A)\, \mid \, N_J(x) \le \alpha\}.
$$
Set $I_{\alpha} = \bigcap_{J \in X_{\alpha}} J$; 
then $I_{\alpha} \subset J_0$ and hence $A/J_0$ is a quotient
of the C*-algebra $A/I_{\alpha}$. Therefore, we have
$$
N_{J_0}(x)
\, = \,
\Vert x + J_0 \Vert \le \Vert x + I_{\alpha} \Vert
\, = \, 
\sup_{J \in X_{\alpha}} \Vert x + J \Vert \, = \, \sup_{J \in X_{\alpha}} N_J(x)
$$
and hence $N_{J_0}(x) \le \alpha$, that is, $J_0 \in X_\alpha$. 
This shows that $X_\alpha$ is closed
for every $\alpha \in \R$; therefore $J \mapsto N_J(x)$ 
is a lower semi-continuous for every $x \in A$.
\par

Next, we show that $F^{-1}$ is even a continuous map. 
Let $S$ be a closed subset of $\Pri(A)$.
Set $I = \bigcap_{J \in S} J$. 
Since $S$ is closed, $I$ is the intersection of the primitive ideals that contain $I$. 
Therefore
$$
F(S) \, = \, \{N_J \mid J \in S \} \, = \, \{N \in \ENorm(A)\, \mid \, N(I) = 0\}.
$$
The set
$$
\{N \in \ENorm(A) \mid N(I) = 0\}
\, = \,
\bigcap_{x \in I} \{N \in \ENorm(A) \mid N(x) = 0 \}
$$
 is closed in $\ENorm(A)$. Therefore $F^{-1}$ is continuous.
\end{proof}

Anticipating on the next section,
and applying Theorem~\ref{EffrosTheoremA} to LC groups, 
we obtain the following result.

\begin{theorem}
% 8.A.24
\label{EffrosTheoremG}
Let $G$ be a second-countable locally compact group.
\par

The primitive dual $\Pri(G)$ is a standard Borel space.
\end{theorem}

\section[C*-algebras and representations of LC groups]
{C*-algebras and representations of LC groups}
% Section 8.B
\label{C*algLCgroup}

\index{Algebras! $3$@convolution algebras, $C^c(G)$, $C^0(G)$, $L^1(G)$}
\index{$j4$@$L^1(G, \mu_G)$ or $L^1(G)$ group algebra of the locally compact group $G$}
Let $G$ be a locally compact group, with left-invariant Haar measure $\mu_G$.
Let $L^1(G, \mu_G)$, or for short $L^1(G)$, denote the Banach $*$-algebra
of integrable functions $G \to \C$ (up to equality $\mu_G$-almost everywhere),
with the convolution product and the usual involution, see Appendix \ref{AppLCG}.
\par

To a representation $\pi \,\colon G \to \U(\Hi_\pi)$,
there is associated an algebra homomorphism $L^1(G) \to \Li (\Hi_\pi)$,
denoted by $\pi$ again, defined by
$$
\pi(f) \, = \, \int_G f(g) \pi(g) d\mu_G(g) 
\hskip.5cm \text{for all} \hskip.2cm 
f \in L^1(G).
$$
Moreover, since $\pi$ is a unitary representation of $G$, we have 
$$
\begin{aligned}
\pi(f^*) & \, = \, \int_G f^*(g) \pi(g) d\mu_G(g)
\, = \, \int_G \Delta(g^{-1}) \overline{f(g^{-1})} \pi(g) d\mu_G(g)
\\
& \, = \, \int_G \overline{f(g)} \pi(g^{-1}) d\mu_G(g)
\, = \, \left(\int_G f(g) \pi(g) d\mu_G(g)\right)^*
\, = \, (\pi(f))^*
\end{aligned}
$$
for all $f \in L^1(G)$, so that
$$
\pi \, \colon \, L^1(G) \to \Li (\Hi_\pi)
$$
is a $*$-representation of the $*$-algebra $L^1(G)$. 
This has an extension to a $*$-representation $M^b(G) \to \Li (\Hi)$ of the measure algebra,
defined by
$$
\pi(\mu) \, = \, \int_G \pi(g) d\mu(g)
\hskip.5cm \text{for all} \hskip.2cm 
\mu \in M^b(G).
$$
\par

We denote by $C^*_\pi(G)$ the C*-subalgebra of $\Li (\Hi_\pi)$
generated by the set of operators $\pi(L^1(G))$.
\par

There are two particularly important cases.

\begin{defn}
% 8.B.1
\label{Def-ReducedMaxC*Alg}
Let $G$ be a LC group.
\par

Let $\lambda$ be left regular representation of $G$ on $L^2(G, \mu_G)$.
The algebra $C^*_\lambda(G)$ is called the \textbf{reduced C*-algebra} of~$G$.
\index{Left regular representation}
\index{Reduced C*-algebra of a LC group}
\index{$j5$@$C^*_\lambda(G)$ reduced C*-algebra
of the locally compact group $G$}
\par

Let $\pi_{\rm max}$ be the direct sum
of all equivalence classes of cyclic representations of~$G$.
The algebra $C^*_{\pi_{\rm max}}(G)$
is called the \textbf{maximal C*-algebra} of $G$ and is denoted by $C^*_{\rm max}(G)$.
\index{Maximal C*-algebra of a LC group} 
\index{$j6$@$C^*_{\rm max}(G)$ maximal C*-algebra of the locally compact group $G$}
\end{defn}

\begin{rem}
% 8.B.2
\label{Rem-DefC*Max}
The maximal C*-algebra $C^*_{\rm max}(G)$ is the completion of $L^1(G)$ with respect to the norm
$$
\Vert f\Vert_{\rm max} \, := \, \sup_\pi \Vert \pi(f) \Vert_{\Li (\Hi_\pi)},
$$
where the supremum is taken
over all equivalence classes of cyclic representations $\pi$ of $G$.
The condition of cyclicity ensures that 
these classes of representations constitute a set.
For the definition of $C^*_{\rm max}(G)$,
one could just as well consider the set of all equivalence classes
of \emph{irreducible} representations of $G$; see \cite[\S~2.7]{Dixm--C*}.
\end{rem}

Let $\pi$ be a representation of $G$. 
The associated $*$-representation $L^1(G) \to \Li (\Hi_\pi)$ extends to 
a non-degenerate C*-representation $C^*_{\rm max}(G) \to \Li (\Hi_\pi)$,
often denoted by~$\pi$ again.
Conversely, every non-degenerate C*-representation of $C^*_{\rm max}(G)$
is associated in this way to exactly one representation of $G$;
see \cite[13.3.4, 13.9.3]{Dixm--C*} or \cite[Sections 3.2 and 7.2]{Foll--16}.
This shows:

\begin{prop}
% 8.B.3
\label{bijrepGrepC*}
For a locally compact group $G$,
the construction above provides a bijection
between group representations of $G$
and non-degenerate C*-representations of $C^*_{\rm max}(G)$.
\par

This bijection respects equivalence and irreducibility.
In particular, it induces a homeomorphism
$$
\widehat G \, \overset{\approx}{\longrightarrow} \, \widehat{ C^*_{\rm max}(G) }
$$
from the dual $\widehat G$ of the group $G$,
with the Fell topology as defined in Chapter \ref{ChapterUnitaryDual},
onto the spectrum of the C*-algebra $C^*_{\rm max}(G)$,
with the topology of Definition \ref{Def-TopSpectre}.
\end{prop}

Note that Propositions \ref{bijrepGrepC*} and Proposition~\ref{P(A)Baire}
imply that $\widehat G$ is a Baire space.
\par

For $\pi \,\colon G \to \U(\Hi_\pi)$ as above, 
the image of $\pi \,\colon C^*_{\rm max}(G) \to \Li (\Hi_\pi)$ is dense in $C^*_\pi(G)$ 
(because the image of $L^1(G)$ is already dense in $C^*_\pi(G)$, by definition),
and therefore is $C^*_\pi(G)$ itself
(because images of morphisms of C*-algebras are closed \cite[1.8.3]{Dixm--C*}),
so that the representation $\pi$ of $G$ gives rise to
a surjective morphism of C*-algebras 
$$
C^*_{\rm max}(G) \twoheadrightarrow C^*_\pi(G).
$$
The kernel of this morphism, denoted by $\textnormal{C*ker}(\pi)$,
is the \textbf{C*-kernel} of $\pi$.
C*-kernels provide for representations of LC groups 
an equivalent definition of weak containment and weak equivalence
\cite[3.4.4 \& 18.1.4]{Dixm--C*}:
\index{Representation! C*-kernel}
\index{C*-kernel}
\index{Representation! weakly contained $\preceq$}
\index{Weakly! contained (representation)}
\index{Representation! weakly equivalent $\sim$}
\index{Weakly! equivalent (representation)}

\begin{prop}
% 8.B.4
\label{tradweqGC*}
Let $G$ be a locally compact group. 
\par

Let $\pi_1, \pi_2$ be two representations of $G$.
The representation $\pi_1$ is weakly contained in $\pi_2$ if and only if
$\textnormal{C*ker}(\pi_2) \subset \textnormal{C*ker}(\pi_1)$,
and therefore $\pi_1$ and $\pi_2$ are weakly equivalent if and only if
$\textnormal{C*ker}(\pi_1) = \textnormal{C*ker}(\pi_2)$.
\par

The bijection of Proposition \ref{bijrepGrepC*}
induces a homeomorphism
$$
\Pri(G) \, \overset{\approx}{\longrightarrow} \, \Pri( C^*_{\rm max}(G) )
$$
from the primitive dual $\Pri(G)$ of the group $G$,
with its topology as defined in Section \ref{PrimIdealSpace},
onto the primitive ideal space of the C*-algebra $C^*_{\rm max}(G)$.
\end{prop}

\begin{rem}
% 8.B.5
\label{Rem-tradweqGC*}
Let $\pi_1, \pi_2$ be two representations of a LC group $G$
and let $B$ be a dense subspace of $C^*_{\rm max}(G)$,
for example $B=C^c(G)$ or $B=L^1(G)$.
The following properties are equivalent: 
\begin{enumerate}[label=(\roman*)]
\item\label{iDERem-tradweqGC*}
$\pi_1$ is weakly contained in $\pi_2$; 
\item\label{iiDERem-tradweqGC*}
$\Vert \pi_1(x)\Vert \le \Vert \pi_2(x)\Vert$ for every $x \in B$;
\item\label{iiiDERem-tradweqGC*}
$\Vert \pi_1(\mu) \Vert \le \Vert \pi_2(\mu) \Vert$ for every $\mu \in M^b(G)$.
\end{enumerate}
Indeed, \ref{iiDERem-tradweqGC*} is equivalent to $\Vert \pi_1(x) \Vert \le \Vert \pi_2(x) \Vert$ 
for every $x \in C^*_{\rm max}(G)$, by density of $B$;
this last property is equivalent to $\textnormal{C*ker}(\pi_2) \subset \textnormal{C*ker}(\pi_1)$
and hence to \ref{iDERem-tradweqGC*}, by Proposition~\ref{tradweqGC*}.
The equivalence of \ref{iiiDERem-tradweqGC*} and \ref{iiDERem-tradweqGC*} 
follows from the fact that $L^1(G)$ is dense in $M^b(G)$ for the weak* topology.
\end{rem}

\index{Primitive dual! $1$@of a topological group}
\index{$c1$@$\Pri(G)$ primitive dual of the topological group $G$}
The canonical surjective map (\ref{candualprim}) of Section \ref{PrimIdealSpace}
can now be written
$$
\kappa^{\rm d}_{\rm prim} \, \colon \, \widehat G
\twoheadrightarrow \Pri(G) = \Pri(C^*_{\rm max}(G)) ,
\hskip.5cm
\pi \mapsto \textnormal{C*ker}(\pi) .
\leqno(\kappa 1)
$$

\begin{rem}
% 8.B.6
\label{remonC*algebraofG}

\index{Fell topology}
(1)
The Fell topology on the dual $\widehat G$ of a LC group $G$
can be defined in at least two ways:
in terms of functions of positive type, as in Section \ref{S-DefUnitD},
or in terms of the Jacobson topology on $\Pri(G) = \Pri(C^*_{\rm max}(G))$
and the map ($\kappa 1$). 
For more on the equivalence of several definitions of this topology,
see \cite[\S~3]{Dixm--C*}.

\vskip.2cm

(2)
The map of ($\kappa 1$) is injective (and therefore bijective) in important cases, 
see Glimm Theorem \ref{ThGlimm}, 
but not in general, see Section \ref{noninjectivity}.

\vskip.2cm

(3)
The C*-kernel of $\pi$ should not be confused 
with the kernel of the group homomorphism $\pi \,\colon G \to \U(\Hi_\pi)$.
For example, if $\lambda$ is the left regular representation of $G$ as above,
so that $\lambda \,\colon G \to \U(L^2(G))$ is injective in all cases, 
we have $\textnormal{C*ker}(\lambda) = \{0\}$
if and only if $G$ is amenable
(Theorem \ref{HulanickiReiter}).

\vskip.2cm

(4) The bijection of Proposition \ref{bijrepGrepC*}
respects also quasi-equivalence, as introduced in Remark~\ref{Rem-Cor-FacRepPrimitive}.

\vskip.2cm

(5)
There are some topological groups $G$ that are not locally compact
for which a C*-algebra $C^*_{\rm max}(G)$ can nevertheless be defined,
and be of some use.
See \cite{StVo--75} for $\U(\infty)$, 
an appropriate inductive limit of the compact Lie groups $\U(n)$.
C*-algebras of groups that are not locally compact
will not be mentioned again in the present book.
\end{rem}

\begin{exe}
% 8.B.7
\label{Exa-C*-AlgAmenableGroup}
Let $G$ be a locally compact group that is amenable.
We repeat that the left regular representation $\lambda$ of $G$
establishes an isomorphism from $C^*_{\rm max}(G)$ onto $C^*_{\lambda}(G)$;
in other words, every representation of $G$
is weakly contained in the regular representation.

\vskip.2cm

(1) 
In case $G$ is an abelian LC group,
the Fourier transform $\mathcal F \,\colon L^2(G) \to L^2(\widehat G)$ 
is a unitary isomorphism of Hilbert spaces
(the $L^2$ spaces are defined in terms of Haar measures
on $G$ and $\widehat G$
that are associated with each other),
and $\mathcal F C^*_{\lambda}(G) \mathcal F^{-1} = C^0(\widehat G)$;
recall that $C^0(\widehat G)$ stands for
the convolution algebra of complex-valued continuous functions
on $\widehat G$ vanishing at infinity.
In particular, the C*-algebras $C^*_{\rm max}(G)$ and $C^0(\widehat G)$ 
are isomorphic.

\vskip.2cm

(2)
When $G$ is compact, 
its C*-algebra $C^*_{\rm max}(G)$ is isomorphic to the restricted product
$\prod^{c_0}_{\pi \in \widehat G} {\rm End}(\Hi_\pi)$
of the family of full matrix algebras,
indexed by the dual of $G$.
Recall that $1 \le \dim \Hi_\pi < \infty$ for all $\pi \in \widehat G$,
so that ${\rm End}(\Hi_\pi)$ is isomorphic to
the algebra ${\rm M}_n(\C)$ of complex $n$-by-$n$ matrices
for $n = \dim \Hi_\pi$;
see for instance Proposition 3.4 in \cite{Will--07}.
\par

\index{Cartesian product of C*-algebras}
\index{C*-algebra! $2$@cartesian product}
Recall also that,
given a family $(A_\iota)_{\iota \in I}$ of C*-algebras,
their Cartesian product $\prod_{\iota \in I} A_\iota$
is the C*-algebra constituted of
$(x_\iota)_{\iota \in I}$ with $x_\iota \in A_\iota$ for all $\iota \in I$
and $\sup_{\iota \in I} \Vert x_\iota \Vert < \infty$,
and their \textbf{restricted product} $\prod^{c_0}_{\iota \in I}A_\iota$
is the C*-algebra constituted of
$(x_\iota)_{\iota \in I}$ such that, for all $\varepsilon > 0$,
there is a cofinite subset $J_\varepsilon$ in $I$ such that
$\Vert x_\iota \Vert < \varepsilon$ for all $\iota \in J_\varepsilon$.
\index{Restricted product of C*-algebras}
\index{C*-algebra! $3$@restricted product}
\par

For example, $C^*_{\rm max}(\T)$ is isomorphic to the abelian C*-algebra
$c_0(\Z)$ of sequences of complex numbers $(c_n)_{n \in \Z}$ such that
$\lim_{\vert n \vert \to \infty} c_n = 0$.
Note that $c_0(\Z)$ is a restricted product of a countable infinite number
of copies of the C*-algebra $\C$.
\end{exe}

\begin{rem}
% 8.B.8
\label{Rem-C*-AlgebraDiscreteGroup}
(1)
Let $G$ be a discrete group.
Then $C^*_\pi(G)$ is a C*-algebra with unit for every representation $\pi$ of $G$. 
There is some kind of converse:
a LC group $G$ is discrete if and only if
its reduced C*-algebra $C^*_{\lambda}(G)$ has a unit,
if and only if its maximal C*-algebra $C^*_{\rm max}(G)$
has a unit \cite{Miln--71}.

\vskip.2cm 

(2)
Specific descriptions of maximal C*-algebras have been worked out
for a few groups,
often as spaces of sections of bundles or continuous fields.
Besides abelian and compact groups, 
here is a small sample of groups for which such descriptions are known, 
with some of the relevant references:
\begin{enumerate}
\item[--]
the special linear groups
$\SL_2(\C)$ \cite{Fell--61} and $\SL_2(\R)$ \cite{Mili--71, Vale--84};
\item[--]
the group ${\rm Spin} (4, 1)$, 
universal cover of the de Sitter group $\SO(4, 1)_0$ \cite{BoMa--77};
\item[--]
the motion group $\R^n \rtimes \SO(n)$ \cite{AbEL--11};
\item[--]
the Heisenberg group $H(\R)$, as in Remark \ref{dualHeisLie};
see \cite{LuTu--11}, and also \cite{Gorb--80};
\index{Heisenberg group! $4$@$H(\R)$}
\item[--]
the affine group
$\begin{pmatrix} \K^\times & \K \\ 0 & 1 \end{pmatrix}$
over a non-discrete locally compact field $\K$,
for example over $\Q_p$ \cite{Rose--76}.
\end{enumerate}
\end{rem}

\section{Functorial properties of group C*-algebras}
% Section 8.C
\label{SS:Functioriality}

For a locally compact group $G$,
the assignments $G \rightsquigarrow C^*_{\rm max}(G)$
and $G\rightsquigarrow C^*_\lambda(G)$
do not have good functorial properties with respect to general continuous group homomorphisms.
Here is an example which illustrates this fact.

\begin{exe}
% 8.C.1
\label{Exa-NonFunctorial}
Let $G$ be a non discrete LC group
and let $G_{\rm disc}$ be the group $G$ equipped with the discrete topology.
The identity map $G_{\rm disc} \to G$ is a continuous isomorphism. 
However, $C^*_{\rm max}(G)$ is not isomorphic to a quotient of $C^*_{\rm max}(G_{\rm disc})$.
Indeed, $C^*_{\rm max}(G_{\rm disc})$ is a unital C*-algebra whereas,
as mentioned above (Remark~\ref{Rem-C*-AlgebraDiscreteGroup}),
$C^*_{\rm max}(G)$ has no unit.
\par

Similarly, if $\Gamma$ is a discrete subgroup of a non-discrete LC group $G$,
then $C^*_{\rm max}(\Gamma)$, a C*-algebra that has a unit,
does not embed in $C^*_{\rm max}(G)$,
a C*-algebra that does not have one.
\end{exe}
 
The reason for the lack of functoriality of $G \rightsquigarrow C^*_{\rm max}(G)$ 
lies in the bad functorial properties of the assignment $G \rightsquigarrow L^1(G)$;
more precisely, given a continuous group homomorphism $G_1 \to G_2$,
there is in general no naturally associated homomorphism
between the convolution algebras $L^1(G_1)$ and $L^1(G_2)$.
The same holds for $C^c(G_1)$ and $C^c(G_2)$.
\par
 
We will see that $G \rightsquigarrow C^*_{\rm max}(G)$ and, to a lesser extend,
$G \rightsquigarrow C^*_\lambda(G)$, do have satisfactory functorial properties
with respect to some special continuous group homomorphisms.

\subsection
{On the functor $G \rightsquigarrow M^b(G)$}
% subsection 8.C.a
\label{SSS:FunctiorialityMeasure}

Let $G$ be a LC group.
Recall from Section~\ref{C*algLCgroup}
that the measure algebra $M^b(G)$ of complex Radon measures on $G$
is a Banach $*$-algebra
and that every representation $\pi$ of $G$ induces a $*$-representation 
$\pi \,\colon M^b(G) \to \Li (\Hi_\pi)$.
\par

As a preliminary step, we show that the assignment $G \rightsquigarrow M^b(G)$
has optimal functorial properties. 

\begin{prop}
% 8.C.2
\label{Prop-GroupHomMeasure}
Let $G_1, G_2$ be locally compact groups
and $\theta \,\colon G_1 \to G_2$ a continuous homomorphism.
For $\mu \in M^b(G_1)$,
denote by $\theta_*(\mu) \in M^b(G_2)$ the image of $\mu$ under $\theta$.
\par

The map $\theta_* \,\colon M^b(G_1) \to M^b(G_2)$
is a continuous homomorphism of $*$-algebras such that $\Vert \theta_*\Vert \le 1$.
\end{prop}

\begin{proof}
Observe that, for $\mu\in M^b(G_1)$,
the measure $\theta_*(\mu)\in M^b(G_2)$ is determined by 
$$
\int_{G_2} \varphi(y) d(\theta_*(\mu)) (y) \, = \, \int_{G_1} \varphi(\theta(x)) d\mu(x)
\hskip.5cm \text{for all} \hskip.2cm
\varphi \in C^b(G_2).
$$
It is clear that the map $\theta_*$ is a continuous linear map and that $\Vert \theta_*\Vert \le 1$.
\par

For $\mu_1, \mu_2 \in M^b(G_1)$ and $\varphi\in C^b(G_2)$, we have
$$
\begin{aligned}
\int_{G_2} \varphi(y) d(\theta_*(\mu_1\ast \mu_2)) (y)
& \, = \, \int_{G_1} \varphi(\theta(x)) d(\mu_1\ast \mu_2) (x) \\
& \, = \, \int_{G_1}\int_{G_1} \varphi(\theta(xz)) d\mu_{1}(x)d\mu_{2}(z) \\
& \, = \, \int_{G_1}\int_{G_1} \varphi(\theta(x)\theta(z)) d\mu_{1}(x)d\mu_{2}(z) \\
& \, = \, \int_{G_2}\int_{G_2} \varphi(uv)d(\theta_*(\mu_1)) (u) d(\theta_*(\mu_2)) (v) \\
& \, = \, \int_{G_2} \varphi (y) d(\theta_*(\mu_1)\ast \theta_*(\mu_2)) (y).
\end{aligned}
$$
Therefore $\theta_*(\mu_1\ast \mu_2) = \theta_*(\mu_1)\ast \theta_*(\mu_2)$.
\par

For $\mu\in M^b(G_1)$ and $\varphi\in C^b(G_2)$, we have 
$$
\begin{aligned}
\int_{G_2} \varphi(y) d(\theta_*(\mu^*)) (y)
& \, = \, \int_{G_1} \varphi(\theta(x)) d\mu^*(x) \\
& \, = \, \int_{G_1} \overline{\varphi(\theta(x^{-1}))} d\mu(x) \\
& \, = \, \int_{G_2} \overline{\varphi(z^{-1})} d\theta_*(\mu) (z) \\
& \, = \, \int_{G_2} \varphi(y) d( \theta_*(\mu))^*(y).
\end{aligned}
$$
Therefore $\theta_*(\mu^*) = (\theta_*(\mu))^*$.
This proves the claim.
\end{proof}
 
\begin{rem}
% 8.C.3
\label{Rem-MeasureAlgebra}
Proposition~\ref{Prop-GroupHomMeasure} shows that the assignment
$G \rightsquigarrow M^b(G)$,
is a functor from the category of LC groups 
to the category of Banach $*$-algebras.
Since $(\beta \circ \alpha)_* = \beta_* \circ \alpha_*$,
when $\alpha \,\colon G_1 \to G_2$ and $\beta \,\colon G_2 \to G_3$
are continuous homomorphisms between LC groups $G_1, G_2, G_3$,
this functor is covariant.
\end{rem}

We need to study the action of representations of $G_2$
on measures in $\theta_*(M_b(G_1))$.

\begin{prop}
% 8.C.4
\label{Prop-RepMeasureGroupAlg}
Let $\theta \,\colon G_1 \to G_2$ and $\theta_* \,\colon M^b(G_1) \to M^b(G_2)$
be as in Proposition~\ref{Prop-GroupHomMeasure}.
\par

For every $\mu \in M^b(G_1)$ and every representation $\pi$ of $G_2$, we have
$$
\pi(\theta_*(\mu)) \, = \, (\pi \circ \theta)_* (\mu).
$$
\end{prop}

\begin{proof}
Let $\pi$ be a representation of $G_2$ and $\mu \in M^b(G_1)$. Then
$$
\pi(\theta_*(\mu_1))
\, = \, \int_{G_2} \pi(y) d\theta_*(\mu) (y)
\, = \, \int_{G_1} \pi(\theta(x)) d\mu(x)\\
\, = \, (\pi \circ \theta)_* (\mu).
$$
This proves the claim. 
\end{proof}

For our study of the assignment $G \rightsquigarrow C^*_{\rm max}(G)$,
it will be helpful to introduce C*-algebras associated to the Banach $*$-algebra $M^b(G)$,
as we did for the Banach *-algebra $L^1(G)$ in Section~\ref{C*algLCgroup}.

\begin{defn}
% 8.C.5
\label{Def-EnvC*Measure}
Let $G$ be a LC group.
The \textbf{enveloping C*-algebra} of $M^b(G)$, denoted by $M^*_{\rm max}(G)$,
is the completion of $M^b(G)$ for the norm
$$
\Vert \mu\Vert_{\rm max} \, := \, \sup_\pi \Vert \pi(\mu) \Vert_{\Li (\Hi_\pi)},
$$
where the supremum is taken over
all equivalence classes of cyclic (or, equivalently, irreducible) representations $\pi$ of $G$.
\par

The \textbf{reduced enveloping C*-algebra} of $M^b(G)$, denoted by $M^*_\lambda(G)$,
is the completion of $M^b(G)$ for the norm 
$$
\Vert \mu\Vert_{\lambda} \, := \, \Vert \lambda(\mu) \Vert.
$$
Equivalently, $M^*_\lambda(G)$ could be defined as the norm closure 
of $\lambda(M^b(G))$ in $\Li (L^2(G))$.
\index{Enveloping C*-algebra of measure algebra}
\index{Reduced Enveloping C*-algebra of measure algebra}
\end{defn}

The next proposition shows in particular
that the assignment $G \rightsquigarrow M^*_{\rm max}(G)$ is a covariant functor
from the category of LC groups to the category of C*-algebras.

\begin{theorem}
% 8.C.6
\label{Theo-C*FunctMeasureAlg}
Let $G_1, G_2$ be locally compact groups
and $\theta \,\colon G_1 \to G_2$ a continuous homomorphism.
\par

The map $\theta_* \,\colon M^b(G_1) \to M^b(G_2)$
of Proposition~\ref{Prop-GroupHomMeasure}
extends to a morphism $\theta_* \,\colon M^*_{\rm max}(G_1) \to M^*_{\rm max}(G_2)$
of C*-algebras.
\end{theorem}

\begin{proof}
In the sequel, $\Vert \cdot \Vert_{\rm max}$ will denote
the norm on $M^*_{\rm max}(G_1)$ 
as well as the norm on $M^*_{\rm max}(G_2)$.
\par

Let $\mu \in M^b(G_1)$ and $\pi$ a representation of $G_2$.
By Proposition~\ref{Prop-RepMeasureGroupAlg}, we have
$\pi(\theta_*(\mu)) = (\pi \circ \theta)_* (\mu)$.
It follows from the definition of the norm $\Vert \cdot \Vert_{\rm max}$ that 
$$
\Vert \pi(\theta_*(\mu))\Vert \, = \, \Vert (\pi \circ \theta)_* (\mu)\Vert
\, \le \, \Vert \mu \Vert_{\rm max}.
$$
Therefore we have 
$$
\Vert \theta_*(\mu) \Vert_{\rm max} \, \le \, \Vert \mu \Vert_{\rm max}
\hskip.5cm \text{for every} \hskip.2cm
\mu\in M^b(G_1).
$$
Therefore, the map $\theta_* \,\colon M^b(G_1) \to M^b(G_2)$,
which is a homomorphism of $*$-algebras, 
extends to a morphism of C*-algebras $M^*_{\rm max}(G_1) \to M^*_{\rm max}(G_2)$.
\end{proof}

 \subsection
 {On the assignment $G \rightsquigarrow C^*_{\rm max}(G)$}
% subsection 8.C.b
\label{SSS:FunctiorialityMax}

Let $G$ be a LC group.
Recall that $C^c(G)$ is a $*$-algebra for the convolution product 
and the involution given by $f^*(x) = \Delta(x^{-1}) \overline{f(x^{-1})}$ for $f\in C^c(G)$.
Recall also that $C^c(G)$ may be viewed a $*$-subalgebra of the Banach $*$-algebra $M^b(G)$ 
by means of the map
$$
f \mapsto f(x) d\mu_G(x)
\hskip.5cm \text{for} \hskip.2cm
f \in C^c(G).
$$ 
\par

Let $G_1, G_2$ be LC groups
and $\theta \,\colon G_1 \to G_2$ a continuous group homomorphism.
Let $\theta_* \,\colon M^b(G_1) \to M^b(G_2)$ be the associated map
as in Proposition~\ref{Prop-GroupHomMeasure}.
Assume that $\theta$ is an \emph{open} homomorphism,
that is, $\theta(U)$ is open in $G_2$ for every open subset $U$ of $G_1$.
In Theorem \ref{Theo-C*-GroupHom},
we will show that, in this case, $\theta_*(C^c(G_1))$ is contained in $C^c(G_2)$ 
and that $\theta_* \,\colon C^c(G_1) \to C^c(G_2)$
extends to a morphism $C^*_{\rm max}(G_1) \rightarrow C^*_{\rm max}(G_2)$.
To this end, we need to treat separately three types of homomorphisms:
topological isomorphisms,
quotient maps,
and open injective homomorphisms.

\begin{prop}
% 8.C.7
\label{Prop-C*-algIsoGroup}
Let $G_1, G_2$ be LC groups with Haar measures $\mu_{G_1}, \mu_{G_2}$,
and let $\theta \,\colon G_1 \to G_2$ be an isomorphism of topological groups.
\begin{enumerate}[label=(\arabic*)]
\item\label{iDEProp-C*-algIsoGroup}
For an appropriate normalization of $\mu_{G_1}$ and $\mu_{G_2}$, 
we have $\theta_*(C^c(G_1)) = C^c(G_2)$ 
and $\theta_* \,\colon C^c(G_1) \to C^c(G_2)$ is an isomorphism of $*$-algebras.
\item\label{iiDEProp-C*-algIsoGroup}
The map $\theta_* \,\colon C^c(G_1) \to C^c(G_2)$
extends to an isomorphism $C^*_{\rm max}(G_1) \rightarrow C^*_{\rm max}(G_2)$
of C*-algebras.
\end{enumerate}
\end{prop}

\begin{proof}
Since $\theta$ is an isomorphism of topological groups,
it follows from Proposition~\ref{Prop-GroupHomMeasure}
that $\theta_* \,\colon M^b(G_1) \to M^b(G_2)$ is an isomorphism of $*$-algebras.
It is obvious that $\theta_*(\mu_{G_1})$ is a multiple of $\mu_{G_2}$
that $\theta_*(C^c(G_1)) = C^c(G_2)$,
and that the restriction of $\theta_*$ is a linear bijection $C^c(G_1) \to C^c(G_2)$.
We may normalize $\mu_{G_1}$ and $\mu_{G_2}$
so that $\mu_{G_2} = \theta_*(\mu_{G_1})$.
With this normalization, $\theta_* \,\colon C^c(G_1) \to C^c(G_2)$
is a morphism of $*$-algebras. This proves Item \ref{iDEProp-C*-algIsoGroup}.

\vskip.2cm

Let $f \in C^c(G_1)$.
Denoting by $\Vert \cdot \Vert_{\rm max}$ the norms on
$M^*_{\rm max}(G_1)$ and $M^*_{\rm max}(G_2)$,
we have, by Theorem~\ref{Theo-C*FunctMeasureAlg},
$$
\Vert \theta_*(f) \Vert_{\rm max} \, \le \, \Vert f \Vert_{\rm max}.
$$
Since the same arguments apply also to $\theta^{-1}$, it follows that 
$$
\Vert \theta_*(f) \Vert_{\rm max} \, = \, \Vert f \Vert_{\rm max},
\hskip.5cm \text{for every}\hskip.2cm
f \in C^c(G_1).
$$
Therefore the $*$-isomorphism $\theta_* \,\colon C^c(G_1) \rightarrow C^c(G_2)$ 
extends to an isomorphism
$$
C^*_{\rm max}(G_1) \rightarrow C^*_{\rm max}(G_2)
$$
of C*-algebras.
\end{proof}

Next, we examine the case of quotient maps. 

\begin{prop}
% 8.C.8
\label{Prop-C*-algQuotientGroup}
Let $G$ be a locally compact group, $N$ a closed normal subgroup of $G$,
and $\theta \,\colon G \to G/N$ be the canonical epimorphism.
\begin{enumerate}[label=(\arabic*)]
\item\label{iDEProp-C*-algQuotientGroup}
For a proper normalization of the Haar measures $\mu_G, \mu_N$ and $\mu_{G/N}$,
we have $\theta_*(C^c(G)) = C^c(G/N)$ 
and $\theta_* \,\colon C^c(G) \to C^c(G/N)$ is a surjective homomorphism of $*$-algebras.
\item\label{iiDEProp-C*-algQuotientGroup}
The map $\theta_* \,\colon C^c(G )\twoheadrightarrow C^c(G/N)$
extends to a surjective morphism $C^*_{\rm max}(G) \twoheadrightarrow C^*_{\rm max}(G/N)$
of C*-algebras. 
\end{enumerate}
\end{prop}

\begin{proof}
Fix a Haar measure $\mu_N$ on $N$ and let $P \,\colon C^c(G) \to C^c(G/N)$ be the linear 
map defined by 
$$
P(f) (xN) = \int_N f(xy) d\mu_N(y)
\hskip.5cm \text{for} \hskip.2cm
f \in C^c(G), x \in G.
$$
It is a surjective linear map \cite[Lemma B.1.2]{BeHV--08}.
For a proper normalization of the Haar measures $\mu_G$ and $\mu_{G/N}$ of $G$ and $G/N$,
we have the so-called Weil's formula
$$
\int_G f(x) d\mu_G(x) = \int_{G/N} P(f) (xN) d\mu_{G/N}(xN)
\hskip.5cm \text{for all}\hskip.2cm
f \in C^c(G)
$$
(see \cite[Chap.~II, \S 9]{Weil--40}).
% or \cite[Appendix B]{BeHV--08}).
\par

Let $f \in C^c(G)$ and $\varphi\in C^c(G/N)$. 
Since $P((\varphi \circ \theta)f) = \varphi P(f)$, we have 
$$
\int_{G/N} \varphi(xN)P(f) (xN) d\mu_{G/N} \, = \,
\int_{G} \varphi (\theta(x)) f(x) d\mu_G(x),
$$
by Weil's formula. This shows that 
$$
\theta_*(f) \, = \, P(f)
\hskip.5cm \text{for all} \hskip.2cm
f \in C^c(G),
$$
that is, the restriction of $\theta_*$ to $C^c(G)$ coincides with $P$.
This implies that 
$$
\theta_*(C^c(G)) \, = \, P(C^c(G)) = C^c(G/N).
$$
Moreover, $\theta_* \,\colon C^c(G) \to M^b(G/N)$
is a homomorphism of $*$-algebras, by Proposition~\ref{Prop-GroupHomMeasure}.
Therefore $\theta_* \,\colon C^c(G) \to C^c(G/N)$ is a surjective homomorphism of $*$-algebras.
This proves Item \ref{iDEProp-C*-algQuotientGroup}.

\vskip.2cm

Let $f \in C^c(G)$. By Theorem~\ref{Theo-C*FunctMeasureAlg}, we have
$$
\Vert \theta_*(f) \Vert_{\rm max} \, \le \, \Vert f \Vert_{\rm max},
$$
where, as before, $\Vert \cdot \Vert_{\rm max}$
denotes the norm on the enveloping C*-algebras of the measure algebras.
Therefore $\theta_* \,\colon C^c(G) \to C^c(G/N)$ extends to a morphism
$C^*_{\rm max}(G) \rightarrow C^*_{\rm max}(G/H)$ of C*-algebras;
this morphism is surjective, since $ \theta_*(C^c(G)) = C^c(G/N)$
and since $ C^c(G/N)$ is dense in $C^*_{\rm max}(G/H)$.
\end{proof}

Finally, we treat the inclusion of open subgroups.
We give first a preliminary result on the inclusion homomorphism $\theta$ of a closed,
not necessarily open, subgroup $H$ of a LC group $G$.
For $f\in C_c(H)$, observe that 
$\theta_*(f)$ is the measure $f(x) d\mu_H(x)$ on $G$.
Observe also that, for every representation $\pi$ of $G$, we have 
$$
\pi(\theta_*( \mu)) \, = \, (\pi \vert_H) (\mu)
\hskip.5cm\text{for all} \hskip.2cm
\mu \in M^b(H).
$$

\begin{prop}
% 8.C.9
\label{Pro-InjectionC*-Max}
Let $G$ be a LC group, $H$ a closed subgroup of $G$,
and $\theta \,\colon H \to G$ the canonical inclusion.
\par

The injective $*$-homomorphism $\theta_* \,\colon C^c(H) \to M^b(G)$ extends to a morphism
$\theta_* \,\colon C^*_{\rm max}(H) \rightarrow M^*_{\rm max}(G)$ of C*-algebras.
Moreover, the following properties are equivalent:
\begin{enumerate}[label=(\roman*)]
\item\label{iDEPro-InjectionC*-Max}
$\theta_* \,\colon C^*_{\rm max}(H) \rightarrow M^*_{\rm max}(G)$ is injective;
\item\label{iiDEPro-InjectionC*-Max}
for every representation $\sigma$ of $H$, there exists a representation $\pi$ of $G$ such that 
$\sigma$ is weakly contained in $\pi \vert_H$.
\end{enumerate}
\end{prop}

\begin{proof}
By Theorem~\ref{Theo-C*FunctMeasureAlg},
the map $\theta_* \,\colon C^c(H) \to M^b(G)$ extends to a morphism
$\theta_* \,\colon M^*_{\rm max}(H) \rightarrow M^*_{\rm max}(G)$ of C*-algebras.
Since obviously $C^*_{\rm max}(H)$ may be viewed as C*-subalgebra of $M^*_{\rm max}(H)$,
we obtain by restriction a morphism
$\theta_* \,\colon C^*_{\rm max}(H) \rightarrow M^*_{\rm max}(G)$ of C*-algebras.

\vskip.2cm

Assume that \ref{iDEPro-InjectionC*-Max} holds.
Then $\theta_*$ is isometric, that is, 
$$
\Vert \theta_*(f)\Vert_{\rm max} \, = \, \Vert f \Vert_{\rm max}
\hskip.5cm \text{for all} \hskip.2cm
f \in C^c(H).
\leqno{(*)}
$$
Let $\sigma$ be a representation of $H$.
Set $\pi_{\text max} := \bigoplus_{\pi \in \widehat G} \pi$.
Since
$$
\Vert \mu \Vert_{\rm max} \, = \, \Vert \pi_{\text max}(\mu) \Vert
\hskip.5cm \text{for all} \hskip.2cm
\mu \in M^b(G),
$$
it follows from $(*)$ that
$$
\Vert \sigma(f) \Vert \le \Vert \pi_{\text max}(\theta_*(f)) \Vert
\, = \, \Vert (\pi_{\text max} \vert_H) (f) \Vert 
\hskip.5cm \text{for all} \hskip.2cm
f \in C^c(H);
$$
this shows that $\sigma$ is weakly contained in $\pi_{\text max} \vert_H$.

\vskip.2cm

Conversely, assume that \ref{iiDEPro-InjectionC*-Max} holds. Let 
$$
\sigma_{\text max} \, := \, \bigoplus_{\sigma \in \widehat{H}} \sigma.
$$
There exists a representation $\pi$ of $G$ such that 
$\sigma_{\text max}$ is weakly contained in $\pi \vert_H$. 
Therefore, we have 
$$
\Vert f\Vert_{\rm max}
\, = \, \Vert \sigma_{\text max}(f) \Vert
\, \le \, \Vert (\pi \vert_H) (f) \Vert
\, = \, \Vert \pi(\theta_*(f)) \Vert
\, \le \, \Vert \theta_*(f) \Vert_{\rm max},
$$
and hence $\Vert f\Vert_{\rm max} = \Vert \theta_*(f)\Vert_{\rm max}$
for all $f\in C^c(H)$. 
So, $\theta_* \,\colon C^*_{\rm max}(H) \rightarrow M^*_{\rm max}(G)$ is injective.
\end{proof}
 
We give two classes of subgroups $H$ 
for which the map $\theta_* \,\colon C^*_{\rm max}(H) \rightarrow M^*_{\rm max}(G)$
as in Proposition~\ref{Pro-InjectionC*-Max} is injective.

\begin{cor}
% 8.C.10
\label{Cor-Pro-InjectionC*-Max}
Let $G$ be a LC group, $H$ a closed subgroup of $G$,
and $\theta \,\colon H \to G$ the canonical injection.
\begin{enumerate}[label=(\arabic*)]
\item\label{iDECor-Pro-InjectionC*-Max}
Assume either that $H$ is open or that $H$ is amenable.
Then the morphism $\theta_* \,\colon C^*_{\rm max}(H) \rightarrow M^*_{\rm max}(G)$ 
is injective.
\item\label{iiDECor-Pro-InjectionC*-Max}
Assume that $H$ is open.
For an appropriate normalization of $\mu_G$ and $\mu_H$,
the image of $C^*_{\rm max}(H)$ under $\theta_*$
is a subalgebra of $C^*_{\rm max}(G)$.
\end{enumerate}
\end{cor}

\begin{proof}
To show Item \ref{iDECor-Pro-InjectionC*-Max},
it suffices to check that Condition \ref{iiDEPro-InjectionC*-Max} of Proposition~\ref{Pro-InjectionC*-Max} is satisfied,
when $H$ is either open or amenable.

\vskip.2cm

$\bullet$ \emph{First case.}
Assume that $H$ is open. Let $\sigma$ be a representation of $H$. 
Consider the induced representation $\pi := \Ind_H ^G \sigma$. 
Then $\sigma$ is contained in $\pi \vert_H$ (see Definition~\ref{openeasierthanclosed}).

\vskip.2cm

$\bullet$ \emph{Second case.} 
Assume that $H$ is amenable. Let $\sigma$ be a representation of $H$. 
Then $\sigma$ is weakly contained in 
$\lambda_H$. Since $\lambda_G \vert_H$ is weakly equivalent to $\lambda_H$, 
it follows that $\sigma$ is weakly contained in $\lambda_G \vert_H$.

\vskip.2cm

To show Item \ref{iiDECor-Pro-InjectionC*-Max},
note first that $C^c(H)$ is naturally a subspace of $C^c(G)$, since $H$ is open in $G$.
Choose then $\mu_G$ and $\mu_H$ such that $\mu_G \vert_H = \mu_H$.
% so that $C^c(H)$ is a subalgebra of $C^c(G)$, for the convolution products.
Then, for every $f \in C^c(H)$, we have 
$$
\theta_*(f) \, = \, f(x) d\mu_H(x) \, = \, f(x) d\mu_G(x) ,
$$
hence $\theta_*$ is an inclusion of $C^c(H)$ as a subalgebra of $C^c(G)$.
% This shows that, under the usual identification of $C^c(G)$ with a subalgebra of $M^b(G)$, 
% we have $\theta_* (C^c(H)) \subset C^c(G)$. 
Since the closure of $C^c(G)$ in $M^*_{\rm max}(G)$ coincides with 
$C^*_{\rm max}(G)$, the claim follows.
\end{proof}

\begin{rem}
% 8.C.11
\label{Rem-Pro-InjectionC*-Max}
(1)
Condition \ref{iiDEPro-InjectionC*-Max} of Proposition~\ref{Pro-InjectionC*-Max}
was introduced in \cite[Page 442]{Fell--64}
as ``Property (WP3)" and studied as a weak version of Frobenius reciprocity for finite groups.
It also appeared in \cite[page 209]{Rief--74} and \cite{BeVa--95} in a context similar to ours,
with the multiplier C*-algebra $M(C^*_{\rm max}(G))$ of $C^*_{\rm max}(G)$
replacing the algebra $M^*_{\rm max}(G)$.

\vskip.2cm

(2)
Condition \ref{iiDEPro-InjectionC*-Max} from Proposition~\ref{Pro-InjectionC*-Max}
does not always hold.
Indeed, this condition fails for $G = \SL_3(\C)$
and a subgroup $H$ isomorphic to $\SL_2(\C)$,
as shown in \cite[Theorem 6.1]{Fell--64} and \cite[Remark 1.13 (i)]{BeLS--92}
(see our Remark \ref{remdefind}).
% with an incomplete proof which was corrected in \cite{BeLS--92}.
The failure of Condition \ref{iiDEPro-InjectionC*-Max} was established in \cite{BeVa--95}
for any lattice $H$ in a higher rank simple Lie group
$G$ as well as for many lattices in a rank one simple Lie group.
\end{rem}

We are now ready to state the final result
concerning the assignment $G \rightsquigarrow C^*_{\rm max}(G)$.

\begin{theorem}
% 8.C.12
\label{Theo-C*-GroupHom}
Let $G_1, G_2$ be locally compact groups with Haar measures $\mu_{G_1}, \mu_{G_2}$
and let $\theta \,\colon G_1 \to G_2$ be an \emph{open} continuous homomorphism.
\begin{enumerate}[label=(\arabic*)]
\item\label{iDETheo-C*-GroupHom}
We have $\theta_*(C^c(G_1)) \subset C^c(G_2)$ and,
for appropriate normalizations of $\mu_{G_1}$ and $\mu_{G_2}$, 
the map $\theta_* \,\colon C^c(G_1) \to C^c(G_2)$ is a morphism of $*$-algebras.
\item\label{iiDETheo-C*-GroupHom}
This map $\theta_* \,\colon C^c(G_1) \to C^c(G_2)$
extends to a morphism of C*-algebras
\hfill\par\noindent
$\theta_* \,\colon C^*_{\rm max}(G_1) \rightarrow C^*_{\rm max}(G_2)$.
\item\label{iiiDETheo-C*-GroupHom}
$\theta_* \,\colon C^*_{\rm max}(G_1) \rightarrow C^*_{\rm max}(G_2)$ is surjective
if and only if $\theta \,\colon G_1 \to G_2$ is surjective.
\item\label{ivDETheo-C*-GroupHom}
$\theta_* \,\colon C^*_{\rm max}(G_1) \rightarrow C^*_{\rm max}(G_2)$ is injective
if and only if $\theta \,\colon G_1 \to G_2$ is injective
\end{enumerate}
\end{theorem}

\begin{proof}
Let $N\subset G_1$ be the kernel and $H\subset G_2$ the image of $\theta$.
We can write 
$$
\theta \, = \, \theta_3 \circ \theta_2 \circ \theta_1,
$$
where
\begin{itemize}
\setlength\itemsep{0em}
\item
$\theta_1 \,\colon G_1 \to G_1/ N$ is the canonical epimorphism, 
\item
$\theta_2 \,\colon G_1/N \to H$ is the group isomorphism induced by $\theta$,
\item
$\theta_3 \,\colon H \to G_2$ is the canonical injection.
\end{itemize}

Since $\theta$ is an open map, $H$ is an open subgroup of $G_2$ and 
$\theta_2 \,\colon G_1/N \to H$ is an isomorphism of topological groups. 
Therefore Items \ref{iDETheo-C*-GroupHom} and \ref{iiDETheo-C*-GroupHom}
follow from Proposition~\ref{Prop-C*-algIsoGroup}, 
Proposition~\ref{Prop-C*-algQuotientGroup}, and Corollary~\ref{Cor-Pro-InjectionC*-Max}.
\par

Items \ref{iiiDETheo-C*-GroupHom} and \ref{ivDETheo-C*-GroupHom}
are immediate consequences of the following observations:
\begin{itemize}
\setlength\itemsep{0em}
\item
the map $C^*_{\rm max}(G_1) \rightarrow C^*_{\rm max}(G_1/N)$ 
induced by $\theta_1$ is surjective;
\item
the maps $C^*_{\rm max}(G_1/N) \rightarrow C^*_{\rm max}(H)$
and $C^*_{\rm max}(H) \rightarrow C^*_{\rm max}(G)$
induced by $\theta_2$ and $\theta_3$ 
are injective;
\item
the map $C^*_{\rm max}(G_1) \rightarrow C^*_{\rm max}(G_1/N)$ 
is injective if and only if $N = \{e\}$, that is, if and only if $\theta$ is injective;
\item
the map $C^*_{\rm max}(H) \rightarrow C^*_{\rm max}(G_2)$ is surjective 
if and only if $H = G_2$, that is, if and only if $\theta$ is surjective.
\end{itemize}
\vskip-.8cm
\end{proof}
 
If $\Gamma$ is a discrete group, observe that
$C^c(\Gamma)$ coincides with the group algebra $\C[\Gamma]$.
We will identify $\Gamma$ in the obvious way as a subset of $\C[\Gamma]$.
Under this identification, if $\theta \,\colon \Gamma_1 \to \Gamma_2$ is a homomorphism
between discrete groups $\Gamma_1$ and $\Gamma_2$,
then $\theta_* \,\colon \C[\Gamma_1] \to \C[\Gamma_2]$ coincides with
the linear extension of $\theta$ to $\C[\Gamma_1]$.
\par

The following immediate consequence of Theorem~\ref{Theo-C*-GroupHom} 
shows in particular that the assignment
$\Gamma \rightsquigarrow C^*_{\rm max}(\Gamma)$
is a covariant functor from the category of discrete groups to the category of C*-algebras.

\begin{cor}
% 8.C.13
\label{Cor-Theo-C*-GroupHom}
Let $\Gamma_1, \Gamma_2$ be discrete groups.
Every group homomorphism $\theta \,\colon \Gamma_1 \to \Gamma_2$
extends to a morphism
$\theta_* \,\colon C^*_{\rm max}(\Gamma_1) \to C^*_{\rm max}(\Gamma_2)$
of C*-algebras.
\end{cor}

\subsection
{On the assignment $G \rightsquigarrow C^*_{\lambda}(G)$}
% 8.C.c
\label{SSS:FunctiorialityRed}

As we now discuss, reduced group C*-algebras have less functorial properties
than maximal group C*-algebras.
\par

Let $G_1, G_2$ be LC groups and $\theta \,\colon G_1 \to G_2$ an open continuous
homomorphism. Recall from Theorem~\ref{Theo-C*-GroupHom} that $\theta$ induces 
a homomorphism $\theta_* \,\colon C^c(G_1) \to C^c(G_2)$ of $*$-algebras,
upon an appropriate normalization of $\mu_{G_1}$ and $\mu_{G_2}$.
The following result gives a characterization of the homomorphisms $\theta$
inducing a morphism between the respective reduced C*-algebras.

\begin{theorem}
% 8.C.14
\label{Theo-RedC*-Group}
Let $G_1,G_2$ be LC groups 
and $\theta \,\colon G_1 \to G_2$ a continuous \emph{open} homomorphism.
Assume that the Haar measures $\mu_{G_1}, \mu_{G_2}$ are normalized
as in Theorem~\ref{Theo-C*-GroupHom}~\ref{iDETheo-C*-GroupHom}. 
The following properties are equivalent:
\begin{enumerate}[label=(\roman*)]
\item\label{iDETheo-RedC*-Group}
the induced homomorphism $\theta_* \,\colon C^c(G_1) \to C^c(G_2)$ 
extends to a morphism $C^*_{\lambda}(G_1) \to C^*_{\lambda}(G_2)$ of C*-algebras;
\item\label{iiDETheo-RedC*-Group}
the kernel of $\theta$ is an amenable subgroup of $G_1$.
\end{enumerate}
\end{theorem}

\begin{proof} 
Recall that $\lambda_{G}$ denotes the left regular representation of a LC group $G$.

\vskip.2cm

$\bullet$ \emph{First step.}
We claim that $\theta_* \,\colon C^c(G_1) \to C^c(G_2)$ 
extends to $C^*_{\lambda}(G_1)$ if and only if
the representation $\lambda_{G_2} \circ \theta$ of $G_1$
is weakly contained in $\lambda_{G_1}$.
\par
 
Indeed, let $f\in C^c(G_1)$. By Proposition~\ref{Prop-RepMeasureGroupAlg}, we have
$$
\lambda_{G_2}(\theta_*(f) \,) = \, (\lambda_{G_2} \circ \theta) (f).
$$
This shows that the inequality 
$$
\Vert \lambda_{G_2}(\theta_*(f)) \Vert \, \le \, \Vert \lambda_{G_1}(f) \Vert
$$
is equivalent to the inequality
$$
\Vert (\lambda_{G_2} \circ \theta) (f) \Vert \, \le \, \Vert \lambda_{G_1}(f) \Vert.
$$
Therefore $\theta_*$ extends from $C^c(G_1)$ to $C^*_{\lambda}(G_1)$
if and only if 
$$
\Vert (\lambda_{G_2} \circ \theta) (f) \Vert \, \le \, \Vert \lambda_{G_1}(f) \Vert
\hskip.5cm\text{for every} \hskip.2cm
f \in C^c(G_1),
$$
that is, if and only if $\lambda_{G_2} \circ \theta$ is weakly contained in $\lambda_{G_1}$.

\vskip.2cm

$\bullet$ \emph{Second step.}
Assume that $\theta \,\colon G_1 \to G_2$ is an isomorphism.
We claim that $\lambda_{G_2} \circ \theta$ equivalent to $\lambda_{G_1}$.
\par

Indeed, define a linear map $U \,\colon L ^2(G_1) \to L^2(G_2)$ by 
$$
(U\xi) (y) \, = \, \xi(\theta^{-1}(y))
\hskip.5cm\text{for } \hskip.2cm
\xi \in L^2(G_1), \ y \in G_2.
$$
By the normalization of $\mu_{G_1}$ and $\mu_{G_2}$, we have 
$\theta_*(\mu_{G_1}) = \mu_{G_2}$ (see the proof of Proposition~\ref{Prop-C*-algIsoGroup}).
It follows, as is easily checked, that $U$ is an isomorphism of Hilbert spaces.
Moreover, we have
$$
\begin{aligned}
(\lambda_{G_2}(\theta(g)) U\xi) (y)
\, &= \, (U\xi) (\theta(g^{-1}) y)\\
\, &= \, \xi (g^{-1} \theta^{-1}(y))\\ 
\, &= \, (\lambda_{G_1}(g) \xi) (\theta^{-1}(y))\\
\, &= \, (U\lambda_{G_1}(g)\xi) (y),
\end{aligned} 
$$
for every $\xi \in L^2(G_1), g\in G_1$ and $y \in G_2$. 
So, $U$ intertwines $\lambda_{G_2} \circ \theta$ and $\lambda_{G_1}$.

\vskip.2cm

$\bullet$ \emph{Third step.}
Assume that $G_1$ is an open subgroup of $G_2$ and that 
$\theta \,\colon G_1 \to G_2$ is the canonical injection. 
We claim that $\lambda_{G_2} \circ \theta$ is weakly equivalent to $\lambda_{G_1}$.
\par

Indeed, this holds since we have
$\lambda_{G_2} \circ \theta = \lambda_{G_2} \vert_{G_1}$ in this case.

\vskip.2cm

$\bullet$ \emph{Fourth step.}
Assume that $G_2 = G_1/N$ for a closed
normal subgroup $N$ of $G_1$ and that 
$\theta \,\colon G_1 \to G_1/N$ is the canonical epimorphism.
We claim that $\lambda_{G_2} \circ \theta$ is weakly contained in $\lambda_{G_1}$
if and only if $N$ is amenable.
\par

Indeed, observe first that $\lambda_{G_2} \circ \theta$
is equivalent to the induced representation
$\Ind_{N}^{G_1}1_{N}$. 
\par

Assume that $N$ is amenable. Then $1_N$ is weakly contained
in $\lambda_{N}$ and hence, by continuity of induction, 
$\Ind_{N}^{G_1}1_{N}$ is weakly contained in $\Ind_{N}^{G_1}\lambda_{N}$,
which is equivalent to $\lambda_{G_1}$.
Therefore $\lambda_{G_2} \circ \theta$ is weakly contained in $\lambda_{G_1}$.
\par

Conversely, assume that $\lambda_{G_2} \circ \theta$ is weakly contained in $\lambda_{G_1}$.
Then, by continuity of restriction, the representation $(\lambda_{G_2} \circ \theta) \vert_N$ 
is weakly contained in $\lambda_{G_1} \vert_N$. 
Observe that $(\lambda_{G_2} \circ \theta) \vert_N$ is 
equivalent to a multiple of the trivial representation $1_N$. Since $\lambda_{G_1} \vert_N$
is weakly equivalent to $\lambda_{N}$,
it follows that $1_N$ is weakly contained in $\lambda_{N}$.
Therefore $N$ is amenable.

\vskip.2cm

$\bullet$ \emph{Fifth step.}
Let $\theta \,\colon G_1 \to G_2$ be an arbitrary open homo\-morphism.
Denote by $N$ the kernel of $\theta$ and by $H$ the image of $\theta$.
As in the proof of Theorem~\ref{Theo-C*-GroupHom}, 
we can write 
$$
\theta= \theta_3 \circ \theta_2 \circ \theta_1,
$$
where
$\theta_1 \,\colon G_1 \to G_1/ N$ is the canonical epimorphism, 
$\theta_2 \,\colon G_1/N \to H$ is the group isomorphism induced by $\theta$,
and $\theta_3 \,\colon H \to G_2$ is the canonical injection.
\par

By the three first steps, $\theta_2$ and $\theta_3$
induce morphisms $C^*_{\lambda}(G_1/N) \to C^*_{\lambda}(H)$
and $C^*_{\lambda}(H) \to C^*_{\lambda}(G_2)$.
The fourth step shows that $\theta_1$ induces a morphism
$C^*_{\lambda}(G_1) \to C^*_{\lambda}(G_1/N)$
if and only $N$ is amenable.
\end{proof}

The following corollary is an immediate consequence
of Theorem~\ref{Theo-RedC*-Group};
it should be compared with Corollary~\ref{Cor-Theo-C*-GroupHom}, in which we established that 
$$
\Gamma \rightsquigarrow C^*_{\rm max}(\Gamma)
$$
is a functor on the category of discrete groups to that of C*-algebras.
Recall that we view a discrete group $\Gamma$ as a subset of its group algebra $\C[\Gamma]$,
which is itself a subalgebra of $C^*_{\lambda}(\Gamma)$.

\begin{cor}
% 8.C.15
\label{Cor-Theo-RedC*-Group}
 Let $\Gamma_1, \Gamma_2$ be discrete groups and
 $\theta \,\colon \Gamma_1 \to \Gamma_2$ a homomorphism.
The following properties are equivalent:
\begin{enumerate}[label=(\roman*)]
\item\label{iDECor-Theo-RedC*-Group}
$\theta$ extends to a morphism
$\theta_* \,\colon C^*_{\lambda}(\Gamma_1) \to C^*_{\lambda}(\Gamma_2)$ of C*-algebras;
\item\label{iiDECor-Theo-RedC*-Group}
the kernel of $\theta$ is an amenable subgroup of $\Gamma_1$.
\end{enumerate}
\end{cor}

\section
{Second-countable and $\sigma$-compact LC groups}
% Section 8.D
\label{C*kernelsprimitive}

Some properties of representations of $\sigma$-compact locally compact groups
are well-documented in the literature for second-countable LC groups only.
For example, the fact that C*-kernels of factor representations of $\sigma$-compact LC groups
are primitive is apparently known to experts (see Line 5 of \cite{Pogu--83}); 
however, the only references of which we are aware deal with second-countable LC groups
(see \cite[Corollaire 3]{Dixm--60a},
and also Corollaries~\ref{Cor-FacRepPrimitive} and \ref{Cor-FacRepPrimitiveG}).
Observe that a $\sigma$-compact LC group need not be second-countable.
For example, an uncountable Cartesian product of non-trivial compact groups
is compact and is not second-countable.
\par

We indicate in this section how, in several situations,
the case of $\sigma$-compact LC group can be reduced to that of second-countable LC groups.
The first tool for this reduction is the Kakutani--Kodaira theorem;
for a proof, see the original article \cite{KaKo--44}, or \cite[2.B.6]{CoHa--16}.

\begin{theorem}[Kakutani--Kodaira]
% 8.D.1
\label{Kakutani--Kodaira}
Let $G$ be a $\sigma$-compact locally compact group.
\par

For every neighbourhood $U$ of $e$ in $G$,
there exists a compact normal subgroup $K$ of $G$ 
such that $K \subset U$ and $G/K$ is second-countable.
\end{theorem}

Using the notion of projective limits of topological groups
(see \cite[Chap.~III, \S~7]{BTG1--4}), 
Theorem~\ref{Kakutani--Kodaira} shows that every $\sigma$-compact LC group
is the projective limit of second-countable LC groups.
\par

The second tool for the reduction 
is the following crucial lemma; a similar result appears in \cite[Proposition 2.2]{Moor--72}, under more restrictive assumptions and with a different proof.

\begin{lem}
% 8.D.2
\label{Lem-FactQuot}
Let $G$ be a topological group 
and $\mathcal{K}$ a family of compact normal subgroups of $G$ 
with the following property:
for every neighbourhood $U$ of $e$ in $G$, 
there exists $K \in \mathcal K$ with $K \subset U$.
Let $\pi$ be a factor representation of $G$. 
\par

Then $\pi$ factorizes through a quotient $G/K$ for some $K \in \mathcal{K}$ .
\end{lem}

\begin{proof}
Let $\xi$ be a unit vector in $\Hi_\pi$. 
Since $\pi$ is continuous, there exists a neighbourhood $U$ of $e$ in $G$ such
that 
$$
\Vert \pi(u) \xi-\xi\Vert <\sqrt{2}
\hskip.5cm \text{for all} \hskip.2cm
u \in U.
$$
Therefore there exists $K \in \mathcal{K}$ such that 
$$
\Vert \pi(g) \xi-\xi \Vert<\sqrt{2}
\hskip.5cm \text{for all} \hskip.2cm
g \in K.
$$
Since $K$ is a group, this implies that 
the space $\Hi_\pi^K$ of $\pi(K)$-invariant vectors in $\Hi_\pi$
is no7 zero (see \cite[Proposition 1.1.5]{BeHV--08}). 
Since $K$ is normal in $G$, 
the subspace $\Hi_\pi^K$ is $\pi(G)$-invariant.
Moreover, $\Hi_\pi^K$ is obviously invariant under the commutant
$\pi(G)'$ of $\pi(G)$. 
Therefore the orthogonal projection 
$p \,\colon \Hi_\pi \to \Hi_\pi^K$ belongs to $\pi(G)' \cap \pi(G)''$,
that is, to the centre of the von Neumann algebra $ \pi(G)''$.
Since $\pi$ is factorial and $\Hi_\pi^K \ne \{0\}$, we have $p = \operatorname{id}_{\Hi_\pi}$,
that is, $\Hi_\pi = \Hi_\pi^K$; this means that $\pi$ factorizes through $G/K$.
\end{proof}

Lemma~\ref{Lem-FactQuot}, combined with Theorem~\ref{Kakutani--Kodaira},
allows us to relate the dual and the quasi-dual a $\sigma$-compact LC group
to the duals and quasi-duals of appropriate second-countable LC groups.

\begin{prop}
% 8.D.3
\label{FacKK}
Let $G$ a $\sigma$-compact locally compact group,
and let $\mathcal K_G$ denote the set of compact normal subgroups $K$ in $G$
such that $G/K$ is second-countable, ordered by inclusion.
\begin{enumerate}[label=(\arabic*)]
\item\label{iDEFacKK}
For every factor representation $\pi$ of $G$,
there exists a compact normal subgroup $K \in \mathcal K_G$ 
such that $\pi$ factorizes through $G/K$.
\item\label{iiDEFacKK}
For every factor representation $\pi$ of $G$, 
the C*-kernel $\textnormal{C*ker}(\pi)$ is a primitive ideal of $C^*_{\rm max}(G)$.
\item\label{iiiDEFacKK}
The dual $\widehat G$ is homeomorphic to the inductive limit
$$
\widehat G = \varinjlim_{K \in \mathcal K_G} \hskip.1cm \widehat{G/K} .
$$
\item\label{ivDEFacKK}
Similarly, the quasi-dual $\QD(G)$ 
is Borel-isomorphic to the inductive limit
$$
\QD(G) \, = \, \varinjlim_{K \in \mathcal K_G} \QD(G/K) .
$$
\end{enumerate}
\end{prop}

\begin{proof}
Claim \ref{iDEFacKK} follows from Lemma~\ref{Lem-FactQuot}
and Theorem~\ref{Kakutani--Kodaira}.

\vskip.2cm

\ref{iiDEFacKK}
Let $\pi$ be a factor representation of $G$.
By \ref{iDEFacKK}, there exists $K \in \mathcal K_G$
and a factor representation $\pi'$ of $G/K$
such that $\pi = \pi' \circ p$, where $p \,\colon G \to G/K$ stands for
the canonical epimorphism.
Let $p_* \,\colon C^*_{\rm max}(G) \twoheadrightarrow C^*_{\rm max}(G/K)$
be the morphism induced on the maximal C*-algebras
(see Proposition~\ref{Prop-C*-algQuotientGroup}).
By Corollary~\ref{Cor-FacRepPrimitiveG} and Proposition \ref{tradweqGC*},
there exists an irreducible representation $\rho'$ of $G/K$
such that $\textnormal{C*ker}(\rho') = \textnormal{C*ker}(\pi')$.
Consider the representation given by the composition
$\rho \,\colon C^*_{\rm max}(G) \overset{p_*}{\twoheadrightarrow} C^*_{\rm max}(G/K)
\overset{\rho'}{\to} \Li (\Hi_{\rho'})$.
Then
$$
\textnormal{C*ker}(\pi) \, = \,
(p_*)^{-1} \left( \textnormal{C*ker}(\pi') \right) \, = \,
(p_*)^{-1} \left( \textnormal{C*ker}(\rho') \right) \, = \,
\textnormal{C*ker}(\rho)
$$
and the ideal $\textnormal{C*ker}(\pi)$ of $C^*_{\rm max}(G)$ is primitive.

\vskip.2cm

\ref{iiiDEFacKK} \& \ref{ivDEFacKK}
These claims follow from \ref{iDEFacKK}. Let us spell out the homeomorphism.
For $K \in \mathcal K_G$, the space $\widehat{G/K}$ with its Fell topology
can be seen as a closed subspace of $\widehat G$.
The Fell topology on $\widehat G$ coincides with 
the inductive limit topology, 
that is, the one for which a subset $A$ in $\widehat G$ is open
if and only if $A \cap \widehat{G/K}$ is
open for all $K \in \mathcal K_G$.
\par

\index{Mackey--Borel structure}
Similarly, for each $K \in \mathcal K_G$,
let the quasi-dual of $G/K$
be given its Mackey--Borel structure.
For its own Mackey--Borel structure,
a subset $A$ in $\QD(G)$ is Borel 
if and only if the intersection 
$A \cap \QD(G/K)$
is Borel for all $K \in \mathcal K_G$.
\end{proof}

As a corollary, we show that, for a $\sigma$-compact locally compact group $G$, 
the quotient map $\widehat G \twoheadrightarrow \Pri(G)$ has a natural extension 
to the quasi-dual $\QD(G)$.

\begin{cor}
% 8.D.4
\label{kapp3extendskapp1}
For a $\sigma$-compact locally compact group $G$, the map
\begin{equation}
\label{eqq/qd/spG}
\tag{$\kappa$3}
\kappa^{\rm qd}_{\rm prim} \, \colon \,
\QD(G) \twoheadrightarrow \Pri(G) ,
\hskip.5cm
\rho \mapsto \textnormal{C*ker}(\rho) 
\end{equation}
is an extension of the canonical map
$\kappa^{\rm d}_{\rm prim} \,\colon \widehat G \twoheadrightarrow \Pri(G)$
of Section \ref{PrimIdealSpace}.
\end{cor}

\begin{proof}
Let $\rho$ be a factor representation of $G$.
On the one hand,
every factor representation of $G$ quasi-equivalent to $\rho$
has the same C*-kernel as $\rho$
(Corollary~\ref{QE-WEQ} and Proposition~\ref{tradweqGC*});
on the other hand, 
$\textnormal{C*ker}(\rho)$ is a primitive ideal
(Proposition \ref{FacKK}~\ref{iiDEFacKK}).
\end{proof}

\section[The central character of a representation]
{The central character of a representation, of a primitive ideal}
% Section 8.E
\label{Section-CentralCharacter}

\index{Projective kernel}
Let $G$ be a topological group.
The \textbf{projective kernel} of a representation 
$(\pi, \Hi)$ of $G$ is the normal subgroup $\Pker \pi$ defined by
$$
\Pker \pi \, = \, \{ g \in G \mid \pi(g) = \chi_\pi(g) \mathrm{Id}_{\Hi} 
\hskip.2cm \text{for some number} \hskip.2cm
\chi_\pi(g) \in \T \} .
$$
Observe that $\Pker \pi$ is a normal closed subgroup of $G$ containing $\ker \pi$, 
that the function $\chi_\pi \,\colon g \mapsto \chi_\pi(g)$ is a unitary character of $\Pker \pi$,
and that
$$
\ker( \pi \,\colon G \to \U(\Hi) ) \, = \, \ker( \chi_\pi \,\colon \Pker \pi \to \T ) .
$$

\begin{prop}
% 8.E.1
\label{Prop-ProjKer}
Let $G$ be a topological group and $\pi$ a factor representation of $G$.
Let $\chi_\pi \,\colon \Pker \pi \to \T$ be as above.
\begin{enumerate}[label=(\arabic*)]
\item\label{iDEProp-ProjKer}
The centre $Z$ of $G$ is contained in $\Pker \pi$,
so that there is a unitary character $\psi_\pi := \chi_\pi \vert_Z$
in $\widehat Z$.
\item\label{iiDEProp-ProjKer}
Let $H$ be a closed subgroup
which contains the commutator subgroup of $G$
(and therefore which is normal in $G$)
and $p \,\colon G \to G/(\ker \pi \cap H)$ the canonical projection.
% utile avec $H \ne G$ juste avant \ref{Theo-PrimIdealTwoStepNilpotent}.
\par
Then $\Pker \pi = p^{-1}(\overline Z)$,
where $\overline Z$ is the centre of $G/(\ker \pi \cap H)$. 
\item\label{iiiDEProp-ProjKer}
Let $\rho$ be any representation of $G$ that is weakly equivalent to $\pi$.
\par
Then $\Pker \rho = \Pker \pi$ and $\chi_\rho = \chi_\pi$. 
\end{enumerate}
\end{prop}

\begin{proof}
\ref{iDEProp-ProjKer}
For $g \in Z$, the operator $\pi(g)$ belongs to the centre of the factor $\pi(G)''$ 
and hence is a scalar operator.

\vskip.2cm

\ref{iiDEProp-ProjKer}
Let $g \in G$. Then $\pi(g)$ commutes with $\pi(g')$ for all $g' \in G$
if and only if $\pi(g)$ belongs to the centre of the factor $\pi(G)''$,
that is, if and only if $\pi(g)$ is a scalar operator. 
Therefore, we have $g \in \Pker \pi$ if and only if
$$
\pi([g,g']) \, = \, [\pi(g), \pi(g')] \, = \, \mathrm{Id}_{\Hi}
\hskip.5cm \text{for all} \hskip.2cm
g' \in G , 
$$
that is, if and only if $[g,g'] \in \ker \pi$ for all $g' \in G$.
\par

Since $H$ contains $[G,G]$,
the last condition holds if and only if $[g,g'] \in \ker \pi \cap H$ for all $g' \in G$,
i.e., if and only if $g \in p^{-1}(\overline Z)$.

\vskip.2cm

\ref{iiiDEProp-ProjKer}
The restriction $\rho \vert_{\Pker \pi}$ is weakly contained in $\pi \vert_{\Pker \pi}$,
and therefore in $\chi_\pi$. 
This implies that every normalized function of positive type
associated to $\rho \vert_{\Pker \pi}$ coincides with $\chi_\pi$,
and therefore that
$$
\rho(g) \, = \, \chi_\pi(g) \mathrm{Id}_{\Hi} 
\hskip.5cm \text{for all} \hskip.2cm 
g \in \Pker \pi.
$$
Therefore $\Pker \pi \subset \Pker \rho$ and $\psi_\rho \vert_{\Pker \pi} = \chi_\pi$.
Similarly, we have $\Pker \rho \subset \Pker \pi$
and $\chi_\pi \vert_{\Pker \rho} = \chi_\rho$.
\end{proof}

\index{Central character}
\index{Character! $3$@central}
By Proposition~\ref{Prop-ProjKer}, we can define two maps
$$
\widehat G \, \to \, \widehat Z
\hskip.5cm \text{and} \hskip.5cm
\QD(G) \, \to \, \widehat Z 
$$
by associating to $\pi \in \widehat G$ 
or $\pi \in \QD(G)$
the unitary character $\psi_\pi \in \widehat Z$ such that
$$
\pi(z) \, = \, \psi_\pi(z) \mathrm{Id}_{\Hi_\pi} 
\hskip.5cm \text{for all} \hskip.2cm 
z \in Z .
$$
called the \textbf{central character} of $\pi$.
Note that this central character is the restriction to $Z$
of the unitary character $\chi_\pi \,\colon \Pker \pi \to \T$ previously defined.
Recall that, when $\pi$ is irreducible, the central character $\chi_\pi$
has already been define, see \ref{defcentralchar}.
% Note also that central characters have already appeared in \ref{Section-IrrRepTwoStepNil}.
\par

By Proposition \ref{Prop-ProjKer}~\ref{iiiDEProp-ProjKer}, 
the map $\widehat G \to \widehat Z$ factorizes through the map 
$$
\kappa^{\rm d}_{\rm prim} \, \colon \, \widehat G \twoheadrightarrow \Pri(G)
$$
of (\ref{candualprim}) in Section \ref{C*algLCgroup}.
Similarly, when $G$ is a $\sigma$-compact LC group,
the map $\QD(G) \to \widehat Z$
factorizes through the map
$$
\kappa^{\rm qd}_{\rm prim} \, \colon \, 
\QD(G) \twoheadrightarrow \Pri(G),
$$
of (\ref{eqq/qd/spG}) in \ref{C*kernelsprimitive}.
We can therefore speak of the central character of a primitive ideal.

\section
{Characterization of type I groups: Glimm theorem}
% Section 8.F
\label{SectionGlimm}

This section is devoted to the statement of Glimm theorem,
which gives several characterizations of LC groups of type I. 

\subsection
{CCR and GCR, or liminal and postliminal C*-algebras}
% subsection 8.F.a
\label{SSectionPreGlimm}

We introduce a class of C*-algebras which play a central role
in Glimm theorem.
\par

Given a Hilbert space $\Hi$, recall that $\mathcal K (\Hi)$ denotes
the ideal of compact operators in $\Li (\Hi)$.
\par

Let $A$ be a C*-algebra. For a representation $\pi$ of $A$,
consider the two-sided ideal 
$$
\pi^{-1}(\mathcal K (\Hi_\pi)) \, = \, \{ a \in A \mid \pi(a) \hskip.2cm \text{is compact}\}
$$
of $A$~;
note that it depends on the equivalence class of $\pi$ only.
Define ${\rm CCR} (A)$ as 
the intersection of $\pi^{-1}(\mathcal K (\Hi_\pi))$
over all equivalence classes of irreducible representations $\pi$ of $A$.

\begin{defn}
% 8.F.1
\label{Def-GCR-C*}
Let $A$ be C*-algebra $A$.
\par

$A$ is defined to be \textbf{CCR} or \textbf{liminal} if,
for every irreducible representation $\pi \,\colon A \to \Li (\Hi_\pi)$,
the image $\pi(A)$ consists of compact operators on $\Hi_\pi$. 
\par

$A$ is defined to be \textbf{GCR} or \textbf{postliminal} 
if ${\rm CCR} (A/J) \ne \{0\}$ for every closed two-sided ideal $J$ in $A$.
\index{GCR! $2$@C*-algebra}
\index{Postliminal! C*-algebra}
\index{CCR! $3$@C*-algebra}
\index{Liminal! $2$@C*-algebra}
\end{defn}

\begin{rem}
% 8.F.2
\label{Rem-Def-GCR-C*}
(1)
The acronym CCR refers to ``completely continuous representations'',
the other term ``liminal", indicating that these algebras are the first to consider
(etymologically, ``liminal'' comes from the Latin ``limen'', meaning threshold).
The acronym GCR refers to ``generalized CCR algebra''.

\vskip.2cm

(2)
We comment on the relationship between the notions 
of GCR groups and CCR groups we defined in Section~\ref{S:ClassTypeI-GCRGroups}
and the notions of CCR and GCR C*-algebras we have just introduced.
\begin{itemize}
\setlength\itemsep{0em}
\item It is clear that a LC group $G$ is a CCR group if and only if 
its maximal C*-algebra $C^*_{\rm max}(G)$ is CCR;
\item it is also clear that, if $C^*_{\rm max}(G)$ is GCR, 
then $G$ is a GCR group;
\item
when $G$ is second-countable,
the fact that $C^*_{\rm max}(G)$ is GCR if $G$ is a GCR group holds;
however, this is a deep fact which is part of Glimm theorem~\ref{ThGlimm} below.
\end{itemize}

\vskip.2cm

(3)
The notion of type I group from Section~{SectionTypeI} makes sense
in the context of C*-algebras: a C*-algebra $A$ is of type I if every representation $\pi$ of $A$,
the von Neumann algebra $\pi(A)''$ is of type I. 
Glimm Theorem~\ref{ThGlimm}, in its general formulation,
shows that the class of separable C*-algebras
coincides with the class of separable GCR C*-algebras.
\end{rem}

\subsection
{Glimm theorem}
% subsection 8.F.b

In the fundamental 1961 article of Glimm, 
only one out of nine theorems is formulated 
for second-countable locally compact groups,
and the other theorems are formulated for separable C*-algebras.
Here, we choose a formulation for groups,
more precisely for $\sigma$-compact LC groups.

\begin{theorem}
% 8.F.3
\label{ThGlimm}
\index{Glimm theorem}
For a $\sigma$-compact locally compact group $G$,
the following conditions are equivalent:
\begin{enumerate}[label=(\roman*)]
\item\label{aDEThGlimm}
$G$ is type I, as defined in Section \ref{SectionTypeI};
\index{Type I! $3$@group}
\item\label{bDEThGlimm}
$G$ is a GCR group, as defined in Section~\ref{S:ClassTypeI-GCRGroups};
\item\label{cDEThGlimm}
$C^*_{\rm max}(G)$ is a GCR C*-algebra;
\item\label{dDEThGlimm}
the map $\kappa^{\rm d}_{\rm prim} \,\colon \widehat G \twoheadrightarrow \Pri(G)$
of (\ref{candualprim}) in \ref{PrimIdealSpace} and \ref{C*algLCgroup}
is a homeomorphism;
\item\label{eDEThGlimm}
the Mackey--Borel structure on $\widehat G$ is countably separated;
\item\label{fDEThGlimm}
the Mackey--Borel structure on $\widehat G$ is standard;
\item\label{gDEThGlimm}
the Mackey--Borel structure on $\widehat G$ is the same as 
the Borel structure defined by the Fell topology.
\end{enumerate}
\end{theorem}
% no condition (*) $\widehat G$ is a T$_0$ space -- essentially contained in (d).
% no condition(**) The map $\kappa^{\rm d}_{\rm qd} \,\colon \widehat G \to \QD(G)$
% of (\ref{eqq/d/qd}) in \ref{SectionQuasidual}
% is an isomorphism of Borel space -- essentially contained in (a) = type I.

We give below several comments on Glimm theorem; for the complete
proof, we refer to the original article \cite{Glim--61a} and \cite[Chap.~9]{Dixm--C*}. 

% \vskip.2cm
% 
% \noindent
\newpage

\textbf{Comments on Theorem~\ref{ThGlimm}}

\vskip.2cm

(1)
The equivalence between the properties \ref{aDEThGlimm} to \ref{fDEThGlimm} listed above
is proved in \cite{Glim--61a} in the context of separable C*-algebras;
applied to the C*-algebra $C^*_{\rm max}(G)$, 
this implies Theorem~\ref{ThGlimm} for a second-countable LC group $G$.
Concerning the equivalence of these properties with property \ref{gDEThGlimm},
see Comment (8) below.

\vskip.2cm

(2)
The equivalence of all properties, except \ref{cDEThGlimm},
carries over from second-countable LC groups
to the more general case of a $\sigma$-compact LC group $G$,
by Proposition \ref{FacKK}.
For the equivalence of \ref{cDEThGlimm} with the other properties,
see Comment (11) below.

\vskip.2cm

(3)
The equivalence \ref{aDEThGlimm} $\iff$ \ref{fDEThGlimm}
was conjectured by Mackey \cite[Page 163]{Mack--57}.

\vskip.2cm

(4)
Property \ref{aDEThGlimm} is equivalent to the fact 
that every factor representation of $G$ is a multiple of an irreducible representation.
See Theorem~\ref{Theo-TypIFactorRep} and Proposition~\ref{Prop-TypeIFac}.

\vskip.2cm

(5)
The implication \ref{bDEThGlimm} $\implies$ \ref{aDEThGlimm} is proved
in Theorem~\ref{Theo-GCR-Group}.

\vskip.2cm

(6)
The implication \ref{dDEThGlimm} $\implies$ \ref{aDEThGlimm}
is a consequence of Corollary~\ref{Cor-FacRepPrimitiveG-NonTypI};

\vskip.2cm

(7)
The implications \ref{cDEThGlimm} $\implies$ \ref{bDEThGlimm} 
and \ref{fDEThGlimm} $\implies$ \ref{eDEThGlimm} are obvious.

\vskip.2cm

(8)
The proofs in \cite{Glim--61a} are independent of previous work,
mainly by Dixmier, Fell and Kaplansky,
in which the equivalence of several of the properties listed above were established,
in the context of separable C*-algebras:
\vskip.1cm

the implication \ref{cDEThGlimm} $\implies$ \ref{aDEThGlimm}
is an early result from \cite{Kapl--51b};
\vskip.1cm

the equivalence \ref{cDEThGlimm} $\iff$ \ref{dDEThGlimm}
was proved in \cite{Dixm--60a};
\vskip.1cm

the equivalences
\ref{cDEThGlimm} $\iff$ \ref{eDEThGlimm} $\iff$
\ref{fDEThGlimm} $\iff$ \ref{gDEThGlimm}
was proved in \cite{Dixm--60c};
\vskip.1cm

the implication \ref{dDEThGlimm} $\implies$ \ref{fDEThGlimm}
was established in \cite[Theorem 4.1]{Fell--60b}.

\vskip.2cm

(9)
In view of the previous comments,
the main novelty in \cite{Glim--61a} is the proof
of the implication \ref{aDEThGlimm} $\implies$ \ref{cDEThGlimm},
apart from the independent proofs it
offers for all the other implications.

\vskip.2cm

(10)
A conceptually new proof of \ref{aDEThGlimm} $\implies$ \ref{fDEThGlimm} 
is given in \cite{Effr--65}.

\vskip.2cm

(11)
Let $G$ be a locally compact group, $\sigma$-compact or not.
It is known that Properties
\ref{aDEThGlimm}, \ref{bDEThGlimm}, 
and \ref{cDEThGlimm} of Theorem \ref{ThGlimm}
are equivalent; see
% \cite[3.1.6]{Dixm--C*} and 
\cite{Saka--66, Saka--67}.

\vskip.2cm

(12)
Recall from Theorem \ref{Theo-SeparationPropertiesDual}~\ref{iDETheo-SeparationPropertiesDual}
that the conditions of Theorem \ref{ThGlimm} are moreover equivalent
to the fact that the dual $\widehat G$ is a T$_0$ space.

\vskip.2cm

(13)
Considering that groups of type I enjoy many good properties which other groups lack,
Kirillov has suggested that
a second-countable locally compact group
should be called \textbf{tame} if it is type I
and \text{wild} otherwise \cite[see Section 8.4]{Kiri--76}.
\index{Tame group}

\vskip.2cm

Let $G$ be a locally compact which is not type I.
Then $G$ has a factor representation of type II or a factor representation of type III.
In fact, when $G$ is second-countable, it was shown in \cite{Glim--61a}
that $G$ has factor representations of both types.

\begin{theorem}
% 8.F.4
\label{ThGlimm2}
Let $G$ be a locally compact group which is not of type I.
\begin{enumerate}[label=(\arabic*)]
\item\label{iDEThGlimm2}
The group $G$ has factor representations of type III.
\item\label{iiDEThGlimm2}
If $G$ is $\sigma$-compact, then
$G$ has factor representations of type II$_\infty$.
\item\label{iiiDEThGlimm2}
Let $\Gamma$ be a countable group that is not of type I.
Then $\Gamma$ has a factor representation of type II$_1$.
\end{enumerate}
\end{theorem}

\begin{proof}[On the proof]
For \ref{iDEThGlimm2},
see \cite{Glim--61a} in case $G$ is second-countable
and \cite{Saka--66} in general.
% \cite[9.5.9]{Dixm--C*}
\par

Claim \ref{iiDEThGlimm2}
follows from \cite{Glim--61a} and \cite{Mare--75},
by translating from separable C*-algebras to $\sigma$-compact groups.
% (see also \cite[9.5.4]{Dixm--C*}) 
\par

For \ref{iiiDEThGlimm2},
it follows from Theorem \ref{discretTypes}
that the regular representation of $\Gamma$ has a part of type II, 
say $\lambda^{II}_\Gamma$.
In the central decomposition 
$$
\lambda^{II}_\Gamma \, = \, \int^\oplus_{\QD(\Gamma)} \pi_\omega d\mu(\omega) 
$$
(see Theorem \ref{thmDirectIntFact}),
$\pi_\omega$ is a factor representation of type II$_1$ 
for almost every $\omega \in \QD(\Gamma)$
\cite[Chap.~II, \S~5, no~2]{Dixm--vN}.
\end{proof}

\section
{The von Neumann algebra of a group representation}
% Section 8.G
\label{vNalgLCgroup}

In Section \ref{S:QE-Factor-VNAlgebras},
we have associated two von Neumann algebras, $\pi(G)''$ and $\pi(G)'$,
to any representation $\pi$ of a topological group $G$. 
Let us record another way to define these when $G$ is locally compact.
\par

\index{Algebras! $3$@convolution algebras, $C^c(G)$, $C^0(G)$, $L^1(G)$}
Let $G$ be a locally compact group
and $\pi$ a representation of $G$ in a Hilbert space $\Hi_\pi$.
Recall that $C^c(G)$ stands for
the convolution algebra of complex-valued continuous functions
on $G$ with compact support.
We denote again by $\pi$ each of the corresponding morphisms from 
$C^c(G)$, $L^1(G)$, $M^b(G)$, and $C^*_{\rm max}(G)$
to $\Li (\Hi_\pi)$.
Since $C^c(G)$, $L^1(G)$, $C^*_{\rm max}(G)$ have approximate units,
the following is elementary to check; see \cite[13.3.5]{Dixm--C*}.

\begin{prop}
% 8.G.1
\label{vNpiGallsame}
Let $\pi$ be a representation of a locally compact group $G$.
\par

With the notation as above,
$\pi(G)$, $\pi(C^c(G))$, $\pi(L^1(G))$, and $\pi(C^*_{\rm max}(G))$
all generate the same von Neumann algebra:
$$
\pi(G)'' \, = \, \pi(C^c(G))'' \, = \, \pi(L^1(G))'' \, = \, \pi(M^b(G))''\, = \, \pi(C^*_{\rm max}(G))''.
$$
\end{prop}

Given a second-countable locally compact group $G$,
the second dual of $C^*_{\rm max} (G)$, say $\vN_{\rm max} (G)$,
has a natural structure of von Neumann algebra,
of which the normal $*$-representations are in natural one-to-one correspondence
with the unitary representations of $G$.
The algebra $\vN_{\rm max}(G)$ is the enveloping von Neumann algebra
of $C^*_{\rm max}(G)$ \cite[\S~12.1]{Dixm--C*}.
On this subject, we only quote two articles by Ernest \cite{Erne--64, Erne--65}.
\index{von Neumann algebra! enveloping}

\begin{rem}
% 8.G.2
In contrast to the equality $\pi(C^*_{\rm max}(G))'' = \pi(G)''$
of Proposition \ref{vNpiGallsame}, 
$C^*_\pi(G)$ \emph{need not} coincide 
with the sub involutive algebra $C^*(\pi(G))$
of $\Li (\Hi_\pi)$ generated by $\pi(G)$,
as the following examples illustrate.

\vskip.2cm

(1)
Note that $C^*(\lambda(G))$ is always a C*-algebra with unit,
whereas, when $G$ is not discrete, 
the reduced C*-algebra $C^*_\lambda(G)$ does not have a unit
(see Remark~\ref{Rem-C*-AlgebraDiscreteGroup}).

\vskip.2cm

(2)
Let $G$ be a LCA group.
Recall from Proposition \ref{Prop-BohrAbelian}
that the Bohr compactification $\Bohr(\widehat G)$
of the dual of $G$
can be identified with the dual $\widehat{G_{\rm disc}}$
of the group $G$ viewed as a discrete group.
\par

Then $C^*(\lambda(G))$ is isomorphic to the algebra
$C(\Bohr(\widehat G)) = C(\widehat{G_{\rm disc}})$
of continuous functions on $\widehat{G_{\rm disc}}$.
It is not isomorphic to $C^*_\lambda(G) \approx C^0(\widehat G)$
when $G$ is not discrete \cite{KoKa--43}.
% this has been extended from the left regular representation $\lambda$ of $G$
% to other representations \cite{Arve--66}.
More information on $C^*(\lambda(G))$ can be found in \cite{Bedo--94}.
\index{$a5$@$\approx$ isomorphism! of C*-algebras}
\end{rem}

\section
{Variants}
% Section 8.H
\label{SectionVariants}

Given a LC group $G$, there are other group algebras than $C^*_{\rm max}(G)$
and other primitive duals than $\Pri(G)$ that can be attached to $G$. 
One is the space $\Pri_*( L^1(G))$ consisting of kernels
of irreducible $*$-representations of $L^1(G)$
associated to irreducible representations of $G$.
Another one is the space $\Pri(L^1(G))$ of all primitive ideals of $L^1(G)$, 
i.e., of all kernels of simple $L^1(G)$-modules (in the algebraic sense). 
There are also spaces of maximal ideals,
$\Max (C^*_{\rm max}(G))$ and $\Max(L^1(G))$.
There is a natural injection
of $\Max (C^*_{\rm max}(G))$ into $\Pri(G)$,
and this is a bijection if and only if $\Pri(G)$ is a T$_1$ space.
\par

The only variant about which we add a comment
is the space $\Pri_*(L^1(G))$.
There is a natural map
$$
\Phi \, \colon \, \Pri(G) \to \Pri_*(L^1(G)), \hskip.2cm
\textnormal{C*ker}(\pi) \mapsto
L^1 \ker (\pi) = \textnormal{C*ker}(\pi) \cap L^1(G)
$$
that is surjective, and continuous for the Jacobson topologies.
In \cite[Th\'eo\-r\`eme 4]{Dixm--60b}, it is shown that this map
is injective (and therefore a homeomorphism)
when $G$ satisfies some ``Property (P)'',
in particular when $G$ is a connected Lie group that is either nilpotent or semisimple.
It is also known that $\Phi$ is a homeomorphism when
$G$ is a LC-group with polynomial growth \cite[Satz~2]{BLSV--78}.
\par

We end this section by an example showing that $\Phi$ need not be injective.

\begin{prop}
% 8.H.1
\label{resPSL2}
Let $\Gamma$ be a lattice in $G = \PSL_2(\R)$
and let $\Phi$ be as above.
\par

Then $\Phi \,\colon Prim (\Gamma) \to Prim_* (\ell^1(\Gamma))$
is not injective.
\end{prop}

\begin{proof}
Let $\pi_1$ and $\pi_2$ be non-equivalent representations of $G$ 
in the complementary series.
Then the restrictions $\pi_1 \vert_\Gamma$ and $\pi_2 \vert_\Gamma$ to $\Gamma$
are irreducible, by \cite[Proposition 2.5]{CoSt--91}. 
By \cite[Th\'eor\`eme 2]{BeHa--94}, we have on the one hand
$$
\textnormal{C*ker}\left( \pi_1 \vert_\Gamma \right) \ne 
\textnormal{C*ker}\left( \pi_2 \vert_\Gamma \right) ,
$$
and on the other hand
$$
\textnormal{C*ker}\left( \pi_1 \vert_\Gamma \right) \subset 
\textnormal{C*ker}( \lambda_\Gamma)
\hskip.5cm \text{and} \hskip.5cm
\textnormal{C*ker}\left( \pi_2 \vert_\Gamma \right) \subset 
\textnormal{C*ker}( \lambda_\Gamma),
$$
where $\lambda_\Gamma$ stands for
the left regular representations of $\Gamma$.
It follows that
$$
L^1 \ker (\pi_1 \vert_\Gamma) 
\subset L^1 \ker ( \lambda_\Gamma)
\hskip.5cm \text{and} \hskip.5cm
L^1 \ker (\pi_2 \vert_\Gamma) 
\subset L^1 \ker ( \lambda_\Gamma).
$$
Now $\lambda_\Gamma$ is faithful on $\ell^1(\Gamma)$,
i.e., $L^1 \ker ( \lambda_\Gamma) = \{0\}$. 
Therefore
$$
L^1 \ker (\pi_1 \vert_\Gamma) = 
L^1 \ker (\pi_2 \vert_\Gamma) = \{0\}
$$
and that concludes the proof.
\end{proof}

%-----------------------------------------------------------------------
% End of chapter 8
%-----------------------------------------------------------------------

\chapter{Examples of primitive duals} 
% Chapter 9
\label{ChapterPrimExa}

\emph{In this chapter, the primitive duals are determined
for the examples of groups
considered in Chapter \ref{Chapter-ExamplesIndIrrRep},
that is for $D_\infty$,
for two-step nilpotent discrete groups
and in particular the Heisenberg groups $H(\K)$ and $H(\Z)$, 
for $\Aff(\K)$, and $\BS(1, p)$, the lamplighter group, and $\GL_n(\K)$
when $\K$ is an infinite algebraic extension of a finite field. }
\par

\emph{The main tool we will use,
Proposition~\ref{Pro-WeakContSemiDirect},
is a weak containment result
for induced representations from an abelian normal subgroup.}

\vskip.5cm

The discussion for the first of our test examples is straightforward.
The infinite dihedral group $D_\infty$ is abelian-by-finite,
and is therefore a type I group.
In particular, the canonical projection from its dual to its primitive ideal space
is a homeomorphism.
The dual of $D_\infty$ has been described in Section \ref{SectionInfDiGroup}.

\section[A weak containment result]
{A weak containment result for induced representations}
% Section 9.A
\label{Sect-PrimIdealSemiDirect}

Proposition~\ref{Pro-WeakContSemiDirect} is implicitly contained in
\cite[Chap.~II, \S~2, Lemme 1]{Guic--63} and \cite[\S3, Lemme 2]{Guic--65}.
It is a particular case of the ``generalized Effros--Hahn conjecture",
a result about the structure of the primitive ideal space
of crossed product C*-algebras established in \cite{GoRo--79}.

\begin{prop}
% 9.A.1
\label{Pro-WeakContSemiDirect}
Let $\Gamma$ be a discrete group, 
$(\pi, \Hi)$ a representation of $\Gamma$,
and $N$ an \emph{abelian} normal subgroup of $\Gamma$.
Assume that there exists $\chi_0 \in \widehat N$ 
which is weakly contained in $\pi \vert_N$
and such that $\chi_0^\gamma \ne \chi_0$ 
for every $\gamma \in \Gamma \smallsetminus N$.
\par

Then the induced representation $\Ind_N^\Gamma \chi_0$
is weakly contained in $\pi$.
\end{prop}

\begin{proof} 
Set $\sigma := \Ind_N^\Gamma \chi_0$.
Let $\widetilde \chi_0$ be the trivial extension of $\chi_0$,
equal to $\chi_0$ on $N$ and to $0$ on $\Gamma \smallsetminus N$;
recall from Proposition \ref{diagcoeffinduced}
that the representation $\sigma$ of $\Gamma$ 
has a cyclic vector, say $\xi_\sigma$,
and that $\widetilde \chi_0$ is the corresponding function of positive type:
$$
\widetilde \chi_0 (\gamma) \, = \, 
\langle \sigma(\gamma) \xi_\sigma \mid \xi_\sigma \rangle
\hskip.5cm \text{for all} \hskip.2cm
\gamma \in \Gamma .
$$
(Even it it is not relevant to the present proof,
recall moreover that $\sigma$ is irreducible;
see Corollary \ref{Cor-NormalSubgIrr}.)
\par

Let $\varepsilon > 0$, and let $F$ be a finite subset of $\Gamma$.
By Proposition \ref{oneisenough} it is enough to show that
there exists a function of positive type $\varphi$ associated to $\pi$
such that
$\sup_{\gamma \in F} \vert \varphi(\gamma) - \widetilde \chi_0 (\gamma) \vert
< \varepsilon$.
We write $F_1$ for $F \cap N$ and $F_2$ for $F \cap (\Gamma \smallsetminus N)$,
so that $F = F_1 \sqcup F_2$.
\par

Let $E \,\colon \mathcal B (\widehat N) \twoheadrightarrow \Proj (\Hi)$
be the projection-valued measure
associated to the representation $(\pi \vert_N, \Hi)$ of the abelian group $N$;
we have
$$
\pi(\gamma) \, = \, \int_{\widehat N} \chi (\gamma) dE(\chi)
\hskip.5cm \text{for all} \hskip.2cm
\gamma \in N.
$$
By assumption,
$\chi_0^\gamma \ne \chi_0$ for every $\gamma \in F_2$.
Therefore by continuity of the action $\widehat N \curvearrowleft \Gamma$,
there exists a neighbourhood $U$ of $\chi_0$ in $\widehat N$ such that
$$
U^{\gamma} \cap U \, = \, \emptyset 
\hskip.5cm \text{for all} \hskip.2cm
\gamma \in F_2.
\leqno{(*)}
$$
Upon replacing $U$ by a smaller neighbourhood of $\chi_0$ if necessary,
we can assume that, moreover, we have
$$
\vert \chi(\gamma) - \chi_0(\gamma) \vert \, < \, \varepsilon 
\hskip.5cm \text{for all} \hskip.2cm
\chi \in U 
\hskip.2cm \text{and} \hskip.2cm
\gamma \in F_1. 
\leqno{(**)}
$$
Since $\chi_0$ is weakly contained in $\pi \vert_N$,
we have $E(U) \ne 0$, by Proposition~\ref{Prop-ContAbelian}.
Let $\xi \in \Hi$ be a unit vector in the range of the projection $E(U)$.
We claim that 
$$
\vert \langle \pi(\gamma) \xi \mid \xi \rangle - \widetilde{\chi_0}(\gamma)\vert 
\, < \, \varepsilon 
\hskip.5cm \text{for all} \hskip.2cm 
\gamma \in F.
\leqno{(\sharp)}
$$
\par 

Indeed, we have
$$
\langle \pi(\gamma) \xi \mid \xi \rangle \, = \,
\int_{\widehat N} \chi (\gamma) d \mu_\xi (\chi)
\hskip.5cm \text{for all} \hskip.2cm 
\gamma \in N ,
$$
where $\mu_\xi$ is the probability measure on $\widehat N$ associated to $E$ and $\xi$.
Since
$$
E(B)\xi \, = \, E(B)E(U)\xi \, = \, 0
$$
for every $B \in \mathcal B (\widehat N)$ with $B \cap U = \emptyset$,
the support of $\mu_\xi$ is contained in $U$ and therefore
$$
\langle \pi(\gamma)\xi \mid \xi \rangle \, = \, 
\int_{U}\chi (\gamma)d\mu_\xi(\chi) 
\hskip.5cm \text{for all} \hskip.2cm
\gamma \in N.
$$
It follows from $(**)$ that we have
$$
\begin{aligned}
\vert \langle \pi(\gamma )\xi \mid \xi \rangle -\widetilde{\chi_0}(\gamma) \vert
\, = \, \left\vert \int_{U}(\chi (\gamma)
- \chi_0(\gamma))d\mu_\xi(\chi) \right\vert &
\\ 
\, \le \, \int_{U} \vert \chi (\gamma)-\chi_0(\gamma) \vert d\mu_\xi(\chi) 
\, < \, \varepsilon &
\hskip.8cm \text{for all} \hskip.2cm 
\gamma \in F_1 .
\end{aligned}
\leqno{(\sharp_1)}
$$
\par

Now, let $\gamma \in F_2$. 
We have 
$$
\pi(\gamma)E(B)\pi(\gamma)^{-1} \, = \, E(B^\gamma)
\hskip.5cm \text{for all} \hskip.2cm 
B \in \mathcal B (\widehat N)
$$
(see Proposition~\ref{Prop-RestNormalSub}).
Since $\xi = E(U)\xi$, it follows that
$$
\pi(\gamma)\xi
\, = \, \pi(\gamma) E(U)\xi 
\, = \, \pi(\gamma) E(U) \pi(\gamma)^{-1} \pi(\gamma)\xi 
\, = \, E(U^\gamma) \pi(\gamma) \xi
$$
and hence, using $(*)$,
$$
\langle \pi(\gamma) \xi \mid \xi \rangle
\, = \, \langle E(U^\gamma) \pi(\gamma) \xi \mid E(U)\xi \rangle
\, = \, \langle E(U) E(U^\gamma) \pi(\gamma) \xi \mid \xi \rangle \, = \, 0.
$$
Therefore, we have
$$
\langle \pi(\gamma)\xi \mid \xi \rangle \, = \, \widetilde{\chi_0}(\gamma)
\hskip.5cm \text{for all} \hskip.2cm 
\gamma \in F_2 .
\leqno{(\sharp_2)}
$$
\par

Since ($\sharp_1$) and ($\sharp_2$) imply $(\sharp)$, 
this concludes the proof.
\end{proof}

\begin{rem}
% 9.A.2
\label{Rem-WeakContSemiDirect}
Under the assumption of Proposition~\ref{Pro-WeakContSemiDirect},
the induced representation $\Ind_N^\Gamma \chi_0^\gamma$
is weakly contained in $\pi$ for every $\gamma \in \Gamma$. 
Indeed, $\Ind_N^\Gamma\chi_0^\gamma$
is equivalent to $\Ind_N^\Gamma\chi_0$
(Proposition~\ref{PropConjIndRep}).
\end{rem}

We need at several places the following general result on
CCR representations as introduced in Definition~\ref{Def-CCR-Group}.
\index{CCR! $1$@representation}
\index{Representation! CCR}
% Needed in
% fin d\'emS \ref{Theo-PrimIdealBS} et \ref{Theo-PrimIdealPGLn} (???),
% proof Sof \ref{Theo-ThomaDualTwoStepNilpotent-bis}
% and \ref{Theo-ThomaDualBS},

\begin{prop}
% 9.A.3
% \label{Lem-RepDimFinie}
\label{PropOnWeAndEq}
Let $G$ be a locally compact group
and $(\pi, \Hi)$ an irreducible CCR representation of $G$.
\begin{enumerate}[label=(\arabic*)]
\item\label{iDEPropOnWeAndEq}
Let $\rho$ be a representation of $G$ which is weakly contained in $\pi$.
Then $\rho$ is equivalent to a multiple of $\pi$.
\item\label{iiDEPropOnWeAndEq}
Assume that $G$ is discrete;
let $H$ be a subgroup of $G$ and $\sigma$ a representation of~$H$.
If $\pi$ is contained in the induced representation $\Ind_H^G\sigma$,
then $H$ has finite index in $G$.
\end{enumerate}
\end{prop}

\begin{proof}
\ref{iDEPropOnWeAndEq}
On the one hand,
$\rho$ factorizes through the quotient algebra
\par\noindent
$C^*_{\rm max}(G)/\textnormal{C*ker}(\pi)$, 
since $\rho$ is weakly contained in $\pi$.
\par

On the other hand, 
since $\pi$ is a CCR representation and is irreducible, 
\par\noindent
$C^*_{\rm max}(G) / \textnormal{C*ker}(\pi)$ 
is isomorphic to the algebra $\Ki(\Hi)$ of compact operators on $\Hi$. 
Since every representation of $\Ki(\Hi)$ is equivalent 
to a multiple of the identity representation on $\Hi$, see Proposition~\ref{Pro-FactsAboutCompactOperators},
this implies that $\rho$ is equivalent to a multiple of $\pi$.

\vskip.2cm

\ref{iiDEPropOnWeAndEq}
Since $G$ is discrete and $\pi$ is a CCR representation,
the identity operator $\pi(e)$ is compact.
Therefore $\pi$ is finite-dimensional.
\par

Assume that $\pi$ is contained 
in the induced representation $\lambda := \Ind_H^G\sigma$.
Since $\pi$ is finite-dimensional, $1_G$ is contained in the tensor product 
$\pi \otimes \overline \pi$, 
where $\overline \pi$ is the conjugate representation of $\pi$
\cite[A.1.13]{BeHV--08}. 
Therefore $1_G$ is contained in $\lambda \otimes \overline \lambda$.
However, $\lambda \otimes \overline \lambda$ is equivalent
to the induced representation 
$\Ind_H^G (\sigma \otimes (\lambda \vert_H))$; 
see \cite[E.2.5]{BeHV--08}.
This implies that $H$ has finite index in $G$ \cite[E.3.1]{BeHV--08}.
\end{proof}

Note that, as a particular case of \ref{iDEPropOnWeAndEq},
if $\rho, \pi$ are irreducible representations of $G$ and $\pi$ is finite dimensional,
then $\rho$ is weakly contained in $\pi$ if and only if $\rho$ is contained in $\pi$.
\par

The following corollary is an immediate consequence
of Proposition~\ref{PropOnWeAndEq}~\ref{iDEPropOnWeAndEq}.

\begin{cor}
% 9.A.4
\label{Cor-RepCCR}
Let $G$ be a locally compact group and $(\pi, \Hi)$ an irreducible CCR representation of $G$.
\par

Then the weak equivalence class of $\pi$ in $\widehat G$
consists of $\pi$ itself;
in other words, $\{\pi\}$ is a closed point in $\widehat G$. 
\end{cor}

In particular, if $\widehat G_{\rm fd}$ is the subset of $\widehat G$
of finite-dimensional irreducible representations of $G$,
as in Chapter \ref{Chapter-AllFiniteDimensionalRep},
the canonical projection of $\widehat G$ onto $\Pri (G)$
restricts to an injection $\widehat G_{\rm fd} \to \Pri (G)$.

\section
{Two-step nilpotent groups}
% Section 9.B
\label{Sect-PrimIdealHeisenberg}

Let $\Gamma$ be a two-step nilpotent discrete group, with centre $Z$.
(Recall that, by our definition of two-step nilpotent groups,
$\Gamma$ could also be an abelian group.)
In Sections \ref{Section-IrrRepTwoStepNil}
and \ref{Section-FiniteDimRepForSomeGroups},
we have constructed
families of non-equivalent irreducible representations of $\Gamma$.
In this section, we describe $\Pri(\Gamma)$, using Proposition~\ref{Pro-WeakContSemiDirect}
and some facts established in \ref{Section-IrrRepTwoStepNil}.
\index{Two-step nilpotent group}
\par

For $\psi \in \widehat Z$,
denote as in Section \ref{Section-IrrRepTwoStepNil}
by $\overline \Gamma$ the quotient $\Gamma / \ker \psi$,
by $\overline Z$ its centre,
and by $p \,\colon \Gamma \twoheadrightarrow \overline \Gamma$ the canonical projection.
Set
$$
Z_\psi \, := \, p^{-1}(\overline Z) .
$$
which is a normal subgroup of $\Gamma$ containing $Z$.
Denote by $p'$ the restriction $Z_\psi \twoheadrightarrow \overline Z$ of $p$.
The map $\overline \chi \mapsto \overline \chi \circ p'$
is a homeomorphism from the dual group $\widehat{\overline Z}$ of $\overline Z$
onto the closed subspace $\{ \chi \in \widehat{Z_\psi} \mid \chi \vert_Z = \psi \}$
of $\widehat{Z_\psi}$.
\par

Let $\pi$ be a factor representation of $\Gamma$.
Recall that the projective kernel $\Pker \pi$,
the associated unitary character $\chi_\pi \in \widehat{\Pker \pi}$,
and the central character $\psi_\pi = \chi_\pi \vert_Z \in \widehat Z$ of $\pi$
have been defined in Section \ref{Section-CentralCharacter}.
Since $\ker \pi \cap Z = \ker \psi_\pi$
and since $[\Gamma, \Gamma] \subset Z$, we have
$$
\Pker \pi \, = \, Z_{\psi_\pi},
$$
by Proposition \ref{Prop-ProjKer}~\ref{iiDEProp-ProjKer}.

\begin{theorem}
% 9.B.1
\label{Theo-PrimIdealTwoStepNilpotent}
Let $\Gamma$ be a two-step nilpotent group with centre $Z$.
For every $\psi \in \widehat Z$,
let $Z_\psi$ be the normal subgroup of $\Gamma$ defined as above.
\par

The map
$$
\Phi \, \colon \, 
\begin{cases}
\hskip.2cm
\displaystyle
\bigsqcup_{\psi \in \widehat Z} \left(
\{\psi\} \times \left\{ \chi \in \widehat{Z_\psi} \mid \chi\vert_Z = \psi \right\}
\right)
& \to \hskip.5cm \Pri(\Gamma)
\\
\hskip.2cm \phantom{\{} (\psi, \chi) 
& \mapsto \hskip.2cm \textnormal{C*ker} (\Ind_{Z_\psi}^\Gamma \chi)
\end{cases}
$$
is a bijection. The inverse map $\Phi^{-1}$ is given by 
$$
\Phi^{-1}(\textnormal{C*ker}(\pi)) \, = \, (\psi_\pi, \chi_\pi),
$$
where $\chi_\pi$ is the unitary character of $\Pker \pi$
determined by $\pi \in \widehat \Gamma$,
and $\psi_\pi = \chi_\pi \vert_Z$ is the central character of $\pi$.
\end{theorem}

\begin{proof}
$\bullet$ \emph{First step.}
We first prove that every primitive ideal in $C^*_{\rm max}(\Gamma)$
is in the image of $\Phi$.
Let $\pi$ be an irreducible representation of $\Gamma$,
with central character $\psi = \psi_\pi \in \widehat Z$
and unitary character $\chi = \chi_\pi \in \widehat{\Pker \pi} = \widehat{Z_\psi}$.
We claim that $\pi$ is weakly equivalent to $\Ind_{Z_\psi}^\Gamma \chi$,
that is,
$$
\textnormal{C*ker}\pi
\, = \,
\textnormal{C*ker} \left( \Ind_{Z_\psi}^\Gamma \chi \right) .
$$
The claim implies that every primitive ideal in $C^*_{\rm max}(\Gamma)$
arises as $\Phi(\psi, \chi)$ 
for a pair $(\psi, \chi)$ in the domain of $\Phi$.

\vskip.2cm

$\circ$ \emph{First case.} 
Assume first that $\psi$ is faithful. Then $Z_\psi = Z$
and we have to show that $\pi$ is weakly equivalent to $\Ind_{Z}^\Gamma \psi$.
\par

Let $N$ be a maximal abelian subgroup of $\Gamma$. 
%as in the proof of Theorem~\ref{Theo-PrimIdealAffGr},
Let $\chi_0 \in \widehat N$ be such that 
$\chi_0$ is weakly contained in $\pi \vert_N$.
By continuity of the restriction process,
$\chi_0 \vert_Z$ is weakly contained in $\pi \vert_Z$
and hence $\chi_0 \vert_Z = \psi$,
that is,
$$
\begin{aligned}
\chi_0 \in \widehat N (\psi)
\, &= \,
\{\chi \in \widehat N \mid \text{the restriction of $\chi$ to $Z$ is $\psi$} \}
\\
\, &= \,
\left\{\chi \in \widehat N \hskip.1cm \big\vert \hskip.1cm
\chi = \chi_0 \rho
\hskip.2cm \text{for some} \hskip.2cm
\rho \in \widehat{N/Z} \right\}
\end{aligned}
$$
($\widehat N (\psi)$ is here as in Lemma \ref{Lem-FreeActDualHeis}).
\par

By Lemma~\ref{Lem-FreeActDualHeis}~\ref{iiiDELem-FreeActDualHeis}, 
the group $\Gamma/N$ acts freely on $\widehat N (\psi)$.
By Proposition~\ref{Pro-WeakContSemiDirect}
and Remark~\ref{Rem-WeakContSemiDirect},
it follows that $\Ind_N^\Gamma \chi_0^\gamma$ is weakly contained in $\pi$
for every $\gamma \in \Gamma$.
By Lemma~\ref{Lem-FreeActDualHeis}~\ref{ivDELem-FreeActDualHeis}, 
the set $\{\chi_0^\gamma \mid \gamma \in \Gamma \}$
is dense in $\widehat N (\psi)$.
It follows that $\Ind_N^\Gamma \chi$ 
is weakly contained in $\pi$ for every $\chi \in \widehat N (\psi)$,
and therefore that $\bigoplus_{\chi \in \widehat N (\psi)} \Ind_N^\Gamma \chi$
is weakly contained in $\pi$.
\par

On the one hand, by Proposition
\ref{Prop-WeakContainmentInduced}~\ref{iiDEProp-WeakContainmentInduced},
we have a weak equivalence
$$
\Ind_Z^N \psi \, \sim \,
\bigoplus_{\rho \in \widehat{N/Z}} \chi_0 \rho \, = \, 
\bigoplus_{\chi \in \widehat N_\psi} \chi .
$$
By continuity of the induction process (see \cite[F.3.5]{BeHV--08}),
it follows that the representation
$
\Ind_Z^\Gamma \psi \simeq \Ind_N^\Gamma( \Ind_Z^N \psi)
$
is weakly equivalent to
$$
\Ind_N^\Gamma \bigg( \bigoplus_{\chi \in \widehat N_\psi} \chi \bigg)
\, \simeq \,
\bigoplus_{\chi \in \widehat N_\psi} \Ind_N^\Gamma \chi,
$$
and therefore that $\Ind_Z^\Gamma \psi$ is weakly contained in~$\pi$. 
\par

On the other hand, as $\Gamma/Z$ is amenable, 
$\pi$ is weakly contained in $\Ind_Z^\Gamma(\pi \vert_Z)$, by Proposition
\ref{Prop-WeakContainmentInduced}~\ref{iDEProp-WeakContainmentInduced}.
Since $\Ind_Z^\Gamma (\pi \vert_Z)$ is weakly equivalent to $\Ind_Z^\Gamma \psi$,
it follows that $\pi$ is weakly contained in $\Ind_Z^\Gamma \psi$. 
\par

Therefore, $\pi$ is weakly equivalent to $\Ind_Z^\Gamma \psi$, 
and the claim is proved in case $\psi$ faithful.

\vskip.2cm

$\circ$ \emph{Second case.} 
Assume now that $\psi$ is not faithful. We have 
$$
\Pker \pi \, = \, Z_\psi \, = \, p^{-1}(\overline Z),
$$
by Proposition~\ref{Prop-ProjKer}. 
\par

The representation $\pi$ factorizes
by a representation $\overline \pi$ of $\overline \Gamma$
such that $\pi = \overline \pi \circ p$.
The unitary character $\chi$ of $\Pker \pi$
factorizes by a faithful unitary character $\overline \chi$ of $\overline Z$,
 which is also the central character of $\overline \pi$.
\par

By what was proved in the first case, 
$\overline \pi$ is weakly equivalent to 
$\Ind_{\overline Z}^{\overline \Gamma} \overline \chi$.
Since 
$$
\Ind_{Z_\psi}^\Gamma \chi \, = \,
\left( \Ind_{\overline Z}^{\overline \Gamma} \overline \chi \right) \circ p,
$$
it follows that $\pi = \overline \pi \circ p$
is weakly equivalent to $\Ind_{Z_\psi}^\Gamma \chi$.

\vskip.2cm 

$\bullet$ \emph{Second step.} 
We claim that every ideal $\textnormal{C*ker} (\Ind_{Z_\psi}^\Gamma \chi)$
in the image $\Phi$ is a primitive ideal of $C^*_{\rm max}(\Gamma)$.
More precisely, given $\psi \in \widehat Z$
and $\chi \in \widehat Z_\psi$ extending $\chi$,
we claim that there exists $\pi \in \widehat \Gamma$ such that 
$$
\textnormal{C*ker} (\Ind_{Z_\psi}^\Gamma \chi)
\, = \, \textnormal{C*ker}(\pi)
$$
and that, moreover, $\psi$ coincides with the central character $\psi_\pi$ of $\pi$
and $\chi$ with the unitary character $\chi_\pi$
of $\widehat{\Pker \pi} = \widehat{Z_\psi}$.
\par

The unitary character $\chi$ factorizes by 
a unitary character of $\overline Z$, say $\overline \chi$.
Let $\overline \pi$ be an irreducible representation of $\overline \Gamma$ 
which is weakly contained in $\Ind_{\overline Z}^{\overline \Gamma} \overline \chi$.
Then $\overline \pi$ has $\overline \chi$ as central character,
since the restriction of $\Ind_{\overline Z}^{\overline \Gamma}\overline \chi$
to $\overline Z$ is a multiple of $\overline \chi$
(see Proposition~\ref{PropConjIndRep}~\ref{iiDEPropConjIndRep}).
Therefore, by the first step,
$$
\textnormal{C*ker}(\overline \pi)
\, = \, 
\textnormal{C*ker}(\Ind_{\overline Z}^{\overline \Gamma} \overline \chi);
$$
lifting $\overline \pi$ to a representation $\pi$ of $\Gamma$, 
we have
$$
\textnormal{C*ker}(\pi) \, = \,
\textnormal{C*ker}(\Ind_{Z_\psi}^\Gamma \chi).
$$
Observe that we also have $\chi_\pi = \chi$ and $\psi_\pi = \psi$, as claimed.

\vskip.2cm

$\bullet$ \emph{Third step.} 
We claim that the map $\Phi$ is injective.
\par

Indeed, for $i = 1, 2$, let $\psi_i \in \widehat Z$
and $\chi_i \in \widehat Z_{\psi_i}$ extending $\psi_i$.
By the second step, there exists $\pi_i \in \widehat \Gamma$
with the following properties:
$\Phi(\psi_i, \chi_i) = \textnormal{C*ker}(\pi_i)$,
the central character of $\pi_i$ is $\psi_i$,
and the unitary character of $\Pker \pi_i=Z_{\psi_i}$ is $\chi_i$.
\par

Assume now that $\Phi(\psi_1, \chi_1) = \Phi(\psi_2, \chi_2)$.
Then $\pi_1$ and $\pi_2$ are weakly equivalent.
Therefore $\Pker \pi_1 = \Pker \pi_2$ and $\chi_1 = \chi_2$,
by Proposition \ref{Prop-ProjKer}~\ref{iiiDEProp-ProjKer}.
It is then obvious that $\psi_1 = \psi_2$.
\end{proof}

\begin{rem}
% 9.B.2
\label{Rem-PrimIdealTwoStepNil}
(1)
The description in Theorem~\ref{Theo-PrimIdealTwoStepNilpotent}
of $\Pri(\Gamma)$ for a two-step nilpotent group $\Gamma$
appears in \cite[Lemma 2]{Kani--82},
which follows ideas from \cite[Proposition 5]{Howe--77b}.
The work \cite{Howe--77b} provides a description of $\Pri(\Gamma)$
for $\Gamma$ a finitely generated, torsion-free, discrete nilpotent group.
For an extension of this description to a large class of nilpotent LC groups, 
see \cite{EcKl--12}.

\vskip.2cm

(2)
As we will see later (Theorem~\ref{Theo-ThomaDualTwoStepNilpotent-bis}),
every representation $\Ind_{Z_\psi}^\Gamma \chi$ of $\Gamma$
appearing in Theorem~\ref{Theo-PrimIdealTwoStepNilpotent}
is a \textbf{factor} representation of finite type.
So, except in some special cases, 
$\Ind_{Z_\psi}^\Gamma \chi$ will not be irreducible
(see Corollaries~\ref{Cor-ThomaDualHeisField-Bis} 
and \ref{Cor-ThomaDualHeisIntegers-Bis}
concerning of the case of Heisenberg groups).
\end{rem}

We are going apply Theorem~\ref{Theo-PrimIdealTwoStepNilpotent}
to Heisenberg groups.

\subsection*{Heisenberg group over a ring}

Let $R$ be a unital commutative ring.
As in Section \ref{Section-IrrRepTwoStepNil},
consider the Heisenberg group $H(R)$, with centre $Z \approx R$,
Recall that, for $\psi \in \widehat Z$, we have
\begin{enumerate}
\item[--]
an ideal $I_\psi \, = \, \{a \in R \mid aR \subset \ker \psi \}$ of $R$;
\item[--]
the normal subgroup $Z_\psi$ of $H(R)$ defined earlier in this section
\item[]
which is also
$Z_\psi \, = \, \{ (a,b,c) \in \Gamma \mid a, b \in I_\psi, \hskip.1cm c \in R \}$
by Lemma~\ref{Lem-IrredHeisRing};
\item[--]
for every unitary character $\chi$ of $Z_\psi$ with $\chi\vert_Z = \psi$,
\item[]
there is a pair $(a, b) \in I_\psi \times I_\psi$ such that $\chi$ is equal to
$$
\chi_{\psi, \alpha, \beta} \, \colon \, Z_\psi \to \T
\hskip.2cm
(a,b,c) \mapsto \alpha(a)\beta(b)\psi(c) .
$$
\end{enumerate}
The following result is a direct consequence
of Theorem~\ref{Theo-PrimIdealTwoStepNilpotent}.

\begin{cor}
% 9.B.3
\label{Cor-PrimIdealHeisRing}
Let $\Gamma = H(R)$ for a unital commutative ring $R$.
For every $\psi \in \widehat Z$,
let $I_\psi$ and $Z_\psi$ be as above.
\par

The map
$$
\Phi \, \colon \, 
\begin{cases}
\hskip.2cm
\displaystyle
\bigsqcup_{\psi \in \widehat Z} \{\psi\} \times \widehat{I_\psi}^2
& \to
\hskip.5cm
\Pri(\Gamma)
\\
\hskip.2cm \phantom{\{} (\psi,(\alpha, \beta))
& \mapsto \hskip.2cm
\textnormal{C*ker} (\Ind_{Z_\psi}^\Gamma \chi_{\psi, \alpha, \beta}),
\end{cases}
$$
is a bijection
\end{cor}

\index{Heisenberg group! $2$@$H(\K)$ over an infinite field $\K$}
We consider the particular case $\Gamma = H(\K)$, where $R = \K$ is a field.
Recall that $Z \approx \K$, and $\Gamma/Z \approx \K^2$.
\par

For $\psi \in \widehat{\K} \smallsetminus \{1_\K\}$,
we have $I_\psi = \{0\}$ and hence $Z_\psi = Z$;
for $\psi =1_\K$, we have $I_\psi = \K$ and hence $Z_\psi= H(\K)$.
\par

We obtain therefore the following consequence of Corollary~\ref{Cor-PrimIdealHeisRing}.

\begin{cor}
% 9.B.4
\label{Cor-PrimIdealHeisField}
Let $\K$ be a field, $\Gamma = H(\K)$ the Heisenberg group over $\K$,
$Z$~the centre of $H(\K)$,
and $p \,\colon \Gamma \twoheadrightarrow \Gamma/Z$ the canonical epimorphism.
\par

The map
$$
\Phi \, \colon \,
\big( \widehat Z \smallsetminus \{1_Z\} \big) \sqcup \widehat{\Gamma/Z} 
\, \to \,
\Pri(\Gamma)
$$
defined by
$$
\begin{aligned}
\Phi(\psi) \, &= \, \textnormal{C*ker}(\Ind_Z^\Gamma \psi)
& \hskip.5cm \text{for} \hskip.3cm
\psi
&\in \widehat Z \smallsetminus \{1_Z\} \approx \widehat K \smallsetminus \{1_K\} 
\\
\Phi(\chi) \, &= \, \textnormal{C*ker}(\chi \circ p) 
& \hskip.5cm \text{for} \hskip.3cm
\chi
&\in \widehat{\Gamma/Z} \approx \widehat K^2,
\end{aligned}
$$
is a bijection.
\end{cor}

We will see later (Corollary~\ref{Cor-ThomaDualHeisField-Bis})
that $\Ind_Z^\Gamma \psi$ is a \textbf{factor} representation of $\Gamma=H(\K)$,
which moreover is \textbf{not} irreducible. 
The following version of Corollary~\ref{Cor-PrimIdealHeisField}
provides an explicit parametrization of $\Pri(\Gamma)$
in terms of weak equivalence classes of irreducible representations.
\par

Recall from Proposition~\ref{Prop-ProjKer} that the central character
of an irreducible representation only depends on its weak equivalence class.

\begin{cor}
% 9.B.5
\label{Cor-PrimIdealHeisField-bis} 
Let $\K$ be a field, $\Gamma = H(\K)$ the Heisenberg group over $\K$,
and $Z$ the centre of $\Gamma$.
Let $N$ the maximal abelian subgroup of $\Gamma$,
defined by $N = (0, b, c) \mid b, c \in \K\}$.
For every $\psi \in \widehat Z\smallsetminus \{1_Z\}$,
let $\chi_\psi$ be the unitary character of $N$
defined by $\chi_\psi(0,b,c) = \psi(c)$ for $b, c \in \K$. 
\begin{enumerate}[label=(\arabic*)]
\item\label{iDECor-PrimIdealHeisField-bis}
For $\psi \in \widehat Z\smallsetminus \{1_Z\}$,
the representation $\Ind_{N}^\Gamma \chi_\psi$ is irreducible;
its weak equivalence class
consists of all irreducible representations of $\Gamma$
with central character~$\psi$.
\item\label{iiDECor-PrimIdealHeisField-bis}
For $\chi \in \widehat{\Gamma/Z}$,
the weak equivalence class of $\chi \circ p$ consists of $\chi \circ p$ itself. 
\end{enumerate}
\end{cor}

\begin{proof}
By Corollary~\ref{Cor-IrredHeisField}),
for every $\psi \in \widehat Z\smallsetminus \{1_Z\}$,
the representation $\Ind_{N}^\Gamma \chi_\psi$ is irreducible
and has $\psi$ as central character.
So, the claim follows from Corollary~\ref{Cor-PrimIdealHeisField}.
\end{proof}

\begin{rem}
% 9.B.6
\label{Rem-Cor-PrimIdealHeisField-bis}
Corollary~\ref{Cor-PrimIdealHeisField-bis} shows that the map 
$$
\Pri(H(\K)) \, \to \, \widehat \Z \approx \widehat \K,
$$
that associates to a weak equivalence class its central character
(see Section~\ref{Section-CentralCharacter}) is a surjective map
with the property that the inverse image of 
an element in $\widehat \K \smallsetminus \{1_\K\}$ is one point 
and the inverse image of $\{1_\K\}$ is a copy of $\widehat{\K}^2$.
\end{rem}

Next, we consider the case where $R$ is the ring $\Z$ of rational integers.

\vskip.2cm

Recall that 
$$
\widehat Z \approx \widehat \Z \, = \,
\{\psi_\theta \mid \theta \in \mathopen[ 0,1 \mathclose[ \},
$$
with $\psi_\theta(c) = e^{2 \pi i \theta c}$ for all $c \in \Z$.
We denote by $\widehat Z_\infty$ the elements of infinite order 
and by $\widehat Z_n$ the elements of order $n \ge 1$ in $\widehat Z$.
\par

Let $\psi = \psi_\theta \in \widehat Z$. Two cases can occur:
\begin{enumerate}
\item[$\bullet$]
$\psi \in \widehat Z_\infty$;
then $I_\psi = \{0\}$ and hence $Z_\psi = Z$.
\item[$\bullet$]
$\psi \in \widehat Z_n$ for some $n \ge 1$,
that is, $\theta = p/n$ for an integer $p \in \{0, 1, \hdots, n-1 \}$, 
with $p$ and $n$ coprime;
then $I_\psi = n\Z$ and hence $Z_\psi$ coincides with the normal subgroup
$$
\Lambda(n) \, = \,
\{(a, b, c) \in \Gamma \mid a \in n\Z, \hskip.1cm b \in n\Z, \hskip.1cm c \in \Z\},
$$
which only depends on $n$.
\end{enumerate}
\par

Fix $\psi = \psi_\theta \in \widehat{Z}_n$.
Every unitary character of $\Lambda(n)$ which coincides with $\psi$ on $Z$ 
is of the form $\chi_{\psi, \alpha, \beta}$,
where 
$$
\chi_{\psi, \alpha, \beta}(a, b, c) \, = \, 
e^{2 \pi i (\alpha a + \beta b)} \psi (c) 
\hskip.5cm \text{for} \hskip.2cm
a, b \in n\Z, \hskip.1cm c \in \Z
$$
for a uniquely determined pair $(\alpha, \beta) \in \mathopen[ 0,1/n \mathclose[^2$. 
\par

In view of these remarks, the following result is an immediate consequence of Corollary~\ref{Cor-PrimIdealHeisRing}.

\begin{cor}
% 9.B.7
\label{Cor-PrimIdealHeisIntegers}
\index{Heisenberg group! $3$@$H(\Z)$}
Let $\Gamma = H(\Z)$ be the Heisenberg group over the integers
and $Z$ its centre.
The map
$$
\Phi \, \colon \, 
\widehat Z_\infty \hskip.1cm \bigsqcup \hskip.1cm
\Big( \bigsqcup_{n \ge 1}
\bigsqcup_{\psi \in \widehat Z_n}
\{\psi\} \times \mathopen[ 0,1/n \mathclose[^2 \Big) 
\, \to \, \Pri(\Gamma)
$$
defined by
$$
\begin{aligned}
\Phi(\psi) \, &= \,
\textnormal{C*ker} (\Ind_Z^\Gamma \psi) 
& \hskip.5cm \text{for} \hskip.3cm
&\psi \, \in \, \widehat Z_\infty \hskip2cm
\\
\Phi(\psi, \alpha, \beta) \, &= \,
\textnormal{C*ker} (\Ind_{\Lambda(n)}^\Gamma \chi_{\psi, \alpha, \beta})
& \hskip.5cm \text{for} \hskip.3cm
&(\psi, \alpha, \beta) \, \in \, \widehat Z_n \times \mathopen[ 0,1/n \mathclose[^2
\end{aligned}
$$
is a bijection.
\end{cor}

As in the case of the Heisenberg group over a field,
we give an explicit parametrization of $\Pri(\Gamma)$
in terms of weak equivalence classes of irreducible representations.
\par

For $\psi \in \widehat Z_\infty$,
denote as before by $\chi_\psi$ the unitary character
of the maximal abelian subgroup $N = \{(0, b, c) \mid b, c \in \Z\}$
defined by $\chi_\psi(0,b,c) = \psi(c)$ for $b,c \in \Z$.
For an integer $n \ge 1$, denote by $\Gamma(n)$ the normal subgroup
$$
\Gamma(n) \, = \,
\{(a, b, c) \in \Gamma \mid a \in n\Z, \hskip.1cm b, c \in \Z \}
$$
of $\Gamma$ (which appeared already
in Corollaries \ref{Cor-IrredFiniHeisIntegers} and \ref{Cor-FinDimRepHeis-Integers}).
For $\psi \in \widehat Z_n$ and $(\alpha, \beta) \in \mathopen[ 0,1/n \mathclose[^2$, 
the previously defined unitary character 
$\chi_{\psi, \alpha, \beta}$ of $\Lambda(n)$
extends to a unitary character of $\Gamma(n)$,
again denoted by $\chi_{\psi, \alpha, \beta}$ and defined by the same formula:
$$
\chi_{\psi, \alpha, \beta}(a,b,c) \, = \, 
e^{2 \pi i (\alpha a+ \beta b)} \psi (c) 
\hskip.5cm \text{for} \hskip.2cm
a \in n\Z, \hskip.1cm b \in \Z, \hskip.1cm c \in \Z.
$$

\begin{cor}
% 9.B.8
\label{Cor-PrimIdealHeisIntegers-bis}
Let $\Gamma = H(\Z)$ and $Z$ its centre. 
\begin{enumerate}[label=(\arabic*)]
\item\label{iDECor-PrimIdealHeisIntegers-bis}
For $\psi \in \widehat Z_\infty$,
the representation $\Ind_{N}^\Gamma \chi_\psi$ is irreducible;
its weak equivalence class
consists of all irreducible representations of $\Gamma$
with central character $\psi$.
\item\label{iiDECor-PrimIdealHeisIntegers-bis}
For every $n \ge 1$
and $(\psi, \alpha, \beta) \in Z_n \times \mathopen[ 0,1/n \mathclose[^2$,
the representation $\pi_{\psi, \alpha, \beta} =
\Ind_{\Gamma(n)}^\Gamma \chi_{\psi, \alpha, \beta}$
is irreducible.
Moreover, every irreducible representation $\pi$ of $\Gamma$
with central character $\psi \in Z_n$
is weakly equivalent to $\pi_{\psi, \alpha, \beta}$
for the unique pair $(\alpha, \beta) \in \mathopen[ 0,1/n \mathclose[^2$
such that $\pi \vert_{\Lambda(n)}$ is a multiple of $\chi_{\psi, \alpha, \beta}$.
\end{enumerate}
\end{cor}

\begin{proof}
This is a direct consequence of Corollaries
\ref{Cor-IrredRepHeisInteger}, \ref{Cor-IrredFiniHeisIntegers},
and \ref{Cor-PrimIdealHeisIntegers}.
\end{proof}

\begin{rem}
% 9.B.9
\label{Rem-Cor-PrimIdealHeisIntegers-bis}

Corollary~\ref{Cor-PrimIdealHeisIntegers-bis} shows that the map 
$$
\Pri(H(\Z)) \, \to \, \widehat \Z \approx \T,
$$
that associates to a weak equivalence class its central character
is a surjective map with the property that
the inverse image of an element of infinite order in $\T$ is one point 
and the inverse image of an element of finite order is a $2$-torus.
% \cite{Howe--77}; see also \cite[VII.5.2]{Davi--96}.
Compare with Example \ref{ExRepHeisZ}.
% Convenu qu'on ne cite pas Kawakami \cite{Kawa--82}.
\par

Let $\psi \in \widehat Z$ be of infinite order,
say $\psi (c) = e^{2 \pi i \theta c}$ for all $c \in Z$
and some $\theta \in \mathopen[ 0,1 \mathclose[$, $\theta \notin \Q$.
Denote by $\widehat \Gamma (\psi)$
the subset of the dual $\widehat \Gamma$
of classes of irreducible representations with central character $\psi$. 
\par

Fix $\psi = \psi_\theta \in \widehat Z$ with infinite order. 
On the one hand, all irreducible representations in $\widehat \Gamma (\psi_\theta)$ 
correspond to the same primitive ideal of $\Gamma$,
that is, the image of $\widehat \Gamma (\psi_\theta)$ under the canonical map
$\kappa^{\rm d}_{\rm prim} \,\colon \widehat \Gamma \twoheadrightarrow \Pri(\Gamma)$ 
consists of a unique point
(see Corollary~\ref{Cor-PrimIdealHeisIntegers}).
On the other hand, $\widehat \Gamma (\psi_\theta)$
contains uncountably many elements 
(see Corollary~\ref{Cor-IrredRepHeisInteger}
and Remark \ref{Rem-Cor-IrredHeisIntegers}).
\par

For every $\pi \in \widehat \Gamma (\psi_\theta)$ 
the C*-subalgebra $\pi(C^*_{\rm max}(\Gamma))$ of $\Li (\Hi_\pi)$
is isomorphic to the C*-algebra $A_\theta$
of the irrational rotation of angle $2 \pi \theta$ \cite{Rief--81}.
Therefore, $\widehat \Gamma (\psi_\theta)$ is in natural bijection
with the spectrum $\widehat{A_\theta}$ of $A_\theta$.
Such a C*-algebra has uncountably many pairwise non-equivalent
irreducible representations with the same kernel,
and uncountably many pairwise non-quasi-equivalent.
For more on the representation theory of $A_\theta$, see \cite{O'Do--81}.
\end{rem}

\section
{Affine groups of infinite fields}
% Section 9.C
\label{Sect-PrimIdealAffGr}

\index{Affine group! $4$@$\Aff(\K)$ of an infinite field $\K$}
Let $\K$ be an \emph{infinite} field,
$\Aff(\K)$ the group of affine transformations of $\K$
and $N$ its derived group,
as in Section \ref{Section-IrrRepAff}.

\begin{theorem}
% 9.C.1
\label{Theo-PrimIdealAffGr}
For $\Gamma = \Aff(\K)$ and $N \approx \K$ as above, 
we have
$$
\Pri(\Gamma) \, = \, \{0\} \sqcup 
\left\{\textnormal{C*ker}(\chi \circ p) \mid \chi \in \widehat{\Gamma/N} \right\}
\, \approx \, \{0\} \sqcup \widehat{\Gamma/N}
$$
where $p$ is the quotient map $\Gamma \to \Gamma/N$.
\end{theorem}

\begin{proof}
Let $\pi$ be an irreducible representation of $\Gamma$.
\par

Assume first that $\pi$ is trivial on $N$.
Then $\pi = \chi \circ p$ for some unitary character $\chi \in \widehat{\Gamma/N}$ 
and hence $\textnormal{C*ker}(\pi) = \textnormal{C*ker}(\chi \circ p)$. 
\par

Assume now that $\pi$ is not trivial on $N$.
We claim that $\textnormal{C*ker}(\pi) = \{0\}$.
\par

Indeed, in this case, there exists $\chi_0 \in \widehat N$ with $\chi_0 \ne 1_N$ 
which is weakly contained in $\pi \vert_N$. 
\par

Since $\Gamma/N$ acts freely on $\widehat N \smallsetminus \{1_N\}$
(Lemma~\ref{Lem-FreeActDualAff}),
$\Ind_N^\Gamma \chi_0^\gamma$ is weakly contained in $\pi$ 
for every $\gamma \in \Gamma$
(Proposition~\ref{Pro-WeakContSemiDirect} and Remark~\ref{Rem-WeakContSemiDirect}).
As $\{\chi_0^\gamma \mid \gamma \in \Gamma\}$ is dense in $\widehat N$ 
(Lemma~\ref{Lem-FreeActDualAff}),
it follows that $\Ind_N^\Gamma \chi$ is weakly contained in $\pi$ 
for every $\chi \in \widehat N$.
\par

On the one hand, the direct sum $\bigoplus_{ \chi \in \widehat N}\chi$
is weakly equivalent to the regular representation $\lambda_N$ of $N$, 
by amenability of $N$. Using continuity of induction,
it follows that
the regular representation $\lambda_\Gamma = \Ind_N^\Gamma \lambda_N$
of $\Gamma$ is weakly contained in $\pi$.
\par

On the other hand, since $\Gamma$ is amenable, 
$\pi$ is weakly contained in $\lambda_\Gamma$. 
Therefore $\pi$ is weakly equivalent to $\lambda_\Gamma$, that is, 
$$
\textnormal{C*ker}(\pi) \, = \, \textnormal{C*ker}(\lambda_\Gamma) \, = \, \{0\},
$$
as claimed.

\vskip.2cm

Conversely, observe that $\textnormal{C*ker}(\chi \circ p) \in \Pri(\Gamma)$ 
for every $\chi \in \widehat{\Gamma/N}$.
Moreover, there exist irreducible representations of $\Gamma$ 
which are non-trivial on $N$; let $\pi$ be one of these.
By what was shown above, 
$\textnormal{C*ker}(\pi) = \{0\}$,
and hence $\{0\} \in \Pri(\Gamma)$.
\end{proof}

\begin{rem}
% 9.C.2
\label{Rem-Theo-PrimIdealAffGr}
Let $\Gamma = \Aff(\K)$ and $N \approx \K$ be as above. 
As in the case of the Heisenberg groups,
we can give an explicit parametrization of $\Pri(\Gamma)$
in terms of weak equivalence classes of irreducible representations.
Indeed, choose any $\psi \in \widehat \K \smallsetminus \{1_\K\}$
and recall from Theorem \ref{Prop-IrredRepAffine}
that $\pi_{\psi} = \Ind_{N}^\Gamma \psi$ is irreducible.
\par

For every $\chi \in \widehat{\Gamma/N}$, the weak equivalence class
of the unitary character $\chi \circ p$
contains only $\chi \circ p$ itself.
By contrast, the weak equivalence class of $\pi_{\psi}$
contains all other irreducible representations of $\Gamma$,
which all have C*-kernel $\{0\}$.
\end{rem}

\section
{Solvable Baumslag--Solitar groups}
% Section 9.D
\label{Section-PrimIdealBS}

Let $p$ be a prime.
Recall that from Section~\ref{Section-IrrRepBS}
that the Baumslag--Solitar $\Gamma = \BS(1, p)$ is the semi-direct product 
$\Gamma = A \ltimes N \approx \Z \ltimes \Z[1/p]$, where
$$
A \, = \, \left\{ \begin{pmatrix} p^k & 0 \\ 0 & 1 \end{pmatrix} 
\mid k \in \Z \right\}\, \approx \, \Z 
\hskip.5cm \text{and} \hskip.5cm
N \, = \, \left\{ \begin{pmatrix} 1& b \\ 0 & 1 \end{pmatrix} 
\mid b \in \Z[1/p] \right\} \, \approx \, \Z[1/p] ,
$$
and the generator $\begin{pmatrix} p & 0 \\ 0 & 1 \end{pmatrix}$ of $A$
acts on $N$ by multiplication by $p$.
We often denote the elements of $\Gamma$ 
simply by $(k, b)$, with $k \in \Z$ and $b \in \Z[1/p]$. 
\index{Baumslag--Solitar group $\BS(1, p)$}
\par

We proceed to describe $\Pri(\Gamma)$. 
The result (Theorem~\ref{Theo-PrimIdealBS})
can be derived from \cite{Guic--65} where a more general situation is considered
and where the case $p=2$ is treated with details; 
our exposition is in many respects different from the one in \cite{Guic--65}.
\par

Recall from Section \ref{Section-IrrRepBS}
that we can identify $\widehat N$ with the $p$-adic solenoid 
$$
\So_p \, = \, (\Q_p \times \R) / \Delta, 
\hskip.5cm \text{where} \hskip.2cm
\Delta \, = \, \{(a,-a) \in \Q_p \times \R \mid a \in \Z[1/p]\} ,
$$
and that we write $[(x,y)]$ for the class in $\So_p$
of $(x,y)$ in $\Q_p \times \R$.
The action of $A$ on $\widehat N$ corresponds to 
the action of $\Z$ on $\So_p$ for which the generator $1 \in \Z$
acts by the map
$$
T_p \, \colon \, \So_p \to \So_p , \hskip.2cm
[(x,y)] \mapsto [px,py]
$$
induced by the multiplication by $p$.
We denote by $\chi_s$ the element in $\widehat N$ 
corresponding to $s \in \So_p$.
For $n \in \N^*$, let $\So_p(n)$
be the set of elements of $\So_p$ with period $n$. Let
$$
\Per(T_p) \, = \, \bigsqcup_{n \ge 1}\So_p(n),
$$
denote the set of periodic elements, and
$$
\So_p(\infty) \, := \, \So_p \smallsetminus \Per(T_p)
$$
the set of elements with an infinite $T_p$-orbit.
\par

We have determined in Corollary~\ref{Cor-FinDimRepBS}
the finite-dimensional representations of $\Gamma$
which are associated to $T_p$-periodic elements.
More precisely,
given $n \ge 1, s \in \So_p(n)$, and $\theta \in \mathopen[0,1/n \mathclose[$,
the induced representation 
$\pi_{s, \theta} = \Ind_{\Gamma(n)}^\Gamma \chi_{s, \theta}$
is an irreducible finite-dimensional representation, 
where $\chi_{s, \theta}$ is a unitary character
of the finite index normal subgroup 
$\Gamma(n) = \left\{(k, b) \in \Gamma \mid k \in n\Z, \hskip.1cm b \in \Z[1/p] \right\}$, 
with $\chi_{s, \theta} = \chi_s$ on $N$.
\par

As in Corollary~\ref{Cor-FinDimRepBS}, we set
$$
X_{\rm fd} \, = \, 
\bigsqcup_{n \ge 1} \hskip.1cm
\big( \So_p(n) \times \mathopen[0,1/n \mathclose[ \big) ,
$$
and $\widetilde T_p \,\colon X_{\rm fd} \to X_{\rm fd}$
is defined by $\widetilde T_p (s, \theta) = (T_ps, \theta)$.
We will need to consider the quasi-orbit space
%(see Subsection~\ref{SNAG-NonAbelian}) 
of the $T_p$-action on $X_{\rm fd} \sqcup \So_p(\infty)$.
\par

Given an action of a group $G$ on a topological $X$,
the corresponding \textbf{quasi-orbit space}
is the quotient space $X/\sim$ for the equivalence relation $\sim$
defined on $X$ by
$$
x \, \sim \, x'
\hskip.5cm \text{if} \hskip.5cm
\overline{Gx} \, = \, \overline{Gx'} ,
$$
or equivalently the quotient of the space of orbits $X/G$
by the equivalence relation defined on $X/G$ by $Gx = G'$ if $\overline{Gx} = \overline{Gx'}$.
See Appendix \ref{AppTop}.

\begin{theorem}
% 9.D.1
\label{Theo-PrimIdealBS}
Let $\Gamma = \BS(1, p) = A \ltimes N$ be the Baumslag--Solitar group;
let $X_{fd}$ and $\So_p(\infty)$ be as above.
The map
$$
\Phi \, \colon \, X_{\rm fd} \sqcup \So_p(\infty) \to \Pri(\Gamma),
$$
defined by
$$
\Phi(s) \, = \, 
\textnormal{C*ker}( \Ind_{\Gamma(n)}^\Gamma\chi_{s, \theta})
\hskip.5cm \text{for} \hskip.2cm 
n \ge 1, \hskip.1cm (s, \theta) \in \So_p(n) \times \mathopen[0,1/n \mathclose[ 
$$
and
$$
\Phi(s) \, = \, 
\textnormal{C*ker}(\Ind_N^\Gamma \chi_s) 
\hskip.5cm \text{for} \hskip.2cm 
s \in \So_p(\infty) ,
$$
factorizes to a bijection between
the quasi-orbit space of the $\widetilde T_p$-action 
on $X_{\rm fd} \sqcup \So_p(\infty)$
and $\Pri(\Gamma)$.
\end{theorem}

\begin{proof}
$\bullet$ \emph{First step.} 
We first prove that every primitive ideal is contained in the range of $\Phi$.
\par

Let $\pi$ be an irreducible representation of $\Gamma$.
Let $E \,\colon \mathcal B (\widehat N) \to \Proj (\Hi)$ be the projection-valued measure
associated to the restriction $\pi \vert_N$ of $\pi$ to $N$
and let $S$ be the support of $E$.
\par

By Proposition~\ref{Prop-RestNormalSub}~\ref{iiDEProp-RestNormalSub}, 
there exists $s \in \So_p$ such that $S$ is the closure
of the $\Gamma$-orbit $\{\chi_s^\gamma \mid \gamma \in A\}$. 
Two cases may occur: $s \in \So_p(\infty)$ or $s \in \Per(T_p)$.
\par

Assume first that $s \in \So_p(\infty)$. We claim that 
$$
\textnormal{C*ker}\pi \, = \, \textnormal{C*ker}(\Ind_N^\Gamma \chi_s) .
$$
\par

Indeed, since $A$ acts freely on the orbit of $\chi_s$, the induced representation
$\Ind_N^\Gamma \chi_s^\gamma$ is weakly contained in $\pi$ 
for every $\gamma \in \Gamma$, by Proposition~\ref{Pro-WeakContSemiDirect}.
It follows that $\Ind_N^\Gamma \chi$ is weakly contained in $\pi$ 
for every $\chi \in S$.
Therefore the representation 
$$
\rho \, := \, \Ind_N^\Gamma (\bigoplus_{\chi \in S} \chi)
\cong \bigoplus_{\chi \in S}(\Ind_N^\Gamma \chi)
$$
is weakly contained in $\pi$. 
\par

On the one hand, observe that $\Ind_N^\Gamma(\pi \vert_N)$ 
is weakly equivalent to $\rho$.
Moreover, since
$\Ind_N^\Gamma(\pi \vert_N) \cong \pi \otimes \Ind_N^\Gamma 1_N$ 
and since $\Gamma/N$ is amenable, 
$\pi$ is weakly contained in $\Ind_N^\Gamma (\pi \vert_N)$; 
compare also the proofs of Theorems
\ref{Theo-PrimIdealTwoStepNilpotent} and \ref{Theo-PrimIdealAffGr}.
Therefore, $\pi$ is weakly contained in $\rho$. 
Therefore $\pi$ is weakly equivalent to $\rho$. 
\par

On the other hand, since $\Ind_N^\Gamma \chi_s^\gamma$ 
is equivalent to $\Ind_N^\Gamma \chi_s$ for every $\gamma \in A$ 
and since $S$ is the closure of 
the $A$-orbit of $\chi_s$, continuity of induction shows that
$\rho$ is weakly equivalent to $\Ind_N^\Gamma \chi_s$.
Therefore $\pi$ is weakly equivalent to $\Ind_N^\Gamma \chi_s$, as claimed.

\vskip.2cm

Assume now that $s$ is a periodic point for $T_p$, 
that is, $s \in \So_p(n)$ for some $n \ge 1$.
We claim that $\pi$ is equivalent to $\Ind_{\Gamma(n)}^\Gamma \chi_{s, \theta}$
for some $\theta \in \mathopen[0,1/n \mathclose[$.
\par
 
Indeed, the support $S$ of $E$ consists of the finite $A$-orbit of $\chi_s$.
Therefore, by Corollary~\ref{Cor-ContAbelian},
we have $\Hi^{\chi} \ne \{0\}$ for every $\chi \in S$ 
and $\Hi = \bigoplus_{\chi \in S} \Hi^{\chi}$, where
$$
\Hi^{\chi} \, = \, \left\{ \xi \in \Hi \mid \pi(n) \xi = \chi(n) \xi 
\hskip.5cm \text{for all} \hskip.2cm
n \in N \right\}.
$$
\par

We now proceed as in proof of the Third Step of
Theorem~\ref{Theo-FiniteDimRepSemiDirect}.
Since $N$ is a normal subgroup of $\Gamma$,
we have $\pi(\gamma^{-1}) \Hi^{\chi} = \Hi^{\chi^\gamma}$
for every $\gamma \in \Gamma$ and $\chi \in S$.
In particular, we have $\pi(\gamma) \Hi^{\chi_s} = \Hi^{\chi_s}$ 
for every $\gamma$ in the stabilizer $\Gamma(n)$ of $\chi_s$ in $\Gamma$.
It follows that $\pi$ is equivalent to the induced representation 
$\Ind_{\Gamma(n)}^\Gamma \sigma$, 
where $\sigma$ is the subrepresentation of $\pi \vert_{\Gamma(n)}$ 
defined on the $\Gamma(n)$-invariant subspace $\Hi^{\chi_s}$.
\par

The representation $\sigma$ is an irreducible representation 
of $\Gamma(n) = A(n) \ltimes N$ 
which is trivial on the commutator subgroup of $\Gamma(n)$.
It follows that $\sigma$ is a unitary character of $\Gamma(n)$.
Therefore, $\pi$ is equivalent to $\Ind_{\Gamma(n)}^\Gamma \chi_{s, \theta}$. 

\vskip.2cm

$\bullet$ \emph{Second step.} 
We claim that every ideal in the image $\Phi$ is a primitive ideal.
\par

Indeed, let $s \in \So_p(\infty)$. The representation 
$\Ind_N^\Gamma \chi_{s}$ is irreducible
(Proposition~\ref{Prop-IrredRepBS-bis}).
So, $\Phi(s) \in \Pri(\Gamma)$.
\par

Similarly, for $n \ge 1$ and $(s, \theta) \in \So_p(n)$, 
the induced representation $\Ind_{\Gamma(n)}^\Gamma \chi_{s, \theta}$ 
is irreducible
(Proposition~\ref{Prop-IrredRepBS-FiniteDim} or Corollary~\ref{Cor-FinDimRepBS})
and so $\Phi(s, \theta) \in \Pri(\Gamma)$.

\vskip.2cm

$\bullet$ \emph{Third step.} 
We claim that $\Phi$ factorizes to a map between 
the quasi-orbit space of $X_{\rm fd} \sqcup \So_p(\infty)$
and $\Pri(\Gamma)$.
\par

Indeed, let $s, s' \in \So_p(\infty)$ be two points with the same quasi-orbit.
Then, the first step applied to $\pi = \Ind_N^\Gamma \chi_{s'}$ shows that 
$$
\textnormal{C*ker}(\Ind_N^\Gamma \chi_{s'}) 
\, = \,
\textnormal{C*ker}(\Ind_N^\Gamma \chi_s) .
$$
\par

If $s,s' \in \So_p(n)$ are in the same quasi-orbit,
then they belong to the same $T_p$-orbit. 
Therefore,
$\Ind_{\Gamma(n)}^\Gamma \chi_{s, \theta}$ 
and $\Ind_{\Gamma(n)}^\Gamma \chi_{s', \theta}$ 
are equivalent for every $\theta \in \mathopen[0,1/n \mathclose[$.

\vskip.2cm

$\bullet$ \emph{Fourth step.} 
We claim that the map induced by $\Phi$ between
the quasi-orbit space of $X_{\rm fd} \sqcup \So_p(\infty)$ 
and $\Pri(\Gamma)$ is injective.
\par

Indeed, let $s \in \So_p(\infty)$ and set $\pi = \Ind_{N}^\Gamma \chi_{s}$.
Then the support of the projection-valued measure associated to $\pi \vert_N$ 
is the closure of the $A$-orbit of $\chi_s$. 
Let $s' \in \So_p(\infty)$ be such that $\pi' = \Ind_{N}^\Gamma \chi_{s'}$ 
is weakly equivalent to $\pi$.
Then $\pi' \vert_N$ is weakly equivalent to $\pi' \vert_N$.
Therefore, the closures of the orbits of $s'$ and $s$ coincide,
that is, $s'$ and $s$ belong to the same quasi-orbit.
\par

Let $s \in \So_p(n)$ and $\theta \in \mathopen[0,1/n \mathclose[$;
set $\pi = \Ind_{\Gamma(n)}^\Gamma \chi_{s, \theta}$.
Let $s' \in \So_p(m)$ and $\theta' \in \mathopen[0,1/m \mathclose[$ 
be such that $\pi' = \Ind_{\Gamma(m)}^\Gamma \chi_{s', \theta'}$
is weakly equivalent to $\pi$.
Then the irreducible representations $\pi$ and $\pi'$ are equivalent,
since they are finite-dimensional (see Proposition~\ref{PropOnWeAndEq}). 
Therefore, by Corollary~\ref{Cor-FinDimRepBS},
$(s, \theta)$ and $(s', \theta')$
belong to the same $\widetilde T_p$-orbit in $X_{\rm fd}$.
\end{proof}

\begin{rem}
% 9.D.2
\label{Rem-PrimIdealBS}
Let $\lambda$ be the normalized Haar measure on $\So_p$.
As was observed in Remark~\ref{Rem-Prop-IrredBS-bis}, 
the action of $T_p$ on $(\So_p, \lambda)$ is ergodic;
it follows that the $T_p$-orbit of $\lambda$-almost every element in $\So_p$ is dense
(see \cite[Chap IV, Proposition 1.5]{BeMa--00}). 
However, besides the periodic points, there are points in $\So_p$ 
with a non-dense orbit (see Lemma~\ref{Lem-NormalLamplighter}).
\end{rem}

\section
{Lamplighter group}
% Section 9.E
\label{Section-PrimIdealLamplighter}

\index{Lamplighter group}
Recall from Section~\ref{Section-IrrRepLamplighter}
that the lamplighter group is the semi-direct product $\Gamma = A \ltimes N$
of $A = \Z$ with $N = \bigoplus_{k \in \Z} \Z / 2 \Z$,
where the action of $\Z$ on $\bigoplus_{k \in \Z} \Z / 2 \Z$
is given by shifting the coordinates.
Recall also that $\widehat N$ can be identified
with $X = \prod_{k \in \Z} \{ 0,1 \}$,
the dual action $\Z$ on $\widehat N$ being given
by the shift transformation $T$ on $X$.
\par

As in Section~\ref{Section-IrrRepLamplighter},
$\Per(T)$ denotes the set of $T$-periodic points
and $X(n)$ the set of points in $X$ with $T$-period $n \ge 1$.
\par

For $x \in X$, we denote by $\chi_x$ the corresponding element in $\widehat N$.
Recall that, for $x \in X(n)$ and $\theta \in \mathopen[ 0,1/n \mathclose[$,
the induced representation 
$\pi_{x, \theta} = \Ind_{\Gamma(n)}^\Gamma \chi_{x, \theta}$
is an irreducible finite-dimensional representation of $\Gamma$,
where $\chi_{x, \theta}$ is a unitary character
of the finite index normal subgroup 
$\Gamma(n) = n\Z \ltimes N$
with $\chi_{x, \theta} = \chi_x$ on $N$.
\par

As in Proposition~\ref{Prop-IrredRepLamplighter-FiniteDim}, we set
$$
X_{\rm fd} \, = \, 
\bigsqcup_{n \ge 1} \hskip.1cm
\big(X(n) \times \mathopen[ 0,1/n \mathclose[ \big) ,
$$
and $T \,\colon X_{\rm fd} \to X_{\rm fd}$ is defined by $T(x, \theta) = (Tx, \theta)$.
We set also $X(\infty) = X \smallsetminus \Per(T)$.
\par

The proof of Theorem~\ref{Theo-PrimIdealBS} carries over \emph{mutatis mutandis}
and yields the following result.

\begin{theorem}
% 9.E.1
\label{Theo-PrimIdealLamplighter}
Let $\Gamma = A \ltimes N$ be the lamplighter group;
we keep the notation above.
The map 
$$
\Phi \, \colon \, X_{\rm fd} \sqcup X(\infty) \to \Pri(\Gamma),
$$
defined by
$$
\Phi(x) \, = \, 
\textnormal{C*ker} (\Ind_{\Gamma(n)}^\Gamma \chi_{x, \theta})
\hskip.5cm \text{for} \hskip.2cm 
n \ge 1, \hskip.1cm (x, \theta) \in X(n) \times \mathopen[ 0,1/n \mathclose[
$$
and
$$
\Phi(x) \, = \, 
\textnormal{C*ker} (\Ind_N^\Gamma \chi_x) 
\hskip.5cm \text{for} \hskip.2cm 
x \in X(\infty) ,
$$
factorizes to a bijection between the quasi-orbit space of the $T$-action 
on $X_{\rm fd} \sqcup X(\infty)$ and $\Pri(\Gamma)$.
\end{theorem}

\section
{General linear groups}
% Section 9.F
\label{Sect-PrimIdealGLn}

Let $\K$ be an infinite field and $n \ge 2$.
Following \cite{HoRo--89}, we will determine the maximal ideals 
in the C*-algebra $C^*_{\rm max} (\Gamma)$ 
for the general linear group $\Gamma = \GL_n(\K)$,
viewed as discrete group.
When $\GL_n(\K)$ is amenable (equivalently, when
$\K$ is an algebraic extension of a finite field, see Proposition~\ref{Pro-AmenableGLn}),
we will describe the primitive dual $\Pri( \GL_n(\K))$.
\par

As preparation, will need several lemmas in which $\K$ will be an arbitrary infinite field.
\par
 
The group $\Gamma = \GL_n(\K)$ acts naturally on $\K^n$,
and hence on $\widehat \K^n \approx \K^n$ by duality,
through the formula
$$
\chi^\gamma(x) \, = \, \chi(\gamma^{-1}x) 
\hskip.5cm \text{for} \hskip.2cm
\gamma \in \GL_n(\K), \hskip.1cm \chi \in \widehat \K^n, \hskip.1cm x \in \K^n.
$$

\begin{lem}
% 9.F.1
\label{Lem-FreeActDualGL_n}
We keep the notation above and assume that $\K$ is \emph{infinite}.
\begin{enumerate}[label=(\arabic*)]
\item\label{iDELem-FreeActDualGL_n}
The $\Gamma$-orbit of every $\chi \in \widehat \K^n \smallsetminus \{1 \}$ 
is dense in $\widehat \K^n$.
\item\label{iiDELem-FreeActDualGL_n}
Assume moreover that $\K$ is countable. The set of $\chi \in \widehat \K^n$
with trivial stabilizer in $\Gamma$ is a dense G$_\delta$ set in $\widehat \K^n$.
\end{enumerate}
\end{lem}

\begin{proof}
\ref{iDELem-FreeActDualGL_n}
Let $\chi = (\chi_1, \hdots, \chi_n) \in \widehat \K^n \smallsetminus \{1 \}$.
We first claim that there exists $\chi' = (\chi_1', \hdots, \chi_n')$
in the $\Gamma$-orbit of $\chi$
with $\chi_i' \ne 1$ for every $i \in\{1, \hdots, n\}$. 
\par

Indeed, choose $i_0\in \{1, \hdots, n\}$ such that $\chi_{i_0} \ne 1$,
and let $i \in\{1, \hdots, n\}$ be such that $\chi_i = 1$. 
Let $\gamma$ be the elementary matrix $E_{i, i_{0}}(1)$.
Then $\chi^\gamma = (\chi_1', \hdots, \chi_n')$, where $\chi'_{i} = \chi_{i_0}$
and $\chi'_j = \chi_j$ for $j \ne i$. Therefore $\chi'_{i} \ne 1$. 
The claim follows by iteration of this procedure.
\par

So, it suffices to prove \ref{iDELem-FreeActDualGL_n}
in case $\chi_i \ne 1$ for every $i \in\{1, \hdots, n\}$.
For a diagonal matrix $\gamma \in A$
with diagonal entries $a_1, \hdots, a_n \in \K^\times$, 
we have
$$
\chi^\gamma \, = \, (\chi_1^{a_1}, \hdots, \chi_n^{a_n}),
$$
where $(a, \chi_i) \to \chi_i^{a}$ denotes
the natural action of $\K^\times$ on $\widehat \K$ 
Since $\chi_i \ne 1$, the $\K^\times$-orbit of $\chi_i$ is dense in $\widehat \K$
for every $i \in\{1, \hdots, n\}$, by Lemma~\ref{Lem-FreeActDualAff}.
This proves the claim.

\vskip.2cm

\ref{iiDELem-FreeActDualGL_n}
Let $\gamma \in \Gamma$ be distinct from the identity matrix $I_n$. 
Denote by ${\rm Fix}(\gamma)$ the set of fixed points of $\gamma$ in 
$\widehat \K^n$. 
\par

For $\chi \in \widehat \K^n$, we have $\chi \in {\rm Fix}(\gamma)$ if and only if
the subspace $V := (\gamma-I_n) \K^n$ of $\K^n$ is contained in $\ker\chi$. 
Therefore, ${\rm Fix}(\gamma)$ coincides with 
the annihilator $V^\perp$ of $V$ in $\widehat \K^n$.
In particular, ${\rm Fix}(\gamma)$ is a closed subgroup of $\widehat \K^n$.
Observe that $V \ne \{0\}$, since $g \ne I_n$.
%Moreover, since $g \ne I_n$, we have $(\gamma-I_n) \K^n \ne \{0\}$
%and hence ${\rm Fix}(\gamma) \ne \widehat \K^n$.
\par

We claim that ${\rm Fix}(\gamma)$ has empty interior. 
Indeed, assume by contradiction that this is not the case. 
Then ${\rm Fix}(\gamma)$ is an open subgroup of $\widehat \K^n$.
Since $\widehat \K^n$ is compact, it follows that ${\rm Fix}(\gamma)$ 
has finite index in $\widehat \K^n$.
However, by Pontrjagin duality, 
the dual group $\widehat V$ of the subgroup $V$ of $\K^n$
can be identified with $\widehat \K^n / V^{\perp}$.
Since ${\rm Fix}(\gamma) = V^{\perp}$, 
it follows that $\widehat V$ and hence $V$ are finite.
This contradicts the fact that $V \ne \{0\}$ and that $\K$ is infinite.
Therefore ${\rm Fix}(\gamma)$ has empty interior for every $\gamma \in \Gamma$ 
with $\gamma \ne I_n$.
\par

Since $\K$ is countable, $\Gamma$ is countable
and it follows from the Baire Category Theorem 
that the $F_\sigma$ subset
$$
X \, := \, \bigcup_{\gamma \in \Gamma, \gamma \ne I_n} {\rm Fix}(\gamma)
$$
of $\widehat \K^n$ has empty interior. 
The set of $\chi \in \widehat \K^n$ with trivial stabilizer in $\Gamma$ 
coincides with the complement of $X$
in $\widehat \K^n$ and is therefore a dense G$_\delta$ set.
\end{proof}

As in \ref{Section-IrrRepGLN}, 
we denote by $A$ and $B$ respectively the subgroups of diagonal 
and upper triangular matrices in $\Gamma$.
Let $W \approx {\mathcal S}_n$ be the subgroup 
of the permutation matrices in $\Gamma$, that is, 
the matrices with exactly one $1$ in each row and column and with $0$~'s elsewhere.
Observe that $W$ acts on $A \approx (\K^\times)^n$ by permuting coordinates.
Recall the \textbf{Bruhat decomposition} of $\Gamma$
(see for instance \cite[Chap.~2]{AlBe--95}): 
$$
\GL_n(\K) \, = \, \bigsqcup_{w \in W} BwB.
$$
\index{Bruhat decomposition of $\GL_n(\K)$} 
\index{General linear group! Bruhat decomposition of $\GL_n(\K)$} 

\begin{lem}
% 9.F.2
\label{Lem-PrincipalSerieRep}
Let $(\pi, \Hi)$ be a cyclic representation of $\Gamma$,
with cyclic unit vector $\xi \in \Hi$.
Assume that there exists a unitary character
$\chi = (\chi_1, \hdots, \chi_n)$ of $B$
with $\chi_1, \hdots, \chi_n$ pairwise distinct such that 
$$
\langle \pi(b) \xi \mid \xi \rangle \, = \, \chi(b)
\hskip.5cm \text{for every} \hskip.2cm
b \in B.
$$
Then $\pi$ is equivalent to
the principal series representation $\Ind_B^\Gamma \chi$.
\end{lem}
\index{Principal series! of $\GL_n(\K)$}
\index{General linear group! principal series of $\GL_n(\K)$}

\begin{proof}
Let $T$ be a set of representatives
for the classes in $B \backslash \Gamma$ with $e \in \Gamma$.
Then, with the notation as in Section~\ref{Section-IrrIndRep}, 
$\rho = \Ind_B^\Gamma \chi$ is realized on $\ell^2(T)$ by 
$$
(\rho(\gamma) f) (t) \, = \, \chi (\alpha(t, \gamma)) f(t \cdot \gamma)
\hskip.5cm \text{for all} \hskip.2cm
\gamma \in \Gamma, \hskip.1cm f \in \ell^2(T),
\hskip.2cm \text{and} \hskip.2cm
t \in T .
$$
The delta function $\delta_e$ is a cyclic vector for $\rho$.
To show that $\pi$ are $\rho$ are equivalent,
it suffices, by Proposition~\ref{GNSbijP(G)cyclic},
to prove that the two functions of positive type 
$$
\varphi_\pi \, \colon \,
\gamma \mapsto \langle \pi(\gamma) \xi \mid \xi \rangle
$$
and
$$
\varphi_\rho \, \colon \,
\gamma \mapsto \langle\rho(\gamma) \delta_e \mid \delta_e \rangle
$$
coincide on $\Gamma$.
\par

We have
$$
\varphi_\rho(b)
\, = \, \langle \rho(b) \delta_e \mid \delta_e \rangle
\, = \, \chi(b)
\hskip.5cm \text{for every} \hskip.2cm
b \in B
$$
and $\varphi_\rho(\gamma) = 0$ for $\gamma \in \Gamma \smallsetminus B$. 
This means that $\varphi_\rho$ is the trivial extension 
$\widetilde \psi$ of $\psi$ to $\Gamma$ (see also Proposition~\ref{diagcoeffinduced}).
It follows from the assumption that $\varphi_\pi = \varphi_\rho$ on $B$; 
it remains therefore to show that
$\varphi_\pi = 0$ on $\Gamma \smallsetminus B$.
\par

Let $b \in B$. We have, by assumption, 
$\langle \pi(b)\xi \mid \xi \rangle = \chi(b)$.
It follows from the equality case of the Cauchy--Schwarz inequality
that $\pi(b)\xi = \xi$.
\par

Let $w \in W$ and $b_1, b_2 \in B$. Then, we have 
$$
\begin{aligned}
\varphi_\pi(b_1 w b_2)
\, &= \, \langle \pi(b_1 wb_2)\xi \mid \xi \rangle
\, = \, \langle \pi(b_1)\pi(w)\pi(b_2)\xi \mid \xi \rangle
\\
\, &= \, \langle \pi(w)\pi(b_2)\xi \mid \pi(b_1^{-1})\xi \rangle
\, = \, \chi(b_2)\chi(b_1) \langle \pi(w)\xi \mid \xi \rangle
\, = \, \chi(b_2)\chi(b_1) \varphi_\pi(w).
\end{aligned} 
$$
In view of the Bruhat decomposition of $\Gamma$, it suffices therefore to prove that
$\varphi_\pi(w) = 0$ for every $w \in W \smallsetminus \{e\}$.
\par

Assume by contradiction that $\varphi_\pi(w) \ne 0$
for some $w \in W \smallsetminus \{e\}$.
Then, for every $a \in A$ we have $waw^{-1} \in A \subset B$ and hence
$$
\begin{aligned}
\chi(a) \varphi_\pi( w )
\, &= \, \langle \pi( wa)\xi \mid \xi \rangle
\\
\, &= \, \langle \pi(waw^{-1})\pi(w)\xi \mid \xi \rangle
\\
\, &= \, \langle \pi(w)\xi \mid \pi(wa^{-1}w^{-1}) \xi \rangle
\\
\, &= \, \overline{\chi(wa^{-1}w^{-1})} \langle \pi(w)\xi \mid \xi \rangle
\\
\, &= \, \chi(waw^{-1})\varphi_\pi(w).
\end{aligned} 
$$
Since $\varphi_\pi(w) \ne 0$,
it follows that $\chi(a) = \chi(waw^{-1})$ for every $a \in A$.
Let $\sigma \in {\mathcal S}$ be the permutation corresponding to $w$, so that
$$
waw^{-1} \, = \, (a_{\sigma(1)}, \hdots, a_{\sigma(n)})
\hskip.5cm \text{for all} \hskip.2cm
a = (a_1, \hdots, a_n) \in A.
$$
Then, for every $a = (a_1, \hdots, a_n) \in A$, we have
$$
\begin{aligned}
\chi_1(a_1) \cdots \chi_n(a_n)
\, &= \, \chi(a) \, = \, \chi(waw^{-1})
\\
\, &= \, \chi_1(a_{\sigma(1)}) \cdots \chi_n(a_{\sigma(n)})
\end{aligned}
$$
and hence $\chi_{\sigma(i)} = \chi_i$ for all $i \in \{1, \hdots, n\}$. 
Since $w \ne e$, there exists $i \in\{1, \hdots, n\}$ with $\sigma(i) \ne i$ 
and this is a contradiction to the assumption
that the $\chi_j$~'s are all pairwise distinct.
\end{proof}

Recall that the centre $Z$ of $\Gamma$ is the subgroup of scalar matrices.
For a unitary character $\psi \in \widehat Z \approx \widehat{\K^\times}$,
we denote by $\lambda_\psi$ the induced representation $\Ind_Z^\Gamma \psi$.

\begin{lem}
% 9.F.3
\label{Lem-PrincipalSerieRepWeakEqui}
Let $\chi = (\chi_1, \hdots, \chi_n)$ be a unitary character of $B$ and 
set $\psi := \chi \vert_Z \in \widehat Z$. 
The principal series representation $\Ind_B^\Gamma \chi$ 
 is weakly equivalent to $\lambda_\psi$.
\end{lem}

\begin{proof}
$\bullet$ {\it First step.}
We claim that $\lambda_\psi$ is weakly contained in 
$\rho := \Ind_B^\Gamma \chi$. 
\par

Let $\widetilde \psi$ be the trivial extension of $\psi$ to $\Gamma$. 
Since $\delta_e$ is a cyclic vector for $\lambda_\psi$
and since $\widetilde \psi$ is the function of positive type 
associated to $\lambda_\psi$ and $\delta_e$ 
(see proof of Lemma~\ref{Lem-PrincipalSerieRep}),
it suffices to show that $\widetilde \psi$ is uniform limit over finite subsets of
normalized functions of positive type associated to $\rho$.
\par

Since $Z \subset B$, it follows from the definition of induced representations that 
$$
\rho(z) \, = \, 
(\Ind_B^\Gamma\chi) (z) \, = \, \chi(z)I \, = \, \psi(z) I 
\hskip.5cm \text{for all} \hskip.2cm 
z \in Z.
$$
Let $F$ be a finite subset of $\Gamma$ and write $F = F_1 \cup F_2$ with
$F_1 \subset Z$ and $F_2 \subset \Gamma \smallsetminus Z$. 
Let $\gamma \in F_2$. Set 
$$
X^\gamma \, := \, \{ x \in \Gamma \mid Bx \gamma = Bx \} 
\, = \, \{ x \in \Gamma \mid x \gamma x^{-1} \in B \}.
$$
Since $B$ is an algebraic subvariety of $\Gamma = \GL_n(\K)$,
we see that $X^\gamma$ is an algebraic subvariety of $\Gamma$.
Observe that
$$
\bigcap_{x \in \Gamma} x^{-1} Bx= Z.
$$
Therefore $X^\gamma$ has positive codimension in $\GL_n(\K)$,
because $\gamma \notin Z$. 
Since $\K$ is infinite,
it follows that $\bigcup_{\gamma \in F_2} X^\gamma$ is a proper subset of $\Gamma$. 
Therefore there exists $x \in \Gamma$ such that $x \notin X^\gamma$,
that is, such that $Bx \gamma \ne Bx$ for every $\gamma \in F_2$. 
\par

Let $T$ be a set of representatives for $B \backslash G$ with $x \in T$. 
For all $\gamma \in F_2$, we have $x \cdot \gamma \ne x$ and hence
$$
\begin{aligned}
\langle \rho(\gamma) \delta_x \mid \delta_x \rangle
\, &= \, \langle (\Ind_B^\Gamma\chi) (\gamma) \delta_x \mid \delta_x \rangle
\\
\, &= \, \chi(\alpha(x, \gamma)) \langle \delta_{x \cdot \gamma} \mid \delta_x \rangle
\\
\, &= \, 0.
\end{aligned}
$$
Moreover, for $\gamma \in F_1$ we have 
$$
\begin{aligned}
\langle \rho(\gamma) \delta_x \mid \delta_x \rangle
\, &= \, \psi(\gamma) \langle \delta_{x} \mid \delta_x \rangle
\\
\, &= \, \psi(\gamma),
\end{aligned}
$$
since $F_1 \subset Z$.
The claim is therefore proved.

\vskip.2cm.

$\bullet$ {\it Second step.}
We claim that $\rho = \Ind_B^\Gamma \chi$ is weakly contained in $\lambda_\psi$. 
\par

Indeed, $B$ is solvable and hence amenable.
Therefore, $1_B$ is weakly contained 
in the quasi-regular representation $\lambda_{B/Z} = \Ind_Z^B 1_Z$.
Using Proposition \ref{InductionQqPropr}~\ref{ivDEInductionQqPropr},
we have
$$
\Ind_Z^B \psi \, = \, \Ind_Z^B (\chi \vert_Z) \, \simeq \, \chi \otimes \Ind_Z^B 1_Z
\, = \, \chi \otimes \lambda_{B/Z},
$$
it follows that $\chi \simeq \chi \otimes 1_B$
is weakly contained in $\chi \otimes \lambda_{B/Z} \simeq \Ind_Z^B \psi$.
Therefore $\rho \, = \, \Ind_B^\Gamma \chi$ is weakly contained in
$\Ind_B^\Gamma(\Ind_Z^B \psi) \simeq \Ind_Z^\Gamma \psi
= \lambda_\psi$,
as was to be shown.
\end{proof}

\begin{rem}
% 9.F.4
\label{Rem-PrincipalSerieRepWeakEqui}
Let $\Ind_B^\Gamma\chi_1$ and $\Ind_B^\Gamma \chi_2$ 
be two principal series representations of $\Gamma$ 
with the same central character, that is, such that $\chi_1 \vert_Z = \chi_2 \vert_Z$. 
It follows from Lemma~\ref{Lem-PrincipalSerieRepWeakEqui} that 
$\Ind_B^\Gamma\chi_1$ and $\Ind_B^\Gamma \chi_2$ are weakly equivalent.
However, observe that $\Ind_B^\Gamma\chi_1$ and $\Ind_B^\Gamma \chi_2$ 
are not equivalent if $\chi_1 \ne \chi_2$, by Theorem~\ref{Theo-RepIrrGLn}.
\end{rem}

Recall that the commutator subgroup of $\Gamma = \GL_n(\K)$ is $\SL_n(\K)$; 
hence, the one-dimensional representations of 
$\Gamma$ correspond to the unitary characters of 
the abelianized group
$$
\Gamma_{\rm ab} \, = \, \GL_n(\K)/\SL_n(\K) \, \approx \, \K^\times.
$$
One checks that $\psi \in \widehat Z \approx \K^\times$ arises as the central character 
of a one-dimensional representation of $\Gamma$ if and only if 
$\psi = 1$ on the group $\mu_n^{\K}$ of the $n$-th roots of unity in $\K$. 
\par

Here is the first main result of this subsection.

\begin{theorem}
% 9.F.5
\label{Theo-WeakContRepGLn}
Let $\K$ be an infinite field and $n \ge 2$ an integer. 
Let $\pi$ be a representation of $\Gamma = \GL_n(\K)$
which is not trivial on $\SL_n(\K)$.
Assume that the centre $Z$ of $\Gamma$
is contained in the projective kernel of $\pi$. 
Then $\lambda_\psi = \Ind_Z^\Gamma\psi$ is weakly contained in $\pi$,
where $\psi \in \widehat Z$ is the central character of $\pi$.
\end{theorem}

\begin{proof}
$\bullet$ {\it First step.} 
We claim that it suffices to prove the theorem 
in the case of a countable infinite field $\K$.

Indeed, assume that the theorem has been proved in the case of a countable field. 
Since $\delta_e$ is a cyclic vector for $\lambda_\psi$, 
we have to show that the trivial extension $\widetilde \psi$ 
of $\psi$ to $\Gamma$ is uniform limit over finite subsets of
normalized functions of positive type associated to $\pi$ 
(see the proof of Lemma~\ref{Lem-PrincipalSerieRepWeakEqui}).

Let $F$ be a finite subset of $\Gamma$. Since $\pi$ is not trivial on $\SL_n(\K)$,
upon adding to $F$ a suitable element of $\SL_n(\K)$, we can assume
that $\pi(\gamma)$ is not the identity for some $\gamma \in F$. 
Let $\mathbf L$ be the subfield of $\K$ generated by 
the coefficients of all matrices in $F$
(in case $\mathbf L$ is finite, we replace $\mathbf L$ 
by a countable infinite subfield of $\K$ containing $\mathbf L$). 
Then $\mathbf L$ is countable and $\pi$ is not trivial on $\SL_n(\mathbf L)$. 
Set $\Gamma' = \GL_n(\mathbf L)$.
Then $\Ind_{Z'}^{\Gamma'}\psi'$ is weakly contained
in $\pi \vert_{\Gamma'}$, where $Z' = Z\cap \Gamma$ is the centre of $\Gamma'$ 
and $\psi' = \psi \vert_{Z'}$.
Since $F\subset \Gamma'$, it follows that $\widetilde \psi'$ and hence $\widetilde \psi$ 
is uniform limit over $F$ of normalized functions of positive type associated to $\pi$. 

The claim is proved and so we may assume from now on that $\K$ is countable. 

Next, we will need to consider the following subgroup 
$Q = H \ltimes U \approx \GL_{n-1}(\K) \ltimes \K^{n-1}$ of $\Gamma$, where
$$
H \, = \, \left\{ \begin{pmatrix} A&0 \\
0 & 1 \\
\end{pmatrix}
\mid A \in \GL_{n-1}(\K) \right\} 
\, \approx \, \GL_{n-1}(\K)
$$
and 
$$
U \, = \, \left\{ \begin{pmatrix} I_{n-1} & x \\
0 & 1 \\
\end{pmatrix}
\mid x \in \K^{n-1}\right\}
\, \approx \, \K^{n-1}.
$$
Observe that $Q \approx \Aff(\K)$ in case $n = 2$.

\vskip.2cm

$\bullet$ {\it Second step.} 
We claim that the restriction $\pi \vert_U$ of $\pi$ to $U$ is not trivial.
\par

Indeed, assume that $\pi \vert_U$ is trivial, that is, $U\subset \ker \pi$. 
Every elementary matrix $E_{i, j}(\alpha)$ 
for $\alpha \in \K$ and $i, j \in\{1, \hdots, n\}, i \ne j$, is conjugate
to an elementary matrix contained in $U$. It follows that $\ker \pi$ contains 
all elementary matrices $E_{i, j}(\alpha)$ and hence $\SL_n(\K)$, 
since $\SL_n(\K)$ is generated by these matrices. 
This is a contradiction, as $\pi$ is not trivial on $\SL_n(\K)$.

\vskip.2cm

$\bullet$ {\it Third step.} 
We claim that the regular representation $\lambda_Q$ of $Q$ is 
weakly contained in $\pi \vert_Q$.
\par

Indeed, by the second step, there exists $\chi \in \widehat U \approx \widehat \K^n$
with $\chi \ne 1$ which is weakly contained in $\pi \vert_U$. 
Then $\chi^\gamma$ is weakly contained in $\pi \vert_U$
for every $\gamma \in H$ (see Remark~\ref{Rem-WeakContSemiDirect})
It follows from Lemma~\ref{Lem-FreeActDualGL_n} that every $\chi \in \widehat U$
is weakly contained $\pi \vert_U$. 
Again by Lemma~\ref{Lem-FreeActDualGL_n} (here we use the countability of $\K$),
there exists $\chi_0\in \widehat U$ with trivial stabilizer in $H$ which. 
is weakly contained in $\pi \vert_U$. 
Proposition~\ref{Pro-WeakContSemiDirect}
implies therefore that 
$\Ind_U^Q \chi_0^\gamma$ is weakly contained in $\pi \vert_Q$ for every $\gamma \in H$.
Again by Lemma~\ref{Lem-FreeActDualGL_n},
the set $\{\chi_0^\gamma \mid \gamma \in H \}$ is dense in $\widehat U$.
It follows that $\Ind_U^Q \chi$ 
is weakly contained in $\pi \vert_Q$ for every $\chi \in \widehat U$.
Since the direct sum $\bigoplus_{\chi \in \widehat U} \chi$
is weakly equivalent to the regular representation $\lambda_U$, 
this implies that $\lambda_Q = \Ind_U^Q \lambda_U$ is weakly contained 
$\pi \vert_Q$. 

\vskip.2cm

$\bullet$ {\it Fourth step.} 
We claim that $\pi \vert_B$ weakly contains
a unitary character $\chi = (\chi_1, \hdots, \chi_n)$ of $B$
with $\chi \vert_Z = \psi$ and such that the $\chi_i$~'s
are pairwise distinct.
\par

Indeed, let $B(n-1)$ be the subgroup of the upper-triangular matrices in
$H \approx \GL_{n-1}(\K)$
and consider the subgroup $S=B(n-1)\ltimes U$ of $Q$.
By the third step, $\lambda_Q$ is weakly contained in $\pi \vert_Q$. 
Therefore $\lambda_{S}$ is weakly contained in $\pi \vert_{S}$.
Since $S$ is amenable, it follows that every one-dimensional representation of $S$
is weakly contained in $\pi \vert_{S}$.
\par

Observe that $B$ is the Cartesian product of the subgroup $S$ with the centre $Z$. 
Since $\psi$ is the central character of $\pi$, 
it follows that every one-dimensional representation of $B$ 
with central character $\psi$ is weakly contained in $\pi \vert_{B}$.
As $\K$ is infinite, we can find a unitary character 
$\chi = (\chi_1, \hdots, \chi_n) \in \widehat{(\K^\times)^n}$ 
of $B$ with $\chi \vert_Z = \psi$ and all $\chi_i$ pairwise distinct. 

\vskip.2cm

$\bullet$ {\it Fifth step.} 
We claim that $\pi$ weakly contains a principal series representations
$\Ind_B^\Gamma \chi$ for $\chi = (\chi_1, \hdots, \chi_n) \in \widehat A$
with $\chi \vert_Z = \psi$ and such that the $\chi_i$~'s
are pairwise distinct.
\par

Indeed, let $\chi = (\chi_1, \hdots, \chi_n)$ be a unitary character of $B$
as in the fourth step.
Since $\chi$ is weakly contained in $\pi \vert_B$, there exists a sequence
$\varphi_n$ of normalized functions of positive type associated to $\pi$ such that 
$$
\lim_n\varphi_n(b) \, = \, \chi(b) 
\hskip.5cm \text{for all} \hskip.2cm
b \in B.
$$
By compactness of the set of normalized functions of positive type on $\Gamma$
for the topology of pointwise convergence, upon passing to a subsequence,
we can assume that $\varphi_n$ converges
to a normalized function of positive type $\varphi$ on $\Gamma$.
Let $(\rho, \Hi, \xi)$ be the associated GNS triple
(see Construction~\ref{constructionGNS2}).
Since $\varphi \vert_B = \chi$ and since the $\chi_i$~'s
are pairwise distinct, Lemma~\ref{Lem-PrincipalSerieRep} 
shows that $\rho$ is equivalent to the
principal series representations $\Ind_B^\Gamma \chi$.
As $\lim_n \varphi_n = \varphi$, 
the representation $\rho$ and hence $\Ind_B^\Gamma \chi$
is weakly contained in $\pi$.

\vskip.2cm

$\bullet$ {\it Sixth step.} 
We claim that $\pi$ weakly contains the representation $\lambda_\psi = \Ind_Z^\Gamma\psi$.
Indeed, this follows from the fifth step and from Lemma~\ref{Lem-PrincipalSerieRepWeakEqui}.
\end{proof}

As a consequence of Theorem~\ref{Theo-WeakContRepGLn},
we obtain a description of the maximal ideals
in the C*-algebra $C^*_{\rm max} (\GL_n(\K))$
for an arbitrary infinite field $\K$.

\begin{cor}
% 9.F.6
\label{Cor-MaxIdealGLn}
Let $\K$ be an infinite field and $\Gamma = \GL_n(\K)$ for $n \ge 2$.
The maximal two-sided ideals of $C^*_{\rm max} (\Gamma)$ are: 
\begin{enumerate}
\item[$\bullet$]
either C*-kernels of one-dimensional representations of $\Gamma$, that is, 
ideals of the form $\textnormal{C*ker}(\chi)$ for a unitary character $\chi$ of
$\Gamma_{\rm ab} = \GL_n(\K) / \SL_n(\K)$;
\item[$\bullet$]
or ideals of the form $\textnormal{C*ker}(\Ind_Z^\Gamma \psi)$, 
for some $\psi \in \widehat Z$.
\end{enumerate}
\end{cor}

\begin{proof} 
Let $J$ be a maximal ideal of $C^*_{\rm max} (\Gamma)$; then 
$J$ is a primitive ideal of $\Gamma$ and hence $J = \textnormal{C*ker}(\pi)$
for an irreducible representation $\pi$ of $\Gamma$. 

If $\pi$ is one-dimensional, then $J$ is an ideal of codimension 1 in $C^*_{\rm max}(\Gamma)$ 
and is hence a maximal ideal.

Assume now that $\pi$ is not one-dimensional. 
Then, by Theorem~\ref{Theo-WeakContRepGLn}, 
$\Ind_Z^\Gamma \psi$ is weakly contained in $\pi$, that is, 
$$
J \, = \, \textnormal{C*ker}(\pi)
\, \subset \, \textnormal{C*ker}(\Ind_Z^\Gamma \psi).
$$
Therefore $J = \textnormal{C*ker}(\Ind_Z^\Gamma \psi)$, by maximality of $J$.
\end{proof}

\begin{rem}
% 9.F.7
\label{Rem-MaxIdealsGLn}
In Corollary~\ref{Cor-MaxIdealGLn}, 
it may happen that an ideal of the form 
$\textnormal{C*ker}(\Ind_Z^\Gamma \psi)$ for $\psi \in \widehat Z$ is not maximal. 
Indeed, assume that $\Gamma$ is amenable
(equivalently, that $\K$ be an infinite algebraic extension of a finite field;
see Proposition~\ref{Pro-AmenableGLn} below). 
Then the trivial representation $1_\Gamma$ 
is weakly contained in $\Ind_Z^\Gamma 1_Z$, that is, 
$\textnormal{C*ker}(\Ind_Z^\Gamma 1_Z)
\subset \textnormal{C*ker}(1_\Gamma)$.
However, $\Ind_Z^\Gamma 1_Z$ is obviously not weakly contained in $1_\Gamma$.
Therefore $\textnormal{C*ker}(\Ind_Z^\Gamma 1_Z)$ is not a maximal ideal.
\end{rem}

When $\GL_n(\K)$ is amenable, Theorem~\ref{Theo-WeakContRepGLn}
enables us to give a complete description of the primitive dual of $\GL_n(\K)$.
The fields $\K$ for which $\GL_n(\K)$ is amenable can easily be characterized.

\begin{prop}
% 9.F.8
\label{Pro-AmenableGLn} 
Let $\K$ be a field.
The following properties are equivalent:
\begin{enumerate}[label=(\roman*)]
\item\label{iDEPro-AmenableGLn}
the group $\GL_n(\K)$ is amenable;
\item\label{iiDEPro-AmenableGLn}
the field $\K$ is an algebraic extension of a finite field.
\end{enumerate}
\end{prop}

\begin{proof}
Assume that $\K$ is an algebraic extension of a finite field $\F_q$.
Then $\K = \bigcup_{m} \K_m$
for an increasing family of finite extensions $\K_m$ of $\F_q$,
and therefore $\GL_n(\K) = \bigcup_m \GL_n(\K_m)$ is the
inductive limit of the finite and hence amenable groups $\GL_n(\K_m)$.
It follows that $\GL_n(\K)$ is amenable.

\vskip.2cm

Assume that $\K$ is not an algebraic extension of a finite field. 
If the characteristic of $\K$ is $0$, then $\K$ contains $\Q$;
if $\K$ is of characteristic $p$,
then $\K$ contains the purely transcendental extension $\F_p(T)$ of $\F_p$.
Therefore $\GL_n(\K)$ contains a copy of $\GL_2(\Q)$
or a copy of $\GL_2(\F_p(T))$.
\par

As is well-known, each of the groups $\GL_2(\Q)$ or $\GL_2(\F_p(T))$
contains a free group on two generators.
\par

Indeed, denote by $K$ either $\Q$ or $\F_p(T)$.
Let $\vert \cdot \vert$ be an absolute value on $K$
with respect to which the completion $k$ of $K$ is a local field.
In the case $K = \Q,$ we could take 
for $\vert \cdot \vert$ the usual Archimedean absolue value for which $k = \R$;
however, in the case $K = \F_p(T)$,
there is essentially only one choice for 
$\vert \cdot\ vert$ and then $k$
is the field of Laurent series $\F_p((T^{-1}))$ over $\F_p$.
\par

Choose $\lambda \in K^*$ with $\vert \lambda \vert \ne 1$
and $a \in K \smallsetminus \{0, 1\}$.
Consider the following matrices $A, B \in \GL_2(K)$
$$
A \, = \,
\begin{pmatrix} \lambda & 0 \\ 0 &1 \end{pmatrix} 
\hskip.5cm \text{and} \hskip.5cm
B \, = \, CAC^{-1},
\hskip.5cm \text{where} \hskip.5cm
 C \, = \,
\begin{pmatrix} 1 & a \\ 1 & 1 \end{pmatrix},
$$
and their action on the projective line ${\mathbf P}(k) = k \cup \{\infty\}$.
The set of fixed points of $A$ is $\{0,\infty\}$ 
and is disjoint from the set $\{a, 1\}$ of fixed points of $B$.
Observe that $\lim_{n \to +\infty} A^{\pm n} x \in \{0, \infty\}$
and $\lim_{n \to +\infty} B^{\pm n} x \in \{1, a\}$ 
for every $x \in {\mathbf P}(k)$.
Choose compact and disjoint subsets $X_1$ and $X_2$ of ${\mathbf P}(k)$ containing respectively $\{0, \infty\}$ and $\{1, a\}$.
There exists $N \ge 1$ such that
we have 
$$
A^{\pm m}(X_2) \subset X_1
\hskip.5cm \text{and} \hskip.5cm
B^{\pm m}(X_1) \subset X_2,
$$
for every $m\ge N$.
It follows from the classical ping-pong lemma
that $A^N$ and $B^N$ generate a free subgroup of $\GL_2(K).$
\par

Since amenability is inherited by subgroups, it follows that $\GL_n(\K)$ is not amenable. 
\end{proof}

%
%
\iffalse
\begin{proof}
% first alinea did not change
\par

% what was written up to end of October 2019
Assume that $\K$ is not an algebraic extension of a finite field. 
Then $\K$ contains $\Q$, if the characteristic of $\K$ is $0$,
and $\K$ contains the purely transcendental extension $\F_p(T)$
if $\K$ is of characteristic $p$.
Therefore $\GL_n(\K)$ contains a copy of $\GL_2(\Q)$ or a copy of $\GL_2(\F_p(T))$.
As is well-known, each of the groups $\GL_2(\Q)$ or $\GL_2(\F_p(T))$
contains a free group on two generators; 
indeed, the matrices
$$
A \, = \, \begin{pmatrix} 1 & 1 \\ 0 & 1 \end{pmatrix} 
\hskip.5cm \text{and} \hskip.5cm
B \, = \, \begin{pmatrix} 1 & 0 \\ 1 & 1 \end{pmatrix},
$$
respectively the matrices 
$$
A = \begin{pmatrix} 1 & T \\ 0 & 1 \end{pmatrix} 
\hskip.5cm \text{and} \hskip.2cm
B = \begin{pmatrix} 1 &0\\ T & 1 \end{pmatrix}
$$
have distinct fixed points in the projective line
${\mathbf P}(\R)$, respectively ${\mathbf P}(\F_p((T^{-1})))$, 
where $\F_p((T^{-1}))$ is the local field of Laurent series over $\F_p$; 
this implies that the subgroup of $\GL_2(\Q)$, respectively of $\GL_2(\F_p(T))$, 
generated by $A$ and $B$ is free on $A$ and $B$.
\par
 
Since amenability is inherited by subgroups, it follows that $\GL_n(\K)$ is not amen\-able. 
\end{proof}
\fi
%
%

\begin{theorem}
% 9.F.9
\label{Theo-PrimIdealPGLn}
Let $\Gamma = \GL_n(\K)$, 
where $n \ge 2$ and $\K$ is an infinite algebraic extension of a finite field.
The map 
$$
\Phi \, \colon \, \widehat{\Gamma_{\rm ab}} \sqcup \widehat Z 
\to 
\Pri(\Gamma),
$$
defined by
$$
\Phi(\chi) = \textnormal{C*ker}(\chi)
\hskip.5cm \text{for} \hskip.2cm
\chi \in \widehat \Gamma_{\rm ab}
$$
and
$$
\Phi(\psi) = \textnormal{C*ker}(\Ind_Z^\Gamma \psi)
\hskip.5cm \text{for} \hskip.2cm
\psi \in \widehat Z \approx \widehat{\K^\times}
$$
is a bijection.
\end{theorem}

\begin{proof}
\par
\noindent
$\bullet$ {\it First step.} 
We claim that the image of the map $\Phi$ is indeed contained in $\Pri(\Gamma)$.
\par

Observe first that, for every $\chi \in \widehat{\Gamma_{\rm ab}}$, the ideal 
$\textnormal{C*ker}(\chi)$ is of course a primitive ideal.
\par

Let $\psi \in \widehat Z$ and set $J = \textnormal{C*ker}(\Ind_Z^\Gamma \psi)$.
We can choose a unitary character 
$\chi = (\chi_1, \hdots, \chi_n) \in \widehat{\K^\times}^n$ of $B$ such that 
$\chi \vert_Z = \psi$
(for instance, by taking $\chi_1 = \psi$ and $\chi_i = 1$ for all $i \ne 1$).
Let $\pi = \Ind_B^\Gamma \chi$ be the corresponding principal series representation.
Then $\pi$ is irreducible (Theorem~\ref{Theo-RepIrrGLn}) 
and the central character of $\pi$ is $\psi$.
On the one hand, since $\pi$ is not one-dimensional, $\Ind_Z^\Gamma \psi$ is weakly contained
in $\pi$, by Theorem~\ref{Theo-WeakContRepGLn}. 
On the other hand, since $\Gamma$ is amenable and since 
$\Ind_Z^\Gamma(\pi \vert_Z)$ and $\pi \otimes \Ind_Z^\Gamma 1_Z$ are equivalent,
$\pi$ is weakly contained $\Ind_Z^G(\pi \vert_Z)$
which is weakly equivalent to $\Ind_Z^\Gamma \psi$.
Therefore $\pi$ and $\Ind_Z^\Gamma \psi$ are weakly equivalent.
Therefore, $J = \textnormal{C*ker}(\pi)$ and this shows that $J\in \Pri(\Gamma)$.
So, the image of the map $\Phi$ is indeed contained in $\Pri(\Gamma)$.

\vskip.2cm

$\bullet$ {\it Second step.} 
We claim that the map $\Phi$ is surjective. 
\par

Indeed, let $\pi$ be an irreducible representation of $\Gamma$, 
with central character $\psi$. 
Assume that $\pi$ is not one-dimensional.
On the one hand, $\Ind_Z^\Gamma \psi$ is weakly contained
in $\pi$, by Theorem~\ref{Theo-WeakContRepGLn}.
 On the other hand, as in the first step, 
$\pi$ is weakly contained in $\Ind_Z^\Gamma \psi$, since $\Gamma$ is amenable.
This shows that 
$$
\textnormal{C*ker}(\pi) \, = \, \textnormal{C*ker}(\Ind_Z^\Gamma \psi).
$$

\vskip.2cm

$\bullet$ {\it Third step.} 
We claim that the map $\Phi$ is injective. 
\par

Indeed, let $\chi_1, \chi_2 \in \widehat{\Gamma_{\rm ab}}$ be such 
that $\Phi(\chi_1) = \Phi(\chi_2)$, that is, such that $\chi_1$ and $\chi_2$
are weakly equivalent. Then $\chi_1 = \chi_2$.
\par

Next, let $\psi_1, \psi_2 \in \widehat Z$ be such that $\Phi(\psi_1) = \Phi(\psi_2)$, 
that is, such that $\lambda_{\psi_1} = \Ind_Z^\Gamma \psi_1)$ 
and $\lambda_{\psi_2} = \Ind_Z^\Gamma \psi_2)$ are weakly equivalent.
Since $\lambda_{\psi_1}(z) = \psi_1(z)I$
and $\lambda_{\psi_2}(z) = \psi_2(z)I$ for every $z\in Z$,
the unitary characters $\psi_1, \psi_2$ of $Z$ are weakly equivalent and hence equal.
\par

Finally, let $\chi \in \widehat{\Gamma_{\rm ab}}$ and $\psi \in \widehat Z$. 
As in the first step, there exists a principal series representation $\pi$ of $\Gamma$ such that 
$\pi$ is weakly equivalent to $\Ind_Z^\Gamma \psi$. 
Since $\pi$ is an infinite-dimensional irreducible representation, $\pi$ cannot be 
weakly equivalent to the finite-dimensional representation $\chi$ 
% (see \cite[F.2.9]{BeHV--08} and hence $\Phi(\chi) \ne \Phi(\psi)$.
(see Proposition~\ref{PropOnWeAndEq}~\ref{iDEPropOnWeAndEq})
and hence $\Phi(\chi) \ne \Phi(\psi)$.
\end{proof}
 
Theorem~\ref{Theo-PrimIdealPGLn} takes a striking form, when we go over to the projective linear group 
$\PGL_n(\K) = \GL_n(\K)/Z$. 

\begin{cor}
% 9.F.10
\label{Cor-PrimIdealPGLn}
Let $\Gamma = \PGL_n(\K)$, 
where $n \ge 2$ and $\K$ is an infinite algebraic extension of a finite field.
Then 
$$
\Pri(\Gamma) \, = \, 
\{0\} \cup
\left\{\textnormal{C*ker}(\chi) \mid \chi \in \widehat{\Gamma_{\rm ab}}\right\}
$$
\end{cor}
 
\begin{proof}
Let $G= \GL_n(\K)$. For $\psi \in \widehat Z$, the
representation $\Ind_Z^G \psi$ factorizes to $\Gamma$ if and only if $\psi = 1_Z$. 
Moreover, $\Ind_Z^G 1_Z$ factorizes to 
the regular representation $\lambda_\Gamma$ of $\Gamma$.
Since $\Gamma$ is amenable, $\textnormal{C*ker}(\lambda_\Gamma) = \{0\}$.
Therefore, the result is a consequence of Theorem~\ref{Theo-PrimIdealPGLn}. 
\end{proof}

The proof of Theorem~\ref{Theo-PrimIdealPGLn} 
carries over to $\Gamma = \SL_n(\K)$ for $n \ge 3$,
with the appropriate modifications. 
However, for $\Gamma = \SL_2(\K)$, the proof of Lemma~\ref{Lem-FreeActDualGL_n}
breaks down. Here is the result
(for more details, see the final remarks in \cite{HoRo--89}).

Recall that the centre $Z$ of $\SL_n(\K)$ consists of the matrices $\lambda I_n$
with $\lambda$ a $n$-th root of unity in $\K$, that is,
$Z \approx \mu_n^\K := \{ \lambda \in \K^* \mid \lambda^n = 1 \}$. 

\begin{theorem}
% 9.F.11
\label{Theo-PrimIdealSLn}
Let $\Gamma = \SL_n(\K)$, 
where $n \ge 3$
and $\K$ is an infinite algebraic extension of a finite field.
The map 
$$
\Phi \, \colon \, \{1_\Gamma\} \sqcup \widehat Z 
\to 
\Pri(\Gamma),
$$
defined by
$\Phi(1_\Gamma) = \textnormal{C*ker}(1_\Gamma)$ and
$$
\Phi(\psi) \, = \, \textnormal{C*ker}(\Ind_Z^\Gamma \psi)
\hskip.5cm \text{for} \hskip.2cm 
\psi \in \widehat Z \approx \widehat{\mu_n^\K}
$$
is a bijection.
\end{theorem}

A consequence of Theorem~\ref{Theo-PrimIdealSLn}
is that $\Pri(\SL_n(\K))$ is finite when $n \ge 3$
and $\K$ is an infinite algebraic extension of a finite field.
The result for $\Gamma = \PSL_n(\K)$ is even more striking. 

\begin{cor}
% 9.F.12
\label{Cor-PrimIdealSLn}
Let $\Gamma = \PSL_n(\K)$, where $n \ge 3$
and $\K$ is an infinite algebraic extension of a finite field.
The primitive dual $\Pri(\Gamma)$ consists of exactly two ideals:
$\{0\}$ and $\textnormal{C*ker}(1_\Gamma)$. 
\end{cor}

\begin{rem}
% 9.F.13
\label{Rem-PrimIdealSLn}
(1)
Let $\Gamma = \PSL_n(\K)$,
where $\K$ be an infinite algebraic extension of a finite field and $n \ge 3$.
 It is worth rephrasing Corollary~\ref{Cor-PrimIdealSLn} in two ways:
\begin{itemize}
\setlength\itemsep{0em}
\item 
In terms of weak containment: every irreducible representation of $\Gamma$ 
which is not the trivial representation $1_\Gamma$ 
is weakly equivalent to the regular representation $\lambda_\Gamma$.
\item
In terms of the ideal structure of the C*-algebra of $\Gamma$:
the only non-trivial quotient of
$C^*_{\rm max} (\Gamma) = C^*_\lambda (\Gamma)$ is 
the copy of $\C$, given by the trivial representation $1_\Gamma$. 
\end{itemize}

\vskip.2cm

(2)
Let $\Gamma$ be a group, with $\Gamma \ne \{e\}$. 
Then $\Pri(\Gamma)$ has at least two elements:
indeed, $\textnormal{C*ker}(1_\Gamma)$ is a closed point in $\Pri(\Gamma)$; 
moreover, there exists an irreducible representation 
$\pi$ of $\Gamma$ which is distinct from $1_\Gamma$ 
and then $\textnormal{C*ker}(\pi) \in \Pri(\Gamma)$ and 
$\textnormal{C*ker}(\pi) \ne \textnormal{C*ker}(1_\Gamma)$.
The groups $\PSL_n(\K)$ as in Corollary~\ref{Cor-PrimIdealSLn}
or $\SL_n(\K)$ as in Remark~\ref{Rem-Cor-PrimIdealSLn}
are the only infinite groups $\Gamma$
for which we know that $\Pri(\Gamma)$ consists of exactly two points. 

\vskip.2cm

(3)
\textbf{Question:} Does Corollary \ref{Cor-PrimIdealSLn} also holds for $n=2$ ?
\end{rem}

\begin{rem}
% 9.B.14
\label{Rem-Cor-PrimIdealSLn}
Note that there are cases for which already $\Gamma = \SL_n(\K)$
has a trivial centre for $n \ge 3$; 
in this case, the same conclusion as in Corollary~\ref{Cor-PrimIdealSLn} holds:
the primitive dual $\Pri(\Gamma)$ consists of exactly the two ideals 
$\{0\}$ and $\textnormal{C*ker}(1_\Gamma)$. 
This is true in particular for $\SL_3(\overline{\F_3})$,
since $\overline{\F_3}$ has no non-trivial cubic root of $1$,
and more generally for $\SL_p(\K)$,
whenever $p$ is an odd prime and $\K$ an infinite algebraic extension of $\F_p$.
\end{rem}

\section[Map from the dual to the primitive dual]
{On the non-injectivity of the map from the dual to the primitive dual}
% Section 9.G
\label{noninjectivity}

Let $G$ be a second-countable LC group which \emph{is not} of type I.
Then Corollary~\ref{Cor-FacRepPrimitiveG-NonTypI} shows that the map
% Then Proposition ???\footnote{Est-ce que \c ca r\'esulte de
% quelque chose qu'on a d\'ej\`a dit~?
% Sinon, le th\'eor\`eme de Glimm dit au moins que cette application
% N'EST PAS injective.}
% \marginpar{Footnote} 
$\kappa^{\rm d}_{\rm prim} \, \colon \widehat G \twoheadrightarrow \Pri(G)$
of Sections \ref{PrimIdealSpace} and \ref{C*algLCgroup}
has fibers of uncountable cardinality.
For several of the discrete groups $\Gamma$ we have met so far, 
this can be checked with specific representations, as shown by the next examples.

\begin{exe}
% 9.G.1
\label{ExRepHeisZ}

(1)
For the Heisenberg group $\Gamma = H(\Z)$ over the integers,
% as in \ref{Section-IrrRepTwoStepNil},
% and \ref{Sect-PrimIdealHeisenberg},
this was observed in Remark \ref{Rem-Cor-PrimIdealHeisIntegers-bis}.

\vskip.2cm

(2) 
For the affine group $\Aff (\K)$ of an infinite field $\K$,
% as in \ref{Section-IrrRepAff},
see Remark \ref{Rem-Theo-PrimIdealAffGr}.

\vskip.2cm

(3)
\index{Baumslag--Solitar group $\BS(1, p)$} 
Let $p$ be a prime and $\Gamma = \BS(1, p)$ be the Baumslag--Solitar group,
as in \ref{Section-IrrRepBS} and \ref{Section-PrimIdealBS}.
Proposition~\ref{Prop-IrredRepBS-bis} shows that 
there is an injective map $\Phi$
from the set of non-periodic $T_p$-orbits in the solenoid $\So_p$
to $\widehat \Gamma$, where $T_p$ is the transformation given by 
multiplication by $p$. 
There are uncountably many dense $T_p$-orbits
(see Remark~\ref{Rem-Prop-IrredBS-bis}).
By Theorem~\ref{Theo-PrimIdealBS} (see First Step of the proof),
every irreducible representation of $\Gamma$ 
corresponding to a dense $T_p$-orbit 
maps under $\Phi$ to the primitive ideal $\{0\}$. 

\vskip.2cm

(4) Let $\Gamma = \Z \rtimes \bigoplus_{k \in \Z} \Z / 2 \Z$ be the lamplighter group,
as in \ref{Section-IrrRepLamplighter} and \ref{Section-PrimIdealLamplighter}.
There are uncountably many dense $T$-orbits
for the shift map $T$ on $X = \prod_{k\in \Z} \Z / 2 \Z$.
As in the case of the Baumslag--Solitar group, 
every irreducible representation of $\Gamma$ corresponding to a dense $T$-orbit 
maps under $\Phi$ to the primitive ideal $\{0\}$,
by Theorem~\ref{Theo-PrimIdealLamplighter}.

\vskip.2cm

(5)
Let $\Gamma = \PGL_n(\K)$, 
where $\K$ is an infinite algebraic extension of a finite field.
% as in \ref{Theo-PrimIdealPGLn}.
Every principal series representation of $\GL_n(\K)$,
see Subsection~\ref{Section-IrrRepGLN},
which factorizes through $\Gamma$,
maps to the primitive ideal $\{0\}$, by Theorem~\ref{Theo-PrimIdealPGLn}. 
There are uncountably many such representations.
\end{exe}

\begin{exe}
% 9.G.2
\label{ExRepC*simple1}
Define a locally compact group $G$ to
\textbf{have a large simple C*-quotient}
if there exists a quotient of $C^*_{\rm max}(G)$
which is a simple infinite-dimensional C*-algebra with unit.
The group $H(\Z)$ of \ref{ExRepHeisZ}(1)
% of the previous example
has a large simple quotient.
Another class of LC groups having large simple C*-quotients 
is that of second-countable \textbf{C*-simple groups},
i.e., of the second-countable LC groups $G$, not $\{1 \}$, 
of which the reduced C*-algebra $C^*_\lambda(G)$ is simple \cite{Harp--07, Kenn}.
\par

The previous example generalizes as follows.
Let $G$ be a second-countable LC group.
Let $A$ be a simple infinite-dimensional C*-algebra;
assume that there exists a surjective C*-morphism
$p \,\colon \pi(C^*_{\rm max}(G)) \twoheadrightarrow A$.
Let $\widehat G_A$ denote the image of the corresponding map
from $\widehat A$ to $\widehat G$.
Then, by the map 
$\kappa^{\rm d}_{\rm prim} \,\colon \widehat G \twoheadrightarrow \Pri(G)$, 
the image of the uncountable subset $\widehat G_A$
consists of a unique point, i.e., of $\ker (p)$.
\end{exe}

\begin{exe}
% 9.G.3
\label{ExRepC*simple2}
Let $\Gamma = \PSL_2(\Z)$ 
and $H \approx \Z$ the subgroup which is the image
of $\begin{pmatrix} 1 & \Z \\ 0 & 1 \end{pmatrix}$ in $\Gamma$.
For every $\chi_\theta \in \widehat H \approx \T$, let 
$\pi_\theta$ denote the induced representation 
$\Ind_H^\Gamma \chi_\theta$. 
We claim that 
\begin{enumerate}[label=(\arabic*)]
\item\label{iDEExRepC*simple2}
the regular representation $\lambda$ of $\Gamma$ is a direct integral 
$
\lambda = \int_{\widehat H}^\oplus \pi_\theta d\theta;
$
\item\label{iiDEExRepC*simple2}
all $\pi_\theta$ are irreducible and pairwise not equivalent to each other;
\item\label{iiiDEExRepC*simple2}
all $\pi_\theta$ are weakly equivalent to each other and to $\lambda$ as well; 
in particular, all $\pi_\theta$ have the same image in $\Pri(\Gamma)$.
\end{enumerate}
\par

To show \ref{iDEExRepC*simple2}, 
observe that, by induction in stages, 
$\lambda = \Ind_H^\Gamma \lambda_H$ 
for the regular representation $\lambda_H$ of $H$.
Let $d\theta$ denote the normalized Haar measure on $\widehat H \approx \T$;
by the Plancherel theorem \ref{PlancherelTh}, 
$\lambda_H = \int_{\widehat H}^\oplus \chi_\theta d\theta$.
It follows that
$$
\lambda
\, = \, 
\Ind_H^\Gamma \left(\int_{\widehat H}^\oplus \chi_\theta\right) d\theta
\, = \,
\int_{\widehat H}^\oplus \left( \Ind_H^\Gamma \chi_\theta \right) d\theta 
\, = \,
\int_{\widehat H}^\oplus \pi_\theta d\theta.
$$
\par

To show \ref{iiDEExRepC*simple2},
it suffices to prove that $H$ is equal to its commensurator
$\Comm_\Gamma (H)$ in $\Gamma$,
by Mackey--Shoda irreducibility and equivalence criteria 
(Theorems~\ref{Theo-IrredInducedRep} and \ref{Theo-EquiInducedRep}).
\par

Let $X$ be the set of equivalence classes of primitive vectors in $\Z^2$,
where two primitive vectors $x, y$ are equivalent if $y = \pm x$; 
the group $\Gamma$ acts transitively on $X$ 
and $H$ is the stabilizer of the class $x_0$ of $(1,0)^t$. 
Moreover, $\{x_0\}$ is the unique finite orbit of $H$ in~$X$; 
hence, $H = \Comm_\Gamma (H)$ (see also \cite[Example 3]{BuHa--97}).
\index{Commensurator}
\par

Finally, since $H$ is amenable,
$\pi_\theta = \Ind_H^\Gamma \theta$ 
is weakly contained in $\lambda$ \cite[F.3.5]{BeHV--08}.
As $\Gamma$ is C*-simple, 
it follows that $\pi_\theta$ is weakly equivalent to $\lambda$
for every $\theta \in \widehat H$ and Item \ref{iiiDEExRepC*simple2} is proved.
\end{exe}

\begin{exe}
% 9.G.4
Let $G$ be a connected real Lie group that is simple, non-compact,
with centre $\{1 \}$, and let $\Gamma$ be a lattice in $G$.
Let $\pi, \pi'$ be non-equivalent irreducible representations of $G$
that are in the principal series of $G$;
denote by $\pi \vert_\Gamma, \pi' \vert_\Gamma$ their restrictions to $\Gamma$.
Then:
\par

$\pi \vert_\Gamma$ and $\pi' \vert_\Gamma$ are irreducible and non-equivalent
\cite{CoSt--91};
\par
 
$\pi \vert_\Gamma$ and $\pi' \vert_\Gamma$ are weakly equivalent to each other
\cite{BeHa--94}.
\par\noindent
(This was already cited for the particular case of $G = \SL_2(\R)$
in the proof of Proposition \ref{resPSL2}.)
In particular, there are fibres of uncountable cardinality
in the projection
$\widehat \Gamma \twoheadrightarrow \Pri(\Gamma)$.
\end{exe}

\begin{exe}
% 9.G.5
Let $\Gamma$ be a weakly branch group acting on the $d$-regular rooted tree,
for some $d \ge 2$; see \cite{Grig--11}.
\par

Assume that $\Gamma$ is subexponentially bounded.
Then the projection $\widehat \Gamma \twoheadrightarrow \Pri(\Gamma)$
has uncountable fibres.
This covers several of the basic examples of groups acting on rooted trees,
such as the first Grigorchuk group, Gupta--Sidki $p$-groups, and the Basilica group.
See \cite{DuGr--17a}, and Corollary 1 in \cite{DuGr--17b}.
\end{exe}

\section
{Borel comparison for duals of groups}
% Section 9.H
\label{SectionEffrosThomas}

It is a natural question to ask what are the possible spaces
that can arise as duals of groups. 
We address now this question for duals of countable discrete groups,
viewed as Borel spaces,
following part of \cite{Thos--15}.

\begin{prop}
% 9.H.1
\label{GpeDenDualNonden}
The dual of a countable infinite discrete group
is un uncountable Borel space.
\end{prop}

\begin{proof}
Let $\Gamma$ be a countable discrete group.
By Theorem \ref{thmDirectIntIrreps}, there exist
a standard Borel space $X$,
a $\sigma$-finite measure $\mu$ on $X$,
and a measurable field $x \mapsto \pi_x$ of irreducible representations of $G$ over $X$,
such that the left-regular representation $\lambda_\Gamma$
is a direct integral $\int^\oplus_X \pi_x d\mu(x)$.
\par

If the dual $\widehat \Gamma$ is countable,
% ??? \footnotesize some argument missing \normalsize ???
% \marginpar{Argument missing}
this direct integral is a direct sum,
so that $\lambda_\Gamma$ is a direct sum of irreducible representations.
It follows from Proposition \ref{Pro-NoIrredSubRegRep} that $\Gamma$ is finite.
\par

By contraposition, if $\Gamma$ is infinite, $\widehat \Gamma$ is uncountable.
% For the details, we refer to \cite[1.5]{Bagg--72}.
\end{proof}

The easiest case to deal with is that of the abelian-by-finite groups.
As a consequence of Proposition \ref{GpeDenDualNonden}
and Theorem \ref{ThGlimm}~\ref{fDEThGlimm}, we have:

\begin{cor}
% 9.H.2
\label{boreldualZ}
Let $\Gamma$ be a countable infinite abelian-by-finite group.
\par

The dual $\widehat \Gamma$ is an uncountable standard Borel space,
and in particular is Borel isomorphic to the dual $\widehat \Z = \T$.
\end{cor}

\index{Borel equivalence relation}
\index{Borel equivalence relation! smooth}
Let $X$ be a standard Borel space.
An equivalence relation $\mathcal R$ on $X$ is \textbf{Borel}
if $R$ is a Borel subset of $X \times X$,
and \textbf{smooth} 
if there exists a Borel map $\varphi \,\colon X \to \R$
such that, for $(x,x') \in X \times X$, 
we have $(x,x') \in \mathcal R$ if and only if $\varphi (x) = \varphi (x')$.
\par

Consider a countable infinite group $\Gamma$.
Let $\Hi$ be an infinite-dimensional separable Hilbert space,
and $\U(\Hi)$ its unitary group, with the strong topology; 
recall that $\U(\Hi)$ is a Polish group.
Denote by ${\rm Irr}_\infty(\Gamma)$
the space of irreducible representations 
of $\Gamma$ in $\Hi$, with its Polish topology defined here via
the natural embedding of ${\rm Irr}_\infty(\Gamma)$
into the product space $\U(\Hi)^\Gamma$.
(This provides the same topology as 
that defined Section \ref{SectionQuasidual} in a more general situation.)
In ${\rm Irr}_\infty(\Gamma)$,
equivalence of representations
is a Borel equivalence relation;
we denote it here by $\simeq_\Gamma$.
Observe that the quotient space
${\rm Irr}_\infty(\Gamma) / \simeq_\Gamma$
is the complement of the finite-dimensional part in the dual $\widehat \Gamma$,
viewed here with its Mackey--Borel structure.
\index{$a1$@$\simeq$ equivalence of representations!
occasionally: $\simeq_\Gamma$}
\par

For the other cases, we know from \cite{Glim--61b} and \cite{Thom--64a}
that the following conditions are equivalent:
\begin{enumerate}[label=(\roman*)]
\item\label{iDEbla}
$\Gamma$ is not abelian-by-finite,
\item\label{iiDEbla}
$\Gamma$ has an infinite-dimensional irreducible representation,
i.e., ${\rm Irr}_\infty(\Gamma) \ne \emptyset$,
\item\label{iiiDEbla}
the equivalence relation $\simeq_\Gamma$ on the space 
${\rm Irr}_\infty(\Gamma)$ is not smooth,
\item\label{ivDEbla}
the dual $\widehat \Gamma$ is not countably separated.
\end{enumerate}
\par

\index{Borel equivalence relation! Borel reducible}
Let $\mathcal R, \mathcal S$ be two Borel equivalence relations 
on Polish spaces $X, Y$ respectively.
Then $\mathcal R$ is \textbf{Borel reducible to} $\mathcal S$,
and this is denoted by $\mathcal R \le_B \mathcal S$, 
if there exists a Borel map $\varphi \,\colon X \to Y$
such that $(x,x') \in \mathcal R$ if and only if $(\varphi(x), \varphi(x')) \in \mathcal S$;
intuitively, this means that $\mathcal S$ is not less complicated than $\mathcal R$.
For example, a Borel equivalence relation $\mathcal R$ is smooth
if and only if it is Borel reducible to the equality equivalence relation
$\{ (x,y) \in \R \times \R \mid x=y \}$ on $\R$.
The relation $\mathcal R$ is \textbf{Borel equivalent to} $\mathcal S$,
denoted by $\mathcal R \hskip.1cm\sim_B \hskip.1cm \mathcal S$, 
if $\mathcal R \le_B \mathcal S$ and $\mathcal S \le_B \mathcal R$.
\index{Borel equivalence relation! Borel equivalent}
\par

Consider the equivalence relation $E_0$ on $2^\N$
defined for $(\omega, \omega') \in 2^\N \times 2^\N$ by
$(\omega, \omega') \in E_0$ if $\omega(n) = \omega'(n)$ for all $n \in \N$ large enough;
this relation $E_0$ is not smooth.
For a Borel equivalence relation $\mathcal R$ on a Polish space $X$,
it is a fundamental result of \cite{HaKL--90},
building up on previous results by Glimm \cite{Glim--61a, Glim--61b} and Effros, that
\begin{enumerate}[label=(\arabic*)]
% \item\label{iDEHaKL}
\item[--]
either $\mathcal R$ is smooth,
% \item\label{iiDEHaKL}
\item[--]
or $E_0 \le_B \mathcal R$.
\end{enumerate}
In particular, $E_0 \le_B \hskip.1cm \simeq_\Gamma$
for any countable group $\Gamma$ that is not abelian-by-finite.

\vskip.2cm

Let $\Gamma, \Delta$ be two countable groups which are not abelian-by-finite.
Then $\simeq_\Gamma \sim_B \simeq_\Delta$
if and only if the duals $\widehat \Gamma$ and $\widehat \Delta$
are Borel isomorphic
\cite[Corollary 1.7]{Thos--15}.
\par

Let $F_\infty$ be the free group on a countable infinite set of generators.
For any countable group $\Gamma$, it is a simple observation that
$\simeq_\Gamma \le_B \simeq_{F_\infty}$;
indeed, if $p \,\colon F_\infty \to \Gamma$ is a surjective homomorphism,
then the natural map $\pi \mapsto \pi \circ p$
from ${\rm Irr}_\infty(\Gamma)$ to ${\rm Irr}_\infty(F_\infty)$
is a Borel reduction from $\simeq_\Gamma$ to $\simeq_{F_\infty}$.
Let $F_2$ denote the free group of rank $2$;
then $\simeq_{F_2} \hskip.1cm \sim_B \hskip.1cm \simeq_{F_\infty}$ \cite{Thos--15}.
\par

Denote by $S_3^{(\N)}$ the restricted product of a countable infinite number
of copies of the non-abelian finite group of order $6$;
this is a countable infinite amenable group that is not abelian-by-finite.
Building on results by Elliot and Sutherland, 
Thomas has also shown:

\begin{theorem}
% 9.H.3
\label{boreldualnext}
For any countable infinite amenable group $\Gamma$ that is not abelian-by-finite,
we have
$$
\simeq_\Gamma \, \hskip.1cm \sim_B \hskip.1cm \,
\simeq_{S_3^{(\N)}} .
% \, \hskip.1cm \sim_B \hskip.1cm \, E_0. NON C'EST FAUX !!!
$$
Equivalently, 
the dual $\widehat \Gamma$ is Borel isomorphic to $\widehat{S_3^{(\N)}}$.
% and also to $2^\N / E_0$. NON C'EST FAUX !!!
\end{theorem}

For groups which are not abelian-by-finite, the following question is open:

\begin{ques}
% 9.H.4
\label{DistinctDuals?}
Does there exist a countable infinite group $\Gamma$
which is not abelian-by-finite and such that 
$\simeq_\Gamma \hskip.1cm \nsim_B \hskip.1cm \simeq_{F_\infty}$ ?
% Dans d'autres termes, voir \cite[Corollary 1.7]{Thos--15}:
Equivalently which is not abelian-by-finite and such that
$\widehat \Gamma$ is not Borel isomorphic to $\widehat{F_\infty}$ ?
\end{ques}

%-----------------------------------------------------------------------
% End of chapter 9
%-----------------------------------------------------------------------

\chapter{Normal quasi-dual and characters}
% Chapter 10
\label{ChapterCharacters} 

\emph{
Let $G$ be a second-countable LC group.
In general, the primitive dual $\Pri(G)$ is a far more accessible object than $\widehat G$. 
Indeed,
$\Pri(G)$ is a standard Borel space (Theorem \ref{EffrosTheoremG}),
whereas the dual $\widehat G$ is standard
if and only if $G$ is of type I (Theorem~\ref{ThGlimm}).
}
\par

\emph{
With the ultimate goal to parametrize the primitive dual $\Pri(G)$,
we introduce the notions of a traceable representation
(Section~\ref{Section-TraceRep})
and of a normal factor representation of $G$
(Section~\ref{Section-NormRepChar}).
The normal quasi-dual $\QD(G)_{\rm norm}$ of $G$,
which is defined as the set of quasi-equivalence classes of normal factor representations,
is a standard Borel subspace of the quasi-dual $\QD(G)$; see Theorem~\ref{ThHalpern}.
We will see in \ref{Section-CharactersPrimitive} that, for some classes of groups $G$, 
the primitive dual $\Pri(G)$ is naturally in bijection with $\QD(G)_{\rm norm}$.
}

\emph{
Normal factor representations of $G$ can be described in terms of characters.
Characters of $G$ are (usually not everywhere defined) lower semi-continuous traces
on the maximal C*-algebra $C^*_{\rm max}(G)$ of $G$
with an appropriate extremality property.
We will show in Section~\ref{Section-NormRepChar}
that there is a natural bijective correspondence $\QD(G)_{\rm norm} \to \Char(G)$
between $\QD(G)_{\rm norm}$ and the space $\Char(G)$
of equivalence classes of characters of $G$,
where two characters are equivalent if they are proportional to each other.
}

\emph{
A detailed account on Hilbert algebras is given in Section~\ref{Section-RepTraceC*}.
Such an algebra is associated to every lower semi-continuous trace on a C*-algebra
(Proposition~\ref{Prop-TraceHilbertAlg}),
a fact which is crucial in order to show the surjectivity of the map
$\QD(G)_{\rm norm} \to \Char(G)$.
Hilbert algebras appear also
in the construction of the standard representation of a semi-finite von Neumann algebra (Section~\ref{Section-StandardRep})
as well as in Chapter~\ref{Thomadual} about traces on general topological groups
(Section~\ref{Section-GNS-Traces}).
}

\section
{Traces and Hilbert algebras}
% Section 10.A
\label{Section-RepTraceC*}

\index{Trace! $1$@on a C*-algebra}
\index{Trace! $4$@ideal of definition}
\index{Ideal of definition of a trace}
Let $A$ be a C*-algebra 
and $t \,\colon A_+ \to \mathopen[ 0, \infty \mathclose]$ a trace on $A$
(as defined in \ref{SectionvN}).
The set
$$
\mathfrak n_t \, := \, \{ x \in A \mid t(x^*x) < \infty \}
$$
is a selfadjoint two-sided ideal of $A$.
The set $\mathfrak m_t$ of linear combinations of elements in
$\mathfrak m_{t,+} := \{x \in A_+ \mid t(x) < \infty \}$
is a selfadjoint two-sided ideal of $A$
which coincides with the ideal $\mathfrak n_t^2$, i.e.,
with the set of elements of the form
$x_1y_1 + \cdots + x_ky_k$ for $x_1, y_1, \hdots, x_k, y_k \in \mathfrak n_t$.
The two-sided ideal $\mathfrak m_t$ is the \textbf{ideal of definition} of $t$. 
The trace $t$ extends from $\mathfrak m_{t,+}$
to a linear form $\mathfrak m_t \to \C$,
denoted by $t$ again,
and
$$
\begin{aligned}
t(x^*)&\, = \, \overline{t(x)}\phantom{y} 
\hskip.5cm \text{for all} \hskip.2cm 
x \in \mathfrak m_t
\\
t(xy) &\, = \, t(yx) 
\hskip.5cm \text{for all} \hskip.2cm 
x \in \mathfrak m_t
\hskip.2cm \text{and} \hskip.2cm 
y \in A
\\
t(xy) &\, = \, t(yx) 
\hskip.5cm \text{for all} \hskip.2cm 
x,y \in \mathfrak n_t.
\end{aligned}
$$
(See \cite[6.1.1]{Dixm--C*}.)
\par

Let $A$ be a C*-algebra.
There is a construction that associates
to every lower semi-continuous semi-finite trace $t$ on $A$
a representation $\pi_t$ of $A$ such that the generated von Neumann algebra
$\pi_t(A)''$ is semi-finite; see Proposition \ref{Prop-TraceHilbertAlg}.
The construction involves basic objects that are Hilbert algebras.
We give a complete account on two main facts concerning these algebras,
i.e., the Commutation Theorem and the existence of a canonical trace, 
following \cite[Chap.~I, \S~5]{Dixm--vN};
see Theorems \ref{TheoHilbertAlg-Commutation} and \ref{TheoHilbertAlg-Trace}.

\vskip.2cm

\index{Hilbert algebra}
\index{Algebras! $4$@Hilbert algebra} 
A \textbf{Hilbert algebra}
is a complex $*$-algebra $\AC$,
equipped
with a scalar product $\langle \cdot \mid \cdot \rangle$,
with the following properties:
\begin{enumerate}[label=(\arabic*)]
\item 
$\langle x \mid y \rangle = \langle y^* \mid x^* \rangle$ for all $x,y \in \AC$;
\item 
$\langle xy \mid z \rangle = \langle y \mid x^*z \rangle$ for all $x,y,z \in \AC$;
\item 
for every $x \in \AC$, the map $\AC \to \AC, \, y \mapsto xy$ is continuous;
\item 
the set $\{xy \mid x,y \in \AC\}$ is total in $\AC$.
\end{enumerate}
Let $\Hi$ be the Hilbert space completion of $\AC$.
It follows from (1) to (4) that:
\begin{enumerate}[label=(\arabic*)]
\item[(5)]
the map $x \mapsto x^*$ from $\AC$ to $\AC$ 
extends to an involutive antilinear isometry $J \,\colon \Hi \to \Hi$;
\item[(6)]
for every $x \in \AC$, the maps defined on $\AC$ 
by $y \mapsto xy$ and $y \mapsto yx$ 
extend to bounded linear operators $\lambda(x)$ and $\rho(x)$ on $\Hi$,
and $x \mapsto \lambda(x)$ [respectively $x \mapsto \rho(x)$] 
is a non-degenerate $*$-representation of the $*$-algebra $\AC$ 
[respectively of the opposite algebra of $\AC$];
\item[(7)]
$J\lambda(x)J = \rho(x^*)$,
and $J\rho(x)J = \lambda(x^*)$,
and $\lambda(x) \rho(y) = \rho(y) \lambda(x)$ for every $x, y \in \AC$.
\end{enumerate}
Let $\lambda(\AC)''$ and $\rho(\AC)''$
be the von Neumann subalgebras of $\Li (\Hi)$ 
generated by $\{\lambda(x) \mid x \in \AC\}$ and $\{\rho(x) \mid x \in \AC\}$ respectively. 
It is clear that $\lambda(\AC)'' \subset \rho(\AC)'$ and $\rho(\AC)'' \subset \lambda(\AC)'$.
Theorem~\ref{TheoHilbertAlg-Commutation} below shows that 
these von Neumann algebras are in fact each other's commutants.
The proof is based on properties of the so-called bounded elements 
in the Hilbert space $\Hi$, that we introduce now.
\par

\index{Bounded element of a Hilbert algebra} 
\index{Hilbert algebra! bounded element}
An element $a \in \Hi$ is said to be \textbf{bounded} 
if the map $\AC \to \Hi, \hskip.1cm x \mapsto \rho(x) a$
extends from $\AC$ to a bounded operator on $\Hi$;
when this is the case,
this extension is a uniquely determined operator, denoted by $\lambda(a)$. 
%in this case, the map $x \mapsto \lambda(x)a$ also extends from $\AC$ 
% to a bounded operator on $\Hi$, denoted $\rho(a)$.
% \par
Denote by $\AC_b$ the set of bounded elements in $\Hi$; 
it is clear that $\AC_b$ is a linear subspace of $\Hi$.
Observe that $\AC \subset \AC_b$ 
and that the map $\lambda \,\colon \AC_b \to \Li (\Hi)$ coincides on $\AC$ 
with the representation $\lambda$ of $\AC$ introduced above.
More on $\AC_b$ in Remark \ref{Rem-HilbertAlgBoundedEl}.
 
\begin{lem}
% 10.A.1
\label{LemHilbertAlg}
Let $\AC$ be a Hilbert algebra,
and let $\Hi, J, \lambda, \rho, \AC_b$ be as above.
\begin{enumerate}[label=(\arabic*)]
\item\label{iDELemHilbertAlg}
The identity operator $i \in \Li (\Hi)$ is in the closure 
of both the $*$-algebras $\lambda(\AC)$ and $\rho(\AC)$, 
for the strong operator topology;
\item\label{iiDELemHilbertAlg}
the linear map $a \mapsto \lambda(a)$ is injective on $\AC_b$;
\item\label{iiiDELemHilbertAlg}
for $a \in \AC_b$, we have $Ja \in \AC_b$ and $\lambda(Ja) = \lambda(a)^*$;
\item\label{ivDELemHilbertAlg} 
for $a \in \AC_b$, there exists a unique bounded operator on $\Hi$, denoted $\rho(a)$, 
which extends the map
$x \mapsto \lambda(x) a$ from $\AC$ to $\Hi$; 
moreover, $\rho(a) = J\lambda(a)^*J$;
\item\label{vDELemHilbertAlg}
the set $\lambda(\AC_b)$ [respectively $\rho(\AC_b)$] 
is a two-sided ideal of $\rho(\AC)'$ [respectively $\lambda(\AC)'$]; 
more precisely, for $a \in \AC_b$ and $T \in \rho(\AC)'$ [respectively $\lambda(\AC)'$] 
we have $Ta \in \AC_b$, $JT^*Ja \in \AC_b$, and
$$
T\lambda(a) \, = \, \lambda(Ta)
\hskip.5cm \text{and} \hskip.5cm 
\lambda(a)T \, = \, \lambda(JT^*Ja)
$$
[respectively, $T\rho(a) = \rho(Ta)$ and $\rho(a)T = \rho(JT^*Ja)$].
\end{enumerate}
\end{lem}

\begin{proof}
\ref{iDELemHilbertAlg}
By Property (4) in the definition of a Hilbert algebra, 
the linear span of $\lambda(\AC)\Hi$ is dense in $\Hi$. 
By von Neumann's bicommutant theorem,
$\lambda(\AC)''$ is the closure of $\lambda(\AC)$ for the strong operator topology 
(see Corollaire 1 of Th\'eor\`eme 2 in \cite[Chap.~I, \S~3]{Dixm--vN});
since $I \in \lambda(\AC)''$, this proves the claim for $\lambda(\AC)$. 
The proof for $\rho(\AC)$
is similar.

\vskip.2cm

\ref{iiDELemHilbertAlg}
Let $a \in \AC_b$ be such that $\lambda(a) = 0$. 
Then $\rho(x)a = 0$ for all $x \in \AC$
and it follows from \ref{iDELemHilbertAlg} that $a = 0$.
 
\vskip.2cm
 
 \ref{iiiDELemHilbertAlg}
 For $a \in \AC_b$ and $x,y \in \AC$, 
 %using Property (1) in the definition of a Hilbert algebra 
 we have
 $$
\begin{aligned}
 \langle \lambda(a) x \mid y \rangle
 &\, = \, \langle \rho(x) a \mid y \rangle \, = \, \langle a \mid \rho(x)^* y \rangle
 \\
 &\, = \, \langle a \mid yx^* \rangle \, = \, \langle xy^* \mid Ja \rangle
 \\
 &\, = \, \langle \rho(y^*)x \mid Ja \rangle,
 \end{aligned}
$$
and hence $\lambda(a)^*y = \rho(y)Ja$, 
showing that $Ja \in \AC_b$ and $\lambda(Ja) = \lambda(a)^*$.

\vskip.2cm

\ref{ivDELemHilbertAlg}
Using \ref{iiiDELemHilbertAlg}, we have, for $a \in \AC_b$ and $x \in \AC$,
$$
\lambda(x) a \, = \, J\rho(x^*) J a \, = \, J \lambda(a)^*x^* \, = \, J\lambda(a)^*J x,
$$
showing that $J\lambda(a)^*J$ coincides with $x \mapsto \lambda(x)a$ on $\AC$.
 
\vskip.2cm
 
 \ref{vDELemHilbertAlg} 
For $a \in \AC_b$ and $x,y \in \AC$, we have
$$
\lambda(a)\rho(y)x \, = \,
\lambda(a) (yx) \, = \, \rho(yx) a \, = \, 
\rho(y)\rho(x)a \, = \, \rho(y) \lambda(a) x,
$$
showing that $\lambda(\AC_b) \subset \rho(\AC)'$.
\par 

Let $T \in \rho(\AC)'$, $a \in \AC_b$. For $x \in \AC$, we have 
$$
\rho(x) Ta \, = \, T \rho(x) a \, = \, T\lambda(a)x;
$$
this shows that $Ta \in \AC_b$ and that $\lambda(Ta) = T\lambda(a)$.
Since $T^*\in \rho(\AC)'$, using \ref{iiiDELemHilbertAlg} twice,
we have therefore
$$
\lambda(a)T \, = \, (T^* \lambda(a)^*)^* \, = \, 
(T^* \lambda(Ja))^* \, = \, \lambda(T^*Ja)^* \, = \, \lambda(JT^*Ja).
$$
% So , $\lambda(\AC_b)$ is a left ideal of $\rho(\AC)'$; since $\lambda(AC')$ is 
% selfadjoint, it follows that $\lambda(AC')$ is a two-sided ideal.
The claim for $\rho(\AC_b)$ is proved in a similar way.
\end{proof}

We are ready to prove the Commutation Theorem for a Hilbert algebra.
It is also Theorem 1 in \cite[Chap.~I, \S~5, no~2]{Dixm--vN}.
 
\begin{theorem}[\textbf{Commutation Theorem}]
% 10.A.2
\label{TheoHilbertAlg-Commutation}

\index{Commutation Theorem for Hilbert algebras}
Let $\AC$ be a Hilbert algebra, $\Hi$ its completion,
and $\lambda(\AC)''$ and $\rho(\AC)''$ the von Neumann subalgebras of $\Li (\Hi)$ 
generated by the left and right multiplication by elements from $\AC$.
We have 
$$
\lambda(\AC)' \, = \, \rho(\AC)''
\hskip.5cm \text{and} \hskip.5cm
\rho(\AC)' \, = \, \lambda(\AC)''.
$$
\end{theorem}

\begin{proof}
$\bullet$ {\it First step.} 
We claim that $\lambda(\AC_b)$ is dense in $\rho(\AC)'$ 
and that $\rho(\AC_b)$ is dense in $\lambda(\AC)'$
for the strong operator topology. 
\par

Indeed, let $T \in \rho(\AC)'$.
Since $I$ is in the strong closure of $\lambda(\AC)'$ 
by Lemma \ref{LemHilbertAlg}~\ref{iDELemHilbertAlg},
there exists in $\lambda(\AC)'$ a sequence $(T_n)_{n \ge 1}$ 
converging strongly to $I$;
hence the sequence $(TT_n)_{n \ge 1}$ converges strongly to $T$.
By Lemma \ref{LemHilbertAlg}~\ref{vDELemHilbertAlg},
$TT_n \in \lambda(\AC_b)$ for all $n \ge 1$,
and it follows that $T$ is in the strong closure of $\lambda(\AC_b)$.
\par

The claim for $\rho(\AC_b)$ is proved similarly.

\vskip.2cm

$\bullet$ {\it Second step.} 
We claim that, for $a, b \in \AC_b$, we have
$\lambda(a)b = \rho(b)a$. Indeed, for $x \in \AC$, we have using 
Lemma~\ref{LemHilbertAlg}, \ref{iiiDELemHilbertAlg} and \ref{ivDELemHilbertAlg},
$$
\begin{aligned}
\langle \lambda(a)b \mid x \rangle
&\, = \, \langle b \mid \lambda(a)^* x \rangle
\, = \, \langle b \mid \lambda(Ja)x \rangle
\, = \, \langle b \mid \rho(x) Ja \rangle
\\
&\, = \, \langle \rho(x^*) b \mid Ja \rangle
\, = \, \langle J\lambda(x)J b \mid Ja \rangle 
\, = \, \langle a \mid \lambda(x) Jb \rangle
\\
&\, = \, \langle a \mid \lambda(x) b^* \rangle
\, = \, \langle a \mid \rho(b)^*x \rangle 
\, = \, \langle \rho(b) a \mid x \rangle.
\end{aligned}
$$

\vskip.2cm

$\bullet$ {\it Third step.} 
We claim that $\lambda(\AC)' = \rho(\AC)''$ and $\rho(\AC)' = \lambda(\AC)''$.
 
As mentioned earlier, we have $\rho(\AC)'' \subset \lambda(\AC)'$.
In order to show that $\lambda(\AC)' \subset \rho(\AC)''$, 
we have to prove that $ST = TS$ for $S \in \lambda(\AC)'$ and $T \in \rho(\AC)'$.
By the first step, it suffices to show that $ST=TS$ for $S = \rho(a)$ and $T = \lambda(b)$,
with $a,b \in \AC_b$.
Using Lemma~\ref{LemHilbertAlg}~\ref{vDELemHilbertAlg}
and the second step, we have for every $x \in \AC$,
$$
ST x \, = \, \rho(a) \lambda(b) x \, = \, \rho(a) \rho(x) b 
\, = \, \rho (\rho(a)x) b \, = \, \lambda(b) \rho(a) x \, = \, TS x.
$$
The proof that $\rho(\AC)' = \lambda(\AC)''$ is similar.
\end{proof}

\begin{rem}
% 10.A.3
\label{Rem-HilbertAlgBoundedEl}
The set $\AC_b$ of bounded elements of a Hilbert algebra $\AC$ 
is a vector subspace of $\AC$; 
moreover, for $a,b \in \AC_b$, 
we have $\lambda(a)b \in \AC_b$ and $Ja \in \AC_b$
(see Lemma~\ref{LemHilbertAlg}). 
We define maps
$$
\AC_b \times \AC_b \to \AC_b, \hskip.2cm (a,b) \mapsto ab 
\hskip.5cm \text{and} \hskip.5cm
\AC_b \to \AC_b, \hskip.2cm a \mapsto a^*
$$
by 
$ab := \lambda(a)b$ and $a^* := Ja$.
These maps extend the multiplication law
and the involution of $\AC$ and turn $\AC_b$ into a complex $*$-algebra.
Indeed, $(a,b) \mapsto ab$ is clearly bilinear 
and $a\mapsto a^*$ is an anti-linear involution;
moreover, we have $\lambda(a)b = \rho(b)a$ for $a,b \in \AC_b$
(see the second step of the proof of 
Theorem~\ref{TheoHilbertAlg-Commutation}).
Therefore for $a,b,c \in \AC_b$, we have 
$$
a(bc) \, = \, \lambda(a) (\rho(c)b) \, = \, \rho(c) (\lambda(a)b) \, = \, (ab)c
$$
and 
$$
a^*b^* \, = \, \lambda(Ja) (Jb) \, = \, J(J\lambda(Ja)J) (b)
\, = \, J \rho(a) (b) \, = \, J(ba) \, = \, (ba)^*.
$$
One checks that $\AC_b$, endowed with the scalar product 
induced from the Hilbert space completion $\Hi$ of $\AC$,
is a Hilbert algebra which contains $\AC$ as a $*$-subalgebra.
Every element in $\Hi$ which is bounded with respect to the Hilbert algebra $\AC_b$ 
is bounded with respect to the Hilbert algebra $\AC$ 
and therefore belongs to $\AC_b$, that is, 
we have $(\AC_b)_b = \AC_b$.
A Hilbert algebra with this property is called a 
\textbf{full Hilbert algebra}.
\index{Full Hilbert algebra}
\index{Hilbert algebra! full}
\end{rem}

The next theorem establishes that
the von Neumann algebra $\lambda(\AC)''$ 
associated to a Hilbert algebra $\AC$ carries a canonical trace,
called the \textbf{natural trace} associated to $\AC$;
and similarly for $\rho(\AC)''$.
for a proof, we will need the following two elementary lemmas.
\index{Hilbert algebra! natural trace}
\index{Trace! $5$@natural trace on the von Neumann algebra
associated to a Hilbert algebra}
 
\begin{lem}
% 10.A.4
\label{Lem-OperVN}
Let $\M\subset \Li (\Hi)$ be a von Neumann algebra and $S,T \in \M_+$ with $S \le T$.
There exists a unique operator $X \in \M$ with $X \vert_{T(\Hi)^{\perp}} = 0$ 
such that $S^{1/2} = X T^{1/2}$. 
Moreover, we have $\Vert X \Vert \le 1$.
\end{lem}

\begin{proof}
For every $\xi \in \Hi$, we have
$$
\Vert S^{1/2}\xi\Vert^2 
\, = \, \langle S \xi \mid \xi \rangle 
\, \le \, \langle T \xi \mid \xi \rangle
\, = \, \Vert T^{1/2}\xi\Vert^2;
$$
hence the map $T^{1/2}(\Hi) \to \Hi, \hskip.2cm T^{1/2}\xi \to S^{1/2}\xi$, 
is well defined and continuous, of norm at most $1$.
It has a linear continuous extension 
$Y \,\colon \overline{T^{1/2}(\Hi)} \to \Hi$.
Then $S^{1/2} = Y T^{1/2}$.
\par

Observe that we have $\ker T= \ker T^{1/2}$,
and therefore
$$
\overline{T(\Hi)} \, = \, (\ker T)^\perp \, = \, (\ker T^{1/2})^\perp \, = \, \overline{T^{1/2}(\Hi)}.
$$
Define $X \in \Li (\Hi)$ by $X = Y$ on $\overline{T(\Hi)}$
and $X = 0$ on $T(\Hi)^{\perp}$. 
Then $S^{1/2} = XT^{1/2}$.
\par

The uniqueness of $X \in \Li (\Hi)$ with the previous properties is clear; 
it follows that $U X U^{-1} = X$ for every unitary operator $U \in \M'$. 
Therefore $X \in \M'' = \M$.
Since $\Vert Y \Vert \le 1$, it is clear that $\Vert X \Vert \le 1$.
\end{proof} 

\begin{lem}
% 10.A.5
\label{Lem-SquareRoot}
Let $(T_i)_i$ in $\Li (\Hi)_+$ be a net of positive operators 
converging to $T \in \Li (\Hi)_+$ in the strong operator topology.
Assume that $\sup_{i} \Vert T_i\Vert \le \Vert T \Vert$.
\par
Then $(T_i^{1/2})_i$ converges to $T^{1/2}$ in the strong operator topology.
\end{lem}

\begin{proof}
We claim first that$(P(T_i))_i$ converges to $P(T)$ in the strong operator topology
for every polynomial $P \in \C[X]$.
To show this, it suffices to prove that $(T_i^n)_i$ converges to $T^n$ 
for every integer $n \ge 0$.
\par

We proceed by induction on $n$; 
the case $n = 0$ is trivial and the case$n=1$ is true by assumption. 
Assume that the claim is true for $n \ge 1$. 
Let $\xi \in \Hi$. Then $\lim_iT_i^n \xi = T^n\xi$
and $\lim_i T_i (T^n \xi ) = T(T^n \xi )$.
Since
$$
\begin{aligned}
\Vert T_i^{n+1} \xi - T^{n+1}\xi \Vert
&\le \, \Vert T_i(T_i^n \xi -T^n \xi )\Vert + \Vert T_i (T^n \xi )- T(T^n \xi )\Vert
\\
&\le \, \Vert T \Vert \Vert T_i^n \xi -T^n \xi \Vert + \Vert T_i(T^n \xi )- T(T^n \xi )\Vert,
\end{aligned}
$$
the claim is true for $n+1$.
\par

Observe that, for every $S \in \Li (\Hi)_+$ with $\Vert S \Vert \le \Vert T \Vert$,
the C*-algebra homomorphism
$$
C(\mathopen[0, \Vert T \Vert \mathclose]) \to \Li (\Hi),
\hskip.5cm
f \mapsto f(S)
$$
defined by functional calculus \cite[\S~1.5]{Dixm--C*} is norm decreasing.
\par
 
Fix $\varepsilon > 0$. Let $P \in \C[X]$ be such that 
$\sup_{t \in \mathopen[0, \Vert T \Vert \mathclose]} \vert P(t) - t^{1/2} \vert \le \varepsilon$; 
then 
$$
\Vert P(T) - T^{1/2} \Vert \, \le \, \varepsilon
\hskip.5cm \text{and} \hskip.5cm 
\Vert P(T_i)-T_i^{1/2}\Vert \, \le \, \varepsilon
\hskip.5cm \text{for all} \hskip.2cm
i.
$$
Let $\xi \in \Hi$. Since $\lim_i P(T_i) \xi = P(T) \xi$, 
we have $\Vert T_i^{1/2}\xi -T^{1/2} \xi \Vert \le 3\varepsilon$,
for $i$ large enough.
\end{proof}

The following theorem is also Theorem 1 in \cite[Chap.~I, \S~6, no~2]{Dixm--vN}.

\begin{theorem}[\textbf{Natural trace associated to a Hilbert algebra}]
% 10.A.6
\label{TheoHilbertAlg-Trace}

\index{Natural trace of a Hilbert algebra}
Let $\AC$ be a Hilbert algebra, $\Hi$ its completion, 
$\AC_b$ the set of its bounded elements, and 
$\lambda(\AC)''$ and $\rho(\AC)''$ the von Neumann subalgebras of $\Li (\Hi)$ 
generated by the left and right multiplications by elements of $\AC$.
For $T \in \lambda(\AC)''_+$, set 
$$
\begin{aligned}
t(T)
\, &= \, \langle a \mid a \rangle 
\hskip.2cm \text{if} \hskip.2cm 
 T^{1/2} = \lambda(a) 
\hskip.2cm \text{for some} \hskip.2cm 
a \in \AC_b ,
\\
t(T)
\, &= \, +\infty \,
\hskip.5cm \text{otherwise}.
\end{aligned}
$$
\par

Then $t$ is a semi-finite faithful normal trace on $\lambda(\AC)''$. 
Moreover, we have
$$
\mathfrak n_t \, = \, \{\lambda(a) \mid a \in \AC_b\}.
$$
% and (denoting the extension of $t$ to $\mathfrak m_t = \mathfrak m_t^2$
% again by $t$)
% $$
% t(\lambda(b)^*\lambda(a)) = \langle a\mid b \rangle 
% \hskip.5cm \text{for all} \hskip.2cm a,b \in \AC_b.
% $$
A similar trace is defined on $\rho(\AC)'' = \lambda(\AC)'$.
\par

In particular, the von Neumann algebras 
$\lambda(\AC)''$ and $\rho(\AC)''$ are semi-finite.
\end{theorem}

\begin{proof}
$\bullet$ {\it First step.} 
It is clear that $t$ takes its values in $\mathopen[ 0, \infty \mathclose]$
and that $t(\lambda T) = \lambda t(T)$ 
for all $T \in \lambda(\AC)''_+$ and $\lambda \ge 0$.
 
 \vskip.2cm

$\bullet$ {\it Second step.} 
Let $S,T \in \lambda(\AC)''_+$ be such $S \le T$. We claim that
$t(S) \le t(T)$. 
\par

Indeed, since this is obvious if $t(T) = +\infty$, 
we can assume that $T^{1/2} = \lambda(a)$ for some $a \in \AC_b$.
By Lemma~\ref{Lem-OperVN}, there exist $A \in \lambda(\AC)''$ 
with $\Vert A\Vert \le 1$ such that $S^{1/2}= A T^{1/2}$. 
It follows from Lemma~\ref{LemHilbertAlg}~\ref{vDELemHilbertAlg}
that $Aa \in \AC_b$ and $S^{1/2} = \lambda(Aa)$. 
We then have 
$$
t(S) \, = \, \Vert Aa\Vert^2 \, \le \, \Vert a \Vert^2 \, = \, t(T).
$$
\vskip.2cm

$\bullet$ {\it Third step.} 
Let $T_1,T_2 \in \lambda(\AC)''_+$; we claim that $t(T_1+T_2) = t(T_1)+t(T_2)$.
\par

Indeed, set $T := T_1+T_2$.
By Lemma~\ref{Lem-OperVN}, for each $j \in \{1, 2 \}$, 
there exists $A_j \in \lambda(\AC)''$ such that 
$T_j^{1/2}= A_j T^{1/2}$ and such that $A_j = 0$ on $T(\Hi)^{\perp}$. 
Set $S := A_1^* A_1 + A_2^* A_2$. We have
$$
T \, = \, 
T_1^{1/2}T_1^{1/2} + T_2^{1/2}T_2^{1/2}
\, = \, 
T^{1/2}A_1^* A_1T^{1/2} + T^{1/2}A_2^* A_2T^{1/2}
\, = \, 
T^{1/2}ST^{1/2}
$$
and hence 
$$
\langle ST^{1/2} \xi \mid T^{1/2}\xi \rangle
\, = \, 
\langle T^{1/2} \xi \mid T^{1/2}\xi \rangle 
\hskip.5cm \text{for all} \hskip.5cm 
\xi \in \Hi;
$$
it follows that $\langle S \eta \mid \eta \rangle = \langle \eta \mid \eta \rangle$ 
for all $\eta \in \overline{T(\Hi)}$.
Since $S = 0$ on $T(\Hi)^{\perp}$, we see that $S$ is 
the orthogonal projection on $\overline{T(\Hi)} = \overline{T^{1/2}(\Hi)}$
and therefore $T^{1/2} = ST^{1/2}$.
\par
 
If $t(T_1) = +\infty$ or $t(T_2) = +\infty$, 
then $t(T) = +\infty$, by the second step.
So, we can assume that $T_1^{1/2} = \lambda(a_1)$ and $T_2^{1/2} = \lambda(a_2)$ 
for $a_1, a_2 \in \AC_b$.
Then 
$$
T^{1/2}
\, = \,
ST^{1/2}= A_1^* A_1T^{1/2} + A_2^* A_2 T^{1/2}
\, = \, 
A_1^* T_1^{1/2}+ A_2^* T_2^{1/2}
$$
and hence $T^{1/2} = \lambda(a)$ for $a = A_1^*a_1 + A_2^*a_2 \in \AC_b$.
\par

Observe that $A_j a= a_j$, since 
$$
\lambda(a_j) \, = \, T_j^{1/2} \, = \, A_j T^{1/2} \, = \, A_j\lambda(a) \, = \, \lambda(A_j a).
$$
Observe also that $a \in \overline{T(\Hi)} = \overline{T^{1/2}(\Hi)}$; 
indeed, by Lemma~\ref{LemHilbertAlg}~\ref{iDELemHilbertAlg},
there exists $x_i \in \AC$ with $\lim_i \rho(x_i) = I$ in the strong operator topology
and hence
$$
a \, = \, \lim_i \rho(x_i)a \, = \, \lim_i \lambda(a) x_i \, = \, \lim_i T^{1/2} x_i 
$$
It follows that $Sa = a$ and therefore
$$
\begin{aligned}
t(T)& \, = \, \langle a \mid a \rangle \, = \, \langle S a \mid a \rangle
\\
& \, = \, \langle A_1 a \mid A_1 a \rangle + \langle A_2 a \mid A_2 a \rangle
\\
& \, = \, \langle a_1 \mid a_1 \rangle + \langle a_2 \mid a_2 \rangle
\, = \, t(T_1) + t(T_2).
\end{aligned}
$$

\vskip.2cm

$\bullet$ {\it Fourth step.} 
Let $T \in \lambda(\AC)$; we claim that $t(TT^*) = t(T^*T)$.
\par

Let $T = W \vert T \vert$ be the polar decomposition of $T$. 
Recall $\vert T \vert = (T^*T)^{1/2}$ and that $W$ is partial isometry with initial space
$(\ker T)^\perp = \overline{\vert T \vert(\Hi)}$ and final space $\overline{T(\Hi)}$;
recall also that $W$ and $\vert T \vert$ belong to the von Neumann algebra 
generated by $T$ and hence to $\lambda(\AC)''$.
 \par
 
Assume that $t(T^*T) < +\infty$.
Then $\vert T \vert = \lambda(a)$ for some $a \in \AC$ and 
$$
t(T^*T) \, = \, \Vert a \Vert^2.
$$
\par

Since $a \in \overline{ \vert T \vert(\Hi) }$ (see the third step), we have, on the one hand,
$$
\Vert Wa\Vert \, = \, \Vert a\Vert.
$$
On the other hand,
using Lemma~\ref{LemHilbertAlg}~\ref{iiiDELemHilbertAlg}
and \ref{ivDELemHilbertAlg}, we have
$$
T^* \, = \, (W \lambda(a))^* \, = \, \lambda(a)^*W^* 
\, = \, \lambda(Ja) W^* \, = \, \lambda(JWa)
$$
and it follows (by an argument from the third step applied to $JWa$) 
that $JW a \in \overline{T^*(\Hi)} = \overline{\vert T \vert(\Hi)}$; 
therefore 
$$
\Vert WJW a \Vert \, = \, \Vert JW a\Vert \, = \, \Vert W a\Vert \, = \, \Vert a \Vert
$$
As $\vert T^* \vert = WT^* = \lambda(WJW a)$, it follows that 
$$
t(TT^*) \, = \, \Vert WJW a\Vert^2 \, = \, \Vert a\Vert^2 \, = \, t(T^*T).
$$
We have shown, for any $T \in \lambda(\AC)$, that, 
if $t(TT^*) < +\infty$, then $\tau(T^*T) < +\infty$ and 
$t(T^*T) = t(TT^*)$; 
it follows that, for any $T \in \lambda(\AC)$, 
we have $t(T^*T) < +\infty$ if and only if $t(TT^*) < +\infty$
and then $t(T^*T) = t(TT^*)$. 
This shows that $t(T^*T) = t(TT^*)$ for every $T \in \lambda(\AC)''$.

\vskip.2cm

$\bullet$ {\it Fifth step.} 
The preceding steps show that $t$ is a trace on $\lambda(\AC)''$.
We claim that $t$ is faithful. 
\par

Indeed, let $T \in \lambda(\AC)''_+$ be such that $t(T) = 0$. 
Then $T^{1/2} = \lambda(a)$ for some $a \in \AC_b$ 
and we have $t(T) = \Vert a \Vert^2 = 0$, that is, $a = 0$;
hence $T^{1/2} = 0$ and $T = 0$.

\vskip.2cm

$\bullet$ {\it Sixth step.} 
We claim that $t$ is normal. 
\par

Indeed, let $T \in \lambda(\AC)''_+$ and let $(T_i)_{i}$ 
be an increasing filtering subset of $\lambda(\AC)''_+$
with $T = \sup_{i} T_i$. 
By the second step, we have $\sup_i t(T_i) \le t(T)$. 
So, it suffices to show that $t(T) \le \sup_i t(T_i)$. 
Since this inequality is obvious if $\sup_i t(T_i) = +\infty$, 
we may assume that $\alpha := \sup_i t(T_i) < +\infty$. 
\par

We have $T_i = \lambda(a_i)$ for some $a_i \in \AC_b$ and
$$
\langle a_i \mid a_i \rangle \, = \, t(T_i) \le \alpha 
\hskip.5cm \text{for all} \hskip.2cm 
i.
$$
Since $\lim_i T_i = T$ in the strong operator topology and since $T_i \le T$, 
we have $\lim T_i^{1/2} = T^{1/2}$ in the strong operator topology 
(Lemma~\ref{Lem-SquareRoot}).
Therefore we have for $x,y \in \AC$,
$$
\langle T^{1/2}y \mid x \rangle \, = \, 
\lim_i \langle T_i^{1/2}y \mid x \rangle \, = \, 
\lim_i \langle \lambda(a_i)y \mid x \rangle \, = \,
\lim_i \langle \rho(y)a_i \mid x \rangle \, = \, 
\lim_i \langle a_i \mid xy^* \rangle.
$$
Since the set $\{xy^* \mid x,y \in \AC\}$ is total in $\Hi$ 
and since $\sup_i \Vert a_i \Vert \le \alpha$, 
the limit $\lim_i a_i = a \in \Hi$
exists in the weak topology on $\Hi$;
we then have $\Vert a\Vert \le \alpha$ and
$$
\langle T^{1/2}y \mid x \rangle \, = \,
\langle a \mid xy^* \rangle \, = \, 
\langle \rho(y)a \mid x \rangle
$$
for all $x,y \in \AC$. Therefore $a \in \AC_b$, $T^{1/2} = \lambda(a)$, and 
$$
t(T) \, = \, \langle a \mid a \rangle \, \le \, \alpha \, = \, \sup_i t(T_i).
$$

\vskip.2cm

$\bullet$ {\it Seventh step.} 
We claim that $\mathfrak n_t = \lambda(\AC_b)$. 
\par

Let $T \in \lambda(\AC)''_+$. We have to show $t(T^*T) < +\infty$ 
if and only if $T \in \lambda(\AC_b)$.
Let $T = W \vert T \vert$ be the polar decomposition of $T$.
Assume that $t(T^*T) = t(\vert T \vert^2) < +\infty$;
then $\vert T \vert \in \lambda(\AC_b)$, by definition of $t$.
Since $\lambda(\AC_b)$ is an ideal of $\lambda(\AC)'' = \rho(\AC)'$, 
it follows that $T = W \vert T \vert \in \lambda(\AC_b)$.
Conversely, if $T = \lambda(\AC_b)$, then $\vert T \vert = W^*T \in \lambda(\AC_b)$
and hence $t(T^*T) = t(\vert T \vert^2) < +\infty$.

\vskip.2cm

$\bullet$ {\it Eighth step.} 
We claim that $t$ is semi-finite. 
\par

Since $t$ is normal, it suffices to prove the following claim 
(see Corollaire 3 in \cite[Chap.~I, \S~6]{Dixm--vN}).
Let $T \in \lambda(\AC)''_+$ with $T \ne 0$. 
Then there exists $S \in \lambda(\AC)''_+$ with $0 \ne S \le T$ and 
$t(S) < +\infty$.
\par
 
Let $a \in \AC_b$. Since $\lambda(\AC_b)$ is a right ideal in $\lambda(\AC)''$ 
we have $\lambda(a)T^{1/2} \in \lambda(\AC_b)$
and hence, by the seventh step,
$$
t T^{1/2} \lambda(a)^* \lambda(a) T^{1/2}) \, < \, +\infty.
$$
Observe that, for every $\xi \in \Hi$, we have 
$$
\begin{aligned}
\langle T^{1/2} \lambda(a)^* \lambda(a) T^{1/2}\xi \mid \xi \rangle
\, &= \, \Vert \lambda(a) T^{1/2}\xi \Vert ^2 
\\
\, &\le \, \Vert \lambda(a)\Vert^2 \Vert T^{1/2}\xi \Vert ^2
\, = \, \Vert \lambda(a)\Vert^2 \langle T\xi \mid \xi \rangle,
\end{aligned}
$$
so that 
$$
T^{1/2} \lambda(a)^* \lambda(a) T^{1/2} \, \le \, \Vert \lambda(a)\Vert^2 T.
$$
We claim that there exists $a \in \AC_b$ with $\lambda(a)T^{1/2} \ne 0$. 
\par

Indeed, otherwise we would have $T^{1/2} = 0$ and hence $T = 0$, 
since $I$ belongs to the closure of $\lambda(\AC)'$ 
in the strong operator topology.
\par

Let $a \in \AC_b$ be such that $\lambda(a)T^{1/2} \ne 0$. 
Then $T^{1/2} \lambda(a)^* \lambda(a) T^{1/2} \ne 0$.
Set 
$$
S \, := \, T^{1/2} \lambda(b)^* \lambda(b) T^{1/2} 
\hskip.5cm \text{for} \hskip.2cm
b \, = \, \frac{1}{\Vert \lambda(a) \Vert} a \in \AC_b.
$$
We have
$$
S \ne 0,
\hskip.5cm
t(S) < +\infty, 
\hskip.5cm \text{and} \hskip.5cm 
S \le T.
$$
This concludes the proof.
\end{proof}

\begin{exe}
% 10.A.7
\label{Ex-HilbertAlgebra}
(1)
Let $\Hi$ be a separable Hilbert space; the $*$-algebra ${\rm TC}(\Hi)$
of trace-class operators on $\Hi$ is a Hilbert algebra for the scalar product
$$
\langle S \mid T \rangle \, = \, \Tr(T^*S);
$$
the completion of ${\rm TC} (\Hi)$
is the Hilbert space ${\rm HS}(\Hi)$ of Hilbert-Schmidt operators on $\Hi$.
Every element in ${\rm HS}(\Hi)$ is bounded.
The associated von Neumann algebras are the isomorphic images of $\Li (\Hi)$ acting 
on ${\rm HS}(\Hi)$ by left and right composition respectively.
Identifying these algebras with $\Li (\Hi)$, 
the associated natural trace $t$ is the standard trace $\Tr$ on $\Li (\Hi)$.

\vskip.2cm

(2)
Let $G$ be a unimodular LC group, with Haar measure $dx$.
The convolution algebra $C^c(G)$, equipped with the
involution $f^*(x) = \overline{f(x^{-1})}$ and the
scalar product
\index{Algebras! $3$@convolution algebras, $C^c(G)$, $C^0(G)$, $L^1(G)$}
$$
\langle f_1 \mid f_2 \rangle \, = \, \int_G f_1(x) \overline{f_2(x)} dx
$$
is a Hilbert algebra;
the completion of $C^c(G)$ is $L^2(G)$.
A function $\xi \in L^2(G)$ is bounded if and only if
$\xi \ast \eta \in L^2(G)$ for every $\eta \in L^2(G)$; 
in this case $\lambda(\xi)$ is defined by $\lambda(\xi) (\eta) = \xi \ast \eta$ 
for all $\eta \in L^2(G)$. 
The associated von Neumann algebras
$\lambda(C^c(G))''$ and $\rho(C^c(G))''$ are the von Neumann algebras
generated by the left and the right regular representations of $G$.
The associated natural trace $t$ is given by 
$$
t(f^* \ast f) \, = \, \Vert f \Vert_2^2 
\hskip.5cm \text{for all} \hskip.2cm 
f \in C^c(G).
$$
\par

In particular, when $G$ is compact, 
$L^2(G)$ is itself a Hilbert algebra.

\vskip.2cm

(3)
Let $(X, \mu)$ be a measure space.
Then $L^\infty(X, \mu) \cap L^2(X, \mu)$ is a commutative Hilbert algebra
with Hilbert space completion $L^2(X, \mu)$.

\vskip.2cm

(4)
Let $\mathcal M$ be a von Neumann algebra,
$t$ a semi-finite faithful normal trace on $\mathcal M$,
and $\mathfrak n_t$ the two-sided ideal defined 
in the beginning of Section \ref{Section-RepTraceC*}.
Then $\mathfrak n_t$, with the scalar product defined by
$\langle x \mid y \rangle = t(xy^*)$,
is a Hilbert algebra.
This is a special case of the more general construction
from Proposition~\ref{Prop-TraceHilbertAlg} below.
% [For more on this, see Theorem 2 in \cite[Chap.~I, \S~6, no~2]{Dixm--vN}.]
We come back to this example in Section \ref{Section-StandardRep}.
\par

Note that (1) above is a particular case of (4).

\vskip.2cm

(5)
Every dense $*$-subalgebra of a Hilbert algebra is a Hilbert algebra.
For example, in the situation of (1), 
the algebra of operators of finite rank on $\Hi$ is a Hilbert algebra.

\vskip.2cm

(6)
Proposition \ref{Lem-HilbertAlgebra} below
provides a family of examples of Hilbert algebras.
\end{exe}

In combination with Example~\ref{Ex-HilbertAlgebra} (2),
we record the following important consequence of 
Theorem~\ref{TheoHilbertAlg-Commutation} and Theorem~\ref{TheoHilbertAlg-Trace}.

\begin{theorem}
% 10.A.8
\label{Theo-RegRepHilbertAlg}
Let $G$ be a unimodular locally compact group.
Let $\Li (G) := \lambda(C^c(G))''$ and $\Ri (G) := \rho(C^c(G))''$
be the von Neumann algebras
generated by the left and the right regular representations of $G$.
Then $\Li (G)$ and $\Ri (G)$ are semi-finite von Neumann algebras and we have
$$
\Li (G)' \, = \, \Ri (G)
\hskip.5cm \text{and} \hskip.5cm
\Ri (G)' \, = \, \Li (G).
$$
\end{theorem}

\begin{rem}
% 10.A.9
\label{Rem-Theo-RegRepHilbertAlg}
Let $G$ be a \emph{non unimodular} LC group, and let $\Li (G)$ and $\Ri (G)$ be as above.
Then it is still true that $\Li (G)' = \Ri (G)$;
see \cite[Chap.~I, \S~5, Th.\ 1 and Exercice 5]{Dixm--vN}.
However, $\Li (G)'$ may no longer be semi-finite;
indeed, we will give in Section \ref{Section-VNRegular} below
examples of LC groups $G$ for which $\Li (G)'$ is a von Neumann algebra of type III.
\end{rem}

The following proposition shows that a lower semi-continuous trace on a C*-algebra
gives rise to a Hilbert algebra (compare \cite[6.4 \& 6.2]{Dixm--C*}).

\begin{prop}
% 10.A.10
\label{Prop-TraceHilbertAlg}
Let $A$ be a C*-algebra and $t$ a lower semi-continuous trace on $A$,
with $t \ne 0$. 
% \begin{enumerate}[label=(\roman*)]
\begin{enumerate}[label=(\alph*)]
\item\label{iDEProp-TraceHilbertAlg}
The map $(x,y) \mapsto t(y^*x)$ 
is a sesquilinear positive Hermitian form on $\mathfrak n_t$.
\item\label{iiDEProp-TraceHilbertAlg}
The set $N_t$ of all $x$ in $\mathfrak n_t$ with $t(x^*x) = 0$ 
is a selfadjoint two-sided ideal of $A$; 
moreover, we have
$$
N_t \, = \,
\{x \in \mathfrak n_t \mid t(y^*x) = t(xy^*) = 0 
\hskip.2cm \text{for all} \hskip.2cm
y \in \mathfrak n_t \}.
$$
\item\label{iiiDEProp-TraceHilbertAlg}
The algebra $\AC_t := \mathfrak n_t / N_t$, equipped with the map 
$$
(x+N_t,y+N_t) \, \mapsto \, t(y^*x)
$$
and the involution $x+N_t \mapsto x^*+N_t$,
\emph{is a Hilbert algebra.}
\item\label{ivDEProp-TraceHilbertAlg}
Let $\lambda_t$ and $\rho_t$ be the
representations of $A$ on the Hilbert space completion $\Hi_t$ of $\AC_t$,
defined by $\lambda_t(a) : x+N_t\mapsto ax+N_t$ 
and $\rho_t(a) \,\colon x + N_t \mapsto xa^* + N_t$
% sic ! corrigŽ la correction le 14 mai
for $a \in A$.
Then the von Neumann subalgebras
$\lambda_t(A)''$ and $\rho_t(A)''$ of $\Li (\Hi_t)$ 
coincide with $\lambda(\AC_t)''$ and $\rho(\AC_t)''$,
respectively, and are each other commutants:
$$
\lambda_t(A)' \, = \, \rho_t(A)''
\hskip.5cm \text{and} \hskip.5cm
\rho_t(A)' \, = \, \lambda_t(A)''.
$$
\item\label{vDEProp-TraceHilbertAlg}
Let $t_t$ denotes the natural trace on $\lambda_t (A))''$;
then
$$
t_t \circ \lambda_t(x) \, = \, t(x)
\hskip.5cm \text{for every} \hskip.2cm 
x \in \mathfrak m_{t,+}.
$$
\end{enumerate}
\end{prop}
 
 \begin{proof}
 \ref{iDEProp-TraceHilbertAlg}
 It is straightforward to check that $(x,y) \mapsto t(y^*x)$ 
is a sesquilinear positive Hermitian form on $\mathfrak n_t$.

\vskip.2cm

\ref{iiDEProp-TraceHilbertAlg}
Let $x \in N_t$, that is, $t(x^*x) = 0$. Cauchy--Schwarz inequality shows that 
$t(y^*x) = 0$ for every $y \in \mathfrak n_t$. 
Therefore, $N_t$ is a two-sided of $\mathfrak n_t$.
Since $t$ is a trace, we have $t(yx^*) = t(x^*y)$ and hence, by sesquilinearity,
$$
t(yx^*) \, = \, t((y^*x)^*) \, = \, \overline{t(y^*x)} \, = \,0
$$
for every $y \in \mathfrak n_t$. This shows that $N_t$ is a $*$-ideal and hence 
a two-sided of $\mathfrak n_t$.

\vskip.2cm

\ref{iiiDEProp-TraceHilbertAlg}
The map on $(x+N_t,y+N_t) \mapsto t(y^*x)$ is well-defined on $\AC_t$, since $N_t$ is a left 
ideal of $\mathfrak n_t$. Also, since $N_t$ is a $*$-ideal, the involution
$x+N_t\mapsto x^*+N_t$ is well-defined on $\AC_t$.
\par

Property (1) in the definition of a Hilbert algebra is satisfied, since $t$ is a trace;
Property (2) is obvious. 
\par.

Let us check Property (3). 
The selfadjoint ideal $\mathfrak m_t = \mathfrak n_t^2$
has an approximate identity (see \ref{AppAlgC*}) and,
for every $y \in \mathfrak n_t$, 
the map $x \mapsto t(y^*xy)$ is a continuous positive linear form
on $\mathfrak m_t$; hence
$$
t(y^*x^*xy) \, \le \, \Vert x^*x \Vert t(y^*y) 
\hskip.5cm \text{for all} \hskip.2cm 
x,y \in \mathfrak n_t ,
$$
see \cite[Proposition 2.1.5]{Dixm--C*}.
It follows that, for all $x,y \in \mathfrak n_t$, we have
$$
\begin{aligned}
\Vert xy+N_t\Vert ^2
& \, = \, t((xy)^* yx) \, = \, t(y^*x^*xy)
\\
&\, \le \, \Vert x^*x \Vert t(y^*y) \, = \, \Vert x^*x \Vert \Vert y+N_t \Vert^2
\end{aligned}
$$
and this shows that $y+N_t\mapsto xy+N_t$ is continuous 
for every $x \in \mathfrak n_t $.
\par

To show Property (4), let $(u_i)_i$ be an increasing approximate identity
for $\mathfrak n_t$ in $\mathfrak n_t^+$
(see \ref{AppAlgC*}).
Let $x \in \mathfrak n_t $. Then 
$$
\begin{aligned}
\Vert (u_i x- x)+N_t \Vert ^2
& \, = \, t((u_i x-x)^*(u_i x-x))
\\
& \, = \, t(x^*u_i^2x )- 2t(x^*u_i x)+ t(x^*x) 
\\
& \, \le \, \Vert u_i^2 \Vert t(x^*x) - 2t(x^*u_ix) +t(x^*x)
\\
& \, \le \, 2 \left(t(x^*x) - t(x^*u_ix)\right).
\end{aligned}
$$
Since $(x^* u_i x)_i$ is an increasing family with $\lim_i x^*u_ix = x^*x$ 
and since $t$ is lower semi-continuous, we have $\lim_i t(x^*u_ix) = t(x^*x)$.
This shows that $\lim_i \Vert (u_i x- x)+N_t \Vert = 0$. 

\vskip.2cm

 \ref{ivDEProp-TraceHilbertAlg}
This is a direct consequence of Theorem~\ref{TheoHilbertAlg-Commutation}

\vskip.2cm

 \ref{vDEProp-TraceHilbertAlg}
Let $x \in \mathfrak m_{t,+}$. Then $a := x^{1/2}\in \mathfrak m_{t,+}$ and 
$\lambda_t(x)^{1/2} = \lambda (a+ N_t) \in \lambda(\AC)$; 
hence
$$
\tau_t(\lambda_t(x)) \, = \, \langle a+N_t \mid a+N_t \rangle \, = \, t(a^*a) \, = \, t(x) ,
$$
and the proof is complete.
\end{proof}

\section
[Standard representation]
{The standard representation of a semi-finite von Neumann algebra}
% Section 10.B
\label{Section-StandardRep}

We will often need to consider the so-called 
standard representation of a semi-finite von Neumann algebra.
Our exposition rests on \cite[Chap.~I, \S~6, no~2]{Dixm--vN}.
\par

Let $\M$ be a von Neumann algebra,
equipped with a faithful semi-finite normal trace $\tau$.
Let $\mathfrak m = \mathfrak m_\tau$ be the ideal of definition of $\tau$,
and $\mathfrak n = \mathfrak n_\tau$
the ideal of all $x \in \M$ with $\tau(x^*x) < \infty$;
recall that $\mathfrak m = \mathfrak n^2$.
\par

Since $\tau$ is normal, $\tau$ is lower-continuous for the weak topology 
and hence for the norm topology. 
By Proposition~\ref{Prop-TraceHilbertAlg}, the space $\mathfrak n$,
equipped with the scalar product $(x,y)\mapsto \tau (y^*x)$, 
is therefore a Hilbert algebra.
\par

Let $L^2(\mathfrak n, \tau)$ denote the Hilbert space completion 
of $\mathfrak n$, with respect to this scalar product. 
Instead of $\lambda$, we denote by $x \mapsto \widehat x$ 
the associated representation of $\mathcal M$ on $L^2(\mathfrak n, \tau)$,
where $\widehat x$ is the continuous extension to $L^2(\mathfrak n, \tau)$
of the linear map
$$
\mathfrak n \to \mathfrak n,
\hskip.5cm
y \mapsto xy .
$$
Moreover, if $\widehat \tau$ denotes the natural trace on $\mathfrak n$,
then 
$$
\widehat \tau (\widehat x) \, = \, \tau(x)
\hskip.5cm \text{for all}\hskip.5cm
x \in \mathcal M.
$$

\begin{prop}
% 10.B.1
\label{Prop-StandRep}
The algebra
$$
\widehat{\mathcal M} \, := \, \{\widehat x \mid x \in \mathcal M\}
$$
is a von Neumann subalgebra of $\Li (L^2(\mathcal M, \tau))$.
\end{prop}

\begin{proof}
Since the map
$$
\mathcal M \, \to \, \Li (L^2(\mathcal M, \tau)), 
\hskip.5cm
x \, \mapsto \, \widehat x
$$
is a unital $*$-algebra homomorphism, it suffices to show that 
this map is normal (see Corollaire 2 \cite[Chap.~I, \S~4, 39]{Dixm--vN}).
\par

Let $(x_i)_i$ be an increasing family in ${\mathcal M}_+$ with 
$x = \sup_{i} x_i \in {\mathcal M}_+$.
Then $(\widehat{x_i})_i$ is an increasing family in $\Li (L^2(\mathcal M, \tau))_+$ 
with $\widehat{x_i} \le \widehat x$ for all $i$ 
and hence $\sup_i \widehat{x_i} \le \widehat x$.
\par

Let $y \in \mathfrak n$. Then 
$$
\langle \widehat{x_i} y \mid y \rangle \, = \,
\langle x_i y \mid y \rangle
\, = \,
\tau(y^*x_iy)
$$
for all $i$. Since $\sup_i y^*xy_i = y^*xy$ and since
$\tau$ is normal, we have 
$\sup_i \tau(y^*xy_i) = \tau(y^*xy)$ and hence 
$$
\sup_i \langle \widehat{x_i} y \mid y \rangle \, = \,
\tau(y^*xy) \, = \, \langle \widehat x y \mid y \rangle 
$$
for all $y \in \mathfrak n$. 
As $\mathfrak n$ is dense in $L^2(\mathcal M, \tau)$,
it follows that $\sup_i \widehat{x_i} = \widehat x$.
\end{proof}

It follows from the previous lemma that $\widehat{\mathcal M}$
coincides with the von Neumann $\lambda(\mathfrak n)''$
associated with the Hilbert algebra
$\mathfrak n$ as in Section~\ref{Section-RepTraceC*}.
\par

\index{Standard! $1$@representation of $(\mathcal M, \tau)$}
\index{Representation! standard}
The representation 
$$
\mathcal M \, \to \, \widehat{\mathcal M}, 
\hskip.5cm
x \, \mapsto \, \widehat x
$$
is called the \textbf{standard representation} 
of the pair $(\mathcal M, \tau)$.
Observe that the standard representation is an isomorphism between $\mathcal M$ and $\widehat{\mathcal M}$;
indeed, if $\widehat x = 0$ for $x \in \mathcal M$, then $x \in \mathfrak n$ and 
$\tau(x^*x) = 0$ and hence $x = 0$, since $\tau$ is faithful. 
The von Neumann algebra $\widehat{\mathcal M}$
is often called the \textbf{standard form} of $\mathcal M$.
%For all this, see Th\'eor\`eme 2 in \cite[I, \S~6, 2]{Dixm--vN}. 

\begin{exe}
% 10.B.2
\label{Ex-StandardRep}
(1)
Let $\Hi$ be a separable Hilbert space and $\mathcal M = \Li (\Hi)$.
The standard representation of the von Neumann algebra $\Li (\Hi)$ is 
the representation 
$T \mapsto \widehat T$ of $\Li (\Hi)$ 
on the Hilbert space ${\rm HS}(\Hi)$ of Hilbert-Schmidt operators on $\Hi$
given by 
$$
\widehat T (S) \, = \, TS, 
\hskip.5cm \text{for all} \hskip.2cm
S \in {\rm HS}(\Hi).
$$

\vskip.2cm

(2) Let $G$ be a unimodular LC group.
The standard representation of the von Neumann algebra $\lambda_G(G)''$ 
of the left regular representation of $G$ 
coincides with $\lambda_G(G)''$ acting on $L^2(G)$. 
\end{exe}

\section
{Group representations associated to traces}
% Section 10.C
\label{Section-TraceRep}

Let $G$ be a LC group.
Let $t$ be a non-zero lower semi-continuous semi-finite trace 
on the maximal C*-algebra $C^*_{\rm max}(G)$ of $G$. 
Let $\lambda_t$ be
% \footnote{CHOIX : $\pi_t$ ou $\lambda_t$;
% \textbf{Reponse:} il faut prendre $\lambda_t$ pour rester coh\'erent avec la notation}
the representation of $C^*_{\rm max}(G)$ associated to $t$,
as in Section~\ref{Section-RepTraceC*}; 
then $\lambda_t$ is a traceable representation of $G$, in the following sense.
\par

\index{Trace representation}
\index{Representation! trace}
A \textbf{trace representation} of a LC group $G$ is a pair $(\pi, \tau)$,
with $\pi$ a representation of $G$
and $\tau$ a normal faithful semi-finite trace on the von Neumann algebra $\pi(G)''$
such that the following condition holds:
\begin{enumerate}
\item[(*)]\label{conditiontracerep}
the set of all $\pi(x)$ with $x \in C^*_{\rm max}(G)$ and 
$0 < \tau(\pi(x^*x)) < +\infty$ is weakly dense in $\pi(G)''$.
% (that is, such that $\pi(C^*_{\rm max}(G)) \cap \mathfrak n_\tau$ 
% is weakly dense in $\pi(G)''$).
\end{enumerate}
(Recall from \ref{vNalgLCgroup} that $\pi(G)'' = \pi(C^*_{\rm max}(G))''$.) 
\par

Note that the existence of a normal faithful semi-finite trace on $\pi(G)''$
\emph{does not} imply Condition (*), 
even when $\pi$ is factorial;
see Example \ref{ExaNormalRep}(3) below.
\par

\index{Quasi-equivalent! trace representations}
\index{Trace representation! quasi-equivalent}
We say that two trace representations $(\pi_1, \tau_1)$ and $(\pi_1, \tau_1)$
are \textbf{quasi-equivalent} 
if there exists an isomorphism $\Phi \,\colon \pi_1(G)'' \to \pi_2(G)''$ 
such that $\Phi(\pi_1(g)) = \pi_2(g)$ for all $g \in G$, 
and such that $\tau_1 = \tau_2 \circ \Phi$.
\par

\index{Traceable representation}
\index{Representation! traceable}
A representation $\pi$ of $G$ is 
\textbf{traceable} 
if there exists a normal faithful semi-finite trace $\tau$
on the von Neumann algebra $\pi(G)''$
such that the pair $(\pi, \tau)$ is a trace representation.
\par

Trace and traceable representations appear for C*-algebras 
in \cite[6.6 and 6.7]{Dixm--C*},
and for LC groups in \cite[17.1]{Dixm--C*}.

\vskip.2cm

Being traceable only depends on the quasi-equivalence class of the representation:

\begin{prop}
% 10.C.1
\label{quasieqtraceabletraceable}
Let $\pi$ and $\rho$ be quasi-equivalent representations of $G$. 
If $\pi$ is traceable, then $\rho$ is traceable.
\end{prop}

\begin{proof}
Let $\tau$ be a normal faithful semi-finite trace on $\pi(G)''$
such that $(\pi, \tau)$ is a trace representation and let 
$\Phi \,\colon \pi(G)'' \to \rho(G)''$ be an isomorphism 
with $\Phi(\pi(g)) = \rho(g)$ for all $g \in G$. 
Since $\Phi$ is an isomorphism of von Neumann algebras
$\Phi$ and $\Phi^{-1}$ are normal maps;
see \cite[Chap.~I, \S~4, Corollaire 1]{Dixm--vN}.
Therefore $\tau' := \tau \circ \Phi^{-1}$ is a normal faithful
semi-finite trace on $\rho(G)''$ 
and we have $\Phi(\pi(x)) = \rho(x)$ for all $x \in C^*_{\rm max}(G)$.
It follows that Condition (*) above is satisfied for $\tau'$ and hence $(\rho, \tau')$ is a trace representation of $G$.
\end{proof}

\begin{exe}
% 10.C.2
\label{Exemples-RepTracables}

(1) 
Let $\pi \,\colon G \to \U(n)$ be a finite-dimensional representation of a LC group $G$.
Then $(\pi, \tau)$ is a trace representation, where $\tau$ is the restriction of the 
standard trace on $M_n(\C)$ to the algebra $\pi(G)''$, 
which is the linear span of $\pi(G)$.
Assume that $\pi$ is not a factor representation, 
that is, $\pi$ is not a multiple of an irreducible representation.
Then there are more faithful traces on $\pi(G)''$
than positive multiples of $\tau$.
For instance,
assume that $\pi = \pi_1 \oplus \pi_2$ is the direct sum 
of two non-equivalent irreducible representations
$\pi_1$ and $\pi_2$ on subspaces $V_1$ and $V_2$ of $\C^n$. 
Then, for every $a > 0, b > 0$, a faithful trace $\tau_{a,b}$ is defined on
$$
\pi(G)'' \, = \, \pi_1(G)'' \oplus \pi_2(G)'' \, \approx \, \Li (V_1) \oplus \Li (V_2)
$$
by 
$$
\tau_{a,b}(T \oplus S) \, = \, a \tau_1(T) + b\tau_2(S)
$$
where $\tau_i$ is the standard trace on $\Li (V_i)$.

\vskip.2cm

(2) 
Examples of traceable representations of a unimodular 
LC group $G$ arise often as follows. 
Let $\nu$ be complex-valued Radon measure on $G$. 
Assume that $\nu$ is of positive type, that is,
% pas comme dans \cite[13.7.1]{Dixm--C*} qui utilise $\widetilde f$ au lieu de $f^*$
$$
\int_G (f^* \ast f) (x) d\nu(x) \, \ge \, 0 
\hskip.5cm \text{for all} \hskip.2cm 
f \in C^c(G),
$$
where $f^*$ is defined by $f^*(x) = \overline{f(x^{-1})}$.
% $wide\tilde f$ d\'efini in \cite[13.2.3]{Dixm--C*}.
% Denote by $\pi_\nu$ the GNS representation associated to $\nu$.
Assume also that $\nu$ is central, that is, 
$$
\int_G (f_1 \ast f_2) (x) d\nu(x) \, = \, \int_G (f_2 \ast f_1) (x)d\nu(x) 
\hskip.5cm \text{for all} \hskip.2cm 
f_1, f_2 \in C^c(G).
$$
Then $\nu$ defines a lower semi-continuous, semi-finite trace $t_\nu$ on $C^*_{\rm max}(G)$,
with ideal of definition containing $C^c(G)$, and given there by 
$$
t_\nu(f \ast f^*) = \int_G (f \ast f^*) (x) d\nu(x) 
\hskip.5cm \text{for all} \hskip.2cm 
f \in C^c(G).
$$
The representation $\pi_\nu$ associated to $t_\nu$ 
as in Proposition \ref{Prop-TraceHilbertAlg}
is traceable.
\par

For instance, when $\nu = \delta_e$ is the Dirac measure at the group unit, 
$t_\nu$ is the semi-finite trace given by
$$
t_\nu (f \ast f^*) \, = \, (f \ast f^*) (e) \, = \, \Vert f \Vert_2^2 
\hskip.5cm \text{for all} \hskip.2cm 
f \in C^c(G).
$$
and $\pi_\nu$ is the left regular representation
(see also Example~\ref{Ex-HilbertAlgebra}(2)).
The trace $t_\nu$ is finite only when $G$ is discrete.
For all this, see \cite{Gode--54}.

\vskip.2cm

(3)
When $G$ is a real Lie group, 
the construction from (2) can be --- and has to be --- extended to
central distributions of positive type on $G$.
The resulting lower semi-continuous, semi-finite traces on $C^*_{\rm max}(G)$ 
have ideals of definition which contain
the space of smooth functions with compact support on $G$.
\end{exe}

As shown in Proposition~\ref{Prop-TraceHilbertAlg},
for every lower semi-continuous semi-finite trace $t$ on $C^*_{\rm max}(G)$,
there exists a representation $\lambda_t$ of $G$ 
(equivalently of $C^*_{\rm max}(G)$)
and a trace $\tau_t$ on $\lambda_t(G)''$ 
such that $(\lambda_t, \tau_t)$ is a trace representation of $G$,
with $t = \tau_t \circ \lambda_t$.
\par

Proposition \ref{traceablerepbij} will show that 
every trace representation of $G$ is of this kind.
It is convenient to state two preliminary lemmas.

\begin{lem}
% 10.C.3
\label{Lemme-Trace-Representation}
Let $(\pi, \tau)$ be a trace representation of $G$. 
Then $\tau \circ \pi$ is a lower semi-continuous semi-finite trace 
on $C^*_{\rm max}(G)$.
\end{lem} 

\begin{proof}
It is clear that $t := \tau \circ \pi$ is a trace on $C^*_{\rm max}(G)$.
Since $(\pi, \tau)$ is a trace representation, $\pi(\mathfrak m_t)$
is a weakly dense in $\pi(G)''$ and hence the restriction of $\pi$ 
to the ideal $\overline{\mathfrak m_t}$ is non-degenerate. 
Let $(u_i)_i$ be an increasing approximate identity for 
$\overline{\mathfrak m_t}$ with $u_i \in \mathfrak m_+$.
Then $(\pi(u_i))_i$ converges to $I$ for the strong operator topology.

Let $x \in C^*_{\rm max}(G)_+$. 
Since $(\pi(x^{1/2}u_i x^{1/2}))_i$ is an increasing family in $\pi(G)''_+$ with 
$$
\lim_i \pi(x^{1/2}u_i x^{1/2}) \, = \, \pi(x)
$$
and since $\tau$ is normal, we have
$$
\lim_i \tau(\pi(x^{1/2}u_i x^{1/2})) \, = \, \tau(\pi(x)).
$$
This shows that $\lim_i t(x^{1/2}u_i x^{1/2}) = t(x)$.
Since $\pi(x^{1/2}u_i x^{1/2}) \in \mathfrak m_t$, we have
$$
t(x^{1/2}u_i x^{1/2}) \, < \, +\infty
$$
and hence $t$ is semi-finite. 
Moreover, $\tau$ being normal is lower semi-continuous for the weak topology
and therefore for the norm topology on $\pi(G)''$. 
It follows that $t = \tau \circ \pi$ is lower semi-continuous.
\end{proof}

\begin{lem}
% [compare with Proposition \ref{GNSbijP(G)cyclic}]
% 10.C.4
\label{Lem-UniqueGNS-TraceRep}
Let $(\pi_1, \tau_1)$ and $(\pi_2, \tau_2)$ be trace representation of $G$ such that
$\tau_1 \circ \pi_1 = \tau_2 \circ \pi_2$.
Then $(\pi_1, \tau_1)$ and $(\pi_2, \tau_2)$ are quasi-equivalent.
\end{lem} 

\begin{proof}
Set 
$$
\mathcal M_j \, = \, \pi_j(G)'' \, = \, \pi_j(C^*_{\rm max}(G))''
$$
and $\mathfrak n_j = \mathfrak n_{\tau_j}$ for $j = 1, 2$. 
Recall that $\pi_j(C^*_{\rm max}(G)) \cap \mathfrak n_j$ 
is weakly dense in $\M_j$,
since $\pi_j$ is a trace representation.
\par

We consider the standard representation 
$x \mapsto \widehat x$ of $\mathcal M_j$ on $L^2(\mathfrak n_j, \tau_j)$. 
The map $\Phi_j \,\colon x \mapsto \widehat x$
is an isomorphism between the von Neumann algebras $\mathcal M_j$ and 
$$
\widehat{\mathcal M}_j \, := \, \{\widehat x \mid x \in \mathcal M_j\}
$$
(see Section \ref{Section-StandardRep}).
Moreover, $\widehat{\tau_j} \circ \Phi_j = \tau_j$,
where $\widehat{\tau_j}$ is the natural trace on $\widehat{\mathcal M_j}$.
\par

Denote by $\widehat{\pi_j}$ the representation of $G$ 
on $L^2(\mathfrak n_j, \tau_j)$ defined by 
$\widehat{\pi_j} = \Phi_j \circ \pi_j$, that is,
$$
\widehat{\pi_j}(g) \, = \, \widehat{\pi_j(g)}
\hskip.5cm \text{for all}\hskip.5cm
g \in G.
$$
We have, for every $a \in C^*_{\rm max}(G)$,
$$
\widehat{\tau_j}(\widehat{\pi_j(a)}) \, = \, \tau_j(\pi_j(a))
$$
and hence $\widehat{\tau_j}(\widehat{\pi_j(a)}) < +\infty$ 
if and only if $\tau_j(\pi_j(a)) < +\infty$.
Since $\Phi_j \,\colon \M_j \to \widehat{\mathcal M}_j$ is an isomorphism, 
this shows that $(\widehat{\pi_j}, \widehat{\tau_j})$ is a trace representation
and that $(\widehat{\pi_j}, \widehat{\tau_j})$ 
is quasi-equivalent to the trace representation $(\pi_j, \tau_j)$. 
\par

As a consequence, it suffices to show that the trace representations 
$(\widehat{\pi_1}, \widehat{\tau_1})$ and $(\widehat{\pi_2}, \widehat{\tau_2})$
are quasi-equivalent.
\par

Set 
$$
\A_j \, := \, 
\widehat{\pi_j}(C^*_{\rm max}(G)) \cap \mathfrak n_{\widehat{\tau_j}};
$$
then $\A_j$ is a Hilbert subalgebra of the Hilbert algebra $\mathfrak n_{\widehat{\tau_j}}$
and, since $\A_j$ is weakly dense in $\widehat{\mathcal M}_j$,
the associated left von Neumann algebra 
coincides with $\widehat{\mathcal M}_j$ 
(see Proposition 1 in \cite[Chap.~I, \S~5, 2]{Dixm--vN}).
Therefore $\A_j$ is dense in the Hilbert space completion of 
$\mathfrak n_{\widehat{\tau_j}}$,
which is $L^2(\mathfrak n_j, \tau_j)$. 
\par

Observe that, by assumption, we have 
$$
\widehat{\tau_1}(\widehat{\pi_1(a)}) \, = \, \tau_1(\pi_1(a)) \, = \, 
\tau_2(\pi_2(a)) \, = \, \widehat{\tau_2}(\widehat{\pi_2(a)}),
$$
for all $a \in C^*_{\rm max}(G)$.
In particular, we have $\widehat{\pi_1}^{-1}(\A_1) = \widehat{\pi_2}^{-1}(\A_2)$.
Since $\widehat{\tau_1}$ and $\widehat{\tau_2}$ are faithful, it follows that the map 
$$
\A_1 \to \A_2,
\hskip.5cm
\widehat{\pi_1}(a) \mapsto \widehat{\pi_2}(a), 
\hskip.5cm \text{for all} \hskip.5cm
a \in \widehat{\pi_1}^{-1}(\A_1)
$$
is well-defined and extends to a Hilbert space isomorphism 
$$
U \, \colon \, L^2(\mathcal M_1, \tau_1) \to L^2(\mathcal M_2, \tau_2).
$$
Let 
$$
\Phi \, \colon \, 
\Li (L^2(\mathcal M_1, \tau_1)) \to \Li (L^2(\mathcal M_2, \tau_2)), 
\hskip.5cm 
T \mapsto UTU^{-1}.
$$
Then $\Phi$ restricts to an isomorphism 
from $\widehat{\mathcal M_1}$ to $\widehat{\mathcal M_2}$,
such that $\Phi( \widehat{\pi_1}(a)) = \widehat{\pi_2}(a)$ 
for all $a \in C^*_{\rm max}(G)$.
It follows that $\Phi( \widehat{\pi_1}(g)) = \widehat{\pi_2}(g)$ for all $g \in G$. 
Moreover, we have $\widehat{\tau_1} = \widehat{\tau_2} \circ \Phi$. 
Therefore $(\widehat{\pi_1}, \widehat{\tau_1})$
and $(\widehat{\pi_2}, \widehat{\tau_2})$ 
are quasi-equivalent.
\end{proof}
 
Let $t$ be a lower semi-continuous semi-finite trace 
on $C^*_{\rm max}(G)$, $t \ne 0$; 
recall from Proposition~\ref{Prop-TraceHilbertAlg} that 
there is an associated trace representation $(\lambda_t, \tau_t)$ 
such that $\tau_t \circ \lambda_t = t$.
Conversely, let $(\pi, \tau)$ be a trace representation of $G$; 
by Lemma \ref{Lemme-Trace-Representation}, $\tau \circ \pi$ is
a lower semi-continuous semi-finite trace on $C^*_{\rm max}(G)$,
which by Lemma \ref{Lem-UniqueGNS-TraceRep}
depends on the equivalence class of $(\pi, \tau)$ only.

\begin{prop}
% 10.C.5
\label{traceablerepbij}
Let $G$ be a locally compact group.
\par

For a lower semi-continuous semi-finite trace $t \ne 0$ 
on the C*-algebra $C^*_{\rm max}(G)$, 
let $(\lambda_t, \tau_t)$ be as in 
Proposition \ref{Prop-TraceHilbertAlg}
the corresponding trace representation.
For a trace representation $(\pi, \tau)$ of $G$,
consider as in Lemma \ref{Lemme-Trace-Representation}
the lower semi-continuous semi-finite traces $\tau \circ \pi$ 
on $C^*_{\rm max}(G)$.
Then the maps $t \mapsto (\text{class of} \hskip.2cm (\lambda_t, \tau_t))$ 
and $(\text{class of} \hskip.2cm (\pi, \tau)) \mapsto \tau \circ \pi$
induce bijections inverse to each other between
\begin{enumerate}[label=(\roman*)]
\item\label{iDEtraceablerepbij}
the set of lower semi-continuous semi-finite traces $t$ 
on the maximal C*-algebra $C^*_{\rm max}(G)$,
and
\item\label{iiDEtraceablerepbij}
the set of trace representations $(\pi, \tau)$ of $G$, 
up to quasi-equivalence.
\end{enumerate}
\end{prop}

\begin{proof}
As a preliminary step,
observe that $\tau \circ \pi$ depends on the equivalence class of $(\pi, \tau)$ only
(Lemma \ref{Lem-UniqueGNS-TraceRep}).
\par

We have to check that the composition of these two maps, in each order,
is the appropriate identity.
\par

For $t$ as in \ref{iDEtraceablerepbij}, we have
$\tau_t \circ \lambda_t = t$
by Proposition \ref{Prop-TraceHilbertAlg}~\ref{vDEProp-TraceHilbertAlg}.
\par

For $(\pi, \tau)$ and $t = \tau \circ \pi$ as in \ref{iiDEtraceablerepbij},
the trace representations $(\pi, \tau)$ and $(\lambda_t, \tau_t)$
are quasi-equivalent by Lemma \ref{Lem-UniqueGNS-TraceRep}.
\end{proof}

In the sequel, we will be mostly interested 
in traceable representations which are \emph{factorial}.

\section
{Normal representations and characters}
% Section 10.D
\label{Section-NormRepChar}

\index{Representation! normal}
\index{Normal representation}
Let $G$ be a LC group. 
A representation $\pi$ of $G$ is \textbf{normal}
if it is factorial and traceable.
Normal factor representations were introduced in \cite{Gode--54}
and further studied and extended to the context of Banach $*$-algebras in \cite{Guic--63}.

\begin{prop}
% 10.D.1
\label{NormalRep}
For a factor representation $\pi$ of a locally compact group $G$,
the following properties are equivalent:
\begin{enumerate}[label=(\roman*)]
\item\label{iDENormalRep}
$\pi$ is a normal factor representation;
\item\label{iiDENormalRep}
the factor $\pi(G)''$ is semi-finite (that is, of type I or II) and
there exists an element $x \in C^*_{\rm max}(G)$ such that
$0 < \tau(\pi(x^*x)) < +\infty$, where $\tau$ is a faithful 
normal trace on $\pi(G)''$.
\end{enumerate}
\end{prop}

\begin{proof}
It is clear that \ref{iDENormalRep} implies \ref{iiDENormalRep}. 
Assume that \ref{iiDENormalRep}
and let $\tau$ be a normal semi-finite faithful trace on $\pi(G)''$. 
Then 
$$
J \, := \, \mathfrak n_\tau \cap \pi( C^*_{\rm max}(G))
$$
is a non-zero two-sided ideal of $C^*_{\rm max}(G)$.
Therefore the weak closure $\mathfrak n$ of $J$ 
in a non-zero weakly closed two-sided ideal of
the von Neumann algebra $\pi(C^*_{\rm max}(G))'' = \pi(G)''$.
Since $\pi(G)''$ is a factor, we have $\mathfrak n = \pi(G)''$; 
see \cite[Chap.~I, \S~4, Corollaire 2]{Dixm--vN}.
Therefore $\pi$ is a normal factor representation.
\par

(The proposition and this proof is that of \cite[6.7.2]{Dixm--C*}.)
\end{proof}

We first characterize irreducible representations which are normal.

\begin{prop}
% 10.D.2
\label{Pro-IrredNormalRep}
Let $G$ be a locally compact group.
An irreducible representation $\pi$ of $G$ is normal 
if and only if $\pi(C^*_{\rm max}(G))$
contains the space $\mathcal K (\Hi_\pi)$ of compact operators on $\Hi_\pi$.
\end{prop}

\begin{proof}
It is obvious that $\pi$ is normal if $\pi(C^*_{\rm max}(G))$
contains $\mathcal K (\Hi_\pi)$. Conversely, assume that
$\pi$ is traceable. Then $\pi(C^*_{\rm max}(G))$ contains a non-zero
trace-class (and hence compact) operator. Since $\pi$ is irreducible, this 
implies that $\pi(C^*_{\rm max}(G))$ contains all compact operators
(see \cite[4.1.10]{Dixm--C*}).
\end{proof}

The following corollary is a consequence of Proposition~\ref{Pro-IrredNormalRep}
and Glimm's theorem \ref{ThGlimm}.
% ??? ThGlimm. ???

\begin{cor}
% 10.D.3
\label{Cor-IrredNormalRep}
Let $G$ be a $\sigma$-compact locally compact group.
The following properties are equivalent.
\begin{enumerate}[label=(\roman*)]
\item\label{iDECor-IrredNormalRep}
every irreducible representation of $G$ is normal;
\item\label{iiDECor-IrredNormalRep}
$G$ is of type I.
\end{enumerate}
\end{cor}

\begin{exe}
% 10.D.4
\label{ExaNormalRep}
Let $G$ be a locally compact group.

\vskip.2cm

(1) A finite-dimensional irreducible representation of $G$ is normal.

\vskip.2cm

(2)
A factor representation $\pi$ of finite type is normal. 
Indeed, when $\pi(G)''$ is a factor that is either of type I$_n$
for some finite integer $n$ or of type II$_1$,
there is a faithful finite trace on this factor.

\vskip.2cm

(3)
A factor representation of type I$_\infty$ need not be normal,
as shown by Corollary~\ref{Cor-IrredNormalRep},
and also by the following explicit example.
\par

\index{Free group}
Consider the representation $\pi$ defined below
of the free group $F_2$ on two generators $a$ and $b$.
We will check that $\pi$ is
infinite-dimensional and irreducible, in particular factorial of type I$_\infty$,
and that it is not normal.
\par

Denote by $C$ the infinite cyclic subgroup of $F_2$ generated by $a$.
Define $\pi$ as the quasi-regular representation of $F_2$ 
acting on $\Hi_\pi := \ell^2(F_2 / C)$;
it is one model for the induced representation $\Ind_C^{F_2}(1_C)$.
The representation $\pi$ is clearly infinite-dimensional.
As $C$ is its own commensurator in $F_2$,
the representation $\pi$ is irreducible (Example \ref{Exa-IrredRepFreeGroup}).
% compare with our Example \ref{ExRepC*simple2}),
It remains to check that $\pi$ is not normal;
by Proposition~\ref{Pro-IrredNormalRep}, 
it suffices to check that $\pi(C^*_{\rm max}(F_2))$
does not contain the space of compact operators on $\Hi_\pi$.
\par

Since $C$ is amenable, the unit representation $1_C$ 
is weakly contained in the regular representation $\lambda_C$ of $C$, 
so that $\pi = \Ind_C^{F_2}(1_C)$ is weakly contained 
in the regular representation 
$\lambda = \Ind_C^{F_2} (\lambda_C)$ of $F_2$,
and therefore $\pi \,\colon C^*_{\rm max}(F_2) \to \Li (\Hi_\pi)$ 
is the composition of the canonical surjection 
$\lambda \,\colon C^*_{\rm max}(F_2) \twoheadrightarrow C^*_\lambda (F_2)$
with a C*-algebra morphism $\pi' \,\colon C^*_\lambda (F_2) \to \Li (\Hi_\pi)$. 
Since the reduced C*-algebra $C^*_\lambda (F_2)$ is simple, 
the ideal $(\pi')^{-1}(\mathcal K (\Hi_\pi))$ is either $\{0\}$ or $C^*_\lambda (F_2)$;
but the second case is impossible because $C^*_\lambda(F_2)$ is unital
and $\Hi_\pi$ infinite-dimensional;
hence $\pi(C^*_{\rm max}(F_2)) \cap \mathcal K (\Hi_\pi) = \{0\}$,
and $\pi$ is not normal.

\vskip.2cm

\index{Heisenberg group! $2$@$H(\K)$ over a field $\K$}
(4)
A factor representation of type II$_\infty$ need not be normal.
For example, consider an infinite field $\K$
and the Heisenberg group $H(\K)$ 
as a nilpotent discrete group
(as in \ref{Sect-PrimIdealHeisenberg}).
Then $H(\K)$ has factor representations of type II$_\infty$
(Theorems \ref{discretTypes}~\ref{iDEdiscretTypes} and \ref{ThGlimm2}),
but every factor representation that is normal is of finite type
(Proposition \ref{Prop-NilGr}).

\vskip.2cm

(5)
A factor representation of type III is not normal.

\vskip.2cm

(6) 
Assume that $G$ is unimodular and second-countable.
By Theorem~\ref{Ex-HilbertAlgebra} (compare also Example~{Ex-HilbertAlgebra} (2)),
$\lambda_G(G)''$ carries a faithful, normal, semi-finite trace.
In other terms, the left regular representation $\lambda_G$ is traceable.
It follows that $G$ has a family of normal factor 
representations which separate the points of $L^1(G)$,
and in particular which separate the points of $G$; 
see Corollaire to Th\'eor\`eme 1 in Chap.~I, \S~3 of \cite{Guic--63}.
%(see \cite[Page 20]{Puka--99}).

\vskip.2cm

\index{Purely infinite! $3$@C*-algebra}
(7)
There are separable C*-algebras which are purely infinite, that is, 
which have no non-zero traces with non-zero ideal of definition; 
an example is the Cuntz algebra $O_2$ of \cite{Cuntz--77}
(it appears already in \cite[Example 2.1]{Dixm--64a}).
\par

For a LC group $G$, 
the maximal C*-algebra of $G$ has at least one normal representation,
the unit representation $1_G$.
It has many other ones when $G$ is unimodular and not $\{1 \}$,
see (6) above,
or when $G$ is a connected real Lie group,
see Theorem \ref{ThPukan} below.
The following question seems to be open.
\par

\textbf{Question:} Does there exist
a LC group $G$ such that $1_G$ is the only normal representation of $G$ ?
\end{exe}

Let $G$ be a LC group. 
The \textbf{normal quasi-dual} 
$\QD(G)_{\rm norm}$ of $G$
is the set of quasi-equivalence classes of normal factor representations of $G$.
\index{Normal quasi-dual}
\index{Quasi-dual! $2$@normal}
\index{$c3$@$\QD(G)_{\rm norm}$ 
normal quasi-dual of the LC group $G$}
\par

The following result from \cite{Halp--75} holds for general separable C*-algebras.
In case of a second-countable locally compact group $G$ that is not of type I,
it shows that $\QD(G)_{\rm norm}$ 
is better behaved than the dual $\widehat G$ 
or the quasi-dual $\QD(G)$.

\begin{theorem}
% 10.D.5
\label{ThHalpern}
Let $G$ be a second-countable locally compact group. 
The space $\QD(G)_{\rm norm}$, equipped with 
the Mackey--Borel structure induced by that of $\QD(G)$, 
is a standard Borel space.
\end{theorem}

Normal factor representations of LC group $G$ 
can be described in terms of characters.

\begin{defn}
% 10.D.6
\label{Char-LCG}
\index{Character! $4$@of a LC group}
A \textbf{character} of the LC group $G$ 
is a non-zero lower semi-continuous semi-finite trace $t$ on $C^*_{\rm max}(G)$, 
which satisfies the following condition: 
every lower semi-continuous semi-finite trace on $C^*_{\rm max}(G)$ dominated by $t$ 
is proportional to~$t$ on~$\overline{ \mathfrak m_{t,+}}$.
We denote by $\Char(G)$ 
the set of equivalence classes of characters of $G$, 
where two characters are equivalent if they are proportional to each other.
\index{$c4$@$\Char(G)$ 
equivalence classes of characters of the LC group $G$}
\end{defn}

Recall from Lemma \ref{Lemme-Trace-Representation} that,
if $(\pi, \tau)$ is a trace representation of $G$, 
then $\tau \circ \pi$ is a lower semi-continuous semi-finite trace
on $C^*_{\rm max}(G)$.
\par

Let $t$ be a lower semi-continuous semi-finite trace
on $C^*_{\rm max}(G)$, $t \ne 0$.
Recall from Proposition~\ref{Prop-TraceHilbertAlg} 
that the Hilbert algebra associated to $t$
is $\AC_t= \mathfrak n_t / N_t$, where
$$
N_t \, = \, \{x \in \mathfrak n_t \mid t(x^*x) = 0 \}.
$$
Recall also that there are corresponding commuting representations 
$\lambda_t$ and $\rho_t$ of $G$
on the Hilbert space completion $\Hi_t$ of $\AC_t$
and there exists a faithful normal trace $\tau_t$ on $\lambda_t(G)''$ 
such that $t = \tau_t \circ \lambda_t$. 
For $x \in \mathfrak n_t$, we denote by $[x]$ the image of $x$ in $\Hi_t$.
%\cite[6.7.3]{Dixm--C*}.
% Voir DC* p.127. 
\par

As we now show, every semi-continuous semi-finite trace on $C^*_{\rm max}(G)$
dominated by $t$ on $\overline{ \mathfrak m_{t,+}}$
defines an element in the centre of $\lambda_t(G)''$.

\begin{lem}
% 10.D.7
\label{Lem-CommutantTrace}
Let $t$ and $t'$ be lower semi-continuous semi-finite traces on $C^*_{\rm max}(G)$
such that $t'(x) \le t(x)$ for every $x \in \overline{ \mathfrak m_t }_+$.
\begin{enumerate}[label=(\arabic*)]
\item\label{iDELem-CommutantTrace}
There exists $T \in \lambda_t(G)' \cap \rho_t(G)'$ with $0 \le T \le I$ 
such that
$$
t'(y^*x) \, = \, \langle T[x] \mid [y] \rangle 
\hskip.5cm \text{for all} \hskip.2cm
x,y \in \mathfrak n_t.
$$
\item\label{iiDELem-CommutantTrace}
The function 
$$
t'' \, \colon \, x \mapsto \tau_t(T\lambda_t(x))
$$
is a semi-continuous semi-finite trace on $C^*_{\rm max}(G)$
such that $t''(x) \le t(x)$ and $t''(x) = t'(x)$ for every $x \in \overline{ \mathfrak m_t }_+$.
\end{enumerate}
\end{lem}

\begin{proof}
Since $t'(x) \le t(x)$ for every $x \in \overline{ \mathfrak m_{t,+}}$, the map
$$
\Phi \, \colon \, (x,y) \, \mapsto \, t'(y^*x)
$$
is a well defined positive Hermitian form on $\AC_t$, which is 
separately continuous for the Hilbert space norm (compare with Lemma~\ref{Lem-CommutantCyclic}).
There exists therefore $T \in \Li (\Hi_t)$ with $0 \le T \le I$ such that
$$
\Phi( x, y) \, = \, \langle T [x] \mid [y] \rangle
\hskip.5cm \text{for all} \hskip.2cm 
x,y \in \mathfrak n_t.
$$
Let us check that $T$ is in the centre of $\lambda_t(G)''$. 
For $z\in C^*_{\rm max}(G), x,y \in \mathfrak n_t$,
we have
$$
\begin{aligned}
\langle T\lambda_t(z)[x] \mid [y] \rangle
&\, = \, \langle T[zx] \mid [y] \rangle
\, = \, t'(y^*zx) \, = \, t'((z^*y)^*x)
\\
&\, = \, \langle T[x] \mid [z^*y] \rangle
\, = \, \langle T[x] \mid \lambda_t(z^*) [y] \rangle
\, = \, \langle \lambda_t(z) T[x] \mid \ [y] \rangle;
\end{aligned}
$$
this shows that $T \in \lambda_t(C^*_{\rm max}(G))' = \lambda_t(G)'$.
Since $t'$ is a trace, we also have
$$
\begin{aligned}
\langle T\rho_t(z)[x] \mid [y] \rangle
&\, = \, \langle T[xz^*] \mid [y] \rangle
\, = \, t'(y^*xz^*) \, = \, t'(z^*y^*x) \, = \, t'((yz)^*x)
\\
&\, = \, \langle T[x] \mid [yz] \rangle
\, = \, \langle T[x] \mid \rho_t(z^*) [y] \rangle
\, = \, \langle \rho_t(z) T[x] \mid \ [y] \rangle,
\end{aligned}
$$
and hence $T \in \rho_t(C^*_{\rm max}(G))' = \rho_t(G)'$ 
and this proves \ref{iDELem-CommutantTrace}.

\vskip.2cm

To show \ref{iiDELem-CommutantTrace}, 
let $P$ be the orthogonal projection on $\Ki = \overline{T(\Hi)}$. 
% \marginpar{?~$P$ ou $p$~?} 24/4 : d\'ecid\'e de garder $P$.
Then $\Ki$ is invariant under $\lambda_t(G)$ and $\rho_t(G)$ 
and hence $P$ is in the centre of $\lambda_t(G)''$.
\par

Let $\lambda_t^{(P)}$ be the representation of $G$ 
defined by restriction of $\lambda_t$ to $\Ki = P(\Hi)$.
Then
$$
\tau^{(P)} \, \colon \, S \mapsto \tau_t(TS) = \tau_t(ST)
$$
is a normal trace on the von Neumann algebra $\lambda_t^{(P)}(G)''$.
Moreover, $\tau^{(P)}$ is faithful; indeed, let $S \in \Li (\K)$ be such that $\tau(ST) = 0$; 
then $ST = 0$ since $\tau$ is faithful and hence $S = 0$.
\par

Let $x \in C^*_{\rm max}(G)_+$ be such that $\lambda_t(x^*x) \in \mathfrak n_{\tau_t}$. 
Then $\tau_t(\lambda_t(x^*x))<+\infty$ and hence 
$$
\tau^{(P)}(\lambda_t^{(P)}(x^*x)) \, = \, 
\tau_t(T\lambda_t(x^*x)) \, \le \,
\Vert T \Vert \tau_t(\lambda_t(x^*x)) \, < \,
+\infty.
$$
This shows that $(\lambda_t^{(P)}, \tau^{(P)})$ is a trace representation of $G$. 
Since $t'' = \tau^{(P)} \circ \lambda_t^{(P)}$, 
it follows from Lemma~\ref{Lemme-Trace-Representation}
that $t''$ is a lower semi-continuous semi-finite trace on $C^*_{\rm max}(G)$.
It is clear that $t'' \le t$ on $\overline{ \mathfrak m_t }_+$. 
\par

Let $x \in \mathfrak n_t$. 
Then $T[x]\in \Hi_t$ is a bounded element for the Hilbert algebra $\AC_t$,
by Lemma~\ref{LemHilbertAlg}~\ref{vDELemHilbertAlg}; 
hence, by Theorem~\ref{TheoHilbertAlg-Trace}, we have
$$
\tau_t(T \lambda_t(x)\lambda_t(x^*)) \, = \, 
\langle T^{1/2}[x] \mid T^{1/2}[x] \rangle \, = \,
\langle T[x] \mid [x] \rangle \, = \, t'(x^*x).
$$
This shows that $t''$ coincides with $t'$ on 
$\overline{\mathfrak n_t}= \overline{ \mathfrak m_t }$.
\end{proof}

Let $\pi_1$ and $\pi_2$ be normal factor representations of $G$
 and $\tau_1$ and $\tau_2$ semi-finite normal faithful traces on 
 $\pi_1(G)''$ and $\pi_2(G)''$. 
Observe that, if $\pi_1$ and $\pi_2$ are quasi-equivalent, 
then the trace $\tau_2 \circ \Phi$ is proportional to $\tau_1$.
Indeed, let $\Phi$ is an isomorphism between $\pi_1(G)''$ and $\pi_2(G)''$;
since $\pi_1(G)''$ is a factor, $\tau_2 \circ \Phi$ is proportional to 
$\tau_1$ (see Corollaire in \cite[Chap.~I, \S~6, 4]{Dixm--vN}).
\par

As a consequence, we see that $\pi \mapsto \tau \circ \pi$ is a well-defined
 map from $\QD(G)_{\rm norm}$ to 
the set of lower semi-continuous semi-finite trace on $C^*_{\rm max}(G)$, 
where two such traces are identified if they are proportional to each other.
The following proposition 
% from \cite[6.7.4]{Dixm--C*}, 
shows that this correspondence is in fact a bijection 
between quasi-equivalence classes of normal factor representations 
and characters of $G$.

\begin{prop}
% 10.D.8
\label{Char-NormalRep}
Let $G$ be a locally compact group.
The map $\pi \mapsto \tau \circ \pi$ described above induces a bijection 
$$
\QD(G)_{\rm norm} 
\hskip.2cm \overset{\approx}{\longrightarrow} \hskip.2cm
\Char(G) ,
$$
with inverse bijection induced by $t \mapsto \lambda_t$
of Proposition \ref{Prop-TraceHilbertAlg}.
\end{prop}

\begin{proof}
Let $\pi$ be a normal factor representation of $G$; 
let $\tau$ be a faithful normal trace on $\pi(G)''$ and $t = \tau \circ \pi$. 
Let $t'$ be a lower semi-continuous semi-finite trace on $C^*_{\rm max}(G)$
such that $t'(x) \le t(x)$ for every $x \in \overline{ \mathfrak m_t }_+$.
By Lemma~\ref{Lem-CommutantTrace}, 
there exists an operator $T$ in the centre of $\pi(G)''$
such that $t'(y^*x) = \langle T[x] \mid [y] \rangle$ for all $x,y \in \mathfrak n_t$.
Since $\pi(G)''$ is a factor, $T$ is a scalar multiple of $I$ and hence 
$t'$ is proportional to $t$ on $\overline{ \mathfrak m_t }$.
Therefore, $t$ is a character of $G$.
\par

Conversely, let $t$ be a character of $G$ and let $(\lambda_t, \tau_t)$ 
be the associated trace representation
on the Hilbert space $\Hi_t$.
We claim that $\lambda_t$ is factorial. 
Indeed, let $T$ be an operator in the centre of $\lambda_t(G)''$ with $0 \le T \le I$.
The map 
$$
t' \, \colon \, x \mapsto \tau_t(T\lambda_t(x))
$$
is a trace on $C^*_{\rm max}(G)$
By Lemma~\ref{Lemme-Trace-Representation}, 
$t'$ is lower semi-continuous and semi-finite.
Since $t = \tau_t \circ \lambda_t$, we have $t'(x) \le t(x)$ for every 
$x \in \overline{ \mathfrak m_t }_+$.
Therefore there exists $\alpha \in \mathopen[ 0,1 \mathclose]$ such that 
$t' = \alpha t$ on $\overline{ \mathfrak m_t }_+$.

Let $x \in \mathfrak n_t$. Then $T[x] \in \Hi_t$ is a bounded element
for the Hilbert algebra associated to $t$,
by Lemma~\ref{LemHilbertAlg}~\ref{vDELemHilbertAlg},
where $[x]$ denotes the image of $x$ in $\Hi_t$.
By Theorem~\ref{TheoHilbertAlg-Trace}, we have therefore
$$
\langle T[x] \mid [x] \rangle= t(T\lambda_t(x^*)\lambda_t(x))
\, = \,
t'(x^*x) \, = \, \alpha t(x^*x) \, = \, \alpha \langle [x], [x] \rangle.
$$
It follows that $T = \alpha I$ and hence that $\lambda_t(G)''$ is a factor.
\end{proof}

For the subclass of normal irreducible representations 
(see Proposition~\ref{Pro-IrredNormalRep}), we obtain the following corollary.

\begin{cor}
% 10.D.9
\label{Cor-Char-NormalRep}
Let $G$ be a locally compact group.
\begin{enumerate}[label=(\arabic*)]
\item\label{iDECor-Char-NormalRep}
Let $\pi$ be a normal irreducible representation of $G$.
The character of $\pi$ can be normalized to be $t_\pi := {\mathrm Tr} \circ \pi$ 
where ${\mathrm Tr}$ is the canonical trace on $\Li (\Hi_\pi)$;
the ideal of definition of $t_\pi$ is the set of $x \in C^*_{\rm max}(G)$
for which $\pi(x)$ is a trace class operator.
\item\label{iiDECor-Char-NormalRep}
Let $\pi$ and $\pi'$ be normal irreducible representations of $G$.
 Then $t_{\pi}$ and $t_{\pi'}$ are proportional 
if and only $\pi$ and $\pi'$ are equivalent.
\end{enumerate}
\end{cor}

\begin{proof}
Item \ref{iDECor-Char-NormalRep} 
follows from Proposition~\ref{Pro-IrredNormalRep} 
and from the fact that ${\mathrm Tr}$ 
is the unique normal semi-finite trace on $\Li (\Hi_\pi)$, up to positive scalar multiples.
\par

Item \ref{iiDECor-Char-NormalRep}
is a consequence of \ref{iDECor-Char-NormalRep}, 
of Proposition~\ref{Char-NormalRep}, 
and of the fact that two irreducible representations are equivalent 
if and only if they are quasi-equivalent.
\end{proof}

\begin{exe}
% 10.D.10
\label{ExCharacters}

(1)
Let $G$ be a compact group.
Characters of $G$ are finite.
It is common to normalize 
the character $t$ of a finite-dimensional representation $\pi$ of $G$
by $\Vert t \Vert = \dim \pi$.
where $\Vert \cdot \Vert$ denotes the norm 
on bounded linear forms on $C^*_{\rm max}(G)$.
Corollary \ref{Cor-Char-NormalRep} establishes the classical fact
that the map $\pi \mapsto \tau \circ \pi$ provides a bijection
from the dual $\widehat G$ of $G$ 
to the set $\Char(G)$ of its irreducible characters.
\par

Note that characters are defined here for factor representations only
whereas, in the setting of finite and compact groups, 
it is usual to define a character 
for every finite-dimensional representation of the group.
 
 \vskip.2cm

(2)
More generally,
let $G$ be a $\sigma$-compact locally compact group of type I.
Every irreducible representation is normal (Corollary~\ref{Cor-IrredNormalRep})
and so every character on $G$ is proportional to 
$t_\pi$ for a unique irreducible representation $\pi$ on $G$.
In this case again, the map $\pi\mapsto t_\pi$ provides a bijection 
$\widehat G \overset{\approx}{\longrightarrow} \Char(G)$.

\vskip.2cm

(3)
Let $\Gamma$ be a discrete group, not the one-element group.
The Dirac function $\delta_e$ at the group unit is of positive type and central.
The associated representation 
is the left regular representation $\lambda$ of $\Gamma$ in $\ell^2(\Gamma)$;
it is of finite type and traceable. 
The group $\Gamma$ is icc if and only if $\lambda$ is factorial
(Proposition \ref{iccfactorII1}).
When this is the case, $\lambda$ is of type II$_1$
and its character $C^*_{\rm max}(\Gamma) \to \C$
is given by $f \mapsto \langle \lambda(f) \delta_e \mid \delta_e \rangle$;
anticipating on the notation of Section \ref{Section-NormalRepFiniteType},
we have $\delta_e \in E(\Gamma)$. 

\vskip.2cm

(4) Let $G$ be a unimodular LC group,
and $\nu$ a complex-valued Radon measure on $G$
which is of positive type and central. 
As discussed in Example~\ref{Exemples-RepTracables},
$\nu$ defines a lower semi-continuous, semi-finite trace 
$t_\nu$ on $C^*_{\rm max}(G)$;
the associated representation $\pi_\nu$ is factorial (and hence normal) 
if and only if $t_\nu$ is a character.
This can be repeated for
a unimodular real Lie group $G$
and a distribution $\nu$ on $G$ 
which is of positive type and central.

\vskip.2cm

(5) Let $G$ be a semisimple or nilpotent connected Lie group. 
Every irreducible representation $\pi$ of $G$ is normal 
(since $G$ is of type I) and the corresponding character $t_\pi$
is given by a central distribution of positive type on $G$.
See \cite[\S~5]{Hari--54} and \cite[\S~7]{Kiri--62}.

\vskip.2cm

(6) 
Let $\Gamma$ be a countable group
given with a measure preserving action
on a standard probability space $(X, {\mathcal F}, \mu)$.
For every $\gamma \in \Gamma$, let
$$
X^\gamma \, = \, \{ x \in X \mid \gamma x = x \} 
$$
denote the corresponding fixed point set.
As we will show (Theorem~\ref{Theo-Vershik}
in Section~\ref{Section:IRS}), the function
%\cite[Theorem 9]{Vers--11} 
$$
\varphi \, \colon \, \Gamma \to \R, \hskip.2cm \gamma \mapsto \mu(X^\gamma).
$$
is central and of positive type;
therefore, $\varphi$ defines a finite trace $t_\varphi$ on $C^*_{\rm max}(\Gamma)$; 
see Section~\ref{Section-NormalRepFiniteType}.
\end{exe}

%-----------------------------------------------------------------------
% End of chapter 10
%-----------------------------------------------------------------------
\chapter{Finite characters and Thoma's dual}
% Chapter 11
\label{Thomadual}

\emph{
In this chapter, we study the finite part $\QD(G)_{\rm fin}$
of the quasi-dual $\QD(G)$ of a topological group $G$,
that is the space of quasi-equivalence classes
of factor representations of $G$ which are of finite type. 
}
\par

\emph{
When $G$ is locally compact,
$\QD(G)_{\rm fin}$ is in a bijective correspondence
with the set $\Char(G)_{\rm fin}$ of finite characters, that is,
characters which are everywhere defined on the maximal C*-algebra of $G$
(see Section~\ref{Section-NormalRepFiniteType}).
}
\par

\emph{
For an arbitrary topological group $G$, 
let $\Tr(G)$ be the set of traces of $G$, that is,
the set of continuous functions on $G$ which are central and of positive type. 
In Section~\ref{Section-GNS-Traces},
we show that the GNS-construction provides a surjective map from $\Tr(G)$
to the set of quasi-equivalence classes of representations $\pi$ of $G$
such that $\pi(G)$'' is a finite von Neumann algebra.
A crucial role is played in this context by the fact that
a Hilbert algebra is associated to every trace on $G$
(see Chapter~\ref{Section-RepTraceC*}).
As a consequence, we will see in Section~\ref{SThoma'sdual} that 
$\QD(G)_{\rm fin}$ is parametrized by the Thoma dual $E(G)$,
which is the set of extreme rays of $\Tr(G)$.
}
\par

\emph{
In Section~\ref{Section-CharactersPrimitive},
we review for a second countable locally compact group $G$
the relationship between
the quasi-dual spaces $\QD(G)_{\rm norm}$ and $\QD(G)_{\rm fin}$, on the one hand,
and the primitive dual $\Pri(G)$, on the other hand,
}

\section[Finite type and finite characters]
{Factor representations of finite type and finite characters}
% Section 11.A
\label{Section-NormalRepFiniteType}

\index{Finite part of the quasi-dual}
\index{Quasi-dual! $3$@finite part}
\index{$c5$@$\QD(G)_{\rm fin}$ finite part of $\QD(G)$}
\index{Representation! finite type}
Let $G$ be a topological group. 
Recall that the type of a representation $\pi$ of $G$ only depends
on the quasi-equivalence class of $G$ (see Remark~\ref{Rem-Def-TypesRep}).

\begin{defn}
% 11.A.1
\label{Def-FiniteQD}
The \textbf{finite part of the quasi-dual} $\QD(G)_{\rm fin}$
is the space of quasi-equivalence classes of factor representations $\pi$ of $G$
which are of finite type, that is, 
such that $\pi(G)''$ is either a finite-dimensional factor or a factor of type II$_1$.
\end{defn}
\par

Let $G$ a locally compact group. 
The set $\Char(G)$ of characters
of $G$ has been introduced in Definition~\ref{Char-LCG}.
A character $t \in \Char(G)$ is \textbf{finite} if $t$ is a finite trace 
on the maximal C*-algebra $C^*_{\rm max}(G)$ of $G$.
We write $\Char(G)_{\rm fin}$ for the subset of $\Char(G)$
consisting of the finite characters $t$ 
normalized by $\Vert t \Vert = 1$. 
\par

The following fact is a consequence of Proposition \ref{Char-NormalRep}.
\index{Character! $5$@finite}
\index{$c6$@$\Char(G)_{\rm fin}$ finite part of $\Char(G)$}

\begin{prop}
% 11.A.2
\label{bijquasidualfinicharfini}
Let $G$ be a locally compact group.
\par

The map $\pi \mapsto \tau \circ \pi$ of Proposition \ref{Char-NormalRep} 
induces a bijection
$$
\QD(G)_{\rm fin}
\hskip.2cm \overset{\approx}{\longrightarrow} \hskip.2cm
\Char(G)_{\rm fin}.
$$
\end{prop}

Let $\mu$ be a Haar measure on the LC group $G$.
Let $t$ be a finite trace on $C^*_{\rm max}(G)$.
The restriction $t \vert_{L^1(G, \mu)}$ is a positive linear functional on $L^1(G,\mu)$
and is therefore defined by a unique continuous function of positive type $\varphi_t$ on $G$,
such that $\varphi_t(e) = \Vert t \Vert$
(see Remark \ref{bijftp} (2) and Example \ref{exGNSpourC*} (2)).
Since $t$ is a trace, $\varphi_t$ is central that is,
$\varphi_t$ is invariant under conjugation by elements of $G$.
\par

Conversely, let $\varphi$ be a central continuous function of positive type on $G$.
Then $\varphi$ defines a positive linear functional on $L^1(G,\mu)$
which extends uniquely to a positive linear functional $t_\varphi$ of $C^*_{\rm max}(G)$
with $\Vert t_\varphi \Vert = \varphi(e)$.
Since $\varphi$ is central, $t_\varphi$ is a trace on $C^*_{\rm max}(G)$.
\par

Observe that the maps $t \mapsto \varphi_t$ and $\varphi \mapsto t_\varphi$ 
are continuous affine maps which are inverse to each other.
\par

Observe also that the set of central continuous functions of positive type $\varphi$ on $G$
with $\Vert \varphi \Vert_\infty = \varphi(e) \le 1$ 
is a convex subset of $L^\infty (G,\mu)$, which is compact for the weak$^*$-topology.
As a result of the discussion above, we obtain the following result.

\begin{prop}
% 11.A.3
\label{Prop-CharFini-CentralPD}
Let $G$ be a locally compact group.
\begin{enumerate}[label=(\arabic*)]
\item\label{iDEProp-CharFini-CentralPD}
The map $t \mapsto \varphi_t$ described above is a bijection
between the set of finite traces on $C^*_{\rm max}(G)$
and the set of central continuous functions of positive type on $G$.
\item\label{iiDEProp-CharFini-CentralPD}
The map $t \mapsto \varphi_t$ restricted to $\Char(G)_{\rm fin}$ is a bijection
between $\Char(G)_{\rm fin}$
and the set of central continuous functions of positive type $\varphi$ on $G$
such that $\varphi(e) = 1$, and which are extremal with these properties.
\end{enumerate}
\end{prop}

Propositions~\ref{bijquasidualfinicharfini} and \ref{Prop-CharFini-CentralPD}
yield a parametrization of $\QD(G)_{\rm fin}$
in terms of central continuous functions of positive type on $G$,
when $G$ is locally compact.
\par

When $G$ is not locally compact, the space $\Char(G)_{\rm fin}$ is not defined.
Nevertheless, as we will see in the next sections, the parametrization of $\QD(G)_{\rm fin}$
in terms of central continuous functions of positive type on $G$ 
carries over to an arbitrary topological group $G$.

\section
{GNS Construction for traces on groups}
% Section 11.B
\label{Section-GNS-Traces}

Let us return to functions of positive type and the GNS construction,
keeping the notation of Section \ref{S-FPosType}.
\par

Let $G$ be a topological group (not necessarily a locally compact group).
Recall that a function $\varphi$ on $G$ is \textbf{central}
if $\varphi(gh) = \varphi(hg)$ for all $g, h \in G$.

\begin{defn}
% 11.B.1
\label{Def-Trace-FiniteRep}
\index{Trace! $6$@on a topological group}
(1)
A \textbf{trace} on $G$ is
a continuous function $\varphi \,\colon G \to \C$ 
that is central and of positive type.
\par

Denote by $\Tr(G)$ the space of traces on $G$.
Let also $\Tr_{\le 1}(G)$ [respectively $\Tr_{1}(G)$]
be the subset of $\Tr(G)$ of functions
such that $\Vert \varphi \Vert_\infty = \varphi(e) \le 1$
[resp.\ $\varphi(e) = 1$].
\index{$d2$@$\Tr(G), \Tr_{\le 1}(G), \Tr_1(G)$ traces on $G$}

\vskip.2cm

(2)
A representation $\pi$ of $G$ is of \textbf{ finite type}
if the von Neumann algebra $\pi(G)''$ has a faithful finite normal trace.
\index{Representation of finite type}
\end{defn}

The term ``trace" as defined above was coined in \cite{CaMo--84};
traces are often called ``characters" in the literature
(see, for instance, \cite{DuMe--14} and \cite{PeTh--16}).

\begin{rem}
% 11.B.2
\label{Rem-FiniteTypeRep}
(1)
Let $\pi$ be representation of $G$.
Assume that $\pi(G)''$ is countably decomposable
(this is for instance the case if the Hilbert space of $\pi$ is separable).
Then $\pi(G)''$ is of finite type if and only if $\pi(G)''$ is a finite von Neumann algebra
(see Remark~\ref{Rem-TypI-Semifinite}).

\vskip.2cm

(2)
Let $\pi$ be a finite type representation of $G$ and assume that $\pi(G)''$ is 
a factor. Then $\pi(G)''$ 
has a \emph{unique} normalized finite normal trace;
in fact, there is a unique linear functional $\tau$ on $\pi(G)''$ such 
that $\tau(xy)=\tau(yx)$ for all $x,y\in \pi(G)''$ and $\tau(I)=1$.
(see Corollaire in \cite [Chap.~III, \S~5, no~1]{Dixm--vN})
\end{rem}

\begin{exe}[compare with Definition \ref{exposfunctdiagcoeff}]
% 11.B.3
\label{excentralposfuncttrace}
Let $\pi$ be a finite type representation of $G$
and $\tau$ a faithful finite normal trace on $\pi(G)''$.
Then the continuous function $\varphi_{\pi, \tau}$ defined by 
$$
\varphi_{\pi, \tau} (g) \, = \, \tau ( \pi(g))
\hskip.5cm \text{for all} \hskip.2cm
g \in G
$$
is clearly central and non-zero.
It is also of positive type; indeed,
for $n \ge 1$, $c_1, \hdots, c_n \in \C$, and $g_1, \hdots, g_n \in G$, we have
$$
\sum_{i, j = 1}^n c_i \overline{c_j} \varphi_{\pi, \tau}(g_j^{-1}g_i) \, = \, 
\tau \Big( 
\big( \sum_{i = 1}^n c_i \pi(g_i) \big) \big( \sum_{j = 1}^n c_j \pi(g_j) \big)^*
\Big) \, \ge \, 0 .
$$
Therefore $\varphi_{\pi, \tau} \in \Tr(G)$.
If $\tau$ is normalized, then $\varphi_{\pi, \tau} \in \Tr_1(G)$.
\end{exe}

\begin{constr}[Compare with Construction \ref{constructionGNS2}]
% 11.B.4
\label{constructionGNS12}

The present construction shows
that every non-zero element in $\Tr(G)$ arises as in the previous example.
More precisely, it shows that,
given $\varphi \in \Tr(G)$, $\varphi \ne 0$, 
there exists a representation $\pi_\varphi$ of $G$
(the GNS representation of \ref{constructionGNS2})
and a faithful finite normal trace $\tau_\varphi$ on $\pi_\varphi(G)''$
such that $\varphi (g) = \tau_\varphi(\pi_\varphi(g))$ for all $g \in G$.
\index{Gel'fand--Naimark--Segal representation! $3$@with $\varphi \in \Tr(G)$}
\par

\index{Gel'fand--Naimark--Segal representation!}
Let $\varphi \in \Tr(G)$, $\varphi \ne 0$. 
Recall from Section~\ref{S-FPosType} how
the corresponding GNS triple $(\pi_\varphi, \Hi_\varphi, \xi_\varphi)$ 
is constructed. 
The set $\C[G]$ of complex-valued functions on $G$ with finite support
is a $*$-algebra endowed with a positive Hermitian form $\Phi_\varphi$,
and the set $J_\varphi$ of functions $f$ such that $\Phi_\varphi(f,f) = 0$
is a two-sided $*$-ideal.
We denote by $\AC_\varphi$ the quotient $\C[G] / J_\varphi$.
We write $\delta_g \in \C[G]$ the characteristic function of $\{g\}$
for each $g \in G$, and $[f]$ for the class in $\AC_\varphi$
of a function $f \in \C[G]$.
We denote by $\Hi_\varphi$ the Hilbert space
obtained by completion of $\AC_\varphi$,
with respect to the scalar product induced on $\AC_\varphi$ by $\Phi_\varphi$,
and by $\xi_\varphi$ the vector $[\delta_e]$ viewed in $\Hi_\varphi$.
For $g \in G$, the operator $\pi_\varphi(g)$ on $\Hi_\varphi$
is defined as in \ref{constructionGNS2} 
to be the continuous linear extension to $\Hi_\varphi$
of the left translation $[f] \mapsto [g \cdot f]$,
where $(g \cdot f) (x) = f(g^{-1}x)$ for $g,x \in G$ and $f \in \AC_\varphi$.
\par

As $\varphi$ is central, $J_\varphi$ is now a \emph{two-sided} $*$-ideal,
the quotient $\AC_\varphi = \C[G] / J_\varphi$ is an algebra,
and there is for all $g \in G$ an operator $\rho_\varphi(g)$ on $\Hi_\varphi$
which is the continuous linear extension
of the \emph{right translation} $[f] \mapsto [f \cdot g]$,
where $(f \cdot g) (x) = f(xg)$ for $g,x \in G$ and $f \in \AC_\varphi$. 
One checks, as for $\pi_\varphi$ in Section~\ref{S-FPosType}, 
that the assignment $g \mapsto \rho_ \varphi(g)$
is continuous, and is a representation of $G$ in $\Hi_ \varphi$
with cyclic vector $\xi_\varphi = [\delta_e]$.
We have $\pi_\varphi(g)\xi_\varphi = [\delta_g]$ and 
$$
\pi_\varphi(g)\rho_\varphi(h) \, = \, \rho_\varphi(h)\pi_\varphi(g)
\hskip.5cm \text{for all} \hskip.2cm
g, h \in G.
$$
Moreover, the map 
$$
\tau_\varphi \, \colon \, T \mapsto \langle T\xi_\varphi \mid \xi_\varphi \rangle
$$
is a finite normal trace on $\pi_\varphi(G)''$ and 
$\varphi = \tau_\varphi \circ \pi_\varphi$.
\par

To check that $\tau_\varphi$ is faithful,
let $T\in \pi_\varphi(G)''$ be such that $\tau_\varphi (T^*T) = 0$. Then 
$$
\Vert T\xi_\varphi \Vert^2 \, = \, \langle T^*T\xi_\varphi \mid \xi_\varphi \rangle \, = \, 0
$$
and so $T\xi_\varphi=0$.
Since obviously $\pi_\varphi(G)''\subset \rho_ \varphi(G)'$
and since $\xi_\varphi$ is a cyclic vector for $\rho_ \varphi$,
it follows that $T (\rho_ \varphi(g) \xi_\varphi) = 0$ for all $g \in G$
and hence that $T = 0$. 
\par

Note: For coherence with Section \ref{Section-RepTraceC*},
we could have written $\lambda_\varphi$ instead of $\pi_\varphi$.
We have chosen the latter for coherence with Section \ref{S-FPosType}.
\end{constr}

We are going to show
that the von Neumann algebras $\pi_\varphi(G)''$ and $\rho_\varphi(G)''$
are each other commutants.

\begin{prop}
% 11.B.5
\label{Lem-HilbertAlgebra}
Let $G$ be a topological group and $\varphi \in \Tr(G)$, $\varphi \ne 0$.
Let $\AC_\varphi, \Phi_\varphi, \pi_\varphi$ and $\rho_\varphi$ 
be as in \ref{constructionGNS12}.
% \begin{enumerate}[label=(\roman*)]
\begin{enumerate}[label=(\alph*)]
\item\label{iDELem-HilbertAlgebra}
The $*$-algebra $\AC_\varphi$, 
equipped with the scalar product induced by $\Phi_\varphi$, is a Hilbert algebra.
\item\label{iiDELem-HilbertAlgebra}
We have $\pi_\varphi(G)'' = \rho_\varphi(G)'$ 
and $\rho_\varphi(G)'' = \pi_\varphi(G)'$.
Consequently,
$\pi_\varphi(G)' \cap \rho_\varphi(G)'$ is the common centre 
of the von Neumann algebras $\pi_\varphi(G)''$ and $\rho_\varphi(G)''$.
\end{enumerate}
\end{prop} 

\begin{proof}
\ref{iDELem-HilbertAlgebra} 
Property (1) of the definition of a Hilbert algebra 
(Section~\ref{Section-RepTraceC*}),
follows from the fact that $\varphi$ is central.
It is straightforward to check Properties (2) and (3),
and Property (4) is obvious since $\AC_\varphi$ is unital.
\par

Claim \ref{iiDELem-HilbertAlgebra}
follows from \ref{iDELem-HilbertAlgebra} 
and from the Commutation Theorem \ref{TheoHilbertAlg-Commutation}.
\end{proof}

Let $\pi$ be a finite type representation of $G$ 
and $\tau$ a normalized faithful finite normal trace on $\pi(G)''$. 
Recall from Example \ref{excentralposfuncttrace}
that $\varphi_{\pi, \tau} = \tau \circ \pi$ belong to $\Tr_1(G)$.

\vskip.2cm

An equivalence relation was defined in Section\ref{Section-TraceRep}
for arbitrary trace representations in the case where $G$ is \emph{locally compact}. 
We extend this definition as follows for finite type representations
of a general topological group.

\begin{defn}
% 11.B.6
\label{Def-QE-FiniteTypeRep}
For $i = 1, 2$, let $\pi_i$ be a finite type representation of the topological group $G$
and let $\tau_i$ be a normalized faithful finite normal trace on the von Neumann algebra $\pi_i(G)''$.
We say that the pairs $(\pi_1, \tau_1)$ and $(\pi_2, \tau_2)$ are \textbf{quasi-equivalent} 
if there exists an isomorphism $\Phi \,\colon \pi_1(G)'' \to \pi_2(G)''$ 
such that $\Phi(\pi_1(g)) = \pi_2(g)$ for all $g \in G$, 
and such that $\tau_1 = \tau_2 \circ \Phi$.
It is clear in this case that the associated functions 
$\varphi_{\pi_1, \tau_1} = \tau_1 \circ \pi_1$
and $\varphi_{\pi_2, \tau_2} = \tau_2 \circ \pi_2$ 
in $\Tr_1(G)$ are equal.
\index{Quasi-equivalent finite type representations}
\end{defn}

The next proposition, which appears in \cite{HiHi--05},
shows that the converse statement is also true.
Observe that, when $G$ is locally compact, this proposition 
follows easily from Proposition~\ref{Prop-CharFini-CentralPD}
combined with Proposition \ref{traceablerepbij}.

\begin{prop}[compare with Proposition \ref{GNSbijP(G)cyclic}]
% 11.B.7
\label{GNSbijTr(G)trace}
Let $G$ be a topological group.
\par

The assignments $\varphi \rightsquigarrow ( \pi_\varphi, \tau_\varphi)$ 
given by the GNS construction and 
$(\pi, \tau) \mapsto \tau \circ \pi$
induce bijections inverse to each other between
\begin{enumerate}[label=(\roman*)]
\item\label{iDEGNSbijTr(G)trace}
the set $\Tr_1(G)$ of normalized central functions of positive type $\varphi$,
and 
\item\label{iiDEGNSbijTr(G)trace}
the set of quasi-equivalence classes of pairs $(\pi, \tau)$, 
where $\pi$ is a finite type representation of $G$ 
and $\tau$ a normalized faithful finite normal trace on $\pi(G)''$.
\end{enumerate}
\end{prop}

\begin{proof}
We first recall some facts about the standard representation 
(see Section~\ref{Section-StandardRep})
in the special case of a \textbf{finite} von Neumann algebra. 
\par

Let $\mathcal M$ be a finite von Neumann algebra,
with a faithful normalized finite normal trace $\tau$.
Equipped with the scalar product $(x, y) \mapsto \tau (y^*x)$, 
the space $\mathcal M$ is a Hilbert algebra. 
Denote by $L^2(\mathcal M, \tau)$ the Hilbert space completion of $\mathcal M$
and let
$$
\mathcal M \, \to \, \widehat{\mathcal M},
\hskip.5cm
x \, \mapsto \, \widehat x,
$$
be the standard representation of $(\mathcal M, \tau)$ on $L^2(\mathcal M, \tau)$.
It is clear that the image $\xi$ of $I$ in $L^2(\mathcal M, \tau)$ 
is a cyclic vector for the representation $x \mapsto \widehat x$;
the natural trace $\widehat \tau$ on $\widehat{\mathcal M}$ is given by 
$$
\widehat \tau (\widehat x) \, = \, \langle \widehat x \xi \mid \xi \rangle \, = \, \tau(x), 
\hskip.5cm \text{for all} \hskip.5cm
x \in \mathcal M.
$$

\vskip.2cm 

For $j = 1, 2$, let $\pi_j \,\colon G \to \U(\Hi_j)$ be a representation 
such that the von Neumann algebra $\mathcal M_j := \pi_j(G)''$ 
has a faithful normalized finite normal trace $\tau_j$.
Assume that $\varphi_{\pi_1, \tau_1} = \varphi_{\pi_2, \tau_2}$.
We have to show that
the pairs $(\pi_1, \tau_1)$ and $(\pi_2, \tau_2)$ are quasi-equivalent.
\par

Still for $j = 1, 2$,
consider the standard representation $x \mapsto \widehat x$ 
of the Neumann algebra $\mathcal M_j$ on $L^2(\mathcal M_j, \tau_j)$.
The map $\Phi_j \,\colon x \mapsto \widehat x$
is an isomorphism between the von Neumann algebras $\mathcal M_j$ 
and $\widehat{\mathcal M_j} := \{\widehat x \, \colon \, x \in \mathcal{M_j} \}$,
and $\widehat{\tau_j} \circ \Phi_j = \tau_j$,
where $\widehat{\tau_j}$ is the natural trace on $\widehat{\mathcal M_j}$.
Therefore the pair $(\pi_j, \tau_j)$ is quasi-equivalent to the 
pair $(\widehat{\pi_j}, \widehat{\tau_j})$,
where $\widehat{\pi_j}$ is the representation of $G$ on $L^2(\mathcal M_j, \tau_j)$
defined by $\widehat{\pi_j}(g) = \widehat{\pi_j(g)}$ for $g \in G$.
As a consequence, it suffices to show that the pairs
$(\widehat{\pi_1}, \widehat{\tau_1})$ and $(\widehat{\pi_2}, \widehat{\tau_2})$
are quasi-equivalent.
\par

Observe that, by assumption, we have 
$$
\widehat{\tau_1}(g) \, = \, 
\varphi_{\pi_1, \tau_1}(g) \, = \, 
\varphi_{\pi_2, \tau_2}(g) \, = \,
\widehat{\tau_2}(g),
$$
that is, 
$$
\langle \widehat{\pi_1}(g)\xi_1 \mid \xi_1 \rangle
\, = \,
\langle \widehat{\pi_2}(g)\xi_2 \mid \xi_2 \rangle ,
$$
for all $g \in G$.
Since, for all $n \ge 1$, $g_1, \hdots, g_n \in G$ and $c_1, \hdots, c_n \in \C$,
we have
$$
\Big\Vert \sum_{i = 1}^n c_i \widehat{\pi_2}(g_i) \xi_2 \Big\Vert^2 
\, = \, 
\sum_{i, j = 1}^n c_i \overline{c_j} \widehat{\tau_2}(g_j^{-1} g_i)
\, = \, 
\sum_{i, j = 1}^n c_i \overline{c_j} \widehat{\tau_1}(g_j^{-1} g_i)
\, = \,
\Big\Vert \sum_{i = 1}^n c_i \widehat{\pi_1}(g_i) \xi_1 \Big\Vert^2,
$$
the map 
$\sum_{i = 1}^n c_i \widehat{\pi_1}(g_i) \xi_1 
\mapsto 
\sum_{i = 1}^n c_i \widehat{\pi_1}(g_i) \xi_2$ 
is well defined on the linear span of $\widehat{\pi_1}(G)\xi_1$ 
and extends to a Hilbert space isomorphism 
$U \,\colon L^2(\mathcal M_1, \tau_1) \to L^2(\mathcal M_2, \tau_2)$.
Let 
$$
\Phi \, \colon \, 
\Li (L^2(\mathcal M_1, \tau_1)) \to \Li (L^2(\mathcal M_2, \tau_2)), 
\hskip.5cm 
T \, \mapsto \, UTU^{-1}.
$$
Then $\Phi$ restricts to an isomorphism 
from $\widehat{\mathcal M_1}$ to $\widehat{\mathcal M_2}$,
such that $\Phi( \widehat{\pi_1}(g)) = \widehat{\pi_2}(g)$ for every $g \in G$.
Moreover, we have $\widehat{\tau_1}= \widehat{\tau_2}\circ \Phi$. 
Therefore $(\widehat{\pi_1}, \widehat{\tau_1})$ and $(\widehat{\pi_2}, \widehat{\tau_2})$ 
are quasi-equivalent, and this concludes the proof.
\end{proof}

\section
{Thoma's dual}
% Section 11.C
\label{SThoma'sdual}

\index{Indecomposable! $3$@for $\varphi \in \Tr_{\le 1}(G)$}
\index{$d3$@${\rm Extr} (\cdot)$ indecomposable elements}
Let $G$ be a topological group. 
Define ${\rm Extr} (\Tr_{\le 1}(G))$ 
as the set of functions $\varphi$ in $\Tr_{\le 1}(G)$ 
that are \textbf{indecomposable} 
in the following sense: 
if $\varphi = t \varphi_1+ (1-t)\varphi_2$ 
for $\varphi_1, \varphi_2 \in \Tr_{\le 1}(G)$ and $t \in \mathopen] 0,1 \mathclose[$, 
then $\varphi = \varphi_1 = \varphi_2$.
Define similarly ${\rm Extr} (\Tr_1(G))$, and call it
the \textbf{Thoma dual} of $G$.
From now on, we use Thoma's notation, $E(G)$, rather than ${\rm Extr} (\Tr_1(G))$.
\index{Thoma's dual} 
\index{$d4$@$E(G)$ Thoma's dual of the topological group $G$}
\par

Thoma used properties of the space $E(G)$ 
in order to achieve his characterization of discrete groups which are of type I,
and he advocated its use as a substitute
for the dual of discrete groups \cite{Thom--64a, Thom--64b}.
\par

The following lemma, 
which relates $E(G)$ and ${\rm Extr} (\Tr_{\le 1}(G))$,
is proved exactly in the same way as Lemma \ref{Lem-ExtP}, which relates
for ${\rm Extr} (P_{\le 1}(G))$ and ${\rm Extr} (P_1(G)) \cup \{0\}$.

\begin{lem}[compare with Lemma~\ref{Lem-ExtP}]
% 11.C.1
\label{Lem-ExtTr}
For every topological group $G$, we have
$$
{\rm Extr} (\Tr_{\le 1}(G)) \, = \, E(G)\cup \{0\}.
$$
\end{lem}

Proposition \ref{encoreblablacentralprop} provides
a necessary and sufficient condition on $\varphi$ 
for the GNS representation $\pi_\varphi$ 
of Sections \ref{S-FPosType} and \ref{Section-GNS-Traces}
to be factorial.
\par
 
Let $\psi \in \Tr(G)$ be such that $\psi \le \varphi$, that is, 
$\varphi - \psi \in P(G)$;
it is then obvious that $\varphi - \psi$ belongs to $\Tr(G)$. 
For $T \in \pi_ \varphi(G)'$ with $0 \le T \le I$, define $\varphi_T \in P(G)$ 
as in Section~\ref{S-FPosType} by
$$
\varphi_T(g) 
\, = \, \langle \pi_\varphi(g)T\xi_\varphi \mid \xi_\varphi \rangle
\, = \, \langle \pi_ \varphi(g)T^{1/2} \xi_ \varphi \mid T^{1/2} \xi_ \varphi \rangle .
$$

\begin{lem}[compare with Lemma \ref{Lem-CommutantCyclic}]
% 11.C.2
\label{Lem-CommutantCyclic-Central}
Let $\varphi \in {\rm Extr} (\Tr(G))$;
let $\pi_\varphi$ and $\rho_\varphi$ be as in \ref{constructionGNS12}.
\par

The map $T \to \varphi_T$ is a bijective correspondence 
from the set of $T \in \pi_ \varphi(G)' \cap \rho_ \varphi(G)'$ with $0 \le T \le I$ 
to the set of $\psi \in \Tr(G)$ with $\psi \le \varphi$. 
\end{lem}

\begin{proof}
For $T \in \pi_ \varphi(G)' \cap \rho_ \varphi(G)'$ with $0 \le T \le I$,
the function $\varphi_T$ is central;
indeed, for $g,x \in G$, we have
$$
\begin{aligned}
\varphi_T(x^{-1}gx)
\, &= \, 
\langle \pi_ \varphi(g)\pi_ \varphi(x)T\xi_ \varphi 
 \mid \pi_ \varphi(x)\xi_ \varphi \rangle
\, = \, 
\langle \pi_ \varphi(g)T\pi_ \varphi(x)\xi_ \varphi 
 \mid \pi_ \varphi(x)\xi_ \varphi \rangle
\\
\, &= \, 
\langle \pi_ \varphi(g)T\pi_ \varphi(x)\xi_ \varphi 
 \mid \rho_ \varphi(x^{-1})\xi_ \varphi \rangle
\, = \, 
\langle \pi_ \varphi(g)T\pi_ \varphi(x)\rho_ \varphi(x)\xi_ \varphi 
 \mid \xi_ \varphi \rangle
\\
\, &= \, 
\langle \pi_ \varphi(g)T\xi_ \varphi 
 \mid \xi_ \varphi \rangle 
\, = \,
\varphi_T(g) .
\end{aligned}
$$
Moreover, since $\varphi - \varphi_T = \varphi_{I - T}$, 
we have $\varphi_T \le \varphi$.
\par

It follows from Lemma~\ref{Lem-CommutantCyclic}
that the map $T \mapsto \varphi_T$ is injective.
\par

Let $\psi \in \Tr(G)$ with $\psi \le \varphi$.
The proof of Lemma~\ref{Lem-CommutantCyclic} shows that 
there exists $T \in \pi_\varphi(G)'$ with $0 \le T \le I$ such that $\psi = \varphi_T$.
It remains to show that $T \in \rho_\varphi(g)'$.
\par

We have, for $x, y, z \in G$,
$$
\begin{aligned}
\psi(y(z^{-1}x)y^{-1})
&\, = \, \langle \pi_\varphi(y(z^{-1}x)y^{-1})\xi_ \varphi 
 \mid T\xi_ \varphi \rangle
\\
&\, = \, \langle \pi_\varphi(z^{-1}x) \pi_\xi(y^{-1})\xi_ \varphi 
 \mid T\pi_\varphi(y^{-1})\xi_ \varphi \rangle
\\
&\, = \, \langle \pi_\varphi(z^{-1}x) \rho_\xi(y)\xi_ \varphi 
 \mid T\rho_\varphi(y)\xi_ \varphi \rangle
\\
&\, = \, \langle \pi_\varphi(x) \xi_ \varphi 
 \mid \rho_\xi(y^{-1})T\rho_\varphi(y)\pi_\varphi(z)\xi_ \varphi \rangle .
\end{aligned}
$$
However, since $\psi$ is central, we also have 
$$
\psi(y(z^{-1}x)y^{-1}) \, = \,
\psi(z^{-1}x) \, = \,
\langle \pi_\varphi(x) \xi_ \varphi \mid T\pi_\varphi(z)\xi_ \varphi \rangle.
$$
Because $\xi_ \varphi$ is a cyclic vector, 
it follows that $\rho_\varphi(y^{-1})T\rho_\varphi(y) =T$
for all $y \in G$, that is, $T \in \rho_\varphi(G)'$.
\end{proof}

\begin{prop}[compare with Proposition \ref{P1=cyclicmodeqIRR}]
% 11.C.3
\label{encoreblablacentralprop}
Let $G$ be a topological group.
Let $\varphi \in \Tr_1(G)$, and $\pi_\varphi$ be the corresponding GNS representation,
as in \ref{constructionGNS12} and \ref{GNSbijTr(G)trace}.
\par

Then $\pi_\varphi$ is factorial if and only if $\varphi \in E(G)$. 
\end{prop}

\begin{proof}
Assume that $\pi_\varphi$ is factorial;
let $\varphi_1, \varphi_2 \in \Tr_1(G)$ and $t \in \mathopen] 0,1 \mathclose[$
be such that $\varphi = t\varphi_1 + (1-t)\varphi_2$.
Let $j \in \{1, 2 \}$; since $\varphi_j \le \varphi$,
there exists by Lemma \ref{Lem-CommutantCyclic-Central}
$T_j \in \pi_\varphi(G)' \cap \rho_\varphi(G)'$ with $0 \le T_j \le I$
such that $\varphi_j = \varphi_{T_j}$.
Since $\pi_\varphi$ is factorial, 
$\pi_\varphi(G)' \cap \rho_\varphi(G)' = \pi_\varphi(G)' \cap \pi_\varphi(G)''
= \mathrm{Id}_{\Hi_\varphi}$
(see Proposition \ref{Lem-HilbertAlgebra}),
so that both $T_1$ and $T_2$ are scalar multiples of the identity.
It follows that $\varphi_1 = \varphi_2 = \varphi$, 
and consequently that $\varphi \in E(G)$.
\par

For the converse implication, we leave it to the reader to adapt
the proof of Proposition \ref{P1=cyclicmodeqIRR}.
\end{proof}

For $i = 1, 2$, let $\pi_i$ be a \emph{factor} representation of finite type of $G$
and let $\tau_i$ be the unique normalized trace on $\pi_i(G)''$
(see Remark~\ref{Rem-FiniteTypeRep}.
If the representations $\pi_1$ and $\pi_2$ are quasi-equivalent, 
then the pairs $(\pi_1, \tau_1)$ and $(\pi_2, \tau_2)$ are quasi-equivalent.
Indeed, let $\Phi$ be a quasi-equivalence between $\pi_1$ and $\pi_2$; 
then $\tau_2 \circ \Phi$ is a normalized finite normal trace on $\pi_1(G)''$,
and hence $\tau_2 \circ \Phi = \tau_1$.
\par

The following result is therefore a consequence
of Proposition~\ref{encoreblablacentralprop}.
Observe that Item \ref{iiDEThm_NormalRepFiniteType}
was already established in Proposition~\ref{Prop-CharFini-CentralPD}.

\begin{theorem}
% 11.C.4
\label{Thm_NormalRepFiniteType}
Let $G$ be a topological group.
For $\varphi \in E(G)$, let $\pi_\varphi$ and $\tau_\varphi$
be as in \ref{constructionGNS12} and \ref{GNSbijTr(G)trace}.
\begin{enumerate}[label=(\arabic*)]
\item\label{iDEThm_NormalRepFiniteType}
The map 
$$
E(G) \to \QD(G)_{\rm fin},
\hskip.2cm
\varphi \mapsto \pi_\varphi
$$
establishes a bijection from the Thoma dual 
onto the set of quasi-equivalence classes
of factor representations of finite type of $G$.
\item\label{iiDEThm_NormalRepFiniteType}
Assume that $G$ is locally compact. 
The map 
$$
E(G) \to \Char(G)_{\rm fin}, 
\hskip.2cm
\varphi \mapsto t_\varphi = \tau_\varphi \circ \pi_\varphi
$$
is a bijection from the Thoma dual 
onto the set of finite normalized characters of $G$.
\end{enumerate}
\end{theorem}
\index{Character! $5$@finite}

\begin{exe}
% 11.C.5
\label{exThomaNonDiscrete}
Let $G$ be a topological group.

\vskip.2cm

(1)
When $G$ is locally compact and abelian,
the Thoma dual $E(G)$ is naturally identified with the Pontrjagin dual $\widehat G$.
When $G$ is compact, $E(G)$ is the set of normalized traces
of irreducible representations of $G$;
these are particular cases of Theorem \ref{Thm_NormalRepFiniteType}.

\vskip.2cm

(2)
When $G$ is a $\sigma$-compact LC group of type I, 
$E(G)$ is the set of normalized traces of 
finite-dimensional irreducible representations of $G$.
See Theorem \ref{Thm_NormalRepFiniteType}\ref{iDEThm_NormalRepFiniteType}
and Corollary~\ref{Cor-IrredNormalRep}.

\vskip.2cm

(3)
When $G$ is a connected LC group, 
then $E(G)$ consists exactly of the normalized traces
of finite-dimensional irreducible unitary representations of $G$; 
indeed, such a group has no representation of type II$_1$.
In particular, if $G$ is a semisimple connected Lie group 
without compact simple constituents,
then $E(G) = \{1_G\}$ by \cite{SevN--50}.
This has already been mentioned in \ref{SectionTypesI+IILC}.

\vskip.2cm

(4)
From the very few results concerning non-locally compact groups we are aware of, 
we mention the description of $E(G)$ in \cite{Dudk--11} 
for $G$ the full group of an ergodic hyperfinite equivalence relation.
\end{exe}

\begin{exe}
% 11.C.6
\label{encoreatrouver}
The Thoma duals of the discrete groups
appearing as the main examples of Chapter \ref{Chapter-ExamplesIndIrrRep}
will be addressed below.

\vskip.2cm

(0)
The Thoma dual of the infinite dihedral group $D_\infty$
consists of the normalized characters of the irreducible representations
$\rho_{0,1}, \rho_{0, -1}, \rho_{\pi, 1}, \rho_{\pi, -1}$ and $\{\rho_\theta\}_{0 < \theta < \pi}$
explicitly written in Section \ref{SectionInfDiGroup}.

\vskip.2cm

\index{Heisenberg group! $2$@$H(\K)$ over a field $\K$}
(1)
Let $H(\K)$ be the Heisenberg group over a field $\K$.
For $E(\Gamma)$ when $\Gamma = H(\K)$ 
is viewed as a discrete group,
see Subsection \ref{ThomaDualTwoStepNil} and
Corollary \ref{Cor-ThomaDualHeisField}.
For $E(H(\K))$ when $\K$ is a non-discrete topological field
and $H(\K)$ is viewed as a topological group homeomorphic with $\K^3$,
see Corollary \ref{ThomaHeisTop}.

\vskip.2cm

(2)
Let $\Aff(\K)$ be the group of affine transformations 
over an infinite field $\K$.
For $E(\Gamma)$ when $\Gamma = \Aff(\K)$ 
is viewed as a discrete group,
see Subsection \ref{ThomaDualAf}
and Theorem \ref{ThomaAffDiscret}.
For $E(\Aff(\K))$ when $\K$ is a non-discrete topological field
and $\Aff(\K)$ is viewed as a topological group 
homeomorphic with $\K^\times \times \K$,
see Corollary \ref{ThomaAffTop}.

\vskip.2cm

(3)
For $\BS(1, p)$, see Subsection \ref{ThomaDualBS}.

\vskip.2cm

(4) For the lamplighter group 
$\Z \rtimes \bigoplus_{k \in \Z} \Z / 2 \Z$
(Section \ref{ThomaDualLamplighter}),

\vskip.2cm

(5)
For $E(\GL_n(\K))$ and $E(\SL_n(\K))$ when $\K$ is an infinite field, 
see Subsection \ref{ThomaDualGL}, 
Theorems \ref{Theo: GLn} and \ref{Theo-SLn}.
\end{exe}

\begin{exe}
% 11.C.7
\label{exThomaDiscrete}
The problem of the description of $E(\Gamma)$ for other groups $\Gamma$
has been considered by several authors, as the following list shows.
\begin{enumerate}[label=(\arabic*)]
\item\label{iDEexThomaDiscrete}
For the infinite symmetric group $\Gamma = \Sym_{\rm fin}(\N)$, 
the group of permutations with finite supports of a countable infinite set,
see \cite{Thom--64c} and, for other methods, \cite{VeKe--81}.
\item\label{iiDEexThomaDiscrete}
For the group $\Gamma = \GL(\infty, \F)$, where $\F$ 
is a finite field, see \cite{Skud--76}, and also \cite{GoKV--14}.
\item\label{iiiDEexThomaDiscrete}
For $\Gamma = \U(\infty)$, the inductive limit of the
unitary groups $\U(n)$ equipped with the inductive limit topology,
see \cite{Voiculescu}, as well as \cite{StVo--75} and \cite{Boyer}.
\item\label{ivDEexThomaDiscrete}
For the group $\Gamma = \SL_n( \Z), n \ge 3$, see \cite{Bekk--07}.
\item\label{vDEexThomaDiscrete}
For $\Gamma$ a lattice in a simple connected Lie group of real rank $\ge 2$, 
see \cite{Pete--15}.
\item\label{viDEexThomaDiscrete}
For $\Gamma$ from the Higman--Thompson families $F_{n, r}$ or $G_{n, r}$,
see \cite{DuMe--14}.
\end{enumerate}
\end{exe}

\section
{Characters and primitive duals}
% Section 11.D
\label{Section-CharactersPrimitive}

Let $G$ be a second-countable LC group.
As Theorem~\ref{ThHalpern} shows, $\QD(G)_{\rm norm}$ is a standard Borel set.
One may thus seek to parametrize $\Pri(G)$ 
by the elements of $\QD(G)_{\rm norm}$
and consider the restriction to $\QD(G)_{\rm norm}$
of the map (\ref{eqq/qd/spG}) of \ref{kapp3extendskapp1}:
\begin{equation}
\label{eqq/ch/spG}
\tag{$\kappa$4}
\kappa^{\rm norm}_{\rm prim} \, \colon \,
\QD(G)_{\rm norm} \rightarrow\Pri(G) ,
\hskip.5cm
\pi \mapsto \textnormal{C*ker}(\pi).
\end{equation}
In contrast to the previous maps
$\kappa^{d}_{\rm prim}$ and $\kappa^{\rm qd}_{\rm prim}$, 
it is not known whether $\kappa^{\rm norm}_{\rm prim}$ is always surjective.
\par 

When $\pi_1$ and $\pi_2$ are \emph{normal} irreducible representations 
such that $\kappa^{\rm norm}_{\rm prim}(\pi_1) = \kappa^{\rm norm}_{\rm prim}(\pi_2)$, 
then $\pi_1$ and $\pi_2$ are equivalent 
(this follows from Proposition~\ref{Pro-IrredNormalRep}
and \cite[4.1.10]{Dixm--C*}).
So, $\kappa^{\rm norm}_{\rm prim}$ is injective on the type I part of 
$\QD(G)_{\rm norm}$.
However, $\kappa^{\rm norm}_{\rm prim}$ is not injective in general 
on the type II part of $\QD(G)_{\rm norm}$. 
\par

\index{Baumslag--Solitar group $\BS(1, p)$}
Indeed, for a prime $p$, consider the Baumslag--Solitar group 
$\Gamma = \Z \ltimes \Z[1/p]$,
where the generator of $\Z$ acts by multiplication by $p$ on $\Z[1/p]$.
As we will see in Corollary~\ref{Cor-Cor-NormalBS}, 
there exist uncountably many normal representations $\pi_i$
of $\Gamma$ of type II$_1$ 
such that $\kappa^{\rm norm}_{\rm prim}(\pi_i) = \kappa^{\rm norm}_{\rm prim}(\pi_j)$ for all $i$ and $j$ 
and such that $\pi_i$ and $\pi_j$ are not quasi-equivalent for $i \ne j$.
\par

The map $\kappa^{\rm norm}_{\rm prim}$ fails to be injective also in the case of
$\Gamma$ the group $\Sym_{\rm fin}(\N)$ of bijections of $\N$ with finite support,
as shown in \cite{Haue--86}.
\par

On the positive side, here is a remarkable result of Puk\'anszky for
the class of connected Lie groups.

\begin{theorem}[\textbf{Puk\'anszky}]
% 11.D.1
\label{ThPukan}
Let $G$ be a connected real Lie group. The map 
$$
\kappa^{\rm norm}_{\rm prim} \, \colon \,
\QD(G)_{\rm norm} \rightarrow\Pri(G) 
$$
is a bijection. 
\end{theorem}

\begin{proof}[References for the proof]
Puk\'anszky's theorem was first proved in \cite[Theorem 1]{Puka--74}.
\par

An exposition of the proof is given in Chapter~III of the monograph \cite{Puka--99}.
Puk\'anszky's theorem can now be deduced from results by Poguntke: 
see Corollary to Theorem 2 in \cite{Pogu--83}.
\end{proof}

Moreover, in the case of $G$ 
a solvable connected and simply connected Lie group, 
\cite{Puka--74} and \cite{Puka--99} provide an explicit description 
of the (usually infinite) characters associated
to the normal factor representations of $G$ 
in terms of a ``generalized Kirillov theory", involving orbits
of the co-adjoint action of $G$ on its Lie algebra.
\par

There are only few non-connected LC groups $G$ for which it is known that 
the map $\kappa^{\rm norm}_{\rm prim}$ is a bijection. 
There is some evidence that this might be true for nilpotent groups;
for what we know:

\begin{prop}
% 11.D.2
\label{Prop-NilGr}
Let $\Gamma$ be a nilpotent-by-finite discrete group. 
\begin{enumerate}[label=(\arabic*)]
\item\label{iDEProp-NilGr}
Every character of $\Gamma$ is finite, that is, 
$\Char(\Gamma) = \Char(\Gamma)_{\rm fin}$;
\item\label{iiDEProp-NilGr}
the map $\kappa^{\rm norm}_{\rm prim}$ is surjective.
\end{enumerate}
\end{prop}

\noindent \emph{Note.}
By Proposition \ref{Char-NormalRep} and Theorem \ref{Thm_NormalRepFiniteType},
the equality in \ref{iDEProp-NilGr} could equally well be written 
$
\QD(\Gamma)_{\rm norm}
=
\QD(\Gamma)_{\rm fin} 
$,
and the sets
$$
\QD(\Gamma)_{\rm norm}, \hskip.2cm
\Char(\Gamma), \hskip.2cm
\text{and} \hskip.2cm E(\Gamma)
$$
are in natural bijection with each other.

\begin{proof}[References for the proof]
The proof was given in \cite{CaMo--84} for nilpotent groups only;
however, the arguments there carry over to nilpotent-by-finite groups.
\par

Indeed, the proof of \ref{iDEProp-NilGr}
is an easy consequence of the fact, 
which was already mentioned in \ref{Sect-PrimIdealHeisenberg}, 
that every primitive ideal of a nilpotent group is maximal, 
and this is still true for a nilpotent-by-finite group (see \cite{Pogu--82}). 
Claim \ref{iiDEProp-NilGr} is a consequence of the amenability of $\Gamma$ and, again,
of the fact that primitive ideals of $C^*_{\rm max}(\Gamma)$ are maximal.
\end{proof}

\begin{theorem}
% 11.D.3
\label{Th-HoweKa}
The map $\kappa^{\rm norm}_{\rm prim} \,\colon E(\Gamma) \rightarrow \Pri(\Gamma)$
is a bijection in the following cases:
\begin{enumerate}[label=(\arabic*)]
\item\label{iDETh-HoweKa}
$\Gamma$ is a finitely generated group with polynomial growth;
\item\label{iiDETh-HoweKa}
$\Gamma = \G(\Q)$ is the group of rational point 
of a unipotent algebraic group over $\Q$.
\item\label{iiiDETh-HoweKa}
$\Gamma = \GL_n(\K)$, 
where $\K$ is an infinite algebraic extension of a finite field.
\end{enumerate}
\end{theorem}

\begin{proof}[References for the proof]
Result \ref{iDETh-HoweKa} was established in \cite{Howe--77a} 
for finitely generated torsion free nilpotent groups,
and extended in \cite{Kani--80}
to finite extensions of finitely generated nilpotent groups, 
which are the groups with polynomial growth by Gromov's theorem.
\par

Part \ref{iiDETh-HoweKa} is from \cite{Pfef--95}.
Concerning \ref{iiiDETh-HoweKa}, see Corollary~\ref{Cor-PrimGLn} below.
\end{proof}

\section
{Summing up diagram: dual and variants}
% Section 11.E
\label{Section-SummingUpDiagram}

Let $G$ be a second-countable locally compact group.
The duals and sets of characters introduced in this book
are organized as in the following diagram.
$$
\begin{matrix}
% row 1
&&&&&&&& \Tr_1(G)
\\
% row 2
&&&&&&&& \cup
\\
% row 3
\widehat G_{\rm fd}
&&
\hookrightarrow
&&
\QD(G)_{\rm fin}
&
\simeq
&
\Char(G)_{\rm fin}
&
\simeq
&
E(G)
\\
%% more space -- but not convenient !
%% &&&&&&&&
%% \\
% row 4
&&&&
\cap
&&
\cap
&&
\\
% row 5
\downarrow
&&&& 
\QD(G)_{\rm norm}
&
\simeq
&
\Char(G)
&&
\\
% row 6
&&&&
\cap
&&&&.
\\
% row 7
\widehat G
&&
\overset{\kappa^{\rm d}_{\rm qd}}{\hookrightarrow}
&&
\QD(G)
&&&&
\\
% row 8
&&&&&&&&
\\
% row 9
& \searrow\kappa^{\rm d}_{\rm prim} && \swarrow\kappa^{\rm qd}_{\rm prim} &&&&&
\\
% row 10
&&&&&&&&
\\
% row 11
&& \Pri(G) &&&&&&
\end{matrix}
$$

%-----------------------------------------------------------------------
% End of chapter 11
%-----------------------------------------------------------------------
\chapter{Examples of Thoma's duals}
% Chapter 12
\label{ChapterThomadualExamples}

\emph{After a short preparation section
on traces (Section \ref{Section-FactTraces}),
we determine Thoma's duals $E(\Gamma)$ 
for our four main test classes of discrete groups $\Gamma$:
\begin{enumerate}[noitemsep]
\item[$\bullet$]
two-step nilpotent discrete groups,
in particular the Heisenberg group over an infinite field $\K$
or over the integers
(Section \ref{ThomaDualTwoStepNil}), 
\item[$\bullet$]
the affine group $\Aff(\K)$ over $\K$
(Section \ref{ThomaDualAf}), 
\item[$\bullet$]
the Baumslag--Solitar group $\BS(1, p)$ for a prime $p$
(Section \ref{ThomaDualBS}),
\item[$\bullet$]
the lamplighter group $\Z \rtimes \bigoplus_{k \in \Z} \Z / 2 \Z$
(Section \ref{ThomaDualLamplighter}),
\item[$\bullet$]
and the general linear group $\GL_n(\K)$ over $\K$
(Section \ref{ThomaDualGL}).
\end{enumerate}
The results are well-known:
see \cite[Example 2.4]{ArKa--97} for the Heisenberg group,
\cite[\S~1, Beispiel 2]{Thom--64b} for the affine group, 
\cite[\S~3, Prop. 1]{Guic--63} for $\BS(1, 2)$,
and \cite{Kiri--65} for the general linear group.
However, apart from the case of Heisenberg group, 
the proofs we give are modelled after the approach 
used in \cite{Bekk--07} for the case $\SL_n(\Z), n \ge 3$.
}

\section
{Elementary lemmas about traces on groups}
% Section 12.A
\label{Section-FactTraces}

We first establish two elementary lemmas to be used in the sequel.
\par

Let $G$ be a topological group.
Let $\varphi \in P_1(G)$ be
a normalized continuous function of positive type on $G$.
(Functions of positive type have been defined in \ref{S-FPosType}.)
Set
$$
P_\varphi \, = \, \{ g \in G \mid \vert \varphi (g) \vert = 1 \} 
\hskip.5cm \text{and} \hskip.5cm
K_\varphi \, = \, \{ g \in G \mid \varphi (g) = 1 \} ,
$$
which are closed subsets of $G$ containing $e$,
and $P_\varphi \supset K_\varphi$.
Let $(\Hi_\varphi, \xi_\varphi, \pi_\varphi)$ be the GNS triple defined by $\varphi$,
as in Construction \ref{constructionGNS2}.
Let $\chi_\varphi = \chi_{\pi_\varphi}$ be the associated unitary character of 
the projective kernel $\Pker \pi_\varphi$ of $\pi_\varphi$ 
as in Section \ref{Section-CentralCharacter}.
\par

Recall from Proposition~\ref{Prop-ProjKer} 
that every factor representation has a central character.

\begin{lem}
% 12.A.1
\label{Lem-CharHeis}
Let $G$ be a topological group, $\varphi \in P_1(G)$, 
and $\pi_\varphi, P_\varphi, K_\varphi$ as above. 
\begin{enumerate}[label=(\arabic*)]
\item\label{iDELem-CharHeis}
We have
$$
P_\varphi \, = \,
\left\{ g \in G \mid \pi_\varphi (g)\xi_\varphi = \varphi(g) \xi_\varphi \right\}
$$
and
$$
K_\varphi \, = \,
\left\{ g \in G \mid \pi_\varphi (g)\xi_\varphi = \xi_\varphi \right\} .
$$
In particular, $P_\varphi$ and $K_\varphi$ are closed subgroups of $G$,
and the restriction $\varphi \vert_{P_\varphi}$
is a unitary character $P_\varphi \to \T$.
\item\label{iiDELem-CharHeis}
We have 
$\varphi(g_1 g g_2) = \varphi(g_1) \varphi(g) \varphi(g_2)$
for all $g_1, g_2 \in P_\varphi$ and $g \in G$.
\end{enumerate}
Assume moreover that $\varphi$ is central, i.e., that $\varphi \in \Tr_1(G)$. 
\begin{enumerate}[label=(\arabic*)]
\addtocounter{enumi}{2}
\item\label{iiiDELem-CharHeis}
Then $P_\varphi = \Pker \pi_\varphi$ and $K_\varphi = \ker \pi_\varphi$; 
moreover, $\varphi$ coincides with the unitary character $\chi_{\varphi}$
of $P_\varphi$ associated to $\pi_\varphi$.
\end{enumerate}
Assume moreover that $\varphi$ is indecomposable, i.e., that $\varphi \in E(G)$.
\begin{enumerate}[label=(\arabic*)]
\addtocounter{enumi}{3}
\item\label{ivDELem-CharHeis} 
The restriction $\varphi \vert_Z$ of $\varphi$ to the centre $Z$ 
coincides with the central character of $\pi_\varphi$.
\end{enumerate}
\end{lem}

\begin{proof}
\ref{iDELem-CharHeis}
For every $g \in G$, we have
$\vert \varphi (g) \vert
= \vert \langle \pi_\varphi(g) \xi_\varphi \mid \xi_\varphi \rangle \vert
\le 1$
by the Cauchy--Schwarz inequality.
By the equality case, if $g \in P_\varphi$,
then $\pi_\varphi (g) \xi_\varphi$ is a scalar multiple of $\xi_\varphi$,
and this scalar multiple is precisely $\varphi (g)$.
Consequently,
$$
P_\varphi \, = \, 
\left\{ g \in G \mid \ \pi_\varphi (g)\xi_\varphi = \varphi(g) \xi_\varphi \right\} ,
$$
and $P_\varphi$ is a subgroup of $G$, 
and the map $P_\varphi \to \T$, $g \mapsto \varphi(g)$ is a unitary character,
with kernel $K_\varphi$.

\vskip.2cm

\ref{iiDELem-CharHeis}
For $g_1, g_2 \in P_\varphi$ and $g \in G$, 
it follows from \ref{iDELem-CharHeis} that
$$
\begin{aligned}
\varphi(g_1 g g_2) 
\, &= \, \langle \pi_\varphi(g_1 g g_2) \xi_\varphi \mid \xi_\varphi \rangle
\, = \, \langle \pi_\varphi(g) \pi_\varphi(g_2) \xi_\varphi 
 \mid \pi_\varphi(g_1^{-1})\xi_\varphi \rangle
\\
\, &= \,
\langle \pi_\varphi(g) \varphi(g_2) \xi_\varphi
\mid \varphi(g_1^{-1})\xi_\varphi \rangle
\, = \, \varphi(g_1) \varphi(g) \varphi(g_2) .
\end{aligned}
$$

\vskip.2cm
 
\ref{iiiDELem-CharHeis} 
Let $g \in G$; the triple $(\Hi_\varphi, \pi(g)\xi_\varphi, \pi_\varphi)$
is a GNS triple for $\varphi^g$, where $\varphi$ is the function of positive type
on $G$ defined by $\varphi^g(x) = \varphi(g^{-1} xg)$.
\par

Assume now that $\varphi$ is central and let $x \in P_\varphi$.
Since $P_\varphi = P_{\varphi^g}$, it follows from \ref{iDELem-CharHeis} that 
$$
\pi_\varphi (x)\pi_\varphi (g)\xi_\varphi 
\, = \,
\varphi^g(x)\pi_\varphi (g)\xi_\varphi
\, = \, 
\varphi(x)\pi_\varphi (g)\xi_\varphi
\hskip.5cm \text{for all} \hskip.2cm
g \in G
$$
 and hence that $\pi_\varphi (x) = \varphi(x)I$, since $\xi_\varphi$
is a cyclic vector for $\pi_\varphi$.
This shows that $P_\varphi$ is contained in $\Pker \pi$. 
Since $\Pker \pi$ is obviously contained in $P_\varphi$, it follows that 
$P_\varphi = \Pker \pi$; it is obvious that $\varphi \vert_{P_\varphi}$ coincides
with $\chi_{\varphi}$ and that 
$K_\varphi = \ker_{\pi_\varphi}$.

\vskip.2cm 

\ref{ivDELem-CharHeis}
Since $\varphi \in E(\Gamma)$, the representation $\pi_\varphi$ is factorial and
the claim follows from Proposition~\ref{Prop-ProjKer}.
\end{proof}

The next elementary lemma
will be the crucial tool in our study of traces on two-step nilpotent groups.
Recall that $[g, h]$ denotes the commutator $g^{-1}h^{-1}gh$ 
of two elements $g, h \in G$.

\begin{lem}
% 12.A.2
\label{Lemmaphi(g)=0}
Let $G$ a topological group, $\varphi \in \Tr(G)$, and $g \in G$. 
Assume that there exists
$h \in G$ such that $[g, h] \in P_\varphi \smallsetminus K_\varphi$. 
Then $\varphi(g) = 0$.
\end{lem}

\begin{proof}
Using the fact $\varphi$ is central and Lemma~\ref{Lem-CharHeis}, we have
$$
\varphi(g) \, = \, \varphi(h^{-1}gh) 
\, = \, \varphi( g [g, h]) 
\, = \, \varphi( g ) \varphi( [g, h] ) ,
$$
and the claim follows.
\end{proof}

\section
{Two-step nilpotent groups}
% Section 12.B
\label{ThomaDualTwoStepNil}

Let $\Gamma$ be a two-step nilpotent discrete group, and $Z$ its centre,
as in Sections \ref{Section-IrrRepTwoStepNil} and \ref{Sect-PrimIdealHeisenberg}.
Recall that, for $\psi \in \widehat Z$, we denote by $Z_\psi$
the inverse image of the centre of $\Gamma/\ker \psi$
under the canonical epimorphism $\Gamma \twoheadrightarrow \Gamma/\ker \psi$.
\par

Recall also that, if $\varphi$ 
is a function of positive type on a subgroup $H$ of $\Gamma$, 
then its trivial extension $\widetilde \varphi$ to $\Gamma$ 
is a function of positive type on $\Gamma$ (Proposition~\ref{diagcoeffinduced}).
\index{Two-step nilpotent group}

\begin{theorem}
% 12.B.1
\label{Theo-ThomaDualTwoStepNilpotent}
Let $\Gamma$ be a two-step nilpotent discrete group, with centre $Z$. Set
$$
X \, := \,
\displaystyle
\bigsqcup_{\psi \in \widehat Z} \left(
\{\psi\} \times \left\{\chi \in \widehat{Z_\psi} \mid \chi \vert_Z = \psi \right\}
\right)
$$
(this $X$ appears in Theorems \ref{Theo-PrimIdealTwoStepNilpotent}
and \ref{Theo-ThomaDualTwoStepNilpotent-bis}).
\par

Then the map 
$$
\Phi \, \colon \, X \, \to \, E(\Gamma),
\hskip.5cm
(\psi, \chi) \mapsto \widetilde \chi
$$
is a bijection.
The inverse map $\Phi^{-1} \,\colon E(\Gamma) \to X$ is given by 
$$
\Phi^{-1}(\varphi) \, = \, (\psi_{\pi_\varphi}, \lambda_{\pi_\varphi}),
$$
where $\psi_{\pi_\varphi}$ is the central character 
of the GNS representation $\pi_\varphi$ associated to $\varphi$ and
$\chi_{\pi_\varphi}$ the unitary character of $\Pker {\pi_\varphi}$
determined by $\pi_\varphi$.
\end{theorem}

\begin{proof}
$\bullet$ \emph{First step.}
We claim that every element from $E(\Gamma)$ 
is of the form $\Phi(\psi, \chi)$ for a pair $(\psi, \chi)$ in the set $X$.
\par

Let $\varphi \in E(\Gamma)$;
let $\pi = \pi_\varphi$ be the associated GNS representation.
The quotient group $\overline \Gamma = \Gamma/\ker \pi$
is two-step nilpotent, and its centre is $\Pker \pi/\ker \pi$
(see Proposition~\ref{Prop-ProjKer}.ii).
\par

A function $\varphi \in E(\Gamma)$ is a unitary character, say $\chi$,
if and only if $\pi_\varphi(g) = \chi(g) \mathrm{Id}_{\Hi_\varphi}$
for all $g \in G$,
if and only if $\Pker \pi_\varphi = \Gamma$.
If these conditions hold, then $\widetilde \chi = \chi$
is in the image of $\Phi$.
We can therefore assume from now on that $\Pker \varphi \ne \Gamma$.
\par

Let $\gamma \in \Gamma$ with $\gamma \notin \Pker \pi$.
Then there exists $\delta \in \Gamma$ with $[\gamma, \delta] \notin \ker \pi$. 
Since $\overline \Gamma$ is two-step nilpotent,
the image of $[\gamma, \delta]$ in $\overline \Gamma$
belongs to the centre of $\overline \Gamma$.
Therefore, $[\gamma, \delta] \in \Pker \pi$. 
Therefore $\varphi(\gamma) = 0$ by Lemma~\ref{Lemmaphi(g)=0}.
This shows that $\varphi$ is the trivial extension to $\Gamma$
of the unitary character $\chi_\pi$ of $\Pker \pi$
determined by $\pi$.
\par

Let $\psi \in \widehat Z$ be the central character of the factor representation $\pi$.
By Proposition~\ref{Prop-ProjKer}.ii,
we have $\Pker \pi = Z_\psi$, where $Z_\psi$ 
is such that $Z_\psi/\ker \psi$ is the centre of $\Gamma/\ker \psi$.
\par

Since $\chi_\pi \vert_Z$ coincides with $\psi$, we have 
$(\psi, \chi) \in X$ and 
$\varphi = \Phi(\psi, \chi)$ for $\chi = \chi_\pi$.

\vskip.2cm

$\bullet$ \emph{Second step.} 
Let $\psi \in \widehat Z$
and $\varphi \in \Tr_1(\Gamma)$ be such that $\varphi \vert_Z = \psi$. 
We claim that $\varphi = 0$ on $\Gamma \smallsetminus Z_\psi$. 
\par

Indeed, since the convex hull of $E(\Gamma)$ is dense $\Tr_1(\Gamma)$
for the topology of pointwise convergence, we can assume that
$\varphi = \sum_{i = 1}^n t_i \varphi_i$ is a convex combination
of functions $\varphi_i \in E(\Gamma)$. 
\par

On the one hand, we have for every $z \in Z$ 
$$
\psi(z) \, = \, \varphi(z) \, = \, \sum_{i = 1}^n t_i\varphi_i(z);
$$
since $\vert \psi (z) \vert = 1$ 
and $\vert \varphi_i (z) \vert \le 1$,
it follows that $\varphi_i\vert_Z = \psi$ for every $i = 1, \hdots, n$.
Therefore we have, by the first step, 
$\varphi_i = 0$ on $\Gamma \smallsetminus Z_\psi$ for every $i = 1, \hdots, n$. 
Therefore $\varphi = 0$ on $\Gamma \smallsetminus Z_\psi$.

\vskip.2cm
 
$\bullet$ \emph{Third step.} 
We claim that every function in the image of $\Phi$ 
belongs to $E(\Gamma)$. 
\par
 
Let $(\psi, \chi) \in X$; 
we have to show that the trivial extension $\varphi := \widetilde \chi$ 
of $\chi$ belongs to $E(\Gamma)$. 
Since $\varphi$ is in $\Tr_1(\Gamma)$ 
by Proposition \ref{diagcoeffinduced},
we have only to check $\varphi$ that it is indecomposable.
\par

Let $\varphi_1, \varphi_2 \in \Tr_1(\Gamma)$ and $t \in \mathopen] 0,1 \mathclose[$
be such that $\widetilde \chi = t \varphi_1+ (1-t) \varphi_2$. 
We have for every $\gamma \in Z_\psi$
$$
\chi(\gamma) \, = \, \varphi(\gamma) \, = \, t \varphi_1(\gamma) + (1-t) \varphi_2(\gamma) ;
$$
since $\vert \chi (\gamma) \vert = 1$ 
and $\vert \varphi_1 (\gamma) \vert \le 1, \vert\varphi_2(\gamma) \vert \le 1$,
it follows that $\varphi_1 \vert_{ Z_\psi} = \varphi_2 \vert_ {Z_\psi} = \chi$.
In particular, $\varphi_1 \vert_{ Z} = \varphi_2 \vert_ Z = \chi \vert_Z = \psi$.
\par

By the second step, we have therefore
$\varphi_1 = 0$ and $\varphi_2 = 0$ on $\Gamma \smallsetminus Z_\psi$.
Therefore $\varphi_1 = \varphi_2= \widetilde \chi$.

\vskip.2cm

$\bullet$ \emph{Fourth step.} 
The map $\Phi$ is injective. 
Indeed, let $(\psi_1, \chi_1), (\psi_2, \chi_2) \in X$ 
be such that $\Phi(\psi_1, \chi_1) = \Phi(\psi_2, \chi_2)$,
that is, $\widetilde \chi_1 = \widetilde \chi_2$. 
Then 
$$
\psi_1 \, = \, \chi_1 \vert_Z \, = \, \widetilde \chi_1 \vert_Z
\, = \, \widetilde \chi_2 \vert_Z \, = \, \chi_2 \vert_Z
 \, = \, \psi_2
$$
and 
$$
\chi_1 \, = \, \widetilde \chi_1 \vert_{Z_{\psi_1}}
 \, = \, \widetilde \chi_2 \vert_{Z_{\psi_2}} \, = \, \chi_2.
$$
Therefore, $(\psi_1, \chi_1) = (\psi_2, \chi_2)$.
\end{proof}

We rephrase Theorem~\ref{Theo-ThomaDualTwoStepNilpotent}
in terms of the space $\QD(\Gamma)_{\rm fin}$
of quasi-equivalence classes of factor representations of finite type,
specifying moreover their possible types.

\begin{theorem}
% 12.B.2
\label{Theo-ThomaDualTwoStepNilpotent-bis}
Let $\Gamma$ be a two-step nilpotent discrete group, with centre $Z$,
and let $X$ be as in Theorem \ref{Theo-ThomaDualTwoStepNilpotent}.
\par

Then the map 
$$
\Phi \, \colon \, X \, \to \, \QD(\Gamma)_{\rm fin},
\hskip.5cm
(\psi, \chi) \, \mapsto \, \Ind_{Z_\psi}^\Gamma \chi
$$
is a bijection.
\par

Moreover, for $(\psi, \chi)\in X$,
the factor representation $\Ind_{Z_\psi}^\Gamma \chi$
is of type II$_1$ if $Z_\psi$ has infinite index in $\Gamma$,
and of type I$_{\rm finite}$ otherwise.
\end{theorem}

\begin{proof}
The fact that $\Phi$ is a bijection between $X$ and $\QD(G)_{\rm fin}$
follows from Theorem~\ref{Theo-ThomaDualTwoStepNilpotent},
in combination with Theorem~\ref{Thm_NormalRepFiniteType}
and Proposition~\ref{diagcoeffinduced}.
\par

Let $(\psi, \chi)\in X$.
If $Z_\psi$ has finite index in $\Gamma$, then $\Ind_{Z_\psi}^\Gamma \chi$
is finite dimensional and hence $\Ind_{Z_\psi}^\Gamma \chi$
is of type I$_{n}$ for an integer $n$.
If $Z_\psi$ has infinite index, then $\Ind_{Z_\psi}^\Gamma \chi$ is of type II$_{1}$;
indeed, otherwise, $\Ind_{Z_\psi}^\Gamma \chi$
would be a multiple of a finite dimensional irreducible representation
and this is impossible by Proposition~\ref{PropOnWeAndEq}~\ref{iiDEPropOnWeAndEq}.
\end{proof}

The following corollary is an immediate consequence
of Theorem~\ref{Theo-ThomaDualTwoStepNilpotent}
and the description of $\Pri(\Gamma)$
from Theorem~\ref{Theo-PrimIdealTwoStepNilpotent}.
Compare also Theorem~\ref{Th-HoweKa}.

\begin{cor}
% 12.B.3
\label{ThomaHeis-PrimIdeal}
For a two-step nilpotent discrete group $\Gamma$, the natural map
$$
\kappa^{\rm norm}_{\rm prim} \, \colon \, E(\Gamma) \rightarrow \Pri(\Gamma)
$$
is a bijection
\end{cor}

We are going to apply Theorem \ref{Theo-ThomaDualTwoStepNilpotent}
to Heisenberg groups.

\subsection*{Thoma's dual for Heisenberg groups}

As in Sections \ref{Section-IrrRepTwoStepNil} and \ref{Sect-PrimIdealHeisenberg},
$H(R)$ denotes here the Heisenberg group
over a unital commutative ring $R$.
\par

Fix a unitary character $\psi$ of the centre $Z$ of $H(R)$.
Let $I_\psi$ be the largest ideal of $R$ contained in $\ker \psi$ 
(which already appeared in Sections
\ref{Section-IrrRepTwoStepNil}, \ref{Section-FiniteDimRepForSomeGroups},
and \ref{Sect-PrimIdealHeisenberg}).
By Lemma~\ref{Lem-IrredHeisRing}, we have 
$$
Z_\psi \, = \, \{(a, b, c) \in \Gamma \mid (a, b) \in I_\psi^2,c \in R\} .
$$
Every unitary character of $Z_\psi$ with $\chi\vert_Z = \psi$ is of the form
$\chi_{\psi, \alpha, \beta}$ for a unique pair $(\alpha, \beta) \in \widehat{I_\psi}^2$,
where $\chi_{\psi, \alpha, \beta} \,\colon Z_\psi \to \T$ is defined by the formula
$$
\chi_{\psi, \alpha, \beta}(a, b, c) \, = \, \alpha(a)\beta(b)\psi(c)
\hskip.5cm \text{for} \hskip.2cm
(a, b) \in I_\psi^2, \hskip.1cm b, c \in R.
$$
In view of these facts,
the following result is a direct consequence of
Theorem~\ref{Theo-ThomaDualTwoStepNilpotent}.

\begin{cor}
% 12.B.4
\label{Cor-ThomaDualHeisRing}
Let $\Gamma=H(R)$ be the Heisenberg group
over a unital commutative ring $R$. Set 
$$
X \, := \, \bigsqcup_{\psi \in \widehat Z} \{\psi\} \times \widehat{I_\psi}^2.
$$
The map $\Phi \,\colon X \to E(\Gamma)$, defined by 
$$
\Phi(\psi, (\alpha, \beta)) \, = \, \widetilde{\chi_{\psi, \alpha, \beta}}
$$
(where $\widetilde{\chi_{\psi, \alpha, \beta}}$
is the trivial extension of $\chi_{\psi, \alpha, \beta}$ to $\Gamma$)
is a bijection.
\end{cor}

We first consider the particular case where $R = \K$ is a field.
For $\psi \in \widehat{\K} \smallsetminus \{1_\K\}$,
we have $I_\psi = \{0\}$ and hence $Z_\psi = Z$;
for $\psi = 1_\K$, we have $I_\psi = \K$ and hence $Z_\psi = H(\K)$.
The following result is therefore an immediate consequence
of Corollary~\ref{Cor-ThomaDualHeisRing}.

\begin{cor}
% 12.B.5
\label{Cor-ThomaDualHeisField}
\index{Heisenberg group! $2$@$H(\K)$ over an infinite field $\K$}
Let $\Gamma = H(\K)$ for a field $\K$ and $Z$ the centre of $\Gamma$. 
Let $p \,\colon \Gamma \twoheadrightarrow \Gamma/Z$ be the canonical projection.
Set
$$
X \, = \, \left( \widehat Z \smallsetminus \{1_Z\} \right) \sqcup \widehat{\Gamma/Z}
\, \approx \, \left( \widehat \K \smallsetminus \{1_{\K}\} \right) \sqcup \widehat \K^2 
$$
(this $X$ appears in Corollary \ref{Cor-PrimIdealHeisField}).
\par

The map $\Phi \,\colon X \to E(\Gamma)$, defined by 
$$
\Phi(\psi) \, = \, \widetilde \psi
\hskip.5cm \text{for} \hskip.2cm
\psi \in \widehat Z \smallsetminus \{1_Z\}
$$
and
$$
\Phi(\chi) \, = \, \chi \circ p 
\hskip.5cm \text{for} \hskip.2cm
\chi \in \widehat{\Gamma/Z},
$$
is a bijection.
\end{cor}

We rephrase Corollary~\ref{Cor-ThomaDualHeisField}
in terms of factor representations,
specifying moreover their possible types.

\begin{cor}
% 12.B.6
\label{Cor-ThomaDualHeisField-Bis}
Let $\Gamma = H(\K)$ for a field $\K$
and $p \,\colon \Gamma \twoheadrightarrow \Gamma/Z$ the canonical projection;
let $X$ be as in Corollary~\ref{Cor-ThomaDualHeisField}.
The map $\Phi$ defined by 
$$
\Phi(\psi) \, = \, \Ind_{Z}^\Gamma \psi
\hskip.5cm \text{for} \hskip.2cm
\psi \in \widehat Z \smallsetminus \{1_Z\}
$$
and 
$$
\Phi(\chi) \, = \, \chi \circ p 
\hskip.5cm \text{for} \hskip.2cm
\chi \in \widehat{\Gamma/Z}
$$
is a bijection between $X$ and the set $\QD(\Gamma)_{\rm fin}$.
Moreover, for $\psi \in \widehat Z \smallsetminus \{1_Z\}$,
the factor representation $\Ind_{Z}^\Gamma \psi$
is of type II$_1$ if $\K$ is infinite
and of type I$_q$ if $\K$ is finite of cardinality $q$.
\end{cor}

\begin{proof}
The fact that $\Phi$ is a bijection between $X$ and $\QD(G)_{\rm fin}$
is a straightforward consequence
of Theorem~\ref{Theo-ThomaDualTwoStepNilpotent-bis}
and Corollary~\ref{Cor-ThomaDualHeisField}.
That $\Ind_{Z}^\Gamma \psi$ is of type II$_1$ if $\K$ is infinite follows
also from Theorem~\ref{Theo-ThomaDualTwoStepNilpotent-bis}. 
\par

Assume that $\K$ is finite of cardinality $q$
and let $\psi \in \widehat Z \smallsetminus \{1_Z\}$. 
Since $\Ind_{Z}^\Gamma \psi$ is a factor representation,
it is a multiple of an irreducible representation $\pi$ of $\Gamma$. 
Since $\pi$ has $\psi$ as central character, $\pi$ is not one-dimensional
and so has dimension $q$, by Corollary~\ref{Cor-FinDimRepHeis-finiteFields}.
It follows that $\Ind_{Z}^\Gamma \psi$ is of type I$_q$. 
\end{proof}

The following corollary shows that Heisenberg groups
over non discrete topological fields have no interesting finite factor representations.

\begin{cor}
% 12.B.7
\label{ThomaHeisTop}
Let $\K$ be a non-discrete topological field.
Consider now $H(\K)$ as a topological group,
with its topology making it homeomorphic to the topological space $\K^3$,
and therefore $H(\K)/Z$ topologically isomorphic to $\K^2$.
\par

Every function in $E(H(\K))$ is the lift of a unitary character of
the quotient group $H(\K)/Z$.
\par

Equivalently, the factor representations of finite type of $H(\K)$
are the multiples of one-dimensional representations. 
\end{cor}

\begin{proof}
Denote here by $H(\K)_{\rm disc}$
the group $\Gamma$ of \ref{Cor-ThomaDualHeisField},
where ``disc'' stands for ``discrete''.
A function $\varphi \,\colon H(\K) \to \C$ that is in $E(H(\K)_{\rm disc})$
is in $E(H(\K))$ if and only if it is continuous.
In particular, a function of the type $\widetilde \psi$ for $\psi \in \widehat Z$
is not in $E(H(\K))$, because $Z$ is not open in $H(\K)$.
\end{proof}

Next, we consider the case where $R = \Z$ is the ring of rational integers.
\par

We recall some notation introduced before Corollary~\ref{Cor-PrimIdealHeisIntegers}.
We denote by $\widehat Z_\infty$ the set of elements of infinite order 
and $\widehat Z_n$ the elements of order $n \ge 1$ in $\widehat Z$.
\par

Let $\psi \in \widehat Z$. If $\psi \in \widehat Z_\infty$,
then $Z_\psi= Z$.
If $\psi \in \widehat Z_n$ for some $n \ge 1$,
then $Z_\psi = \Lambda(n)$, where
$$
\Lambda(n) \, = \, \{(a, b, c) \in \Gamma \mid a, b \in n\Z, \hskip.1cm c \in \Z\}.
$$
Recall that 
$$
\widehat Z \approx \widehat \Z \, = \,
\{\psi_\theta \mid \theta \in \mathopen[ 0,1 \mathclose[ \},
$$
with $\psi_\theta(n) = e^{2 \pi i n \theta}$ for $n \in \Z$.
\par

Fix $\psi = \psi_\theta \in \widehat{Z}_n$.
For $(\alpha, \beta) \in \mathopen[ 0,1/n \mathclose[^2$,
we denote by $\chi_{\psi, \alpha, \beta}$ the unitary character 
of $\Lambda(n)$ which coincides with $\psi$ on $Z$, 
defined by
$$
\chi_{\psi, \alpha, \beta}(a, b, c) \, = \, 
e^{2 \pi i (\alpha a + \beta b)} \psi (c) 
\hskip.5cm \text{for} \hskip.2cm
a, b \in n\Z, \hskip.1cm c \in \Z.
$$
\par

The following result is an immediate consequence
of Corollary~\ref{Cor-ThomaDualHeisRing}.

\begin{cor}
% 12.B.8
\label{Cor-ThomaDualHeisIntegers}
\index{Heisenberg group! $3$@$H(\Z)$}
Let $\Gamma = H(\Z)$ be the Heisenberg group over the integers
and $Z$ its centre.
Set
$$
X \, = \, \widehat Z_\infty \bigsqcup 
\bigg(
\bigsqcup_{n \ge 1} \bigsqcup_{\psi \in \widehat Z_n}
\{\psi\} \times \mathopen[ 0,1/n \mathclose[^2
\bigg)
$$
(this $X$ appears in Corollary \ref{Cor-PrimIdealHeisIntegers}).
\par

The map $\Phi \,\colon X \to E(\Gamma)$, defined by 
$$
\Phi(\psi) \, = \, \widetilde \psi 
\hskip.5cm \text{for} \hskip.2cm
\psi \in \widehat Z_\infty
$$
and 
$$
\Phi(\psi, \alpha, \beta) \, = \,
\widetilde{\chi_{\psi, \alpha, \beta}}
\hskip.5cm \text{for} \hskip.2cm
(\psi, \alpha, \beta) \in \widehat Z_n \times \mathopen[ 0,1/n \mathclose[^2,
$$
is a bijection.
\end{cor}

We rephrase Corollary~\ref{Cor-ThomaDualHeisIntegers}
in terms of factors representations,
specifying moreover their possible types.

\begin{cor}
% 12.B.9
\label{Cor-ThomaDualHeisIntegers-Bis}
Let $\Gamma = H(\Z)$ 
and $X$ as in Corollary~\ref{Cor-ThomaDualHeisIntegers}.
The map $\Phi$ defined by 
$$
\Phi(\psi) \, = \, \Ind_{Z}^\Gamma \psi
\hskip.5cm \text{for} \hskip.2cm
\psi \in \widehat Z_\infty
$$
and 
$$
\Phi(\psi, \alpha, \beta) \, = \, 
\Ind_{\Lambda(n)}^\Gamma \chi_{\psi, \alpha, \beta}
\hskip.5cm \text{for} \hskip.2cm
(\psi, \alpha, \beta) \in \widehat Z_n \times \mathopen[ 0,1/n \mathclose[^2,
$$
is a bijection between $X$ and the set $\QD(\Gamma)_{\rm fin}$.
\par

Moreover, $\Ind_{Z}^\Gamma \psi$
is of type II$_1$ for $\psi \in \widehat Z_\infty$
and $\Ind_{\Lambda(n)}^\Gamma \chi_{\psi, \alpha, \beta}$
is of type I$_n$ for
$(\psi, \alpha, \beta) \in \widehat Z_n \times \mathopen[ 0,1/n \mathclose[^2$.
\end{cor}

\begin{proof}
The fact that $\Phi$ is a bijection between $X$ and $\QD(\Gamma)_{\rm fin}$
is a consequence of Theorem~\ref{Theo-ThomaDualTwoStepNilpotent-bis}
as well as the fact that $\Ind_{Z}^\Gamma \psi$ is of type II$_1$
for $\psi \in \widehat Z_\infty$.
\par

Let $(\psi, \alpha, \beta) \in \widehat Z_n \times \mathopen[ 0,1/n \mathclose[^2$.
The factor representation
$\Ind_{\Lambda_n}^\Gamma \chi_{\psi, \alpha, \beta}$
is a multiple of the irreducible representation
$\Ind_{\Gamma(n)}^\Gamma \chi_{\psi, \alpha, \beta}$,
of dimension $n$ from Corollary~\ref{Cor-IrredFiniHeisIntegers}.
Therefore $\Ind_{\Lambda_n}^\Gamma \chi_{\psi, \alpha, \beta}$
is of type I$_n$.
\end{proof}

\section
{Affine groups of infinite fields}
% Section 12.C
\label{ThomaDualAf}

As in Sections \ref{Section-IrrRepAff} and \ref{Sect-PrimIdealAffGr},
let $\K$ be an \emph{infinite} field,
$\Gamma = \Aff(\K)$ the group of affine transformations of $\K$,
and $N$ its derived group.
\index{Affine group! $4$@$\Aff(\K)$ of an infinite field $\K$}
\par

We will describe below (Theorem~\ref{ThomaAffDiscret})
the Thoma dual of $\Gamma$;
the main ingredients for the proof are the following two lemmas that deal with restrictions of characters to normal subgroups.

\begin{lem}
% 12.C.1
\label{Lem-ExtNorm}
Let $\Gamma$ be a group and $N$ a normal subgroup of $\Gamma$.
Let $\varphi \in E(\Gamma)$ be such that $\varphi \vert_N = \delta_e$.
\par

Then $\varphi(\gamma) = 0$ for every $\gamma \in \Gamma$ 
for which the set of commutators
$\{[\gamma, x] \mid x \in N\}$ is infinite. 
\end{lem}

\begin{proof}
By assumption, there exists an infinite sequence
of $(x_k)_k$ of elements in $N$ such that,
for $y_k := [\gamma,x_k]$, we have
$$
y_k \ne y_l 
\hskip.5cm \text{for all} \hskip.2cm 
k \ne l.
$$
Let $(\pi, \Hi, \xi)$ be a GNS triple associated to $\varphi$.
Since $y_l^{-1}y_k\in N \smallsetminus \{e\}$, we have $\varphi(y_l^{-1}y_k) = 0$,
that is, $\langle \pi(y_k)\xi \mid \pi(y_l)\xi \rangle = 0$ for all $k \ne l$.
This means that $(\pi(y_k) \xi)_k$ is an orthonormal sequence in $\Hi$,
so that $(\pi(y_k) \xi)_k$ is weakly convergent to $0$ in $\Hi$.
\par

For every $k$, we have
$$
\varphi(\gamma) \, = \, \varphi(x_k \gamma x_k^{-1})
\, = \, \varphi(\gamma [\gamma,x_k])
\, = \, \langle \pi(y_k)\xi \mid \pi(\gamma^{-1})\xi \rangle.
$$
Therefore
$
\varphi(\gamma) = 
\lim_k \langle \pi(y_k)\xi \mid \pi(\gamma^{-1})\xi \rangle = 0. 
$
\end{proof}

\begin{lem}
% 12.C.2
\label{Lem-Aff-Norm}
Let $N$ be an \emph{infinite} abelian normal subgroup of the group $\Gamma$.
Assume that the $\Gamma$-conjugacy classes
contained in $N$ are $\{e\}$ and $N \smallsetminus \{e\}$. 
Let $\varphi$ be a normalized function of positive type on $N$ 
which is $\Gamma$-invariant.
\par

Then $\varphi = t 1_N+(1-t)\delta_e$ 
for some $t \in \mathopen[ 0,1 \mathclose]$.
\end{lem}

\begin{proof}
Since $\varphi$ is $\Gamma$-invariant, we have
$\varphi = t 1_N + s \delta_e$ for some $s, t \in \C$.
Moreover, $s + t = 1$, as $\varphi(e) = 1$. 
It remains to check that $t \in \mathopen[ 0,1 \mathclose]$
\par
 
Since $\varphi$ is a function of positive type on the abelian group $N$, 
by Bochner's theorem \ref{TheoremBochner},
$\varphi$ is the Fourier--Stieltjes transform $\mathcal{F}(\mu)$
of a probability measure $\mu$ on the compact dual group 
$\widehat N$ of $N$. 
Then, $\mu = t \delta_{e} + s \lambda_{\widehat N}$, 
where $\lambda_{\widehat N}$ is the
normalized Haar measure on $\widehat N$ 
and $\delta_{e}$ the Dirac measure at the group unit of $\widehat N$.
As $N$ is infinite, $\widehat N$ is infinite and $\lambda_{\widehat N}$ has no atoms. 
So, $\mu(\{e\}) = t \in \mathopen[ 0,1 \mathclose]$. 
\end{proof} 

We are now in position to describe the Thoma dual of $\Gamma$.

\begin{theorem}
% 12.C.3
\label{ThomaAffDiscret}
Let $\K$ be an infinite field and $\Aff(\K)$ the corresponding affine group.
A function $\varphi \,\colon \Gamma \to \C$ belongs to $E(\Gamma)$ if and only if 
$\varphi$ has one of the following forms:
\begin{enumerate}[label=(\arabic*)]
\item\label{iDEThomaAffDiscret}
the lift to $\Gamma$ of a unitary character 
of the quotient group $\Gamma/N$
(the latter isomorphic to the discrete group $\K^\times$); 
\item\label{iiDEThomaAffDiscret}
the Dirac function $\delta_e$ at the group unit.
\end{enumerate}
In other terms, $E(\Gamma)$ is in bijection with
$\widehat{\K^\times} \sqcup \{\delta_e\}$.
\par

Equivalently, every factor representation of finite type of $\Gamma$
is quasi-equivalent to one of:
\begin{enumerate}[label=(\arabic*)]
\addtocounter{enumi}{2}
\item\label{iiiDEThomaAffDiscret}
the composition of the quotient map $\Gamma \to \Gamma/N$
with a unitary character $\Gamma/N \to \T$, which is of type I$_1$;
\item\label{ivDEThomaAffDiscret}
the left regular representation of $\Gamma$, which is of type II$_1$.
\end{enumerate}
\end{theorem} 

\begin{proof}
A function $\varphi$ as in \ref{iDEThomaAffDiscret} of Theorem \ref{ThomaAffDiscret}
is clearly in $E(\Gamma)$.
For $\varphi = \delta_e$ as in \ref{iiDEThomaAffDiscret},
we have $\varphi \in \Tr_1(\Gamma)$ by Proposition \ref{diagcoeffinduced}, 
and indeed $\varphi \in E(\Gamma)$ because the group $\Gamma$ is icc. 
We have to show that, conversely, every $\varphi \in E(\Gamma)$
is either as in \ref{iDEThomaAffDiscret} or as in \ref{iiDEThomaAffDiscret}.
\par

We identify $\Gamma = \Aff(\K)$ with $\K^\times \ltimes \K$ 
and write $(a,s)$ for the matrix 
$$
\begin{pmatrix}
a & s
\\
0 & 1
\end{pmatrix}.
$$
The conjugacy class of an element $(a, s) \in \Gamma$ is 
$$
\left\{ (a, (1-a)t + \beta s) \mid \ t \in \K, b \in \K^\times \right\}.
$$
In particular, the conjugacy classes in $\Gamma$ 
contained in the normal subgroup
$N = \{(1,s) \mid \ s \in \K\} \approx \K$ 
are $\{0\}$ and $N \smallsetminus \{0\}$.
\par

Let $\varphi \in E(\Gamma)$. 
Consider the restriction $\psi := \varphi \vert_N$ of $\varphi$ to $N$.
Since $\psi$ is constant on $\Gamma$-conjugacy classes in $N$, it follows from
Lemma~\ref{Lem-Aff-Norm} that $\psi = t 1_N+ (1-t) \delta_e$ 
for some $t \in \mathopen[ 0,1 \mathclose]$.
We claim that $t = 0$ or $t = 1$.
\par

Indeed, let $(\pi, \Hi, \xi)$ be the GNS triple associated to $\varphi$.
Then, corresponding to the decomposition $\psi = t 1_N+ (1-t)\delta_e$,
we have an orthogonal decomposition 
$\Hi = \Hi^N \oplus \Hi_\infty$ into $\pi(N)$-invariant subspaces, where
$\Hi^N$ is the space of $N$-invariant vectors in $\Hi$.
On the one hand,
since $N$ is normal in $\Gamma$, the space $\Hi^N$ is $\pi(\Gamma)$-invariant.
On the other hand, 
$\Hi^N$ is clearly invariant under 
the commutant $\pi(\Gamma)'$ of $\pi(\Gamma)$ in ${\Li}(\Hi)$.
So, the orthogonal projection $p \,\colon \Hi \twoheadrightarrow \Hi^N$ onto $\Hi^N$
belongs to $\pi(\Gamma)' \cap \pi(\Gamma)''$, which is the centre of $\pi(\Gamma)''$.
Since $\pi(\Gamma)''$ is a factor, we have
either $p = 0$ or $p = \mathrm{Id}_\Hi$.
Therefore, $t = 0$ or $t = 1$, as claimed.
There are two cases to consider.

\vskip.2cm

$\circ$ \emph{First case.}
Assume that $t=1$, that is, $\varphi \vert_N = \psi = 1_N$. 
Then $\pi(s) = \mathrm{Id}_\Hi$ for all $s \in N$;
hence, $N$ is in the kernel of $\pi$.
So, $\varphi$ is the lift to $\Gamma$ 
of a unitary character of $\Gamma/N \approx \K^\times$.

\vskip.2cm

$\circ$ \emph{Second case.}
Assume that $t = 0$, that is, $\varphi \vert_N = \psi = \delta_e$.
We claim that $\varphi = \delta_e$.
Indeed, let $\gamma = (a, s) \in \Gamma \smallsetminus N$; then 
$a \ne 1$ and hence 
$$
\{[\gamma, t] \mid t \in N\} \, = \, \{(a-1)t \mid t \in \K\}
$$
is infinite. 
So, $\varphi(\gamma) = 0$, as claimed, by Lemma~\ref{Lem-ExtNorm}.
This ends the proof.
\end{proof}

The following corollary is an immediate consequence of 
Theorem~\ref{ThomaAffDiscret} and the description of
$\Pri(\Aff(\K))$ from Theorem~\ref{Theo-PrimIdealAffGr}.

\begin{cor}
% 12.C.4
\label{ThomaAff-PrimIdeal}
For $\Gamma = \Aff(\K)$, the natural map
$\kappa^{\rm norm}_{\rm prim} \,\colon E(\Gamma) \rightarrow \Pri(\Gamma)$ is a bijection
\end{cor}

Observe that, when $\K$ be a non-discrete topological field, the functions
$\varphi$ from Theorem~\ref{ThomaAffDiscret} are continuous 
on the topological group $\Aff(\K)$
only if $\varphi$ is the lift of a continuous unitary character
of $\Aff (\K) / N$. 

\begin{cor}
% 12.C.5
\label{ThomaAffTop}
Let $\K$ be a non-discrete topological field.
Consider now $\Aff(\K)$ as a topological group,
with its topology making it homeomorphic to $\K^\times \times \K$,
and therefore $\Aff(\K) / N$ topologically isomorphic to $\K^\times$.
\par

Every function in $E(\Aff(\K))$ is the lift of a unitary character
of the quotient group $\Aff(\K) / N$.
\par

Equivalently, every factor representation of finite type of $\Aff(\K)$
is the composition of the quotient map 
$\Aff(\K) \to \Aff(\K) / N \approx \K^\times$
with a unitary character $\K^\times \to \T$.
\end{cor}

\section
{Solvable Baumslag--Solitar groups}
% Section 12.D
\label{ThomaDualBS}

\index{Baumslag--Solitar group $\BS(1, p)$}
Let $p$ be a prime.
Consider the Baumslag--Solitar group
$\Gamma = \BS(1, p) = A \ltimes N \approx \Z \ltimes \Z[1/p]$,
as in Sections \ref{Section-IrrRepBS} and \ref{Section-PrimIdealBS}.
Recall that we can identify $\widehat N$ with the $p$-adic solenoid 
$$
\So_p \, = \, (\Q_p \times \R) / \Delta
\hskip.5cm \text{for} \hskip.2cm 
\Delta \, = \, \{(a,-a)\mid a \in \Z[1/p]\} ,
$$
and that the action of $A$ on $\widehat N$ corresponds to 
the action of $\Z$ on $\So_p$ given by the transformation $T_p \,\colon \So_p \to \So_p$
induced by the multiplication by $p$ (see Lemma~\ref{Lem-DualAdicInteger}).
\par 

Recall also that we parametrized in Corollary~\ref{Cor-FinDimRepBS}
the space $\widehat \Gamma_{\rm fd}$ 
of finite-dimensional irreducible representations of $\Gamma$
by the set $X_{\rm fd}/\sim$ of $T_p$-orbits in $X_{\rm fd}$, where
$$
X_{\rm fd} \, = \, \bigsqcup_{n \ge 1}
\big( \So_p(n) \times \mathopen[0,1/n \mathclose[ \big) .
$$
For $(s, \theta) \in \So_p(n) \times \mathopen[0,1/n \mathclose[$, we denote by 
$\pi_{s, \theta} = \Ind_{\Gamma(n)}^\Gamma \chi_{s, \theta}$ 
the corresponding finite-dimensional irreducible representation of $\Gamma$.
The associated character 
$\varphi_{s, \theta} \in E(\Gamma)$ of $\Gamma$
is the function defined by
$$
\varphi_{s, \theta}(\gamma) \, = \, \Tr (\pi_{s, \theta}(\gamma))
\hskip.5cm \text{for} \hskip.2cm 
\gamma \in \Gamma,
$$
where $\Tr(z)$ denotes the normalized trace of a matrix $z$.
\par

One can give an explicit formula for $\varphi_{s, \theta}$: 
for $n \ge 1$,
$(s, \theta) \in \So_p(n) \times \mathopen[0,1/n \mathclose[$, 
and $(k, b) \in \Gamma$, we have
$$
\varphi_{s, \theta}(k, b) \, = \, 
\begin{cases}
\hskip.2cm
\frac{1}{n} e^{2 \pi i k \theta /n} \sum_{\chi \in \mathcal O_s} \chi(b)
& \hskip.2cm \text{if} \hskip.2cm k \in n\Z
\\
\hskip.2cm
0
& \hskip.2cm \text{otherwise,}
\end{cases}
$$
where $\mathcal O_s$ denotes the $A$-orbit of $\chi_{s, \theta}$ in $\widehat N$.
\par

Denote by ${\rm EP}(\So_p)$ the set of 
ergodic $T_p$-invariant probability measures $\mu$ on $\So_p$, 
which are \emph{non-atomic} 
(that is, such that $\mu(\{s\}) = 0$ for every $s \in \So_p$).
We view the elements of ${\rm EP}(\So_p)$ 
as ergodic $A$-invariant probability measures on $\widehat N$, 
under the identification $\widehat N \approx \So_p$.
\par

Observe that, by Bochner Theorem \ref{TheoremBochner},
the (inverse) Fourier--Stieltjes transform 
$\overline{\mathcal{F}}(\mu)$ of a probability measure $\mu$ on $\widehat N$ 
is a function of positive type on $N$.
The trivial extension $\widetilde{\overline{\mathcal{F}}(\mu)}$
of $\overline{\mathcal{F}}(\mu)$ to $\Gamma$
is then a function of positive type on $\Gamma$ 
(Proposition~\ref{diagcoeffinduced}).

\begin{theorem}
% 12.D.1
\label{Theo-ThomaDualBS}
Let $\Gamma = A \ltimes N \approx \Z \ltimes \Z[1/p]$ 
be the Baumslag--Solitar group $\BS(1, p)$.
Set
$$
X \, = \, {\rm EP}(\So_p) \sqcup X_{\rm fd}
\, = \, {\rm EP}(\So_p) \sqcup
\bigsqcup_{n \ge 1}
\big( \So_p(n) \times \mathopen[0,1/n \mathclose[ \big) ,
$$
and identify $\widehat N$ with $\So_p$. 
\par

The map $\Phi \,\colon X \to E(\Gamma)$, defined by 
$$
\Phi(\mu) \, = \, \widetilde{\overline{\mathcal{F}}(\mu)} 
\hskip.5cm \text{for} \hskip.2cm
\mu \in {\rm EP}(\So_p)
$$
and
$$
\Phi(s, \theta) \, = \, \varphi_{s, \theta}
\hskip.5cm \text{for} \hskip.2cm
(s, \theta) \in X_{\rm fd} ,
$$
is a bijection.
\end{theorem}

\begin{proof}
$\bullet$ \emph{First step.}
We claim that every element from $E(\Gamma)$ is in the image of $\Phi$.
\par

Let $\varphi \in E(\Gamma)$; let $(\pi, \Hi, \xi)$ be a GNS triple
associated to $\varphi$.
Consider the projection-valued measure
$E \,\colon {\mathcal B}(\widehat N)$ on $\widehat N$
associated to $\pi \vert_N$.
Let $\mu = \mu_\xi$ be the probability measure on $\widehat N$
corresponding to $(E, \xi)$, which is defined by
$$
\mu(B) \, = \, \langle E(B)\xi \mid \xi \rangle
\hskip.5cm \text{for} \hskip.2cm 
B \in {\mathcal B}(\widehat N).
$$
Recall that, for every $B \in {\mathcal B}(\widehat N)$,
the projection $E(B)$ belongs to the von Neumann algebra $\pi(N)''$
and hence to $\pi(\Gamma)''$.
\par

Observe that the support of $\mu$ coincides with the support of $E$.
Indeed, let $B \in {\mathcal B}(\widehat N)$; if $E(B) = 0$, then obviously $\mu(B) = 0$.
Conversely, assume that $\mu(B) = 0$. 
Then 
$$
\Tr_\varphi(E(B)) \, = \, \langle E(B)\xi \mid \xi \rangle \, = \, 0,
$$
where $\Tr_\varphi$ is the trace on $\pi(\Gamma)''$ associated to $\varphi$
as in Construction~\ref{constructionGNS12}.
Therefore $E(B) = 0$ since $\Tr_\varphi$ is faithful.
\par

Since
$$
\pi(\gamma)E(B)\pi(\gamma)^{-1} \, = \, E(B^{\gamma^{-1}})
\hskip.5cm \text{for all} \hskip.2cm
B \in \mathcal B (\widehat N), \hskip.2cm \gamma \in \Gamma,
\leqno (*) 
$$
(see Proposition~\ref {Prop-RestNormalSub})
and since $\varphi$ is conjugation-invariant,
the measure $\mu$ is $\Gamma$-invariant.
\par

We claim that $\mu$ is ergodic under the $\Gamma$-action.
Indeed, let $B \in \mathcal B (\widehat N)$ be a $\Gamma$-invariant subset
with $\mu(B) > 0$.
Then $E(B)$ is a non-zero projection 
which belongs to the von Neumann algebra $\pi(\Gamma)''$. 
Moreover, relation $(*)$ above shows that $E(B)$ 
commutes with $\pi(\gamma)$ for every $\gamma \in \Gamma$.
So, $E(B)$ is in the centre of $\pi(\Gamma)''$. 
Therefore $E(B) = I$, since $\pi$ is factorial and $E(B) \ne 0$.
It follows that 
$$
E(\widehat N \smallsetminus B) \, = \, I - E(B) \, = \, 0
$$
and hence $\mu(\widehat N \smallsetminus B) = 0$; this proves the claim.
\par

Two cases may now occur.

\vskip.2cm

$\circ$ \emph{First case.} 
$\mu$ is atomic, that is, there exists $s \in \widehat N \approx \So_p$
such that $\mu(\{s\}) > 0$. 
Since $\mu$ is a $\Gamma$-invariant probability measure, 
it follows that $s$ is a periodic point for $T_p$ and the support of $\mu$, 
which is the is the support of $E$, 
is the $T_p$-orbit of $s$.
Now, since $\pi$ is factorial,
$\pi$ is weakly equivalent to an \emph{irreducible} representation $\rho$ of $\Gamma$
(Proposition~\ref{FacKK}~\ref{iiDEFacKK}).
Then $\pi \vert_N$ and $\rho \vert_N$ are weakly equivalent
and hence the projection-valued measure $E'$ associated to $\rho \vert_N$
has the same support as $E$
(see Proposition~\ref{Prop-ContAbelian}~\ref{iiDEProp-ContAbelian}).
So, the support of $E'$ is a finite $T_p$-orbit. By Theorem~\ref{Theo-PrimIdealBS},
$\rho$ has to be a finite-dimensional representation
of the form 
$$
\pi_{s, \theta} \, = \,
\Ind_{\Gamma(n)}^\Gamma \chi_{s, \theta}.
$$
It follows that $\pi$ is a multiple of $\rho$, by Proposition~\ref{PropOnWeAndEq}.
Therefore $\varphi$ coincides with 
the normalized trace $\varphi_{s, \theta}$ of $\rho = \pi_{s, \theta}$.

\vskip.2cm

$\circ$ \emph{Second case.}
The measure $\mu$ is non-atomic.
We claim that $\varphi = \widetilde{\overline{\mathcal{F}}(\mu)}$. 
\par

Indeed, on the one hand, for every $n \in N$, we have, by the SNAG Theorem,
$$
\varphi(n) \, = \, \langle \pi(n)\xi \mid \xi \rangle
\, = \, \int_{\widehat N} \chi(n) d\mu(n)
\, = \, \overline{\mathcal{F}}(\mu) (n),
$$
and so $\varphi = \overline{\mathcal{F}}(\mu)$ on $N$.
\par

On the other hand, let $\gamma \in \Gamma \smallsetminus N$.
Since $\mu(\Per(T_p)) = 0$, the action of the quotient group $\Gamma/N$
on the measure space $(\widehat N, \mu)$ is essentially free 
(see Section~\ref{SectionMSC} for this notion). 
Therefore by Proposition \ref{prop-EssentialFree}, there exists
 a sequence $(B_n)_{n}$ in $\mathcal B (\widehat N)$ with 
$$
\begin{aligned}
\Un_{B_n} \Un_{B_m}
\, &= \, 0 
\hskip.5cm \text{for} \hskip.2cm
n \ne m
\\
\sum_n \Un_{B_n}
\, &= \, \Un_{\widehat N} \
\hskip.2cm \text{and}
\\
\Un_{B_n}\theta_\gamma (\Un_{B_n})
\, &= \, \Un_{B_n}\Un_{B_n^\gamma} \, = \, 0
\hskip.5cm \text{for every} \hskip.2cm
n.
\end{aligned}
$$
It follows from the general properties of a projection-valued measure that 
$$
\begin{aligned}
E(B_n)E(B_m) &\, = \, 0
\hskip.5cm \text{for} \hskip.2cm
n \ne m
\\
\sum_n E(B_n) &\, = \, I 
\hskip.5cm \text{(in the strong operator topology)} \hskip.2cm
\text{and}
\\
E(B_n) E(B_n^{\gamma^{-1}})&\, = \, 0
\hskip.5cm \text{for every} \hskip.2cm 
n.
\end{aligned}
$$
Using the trace property of $\Tr_\varphi$, it follows that
$$
\langle \pi(\gamma)E(B_n)\xi \mid E(B_m)\xi \rangle \, = \,
\Tr_\varphi (E(B_n) \pi(\gamma) E(B_m)) \, = \,
\Tr_\varphi (E(B_m)E(B_n) \pi(\gamma)) \, = \, 0,
$$
for $n \ne m$. Therefore we have
$$
\begin{aligned}
\varphi(\gamma)&\, = \, \langle \pi(\gamma)\xi \mid \xi \rangle
\\
&\, = \, \langle \pi(\gamma) (\sum_n E(B_n) \xi) \mid \sum_n E(B_n) \xi \rangle
\\
&\, = \, \sum_{n, m} \langle \pi(\gamma) E(B_n) \xi \mid E(B_m) \xi \rangle
\\
&\, = \, \sum_{n} \langle \pi(\gamma) E(B_n) \xi \mid E(B_n) \xi \rangle
\\
&\, = \, \sum_{n} \langle E(B_n^{\gamma^{-1}})\pi(\gamma) \xi \mid E(B_n) \xi \rangle
\\
&\, = \, \sum_{n} \langle E(B_n) E(B_n^{\gamma^{-1}})\pi(\gamma) \xi \mid \xi \rangle
\\
&\, = \, 0.
\end{aligned}
$$
Summarizing, we have shown that $\varphi$ coincides with the trivial 
extension of $\overline{\mathcal{F}}(\mu)$ to~$\Gamma$.

\vskip.2cm

$\bullet$ \emph{Second step.} 
We claim that every function in the image of $\Phi$ 
belongs to $E(\Gamma)$. 
\par
 
This is clear for the functions of the form $\varphi_{s, \theta}$.
\par

Let now $\mu$ be an ergodic $\Gamma$-invariant non-atomic probability measure
on $\widehat N$ and set $\varphi = \widetilde{\overline{\mathcal{F}}(\mu)}$.
It is clear (compare Proposition \ref{diagcoeffinduced})
that $\varphi$ is in $\Tr_1(\Gamma)$; 
so, we have only to check that $\varphi$ is indecomposable.
\par

Let $\varphi_1, \varphi_2 \in \Tr_1(\Gamma)$ and $t \in \mathopen] 0,1 \mathclose[$
be such that $\varphi = t \varphi_1 + (1-t) \varphi_2$. 
As above, we have $\varphi_1 \vert_N = \overline{\mathcal{F}}(\mu_1)$
and $\varphi_2 \vert_N = \overline{\mathcal{F}}(\mu_1)$
for $\Gamma$-invariant probability measures $\mu_1$ and $\mu_2$ on $\widehat N$.
By injectivity of the Fourier--Stieltjes transform,
we have then $\mu = t \mu_1 + (1-t) \mu_2$.
Since $\mu$ is ergodic, it follows that $\mu_1 = \mu_2 = \mu$.
In particular, $\mu_1$ and $\mu_2$ are non-atomic.
\par

As in the second case of the first step, using the essential freeness
of the action of $\Gamma$ on $(\widehat N, \mu_i) = (\widehat N, \mu)$
and the trace property of $\Tr_{\varphi_i}$,
it follows that $\varphi_i = 0$ on $\Gamma \smallsetminus N$ for $i = 1, 2$.
Therefore $\varphi_1 = \varphi_2= \widetilde{\overline{\mathcal{F}}(\mu)}$.
\end{proof}

\begin{rem}
% 12.D.2
\label{Rem-Theo-ThomaDualBS}
An ergodic $T_p$-invariant probability measure on $\So_p$ which is atomic 
is necessarily the uniform measure on a periodic $T_p$-orbit; 
so, there are countably many such measures. 
\par

Observe that the normalized Haar measure on $\So_p$
is ergodic $T_p$-invariant and non-atomic 
(see Remark~\ref{Rem-Prop-IrredBS-bis}),
so that ${\rm EP}(\So_p)$ is non-empty. 
In fact, as we will show (see Corollary~\ref{Cor-InvMeasureBS}),
the set ${\rm EP}(\So_p)$ is uncountable.
\end{rem}

\begin{rem}
% 12.D.3
\label{Rem-Theo-ThomaDualLamplighter}
The atomic ergodic $T$-invariant probability measures on $X$ 
are the uniform measures on periodic $T$-orbits.
So, there are countably many such measures.
By contrast, the set ${\rm EP}(X)$ is uncountable,
as we now see (Proposition~\ref{Pro-InvMeasureLamplighter}).
There is no known classification of the measures in ${\rm EP}(X)$.
\end{rem}

For later use (see Corollary~\ref{Cor-InvMeasureBS}),
we need to work on a shift space $(X,T)$ over possibly more than two symbols, 
that is, $X = \{0, \hdots, \beta-1 \}^\Z$ for an integer $\beta \ge 2$,
equipped with the product topology
and the shift transformation $T \,\colon X \to X$.

\begin{prop}
% 12.D.4
\label{Pro-InvMeasureLamplighter}
Let $(X, T) $ be the two-sided shift over $\{0, \hdots, \beta-1 \}$
for an integer $\beta \ge 2$.
\par

There exist uncountably many
mutually singular ergodic $T$-invariant non-atomic probability measures on $X$,
all with (topological) support equal to $X$.
\end{prop}

\begin{proof}
For a real number $t \in \mathopen] 0,1 \mathclose[$,
let $\nu_t$ be a probability measure on $\{0, \hdots, \beta-1 \}$
such that $\nu_t(\{i\}) > 0$ for every $i \in \{0, \hdots, \beta-1 \}$;
below, we will impose a further condition depending on $t$.
Let $X$ be the space of sequences $(x_n)_{n \in \Z}$ 
with $x_n \in \{0, \hdots, \beta-1 \}$ for all $n \in \Z$..
This is a compact space for the product topology,
and a non-atomic measure space for the product measure
$\mu_t = \otimes^\Z \nu_t$ of copies of $\nu_t$.
\par

We have $\mu_t(U) > 0$ for ever non-empty open subset $U$ of $X$.
Indeed, such a set $U$ contains a cylinder of the form
$$
C_{a_1, \hdots, a_k}^{n_1, \hdots, n_k} \, = \,
\left\{ (x_n)_{n \in \Z} \in X \mid x_{n_i} = a_i
\hskip.2cm \text{for} \hskip.2cm
i =1, \hdots, k \right\},
$$
for integers $n_1 < \cdots < n_k$ in $\Z$
and elements $a_1, a_2, \hdots, a_k$ in $\{0, \hdots, \beta-1 \}$,
and 
$$
\mu_t(C_{a_1, \hdots, a_k}^{n_1, \hdots, n_k}) \, = \, \prod_{i = 1}^k\nu_t(a_i) > 0 ,
$$
In particular, the support of $\mu_t$ is $X$.
Denote by $T$ the shift map: $(Tx)_n = x_{n+1}$
for all $x = (x_m)_{m \in \Z} \in X$ and $n \in \Z$.
Then $\mu_t$ is $T$-invariant and ergodic (see \cite[Theorem 1.12]{Walt--82}).
\par

We assume now one more condition on $\nu_t$.
Its mean is equal to $t$, that is:
$$
\sum_{i = 0}^{\beta - 1} i\nu_t(i) \, = \, t .
$$
For definiteness, we could define $\nu_t$ by
$$
\nu_t(\{0\}) \, = \, 1 - \frac{2t}{\beta}
\hskip.5cm \text{and} \hskip.5cm
\nu_t(\{i\}) \, = \, \frac{2t}{\beta(\beta - 1)}
\hskip.2cm \text{for} \hskip.2cm
i \in \{1, \hdots, \beta - 1 \} .
$$
\par

It remains to show that the $\mu_t$~'s are mutually singular.
For every $t \in \mathopen] 0,1 \mathclose[$,
there exists by the law of large numbers
a measurable subset $A_t$ of $X$ with $\mu_t(A_t) = 1$ such that 
$$
\lim_{N \to + \infty} \frac{1}{N} \sum_{n = 1}^N x_n \, = \, t
\hskip.5cm \text{for all} \hskip.2cm
x = (x_n)_{n \in \Z} \in A_t.
$$
It follows that $A_s \cap A_t = \emptyset$
for $s, t \in \mathopen] 0,1 \mathclose[)$ with $s \ne t$,
and this ends the proof.
\end{proof}

Recall that every real number $x \in \mathopen[ 0,1 \mathclose[$
has a \emph{unique} $\beta$-expansion 
$$
x \, = \, 0.x_1 x_2 x_3 \hdots \, = \, \sum_{i = 1}^{+\infty} \frac{x_i}{\beta^i}
$$
with $x_i \in \{0, \hdots, \beta-1 \}$, provided we require that
$x_i \ne \beta-1$ for infinitely many $i$~'s. 
\par

Let again $p$ be a prime, $p \ge 2$.
Let $\So_p$ be the $p$-adic solenoid,
and $T_p \,\colon \So_p \to \So_p$ the multiplication by $p$,
as earlier in this section.
The following proposition shows that $(\So_p,T_p)$ is strongly related 
to the two-sided shift $(X, T)$ over $\{0, \hdots, p-1 \}$,
as a Borel space dynamical system
(that is, as a Borel space equipped with a Borel transformation).
\par

Denote by $X_1$ the measurable $T$-invariant subset of $X$
consisting of sequences $(x_n)_{n \in \Z}$ in $X$
which are eventually constant equal to $p-1$, that is, 
for which there exists an $N \ge 0$ such that $x_n = p-1$ for every $n \ge N$.

\begin{prop}
% 12.D.5
\label{Pro-BS-Lamplighter}
There is a Borel isomorphism
$$
\Phi \, \colon \, X \smallsetminus X_1 \to \So_p
$$
such that $\Phi(Tx) = T_p(\Phi(x))$ for every $X \smallsetminus X_1$.
\end{prop}

\begin{proof}
Recall from Section~\ref{Section-IrrRepBS}
that $\So_p = (\Q_p \times \R) / \Z[1/p]$,
with $\Z[1/p]$ identified with the subgroup
$\{(a,-a) \in \Q_p \times \R \mid a \in \Z[1/p]\}$. 
\par

Let $\Z_p$ is the ring of $p$-adic integers. The obvious map 
$$
\Z_p \times \R \to (\Q_p \times \R) / \Z[1/p]
$$
has $\Z$ as kernel (identified as above as subgroup of $\Q_p \times \R$)
and is clearly surjective.
It follows that $\So_p$ is isomorphic (as topological group) to $\Z_p \times \R/\Z$.
Therefore $\Z_p \times \mathopen[ 0,1 \mathclose[$
is a transversal for the $\Z$-cosets in $\Z_p \times \R/\Z$,
and so we have an identification of Borel spaces
$$
\So_p \, \cong \, \Z_p\times \mathopen[ 0,1 \mathclose[ .
$$
Under this identification, the map $T_p$ is given on
$\Z_p \times \mathopen[ 0,1 \mathclose[$ by 
$$
T_p(a, \lambda) \, = \, (pa + [p \lambda], \{p \lambda\})
$$
for $a \in \Z_2$ and $\lambda \in \mathopen[ 0,1 \mathclose[$,
where $[\alpha] \in \N$ and $\{\alpha\} \in \mathopen[ 0,1 \mathclose[$ denote the 
integer part and the fractional part of a real number $\alpha$. 
\par

We define a measurable map 
$$
\Phi \, \colon \, X \smallsetminus X_1 \to \Z_p\times \mathopen[ 0,1 \mathclose[
$$
as follows. 
For $x=(x_n)_{n \in \Z}\in X \smallsetminus X_1$, set 
$$
x^+ \, = \, (x_n)_{n \in \N}
\hskip.5cm \text{and} \hskip.5cm
x^- \, = \, (x_{-n})_{n \in \N^*}.
$$
Let $\lambda(x^+)$ be the real number in $\mathopen[ 0,1 \mathclose[$
with $p$-expansion given by the sequence $x^+$, that is, 
$$
\lambda(x^+) \, = \, 0.x_0x_1x_2 \cdots
\, = \, \sum_{n=0}^{+\infty} \frac{x_n}{p^{n+1}}
$$
and let $a(x^-)$ be the $p$-adic integer 
$$
a(x^-) \, = \, \sum_{n = 1}^{+\infty} x_{-n}p^{n-1}
\, = \, x_{-1}+ x_{-2} p+ x_{-3} p^2 + \cdots .
$$
Define 
$$
\Phi(x) \, := \, (a(x^-), \lambda(x^+))
\hskip.5cm \text{for} \hskip.5cm
x \in X \smallsetminus X_1.
$$
Since every element $\mathopen[ 0,1 \mathclose[$ has a unique $p$-expansion
by a sequence which is not eventually constant equal to $p-1$,
the map $\Phi$ is a bijection.
\par

The inverse map of $\Phi$ is given as follows.
For $(a,x) \in \Z_2 \times \mathopen[ 0,1 \mathclose[$, 
let 
$$
x \, = \, 0.x_0x_1x_2 \cdots \, = \, \sum_{i = 0}^{+\infty} \frac{x_i}{p^{i+1}}
$$
be the unique $p$-expansion of $x$ by a sequence
which is not eventually constant equal to $p-1$ and let 
$$
a \, = \, x_{-1} + x_{-2} p+ x_{-3} p^2 + \cdots
$$
be the $p$-adic expansion of $a$.
Then $\Phi^{-1}(a, \lambda) = (x^-, x^+)$, where
$$
x^+ \, = \, (x_0, x_1, x_2, \hdots)
\hskip.5cm \text{and} \hskip.5cm
x^- \, = \, (\hdots, x_{-3}, x_{-2}, x_{-1}).
$$
Since $\Phi^{-1}$ is obviously measurable,
$\Phi \,\colon X \smallsetminus X_1 \to \Z_2 \times \mathopen[ 0,1 \mathclose[$
is a Borel isomorphism.
\par

We claim that 
$$
\Phi(Tx) \, = \, T_p (\Phi(x))
\hskip.5cm \text{for} \hskip.5cm
x \, = \, (x_n)_{n \in \Z} \in X \smallsetminus X_1.
$$
Indeed, we have 
$$
(Tx)^+ \, = \, (x_{n+1})_{n \in \N}
\hskip.5cm \text{and} \hskip.5cm
(Tx)^- \, = \, (x_{-n+1})_{n \in \N^*}
$$
and so 
$$
\lambda((Tx)^+) \, = \, 0.x_1x_2x_3 \cdots
\hskip.5cm \text{and} \hskip.5cm
a((Tx)^-) \, = \, x_0 + x_{-1}p+ x_{-2} p^2 + \cdots.
$$
Therefore 
$$
\lambda((Tx)^+)\, = \, \{p\lambda(x^+)\}
\hskip.5cm \text{and} \hskip.5cm
a((Tx)^-) \, = \, p a(x^-) + x_0 \, = \, p a(x^-) + [p\lambda(x^+)]
$$
and it follows that
$$
\Phi(Tx) \, = \, (a((Tx)^-), \lambda ((Tx)^+)) \, = \, T_p(a(x^-), \lambda(x^+)) \, = \, T_p(\Phi(x)) ,
$$
as was to be shown.
\end{proof}

We now use
Propositions~\ref{Pro-InvMeasureLamplighter} and \ref{Pro-BS-Lamplighter}
in order to show that $\So_p$ has uncountable many
ergodic $T_2$-invariant non-atomic probability measures.

\begin{cor}
% 12.D.6
\label{Cor-InvMeasureBS}
Let $(\So_p,T_p)$ be the $p$-adic solenoid.
There exist uncountably many mutually singular
ergodic $T_p$-invariant non-atomic probability measures on $\So_p$,
all with (topological) support equal to $\So_p$.
\end{cor}

\begin{proof}
For $t \in (0,t)$, let $\nu_t$ and $\mu_t$
be the probability measures on $\{0, \hdots, p-1 \}$ and $X$
defined in the proof of Proposition~\ref{Pro-InvMeasureLamplighter}.
\par

We claim that $\mu_t(X_1) = 0$,
where $X_1$ is the measurable subset of $X$
as in Proposition~\ref{Pro-BS-Lamplighter}.
Indeed, for $N \ge 0$, let $X_1^N$ be the set of sequences
$(x_n)_{n \in \Z}$ in $\{0, \hdots, p-1 \}^\Z$ such that
$x_n = p-1$ for every $n \ge N$.
Then 
$$
\mu_t(X_1^N) \, = \, \lim_{n\to +\infty} (\nu_t(\{p-1 \}))^{n-N} \, = \, 0
$$
and hence $\mu_t(X_1) = 0$, since $X_1 = \bigcup_{N \ge 0} X_1^N$.
As a result, we may view $\mu_t$ as probability measure on $X \smallsetminus X_1$.
\par

Let $\nu_t = \Phi_*(\mu_t)$ be the image of $\mu_t$ under the map 
$$
\Phi \, \colon \, X \smallsetminus X_1 \to \So_p
$$
from Proposition~\ref{Pro-BS-Lamplighter}.
The probability measure $\nu_t$ is $T_p$-invariant, ergodic and non-atomic,
since $\mu_t$ is ergodic $T$-invariant and non-atomic
and since $\Phi$ is an equivariant Borel isomorphism. 
\par

Moreover, the $\nu_t$~'s are mutually singular measures on $\So_p$,
since the $\mu_t$~'s are mutually singular measures on 
$X \smallsetminus X_1$. 
\par

We claim that $\mathrm{supp}(\nu_t) = \So_p$.
Indeed, let $U$ be a non-empty open subset of $\So_p$. 
Then there exists $N \ge 2$ such that
image $U'$ of $U$ in $\Z_p \times [\mathopen[ 0,1 \mathclose[$ 
contains $p^N \Z_2 \times [1/p^N, 1/p^{N - 1})$.
Therefore $\Phi^{-1}(U')$ contains the cylinder 
$$
C_N\, = \, \left\{
(x_n)_{n \in \Z} \in X \smallsetminus X_1
\hskip.2cm \Big\vert \hskip.2cm
\begin{aligned}
&x_{-N} = \cdots = x_{-1} = 0
\hskip.2cm \text{and}
\\
&x_0 = \cdots = x_{N - 2} = 0, \hskip.2cm x_{N - 1} \ne 0
\end{aligned}
\right\} .
$$
Since $\mu_t(C_N) > 0$, it follows that $\nu_t(U) > 0$.
\end{proof}

\section
{Lamplighter group}
% Section 12.E
\label{ThomaDualLamplighter}

\index{Lamplighter group}
Recall from Sections \ref{Section-IrrRepLamplighter}
and \ref{Section-PrimIdealLamplighter}
that the lamplighter group is the semi-direct product $\Gamma = A \ltimes N$
of $A = \Z$ with $N = \bigoplus_{k \in \Z} \Z / 2 \Z$,
where the action of $\Z$ on $\bigoplus_{k \in \Z} \Z / 2 \Z$
is given by shifting the coordinates.
Recall also that $\widehat N$ can be identified
with $X = \{0,1 \}^\Z$, the dual action $\Z$ on $\widehat N$
being given by the shift transformation $T$ on $X$.
\par

As in Section~\ref{Section-IrrRepLamplighter},
$X(n)$ denotes the set of points in $X$
with $T$-period $n \ge 1$. 
Recall from Proposition~\ref{Prop-IrredRepLamplighter-FiniteDim}
that the space $\widehat \Gamma_{\rm fd}$ 
of finite-dimensional irreducible representations of $\Gamma$ is parametrized
by the set $X_{\rm fd} / \sim$ of $T$-orbits in $X_{\rm fd}$, where
$$
X_{\rm fd} \, = \, \bigsqcup_{n \ge 1}
\big(X(n) \times \mathopen[ 0,1/n \mathclose[ \big) .
$$
For $(x, \theta) \in X(n) \times \mathopen[ 0,1/n \mathclose[$,
we denote by $\pi_{x, \theta} = \Ind_{\Gamma(n)}^\Gamma \chi_{x, \theta}$ 
the corresponding finite-dimensional irreducible representation of $\Gamma$
and by $\varphi_{x, \theta}$ the associated character in $E(\Gamma)$. 
\par

Denote by ${\rm EP}(X)$
the set of ergodic $T$-invariant non-atomic probability measures on $X$,
which we view as ergodic $A$-invariant probability measures on $\widehat N$, 
under the identification $\widehat N \approx X$.
\par

The following result can be proved in exactly the same way
as Theorem~\ref{Theo-ThomaDualBS}.

\begin{theorem}
% 12.E.1
\label{Theo-ThomaDualLamplighter}
Let $\Gamma = A \ltimes N$ 
be the lamplighter group.
Set
$$
X \, = \, {\rm EP}(X) \sqcup X_{\rm fd}
\, = \, {\rm EP}(X) \sqcup
\bigsqcup_{n \ge 1}
\big( X(n) \times \mathopen[ 0,1/n \mathclose[ \big) ,
$$
and identify $\widehat N$ with $X$. 
\par

The map $\Phi \,\colon X \to E(\Gamma)$, defined by 
$$
\Phi(\mu) \, = \, \widetilde{\overline{\mathcal{F}}(\mu)} 
\hskip.5cm \text{for} \hskip.2cm
\mu \in {\rm EP}(X)
$$
and
$$
\Phi(x, \theta) \, = \, \varphi_{x, \theta}
\hskip.5cm \text{for} \hskip.2cm
(X, \theta) \in X_{\rm fd} ,
$$
is a bijection.
\end{theorem}

\section
{General linear groups}
% Section 12.F
\label{ThomaDualGL}

Let $\K$ be an infinite field and $n$ an integer, $n \ge 2$.
Recall that the commutator group 
of the general linear group $\GL_n(\K)$
is the special linear group $\SL_n(\K)$,
and the centre $Z$ of $\GL_n(\K)$ consists of 
the scalar matrices $\lambda I_n$, for $\lambda \in \K^\times$.
The following result is from \cite{Kiri--65}.
\index{General linear group! $\GL_n(\K)$ with $\K$ a field} 

\begin{theorem}[\textbf{Kirillov}]
% 12.F.1
\label{Theo: GLn}
Let $\GL_n(\K)$ be as above.
A function $\varphi \,\colon \GL_n(\K) \to \C$ belongs to $E(\GL_n(\K))$
if and only if $\varphi$ has one of the following forms:
\begin{enumerate}[label=(\arabic*)]
\item\label{iDETheo: GLn}
the lift to $\GL_n(\K)$ of a unitary character of 
$\GL_n(\K)/\SL_n(\K) \approx \K^\times$,
\item\label{iiDETheo: GLn}
the trivial extension $\widetilde \chi$ 
of a unitary character $\chi$ of the centre $Z \approx \K^\times$ of $\GL_n(\K)$. 
\end{enumerate}
In other terms, $E(\GL_n(\K))$ is in natural bijection 
with $\widehat{\K^\times} \sqcup \widehat Z$.
\par

Equivalently, 
every factor representation of finite type of $\GL_n(\K)$
is quasi-equivalent to one of:
\begin{enumerate}[label=(\arabic*)]
\addtocounter{enumi}{2}
\item\label{iiiDETheo: GLn}
the composition of the quotient map 
$\GL_n(\K) \twoheadrightarrow \GL_n(\K)/\SL_n(\K) \approx \K^\times$
with a unitary character $\K^\times \to \T$;
\item\label{ivDETheo: GLn} 
the induced representation $\Ind_Z^{\GL_n(\K)} \chi$ 
for a unitary character $\chi$ of $Z$.
\end{enumerate}
\end{theorem}

\noindent
\textbf{Proof of Theorem~\ref{Theo: GLn} for $n = 2$.}
% 12.F.1 Proof of $n=2$
We write $\Gamma$ for $\GL_2(\K)$,
and we consider the following subgroups of $\Gamma$ :
$$
\begin{aligned}
P&\, = \, \left\{
\begin{pmatrix}
a & b \\ 0 & 1
\end{pmatrix}
\hskip.1cm \Big\vert \hskip.1cm
a \in \K^\times, \hskip.1cm b \in \K 
\right\}
\\
P^{t}&\, = \, \left\{
\begin{pmatrix}
a & 0 \\ c & 1
\end{pmatrix}
\hskip.1cm \Big\vert \hskip.1cm
a \in \K^\times, \hskip.1cm c \in \K 
\right\}
\\
N&\, = \, \left\{
\begin{pmatrix}
1 & b \\ 0 & 1
\end{pmatrix}
\hskip.1cm \Big\vert \hskip.1cm
b \in \K 
\right\}
\\
N^{t}&\, = \, \left\{
\begin{pmatrix}
1 & 0 \\ c & 1
\end{pmatrix}
\hskip.1cm \Big\vert \hskip.1cm
c \in \K 
\right\}
\\
A&\, = \, P\cap P^{t} = \left\{
\begin{pmatrix}
a & 0 \\ 0 & 1
\end{pmatrix}
\hskip.1cm \Big\vert \hskip.1cm 
a \in \K^\times 
\right\}.
\end{aligned}
$$
Observe that $P = A \ltimes N$ and $P^t = A \ltimes N^t$ are isomorphic 
to the group $\Aff(\K)$ of affine transformations of $\K$.
\par

Let $\varphi \in E(\Gamma)$.
The restrictions $\varphi\vert_P$ and $\varphi\vert_{P^t}$
to $P$ and $P^t$ are in $\Tr_1(P)$ and $\Tr_1(P^t)$, respectively.
By Lemma \ref{Lem-Aff-Norm},
there exist $\varphi_1 \in \Tr_1(P), \varphi_2 \in \Tr_1(P^t)$ 
with $\varphi_1 \vert_N = 1, \varphi_2 \vert_{N^t} = 1$ 
and real numbers $t_1,t_2 \in \mathopen[ 0,1 \mathclose]$ such that
$$
\varphi \vert_P \, = \, t_1 \varphi_1+ (1-t_1) \delta_e
\hskip.5cm \text{and} \hskip.5cm
\varphi \vert_{P^{t}} \, = \, t_2 \varphi_2+ (1-t_2) \delta_e. 
$$
Let $(\pi, \Hi, \xi)$ be the GNS triple associated to $\varphi$.
Then, corresponding to the two previous decompositions of $\varphi$,
we have two decompositions of $\Hi$ under the action of $P$ and $P^t$
respectively:
$$
\Hi \, = \, \Hi_{(1)} \bigoplus \Hi^\infty_{(1)}
\hskip.2cm \text{and} \hskip.2cm
\Hi \, = \, \Hi_{(2)} \bigoplus \Hi^\infty_{(2)};
$$
The group $P$ acts on $\Hi^\infty_{(1)}$ 
as a multiple of the regular representation $\lambda_P$ of $P$
and, similarly, $P^t$ acts on $\Hi^\infty_{(2)}$ as a multiple of $\lambda_{P^t}$.
Moreover, $N$ acts as the identity on $\Hi_{(1)}$ 
and $N^t$ acts as the identity on $\Hi_{(2)}$.
(Some of the subspaces $\Hi_{(i)}$ or $\Hi^\infty_{(i)}$ may be trivial.)
\par

Since the restrictions $\lambda_{P} \vert_{N}$ and $\lambda_{P^t} \vert_{N^t}$
are multiples of $\lambda_{N}$ and $\lambda_{N^t}$
and since $N$ and $N^t$ are infinite, 
$\Hi^\infty_{(1)}$ [respectively $\Hi^\infty_{(2)}$]
contains no non-zero vector invariant under $N$ [respectively $N^t$]. 
It follows that $\Hi_{(1)}$ is the space $\Hi^N$ of $\pi(N)$-invariant vectors 
and $\Hi_{(2)}$ is the space $\Hi^{N^t}$ of $\pi(N^t)$-invariant vectors in $\Hi$.

\vskip.2cm
We claim that $\Hi^N \cap \Hi^{N^t}$ 
is invariant by $\pi(\Gamma)'' \cup \pi(\Gamma)'$.
\par

Indeed, $\Hi^N \cap \Hi^{N^t}$ is a closed subspace of $\Hi$ 
which is invariant under $\pi(A)$,
as $A$ normalizes $N$ and $N^t$. 
So, $\Hi^N\cap \Hi^{N^t}$ is invariant under $\pi(P)$ and $\pi(P^{t})$,
hence under $\pi(\Gamma)$
since $P \cup P^t$ generates $\Gamma = \GL_2(\K)$,
hence under $\pi(\Gamma)''$.
Moreover, $\Hi^N$ and $\Hi^{N^t}$ are invariant under
the commutants $\pi(N)'$ and $\pi(N^t)'$ of $\pi(N)$ and $\pi(N^t)$, 
respectively;
% \marginpar{Argument ou r\'ef !!!}
% \footnote{c'est \'evident: si $x$ est un point fixe d'une application $f$ 
% d'un ensemble quelconque dans lui-m\^eme et $g$
% est une application qui commute avec $f$, alors $g(x)$ est un point fixe de $f$.}
therefore $\Hi^N \cap \Hi^{N^t}$
is invariant under $\pi(\Gamma)'$.

\vskip.2cm

\noindent
$\bullet$ \emph{First case.} 
Assume that $\Hi^N \cap \Hi^{N^t} \ne \{0\}$.
We claim that $\varphi(\gamma) = 1$ for all $\gamma \in \SL_2(\K)$.
\par

Indeed, since $\Hi^N \cap \Hi^{N^t}$ 
is invariant by $\pi(\Gamma)'' \cup \pi(\Gamma)'$,
the orthogonal projection $p$ onto $\Hi^N \cap \Hi^{N^t}$
belongs to the centre of $\pi(\Gamma)''$. 
However, $\pi$ is a factor representation, as $\varphi \in E(\Gamma)$.
It follows that $p=I$, that is, $\Hi = \Hi^N\cap \Hi^{N^t}$.
This shows that $\varphi(\gamma) = 1$ for all
$\gamma$ either in $N$ or in $N^t$, 
and therefore also for all $\gamma$ in the subgroup $\SL_2(\K)$
generated by $N \cup N^t$.

\vskip.2cm

\noindent
$\bullet$ \emph{Second case.} 
Assume that $\Hi^N\cap \Hi^{N^t}= \{0\}$.
Set 
$$
\Ki \, := \, \Hi^\infty_{(1)} + \Hi^\infty_{(2)}.
$$
Then $\Ki$ is a dense subspace of $\Hi$, 
because $\Ki^\perp = \Hi^N\cap \Hi^{N^t} = \{0\}$.
\par

Since $P$ acts on $\Hi^\infty_{(1)}$ 
as a multiple of the regular representation $\lambda_P$,
all matrix coefficients of $\pi$ restricted to $\Hi^\infty_{(1)}$ vanish at infinity on $P$,
that is, all functions
$$
P \, \to \, \C, \hskip.2cm 
p \, \mapsto \, \langle \pi(p) \eta \mid \eta' \rangle 
$$
for $\eta, \eta' \in \Hi^\infty_{(1)}$ vanish at infinity;
a similar fact is true for the restriction of $\pi \vert_{P^t}$ to $\Hi^\infty_{(2)}$.
Since $A= P\cap P^t$, 
it follows that \emph{all} matrix coefficients of $\pi$ vanish at infinity on $A$.
\par

Let $\gamma \in \Gamma \smallsetminus Z$. We claim that $\varphi(\gamma) = 0$.
\par

To show this, observe first that $\gamma$ is conjugated inside $\Gamma$ 
to a matrix of the form
$\begin{pmatrix}
0 & b \\ 1& d
\end{pmatrix}$,
for $b \in \K^\times$ and $d \in \K$.
Indeed, choose $v \in \K^2$ such that $\gamma(v)$ is not a multiple to $v$ 
(this is possible, since $\gamma \notin Z$);
if $\delta \in \Gamma$ maps the canonical basis of $\K^2$
onto $\{v, \gamma(v)\}$, then
$\delta^{-1} \gamma \delta =
\begin{pmatrix} 0 & b \\ 1& d \end{pmatrix}$.
So, it suffices to show that $\varphi(\gamma) = 0$ for 
$\gamma = \begin{pmatrix} 
0 & b \\ 1& d
\end{pmatrix}$.
Set $\varepsilon = \begin{pmatrix} 0 & 1 \\ 1 & 0 \end{pmatrix}$.
\par

Let $(a_k)_k$ be a sequence of pairwise distinct elements in $\K^\times$; 
for all $k$, set 
$$
\delta_k \, = \, \begin{pmatrix} 
a_k & 0 \\ 0 & 1
\end{pmatrix} 
\hskip.5cm \text{and} \hskip.5cm
\gamma_k = \begin{pmatrix} 
a_k^{-1} & d \\ 0 & a_k b
\end{pmatrix} .
$$
Then
$$
\delta_k \gamma \delta_k^{-1} 
\, = \, 
\begin{pmatrix} 
0 & a_kb \\ a_k^{-1} & d
\end{pmatrix}
\, = \,
\begin{pmatrix} 
0 & 1 \\ 1 & 0
\end{pmatrix}
\begin{pmatrix} 
a_k^{-1} & d \\ 0 & a_k b
\end{pmatrix}
\, = \,
\varepsilon \gamma_k ,
$$
and we have $\varphi(\gamma) = \varphi(\varepsilon \gamma_k)$.
Write $\gamma_k = z_k \alpha_k n_k$ with
$$
z_k \, = \, 
\begin{pmatrix} 
a_k b & 0 \\ 0 & a_k b
\end{pmatrix} \, \in \, Z,
\hskip.2cm
\alpha_k \, = \, 
\begin{pmatrix} 
a_k^{-2} b^{-1} & \\ 0 & 1
\end{pmatrix} \, \in \, A,
\hskip.2cm \text{and} \hskip.2cm
n_k \, = \, 
\begin{pmatrix} 
1 & da_k \\ 0 & 1
\end{pmatrix} \, \in \, N.
$$
Since $z_k\in Z$, we have, using Lemma \ref{Lem-CharHeis},
$$
\varphi(\gamma) 
\, = \, \varphi(\varepsilon \gamma_k)
\, = \, \varphi(\varepsilon z_k \alpha_k n_k)
\, = \, \varphi( z_k) \varphi(\varepsilon \alpha_k n_k) .
$$
Set $\eta := \pi(\varepsilon)\xi$ and
write the orthogonal decompositions 
$$
\xi \, = \, \xi_1+ \xi_2
\hskip.2cm \text{and} \hskip.2cm
\eta \, = \, \eta_1+ \eta_2
$$
of $\xi$ and $\eta$, 
with $\xi_1, \eta_1 \in \Hi^N$ and $\xi_2, \eta_2 \in \Hi^{\infty}_{(1)}$; 
for every $k$, we have
$$
\begin{aligned}
\varphi(\varepsilon \alpha_k n_k)
&\, = \, \langle \pi(\alpha_k n_k) \xi \mid \eta \rangle 
\\
&\, = \, \langle \pi(\alpha_k n_k)\xi_1 \mid \eta_1 \rangle
 + \langle \pi(\alpha_kn_k) \xi_2 \mid \eta_2 \rangle
\\
&\, = \, \langle \pi(\alpha_k) \xi_1 \mid \eta_1 \rangle
 + \langle \pi(\alpha_kn_k) \xi_2 \mid \eta_2 \rangle.
\end{aligned}
$$
On the one hand, 
since the matrix coefficients of $\pi$ vanish at infinity on $A$,
we have
$$
\lim_k \langle \pi(\alpha_k) \xi_1, \mid \eta_1 \rangle \, = \, 0.
$$
On the other hand, 
since the matrix coefficients of $\pi$ restricted to $\Hi^\infty_{(1)}$ 
vanish at infinity on $P$, 
we have 
$$
\lim_k \langle \pi(\alpha_kn_k) \xi_2 \mid \eta_2 \rangle \, = \, 0.
$$
Therefore,
$$
\varphi(\gamma) \, = \, 
\lim_k \varphi(z_k) \big(\langle \pi(\alpha_k) \xi_1 \mid \eta_1 \rangle
 + \langle \pi(\alpha_kn_k) \xi_2 \mid \eta_2 \rangle \big) 
 \, = \, 0. 
$$
\par

It remains to show that, if a function $\varphi$ on $\Gamma$ 
has one of the forms stated in Theorem~\ref{Theo: GLn},
then $\varphi \in E(\Gamma)$. 
This is obvious if 
$\varphi$ is the lift of a unitary character of $\GL_2(\K)/SL_2(\K)$.
Assume that $\varphi = \widetilde \chi$ is the trivial extension 
of a unitary character $\chi$ of the centre $Z$. 
\par

Since the function $\widetilde \chi$ is in $\Tr_1(\Gamma)$ 
by Proposition \ref{diagcoeffinduced},
we have to check that it is indecomposable.
Let $\varphi_1, \varphi_2 \in \Tr_1(\Gamma)$ 
and $t \in \mathopen] 0,1 \mathclose[$
be such that 
$$
\widetilde \chi \, = \, t \varphi_1+ (1-t) \varphi_2.
$$
By what we have just seen, for $j = 1, 2$, we can write
$$
\varphi_j \, = \, t_j \psi^{(j)}_1 +(1-t_j)\psi^{(j)}_2 ,
$$
where $t_j \in \mathopen[ 0,1 \mathclose]$, 
$\psi^{(j)}_1 \in \Tr_1(\Gamma)$ is equal to $1$ on $\SL_2(\K)$,
and $\psi^{(j)}_2 \in \Tr_1(\Gamma)$ is equal to $0$ on $\Gamma \smallsetminus Z$.
It follows that can assume that each $\varphi_j$ is 
either equal to $1$ on $\SL_2(\K)$ or equal to $0$ on $\Gamma \smallsetminus Z$.
Since $\varphi$ vanishes outside $Z$, 
we necessarily have that both $\varphi_1$ and $\varphi_2$ 
are $0$ on $\Gamma \smallsetminus Z$.
Then, as in the proof of Theorem \ref{ThomaAffDiscret},
% \marginpar{Compl\'eter !!! Ancien~: 9.F.6.}
$\varphi_1$ and $\varphi_2$ have to agree with $\chi$ on $Z$.
So, $\varphi_1 = \varphi_2 = \varphi$ and hence
$\varphi \in E(\GL_n(\K))$.
$\square$

\vskip.2cm

\noindent
\textbf{Proof of Theorem~\ref{Theo: GLn} for $n \ge 3$.}
% 12.F.1 Proof of $n \ge 3$
We denote by $e_1, \hdots, e_n$ the vectors of the canonical basis of $\K^n$.
For every $i \in \{1, \hdots, n\}$, 
we consider the copy $V_i$ of $\K^{n-1}$ inside $\Gamma = \GL_n(\K)$
consisting of the matrices $\gamma \in \SL_n(\K)$ with 
$$
\gamma(e_j) \, = \, e_j 
\hskip.5cm \text{for all} \hskip.2cm 
j \ne i.
$$
Observe that $V_i$ is normalized by 
a copy $G_i$ of $\GL_{n-1}(\K) \ltimes \K^{n-1}$ inside $\Gamma$.
We will have also to consider the transpose subgroups $V_i^t$; 
of course, $V_i^t$ is normalized by $G_i^t$.
Observe that $V_j \cap V_i^t$ is a non-trivial subgroup
of $\Gamma$ for all $i \ne j$.
We will refer to subgroups of the form $V_i$ or $V_i^t$
as to the copies of $\K^{n-1}$ inside $\Gamma$.
\par

Let $\varphi \in E(\Gamma)$. 
Fix a copy $V$ of $\K^{n-1}$ inside $\Gamma$.
Lemma~\ref{Lem-Aff-Norm} applied to $\GL_{n-1}(\K) \ltimes \K^{n-1}$ 
(with $N=V$) shows that 
$$
\varphi \vert_V \, = \, t 1_V+ (1-t) \delta_e,
$$
for some $t \in \mathopen[ 0,1 \mathclose]$.
Let $(\pi, \Hi, \xi)$ be the GNS triple associated to $\varphi$.
We have a corresponding 
orthogonal decomposition of $\Hi$ under the action of $V$ 
$$
\Hi \, = \, \Hi^{V} \bigoplus \Hi^\infty_{V},
$$
where $V$ acts trivially on $\Hi^V$ and
as multiple of the regular representation on $\Hi^\infty_{V}$. 

\vskip.2cm

$\bullet$ {\it First step.} 
Let $W$ be another copy of $\K^{n-1}$
with $V \cap W \ne \{0\}$. We claim that $\Hi^{V} = \Hi^W$ and
$\Hi_{\infty}^{V} = \Hi_{\infty}^{W}$.
Indeed, since the restriction to $A = V\cap W$ of $\pi$ to $\Hi_{\infty}^{V}$ 
(or $\Hi_{\infty}^{W}$)
is a multiple of the regular representation of $A$, we have
$\Hi^A = \Hi^V = \Hi^W$.

\vskip.2cm

$\bullet$ {\it Second step.} 
Let $W$ be an arbitrary copy of $\K^{n-1}$.
We claim that $\Hi^{V} = \Hi^W$ and $\Hi_{\infty}^{V} = \Hi_{\infty}^{W}$.
Indeed, we can find two copies $W^1$ and $W^2$ of $\K^{n-1}$ with
$$
V \cap W^1 \, \ne \, \{0\}, 
\hskip.5cm 
W^1 \cap W^2 \, \ne \, \{0\}, 
\hskip.5cm \text{and} \hskip.5cm
W^2 \cap W \, \ne \, \{0\}.
$$
The claim follows from the first step.

\vskip.2cm

$\bullet$ {\it Third step.} 
Assume that $\Hi^V \ne \{0\}$. 
We claim that $\varphi =1$ on $\SL_n(\K)$.
Indeed, by the second step, $\Hi^{V} = \Hi^{W}$ 
for every copy $W$ of $\K^{n-1}$ 
and hence $\Hi^V = \Hi^{\SL_n(\K)}$, since $\SL_n(\K)$ 
is generated by the various $W$~'s. 
In particular, $\Hi^V$ is invariant under $\pi(\Gamma)$, 
as $\SL_n(\K)$ is normal in $\Gamma$.
Moreover, $\Hi^V$ is also invariant under
the commutant $\pi(\Gamma)'$ of $\pi(\Gamma)$.
Therefore $\Hi^V = \Hi$, as $\pi$ is a factor representation.

\vskip.2cm

$\bullet$ {\it Fourth step.} 
Assume that $\Hi^V = \{0\}$. 
We claim that $\varphi = 0$ on $\Gamma \smallsetminus Z$.
\par

Indeed, by the second step, $\Hi = \Hi^\infty_V$ 
and hence $\varphi \vert_V = \delta_e$.
It follows from Lemma~\ref{Lem-ExtNorm}, that $\varphi \vert_{G_V} = \delta_e$ for 
the subgroup $G_V \approx \GL_{n-1}(\K) \ltimes \K^{n-1}$ normalizing $V$.
\par

Let $\gamma \in \Gamma \smallsetminus Z$. 
As in the last part of the proof of the case $n=2$,
upon passing to a conjugate of $\gamma$ in $\Gamma$, 
we can find $\varepsilon \in \Gamma$, 
a sequence $(z_k)_k$ in $Z$, 
and a sequence $(\gamma_k)_k$ of pairwise distinct elements from $G_V$ 
such that $\gamma = z_k \varepsilon \gamma_k$.
Since 
$$
\langle \pi(\gamma_k)\xi \mid \pi(\gamma_l)\xi \rangle
\, = \, \varphi(\gamma_l^{-1} \gamma_k) \, = \, 0
$$
for $k \ne l$, 
the sequence $(\pi(\gamma_k)\xi)_k$ is weakly convergent to $0$ in $\Hi$ 
and therefore $\varphi(\gamma) = 0$.
\par

Finally, one checks, as in the case $n=2$, that every function 
on $\Gamma$ which has one of the forms stated in Theorem~\ref{Theo: GLn}
is in $E(\Gamma)$.
\hfill $\square$

\vskip.2cm

Observe that, for $\psi \in \widehat Z$,
the induced representation $\Ind_Z^\Gamma \psi$
is a factor representation of $\Gamma$ with character $\widetilde \chi$
(the extension of $\chi$ to $\Gamma$ by $0$).
The following corollary is therefore an immediate consequence of 
Theorem~\ref{Theo: GLn} and Theorem~\ref{Theo-PrimIdealPGLn}.

\begin{cor}
% 12.F.2
\label{Cor-PrimGLn}
Let $\K$ be an infinite algebraic extension of a finite field,
$n \ge 2$, and $\Gamma= \GL_n(\K)$.
The map $\kappa^{\rm norm}_{\rm prim} \,\colon E(\Gamma) \rightarrow \Pri(\Gamma)$ is a bijection.
\end{cor}

Consider now the group $\Gamma = \SL_n(\K)$.
Its centre $Z$ is finite, it consists of the matrices $\lambda I_n$
with $\lambda$ a $n$th root of unity in $\K$, that is,
$Z \approx \mu_n^\K := \{ \lambda \in \K^* \mid \lambda^n = 1 \}$. 
The quotient $\PSL_n(K)$ is an infinite simple group.
The proof of Theorem~\ref{Theo: GLn} 
carries over without change to $\Gamma = \SL_n(\K)$ for $n \ge 3$,
after replacing there the copies of $\GL_{n-1}(\K) \ltimes \K^{n-1}$
by copies of $\SL_{n-1}(\K) \ltimes \K^{n-1}$.
However, for $\Gamma = \SL_2(\K)$, the proof involves
further arguments; see \cite{PeTh--16}.
The final result is that, for $n \ge 2$,

\begin{theorem}
% 12.F.3
\label{Theo-SLn}
Let $\K$ be an infinite field and $n$ an integer, $n \ge 2$. Then
$$
E(\SL_n(\K)) \, = \, 
\{\Un_{\SL_n(\K)}\} \sqcup \{ \widetilde \chi \mid \chi \in \widehat Z \},
$$
where $\widetilde \chi$ denotes the trivial extension 
of a unitary character $\chi$ of the centre $Z \approx \mu_n^\K$ 
of $\SL_n(\K)$. 
\par

Consequently, we have
$$
E(\PSL_n(\K)) \, = \, \{\Un_{\PSL_n(\K)}, \ \delta_e\}.
$$
\end{theorem}

%-----------------------------------------------------------------------
% End of chapter 12
%-----------------------------------------------------------------------
\chapter
{The group measure space construction}
% Chapter 13
\label{ChapterGroupMeasureSpace}

\emph{
The previous chapters have shown the importance
of finding representations of a discrete group $\Gamma$
which generate factors of various types. 
As a preliminary step, we discuss in this chapter
the group measure space construction of Murray and von Neumann \cite{MuvN--36},
which plays a fundamental role in the theory of operator algebras.
This construction associates to an action of $\Gamma$
on a standard measure space $(X, \mathcal B)$,
equipped with a quasi-invariant measure $\mu$,
a von Neumann algebra $L^\infty(X, \mu) \rtimes \Gamma$.
The group measure space construction is described in detail
in Section \ref{sectionNoinvariantmeasure}.
}

\emph{
Assume that the $\Gamma$-action on $(X,\mu)$ is ergodic and essentially free. 
We will show that $L^\infty(X, \mu) \rtimes \Gamma$ is a factor 
and that its type depends on properties of the measure $\mu$.
In particular, this factor is of type $III$
when there is no $\Gamma$-invariant on $(X, \mathcal B)$ which is equivalent to $\mu$. 
In order to check this last condition, we discuss in Section \ref{sectionNoinvariantmeasure}
a criterion based on the essential range of the action of $\Gamma$ on $(X,\mu)$.
This allows us to show that
$L^\infty(\mathbf{P}^1(\R), \mu) \rtimes \Gamma$ is a factor of type III$_1$,
for any lattice $\Gamma$ in $\PSL_2(\R)$. 
}

\emph{
The results in this chapter will be used in Chapter \ref{Chap:NormalInfiniteRep}
in order to construct factor representations of certain countable groups.
}

\section
{Construction of factors}
% Section 13.A
\label{SectionMSC}

The main result of this section is Theorem \ref{ROIV}.
Our exposition follows
\cite[Chap.~I, \S~9]{Dixm--vN} and \cite[Chapter V, Section 7]{Take--79}.

\subsection*{Countable group actions on standard Borel spaces}

Consider a standard measure space $(X, \mathcal B, \mu)$,
i.e., a standard Borel space $(X, \mathcal B)$
together with a $\sigma$-finite measure $\mu$.
To this is naturally associated
a commutative von Neumann algebra $L^\infty(X, \mu)$
acting on the Hilbert space $L^2(X, \mu)$.
\index{Standard! $2$@measure}
\par

Consider moreover a countable group $\Gamma$ acting on $X$
by Borel isomorphisms $X \to X$ such that 
$\mu$ is quasi-invariant by the $\Gamma$-action.
We assume that this action is \emph{on the right},
and we write $(X, \mathcal B, \mu) \curvearrowleft \Gamma$
or more simply $(X, \mu) \curvearrowleft \Gamma$.
\par

Such an action $(X, \mu) \curvearrowleft \Gamma$ is
\begin{itemize}
\setlength\itemsep{0em}
\item[$\bullet$]
\textbf{essentially free} 
if $\mu(F) = 0$, 
where 
$$
F \, := \, \{x \in X \mid x\gamma = x
\hskip.2cm \text{for some} \hskip.2cm
\gamma \ne e
\hskip.2cm \text{in} \hskip.2cm
\Gamma \},
$$
\index{Essentially free action}
\index{Action! essentially free}
\item[$\bullet$]
\textbf{ergodic}
if either $\mu(Y) = 0$ or $\mu(X \smallsetminus Y) = 0$ 
for every $\Gamma$-invariant subset $Y \in \mathcal B$
(ergodic actions are defined in Section \ref{AppGpactions}
for general topological groups),
\index{Ergodic! action}
\index{Action! ergodic}
\item[$\bullet$]
\textbf{essentially transitive} 
if there exists an orbit $Y = x\Gamma $
in $X$ such that 
\hfill\par\noindent
$\mu(X \smallsetminus Y) = 0$.
\index{Essentially transitive action}
\index{Action! essentially transitive}
\end{itemize}
The action of $\Gamma$ on $(X, \mathcal B, \mu)$ 
provides a \emph{left action} $\gamma \mapsto \theta_\gamma$
of $\Gamma$ on $L^\infty(X, \mu)$
by automorphisms of von Neumann algebras, 
defined by
$$
\theta_\gamma(\varphi) (x) \, = \,
\varphi ( x \gamma)
\hskip.5cm \text{for all} \hskip.2cm
\gamma \in \Gamma, \hskip.1cm \varphi \in L^\infty(X, \mu) ,
\hskip.2cm \text{and} \hskip.2cm x \in X .
\leqno (*)
$$
It will be useful to characterize the essential freeness 
of an action of $\Gamma$ on $(X, \mathcal B, \mu)$
in terms of the associated action
on the von Neumann algebra $L^\infty(X, \mu)$.
In this section, we will only use the equivalence
of \ref{iDEprop-EssentialFree} and \ref{iiiDEprop-EssentialFree}; 
the equivalence with \ref{iiDEprop-EssentialFree} 
has been used earlier,
in the proof of Theorem \ref{Theo-ThomaDualBS}.

\begin{prop}
% 13.A.1
\label{prop-EssentialFree}
Let $\Gamma$ be a countable group acting
on a standard measure space $(X, \mathcal B, \mu)$
by measurable transformations for which the class of $\mu$ is invariant.
The following properties are equivalent:
\begin{enumerate}[label=(\roman*)]
\item\label{iDEprop-EssentialFree}
the action of $\Gamma$ on $(X, \mathcal B, \mu)$ is essentially free;
\item\label{iiDEprop-EssentialFree}
for every $\gamma \in \Gamma \smallsetminus \{e\}$, 
there exists a sequence $(Y_n)_{n \ge 1}$ in $\mathcal B$
with the following properties:
$$
\begin{aligned}
\Un_{Y_n} \Un_{Y_m} \, &= \, 0
\hskip1cm \text{for} \hskip.2cm 
m, n \ge 1 \hskip.2cm \text{such that} \hskip.2cm m \ne n,
\\
\sum_{n \ge 1} \Un_{Y_n} \, &= \, \Un_{X} 
\\
\Un_{Y_n} \theta_\gamma(\Un_{Y_n}) \, &= \, 0
\hskip1cm \text{for every} \hskip.2cm 
n \ge 1 
\end{aligned}
$$
(where the characteristic function $\Un_Y$ of $Y \in \mathcal B$
is viewed as an element of the von Neumann algebra $L^\infty(X, \mu)$);
\item\label{iiiDEprop-EssentialFree}
if $\varphi \in L^\infty(X, \mu)$ 
and $\gamma \in \Gamma \smallsetminus\{e\}$
are such that $\varphi \theta_\gamma(\psi) = \varphi \psi$
for every $\psi \in L^\infty(X, \mu)$, then $\varphi = 0$.
% \item[(iv)]
% for every $\gamma \in \Gamma \smallsetminus\{e\}$ 
% and for every $Y \in \mathcal B$ with $\mu(Y) > 0$,
% there exists $Z\in \mathcal B$ such that $Z \subset Y$, $\mu(Z) > 0$, and 
% $\gamma \cdot Z \cap Z = \emptyset$. 
\end{enumerate}
\end{prop}

\begin{proof}
\ref{iDEprop-EssentialFree} $\Rightarrow$ \ref{iiDEprop-EssentialFree}
Since $X$ is a standard Borel space, 
there exists a sequence $(B_n)_{n \ge 1}$ in $\mathcal B$
which separates the points of $X$.
\par

Fix $\gamma \in \Gamma \smallsetminus \{e\}$.
Set $X^\gamma := \{ x \in X \mid x\gamma = x \}$
and $X'^\gamma = X \smallsetminus X^\gamma$.
For every $x \in X'^\gamma$, we have $x \ne x\gamma$,
and so there exists $n \ge 1$ such that $x \in B_n$ 
and $x\gamma \notin B_n$.
Therefore 
$$
X'^\gamma \, = \,
\bigcup_{n \ge 1} (X'^\gamma\cap (B_n \smallsetminus B_n \gamma^{-1} )).
$$
For every $n \ge 1$,
set $Z_n := X'^\gamma \cap (B_n \smallsetminus B_n \gamma^{-1} )$.
We have $Z_n \cap Z_n\gamma = \emptyset$.
Define recursively a sequence of Borel sets $(Y_n)_{n \ge 1}$ by 
$$
% Y_0 \, = \, \emptyset
% \hskip.5cm \text{and} \hskip.5cm
Y_n \, = \, Z_n \smallsetminus \bigcup_{k < n} Z_k
\hskip.5cm \text{for} \hskip.2cm
n \ge 1 .
$$
Then 
$$
\begin{aligned}
Y_m \cap Y_n \, &= \, \emptyset
\hskip1cm \text{for} \hskip.2cm
m, n \ge 1 \hskip.2cm \text{such that} \hskip.2cm m \ne n,
\\
X'^\gamma
\, &= \, \bigcup_{n \ge 0} Y_n 
\\
Y_n \cap Y_n \gamma
\, &= \, \emptyset
\hskip1cm \text{for every} \hskip.2cm 
n \ge 1 .
\end{aligned}
\leqno (\sharp)
$$
\par

Assume now that \ref{iDEprop-EssentialFree} holds.
Then $\mu(X \smallsetminus X'^\gamma) = \mu(X^\gamma) = 0$,
so that $\Un_{X'^\gamma} = \Un_{X} \in L^\infty(X, \mu)$,
and ($\sharp$) shows that \ref{iiDEprop-EssentialFree} holds.

\vskip.2cm

\ref{iiDEprop-EssentialFree} $\Rightarrow$ \ref{iiiDEprop-EssentialFree}
Assume that there exists a sequence $(Y_n)_{n \ge 1}$
as in \ref{iiDEprop-EssentialFree}.
Let $\varphi \in L^\infty(X, \mu)$
and $\gamma \in \Gamma \smallsetminus \{e\}$ 
be such that 
$\varphi \theta_\gamma(\psi) = \varphi \psi$ for all $\psi \in L^\infty(X, \mu)$.
Then, for all $n \ge 1$, we have
$\varphi \theta_\gamma(\Un_{Y_n}) = \varphi \Un_{Y_n}$ and $$
0 \, = \,
\varphi \left( \Un_{Y_n} \theta_*(\Un_{Y_n}) \right) \, = \,
\left( \varphi \theta_*(\Un_{Y_n}) \right) \Un_{Y_n} \, = \,
\left( \varphi \Un_{Y_n} \right) \Un_{Y_n} \, = \,
\varphi \Un_{Y_n} ,
$$
so that $\Un_{Y_n}\varphi = 0$. Therefore
$$
\varphi \, = \, \sum_{n \ge 1} \varphi \Un_{Y_n} \, = \, 0
$$
and this shows that \ref{iiDEprop-EssentialFree}
implies \ref{iiiDEprop-EssentialFree} .

\vskip.2cm

\ref{iiiDEprop-EssentialFree} $\Rightarrow$ \ref{iDEprop-EssentialFree}
Let $\gamma \in \Gamma \smallsetminus \{e\}$.
For $x \in X^\gamma$ and $\psi \in L^\infty(X, \mu)$, we have 
$$
\theta_\gamma(\psi) (x) \, = \, \psi(x\gamma) \, = \, \psi(x).
$$
Therefore for $\varphi := 1_{X^\gamma}$, we have
$\varphi\theta_\gamma(\psi) = \varphi \psi$.
By assumption \ref{iiiDEprop-EssentialFree}, 
we have therefore $\varphi = 0$, that is, 
$\mu(X^\gamma) = 0$.
This shows that \ref{iiiDEprop-EssentialFree}
implies \ref{iDEprop-EssentialFree}.
\end{proof}

We also translate the ergodicity of a group action
on a measure space $(X, \mu)$ in terms of
the associated action on $L^\infty(X, \mu)$.

\begin{prop}
% 13.A.2
\label{prop-ergodicity}
Let $\Gamma$ be a countable group acting 
on a standard measure space $(X, \mathcal B, \mu)$
by measurable transformations for which the class of $\mu$ is invariant.
The following properties are equivalent:
\begin{enumerate}[label=(\roman*)]
\item\label{iDEprop-ergodicity}
the action of $\Gamma$ on $(X, \mathcal B, \mu)$ is ergodic;
\item\label{iiDEprop-ergodicity}
every $\varphi \in L^\infty(X, \mu)$
such that $\theta_\gamma(\varphi) = \varphi$
for all $\gamma \in \Gamma$ is constant $\mu$-almost everywhere.
\end{enumerate}
\end{prop}

\noindent
\emph{Note:} in \ref{iiDEprop-ergodicity},
one can equivalently replace $L^\infty(X, \mu)$
by the space of complex-valued measurable functions on $X$
up to equality $\mu$-almost everywhere.

\begin{proof}
Assume that \ref{iiDEprop-ergodicity} holds
and let $Y \in \mathcal B$ be $\Gamma$-invariant.
Then $\theta_\gamma(\Un_Y) = \Un_Y$ for all $\gamma \in \Gamma$,
and hence $\Un_Y$ is constant $\mu$-almost everywhere.
This implies that $\mu(Y) = 0$ or $\mu(X \smallsetminus Y) = 0$.

\vskip.2cm

Conversely, assume that 
the action of $\Gamma$ on $(X, \mathcal B, \mu)$ is ergodic.
Let $\varphi \in L^\infty(X, \mu)$ 
be such that $\theta_\gamma(\varphi) = \varphi$,
for every $\gamma \in \Gamma$. 
Upon considering the real and imaginary part of $\varphi$,
we may assume that $\varphi$ is real-valued.
\par

For every $\gamma \in \Gamma$,
we have $\mu(N_\gamma) = 0$, where 
$$
N_\gamma \, = \, 
\{ x \in X \mid \varphi(x\gamma) \ne \varphi(x)\}.
$$
Since $\Gamma$ is countable, we have $\mu(N) = 0$,
for $N := \bigcup_{\gamma \in \Gamma \smallsetminus \{e\}} N_\gamma$.
Let $\varphi' \,\colon X \to \R$ be defined by $\varphi' = \varphi$ on 
$X \smallsetminus N$ and $\varphi' = 0$ on $N$;
then $\varphi'$ coincides with $\varphi$
almost everywhere and $\varphi'$ is $\Gamma$-invariant: 
$\varphi'(x\gamma) = \varphi'(x)$ for every $x \in X$
and $\gamma \in \Gamma$.
So, me may assume that $\varphi$ is $\Gamma$-invariant.
\par

Consider the ``distribution function" 
$F \,\colon \R \to \mathopen[ 0, +\infty \mathclose[ \cup \{+\infty \}$ of $\varphi$,
defined by 
$$
F(t) \, = \, \mu ( \{ x \in X \mid \varphi(x) \le t \} )
\hskip.5cm \text{for} \hskip.2cm
t \in \R.
$$
On the one hand,
$F$ is an increasing function with 
$$
\lim_{t \to -\infty}F(t) \, = \, 0
\hskip.5cm \text{and} \hskip.5cm 
\lim_{t \to +\infty}F(t) \, = \, \mu(X).
$$
On the other hand, 
since $\{x \in X \mid \varphi(x) \le t \}$ is $\Gamma$-invariant,
we have $F(t) = 0$ or $F(t) = \mu(X)$ for every $t \in \R$,
by ergodicity of the $\Gamma$-action.
Then 
$$
t_0 \, := \, \inf \{ t \in \R \mid F(t) = \mu(X) \}
$$
exists in $\R$ and 
$\mu ( \{ x \in X \mid \varphi(x) = t_0 \} ) = \mu(X)$.
Therefore $\varphi$ is constant $\mu$-almost everywhere.
\end{proof}

\begin{rem}
% 13.A.3
\label{caracterisationesstrans}
Let $\Gamma$ be a countable group acting 
on a standard measure space $(X, \mathcal B, \mu)$.
The following properties are equivalent:
\begin{enumerate}[label=(\roman*)]
\item\label{iDEcaracterisationesstrans}
the action of $\Gamma$ on $(X, \mathcal B, \mu)$ is essentially transitive;
\item\label{iiDEcaracterisationesstrans}
the action is ergodic and the measure $\mu$ is atomic.
\end{enumerate}
\end{rem}

\begin{defn}
% 13.A.4
\label{defKoopman}
\index{Koopman representation}
\index{Representation! Koopman}
Let $\Gamma$ be a group 
acting on a standard measure space $(X, \mathcal B, \mu)$, on the left.
Denote by $\gamma \mapsto u_\gamma$
the associated \textbf{Koopman representation}
of $\Gamma$ on $L^2(X, \mu)$, which is defined by 
$$
(u_\gamma f) (x) \, = \, d(\gamma^{-1}, x)^{1/2} f(\gamma^{-1} x) 
\hskip.5cm \text{for all} \hskip.2cm
\gamma \in \Gamma, \hskip.1cm
f \in L^2(X, \mu), \hskip.1cm
x \in X ,
$$
where $d(\gamma, x) = \frac{d (h^{-1})_*(\mu)}{d \mu} (x)$ denotes
the Radon--Nikodym cocycle of the action
$\Gamma \curvearrowright (X, \mathcal B, \mu)$,
as in Appendix \ref{AppMeasureB} and in Definition \ref{defnRN}.
% see for example in \cite{KuSt--96} and \cite{Dudk--18}.
\par

In case of a \textbf{right} group action $(X, \mathcal B, \mu) \curvearrowleft \Gamma$
by measurable transformations preserving the class of $\mu$,
the relevant Radon--Nikodym cocycle is given by
$c(x, \gamma) = \frac{d \mu_*(\gamma^{-1}) }{d \mu} (x)$,
and the Koopman representation $\gamma \mapsto u_\gamma$
of $\Gamma$ on $L^2(X, \mu)$ is defined by
$$
(u_{\gamma}f) (x) \, = \, c(x, \gamma)^{1/2}f(x\gamma) 
\hskip.5cm \text{for all} \hskip.2cm
\gamma \in \Gamma, \hskip.1cm
f \in L^2(X, \mu), \hskip.1cm
x \in X,
$$
where $c(x, \gamma) = \frac{d \mu_*(\gamma^{-1})}{d \mu} (x)$.
\par

The definition of Koopman representations can be extended
to $\sigma$-compact locally compact groups acting by non-singular automorphisms
on measure spaces $(X, \mathcal B, \mu)$ such that $L^2(X, \mu)$ is a separable Hilbert space;
these technical conditions, $\sigma$-compactness and separability,
are needed for the associated representation to be continuous.
See Proposition A.6.1 in \cite{BeHV--08}.
\par

The terminology refers to an article by Koopman \cite{Koop--31},
where the author considers the time evolution
on an energy level in the $(p, q)$-space of an Hamiltonian system
and associates to this situation a ``one-parameter group of unitary transformations''
in the appropriate Hilbert space, i.e., a representation of $\R$.
\end{defn}

\subsection*{Associated representations and the group measure space construction}

There is also a representation $\varphi \mapsto m(\varphi)$
of the abelian von Neumann algebra $\AC := L^\infty(X, \mu)$
on $L^2(X, \mu)$ by multiplication operators: 
$$
m(\varphi) (f) (x) \, = \, \varphi(x) f(x)
\hskip.5cm \text{for all} \hskip.2cm
\varphi \in \AC, \hskip.1cm f \in L^2(X, \mu), \hskip.1cm x \in X .
$$
The following relation holds between
the representations $\gamma \mapsto u_\gamma$ of $\Gamma$
and $\varphi\mapsto m(\varphi)$ of $L^\infty(X, \mu)$:
$$
u_\gamma m(\varphi) u_{\gamma}^{-1}
\, = \, m(\theta_\gamma(\varphi)) 
\hskip.5cm \text{for all} \hskip.2cm
\gamma \in \Gamma, \hskip.1cm \varphi \in \AC.
\leqno (**)
$$
Denote by $\Li$ the von Neumann algebra
acting on the Hilbert space $L^2(X, \mu)$ and generated by
$$
\{u_\gamma \in \U(L^2(X, \mu)) \mid \gamma \in \Gamma\}
\, \cup \,
\{m(\varphi) \in \Li (L^2(X, \mu)) \mid \varphi \in \AC)\} .
$$

\begin{prop}
% 13.A.5
Let $\Gamma$ be a countable group acting on the right
on a standard measure space $(X, \mathcal B, \mu)$
and let $\mathcal L$ be the von Neumann algebra on $L^2(X, \mu)$ defined as above.
Assume that the action $(X, \mathcal B, \mu) \curvearrowleft \Gamma$
is ergodic.
\par

Then $\mathcal L$ is the algebra $\Li (L^2(X, \mathcal B, \mu))$
of all bounded operators on $L^2(X, \mathcal B, \mu)$.
In particular, $\mathcal L$ is a factor of type I.
\end{prop}

\begin{proof}
The subset
$$
\{ m(\varphi) \in \Li (L^2(X, \mu)) \mid \varphi \in \AC \}
$$
is a maximal abelian subalgebra of $\Li (L^2(X, \mathcal B, \mu))$,
by Corollary~\ref{Cor-AbelianVN}.
\par

Let $T \in \Li'$.
On the one hand, $T$ commutes with $m(\varphi)$ for all $\varphi \in \AC$,
hence $T = m(\psi)$ for some $\psi \in \AC$.
On the other hand, $T = m(\psi)$ commutes with
$\{u_\gamma \in \U(L^2(X, \mu)) \mid \gamma \in \Gamma\}$,
and it follows from relation $(**)$ that $\psi$ is $\Gamma$-invariant.
By Proposition \ref{prop-ergodicity},
$\psi$ is constant almost everywhere, that is, $T$ is a scalar operator.
Therefore $\Li' = \C \mathrm{Id}_{L^2(X, \mu)}$, and $\Li = \Li (L^2(X, \mu))$.
\end{proof}

In order to construct factors of other types,
we consider now representations of $\Gamma$ and $\AC$ 
on the Hilbert space $\ell^2(\Gamma, L^2(X, \mu))$,
namely:
\begin{itemize}
\setlength\itemsep{0em}
\item[$\bullet$]
the representation $\gamma \mapsto U_\gamma$
of the group $\Gamma$ defined by
$$
U_\gamma(F) (\gamma') \, = \,
u_\gamma(F(\gamma' \gamma))
$$
for all $\gamma, \gamma' \in \Gamma$
and $F \in \ell^2(\Gamma, L^2(X, \mu))$,
\item[$\bullet$]
the representation $\varphi \mapsto M(\varphi)$ of $\AC$ defined by
$$
M(\varphi) (F) (\gamma) \, = \, m(\varphi) (F(\gamma))
$$
for all $\varphi \in \AC$, $F \in \ell^2(\Gamma, L^2(X, \mu))$,
and $\gamma \in \Gamma$.
\end{itemize}
The following relations hold:
$$
U_\gamma M(\varphi) U_{\gamma}^{-1} \, = \, M(\theta_\gamma(\varphi)) 
\hskip.5cm \text{for all} \hskip.2cm 
\gamma \in \Gamma
\hskip.2cm \text{and} \hskip.2cm
\varphi \in \AC.
\leqno(***)
$$

\index{Crossed product} 
\index{Group measure space construction, see crossed product}
\index{von Neumann algebra! crossed product}
\index{von Neumann algebra! group measure space construction}
By definition, the \textbf{group measure space von Neumann algebra} 
or the \textbf{crossed product von Neumann algebra} 
associated with the $\Gamma$-action on $(X, \mathcal B, \mu)$ 
is the von Neumann algebra denoted by\footnote{We
feel obliged to respect tradition and write $\Gamma$ \emph{on the right}
in $L^\infty(X, \mu) \rtimes \Gamma$,
even though the action of $\Gamma$ on $L^\infty(X, \mu)$ is here a \emph{left action}.}
$$
L^\infty(X, \mu) \rtimes \Gamma ,
\hskip.5cm \text{or} \hskip.5cm
\AC \rtimes \Gamma ,
\hskip.5cm \text{or below} \hskip.5cm
\M ,
$$
of operators on $\ell^2(\Gamma, L^2(X, \mu))$ generated by 
$$
\{ M(\varphi) \mid \varphi \in \AC \} 
\, \cup \,
\{ U_\gamma \mid \gamma \in \Gamma \}.
$$
Observe that, as relation $(***)$ shows, 
the linear subspace 
$$
{\M}_0 \, := \, \{\sum_{k = 1}^n M(\varphi_k) U_{\gamma_k}
\mid
n \in \N, \hskip.1cm \varphi_k \in \AC, \hskip.1cm \gamma_k \in \Gamma\}
$$
is an algebra, which is moreover selfadjoint
and weakly dense in $\M$.
\par

We identify $\ell^2(\Gamma, L^2(X, \mu))$ with the tensor product
$\ell^2(\Gamma) \otimes L^2(X, \mu)$.
For $T \in \Li (\ell^2(\Gamma, L^2(X, \mu)))$, we denote by 
$(T_{s, t})_{s, t \in \Gamma}$ the matrix associated to $T$
with respect to the decomposition
$$
\ell^2(\Gamma) \otimes L^2(X, \mu)
\, = \, 
\bigoplus_{\gamma \in \Gamma} (\C \delta_\gamma \otimes L^2(X, \mu)) ,
$$
where $(\delta_\gamma)_{\gamma \in \Gamma}$
is the orthonormal basis of $\ell^2(\Gamma)$
constituted of Dirac delta functions.
Thus, $T_{s, t} \in \Li (L^2(X, \mu))$ is defined by 
$$
\langle T_{s, t} f_1 \mid f_2 \rangle
\, = \,
\langle T(\delta_t \otimes f_1) \mid \delta_s \otimes f_2 \rangle
\hskip.5cm \text{for} \hskip.2cm
s, t \in \Gamma
\hskip.2cm \text{and} \hskip.2cm
f_1, f_2 \in L^2(X, \mu).
$$
\par

Let $\varphi \in \AC$ and $\gamma \in \Gamma$. 
The matrix $(M_{s, t})_{s, t \in \Gamma}$ of $M(\varphi)$ is given by 
$$
M_{s, t} \, = \,
\begin{cases} 
m(\varphi) & \text{for} \hskip.2cm s = t
\\
0 & \text{otherwise,}
\end{cases}
$$
and the matrix $(U_{s, t})_{s, t \in \Gamma}$
of $U_\gamma$ is given by 
$$
U_{s, t} \, = \,
\begin{cases} 
u_\gamma & \text{if} \hskip.2cm s t^{-1} = \gamma
\\
0 & \text{otherwise}.
\end{cases}
$$
Let $T = \sum_{k = 1}^n M(\varphi_{\gamma_k}) U_{\gamma_k} \in \M_0$
with $\gamma_k \ne \gamma_l$ for $k \ne l$.
The matrix $(T_{s, t})_{s, t \in \Gamma}$ of $T$ is given by 
$$
T_{s, t} \, = \,
\begin{cases} 
m(\varphi_{\gamma_k}) u_{\gamma_k} & \text{if} \hskip.2cm s t^{-1} = \gamma_k 
\in \{\gamma_1, \hdots, \gamma_n\}
\\
0 & \text{otherwise}.
\end{cases}
$$
This last property is preserved by going over to weak limits: 
for every $T \in \M$, there exists a family 
$\{\varphi_\gamma\}_{\gamma \in \Gamma}$ in $\AC$
such that the matrix of $T$ is
$(m(\varphi_{st^{-1}}) u_{st^{-1}})_{s, t \in \Gamma}$.
\par

Henceforth, for $T \in \M$, we will write 
$T = \sum_{\gamma \in \Gamma} \varphi_{\gamma} u_{\gamma}$
instead of writing the matrix $(m(\varphi_{st^{-1}}) u_{st^{-1}})_{s, t \in \Gamma}$.
We will also often consider $\AC$ as a subalgebra of $\M$,
by identifying a function $\varphi \in \AC$
with the corresponding operator $M(\varphi) \in \M$.
\par

Observe that,
for $S = \sum_{\gamma \in \Gamma} \varphi_{\gamma} u_{\gamma}$
and $T = \sum_{\gamma \in \Gamma} \psi_{\gamma} u_{\gamma}$ in $\M$,
we have, by matrix multiplication in combination with relation $(**)$ above:
$$
S^* \, = \, 
\sum_{\gamma \in \Gamma} 
\theta_\gamma (\overline{\varphi_{\gamma^{-1}}}) u_{\gamma} 
$$
and
$$
ST \, = \,
\sum_{\gamma \in \Gamma}
\bigg( \sum_{\delta \in \Gamma}\varphi_{\delta} 
%\theta_\delta(\psi_{\delta\gamma}) \bigg) u_\gamma,
\theta_\delta(\psi_{\delta^{-1}\gamma}) \bigg) u_\gamma,
$$
where the infinite sum in parentheses above is weakly convergent in $\AC$.

\begin{lem}
% 13.A.6
\label{Lem-FreeActionMaxAbelian}
Assume that the action of the countable group $\Gamma$ 
on the standard measure space $(X, \mathcal B, \mu)$
is essentially free.
\par

Then the von Neumann subalgebra $M(\AC) \approx \AC$
of $\M = L^\infty(X, \mu) \rtimes \Gamma$
is maximal abelian.
\end{lem}

\begin{proof}
Let $T = \sum_{\gamma \in \Gamma} \varphi_{\gamma} u_{\gamma} \in \M$
be in the commutant of $M(\AC)$.
For every $\varphi \in \AC$, we have 
$TM(\varphi) = M(\varphi)T$ and therefore
$\varphi_{\gamma} \theta_\gamma (\varphi)u_{\gamma} \, = \,
\varphi \varphi_{\gamma} u_{\gamma}$,
that is
$$
\varphi_{\gamma} \theta_\gamma (\varphi) \, = \, \varphi \varphi_{\gamma}
\hskip.5cm \text{for all} \hskip.2cm
\varphi \in \AC \hskip.2cm \text{and} \hskip.2cm \gamma \in \Gamma .
$$
Since the action of $\Gamma$ is essentially free,
it follows from Proposition \ref{prop-EssentialFree}
that $\varphi_{\gamma } = 0$ for every 
$\gamma \in \Gamma \smallsetminus \{e\}$.
This shows that $T = M(\varphi_{e})$
and the claim is proved.
\end{proof}

\index{Conditional expectation}
We define a linear map 
$$
E \, \colon \, \M = L^\infty(X, \mu) \rtimes \Gamma
\, \longrightarrow \,
\AC = L^\infty (X, \mu),
$$
called the \textbf{conditional expectation} onto $\AC$, by
$$
E(T) \, = \, \varphi_e 
\hskip.5cm \text{for} \hskip.5cm 
T \, = \, \sum_{\gamma \in \Gamma} \varphi_\gamma u_\gamma \in \M.
$$
The following properties of $E$ are easily checked.

\begin{lem}
% 13.A.7
\label{Lem-CondExpCrossProduct}
Let $E \,\colon \M \to \AC$ be as above.
\begin{enumerate}[label=(\arabic*)]
\item\label{iDELem-CondExpCrossProduct}
$E(I) = I$.
For $T \in \M_+$, we have $E(T) \in \AC_+$.
\item\label{iiDELem-CondExpCrossProduct}
For $T = \sum_{\gamma \in \Gamma} \varphi_\gamma u_\gamma \in \M$,
we have 
$$
E(T^*T) \, = \, \sum_{\gamma \in \Gamma}
\vert \theta_\gamma(\varphi_{\gamma^{-1}}) \vert^2
\hskip.5cm \text{and} \hskip.5cm 
E(TT^*) \, = \, \sum_{\gamma \in \Gamma} \vert \varphi_{\gamma} \vert^2 .
$$
\item\label{iiiDELem-CondExpCrossProduct}
For $T \in \M$, we have $E(T^*T) = 0$ if and only if $T = 0$.
\item\label{ivDELem-CondExpCrossProduct}
$E$ is a normal map between the von Neumann algebras $\M$ and $\AC$.
\item\label{vDELem-CondExpCrossProduct}
for $\varphi, \psi \in \AC$ and $T \in \M$, 
we have $E(\varphi T\psi) = \varphi E(T)\psi$;
in particular, $E$ is the identity map on $\AC$.
\end{enumerate}
\end{lem}

In the following proposition, we relate traces on the crossed product 
algebra $L^\infty(X, \mu) \rtimes \Gamma$ and invariant measures
on $X$. 

\begin{prop}
% 13.A.8
\label{Prop-CondExpTrace} 
Let $(X, \mathcal B, \mu) \curvearrowleft \Gamma$ be
an action of a countable group $\Gamma$
on the standard measure space $(X, \mathcal B, \mu)$. Let 
$\M, \AC$ and $E$ be as above.
\begin{enumerate}[label=(\arabic*)]
\item\label{iDEProp-CondExpTrace}
Assume that there exists 
a $\sigma$-finite $\Gamma$-invariant measure $\nu$ on $(X, \mathcal B)$ 
which is equivalent to $\mu$.
Then 
$$
\tau \, \colon \, \M_+ \to \mathopen[ 0, +\infty \mathclose] , 
\hskip.2cm 
T \mapsto \nu \circ E(T) = \int_X E(T) (x) d\nu(x)
$$
is a faithful semi-finite normal trace on $\M$.
Moreover, $\tau$ is finite if and only if $\nu$ is a finite.
\item\label{iiDEProp-CondExpTrace}
Assume that there exists a faithful semi-finite normal trace
$\tau \,\colon \M_+ \to \mathopen[ 0, +\infty \mathclose]$ on $\M$
and, moreover, that the action of $\Gamma$ is ergodic and essentially free.
Then there exists a $\sigma$-finite $\Gamma$-invariant measure
$\nu$ on $(X, \mathcal B)$ 
which is equivalent to $\mu$ and such that $\tau \vert_{\AC} = \nu$, that is, 
$$
\tau (M(\varphi)) \, = \, \int_X \varphi(x)d\nu(x)
\hskip.5cm \text{for all} \hskip.2cm
\varphi \in \AC_+.
$$
\end{enumerate}
\end{prop}

\begin{proof}
\ref{iDEProp-CondExpTrace}
For $S,T \in \M_+$ and $\lambda \ge 0$, it is straightforward to check that
$$
\tau(S+T) \, = \, \tau(S)+\tau(T), 
\hskip.5cm \text{and} \hskip.2cm
\tau(\lambda S) \, = \, \lambda \tau(S).
$$
For $T = \sum_{\gamma \in \Gamma} \varphi_\gamma u_\gamma$,
we have, using the formulae in Lemma
\ref{Lem-CondExpCrossProduct}~\ref{iiDELem-CondExpCrossProduct}
and the $\Gamma$-invariance of the measure $\nu$,
$$
\begin{aligned}
\tau(T^*T)
& = \int_X \left(\sum_{\gamma \in \Gamma}
\vert \theta_\gamma(\varphi_{\gamma^{-1}}) (x) \vert^2 \right)d\nu(x)
\\
& = \sum_{\gamma \in \Gamma} \int_X 
\vert\theta_\gamma (\varphi_{\gamma^{-1}}) (x) \vert^2 d\nu(x)
\\
& = \sum_{\gamma \in \Gamma} \int_X \vert \varphi_{\gamma^{-1}}(x)\vert^2 d\nu(x)
\\
& = \int_X \left(\sum_{\gamma \in \Gamma}
\vert \varphi_{\gamma^{-1}}(x)\vert^2 \right)d\nu(x)\\
& = \int_X \left(\sum_{\gamma \in \Gamma}
\vert \varphi_{\gamma}(x)\vert^2 \right)d\nu(x)
\\
& = \tau(TT^*)
\end{aligned}
$$
Thus, $\tau(T^*T) = \tau(TT^*)$ for every $T \in \M$.
The above formula shows also that $\tau$ is faithful on $\M_+$
since $\nu$ is faithful on $\AC_+$.
Moreover, $\tau$ is normal since $E$ and $\nu$ are normal.
It remains to show that $\tau$ is semi-finite.
\par

Since $\nu$ is semi-finite,
there exists an increasing sequence $(\varphi_n)_{n \ge 1}$ in $\AC_+$
with $\sup_n \varphi_n = \Un_X$ and $\int_X f_n(x) d\nu(x) < +\infty$
for all $n \ge 1$.
Then $\tau(\varphi_n) = \int_X f_n(x) d\nu(x) < +\infty$ and 
$\sup_n \varphi_n = I$ in $\M$ for all $n \ge 1$.
This implies that $\tau$ is semi-finite;
see for instance
\cite[Chap.~I, \S~6, Prop.\ 1, Cor.\ 2]{Dixm--vN}.

\vskip.2cm

\ref{iiDEProp-CondExpTrace}
We define $\nu \,\colon \mathcal B \to \mathopen[ 0, +\infty \mathclose]$ by 
$$
\nu(B) \, = \, \tau(M(\Un_B))
\hskip.5cm \text{for all} \hskip.2cm
B \in \mathcal B.
$$
Since $\tau$ is normal, $\nu$ is $\sigma$-additive.
Also, since $\tau$ is faithful, $\nu$ is equivalent to $\mu$. 
Moreover, $\nu$ is $\Gamma$-invariant; indeed, for 
$B \in \mathcal B$ and $\gamma \in \Gamma$, we have
$$
\begin{aligned}
\nu(\gamma B)
\, &= \, \tau(M(\Un_{\gamma B}))
\, = \, \tau (M(\theta_\gamma(\Un_B)))
\\
\, &= \, \tau (U_\gamma M(\Un_B) U_\gamma^{-1})
\, = \, \tau (M(\Un_B)) 
\, = \, \nu(B).
\end{aligned}
$$
It remains to show that $\nu$ is $\sigma$-finite. 
To prove this, it suffices to show the following 

\vskip.2cm

\noindent
\textbf{Claim:} there exists $B \in \mathcal B$ with $0 < \nu(B) < +\infty$.
\vskip.2cm

Indeed, assume that the claim is proved. 
Since the action of $\Gamma$ on $X$ is $\mu$-ergodic 
and since $\nu$ is equivalent to $\mu$, 
this action is $\nu$-ergodic.
Let $X' := \bigcup_{\gamma} \gamma B$;
then $X' \in \mathcal B$, $X$ is $\Gamma$-invariant, and $\nu(X') > 0$.
Therefore $\nu(X \smallsetminus X') = 0$.
Moreover, $\nu(\gamma B) = \nu(B) < +\infty$
for every $\gamma \in \Gamma$.
Since $\Gamma$ is countable, this shows that $\nu$ is $\sigma$-finite.

\vskip.2cm

To prove the claim, let $T \in \M_+$ with $0 < \tau(T) < +\infty$.
The first step is to show that $\tau(E(T)) < +\infty$. 
\par

Indeed, denote by $K(T)$ the convex hull of 
$\{UTU^* \mid U \in \U(\AC)\}$,
where $\U(\AC)$ is the unitary group of $\AC$.
Let $K'(T)$ be the closure of $K(T)$
for the ultra-weak topology on $\Li (\Hi)$,
where $\Hi = \ell^2(\Gamma, L^2(X, \mu))$. 
Since $K(T)$ is norm bounded, $K'(T)$ is a compact convex
subset for the ultra-weak topology.
(On the ultra-weak topology, see Appendix \ref{AppHspacesop}).
The group $\U(\AC)$ acts on $K'(T)$ by the affine
ultra-weakly continuous transformations given by 
$$
K'(T) \ni S \, \mapsto \, USU^*\in K'(T).
$$
for $U \in \U(\AC)$.
Since $\U(\AC)$ is abelian,
it follows from the Markov-Kakutani theorem that there exists
a common fixed point $T_0 \in K'(T)$ for all $U \in \U(\AC)$.
(For this theorem, see for example \cite[Chap.\ IV, Appendice]{BEVT1--5}.)
Thus, we have $UT_0 = T_0U$ for every $U \in \U(\AC)$. 
Since $\AC$ coincides with the linear span of $\U(\AC)$,
it follows that $S T_0 =T_0 S$ for every $S \in \AC$.
Since the $\Gamma$-action is essentially free,
$\AC$ is a maximal abelian subalgebra of $\M$
(Lemma~\ref{Lem-FreeActionMaxAbelian}) and it follows that $T_0 \in \AC$.
\par

We have 
$$
E(UTU^*) \, = \, UE(T)U^* \, = \, E(T)
\hskip.5cm \text{for all} \hskip.2cm
U \in \U(\AC).
$$
Since $T_0\in K'(T)$ and since $E \,\colon \M \to \AC$ is normal
and therefore continuous for the ultra-weak topology,
it follows that $T_0 = E(T)$.
\par

We have $\tau(S) = \tau (T)$ for every $S \in K(T)$, by the trace property of $\tau$.
Moreover, since $\tau$ is normal, $\tau$ is lower continuous
for the weak topology (as mentioned before Example~\ref{ExampleTracevN})
and hence for the ultra-weak topology.
It follows that $\tau(S) \le \tau (T)$ for every $S \in K'(T)$; in particular,
$$
\tau(T_0) \, = \, \tau(E(T)) \, \le \, \tau (T) \, < \, +\infty.
$$
\par

Let now $\varphi$ be the function in $\AC=L^\infty(X, \mu)$
such that $M(\varphi) =T_0 = E(T)$. 
By Items \ref{iDELem-CondExpCrossProduct} 
and \ref{iiiDELem-CondExpCrossProduct}
of Lemma~\ref{Lem-CondExpCrossProduct}, 
we have $\varphi \ne 0$ and $\varphi \ge 0$, 
since $T \ne 0$ and $T \in \M_+$. 
So, $\lambda := \Vert \varphi \Vert_\infty > 0$.
Therefore there exists $B \in \mathcal B$ with $0 < \mu(B) < +\infty$
with $\varphi(x) \ge \lambda/2$ for every $x \in B$.
\par

Since $\nu$ is equivalent to $\nu$, we have $\nu(B) > 0$. 
Moreover, we have $\Un_B \le \dfrac{2}{\lambda} \varphi$ and hence 
$$
\nu(B) \, = \,
\tau(M(\Un_B)) \, \le \,
\frac{2}{\lambda}\tau(M(\varphi)) \, = \,
\frac{2}{\lambda}\tau(T_0) \, < \,
\infty. 
$$
The proof of the Claim is now complete.
\end{proof}

\begin{rem}
% 13.A.9
\label{Rem-Prop-CondExpTrace} 
Under the assumptions of 
Proposition~\ref{Prop-CondExpTrace}~\ref{iiDEProp-CondExpTrace},
a stronger result is true: 
if $\tau$ is a faithful semi-finite normal trace on
$\M = L^\infty(X, \mu) \rtimes \Gamma$,
then $\tau$ is of the form $\nu \circ E$
for a $\sigma$-finite $\Gamma$-invariant measure $\nu$ on $(X, \mathcal B)$.
\par

Indeed, as we will see below (Theorem~\ref{ROIV}),
the crossed product algebra $\M$ is a factor.
On the one hand,
% by Proposition~\ref{Prop-CondExpTrace}~\ref{iiDEProp-CondExpTrace},
$\tau \vert_{\AC}$ is given by a 
$\sigma$-finite $\Gamma$-invariant measure $\nu$ on $X$.
On the other hand, 
by Proposition~\ref{Prop-CondExpTrace}~\ref{iDEProp-CondExpTrace},
$\nu$ defines a faithful semi-finite normal trace $\tau'$ on $\M$
such that $\tau' = \nu' \circ E$.
However, since $\M$ is a factor,
$\tau$ and $\tau'$ have to be proportional to each other 
(Corollaire to Th\'eor\`eme 3 in \cite[Chap.~I, \S~6]{Dixm--vN});
hence, we have $\tau = \tau'$, 
since $\tau$ and $\tau'$ coincide on $\AC$.
\end{rem}

The following result is due to Murray and von Neumann;
see Lemmas 12.3.4, 13.1.1, and 13.1.2 in \cite{MuvN--36}
and Theorem IX in \cite{vNeu--40}. 

\begin{theorem}
% 13.A.10
\label{ROIV}
Let $\Gamma$ be a countable group acting 
on a standard measure space $(X, \mathcal B, \mu)$, as above. 
Denote by $\mathcal M$ the corresponding crossed product
$L^\infty(X, \mu) \rtimes \Gamma$.
Assume that the action of $\Gamma$ is essentially free and ergodic. 
\begin{enumerate}[label=(\arabic*)]
\item\label{iDEROIV}
% i
The von Neumann algebra $\mathcal M$ is a factor.
\item\label{IDEROIV}
% ii
Assume that the action of $\Gamma$ on $(X, \mathcal B)$ 
is essentially transitive; then 
$\mathcal M$ is a factor of type I.
\item\label{II1DEROIV}
% iii
Assume that there exists 
a $\Gamma$-invariant non-atomic probability measure on $(X, \mathcal B)$ 
which is equivalent to $\mu$; then $\mathcal M$ is a factor of type II$_1$.
\item\label{IIDEROIV}
% iv
Assume that there exists an infinite $\sigma$-finite and $\Gamma$-invariant
non-atomic measure on $(X, \mathcal B)$ which is equivalent to $\mu$; 
then $\mathcal M$ is a factor of type II$_\infty$.
\item\label{IIIDEROIV}
% v
Assume that there does not exist any
$\Gamma$-invariant $\sigma$-finite measure on $X$ which is equivalent to $\mu$; then
$\mathcal M$ is a factor of type III.
\end{enumerate}
\end{theorem}

\begin{proof}
\ref{iDEROIV}
Let $T \in \M \cap \M'$. In particular, $T \in M(\AC)'$. 
Since the action of $\Gamma$ is essentially free, 
$M(\AC)$ is a maximal abelian von Neumann subalgebra of $\M$
(Lemma~\ref{Lem-FreeActionMaxAbelian})
and hence $T = M(\varphi)$ for some $\varphi \in \AC$.
As $T$ commutes with $U_\gamma$, we have then, by relation $(***)$ above,
$\theta_\gamma(\varphi) = \varphi$ for every $\gamma \in \Gamma$.
By ergodicity of the $\Gamma$-action,
it follows that $\varphi$ is constant almost everywhere,
that is, $T$ is a multiple of the identity operator.
This shows that $\M$ is a factor.

\vskip.2cm

\ref{IDEROIV}
Upon removing a set of $\mu$-measure $0$, we may
assume that $X$ is a $\Gamma$-orbit.
Since $\Gamma$ is countable,
$X$ is countable and $\mu(\{x\}) > 0$ for every $x \in X$.
Let $x \in X$. 
Then $p := M(1_{\{x}\})$ is a non-zero minimal projection
in the von Neumann algebra $M(\AC) \approx \AC$.
We claim that $p$ is a non-zero minimal projection in $\M$.
Once proved, this will show that $\M$ is a factor of type I.
\par

Let $q \in \M$ be a non-zero projection with $q \le p$,
that is, such that $pq = qp = q$.
Then, for every $\varphi \in \AC$, we have
$$
q M(\varphi) \, = \, qp M(\varphi) \, = \, q M(1_{\{x\}} \varphi) \, = \, 
\varphi(x) q M(1_{\{x\}}) \, = \, \varphi(x) qp \, = \, \varphi(x) q
$$
and 
$$
M(\varphi) q \, = \, M( \varphi)pq \, = \, M(\varphi 1_{\{x\}}) q \, = \,
\varphi(x) M(1_{\{x\}}) q \, = \, \varphi(x) pq \, = \, \varphi(x) q ,
$$
showing that $q$ commutes with $M(\AC)$.
Since $M(\AC)$ is a maximal abelian subalgebra of $\M$
(Lemma~\ref{Lem-FreeActionMaxAbelian}), it follows that $q \in M(\AC)$ 
and hence $q = p$. So, $p$ is a minimal projection in $\M$.

\vskip.2cm

\ref{II1DEROIV} and \ref{IIDEROIV}
Let $\nu$ be a $\Gamma$-invariant $\sigma$-finite
non-atomic measure on $(X, \mathcal B)$ 
which is equivalent to $\mu$. It follows from \ref{iDEROIV} and
Proposition~\ref{Prop-CondExpTrace}~\ref{iDEProp-CondExpTrace}
that $\M$ is a semi-finite factor, with trace $\tau$ such that 
$$
\tau (M(\varphi)) \, = \, \int_X \varphi(x)d\nu(x)
\hskip.5cm \text{for all} \hskip.2cm
\varphi \in \AC_+.
$$
\par

Assume that $\mu$ is a probability measure.
Since $\mu$ is non-atomic, for every $\alpha \in\mathopen[ 0,1 \mathclose]$, 
there exists $B \in \mathcal B$ such that 
$\nu(B) = \alpha$ and hence $\tau(M(\mathbf{1}_A)) = \alpha$. 
This shows that $\tau$ takes on projections in $\M$
every value in $\mathopen[ 0,1 \mathclose]$.
Therefore, $\M$ is a factor of type II$_1$. 
\par

Assume that $\mu$ is an infinite measure.
Then, similarly, $\tau$ takes on projections in $\M$
every value in $\mathopen[ 0, +\infty \mathclose]$.
Therefore, $\M$ is a factor of type II$_\infty$. 

\vskip.2cm

\ref{IIIDEROIV} follows from
Proposition~\ref{Prop-CondExpTrace}~\ref{iiDEProp-CondExpTrace}.
\end{proof}

Here is an addendum to Theorem~\ref{ROIV},
concerning the assumptions of essential freeness and ergodicity
of the action $\Gamma \curvearrowright (X, \mu)$.
 
\begin{prop}
% 13.A.11
\label{Prop-NecCondROIV}
Let $\Gamma$ be a countable group acting 
on a standard measure space $(X, \mathcal B, \mu)$.
Let $\M = L^\infty(X, \mu) \rtimes \Gamma$
be the associated crossed product von Neumann algebra.
\begin{enumerate}[label=(\arabic*)]
\item\label{iDEProp-NecCondROIV}
The algebra $M(\AC) \approx \AC$
is a maximal abelian von Neumann subalgebra of $\M$
if and only the action $(X, \mu) \curvearrowleft \Gamma$ is essentially free. 
\item\label{iiDEProp-NecCondROIV}
If $\M$ is a factor, then $(X, \mu) \curvearrowleft \Gamma$is ergodic. 
\end{enumerate}
\end{prop}

\begin{proof}
\ref{iDEProp-NecCondROIV}
We have already proved (Lemma~\ref{Lem-FreeActionMaxAbelian}) that,
if $(X, \mu) \curvearrowleft \Gamma$ is essentially free,
then $M(\AC)$ is a maximal abelian von Neumann subalgebra of $\M$.
\par

Assume that $(X, \mu) \curvearrowleft \Gamma$ is not essentially free.
Then, by Proposition~\ref{prop-EssentialFree},
there exists $\varphi \in L^\infty(X, \mu)$ with $\varphi \ne 0$
and $\gamma_0 \in \Gamma \smallsetminus \{e\}$
such that $\varphi \theta_{\gamma_0}(\psi) = \psi\varphi$
for every $\psi \in L^\infty(X, \mu)$. 
\par

Let $T := \varphi u_{\gamma_0} \in \M$. We claim that $T \in M(\AC)'$. 
Indeed, let $S = \psi \in \AC$.
Then 
$$
\begin{aligned}
TS 
\, &= \, \varphi u_{\gamma_0}\psi
\, = \, \varphi (u_{\gamma_0}\psi u_{\gamma_0^{-1}})u_{\gamma_0}
\\
\, &= \, \varphi \theta_{\gamma_0} (\psi) u_{\gamma_0}
\, = \, \psi \varphi u_{\gamma_0}
\, = \, ST.
\end{aligned}
$$
Since $T \notin M(\AC)$, it follows that $M(\AC)$ is not maximal abelian.

\vskip.2cm

\ref{iiDEProp-NecCondROIV}
Assume that $(X, \mu) \curvearrowleft \Gamma$ is not ergodic. 
Then, there exists (Proposition~\ref{prop-ergodicity})
a non-constant $\varphi \in L^\infty(X, \mu)$
such that $\theta_\gamma(\varphi) = \varphi$ for all $\gamma \in \Gamma$.
Let $T := \varphi \in M(\AC)$. We claim that $T \in \M'$. 
Indeed, let
$S = \sum_{\gamma \in \Gamma} \psi_{\gamma} u_{\gamma} \in \M$.
Then 
$$
\begin{aligned}
TS
\, &= \, \sum_{\gamma \in \Gamma} \varphi \psi_{\gamma} u_{\gamma}
\, = \, \sum_{\gamma \in \Gamma} \psi_{\gamma} \varphi u_{\gamma}
\, = \, \sum_{\gamma \in \Gamma} \psi_{\gamma} u_{\gamma}
(u_{\gamma^{-1}} \varphi u_{\gamma})
\\
\, &= \, \sum_{\gamma \in \Gamma} \psi_{\gamma}u_{\gamma}
 \theta_{\gamma^{-1}}(\varphi)
\, = \, \sum_{\gamma \in \Gamma} \psi_{\gamma}u_{\gamma} \varphi
\, = \, ST.
\end{aligned}
$$
So, $T$ belongs to the centre of $\M$ and is not a scalar multiple of the identity.
Therefore $\M$ is not a factor.
\end{proof}

\begin{rem}
% 13.A.12
\label{Rem-NecCondROIV}
Unlike the ergodicity assumption, the essential freeness assumption
for the action $(X, \mu) \curvearrowleft \Gamma$ is not necessary
for $L^\infty(X, \mu) \rtimes \Gamma$ to be a factor.
\par

An easy example is given by the action of an icc countable group $\Gamma$
on a one-point space $X = \{\ast\}$.
In this case, $L^\infty(X, \mu) \rtimes \Gamma$
coincides with the von Neumann algebra $\Li (\Gamma)$
of the regular representation of $\Gamma$
and is therefore a factor (Proposition~\ref{iccfactorII1}).
\end{rem}

Using Theorem~\ref{ROIV},
we give various examples of group actions $(X, \mu) \curvearrowleft \Gamma$ 
on measure spaces $(X, \mu)$,
for which the crossed product von Neumann algebra $L^\infty(X, \mu) \rtimes \Gamma$
is a factor of type II$_1$, II$_\infty$, or III.
For other examples of the use of Theorem~\ref{ROIV},
in particular in connection with the regular representation of a group, 
see Chapter~\ref{Chap:NormalInfiniteRep}.

\begin{exe}[\textbf{type II$_1$}]
% 13.A.13
\label{Examples-ROIV-12}
Let $G$ be an infinite metrizable compact group,
with normalized Haar measure $\mu_G$. 
Let $\Gamma$ be a countable dense subgroup of $G$.
The action of $\Gamma$ on $(G, \mu_G)$ by right translations
is measure preserving, free, and ergodic.
Since $G$ is infinite, $\mu_G$ is non-atomic.
Therefore by Theorem~\ref{ROIV}, 
$L^\infty(G, \mu_G) \rtimes \Gamma$ is a factor of type II$_1$.
\par
 
Examples of pairs $(G, \Gamma)$ as above are for instance:
the pair $(\T, \Gamma_\theta)$ with
$\Gamma_\theta = \{e^{2 \pi i n \theta} \mid n \in \Z\}$
for an irrational number $\theta$;
the pair $(\SO(3), \Gamma)$, where $\Gamma$ is one of the subgroups of $\SO(3)$ 
isomorphic to $(\Z / 2 \Z) \ast (\Z/3\Z)$ mentioned
in Example~\ref{Exa-TotallyDisconnectedBohr}(3).
\end{exe}

\begin{exe}[\textbf{type II$_1$}]
% 13.A.14
\label{Examples-ROIV-13}
Let $G$ be a simple non-compact connected real Lie group with trivial centre
and let $\Lambda$ be a lattice in $G$.
Let $\mu$ be the unique $G$-invariant probability measure
on the Borel subsets of $\Lambda\backslash G$.
Let $\Gamma$ be a countable subgroup of $G$
which is not relatively compact in $G$. 
The action of $\Gamma$ on $(\Lambda\backslash G, \mu)$ by right translations
is measure preserving.
Moreover, this action is essentially free (see \cite[Lemma 4.1]{Boul--16})
and ergodic by Moore's ergodicity theorem (see \cite[Theorem 2.2.6]{Zimm--84}).
Therefore, by Theorem~\ref{ROIV}, 
$L^\infty(\Lambda\backslash G, \mu) \rtimes \Gamma$ is a factor of type II$_1$.
\par
 
An example of a triple $(G, \Lambda, \Gamma)$ as above is
for instance $(PSL_n(\R), PSL_n(\Z), \Gamma)$ for $n \ge 2$,
where $\Gamma$ is any countable subgroup of $\PSL_n(\R)$
which is not relatively compact. 
\end{exe}
 
\begin{exe}[\textbf{type II$_\infty$}]
% 13.A.15
\label{Examples-ROIV-14}
Let $G$ be a non-compact and non-discrete second-countable LC group,
with right Haar measure $\mu_G$. 
Let $\Gamma$ be a countable dense subgroup of $G$.
The action of $\Gamma$ on $(G, \mu_G)$ by right translations
is measure preserving, free, and ergodic.
Since $G$ is infinite, $\mu_G$ is non-atomic.
Therefore by Theorem~\ref{ROIV}, 
$L^\infty(G, \mu_G) \rtimes \Gamma$ is a factor of type II$_\infty$.
\par

Examples of pairs $(G, \Gamma)$ as above are for instance the pair $(\R, \Q)$ and
the pair $(SL_n(\R), \SL_n(\Q))$.
\end{exe}
 
\begin{exe}[\textbf{type III}]
% 13.A.16
\label{Examples-ROIV-15}
The group $G = PSL_2(\R)$ acts (on the right) on the real projective line $\mathbf{P}^1(\R)$
by M\"obius transformations.
The Lebesgue measure $\mu$ on $\mathbf{P}^1(\R) \cong \R \cup \{\infty \}$
is quasi-invariant (but not invariant) by the action of $G$.
Let $\Gamma$ be a lattice in $G$ (as, for instance, $\Gamma = PSL_2(\Z)$).

\begin{itemize}
\setlength\itemsep{0em}
\item[$\bullet$]
The action $\Gamma \curvearrowright (\mathbf{P}^1(\R), \mu)$
is essentially free.
Indeed, every $g \in G, g \ne e$,
has at most two fixed points in $\mathbf{P}^1(\R)$.
\par

\item[$\bullet$]
The action $\Gamma \curvearrowright (\mathbf{P}^1(\R), \mu)$ is ergodic.
Indeed, observe that $G$ acts transitively on $\mathbf{P}^1(\R)$
and that the stabilizer of the point $\infty$
is the image $P$ in $G = PSL_2(\R)$ of the group
$$
\left\{ \begin{pmatrix} a & 0 \\ b& c \end{pmatrix} 
\hskip.1cm \Big\vert \hskip.1cm
a,b,c \in \R, \hskip.1cm ac = 1 \right\}
$$
of lower triangular matrices.
So, $\mathbf{P}^1(\R)$ can be identified, as a $G$-space, with $P\backslash G$.
By Moore's duality theorem (see \cite[Corollary 2.2.3]{Zimm--84}),
it suffices to show that the measure preserving action
of $P$ on the probability space $(G/\Gamma, \nu)$ is ergodic,
where $\nu$ is the unique $G$-invariant probability measure
on the Borel subset of $G/\Gamma$.
Since $P$ is not relatively compact in $G$,
this is the case by Moore ergodicity theorem again
(see the reference in \ref{Examples-ROIV-13} above).
\par
 
\item[$\bullet$]
There exists no $\Gamma$-invariant $\sigma$-finite measure
on $\mathbf{P}^1(\R)$ which is equivalent to $\mu$.
This is shown below, see Corollary
\ref{Cor-NoInvMeasureProjectiveSpace}.
\end{itemize}

\noindent
Therefore, by Theorem~\ref{ROIV},
$L^\infty(\mathbf{P}^1(\R), \mu) \rtimes \Gamma$ is a factor of type III.
\end{exe}

\section
{Ergodic group actions without invariant measure}
% Section 13.B
\label{sectionNoinvariantmeasure}

Theorem \ref{ROIV}~\ref{IIIDEROIV} and Example \ref{Examples-ROIV-15}
are motivations for the study of actions of groups on standard measure spaces
without invariant measures.

\vskip.2cm

Let $\Gamma$ be a countable group, 
$(X, \mathcal B, \mu)$ a standard measure space,
and $(X, \mu) \curvearrowleft \Gamma$ an action of $\Gamma$
on $X$ by measurable transformations preserving the class of $\mu$.
Recall from Section \ref{Section-MoreIrredRep}
the Radon--Nikodym cocycle $c_\mu \,\colon X \times \Gamma \to \R^\times_+$,
given by
$$
c_\mu(x, \gamma) \, = \, \frac{d\gamma_*(\mu)}{d\mu}(x)
\hskip.5cm \text{for all} \hskip.2cm
x \in X, \hskip.1cm \gamma \in \Gamma .
$$
\index{Radon--Nikodym! cocycle}
\index{Cocycle! $1$@Radon--Nikodym}
\par

Assume that the action $(X, \mu) \curvearrowleft \Gamma$ is ergodic
and $\mu$ non-atomic. 
Following \cite[\S~3]{Schm--77a} and \cite[\S~8]{FeMo--77},
we define the essential range of $c_\mu$.
\par

\index{$l 8$@$\overline{\R^\times_+} = \mathopen] 0, \infty \mathclose]
= \R^\times_+ \cup \{\infty \}$ extended positive real numbers}
We denote by $\overline{\R^\times_+} = \R^\times_+ \cup \{\infty \}$
the one-point compactification of $\R^\times_+$.
A number $\lambda\in \R^\times_+$ is an essential value of $c_\mu$ if,
for every $\varepsilon > 0$ and 
every measurable subset $A$ of $X$ with $\mu(A) > 0$,
there exists a measurable subset $B \subset A$ with $\mu(B) > 0$
and $\gamma \in \Gamma$ such that $B \gamma \subset A$ and
$$
\vert c_\mu(x, \gamma)- \lambda \vert \, \le \, \varepsilon
\hskip.5cm \text{for all} \hskip.2cm
x \in B;
$$
the point $\infty$ is an essential value of $c_\mu$ if,
for every $n \in \N^*$ and every measurable subset $A$ of $X$ with $\mu(A) > 0$,
there exists a measurable subset $B \subset A$ with $\mu(B) > 0$
and $\gamma \in \Gamma$ such that $B \gamma \subset A$ and
$$
\vert c_\mu(x, \gamma) \vert \, \notin \, \mathopen[ 1/n, n \mathclose]
\hskip.5cm \text{for all} \hskip.2cm
x \in B.
$$
The set of essential values of $c_\mu$ in $\overline{\R^\times_+}$
is called the \textbf{essential range}
of the action $(X, \mu) \curvearrowleft \Gamma$ or of the cocycle $c_\mu$,
and will be denoted by $r(X, \Gamma, \mu)$;
it is also called the asymptotic range or the ratio range of Krieger--Araki--Woods
of this action.
\index{Essential range! of a cocycle}
\par

It is easy to check that the essential range $r(X, \Gamma, \mu)$
only depends on the equivalence class of $\mu$.
 and on the orbit equivalence class of $(X, \mu) \curvearrowleft \Gamma$
Moreover, $r(X, \Gamma, \mu)$ is obviously
a closed subset of $\overline{\R^\times_+}$
and $r(X, \Gamma, \mu) \cap \R^\times_+$ 
is a closed subgroup of the multiplicative group $\R^\times_+$
(see Lemma~3.3 in \cite{Schm--77a}).
Therefore the possible values of the asymptotic range $r(X, \Gamma, \mu)$ are
\begin{itemize}
\setlength\itemsep{0em}
\item[(1)]
$r(X, \Gamma, \mu) = \{1 \}$;
\item[(2)]
$r(X, \Gamma, \mu) = \{1, \infty \}$;
\item[(3)]
$r(X, \Gamma, \mu) = \overline{\R^\times_+}$;
\item[(4)]
$r(X, \Gamma, \mu) = \lambda^\Z\cup\{ \infty \}$
for some $\lambda \in \R^\times_+$ with $0 < \lambda < 1$.
\end{itemize} 
The action $(X, \mu) \curvearrowleft \Gamma$ is said to be
of \textbf{type II} in case (1),
of \textbf{type III$_0$} in case (2),
of \textbf{type III$_1$} in case (3),
and of \textbf{type III$_\lambda$} in case (4).
\index{Type II, III! $3$@action}
\index{Action! type II, III$_\lambda$}

\begin{rem}
% 13.B.1
\label{Rem-EssentialRange}
Actually, given any Borel cocycle $c \,\colon X \times \Gamma \to G$
with values in a second-countable locally compact group $G$,
one can define the essential range of $c$ in exactly the same way
(see \cite[Definition 3.1]{Schm--77a}).
\end{rem}

Next, we show that if $(X, \mu) \curvearrowleft \Gamma$ is of type III,
then there exists no $\sigma$-finite $\Gamma$-invariant measure on $X$
which is equivalent to $\mu$.

\begin{prop}
% 13.B.2
\label{Pro-TypeIIINoInvMeasure}
Assume that there exists a $\Gamma$-invariant $\sigma$-finite measure on $X$
which is equivalent to $\mu$.
Then $(X, \mu) \curvearrowleft \Gamma$ is of type II.
\end{prop}

\begin{proof}
We claim that $r(X, \Gamma, \mu) = \{1 \}$.
Indeed, it follows from the assumption that the cocycle $c_\mu$ is a coboundary,
that is, we have
$$
c_\mu(x, \gamma) \, = \, \frac{\varphi(x)}{\varphi(x\gamma)}
\hskip.5cm \text{for all} \hskip.2cm
x \in X, \hskip.1cm \gamma \in \Gamma
$$
for some measurable function $\varphi \,\colon X \to \R^\times_+$
(see Example~\ref{examplecocycle}).
\par

Since $X$ is a non-atomic standard Borel space,
we can assume that $X = \mathopen[ 0,1 \mathclose]$,
equipped with the Lebesgue measure on its Borel subsets.
\par

Fix $\varepsilon > 0$.
Since $\varphi$ is measurable, there exists by Lusin's theorem (see \cite[2.23]{Rudi--66})
a continuous function $\psi \,\colon X \to \R^\times_+$ such that, for the measurable subset 
$$
X_0 \, := \, \{x \in X\mid \varphi(x) \, = \, \psi(x) \},
$$
we have $\mu(X_0) > 0$.
By uniform continuity of $\psi$, there exists $N > 0$
such that, for every $i \in \{0, \hdots, N - 1 \}$, 
$$
\left\vert \frac{\psi(x)}{\psi(y)} - 1 \right\vert \, \le \, \varepsilon
\hskip.5cm \text{for all} \hskip.5cm
x, y \in \left[ \frac{i}{N}, \frac{i + 1}{N} \right].
$$
Since $\mu(X_0) > 0$, we have $\mu(X_0 \cap \left[ \frac{i}{N}, \frac{i+1}{N} \right]) > 0$
for some $i \in \{0, \hdots, N - 1 \}$.
Set 
$$
A \, := \, X_0 \cap \left[ \frac{i}{N}, \frac{i+1}{N} \right]
%\{x \in X \mid \left\vert \frac{\varphi(x)}{\varphi(\gamma x)} - 1 \right\vert \le \varepsilon \}.
$$
and observe that 
$$
\left \vert \frac{\varphi(x)}{\varphi(y)} - 1 \right\vert \, \le \, \varepsilon
\hskip.5cm \text{for all} \hskip.2cm
x, y \in A.
$$
\par

Let $B$ be measurable subset with $B \subset A$
and let $\gamma \in \Gamma$ be such that $ B \gamma \subset A$. 
Let $x \in B$.
Since $x \in A$ and $x \gamma \in A$, we have 
$$
\vert c_\mu(x, \gamma) - 1 \vert \, = \,
\left\vert \frac{\varphi(x)}{\varphi(x \gamma)} - 1 \right\vert \, \le \, \varepsilon
\hskip.5cm \text{for all} \hskip.2cm
x \in B.
$$
As $\varepsilon$ was arbitrary,
it follows from the definition of an essential value of $c_\mu$ that
$r(X, \Gamma, \mu) =1$.
\end{proof}

\begin{rem}
% 13.B.3
\label{Rem-TypeIIINoInvMeasure}
(1)
The converse statement in Proposition~\ref{Pro-TypeIIINoInvMeasure} is also true:
if $r(X, \Gamma, \mu) = \{1 \}$ (that is, if $(X, \mu) \curvearrowleft \Gamma$ is of type II),
then there exists a $\Gamma$-invariant $\sigma$-finite measure on $X$
which is equivalent to $\mu$; see Theorem~3.9 in \cite{Schm--77a}.

\vskip.2cm

(2)
There is a subdivision of von Neumann factors of type III into further subclasses:
types III$_0$, III$_1$, and III$_\lambda$ for $\lambda \in (0,1)$;
see Chapter 3 in \cite{Sund--87}.
It is known that, for a group action $(X, \mu) \curvearrowleft \Gamma$ 
as in Theorem~\ref{ROIV}, the crossed product
$L^\infty(X, \mu) \rtimes \Gamma$ is of type III$_\lambda$ if and only if 
$(X, \mu) \curvearrowleft \Gamma$ is of type III$_\lambda$
for $\lambda \in \mathopen[ 0,1 \mathclose]$
(see Proposition 4.3.18 in \cite{Sund--87}).
\index{Type II, III! $4$@III$_\lambda$}
\end{rem}

Let $(X, \mu) \curvearrowleft \Gamma$ be an ergodic action as above.
Consider the product space
$$
\widetilde{X} \, := \, X \times \R
$$
equipped with the product $\sigma$-algebra and with the measure
$$
d\widetilde{\mu}(x, t) \, = \, d\mu(x) e^t dt,
$$
where $dt$ is the Lebesgue measure on $\R$.
We are going to associate to the Radon--Nikodym cocycle $c_\mu$
a \emph{measure preserving} action of 
$\Gamma$ on $(\widetilde{X}, \widetilde{\mu})$.
\par

For $\gamma \in \Gamma$ and $(x, t) \in \widetilde{X}$, define
$$
(x, t) \gamma \, = \, (x \gamma, t - \log c_\mu(x, \gamma)).
$$
Using the cocycle identity, we check that this defines indeed
an action $\widetilde{X} \curvearrowleft \Gamma$:
for $\gamma, \delta \in \Gamma$ and $(x, t) \in \widetilde{X}$, we have
$$
\begin{aligned}
((x, t) \gamma)\delta
\, &= \, (x\gamma, t - \log c_\mu(x, \gamma)) \delta
\\
\, &= \, ((x\gamma)\delta, t - \log c_\mu(x, \gamma) -\log c_\mu(x \gamma, \delta))
\\
\, &= \, (x(\gamma\delta), t - \log c_\mu(x, \gamma \delta))
\\
\, &= \, (x, t) \gamma\delta.
\end{aligned}
$$
The action $\widetilde{X} \curvearrowleft \Gamma$ preserves $\widetilde{\mu}$.
for $\gamma \in \Gamma$ and a measurable function $f \,\colon \widetilde{X} \to \R_+$,
we have 
$$
\begin{aligned}
\int_{\widetilde{X}} f((x, t)\gamma) d\widetilde{\mu}(x, t)
\, &= \, \int_{X} \int_\R f(x\gamma, t - \log c_\mu(x, \gamma)) e^t d\mu(x) dt
\\
\, &= \, \int_{X} \int_\R f(x, t - \log c_\mu(x\gamma^{-1}, \gamma)) c(x, \gamma^{-1}) e^t
d\mu(x) dt
\\
\, &= \, \int_{X} \int_\R f(x, t) c(x, \gamma^{-1})c_\mu(x\gamma^{-1}, \gamma)) e^t d\mu(x) dt
\\
\, &= \, \int_{X} \int_\R f(x, t) e^t d\mu(x) dt
\\
\, &= \, \int_{\widetilde{X}} f(x, t) d\widetilde{\mu}(x, t).
\end{aligned}
$$
The action $(\widetilde{X}, \widetilde{\mu}) \curvearrowleft \Gamma$
is called the \textbf{Maharam extension} or Maharam skew product extension
of the action $(X, \mu) \curvearrowleft \Gamma$
(see \cite[\S~3.4]{Aaro--97} and \cite[\S~5]{Schm--77a}).
\index{Maharam extension}
\par

The type of the action $(X, \mu) \curvearrowleft \Gamma$ can be read off 
properties of the Maharam extension (see Theorem~5.2 in \cite{Schm--77a}).
We will only need the following very special result.

\begin{prop}
% 13.B.4
\label{Pro-MaharamErgodic} 
Assume that the Maharam extension
$(\widetilde{X}, \widetilde{\mu}) \curvearrowleft \Gamma$ is ergodic.
Then the action $(X, \mu) \curvearrowleft \Gamma$ is of type III$_1$.
\end{prop}

\begin{proof}
Since $r(X, \Gamma, \mu)$ is closed in $\overline{\R^\times_+}$, it suffices to show
that $\R^\times_+ \subset r(X, \Gamma, \mu)$. 
\par

Let $\lambda \in \R^\times_+$ and set $a = \log \lambda$.
Let $\varepsilon > 0$ and let $A$ be a measurable subset of $X$ with $\mu(A) > 0$. 
Then
$$
\widetilde{A} \, := \, A \times (a - \varepsilon/2, a + \varepsilon/2)
$$
is a measurable subset of $\widetilde{X}$ with $\widetilde{\mu}(\widetilde{A}) > 0$.
The set 
$$
\widetilde{A_0} \, := \, \bigcup_{\gamma\in \Gamma}
(A \times (-\varepsilon/2, \varepsilon/2))\gamma
$$
is a $\Gamma$-invariant measurable subset of $\widetilde{X}$
with $\widetilde{\mu}(\widetilde{A_0}) > 0$.
It follows from the invariance of $\widetilde{\mu}$
and the ergodicity of $(\widetilde{X}, \widetilde{\mu}) \curvearrowleft \Gamma$
that $\widetilde{\mu}(\widetilde{X} \smallsetminus \widetilde{A_0}) = 0$.
Since $\Gamma$ is countable and $\widetilde{\mu}(\widetilde{A_0}) > 0$,
there exists an element $\gamma \in \Gamma$ such that
$$
\widetilde{\mu}(\widetilde{A} \cap
(A\times (-\varepsilon/2, \varepsilon/2))\gamma^{-1})
\, > \, 0.
$$
Choose $t \in (a-\varepsilon/2, a+\varepsilon/2)$ such that, for the $t$-section
$$
B \, := \, \left\{x \in A \mid (x, t) \in
\widetilde{A} \cap (A \times (-\varepsilon/2, \varepsilon/2))\gamma^{-1} \right\}
$$
of $\widetilde{A} \cap (A \times (-\varepsilon/2, \varepsilon/2))\gamma^{-1}$,
we have $\mu(B) > 0$.
\par 

Let $x \in B$.
There exists $(y, s) \in A\times (-\varepsilon/2, \varepsilon/2)$ such that 
$$
x \, = \, y \gamma^{-1}
\hskip.5cm \text{and} \hskip.5cm
t \, = \, s - \log c_\mu(y, \gamma^{-1}).
$$
Since 
$$
\log c_\mu(y, \gamma^{-1}) \, = \,
\log c_\mu(x \gamma, \gamma^{-1}) \, = \,
\log c_\mu(x, \gamma)^{-1} \,= \,
-\log c_\mu(x, \gamma),
$$
we have 
$$
\vert \log c_\mu(x, \gamma) - a \vert \, \le \,
\vert \log c_\mu(x, \gamma) - t\vert + \vert t - a \vert \, = \,
\vert s \vert + \vert t - a \vert \, \le \, \varepsilon.
$$
Observe that $x \gamma = y \in A$. So, we have $B\gamma \subset A$ and
$$
\vert \log c_\mu(x, \gamma) - \log \lambda \vert \, \le \, \varepsilon
\hskip.5cm \text{for all} \hskip.2cm
x \in B.
$$
This shows that $\lambda\in r(X, \Gamma, \mu)$.
\end{proof}

We return to the situation of Example~\ref{Examples-ROIV-15}.

\begin{prop}
% 13.B.5
\label{Pro-ProjectiveSpaceTypeIII} 
Let $\Gamma$ be a lattice in $G = \SL_2(\R)$
and $\overline \Gamma$ its image in $\PSL_2(\R)$.
Let $X = \mathbf{P}^1(\R)$ be equipped with the Lebesgue measure $\mu$.
\par

The action of $\Gamma$ on $(X, \mu)$ by M\"obius transformations is of type III$_1$.
Therefore $L^\infty(X, \mu) \rtimes \overline \Gamma$ is a factor of type III$_1$.
\end{prop}

\begin{proof}
In view of Proposition~\ref{Pro-MaharamErgodic}, it suffices to show
that the Maharam extension $(\widetilde{X}, \widetilde{\mu}) \curvearrowleft \Gamma$
of $(X, \mu) \curvearrowleft \Gamma$ is ergodic.
\par

\index{$l 4$@$\overline \R = \R \cup \{\infty \} = \mathbf{P}^1(\R)$
real projective line}
We identify $X = \mathbf{P}^1(\R)$ with $\overline \R = \R \cup \{\infty \}$.
The (right) action of $\Gamma$ on $\widetilde{X} = \overline \R \times \R$ is given by 
$$
(x, t) \gamma \, = \,
\left( \frac{ax + c}{bx + d}, t - \log (bx + d)^2 \right)
\hskip.5cm \text{for} \hskip.2cm
(x, t) \in \overline \R \times \R, \hskip.1cm
\gamma = \begin{pmatrix} a & b \\ c & d \end{pmatrix} \in \Gamma .
$$
(Observe that this action extends to a measure preserving action
of $G$ on $(\widetilde{X}, \widetilde{\mu})$,
given by the same formula.)
\par

Consider the Borel isomorphism 
$\Phi \,\colon \R^2 \smallsetminus \{(0,0)\} \to \overline \R \times \R$
defined by 
$$
\Phi(u,v)\, = \, \begin{cases} 
(u/v, -\log v^2) & \text{if } v \ne 0
\\
(0, -\log u^2) & \text{if } v = 0 .
\end{cases}
$$
The image under $\Phi$ of the Lebesgue measure $du dv$ on $\R^2$
is $d\widetilde{\mu}(x, t) = 2e^t dx dt$.
Indeed, since the Jacobian of $\Phi$ at $(u,v)$ is $2/v^2$,
we have, for every measurable function $f \,\colon \widetilde{X} \to \R_+$,
$$
\begin{aligned}
\int_{\R^2} f(\Phi (u,v)) du dv
\, &= \, \int_{\R^2} f(u/v, -\log v^2) du dv
\\
\, &= \, \int_{\R^2} f(u/v, -\log v^2) du dv
\\
\, &= \, 2 \int_{\R^2} f(x, t) e^t dx dt.
\end{aligned}
$$
We claim that $\Phi$ is $\Gamma$-equivariant,
where the action of $\Gamma$ on $\R^2$ is the usual linear (right) action.
Indeed, let $\gamma = \begin{pmatrix} a & b \\ c & d \end{pmatrix} \in \Gamma$.
On the subset of $\R^2 \smallsetminus \{(0,0)\}$ of full Lebesgue measure
consisting of vectors $(u,v) \in \R^2 \smallsetminus\{(0,0)\}$
with $v \ne 0$ and $bu + dv \ne 0$, we have
$$
\begin{aligned}
\Phi((u,v) \gamma)
\, &= \, \Phi(au+cv, bu+dv)
\\
\, &= \, \left(\frac{au+cv}{bu+dv}, -\log (bu+dv)^2 \right)
\\
\, &= \, \left(\frac{a(u/v)+c}{b(u/v)+d}, -\log v^2 -\log (b(u/v)+d)^2 \right)
\\
\, &= \, (u/v, -\log v^2)\gamma
\\
\, &= \, \Phi(u,v) \gamma
\end{aligned} 
$$
As a consequence, it suffices to show that the linear action 
of $\Gamma$ on $(\R^2 \smallsetminus \{(0,0)\}, du dv)$ is ergodic.
This is indeed the case:
$G$ acts transitively on $\R^2 \smallsetminus \{(0,0)\}$
and the stabilizer of the point $(1,0) $ is the subgroup
$$
N \, = \, \left\{ \begin{pmatrix} 1& 0 \\ c& 1 \end{pmatrix} 
\hskip.1cm \Big\vert \hskip.1cm
c \in \R \right\}.
$$
Therefore $\R^2 \smallsetminus \{(0,0\}$ can be identified,
as a $G$-space, with $N\backslash G$.
Since $\Gamma$ is a lattice in $G$ and since $N$ is not relatively compact in $G$,
the claim follows from Moore's duality theorem and Moore ergodicity theorem
(compare with Example~\ref{Examples-ROIV-15}).
\end{proof}

The following corollary is an immediate consequence
of Propositions~\ref{Pro-ProjectiveSpaceTypeIII} and \ref{Pro-TypeIIINoInvMeasure}.

\begin{cor}
% 13.B.6
\label{Cor-NoInvMeasureProjectiveSpace} 
Let $\Gamma$ be a lattice in $SL_2(\R)$.
There exists no $\Gamma$-invariant $\sigma$-finite measure
on $\mathbf{P}^1(\R)$ which is equivalent to the Lebesgue measure.
\end{cor}

\begin{rem}
% 13.B.7
\label{Rem-ProjectiveSpaceTypeIII}
The proof of Proposition~\ref{Pro-ProjectiveSpaceTypeIII}
can easily be extended to show the following more general result.
Let $G$ be a simple non-compact connected Lie group,
$P$ a minimal parabolic subgroup in $G$,
and $\Gamma$ a lattice in $G$.
Then the action of $\Gamma $ on $(G/P, \mu)$ is of type III$_1$,
where $\mu$ is a $G$-quasi-invariant Radon measure on $G/P$
(see Example 4.3.15 in \cite{Zimm--84}). 
For further results establishing that such ``boundary actions" are of type III$_1$,
see \cite{Kaim--00} and \cite{Spat--87}.
\end{rem}

%-----------------------------------------------------------------------
% End of chapter 13
%-----------------------------------------------------------------------
\chapter[Constructing factor representations]
{Constructing factor representations for some semi-direct products}
% Chapter 14
\label{Chap:NormalInfiniteRep}

\emph{
Following \cite{Guic--63}, we first describe in this chapter
factor representations of a semi-direct product $G = H \ltimes N$, where
\begin{itemize}
\setlength\itemsep{0em}
\item
$H$ is discrete countable subgroup;
\item
$N$ is a second-countable locally compact abelian normal subgroup of $G$.
\end{itemize}
The factor representations we construct are induced representations
$\widetilde \pi_\mu = \Ind_N^\Gamma \pi_\mu$,
where $\mu$ is a convenient $\sigma$-finite $H$-quasi-invariant measure
on $\widehat N$
and $\pi_\mu$ is the canonical representation of $N$
on $L^2(\widehat N, \mu)$. 
The proof that $\widetilde \pi_\mu$
is indeed factorial as well as the identification of its type
rely on results on the measure space construction of Murray and von Neumann
from Section~\ref{SectionMSC}.
}
\par

\emph{
The particular case where $\mu$ is the Haar measure on $\widehat N$
is of special interest as the associated representation $\widetilde \pi_\mu$
is equivalent to the regular representation $\lambda_G$.
Criteria are given ensuring that $\lambda_G$ is a factor representation.
We give examples showing that $\widetilde \pi_\mu$
can be a factor representation of all possible types, with the exception of type I.
}
\par

\emph{
The representation $\widetilde \pi_\mu$ of a group $G = H \ltimes N$ as above
is a normal factor representation of type I$_\infty$ or of type II$_\infty$,
provided $\mu$ is a Radon measure
on an open and $H$-invariant subset $U$ of $\widehat N$ 
with the following properties:
\begin{itemize}
\setlength\itemsep{0em}
\item
$\mu$ is $H$-invariant;
\item
$\mu$ is infinite and non-atomic;
\item
the action of $H$ on $(U, \mu)$ is essentially free and ergodic.
\end{itemize}
In this case, the character associated to $\widetilde \pi_\mu$
admits an expression in terms of $\mu$. 
In Section \ref{Section-ExamplesInfiniteChar},
we will give examples of groups for which we can show
that measures $\mu$ as above exist;
we obtain in this way normal factor representations of type I$_\infty$
for the Baumslag--Solitar group $\BS(1, p)$,
Proposition \ref{Cor-NormalBS},
and normal factor representations of type I$_\infty$ and of type II$_\infty$
for the lamplighter group $\Z \wr (\Z / 2 \Z)$,
Propositions \ref{Cor-NormalLamplighter} and \ref{Cor-NormalLamplighter2}.
}
\par

\emph{
As a consequence of the construction
of normal factor representations of type I$_\infty$ for $\Gamma = \BS(1, p)$
or $\Gamma=\Z\wr (\Z /2 \Z)$
we will see that the natural map from $E(\Gamma)$ to $\Pri(\Gamma)$
is neither injective nor surjective,
in contrast to the case of nilpotent groups 
(Proposition~\ref{Prop-NilGr} and Corollary~\ref{ThomaHeis-PrimIdeal}),
the affine group over a field 
(Corollary~\ref{ThomaAff-PrimIdeal}), 
or the general linear group over an infinite algebraic extension 
of a finite field (Corollary~\ref{Cor-PrimGLn}).
}

\section
[Crossed product for semi-direct products]
{Crossed product von Neumann algebras for semi-direct products}
% Section 14.A
\label{Section-IndRepGMS}

Let $G = H \ltimes N$ be a semi-direct product of a countable discrete subgroup $H$
with a second-countable locally compact abelian normal subgroup $N$.

Recall (see Section~\ref{Section-IrrIndRep}) that 
$H$ acts on the right on $\widehat N$ by $(\chi, h) \mapsto \chi^h$,
where $\chi^h(n) = \chi(hnh^{-1})$ for
$\chi \in \widehat N, \hskip.1cm h \in H$, and $n \in N$. 

\begin{constr}
% 14.A.1
\label{constructionInd+mu}
Consider a group $G$ as above,
and a $\sigma$-finite Borel measure $\mu$
on the second-countable LC space $\widehat N$
(we do \emph{not} assume that $\mu$ is finite). 
We assume that $\mu$ is quasi-invariant by the action of $H$. 
\par

Let $\pi_\mu$ be the canonical representation of $N$ associated to $\mu$,
as in Construction~\ref{defpimupourGab};
recall that $\pi_\mu$ is defined on the Hilbert space $L^2(\widehat N, \mu)$ by 
$$
(\pi_\mu(n)f) (\chi) \, = \, \chi(n) f( \chi) \, = \, \widehat{n}(\chi) f(\chi)
\hskip.5cm \text{for} \hskip.2cm
n \in N, \hskip.1cm f \in L^2(\widehat N, \mu), \hskip.1cm \chi \in \widehat N. 
$$
Denote by $\widetilde \pi_\mu$
the induced representation $\Ind_N^G \pi_\mu$.
Since $H$ is a transversal for $N$ in $G$, we can realize $\widetilde \pi_\mu$
on the Hilbert space $\Hi_\mu := \ell^2(H, L^2(\widehat N, \mu))$ by 
$$
\begin{aligned}
(\widetilde \pi_\mu (h,n) F) (h') (\chi)
\, &= \, \left( \pi_\mu(n^{h'}) F(h' h) \right) (\chi)
\\
\, &= \, \chi^{h'}(n) F(h' h) (\chi),
\end{aligned}
$$
for $h, h' \in H, \hskip.1cm n \in N, \hskip.1cm F \in \Hi_\mu$, and $\chi \in \widehat N$
(see Construction~\ref{constructionInd}(2)).
\par

We denote by 
$$
\M \, := \, \widetilde \pi_\mu (G)'' \, \subset \, \Li (\Hi_\mu)
$$
the von Neumann algebra generated by $\widetilde \pi_\mu (G)$.
\end{constr}

We first identify $\M$ with the crossed product von Neumann algebra
$L^\infty(\widehat N, \mu) \rtimes H$
associated to the action $(\widehat N, \mu) \curvearrowleft H$,
as in Section~\ref{SectionMSC}. 
\par

Recall that, for $h \in H$ and $\chi \in \widehat N$, we set
$$
c(\chi, h) \, = \, \frac{dh_*(\mu)}{d\mu}(\chi),
$$
where $\dfrac{dh_*(\mu)}{d\mu}$
is the Radon--Nikodym derivate of $h_*(\mu)$ with respect to $\mu$.

\begin{prop}
% 14.A.2
\label{Pro-IndRepGMS}

Let $G = H \ltimes N$, $\mu$, $\Hi_\mu$, $\widetilde \pi_\mu$, and $\mathcal M$
be as in the previous construction.
\par

The von Neumann algebras $\M$ and $L^\infty(\widehat N, \mu) \rtimes H$
inside $\Li (\Hi_\mu)$ are unitarily equivalent.
More precisely, for the unitary operator $U \,\colon \Hi_\mu \to \Hi_\mu$ defined by 
$$
(UF) (h') (\chi) \, = \, c( \chi, h')^{1/2} F(h') (\chi^{h'})
\hskip.5cm \text{for} \hskip.2cm
h' \in H, \hskip.2cm F \in \Hi_\mu, \hskip.2cm \chi \in \widehat N ,
$$
we have
$$
\M \, = \, U^{-1}\left( L^\infty(\widehat N, \mu) \rtimes H) \right) U.
$$
\par

Moreover, if $\rho_\mu$ denotes the conjugate representation
$(h,n)\mapsto U^{-1} \widetilde \pi_\mu (h,n) U$ of $G$, we have 
$$
(\rho_\mu(h,n) F) (h') (\chi) \, = \, c( \chi, h)^{1/2}\chi(n) F(h'h) ( \chi^{h}),
$$
for $h, h' \in H, \hskip.1cm n \in N, \hskip.1cm F \in \Hi_\mu$, and $\chi \in \widehat N$.
\end{prop}

\begin{proof}
Let us first check that $\rho_\mu$ is given by the formula stated in the proposition.
Indeed, using the cocycle relation for $c$, one checks that
$$
(U^{-1}F) (h') (\chi) \, = \, c( \chi, h'^{-1})^{1/2} F(h') (\chi^{h'^{-1}}).
$$
Recalling that
$$
(\widetilde \pi_\mu (h,n) F) (h') (\chi) \, = \, \chi^{h'}(n) F(h' h) (\chi),
$$
we have therefore
$$
\begin{aligned}
(U^{-1}\widetilde \pi_\mu (h,n) UF) (h') (\chi)
\, &= \, c(\chi, h'^{-1})^{1/2} (\widetilde \pi_\mu (h,n) UF) (h')) (\chi^{h'^{-1}}))
\\
\, &= \, c(\chi, h'^{-1})^{1/2} \chi(n) (UF(h'h)) ( \chi^{h'^{-1}})
\\
\, &= \, c(\chi, h'^{-1})^{1/2} \chi(n) c( \chi^{h'^{-1}}, h'h)^{1/2}F(h'h) ( \chi^{h})
\\
\, &= \, c(\chi, h)^{1/2} \chi(n)F(h'h) (\chi^{h})
\\
\, &= \, (\rho_\mu(h,n) F) (h') (\chi),
\end{aligned}
$$
for $h, h' \in H, \hskip.1cm n \in N, \hskip.1cm F \in \Hi_\mu$, and $\chi \in \widehat N$.
\par

Next, we claim
% \footnote{Ou $U \M U^{-1}$ ?}
% \marginpar{Footnote}
that $U^{-1} \M U$ and $L^\infty (\widehat N, \mu) \rtimes H$
coincide, that is, 
$$
\rho_\mu(G)'' \, = \, L^\infty(\widehat N, \mu) \rtimes H.
$$
\par
 
The Koopman representation $h \mapsto u_h$
of $H$ on $L^2(\widehat N, \mu)$ is given by 
$$
u_h(f) (\chi) \, = \, c(\chi, h)^{1/2} f(\chi^{h}) 
\hskip.5cm \text{for} \hskip.2cm
f \in L^2(\widehat N, \mu), \hskip.1cm \chi \in \widehat N.
$$
Let $h \in H$. We have 
$$
(\rho_\mu(h, e) F) (h') (\chi) \, = \, F(h' h) ( \chi^{h})
\hskip.5cm \text{for} \hskip.2cm
h' \in H, \hskip.1cm F \in \Hi_\mu, \hskip.1cm \chi \in \widehat N
$$
and hence
$$
(\rho_\mu(h, e) F) (h') \, = \, u_h(F( h' h))
\hskip.5cm \text{for} \hskip.2cm
h' \in H, \hskip.1cm F \in \Hi_\mu.
$$
It follows that $\rho_\mu(h, e) = U_h$,
where $h \mapsto U_h$ is the representation of $H$
on $\Hi_\mu = \ell^2(H, L^2(\widehat N, \mu))$ involved in
the definition of the crossed product von Neumann algebra
$L^\infty(\widehat N, \mu) \rtimes H$. 
\par

Let $n \in N$. We have
$$
\begin{aligned}
(\rho_\mu(e,n) F) (h') (\chi) 
\, = \, &\chi(n) F(h') (\chi) \, = \, \widehat n (\chi) F(h') (\chi)
\\
&\text{for} \hskip.2cm
F \in \Hi_\mu, \hskip.1cm h' \in H, \hskip.1cm \chi \in \widehat N,
\end{aligned}
$$
where $\widehat n \in C^b(\widehat N)$ is defined by $\widehat n (\chi) = \chi(n)$.
It follows that $\rho_\mu(e,n) = M({\widehat n})$,
where $\varphi\mapsto M({\varphi})$
is the representation of $L^\infty(\widehat N, \mu)$
involved in the definition of $L^\infty(\widehat N, \mu) \rtimes H$. 
As a consequence, we see that $\M$
is contained in $L^\infty(\widehat N, \mu) \rtimes H$.
\par

It remains to show
% \footnote{Ou $U^? \mathcal M U^?$
% voir Footnote pr\'ec\'edente}
% \marginpar{Footnote}
that $L^\infty(\widehat N, \mu) \rtimes H$ is contained in $\M$.
Since $L^\infty(\widehat N, \mu) \rtimes H$ is generated by
$$
\{U_h \mid h \in H\} \, \cup \, \{M({\varphi}) \mid \varphi \in L^\infty(\widehat N, \mu) \}
$$
and since we have shown that $\{U_h\mid h \in H\} \subset \M$,
it suffices to show that 
$$
\{M({\varphi}) \mid \varphi \in L^\infty(\widehat N, \mu) \}
$$
is contained in $\M$.
\par

\index{Algebras! $2$@trigonometric polynomials on $\widehat G$}
\index{Trigonometric polynomials on $\widehat G$}
Let $\mathrm{Trig} (\widehat N)$ be the algebra
of trigonometric polynomials on $\widehat N$, that is, 
the linear span of $\{\widehat n \mid n \in N\}$.
% (see Appendix \ref{AppLCA+Pont}).
By Lemma~\ref{Lemma-DensityFourierTransform},
$\mathrm{Trig} (\widehat N)$ is dense 
in $L^\infty(\widehat N, \mu)$ for the weak *-topology.
Moreover, the map 
$$
L^\infty(\widehat N, \mu) \, \to \, \Li (\Hi_\mu),
\hskip.2cm
\varphi \mapsto M({\varphi})
$$
is continuous,
when $L^\infty(\widehat N, \mu)$ is equipped with the weak *-topology
and $\Li (\Hi_\mu)$ with the weak operator topology,
by Proposition~\ref{Prop-AbelianVN}.
\par

It follows that $\{M(\varphi) \mid \varphi \in \mathcal A \}$
is dense in $\{M({\varphi}) \mid \varphi \in L^\infty(\widehat N, \mu) \}$
for the weak operator topology.
Therefore, $\{M({\varphi})\mid\varphi \in L^\infty(\widehat N, \mu) \}$
is contained in $\M$.
\end{proof}

\section[Factor representations of semi-direct products]
{Constructing factor representations of some semi-direct products}
% Section 14.B
 \label{SectionConstrFactorRep}
 
As in Section~\ref{Section-IndRepGMS},
let $G = H \ltimes N$ be a semi-direct product
of a countable discrete subgroup $H$
with a second-countable locally compact abelian normal subgroup $N$.
We will use Proposition~\ref{Pro-IndRepGMS},
in combination with results from Section~\ref{SectionMSC}
on the crossed product von Neumann algebras,
in order to construct various factor representations of $G$.
\par

Recall from Section~\ref{Section-IndRepGMS} that we have associated
a representation $\widetilde \pi_\mu$ of $G$
in $\ell^2(H, L^2(\widehat N, \mu))$
to every $\sigma$-finite Borel measure $\mu$
on $\widehat N$ which is quasi-invariant by the action of $H$.
 
\begin{theorem}
% 14.B.1
\label{Theo-FacRepSemiDirect}
Let $G = H \ltimes N$,
where $H$ is a countable discrete subgroup
and $N$ a second-countable locally compact abelian normal subgroup.
Let $\mu$ be a $\sigma$-finite Borel measure $\mu$ on $\widehat N$
which is quasi-invariant by the action of $H$
and $\widetilde \pi_\mu$ the associated representation of $G$
as in Section~\ref{Section-IndRepGMS}.
Assume that the action $H \curvearrowright (\widehat N, \mu)$
is essentially free and ergodic.
\begin{enumerate}[label=(\arabic*)]
\item\label{iDETheo-FacRepSemiDirect}
The representation $\widetilde \pi_\mu$ is a factor representation of $G$.
\item\label{iiDETheo-FacRepSemiDirect}
If $\mu$ is an atomic measure,
then $\widetilde \pi_\mu$ is a factor representation of type I.
\item\label{iiiDETheo-FacRepSemiDirect}
If $\mu$ is $H$-invariant non-atomic finite measure,
then $\widetilde \pi_\mu$ is a factor representation of type II$_1$.
\item\label{ivDETheo-FacRepSemiDirect}
If $\mu$ is a $H$-invariant non-atomic infinite measure,
then $\widetilde \pi_\mu$ is a factor representation of type II$_\infty$.
\item\label{vDETheo-FacRepSemiDirect}
Assume that there does not exist any $H$-invariant $\sigma$-finite measure
on $\widehat N$ which is equivalent to $\mu$;
then $\widetilde \pi_\mu$ is a factor representation of type III.
\end{enumerate}
\end{theorem}

\begin{proof}
By Proposition~\ref{Pro-IndRepGMS}, $\M := \widetilde \pi_\mu (G)''$
is unitarily equivalent to the crossed product von Neumann algebra
$L^\infty(\widehat N, \mu) \rtimes H$.
Since the action $(\widehat N, \mu) \curvearrowleft H$
is essentially free and ergodic,
Items \ref{iDETheo-FacRepSemiDirect} to \ref{vDETheo-FacRepSemiDirect}
follow from Items \ref{iDEROIV} to \ref{IIIDEROIV} in Theorem~\ref{ROIV},
respectively.
\end{proof}

\begin{exe}
% 14.B.2
\label{Exa-FacRepSemiDirect}
We can use Theorem \ref{Theo-FacRepSemiDirect} in order to produce 
factor representations of type II$_\infty$ and of type III
for any one of our favorite groups $\Gamma$
(with the exception of the general group $GL_n(\K)$)
as well as recover some of the irreducible 
representations and representations of type II$_1$ of $\Gamma$
considered in previous sections.
We will consider here only the case where $\Gamma$
is the Heisenberg group over $\Z$.
Similar results hold for the Heisenberg group and the affine group
over a countable infinite field
and for the Baumslag--Solitar groups $\BS(1, p)$.
\par

Recall that $\Gamma = H(\Z)$ is the semi-direct product $H \ltimes N$, where 
$$
H \, = \, \{(a, 0, 0) \in \Gamma \mid a \in \Z\}
\hskip .5cm \text{and} \hskip .5cm
N \, = \, \{(0, b, c) \in \Gamma \mid b,c \in \Z\}
$$
Fix $\theta \in \mathopen[ 0,1 \mathclose[$ irrational
and let $X \subset \widehat N$ be the $H$-invariant set
of unitary characters of $N$ with $\chi\vert_Z = \psi_{\theta}$,
where $\psi_\theta(0,0,c) = e^{2 \pi i \theta c}$. 
Then $X$ is homeomorphic to $\T$
and the action of $H$ on $X$ corresponds to
the action of the irrational rotation $R_\theta$ on the circle group $\T$
(see Corollary~\ref{Cor-IrredRepHeisInteger}).
In particular, the action of $H$ on $X$ is free. 
So, the action of $H$ on $(X, \mu)$ is automatically
essentially free for any $H$-quasi-invariant measure $\mu$ on $X$.
\begin{itemize}
\setlength\itemsep{0em} 
\item[$\bullet$]
Let $\mu$ be the counting measure on a $R_\theta$-orbit $\OO$ in $X$.
Then the representation $\pi_\mu$ 
is equivalent to the direct sum $\bigoplus_ {\chi \in \OO} \chi$ and hence
$\widetilde \pi_\mu = \Ind_N^\Gamma \pi_\mu$ is equivalent to 
$$
\bigoplus_{\chi \in\OO} \Ind_N^\Gamma \chi.
$$
As seen in Corollary~\ref{Cor-IrredRepHeisInteger},
the representations $\Ind_N^\Gamma \chi$ are irreducible and mutually equivalent.
As a result, $\widetilde \pi_\mu$ is a multiple of the irreducible representation 
$\Ind_N^\Gamma \chi$ for any $\chi \in \OO$. 
Thus, we recover the fact proved in Theorem~\ref{Theo-FacRepSemiDirect} 
that $\widetilde \pi_\mu$ is a factor representation of type I;
in this case, $\widetilde \pi_\mu$ is of case I$_\infty$.
\par
 
\item[$\bullet$]
Take as $\mu$ the Lebesgue measure on $X$.
Then the representation $\pi_\mu$ is equivalent to
the induced representation $\Ind_Z^N \psi_\theta$
and hence $\widetilde \pi_\mu$ is equivalent to $\Ind_Z^\Gamma \psi_\theta$. 
This is the representation of type II$_1$ associated to the 
character $\widetilde \psi_\theta \in E(\Gamma)$
from Corollary~\ref{Cor-ThomaDualHeisIntegers}.
\par

\item[$\bullet$]
By Theorem~\ref{Theo-Schmidt},
there exist uncountably many
non-equivalent infinite and $\sigma$-finite non-atomic measures
on $X$ which are $T_\theta$-invariant and ergodic. 
For every such measure $\mu$, the representation $\widetilde \pi_\mu$
of $\Gamma$ is a factor representation of type II$_\infty$.
\par

Observe that we do not know
whether non-equivalent measures $\mu$ as above
yield non-quasi-equivalent factor representations $\widetilde \pi_\mu$ 
(see however Proposition~\ref{Prop-CharactersSemiDirect} below).

\item[$\bullet$]
Theorem~\ref{Theo-Schmidt} has been extended in \cite{Krie--76}
to cover the existence of quasi-invariant but non-invariant measures,
at least in the case of $\Z$-actions.
This result implies that there exist uncountably many non-equivalent
non-atomic probability measures on $X$ which are quasi-invariant by $R_\theta$
and for which there is no equivalent $\sigma$-finite $R_\theta$-invariant measure.
For every such measure $\mu$,
the representation $\widetilde \pi_\mu$ of $\Gamma$
is a factor representation of type III.
\end{itemize}
\end{exe}

\section[Some normal factor representations]
{Some normal factor representations}
% Section 14.C
\label{SectionCharactersNormalRep}

As in Section~\ref{Section-IndRepGMS},
let $G = H \ltimes N$ be a semi-direct product
% \footnote{D\'ecider
% si $H$ agit \`a gauche ou \`a droite sur $N$.}
% \marginpar{D\'ecider~!}
of a countable discrete subgroup $H$
with a second-countable locally compact abelian normal subgroup~$N$. 
\par
 
Let $\mu$ be a $\sigma$-finite Borel measure on $\widehat N$ 
which is quasi-invariant, essentially free, and ergodic by the action of $H$.
As seen in Theorem~\ref{Theo-FacRepSemiDirect},
the associated representation $\widetilde \pi_\mu$ of $G$
is then a factor representation.
We will show below that, under a further assumption,
$\widetilde \pi_\mu$ is a \emph{normal} factor representation.
\par

Let $\nu_N$ denote a Haar measure on $N$. 
Observe that $\nu_N$ is relatively invariant under the action of $H$:
for every $h \in H$, there exists $\Delta(h) > 0$
such that $h_*(\nu_N) = \Delta(h) \nu_N$;
indeed,
the Borel measure $h_*(\nu_N)$ on $N$ is invariant under $N$
and is therefore equal to $\Delta(h)\nu_N$ for a constant $\Delta(h) > 0$,
by uniqueness of Haar measures.
A Haar measure $\nu_G$ on $G$
is given by $\nu_G = (\Delta\nu_H) \otimes \nu_N$,
where $\nu_H$ is the counting measure on $H$; so, we have
$$
\int_G f(g) d\nu_G(g) \, = \, \sum_{h \in H} \Delta(h) \int_N f(h,n) d\nu_N(n)
$$
for every measurable function $f \,\colon G \to \mathopen[ 0, \infty \mathclose]$.
\par

For a function $f \,\colon G \to \C$, we denote by $f_e$ the function on 
$N$ defined by
$$
f_e(n) \, = \, f(e,n) 
\hskip.5cm \text{for all} \hskip.2cm
n \in N.
$$
Observe that $f_e \in L^1(N, \nu_N)$, if $f \in L^1(G, \nu_G)$.
We set 
$$
L^1(G, \nu_G)_+ \, = \, \{f^* \ast f \mid f \in L^1(G, \nu_G)\} .
$$
Observe that every $f \in L^1(G, \nu_G)_+$
is positive as element of the C*-algebra $C^*_\lambda (G)$.
Recall that the involution on $L^1(G, \nu_G)$
is given here by $f^*(g) = \overline{f(g^{-1})}$ for $f \in L^1(G, \nu_G)$ and $g \in G$. 
\par

The following technical lemma will be used
in the proofs of Theorem~\ref{Theo-CharactersSemiDirect}
and Proposition~\ref{Prop-CharactersSemiDirect} below.
As in Appendix \ref{AppTop}, we denote by $C^{c, K}(\widehat N)$
the space of complex-valued continuous functions on $\widehat N$
with support contained in some compact subspace $K$ of $\widehat N$.

\begin{lem}
% 14.C.1
\label{Lem-CharactersSemiDirect}
Let $U$ be an open subset of $\widehat N$ and $K$ a compact subset of $U$. 
\par

Then there exists a compact subset $K'$ of $U$ with $K \subset K'$ such that
$C^{c, K}(\widehat N)$ is contained in the closure of the linear span of
$$
\{ \overline{\mathcal F} (f_e) \mid f \in L^1(G, \nu_G)_+ \} \cap C^{c, K'}(\widehat N)
$$
for the topology of uniform convergence on $\widehat N$.
\end{lem}

\begin{proof}
Since $U$ is open and $K$ compact,
we can find a compact neighbourhood $V$ of $\{e\}$
such that $KV$ in contained in $U$.
Let $K' = KV$ and set 
$$
A_{K'} \, := \, \{{\overline{\mathcal F}}(f_e) \mid f \in L^1(G, \nu_G)_+\} \cap C^{c, K'}(\widehat N) .
$$
Let $\varphi \in C^{c, K}(\widehat N)$.
We claim that $\varphi$ is in the closure of the linear span of $A_{K'}$. 
\par

Upon considering the positive and negative parts
of the real and imaginary parts of $\varphi$,
we can assume that $\varphi \ge 0$.
\par

We assume that the Haar measures
$\nu_N$ on $N$ and $\nu_{\widehat N}$ on $\widehat N$
are so normalized that the Fourier transform
$$
\mathcal{F} \, \colon \, L^2(N, \nu_N) \to L^2(\widehat N, \nu_{\widehat N})
$$
is an isomorphism of Hilbert spaces (Plancherel theorem \ref{PlancherelTh}).
Then this is also true for the inverse Fourier transform 
$$
\overline{\mathcal F} \, \colon \, L^2(N, \nu_N)\to L^2(\widehat N, \nu_{\widehat N}).
$$
\par
 
Choose an approximate identity supported in $V$, that is,
a sequence $(\psi_i)_i$ in $C^c(\widehat N)$ such that 
$\mathrm{supp}(\psi_i) \subset V$ and $\lim_i \psi_i \ast \psi = \psi$ 
uniformly on $\widehat N$ for every $\psi \in C^c(\widehat N)$.
\par

Set $\psi := \varphi^{1/2}$ and set 
$$
r \, := \, \overline{\mathcal F}^{-1}(\psi)
\hskip.5cm \text{and} \hskip.5cm
r_i \, := \, \overline{\mathcal F}^{-1}(\psi_i).
$$
Since $r_i \in L^2(N, \nu_N)$ and $r \in L^2(N, \nu_N)$,
we have $r_i r \in L^1(N, \nu_N)$.
Moreover, 
$$
\overline{\mathcal F}(r_ir) \, = \, 
\overline{\mathcal F}(r_i) \ast \overline{\mathcal F}(r) \, = \,
\psi_i \ast \psi
$$
and hence $\lim_i \overline{\mathcal F}(r_ir) = \psi$ uniformly on $\widehat N$.
\par

Set 
$$
s_i \, := \, (r_ir)^* \ast r_i r \in L^1(N, \nu_N).
$$
We claim that $\lim_i \overline{\mathcal F}(s_i) = \varphi$
uniformly on $\widehat N$
and that $\overline{\mathcal F}(s_i) \in A_{K'}$, for every $i$.
\par 

Indeed, we have 
$$
\overline{\mathcal F}(s_i) \, = \, \vert \overline{\mathcal F}(r_ir) \vert^2
$$
and hence
$$
\lim_i \overline{\mathcal F}(s_i) \, = \, \vert \psi \vert^2 \, = \, \varphi
$$
uniformly on $\widehat N$. 
\par
 
Moreover, since $\mathrm{supp}(\psi_i)\subset V$
and $\mathrm{supp}(\psi)\subset K$,
we have $\mathrm{supp}(\psi_i \ast \psi) \subset KV = K'$.
As $\overline{\mathcal F}(s_i) = \vert \psi_i \ast \psi \vert^2$, it follows that
$$
\mathrm{supp}(\overline{\mathcal F}(s_i))
\, = \, 
\mathrm{supp}(\psi_i \ast \psi) \subset K'.
$$
Therefore $\overline{\mathcal F}(s_i) \in C^{c, K'}(\widehat N)$.
\par
 
It remains to show that every $s_i$ is of the form $(f^{(i)})_0$
for some $f^{(i)} \in L^1(G, \nu_G)_+$.
Let $g_i \in L^1(G, \nu_G)$ be the trivial extension of $r_ir$ to $G$,
defined by 
$$
g_i(h,n) \, = \,
\begin{cases}
(r_ir) (n)
\, & \, \text{if} \hskip.2cm h = e
\\
0
\, & \, \text{if} \hskip.2cm h \ne e.
\end{cases}
$$
For $n \in N$, we have
$$
\begin{aligned}
(g_i^* \ast g_i)_0(n)
\, &= \, (g_i^* \ast g_i) (e,n)
\, = \, \sum_{h' \in H}\Delta(h') \int_N g_i^*((h',n')^{-1}) g_i((h',n') (e,n)) d\nu_N(n')
\\
\, &= \, \sum_{h' \in H}\Delta(h') \int_N \overline{g_i(h',n')}g_i(h', n' h' \cdot n)d\nu_N(n')
\\
\, &= \, \int_N \overline{g_i(e,n')}g_i(e, n'n)d\nu_N(n')
\\
\, &= \, \int_N \overline{(r_ir) (n')}(r_ir) (n'n))d\nu_N(n')
\\
\, &= \, \int_N ((r_ir)^*(n'^{-1}) (r_ir) (n'n)d\nu_N(n')
\\
\, &= \, \left( (r_ir)^* \ast (r_ir) \right) (n)
\\
\, &= \, s_i(n)
\end{aligned}
$$
So, for $f^{(i)} := g_i^* \ast g_i \in L^1(G, \nu_G)_+$, we have $s_i = f^{(i)}_0$.
Therefore for every $i$, we have $\overline{\mathcal F}(s_i) \in A_{K'}$. 
\end{proof}

\begin{theorem}
% 14.C.2
\label{Theo-CharactersSemiDirect}
Let $G = H \ltimes N$,
where $H$ is a countable discrete subgroup
and $N$ a second-countable locally compact abelian normal subgroup.
Let $\mu$ be a non-zero $H$-invariant Radon measure
on an open $H$-invariant subset $U$ of $\widehat N$.
Assume moreover that the action $H \curvearrowright (\widehat N, \mu)$
is essentially free and ergodic.
\begin{enumerate}[label=(\arabic*)]
\item\label{iDETheo-CharactersSemiDirect}
The representation $\widetilde \pi_\mu$ associated to $\mu$
is a normal factor representation;
moreover, the corresponding character $t$ is given by 
$$
t(f) \, = \, \int_{\widehat N} {\overline{\mathcal F}}(f_e) (\chi) d\mu(\chi)
\hskip.5cm \text{for} \hskip.2cm
f \in L^1(G)_+,
$$
where ${\overline{\mathcal F}}$ denotes the inverse Fourier--Stieltjes transform.
\item\label{iiDETheo-CharactersSemiDirect}
The factor representation $\widetilde \pi_\mu$ is of type I$_\infty$
if $\mu$ is an infinite atomic measure.
\item\label{iiiDETheo-CharactersSemiDirect}
The factor representation $\widetilde \pi_\mu$ is of type II$_\infty$
if $\mu$ is an infinite non-atomic measure. 
\end{enumerate}
\end{theorem}

\begin{proof}
\ref{iDETheo-CharactersSemiDirect}
By Theorem~\ref{Theo-FacRepSemiDirect},
the von Neumann algebra $\M := \widetilde \pi_\mu (G)''$ is a semi-finite factor.
We have to show that $\widetilde \pi_\mu$ is a trace representation of $G$.
\par
 
Upon conjugating $\M$ by a unitary operator, we can assume that
$\M$ coincides with the crossed product von Neumann algebra
$L^\infty(\widehat N, \mu) \rtimes H$
and that $\widetilde \pi_\mu$ is the representation $\rho_\mu$
from Proposition~\ref{Pro-IndRepGMS};
thus, we can assume that
$$
(\widetilde \pi_\mu (h,n) F) (h') (\chi) \, = \, \chi(n) F(h^{-1}h') ( \chi^{h}),
$$
for $h, h' \in H, \hskip.1cm n \in N, \hskip.1cm F \in \Hi_\mu$, and $\chi \in \widehat N$.
\par

Let $\tau$ be the unique trace on $\M$ (up to positive scalar multiples).
By Proposition~\ref{NormalRep}, it suffices to show
there exists an element $f \in L^1(G)$ such that
$$
0 \, < \, \tau(\pi_\mu(f^* \ast f)) \, < \, +\infty.
$$
\par

Set $\Hi_\mu := \ell^2(H, L^2(\widehat N, \mu))$ and 
let $E \,\colon \M \to L^\infty(\widehat N, \mu)$ be the conditional expectation onto 
$L^\infty(\widehat N, \mu)$. Recall from Section~\ref{SectionMSC} that, for $T \in \M$, 
the function $E(T)$ from $L^\infty(\widehat N, \mu)$ is determined by 
$$
\int_{\widehat N} E(T) (\chi) f_1(\chi) \overline{f_2(\chi)} d\mu(\chi)
\, = \,
\langle T(\delta_e \otimes f_1) \mid \delta_e \otimes f_2 \rangle
\hskip.5cm \text{for all} \hskip.2cm
f_1, f_2 \in L^2(\widehat N, \mu),
$$
where $\delta_e \otimes f_i \in \Hi_\mu$
is defined by $(\delta_e \otimes f_i) (h,n) = f_i(n)$ if $h = e$
and $(\delta_e \otimes f_i) (h,n) = 0$ otherwise.
\par

By Proposition~\ref{Prop-CondExpTrace}~\ref{iiDEProp-CondExpTrace}
(see also Remark~\ref{Rem-Prop-CondExpTrace}),
we have 
$$
\tau(T) \, = \, \mu \circ E(T) \, = \, \int_{\widehat N} E(T) (\chi) d\mu(\chi)
\hskip.5cm \text{for all} \hskip.2cm
T \in \M.
$$
\par

Let $f \in L^1(G)_+$; set 
$$
T \, := \, \pi_\mu(f) \in \Li (\Hi_\mu).
$$
For $h' \in H, \hskip.1cm F, F' \in \Hi_\mu$, and $\chi \in \widehat N$,
we have
$$
\begin{aligned}
\langle T(F) \mid F' \rangle
\, &= \, \sum_{h \in H} \Delta(h) \int_N f(h,n)
\langle \pi_\mu(h,n) F \mid F' \rangle d\nu_N(n)
\\
\, &= \, \sum_{h, h' \in H} \Delta(h) \int_N \int_{\widehat N}
f(h,n) \chi(n) F(h^{-1} h') (\chi^h) \overline{F'(h') (\chi)} d\nu_N(n)d\mu(\chi);
\end{aligned}
$$
in particular, we have
$$
\begin{aligned}
\langle T(\delta_e \otimes f_1)\mid \delta_e \otimes f_1 \rangle
\, &= \,
\int_N \int_{\widehat N} f(e,n) \chi(n) f_1(\chi) \overline{f_2(\chi)} d\nu_N(n)d\mu(\chi)
\\
\, &= \, \int_{\widehat N }\left( \int_N f(e,n)\chi(n) d\nu_N(n) \right)
f_1(\chi) \overline{f_2(\chi)} d\mu(\chi)
\end{aligned}
$$
for all $f_1, f_2 \in L^2(\widehat N, \mu)$.
\par

This shows that 
$$
E(T) (\chi) \, = \, \int_N f(e,n)\chi(n) d\nu_N(n) \, = \, \overline{{\mathcal F}}(f_e) (\chi);
$$
we have therefore
$$
\tau(\pi_\mu(f)) \, = \, \int_{\widehat N}\overline{{\mathcal F}}(f_e) (\chi) d\mu(\chi)
\hskip.5cm \text{for all} \hskip.2cm
f \in L^1(G)_+.
$$
\par

Since $\mu$ is a Radon on $U$, 
it takes finite values on $C^{c, K}(U)$ for every 
compact subspace $K$ of $U$. 
Moreover, since $\mu$ is non-zero,
there exists a compact subset $K$ of $U$ and 
$\varphi \in C^{c,K}(U)$ with 
$$
\int_{\widehat N} \varphi(\chi) d\mu(\chi) \, > \, 0.
$$
It follows then from Lemma~\ref{Lem-CharactersSemiDirect},
that there exist a compact subset $K'$ of $U$ with $K \subset K'$
and a function $f \in L^1(G)_+$ such that
$\overline{{\mathcal F}}(f_e) \in C^{c, K'}(U)$ and such that 
$$
\int_{\widehat N} \overline{{\mathcal F}}(f_e) (\chi) d\mu(\chi) \, > \, 0.
$$
Since $\overline{{\mathcal F}}(f_e) \in C^{c, K'}(U)$, we have also 
$$
\int_{\widehat N} \overline{{\mathcal F}}(f_e) (\chi) d\mu(\chi) \, < \, +\infty.
$$
So, $\widetilde \pi_\mu$ is a normal factor representation
and we have identified its trace. 

\vskip.2cm

\ref{iiDETheo-CharactersSemiDirect}
Assume that $\mu$ is an infinite atomic measure. 
Then, by Theorem~\ref{Theo-FacRepSemiDirect}, $\widetilde \pi_\mu$ is of type I.
We have to check that $\M$ is an infinite factor. 
\par

By invariance and ergodicity,
$\mu$ takes a constant value $\alpha > 0$
on every $\chi$ in the support $S$ of $\mu$.
For every subset $A \subset S$, we have,
viewing $\Un_{A}$ as a projection in $L^{\infty}(\widehat N, \mu) \subset \M$,
$$
\tau(\Un_A) \, = \, \int_{\widehat N} \Un_A(\chi) d\mu(\chi) \, = \, \alpha \vert A \vert
$$
Since $S$ is infinite, this shows that $\tau$ takes infinitely many values on projections.
Therefore $\M$ is a factor of infinite type.

\vskip.2cm

Item \ref{iiiDETheo-CharactersSemiDirect}
follows from~\ref{Theo-FacRepSemiDirect}.
\end{proof}

\begin{rem}
% 14.C.3
\label{Rem--CharactersSemiDirect1}
Let $\mu$ be a measure as in Theorem~\ref{Theo-CharactersSemiDirect};
assume that $\mu$ is atomic.
Then the support $S$ of $\mu$ is an $H$-orbit, by ergodicity of $\mu$.
Since $\mu$ is a regular measure, $S$ has to be a discrete subset of $\widehat N$.
\par

For instance, let $S$ be an orbit of an irrational rotation on $\T$.
Then $S$ is dense in $\T$ and so the corresponding counting measure $\mu$
does not fulfill the assumption of Theorem~\ref{Theo-CharactersSemiDirect}.
The associated representation $\widetilde \pi_\mu$, 
which is a factor representation of type I$_\infty$
of the Heisenberg group $H(\Z)$ as in Example~\ref{Exa-FacRepSemiDirect},
is \emph{not} normal.
Indeed, $H(\Z)$ has no normal factor representation of infinite type at all
(Proposition~\ref{Prop-NilGr}).
\end{rem}

Next, we show that non-equivalent measures $\mu$
as in Theorem~\ref{Theo-CharactersSemiDirect}
yield non-quasi-equivalent factor representations $\widetilde \pi_\mu$.

\begin{prop}
% 14.C.4
\label{Prop-CharactersSemiDirect}
Let $G$ be as in Theorem~\ref{Theo-CharactersSemiDirect}.
Let $\mu$ and $\mu'$ be $H$-invariant Radon measures
on open $H$-invariant subsets $U$ and $U'$ of $\widehat N$.
Assume moreover that the actions $(\widehat N, \mu) \curvearrowleft H$
and $(\widehat N, \mu') \curvearrowleft H$
are essentially free and ergodic.
The following properties are equivalent:
\begin{enumerate}[label=(\roman*)]
\item\label{iDEProp-CharactersSemiDirect}
the factor representations $\widetilde \pi_\mu$ and $\widetilde \pi_{\mu'}$
are quasi-equivalent;
\item\label{iiDEProp-CharactersSemiDirect}
the restrictions of $\mu$ and $\mu'$ to $U\cap U'$
are proportional to each other and non-zero.
\end{enumerate}
\end{prop}

\begin{proof}
Set $U'' := U \cap U'$. Since $\mu$ is ergodic and $U''$ is $H$-invariant, we have 
either $\mu(U'') = 0$ or $\mu(\widehat N \smallsetminus U'') = 0$.
Similarly, we have 
either $\mu'(U'') = 0$ or $\mu'(\widehat N \smallsetminus U'') = 0$.
\par
Assume that $\mu\vert_{U''}$ and $\mu' \vert_{U''}$
are proportional to each other and non-zero.
Then, we have necessarily
$$
\mu(\widehat N \smallsetminus U'') \, = \, \mu'(\widehat N \smallsetminus U'') \, = \, 0.
$$
It follows from Theorem~\ref{Theo-CharactersSemiDirect}
that the characters $t$ and $t'$
associated to $\widetilde \pi_\mu$ and $\widetilde \pi_{\mu'}$
are proportional to each other.
Therefore $\widetilde \pi_\mu$ and $\widetilde \pi_{\mu'}$
are quasi-equivalent, by Proposition~\ref{Char-NormalRep}.

\vskip.2cm

Conversely, assume that $\widetilde \pi_\mu$ and $\widetilde \pi_{\mu'}$
are quasi-equivalent.
Then the characters $t$ and $t'$ are proportional to each other
(Proposition~\ref{Char-NormalRep}),
that is, there exists $\alpha > 0$ such that 
$$
\hskip.5cm
\int_{\widehat N} {\overline{\mathcal F}} (f_e) (\chi) d\mu(\chi)
\, = \,
\alpha \int_{\widehat N} {\overline{\mathcal F}}(f_e) (\chi) d\mu'(\chi)
\hskip.5cm \text{for} \hskip.2cm
f \in L^1(G)_+,
\leqno{(*)}
$$
Assume, by contradiction, that $\mu(U'') = 0$. Then $\mu(U') = 0$,
since $\mu(\widehat N \smallsetminus U) = 0$ and $U'' = U \cap U'$.
\par

Lemma~\ref{Lem-CharactersSemiDirect},
applied to the Radon measure $\mu'$ on $U'$, shows that we can find 
$f \in L^1(G)_+$ such that $\overline{\mathcal F}(f_e) = 0$
on $\widehat N \smallsetminus U'$ and such that 
$$
\int_{\widehat N}{\overline{\mathcal F}}(f_e) (\chi) d\mu'(\chi) \, > \, 0.
$$
On the one hand, we then have, by $(*)$,
$$
\begin{aligned}
\int_{U'}{\overline{\mathcal F}}(f_e) (\chi) d\mu(\chi)
\, &= \, \int_{\widehat N}{\overline{\mathcal F}}(f_e) (\chi) d\mu(\chi)
\\
\, &= \, \alpha \int_{\widehat N}{\overline{\mathcal F}}(f_e) (\chi) d\mu'(\chi)
\, > \, 0 .
\end{aligned}
$$
On the other hand, we have
$$
\int_{U'}{\overline{\mathcal F}} (f_e) (\chi) d\mu(\chi) \, = \, 0,
$$
since $\mu(U') = 0$.
This contradiction shows that $\mu(\widehat N \smallsetminus U'') = 0$.
Similarly, we have $\mu'(\widehat N \smallsetminus U'') = 0$.
As a result, we may assume that both $\mu$ and $\mu'$
are Radon measures on the open set $U''$.
\par

In view of Lemma~\ref{Lem-CharactersSemiDirect},
it follows from $(*)$ that $\mu$ agrees with $\alpha\mu'$ on $C^{c, K}(U'')$,
for every compact subset $K$ of $U''$.
Therefore $\mu = \alpha\mu'$ on $U''$.
\end{proof}

\begin{rem}
% 14.C.5
\label{Rem--CharactersSemiDirect}
One cannot expect that every normal factor representation of $G = H \ltimes N$ 
is quasi-equivalent to a representation $\widetilde \pi_\mu$
for a measure $\mu$ as in Theorem~\ref{Theo-CharactersSemiDirect}. 
\par

Indeed, as is easily checked,
the equivalence class of measures attached to $\widetilde \pi_\mu\vert_N$
as in Proposition~\ref{Prop-RestNormalSub-Bis} is the class $[\mu]$.
If $\pi$ is a factor representation which is quasi-equivalent to $\widetilde \pi_\mu$,
then the class attached to $\pi \vert_N$ is $[\mu]$. 
However, there are examples of factor representations,
for which the action $H$ on $(\widehat N, \mu')$ 
is not essentially free, for some (equivalently, for any) measure $\mu'$
in the class attached to $\pi \vert_N$.
Indeed, if $\pi$ is a normal factor representation of $H$ lifted to $G$,
then the Dirac measure at the group unit of $N$ is attached to $\pi \vert_N$;
however, $(\widehat N, \mu')$
is not essentially free (when $H$ and $N$ are non-trivial). 
\par

The answer to the following question, which is more sensible to ask,
seems to be unknown;
for some partial results, see \cite[Proposition 1, \S~2]{Guic--63}.

\vskip.2cm

\textbf{Question:}
Let $\pi$ be a normal factor representation of $G$
and $\mu'$ an ergodic and quasi-invariant probability measure on $\widehat N$
attached to $\pi \vert_N$.
Assume that the action $H$ on $(\widehat N, \mu')$ is essentially free.
Is $\pi$ quasi-equivalent to $\widetilde \pi_\mu$ for a measure $\mu$
as in Theorem~\ref{Theo-CharactersSemiDirect}?
\end{rem}

\section[Regular representation]
{The von Neumann algebra of the regular representation as a crossed product}
% Section 14.D
\label{Section-VNRegular}

As in the previous sections, $G = H \ltimes N$
is a semi-direct product of a discrete countable subgroup $H$
with a second-countable locally compact abelian normal subgroup $N$.
\par
 
We apply the results from Sections~\ref{Section-IndRepGMS}, \ref{SectionConstrFactorRep}
and \ref{SectionCharactersNormalRep}
to the special case where $\mu$ is a Haar measure on $\widehat N$. 
As already mentioned in Section~\ref{SectionCharactersNormalRep}
(for $N$ instead of $\widehat N$),
$\mu$ is relatively invariant (and hence quasi-invariant) by the action of $H$:
for every $h \in H$, there exists a constant $\Delta(h) > 0$
such that $h_*(\mu) = \Delta(h) \mu$.

\begin{prop}
% 14.D.1
\label{Prop-VNRegular}
Let $G$ be a second-countable locally compact group;
assume that $G = H \ltimes N$
is a semi-direct product of a discrete subgroup $H$
with a locally compact abelian normal subgroup $N$.
Let $\mu$ be a Haar measure on $\widehat N$ and let $\widetilde \pi_\mu$ 
be the associated representation as in Section~\ref{Section-IndRepGMS}.
\begin{enumerate}[label=(\arabic*)]
\item\label{iDEProp-VNRegular}
The representation $\widetilde \pi_\mu$
is equivalent to the regular representation $\lambda_G$ of $G$.
\item\label{iiDEProp-VNRegular}
The crossed product von Neumann algebra $L^\infty(\widehat N, \mu) \rtimes H$
is unitarily equivalent
to the von Neumann algebra $\Li (G) = \lambda_G(G)''$
generated by the regular representation of $G$.
\end{enumerate}
\end{prop} 

\begin{proof}
\ref{iDEProp-VNRegular}
With $\pi_\mu$ as in Section~\ref{Section-IndRepGMS},
we have $\widetilde \pi_\mu = \Ind_N^G \pi_\mu$.
Since $\lambda_G$ is equivalent to $\Ind_N^G \lambda_N$,
it suffices to show that $\pi_\mu$ is equivalent to $\lambda_N$. 
\par

We can assume that the Haar measures $\nu$ on $N$ and $\mu$ on $\widehat N$
are so normalized that the inverse Fourier transform
$\overline{\mathcal{F}} \,\colon L^2(N, \nu) \to L^2(\widehat N, \mu)$
is a Hilbert space isomorphism.
Moreover, $\overline{\mathcal{F}}$ intertwines $\lambda_N$ and $\pi_\mu$:
$$
\begin{aligned}
\overline{\mathcal{F}}(\lambda_N(n)f) (\chi)
\, &= \, \int_N (\lambda_N(n)f) (x)\chi(x) d\nu(x)
% \\
\, = \, \int_N f(n^{-1} x)\chi(x) d\nu(x)
\\
\, &= \, \int_N f(x)\chi(nx) d\nu(x)
% \\
\, = \, \chi(n) \int_N f(x)\chi(x)) d\nu(x)
\\
\, &= \, \chi(n) \overline{\mathcal{F}}(f) (\chi)
% \\
\, = \, (\pi_\mu(n) (\overline{\mathcal{F}}(f))) (\chi),
\end{aligned}
$$
for $f \in L^2(N, \nu)$ and $n \in N$.
So, $\pi_\mu$ is equivalent to $\lambda_N$.
\vskip.2cm

\ref{iiDEProp-VNRegular}
By \ref{iDEProp-VNRegular}, the von Neumann algebra $\widetilde \pi_\mu (G)''$
is unitarily equivalent to $\lambda_G(G)''$.
The claim follows now from Proposition~\ref{Pro-IndRepGMS}.
\end{proof}

We can now give sufficient conditions for the regular representation of $G$
to be a factor representation and determine the corresponding type.

\begin{cor}
% 14.D.2
\label{Cor-VNRegular}
Let $G = H \ltimes N$ be as in Proposition~\ref{Prop-VNRegular}
and let $\mu$ a Haar measure on $\widehat N$.
\begin{enumerate}[label=(\arabic*)]
\item\label{iDECor-VNRegular}
If $N$ is compact and $N \ne \{e\}$, then $\lambda_G$ is not a factor representation.
\end{enumerate}
Assume for the following Items that the action $H \curvearrowright (\widehat N, \mu)$
is essentially free and ergodic.
\begin{enumerate}[label=(\arabic*)]
\addtocounter{enumi}{1}
\item\label{iiDECor-VNRegular}
If $N$ is discrete and infinite,
then $\lambda_G$ is a factor representation of type II$_1$.
\item\label{iiiDECor-VNRegular}
If $N$ is neither compact nor discrete and $\mu$ is $H$-invariant, then $\lambda_G$ is a factor representation of type II$_\infty$. 
\item\label{ivDECor-VNRegular}
If $\mu$ is not $H$-invariant, then $\lambda_G$ is a factor representation of type III. 
\end{enumerate}
\end{cor}

\begin{proof}
\ref{iDECor-VNRegular}
Assume that $N$ is compact. Then $\widehat N$ is discrete
and so $\mu$ is the counting measure on $\widehat N$.
In particular, we have $\mu(\{(e\}) > 0$.
Since the action of $H$ on $\widehat N$ fixes $e$, 
it cannot be ergodic (as $N \ne \{e\}$) and so $\lambda_G(G)''$ is not a factor,
by Proposition~\ref{Prop-NecCondROIV}.

\vskip.2cm

\ref{iiDECor-VNRegular}
Assume that $N$ is discrete and infinite. Then $\widehat N$ is compact and infinite.
Therefore $\mu$ is a finite non-atomic measure,
which can be normalized to be a probability measure. 
It follows that the action of $H$ on $\widehat N$
preserves $\mu$, by uniqueness of $\mu$.
The claim follows then from Theorem~\ref{Theo-FacRepSemiDirect}.

\vskip.2cm

\ref{iiiDECor-VNRegular}
Assume that $N$ is neither compact nor discrete and $\mu$ is $H$-invariant.
Then $\widehat N$ is neither discrete nor compact.
Therefore $\mu$ is an infinite non-atomic measure.
The claim follows then from Theorem~\ref{Theo-FacRepSemiDirect}.

\vskip.2cm

\ref{ivDECor-VNRegular}
Assume $\mu$ is not $H$-invariant.
So, the modular function $\Delta \,\colon H \to \R^\times_*$
defined by $h_*(\mu) = \Delta(h) \mu$ is not trivial. 
We claim that there exists no $\sigma$-finite $H$-invariant measure on $\widehat N$
which is equivalent to $\mu$.
Once this is proved, Item \ref{ivDECor-VNRegular} will follow
from Theorem~\ref{Theo-FacRepSemiDirect}.
\par

To prove the claim, assume, by contradiction,
that there exists a $\sigma$-finite $H$-invariant measure $\nu$ on $\widehat N$
which is equivalent to $\mu$.
Then $\nu_H \otimes \nu$ is a 
left (and right) invariant $\sigma$-finite measure
on the LC group $H \ltimes \widehat N$, 
where $\nu_H$ is the counting measure on $H$.
Moreover, $\nu_H \otimes \nu$ is equivalent to $\nu_H \otimes \mu$,
which is itself equivalent to the left Haar measure $(\Delta\nu_H) \otimes \mu$
on $H \ltimes \widehat N$.
Therefore $(\Delta\nu_H) \otimes \mu = f(\nu_H \otimes \nu)$
for a measurable function $f \,\colon \widehat N \to \R^\times_+$.
Since $(\Delta\nu_H) \otimes \mu$ and $\nu_H \otimes \nu$ are both invariant,
$f$ is an invariant function.
It follows that $f$ is constant $\nu_H \otimes \nu$-almost everywhere.
This shows that $\Delta$ is constant 
and hence trivial on $H$ and this is a contradiction. 
\end{proof}

\begin{rem}
% 14.D.3
\label{Rem-Cor-VNRegular}
The condition that $H \curvearrowright (\widehat N, \mu)$ is ergodic is necessary 
for $\lambda_G$ to be a factor representation, by Proposition~\ref{Prop-NecCondROIV}.
However, it can happen that $\lambda_G$ is a factor representation
even when $H \curvearrowright (\widehat N, \mu)$ is not essentially free
(see the example in Remark~\ref{Rem-NecCondROIV}).
\end{rem}

Assuming that $\lambda_G$ is factorial, we ask now:
when is $\lambda_G$ a normal factor representation?

\begin{cor}
% 14.D.4
\label{Cor-CharVNRegular}
Let $G = H \ltimes N$ be as in Proposition~\ref{Prop-VNRegular}
and let $\mu$ be a Haar measure on $\widehat N$.
Assume that $\mu$ is non-atomic
and that $H \curvearrowright (\widehat N, \mu)$ is essentially free and ergodic.
\begin{enumerate}[label=(\arabic*)]
\item\label{iDECor-CharVNRegular}
The regular representation $\lambda_G$ is a normal factor representation
if and only if $\mu$ is $H$-invariant.
\item\label{iiDECor-CharVNRegular}
Assume that $\mu$ is $H$-invariant.
The character $t_G$ of $\lambda_G$ is given by 
$$
t_G(f^* \ast f) \, = \, (f^* \ast f) (e)
\hskip.5cm \text{for} \hskip.2cm
f \in C^c(G) ,
$$
for a convenient normalization of $\mu$.
\end{enumerate}
\end{cor}

\begin{proof}
\ref{iDECor-CharVNRegular}
By Corollary~\ref{Cor-VNRegular}, $\lambda_G$ is a factor representation.
Again, by Corollary~\ref{Cor-VNRegular}, $\lambda_G(G)''$ is semi-finite
if and only if $\mu$ is $H$-invariant.

\vskip.2cm

\ref{iiDECor-CharVNRegular}
Assume that $\mu$ is $H$-invariant.
Since $\mu$ is a Radon measure on $\widehat N$,
it follows from Theorem~\ref{Theo-CharactersSemiDirect}
that $\lambda_G$ is a normal factor representation
and that its character $t_G$ is given by 
$$
t_G(f) \, = \, \int_{\widehat N}{\overline{\mathcal F}}(f_e) (\chi) d\mu(\chi)
\hskip.5cm \text{for} \hskip.2cm
f \in L^1(G)_+ .
$$
We can choose normalizations of $\mu$ and of the Haar measure $\nu$ on $N$
so that $\overline{\mathcal F}$ is a Hilbert space isomorphism between $L^2(N, \nu)$
and $L^2(\widehat N, \mu)$.
Then the inversion formula holds (see Proposition 4 in \cite[II, \S~4]{BTS1--2}):
$$
\varphi(e) \, = \, \int_{\widehat N} {\overline{\mathcal F}}(\varphi) (\chi) d\mu(\chi)
$$
for every $\varphi \in A(N)$, where $A(N)$ is the linear span of 
$$
\{\varphi_1 \ast \varphi_2 \mid \varphi_1, \varphi_2 \in C^c(N)\}.
$$
One checks that $(f^* \ast f)_0 \in A(N)$ for $f \in C^c(G)$ and the claim follows.
\end{proof}

\begin{rem}
% 14.D.5
\label{Rem-Cor-CharVNRegular}
Assume that $G$ is discrete and that $\lambda_G$ is factorial. 
As it has to be the case, the trace $t_G$ from Corollary~\ref{Cor-CharVNRegular} 
coincides with the canonical trace on $\Li (G)$ as in Example~\ref{ExampleTracevN}.
\end{rem}

In order to apply the results from the previous sections, we have to ensure that 
$H \curvearrowright (\widehat N, \mu)$ is essentially free and ergodic.
\par

For $h \in H$, we set 
$$
N_h \, := \, \{(h\cdot n) n^{-1} \mid n \in N\}.
$$
Observe that, since $N$ is abelian, $N_h$ is a subgroup of $N$. 

\begin{prop}
% 14.D.6
\label{Pro-VNRegular-FreeErgodic}
Let $\mu$ be a Haar measure on $\widehat N$. 
\begin{enumerate}[label=(\arabic*)]
\item\label{iDEPro-VNRegular-FreeErgodic}
For $h \in H$, the set $\mathrm{Fix}(h)$ of fixed points of $h$ in $\widehat N$
coincides with the annihilator $N_h^\perp$ of $N_h$ in $\widehat N$.
\item\label{iiDEPro-VNRegular-FreeErgodic}
The action $H \curvearrowright (\widehat N, \mu)$ is essentially free
if and only if $\mu(N_h^\perp) = 0$, for every $h \in H \smallsetminus \{e\}$.
\item\label{iiiDEPro-VNRegular-FreeErgodic}
Assume that $N$ is discrete.
The action $H \curvearrowright (\widehat N, \mu)$ is essentially free
if and only if $N_h$ is infinite, for every $h \in H \smallsetminus \{e\}$. 
\item\label{ivDEPro-VNRegular-FreeErgodic}
Assume that $N$ is discrete.
The action $H \curvearrowright (\widehat N, \mu)$ is ergodic
if and only if the $H$-orbit $\{h \cdot n \mid h \in H\}$
of every $n \in N \smallsetminus \{e\}$ is infinite.
\end{enumerate}
\end{prop}

\begin{proof}
\ref{iDEPro-VNRegular-FreeErgodic}
For $\chi \in \widehat N$, we have 
$$
\begin{aligned}
\chi \in \mathrm{Fix}(h) 
\, & \, \Longleftrightarrow \chi^h \, = \, \chi
\\
\, & \, \Longleftrightarrow \chi(h\cdot n n^{-1}) \, = \, 1 
\hskip.2cm \text {for all} \hskip.2cm n \in N
\\
\, & \, \Longleftrightarrow \chi \in N_h^\perp.
\end{aligned}
$$

\vskip.2cm

\ref{iiDEPro-VNRegular-FreeErgodic}
follows immediately from Item \ref{iDEPro-VNRegular-FreeErgodic}.

\vskip.2cm

\ref{iiiDEPro-VNRegular-FreeErgodic}
Assume that $N$ is discrete and let $h \in H$.
By Pontrjagin duality, the dual group of $N_h$ 
can be identified with $\widehat N / N_h^{\perp}$.
Therefore, in view of \ref{iiDEPro-VNRegular-FreeErgodic},
we have to prove that $\mu(N_h^\perp) > 0$ if and only if 
$N_h^{\perp}$ has finite index in $\widehat N$.
\par

Assume that $\mu(N_h^\perp) > 0$.
Then $N_h^\perp$ is an open subgroup of $\widehat N$.
So, $\widehat N / N_h^{\perp}$ is a discrete group.
Since $\widehat N$ and hence $\widehat N / N_h^{\perp}$ is compact,
it follows that $N_h^{\perp}$ has finite index in $\widehat N$.
\par

Conversely, if $N_h^{\perp}$ has finite index in $\widehat N$,
then $N_h^{\perp}$ is an open subgroup 
of $\widehat N$ and hence $\mu(N_h^\perp) > 0$.

\vskip.2cm

\ref{ivDEPro-VNRegular-FreeErgodic}
Since $N$ is discrete, $\widehat N$ is compact;
hence, the Haar measure $\mu$ on $\widehat N$
can be normalized to be a probability measure.
Moreover, $\mu$ is preserved under the action of $H$. 
\par

Since, by Pontrjagin duality, the dual group of $\widehat N$ is $N$,
the claim follows then from a classical result in ergodic theory
(see \cite[Proposition 1.5]{BeMa--00}).
\end{proof}

\begin{exe}
% 14.D.7
\label{Exa-VNRegular}
Using Corollary~\ref{Cor-VNRegular}
and Proposition~\ref{Pro-VNRegular-FreeErgodic},
we give some examples of groups $G$
for which $\lambda_G$ is a factor of various types.

\vskip.2cm

(1)
Let $H$ be an arbitrary countable infinite group
and let $N = \bigoplus_{h \in H} \Z / 2 \Z$.
Then $H$ acts on $N$ by shifting the coordinates:
$$
h \cdot (\bigoplus_{h' \in H} x_{h'}) \, = \, 
\bigoplus_{h' \in H} x_{h^{-1}h'}
\hskip.5cm \text{for} \hskip.2cm
h \in H, \bigoplus_{h' \in H} x_{h'} \in N.
$$
We identify the dual group $\widehat N$
with the direct product $\prod_{h \in H} \Z / 2 \Z$.
The Haar measure on $\widehat N$ is $\mu = \bigotimes_{h \in H} \nu$,
where $\nu$ is the equidistribution of $\Z / 2 \Z$;
the action of $H$ on $\widehat N$ is given by shifting the coordinates again. 
\par

We claim that $H \curvearrowright (\widehat N, \mu)$ is essentially free.
Indeed, let $h \in H \smallsetminus \{e\}$.
In view of Proposition~\ref{Pro-VNRegular-FreeErgodic},
we have to prove that the subgroup $N_h$ 
appearing there is infinite.
\par

We can choose an infinite subset $L$ of $H$ such that 
$$
h^{-1}L \cap L \, = \, \emptyset.
$$
For every $l \in L$, let $a^{(l)} := \bigoplus_{h' \in H} x_{h'} \in N$ be defined by 
$x_{l} = x_{h^{-1}l} = 1$ and $x_{h'} = 0$ if $h' \notin \{h^{-1}l, l\}$.
Then clearly, 
$$
h \cdot a^{(l)} + a^{(l)} \ne h \cdot a^{(l')} + a^{(l')}
\hskip.5cm \text{for} \hskip.2cm
l \ne l'.
$$
Therefore $N_h$ is infinite.
\par

We claim that $H \curvearrowright (\widehat N, \mu)$ is ergodic.
Indeed, let $\bigoplus_{h' \in H} x_{h'} \in N$, with $x_{h'}$ not all $0$.
Then, obviously, the $H$-orbit of $\bigoplus_{h' \in H} x_{h'}$ is infinite.
The claim follows now from Proposition~\ref{Pro-VNRegular-FreeErgodic}.
\par

Therefore $\lambda_G$ is a factor representation of type II$_1$. 
Equivalently, the crossed product algebra
$L^\infty(\prod_{h \in H} \Z / 2 \Z, \mu) \rtimes H$ is a factor of type II$_1$. 

\vskip.2cm

(2)
The group $H = \SL_2(\Z)$ acts naturally on $N = \Z^2$.
The dual group $\widehat N$ is the torus $\T^2$,
with Lebesgue measure $\mu$ as Haar measure.
\par
 
The action $H \curvearrowright (\widehat N, \mu)$ is essentially free.
Indeed, for every $h \in \SL_2(\Z) \smallsetminus \{I\}$, 
the corresponding subgroup $N_h = \{ (h- I) (x) \mid x \in \Z^2 \}$
as in Proposition~\ref{Pro-VNRegular-FreeErgodic} is clearly infinite.
\par

The action $H \curvearrowright (\widehat N, \mu)$ is ergodic.
Indeed the $\SL_2(\Z)$-orbit of every vector $n \in \Z^2 \smallsetminus \{0\}$ is infinite.
\par

Therefore $\lambda_G$ is a factor representation of type II$_1$. 
Equivalently, the crossed product algebra
$L^\infty(\T^2, \mu) \rtimes \SL_2(\Z)$ is a factor of type II$_1$. 

\vskip.2cm

(3)
Consider the group $H = \SL_2(\Q)$ acting on $N = \R^2$ in the usual way.
Let $\mu$ denote the Haar measure on $\widehat N = \R^2$,
i.e., the Lebesgue measure on $\R^2$.
As in the previous examples,
one checks that the action $H \curvearrowright (\widehat N, \mu)$
is essentially free and ergodic. Since $\mu$ is an infinite measure, 
$\lambda_G$ is a factor representation of type II$_\infty$.

\vskip.2cm

(4)
Take Example (3) with $H = \GL_2(\Q)$ instead of $\SL_2(\Q)$.
Similarly, one shows that $H \curvearrowright (\widehat N, \mu)$
is essentially free and ergodic.
Since the Lebesgue measure $\mu$ on $\R^2$ is not $\GL_2(\Q)$-invariant,
it follows that $\lambda_G$ is a factor representation of type III.
\end{exe}

\begin{rem}
% 14.D.8
\label{Rem-Exa-VNRegular}

(1)
Note that Examples \ref{Exa-VNRegular}, (1) and (2),
are also an illustration of Proposition \ref{iccfactorII1}.

\vskip.2cm

(2)
Example \ref{Exa-VNRegular}(4) is attributed to Godement
(for example in \cite[page 247]{Suth--78}).

\vskip.2cm

(3)
For examples of subtypes of type III (see Remark~\ref{Rem-TypeIIINoInvMeasure}),
consider moreover the triangular group
${\rm T}_2(\Q) = 
\begin{pmatrix}
\Q^\times & \Q\phantom{^\times} \\ 
0\phantom{^\times} & \Q^\times
\end{pmatrix}$,
its subgroup
$$
{\rm ST}_2(\Q) \, = \,
\left\{ \begin{pmatrix}
a & b \\
0 & d
\end{pmatrix}
\in {\rm T}_2(\Q)
\hskip.1cm \Big\vert \hskip.1cm ad = 1 \right\}
$$
and, for every $\lambda \in \mathclose] 0,1 \mathopen[$, 
its subgroup ${\rm S}_\lambda \hskip-.1cm {\rm T}_2(\Q)$
generated by ${\rm ST}_2(\Q) \cup 
\begin{pmatrix}
\sqrt{\lambda} & 0 \\ 0 & \sqrt{\lambda}
\end{pmatrix}$~;
each of these groups acts naturally on $\R^2$. Then,
as shown in \cite{Suth--78}:
\begin{itemize}
\setlength\itemsep{0em}
\item[$\bullet$]
the regular representation of 
${\rm ST}_2(\Q) \ltimes \R^2$
is a factor representation of type II$_\infty$;
\item[$\bullet$]
the regular representation of 
${\rm T}_2(\Q) \ltimes \R^2 $
is a factor representation of type III$_1$;
\item[$\bullet$]
the regular representation of
${\rm S}_\lambda\hskip-.1cm {\rm T}_2(\Q) \ltimes \R^2$
is a factor representation of type III$_\lambda$.
\index{Type II, III! $4$@III$_\lambda$}
\end{itemize}

\vskip.2cm

(4)
In \cite{Blac--77},
there are other constructions of groups $G$ for which
$\lambda_G(G)''$ are factors of type III$_\lambda$,
that are moreover infinite tensor products of factors of type~I.

\vskip.2cm

\index{Crossed product}
(5)
For a countable group $H$ and 
an action of $H$ by homeomorphism on a locally compact space $X$,
there is a related notion of reduced crossed product C*-algebra 
$C^0(X) \rtimes_{\rm red} H$.
When $H$ acts on $N$, and therefore on $\widehat N$, as above,
the restriction of the isomorphism between $L^\infty(\widehat N, \mu) \rtimes H$
and $\Li (G)$ that appears in Proposition~\ref{Prop-VNRegular}
provides an isomorphism from the crossed product 
$C^0(X) \rtimes_{\rm red} H$
to the reduced C*-algebra of $G = H \ltimes N$.
When moreover $G$ is amenable, i.e., when $H$ is amenable,
so that the reduced C*-algebra of $G$ coincides with its maximal C*-algebra,
this provides examples of descriptions of maximal C*-algebras
that can be added to those of Section \ref{C*algLCgroup}.
\end{rem}

\section
{Examples of normal factor representations}
% Section 14.E
\label{Section-ExamplesInfiniteChar}

In Chapter~\ref{ChapterThomadualExamples},
we determined the characters of finite type as well as associated factor representations
for our favorite groups (the two-step nilpotent groups, the affine group, 
the Baumslag--Solitar groups $\BS(1, p)$, 
the lamplighter group,
and the general linear group over an infinite field).
So far, we had no example of normal factor representations of type I$_\infty$ 
and of type II$_\infty$ for any one of our examples. 
\par

Using Theorem~\ref{Theo-CharactersSemiDirect},
we will give an example of a normal factor representation of type I$_\infty$
for the Baumslag--Solitar group $\BS(1, p)$
and of a normal factor representation type of II$_\infty$
for the lamplighter group $\Z \wr (\Z / 2 \Z)$. 

\subsection
{Normal representation of solvable Baumslag--Solitar groups}
% subsection 14.E.a
\label{SS:Subsection-NormalBS}

Let $p$ be a prime.
Let $\Gamma = A \ltimes N \approx \Z \ltimes \Z[1/p]$
be the Baumslag--Solitar group $\BS(1, p)$
as in Sections \ref{Section-IrrRepBS}, \ref{Section-PrimIdealBS}, and \ref{ThomaDualBS}.
Recall that we can identify $\widehat N$ with the $p$-adic solenoid 
$$
\So_p \, = \, (\Q_p \times \R) / \Delta
\hskip.5cm \text{for} \hskip.2cm 
\Delta \, = \, \{(a,-a)\mid a \in \Z[1/p]\};
$$
the action of $A$ on $\widehat N$ is given by the transformation
$T \,\colon \So_p \to \So_p$ induced by the multiplication by $p$.
\par 

The crucial tool in our construction is the following lemma
which shows a peculiarity of the structure of $T$-orbits in $\So_p$.

\begin{lem}
% 14.E.1
\label{Lem-NormalBS}
There exists a point $s \in \So_p$ for which the $T$-orbit
$\{T^n s \mid n \in \Z\}$ is an infinite and \emph{discrete} subset of $\So_p$.
\end{lem}.

\begin{proof}
Let $s = \mathopen[ (1,0) \mathclose] \in \So_p = (\Q_p \times \R) / \Delta$. 
We claim that we have
$$
\lim_{n \to + \infty} T^n s \, = \, \mathopen[ (0,0) \mathclose] 
\hskip.5cm \text{and} \hskip.5cm
\lim_{n \to + \infty} T^{-n} s \, = \, \mathopen[ (0,0) \mathclose] .
$$
Indeed, we have for every $n \in \N$, 
$$
T^n s \, = \, [(T^n1, 0)] \, = \, [(p^n, 0)]
$$
and
$$
T^{-n} s
\, = \, [(T^{-n}1, 0)] 
\, = \, [(p^{-n}, 0)]
\, = \, [(p^{-n}, 0)+ (-p^{-n}, p^{-n})]
\, = \, [ (0, p^{-n})] .
$$
Since $\lim_{n \to +\infty} p^n = 0$ in $\Q_p$ 
and $\lim_{n \to +\infty} p^{-n} = 0$ in $\R$, 
the claim is proved. 
It follows that the $T$-orbit of $s$ is an infinite and discrete subset of $\So_p$.
\end{proof}

\begin{cor}
% 14.E.2
\label{Cor-NormalBS}
Let $s$ be a point in $\So_p$, as in Lemma~\ref{Lem-NormalBS}. 
Let $\mu$ be the measure $\sum_{n \in \Z} \delta_{T^n s}$ on 
$\widehat N = \So_p$. 
\par

Then the associated representation $\widetilde \pi_\mu$
as in Section \ref{Section-IndRepGMS}
is a normal factor representation of type I$_\infty$ of $\Gamma = \BS(1, p)$. 
\end{cor}

\begin{proof}
Indeed, since the $T$-orbit of $s$ is a discrete subset of $\So_p$,
the atomic $T$-invariant measure $\mu = \sum_{n \in \Z} \delta_{T^n s}$
is a Radon measure on $\widehat N = \So_p$. 
Moreover, $\mu$ is infinite, since $s$ is not a $T$-periodic point.
Observe that the action of $A$ on $(\widehat N, \So_p)$ is free
(since $s$ is not periodic).
Therefore by Theorem~\ref{Theo-CharactersSemiDirect},
$\widetilde \pi_\mu$ is a normal factor representation of type I$_\infty$. 
\end{proof}

\begin{rem}
% 14.E.3
\label{Rem-NormIrrBS}
Let $s \in \So_p$ with an infinite $T$-orbit
and let $\chi_s$ the associated unitary character of $N$. 
Recall from Proposition~\ref{Prop-IrredRepBS-bis}
that the induced representation $\Ind_N^\Gamma \chi_s$ corresponding to $s$
is irreducible.
Corollary~\ref{Cor-NormalBS} implies that,
if the $T$-orbit of $s$ is a discrete subset of $\So_p$,
the representation $\Ind_N^\Gamma \chi_s$ is normal.
We suspect that the converse statement is true: 
every normal infinite-dimensional irreducible representation of $\Gamma$ 
is equivalent to a representation of the form $\Ind_N^\Gamma \chi_s$,
where the $T$-orbit of $s$ is an infinite discrete subset of $\So_p$.
\end{rem}

The following corollary shows that the Thoma dual $E(\Gamma)$
does not provide a parametrization of $\Pri(\Gamma)$
in the case of $\Gamma = \BS(1, p)$,
in contrast to the case of nilpotent groups 
(Proposition~\ref{Prop-NilGr} and Corollary~\ref{ThomaHeis-PrimIdeal}),
the affine group over a field 
(Corollary~\ref{ThomaAff-PrimIdeal}), 
or the general linear group over an infinite algebraic extension 
of a finite field (Corollary~\ref{Cor-PrimGLn}).

\begin{cor}
% 14.E.4
\label{Cor-Cor-NormalBS}
Let $\Gamma = \BS(1, p)$ be the Baumslag--Solitar group.
The map $\kappa^{\rm norm}_{\rm prim} \,\colon E(\Gamma) \rightarrow \Pri(\Gamma)$
is neither injective nor surjective.
\end{cor}

\begin{proof}
Let us first show that $\kappa^{\rm norm}_{\rm prim}$ is not injective.
\par

Let $\mu_1, \mu_2$ be two ergodic $T$-invariant non-atomic probability measures on $\So_p$,
with $\mu_1 \ne \mu_2$ and with support $\So_p$ 
%(such measures exist, see Remark~\ref{Rem-Theo-ThomaDualBS}).
(such measures exist, see Remark~\ref{Cor-InvMeasureBS}).
For $i = 1, 2$, denote by $\pi_i$ the factor representation of $\Gamma$
associated to $\Phi(\mu_i) \in E(\Gamma)$,
as in Theorem~\ref{Theo-ThomaDualBS}.
The support $S_i$ of the projection-valued measure
associated to the restriction $\pi_i \vert_N$ 
coincides with the support of $\mu_i$.
Therefore $S_1 = S_2 = \So_2$. 
Therefore, by Theorem~\ref{Theo-PrimIdealBS} (see First Step of the proof), 
$\pi_1$ and $\pi_2$ are weakly equivalent, that is, 
$\kappa^{\rm norm}_{\rm prim}(\mu_1) = \kappa^{\rm norm}_{\rm prim}(\mu_2)$. 
However, we have $\Phi(\mu_1) \ne \Phi(\mu_2)$ since $\mu_1 \ne \mu_2$. 

\vskip.2cm

Next, we show that $\kappa^{\rm norm}_{\rm prim}$ is not surjective.
\par

Let $s \in \So_p$ be a point for which the $T$-orbit
$\{Ts \mid n \in \Z\}$ is an infinite and discrete subset of $\So_p$
(such a point exists by Lemma~\ref{Lem-NormalBS}). 
The induced representation $\pi = \Ind_N^\Gamma \chi_s$ corresponding to $s$
is irreducible (see Proposition \ref{Prop-IrredRepBS}).
\par

We claim that $\textnormal{C*ker}(\pi)$
is not in the range of $\kappa^{\rm norm}_{\rm prim}$.
\par

Assume, by contradiction,
that $\textnormal{C*ker}(\pi)$ lies in the image of $\kappa^{\rm norm}_{\rm prim}$.
Since the $T$-orbit of $s$ is infinite, 
$\pi$ is infinite-dimensional and there exists therefore
$\mu \in {\rm EP}(\So_p)$ such that the factor representation $\pi_\mu$
associated to $\Phi(\mu)$ is weakly equivalent to $\pi$.
\par

On the one hand, the support $S$ 
of the projection-valued measure 
associated to the restriction $\pi \vert_N$
coincides with the closure of the $T$-orbit of $s$. 
On the other hand, the support 
of the projection-valued measure 
associated to the restriction $\pi_\mu \vert_N$ 
coincides with the support of the probability measure $\mu$.
Since $\pi_\mu$ is weakly equivalent to $\pi$, 
it follows that the support of $\mu$ 
coincides with $S$, that is, with the closure of the $T$-orbit of $s$. 
\par

However, since the $T$-orbit of $s$ is a discrete subset of $\So_p$,
there exists a neighbourhood $U$ of $s$ such that $U \cap S = \{s\}$.
Since $s$ belongs to the support $S$ of $\mu$, we have $\mu(U) \ne 0$.
It follows that $\mu(\{s\}) \ne 0$. 
This contradicts the fact that $\mu$ is non-atomic.
\end{proof}

\subsection
[Normal representations of lamplighter group]
{Normal factor representations of type I$_\infty$ for the lamplighter group}
% subsection 14.E.b
\label{SS:Subsection-NormalLamp}

\index{Lamplighter group}
Recall from Sections~\ref{Section-IrrRepLamplighter}, Section~\ref{Section-PrimIdealLamplighter}
and Section~\ref{ThomaDualLamplighter}
that the lamplighter group is the semi-direct product $\Gamma = A \ltimes N$
of $A = \Z$ with $N = \bigoplus_{k \in \Z} \Z / 2 \Z$,
where the action of $\Z$ on $\bigoplus_{k \in \Z} \Z / 2 \Z$
is given by shifting the coordinates.
Recall also that $\widehat N$ can be identified
with $X = \prod_{k \in \Z} \{0,1 \}$, the dual action $\Z$ on $\widehat N$ being 
given by the shift transformation $T$ on $X$.
Moreover, the normalized Haar measure on $\widehat N = X$
is the measure $\mu = \otimes^{\Z} \nu$,
where $\nu$ is the uniform probability measure on $\{0,1 \}$.
\par

The following lemma is similar to Lemma~\ref{Lem-NormalBS}
and will allow us construct a countable infinite number
of normal factor representations of type I$_\infty$
for the lamplighter group $\Gamma$.
\par

We denote by $\mathbf{0}$ the element in $X$ with all coordinates equal to $0$.

\begin{lem}
% 14.E.5
\label{Lem-NormalLamplighter}
Let $(x_n)_{n \in \Z} \in X$ be such that $x_n = 1$ for finitely many $n \in \Z$.
Then the $T$-orbit of $(x_n)_{n \in \Z} \in X$
is an infinite and \emph{discrete} subset of $X$.
\end{lem}

\begin{proof}
Let $x := (x_n)_{n \in \Z} \in X$ be such that $x_n = 1$ for finitely many $n \in \Z$.
It is clear that we have
$$
\lim_{n \to +\infty} T^n x \, = \, \mathbf{0}
\hskip.5cm \text{and} \hskip.5cm
\lim_{n \to +\infty} T^{-n} x \, = \, \mathbf{0}.
$$
Therefore the claim follows.
\end{proof}

As in the case of the Baumslag--Solitar group (Corollary \ref{Cor-NormalBS}), 
Lemma \ref{Lem-NormalLamplighter} has the following consequence.

\begin{cor}
% 14.E.6
\label{Cor-NormalLamplighter}
Let $x$ be a point in $X$ as in Lemma~\ref{Lem-NormalLamplighter}. 
Let $\mu$ be the measure $\sum_{n \in \Z} \delta_{T^n x}$
on $\widehat N = X$. 
\par

Then the associated representation $\widetilde \pi_\mu$
as in Section~\ref{Section-IndRepGMS}
is a normal factor representation of type I$_\infty$
of the lamplighter group.
\end{cor}

The following corollary shows that the Thoma dual $E(\Gamma)$
does not provide a parametrization of $\Pri(\Gamma)$. 
Its proof is similar to the proof of the corresponding result
(Corollary~\ref{Cor-Cor-NormalBS})
for the Baumslag--Solitar group $\BS(1, 2)$.

\begin{cor}
% 14.E.7
\label{Cor-Cor-NormalLamplighter}
Let $\Gamma $ be the lamplighter group.
The map $\kappa^{\rm norm}_{\rm prim} \,\colon E(\Gamma) \rightarrow \Pri(\Gamma)$
is neither injective nor surjective.
\end{cor}

Next, we aim to construct a normal factor representation of type II$_\infty$
for the lamplighter group $\Gamma$.
In view of Theorem \ref{Theo-CharactersSemiDirect}, it suffices to exhibit
an ergodic $T$-invariant Radon measure $\mu$
on an open and $T$-invariant subset $U$ of $X$
which is infinite and non-atomic.

\subsection
{Another construction of infinite invariant measures}
% subsection 14.E.c
\label{SS:SubsectionInvariantMeasuresInduced}

Let $X$ be a locally compact space and $T \,\colon X \to X$ a homeomorphism.
The infinite $T$-invariant measures on $X$
obtained via the method from Section~\ref{Section-TheoSchmidt}
seem to always produce \emph{non-regular} measures.
We will use introduce another well-known construction
(see \cite[\S~1.5]{Aaro--97}) of an infinite invariant measure, 
associated to a given probability measure on a subset $A$ of $X$,
which is invariant for the return map to $A$.

\vskip.2cm 

Let $(X, \mathcal B)$ be a Borel space, $T \,\colon X \to X$ a Borel automorphism,
and $A$ a measurable subset of $X$, fixed throughout this subsection.

\begin{defn}
% 14.E.8
\label{Def-InducedTransformation}
The \textbf{first return time} to $A$ is the map 
$$
r_A \, \colon \, A \to \N^* \cup \{+\infty \}
$$
defined by $r_A(x) =: \min\{n \ge 1 \mid T^n x \in A\}$ if there exists $n \ge 1$
such that $T^n x \in A$ and $r_A(x) = +\infty$ otherwise.
\index{First return time} 
\par

The \textbf{induced transformation} by $T$ to $A$
(also known as the first return transformation or the Poincar\'e map)
is the map $T_A \,\colon A \to A$ defined by 
$$
T_A(x) \, = \, \begin{cases}
\hskip.1cm T^{r_A(x)} (x)
\hskip.2cm \text{if} \hskip.2cm
r_A(x) < +\infty
\\
\hskip.1cm x
\hskip.2cm \text{if} \hskip.2cm
r_A(x) = +\infty.
\end{cases}
$$
\index{Induced transformation}
\end{defn}

We endow $A$ with the induced $\sigma$-algebra 
$$
\mathcal F_A \, = \, \{ A\cap B \mid B \in \mathcal F\}.
$$
Then the induced transformation $T_A$ is measurable; indeed, the sets
$$
r_A^{-1} (\{N\}) \, = \,
A \bigcap_{1 \le n < N} (X \smallsetminus T^{-n}(A)) \cap T^{-N}(A)
$$
for $N \in \N^*$ and 
$$
r_A^{-1} (\{+\infty \}) \, = \, A \bigcap_{n \in \N^*} (X \smallsetminus T^{-n}(A))
$$
are measurable. Therefore 
$$
T_A^{-1}(A\cap B) \, = \, \bigcup_{N\in \N^*\cup \{+\infty \}} r_A^{-1} (\{N\})\cap T^{-N}(A\cap B)
$$
is measurable for every measurable subset $B$ of $X$.
\par

Set 
$$
A_N \, := \, r^{-1}(\{N\})
\hskip.2cm \text{for} \hskip.2cm
N \in \N^*
\hskip.5cm \text{and} \hskip.5cm
A_{\infty} \, := \, r^{-1}(\{+\infty \}) ,
$$
and observe that we have a partition 
$$
A\, = \,A_{\infty} \sqcup \bigsqcup_{N^* \in \N} A_N .
$$

\begin{theorem}
% 14.E.9
\label{Theo-InducedInvMeasure}
Let $A$ be a measurable subset of $X$ and $\mu_A$
$T_A$-invariant probability measure on $(A, \mathcal F)$.
Assume that $\mu_A(A_{\infty}) = 0$.
Let $\mu$ be the measure on $X$ defined by 
$$
\mu(B) \, = \, \sum_{n = 1}^\infty \sum_{k = 0}^{N - 1} \mu_A(T^{-k}(B) \cap A_N)
\, = \, \sum_{N = 1}^\infty \sum_{k = 0}^{N - 1} \mu_A(T^{-k}(B \cap T^k(A_N))),
$$
for every $B \in \mathcal F$.
\begin{enumerate}[label=(\arabic*)]
\item\label{iDETheo-InducedInvMeasure}
The measure $\mu$ on $X$ is $\sigma$-finite and $T$-invariant.
\item\label{iiDETheo-InducedInvMeasure}
The measure $\mu$ is infinite if and only if 
$$
\sum_{N = 1}^{\infty} N \mu_A(A_N) \, = \, +\infty.
$$
\item\label{iiiDETheo-InducedInvMeasure}
If $\mu_A$ is $T_A$-ergodic, then $\mu$ is $T$-ergodic.
\end{enumerate}
\end{theorem}

\begin{proof}
\ref{iDETheo-InducedInvMeasure}
Observe that, by definition of the return map,
the $T^k(A_N)$~'s are pairwise disjoint for $N \in \N^*$ and $k \in \{0, \hdots, N - 1 \}$.
Since $\mu_A$ is a measure on $A$,
it follows from the formula for $\mu$ that $\mu$ is $\sigma$-additive. 
Moreover, we have 
$$
\mu(T^k(A_N)) \, = \, \mu_A(A_N) \, < \, +\infty.
$$
Since 
$$
\mu(X) \, = \, \sum_{N = 1}^\infty \sum_{k = 0}^{N - 1}\mu (T^k(A_N)),
$$
we see that $\mu$ is $\sigma$-finite. 
\par
 
Let us show that $\mu$ is $T$-invariant. Let $B \in \mathcal F$.
We have 
$$
\begin{aligned}
\mu(T^{-1}(B))
\, &= \, \sum_{N = 1}^\infty \sum_{k = 0}^{N - 1} \mu_A(T^{-(k+1))}(B) \cap A_N)
\\
\, &= \, \sum_{N = 1}^\infty \sum_{k = 1}^{N} \mu_A(T^{-k}(B) \cap A_N)
\\
\, &= \, \sum_{N = 1}^\infty \sum_{k = 1}^{N - 1} \mu_A(T^{-k}(B) \cap A_N)
 + \sum_{N = 1}^\infty \mu_A(T^{-N}(B) \cap A_N).
\end{aligned}
$$
For every $N\in \N^*$, we have by definition 
of the return map, $T_A = T^N$ on $A_N$ and hence
$$
T_A^{-1}(B \cap A) \, = \,
(B \cap A_{\infty}) \sqcup \bigsqcup_{N^* \in \N} (T^N(B) \cap A_N),
$$
since $A = A_{\infty} \sqcup \bigsqcup_{N^* \in \N} A_N$.
As $\mu_A(A_{\infty}) = 0$, it follows that
$$
\sum_{N = 1}^\infty \mu_A(T^{-N}(B) \cap A_N) \, = \, \mu_A(T_A^{-1}(B \cap A)).
$$
Using the $T_A$-invariance of $\mu_A$, we have therefore 
$$
\sum_{N = 1}^\infty \mu_A(T^{-N}(B) \cap A_N) \, = \, \mu_A(B \cap A) .
$$
Since 
$$
\mu_A(B \cap A) \, = \, \sum_{N = 1}^\infty \mu_A(B \cap A_N),
$$
we obtain 
$$
\begin{aligned}
\mu(T^{-1}(B))
\, &= \, 
\sum_{N = 1}^\infty \sum_{k = 1}^{N - 1} \mu_A(T^{-k}(B) \cap A_N)+ \mu_A(B \cap A)
\\
\, &= \, \sum_{N = 1}^\infty \sum_{k = 0}^{N - 1} \mu_A(T^{-k}(B) \cap A_N)
\\
\, &= \, \mu(B).
\end{aligned}
$$
This proves the $T$-invariance of $\mu$.
 
\vskip.2cm

\ref{iiDETheo-InducedInvMeasure}
The claim follows from the fact that
$$
\mu(X)
\, = \, \sum_{N = 1}^\infty \sum_{k = 0}^{N - 1} \mu_A( A_N)
\, = \, \sum_{N = 1}^\infty N \mu_A( A_N).
$$

\vskip.2cm

\ref{iiiDETheo-InducedInvMeasure}
Let $B$ be a $T$-invariant measurable subset of $X$ such that $\mu(B) > 0$. 
We have to show that $\mu(X \smallsetminus B) = 0$. 
\par

It follows from the definition of $\mu$
that there exists $N \in \N^*$ and $k \in \{0, \hdots, N - 1 \}$ such that 
$\mu_A(T^{-k}(B) \cap A_N) > 0$.
By $T$-invariance of $B$, we have therefore $\mu_A(B \cap A) > 0$.
Since $T_A = T^N$ on $A_N$, we have, again by $T$-invariance of $B$, 
$$
\begin{aligned}
T_A^{-1}(B \cap A)
\, &= \, (B \cap A_{\infty}) \cup \bigcup_{N^*\in \N} (T^{-N} (B) \cap A_N)
\\
\, &= \, (B \cap A_{\infty})\cup\bigcup_{N^*\in \N} (B \cap A_N)
\\
\, &= \, B \cap A.
\end{aligned}
$$
So, $B \cap A$ is $T_A$-invariant. Since $ \mu_A(B \cap A) > 0$,
it follows from the $T_A$-ergodicity of $\mu_A$
that $\mu_A((X \smallsetminus B) \cap A) = 0$.
Since $X \smallsetminus B$ is $T$-invariant, we have therefore
$$
\mu(X \smallsetminus B)
\, = \, \sum_{N = 1}^\infty \sum_{k = 0}^{N - 1} \mu_A((X \smallsetminus B) \cap A_N)
\, = \, 0.
$$
\end{proof}

\subsection
[Further normal representations of the lamplighter group]
{Normal factor representations of type II$_\infty$ for the lamplighter group}
% subsection 14.E.d

Let $X = \prod_{k \in \Z} \{0,1 \}$ and $T$ the shift transformation on $X$,
as in Subsection \ref{SS:Subsection-NormalLamp}.
In order to construct normal factor representation of type II$_\infty$ for the lamplighter group,
we have to find appropriate $T$-invariant infinite measures on $X$
(see Theorem~\ref{Theo-CharactersSemiDirect}). 
For this, we will use Theorem~\ref{Theo-InducedInvMeasure}.
As an intermediate step, we will first consider shift spaces over an infinite alphabet.

\vskip.2cm

Let $S$ be a countable set (a ``state" space or ``alphabet") and consider the product space
$$
\Sigma_S \, := \, S^\Z \, = \, \left\{ (s_n)_{n \in \Z} \mid s_n \in S \right\},
$$
endowed with the shift transformation $\sigma \,\colon \Sigma_S \to \Sigma_S$, given by 
$$
\sigma ((s_n)_{ n \in \Z}) \, = \,
(s_{n+1})_{n \in \Z}
\hskip.5cm \text{for} \hskip.2cm
(s_n)_{n \in \Z}\in \Sigma_S.
$$
Let $\mathcal F$ be the $\sigma$-algebra on $\Sigma_S$ generated by the cylinders
$$
C_{a_1, \hdots, a_k}^{n_1, \hdots, n_k} \, = \,
\left\{ (s_n)_{n \in \Z} \in \Sigma_S \mid s_{n_i} = a_i
\hskip.2cm \text{for} \hskip.2cm
i = 1, \hdots, k \right\},
$$
for integers $n_1 < n_2 < \cdots < n_k$ in $\Z$
and (not necessarily distinct) elements $a_1, a_2, \hdots, a_k$ in $S$.
\par

Let $\mathbf{p}$ be a probability measure on $S$. The product measure 
$$
\nu_{\mathbf{p}} \, := \, \mathbf{p}^{\otimes \N}
$$
is the unique probability measure on $(\Sigma_S, \mathcal F)$ such that,
for every cylinder $C_{s_1, \hdots, s_k}^{n_1, \hdots, n_k}$, we have
$$
\nu_{\mathbf{p}} (C_{a_1, \hdots, a_k}^{n_1, \hdots, n_k}) \, = \,
p_{a_1} p_{a_2} \cdots p_{a_k},
$$
with $p_{a} = \mathbf{p}(\{a\})$.
One checks that the shift transformation
$$
\sigma \, \colon \, \Sigma_S \to \Sigma_S
$$
is measurable with respect to $\mathcal F$. Since
$$
\sigma^{-1} (C_{a_1, \hdots, a_k}^{n_1, \hdots, n_k}) \, = \, 
C_{a_1 \hdots, a_k}^{n_1+1, n_1, \hdots, n_k+1},
$$
we have
$$
\nu_{\mathbf{p}}(\sigma^{-1} (C_{a_1, \hdots, a_k}^{n_1, \hdots, n_k}))
\, = \, p_{a_1} p_{a_2} \cdots p_{a_k}
\, = \, \nu_{\mathbf{p}} (C_{a_1, \hdots, a_k}^{n_1, \hdots, n_k}),
$$
and it follows that $\nu_{\mathbf{p}}$ is $\sigma$-invariant.
\par

For a proof of the following classical fact,
see e.g.\ Theorem 1.12 in \cite{Walt--82};
the proof given there for a finite set $S$ carries over without change
to a countable infinite set $S$.

\begin{prop}
% 14.E.10
\label{Pro-ErgodicShift}
Let $S$ be a countable set and $\mathbf{p}$ a probability measure on $S$. 
\par

The product measure $\nu_{\mathbf{p}}$ on $\Sigma_S$ is ergodic
under the shift transformation $\sigma$.
\hfill $\square$
\end{prop} 

Coming back to our space $X = \prod_{k \in \Z} \{0,1 \}$, we set
$$
U \, := \, X \smallsetminus \{\mathbf{0}\}.
$$

\begin{theorem}
% 14.E.11
\label{Theo-NormalLamplighter-bis}
There exists an ergodic $T$-invariant Radon measure $\mu$
on the open and $T$-invariant subset $U$ of $X$
which is infinite and non-atomic.
\end{theorem}

\begin{proof}
Consider the open and closed subset 
$$
A \, := \, \left\{ (x_n)_{n \in \Z}\mid x_0 \, = \, 1 \right\}
$$
of $U$. 
For $S := \N$, we define a map $\Phi \,\colon \Sigma_S \to A$ as follows:
for $\mathbf{n} = (n_k)_{k \in \Z} \in \Sigma_S$,
we set $\Phi(\mathbf{n}) = (x_n)_{n \in \Z}\in A$, where 
$$
x_n \, = \, \begin{cases} 
&1 \hskip.2cm \text{if} \hskip.2cm
n = n_0 + \cdots + n_k + k + 1
\hskip.1cm \text{for some} \hskip.1cm
k \in \N
\\
&1 \hskip.2cm \text{if} \hskip.2cm
n = -(n_{-1} + \cdots + n_{-k} + k)
\hskip.1cm \text{for some} \hskip.1cm
k \in \N^*
\\
& 0 \hskip.2cm \text{otherwise};
\end{cases}
$$
thus, 
$$
\Phi(\mathbf{n}) \, = \,
(\hdots, 1, [n_{-k}], 1, \hdots, 1, [n_{-1}], \underset{*}{1},
[n_0], 1, [n_1], 1, \hdots, 1,[n_k], 1, \hdots),
$$
where $[n] = \underbrace{(0, \hdots, 0)}_{n-\text{times}}$ for $n \in \N$
and $*$ denotes the $0$-th position in the sequence.

\vskip.2cm

$\bullet$ {\it First step.}
We claim that $\Phi$ is measurable and that 
$\Phi \circ\sigma= T_A \circ \Phi$.
\par

Indeed, to show that $\Phi$ is measurable, it suffices to show that 
$\Phi^{-1} (C^{N}_\varepsilon)$ is a measurable subset of $\Sigma_S$
for every $N\in \Z$, where
$$
C^{N}_\varepsilon \, = \, \{ (x_n)_{n \in \Z} \mid x_{N} = \varepsilon \}
$$
and $\varepsilon \in \{0,1 \}$.
This is the case for $\Phi^{-1} (C^{N}_1)$ for $N \ge 0$, since
$$
\Phi^{-1} (C^{N}_1) \, = \,
\bigcup_{0 \le k \le N} \bigcup_{(n_0, n_1, \hdots, n_k) \in S_{N, k}}
C^{0, 1, \hdots, k}_{n_0, n_1, \hdots, n_k},
$$
where $S_{N, k}$ is the finite subset of $\{0, \hdots, N\}^k$
consisting of all $k$-tuples $(n_0, n_1, \hdots, n_k)$
such that $n_0+n_2+ \cdots n_k + k+1 = N$.
The case $N < 0$ is treated similarly.
It follows that we also have $\Phi^{-1} (C^{N}_0) \in \mathcal F$,
since $C^{N}_0 = A \smallsetminus C^{N}_1$.
\par

Let $\mathbf{n} = (n_k)_{k \in \Z} \in \Sigma_S$. Then 
$$
\begin{aligned}
T_A(\Phi(\mathbf{n}))
\, &= \, ( \hdots, 1 ,[n_{-1}], 1,[n_0], \underset{*}{1}, [n_1], 1,[n_2], 1 \hdots 
\\
\, &= \, \Phi((n_{k+1})_{k \in \Z}) = \Phi(\sigma(\mathbf{n}).
\end{aligned}
$$
This proves the claim. 

\vskip.2cm

Let $\mathbf{p}$ be a probability measure on $\N$ 
with $p_n := \mathbf{p}(n) < 1$ for every $n \in \N$. 
We will make below a more precise choice of $\mathbf{p}$.
\par

Let 
$$
\mu_{A} \, := \, \Phi_*(\mathbf{p})
$$
be the image of $\mathbf{p}$ under $\Phi$. 

\vskip.2cm

$\bullet$ {\it Second step.}
We claim that $\mu_{A}$ is a $T_A$-invariant 
non-atomic probability measure on $A$ which is ergodic.
\par

Indeed, $\mathbf{p}$ is non-atomic: since $\sum_{n \in \N} p_n = 1$ and $p_n < 1$,
there exists $0 \le q < 1$ such that $p_n \le q$ for every $n \in \N$
and hence 
$$
\mathbf{p}(\{(n_k)_{k \in \Z}\})
\, = \, \lim_{k \to + \infty}\prod_{i = -k}^k p_{n_i} \le \lim_{k \to + \infty} q^{2k+1}
\, = \, 0.
$$
Moreover, $\mathbf{p}$ is ergodic by Proposition \ref{Pro-ErgodicShift}).
Since $\Phi \circ \sigma = T_A \circ \Phi$ by the first step, the claim follows.

\vskip.2cm

$\bullet$ {\it Third step.}
We claim that $\mu_{A}(r_A^{-1} (\{+\infty \})) = 0$ 
and that $\mu_{A}(r_A^{-1} (\{N\}) = p_N$ for every $N \in \N^*$.
\par

Indeed, we have
$$
r_A^{-1} (\{+\infty \})) \, = \,
\left\{
(x_n)_{n \in \Z} \in A \mid x_n = 0
\hskip.2cm \text{for all} \hskip.2cm
n \ge 1
\right\};
$$
hence $\Phi^{-1}(r_A^{-1} (\{+\infty \}) = \emptyset$ and so 
$$
\mu_{A}(r_A^{-1} (\{+\infty \})) \, = \, \mathbf{p}(\Phi^{-1}(r_A^{-1} (\{+\infty \})) \, = \, 0.
$$
For $N \in \N^*$, we have, by definition of $\Phi$,
$$
\Phi^{-1}(r_A^{-1} (\{N\})) \, = \, \left\{ (n_k)_{k \in \Z} \in \Sigma_S \mid n_{0} = N \right\}
$$
and hence
$$
\mu_{A}(r_A^{-1} (\{N\})) \, = \, \mathbf{p}(\Phi^{-1}(r_A^{-1} (\{N\})) \, = \, p_N.
$$

\vskip.2cm

Let $\mu$ be the measure on $U$ associated to
the $T_A$-invariant probability measure $\mu_{A}$
as defined in Theorem~\ref{Theo-InducedInvMeasure}.

\vskip.2cm

$\bullet$ {\it Fourth step.}
We claim that $\mu$ is an ergodic $T$-invariant Borel measure on $U$
which is non-atomic and $\sigma$-finite.
\par

Indeed, the fact that $\mu$ is $T$-invariant, ergodic, and $\sigma$-finite
follows from Theorem \ref{Theo-InducedInvMeasure} and the second step
(observe that $\mu_{A}(r_A^{-1} (\{+\infty \})) = 0$, by the third step).
Since $\mu_A$ is non-atomic,
the formula defining $\mu$ shows that $\mu$ is non-atomic.

\vskip.2cm

Choose the probability measure $\mathbf{p}$ on $\N$ so that
$$
\sum_{n \in \N} np_n \, = \, +\infty.
$$

\vskip.2cm 

$\bullet$ {\it Fifth step.}
We claim that the associated measure $\mu$ on $U$ is infinite.
\par

Indeed, this follows from Theorem~\ref{Theo-InducedInvMeasure}.ii and the fact that 
$\mu_{A}(r_A^{-1} (\{N\}) = p_N$ for every $N \in \N^*$ (see the third step). 

\vskip.2cm 

Assume now that the probability measure $\mathbf{p}$ on $\N$
is chosen so that, moreover,
$$
\sum_{n \in \N} \sqrt{n}p_n \, < \, +\infty.
$$
(To fulfill all requirements, we could take
$p_n = \dfrac{6}{\pi^2 (n+1)^2}$ for every $n \in \N$.)

\vskip.2cm

$\bullet$ {\it Sixth step.}
We claim that $\mu$ is a Radon measure on $U$.
\par

Indeed, since $U$ is a second-countable locally compact space,
it suffices to show that $\mu$ is finite on every compact subset of $U$
(see \ref{regmeas2ndc}).
 \par
 
Let $B$ be a compact subset of $U = X \smallsetminus \{\mathbf{0}\}$. 
Then there exists an open neighbourhood $V$ of $\mathbf{0}$ in $X$
such that $B \cap V = \emptyset$. Choose $N_0 \in \N$ so that the cylinder
$$
C \, := \,
%C^{-N_0, \hdots, -1, 0,1, \hdots, N_0}_{1, \hdots, 1, 1,0, \hdots, 0}=
\left\{
(x_n)_{n \in \Z} \in X \mid x_n = 0
\hskip.5cm \text{for all} \hskip.2cm
n
\hskip.2cm \text{such that} \hskip.2cm
-N_0 \le n \le N_0 \right\}
$$
is contained in $V$. 
\par

Let $N > N_0^2$. We claim that 
$$
T^{-k}(B) \cap A_N \, = \, \emptyset
\hskip.5cm \text{for all} \hskip.5cm
\sqrt{N} \le k \le N - \sqrt{N}.
$$
Indeed, let $x = (x_n)_{n \in \Z} \in B$.
Since $B \cap C = \emptyset$, there exists $i \in \{-N_0, \hdots, N_0\}$
such that $x_i = 1$. 
Let $k \in \N$ be such that $\sqrt{N} \le k \le N - \sqrt{N}$.
Then 
$$
i+k \le N_0 + N - \sqrt{N} < N
\hskip.5cm \text{and} \hskip.5cm
i+k \ge -N_0+\sqrt{N} > 0,
$$
that is, $0 < i+k < N$.
For $T^{-k}(x) = (y_n)_{n \in \Z}$, we have $y_n = x_{n-k}$ and hence $y_{i+k} = 1$.
Since $0 < i+k < N$, this implies that $T^{-k}(x)$ does not belong to 
$A_N$ and the claim is proved.
\par

It follows that 
$$
\begin{aligned}
\mu(B)
\, &= \, \sum_{N = 1}^\infty \sum_{k = 0}^{N - 1} \mu_A(T^{-k}(B) \cap A_N)
\\
\, &= \, \sum_{N = 1}^{N_0^2} \sum_{k = 0}^{N - 1} \mu_A(T^{-k}(B) \cap A_N)
\\
\, &+ \, \sum_{N > N_0^2}\left(
\sum_{k = 0}^{\sqrt{N} - 1} \mu_A(T^{-k}(B) \cap A_N)
\, + \, \sum_{k = N - \sqrt{N}}^{N - 1} \mu_A(T^{-k}(B) \cap A_N)
\right)
\\
\, &= \, \le \sum_{N = 1}^{N_0^2} \sum_{k = 0}^{N - 1} \mu_A( A_N)
\, + \, \sum_{N > N_0^2} \left(
\sum_{k = 0}^{\sqrt{N}-1} \mu_A( A_N) + \sum_{k = N - \sqrt{N}}^{N - 1} \mu_A( A_N)
\right)
\\
\, &= \, \sum_{N = 1}^{N_0^2} N \mu_A( A_N)
\, + \, \sum_{N > N_0^2} 2 \sqrt{N} \mu_A( A_N).
\end{aligned}
$$
Since, by assumption, the series $\sum_{n \in \N} \sqrt{n} p_n$ converges,
we see that $\mu(B) < +\infty$.
\end{proof}

The following corollary is a direct consequence of 
Theorem \ref{Theo-NormalLamplighter-bis}
and Theorem~\ref{Theo-CharactersSemiDirect}.

\begin{cor}
% 14.E.12
\label{Cor-NormalLamplighter2}
Let $\Gamma = A \rtimes N$ be the lamplighter group,
with $A = \Z$ and $N = \bigoplus_{k \in \Z} \Z / 2 \Z$.
Identify the dual $\widehat N$ of $N$ with $X = \prod_{k \in \Z} \Z / 2 \Z$
and the action of the generator $1$ of $\Z$ with the two-sided shift $T$ of $X$.
\par

Let $\mu$ be an ergodic $T$-invariant Radon measure on $X$ 
which is infinite and non-atomic, 
as in Theorem \ref{Theo-NormalLamplighter-bis}.
Let $\widetilde \pi_\mu$ be the corresponding representation of $\Gamma$,
as in Construction \ref{constructionInd+mu}.
\par

Then $\widetilde \pi_\mu$ is a normal factor representations of $\Gamma$
of type II$_\infty$.
\end{cor}

%-----------------------------------------------------------------------
% End of chapter 14
%-----------------------------------------------------------------------
\chapter{Separating families of finite type representations}
% Chapter 15
\label{Chap-Applications}

\emph{
Finite type representations of a topological group $G$
have been introduced in Section~\ref{Section-GNS-Traces}
and several examples of such representations were given
in Chapter~\ref{ChapterThomadualExamples},
mainly in the case where $G$ is discrete.
}
\par

\emph{
In this chapter, we will deal with the problem whether $G$ has ``enough" finite type representations;
more precisely, we study the following two questions:
\begin{itemize}
\setlength\itemsep{0em}
\item[$\bullet$]
when does the family of finite type representations of $G$ 
separate the points of $G$?
\item[$\bullet$]
when does the family of finite type representations of $G$ 
separate the points of the maximal C*-algebra $C^*_{\rm max}(G)$ of $G$,
when $G$ is a discrete group?
\end{itemize}
}
\par

\emph{
We will see in Section~\ref{S:FinTypRep-Metric}
that a faithful finite type representation of $G$ gives rise to a bi-invariant metric on $G$. 
Among the class of compactly generated locally compact groups,
the groups with a point separating family of finite type representations
will be characterized in Section~\ref{S:FaithfulFinTypRep} as being the SIN groups.
The class of SIN groups coincides also with the class of locally compact groups
for which the von Neumann algebra generated by the regular representation is finite
(Section~\ref{S:LC-FiniteVN}).
It follows that a connected locally compact group $G$
has a point separating family of finite type representations
if and only if $G$ is the product of a compact group with $\R^n$.
Moreover, the class of compactly generated totally disconnected locally compact groups
with a point separating family of finite type representations
agrees with the class of projective limits of discrete groups.}
\par

\emph{
Concerning the second question,
which is a much deeper problem, only partial results are available.
We show in Section \ref{SS:Faithful}
that the answer to this question is negative for $G = \SL_n(\Q)$ for $n \ge 2$.
In contrast to this, the answer is positive for amenable groups as well as for free groups.
In fact, we will prove that free groups satisfy even two stronger properties,
which we call the C*-MAP property and the C*-residually finiteness property;
this last property means that the representations with finite image
separate the points of the maximal C*-algebra of a free group.
}

\emph{In Section~\ref{Section:IRS}, we show that every so-called random invariant subgroup
of a discrete group $G$ gives rive to a trace of $G$.
As a consequence, if $G$ is character rigid in the sense
that $E(G)$ consists only of $1_G$ and $\delta_e$, then every action of $G$ by measure preserving transformations of a probability space is essentially free.}

\section{Finite type representations and bi-invariant metrics}
% 15.A
\label{S:FinTypRep-Metric}

Let $G$ be a topological group.
Recall from Definition~\ref{Def-Trace-FiniteRep}
that a representation $\pi$ of $G$ is of finite type
if the von Neumann algebra $\pi(G)''$ has a faithful finite normal trace.
Denote by ${\rm Rep}_{\rm FT}(G)$
the family of finite type representations of $G$, up to quasi-equivalence. 
\par

We show that a finite type representation of $G$
defines a bi-invariant metric on an appropriate quotient of $G$.
Recall that a \emph{bi-invariant pseudo-metric} on a group $G$
is a pseudo-metric $d$ with the property 
$$
d(gx, gy) \, = \, d(xg, yg) \, = \, d(x, y)
\hskip.5cm \text{for all} \hskip.2cm
x,y, g \in G .
$$

\begin{prop}
% 15.A.1
\label{Pro-BiInvMetric}
Let $G$ be a topological group and $\pi$ a representation of $G$.
Let $\tau$ be a normalized finite normal trace on $\pi(G)''$, 
let $\varphi = \tau \circ \pi$ be the normalized trace on $G$ associated to $\pi$,
and define
$$
d_{(\pi, \tau)} \, \colon \, G \times G \to \R,
\hskip.5cm
(x,y) \mapsto \tau \left( (\pi(x)-\pi(y))^*(\pi(x)-\pi(y)) \right)^{1/2}.
$$
\begin{enumerate}[label=(\arabic*)]
\item\label{iDEPro-BiInvMetric}
For every $x, y \in G$, we have 
$$
d_{(\pi, \tau)} \, = \, 2 \left(1-{\rm Re}\,\varphi(x^{-1} y)\right)^{1/2}.
$$
\item\label{iiDEPro-BiInvMetric}
$d_{(\pi, \tau)}$ is a bi-invariant continuous pseudo-metric on $G$.
\item\label{iiiDEPro-BiInvMetric}
Assume that $\tau$ is faithful.
Then $d_{(\pi, \tau)}$ factorizes to a bi-invariant metric on $G/\ker \pi$.
\end{enumerate}
\end{prop}

\begin{proof}
Claim \ref{iDEPro-BiInvMetric} follows from a straightforward computation.
Let $(\pi', \Hi', \xi')$ be a GNS triple associated to
the normalized function of positive type $\varphi = \tau \circ \pi$, so that
$$
\tau(\pi(g)) \, = \, \varphi(g) \, = \, \langle \pi'(g)\xi' \mid \xi' \rangle
\hskip.5cm \text{for all} \hskip.2cm
g \in G.
$$
We have
$$
d_{(\pi, \tau)}(x,y) \, = \,
2 \left( 1 - {\rm Re} \hskip.1cm \varphi(x^{-1} y) \right)^{1/2} \, = \,
\Vert \pi'(x) \xi' - \pi'(y) \xi' \Vert 
$$
for all $x, y \in G$.
The fact that $d_{(\pi, \tau)}$ is a left invariant pseudo-metric
follows immediately from this equality. 
Since $\varphi$ is conjugation invariant and continuous,
$d_{(\pi, \tau)}$ is right invariant and continuous.
\par

Assume now that $\tau$ is faithful. 
Then $d_{(\pi, \tau)}(x,y) = 0$ if and only if $\pi(x) = \pi(y)$, that is, $x^{-1} y \in \ker \pi$.
This shows that $d_{(\pi, \tau)}$ factorizes to a metric on $G/\ker \pi$.
\end{proof}

Recall that a family $\mathcal{F}$ of representations of $G$ 
\emph{separates the points of} $G$, or is \emph{point separating},
if, for every $g \in G$ such that $g \ne e$,
there exists $\pi \in \mathcal{F}$ with $\pi(g) \ne I$.
\par

The following corollary is an immediate consequence
of Proposition~\ref{Pro-BiInvMetric}.

\begin{cor}
% 15.A.2
\label{CorPro-BiInvMetric}
Let $G$ be a topological group.
Assume that ${\rm Rep}_{\rm FT}(G)$ contains a point separating sequence $(\pi_n)_{n \ge 1}$.
For every $n \ge 1$, let $\tau_n$ be a normalized faithful finite normal trace on $\pi_n(G)''$. 
Then 
$$
d \, := \, \sum_{n \ge 1} \frac{1}{2^n} d_{(\pi_n, \tau_n)}
$$
is a continuous bi-invariant metric on $G$.
\end{cor}

\begin{rem}
% 15.A.3
\label{Rem-CorPro-BiInvMetric}
Groups as in Corollary~\ref{CorPro-BiInvMetric}
are related to a class of topological groups introduced by Popa in \cite{Popa--07}
and defined as follows. 
\par

\index{Finite type, for a Polish group}
A Polish group $G$ is said to be of \textbf{finite type} if $G$ is isomorphic, as topological group,
to a closed subgroup of the unitary group $\U(\mathcal M)$
of a finite von Neumann algebra $\mathcal M$ acting on a separable Hilbert space,
where $\U(\mathcal M)$ is equipped with the strong operator topology.
\par

It is shown in \cite[Theorem 2.7]{AnMa--12} that a Polish $G$ is of finite type
if and only if ${\rm Rep}_{\rm FT}(G)$ contains a point separating sequence
and the topology of $G$ is generated by the metric $d$
as in Corollary~\ref{CorPro-BiInvMetric}.
\end{rem}

\section{Point separating finite type representations}
% 15.B
\label{S:FaithfulFinTypRep}

We characterize below compactly generated locally compact groups $G$
for which ${\rm Rep}_{\rm FT}(G)$ is point separating.

\begin{defn}
% 15.B.1
\label{defSIN}
\index{SIN group}
A topological group $G$ is a \textbf{SIN group}
if $G$ has a basis of neighbourhoods of $e$
which are invariant under conjugation, that is, 
of neighbourhoods $V$ of $e$ such that $gVg^{-1} = V$ for all $g \in G$.
\end{defn}

The class of SIN topological groups is closed under direct products,
and includes discrete groups, abelian LC groups, and compact groups.

\begin{lem}
% 15.B.2
\label{Lem-HaarSIN}
Let $G$ be a locally compact SIN group.
\par

Then $G$ is unimodular.
\end{lem}

\begin{proof}
Let $\mu_G$ be a left Haar measure on $G$
and let $V$ be a conjugation invariant compact neighbourhood of $e$.
For every $g \in G$, we have, on the one hand,
$\mu_G (g^{-1}Vg) = \mu(V) > 0$
and, on the other hand,
$\mu_G (g^{-1}Vg) = \Delta(g) \mu(V)$,
where $\Delta_G$ is the modular function of $G$ (see Appendix \ref{AppLCG}).
So, $\Delta(g) = 1$.
\end{proof}

The following result is Th\'eor\`eme 10 in \cite{Gode--54}.

\begin{theorem}
% 15.B.3
\label{Theo-SINSepTrace}
Let $G$ be a locally compact SIN group.
\par

Then ${\rm Rep}_{\rm FT}(G)$ separates the points of $G$.
\end{theorem}

\begin{proof}
Let $\mu_G$ be a left Haar measure on $G$.
By Lemma~\ref{Lem-HaarSIN}, $\mu_G$ is also a right Haar measure.
\par

Let $g \in G$ such that $g \ne e$.
We can find a compact neighbourhood $V$ of $e$
which is invariant under conjugation 
and such that $V = V^{-1}$ and $g \notin V^2$.
Set 
$$
f_V \, := \, \frac{\Un_V}{\sqrt{\mu_G(V)}} \in L^2(G,\mu_G)
$$
and let $\varphi_V$ be the normalized continuous function of positive type on $G$ defined by 
$$
\varphi_V(x) \, = \, \langle \lambda_G(x) f_V \mid f_V \rangle \, = \, \mu_G(V)^{-1} \mu_G (x V \cap V)
\hskip.5cm \text{for every} \hskip.2cm
x \in G.
$$
Then, for every $x, y \in G$, we have
$$
\begin{aligned}
\varphi_V(yxy^{-1})
\, &= \,\mu_G(V)^{-1} \mu_G(yxy^{-1} V \cap V)
\, = \, \mu_G(V)^{-1} \mu_G(xV \cap y^{-1}Vy)
\\
\, &= \, \mu_G(V)^{-1} \mu_G(xV \cap V)
\, = \, \varphi_V(x),
\end{aligned}
$$
since $\mu_G$ and $V$ are invariant under conjugation by $y$. 
So, $\varphi_V \in \Tr_1(G)$.
Moreover, since $g \notin V^2 = V^{-1}V$, we have
$gV \cap V = \emptyset$ and hence $\varphi_V(g) = 0$.
\par

Let $\pi_V$ be the cyclic subrepresentation of $\lambda_G$ defined
on the closure of the linear span of $\{ \lambda_G(x)f_V \mid x \in G \}$.
Then $\pi_V(g) \ne 1$ by Proposition~\ref{Pro-BiInvMetric} or Lemma~\ref{Lem-CharHeis}.
\end{proof}

We now show that the converse of Theorem~\ref{Theo-SINSepTrace} is true,
under the additional assumption that $G$ is compactly generated.
This result appears in \cite[Proposition 17.3.5]{Dixm--C*}
and also, in a related version, in \cite[Th\'eor\`eme 6]{Gode--51b}.

\begin{theorem}
% 15.B.4
\label{Theo-SepTraceSIN}
Let $G$ be a compactly generated locally compact group. 
Assume that ${\rm Rep}_{\rm FT}(G)$ separates the points of $G$.
\par

Then $G$ is a SIN group.
\end{theorem}

\begin{proof}
$\bullet$ \emph{First step.}
Let $C$ be a compact subset of $G$ with $e \notin C$.
We claim that there exists a neighbourhood $U$ of $e$
which is disjoint from $C$ and invariant under conjugation.
\par

Indeed, let $g \in C$.
Since $g \ne e$, there exists $\pi \in {\rm Rep}_{\rm FT}(G)$ such that $g \notin \ker \pi$.
Let $\tau$ be a normalized faithful finite normal trace on $\pi(G)''$
and $d = d_{(\pi, \tau)}$ the associated pseudo-metric on $G$
as in Proposition~\ref{Pro-BiInvMetric}.
Since $d$ is continuous, since $d(g,e) > 0$, and since $d$ is bi-invariant, the set 
$$
V_g \, = \, \{ x \in G \mid d(x,e) \ge d(g,e)/2 \}
$$
is a closed neighbourhood of $g$ such that $e \notin V_g$
and $V_g$ is invariant under conjugation.
\par

By compactness of $C$, there exists finitely many elements $g_1, \dots, g_n$ of $C$
such that $C \subset \cup_{i = 1}^n V_{g_i}$.
The set 
$$
U \, := \, G \smallsetminus \bigcup _{i = 1}^n V_{g_i}
$$
is a neighbourhood of $e$ with the required properties.

\vskip.2cm

$\bullet$ \emph{Second step.}
We claim that $G$ is a SIN group.
\par

Indeed, let $K$ be a compact generating set for $G$
and let $V$ be an open neighbourhood of $e$. 
Upon replacing $K$ by $(K \cup \overline{V}) \cup (K \cup \overline{V})^{-1}$,
we can assume that $V \subset K$ and that $K = K^{-1}$.
\par

Set $C := K^3 \cap (G \smallsetminus V)$.
Then $C$ is a compact subset of $G$ with $e \notin C$.
Therefore, by the first step,
there exists a conjugation invariant neighbourhood $U$ of $e$
such that $U \cap C = \emptyset$.
\par

Let $g \in K$ and $x \in V\cap U$.
Then $gxg^{-1} \in K^3$; moreover, $gxg^{-1} \in U$ and therefore $gxg^{-1} \notin C$.
Hence, $gxg^{-1} \in V \cap U$. 
This shows that $V \cap U$ is invariant under conjugation by every element of $K$.
Since $K$ generates $G$,
it follows that $V \cap U$ is invariant under conjugation by every element of $G$.
As $V \cap U$ is a neighbourhood of $e$ contained in $V$, the claim is proved.
\end{proof}

When $G$ is a MAP group (see Subsection~\ref{SS:MAP}),
it is obvious that ${\rm Rep}_{\rm FT}(G)$ separates the points of $G$.
The following corollary is therefore an immediate consequence
of Theorem~\ref{Theo-SepTraceSIN}.

\begin{cor}
% 15.B.5
\label{Cor-Theo-SepTraceSIN}
Let $G$ be a compactly generated locally compact MAP group. 
Then $G$ is a SIN group.
\end{cor}
 
The assumption about the compact generation of $G$ in Corollary~\ref{Cor-Theo-SepTraceSIN}
and \emph{a fortiori} in Theorem~\ref{Theo-SepTraceSIN} is necessary,
as shown by the following example which was communicated to us by P.-E.~Caprace.

\begin{exe}
% 15.B.6
\label{Exa-MAP-NonSIN}
For every $n \in \N$, let $A_n = {\mathrm Alt}(5)$ be the alternating group on $5$ symbols
and let $C_n$ be any proper subgroup of $A_n$ with $C_n \neq \{e\}$.
The product
$$
A \, := \, \prod_{n \in \N} A_n ,
$$
equipped with the product topology, is a compact group
containing $C := \prod_{n \in \N} C_n$ as a closed subgroup.
Let $G$ be the \emph{restricted product} of the $A_n$'s with respect to the $C_n$'s;
recall that $G$ is the subgroup of $A$ defined as follows:
\begin{itemize}
\item
$G = \bigcup_{F} \left( \prod_{n \in F} A_n \times \prod_{n \in \N \smallsetminus F} C_n \right)$,
where $F$ runs over the finite subsets of $\N$;
\item
$G$ is equipped with the unique group topology 
for which $C$ is an open compact subgroup of $G$.
\end{itemize}
Since the canonical injection $G \to A$ is obviously continuous,
$G$ is a MAP group (see Proposition~\ref{Pro-AP-MAP}).
We are going to show that $G$ is not a SIN group.

\vskip.2cm

We claim first that every open normal subgroup of $G$ has finite index in $G$.
Indeed, let $N$ be an open normal subgroup of $G$.
Then $N \cap C$ is an open normal subgroup of the compact group $C$
and hence $N$ contains a subgroup of finite index in $C$.
It follows that $N$ contains a subgroup $H$ of the form 
$$
H \, = \, \prod_{n \in F} \{e\} \hskip.2cm \times \prod_{n \in \N \smallsetminus F} C_n,
$$
for some finite subset $F$ of $\N$.
Since $N$ is normal, $N$ contains all $G$-conjugates of $H$.
As $A_n$ is simple and $C_n \neq \{e\}$,
it follows that $N$ contains the subgroup 
$$
\prod_{n \in F} \{e\} \times \prod_{n \in \N \smallsetminus F} A_n,
$$
which has finite index in $G$.
\par
 
Assume now, by contradiction, that $G$ is a SIN group.
Then $C$, which is a neighbourhood of $\{e\}$ in $G$, 
contains a $G$-invariant neighbourhood $U$ of $\{e\}$ such that $U \neq C$.
It follows that
$$
N \, := \, \bigcap_{g \in G} gCg^{-1}
$$
is a normal subgroup of $G$ containing $U$.
Hence, $N$ is open in $G$.
By what we have seen above, this implies that $N$ has finite index in $G$.
This is a contradiction,
since $N$ is contained in $C$ and since $C$ has infinite index in $G$.
\end{exe}

\section{Locally compact groups with a finite von Neumann algebra}
% 15.C
\label{S:LC-FiniteVN}

In this section, we will show that the class of LC groups 
with a finite group von Neumann algebra
coincides with the class of SIN groups.
\par

We will need the following elementary result.

\begin{lem}
% 15.C.1
\label{Lem-C0-Normal}
Let $G$ be a LC group and $\tau \,\colon \Li(G) \to \C$
a normal positive linear functional on the von Neumann algebra
$\Li(G ) =\lambda_G(G)''$ generated by the regular representation of $G$.
Then the function $\tau \circ \lambda_G \,\colon G \to \C$ belongs to $C^0(G)$.
\end{lem}

\begin{proof}
First, observe that every matrix coefficient of $\lambda_G$ belongs to $C^0(G)$;
indeed, a matrix coefficient
$g \mapsto \langle \lambda_G(g) f_1 \mid f_2 \rangle$
belongs to $C^c(G)$ when $f_1, f_2 \in C^c(G)$.
The claim follows by density of $C^c(G)$ in $L^2(G)$.
\par

Since $\tau$ is normal, $\tau$ is continuous
with respect to the ultraweak topology on $\Li(L^2(G))$
and hence there exists a sequence $(f_i)_{i \ge 1}$ in $L^2(G)$
with $\sum_{i = 1}^\infty \Vert f_i \Vert < +\infty$
such that
$$
\tau(T) \, = \, \sum_{i = 1}^\infty \langle T f_i \mid f_i \rangle
\hskip.5cm \text{for every} \hskip.2cm
T \in \lambda_G(G)''
$$
(see Th\'eor\`eme 1 in \cite[Chap.~I, \S~4, no~2]{Dixm--vN}).
For $N \ge 1$, the function $\varphi_N$, defined by 
$\varphi_N(g) = \sum_{i = 1}^N \langle \lambda(g) f_i \mid f_i \rangle$,
belongs to $C^0(G)$ and we have
$$
\begin{aligned}
\left\vert \varphi_N(g) - \tau(\lambda(g) \right\vert
\, &\le \,
\sum_{i = N+1}^\infty \left\vert \langle \lambda(g) f_i \mid f_i \rangle \right\vert
\\
\, &\le \,
\sum_{i = N+1}^\infty \Vert f_i \Vert^2,
\end{aligned}
$$
for every $ \in G$. 
This implies that $(\varphi_N)_{N \ge 1}$ converges to
$\tau \circ \lambda_G$ uniformly on $G$ and the claim follows.
\end{proof}

\begin{rem}
% 15.C.2
\label{Rem-Lem-C0-Normal}
A more precise version of Lemma~\ref{Lem-C0-Normal} is in fact true:
every normal positive linear functional on $\Li(G)$
is a matrix coefficient of $\lambda_G$ (see \cite[Th\'eor\`eme p.~218]{Eyma--64}).
\end{rem}

The following result is proved in \cite[13.10.5. Proposition]{Dixm--C*}
under the additional assumption that $G$ is unimodular.

\begin{theorem}
% 15.C.3
\label{Theo-SIN-FiniteVN}
Let $G$ be a locally compact group. The following 
properties are equivalent.
\begin{enumerate}[label=(\roman*)]
\item\label{iDETheo-SIN-FiniteVN}
$G$ is a SIN group;
\item\label{iiDETheo-SIN-FiniteVN}
$\Li(G)$ is a finite von Neumann algebra.
\end{enumerate}
\end{theorem}

\begin{proof}
Assume that $G$ is a SIN group.
Let $(V_i)_{i \in I}$ be a basis of compact neighbourhoods of $e$
which are invariant under conjugation. 
For every $i \in I$, set 
$$
f_{i} \, := \, \frac{\Un_{V_i}}{\sqrt{\mu_G(V_i)}}.
$$
The function of positive type $\varphi_i$ on $G$ 
associated to $(\lambda_G, f_i)$ is conjugation invariant
(see the proof of Theorem~\ref{Theo-SINSepTrace}).
\par

Let $\tau_i$ be the linear extension of $\varphi_i$ to $\Li(G)$;
thus, $\tau_i$ is the normal positive linear functional on $\Li(G)$ defined by 
$$
\tau_i(T) \, = \, \langle T f_i \mid f_i \rangle
\hskip.5cm \text{for every} \hskip.2cm
T \in \Li(G).
$$
Since $\varphi_i$ is invariant under conjugation, $\tau_i$ is a trace on $\Li(G)$.
\par

Let $T \in \Li(G)$. Assume that $\tau_{i}(T^*T) = 0$ for every $i \in I$.
Since
$$
\tau_{i}(T^*T) \, = \, \langle T f_{i} \mid T f_{i} \rangle \, = \, \Vert Tf_{i} \Vert^2,
$$
it follows that $Tf_{i} = 0$ for every $i \in I$. 
Observe that, since $(V_i)_{i \in I}$ is a basis of neighbourhoods of $e$,
the family $(f_i)_{i \in I}$ is an approximate identity for the $L^1(G)$-module $L^2(G)$,
that is, 
$$
\lim_i \Vert f_{i} \ast f - f \Vert_2 \, = \, 0
\hskip.5cm \text{for every} \hskip.2cm
f \in L^2(G).
$$
Since $T$ commutes with right convolutions by elements of $C_c(G)$, 
we have therefore, for every $f \in C_c(G)$,
$$
Tf \, = \, \lim_{i} T(f_i \ast f) \, = \, \lim_{i} T(f_i) \ast f \, = \, 0.
$$
Hence $T = 0$. This shows that $\Li(G)$ is a finite von Neumann algebra.

\vskip.2cm

Assume conversely that $\Li(G)$ is a finite von Neumann algebra.
Let $\mathcal{T}$ be the set of finite normal traces $\tau$ on $\Li(G)$
normalized by $\tau(I) = 1$ and set
$$
\Phi \, := \, \{ \tau \circ \lambda_G \mid \tau \in \mathcal{T} \}.
$$
Observe that every $\varphi \in \Phi$
is a conjugation invariant function of positive type on $G$.
\par

We claim that $\Phi$ separates the points of $G$.
Indeed, let $x, y\in G$ with $x\neq y$. For 
$$
T \, := \, \lambda_G(x) - \lambda_G(y) \in \Li(G),
$$
we have $T \neq 0$
and hence there exists $\tau \in \mathcal{T}$ such that $\tau(T^*T) \ne 0$.
It follows from Proposition~\ref{Pro-BiInvMetric}
that $\tau \circ \lambda_G(x) \ne \tau \circ \lambda_G(y)$. 
\par

Fix $\varphi_0 \in \Phi$ and set 
$$
K_0 \, := \, \{g \in G \mid \vert \varphi_0(g)-1 \vert \le 1\}.
$$
Then $K_0$ is closed neighbourhood of $e$ and is conjugation invariant.
Moreover, $K_0$ is compact, since $\varphi_0\in C^0(G)$, by Lemma~\ref{Lem-C0-Normal}. 
\par

Consider the product space $\prod_{\varphi \in \Phi} \C$
equipped with the product topology. The map
$$
\psi \, \colon \, K_0 \to \prod_{\varphi \in \Phi} \C, \hskip.2cm g \mapsto \varphi(g)
$$
is continuous.
Moreover, $\psi$ is injective, since $\Phi$ separates the points of $G$.
It follows from the compactness of $K_0$ that the inverse map 
$$
\psi^{-1} \, \colon \, \psi(K_0) \to K_0
$$
is continuous.
This implies that the family of the sets $U_{F, \varepsilon}$, defined by
$$
U_{F, \varepsilon} \, = \, \{g \in K_0 \mid \vert \varphi(g) - 1 \vert \le \varepsilon
\hskip.2cm \text{for all} \hskip.2cm
\varphi \in F\}
$$
for a finite subset $F$ of $\Phi$ and $\varepsilon > 0$,
is a fundamental family of neighbourhoods of $e$.
Since every set $U_{F, \varepsilon}$ is conjugation invariant,
$G$ is a SIN group.
\end{proof}

\section{Connected and totally disconnected LC groups}
% 15.D
\label{S: ConnectedTotDisconnected}

We apply Theorem~\ref{Theo-SepTraceSIN} to two classes of groups,
the connected LC groups and the totally disconnected LC groups.
For this, we will need to use some structure theorems for SIN groups
contained in either one of these classes.
The first one is due to Iwasawa (\cite[Theorem 3]{Iwas--51}).
For the proof of a more general result, see Theorem 16.4.6 in \cite{Dixm--C*}
(compare also Example~\ref{Ex-MAP-Groups}).
 
\begin{theorem}[\textbf{Iwasawa's theorem}]
% 15.D.1
\label{Theo-Iwasawa}
Let $G$ be a connected topological group.
The following properties are equivalent:
\begin{enumerate}[label=(\roman*)]
\item\label{iDETheo-Iwasawa}
$G$ is a SIN group;
\item\label{iiDETheo-Iwasawa}
$G$ is isomorphic, as topological group, to a direct product
$K \times \R^n$ for a compact group $K$ and some $n \ge 0$.
\end{enumerate} 
\end{theorem} 

Observe that a connected LC group
is generated by any compact neighbourhood of its unit element.
The following result from \cite{KaSi--52}
is an immediate consequence
of Theorem~\ref{Theo-SINSepTrace}, Theorem~\ref{Theo-SepTraceSIN},
and Iwasawa's theorem;
the result was previously proved for semi-simple real Lie groups in \cite{SevN--50}.

\begin{theorem}
% 15.D.2
\label{Theo-KadisonSinger}
Let $G$ be a connected LC group. 
The following properties are equivalent:
\begin{enumerate}[label=(\roman*)]
\item\label{iDETheo-KadisonSinger}
${\rm Rep}_{\rm FT}(G)$ separates the points of $G$;
\item\label{iiDETheo-KadisonSinger}
$G$ is isomorphic, as topological group, to a direct product $K\times \R^n$
for a compact group $K$ and some $n \ge 0$.
\end{enumerate}
\end{theorem}

Theorem~\ref{Theo-KadisonSinger} has the following consequence
for the finite characters of a connected LC group.

\begin{cor}
% 15.D.3
\label{Cor-Theo-KadisonSinger}
Let $G$ be a connected LC group.

Every factor representation of finite type of $G$
is a multiple of a finite dimensional irreducible representation;
equivalently, $E(G)$ consists of the characters
associated to the finite dimensional irreducible representations of $G$.
\end{cor}
 
\begin{proof}
Let $\pi$ be a factor representation of finite type of $G$.
Then $\pi \circ p$ is a faithful factor representation of finite type of $G/\ker \pi$,
where $p \,\colon G \to G/\ker \pi$ is the canonical epimorphism.
Since $G/\ker \pi$ is a connected LC group,
it follows from Theorem~\ref{Theo-KadisonSinger}
that $G/\ker \pi$ is isomorphic to a direct product $K \times \R^n$
for a compact group $K$ and some $n \ge 0$.
Every factor representation of $K \times \R^n$
is a multiple of a finite dimensional irreducible representation.
\end{proof}
 
We now turn to totally disconnected groups.
Here, we will use von Dantzig's theorem;
for a proof, see \cite[Theorem 7.7]{HeRo--63} or \cite[Theorem 2.E.6]{CoHa--16}.
 
\begin{theorem}[\textbf{van Dantzig's theorem}]
% 15.D.4
\label{Theo-Dantzig}
Let $G$ be a totally disconnected LC group.

The compact open subgroups of $G$ form a basis of neighbourhoods 
of the group unit $e$.
\end{theorem} 

The following result was suggested to us by P.-E.\ Caprace.

\begin{theorem}
% 15.D.5
\label{Theo-FiniteRepTotallyDisconnected}
Let $G$ be a \emph{compactly generated} totally disconnected LC group.
The following properties are equivalent:
\begin{enumerate}[label=(\roman*)]
\item\label{iDETheo-FiniteRepTotallyDisconnected}
${\rm Rep}_{\rm FT}(G)$ separates the points of $G$.
\item\label{iiDETheo-FiniteRepTotallyDisconnected}
$G$ is an SIN group.
\item\label{iiiDETheo-FiniteRepTotallyDisconnected}
The compact open \emph{normal}\ subgroups of $G$
form a basis of neighbourhoods of the group unit $e$.
\item\label{ivDETheo-FiniteRepTotallyDisconnected}
$G$ is the projective limit of discrete groups.
\end{enumerate}
\end{theorem}

\begin{proof}
The fact that \ref{iDETheo-FiniteRepTotallyDisconnected}
and \ref{iiDETheo-FiniteRepTotallyDisconnected}
are equivalent follows from Theorem~\ref{Theo-SepTraceSIN}
and Theorem~\ref{Theo-SINSepTrace}.
\par

Assume that $G$ is a SIN group and let $U$ be a conjugation invariant neighbourhood of $e$. 
By van Dantzig's theorem, $U$ contains a compact open subgroup $H$ of $G$.
In turn, $H$ contains a conjugation invariant neighbourhood $V$ of $e$,
since $G$ is a SIN group. Let 
$$
N \, := \, \bigcap_{g \in G} gHg^{-1}.
$$
Then $N$ is a compact normal subgroup of $G$ and we have
$V \subset N \subset U$.
So, $N$ is a neighbourhood of $e$ contained in $U$.
This proves that \ref{iiDETheo-FiniteRepTotallyDisconnected}
implies \ref{iiiDETheo-FiniteRepTotallyDisconnected}
\par

The fact that \ref{iiiDETheo-FiniteRepTotallyDisconnected}
and \ref{ivDETheo-FiniteRepTotallyDisconnected} are equivalent
follows from the definition of a projective limit of LC groups
(see \cite[Chap.~III, \S~7]{BTG1--4}).
\end{proof}

\begin{rem}
% 15.D.6
\label{Rem-SINTD}
(1)
Let $G$ be a (not necessarily compactly generated)
totally disconnected LC group.
Assume that $G$ is a SIN group.
Then $E(G)$ admits a neat description
in terms of the characters of the discrete quotients of $G$.
Indeed, as the proof of Theorem~\ref{Theo-FiniteRepTotallyDisconnected} shows,
$G$ is the projective limit of discrete groups.
Hence, using Lemma~\ref{Lem-FactQuot}, we have
$$
E(G) \, = \, \bigcup_{N \in \mathcal{N}} E(G/N) ,
$$
where $\mathcal{N}$ is the family of compact open normal subgroups of $G$.

\vskip.2cm

(2)
In Theorem~\ref{Theo-FiniteRepTotallyDisconnected},
it is necessary to assume that $G$ is compactly generated; 
indeed, the group $G$ from Example~\ref{Exa-MAP-NonSIN}
is a totally disconnected group which is a MAP group but not a SIN group.
\end{rem}

The following corollary of Theorem~\ref{Theo-FiniteRepTotallyDisconnected}.
applies in particular to compactly generated totally disconnected LC groups
which are topologically simple and non discrete.

\begin{cor}
% 15.D.7
\label{Cor-Theo-FiniteRepTotallyDisconnected}
Let $G$ be a compactly generated totally disconnected LC group.
Assume that the only open normal subgroup of $G$ is $G$ itself.
\par

Then every representation of finite type of $G$
is a multiple for the trivial representation $1_G$;
equivalently, $E(G) = \{1_G\}$.
\end{cor}

\begin{proof}
Let $\pi$ be a representation of finite type of $G$.
We have to show that $\ker \pi = G$. 
\par

Assume, by contradiction, that $\ker \pi$ is a proper normal subgroup of $G$.
Then $G/\ker \pi$ is a non-trivial compactly generated totally disconnected LC group.
Let $p \,\colon G \twoheadrightarrow G/\ker \pi$ be the canonical epimorphism.
Since $\pi \circ p$ is a faithful representation of finite type of $G/\ker \pi$,
it follows from Theorem~\ref{Theo-FiniteRepTotallyDisconnected}
that $G/\ker \pi$ has a proper open normal subgroup $N$.
Then $p^{-1}(N)$ is a proper open normal subgroup of $G$ and this is a contradiction.
\end{proof}

\section{Faithful traces on group C*-algebras}
% Section 15.E
\label{SS:Faithful}

Let $\Gamma$ be a discrete group. Observe that, since $\Gamma$ is a SIN group,
${\rm Rep}_{\rm FT}(\Gamma)$ separates the points of $\Gamma$
(see Theorem~\ref{Theo-SINSepTrace}).
In fact, the regular representation $\lambda_\Gamma$
is a representation of finite type (see Proposition~\ref{Pro-TraceDiscreteGroup})
and is obviously point separating.
\par

Every finite type representation of $\Gamma$
extends to a finite type representation
of the maximal C*-algebra $C^*_{\rm max}(\Gamma)$ of $\Gamma$
and every finite type representation of $C^*_{\rm max}(\Gamma)$ is such an extension.
We will deal in this section with the question
whether the family of finite type representations of $\Gamma$
separates the points of $C^*_{\rm max}(\Gamma)$. 
This question may be rephrased as follows (see Corollary~\ref{Cor-FaithfulTrace}):

\begin{ques}
% 15.E.1
\label{quesfaithfultrace}
Does $C^*_{\rm max}(\Gamma)$ admit a faithful finite trace?
\end{ques}

\index{Connes embedding problem}
The question above is of interest.
For instance, it was shown in \cite[Proposition 8.1]{Kirc--93}
that the existence of such a trace for $\Gamma = \SL_2(\Z) \times \SL_2(\Z)$
is equivalent to a positive solution of the Connes embedding problem.
This problem is an outstanding open question in operator algebras;
its original formulation \cite[Page 105]{Conn--76} asks whether every factor of type II$_1$
acting on a separable Hilbert space
can be embedded in an ultrapower $R^\omega$
of the hyperfinite factor $R$ of type II$_1$.
% page 105 de l'article original de Connes de 1976, en termes de
% "approximate embedding"
The conjecture that the answer is positive is equivalent
to a variety of other important problems in operator algebras,
most notably by the work of Kirchberg;
see also \cite{Ozaw--13}.
% Voir aussi le Springer Lecture Notes 2136 de Capraro et Lupini ? --- non !
\par

We first observe that the existence of a faithful trace 
on the maximal C*-algebra of a group is inherited by each of its subgroups.

\begin{prop}
% 15.E.2
\label{Pro-FaithfulTRaceSubgroup}
Let $\Gamma$ be a discrete group and $\Delta$ a subgroup.
\par

If $C^*_{\rm max}(\Gamma)$ admits a faithful finite trace,
then so does $C^*_{\rm max}(\Delta)$.
\end{prop}

\begin{proof}
Indeed, let $H$ be a subgroup of $\Gamma$. 
As $\Gamma$ is discrete, the canonical map
$\C[H] \hookrightarrow \C[\Gamma]$ extends to an injective homomorphism 
$C^*_{\rm max}(H) \hookrightarrow C^*_{\rm max}(\Gamma)$,
by Theorem~\ref{Theo-C*-GroupHom}.
The claim is now clear.
\end{proof}

Recall that $\Tr(\Gamma)$ is the space of traces on $\Gamma$, 
that is, the space of functions on $\Gamma$ 
which are central of positive type.
Recall from Section~\ref{Section-GNS-Traces} 
that every trace on $\Gamma$ extends to a unique finite trace
of the maximal C*-algebra $C^*_{\rm max}(\Gamma)$ of $\Gamma$ 
and that, conversely,
every finite trace of $C^*_{\rm max}(\Gamma)$ 
defines by restriction a trace on $\Gamma$.
\par

We first relate the faithfulness of a trace $t$ on $C^*_{\rm max}(\Gamma)$
to the faithfulness of the associated GNS representation of $C^*_{\rm max}(\Gamma)$.
This is actually true for any positive functional on $C^*_{\rm max}(\Gamma)$. 
\par

Let $\varphi \in P(\Gamma)$ be a function of positive type on $\Gamma$. 
We will need to use a model for 
the GNS triple $(\pi_\varphi, \Hi_\varphi, \xi_\varphi)$
associated to $\varphi$ which is different from the one
from Construction~\ref{constructionGNS2}.
\par

Denote by $\Phi$ the unique extension of $\varphi$ as a
positive functional on $C^*_{\rm max}(\Gamma)$. 
The set $J = \{x \in C^*_{\rm max}(\Gamma) \mid \Phi(x, x) = 0\}$
is a left ideal of $C^*_{\rm max}(\Gamma)$ and a scalar product is defined on 
$C^*_{\rm max}(\Gamma) / J$ by 
$$
\langle x+J \mid y+J \rangle \, = \, \Phi( x, y) 
\hskip.5cm \text{for all} \hskip.2cm 
x, y \in C^*_{\rm max}(\Gamma).
$$
Let $\Hi_\varphi$ be the Hilbert space completion of $C^*_{\rm max}(\Gamma) / J$; 
then a representation $\pi_\varphi$ of $C^*_{\rm max}(\Gamma)$ on $\Hi_\varphi$ 
is defined on the unital $*$-algebra $C^*_{\rm max}(\Gamma) / J$ by 
$$
\pi_\varphi(a) (x+J) \, = \, ax +J 
\hskip.5cm \text{for all} \hskip.2cm 
x \in C^*_{\rm max}(\Gamma).
$$
One checks that the triple $(\pi_\varphi, \Hi_\varphi, \xi_\varphi)$ 
consisting of the restriction of $\pi_\varphi$ to $\Gamma$ and the cyclic vector
$\xi_\varphi = 1+J$ is equivalent to the GNS triple for $\varphi$ from 
Construction~\ref{constructionGNS2}.

\begin{lem}
% 15.E.3
\label{Lem-FaithfulTrace}
Let $\varphi \in P(\Gamma)$ and $\pi_\varphi$ 
the associated GNS representation.
Then 
$$
\textnormal{C*ker}(\pi_\varphi) \, = \,
\{x \in C^*_{\rm max}(\Gamma) \mid \Phi(x^*x) = 0\},
$$
where $\Phi$ denotes the unique extension of $\varphi$ 
to a positive functional on $C^*_{\rm max}(\Gamma)$.
\end{lem}

\begin{proof}
Let $x \in C^*_{\rm max}(\Gamma)$. 
Then, with the notation as above, 
$x \in \textnormal{C*ker}(\pi_\varphi)$ if and only if 
$\pi_\varphi(x) (y+J) = xy+J = J$ for every $y \in C^*_{\rm max}(\Gamma)$, 
%\marginpar{check $=J$}
that is, if and only if $x \in J$. 
Since $J = \{x \in C^*_{\rm max}(\Gamma) \mid \Phi(x^*x) = 0\}$,
this proves the claim.
\end{proof}

The following corollary is an immediate consequence of Lemma~\ref{Lem-FaithfulTrace}.

\begin{cor}
% 15.E.4
\label{Cor-FaithfulTrace}
Let $\varphi \in \Tr_1(\Gamma)$ and $\Phi$ the unique extension of $\varphi$ 
to a positive functional on $C^*_{\rm max}(\Gamma)$;
let $\pi_\varphi$ be the associated GNS representation.
The following properties are equivalent:
\begin{enumerate}[label=(\roman*)]
\item\label{iDECor-FaithfulTrace}
$\Phi$ is faithful;
\item\label{iiDECor-FaithfulTrace}
$\textnormal{C*ker}(\pi_\varphi) = \{0\}$.
\end{enumerate}
\end{cor}

\subsection{The case of amenable groups}
% 15.E.a
\label{SS:AmenGroupFaithfulTr}

We next show that Question \ref{quesfaithfultrace}
has a positive answer when $\Gamma$ is an amenable group.
Recall first that, for every discrete group $\Gamma$, 
the Dirac function $\delta_e$ is a trace on $\Gamma$.

\begin{prop}
% 15.E.5
\label{Pro-Amenable}
Let $\Gamma$ be a discrete group.
The following properties are equivalent:
\begin{enumerate}[label=(\roman*)]
\item\label{iDEPro-Amenable}
the unique extension of $\delta_e$ to a finite trace $\Phi$ 
on $C^*_{\rm max}(\Gamma)$ is faithful;
\item\label{iiDEPro-Amenable}
$\Gamma$ is amenable
\end{enumerate}
\end{prop}

\begin{proof} 
The GNS representation corresponding to $\delta_e$ 
is the regular representation $\lambda_\Gamma$ of $\Gamma$.
Therefore, by Corollary~\ref{Cor-FaithfulTrace}, 
$\Phi$ is faithful if and only if $\textnormal{C*ker}(\lambda_\Gamma) = \{0\}$.
This happens if and only if $\Gamma$ is amenable
(see Theorem \ref{HulanickiReiter}).
\end{proof}

\subsection{C*-maximally almost periodic groups}
% 15.E.b
\label{SS:C*-MAP}

The situation with regard to Question \ref{quesfaithfultrace}
is less understood in the case of non-amenable groups.
We introduce a class of groups for which the answer is positive.

\begin{defn}
% 15.E.6
\label{Def-CMaxAP}
\index{C*-maximally almost periodic group, or C*-MAP group}
A discrete group $\Gamma$ is said to a \textbf{C*-maximally almost periodic} group,
or a \textbf{C*-MAP} group, 
if the finite-dimensional representations of $\Gamma$ 
separate the points of $C^*_{\rm max}(\Gamma)$, 
that is, for every $x \in C^*_{\rm max}(\Gamma)$ with $x \ne 0$, 
there exists a finite-dimensional representation $\pi$ of $\Gamma$ such that 
$x \notin \textnormal{C*ker}(\pi)$.
\end{defn}

We first give a reformulation of the property of being a C*-MAP group 
in terms of Fell's topology.

\begin{prop}
% 15.E.7
\label{Pro-C*MAP-Reformulation}
Let $\Gamma$ be a discrete group.
Let $\widehat \Gamma_{\rm fd}$ be as in Section \ref{SectionFdrep}
the set of equivalence classes 
of irreducible finite-dimensional representations of $\Gamma$.
The following properties are equivalent:
\begin{enumerate}[label=(\roman*)]
\item\label{iPro-C*MAP-Reformulation}
$\Gamma$ is a C*-MAP group;
\item\label{iiPro-C*MAP-Reformulation}
$\widehat \Gamma_{\rm fd}$ is dense in $\widehat \Gamma$ for the Fell topology;
\item\label{iiiPro-C*MAP-Reformulation}
for every $x \in C^*_{\rm max}(\Gamma)$, we have 
$$
\Vert x \Vert \, = \, \sup_{\pi \in \widehat \Gamma_{\rm fd}} \Vert \pi(x) \Vert.
$$
\end{enumerate}
\end{prop}

\begin{proof}
Let $\pi$ be a finite-dimensional representation of $\Gamma$; 
then we can decompose $\pi$ as a direct sum $\pi = \bigoplus_{i = 1}^n \pi_i$ 
of irreducible representations, which are of course finite-dimensional.
We have
$$
\Vert \pi(x) \Vert \, = \, \sup_{i = 1, \hdots, n} \Vert \pi_i(x) \Vert 
\hskip.5cm \text{for all} \hskip.2cm 
x \in C^*_{\rm max}(\Gamma).
$$
This shows that \ref{iPro-C*MAP-Reformulation}
and \ref{iiiPro-C*MAP-Reformulation} are equivalent.

\vskip.2cm

Observe that \ref{iiPro-C*MAP-Reformulation} holds if and only if 
$$
\bigcap_{\pi \in \widehat \Gamma_{\rm fd}} \textnormal{C*ker}(\pi)
\, = \, 
\bigcap_{\pi \in \widehat \Gamma} \textnormal{C*ker}(\pi),
$$
that is, if and only if,
$$
\bigcap_{\pi \in \widehat \Gamma_{\rm fd}} \textnormal{C*ker}(\pi)
 \, = \, \{0\}.
$$
his shows that \ref{iPro-C*MAP-Reformulation}
and \ref{iiPro-C*MAP-Reformulation} are equivalent.
\end{proof}

Next, we observe that a C*-MAP group is a MAP-group
(see Section~\ref{SectionFdrep})
and that its maximal C*-algebra has a faithful trace.

\begin{prop}
% 15.E.8
\label{Pro-C*MAP}
Let $\Gamma$ be a C*-MAP group. Then
\begin{enumerate}[label=(\arabic*)]
\item\label{iDEPro-C*MAP}
$\Gamma$ is a MAP-group;
\item\label{iiDEPro-C*MAP}
$C^*_{\rm max}(\Gamma)$ has a faithful finite trace if $\Gamma$ is countable.
\end{enumerate}
\end{prop}

\begin{proof}
\ref{iDEPro-C*MAP}
For $\gamma \in \Gamma \smallsetminus \{e\}$, we have 
$\delta_\gamma - \delta_e \in C^*_{\rm max}(\Gamma) \smallsetminus \{0\}$,
hence, there exists a finite-dimensional representation $(\pi, \Hi)$ of $\Gamma$
such that $\delta_\gamma - \delta_e \notin \textnormal{C*ker}(\pi)$,
that is, $\pi(\gamma) \ne {\mathrm{Id}}_{\Hi}$.

\vskip.2cm

\ref{iiDEPro-C*MAP}
Since $\Gamma$ is countable, $C^*_{\rm max}(\Gamma)$ is separable,
hence, there exists a sequence of finite-dimensional representations 
$\pi_n \,\colon C^*_{\rm max}(\Gamma) \to M_{k_n}(\C)$ such that 
$$
\bigcap_{n} \textnormal{C*ker}(\pi_n) \, = \, \{0\}.
$$
For every $n$, let $\tau_n$ be the normalized standard trace on $M_{k_n}(\C)$. 
Choose a sequence $\alpha_n$ of positive real numbers with $\sum_n \alpha_n = 1$
and define $t \,\colon C^*_{\rm max}(\Gamma) \to \C$ by 
$$
t(x) \, = \, 
\sum_{n} \alpha_n tau_n(\pi_n(x)) 
\hskip.5cm \text{for all} \hskip.2cm 
x \in C^*_{\rm max}(\Gamma).
$$
Observe that $t(x)$ is well-defined since 
$$
\vert \tau_n(\pi_n(x)) \vert \, \le \, \Vert \pi_n(x) \Vert \, \le \, \Vert x \Vert 
\hskip.5cm \text{for all} \hskip.2cm
x \in C^*_{\rm max}(\Gamma).
$$
One checks that $t$ is a (finite) trace on $C^*_{\rm max}(\Gamma)$. 
Moreover, if $t(x^*x) = 0$ for $x \in C^*_{\rm max}(\Gamma)$, 
then $\tau_n(\pi_n(x^*x)) = 0$ and hence
$\pi_n(x^*x) = 0$ for all $n$. So, $x = 0$ and therefore $t$ is faithful.
\end{proof}

\begin{rem}
% 15.E.9
\label{Rem-FaithfulTraceC*MAP}
There are groups which are not C*-MAP
but which have a maximal C*-algebra admitting a faithful finite trace.
Indeed, if $\Gamma$ is an amenable group which is not MAP, 
$C^*_{\rm max}(\Gamma)$ has a faithful finite trace (Lemma~\ref{Pro-Amenable}) 
and $\Gamma$ is not C*-MAP.
An example of such a group is $SL_n(\K)$, 
where $n \ge 2$ and $\K$ is an infinite algebraic extension of a finite field.
\end{rem}

The next result will provide us with the first example 
of a non-amenable group $\Gamma$
for which $C^*_{\rm max}(\Gamma)$ admits a faithful finite trace. 
For the proof, we need the following elementary fact on operator theory.

\begin{lem}
% 15.E.10
\label{Lem-ContractionUnitary}
Let $\Hi$ be a Hilbert space and $T \in \Li (\Hi)$ with $\Vert T \Vert \le 1$.
\par

The operator $U$ on $\Hi \oplus \Hi$ defined by the matrix
$$
 \begin{pmatrix} T & (I - TT^*)^{1/2}
 \\
 (I - T^*T)^{1/2} & -T^*
 \end{pmatrix}
 $$
 is unitary. 
 \end{lem}
 
\begin{proof}
One computes that 
$
UU^* \, = \, U^*U \, = \, 
\begin{pmatrix} I & 0 \\ 0 & I \end{pmatrix}.
$
\end{proof}

The next theorem is from \cite{Choi--80}.

\begin{theorem}[\textbf{Choi}]
% 15.E.11
\label{Theo-FreeGroup}
Let $\Gamma = F_2$ be the free group on 2 generators.
\par

Then $\Gamma$ is a C*-MAP group; 
in particular, $C^*_{\rm max}(\Gamma)$ admits a faithful finite trace. 
\end{theorem}

\begin{proof}
Let $\pi \,\colon C^*_{\rm max}(\Gamma) \to \Li (\Hi)$ be a faithful representation
of $C^*_{\rm max}(\Gamma)$ on a separable Hilbert space $\Hi$
(one may take, for instance, $\pi = \bigoplus_n \pi_n$ for a sequence of irreducible
representations $\pi_n$ of $\Gamma$ 
such that $\bigcap_{n} \textnormal{C*ker}(\pi_n) = \{0\}$).
Set $A = \pi(a)$ and $B = \pi(b)$,
where $a$ and $b$ are free generators of $\Gamma$.
\par

Let $(P_n)_{n \ge 1}$ be a sequence of orthogonal projections in $\Li (\Hi)$
such that $\dim P_n(\Hi) = n$ for every $n \ge 1$
and $\lim_n P_n = I$ in the strong operator topology.
Set 
$$
S_n \, := \, P_n A P_n
\hskip.5cm \text{and} \hskip.5cm
T_n \, := \, P_n B P_n.
$$
Then $S_n, T_n$ may be viewed as operators on $\Hi_n := P_n(\Hi)$
and we have $\Vert S_n \Vert, \Vert T_n \Vert \le 1$.
Define $A_n, B_n \in \Li (\Hi_n \oplus \Hi_n)$ by 
$$
\begin{aligned}
A_n \, &= \, \begin{pmatrix} S_n & (I - S_nS_n^*)^{1/2}
\\
(I - S_n^*S_n)^{1/2} & -S_n^* \end{pmatrix} ,
\\
B_n \, &= \, \begin{pmatrix} T_n & (I- T _nT_n^*)^{1/2}
\\ (I - T_n^*T_n)^{1/2} & -T_n^* \end{pmatrix}
\end{aligned} .
$$
By Lemma~\ref{Lem-ContractionUnitary},
$A_n$ and $B_n$ are unitary operators on $\Hi_n \oplus \Hi_n$.
So, a finite-dimensional representation
$\pi_n \,\colon \Gamma \to \Li (\Hi_n \oplus \Hi_n)$ of $\Gamma$ is defined by setting 
$$
\pi_n(a) \, := \, A_n
\hskip.5cm \text{and} \hskip.5cm
\pi_n(b) \, := \, B_n.
$$
We claim that the sequence $(\pi_n)_{n \ge 1}$
separates the points of $C^*_{\rm max}(\Gamma)$.
\par

Indeed, observe that, viewing $A_n$ and $B_n$ as operators on $\Hi \oplus \Hi$
by defining them to be $0$ on $(\Hi_n \oplus \Hi_n)^\perp$, we have
$$
\lim_{n \to +\infty} A_n \, = \, 
\begin{pmatrix} A & 0 \\ 0 & -A^* \end{pmatrix} 
\hskip.5cm \text{and} \hskip.5cm
\lim_{n \to +\infty} B_n \, = \,
\begin{pmatrix} B & 0 \\ 0 & -B^* \end{pmatrix}
$$
and likewise 
$$
\lim_{n \to +\infty} A_n^* \, = \, 
\begin{pmatrix} A^* & 0 \\ 0 & -A \end{pmatrix} 
\hskip.5cm \text{and} \hskip.5cm
\lim_{n \to +\infty} B_n^* \, = \,
\begin{pmatrix} B^* 0 \\ 0 & -B \end{pmatrix}
$$
in the strong operator topology.
Let $x \in \C[\Gamma]$; then $x$ can be written as $x = P(a, a^{-1}, b, b^{-1})$ 
for a non-commutative polynomial $P$ in 4 variables.
Therefore, we have 
$$
\begin{aligned}
\lim_{n \to +\infty}\pi_n(x)& \, = \, \lim_{n \to +\infty} P(A_n, A_n^*, B_n, B_n^*)
\\
& \, = \, 
\begin{pmatrix} P(A,A^*, B, B^*) & 0 \\ 0 & P(-A^*,-A, - B^*, -B) \end{pmatrix} 
\end{aligned}
$$
in the strong operator topology.
It follows that 
$$
\left \Vert
\begin{pmatrix} P(A, A^*, B, B^*) & 0 \\ 0 & P(-A^*,-A, - B^*, -B) \end{pmatrix}
\right\Vert 
\, \le \, 
\liminf_{n \to +\infty} \Vert \pi_n(x)\Vert.
$$
Since $\Vert \pi(x)\Vert = \Vert P(A, A^*,B, B^*) \Vert$, we have
$$
\Vert \pi(x)\Vert \,\le \, 
\left \Vert 
\begin{pmatrix} P(A, A^*, B, B^*) & 0 \\ 0 & P(-A^*,-A, - B^*, -B) \end{pmatrix}
\right\Vert.
$$
We deduce that 
$$
\Vert \pi(x) \Vert \, \le \, 
\liminf_{n \to +\infty} \Vert \pi_n(x) \Vert
$$
and hence $\Vert \pi(x)\Vert = \sup_{n \ge 1} \Vert \pi_n(x)\Vert$
for every $x \in \C[\Gamma]$.
\par
 
Therefore, the representation
$$
C^*_{\rm max}(\Gamma) \, \to \, \Li \Big( \bigoplus_n (\Hi_n \oplus \Hi_n) \Big), 
\hskip.5cm
x \, \mapsto \, \bigoplus_{n \ge 1} \pi_n(x)
$$
is isometric on the dense subspace $\C[\Gamma]$
and hence isometric on $C^*_{\rm max}(\Gamma)$.
\end{proof}

\subsection{C*-residually finite groups}
% 15.E.c
\label{SS:C*-RF}

We introduce a stronger version of the C*-MAP property from the previous subsection
and show that it is shared by free groups.

\begin{defn}
% 15.E.12
\label{Def-CRF}
\index{C*-residually finite group}
We say that a discrete group $\Gamma$ is \textbf{C*-residually finite} 
if the representations of $\Gamma$ with finite images 
separate the points of $C^*_{\rm max}(\Gamma)$.
\end{defn}

We want to show that free groups are C*-residually finite;
this strengthening of Theorem~\ref{Theo-FreeGroup} is Theorem 2.2 in \cite{LuSh--04}, 
where C*-residually finite groups are called groups with property FD.
\par

We will give a variation of the proof from \cite{LuSh--04} using 
the fact that every bijective measure preserving transformation
of the interval $\mathopen[ 0, 1 \mathclose)$
can be approximated by dyadic permutations (\cite[p.65]{Halm--60});
the crucial tool we will use is Proposition~\ref{Pro-Dyadic},
which is based on \cite[Theorem 2.1]{AlEd--79}.
\par

Let $n \in \N$.
A \emph{dyadic interval} of rank $n$
is an interval of the form $\mathopen[ i/2^n, (i+1)/2^n \mathclose)$
for some $i \in \{0, \dots, n-1\}$. 
A \emph{dyadic decomposition} of rank $n$ of $\mathopen[ 0, 1 \mathclose)$
consists of all dyadic intervals of length $1/n$.
A \emph{dyadic permutation} of rank $n$ 
is a bijection from $\mathopen[ 0, 1 \mathclose)$ to $\mathopen[ 0, 1 \mathclose)$
which maps every dyadic interval of rank $n$
into a dyadic interval of rank $n$ by an ordinary translation. 
\par

Observe that every dyadic permutation
preserves the Lebesgue measure on $\mathopen[ 0, 1 \mathclose)$
and that the set $\mathcal{P}_n$ of dyadic permutations of rank $n$ 
is a finite group isomorphic to the symmetric group on $2^n$ letters.
More precisely, every $P \in \mathcal{P}_n$ is uniquely determined
by a permutation $\sigma$ of $\{ 1, \dots, 2^n \}$ for which $P(I_i) = I_{\sigma(i)}$,
where $I_i = \mathopen[ i/2^n, (i+1)/2^n \mathclose)$.

\begin{prop}
% 15.E.13
\label{Pro-Dyadic}
Let $S$ be a finite set of measurable Borel automorphisms
of $\mathopen[ 0, 1 \mathclose)$ 
preserving the Lebesgue measure $\mu$.
\par

For given $\delta > 0$ and $\varepsilon > 0$, there exist an integer $n \ge 1$
and a set $\{P^s \mid s \in S \}$ of dyadic permutations of rank $n$ such that 
$$
\mu \left(
\{ x \in \mathopen[ 0, 1 \mathclose) \mid
\vert s(x) - P^s(x) \vert \ge \delta \}
\right) \, < \, \varepsilon
\hskip.5cm \text{for every} \hskip.2cm
s \in S.
$$
\end{prop}

\begin{proof}
Let $I_1, \dots, I_{2^m}$ be a dyadic decomposition of $\mathopen[ 0, 1 \mathclose)$
such that $\mu(I_i) = \dfrac{1}{2^m} < \delta$.
For every $s \in S$ and every $i \in \{ 1, \dots, 2^m \}$, 
we can choose a compact subset $C_i^s$ of $ s^{-1}(I_i)$ such that
$$
%\mu\left(\bigcup_{i=1}^m C_i^s\right) \ge 1- \frac{\varepsilon}{2}.
\mu (C_i^s) \, \ge \, (1- \varepsilon) \mu (s^{-1}(I_i)) \, = \, \frac{1- \varepsilon}{2^m},
$$
by regularity of $\mu$.
Let
$$
c \, := \, \inf_{s \in S} \inf_{1 \le i \ne j \le m} d(C_i^s, C_j^s),
$$
where $d(A, B) = \inf\{ \vert x - y \vert \mid x \in A, y \in B \}$
for subsets $A, B$ of $\mathopen[ 0, 1 \mathclose)$.
Since the $C_i^s$ 's are compact and since $C_i^s \cap C_j^s = \emptyset$ for $i \ne j$,
we have $c > 0$.
\par

Choose another dyadic decomposition $J_1, \dots, J_{2^n}$ of $\mathopen[ 0, 1 \mathclose)$
such that $\mu(J_i) < c$. 
For every $j \in \{1, \dots, 2^n \}$ and $s \in S$, we have then 
$J_j \cap C_i^s \ne \emptyset$ for at most one $i \in \{1, \dots, 2^m\}$.
Upon choosing $n$ sufficiently large, we may assume that, moreover,
for every $i \in \{1,\dots, 2^m\}$, we have 
$$
\mu \Big( \bigcup_{j \in \tau_i^s} J_j \cap C_i^s \Big)
\, \le \, \mu(s^{-1}(I_i)) \, = \, 1/2^m ,
$$
where $\tau_i^s$ is the set $j \in \{1, \dots, 2^n \}$ with $J_j\cap C_i^s \ne \emptyset$.
\par

Fix $s \in S$.
We are going to define a dyadic permutation $P^s$ of rank $n$ 
with $P^s(J_j) \subset I_i$ for every $i \in \{1, \dots, 2^m\}$ and every $j \in \sigma^s_i$.
\par

Let $i \in \{1, \dots, 2^m\}$.
Let $\tau_i \subset I$ be such that $(J_{i'})_{i' \in \tau_i}$
are the dyadic intervals of length $1/2^n$ contained in $I_i$.
Since $\mu( \bigcup_{j \in \tau_i^s} J_j ) \le 1/2^m = \mu(I_i)$,
we can find an injective map $\sigma^s \,\colon \tau_i^s \to \tau_i$.
We then extend $\sigma^s$ to a permutation of $\{1, \dots, 2^n \}$
and let $P^s$ be the associated dyadic permutation.
\par

Fix $s\in S$. Set 
$$
A^s \, := \, \bigcup_{i=1}^{2^m} C_i^s \, = \,
\bigcup_{i = 1}^{2^m} \bigcup_{j \in \tau_i^s} (J_j \cap C_i^s) .
$$
Since the $C_i^s$ are pairwise disjoint, we have 
$$
\mu(A^s) \, = \, \sum_{i=1}^{2^m} \mu(C_i^s)
\, \ge \, \sum_{i=1}^{2^m} \frac{1- \varepsilon}{2^m} \, = \, 1 - \varepsilon.
$$
Let $x \in A^s$. Then $x \in C_i^s$ for some $i \in \{1,\dots, 2^m\}$.
It follows that $s(x) \in I_i$.
Since we also have $P^s(x) \in I_i$, we have
$$
\vert s(x) - P^s(x) \vert \, \le \, \frac{1}{2^m} \, < \, \delta.
$$
\end{proof}

\begin{theorem}[\textbf{Lubotzky--Shalom}]
% 15.E.14
\label{Theo-FreeGroup2}
Let $r \ge 1$ be an integer.
\par

The free group $F_r$ on $r$ generators is C*-residually finite.
\end{theorem}
 
\begin{proof}
Let $\pi$ be a representation of $\Gamma = F_r$
in a separable Hilbert space.
We have to show that $\pi$ is weakly contained 
in a direct sum of representations of $\Gamma$ with finite images.

%(see Section~\ref{SectionWC+FellTop}).

\vskip.2cm

$\bullet$ \emph{First step.}
We may assume that $\pi$ is the Koopman representation $\pi_X$ on $L^2(X, \mu)$
associated to a measure-preserving action of $\Gamma$
on a standard Borel space $X$, equipped with a probability measure $\mu$.
\par

Indeed, it is well-known (see \cite[A.7.15]{BeHV--08} and \cite[5.2.13]{Zimm--84})
that there exists a measure-preserving action
of $\Gamma$ on a standard probability space $(X,\mu)$
such that $\pi$ is equivalent to a subrepresentation
of the Koopman representation of $\Gamma$ on $L^2(X, \mu)$.
The claim follows immediately.

\vskip.2cm

$\bullet$ \emph{Second step.}
Let $\Gamma \curvearrowright (X,\mu)$ be a 
measure-preserving action of $\Gamma$ on a standard Borel space $X$
equipped with a probability measure $\mu$.
We claim that the Koopman representation $\pi_X$
is weakly contained in a direct sum of representations of $\Gamma$ with finite images.
\par

Assume first that the set $A \subset X$ of atoms of $\mu$ is not empty. 
Since $\mu$ is a probability measure,
the subset $A_\alpha$ of atoms of $\mu$-measure $\alpha$ is finite for every $\alpha > 0$
and $A$ is at most countable.
It follows that there exists a sequence $(\alpha_i)_{i \ge 1}$ of positive numbers
such that $A = \bigcup_{i \ge 1} A_{\alpha_i}$.
Every subset $A_{\alpha_i}$ is $\Gamma$-invariant.
We have therefore a decomposition 
$$
\pi_X \, = \, \Big( \bigoplus_{i \ge 1} \pi_{A_{\alpha_i}} \Big) \oplus \pi_{X \smallsetminus A}.
$$
Since $\Gamma$ acts on $A_{\alpha_i}$ as a permutation, for every $i \ge 1$,
the representation $\pi_{A_{\alpha_i}}$ factorizes through a finite quotient of $\Gamma$.
It follows that $\pi_A = \bigoplus_{i \ge 1} \pi_{A_{\alpha_i}} $ is weakly contained 
in a direct sum of representations of $\Gamma$ with finite images.
\par

It remains to show that $\pi_{X \smallsetminus A}$ is weakly contained 
in a direct sum of representations of $\Gamma$ with finite images.
\par

Since $X \smallsetminus A$ is a standard Borel space
and $\mu$ has no atoms in $X \smallsetminus A$,
we may assume that $X = \mathopen[ 0, 1 \mathclose)$
and $\mu$ is the Lebesgue measure (see \cite[14.41]{Kech--95}). 
\par
 
Let $S$ be a finite set such that $\Gamma$ is the free group over $S$.
Let $f_1, \dots, f_k$ be continuous functions on $\mathopen[ 0, 1 \mathclose]$
and let $\varepsilon > 0$. 
\par

Let $\delta > 0$ be such that $\vert f_i(x) - f_i(y) \vert < \varepsilon$
for all $x, y \in \mathopen[0, 1 \mathclose]$ such that $\vert x - y \vert < \delta$
and all $i \in \{1, \dots, k\}$.
\par

By Proposition~\ref{Pro-Dyadic},
there exist a set $\{P_s \mid s \in S\}$ of dyadic permutations of the \emph{same} rank such that 
$$
\mu \left( \{x \in \mathopen[ 0, 1 \mathclose)
\mid \vert s(x) - P^s(x) \vert \ge \delta \} \right) \, < \, \varepsilon
\hskip.5cm \text{for every} \hskip.2cm
s \in S.
$$
The set $(P_s)_{s \in S}$ defines a measure preserving action 
of the free group $\Gamma$ on $(\mathopen[ 0, 1 \mathclose), \mu)$,
which factorizes through a finite quotient of $\Gamma$. 
Denote by $\rho$ the Koopman representation 
on $L^2(\mathopen[ 0, 1 \mathclose), \mu)$ associated to this action.
\par

Setting $\pi := \pi_{X \smallsetminus A}$, we have,
for every $s \in S$ and every $i \in \{ 1, \dots, k \}$, 
$$
\begin{aligned}
\Vert \pi(s^{-1}) f_i - \rho(s^{-1}) f_i \Vert_2^2
\, &= \, \int_0^1 \vert f_i(s(x)) - f_i(P^s(x)) \vert^2 d\mu(x)
\\
& = \, \int_{\{x \hskip.1cm : \hskip.1cm
\vert s(x)-P^s(x) \vert \ge \delta \}} \vert f_i(s(x))-f_i(P^s(x)) \vert^2 d\mu(x)
\\
& \, + \int_{\{x \hskip.1cm : \hskip.1cm
\vert s(x)-P^s(x) \vert < \delta \}} \vert f_i(s(x))-f_i(P^s(x)) \vert^2d\mu(x)
\\
& \, \le 4 \max_{1 \le j \le k} \Vert f_j\Vert_\infty^2 \varepsilon + \varepsilon^2.
\end{aligned}
$$
It follows that there exists a sequence of representations 
$(\rho_n)_{n \ge 1}$ of $\Gamma$, every one factorizing through a finite quotient,
such that 
$$
\lim_{n \to +\infty} \langle \rho_n(s) f_i \mid f_i \rangle \, = \, \langle \pi(s) f_i \mid f_i \rangle 
\hskip.5cm \text{for all} \hskip.2cm
i \in \{ 1, \dots, k \}.
$$
Since $C(\mathopen[ 0, 1 \mathclose])$ is dense in $L^2(\mathopen[ 0, 1 \mathclose])$,
this implies that $\pi$ is weakly contained in $\oplus_{n \ge 1} \rho_n$.
\end{proof}

\begin{rem}
% 15.E.15
\label{Rem-Lubotzsky--Shalom}
It was shown in \cite{LuSh--04} that 
some other lattices in rank one Lie groups,
as for instance surface groups, are C*-residually finite.
\end{rem}

\subsection{Group C*-algebras with no faithful finite trace}
% 15.E.d
\label{SS:C*-AlgNotSep}

Using Theorem~\ref{Theo-SLn}, we now give examples of groups
for which the maximal C*-algebra has no faithful finite trace.
Observe that $\Gamma = GL_n(\K)$ or $\Gamma = \SL_n(\K)$ is amenable 
if $\K$ is an algebraic extension of a finite field;
so, in this case, $C^*_{\rm max}(\Gamma)$ has a faithful finite trace
by Proposition~\ref{Pro-Amenable}.

\begin{theorem}
% 15.E.16
\label{Thm-SLn-FaithfulTrace}
Let $\K$ be a field which is not an algebraic extension of a finite field
and $\Gamma = GL_n(\K)$ or $\Gamma = \SL_n(\K)$ for $n \ge 2$.
Then $C^*_{\rm max}(\Gamma)$ has no faithful finite trace.
\end{theorem}

\begin{proof}
If the characteristic of $\K$ is $0$, then $\K$ contains $\Q$; 
if the characteristic of $\K$ is $p$, then $\K$ 
contains the purely transcendental extension $\F_p(T)$.
In view of Lemma~\ref{Pro-FaithfulTRaceSubgroup},
it suffices to prove the result for $\Gamma = \SL_2(\Q)$
and for $\Gamma = \SL_2(\F_p(T))$.
\par

Let $t$ be a finite normalized trace on $C^*_{\rm max}(\Gamma)$;
we claim that $t$ is not faithful. 
\par

Let $\varphi \in \Tr_1(\Gamma)$ be the trace on $\Gamma$ corresponding to $t$.
By Corollary~\ref{Cor-FaithfulTrace}, it suffices to show that 
the GNS representation $\pi_\varphi$ corresponding to $\varphi$
is not faithful on $C^*_{\rm max}(\Gamma)$.
\par

It follows from Theorem~\ref{Theo-SLn} 
that there exists $\alpha \in \mathopen[ 0,1 \mathclose]$
and $\psi \in \Tr_1(\Gamma)$ with $\psi = 0$ on $\Gamma \smallsetminus Z$ such that 
$$
\varphi \, = \, \alpha \Un_{\Gamma}+ (1-\alpha) \psi,
$$
where $Z$ is the centre of $\Gamma$. 
\par

Since $Z$ is amenable, the GNS representation $\pi_\psi$ corresponding to $\psi$
is weakly contained in $\lambda_\Gamma$.
As a consequence, $\pi_\varphi$ is weakly contained 
in the direct sum $1_\Gamma \oplus \lambda_\Gamma$.
Therefore to show that $\pi_\varphi$ is not faithful on $C^*_{\rm max}(\Gamma)$,
it suffices to show that the representation $1_\Gamma \oplus \lambda_\Gamma$
is not faithful on $C^*_{\rm max}(\Gamma)$. 
 
\vskip.2cm

Observe that $\SL_2(\Q)$ contains the subgroup $\SL_2(\Z)$ 
which is a lattice in the LC group $\SL_2(\R)$ and that $\SL_2(\F_p[T])$
can be embedded as a lattice in the LC group $\SL_2(\F_p((T^{-1})))$, 
where $\F_p((T^{-1}))$ is the local field of Laurent series over $\F_p$.
\par

Set $\Lambda = \SL_2(\Z)$ when $\Gamma = \SL_2(\Q)$ 
and $\Lambda = \SL_2(\F_p[T])$ when $\Gamma = \SL_2(\F_p(T))$.
We see that in either case $\Gamma$ contains a lattice $\Lambda$ 
in the LC group $G = \SL_2(\mathbf k)$ for a local field $\mathbf k$.
\par

Let $\pi$ be a irreducible representation of $G$ such that $\pi \ne 1_G$
and which is not weakly contained in $\lambda_G$
(that is, $\pi$ is a representation from the complementary series of $G$).
We claim that $\pi \vert_{\Gamma}$
is not weakly contained in $1_\Lambda \oplus \lambda_\Gamma$.
Once proved, this will imply that $1_\Gamma \oplus \lambda_\Gamma$ 
is not faithful on $C^*_{\rm max}(\Gamma)$.
\par

Assume by contradiction that $\pi \vert_{\Gamma}$ 
is weakly contained in $1_\Gamma \oplus \lambda_\Gamma$.
Then $\pi \vert_{\Lambda}$ is weakly contained in $1_\Lambda \oplus \lambda_\Lambda$,
since $\lambda_\Lambda$ is weakly equivalent to $\lambda_\Gamma \vert_\Lambda$.
\par

Let $X$ denote the set of all $\sigma \in \widehat \Lambda$ 
which are weakly contained in $\pi$. 
By assumption, every $\sigma \in X \smallsetminus\{1_\Lambda \}$ 
is weakly contained in $\lambda_\Lambda$.
Since $\Lambda$ is not amenable, 
$1_\Lambda$ is not weakly contained in $\lambda_\Lambda$
and hence $1_\Lambda$ does not belong
to the closure of $X \smallsetminus \{1_\Lambda \}$ in the Fell topology.
It follows that either $1_\Lambda$ is \emph{contained} in $\pi \vert_{\Lambda}$
or $X \subset \widehat \Lambda \smallsetminus \{1_\Lambda \}$.
\par

By Howe--Moore theorem (see \cite[Theorem 5.2]{HoMo--79}),
$\pi$ has vanishing matrix coefficients; 
since $\Lambda$ is a discrete infinite subgroup of $G$, 
this implies that $1_\Lambda$ is not contained in $\pi \vert_{\Lambda}$. 
So, $X \subset \widehat \Lambda \smallsetminus \{1_\Lambda \}$ 
and this means that $\pi \vert_{\Lambda}$ is weakly contained in $\lambda_\Lambda$. 
\par

It follows, by continuity of induction, that $\Ind_{\Lambda}^G\pi \vert_{\Lambda}$
is weakly contained in $\Ind_{\Lambda}^G\lambda_\Lambda \simeq \lambda_G$.
However, 
$$
\Ind_{\Lambda}^G \pi \vert_{\Lambda} 
\, \simeq \,
\pi \otimes \Ind_{\Lambda}^G1_{\Lambda}
$$
and, since $\Lambda$ is a lattice in $G$, the trivial representation $1_G$
is contained in $\Ind_{\Lambda}^G1_{\Lambda}$.
This implies that $\pi$ is weakly contained in $\lambda_G$
and this contradicts the fact that $\pi$ belongs to the complementary series of $G$.
\end{proof}

\begin{rem}
% 15.E.17
The groups from Theorem~\ref{Thm-SLn-FaithfulTrace} are not residually finite 
(they are not even MAP).
In \cite{Bekk--99}, it is shown that the residually finite group $SL_n(\Z)$ 
for $n \ge 3$ is not C*-MAP.
\end{rem}

\section{Traces and Invariant Random Subgroups}
% Section 15.F
\label{Section:IRS}

An invariant random subgroup (IRS for short) of a LC group $G$
is a conjugation invariant probability measure on the space of closed subgroups of $G$.
IRS's generalize both closed normal subgroups and lattices of $G$ (see Example~\ref{Exa-IRS}).
The notion of an IRS was introduced in \cite{AbGV--14} and appeared implicitly in \cite{StZi--94};
see \cite{Gela--18} for a survey on IRS's.
\par

In this section, we establish a link between traces of a discrete group $\Gamma$
and invariant random subgroups of $\Gamma$.
\par
 
Let $G$ be a locally group. Let 
$$
\Sub (G) \, = \, \{ H \subset G \mid
H \hskip.2cm \text{is a closed subgroup of} \hskip.2cm G \}
$$
be equipped with the Chabauty topology;
recall that a basis for this topology is given by sets of the form
$$
\mathcal{U}(K; V_1, \dots, V_n) \, = \, \{ H \in \Sub(G) \mid
H \cap C = \emptyset
\hskip.1cm \text{and} \hskip.1cm 
H \cap V_i \neq \emptyset
\hskip.1cm \text{for} \hskip.1cm
i = 1, \dots, n\},
$$
where $K$ is a compact subset and the $V_i$'s are open subsets of $G$.
The space $\Sub(G)$ is compact (see \cite[Lemma 1]{Fell--62b}) 
and conjugation gives rise to an action of $G$ 
by homeomorphisms on $\Sub(G)$.
\par

\begin{defn}
% 15.F.1
\label{Def-IRS}
\index{Invariant random subgroup, or IRS}
An \textbf{invariant random subgroup} of $G$, abbreviated IRS, 
is a $G$-invariant probability measure on 
the Borel subsets of $\Sub(G)$.
\par

We denote by $\mathrm{IRS}(G)$ the set of IRS's of $G$. 
\end{defn}

The set $\mathrm{IRS}(G)$ is a convex compact subset
of the set of probability measures on $\Sub(G)$, equipped with the weak$^*$-topology.
The extremal points of $\mathrm{IRS}(G)$ are the \emph{ergodic} IRS's,
that is, the $G$-invariant ergodic probability measures on $\Sub(G)$.
Hence, the convex hull of the set of ergodic IRS's
is dense in $\mathrm{IRS}(G)$, by Krein--Milman's theorem.

\begin{exe}
% 15.F.2
\label{Exa-IRS}
(1)
Let $N$ be a closed normal subgroup of $G$. 
The Dirac measure $\delta_{N}$ is an IRS of $G$.
IRS's of this form are exactly the IRS's of $G$ which are Dirac measures.

\vskip.2cm

(2)
Let $\Gamma$ be a lattice in the LC group $G$
and let $\mu_{G/\Gamma}$ be the unique $G$-invariant regular probability measure
on the Borel subsets of $G/\Gamma$.
Let $\Phi_*(\mu_{G/\Gamma})$ be the image of $\mu_{G/\Gamma}$ under the measurable map 
$$
\Phi \, \colon \,G/\Gamma \to \Sub(G), \hskip.2cm g \Gamma \mapsto g\Gamma g^{-1}.
$$
Then $\Phi_*(\mu) \in \mathrm{IRS}(G)$.

\vskip.2cm

(3)
The construction of IRS 's from (2) can be generalized as follows.
Let $G$ be a second countable LC group
and $G \curvearrowright (X, \mathcal{B})$ a measurable action of $G$
on a countably separated Borel space $(X, \mathcal{B})$.
Then, for every $x \in X$, the stabilizer $G_x$ of $x$ in $G$
is a closed subgroup of $G$ (see Corollary 2.1.21 in \cite{Zimm--84}).
\par

Let $\mu$ be a probability measure on $\mathcal{B}$
and assume that $G \curvearrowright (X, \mathcal{B})$ preserves $\mu$.
The image $\Phi_*(\mu)$ of $\mu$ under the measurable map
$$
\Phi \, \colon \, X \to \Sub(G), \hskip.2cm x \mapsto G_x.
$$
is an IRS of $G$.

\vskip.2cm

(4)
Let $G$ be a second countable LC group. 
Every IRS of $G$ arises as in (3) for some measure preserving action 
of $G$ on a countably separated probability space $(X, \mathcal{B},\mu)$;
see \cite[Proposition 14]{AbGV--14} in the discrete case
and \cite[Theorem 2.6]{ABB+--17} in general.

\vskip.2cm

(5)
Let $G$ be a connected simple Lie group with trivial center and with real rank greater than two.
It is a remarkable result from \cite{StZi--94} (see also Theorem~2.4 in \cite{ABB+--17})
that every IRS of $G$ is a convex combination of $\delta_{G}, \delta_{\{e\}}$,
and IRS's associated to a lattice in $G$ as in (2).

\vskip.2cm

(6) 
IRS's of the infinite symmetric group $\Sym_{\rm fin}(\N)$
as in Example \ref{exThomaDiscrete}~\ref{iDEexThomaDiscrete}
have been classified in \cite{Vers--12}.
\end{exe}

The following result was first established in \cite{Vers--11},
with a different proof (see Theorem~\ref{Theo-Vershik}).

\begin{prop}
% 15.F.3
\label{Pro-IRS-Trace}
Let $\Gamma$ be a discrete group and $\nu \in \mathrm{IRS}(\Gamma)$
The function $\varphi_\nu \,\colon \Gamma \to \C$ defined by 
$$
\varphi_\nu(\gamma) \, = \, \nu \left( \{H \in \Sub(G) \mid \gamma \in H \} \right).
$$
is a normalized trace of $\Gamma$,
that is, $\varphi_\nu \in \Tr_1(\Gamma)$.
\end{prop}

\begin{proof}
The function $\varphi_\nu$ is central, since 
for $\gamma, \delta \in \Gamma$, we have, by $\Gamma$-invariance of $\nu$,
$$
\begin{aligned}
\varphi_\nu (\delta^{-1} \gamma \delta) 
&\, = \, \nu \left( \{H \in \Sub(G) \mid \delta^{-1} \gamma \delta \in H \} \right)
\\
&\, = \, \nu \left( \{\delta H \delta^{-1} \in \Sub(G) \mid \delta^{-1} \gamma \delta \in H \} \right) 
\\
&\, = \, \nu \left( \{H \in \Sub(G) \mid \gamma \in H \} \right)
\, = \, \varphi_\nu(\gamma).
\end{aligned}
$$
Moreover, $\varphi$ is normalized: $\varphi_\nu(e) = \nu(\Sub(G)) = 1$.
\par

It remains to check that $\varphi_\nu$ is positive definite.
By definition of $\varphi_\nu$, we have, for every $\gamma \in \Gamma$, 
$$
\varphi_\nu(\gamma) \, = \, \int_{\Sub(G)} \mathbf{1}_{H}(\gamma) d\nu (H)
$$
\par

Let $c_1, \hdots, c_n \in \C$ and $\gamma_1, \hdots, \gamma_n \in \Gamma$.
Since $\mathbf{1}_H$ is a function of positive type (see Proposition~\ref{diagcoeffinduced})
for every $H \in \Sub(G)$ and since
$\nu$ is a positive measure on $\Sub(G)$, we have
$$
\begin{aligned}
\sum_{i, j = 1}^n c_i \overline{c_j} \varphi_\nu(\gamma_j^{-1}\gamma_i)
\, &= \,
\sum_{i, j = 1}^n c_i \overline{c_j} \int_{\Sub(G)} \mathbf{1}_{H}(\gamma_j^{-1}\gamma_i) d\nu (H)
\\
\, &= \,
\int_{\Sub(G)} \left( \sum_{i, j = 1}^n
c_i \overline{c_j} \mathbf{1}_{H}(\gamma_j^{-1}\gamma_i) \right) d\nu (H)
\, \ge \, 0.
\end{aligned}
$$
\end{proof}

Let $N$ be a normal subgroup of the discrete group $\Gamma$.
The trace of $\Gamma$ corresponding to the IRS $\nu = \delta_{N}$ is $\Un_N$.
\par

\vskip.2cm

Let $\Gamma$ be a countable group and let 
$\Gamma \curvearrowright (X, \mathcal B, \mu)$ 
be a measure preserving action of $\Gamma$
on a standard Borel space $(X, \mathcal B)$,
equipped with a probability measure $\mu$ on $\mathcal B$.
As shown in Example~\ref{Exa-IRS} (3),
this action defines an IRS $\nu$ of $\Gamma$;
the associated trace $\varphi_\nu$ as in Proposition~\ref{Pro-IRS-Trace},
which we prefer to denote by $\varphi_\mu$, is given by 
$$
\varphi_\mu(\gamma) \, = \, \mu \left( X^\gamma \right)
\hskip.5cm \text{for all} \hskip.2cm
\gamma \in \Gamma,
$$
where $X^\gamma = \{x \in X \mid \gamma x = x \}$ is the set of fixed points of $\gamma$ in $X$.
\par

We are going to describe a concrete model for the GNS triple
associated to the trace $\varphi_\mu$.
For this, we need to discuss a variant of the group measure space 
construction of Murray and von Neumann from Section~\ref{SectionMSC}.
\index{Koopman representation}
\par

Let $\mathcal R$ be the equivalence relation on $X$, 
with the $\Gamma$-orbits as equivalence classes.
We view $\mathcal R$ as a subset of $X \times X$: 
$$
(x, y ) \in \mathcal R 
\hskip.5cm \text{if and only if} \hskip.5cm 
y \in \Gamma x .
$$
For $A \subset \mathcal R$ and $x \in X$, we set $A_x = A \cap (\{x\} \times X)$,
and let $\vert A_x \vert$ denote its cardinality.
We define a $\sigma$-finite measure $m$ 
on the Borel subsets $A$ of $\mathcal R$ by 
$$
m(A) \, = \, \int_X \vert A_x \vert d\mu(x) .
$$
There are two commuting representations $\pi$ and $\rho$ 
of $\Gamma$ on the Hilbert space $L^2(\mathcal R, m)$,
given by 
$$
(\pi(\gamma) f) (x, y) \, = \, f(\gamma^{-1} x, y)
\hskip.5cm \text{and} \hskip.5cm 
(\rho(\gamma) f) (x, y) \, = \, f(x, \gamma^{-1}y),
$$
for $\gamma \in \Gamma, \hskip.1cm f \in L^2(\mathcal R,m)$, and $(x, y) \in \mathcal R$.

Let $\Un_\Delta$ be the characteristic function
of the diagonal $\Delta = \{(x, x) \mid x \in X \}$ of $\mathcal R$.

\vskip.2cm

The following result is from \cite{Vers--11}.

\begin{theorem}
%[\textbf{Vershik}]
% 15.F.4
\label{Theo-Vershik}
Let $\Gamma \curvearrowright (X, \mu)$ be a measure preserving action 
of the countable group $\Gamma$ on a standard Borel space $X$
equipped with a probability measure $\mu$. Let 
$$
\varphi_\mu \, \colon \, \Gamma \to \C, \hskip.2cm
\gamma \mapsto \mu \left( \{x \in X \mid \gamma x = x \} \right)
$$
be the associated trace of $\Gamma$.
Let $R, m, \pi, \rho, \Delta$ be as above.
\par

Then $\mathbf{1}_\Delta$ belongs to $L^2(\mathcal R,m)$
and $(\rho, \Ki, \Un_\Delta)$
as well as $(\pi, \Ki, \Un_\Delta)$ 
are GNS--triples for $\varphi_\mu$,
where $\Ki$ is the closure of the linear span of $\rho(\Gamma)\mathbf{1}_\Delta$
as well as the closure of the linear span of $\pi(\Gamma)\mathbf{1}_\Delta$.
\end{theorem}

\begin{proof} 
We first show that $\Un_\Delta$
belongs to $L^2(\mathcal R, m)$ and that $\Vert \Un_\Delta \Vert = 1$.
Indeed,
$$
\int_\mathcal R \vert \Un_\Delta(x, y) \vert^2 dm(x, y)
\, = \, 
\int_X \vert \Delta_x \vert d\mu(x) 
\, = \, \mu(X) \, = \, 1.
$$
Moreover, for every $\gamma \in \Gamma$, we have
$$
\begin{aligned}
\langle \rho(\gamma) \Un_\Delta, \Un_\Delta \rangle
&\, = \, \int_\mathcal R \Un_\Delta (x, \gamma^{-1}y) \Un_\Delta (x, y) dm(x, y)
\\
&\, = \, \int_\mathcal R \Un_{ \{(x, x) \in X \times X \, \colon \, \gamma x = x \}}dm
\\
&\, = \, \mu \left( \{x \in X \mid \gamma x = x \} \right)
\, = \, \varphi_\mu(\gamma);
\end{aligned}
$$
a similar computation shows that
$$ 
\langle \pi(\gamma) \Un_\Delta, \Un_\Delta \rangle \, = \, \varphi_\mu(\gamma)
$$
Observe that 
$$
\rho(\gamma) \Un_\Delta \, = \, \Un_{ \{(x, \gamma x) \mid x\in X \}}
\hskip.5cm \text{and} \hskip.5cm 
\rho(\gamma) \Un_\Delta \, = \, \Un_{ \{(x, \gamma^{-1} x) \mid x\in X \}},
$$
so that
$$
\rho(\Gamma)\mathbf{1}_\Delta \, = \, \pi(\Gamma)\mathbf{1}_\Delta.
$$
Therefore, $(\rho, \Ki, \Un_\Delta)$ and $(\pi, \Ki, \Un_\Delta)$
are both GNS--triples for $\varphi_\mu$,
where $\Ki$ is the closure of the linear span of $\rho(\Gamma)\mathbf{1}_\Delta$.
\end{proof}
 
\begin{rem}
% 15.F.5
\label{Rem-RIS-Char}
One may ask which properties of the action $\Gamma \curvearrowright (X, \mu)$ 
imply that the associated function $\varphi \in \Tr_1(\Gamma)$
is indecomposable, that is, $\varphi \in E(\Gamma)$.
For a discussion of this question, see \cite{Vers--11} and \cite{DuGr--18}.
\end{rem}

Next, we introduce a notion which generalizes the notion of simplicity of a discrete group.

\begin{defn}
% 15.F.6
\label{Def-CharacterRigidGroup}
A discrete group $\Gamma$ is \textbf{character rigid}
if $E(\Gamma) = \{\delta_e, \Un_\Gamma \}$.
\end{defn}

\begin{rem}
% 15.F.7
\label{Rem-Def-CharacterRigidGroup}
Let $\Gamma$ be a character rigid group.
Then $\Gamma$ is simple;
indeed, let $N$ be a normal subgroup of $\Gamma$.
Since $\Un_N \in \Tr_1(\Gamma)$,
it follows that $\Un_N = t \delta_e + (1 - t) \Un_\Gamma$ for some $t \in \mathopen[ 0,1 \mathclose]$.
This implies that $N = \{e\}$ or $N = \Gamma$.
\end{rem}

The following result from \cite{DuMe--14} and \cite{PeTh--16}
shows that the ergodic actions of character rigid groups
are essentially free.
 
\begin{theorem}
%[\textbf{Dudko-Medynets, Peterson-Thom}]
% 15.F.8
\label{Theo-CharacterRigid}
Let $\Gamma$ be a countable character rigid group.
\par

Then every measure preserving ergodic action $\Gamma \curvearrowright (X, \mu)$
of $\Gamma$ on a standard Borel space $X$,
equipped with a non-atomic probability measure $\mu$, is essentially free,
that is, $\Gamma_x = \{e\}$ for $\mu$-almost every $x \in X$.
\end{theorem}

\begin{proof}
Since $\Gamma \curvearrowright (X, \mu)$ is ergodic and $\mu$ non-atomic,
the orbit $\Gamma x$ of $\mu$-almost every $x \in X$ is infinite,
by Proposition~\ref{Prop-Aperiodic}.
\par

Assume, by contradiction, that 
$$
\mu \left\{ x \in X \mid \Gamma_x \ne \{e\} \right\} \, > \, 0.
$$
Since $\Gamma$ is countable, 
there exists $\gamma \in \Gamma \smallsetminus \{e\}$ such that
$$
\mu \left( \left\{ x \in X \mid \gamma x \ne x \right\} \right) \, > \, 0.
$$
This means that, for the associated trace $\varphi_\mu \in \Tr(\Gamma)$ 
as in Theorem~\ref{Theo-Vershik}, we have $\varphi_\mu \ne \delta_e$. 
Since $E(\Gamma) = \{\delta_e, \Un_\Gamma \}$,
it follows that 
$$
\varphi_\mu \, = \, t \Un_\Gamma + (1-t) \delta_e
$$
for some $t > 0$.
\par

By Proposition~\ref{Theo-Vershik}, $(\rho, \Ki, \Un_\Delta)$ is a GNS triple for $\varphi$
for some closed subspace $\Ki$ of $L^2(\mathcal R, m)$.
Therefore, there exists a non-zero $\rho(\Gamma)$-invariant function 
in $L^2(\mathcal R,m)$.
It follows that there exists a $\rho(\Gamma)$-invariant Borel subset $A$
of $\mathcal R$ such that $0 < m(A) < \infty$.
\par

Recalling that $A_x = \{ y \in X \mid (x, y) \in \mathcal R \}$, set 
$$
Y \, := \, \left\{ x \in X \mid A_x \ne \emptyset \right\}.
$$
Since $m(A) > 0$, we have $\mu(Y) > 0$. 
Observe that the $\rho(\Gamma)$-invariance of $A$ means
that $A_x = \Gamma x$ for every $x \in Y$.
Therefore $\vert A_x \vert = \infty$ for almost every $x \in Y$, 
as almost every $\Gamma$-orbit is infinite. 
It follows that 
$$
m(A) \, = \, \int_X \vert A_x \vert dm(x) \, = \, \infty \cdot \mu(Y) \, = \, \infty ,
$$ 
and this is a contradiction.
\end{proof}

\begin{rem}
% 15.F.9
\label{Rem-Theo-CharacterRigid}
As proved in \cite{DuMe--14} and \cite{PeTh--16},
the conclusion of Theorem~\ref{Theo-CharacterRigid}
holds more generally if $\Gamma$ is a countable icc group
and if $E(\Gamma)$ consists only of $\delta_e$
and of countably many characters
all associated to finite dimensional irreducible representations.
\end{rem}

We draw the following consequence of Theorem~\ref{Theo-CharacterRigid}.

\begin{cor}
% 15.F.10
\label{Cor-Theo-CharacterRigid}
Let $\Gamma$ be a countable character rigid group.
\par

Then the ergodic IRS's of $\Gamma$ are $\delta_{\Gamma}$ and $\delta_{\{e\}}$.
\end{cor}

\begin{proof}
We may clearly assume that $\Gamma$ is not the group with one element.
\par

Since $\Gamma$ is countable, the compact space $\Sub(\Gamma)$ 
is metrizable and is therefore a standard Borel space.
Let $\mu$ be an ergodic $\Gamma$-invariant probability measure on $\Sub(\Gamma)$.
\par

We claim that $\mu$ is atomic.
Indeed, otherwise, $\Gamma \curvearrowright (\Sub(\Gamma), \mu)$
would be essentially free (Theorem~\ref{Theo-CharacterRigid})
and this is possible only if $\Gamma$ is the group with one element.
\par

Since $\mu$ is atomic, $\mu$ is a Dirac measure
and hence $\mu = \delta_{N}$ for some normal subgroup $N$ of $\Gamma$.
It follows from the simplicity of $\Gamma$ (see Remark~\ref{Rem-Def-CharacterRigidGroup})
that $\mu = \delta_{G}$ or $\mu = \delta_{\{e\}}$.
\end{proof}

\begin{exe}
% 15.F.11
\label{Exa-Theo-CharacterRigid}

 (1)
Let $\Gamma = \PSL_n(\K)$, where $\K$ is an infinite countable field and $n \ge 2$.
Then $\Gamma$ is a countable character rigid group (see Theorem~\ref{Theo-SLn}).
By Corollary~\ref{Cor-Theo-CharacterRigid}, 
the ergodic IRS's of $\Gamma$ are $\delta_{G}$ and $\delta_{\{e\}}$.

\vskip.2cm

(2)
Example (1) has been extended in \cite{Bekk} to more general algebraic groups as follows.
Let $G$ an algebraic group defined and simple over an infinite countable field $\K$,
and let $\Gamma = G(\K)^+$ be the subgroup of $G(\K)$
generated by the unipotent radicals of parabolic subgroups of $G$ defined over $\K$;
then $\Gamma$ is character rigid.

\vskip.2cm

(3)
Let $\Gamma$ be the commutator subgroup of 
a group from the Higman-Thompson families $F_{n,r}$ or $G_{n,r}$.
Then $\Gamma$ is character rigid ; see \cite{DuMe--14}.
\end{exe}

%-----------------------------------------------------------------------
% End of chapter 15
%-----------------------------------------------------------------------

\appendix
\renewcommand{\thechapter}{A\Alph{chapter}}
\renewcommand{\chaptername}{}

\chapter*{Appendix}
% Appendix
\label{Appendix}
\markboth{APPENDIX}{}

\emph{
Reading part of this book requires various notions that the potential reader
either may not know or may have forgotten.
The purpose of the present Appendix is to recall a few definitions
and state some ``well-known facts''
which will hopefully make it easier to follow other parts of the book. 
}

\section
{Topology}
% Appendix A.1
\label{AppTop}

\subsection*{Properties of topological spaces}

A topological space $X$ is T$_0$ if,
for any two points $x, y \in X$, $x \ne y$, there exists an open subset $U$ in $X$
such that $x \in U$ and $y \notin U$, or $x \notin U$ and $y \in U$.
Equivalently: if $\overline{\{x\}} \ne \overline{\{y \}}$
for any points $x, y \in X$ such that $x \ne y$.
\par

The space $X$ is T$_1$ if,
for any two points $x, y \in X$, $x \ne y$,
there exist open subsets $U, V$ in $X$
such that $x \in U$, $y \notin U$, $x \notin V$, $y \in V$.
Equivalently: if $\{x\}$ is closed in $X$ for every $x \in X$.
\par

A T$_2$ space, also known as a \textbf{Hausdorff space},
is a space $X$ satisfying Hausdorff's separation axiom:
for any two points $x, y \in X$, $x \ne y$, 
there exist disjoint open subsets $U, V$ in $X$
such that $x \in U$ and $y \in V$.
\index{Hausdorff topological space, T$_2$}
\par

It follows from the definitions that a T$_1$ space is T$_0$,
and a T$_2$ space is T$_1$.
\par

Bourbaki calls a T$_0$ topological space a \textbf{Kolmogorov space}
and a T$_1$ topological space an \textbf{accessible space}.
See Exercice 2 of Chapter I, \S~1
and Exercice 1 of Chapter I, \S~8 in \cite{BTG1--4}.
We recall the following proposition,
which appears in Exercice 27 of Chapter I, \S~8 of \cite{BTG1--4},
and in \cite[Section 6.2]{Will--07}.
It shows that every topological space admits a T$_0$ topological quotient space
satisfying an appropriate universal property.
\index{Topological space! separation axioms T$_0$, T$_1$, T$_2$}
\index{Separation axioms T$_0$, T$_1$, T$_2$}
\index{Topological space! Kolmogorov}
\index{Kolmogorov topological space, T$_0$}
\index{Topological space! accessible}
\index{Accessible topological space, T$_1$}
\index{Topological space! $\sigma$-compact}

\begin{prop}
% A.1.1
\label{Pro-UniversalKolmogorov}
Let $X$ be topological space.
There exists a pair $(Y,p)$ consisting of a T$_0$ topological space $Y$
and a surjective continuous map $p \,\colon X \twoheadrightarrow Y$ with the following property:
\par

for every pair $(Z,f)$ consisting of a T$_0$ topological space $Z$
and a continuous map $f \,\colon X \to Z$,
there exists a continuous map $g \,\colon Y \to Z$ such that $f = g \circ p$.
\end{prop}

\begin{proof}
We consider the equivalence relation
$\sim$ on $X$ defined by 
$$
x \sim x' 
\hskip.5cm \text{if} \hskip.5cm
\overline{\{x\}} = \overline{\{x' \}} .
$$
Let $Y$ be the quotient space $X/\sim$ equipped with the quotient topology
and
$$
p \, \colon \, X \twoheadrightarrow X/\sim
$$
the canonical projection.
\par

We claim that $X / \sim$ is a T$_0$ space.
Indeed, let $x, x' \in X$ be such that $p(x) \ne p(x')$.
We assume that $\overline{\{x\}} \not\subset \overline{\{x' \}}$
(the other case, $\overline{\{x' \}} \not\subset \overline{\{x\}}$,
can be dealt with similarly).
It is enough to check that there exists an open subset of $X/\sim$
which contains $p(x)$ and not $p(x')$.
\par

Set $U := X \smallsetminus \overline{\{x' \}}$.
Observe that every $v \in X$ such that $\overline{\{v\}} \not\subset \overline{\{x' \}}$
is necessarily a point in $U$.
\par

We claim that $U$
is saturated for the equivalence relation.
Indeed, for $u \in U$ and $u' \in X$ such that $u \sim u'$,
we have $\overline{\{u' \}} = \overline{\{u\}} \not\subset \overline{\{x' \}}$,
so that $u' \in U$.
It follows that $p(U)$ is an open subset of $X / \sim$
such that $p(x) \in p(U)$ and $p(x') \notin p(U)$.
This shows that $X / \sim$ is a T$_0$ space.

\vskip.2cm

Next, let $Z$ be T$_0$ topological space $Z$
and $f \,\colon X \to Z$ a continuous map.
Let $x, x' \in X$ such that $x \sim x' $. Then 
$$
f(x' ) \, \in \, f(\overline{\{x' \}}) \, = \, f(\overline{\{x\}}) \, \subset \, \overline{\{f(x)\}},
$$
by continuity of $f$, and hence $\overline{\{f(x' )\}} \subset \overline{\{f(x)\}}$.
Similarly, we have 
$\overline{\{f(x)\}} \subset \overline{\{f(x' )\}}$.
Therefore, $\overline{\{f(x)\}} = \overline{\{f(x')\}}$.
Since $Z$ is a T$_0$ topological space, it follows that $f(x) = f(x' )$.
Therefore the map $g \,\colon X / \sim \hskip.1cm \to Z, \hskip.2cm p(x) \mapsto f(x)$ is well-defined,
and $g \circ p = f$. 
\par

Moreover, $g$ is continuous.
Indeed, for every open subset $U$ of $Z$,
the set $p^{-1}(g^{-1}(U)) = f^{-1}(U)$ is open in $X$
and hence $g^{-1}(U)$ is open in $Y$.
\end{proof}
 
For a given topological space $X$, 
it follows from the universal property that the pair $(Y,p)$
as in Proposition~\ref{Pro-UniversalKolmogorov}
is unique in the following sense:
if $(Y',p')$ is a pair with the same properties as $(Y,p)$, 
then there exists a homeomorphism $f \,\colon Y \to Y'$ such that $f \circ p = p'$.
We will call $(Y,p)$ or $Y$ the \textbf{universal Kolmogorov quotient} of~$X$.
\par

Such spaces appear in Section \ref{PrimIdealSpace},
where $\Pri(G)$ is identified as the universal Kolmogorov quotient of $\widehat G$,
and in Section~\ref{Section-PrimIdealBS},
as quasi-orbit spaces.
\index{Universal Kolmogorov quotient}

\vskip.2cm

Given a continuous action of a group $G$ on a topological $X$
(see Definition \ref{appendixaction}),
the \textbf{quasi-orbit space} of this action
is the quotient space $X/\sim$ for the equivalence relation $\sim$
defined on $X$ by
$$
x \, \sim \, x'
\hskip.5cm \text{if} \hskip.5cm
\overline{Gx} \, = \, \overline{Gx'}.
$$
An equivalence class for $\sim$ is called a \textbf{quasi-orbit} for $G$.
The space $X/\sim$ equipped with the quotient topology is a T$_0$ topological space.
In fact, let $X/G$ be the space of $G$-orbits in $X$ equipped with the quotient topology.
Then 
\begin{center}
\emph{$X/\sim$ is the universal Kolmogorov quotient of $X/G$, as defined above.}
\end{center}
\index{Quasi-orbit}
\index{Quasi-orbit space}

\vskip.2cm

\index{Quasi-compact topological space}
A topological space is \textbf{quasi-compact} 
if every open covering contains a finite covering,
and is \textbf{compact} if it is both quasi-compact and Hausdorff.
It is \textbf{locally compact} if any of its points has a compact neighbourhood
and if it is Hausdorff.

\vskip.2cm

\index{Baire space}
\index{Topological space! Baire}
A topological space is a \textbf{Baire space}
if countable intersections of open dense subspaces are dense.
Locally compact spaces and completely metrizable spaces are Baire spaces
(Baire Category Theorem \cite[Chap.~9, \S~5]{BTG5--10}).
The quotient of a Baire space by an open equivalence relation
is a Baire space \cite[Chap.~9, \S~5, Exercise 25]{BTG5--10}.

\vskip.2cm

\index{Polish space, Polish topology}
A topological space is \textbf{Polish} if it is 
second-countable and homeomorphic to a complete metric space.
A \textbf{Polish group} is a topological group of which the topology is Polish.
For the following proposition, we refer to \cite[Chap.~9, \S~6, no~1]{BTG5--10};
% Alternative reference: \cite[Theorem 3.11]{Kech--95}.
it is used in the proof of our Proposition \ref{Pro-EtatsStandard}.
Recall that a G$_\delta$ set is the intersection of a countable family of open subspaces.

\begin{prop}
% A.1.2
\label{subspacePolish}
Let $X$ be a Polish space.
\begin{enumerate}[label=(\arabic*)]
\item\label{1DEsubspacePolish}
Every closed subspace of $X$ is Polish,
\item\label{2DEsubspacePolish}
every open subspace of $X$ is Polish,
\item\label{3DEsubspacePolish}
a subspace of $X$ is Polish if and only if it is a G$_\delta$ set.
\end{enumerate}
\end{prop}

Most topological spaces and groups of interest in this book are locally compact.
At several places, they have to satisfy the following restrictive conditions.
\par

\index{$\sigma$-compact topological space}
A topological space is \textbf{$\sigma$-compact}
if it is a countable union of compact subspaces.
In a locally compact space $X$ which is $\sigma$-compact,
there exists a nested sequence of compact subspaces
$K_0 \subset K_1 \subset \hdots \subset K_n \subset K_{n+1} \subset \hdots$
such that $X = \bigcup_{n \ge 0} K_n$,
and $K_n$ is contained in the interior of $K_{n+1}$ for all $n \ge 0$
\cite[Page I.68]{BTG1--4}.
\par

\index{Topological space! second-countable}
\index{Topological space! separable}
A topological space is \textbf{second-countable}
if its topology has a countable basis of open sets,
and \textbf{separable}
if it contains a countable dense subset.
It is straightforward that a second-countable space is separable.
The converse holds for metrizable spaces (see Theorem \ref{Urysohn}),
but it does not hold in general.
\par

For example, let $F$ be a finite group with at least two elements
and $G$ the direct product of copies of $F$ indexed by $\mathopen[0, 1 \mathclose]$;
then $G$ is a compact group which is separable \cite[Corollary 2.3.16]{Enge--89}
and which is not second-countable.
\par

We record the following equivalences,
which contain a particular case of Urysohn's metrization theorem
(see for example \cite[Page IX.21]{BTG5--10}).

\begin{theorem}[\textbf{Second-countable locally compact spaces}]
% A.1.3
\label{Urysohn}
For a locally compact space $X$, 
the following conditions are equivalent:
\begin{enumerate}[label=(\arabic*)]
\item\label{1DEUrysohn}
$X$ is second-countable;
\item\label{2DEUrysohn}
$X$ is metrizable and $\sigma$-compact;
\item\label{3DEUrysohn}
$X$ is metrizable and contains a countable dense subset;
\item\label{4DEUrysohn}
$X$ is Polish.
\end{enumerate}
\end{theorem}

Spaces which satisfy the conclusion of Proposition \ref{hereditarysigmac}
appear in Theorem \ref{regmeas2ndc} below.
Observe first that,
in a $\sigma$-compact locally compact space,
an open subset need not be $\sigma$-compact;
indeed, any non-$\sigma$-compact locally compact space
can be seen as an open subset of its one-point compactification.
However:

\begin{prop}
% A.1.4
\label{hereditarysigmac}
Let $X$ be a second-countable locally compact space.
\par

Then every open subset $U$ of $X$ is $\sigma$-compact.
\end{prop}

\begin{proof}
By Theorem \ref{Urysohn}, there exists a metric $d$ on $X$ which induces the topology of $X$.
For any positive integer $n$, set
$$
F_n \, = \, \{ x \in X \mid d(x, X \smallsetminus U) \ge 1/n \} .
$$
Then $F_n$ is closed in $X$, hence $F_n$ is $\sigma$-compact,
and $U$ is equal to the countable union $\bigcup_n F_n$,
Therefore $U$ is $\sigma$-compact.
\end{proof}

Because of Theorem \ref{Urysohn},
this can also be stated as follows:
\emph{in a second-countable locally compact space,
any open subset is $\sigma$-compact.}

\subsection*{Spaces of complex-valued continuous functions}

\index{$k1$@$C(X)$ space of continuous functions on $X$}
\index{$k2$@$C^b(X)$ space of bounded continuous functions on $X$}
\index{$k3$@$C^c(X)$ space of continuous functions on $X$ with compact support}
\index{$k4$@$C^{c, K}(X)$ space of continuous functions on $X$
with support contained in a compact subspace $K$ of $X$}
\index{$k5$@$C^0(X)$ space of continuous functions on $X$ vanishing at infinity}
For $X$ a topological space, we denote in this book by
\begin{enumerate}
\item[]
$C(X)$ the space of continuous complex-valued functions on $X$,
\item[]
$C^b(X)$ the subspace of bounded functions in $C(X)$,
\item[]
$C^c(X)$ the subspace of compactly supported functions in $C(X)$,
\item[]
$C^{c, K}(X)$ the subspace of $C(X)$
of functions with support
\item[] 
\hskip1.3cm
contained in a compact subspace $K$ of $X$.
\end{enumerate}
When $X$ is locally compact,
\begin{enumerate}
\item[]
$C^0(X)$ denotes the subspace of $C(X)$
of functions vanishing at infinity.
\end{enumerate}
All these spaces are commutative complex $*$-algebras
for the pointwise product and the involution defined by
$f^*(x) = \overline{f(x)}$ for all $x \in X$.
In special situations, they have other relevant structures;
for example, when $G$ is a locally compact group,
with a given left-invariant Haar measure,
$C^c(G)$ is also an algebra for the convolution product.
\index{Algebras! $1$@continuous functions, $C(X)$, $C^b(X)$, $\C^c(X)$, $C^0(X)$}
\par

When $X$ is locally compact, $C^0(X)$ is a commutative C*-algebra.
Conversely, the Gel'fand representation states that every commutative C*-algebra $A$
is isomorphic to $C^0(X)$,
where $X$ is the space of characters of $A$ equipped 
the weak$^*$-topology,
which is a locally compact space \cite[1.4.1]{Dixm--C*}.
(A character of $A$ is a non-zero morphism from $A$ to $\C$.)
\index{C*-algebra! $1$@commutative C*-algebra}

\section
{Borel spaces}
% Appendix A.2
\label{AppBorel}

Let $X$ be a set.
A \textbf{$\sigma$-algebra} of subsets of $X$
is a collection $\mathcal B$ of subsets of $X$
such that $X \in \mathcal B$,
$X \smallsetminus B \in \mathcal B$ for all $B \in \mathcal B$,
and $\bigcup_{n \ge 1} B_n \in \mathcal B$
for all sequences $(B_n)_{n \ge 1}$ of subsets $B_n$ in $\mathcal B$.
A \textbf{measurable space} is a set $X$ endowed with
a $\sigma$-algebra $\mathcal B$ of subsets;
a measurable space is also called a \textbf{Borel space}.
Let $(X, \mathcal B)$ and $(Y, \mathcal C)$ be two measurable spaces;
a map $f \,\colon X \to Y$ is a \textbf{measurable map}
if $f^{-1}(C) \in \mathcal B$ for all $C \in \mathcal C$;
a measurable map is also called a \textbf{Borel map}.
\index{Borel space}
\index{Measurable! $1$@space}
\index{Borel map}
\index{Measurable! $2$@map}

\vskip.2cm

\index{Borel space! countably separated}
\index{Borel space! countably generated}
Let $(X, \mathcal B)$ be a Borel space.
A subset $\mathcal S$ of $\mathcal B$ is \textbf{separating}
if, for every pair $(x, x')$ of distinct points in $X$,
there exists $S \in \mathcal S$ such that $x \in S$ and $x' \notin S$.
The space $(X, \mathcal B)$ is \textbf{countably separated}
if there exists in $\mathcal B$ a countable subset
which is separating,
and \textbf{countably generated}
if there exists in $\mathcal B$ a countable subset
which generates $\mathcal B$ and which is separating.
(We follow \cite{Mack--57} and much of the literature -- note however that
some authors define ``countably generating''
without requiring ``separating'', see for example \cite{Kech--95}.)
Note that
a countably generated Borel space is countably separated, by definition;
moreover,
countable separation and countable generation are inherited by subspaces,
as it is straightforward to check.

\begin{prop}
% A.2.1
\label{propcountablesg}
Let $(X, \mathcal B)$ be a Borel space.
\begin{enumerate}[label=(\arabic*)]
\item\label{1DEpropcountablesg}
$(X, \mathcal B)$ is countably separated if and only if there exists
an injective Borel map $X \to \mathopen[ 0, 1 \mathclose]$;
\item\label{2DEpropcountablesg}
$(X, \mathcal B)$ is countably generated if and only if is is isomorphic
to a subset of $\mathopen[ 0, 1 \mathclose]$.
\end{enumerate}
\end{prop}

\begin{proof}[Remarks, and references for the proof]
First, the reader should mark the difference
between ``subset'' here in \ref{2DEpropcountablesg}
and ``Borel subset'' in \ref{bDEstandard} below.
Then note that, because of Theorem \ref{KuratowskiTH},
the proposition could be stated with any standard Borel space
instead of $\mathopen[ 0, 1 \mathclose]$.
\par

For the proof of \ref{1DEpropcountablesg}, recall that $\{0,1\}^\N$
is Borel isomorphic to the Cantor set, a subspace of $\mathopen[ 0, 1 \mathclose]$;
details of the argument for example in \cite[Appendice A]{Zimm--84}.
For the proof of \ref{2DEpropcountablesg}, see \cite[Theorem 2.1]{Mack--57}.
\end{proof}

\index{Borel space! standard}
\index{Standard! $3$@Borel space}
Let $X$ be a topological space.
The \textbf{Borel $\sigma$-algebra} $\mathcal B_X$,
also denoted by $\mathcal B (X)$,
is the $\sigma$-algebra
generated by the open subsets of $X$.
For a Borel space $(X, \mathcal B)$, the following conditions are equivalent:
% (see \cite[Section 3]{Mack--57}): serait OK, mais doublon avec la suite
% bizarre : \cite[Proposition 12.1]{Kech--95}):
\begin{enumerate}[label=(\alph*)]
\item\label{aDEstandard}
$(X, \mathcal B)$ is Borel isomorphic to $(Y, \mathcal B_Y)$ for some Polish space $Y$;
\item\label{bDEstandard}
$(X, \mathcal B)$ is Borel isomorphic to $(C, \mathcal B_C)$ for some 
Borel subset $C$ of a Polish space~$Y$.
\end{enumerate}
The pair $(X, \mathcal B)$ is a \textbf{standard Borel space}
if these conditions are satisfied.
A standard Borel space is countably generated, 
but the converse does not hold.
For standard Borel spaces, see \cite[Section 3]{Mack--57}.
\par

It is obvious from the equivalence above that 
every Borel subset of  a standard Borel space is a standard Borel space.
\par

Standard Borel spaces can be classified as follows;
see \cite[Page 41]{Kura--66} or \cite[Theorem 15.6]{Kech--95}.
% or \cite[Theorem 3.3.13]{Sriv--98}.

\begin{theorem}[\textbf{Kuratowski}]
% A.2.2
\label{KuratowskiTH}
A standard Borel space is Borel isomorphic to one of the following spaces, 
depending on its cardinal:
\begin{enumerate}[label=(\arabic*)]
\item\label{1DEKuratowskiTH}
the interval $\mathopen[ 0,1 \mathclose]$,
\item\label{2DEKuratowskiTH}
the countable infinite set $\N$,
\item\label{3DEKuratowskiTH}
a finite space.
\end{enumerate}
\end{theorem}

\section
{Measures on Borel spaces and $\sigma$-finiteness}
% Appendix A.3
\label{AppMeasureB}

A \textbf{measure}, more precisely a \textbf{positive measure},
on a measurable space $(X, \mathcal B)$
is a map $\mathcal B \to \overline{\R}_+$ such that
$\mu(\bigcup_{n \ge 1} B_n) = \sum_{n \ge 1} \mu(B_n)$
for all sequences $(B_n)_{n \ge 1}$
of pairwise disjoint subsets $B_n$ in $\mathcal B$.
A \textbf{measure space} is a triple $(X, \mathcal B, \mu)$
where $(X, \mathcal B)$ is a measurable space
and $\mu$ a measure on $(X, \mathcal B)$.
Instead of ``a measure on $(X, \mathcal B)$'',
we also write ``a measure on $X$'', or ``a measure on $\mathcal B$'',
when this is clear enough from the context.
\par

A \textbf{finite measure} is a measure $\mu$ on $(X, \mathcal B)$
such that $\mu(X) < \infty$, and a \textbf{probability measure}
is one such that $\mu(X) = 1$.
\par

\index{Measure! $\sigma$-finite}
A measure space $(X, \mathcal B, \mu)$ is \textbf{$\sigma$-finite}
if there is a countable sequence $(B_n)_{n \ge 1}$ of subsets in $\mathcal B$
such that $X = \bigcup_{n \ge 1} B_n$ 
and $\mu(B_n) < \infty$ for all $n \ge 1$. 
This condition on a measure space is necessary for several important results,
for example for the Fubini theorem and the Radon--Nikodym theorem.
\par

An \textbf{atom} in a measure space $(X, \mathcal B, \mu)$
is a point $x \in X$ such that $\mu(\{x\}) > 0$.
The measure $\mu$ is \textbf{non-atomic} if there does not exist any atom in $X$.
\index{Atom! of a measure}
\index{Measure! non-atomic}

Given a measure $\mu$ on $X$,
two measurable functions $f, g$ defined on $X$ are
\textbf{equal $\mu$-almost everywhere} if 
$\mu\left( \{ x \in X \mid f(x) \ne g(x) \} \right) = 0$.

\vskip.2cm

\index{Integral}
Measurable functions with values either in $\overline{\R}_+$ or in a Banach space
can be integrated, as discussed in standard expositions,
for example in \cite[Chapter 1]{Rudi--66} and \cite[Chapter VI]{Lang--93}.
Integration of functions with values in Banach spaces
are important in this book, for example for the definition
of the representation of the group algebra $L^1(G)$
associated to a representation of a locally compact group $G$,
see Section \ref{C*algLCgroup}.
\par

More precisely, let $(X, \mathcal B, \mu)$ be a measure space
and let $E$ be a Banach space, real or complex;
the two cases $E = \R$ and $E = \C$ are particularly important.
First, any measurable function $f \,\colon X \to \overline{\R}_+$ has an \textbf{integral},
denoted by $\int_X f d\mu$, or $\int_X f(x) d\mu(x)$, or $\mu(f)$,
which is a number in $\overline{\R}_+$.
Then, any measurable function $g \,\colon X \to E$
such that $\int_X \Vert g(x) \Vert d\mu(x) < \infty$
can be integrated, and the integral $\int_X g d\mu$ is a vector in $E$.
% See for example \cite[Chapter VI]{Lang--93}, and in particular Corollary 5.9 there.
The integrals of two measurable functions which are equal $\mu$-almost everywhere
are equal.
It is abusive, but often convenient, to identify a measurable function
to its equivalence class modulo $\mu$-almost everywhere equality.

\vskip.2cm

Let $(X, \mathcal B, \mu)$ be a measure space
and $p$ a real number, $p \ge 1$.
For a measurable function $f \,\colon X \to \C$, set
$\Vert f \Vert_p := \big( \int_X \vert f \vert^p d\mu \big)^{1/p}$.
Denote by $L^p(X, \mu)$, or simply by $L^p(X)$,
the space of measurable functions $f \,\colon X \to \C$
(modulo $\mu$-almost everywhere equality)
such that $\Vert f \Vert_p < \infty$;
then $L^p(X)$ is naturally a Banach space.
For $p = 2$, the space $L^2(X)$ is naturally a Hilbert space.
For $p = 1$, note that $f \mapsto \int_X f d\mu$ defines a continuous linear form on $L^1(X)$.
\par

\index{Essential range! of a function}
\index{Essential supremum}
To define a Banach space $L^\infty(X)$, one may proceed as follows.
Let $(D_i)_i$ be a countable base of the topology of $\C$,
consisting of open discs.
Let $f \,\colon X \to \C$ be a measurable function.
Let $S_f$ be the union of the discs $D_i$
for which $f^{-1}(D_i)$ has zero measure.
The \textbf{essential range} of $f$ is the complement $R_f$ of $S_f$ in $\C$;
it is a closed subset of $\C$.
The \textbf{essential supremum} $\Vert f \Vert_\infty$ of $f$
is the supremum of the moduli $\vert z \vert$ of the points $z \in R_f$.
The function $f$ is \textbf{essentially bounded} if $\Vert f \Vert_\infty < \infty$.
Denote by $L^\infty(X, \mu)$, or simply by $L^\infty(X)$,
the space of (equivalence classes of) measurable functions $f \,\colon X \to \C$
such that $\Vert f \Vert_\infty < \infty$;
then $L^\infty(X)$ is naturally a Banach space.
See \cite[Chapter 3]{Rudi--66}.

\begin{rem}
% A.3.1
\label{DualitySpacesLp}
Let $p,q \in \mathopen[ 1, \infty \mathclose]$ be such that $\frac{1}{p} + \frac{1}{q} = 1$.
For all $f \in L^p(X, \mu)$ and $g \in L^q(X, \mu)$,
the function $fg$ is in $L^1(X, \mu)$,
hence there is a duality
$$
L^p(X, \mu) \times L^q(X, \mu) \to \C, \hskip.5cm (f,g) \mapsto \int_X fg d\mu .
$$
Consider $L^p(X, \mu)$ as a Banach space when $p \in \mathopen[ 1, \infty \mathclose[$,
and with the topology defined by this duality when $p = \infty$
(i.e., a posteriori, with the weak$^*$-topology);
Denote by $L^p(X, \mu)'$ the dual space of $L^p(X, \mu)$.
Assume that the measure $\mu$ is $\sigma$-finite. Then there is an isomorphism
$$
L^q(X, \mu) \, \to \, L^p(X, \mu)' , \hskip.5cm
g \, \mapsto \, \Big(f \mapsto \int_X f g d\mu \Big) .
$$
\par

When $p < \infty$, a proof can be read in \cite[Theorem 6.16]{Rudi--66}.
When $p = \infty$, it follows from general principles
(the dual of a locally convex topological vector space $E$ with a topology $\sigma(E, F)$
can be identified with $F$, see \cite[Chap.~II, \S~6, no~2, Proposition 3]{BEVT1--5}).
% Bourbaki, EVT, chapitre II, \S~6, no 2, Proposition 3, page EVT II.46 de mon \'ed-pdf [Springer 2007].
\end{rem}

Similarly, for $E$ a Banach space,
there is for each $p \in \mathopen[ 1, \infty \mathclose[$
a Banach space $L^p(X, \mu, E)$
of measurable functions $g \,\colon X \to E$ (up to equality $\mu$-almost everywhere)
such that $\Vert g \Vert_p := \big( \int_X \Vert g(x) \Vert_E d\mu(x) \big)^{1/p} < \infty$,
and a Banach space $L^\infty(X, \mu, E)$
of measurable functions $g$ such that the essential supremum $\Vert g \Vert_\infty$
of the function $X \ni x \mapsto \Vert g(x) \Vert_E \in \R_+$ is finite.
In case $E$ is a Hilbert space, $L^2(X, \mu, E)$ is also a Hilbert space;
it appears for example in Construction \ref{defpimupourGab++}.
% Construction 2.A.11 page 78

\subsection*{Radon--Nikodym derivatives and equivalence of measures}

\index{Disjoint! measures}
\index{Measure! disjoint}
\index{Measure! singular}
\index{Measure! absolutely continuous}

Let $\mu_1$ and $\mu_2$ be two positive measures
on a Borel space $(X, \mathcal B)$.
The measures $\mu_1$ and $\mu_2$
are \textbf{singular with each other}, or \textbf{disjoint},
if there exist disjoint sets $A_1, A_2 \in \mathcal B$ such that 
$\mu_1(X \smallsetminus A_1) = 0$ and $\mu_2(X \smallsetminus A_2) = 0$.
The measure $\mu_2$ is \textbf{absolutely continuous} with respect to $\mu_1$
if $\mu_2(B) = 0$ for every $B \in \mathcal B$ such that $\mu_1(B) = 0$.
The measures $\mu_1$ and $\mu_2$ are \textbf{equivalent}
if $\mu_2$ is absolutely continuous with respect to $\mu_1$
and $\mu_1$ absolutely continuous with respect to $\mu_2$.
\index{Measure! equivalent}
\index{Equivalent! measures}

\begin{theorem}[\textbf{Radon--Nikodym}]
% A.3.2
\label{Radon--NikodymTheorem}
Let $\mu_1, \mu_2$ be
two $\sigma$-finite positive measures on $(X, \mathcal B)$.
Suppose that $\mu_2$ is absolutely continuous with respect to $\mu_1$.
\par

Then there exists a measurable function $f \,\colon X \to \R_+$ such that
$$
\mu_2 (B) \, = \, \int_B f(x) d\mu_1(x)
\hskip.5cm \text{for all} \hskip.2cm
B \in \mathcal B .
$$
The function $f$ is uniquely determined by $\mu_1$ and $\mu_2$
up to equality $\mu_1$-almost everywhere.
\par

Moreover, $f$ is locally $\mu_1$-integrable, i.e., 
$\int_B f(x) d\mu_1(x) < \infty$
for all $B \in \mathcal B$ such that $\mu_1(B) < \infty$.
\end{theorem}

\index{Radon--Nikodym! derivative}
We refer to \cite[Chapter 6]{Rudi--66}.
In the situation of the Theorem, the function $f$ is called
the \textbf{Radon--Nikodym derivative} of $\mu_2$ with respect to $\mu_1$,
and it is denoted by $d\mu_2/d\mu_1$.
This is also written $\mu_2 = f \mu_1$, or $d\mu_2(x) = f(x) d\mu_1(x)$.
\par

Observe that $f = d\mu_2/d\mu_1$ is determined by the relation
$$
\int_X \varphi(x) d\mu_2 (x)
\, = \, \int_X \varphi(x) \dfrac{d\mu_2}{d\mu_1}(x) d\mu_1(x)
\leqno{(*)}
$$
for every measurable function $\varphi \,\colon X \to \R_+$.
\par

There is a chain rule: let $\mu_1, \mu_2, \mu_3$ be
three $\sigma$-finite positive measures on $(X, \mathcal B)$
such that $\mu_3$ is absolutely continuous with respect to $\mu_2$
and $\mu_2$ absolutely continuous with respect to $\mu_1$;
then $d\mu_3/d\mu_1 = (d\mu_3/d\mu_2) \hskip.1cm (d\mu_2/d\mu_1)$,
the equality being $\mu_1$-almost everywhere.
It follows that, if $\mu_1$ and $\mu_2$ are equivalent,
the Radon--Nikodym derivative $d\mu_2/d\mu_1$
can be viewed as a measurable function $X \to \R^\times_+$.
\par

\begin{theorem}[\textbf{Lebesgue decomposition theorem}]
% A.3.3
\label{Theo-LebesgueDec}
Let $\mu_1, \mu_2$ be two $\sigma$-finite positive measures on $(X, \mathcal B)$.
There is a unique pair $(\mu_{2,a}, \mu_{2,s})$
of $\sigma$-finite positive measures on $(X, \mathcal B)$
such that $\mu_2 = \mu_{2,a} + \mu_{2, s}$,
with $\mu_{2,a}$ absolutely continuous with respect to $\mu_1$
and $\mu_{2,s}, \hskip.1cm \mu_1$ disjoint.
\end{theorem}

It follows from Theorem~\ref{Theo-LebesgueDec}
that two $\sigma$-finite positive measures $\mu_1$ and $\mu_2$ on $(X, \mathcal B)$
are not singular with each other if and only if
there exists a non-zero measure $\mu$
which is absolutely continuous with respect to both $\mu_1$ and $\mu_2$.

\subsection*{Radon--Nikodym cocycle}

Let $(X, \mathcal B)$ be a Borel space,
equipped with a $\sigma$-finite measure $\mu$. 
The \textbf{direct image} of $\mu$ by a Borel map $\theta \,\colon X \to X$ 
is the measure $\theta_*(\mu)$ on $(X, \mathcal B)$ defined by
$$
\theta_*(\mu)(B) \, = \, \mu(\theta^{-1}(B))
\hskip.5cm \text{for all} \hskip.2cm
B \in \mathcal B .
$$
Equivalently, $\theta_*(\mu)$ is defined by 
$$
\int_X \varphi(x) d\theta_*(\mu) (x) \, = \, \int_X \varphi(\theta(x)) d\mu(x)
\leqno{(**)}
$$
for every measurable function $\varphi \,\colon X \to \R_+$.
\par

Let ${\mathrm Aut}(X, \mu)$ denote
the set of Borel automorphisms $\theta \,\colon X \to X$ such that 
$\theta_*(\mu)$ is equivalent to $\mu$.
It is obvious that $\theta_1 \circ \theta_2^{-1} \in {\mathrm Aut}(X, \mu)$
for $\theta_1, \theta_2 \in {\mathrm Aut}(X,\mu)$,
so that ${\mathrm Aut}(X, \mu)$ is a group.
\par

We define $d \,\colon {\mathrm Aut}(X,\mu) \times X \to \R^\times_+$ by 
$$
d(\theta,x) \, = \, \frac{ d (\theta^{-1})_* (\mu) }{d\mu} (x)
\hskip.5cm \text{for all} \hskip.2cm
x \in X, \hskip.1cm \theta \in {\mathrm Aut}(X,\mu).
$$
We have
$$
\int_X \varphi(\theta(x)) d(\theta,x) d\mu (x) \, = \, \int_X \varphi(x) d\mu(x)
\leqno{(***)}
$$
for every measurable function $\varphi \,\colon X \to \R_+$;
indeed, $\int_X \varphi(\theta(x)) d (\theta^{-1})_*(\mu) (x)$
is equal to the left-hand term by ($*$)
and to the right-hand term by ($**$).
The function $x \mapsto d(\theta, x)$ is determined by the identity ($***$),
as an identity for all measurable function $\varphi \,\colon X \to \R_+$.

\begin{prop}
% A.3.4
\label{cocycleRNtaumu}
Let $(X, \mathcal B)$ be a Borel space
and $\mu$ a $\sigma$-finite measure on $(X, \mathcal B)$.
\par

The function $d \,\colon {\mathrm Aut}(X,\mu) \times X \to \R^\times_+$ 
introduced above satisfies the following cocycle relation
$$
d(\theta_1\circ \theta_2,x) \, = \, d(\theta_1, \theta_2(x)) d(\theta_2, x)
\leqno{(\rm RN)}
$$
for all $\theta_1, \theta_2 \in {\mathrm Aut}(X,\mu)$ and $\mu$-almost all $x \in X$.
\end{prop}

\begin{proof}
Using relation ($***$) above for $\theta_1$ and for $\theta_2,$ we have
$$
\begin{aligned}
\int_X \varphi(x) d\mu(x)
\, &= \, \int_X \varphi(\theta_1(x)) d(\theta_1,x) d\mu(x)
\\
\, &= \, \int_X \varphi(\theta_1(\theta_2(x)) d( \theta_1, \theta_2(x)) d(\theta_2, x) d\mu(x)
\end{aligned}
$$
for every measurable function $\varphi \,\colon X \to \R_+$ and this proves the claim.
\end{proof}

Using the terminology introduced in Definition~\ref{defcocycle},
Proposition~\ref{cocycleRNtaumu} shows that 
$d \,\colon {\mathrm Aut}(X,\mu) \times X \to \R^\times_+$
is a cocycle over $(X,\mu)$ with values in $\R^\times_+$
for the tautological (left) action of the group ${\mathrm Aut}(X,\mu)$ on $X$.

\vskip.2cm

For reference, we record the following elementary fact.

\begin{prop}
% A.3.5
\label{PropEquivProbaMeasure}
Let $\mu$ be a non-zero $\sigma$-finite measure
on a Borel space $(X, \mathcal B)$.
\par

Then there exists a probability measure $\nu$ on $(X, \mathcal B)$
which is equivalent to $\mu$.
\end{prop}

\begin{proof}
Let $(B_n)_n$ be a sequence in $\mathcal B$
such that $X$ is the disjoint union $\bigcup_n B_n$,
and $0 < \mu(B_n) < \infty$ for every $n$. 
For each $n$, let $\mu_{n}$ be the probability measure defined by 
$$
\mu_{n} (B) \, = \, \frac{1}{\mu(B_n)} \mu(B \cap B_n)
\hskip.5cm \text{for all} \hskip.2cm
B \in \mathcal B.
$$
Let $(c_n)_n$ be a sequence of positive real numbers
such that $\sum_n c_n = 1$.
The measure $\nu = \sum_n c_n \mu_{n}$
is a probability measure on $(X, \mathcal B)$.
For $B \in \mathcal B$, we have $\nu(B) = 0$
if and only if $\mu(B \cap B_n) = 0$ for every $n \ge 1$,
that is, if and only if $\mu(B) = 0$.
Therefore $\nu$ is equivalent to $\mu$.
\end{proof}

We state two classical  results of von Neumann.
The first one is a version of his selection theorem \cite[Lemma 5]{vNeu--49};
we refer to the proof of Theorem Z.2 in Chap. 2, \S 2.2 in \cite{Mack--76}.
\par

A $\sigma$-finite measure $\mu$ on a Borel space $(X, \mathcal B)$
is said to be a \textbf{standard measure}
if there exists a Borel subset $N$ of $X$ with $\mu(N) = 0$
such that $X \smallsetminus N$ is a standard Borel space.
\index{Standard measure on a Borel space}

\begin{theorem}
% A.3.6
\label{Theo-vNSelection}
Let $X$ be a Borel space equipped with a standard Borel measure $\mu$
and let $Y$ be a standard Borel space.
Let $\Omega$ be a Borel subset 
of $X \times Y$ such that the restriction to $\Omega$
of the projection $X\times Y\to X$ is surjective.
\par

Then there exists a Borel subset $X_1$ of $X$
with $\mu(X \smallsetminus X_1) = 0$
and a Borel map $\varphi \,\colon X_1 \to Y$
such that $(x, \varphi(x)) \in \Omega$ for every $ x\in X_1$.
\end{theorem}

The second result is \cite[Satz 1]{vNeu--32}
and shows that a standard probability space
is characterized by its algebra of essentially bounded functions; 
we refer to the proof in \cite[Appendix IV]{Dixm--vN}.

\begin{theorem}
% A.3.7
\label{Theo-VonNeumann}
Let $X$ and $Y$ be standard Borel spaces,
equipped with $\sigma$-finite measures $\mu$ and $\nu$.
Assume that there exists an isomorphism of $*$-algebras 
$$
F \, \colon \, L^\infty (X, \mu) \to L^\infty (Y, \nu).
$$
\par

Then there exist Borel subsets $N, M$ of $X, Y$ respectively
such that $\mu(N) = 0$, $\nu(M) = 0$
and a Borel isomorphism $f \,\colon X \smallsetminus N \to Y \smallsetminus M$
such that $f_*(\mu)$ is equivalent to $\nu$
and such that $F(\varphi) (f(x)) = \varphi(x)$
for all $\varphi \in L^\infty (X, \mu), x \in X \smallsetminus N$.
\end{theorem}

\subsection*{Complex measures}

\index{Measure! complex}
\index{Measure! total variation}
Unless explicitly written otherwise, measures in this book are positive measures.
A \textbf{complex measure} on a Borel space $(X, \mathcal B)$
is a map $\mu \,\colon \mathcal B \to \C$ such that
$\mu(B) = \sum_i \mu(B_i)$ for all $B \in \mathcal B$
and all partitions $(B_i)_i$ of $B$ such that $B_i \in \mathcal B$ for all $i$
and $B_i \cap B_j = \emptyset$ for all $i, j$ such that $i \ne j$.
The \textbf{total variation} of such a complex measure $\mu$
is the finite positive measure $\vert \mu \vert$ defined by
$$
\vert \mu \vert (B) \, = \, \sup \sum_i \vert \mu (B_i) \vert
$$
for all $B \in \mathcal B$,
where the supremum is taken over all Borel partitions $(B_i)_i$ of $B$ as above.
The space of complex measures on $(X, \mathcal B)$
is a complex vector space,
with a norm defined by $\Vert \mu \Vert = \vert \mu \vert (X)$.
\par

There is a version of the Radon--Nikodym theorem which implies that
there exists a measurable function $h \,\colon X \to \C$
such that $\vert h(x) \vert = 1$ for all $x \in X$,
and $\mu = h \vert \mu \vert$,
i.e., $\mu(B) = \int_B h d\vert \mu \vert$ for all $B \in \mathcal B$.
For a measurable function $f \in L^1(X, \mathcal B, \vert \mu \vert)$, 
the integral $\int_X f d\mu$ is defined to be $\int_X f h d\vert \mu \vert$.
For complex measures, see \cite[Chapter 6]{Rudi--66}.

\section
{Radon measures on locally compact spaces}
% Appendix A.4
\label{AppMeasureLC}

Let $X$ be a topological space;
recall that the Borel $\sigma$-algebra $\mathcal B_X$
has been defined in Section \ref{AppBorel}.
Let $\mu$ be a measure on $(X, \mathcal B_X)$ such that
$$
\mu(K) \, < \, \infty
\hskip.5cm \text{for every compact subset} \hskip.2cm
K \subset X .
\leqno{(R_1)}
$$
The measure $\mu$ is \textbf{outer regular} if 
$$
\mu(B) \, = \, \inf \{ \mu(U) \mid
B \subset U \hskip.2cm \text{and} \hskip.2cm U \hskip.2cm \text{open} \}
\leqno{(R_2)}
$$
for all $B \in \mathcal B_X$,
and \textbf{inner regular} if
$$
\mu(B) \, = \, \sup \{ \mu(K) \mid B \supset K
\hskip.2cm \text{and} \hskip.2cm
K
\hskip.2cm \text{compact} \}
\leqno{(R_3)}
$$
for all $B \in \mathcal B_X$ such that $\mu(B) < \infty$
and for all $B$ which is open in $X$.
A \textbf{Radon measure} on a locally compact space $X$
is a measure on $(X, \mathcal B_X)$
which satisfies Conditions $(R_1)$, $(R_2)$, and $(R_3)$.
\index{Measure! Radon}
\index{Radon measure}
\par

Note that, if $X$ is a $\sigma$-compact locally compact space,
a Radon measure on $(X, \mathcal B_X)$ is $\sigma$-finite, by definition.

\vskip.2cm

For the hypothesis on $X$ in the following theorem,
see Proposition \ref{hereditarysigmac}.

\begin{theorem}
% A.4.1
\label{regmeas2ndc}
Let $X$ be a second-countable locally compact space.
Let $\mu$ a Borel measure on $(X, \mathcal B_X)$
which is finite on every compact subset of $X$.
Then
\begin{enumerate}[label=(\arabic*)]
\item\label{1DEregmeas2ndc}
the measure $\mu$ is a Radon measure;
\item\label{2DEregmeas2ndc}
the Hilbert space $L^2(X, \mu)$ is separable.
\end{enumerate}
\end{theorem}

\begin{proof}
For Claim \ref{1DEregmeas2ndc},
we refer to \cite[Theorem 2.18]{Rudi--66}.
Here is a proof for Claim \ref{2DEregmeas2ndc}.
\par

Since $X$ is locally compact and second-countable,
there exists a countable family $\{U_n\}_n$ of relatively compact open subsets of $X$
which is a base of the topology.
We are going to show that $\{1_{U_n}\}_{_n}$ is a total sequence in $L^2(X, \mu)$.
For this, consider $f \in L^2(X, \mu)$ such that $\langle f \mid 1_{U_n} \rangle = 0$ for every $n$.
It suffices to show that $f = 0$.
\par

Let $U$ be any open subset of $X$ of finite measure.
Since $U$ is a union of some of the $U_n$~'s,
we have $\int_U f(x)d\mu (x) = \langle f \mid \mathbf{1}_U \rangle = 0$.
Let $B$ a Borel subset of X of finite measure.
By outer regularity of $\mu$, there exists a sequence of open sets $(U_i)_i$
such that $B \subset U_i$ and $\mu(U_i) \le \mu(B) + \frac{1}{i}$ for all $i$.
Then
$\Vert \mathbf{1}_{B} - \mathbf{1}_{U_i} \Vert_2^2 = \mu(U_i \smallsetminus B) \le \frac{1}{i}$,
and so $\lim_i\mathbf{1}_{U_i} = \mathbf{1}_{B}$ in $L^2(X, \mu)$.
It follows that $\int_B f(x)d\mu (x) = 0$.
Since the linear span of
$\{ \mathbf{1}_{B} \mid B \in \mathcal{B}(X), \hskip.1cm \mu(B) < \infty \}$
is dense in $L^2(X, \mu)$, the claim follows.
\end{proof}

\begin{rem}
% A.4.2
\label{importantsecondcountable}
Locally compact spaces which are second-countable are of special importance.
This hypothesis of second-countability
implies that probability measures are Radon measures
and that Radon measures are $\sigma$-finite,
so that various fundamental theorems can be applied (Fubini, Radon--Nikodym).
Note also that the results of our Section \ref{Section-DecomposableOperators}
require separable Hilbert spaces,
so that applications to $L^2(\widehat G, \mu)$ for a LCA group $G$
require $\widehat G$ (and therefore $\G$) to be second-countable.
\par

There are known examples of compact spaces, necessarily not second-countable ones,
having probability measures which are not Radon measures
\cite[Volume II, Example 7.1.3]{Boga--07}.
\end{rem}

Let $X$ be a locally compact space, 
$\mu$ a Radon measure on $(X, \mathcal B_X)$,
and $C^c(X)$ the algebra of complex-valued continuous functions on $X$
with compact supports.
Then the map
$$
C^c(X) \, \to \, \C,
\hskip.5cm
f \, \mapsto \, \int_X f d\mu
$$
is well-defined, linear, and positive in the sense that $\int_X f d\mu \ge 0$
whenever $f(x) \ge 0$ for all $x \in X$.
Conversely:
% \cite[Theorem 2.14]{Rudi--66}:

\begin{theorem}[\textbf{Riesz representation theorem}]
% A.4.3
\label{Riesz}
Let $X$ be a locally compact space.
\par

Every positive linear form on $C^c(X)$
is of the form $f \mapsto \int_X f d\mu$
for a Radon measure $\mu$ on $(X, \mathcal B_X)$,
for all $f \in C^c(X)$.
\end{theorem}

For all this, see \cite[Chapter 2]{Rudi--66}.

\subsection*{Complex Radon measures}

\index{Radon measure! complex}
Since ``measure'' means ``positive measure'' in this book,
we also agree that ``Radon measure'' means ``positive Radon measure'',
unless explicitly written otherwise.
A \textbf{complex Radon measure} on a locally compact space $X$
is a complex measure on $(X, \mathcal B)$
of which the total variation is a Radon measure.
We denote by $M^b(X, \mathcal B_X)$, or for short by $M^b(X)$,
the space of complex Radon measures on $X$.
It is a complex vector space.
For the norm defined by $\Vert \mu \Vert = \vert \mu \vert (X)$ as in Appendix \ref{AppBorel},
it is a Banach space \cite[Exercice 6.3]{Rudi--66}.
\par

The following theorem is a version of Riesz representation theorem for complex Radon measures.
Since the total variation of a complex measure is a finite measure,
it is convenient to use
the space $C^0(X)$ of continuous functions on $X$ vanishing at infinity,
rather than the space $C^c(X)$ as in Theorem \ref{Riesz};
recall that $C^c(X)$ is dense in the Banach space $C^0(X)$.
For a proof of this theorem, we refer to \cite[Theorem 6.19]{Rudi--66}.

\begin{theorem}[\textbf{Riesz representation theorem for complex measures}]
% A.4.4
\label{RieszC}
Let $X$ be a locally compact space.
\par

To each bounded linear form $\Phi \,\colon C^0(X) \to \C$
there corresponds a unique complex Radon measure $\mu$ on $(X, \mathcal B_X)$
such that $\Phi (f) = \int_X f d\mu$ for all $f \in C^0(X)$.
Moreover $\Vert \Phi \Vert = \vert \mu \vert (X)$.
\par

The assignment $\Phi \to \mu$ provides an isomorphism
of the Banach space $M^b(X)$ of complex Radon measures on $(X, \mathcal B)$
with the dual of the Banach space $C^0(X)$.
\index{$k6$@$M^b(X)$ space of bounded complex measures on $(X, \mathcal B)$}
\end{theorem}

\section
{Groups and actions}
% Appendix A.5
\label{AppGpactions}

\index{Property! $\mathcal P$-by-$\mathcal Q$}
Let $\mathcal P$ and $\mathcal Q$ be group properties.
A group $\Gamma$ is $\mathcal P$-by-$\mathcal Q$
if it has a normal subgroup $N$ with Property $\mathcal P$
such that the quotient $\Gamma/N$ has Property $\mathcal Q$.
When $\mathcal P$ is a property inherited by subgroups
and $\mathcal Q$ is the property of being finite, 
a group with property $\mathcal P$-by-$\mathcal Q$
is also called a virtually $\mathcal P$ group.
In particular, a \textbf{virtually abelian} group is a group
which has an abelian subgroup of finite index.
\par

For a group $\Gamma$,
we denote by $\FC(\Gamma)$ its \textbf{FC-centre},
that is the normal subgroup of elements in $\Gamma$ 
with a finite conjugacy class.
The following result is classical.
\index{FC-centre $\FC(\Gamma)$ of a group $\Gamma$}

\begin{prop}
% A.5.1
\label{dueToSchur}
Let $\Gamma$ be group.
\begin{enumerate}[label=(\arabic*)]
\item\label{1DEduetoSchur}
If $\Gamma$ is an FC-group, i.e., if $\FC(\Gamma) = \Gamma$,
then the derived group $D(\Gamma)$ is finite.
\item\label{2DEduetoSchur}
If $\Gamma$ is finitely generated
and the FC-centre $\FC(\Gamma)$ is of finite index in $\Gamma$,
then $\Gamma$ is virtually abelian.
\end{enumerate}
\end{prop}

\begin{proof}
The result of Claim \ref{1DEduetoSchur} is due to Schur.
We refer to \cite[Theorem 10.1.4]{Robi--96}.
\par

For Claim \ref{2DEduetoSchur}, note first that
$\FC(\Gamma)$ is finitely generated,
as a finite index subgroup of a finitely generated group.
Choose a finite generating subset $S$.
For every $s \in S$, the centralizer $Z_{\FC(\Gamma)} (s)$ is of finite index in $\FC(\Gamma)$.
Then the subgroup $A := \bigcap_{s \in S} Z_{\FC(\Gamma)} (s)$
is of finite index in $\FC(\Gamma)$.
Moreover $A = \bigcap_{\gamma \in \FC(\Gamma)} Z_{\FC(\Gamma)} (\gamma)$
is the centre of $\FC(\Gamma)$.
It follows that $A$ is a normal abelian subgroup of finite index in $\Gamma$.
\end{proof}

\begin{rem}
% A.5.2
\label{FCetcrem}
(1)
The hypothesis of finite generation is necessary
in Proposition \ref{dueToSchur}~\ref{2DEduetoSchur}:
see Example \ref{HeisenbergFp}.

\vskip.2cm

(2)
A virtually abelian finitely generated group
need not have a centre of finite index.
For example, the infinite dihedral group is virtually abelian and has a trivial centre.
\end{rem}

\begin{defn}
% A.5.3
\label{deficc}
\index{Icc group}
A group $\Gamma$ is \textbf{icc}, or has infinite conjugacy classes,
if it is not $\{1 \}$ and if all its conjugacy classes distinct from $\{1 \}$ are infinite;
equivalently: if $\{1 \} = \FC(\Gamma) \ne \Gamma$.
\end{defn}

\begin{exe}
% A.5.4
\label{examplesicc}
(1)
Let $\Aff(\K)$ be the affine group over an infinite field $\K$,
as in Section \ref{Section-IrrRepAff}.
This group is icc, as it is easy to check
\index{Affine group! $4$@$\Aff(\K)$ of an infinite field $\K$}

\vskip.2cm

(2)
Similarly, for any integer $n \ge 2$,
the Baumslag--Solitar group $\BS(1,n)$ 
of Section \ref{Section-IrrRepBS} is icc.
\index{Baumslag--Solitar group $\BS(1, p)$}

\vskip.2cm

(3)
For any infinite group $H$ and any group $A \ne \{e\}$,
the restricted wreath product $H \wr A$ is icc.
This applies to
the lamplighter group $\Z \wr (\Z / 2 \Z)$ of \ref{Section-IrrRepLamplighter}.
\index{Lamplighter group}

\vskip.2cm

(4)
For an infinite field $\K$ and an integer $n \ge 2$,
the general linear group $\Gamma = \GL_n(\K)$
of Section \ref{Section-IrrRepGLN} is not type I.
The quotient group $\PGL_n(\K)$ can be shown to be icc,
for example because it is Zariski dense
in the appropriate algebraic group with trivial centre
(compare \cite[Proposition 2]{BeHa--94}).
\index{General linear group! $\GL_n(\K)$ with $\K$ a field}

\vskip.2cm

(5)
Examples of icc groups are very abundant.
Besides those occurring above, let us mention
non-abelian free groups and appropriate groups acting on trees
\cite{Corn--09},
lattices in semisimple connected real Lie groups 
without centre and without compact factors
(as a consequence of the Borel density theorem),
and weakly branch groups
\cite[Theorem 9.17]{Grig--11}.
\index{Free group}
\end{exe}

\begin{defn}
% A.5.5
\label{defcommgrps}
\index{Commensurate subgroups}
\index{Commensurator}
Let $\Gamma$ be a group.
Two subgroups $H, H'$ of $\Gamma$ are \textbf{commensurate}
if $H \cap H'$ has finite index in both $H$ and $H'$.
The \textbf{commensurator} of $H$ in $\Gamma$ is the subgroup
$$
\textnormal{Comm}_\Gamma (H) \, = \,
\{
\gamma \in \Gamma \mid H
\hskip.2cm \text{and} \hskip.2cm
\gamma^{-1} H \gamma
\hskip.2cm \text{are commensurate}
\}
$$
of $\Gamma$.
\end{defn}

\begin{prop}
% A.5.6
\label{defScommensurator}
Let $\Gamma$ be a group and $H$ a subgroup.
\begin{enumerate}[label=(\arabic*)]
\item\label{1DEdefScommensurator}
The commensurator $\Comm_\Gamma (H)$ is a subgroup of $\Gamma$ containing $H$.
\end{enumerate}
For $\gamma \in \Gamma$, the following properties are equivalent:
\begin{enumerate}[label=(\arabic*)]
\addtocounter{enumi}{1}
\item\label{2DEdefScommensurator}
$\gamma \in \Comm_\Gamma (H)$, i.e.,
$H \cap \gamma^{-1} H \gamma$ is of finite index
in both $H$ and $\gamma^{-1} H \gamma$;
\item\label{3DEdefScommensurator}
$H \cap \gamma^{-1} H \gamma$ and $H \cap \gamma H \gamma^{-1}$
are both of finite index in $H$;
\item\label{4DEdefScommensurator}
the $H$-orbits of the classes $H\gamma$ and $H \gamma^{-1}$
in $H \backslash \Gamma$ are both finite;
\item\label{5DEdefScommensurator}
there exists a finite index subgroup $K$ of $H$ such that
$H\gamma$ and $H \gamma^{-1}$
\par
are fixed by K.
\end{enumerate}
(Note that \ref{4DEdefScommensurator} and \ref{5DEdefScommensurator}
refer to the canonical right action $H \backslash \Gamma \curvearrowleft H$.)
\par

Consider a right transversal $T$ for $H$ in $\Gamma$
and the restriction to $H$ of the right action $T \curvearrowleft \Gamma$
defined in Construction \ref{constructionInd}(2).
Recall from there that, for $t \in T$,
the point $t^* \in T$ is defined by $Ht^{-1} = Ht^*$.
Then, for $t \in T$~:
\begin{enumerate}[label=(\arabic*)]
\addtocounter{enumi}{5}
\item\label{6DEdefScommensurator}
$t \in \Comm_\Gamma (H)$ if and only if the $H$-orbits of 
$t$ and $t^*$ in $T$ are both finite.
\end{enumerate}
Moreover $\Comm_\Gamma (H)$ depends on the commensurability class of $H$ only:
\begin{enumerate}[label=(\arabic*)]
\addtocounter{enumi}{6}
\item\label{7DEdefScommensurator}
If $H$ and $H'$ are two commensurate subgroups of $\Gamma$,
\par
then $\Comm_\Gamma (H) = \Comm_\Gamma (H')$.
\end{enumerate}
\end{prop}

There are similar equivalent properties
in terms of left actions $\Gamma \curvearrowright \Gamma/H$
and $\Gamma \curvearrowright S$, with $S$ a left transversal for $H$ in $\Gamma$.

\begin{proof}
We show that
\ref{3DEdefScommensurator} $\Leftrightarrow$ \ref{4DEdefScommensurator}
$\Leftrightarrow$ \ref{5DEdefScommensurator},
and leave the verification of the other claims to the reader.

Let $\gamma \in \Gamma$.
For the action $H \backslash \Gamma \curvearrowleft H$,
denote by $\text{Stab}_H (H\gamma)$ the stabilizer in $H$ of $H\gamma$,
i.e., $\text{Stab}_H (H\gamma) = H \cap \gamma^{-1} H \gamma$.
The orbit $(H\gamma)H$ is in canonical bijection with
the coset space
$\text{Stab}_H(H\gamma) \backslash H = (H \cap \gamma^{-1} H \gamma) \backslash H$,
and the orbit $(H \gamma^{-1}) H$ is in canonical bijection with
$(H \cap \gamma H \gamma^{-1}) \backslash H$.
The equivalence between \ref{3DEdefScommensurator} and \ref{4DEdefScommensurator} follows.
\par

A point $H\gamma$ in the space $H \backslash \Gamma$ has a finite $H$-orbit
if and only if there exists a subgroup of finite index $K$ in $H$ which fixes $H\gamma$.
The equivalence between \ref{4DEdefScommensurator} and \ref{5DEdefScommensurator} follows.
\end{proof}

\subsection*{Group actions}

\begin{defn}
% A.5.7
\label{appendixaction}
\index{Action of $G$ on $X$, respectively $G \curvearrowright X$ and $X \curvearrowleft G$}
An \textbf{action from the left} of a group $G$ on a set $X$
is a map $\alpha \,\colon G \times X \to X, \hskip.2cm (g,x) \mapsto gx$
such that $(g_1g_2)x = g_1(g_2x)$ and $ex = x$
for all $g_1, g_2 \in G$ and $x \in X$.
Such an action is equivalently described by the group homomorphism
$G \to \Sym (X)$ mapping $g \in G$
to the permutation $x \mapsto gx$ of $X$.
As agreed in Item \ref{notationaction},
we refer to such an action by the notation
$$
G \curvearrowright X .
$$
An \textbf{action from the right} of $G$ on $X$
is a map $X \times G \to X, \hskip.2cm (x,g) \mapsto xg$
such that $x(g_1g_2) = (xg_1)g_2 = x(g_1g_2)$ and $xe = x$
for all $g_1, g_2 \in G$ and $x \in X$.
Our corresponding notation is
$$
X \curvearrowleft G .
$$
For brevity, we often write ``action'' for ``action from the left''.
\end{defn}

\index{Measurable! $3$@action $G \curvearrowright (X, \mathcal{B})$}
Let $(X, \mathcal B)$ be a Borel space and $G$ a topological group.
A \textbf{measurable action} of $G$ on $(X, \mathcal B)$
is an action of $G$ by measurable automorphisms of $(X, \mathcal B)$
such that the map $G \times X \to X, \hskip.1cm (g, x) \mapsto gx$ is Borel,
where $G$ is equipped with the Borel structure induced by its topology
and $G \times X$ with the product Borel structure.

\begin{defn}
% A.5.8
\label{defqiergap}
Consider a measurable action of a topological group $G$
on a measurable space $(X, \mathcal B)$,
and a measure $\mu$ on $(X, \mathcal B)$.
\begin{itemize}
\setlength\itemsep{0em}
\item
The measure $\mu$ is \textbf{invariant} by $G$
if $\mu(gB) = \mu(B)$ for all $g \in G$ and $B \in \mathcal B$;
equivalently: if the image $g_*(\mu)$ is equal to $\mu$
for all $g \in G$.
\index{Invariant measure}
\index{Measure! invariant}
\item
The measure $\mu$ is \textbf{quasi-invariant} by $G$
if the set of $\mu$-null elements in $\mathcal B$ is $G$-invariant,
that is, if we have $\mu(g B) = 0$
for all $g \in G$
and $B \in \mathcal B$ such that $\mu(B) = 0$;
equivalently: if the image $g_*(\mu)$ is equivalent to $\mu$
for all $g \in G$.
\index{Quasi-invariant measure}
\index{Measure! quasi-invariant}
\item
Suppose moreover that $G$ is locally compact and second-countable,
and that $\mu$ is $\sigma$-finite and quasi-invariant by $G$.
The action of $G$ on $(X, \mathcal B, \mu)$ is \textbf{ergodic},
and the measure $\mu$ is \textbf{ergodic} for the action of $G$,
if, for every $G$-invariant subset $B \in \mathcal B$,
we have either $\mu(B) = 0$ or $\mu(X \smallsetminus B) = 0$.
Equivalently (see Proposition \ref{ergodicity1=2}),
this action is ergodic if, for every $B \in \mathcal B$
such that $\mu(gB \Delta B) = 0$ for all $g \in G$,
we have either $\mu(B) = 0$ or $\mu(X \smallsetminus B) = 0$.
\index{Ergodic! action}
\index{Action! ergodic}
\index{Ergodic! measure}
\index{Measure! ergodic}
\item
Suppose now that $\Gamma$ is countable and that $\mu$ is non-atomic.
The action is \textbf{aperiodic}
if the orbit $\Gamma x$ is infinite for $\mu$-almost all $x \in X$.
\index{Action! aperiodic}
\index{Aperiodic action}
\end{itemize} 
\end{defn}

Given an action of a group $G$ on a measure space $(X, \mathcal B, \mu)$,
a measurable function $f \,\colon X \to \C$ is \textbf{essentially $G$-invariant} if,
for every $g \in G$, the function $x \mapsto f(gx)$ is constant $\mu$-almost everywhere.

\begin{prop}
% A.5.9
\label{ergodicity1=2}
Let $G$ be a second-countable locally compact group
acting measurably on a Borel space $(X, \mathcal{B})$
and let $\mu$ be a quasi-invariant $\sigma$-finite measure on $(X, \mathcal B)$.
The following properties are equivalent:
\begin{enumerate}[label=(\roman*)]
\item\label{iDEergodicity1=2}
the action of $G$ on $(X, \mathcal B, \mu)$ is ergodic;
\item\label{iiDEergodicity1=2}
every essentially $G$- invariant measurable function $f \,\colon X \to \C$
is constant $\mu$-almost everywhere.
\end{enumerate}
\end{prop}

For the proof, see Theorem 1.3 in \cite{BeMa--00},
or Proposition \ref{prop-ergodicity} for the particular case
of a countable group $G$.
\par

The following fact is given here with a proof, for the reader's convenience.

\begin{prop}
% A.5.10
\label{Prop-Aperiodic}
Let $\Gamma$ be a countable group acting measurably
on a standard Borel space $(X, \mathcal B)$
and let $\mu$ be a non-atomic probability measure on $(X, \mathcal B)$
which is invariant and ergodic.
\par

Then the action $\Gamma \curvearrowright (X, \mathcal B, \mu)$
is aperiodic.
\end{prop} 

\begin{proof}
For an integer $n \ge 1$,
let $X_n$ be set of elements $x \in X$ such that $\vert \Gamma x \vert \le n$.
It is clear that $X_n$ is $\Gamma$-invariant.
\par

We claim that $X_n \in \mathcal B$.
Indeed, let $\mathcal P_n$ be the set of finite subsets $\Delta$ of $\Gamma$
such that $\vert \Delta \vert \le n$. We have
$$
X_n \, = \, \bigcup_{\Delta \in \mathcal P_n} \bigcap_{\gamma \in \Gamma}
\bigcup_{\delta \in \Delta} \{x \in X \mid \gamma \cdot x = \delta \cdot x\}.
$$
Since $\Gamma$ and $\mathcal P_n$ are countable and since 
$\{x \in X \mid \gamma \cdot x = \delta \cdot x\} \in \mathcal B$
for every $\gamma, \delta \in \Gamma$, the claim follows.
\par

We have to prove that 
$$
\mu \Big( \bigcup_{n \ge 1} X_n \Big) \, = \, 0.
$$
For this, it suffices to show that 
$\mu( X_n) = 0$ for every $n \ge 1$.
\par

Assume, by contradiction, that $\mu(X_n) \ne 0$ for some $n \ge 1$. 
It follows from the ergodicity of the action
$\Gamma \curvearrowright (X, \mathcal B, \mu)$
that $\mu(X_n) = 1$.
\par

Since $ (X, \mathcal B)$ is a standard Borel space,
we may assume that $X = \mathopen[ 0,1 \mathclose]$ and that $\mathcal B$
is the $\sigma$-algebra of all Borel subsets of $\mathopen[ 0,1 \mathclose]$. 
Define
$$
B \, := \, \{ x \in X_n \mid x \le \gamma \cdot x
\hskip.2cm \text{for all} \hskip.2cm
\gamma \in \Gamma\}
\, = \,
\bigcap_{\gamma \in \Gamma} \{x \in X_n \mid x \le \gamma \cdot x\}.
$$
Then $B \in \mathcal B$ and $X_n = \bigcup_{\gamma \in \Gamma} \gamma B$.
As $\mu(X_n)>0$, we have $\mu(B) > 0$ by invariance of $\mu$. 
Since $\mu$ is non-atomic,
there exists therefore a set $A \in \mathcal B$ such that $A \subset B$
and $0 < \mu(A) < \mu(B)$. Then
$$
A' \, := \, \bigcup_{\gamma \in \Gamma} \gamma A
$$
is a $\Gamma$-invariant set in $\mathcal B$ such that
$$
0 \, < \, \mu(A') \, < \, \mu(X_n) \, = \, 1.
$$
This is a contradiction to the ergodicity of the $\Gamma$-action on $X$.
\end{proof}

\section
{Locally compact groups}
% Appendix A.6
\label{AppLCG}

Let $G$ be a locally compact group.
There exists on $G$ a non-zero Radon measure $\mu_G$
which is left-invariant, that is,
such that $\mu_G(gB) = \mu_G(B)$ for all $g \in G$ and Borel subset $B$ of $G$,
and such a measure is unique up to multiplication by a positive real number;
this is called a left \textbf{Haar measure} on $G$,
or often, abusively, the Haar measure on $G$.
\index{Haar measure}
\par

\index{Modular function of a locally compact group}
For each $g \in G$, there is a positive number $\Delta(g)$ such that,
for each Borel subset $B$ of $G$, we have
$$
\mu_G(Bg) \, = \, \Delta(g) \mu_G(B) .
$$
The map
$$
\Delta \, \colon \, G \, \to \, \R_+^\times , \hskip.5cm g \, \mapsto \ \Delta(g)
$$
is a continuous homomorphism called the \textbf{modulus}, or \textbf{modular function}, of~$G$.
For $\varphi \in C_c(G)$, we have
$$
\begin{aligned}
\int_G \varphi(xg^{-1}) d\mu_G(g)
\, &= \,
\Delta_G(g) \int_G \varphi(x) d\mu_G (x) 
\\
\int_G \varphi(x^{-1}) d\mu_G(g)
\, &= \, \int_G \varphi(x) \Delta_G(x^{-1}) d\mu_G(x) 
\end{aligned}
$$
for all $g \in G$.
\par

The group $G$ is \textbf{unimodular} if $\Delta(g) = 1$ for all $g \in G$,
i.e., if the measure $\mu_G$ on $G$ is right-invariant (and so left and right-invariant).
For example, locally compact abelian groups, compact groups,
and semisimple Lie groups are unimodular.
More on modular functions and unimodular groups in \cite[A.3]{BeHV--08}.
\index{Unimodular group}

\vskip.2cm

\index{Convolution! of two measures on a LC group}
On the Banach space $M^b(G)$ of complex Radon measures on $G$
(see Theorem \ref{RieszC}),
there is a \textbf{convolution product} $(\nu, \sigma) \mapsto \nu \ast \sigma$,
where $\nu \ast \sigma$ is the direct image of the measure $\nu \times \sigma$
by the multiplication $G \times G \to G$;
alternatively, using the Riesz representation theorem, we have
$$
\int_G \varphi(x) d(\nu \ast \sigma) (x) \, = \,
\int_G \int_G \varphi(yz) d\nu(y) d\sigma(z)
\hskip.5cm \text{for all} \hskip.2cm
\varphi \in C^0(G) .
$$
It is easy to check that
$\Vert \nu \ast \sigma \Vert \le \Vert \nu \Vert \hskip.1cm \Vert \sigma \Vert$.
On $M^b(G)$, there is also an involution $\nu \mapsto \nu^*$
defined by
$$
\nu^*(B) \, = \, \overline{\nu(B^{-1})}
\hskip.5cm \text{for all Borel set} \hskip.2cm
B \in \mathcal B_G ;
$$
alternatively:
$$
\int_G \varphi(x) d(\mu^*) (x) \, = \, \int_G \overline{\varphi(x^{-1})} d\mu(x)
\hskip.5cm \text{for all} \hskip.2cm
\varphi \in C^0(G) .
$$
This makes $M^b(G)$ a Banach $*$-algebra called the
\textbf{measure algebra} of the locally compact group $G$.
It is an algebra with unit, which is the Dirac measure supported at the origin of $G$.
\index{Measure! algebra}
\index{Measure algebra}
\par

\index{Convolution! of two functions on a LC group}
Let $\varphi, \psi$ be two complex-valued measurable functions on $G$.
Their \textbf{convolution} is the function defined on $G$ by
$$
(\varphi \ast \psi) (g) \, = \, \int_G \varphi(gh) \psi(h^{-1}) d\mu_G(h)
\, = \, \int_G \varphi(h) \psi( h^{-1}g) d\mu_G(h)
\hskip.5cm \text{for all} \hskip.2cm
h \in G
$$
whenever these integrals make sense.
There is an involution defined on these functions by
$$
\varphi^*(g) \, = \, \Delta(g^{-1}) \overline{\varphi(g^{-1})} .
$$
For $\varphi, \psi \in L^1(G, \mu_G)$, the value $(\varphi \ast \psi) (g)$
is defined for $\mu_G$-almost all $g \in G$;
the resulting (equivalence class of the) function $\varphi \ast \psi$
is again in $L^1(G, \mu_G)$.
Similarly, $\varphi^*$ is in $L^1(G, \mu_G)$.
It is important to note that the convolution product $\ast$
depends on the choice of the Haar measure $\mu_G$.
\par

For $\varphi, \psi \in L^1(G, \mu_G)$,
the measures $\varphi \mu_G$ and $\psi \mu_G$ are in $M^b(G)$.
The convolution $\varphi \ast \psi$ can equivalently be defined by
$$
(\varphi \ast \psi) \mu_G \, = \, (\varphi \mu_G) \ast (\varphi \mu_G) ,
$$
and the involution by
$$
\varphi^* \mu_G \, = \, (\varphi \mu_G)^* .
$$
The Banach space $L^1(G, \mu_G)$ together with this convolution product and this involution
is called the \textbf{group algebra} of the locally compact group $G$.
The map
$$
L^1(G, \mu_G) \to M^b(G) , \hskip.5cm \varphi \mapsto \varphi \mu_G
$$
is an isometric morphism of Banach $*$-algebras;
it identifies $L^1(G, \mu_G)$ as a closed *-ideal of $M^b(G)$.
\index{$j4$@$L^1(G, \mu_G)$ or $L^1(G)$ group algebra of the locally compact group $G$}
\par

For a discrete group $\Gamma$ and the counting measure $\mu_\Gamma$,
we have $L^1(\Gamma, \mu_\Gamma) = M^b(\Gamma)$;
the notation is then $\ell^1(\Gamma)$ more often than $L^1(\Gamma, \mu_\Gamma)$.
For a locally compact group $G$ which is not discrete,
the group algebra $L^1(G, \mu_G)$ does not have a unit;
in particular $L^1(G, \mu_G)$ does not coincide with $M^b(G)$.
For an analysis of the complement $M^b(G) \smallsetminus L^1(G, \mu_G)$,
see Theorems 19.20 and 20.25 in \cite{HeRo--63}.
% pages 273 and 303
\par

For $\varphi, \psi \in C^c(G)$,
note that $\varphi \ast \psi$ and $\varphi^*$ have compact supports,
so that $C^c(G)$ is also a convolution algebra,
and a $*$-subalgebra of $L^1(G, \mu_G)$ and of $M^b(G)$.
Also, for two functions $\varphi, \psi$ with finite support,
the convolution has finite support, so that $\C[G]$ is a convolution algebra.
When $G$ is discrete, $\C[G]$ is a $*$-subalgebra of the group algebra $\ell^1(G)$.
But $\C[G]$ can be of use even when $G$ is not discrete,
see for example our Proposition \ref{constructionGNS2}.

\vskip.2cm

The following result can be found in many references.
For discussion and details, see, e.g., \cite[Appendix B]{BeHV--08}.
Quasi-invariance has been defined in Item \ref{defqiergap}.

\begin{prop}
% A.6.1
\label{qimeasG/H}
Let $G$ be a locally compact group, $H$ a closed subgroup,
and $G/H$ the corresponding $G$-space.
\begin{enumerate}[label=(\arabic*)]
\item\label{1DEqimeasG/H}
There exists a Radon measure $\mu$ on $G/H$ which is quasi-invariant by $G$.
The support of such a measure is the whole of $G/H$.
Two quasi-invariant Radon measures on $G/H$ are equivalent to each other.
\item\label{2DEqimeasG/H}
There exists a Radon measure on $G/H$ which is invariant by $G$
if and only if the modular function of $H$
is the restriction to $H$ of the modular function of~$G$.
\end{enumerate}
\end{prop}

\section
{Locally compact abelian groups and duality}
% Appendix A.7
\label{AppLCA+Pont}

Let $G$ be a locally compact abelian group, for short a \textbf{LCA group}.
Let $\widehat G$ its dual group,
which is the group of all homomorphisms from $G$ to $\T$,
with the compact-open topology;
it is itself a locally compact abelian group. 
Let $\mu_G$ be a Haar measure on~$G$.

\begin{rem}
% A.7.1
\label{topsurdualvN}
The topology on $\widehat G$ admits another description.
Elements of $\widehat G$ are bounded continuous functions from $G$ to $\C$,
so that $\widehat G$ embeds naturally in $L^\infty(G, \mu_G)$.
The restriction to $\widehat G$ of the weak$^*$-topology of $L^\infty(G, \mu_G)$
coincides with the previously defined topology of $\widehat G$
(see Proposition 1 in \cite[Chap.~II, \S~1, No 1]{BTS1--2}).
\end{rem}

For $g \in G$, denote by $\widehat g \in \widehat{\widehat G}$
the unitary character of $\widehat G$ defined by 
$$
\widehat g (\chi) \, = \, \chi(g)
\hskip.5cm \text{for} \hskip.2cm
\chi \in \widehat G.
$$

\begin{theorem}[\textbf{Pontrjagin duality}]
% A.7.2
\label{PontDuality17}
Let $G$ be a locally compact abelian group, and $\widehat G$ its dual.
\par

The map 
$$
G \mapsto \widehat{\widehat G},
\hskip.5cm
g \mapsto \widehat g
$$
is an isomorphism of topological groups.
\end{theorem}

\begin{proof}[References for proofs]
The theorem is essentially due to Pontrjagin.
A proof and historical notes can be found in Chapter VI of \cite{Weil--40}.
We also refer to Th\'eor\`eme 2 in \cite[Chap.~II, \S~1, No 5]{BTS1--2}.
\end{proof}

The theorem suggests that properties of $G$
are reflected by properties of $\widehat G$.
For example, and as shown in Proposition \ref{Prop-AbelianDiscreteCompactGroups},
a locally compact abelian group $G$
is compact if and only if $\widehat G$ is discrete.
Also:

\begin{prop}
% A.7.3
\label{Prop-DualAbelianSecondCountGroup}
Let $G$ be a locally compact abelian group.
Then $G$ is second-countable if and only if 
$\widehat G$ is second-countable.
\end{prop}

\begin{proof}
By Pontrjagin duality, it suffices to show that, if $G$ is second-countable,
then $\widehat G$ is second-countable.
\par

We consider $\widehat G$ as a subspace of $L^\infty(G, \mu_G)$,
and the latter with the weak$^*$-topology; see Remark \ref{topsurdualvN}. 
Fix a countable dense subset $Q$ of $\C$.
A basis for the topology of $\widehat G$ is given by 
finite intersections of sets of the form
$$
U_{f, q, n} \, = \,
\left\{ \chi \in \widehat G \mid \vert \chi(f)-q \vert < 1/n \right\},
$$
with $f \in L^1(G, \mu_G)$, $q \in Q$, and an integer $n \ge 1$,
where we have written $\chi(f)$ for $\int_{G} \chi(g)f(g) d\mu_G(g)$. 
\par

Since $G$ is second-countable, $L^1(G, \mu_G)$ is separable;
so, there exists a dense sequence $(f_k)_{k \ge 1}$ in $L^1(G, \mu_G)$.
The family of finite intersections of the sets
$$
U_{f_k, q, n} \, = \,
\left\{ \chi \in \widehat G \mid \vert \chi(f_k)-q \vert < 1/n \right\}
$$
for $k, \hskip.1cm n \ge 1$ and $q \in Q$, 
is countable and is clearly a basis for the topology of $\widehat G$.
\end{proof}

Let $\mu$ be a complex Radon measure on $G$.
The \textbf{Fourier--Stieltjes transform} of $\mu$
is the function on $\widehat G$, denoted by $\mathcal F (\mu)$ or $\widehat \mu$, defined by 
\index{Fourier--Stieltjes transform}
$$
\mathcal F (\mu) (\chi) \, = \, \int_{G} \overline{\chi(g)} \hskip.1cm d\mu(g)
\hskip.5cm \text{for all} \hskip.2cm
\chi \in \widehat G .
$$
The \textbf{inverse Fourier--Stieltjes transform} of $\mu$
is the function on $G$, denoted by $\overline{\mathcal F}(\nu)$ on $G$, defined by 
$$
\overline{\mathcal F}(\mu) (\chi) \, = \, 
\mathcal F (\mu) (\chi^{-1}) \, = \,
\int_G \widehat g (\chi) d\mu(\chi) \, = \,
\int_G \chi(g) d\mu(\chi).
$$
(Our notation agrees with that of \cite{BTS1--2}.)
Both $\mathcal F (\mu)$ and $\overline{\mathcal F}(\mu)$
are bounded continuous complex-valued functions on $\widehat G$.
If $\check{\mu}$ denotes the direct image of a measure $\mu$ on $G$
by the map $G \to G, \hskip.1cm g \mapsto g^{-1}$,
observe that $\overline{\mathcal F} (\mu) = \mathcal F (\check{\mu})$.
\par

\index{Fourier transform}
Let $f \in L^1(G, \mu_G)$; observe that $f \mu_G$ is a complex measure on $G$.
The \textbf{Fourier transform} of the function $f$,
denoted by $\mathcal F (f)$ or $\widehat f$,
is the Fourier--Stieltjes transform of the measure $f \mu_G$, that is,
$$
\mathcal F (f) (\chi) \, = \, \int_{G} \overline{\chi(g)} f(g) d\mu_G(g)
\hskip.5cm \text{for all} \hskip.2cm
\chi \in \widehat G.
$$
Similarly, the \textbf{inverse Fourier transform}
$\overline{\mathcal F}(f)$ of $f$
is the inverse Fourier transform of $f \mu_G$, that is
$$
\overline{\mathcal F} (f) (\chi) \, = \,
\int_G \chi(g) f(g) d\mu_G (\chi)
\hskip.5cm \text{for all} \hskip.2cm
\chi \in \widehat G.
$$
\par

For the next proposition,
we refer to \cite[Chap.~II, \S~1]{BTS1--2}, or to \cite[Chapter 4]{Foll--16}.
Recall that $M^b(G)$ is a commutative Banach $*$-algebra for convolution,
and $C^b(G)$ a commutative C*-algebra for the pointwise product.

\begin{prop}
% A.7.4
\label{Pro-InjectivityFourier}
Let $G$ be a locally compact abelian group.
\begin{enumerate}[label=(\arabic*)]
\item\label{1DEPro-InjectivityFourier}
For $\mu \in M^b(G)$ a complex Radon measure on $G$,
the Fourier--Stieltjes transform of $\mu$
and the inverse Fourier--Stieltjes transform of $\mu$
are in $C^b(\widehat G)$.
\item\label{2DEPro-InjectivityFourier}
The mappings $\mathcal F \,\colon M^b(G) \to C^b(\widehat G)$
and $\overline{ \mathcal F } \,\colon M^b(G) \to C^b(\widehat G)$
are morphisms of Banach $*$-algebras.
\item\label{3DEPro-InjectivityFourier}
For $f \in L^1(G, \mu_G)$,
the Fourier--Stieltjes transform of $f$
and the inverse Fourier--Stieltjes transform of $f$
are in $C^0(\widehat G)$.
\item\label{4DEPro-InjectivityFourier}
The mappings $\mathcal F \,\colon M^b(G) \to C^b(\widehat G)$
and $\overline{ \mathcal F } \,\colon M^b(G) \to C^b(\widehat G)$
are injective.
\end{enumerate}
\end{prop}

The next theorem characterizes the images of positive measures.
We refer to \cite[4.23]{Foll--16} or \cite[D.2.2]{BeHV--08}.

\begin{theorem}[\textbf{Bochner Theorem}]
% A.7.5
\label{TheoremBochner}
Let $G$ be a locally compact abelian group.
\par

The map $\mu \mapsto \mathcal F (\mu)$
establishes a bijection
from the space $M^b(\widehat G)_+$ of finite positive Radon measures on the dual of $G$
onto the space $P(G)$ of functions of positive type on $G$
(as defined in Section \ref{S-FPosType}).
\end{theorem}

Let $\mu_{\widehat G}$ be a Haar measure on $\widehat G$.
For $f \in L^1(G, \mu_G) \cap L^2(G, \mu_G)$, 
the Fourier transform $\widehat f$ is in 
$L^1(\widehat G, \mu_{\widehat G}) \cap L^2(\widehat G, \mu_{\widehat G})$.
Recall that $L^1(G, \mu_G) \cap L^2(G, \mu_G)$ is dense in $L^2(G, \mu_G)$.
The map $f \mapsto \widehat f$ extends to a continuous linear map
from $L^2(G, \mu_G)$ onto $L^2(\widehat G, \mu_{\widehat G})$,
denoted again by $f \mapsto \widehat f$.
Moreover, there is a choice of the Haar measure $\mu_{\widehat G}$ such that
$\Vert \widehat f \Vert_{L^2(\widehat G, \mu_{\widehat G})} =
\Vert f \Vert_{L^2(G, \mu_G)}$ for all $f \in L^2(G, \mu_G)$;
the Haar measure $\mu_{\widehat G}$ on $\widehat G$ is said to be
\textbf{associated} to the Haar measure $\mu_G$ on $G$.
For this and for the next Theorem, see \cite[Chap.~II, \S~1, no~3]{BTS1--2}.

\begin{theorem}[\textbf{Plancherel Theorem}]
% A.7.6
\label{PlancherelTh}
Let $G$ be a locally compact abelian group, and $\widehat G$ its dual.
Let $\mu_G$ be a Haar measure on $G$
and $\mu_{\widehat G}$ the Haar measure on $\widehat G$
which is associated to $\mu_G$.
\par

Then the Fourier transform
$$
L^2(G, \mu_G) \, \to \, L^2(\widehat G, \mu_{\widehat G}) ,
\hskip.5cm
f \mapsto \widehat f ,
$$
is an isometry of the Hilbert space $L^2(G, \mu_G)$
onto the Hilbert space $L^2(\widehat G, \mu_{\widehat G})$.
\end{theorem}

\index{Algebras! $2$@trigonometric polynomials on $\widehat G$}
\index{Trigonometric polynomials on $\widehat G$}
Let $G$ be a locally compact abelian group
and $C^b(\widehat G)$ the algebra
of bounded continuous complex-valued function on its dual.
The \textbf{algebra of trigonometric polynomials} on $\widehat G$
is the linear span of $\{ \widehat g \mid g \in G \}$;
it is a subalgebra of $C^b(\widehat G)$.
\par

Note that $\mathrm{Trig} (\widehat G)$
is a dense subalgebra of the algebra $AP (\widehat G)$
of almost periodic functions on $\widehat G$
for the topology of uniform convergence.
Indeed, $\mathrm{Trig} (\widehat G)$ is a point separating $*$-subalgebra
of the algebra $C(\mathrm{Bohr} (\widehat G))$ of continuous functions
on the Bohr compactification of $\widehat G$,
and the claim follows from the Stone--Weierstra\ss\ theorem.
(For $AP (\widehat G)$ and $\mathrm{Bohr} (\widehat G)$,
see Section~\ref{SectionFdrep}).
\par

Given a finite Radon measure $\mu$ on $\widehat G$,
there is a canonical map
from the algebra $C^b(\widehat G)$ to the algebra $L^\infty(\widehat G, \mu)$.
The image by this map of $\mathrm{Trig} (\widehat G)$
is naturally a subalgebra of $L^\infty(\widehat G, \mu)$,
and also a subspace of $L^p(\widehat G, \mu)$
for every $p \ge 1$.
\par

The following lemma is used on several occasion
(\ref{Prop-RepAbCyclic}, \ref{Prop-RepAb-Equiv}, \ref{Pro-IndRepGMS}).
% 2.A.6 et 2.B.2 et 14.A.2

\begin{lem}
% A.7.7
\label{Lemma-DensityFourierTransform}
Let $G$ be a second-countable locally compact abelian group
and $\mu$ a probability measure on its dual $\widehat G$.
\begin{enumerate}[label=(\arabic*)]
\item\label{1DELemma-DensityFourierTransform}
The image of the algebra $\mathrm{Trig} (\widehat G)$ is norm dense
in $L^p(\widehat G, \mu)$ for every $p \in \mathopen[1,+\infty \mathclose[$.
\item\label{2DELemma-DensityFourierTransform}
The image of the algebra $\mathrm{Trig} (\widehat G)$ is weak*-dense
in $L^\infty(\widehat G, \mu)$.
\end{enumerate}
\end{lem}

\begin{proof}
Given $p \in \mathopen[ 1,+\infty \mathclose]$,
denote by $q \in \mathopen[ 1,+\infty \mathclose]$
the conjugate exponent, defined by $\frac{1}{p} + \frac{1}{q} = 1$.
\par

Let $L^p(\widehat G, \mu)$ be equipped with the norm topology when $p<+\infty$
and with the weak$^*$-topology when $p = +\infty$.
The topological dual space of $L^p(\widehat G, \mu)$
may be identified in the usual way with $L^q(\widehat G, \mu)$.
Therefore, by Hahn--Banach Theorem, it suffices to show the following
\par

\emph{Claim: If $f \in L^q(\widehat G, \mu)$ is such that}
$$
\int_{\widehat G} \widehat g (\chi) f(\chi) d\mu(\chi) \, = \, 0
\hskip.5cm \text{for all} \hskip.2cm
g \in G ,
\leqno{(\dagger)}
$$
\emph{then $f = 0$.}
\par

Since $\mu$ is a probability measure, 
$f \in L^q(\widehat G, \mu) \subset L^1(\widehat G, \mu)$,
and $(\dagger)$ can be written $\overline{\mathcal F}(f \mu) = 0$.
By Proposition \ref{Pro-InjectivityFourier},
we have therefore $f \mu = 0$,
i.e., $f(\chi) = 0$ for $\mu$-almost all $\chi \in \widehat G$,
hence $f = 0$.
\end{proof}

\section
{Hilbert spaces and operators}
% Appendix A.8
\label{AppHspacesop}

\index{Dimension! of a Hilbert space}
In a Hilbert space, a Schauder basis, or simply a basis,
is a maximal family of orthonormal vectors.
Every Hilbert space has a basis (as a consequence of Zorn's lemma)
and two basis have the same cardinality.
By definition this cardinality is the \textbf{dimension} of the Hilbert space.
% When the dimension is finite, it coincides with the vector space dimension
% defined as the cardinality of a Hamel basis, as in linear algebra.
% When the dimension is infinite, \emph{it does not};
% for example, in a Hilbert space of countable infinite dimension as defined above,
% every Hamel bais in uncountable (this can be shown by using Baire's Category Theorem).

\vskip.2cm

\index{Strong topology on $\Li(\Hi)$ and $\U(\Hi)$}
\index{Weak topology on $\Li(\Hi)$ and $\U(\Hi)$}
\index{Ultra-weak topology on $\Li(\Hi)$ and $\U(\Hi)$}
Let $\Hi$ be a Hilbert space. We denote by $\Li(\Hi)$ the involutive algebra
of bounded linear operators on $\Hi$.
The algebra $\Li(\Hi)$ has several useful topologies.
The \textbf{strong topology}
is the locally convex topology defined by the set of seminorms
$x \mapsto \Vert x \xi \Vert$, for $\xi \in \Hi$.
The \textbf{weak topology}
is the locally convex topology defined by the set of seminorms
$x \mapsto \vert \langle x \xi \mid \eta \rangle \vert$, for $\xi, \eta \in \Hi$.
The \textbf{ultra-weak topology}
is the locally convex topology defined by the set of seminorms
$x \mapsto \left\vert \sum_{i=1}^\infty \langle x\xi \mid \eta_i \rangle \right\vert$,
for sequences $(\xi_i)_{i \ge 1}$ and $(\eta_i)_{i \ge 1}$ of vectors in $\Hi$
such that $\sum_{i=1}^\infty \Vert \xi_i \Vert^2 < \infty$
and $\sum_{i=1}^\infty \Vert \eta_i \Vert^2 < \infty$.
\par

In $\Li(\Hi)$, the unit ball is weakly compact;
more generally, bounded subsets of $\Li(\Hi)$
which are weakly closed are weakly compact.
On bounded subsets of $\Li(\Hi)$,
the weak topology and the ultra-weak topology coincide;
this implies that bounded subsets of $\Li(\Hi)$
which are ultra-weakly closed are ultra-weakly compact.
\par

On $\Li(\Hi)$, there is also the \textbf{norm topology}, defined by the norm
$\Vert x \Vert = \sup_{\xi \in \Hi, \Vert \xi \Vert \le 1} \Vert x\xi \Vert$.
This makes it a Banach algebra, and indeed a C*-algebra.
\par

For all this, see \cite[Chap.~I, \S~3]{Dixm--vN}.

\vskip.2cm

Denote by
$$
\U(\Hi) \, = \, \{ x \in \Li (\Hi) \mid x^*x = xx^* = \mathrm{Id}_{\Hi} \}
$$
the group of unitary operators on $\Hi$.
The restrictions to $\U(\Hi)$ of the strong, weak, and ultra-weak topologies coincide,
and they make $\U(\Hi)$ a topological group;
this is the topology to consider for the theory of unitary representations.
% \par
The restriction to $\U(\Hi)$ of the norm topology
is of little interest for the subject of this book because,
for many topological groups $G$, the only homomorphism
$G \to (\U(\Hi), {\rm norm \hskip.1cm topology})$ is the trivial homomorphism
\cite{Sing--52}.

\vskip.2cm

\index{Positive element in a C*-algebra}
\index{$e3$@$\M_+, A_+$ positive cone}
For a bounded linear operator on a Hilbert space, the following conditions are equivalent:
\begin{enumerate}
\item[(P1)]
$x^* = x$ and the spectrum of $x$ is contained in the set $\R_+$ of non-negative real numbers;
\item[(P2)]
there exists $y \in \Li (\Hi)$ such that $x = y^*y$;
\item[(P3)]
there exists $z \in \Li (\Hi)$ such that $z^* = z$ and $x = z^2$.
\item[(P4)]
$\langle x \xi \vert \xi \rangle \ge 0$ for all $\xi \in \Hi$.
\end{enumerate}
The operator $x$ is called \textbf{positive} if it satisfies these conditions.
The set $\Li (\Hi)_+$ of positive operators is a closed convex cone
such that $\Li (\Hi)_+ \cap (-\Li (\Hi)_+) = \{0\}$.
There is an order on positive operators,
defined by $x \le y$ if $y-x$ is positive.
For all this, see, e.g., \cite[1.6]{Dixm--C*}.
\par

For a positive operator $x$ on $\Hi$,
functional calculus provides a positive operator, denoted by $\sqrt x$, or $x^{\frac{1}{2}}$,
of which the square is $x$; it is the \textbf{square root} of~$x$.
Note that, in $\Li (\Hi)$ with the strong topology,
$\sqrt x$ is the limit of polynomials in $x$.
\index{Square root of a positive operator}
\index{Positive element in a C*-algebra! square root}
\par

\index{Absolute value of an operator on a Hilbert space}
For $x \in \Li (\Hi)$, the square root of $x^*x$
is the \textbf{absolute value} of $x$, denoted by~$\vert x \vert$.

\vskip.2cm

\index{Partial isometry}
\index{Partial isometry! $1$@initial space} 
\index{Partial isometry! $2$@final space}
Let $\Hi_1, \Hi_2$ be two Hilbert spaces.
A bounded linear operator $u \,\colon \Hi_1 \to \Hi_2$
is a \textbf{partial isometry}
if there exists a closed subspace $\Ki$ of $\Hi_1$ such that
the restriction of $u$ to $\Ki$ is an isometry
and such that $u = 0$ on $\Ki^{\perp}$.
The subspace $\Ki = (\ker u)^{\perp}$ is the \textbf{initial space}
and its image $u(\Ki)$ is the \textbf{final space} of the partial isometry~$u$.
Note that $u^*u$ is the orthogonal projection of $\Hi_1$ onto $\Ki$,
and $uu^*$ is the orthogonal projection of $\Hi_1$ onto $u(\Ki)$.

\vskip.2cm

\index{Polar decomposition of an operator}
Every operator $x \,\colon \Hi_1 \to \Hi_2$ has a \textbf{polar decomposition}
$x = u \vert x \vert$, where $u$ is the partial isometry with initial space
$(\ker x)^{\perp}$ and final space the closure of the image of $x$,
defined by
$$
u\xi \, = \, 0
\hskip.2cm \text{if} \hskip.2cm
\xi \in \ker x
\hskip.5cm \text{and} \hskip.5cm
u\xi \, = \, x\xi
\hskip.2cm \text{if} \hskip.2cm
\xi \in (\ker x)^\perp ,
$$
and where $\vert x \vert = (x^*x)^{1/2}$.
Moreover $u^*x = \vert x \vert$.

\section
{Projection-valued measures}
% Appendix A.9
\label{AppProjValMeas}

\index{Projection-valued measure}
\index{$h7$@$\Proj(\Hi)$ projections in $\Li (\Hi)$}
Let $(X, \mathcal B)$ be a Borel space and $\Hi$ be a Hilbert space.
Denote by $\Proj(\Hi)$ the set of orthogonal projections in $\Li (\Hi)$. 
A \textbf{projection-valued measure} on $(X, \mathcal B)$ and $\Hi$
is a map $E \,\colon \mathcal B \to \Proj(\Hi)$ such that
\begin{enumerate}[label=(\arabic*)]
\item\label{1DESNAG}
$E(\emptyset) = 0$ and $E(X) = I$;
\item\label{2DESNAG}
$E(B \cap B') = E(B)E(B')$ for all $B, B' \in \mathcal B$;
\item\label{3DESNAG}
if $(B_i)_{i \in I}$ is a sequence of pairwise disjoint sets
in $\mathcal B $, then
$$
E(\bigcup_{i \in I} B_i) \, = \, \sum_{i \in I} E(B_i) ,
$$
where the sum converges in the strong topology.
\end{enumerate}
We also say that $E$ is on $X$ and $\Hi$ if $\mathcal B$ is clear from the context.
\par

Projection-valued measures are also called spectral measures \cite{Halm--51},
and resolutions of the identity \cite{Rudi--73}.
% , for example in \cite{Rudi--73} or in the books by Dunford and Schwartz,
It can easily be shown that
Condition \ref{2DESNAG} is redundant,
as it follows from \ref{1DESNAG} and \ref{3DESNAG};
see, e.g., \cite[\S~36]{Halm--51}.
\par

Let $E \,\colon \mathcal B \to \Proj(\Hi)$ be a projection-valued measure as above.
For every $\xi, \eta \in \Hi$,
a (signed) bounded Borel measure $\mu_{\xi, \eta}$ is defined on $X$ by
$$
\mu_{\xi, \eta}(B) \, = \, \langle E(B) \xi \mid \eta\rangle
\hskip.5cm \text{for all} \hskip.2cm
B \in \mathcal B.
$$
In the case of a locally compact space $X$
and its $\sigma$-algebra $\mathcal B$ of Borel sets, 
it is customary to define a \textbf{regular} projection-valued measure
by requiring that all these measures $\mu_{\xi, \eta}$ are Radon measures;
this is automatic when $X$ is a second-countable locally compact space,
by Theorem \ref{regmeas2ndc}.
\index{Regular projection-valued measure}
\index{Projection-valued measure! regular}
\par

When $\xi = \eta$, the measure $\mu_{\xi, \eta}$ is denoted by $\mu_\xi$.
If $\Vert \xi \Vert = 1$, then $\mu_\xi$ is a probability measure.

\begin{constr}
% A.9.1
\label{pvmfc}
To a projection-valued measure $E \,\colon \mathcal B \to \Proj(\Hi)$,
there is an associated functional calculus, i.e., 
a morphism of an appropriate algebra of functions on $X$
to the algebra of operators $\Li (\Hi)$,
defined as follows.
\par

Denote by ${\rm Bor}^b(X)$ the commutative $*$-algebra of
complex-valued bounded measurable functions on $X$,
with the norm defined by $\Vert f \Vert = \sup_{x \in X} \vert f(x) \vert$;
it is a Banach $*$-algebra.
For every $f \in {\rm Bor}^b(X)$,
there exists a unique operator $\Psi(f) \in \Li (\Hi)$ such that 
$$
\langle \Psi(f)\xi \mid \eta \rangle \, = \, \int_X f(x)d \mu_{\xi, \eta}(x)
\hskip.5cm \text{for all} \hskip.2cm
\xi, \eta \in \Hi.
$$
Moreover, the map
$$
{\rm Bor}^b(X) \, \to \, \Li (\Hi) ,
\hskip.5cm
f \, \mapsto \, \Psi(f)
$$
is a homomorphism of Banach $*$-algebras.
The operator $\Psi(f)$ is usually denoted by 
$$
\int_X f(x) dE(x).
$$
For all this, see Chapter 12 in \cite{Rudi--73}.
\end{constr}

\begin{exe}
% A.9.2
\label{spectraltheorem}
Let $S$ a normal operator on a Hilbert space $\Hi$.
Let $X$ denote its spectrum, which is a non-empty compact subset of the complex plane.
One formulation of the spectral theorem is that there exists a unique
projection-valued measure $E \,\colon \mathcal B (X) \to \Proj (\Hi)$
such that $S = \int_X z dE(z)$.
\end{exe}

Let $X$ be a topological space
and let $E$ be a projection-valued measure on $(X, \mathcal B_X)$ and some Hilbert space.
The \textbf{support} of $E$ is the set of $x \in X$ 
such that, for every neighbourhood $U$ of $x$, we have $E(U) \ne 0$;
equivalently, the support of $E$ is the complement of the largest open set 
$U$ of $X$ such that $E(U) = 0$.
For $x \in X$, we say that $x$ is an \textbf{atom} of $E$ if $E(\{x\}) \ne 0$.
Obviously, every atom of $E$ lies in the support of $E$.
\index{Support of a projection-valued measure}
\index{Projection-valued measure! support}
\index{Atom! of a projection-valued measure}
\index{Projection-valued measure! atom}

\section
{C*-algebras}
% Appendix A.10
\label{AppAlgC*}

\index{Algebras! $5$@C*-algebra}
\index{C*-algebra!}
A \textbf{C*-algebra} is a complex Banach algebra $A$
with a conjugate linear involution $x \mapsto x^*$
such that $\Vert x^*x \Vert = \Vert x \Vert^2$ for all $x \in A$.
Though we expect from the reader some familiarity with this notion,
we do recall a few of the relevant definitions, as follows.

\vskip.2cm

Let $A$ be a C*-algebra. 
A \textbf{projection}
in $A$ is an element $p \in A$ such that $p = p^* = p^2$.
\index{Projection in a C*-algebra}

\vskip.2cm

\index{Positive element in a C*-algebra}
\index{$e3$@$\M_+, A_+$ positive cone}
\index{Positive element in a C*-algebra! square root}
For $x \in A$, the analogues of Conditions (P1), (P2), (P3) of \ref{AppProjValMeas}
are equivalent
An element $x \in A$ is \textbf{positive} if it satisfies these conditions.
The set $A_+$ of positive elements in $A$ is a closed convex cone
such that $A_+ \cap (-A_+) = \{0\}$,
and every positive element in $A$ has a unique positive square root ion $A$
(see again \cite[1.6]{Dixm--C*}).
There is an order on positive elements of $A$,
defined for $x, y \in A_+$ by $x \le y$ if $y - x \in A_+$.

\vskip.2cm

\index{Approximate identity}
Let $A$ be a complex normed algebra.
An \textbf{approximate identity} in $A$
is a family $(u_\iota)_{\iota \in I}$ of elements of $A$
indexed by a directed set $I$ such that
\begin{enumerate}[label=(\arabic*)]
\item
$\Vert u_\iota \Vert \le 1$ for all $\iota \in I$;
\item
$\lim_{\iota \in I} \Vert u_\iota x - x \Vert =
\lim_{\iota \in I} \Vert x u_\iota - x \Vert = 0$
for all $x \in A$.
\end{enumerate}
For example, let $L^1(G, \mu_G)$ be the group algebra of a LC group $G$,
with Haar measure~$\mu_G$.
Let $\mathcal V$ be a base of neighbourhood of $e$ in $G$.
For each $V \in \mathcal V$,
there exists a continuous function $u_V \,\colon G \to \R_+$
such that $u_V(g^{-1}) = u_V(g)$ for all $g \in G$, 
with support compact and contained in $V$,
and such that $\int_G u_V(g) d\mu_G(g) = 1$.
Then $(u_V)_V$ is an approximate identity in $L^1(G, \mu_G)$
\cite[Proposition 2.44]{Foll--16}.
\par

Let $A$ be moreover a C*-algebra.
An \textbf{increasing approximate identity} in $A$ 
is an approximate identity $(u_\iota)_{\iota \in I}$ such that
\begin{enumerate}[label=(\arabic*)]
\addtocounter{enumi}{2}
\item
$u_\iota \ge 0$ for all $\iota \in I$,
\item
$u_\iota \le u_{\iota '}$ for all $\iota, \iota ' \in I$ with $\iota \le \iota '$.
\end{enumerate}
Every C*-algebra has an increasing approximate identity.
More precisely, given a dense two-sided ideal $\mathfrak m$ in $A$,
there exists an increasing approximate identity
$(u_\iota)_{\iota \in I}$ in $A$ 
such that $u_\iota \in \mathfrak m$ for all $\iota \in I$.
If $A$ is separable, there exists such an approximate identity
indexed by $\N$ \cite[1.7.1]{Dixm--C*}.

\section
{Von Neumann algebras}
% Appendix A.11
\label{AppAlgvN}

\index{Commutant}
Let $\Hi$ be a Hilbert space,
and $\Li (\Hi)$ the involutive algebra of bounded linear operators on $\Hi$,
as in \ref{S-DefUnitD}.
The \textbf{commutant} of a subset $S$ of $\Li (\Hi)$ 
is $S' := \{ x \in \Li (\Hi) \mid xs = sx \hskip.2cm \text{for all} \hskip.2cm s \in S \}$,
and the \textbf{bicommutant} is $S'' := (S')'$.
Note that $S'$ and $S''$ are subalgebras of $\Li (\Hi)$,
and that $S \subset S''$.
For $S, T$ two subsets of $\Li (\Hi)$,
the inclusion $S \subset T$ implies $T' \subset S'$.
In particular, $S \subset S''$ implies $(S'')' \subset S' \subset (S')''$,
so that $S''' = S'$.
\par

\index{von Neumann algebra}
\index{$e1$@$\M$ von Neumann algebra}
\index{Algebras! $6$@von Neumann algebra}
A \textbf{von Neumann algebra}
acting on $\Hi$ is an involutive subalgebra $\mathcal M \subset \Li (\Hi)$
that is equal to its bicommutant $\mathcal M''$.
For example, $S'$ and $S''$ are von Neumann algebras
for any selfadjoint subset $S$ of $\Li(\Hi)$.
In particular, $\Li (\Hi)$ and $\C \mathrm{Id}_{\Hi}$ are von Neumann algebras acting on $\Hi$.
\par

A \textbf{homomorphism} of a von Neumann algebra $\mathcal M$ acting on $\Hi$
to a von Neumann algebra $\mathcal N$ acting on $\Ki$
is a linear map $\Phi \,\colon \mathcal M \to \mathcal N$ such that
$\Phi(ST) = \Phi(S)\Phi(T)$ and $\Phi(S^*) = \Phi(S)^*$ for all $S, T \in \mathcal M$.
For such a homomorphism $\Phi$, it follows that
$\Phi(\mathcal M_+) \subset \Phi(\mathcal N_+)$,
that the image by $\Phi$ of a projection is a projection,
that $\Vert \Phi(S) \Vert \le \Vert S \Vert$ for all $S \in \mathcal M$,
and when $\Phi$ is injective that $\Vert \Phi(S) \Vert = \Vert S \Vert$ for all $S \in \mathcal M$
(Chap.~I, \S~1, no~5 in \cite{Dixm--vN}).
\par

\index{Normal positive linear map from a von Neumann algebra to another}
Let $\mathcal M, \mathcal N$ be two von Neumann algebras.
A linear map $\Phi \,\colon \mathcal M \to \mathcal N$
is positive if $\Phi(\mathcal M_+) \subset \mathcal N_+$,
and \textbf{normal} if, moreover, $\Phi(x)$ is the least upper bound of $(\Phi(x_\iota))_{\iota \in I}$
for every increasing generalized sequence $(x_\iota)_{\iota \in I}$ in $\mathcal M$
with least upper bound $x$.
A normal positive linear map is ultraweakly continuous \cite[Chap~I, \S~4, no~3]{Dixm--vN}.

\vskip.2cm

For the next theorem, see Chap~I, \S~3, no~4 in \cite{Dixm--vN}.

\begin{theorem}[\textbf{von Neumann bicommutant theorem}]
% A.11.1
\label{thmvNbicommutant}
Let $\Hi$ be a Hilbert space
and $\Li (\Hi)$ the $*$-algebra of bounded operators on $\Hi$.
\par

Let $\mathcal M$ a $*$-algebra containing the identity operator of $\Hi$
and $\mathcal M''$ its bicommutant.
The following conditions are equivalent:
\begin{enumerate}[label=(\roman*)]
\item\label{iDEthmvNbicommutant}
$\mathcal M'' = \mathcal M$, i.e., $\mathcal M$ is a von Neumann algebra in $\Li (\Hi)$;
\item\label{iiDEthmvNbicommutant}
$\mathcal M$ is closed in $\Li (\Hi)$ for the weak operator topology;
\item\label{iiiDEthmvNbicommutant}
$\mathcal M$ is closed in $\Li (\Hi)$ for the strong operator topology.
\end{enumerate}
\par

The von Neumann algebra generated by a selfadjoint subset $S$ of $\Li (\Hi)$
is equal to the bicommutant $S''$.
\end{theorem}

\begin{defn}
% A.11.2
\index{von Neumann algebra! $1$@factor}
\index{Factor}
A von Neumann subalgebra $\mathcal M$ of $\Li (\Hi)$ is a \textbf{factor}
if the centre of $\mathcal M$ consists only of the scalar multiples of $\mathrm{Id}_{\Hi}$,
that is, if 
$$
\mathcal M \cap \mathcal M' \, = \, \{\lambda \mathrm{Id}_{\Hi} \mid \lambda \in \C\}.
$$ 
Observe that, $\mathcal M$ is a factor if and only if $\mathcal M'$ is a factor.
\end{defn}

\begin{exe}
% A.11.3
For each positive integer $n$, the algebra $\Li(\C^n)$ is a factor,
by definition a factor of type I$_n$; its dimension is finite, equal to $n^2$.
For an infinite dimensional Hilbert space $\Hi$, the algebra $\Li(\Hi)$ is a factor,
by definition a factor of type I$_\infty$.
\par

It is an insight of Murray and von Neumann that there exist other types of factors,
in particular factors of type II$_1$, which are infinite dimensional factors, say $\mathcal M$,
such that there exists a linear form $t \,\colon \mathcal M \to \C$
such that $t(xy) = t(yx)$ for all $x,y \in \mathcal M$ and $t(1) = 1$.
There is more on types of von Neumann algebras and of factors
in Section~\ref{SectionvN}.
\end{exe}

\begin{prop}
% A.11.4
\label{Prop-ProjvN}
Let $\Hi$ be a Hilbert space and $\mathcal M$ a subalgebra of $\Li (\Hi)$.
Denote by $\mathcal P$ the set of the orthogonal projections $P \in \Li (\Hi)$
which belong to $\mathcal M$. Then:
\begin{enumerate}[label=(\arabic*)]
\item\label{1DEProp-ProjvN}
$\mathcal M$ is the von Neumann algebra generated by $\mathcal P$,
i.e., $\mathcal P'' = \mathcal M$;
\item\label{2DEProp-ProjvN}
$\mathcal M$ is abelian if and only if $PQ = QP$ for all $P, Q \in \mathcal P$;
\item\label{3DEProp-ProjvN}
$\mathcal M$ is a factor if and only if $0$ and $\mathrm{Id}_{\Hi}$ 
are the only projections in the centre of $\mathcal M$.
\end{enumerate}
\end{prop}

\begin{proof}
Since $\mathcal P \subset \mathcal M$,
we have $\mathcal P'' \subset \mathcal M''$,
hence $\mathcal P'' \subset \mathcal M$
by the von Neumann bicommutant theorem \ref{thmvNbicommutant}.
To show \ref{1DEProp-ProjvN}, it remains to show that $\mathcal M \subset \mathcal P''$.
We choose $S \in \mathcal P'$ and $T \in \mathcal M$,
and we have to show that $ST = TS$.
Upon considering the real part and the imaginary part of $T$
(which both belong to $\mathcal M$),
we may assume that $T$ is a selfadjoint operator.
\par
By the spectral theory of selfadjoint operators, we have
$$
\{T\}'' \, = \, \{ \mathbf{1}_{B}(T) \mid B \in \mathcal B (\R) \}'',
$$
where the projection $\mathbf{1}_{B}(T)$ is defined by functional calculus
(see \cite[Theorem 12.23]{Rudi--73};
also, compare with \ref{iiDEthmSNAG} in the SNAG Theorem~\ref{thmSNAG}).
We have $\mathbf{1}_B(T) \in \mathcal P$ for all $B \in \mathcal B )\R)$;
see \cite[Chap.~I, \S~1, no.~1]{Dixm--vN}.
It follows that $S \in \{T\}'$.
\par

Claims \ref{2DEProp-ProjvN} and \ref{3DEProp-ProjvN}
follow from Claim \ref{1DEProp-ProjvN} and the definitions.
\end{proof}

\begin{defn}
% A.11.5
\label{defsupportoperator}
\index{Support of an operator}
Let $\Hi$ be a Hilbert space.
For an operator $T \in \Li (\Hi)$, the \textbf{support} of $T$
is the projection $E$ of $\Hi$ onto the orthogonal of the kernel of~$T$.
\par

Note that $E$ is the smallest of the projections $E_1 \in \Li (\Hi)$
such that $T E_1 = T$
(projections are ordered by $E_1 \le E_2$ if $E_2 - E_1$ is a positive operator).
\par

Let $F$ denote the projection of $\Hi$ onto the closure of the image $T(\Hi)$.
Then $F$ is the smallest of the projections $F_1 \in \Li (\Hi)$
such that $F_1 T = T$.
The last equality is equivalent to $T^* F_1 = F_1$;
it follows that $F$ is the support of $T^*$.
\end{defn}

\begin{defn}
% A.11.6
\label{defcentralsupportoperator}
\index{Central support of an element in a von Neumann algebra}
Let $\mathcal M$ be a von Neumann algebra on $\Hi$, and $\mathcal Z$ its centre.
Let $T$ be in the commutant $\mathcal M'$ of $\mathcal M$.
For a projection $E_1 \in \mathcal Z = \mathcal M' \cap \mathcal M$ such that $T E_1 = T$,
we have $T^* E_1 = E_1 T^* = (T E_1)^* = T^*$.
It follows that the smallest projection $E$ in $\mathcal Z$
which is larger than the support of $T$
is also the smallest projection in $\mathcal Z$
which is larger than the support of $T^*$.
This $E$ is the \textbf{central support} of $T$.
\par

Let $P$ be a projection in $\mathcal M'$; set $\Ki = P(\Hi)$.
The subspace $\Ki$ of $\Hi$ is invariant by $\mathcal M$.
The central support $E$ of $P$ belongs to $\mathcal Z = \mathcal M' \cap \mathcal M$,
and is the smallest of the projections $E_1$ in $\mathcal Z$ such that $P \le E_1$.
The set
$$
\{P S\vert_{\Ki} \mid S \in \mathcal M \} \, \cong \, P \mathcal M P
$$
is a von Neumann subalgebra of $\Li (\Ki)$,
denoted by $\mathcal M_P$ and called the \textbf{induced von Neumann algebra},
or the von Neumann algebra induced by $\mathcal M$ on $\Ki$.
The \textbf{induction}
$$
\mathcal M \to \mathcal M_P ,
\hskip.3cm
S \mapsto PS \vert_\Ki
$$
is a homomorphism of von Neumann algebras.
For another projection $Q \in \mathcal M'$ such that $P \le Q$,
we have inductions $\mathcal M \to \mathcal M_Q$ and $\mathcal M_Q \to \mathcal M_P$,
and their composition is the induction $\mathcal M \to \mathcal M_P$.
\end{defn}

For the next propositions, see in \cite[Chap.~I]{Dixm--vN},
Proposition 2 of \S~2, no~1, and Theorem 3 of \S~4, no~4.
\index{Induced von Neumann algebra}

\begin{prop}
% A.11.7
\label{inductionisiso}
Let $\mathcal M$ be a von Neumann algebra acting on $\Hi$.
Let $P$ be a projection in $\mathcal M'$
\par

If the central support of $P$ is $\mathrm{Id}_\Hi$,
then the induction $\mathcal M \to \mathcal M_P$
is an isomorphism.
\end{prop}

\begin{prop}
% A.11.8
\label{isovNalgebra}
Let $\mathcal M, \mathcal N$ be two von Neumann algebras
acting on Hilbert spaces $\Hi, \Ki$ respectively,
and $\Phi \,\colon \mathcal M \to \mathcal N$ an isomorphism.
\par

Then there exist a Hilbert space $\Hi_2$,
a projection
$P \in \left( \mathcal M \otimes \C \mathrm{Id}_{\Hi_2} \right)' \subset \Li (\Hi \otimes \Hi_2)$,
with range some closed subspace $\Ki_2$ of $\Hi \otimes \Hi_2$,
and an isomorphism of Hilbert spaces $U \,\colon \Ki_2 \to \Ki$, such that
$$
\Phi(T) \, = \, U P \big(T \otimes \mathrm{Id}_{\Hi_2} \big) P U^{-1}
\hskip.2cm \text{for every} \hskip.2cm
T \in \mathcal M .
$$
\end{prop}

[Observe that $\Ki_2$ is invariant under $\mathcal M \otimes \mathrm{Id}_{\Hi_2}$,
because $P \in \left( \mathcal M \otimes \mathrm{Id}_{\Hi_2} \right)'$.]

%-----------------------------------------------------------------------
% End of Appendix
%-----------------------------------------------------------------------

%-----------------------------------------------------------------------------
% Beginning of biblio.tex
%-----------------------------------------------------------------------------

\bibliographystyle{amsalpha}

\backmatter
%    Bibliographies can be prepared with BibTeX using amsplain,
%    amsalpha, or (for "historical" overviews) natbib style.
\bibliographystyle{amsplain}
\bibliography{}

%    See note above about multiple indexes.
\printindex
% \printindex[notation] does not work December 7

\end{document}